# Mathematics for Machine Learning and Data Science:

## Optimization with Mathematica® Applications

**Mohamed M. Hammad**

**Mohamed M. Yahia**

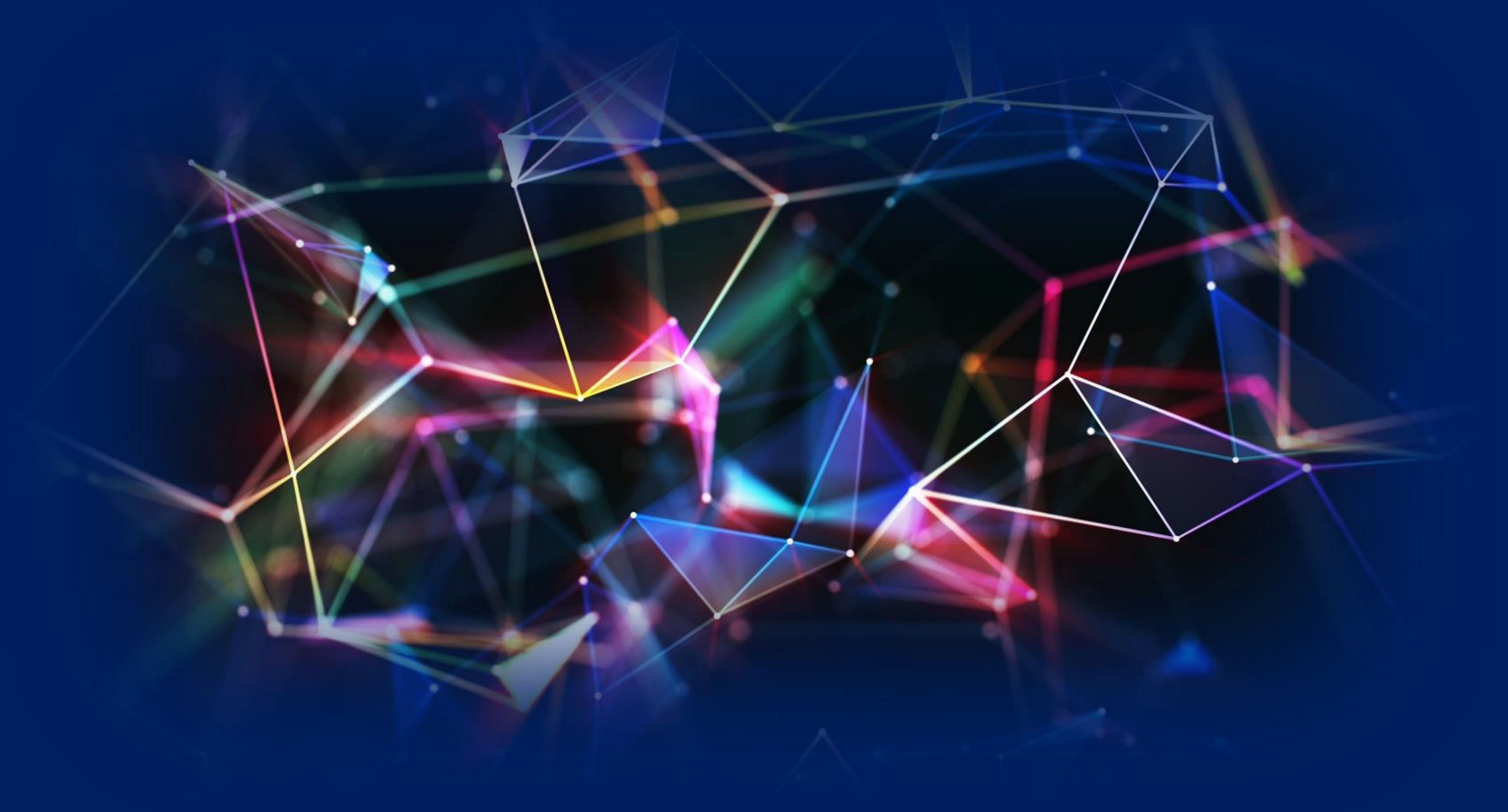

# Mathematics for Machine Learning and Data Science:

## Optimization with Mathematica® Applications

**M. M. Hammad**

Department of Mathematics

Faculty of Science

Damanhour University, Egypt

https://orcid.org/0000-0003-0306-9719

m_hammad@sci.dmu.edu.eg

**M. M. Yahia**

Department of Mathematics

Faculty of Science

Damanhour University, Egypt

https://orcid.org/0000-0001-9274-6915

mohamed.yehia@sci.dmu.edu.eg

*To*

*my*

*mother*

# Preface

Optimization, linear algebra, probability, and statistics are cornerstones of machine learning. The majority of current machine learning textbooks concentrate on machine learning techniques and algorithms and presuppose that readers are proficient in mathematics and statistics. The gap between high school mathematics and the mathematical level necessary to read a typical machine learning textbook is too large for many people who have taught undergraduate and graduate courses at institutions. So, we present the first part of a series of books titled "Mathematics for Machine Learning and Data Science." In this part, we focus on the fundamental optimization theory and contemporary computing techniques that form the foundations of modern machine learning.

Students studying engineering, operations research, data science, and mathematics at the undergraduate and graduate levels might find this book helpful. The book is based on the lectures given to students in the Department of Mathematics, Faculty of Science, Damanhour University, Egypt. We assume that the reader has no prior experience in optimization and scripts. We have provided manually solved examples and examples solved using codes. We tried to provide proofs as simply as possible so that any reader with a background in calculus could easily follow them. It can be used as a textbook for a course spanning two semesters.

Most books on optimization tend to either be overly theoretical or present computational algorithms without enough mathematical background. The present work adopts a strategy that lies somewhere in the middle of these two directions. The computational optimization books present optimization methods to readers to carry out calculations manually or develop these algorithms on their own. This is obviously unrealistic for a typical introductory optimization course, which discusses a wide range of optimization strategies. Recently, some important books have presented optimization programs in built-in computer languages like R, Julia, and Python.

This book offers an introduction to optimization with a focus on practical algorithms using Wolfram Language. We built the optimization algorithms from scratch. For educational goals, we have developed 27 Mathematica codes (functions) on the lines of 18 algorithms to minimize nonlinear objective functions. The outputs of each Mathematica code are tables of the results of each iteration, 3D plots, a contour plot of the results of each iteration, and an interactive application (manipulate the data to make it more visible and understandable). However, a knowledge of Mathematica is not a prerequisite for benefiting from this book. Moreover, we provided examples of built-in functions of Mathematica. Actually, a full range of state-of-the-art local and global optimization techniques, both numeric and symbolic, including constrained nonlinear optimization, interior point methods, and integer programming—as well as original symbolic methods, are integrated into the Wolfram Language. The Wolfram Language symbolic architecture provides seamless access to industrial-strength systems and model optimization, efficiently handling million-variable linear programming and multi-thousand-variable nonlinear problems.

Finally, we express our deep gratitude to Professor Amr R. El Dhaba for his valuable discussions and for supporting this work. Also, we are grateful to many people who looked at early drafts and provided valuable comments: Dr. Fatma El-Safty, Dr. Hamdy El Shamy, Ayman A. Abdelaziz, Eman Farag, Hassan M. Shetawy, Walaa Mansour, Moaz El-Essawey, Aziza Salah, and Eman R. Hendawy.

We hope this book will give data scientists the tools they need to succeed in a data-driven world.

*knowledge itself is power - Sir Francis Bacon 1597*

Egypt 2023

M. M. Hammad
M. M. Yahia

# Abstract

The field of optimization has gotten a lot of interest in recent years owing to significant advances in computer technology. Numerous issues in machine learning, economics, finance, geophysics, molecular modeling, computational systems biology, operations research, and all areas of engineering are now being resolved owing to the rapid growth of optimization methods and algorithms. This monograph presents the main theorems in linear algebra, convex sets, convex functions, single variable optimization, multivariable optimization, and their corresponding algorithms. We also briefly touch upon the constrained nonlinear optimization. We have found the Wolfram language to be ideal for specifying algorithms in human-readable form. To minimize nonlinear objective functions, we have created 27 Mathematica functions (programs) that follow the principles of 18 algorithms. The code examples were carefully designed to demonstrate the purpose of a given algorithm. The code for each algorithm will run as-is with no code from prior algorithms or third parties required beyond the installation of Mathematica.

# Contents

# CHAPTER 1

# MATHEMATICA NOTEBOOKS

## 1.1 Mathematica Notebooks, Palettes, Packages, and Help

Mathematica is a computer algebra system that performs numeric, symbolic, and graphical computations. Although Mathematica can be used as a programming language, its high-level structure is more appropriate for performing sophisticated operations through the use of built-in functions. For example, Mathematica can find limits, derivatives, integrals, and determinants, as well as plot the graph of functions and perform symbolic computations. The number of built-in functions in Mathematica is enormous [1-11]. Our goals in this introductory chapter are modest. Namely, we introduce a small subset of Mathematica commands necessary to explore Mathematica discussed in this book.

### Notebooks

A notebook is a document which allows us to interact with Mathematica. Each notebook is divided up into a sequence of individual units called cells, each containing a specific type of information such as text, graphics, input or output. Text cells contain information to be read by the user but contain no executable Mathematica commands. The following cell, displaying `In[1]` $2^{20}$, is an example of an input cell containing executable Mathematica commands. Mathematica computes the value of $2^{20}$ and the results of the calculation are displayed as `Out[1]:1048576` in an output cell. When we create a new cell, the default cell type is an input cell. Suppose instead, we want to create a text cell. To do this, use the mouse to click in an area where we want to create a new cell and a horizontal line will appear. Then from the Format menu, select Style and then Text. A new text cell will then be created as soon as we begin typing. We can experiment with creating other types of cells by selecting a cell style of our choice, after first choosing Format and Style from the menu.

### Palettes

A palette is similar to a set of calculator buttons, providing shortcuts to entering commands and symbols into a notebook. The name of a useful palette is “Basic Math Assistant Input” and it can be found by selecting Palettes menu and then Basic Math Assistant. After opening Basic Math Assistant, drag it to the right side of the screen and resize the notebook, if necessary, so that both the notebook and palette are visible in non-overlapping windows. To demonstrate the usefulness of palettes, suppose we wish to calculate $\sqrt{804609}$. The Mathematica command for computing the square root of n is `Sqrt[n]`. The following input cell was created by typing in the information exclusively from the keyboard.

```
Input     Sqrt[804609]
Output    897
```

A quicker and more natural way of entering $\sqrt{804609}$ can be accomplished by clicking on the square root button √□ in the palette and then entering 804609.

```
Input      √(804609)
Output     897
```

### Packages

Note that, many of Mathematica functions are available at startup, but additional specialized functions are available from add-in packages. You can load a built-in or installed package in two ways, with the `Needs[ ]` function or with the symbols `<<`. The package name has quotation marks if you use the `Needs[ ]` function, but does not has a mark with `<<`. Package names are always indicated with a backwards apostrophe at the end of the name, `` ` ``.

```
Input       Needs["PackageName`"]
Input       <<PackageName`
```

### Help with Mathematica

There are four important tricks to keep in mind to help with Mathematica:

1- If you want to know something about a Mathematica function or procedure, just type ? followed by a Mathematica command name, and then enter the cell to get information on that command.

```
Input    ?FactorInteger
```

After you press Enter, Mathematica responds:

```
Output   FactorInteger[n] gives a list of the prime factors of the integer n, together with
         their exponents.
```

2- Mathematica also can finish typing a command for you if you provide the first few letters. Here is how it works: After typing a few letters choose Complete Selection from the Edit menu. If more than one completion is possible, you will be presented with a pop-up menu containing all of the options. Just click on the appropriate choice.

3- If you know the name of a command, but have forgotten the syntax for its arguments, type the command name in an input cell, then choose Make Template from the Edit menu. Mathematica will paste a template into the input cell showing the syntax for the simplest form of the command. For example, if you typed `Plot`, and then choose Make Template, the input cell would look like this:

```
Input       Plot[f,{x,xmin, xmax}]
```

4- The Wolfram Documentation is the most useful feature imaginable; learn to use it and use it often. Go to the Help menu and choose Wolfram Documentation. A window will appear displaying the documentation home page.

## 1.2 Basic Concepts

### Document the Code

When you write programs in the Wolfram Language, there are various ways to document your code. As always, by far the best thing is to write clear code, and to name the objects you define as explicitly as possible. Sometimes, however, you may want to add some "commentary text" to your code, to make it easier to understand. You can add such text at any point in your code simply by enclosing it in matching `(*  *)`. Notice that in the Wolfram Language, "comments" enclosed in `(*  *)` can be nested in any way.

| | |
|---|---|
| `(* text *)` | a comment that can be inserted anywhere in Wolfram Language code. |

***Mathematica Examples 1.1***

```
Input      If[a > b,(* then *) p,(* else *) q]
Output     If[a > b, p, q]
```

### Arithmetic Operations

Mathematica can be thought of as a sophisticated calculator, able to perform exact as well as approximate arithmetic computations. You can always control grouping the arithmetic computations by explicitly using parentheses. The following list summarizes the Mathematica symbols used for addition, subtraction, multiplication, division and powers.

| | |
|---|---|
| `x+y+z` | gives the sum of three numbers. |
| `x*y*z, x×y×z, or x y z` | represents a product of terms. |
| `x-y` | is equivalent to $x + (-1 * y)$. |
| `x^y` | gives x to the power y. |
| `x/y` | is equivalent to $x\ y^\wedge - 1$. |

***Mathematica Examples 1.2***

```
Input     2.3+5.63
Output    7.93

Input     2.4/8.9^2
Output    0.0302992

Input     2*3*4
Output    24

Input     (3+4)^2-2(3+1)
Output    41
```

Precedence of common operators is generally defined so that "higher-level" operations are performed first. For simple expressions, operations are typically ordered from highest to lowest in the order: 1. Parenthesization, 2. Factorial, 3. Exponentiation, 4. Multiplication and division, 5. Addition and subtraction. Consider the expression `3×7+2^2`. This expression has value `(3×7)+(2^2)=25`.

Mathematica has several built-in constants. The three most commonly used constants are $\pi$, $e$ and $i$. You can find each of these constants on the Basic Math Assistant palette. Some built-in constants are listed below.

| | |
|---|---|
| `I` | $(i = \sqrt{(-1)})$. |
| `E` | (2.71828). |
| `Pi` | $(\pi = 3.14159)$. |

## Relational and Logical Operators

Relational and logical operators are instrumental in program flow control. They are used in Mathematica to test various conditions involving variables and expressions. The relational operators are listed below.

| | |
|---|---|
| `lhs==rhs` | returns True if lhs and rhs are identical. |
| `lhs!=rhs or lhs≠rhs` | returns False if lhs and rhs are identical. |
| `x>y` | yields True if x is determined to be greater than y. |
| `x>=y or x≥y` | yields True if x is determined to be greater than or equal to y. |
| `x<y` | yields True if x is determined to be less than y. |
| `x<=y or x≤y` | yields True if x is determined to be less than or equal to y. |

Logical operators are used to negate or combine relational expressions. The standard logical operators are listed below.

| | |
|---|---|
| $e_1$`&&`$e_2$`&&...` | is the logical AND function. It evaluates its arguments in order, giving False immediately if any of them are False, and True if they are all True. |
| $e_1$`\|\|`$e_2$`\|\|...` | is the logical OR function. It evaluates its arguments in order, giving True immediately if any of them are True, and False if they are all False. |
| `!expr` | is the logical NOT function. It gives False if expr is True, and True if it is False. |

***Mathematica Examples 1.3***

```
Input       10<7
Output      False

Input       Pi^E < E^Pi
Output      True

Input       2+2 == 4
Output      True

Input       (* Represent an equation: *)
            x^2 == 1+x
Output      x² == 1+x

Input       (* Returns True if elements are guaranteed unequal, and otherwise stays
            unevaluated: *)
            a!=b
Output      a≠b

Input       1!=2
Output      True

Input       1>2 || Pi>3
Output      True

Input       2>1 && Pi>3
Output      True

Input       (3 < 5) || (4 < 5)
Output      True

Input       (3 < 5) && !(4 > 5)
Output      True
```

## 1.3 Elementary Functions

In this section we discuss some of the more commonly used functions Mathematica offers [1,2]. The Wolfram Language has nearly 6000 built-in functions. All have names in which each word starts with a capital letter. Remember that the argument of a function must be contained within square brackets, `[ ]`. Arguments to functions are always separated by commas.

### Common Functions

| | |
|---|---|
| `Log[z]` | gives the natural logarithm of z (logarithm to base *e*). |
| `Log[b,z]` | gives the logarithm to base b. |
| `Exp[z]` | gives the exponential of z. |
| `Sqrt[z]` or `√z` | gives the square root of z. |
| `N[expr]` | gives the numerical value of expr. |
| `Abs[z]` | gives the absolute value of the real or complex number z. |
| `Floor[x]` | gives the greatest integer less than or equal to x. |

### Trigonometric Functions

| | |
|---|---|
| `Sin[z]` | gives the sine of z. |
| `Cos[z]` | gives the cosine of z. |
| `Tan[z]` | gives the tangent of z. |

## Hyperbolic Functions

| | |
|---|---|
| `Sinh[z]` | gives the hyperbolic sine of z. |
| `Cosh[z]` | gives the hyperbolic cosine of z. |
| `Tanh[z]` | gives the hyperbolic tangent of z. |

## Numerical Functions

| | |
|---|---|
| `IntegerPart[x]` | integer part of x. |
| `FractionalPart[x]` | fractional part of x. |
| `Round[x]` | integer x closest to x. |
| `Max[x₁,x₂,...]` | the maximum of $x_1, x_2, \ldots$. |
| `Min[x₁,x₂,...]` | the minimum of $x_1, x_2, \ldots$. |
| `Re[z]` | the real part Re z. |
| `Im[z]` | the imaginary part Im z. |
| `Conjugate[z]` | the complex conjugate $z^*$. |

## Combinatorial Functions

| | |
|---|---|
| `n!` | factorial $n(n-1)(n-2)\ldots \times 2 \times 1$. |
| `n!!` | double factorial $n(n-2)(n-4)\ldots \times 3 \times 1$. |
| `Binomial[n,m]` | binomial coefficient $\binom{n}{m} = \frac{(n!)}{(m!(n-m)!)}$. |
| `Multinomial[n₁,n₂,...]` | multinomial coefficient $(n_1 + n_2 + \ldots)/(n_1!\, n_2! \ldots)$. |

### *Mathematica Examples 1.4*

```
Input     Log[10,1000]
Output    3

Input     Exp[I Pi/5]
Output    E^((I Pi)/5)

Input     Sin[Pi/3]
Output    √3/2

Input     Sinh[1.4]
Output    1.9043

Input     N[1/7]
Output    0.142857

Input     Floor[2.4]
Output    2

Input     30!
Output    265252859812191058636308480000000

Input     Binomial[n,2]
Output    1/2 (-1+n) n

Input     Multinomial[6,5]
Output    462
```

### Sum and Product Functions

Sums and products are of fundamental importance in mathematics, and Mathematica makes their computation simple. Unlike other computer languages, initialization is automatic and the syntax is easy to apply, particularly if the Basic Math Assistant Input palette is used. Any symbol may be used as the index of summation. Negative increments are permitted wherever increment is used.

| | |
|---|---|
| `Sum[f,{i,imax}]` | evaluates the $\sum_{i=1}^{i_{max}} f$. |
| `Sum[f,{i,imin,imax}]` | starts with $i = i_{min}$. |
| `Sum[f,{i,imin,imax,di}]` | uses steps di. |
| `Sum[f,{i,{i₁,i₂,…}}]` | uses successive values $i_1, i_2, \ldots$. |
| `Sum[f,{i,imin,imax},{j,jmin,jmax},…]` | evaluates the multiple sum $\sum_{i=1}^{i_{max}} \sum_{j=1}^{j_{max}} f$. |
| `Product[f,{i,imax}]` | evaluates the $\prod_{i=1}^{i_{max}} f$. |
| `Product[f,{i,imin,imax}]` | starts with $i = i_{min}$. |
| `Product[f,{i,imin,imax,di}]` | uses steps di. |
| `Product[f,{i,{i₁,i₂,…}}]` | uses successive values $i_1, i_2, \ldots$. |
| `Product[f,{i,imin,imax},{j,jmin,jmax},…]` | evaluates the multiple sum $\prod_{i=1}^{i_{max}} \prod_{j=1}^{j_{max}} f$. |

***Mathematica Examples 1.5***

```
Input       (* Numeric sum: *)
            Sum[i^2,{i,10}]
Output      385

Input       (* Symbolic sum: *)
            Sum[i^2,{i,1,n}]
Output      1/6 n (1+n) (1+2n)

Input       Sum[1/i^6, {i, 1, Infinity}]
Output      π^6/945

Input       (* Multiple sum with summation over j performed first: *)
            Sum[1/(j^2 (i + 1)^2), {i, 1, Infinity}, {j, 1, i}]
Output      π^4/120

Input       Product[i^2,{i,1,6}]
Output      518400

Input       Product[i^2,{i,1,n}]
Output      (n!)²

Input       Product[2^j+i,{i,1,p},{j,1,i}]
Output      2^(1/2 p (1+p)²)
```

### Limit and Series Functions

| | |
|---|---|
| `Limit[expr,x->x0]` | finds the limiting value of expr when x approaches $x_0$. |
| `Series[f,{x,x0,n}]` | generates a power series expansion for f about the point $x=x_0$ to order $(x - x_0)^n$. |

***Mathematica Examples 1.6***

```
Input       Limit[(Sin[x])/x,x->0]
Output      1

Input       Limit[(1+x/n)^n, n→Infinity]
Output      E^x
```

```
Input     (* Power series for the exponential function around x=0: *)
          Series[Exp[x],{x,0,10}]
Output    1+ x + x^2/2 + x^3/6 + x^4/24 + x^5/120 + x^6/720 + x^7/5040 + x^8/40320 +
          x^9/362880 + x^10/3628800 + O[x]^11

Input     (* Power series of an arbitrary function around x=a: *)
          Series[f[x],{x,a,3}]
Output    f[a]+f' [a] (x-a)+1/2 f'' [a] (x-a)^2+1/6 f^(3) [a] (x-a)^3+ O[x-a]^4

Input     Series[x^x,{x,0,4}]
Output    1+Log[x]x+1/2 Log[x]^2 x^2+1/6 Log[x]^3 x^3+1/24 Log[x]^4 x^4+O[x]^5
```

### Differentiation Function

| | |
|---|---|
| `D[f,x]` | gives the partial derivative $\frac{\partial}{\partial x}$ f. |
| `D[f,x,y,...]` | gives the derivative $\frac{\partial}{\partial x}\frac{\partial}{\partial y}$ ... f. |
| `D[f,{x,n}]` | gives the multiple derivative $\frac{\partial^n}{\partial x^n}$ f. |

***Mathematica Examples 1.7***

```
Input     (*Derivative with respect to x:*)
          D[x^n,x]
Output    n x^(-1+n)

Input     (* Fourth derivative with respect to x: *)
          D[Sin[x]^10,{x,4}]
Output    5040 Cos[x]^4 Sin[x]^6-4680 Cos[x]^2 Sin[x]^8+280 Sin[x]^10

Input     (* Derivative with respect to x and y: *)
          D[(Sin[x y])/((x^2+y^2)),x,y]
Output    -(2 x^2 Cos[xy])/ (x^2+y^2)^2 -(2 y^2 Cos[x y])/(x^2+y^2)^2 +(Cos[x y])/(x^2+y^2
          )+(8 x y Sin[x y])/(x^2+y^2)^3 -(x y Sin[x y])/(x^2+y^2 )

Input     (* Derivative involving a symbolic function f: *)
          D[x f[x] f'[x],x]
Output    f[x] f' [x] + x f' [x]^2+ x f[x] f'' [x]

Input     D[Sin[x] Cos[x+y], x, y]
Output    -Cos[x+y] Sin[x]-Cos[x] Sin[x+y]

Input     D[ArcCoth[x],{x,2}]
Output    2 x/ (1-x^2)^2
```

### Integration Functions

| | |
|---|---|
| `Integrate[f,x]` | gives the indefinite integral $\int$ f d x. |
| `Integrate[f,{x,xmin,xmax}]` | gives the definite integral $\int_{x_{min}}^{x_{max}}$ f d x. |
| `Integrate[f,{x,xmin,xmax},{y,ymin,ymax},…]` | gives the multiple integral $\int_{x_{min}}^{x_{max}}$ d x $\int_{y_{min}}^{y_{max}}$ dy ... f . |

***Mathematica Examples 1.8***

```
Input     (* Compute an indefinite integral: *)
          Integrate[1/((x^3+1)), x]
Output    (ArcTan[(-1+2 x)/√3])/√3+1/3 Log[1+x]-1/6 Log[1-x+x^2]
```

```
Input       ∫Sqrt[x+ Sqrt[x]]dx
Output      1/12 √(√x+x) (-3+2 √x+8 x)+1/8 Log[1+2√x+2√(√x+x)]

Input       Integrate[1/((x^4+x^2+1)),{x,0,Infinity}]
Output      π/(2√3)

Input       (* Compute an definite integral: *)
            Integrate[x/(Sqrt[1-x]),{x,0,1}]
Output      4/3

Input       Integrate[1/(((2+x^2) Sqrt[4+3 x^2])),{x,-Infinity,Infinity}]
Output      ArcCosh[√(3/2)]

Input       Integrate[x^2+y^2,{x,0,1},{y,0,x}]
Output      1/3
```

## Algebraic Operations

Mathematica has many functions for transforming algebraic expressions [3]. The following list summarizes them.

| | |
|---|---|
| `Simplify[expr]` | performs a sequence of algebraic and other transformations on expr and returns the simplest form it finds. |
| `Expand[expr]` | expands out products and positive integer powers in expr. |
| `Factor[expr]` | factors a polynomial over the integers. |
| `Together[expr]` | puts terms in a sum over a common denominator, and cancels factors in the result. |
| `ExpandAll[expr]` | expands out all products and integer powers in any part of expr. |
| `FunctionExpand[expr]` | tries to expand out special and certain other functions in expr when possible reducing compound arguments to simpler ones. |
| `Reduce[expr,vars]` | reduces the statement expr by solving equations or inequalities for vars and eliminating quantifiers. |

***Mathematica Examples 1.9***

```
Input       Simplify[Sin[x]^2+Cos[x]^2]
Output      1

Input       Expand[(1+x)^10]
Output      1+10 x+45 x²+120 x³+210 x⁴+252 x⁵+210 x⁶+120 x⁷+45 x⁸+10 x⁹+x¹⁰

Input       Factor[x^10-1]
Output      (-1+x)(1+x)(1-x+x²-x³+x⁴)(1+x+x²+x³+x⁴)

Input       Together[x^2/(x^2-1)+x/(x^2-1)]
Output      (x/((-1+x)))

Input       (* Expand polynomials anywhere inside an expression: *)
            ExpandAll[1/(1+x)^3+Sin[(1+x)^3]]
Output      (1/((1+3 x+3 x²+x³)))+Sin[1+3 x+3 x²+x³]

Input       FunctionExpand[Sin[24 Degree]]
Output      -(1/8)Sqrt[3](-1-Sqrt[5])-1/4 Sqrt[1/2 (5-Sqrt[5])]
```

## Solving Equations

Solutions of general algebraic equations may be found using the `Solve` command. `Solve` always tries to give you explicit formulas for the solutions to equations. However, it is a basic mathematical result that, for sufficiently complicated equations, explicit algebraic formulas in terms of radicals cannot be given. If you have an algebraic equation in one variable, and the highest power of the variable is at most four, then the Wolfram Language can always give you formulas for the solutions. However, if the highest power is five or more, it may be mathematically impossible to give explicit algebraic formulas for all the solutions.

You can also use the Wolfram Language to solve sets of simultaneous equations. You simply give the list of equations and specify the list of variables to solve for. Not all algebraic equations are solvable by Mathematica, even if theoretical solutions exist. If Mathematica is unable to solve an equation, it will represent the solution in a symbolic form. For the most part, such solutions are useless, and a numerical approximation is more appropriate. Numerical approximations are obtained with the command `NSolve`.

| | |
|---|---|
| `Solve[lhs==rhs,x]` | solve an equation for x. |
| `Solve[{lhs1==rhs1,lhs2==rhs2,…},{x,y,…}]` | solve a set of simultaneous equations for x, y, .... |
| `Eliminate[{lhs1==rhs1,lhs2==rhs2,…},{x,…}]` | eliminate x, ... in a set of simultaneous equations. |
| `Reduce[{lhs1==rhs1,lhs2==rhs2,…},{x,y,…}]` | give a set of simplified equations, including all possible solutions. |
| `NSolve[expr,vars]` | attempts to find numerical approximations to the solutions of the system expr of equations or inequalities for the variables vars. |
| `FindRoot[f,{x,x₀}]` | searches for a numerical root of f, starting from the point $x = x_0$. |

***Mathematica Examples 1.10***

```
Input     Solve[x^2+ax+1==0,x]
Output    {{x→1/2 (-a-√(-4+a²)},{x→1/2 (-a+√(-4+a²))}}

Input     Solve[a x+y==7&& b x-y==1,{x,y}]
Output    {{x→8/(a+b),y→-(a-7 b)/(a+b)}}

Input     (* Eliminate the variable y between two equations: *)
          Eliminate[{x==2+y,y==z},y]
Output    2+z==x

Input     NSolve[x^5-2x+3==0,x,Reals]
Output    {{x→-1.42361}}

Input     Reduce[x^2-y^3==1,{x,y}]
Output    y==(-1+x²)^{1/3}||y==-(-1)^{1/3}(-1+x²)^{1/3}||y==(-1)^{2/3}(-  1+x²)^{1/3}

Input     FindRoot[Sin[x]+Exp[x],{x,0}]
Output    {x→-0.588533}
```

## Some Notes

1- In doing calculations, you will often need to use previous results that you have got. In the Wolfram Language, `%` always stands for your last result.

| | |
|---|---|
| `%` | the last result generated. |
| `%%` | the next-to-last result. |
| `% n` | the result on output line Out[n]. |
| `Out[n]` | is a global object that is assigned to be the value produced on the $n^{th}$ output line. |

***Mathematica Examples 1.11***

```
Input       77^2
Output      5929

Input       %+1
Output      5930

Input       3 % + %^2 + %%
Output      35188619

Input       % 2+% 3
Output      175943095
```

2- Although Mathematica is a powerful calculating tool, it has its limits. Sometimes it will happen that the calculations you tell Mathematica to do are too complicated or may be the output produced is too long. In these cases, Mathematica could be calculating for too long to get an output so you might want to stop these calculations. To abort a calculation: go to "Kernel" and select "Abort evaluation". It can take long to abort a calculation. If the computer does not respond an alternative is to close down the Kernel. By doing this you do not lose the data displayed in your notebooks, but you do lose all the results obtained so far from the Kernel, so in case you are running a series of calculations, you would have to start again. To close down the Kernel: go to "Kernel" and select "Quit Kernel" and then "Local". Closing down the Kernel is not a practice that is done only when you want to stop a calculation. Sometimes, when you have been using Mathematica for a long time you forget about the definitions and calculations that you have done before (you might have defined values for variables or functions, for example). Those definitions can clash with the calculations you are doing, so you might want to close down the Kernel and start your new calculations from scratch. In general, it is a good idea to close down the Kernel after you have finished with a series of calculations, so that when you move to a different problem your new calculations do not interact with the previous ones.

## 1.4 Variables and Functions

When you perform long calculations, it is often convenient to give names to your intermediate results. Just as in standard mathematics, or in other computer languages, you can do this by introducing named variables [4]. It is very important to realize that values you assign to variables are permanent. Once you have assigned a value to a particular variable, the value will be kept until you explicitly remove it. The value will, of course, disappear if you start a whole new Wolfram Language session.

| | |
|---|---|
| `x=value` | assign a value to the variable x. |
| `x=y=value` | assign a value to both x and y. |
| `x=.or Clear[x]` | remove any value assigned to x. |
| `{x,y}={value1,value2}` | assign different values to x and y. |
| `{x,y}={y,x}` | interchange the values of x and y. |

***Mathematica Examples 1.12***

```
Input       x=5
Output      5

Input       x^2
Output      25

Input       x=7+4
Output      11
```

In Mathematica one can substitute an expression with another using rules. In particular one can substitute a variable with a value without assigning the value to the variable.

| | |
|---|---|
| `lhs:=rhs` | assigns rhs to be the delayed value of lhs. rhs is maintained in an unevaluated form. When lhs appears, it is replaced by rhs, evaluated afresh each time. |
| `expr/.rules` | applies a rule or list of rules in an attempt to transform each subpart of an expression expr. |
| `lhs->rhs or lhs->rhs` | represents a rule that transforms lhs to rhs. |

**Mathematica Examples 1.13**

```
Input       x = 5
Input       y := x + 2
Input       y
Output      7

Input       x = 10
Output      10

Input       y
Output      12
```

**Mathematica Examples 1.14**

```
Input       x+y /. x-> 2
Output      2+y

Input       x + y /. {x -> a, y -> b}
Output      a+b

Input       x^2 + y/.x -> y/.y -> x
Output      x + x^2

Input       x + 2 y/.{x -> y, y -> a}
Output      2a + y
```

The last example reveals that Mathematica goes through the expression only once and replaces the rules. If we need Mathematica to go through the expression again and replace any expression which is possible until no substitution is possible, one uses `//.` . In fact `/.` and `//.` are shorthand for `Replace` and `ReplaceRepeated`, respectively.

**Mathematica Examples 1.15**

```
Input       {x,x^2,a,b} /. x->3
Output      {3,9,a,b}

Input       x + 2y //. {x -> y, y -> a}
Output      3a

Input       x + 2 y //. {x -> b, y -> a, b -> c}
Output      2 a + c

Input       x + 2 y /. {x -> b, y -> a, b -> c}
Output      2 a + b

Input       Sin[x] /. Sin->Cos
Output      Cos[x]
```

There are many functions that are built into the Wolfram Language. Here we discuss how you can add your own simple functions to the Wolfram Language. As a first example, consider adding a function called f which squares its argument. The Wolfram Language command to define this function is `f[x]:=x^2`. The names like f that you use for functions in the Wolfram Language are just symbols. Because of this, you should make sure to avoid using names that begin with capital letters, to prevent confusion with built-in Wolfram Language functions. You should also make sure that you have not used the names for anything else earlier in your session.

| | |
|---|---|
| `f[x]=value` | definition for a specific expression x. |
| `f[x_]=value` | definition for any expression, referred to as x. |
| `Clear[f]` | clear all definitions for f. |
| `Function[x,body]` | is a pure function with a single formal parameter x. |
| `Function[{x_1,x_2,…},body]` | is a pure function with a list of formal parameters. |
| `Map[f,expr] or f/@expr` | applies f to each element on the first level in expr. |
| `Map[f,expr,levelspec]` | applies f to parts of expr specified by levelspec. |

The character _ (referred to as "blank") on the left-hand side is very important.

***Mathematica Examples 1.16***

```
Input     f[x_]:=x^2
Input     f[a+1]
Output    (1+a)^2

Input     f[4]
Output    16

Input     f[3 x + x^2]
Output    (3 x + x^2)^2

Input     Expand[f[(x+1+y)]]
Output    1+2 x+ x^2 +2 y+ 2 x y + y^2

Input     Function[u, 3 + u][x]
Output    3 + x

Input     Function[{u, v}, u^2 + v^4][x, y]
Output    x^2 + y^4

Input     (* Evaluate f on each element of a list: *)
          Map[f,{a,b,c,d,e}]
Output    {f[a],f[b],f[c],f[d],f[e]}

Input     f/@{a,b,c,d,e}
Output    {f[a],f[b],f[c],f[d],f[e]}

Input     (* Map at top level: *)
          Map[f,{{a,b},{c,d,e}}]
Output    {f[{a,b}],f[{c,d,e}]}

Input     (* Map at level 2: *)
          Map[f,{{a,b},{c,d,e}},{2}]
Output    {{f[a],f[b]},{f[c],f[d],f[e]}}

Input     (* Map at levels 1 and 2: *)
          Map[f,{{a,b},{c,d,e}},2]
Output    {f[{f[a],f[b]}],f[{f[c],f[d],f[e]}]}
```

One can define functions of several variables. Here is a simple example defining $f(x,y) = \sqrt{x^2 + y^2}$.

***Mathematica Examples 1.17***

```
Input[1]     f[x_, y_] := Sqrt[x^2 + y^2]
Input[2]     f[3, 4]
Output[2]    5
```

### Some Notes

1- There are four kinds of bracketing used in the Wolfram Language. Each kind of bracketing has a very different meaning.

| | |
|---|---|
| `(term)` | parentheses for grouping. |
| `f[x]` | square brackets for functions. |
| `{a,b,c}` | curly braces for lists. |
| `v[[i]]` | double brackets for indexing (Part[v,i]). |

2- Compound expression

| | |
|---|---|
| `expr₁;expr₂;expr₃` | do several operations and give the result of the last one. |
| `expr₁;expr₂;` | do the operations but print no output. |
| `expr;` | do an operation but display no output. |

***Mathematica Examples 1.18***

```
Input [1]     x=4;y=6;z=y+6
Output[1]     12

Input [2]     a=2;b=3;a+b
Output[2]     5
```

3- Particularly when you write procedural programs in the Wolfram Language, you will often need to modify the value of a particular variable repeatedly. You can always do this by constructing the new value and explicitly performing an assignment such as `x=value`. The Wolfram Language, however, provides special notations for incrementing the values of variables, and for some other common cases.

| | |
|---|---|
| `i++` | increment the value of i, by 1 returning the old value of i. |
| `i--` | decrement the value of i, by 1 returning the old value of i. |
| `++i` | pre-increment i, returning the new value of i. |
| `--i` | pre-decrement i, returning the new value of i. |
| `i+=di` | add di to the value of i and returns the new value of i. |
| `i-=di` | subtract di from i and returns the new value of i. |
| `x*=c` | multiply x by c. |
| `x/=c` | divide x by c. |

***Mathematica Examples 1.19***

```
Input [1]     k = 1; k++
Output[1]     1

Input [2]     k
Output[2]     2

Input [3]     k = x
Output[3]     x

Input [4]     k++
Output[4]     x
```

```
Input [5]     k
Output[5]     1+x

Input [6]     k = 1; ++ k
Output[6]     2

Input [7]     k
Output[7]     2

Input [8]     k = 1; k--
Output[8]     1

Input [9]     k
Output[9]     0

Input [10]    k = 1; k -= 5
Output[10]    -4

Input [11]    k
Output[11]   -4
```

4- Primarily there are three equalities in Mathematica, `=`, `:=`, `==`. There is a fundamental differences between `=` and `:=` explained in the following examples:

***Mathematica Examples 1.20***

```
Input [1]    x = 5; y = x + 2;
Input [2]    y
Output[2]    7

Input [3]    x = 10
Output[3]    10

Input [4]    y
Output[4]    7

Input [5]    x = 15
Output[5]    15

Input [6]    y
Output[6]    7
```

```
Input [1]    x = 5; y := x + 2;
Input [2]    y
Output[2]    7

Input [3]    x = 10
Output[3]    10

Input [4]    y
Output[4]    12

Input [5]    x = 15
Output[5]    15

Input [6]    y
Output[6]    17
```

It is clear that when we defined `y=x+2` then y takes the value of x+2 and this will be assigned to y. No matter if x changes its value, the value of y remains the same. In other words, y is independent of x. But in `y:=x+2`, y is dependent on x, and when x changes, the value of y changes too. Namely using `:=` then y is a function with variable x. Finally, the equality `==` is used to compare:

***Mathematica Examples 1.21***

```
Input [1]   5==5
Output[1]   True

Input [2]   3==5
Output[2]   False
```

## 1.5 Lists

Lists are extremely important objects. In doing calculations, it is often convenient to collect together several objects, and treat them as a single entity. Lists give you a way to make collections of objects [5]. Lists are sequences of Mathematica objects separated by commas and enclosed by curly brackets. A list such as `{3,5,1}` is a collection of three objects. But in many ways, you can treat the whole list as a single object. You can, for example, do arithmetic on the whole list at once, or assign the whole list to be the value of a variable.

Defining your own lists is easy. You can, for example, type them in full, like this:

```
Input     oddList = {81, 3, 5, 7, 9, 11, 13, 15, 17}
Output    {81, 3, 5, 7, 9, 11, 13, 15, 17}
```

Alternatively, if (as here) the list elements correspond to a rule of some kind, the command `Table` can be used, like this:

```
Input     oddList = Table[2 n + 1, {n, 0, 8}]
Output    {1, 3, 5, 7, 9, 11, 13, 15, 17}
```

The functions for obtaining elements of lists are

| | |
|---|---|
| `First[list]` | the first element. |
| `Last[list]` | the last element. |
| `Part[list,n] or list[[n]]` | the nth element. |
| `Part[list,-n] or list[[-n]]` | the nth element from the end. |
| `Part[list,{n1,n2,...}] or list[[{n1,n2,...}]]` | the list of the n1th, n2th, ... elements. |
| `Take[list,n]` | the list of the first n elements. |
| `Take[list,-n]` | the list of the last n elements. |
| `Take[list,{m,n}]` | the list of the mth through nth elements. |
| `Rest[list]` | list without the first element. |
| `Most[list]` | list without the last element. |
| `Drop[list,n]` | list without the first n elements. |
| `Drop[list,-n]` | list without the last n elements. |
| `Drop[list,{m,n}]` | list without the mth through nth elements. |

***Mathematica Examples 1.22***

```
Input     First[{a, b, c}]
Output    a

Input     First[{{a, b}, {c, d}}]
Output    {a, b}

Input     Last[{a, b, c}]
Output    c

Input     {a, b, c, d, e, f}[[3]]
Output    c

Input     {{a, b, c}, {d, e, f}, {g, h, i}}[[2, 3]]
Output    f

Input     Take[{a, b, c, d, e, f}, 4]
Output    {a, b, c, d}

Input     Rest[{a, b, c, d}]
Output    {b, c, d}

Input     Most[{a, b, c, d}]
Output    {a, b, c}

Input     Drop[{a, b, c, d, e, f}, 2]
Output    {c, d, e, f}
```

Some functions for inserting, deleting, and replacing list and sublist elements are

| | |
|---|---|
| `Prepend[list,elem]` | insert elem at the beginning of list. |
| `Append[list,elem]` | insert elem at the end of list. |
| `Insert[list,elem,i]` | insert elem at position i in list. |
| `Insert[list,elem,{i,j,...}]` | insert elem at position {i, j, ...} in list. |
| `Insert[list,elem,{{i1,j1,...},{i2,...},...}]` | insert elem at positions {i1, j1, ...},{i2,...}, ... in list. |
| `Delete[list,i]` | delete the element at position i in list. |
| `Delete[list,{i,j,...}]` | delete the element at position {i, j, ...} in list. |
| `Delete[list,{{i1,j1,...},{i2,...},...}]` | delete elements at positions {i1, j1, ...},{i2,... }, ... in list. |
| `ReplacePart[list,elem,i]` | replace the element at position i in list with elem. |
| `ReplacePart[list,elem,{i,j,...}]` | replace the element at position {i, j, ...} with elem. |
| `ReplacePart[list,elem,{{i1,j1,...},{i2,...},...}]` | replace elements at positions {i1, j1, ...},{i2,...}, ... with elem. |

***Mathematica Examples 1.23***

```
Input     Prepend[{a, b, c, d}, x]
Output    {x, a, b, c, d}

Input     Append[{a, b, c, d}, x]
Output    {a, b, c, d, x}

Input     Insert[{a, b, c, d, e}, x, 3]
Output    {a, b, x, c, d, e}

Input     Insert[{a, b, c, d, e}, x, -2]
Output    {a, b, c, d, x, e}

Input     Delete[{a, b, c, d}, 3]
Output    {a, b, d}

Input     Delete[{a, b, c, d}, {{1}, {3}}]
Output    {b, d}

Input     ReplacePart[{a, b, c, d, e}, 3 -> xxx]
Output    {a, b, xxx, d, e}

Input     ReplacePart[{a, b, c, d, e}, {2 -> xx, 5 -> yy}]
Output    {a, xx, c, d, yy}
```

Some functions for rearranging lists are

| | |
|---|---|
| `Sort[list]` | sort the elements of list into canonical order. |
| `Union[list]` | give a sorted version of list, in which all duplicated elements have been dropped. |
| `Reverse[list]` | reverse the order of the elements in list. |
| `RotateLeft[list]` | cycle the elements in list one position to the left. |
| `RotateLeft[list,n]` | cycle the elements in list n positions to the left. |
| `RotateRight[list]` | cycle the elements in list one position to the right. |
| `RotateRight[list,n]` | cycle the elements in list n positions to the right. |
| `Permutations[list]` | generate a list of all possible permutations of the elements in list. |
| `Partition[list,n]` | partition list into nonoverlapping sublists of length n. |
| `Partition[list,n,d]` | generate sublists with offset d. |
| `Split[list]` | split list into sublists consisting of runs of identical elements. |
| `Transpose[list]` | transpose the first two levels in list. |
| `Flatten[list]` | flatten out nested lists. |
| `Flatten[list,n]` | flatten out the top n levels. |
| `FlattenAt[list,i]` | flatten out a sublist that appears as the ith element of list. |
| `FlattenAt[list,{i,j,...}]` | flatten out the element of list at position {i,j, ...}. |
| `FlattenAt[list,{{i1,j1,...},{i2,j2,...}}]` | flatten out elements of list at several positions. |
| `Join[list1,list2,...]` | concatenate lists together. |
| `Union[list1,list2,...]` | give a sorted list of all the distinct elements that appear in any of the listi. |

***Mathematica Examples 1.24***

```
Input     Sort[{d, b, c, a}]
Output    {a, b, c, d}

Input     Sort[{4, 1, 3, 2, 2}, Greater]
Output    {4, 3, 2, 2, 1}

Input     Union[{1, 2, 1, 3, 6, 2, 2}]
Output    {1, 2, 3, 6}

Input     Union[{a, b, a, c}, {d, a, e, b}, {c, a}]
Output    {a, b, c, d, e}

Input     Reverse[{a, b, c, d}]
Output    {d, c, b, a}

Input     RotateLeft[{a, b, c, d, e}, 2]
Output    {c, d, e, a, b}

Input     RotateRight[{a, b, c, d, e}, 2]
Output    {d, e, a, b, c}

Input     Permutations[{a, b, c}]
Output    {{a, b, c}, {a, c, b}, {b, a, c}, {b, c, a}, {c, a, b}, {c, b, a}}

Input     Partition[{a, b, c, d, e, f}, 2]
Output    {{a, b}, {c, d}, {e, f}}

Input     Flatten[{{a, b}, {c, {d}, e}, {f, {g, h}}}]
Output    {a, b, c, d, e, f, g, h}

Input     Transpose[{{a, b, c}, {x, y, z}}]
Output    {{a, x}, {b, y}, {c, z}}

Input     Join[{a, b, c}, {x, y}, {u, v, w}]
Output    {a, b, c, x, y, u, v, w}
```

Vectors and matrices in the Wolfram Language are simply represented by lists and by lists of lists, respectively. Functions for generating lists are `Range[ ]`, `Table[ ]` and `Array[ ]`.

## Vectors

Mathematica has many functions for generating vectors. The following list summarizes them.

| | |
|---|---|
| `{e1,e2,...}` | is a list of elements. |
| `Range[n]` | create the list $\{1,2,3,\ldots,n\}$. |
| `Range[n1,n2]` | create the list $\{n_1, n_1+1,\ldots,n_2\}$. |
| `Range[n1,n2,dn]` | create the list $\{n_1, n_1+dn,\ldots,n_2\}$. |
| `Table[f,{i,n}]` | build a length-n vector by evaluating f with $i = 1,2,\ldots,n$. |
| `Length[list]` | give the number of elements in list. |
| `List[[i]] or Part[list,i]` | give the $i^{th}$ element in the vector list. |

***Mathematica Examples 1.25***

```
Input       List[a,b,c,d]
Output      {a,b,c,d}

Input       v = {x,y}
Output      {x,y}

Input       Range[4]
Output      {1,2,3,4}

Input       Range[x,x+4]
Output      {x,1+x,2+x,3+x,4+x}

Input       Table[i^2,{i,10}]
Output      {1,4,9,16,25,36,49,64,81,100}

Input       Length[{a,b,c,d}]
Output      4

Input       {5,8,6,9}[[2]]
Output      8
```

| | |
|---|---|
| `c v` | multiply a vector v by a scalar. |
| `a.b` | dot product of two vectors a. b. |
| `Cross[a,b]` | cross product of two vectors (also input as $a \times b$). |
| `Norm[v]` | Euclidean norm of a vector v. |
| `Normalize[v]` | gives the normalized form of a vector v. |
| `Orthogonalize [{v1,v2,…}]` | gives an orthonormal basis found by orthogonalizing the vectors $v_i$. |

***Mathematica Examples 1.26***

```
Input       {a,b,c}.{x,y,z}
Output      a x + b y + c z

Input       Cross[{a,b,c},{x,y,z}]
Output      {-c y + b z,c x - a z,-b x + a y}

Input       Norm[{x,y,z}]
Output      √(Abs[x]²+Abs[y]²+Abs[z]²)

Input       Normalize[{1,5,1}]
Output      {1/(3√3),5/(3√3),1/(3√3)}

Input       Orthogonalize[{{1,0,1},{1,1,1}}]
Output      {{1/√2,0,1/√2},{0,1,0}}
```

### Matrix

Mathematica has many functions for generating matrices. The following list summarizes them.

| | |
|---|---|
| `{{a,b},{c,d}}` | matrix $\begin{pmatrix} a & b \\ c & d \end{pmatrix}$. |
| `Table[f,{i,m},{j,n}]` | build an m × n matrix by evaluating f with i ranging from 1 to m and j ranging. |
| `List[[i,j]] or Part[list,i,j]` | give the i, j th element in the matrix list. |
| `DiagonalMatrix[list]` | generate a square matrix with the elements in list on the main. |
| `Dimensions[list]` | give the dimensions of a matrix represented by list. |
| `Column[list]` | display the elements of list in a column. |
| `c m` | multiply a matrix m by a scalar. |
| `a.b` | dot product of two matrices a. b. |
| `Inverse[m]` | matrix inverse m. |
| `MatrixPower[m,n]` | gives the $n^{th}$ power of a matrix m. |
| `Det[m]` | Determinant m. |
| `Tr[m]` | Trace m. |
| `Transpose[m]` | Transpose m. |

***Mathematica Examples 1.27***

```
Input     m = {{a, b},{c, d}}
Output    {{a, b},{c, d}}

Input     m[[1]]
Output    {a, b}

Input     m[[1, 2]]
Output    b

Input     v = {x, y}
Output    {x, y}

Input     m . v
Output    {a x + b y, c x + d y}

Input     m.m
Output    {{a²+b c, a b + b d},{a c + c d , b c + d²}}

Input     s = Table[i+j, {i,3},{j,3}]
Output    {{2, 3, 4}, {3, 4, 5}, {4, 5, 6}}

Input     DiagonalMatrix[{a, b, c}]
Output    {{a, 0, 0}, {0, b, 0}, {0, 0, c}}

Input     Det[m]
Output    -b c + a d

Input     Transpose[m]
Output    {{a, c}, {b, d}}

Input     h = Table[1/(i+j-1), {i, 3}, {j, 3}]
Output    {{1, 1/2, 1/3}, {1/2, 1/3, 1/4}, {1/3, 1/4, 1/5}}

Input     Inverse[h]
Output    {{9, -36, 30}, {-36, 192, -180}, {30, -180, 180}}
```

## Array

| | |
|---|---|
| `Array[f,n]` | generates a list of length n, with elements f[i]. |
| `Array[f,n,r]` | generates a list using the index origin r. |
| `Array[f,n,{a,b}]` | generates a list using n values from a to b. |
| `Array[f,{n1,n2,...}]` | generates an $n_1 \times n_2 \times \ldots$ array of nested lists, with elements $f[i_1, i_2, \ldots]$. |

***Mathematica Examples 1.28***

```
Input [1]    Array[f, 10]
Output[1]    {f[1], f[2], f[3], f[4], f[5], f[6], f[7], f[8], f[9], f[10]}

Input [2]    Array[f, {3, 2}]
Output[2]    {{f[1, 1], f[1, 2]}, {f[2, 1], f[2, 2]}, {f[3, 1], f[3, 2]}}
```

## Layout & Tables

| | |
|---|---|
| `Print[expr]` | prints expr as output. |
| `MatrixForm[list]` | prints with the elements of list arranged in a regular array. |
| `TableForm[list]` | prints with the elements of list arranged in an array of rectangular cells. |
| `Grid[{{expr11,expr12,...},{expr21,expr22,...},...}]` | is an object that formats with the $expr_{ij}$ arranged in a two-dimensional grid. |
| `Row[expr1,expr2,...]` | is an object that formats with the $expr_i$ arranged in a row, potentially extending over several lines. |
| `Row[list,s]` | inserts s as a separator between successive elements. |
| `Column[expr1,expr2,...]` | is an object that formats with the $expr_i$ arranged in a column, with $expr_1$ above $expr_2$, etc. |
| `Multicolumn[list,cols]` | is an object that formats with the elements of list arranged in a grid with the indicated number of columns. |
| `Multicolumn[list,{rows,Automatic}]` | formats as a grid with the indicated number of rows. |

***Mathematica Examples 1.29***

```
Input [1]    MatrixForm[{{1, 2}, {3, 4}}]
Output[1]
              ( 1  2 )
              ( 3  4 )

Input [2]    MatrixForm[Table[1/(i + j), {i, 4}, {j, 4}]]
Output[2]
              ( 1/2  1/3  1/4  1/5 )
              ( 1/3  1/4  1/5  1/6 )
              ( 1/4  1/5  1/6  1/7 )
              ( 1/5  1/6  1/7  1/8 )

Input [3]    TableForm[Table[1/(i + j), {i, 4}, {j, 4}]]
Output[3]
              1/2  1/3  1/4  1/5
              1/3  1/4  1/5  1/6
              1/4  1/5  1/6  1/7
              1/5  1/6  1/7  1/8

Input [4]    Grid[{{a, b, c}, {x, y, z}}]
Output[4]    a  b  c
             x  y  z

Input [5]    Grid[{{a, b, c}, {x, y^2, z^3}}, Frame -> All]
Output[5]
```

| | | |
|---|---|---|
| a | b | c |
| x | y^2 | z^3 |

```
Input [6]     Row[{aaa, b, cccc}]
Output[6]     aaabcccc

Input [7]     Row[{aaa, b, cccc}, "----"]
Output[7]     aaa----b----cccc

Input [8]     Column[{1, 12, 123, 1234}]
Output[8]     1
              12
              123
              1234

Input [9]     Column[{1, 22, 333, 4444}, Frame -> True]
Output[9]     1
              22
              333
              4444

Input [10]    Multicolumn[Range[50], {6, Automatic}]
Output[10]    1  7   13  19  25  31  37  43  49
              2  8   14  20  26  32  38  44  50
              3  9   15  21  27  33  39  45
              4  10  16  22  28  34  40  46
              5  11  17  23  29  35  41  47
              6  12  18  24  30  36  42  48
```

## 1.6 2D and 3D Graphing

The graph of a function offers tremendous insight into the behavior of the function and can be of a great value in the solution of problems in mathematics. One of the outstanding features of Mathematica is its graphing capabilities [6]. Mathematica contains functions for 2D and 3D graphing of functions, lists, and arrays of data.

### Basic Plotting

| | |
|---|---|
| `Plot[f,{x,xmin,xmax}]` | plot f as a function of x from $x_{min}$ to $x_{max}$. |
| `Plot[{f1,f2,...},{x,xmin,xmax}]` | plot several functions together. |

When the Wolfram Language plots a graph for you, it has to make many choices. It has to work out what the scales should be, where the function should be sampled, how the axes should be drawn, and so on. Most of the time, the Wolfram Language will probably make pretty good choices. However, if you want to get the very best possible pictures for your particular purposes, you may have to help the Wolfram Language in making some of its choices.

There is a general mechanism for specifying "options" in Wolfram Language functions. Each option has a definite name. As the last arguments to a function like `Plot`, you can include a sequence of rules of the form `name->value`, to specify the values for various options. Any option for which you do not give an explicit rule is taken to have its "default" value.

Some options for `Plot` function are

| | |
|---|---|
| `AspectRatio` | the height-to-width ratio for the plot; Automatic sets it from the absolute x and y coordinates |
| `Axes` | whether to include axes |
| `AxesLabel` | labels to be put on the axes; ylabel specifies a label for the y axis, {xlabel,ylabel} for both axes |
| `AxesOrigin` | the point at which axes cross |
| `BaseStyle` | the default style to use for the plot |
| `FormatType` | the default format type to use for text in the plot |
| `Frame` | whether to draw a frame around the plot |
| `FrameLabel` | labels to be put around the frame; give a list in clockwise order starting with the lower x axis |
| `FrameTicks` | what tick marks to draw if there is a frame; None gives no tick marks |
| `GridLines` | what grid lines to include; Automatic includes a grid line for every major tick mark |
| `PlotLabel` | an expression to be printed as a label for the plot |
| `PlotRange` | the range of coordinates to include in the plot; All includes all points |
| `Ticks` | what tick marks to draw if there are axes; None gives no tick marks |
| `PlotStyle` | a list of lists of graphics primitives to use for each curve (see "Graphics Directives and Options") |
| `ClippingStyle` | what to draw when curves are clipped |
| `Filling` | filling to insert under each curve |
| `FillingStyle` | style to use for filling |
| `PlotPoints` | the initial number of points at which to sample the function |
| `MaxRecursion` | the maximum number of recursive subdivisions allowed |

***Mathematica Examples 1.30***

Input
```
Plot[
 Sin[x],
 {x,0,2 Pi}
 ]
```
Output

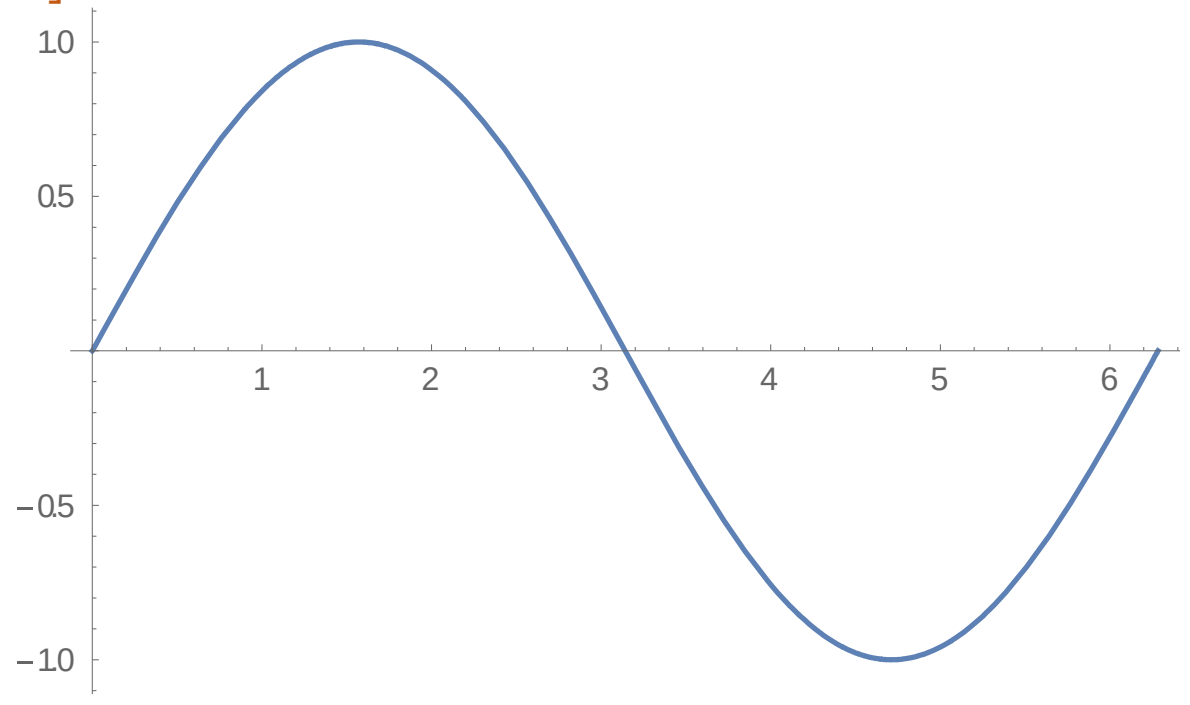

Input
```
Plot[
 {Sin[x],Sin[2 x],Sin[3 x]},
 {x,0,2 Pi}
 ]
```
Output

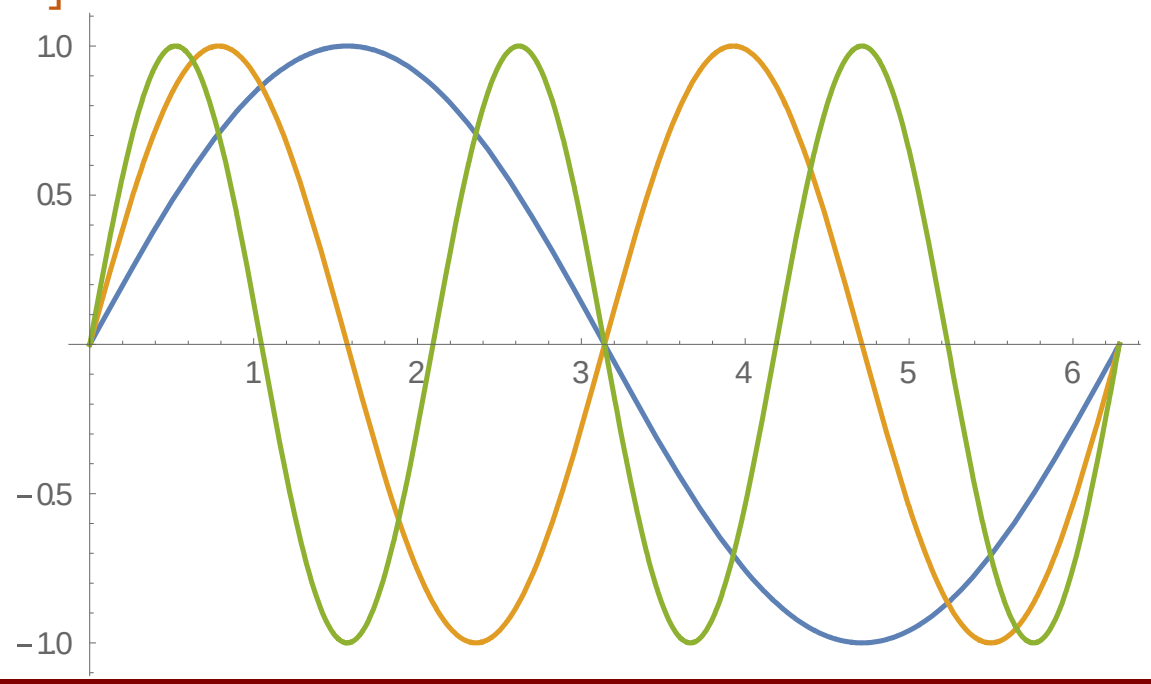

### 3D Plot

| | |
|---|---|
| `Plot3D[f,{x,xmin,xmax},{y,ymin,ymax}]` | make a three-dimensional plot of f as a function of the variables x and y. |

Some options for `Plot3D` function are

| | |
|---|---|
| `Axes` | whether to include axes |
| `AxesLabel` | labels to be put on the axes: zlabel specifies a label for the z axis, {xlabel,ylabel,zlabel} for all axes |
| `BaseStyle` | the default style to use for the plot |
| `Boxed` | whether to draw a three-dimensional box around the surface |
| `FaceGrids` | how to draw grids on faces of the bounding box; All draws a grid on every face |
| `LabelStyle` | style specification for labels |
| `Lighting` | simulated light sources to use |
| `Mesh` | whether an xy mesh should be drawn on the surface |
| `PlotRange` | the range of z or other values to include |
| `SphericalRegion` | whether to make the circumscribing sphere fit in the final display area |
| `ViewAngle` | angle of the field of view |
| `ViewCenter` | point to display at the center |
| `ViewPoint` | the point in space from which to look at the surface |
| `ViewVector` | position and direction of a simulated camera |
| `ViewVertical` | direction to make vertical |
| `BoundaryStyle` | how to draw boundary lines for surfaces |
| `ClippingStyle` | how to draw clipped parts of surfaces |
| `ColorFunction` | how to determine the color of the surfaces |
| `Filling` | filling under each surface |
| `FillingStyle` | style to use for filling |
| `PlotPoints` | the number of points in each direction at which to sample the function; {nx,ny} specifies different numbers in the x and y directions |
| `PlotStyle` | graphics directives for the style of each surface |

***Mathematica Examples 1.31***

Input
```
Plot3D[
 Sin[x+y^2],
 {x,-3,3},
 {y,-2,2}
 ]
```
Output

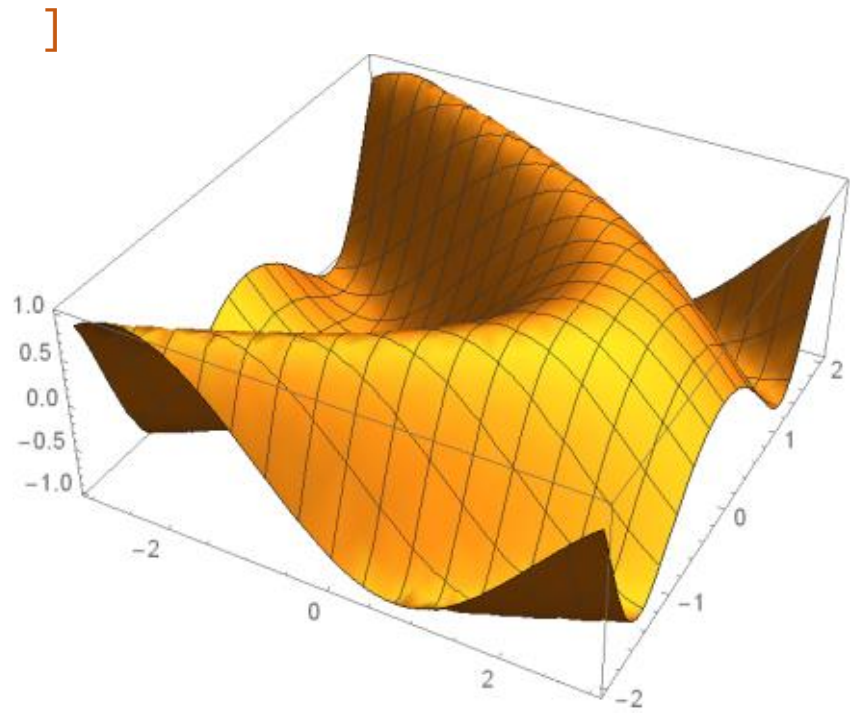

### Plotting Lists of Data

| | |
|---|---|
| `ListPlot[{y₁,y₂,...}]` | plot $y_1, y_2, \ldots$ at x values 1, 2, .... |
| `ListPlot[{{x₁,y₁},{x₂,y₂},...}]` | plot points $(x_1, y_1), \ldots$ |
| `ListLinePlot[list]` | join the points with lines. |
| `ListPlot3D[array]` | generates a three-dimensional plot of a surface representing an array of height values. |
| `ListPlot3D[{{x₁,y₁,z₁},{x₂,y₂,z₂},...}]` | generates a plot of the surface with heights $z_i$ at positions $x_i$, $y_i$. |
| ListPlot3D[{data1,data2,...}] | plots the surfaces corresponding to each of the data$_i$. |
| ListPointPlot3D[array] | generates a 3D scatter plot of points with a 2D array of height values. |
| ListPointPlot3D[{{$x_1$,$y_1$,$z_1$},{$x_2$,$y_2$,$z_2$},...}] | generates a 3D scatter plot of points with coordinates $x_i, y_i, z_i$ |
| ListPointPlot3D[{data1,data2,...}] | plots several collections of points, by default in different colors. |
| DensityPlot[f,{x,xmin,xmax]},{y,ymin,ymax}] | makes a density plot of f as a function of x and y. |
| ContourPlot[f,{x,xmin,xmax},{y,ymin,ymax}] | generates a contour plot of f as a function of x and y. |

***Mathematica Examples 1.32***

Input
```
ListPlot[
 Prime[Range[25]]
 ]
```
Output

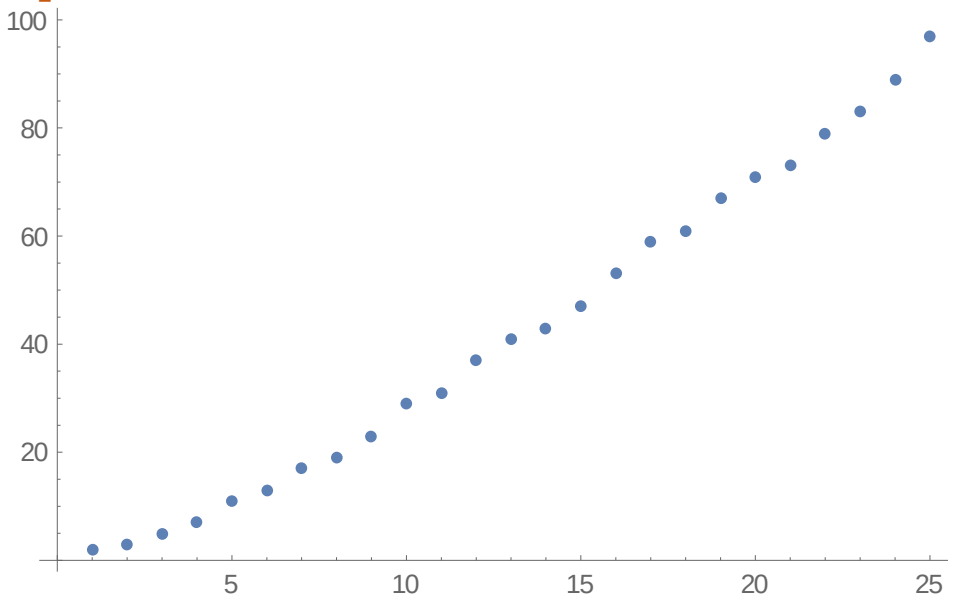

Input
```
ListPlot[
 Table[
  {Sin[n],Sin[2n]},
  {n,50}]
 ]
```
Output

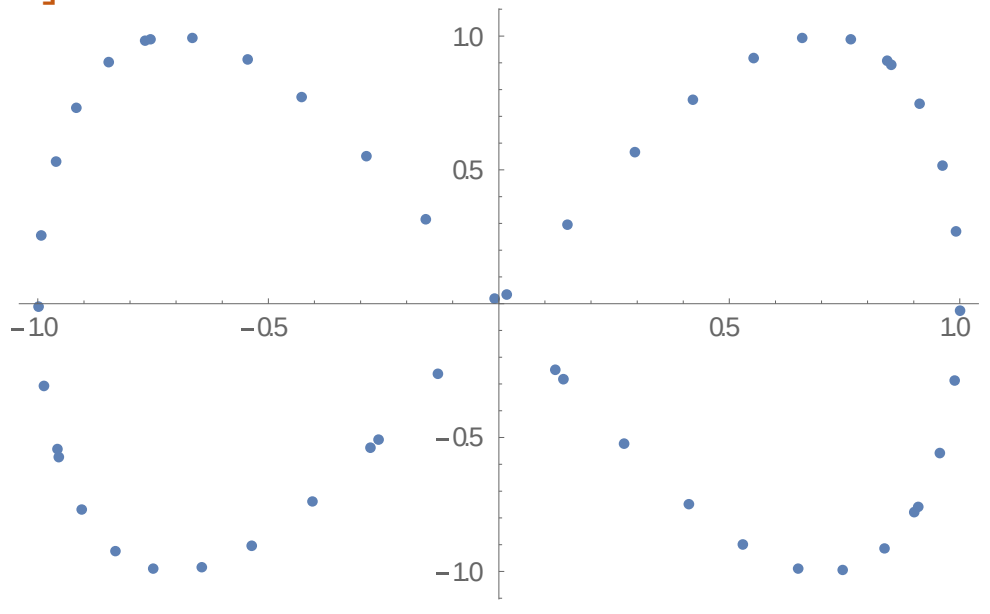

Input
```
ListLinePlot[
 {1,1,2,3,5,8}
 ]
```
Output

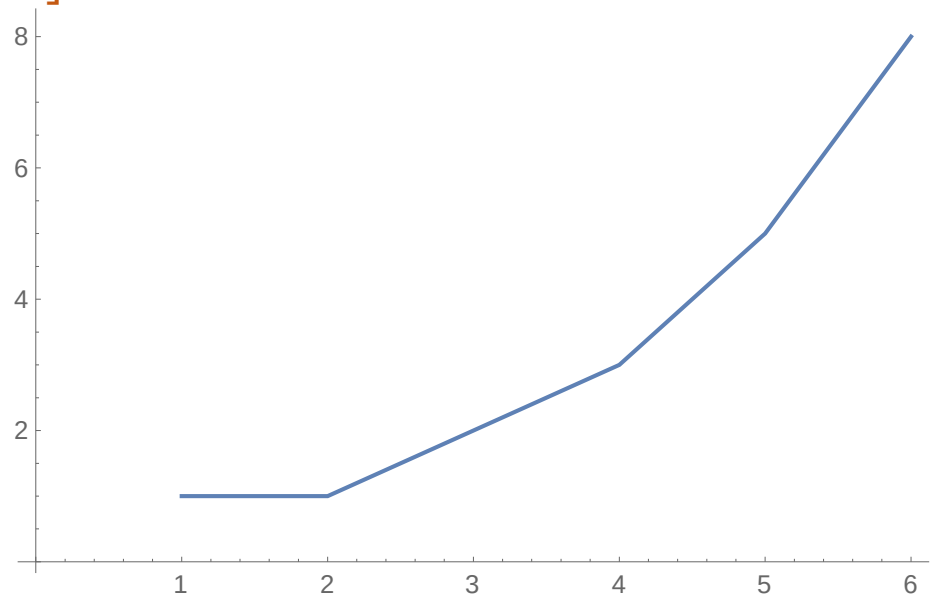

Input
```
ListPlot3D[
 {{1,1,1,1},{1,2,1,2},{1,1,3,1},{1,2,1,4}},
 Mesh->All
```

```
          ]
Output
```

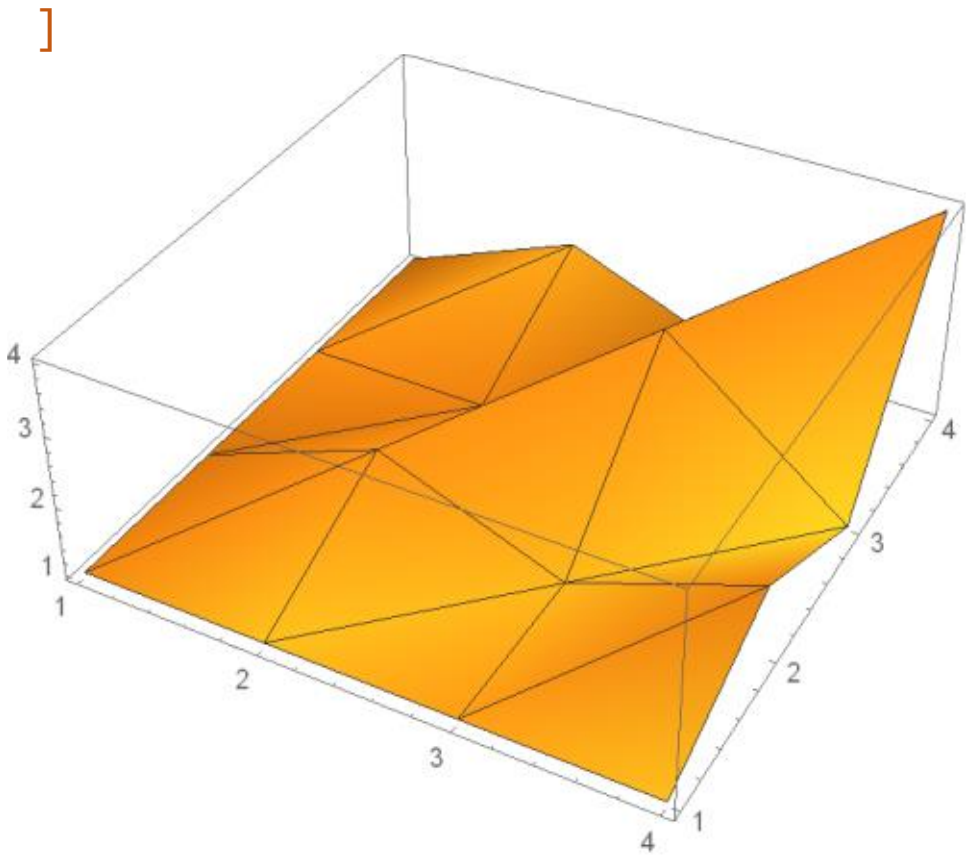

```
Input     data=Table[
             Sin[j^2+i],
             {i,0,Pi,Pi/5},
             {j,0,Pi,Pi/5}
             ];
          ListPlot3D[
           data,
           Mesh→None,
           InterpolationOrder→3,
           ColorFunction→"SouthwestColors"
           ]
Output
```

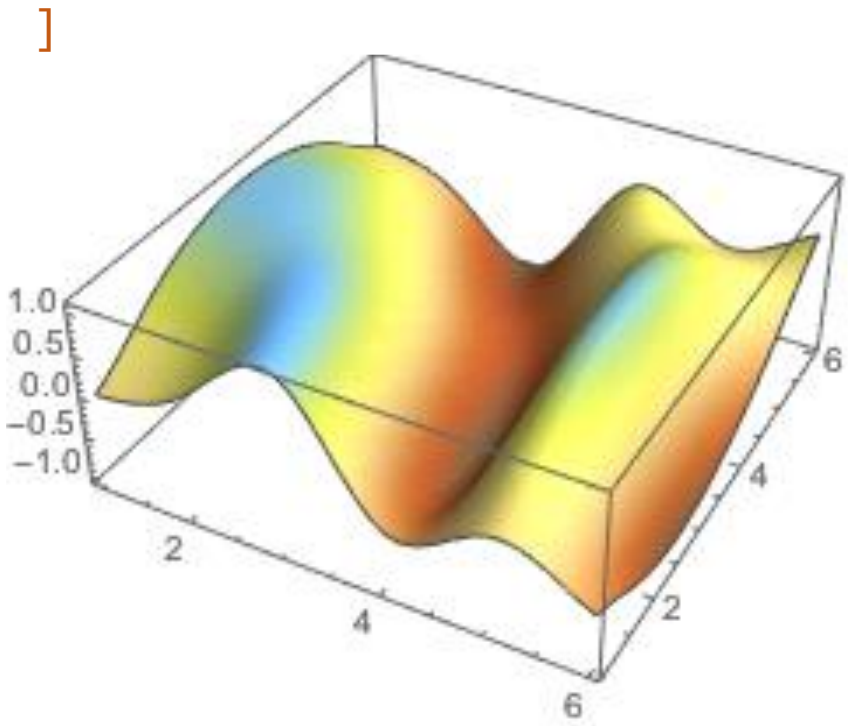

```
Input     ListPointPlot3D[
           Table[
            Sin[j^2+i],
            {i,0,3,0.1},
            {j,0,3,0.1}
            ]
           ]
Output
```

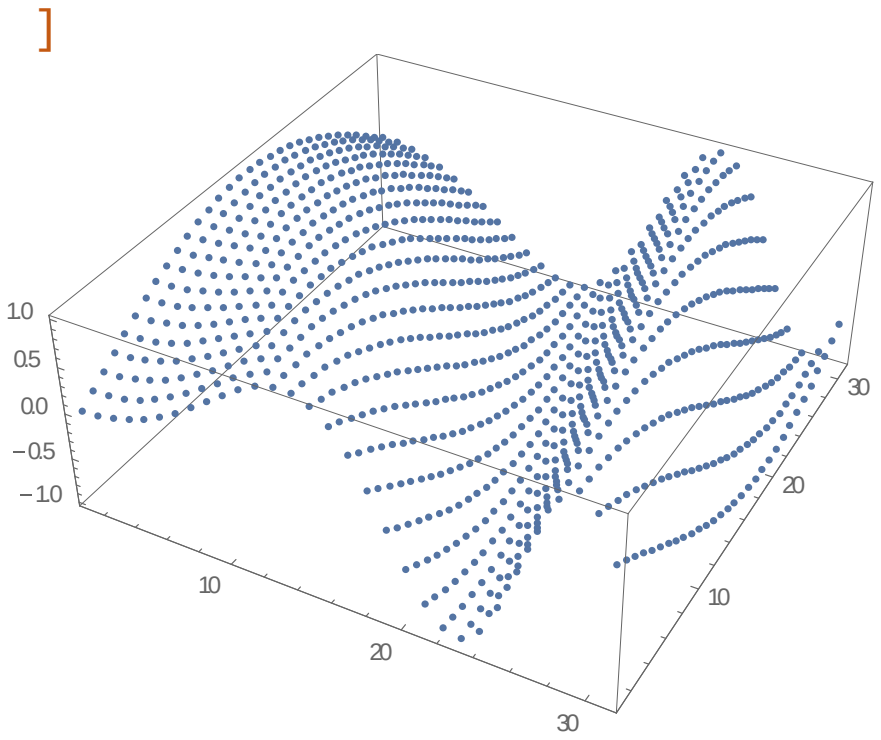

```
Input     ListPointPlot3D[
```

```
             {Table[
               Sin[j^2+i],
               {i,0,3,0.1},
               {j,0,3,0.1}
               ],
              Table[
               Sin[j^2+i]+3,
               {i,0,3,0.1},
               {j,0,3,0.1}
               ]
              }
             ]
Output
```

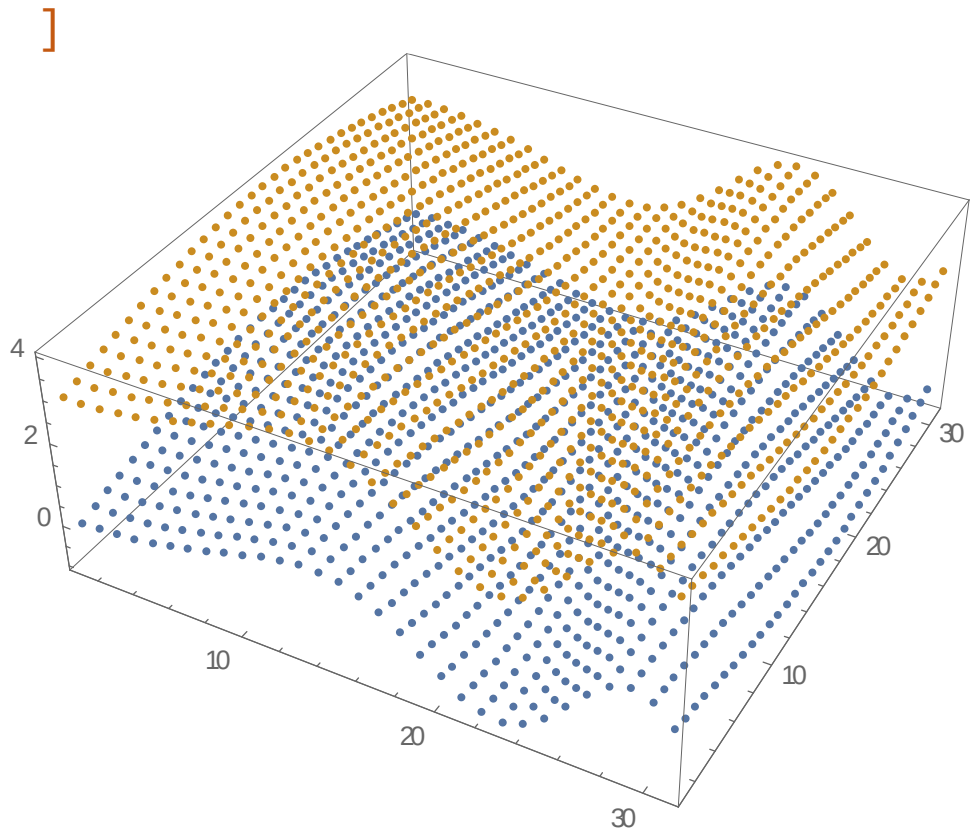

```
Input       DensityPlot[
             Sin[x] Sin[y],
             {x,-4,4},
             {y,-3,3},
             ColorFunction->"SunsetColors",
             PlotLegends->Automatic
             ]
Output
```

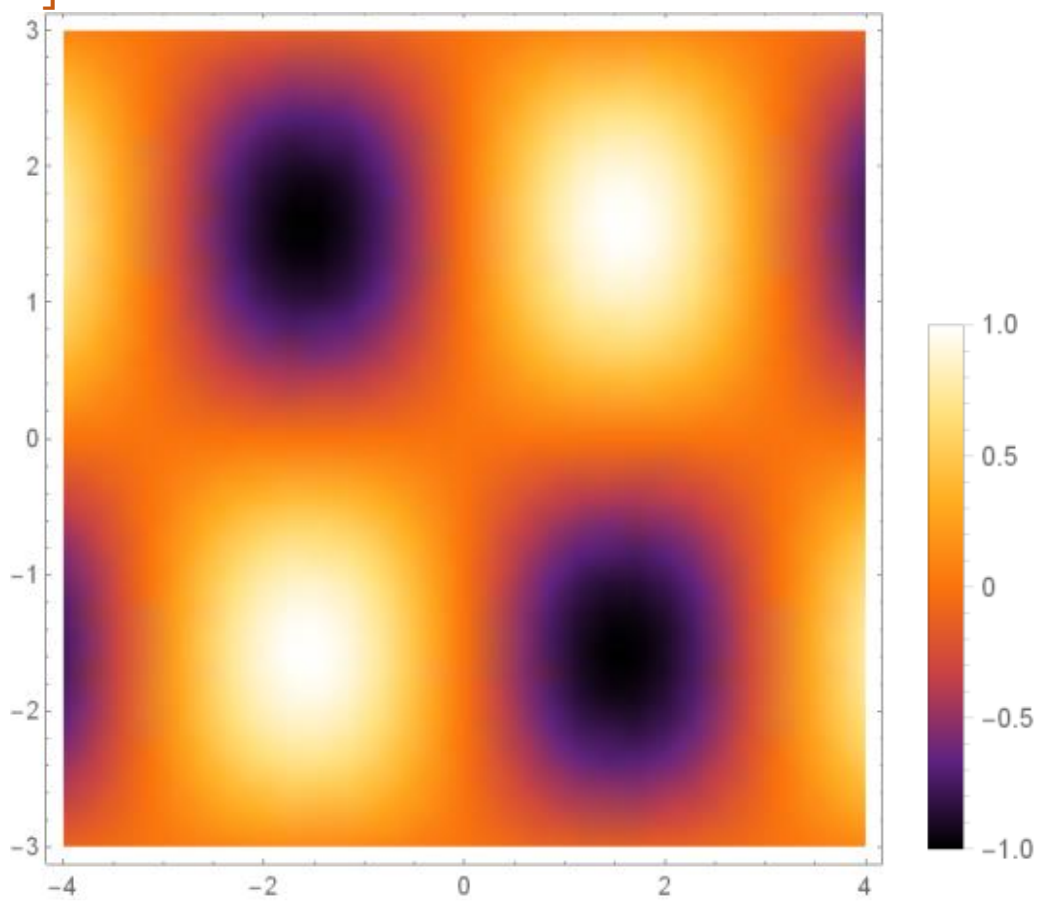

```
Input       ContourPlot[
             Cos[x]+Cos[y],
             {x,0,4 Pi},
             {y,0,4 Pi},
             PlotLegends->Automatic
             ]
```

Output

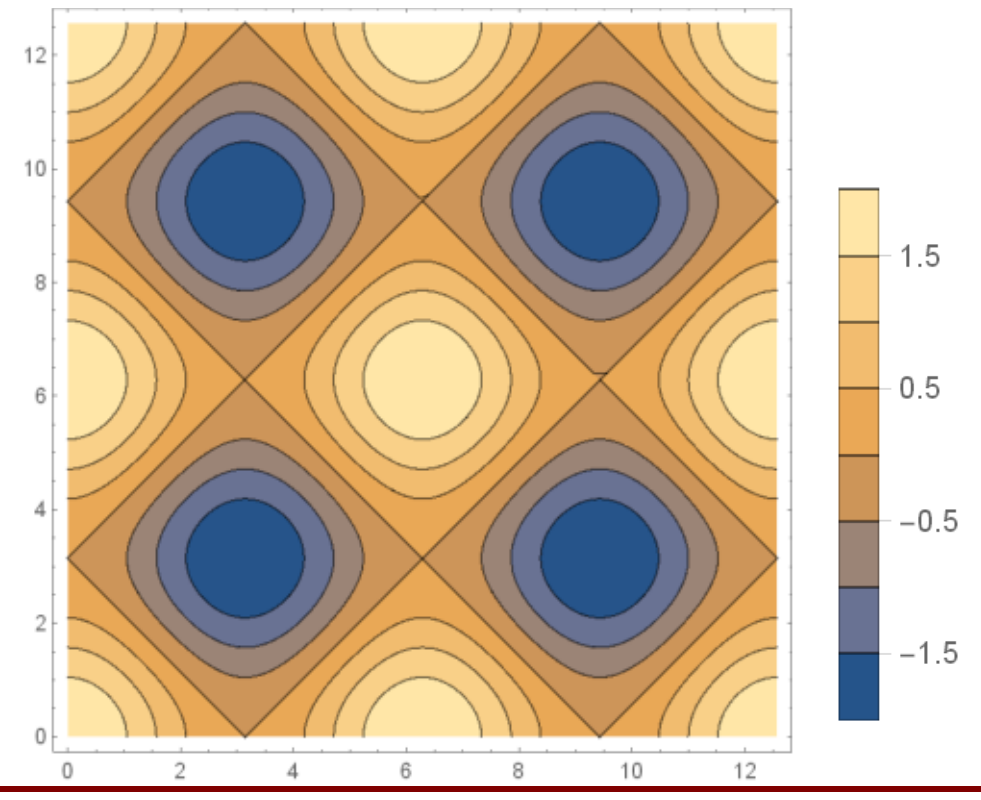

## Combining Plots

| | |
|---|---|
| `Show[plot₁,plot₂,...]` | combine several plots. |
| `GraphicsGrid[{plot₁,plot₂,...},...}]` | draw an array of plots. |
| `GraphicsRow[{plot₁,plot₂,...}]` | draw several plots side by side. |
| `GraphicsColumn[{plot₁,plot₂,...}]` | draw a column of plots. |

### *Mathematica Examples 1.33*

Input

```
Show[Plot[x^2,{x, 0, 3.5}], ListPlot[{1, 4, 9}]]
```

Output

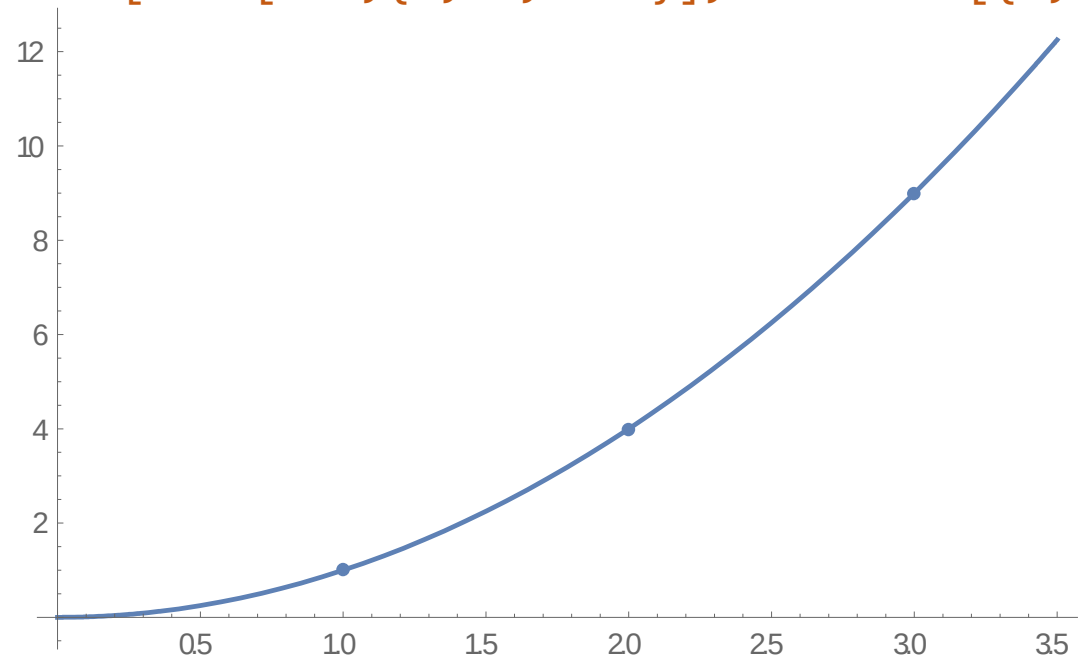

## Vector Field Plots

| | |
|---|---|
| `VectorPlot[{vx,vy },{x,xmin,xmax},{y,ymin,ymax }]` | generates a vector plot of the vector field $\{v_x, v_y\}$ as a function of x and y. |
| `VectorPlot3D[{vx,vy,vz},{x,xmin,xmax},{y,ymin,ymax },{z,zmin,zmax}]` | generates a 3D vector plot of the vector field $\{v_x, v_y, v_z\}$ as a function of x, y, and z. |
| `VectorDensityPlot[{{vx,vy},s},{x,xmin,xmax},{y,ymin,ymax}]` | generates a vector plot of the vector field $\{v_x, v_y\}$ as a function of x and y, superimposed on a density plot of the scalar field s. |

### *Mathematica Examples 1.34*

Input

```
VectorPlot[
 {x,-y},
 {x,-3,3},
 {y,-3,3}
 ]
```

Output

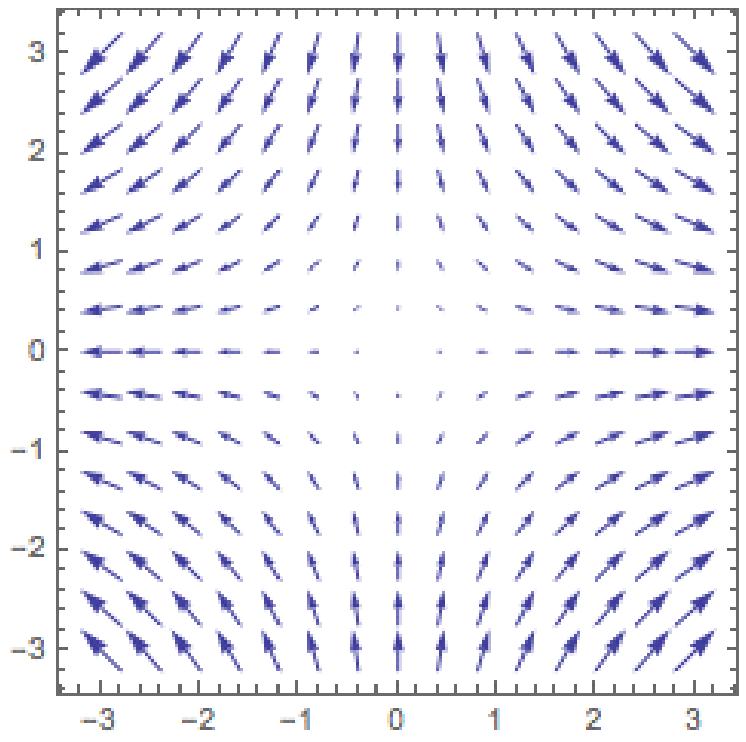

Input

```
VectorPlot3D[
 {x,y,z},
 {x,-1,1},
 {y,-1,1},
 {z,-1,1},
 VectorPoints->points,
 PlotRange->All,
 VectorScale->{Automatic,Scaled[0.5]},
 VectorColorFunction->Hue
 ]
```

Output

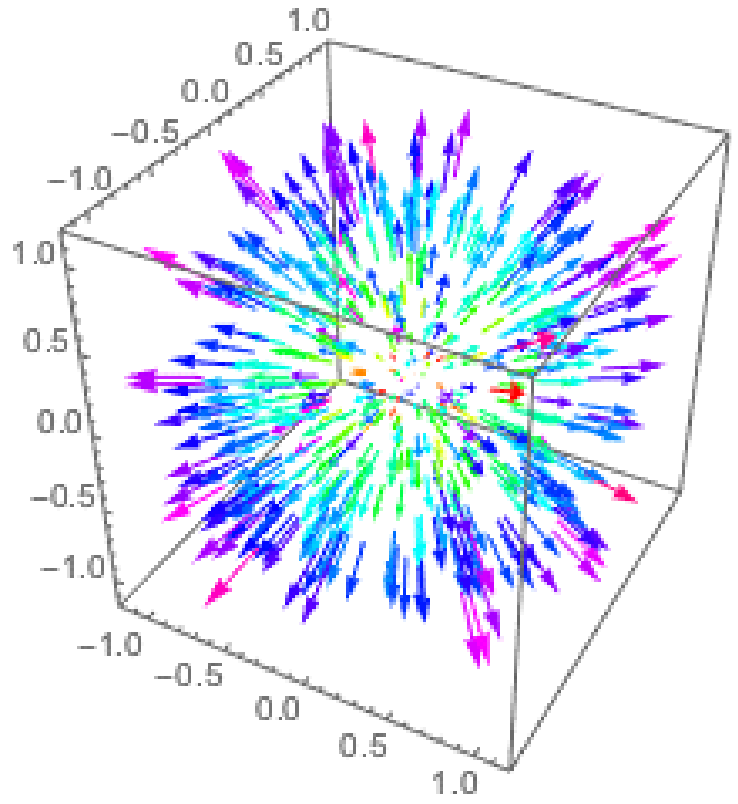

Input

```
VectorDensityPlot[
 {x,-y},
 {x,-3,3},
 {y,-3,3}
 ]
```

Output

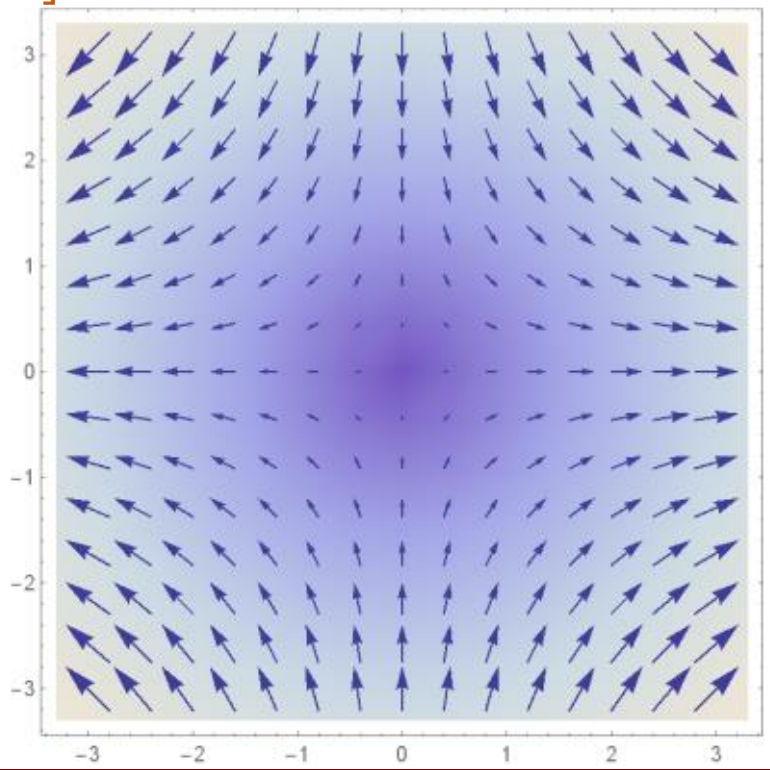

## 1.7 Manipulate

The single command `Manipulate` lets you create an astonishing range of interactive applications with just a few lines of input [7]. The output you get from evaluating a `Manipulate` command is an interactive object containing one or more controls (sliders, etc.) that you can use to vary the value of one or more parameters. The output is very much like a small applet or widget: it is not just a static result, it is a running program you can interact with.

| | |
|---|---|
| `Manipulate[expr,{u,umin,umax}]` | generates a version of expr with controls added to allow interactive manipulation of the value of u. |
| `Manipulate[expr,{u,umin,umax,du}]` | allows the value of u to vary between $u_{min}$ and $u_{max}$ in steps du. |
| `Manipulate[expr,{{u,uinit},umin,umax,...}]` | takes the initial value of u to be $u_{init}$. |

***Mathematica Examples 1.35***

Input
```
Manipulate[
 Plot[
  Sin[x (1+a x)],
  {x,0,6}
  ],
 {a,0,2}
 ]
```
Output

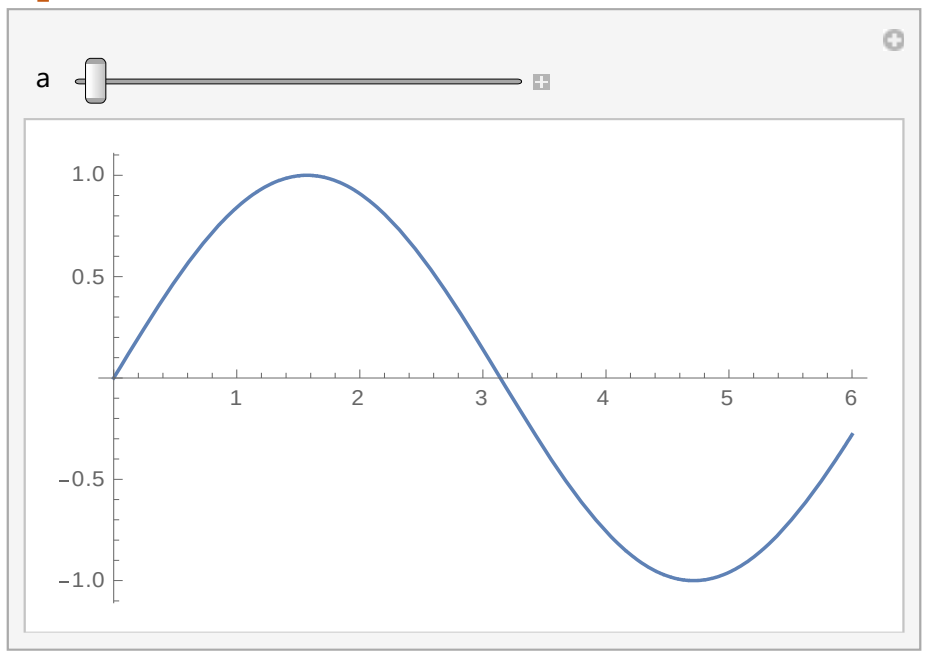

Input
```
Manipulate[
 Plot[
  Sin[a x+b],
  {x,0,6}
  ],
 {a,1,4},
 {b,0,10}
 ]
```
Output :

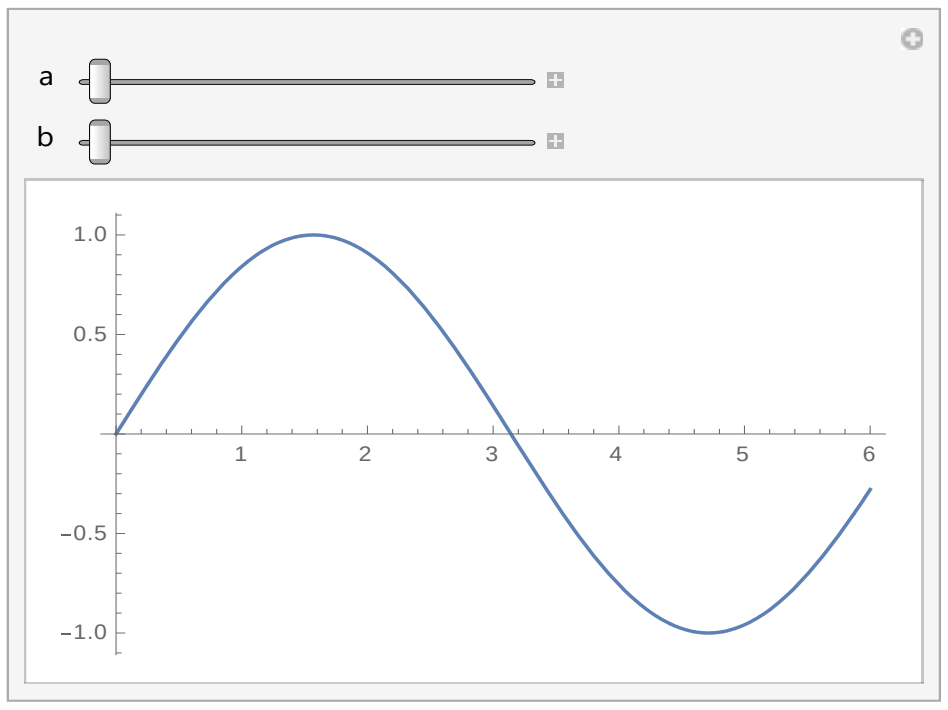

Input
```
Manipulate[
 ContourPlot3D[
```

```
  x^2+y^2+a z^3==1,
  {x,-2,2},
  {y,-2,2},
  {z,-2,2},
  Mesh->None
  ],
 {a,-2,2}
 ]
```

Output

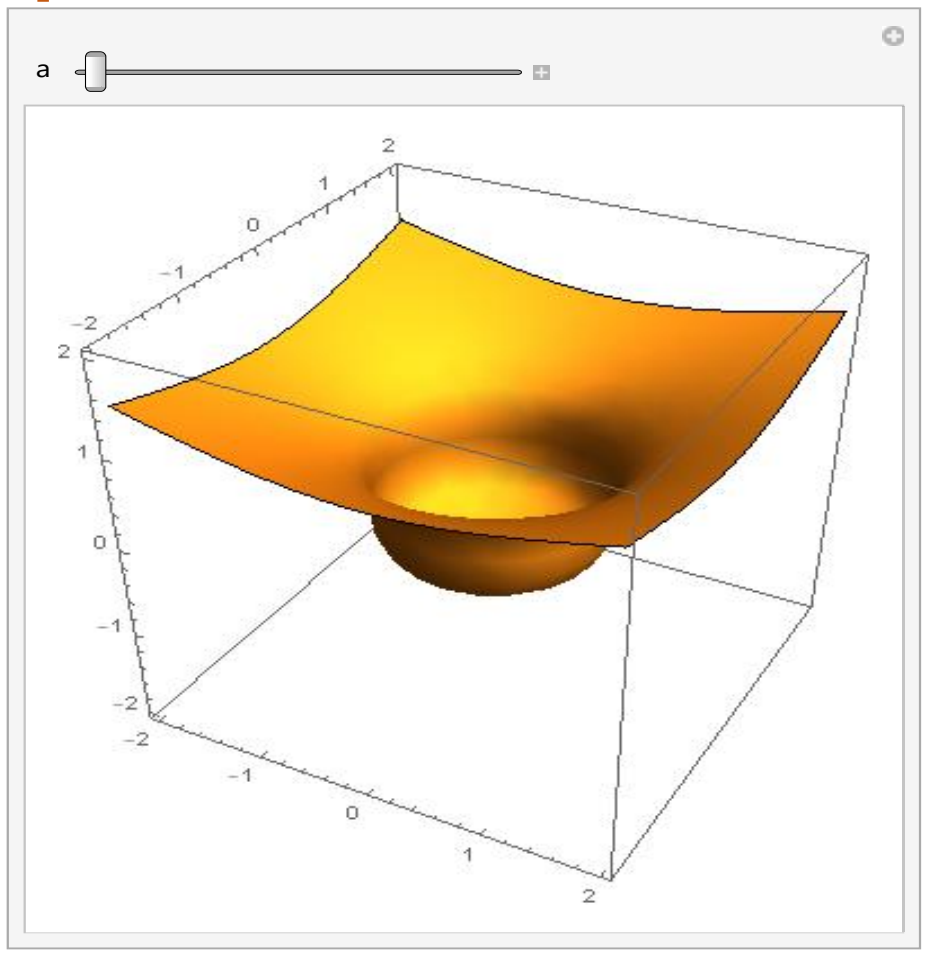

Input

```
Manipulate[
 Plot[
  If[
   t,
   Sin[x],
   Cos[x]
   ],
  {x,0,10}
  ],
 {t,{True,False}}
 ]
```

Output

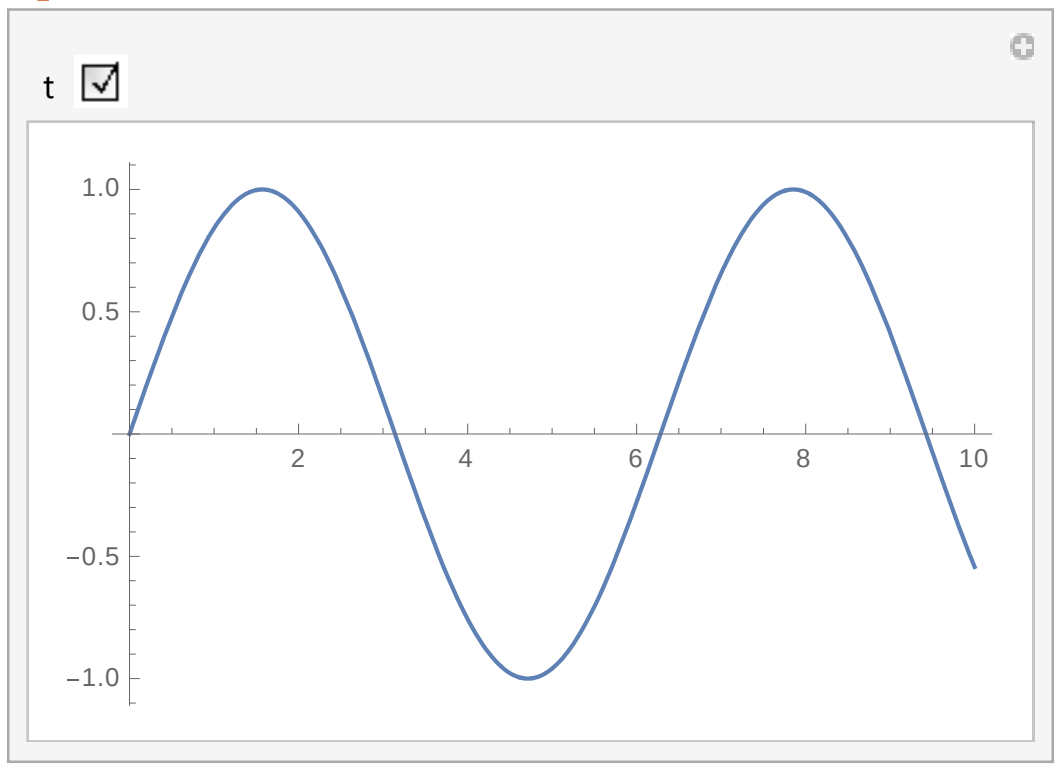

Input

```
Manipulate[
 Plot[
  f[x],
  {x,0,2 Pi}
  ],
 {f,{Sin->"sine",Cos->"cosine",Tan->"tangent"}}
 ]
```

Output

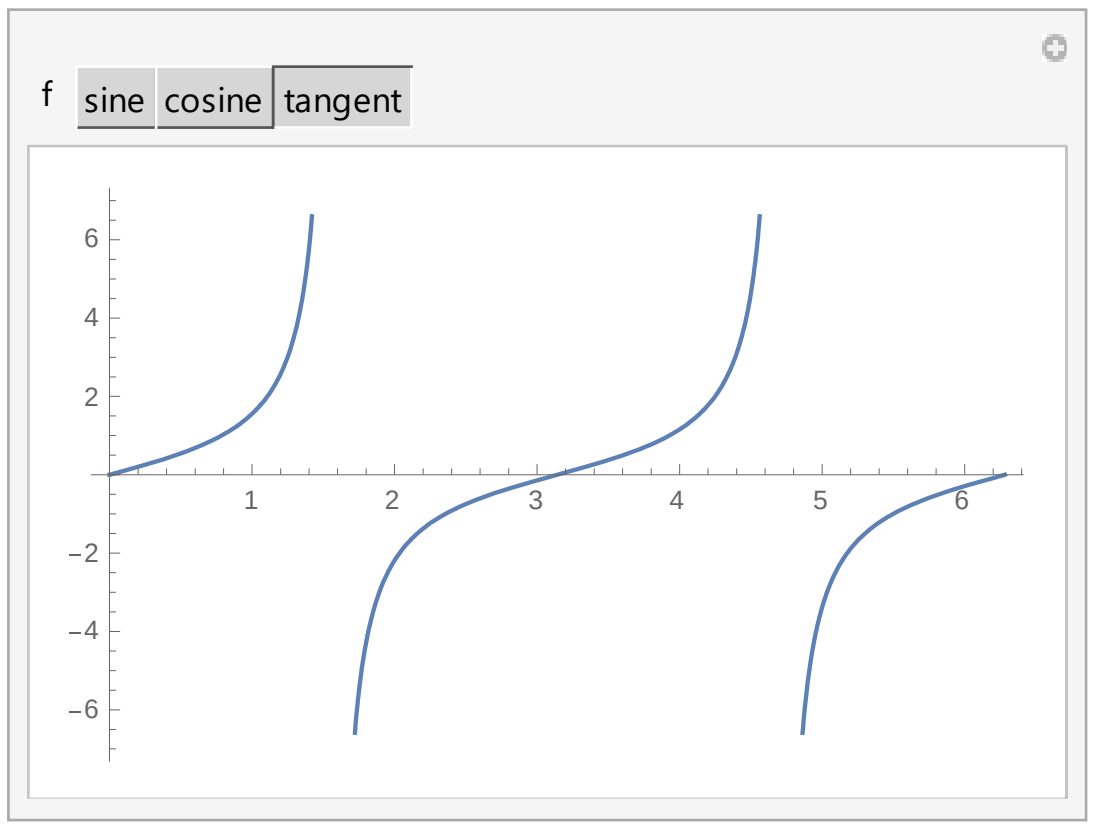

Input

```
Manipulate[
 ParametricPlot[
  {a1 Sin[n1 (x+p1)],a2 Cos[n2 (x+p2)]},
  {x,0,20 Pi},
  PlotRange->1,
  PerformanceGoal->"Quality"
  ],
 Style["Vertical",Bold,Medium],
 {{n1,1,"Frequency"},1,4},
 {{a1,1,"Amplitude"},0,1},
 {{p1,0,"Phase"},0,2 Pi},
 Delimiter,
 Style["Horizontal",Bold,Medium],
 {{n2,5/4,"Frequency"},1,4},
 {{a2,1,"Amplitude"},0,1},
 {{p2,0,"Phase"},0,2 Pi},
 ControlPlacement->Left
 ]
```

Output

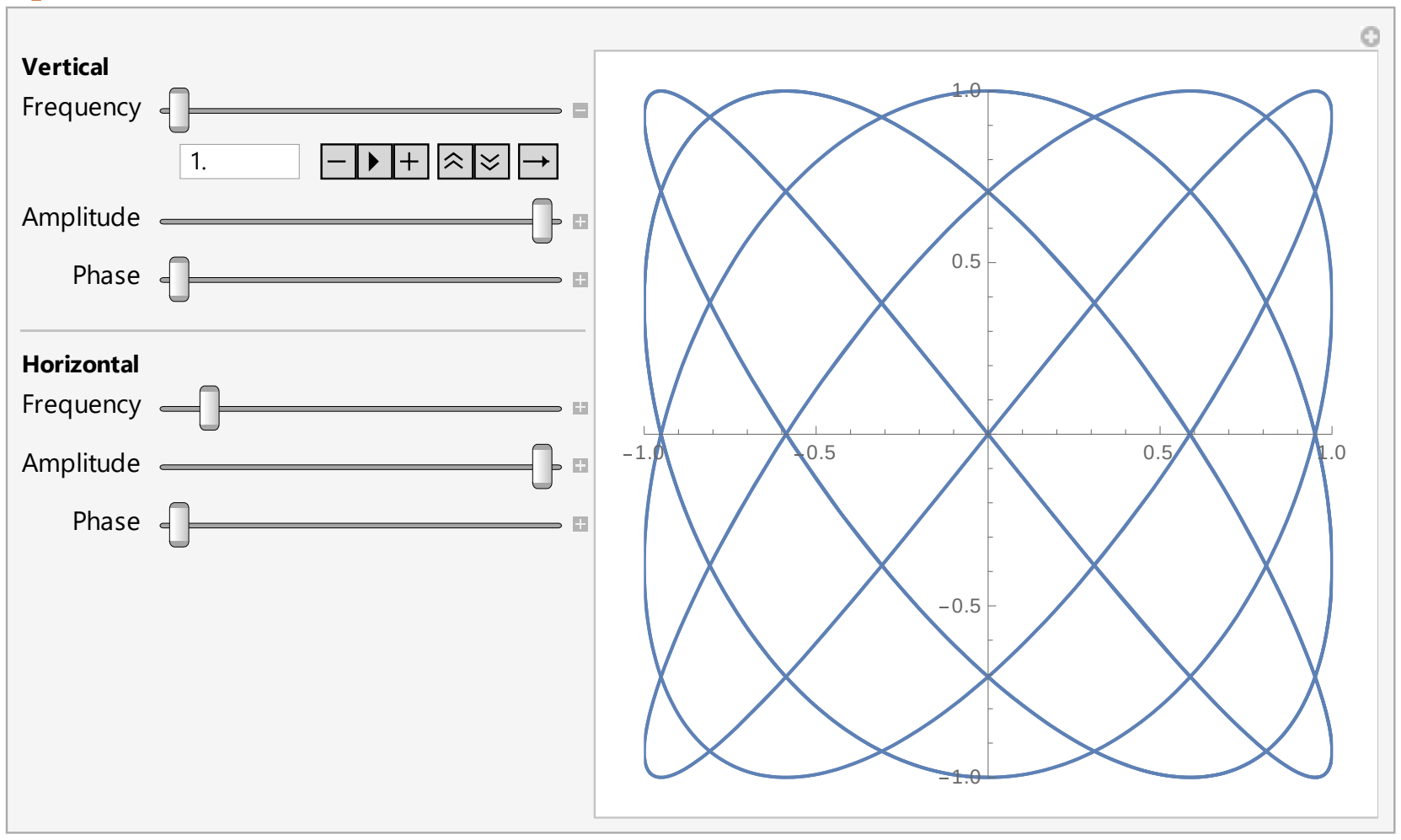

Input

```
Manipulate[
 Plot3D[
  Sin[x y+a],
  {x,0,3},
  {y,0,3}
  ],
 {a,0,1}
 ]
```

Output

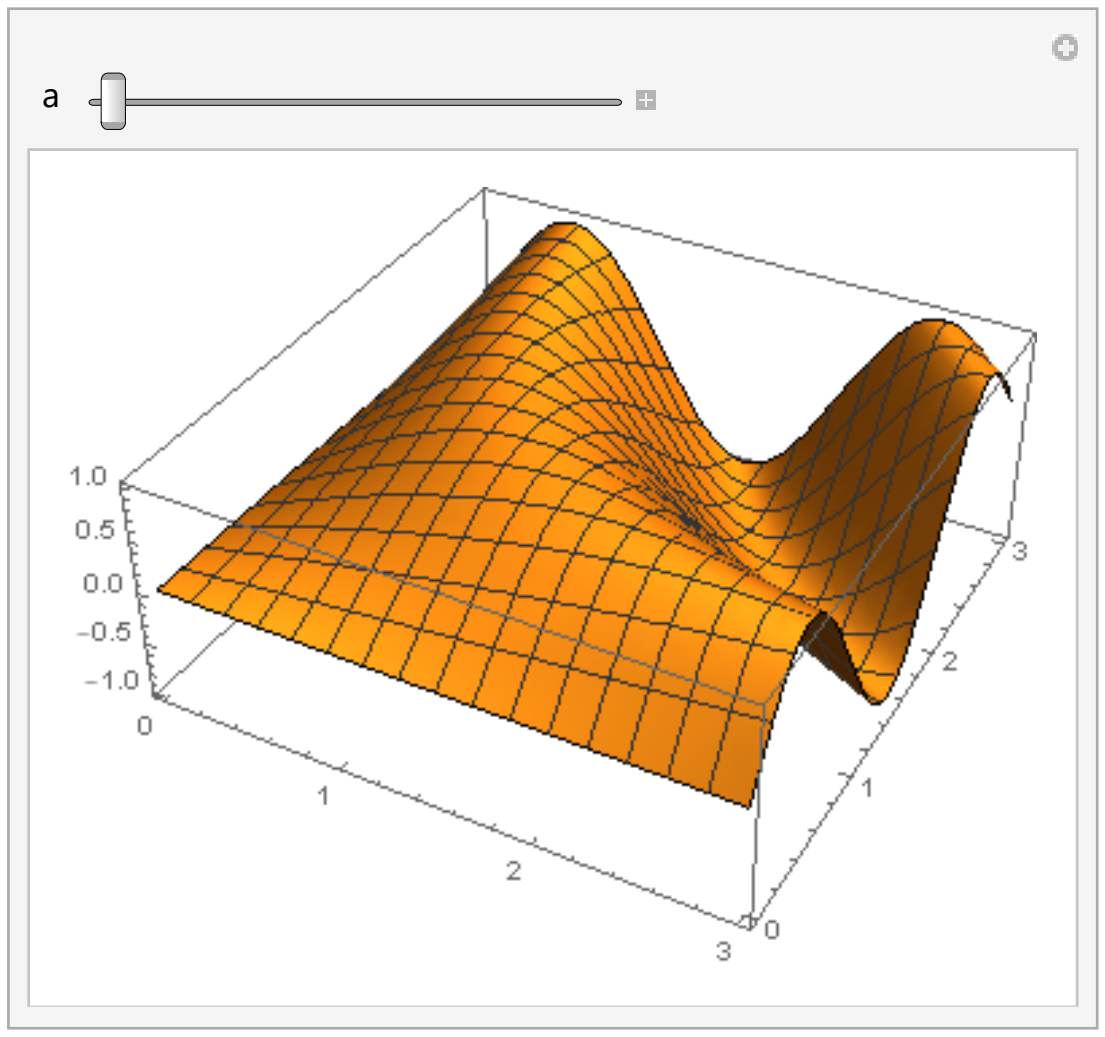

## 1.8 Control Structure

Most of the programming languages use control structures to control the flow of a program. The control structures include decision-making and loops [8]. Decision-making is done by applying different conditions in the program. If the conditions are true, the statements following the condition are executed. The values in a condition are compared by using the comparison operators. The loops are used to run a set of statements several times until a condition is met. If the condition is true, the loop is executed. If the condition becomes false, the loop is terminated, and the control passes to the next statement that follows the loop block.

### Conditional Statements

Programmers often need to check the status of a computed intermediate result to branch the program to such or another block of instructions to pursue the computation. Several examples of the branching condition structures [9] are next.

| | |
|---|---|
| `lhs:=rhs/;test` | is a definition to be used only if test yields True. |
| `If[test,then,else]` | evaluate then if test is True, and else if it is False. |
| `Which[test₁,value₁,test₂,...]` | evaluate the test$_1$ in turn, giving the value associated with the first one that is True. |
| `Switch[expr,form₁,value₁,form₂,...]` | compare expr with each of the form$_i$, giving the value associated with the first form it matches. |
| `Switch[expr,form₁,value₁,form₂,...,_,def]` | use def as a default value. |
| `Piecewise[{{value1,test1 },{value1,test1 },...}]` | represents a piecewise function with values value$_i$ in the regions defined by the conditions test$_i$. |
| `Piecewise[{{value₁,test₁},...},def]` | give the value corresponding to the first test$_1$ which yields True. |

Note that,

1- `If[condition, t, f]` is left unevaluated if condition evaluates to neither `True` nor `False`.
2- `If[condition, t]` gives `Null` if condition evaluates to `False`.

***Mathematica Examples 1.36***

Input
```
(* If can be used as a statement: *)
x=-2;
If[
  x<0,
  y=-x,
```

```
              y=x
              ];
              y
Output      2

Input       (* If can also be used as an expression returning a value: *)
            x=-2;
            y=If[
              x<0,
              -x,
              x]
Output      2

Input       If[
             7>8,
             x,
             y
             ]
Output      y

Input       x=2;
            If[
             x==0,
             Print["x is 0"],
             Print["x is different from 0"]
             ]
Output      x is different from 0

Input       x=3;
            y=0;
            If[
              x>1,
              y=Sqrt[x],
              y=x^2
              ];
            Print[y]
            m:=If[
               x>5,
               1,
               0
               ];
            Print[m]
Output      √3
Output      0

Input       a=2;
            Which[
             a==1,x,
             a==2,b
             ]
Output      b

Input       (* Use True for an else clause that always matches: *)
            sign[x_]:=Which[
              x<0,-1,
              x>0,1,
              True,Indeterminate
              ]
            {sign[-2],sign[0],sign[3]}
Output      {-1,Indeterminate,1}
```

```
Input       (* Use PiecewiseExpand to convert Which to Piecewise: *)
            PiecewiseExpand[
             Which[
              c1,a1,
              c2,a2,
              True,a3
              ]
             ]
Output          a1, c1
             {  a2, !c1&&c2
                a3, True

Input       expr=3;
            Switch[
             expr,
             1,Print["expr is 1"],
             2,Print["expr is 2"],
             3,Print["expr is 3"],
             _,Print["expr has some other value"]
             ]
Output      expr is 3

Input       k=2;
            n=0;
            Switch[
              k,
              1,n=k+10,
              2,n=k^2+3,
              _,n=-1];
            Print[n]
            k=5;
            n:=Switch[
               k,
               1,k+10,
               2,k^2+3,
               _,-1
               ];
            Print[n]
Output      7
Output      -1

Input       (* Find the derivative of a piecewise function: *)
            D[
             Piecewise[
              {
               {x^2,x<0},
               {x,x>0}}
              ],
             x
             ]
Output          2 x,             x<0
             {  1,               x>0
                Indeterminate,   True

Input       (* Define a piecewise function: *)
            pw=Piecewise[
              {
               {Sin[x]/x,x<0},
               {1,x==0}
               },
              -x^2/100+1
```

```
              ]
           (* Evaluate it at specific points: *)
           pw/. {{x->-5},{x->0},{x->5}}
Output
               Sin[x]/x,   x<0
            {  1,          x==0
               1-x2/100,   True
Output     {Sin[5]/5,1,3/4}

Input      Piecewise[
            {
             {x^2,x<0},
             {x,x>0}
             }
            ]
Output
               x2, x<0
            {  x, x>0
               0, True


Input      Piecewise[
            {
             {Sin[x]/x,x<0},
             {1,x==0}
             },
            -x^2/100+1
            ]
Output
               Sin[x]/x,   x<0
            {  1,          x==0
               1-x2/100,   True

Input      (* Remove unreachable cases: *)
           Piecewise[
            {
             {e1,d1},
             {e2,d2},
             {e3,True},
             {e4,d4},
             {e5,d5}},
            e
            ]
Output
               e1, d1
            {  e2, d2
               e3, True

Input      Piecewise[
            {
             {e1,d1},
             {e2,d2},
             {e3,d2&&d3},
             {e4,d4}
             },
            e
            ]
Output         e1, d1
            {  e2, d2
               e4, d4
               e, True
```

```
Input     x=4;
          If[
           x>0,
           y=Sqrt[x],
           y=0
           ]
Output    2

Input     (* PiecewiseExpand converts nested piecewise functions into a single piecewise
          function: *)
          pw=Piecewise[{{Piecewise[{{1,x>=0}},2],Piecewise[{{x,x<=1}},x/2]^2>=1/2}},3]
          PiecewiseExpand[pw]
Output    {
           {\[Piecewise], {
             {{
               {\[Piecewise], {
                 {1, x>=0},
                 {2, True}
                }}
              }, ({
                 {\[Piecewise], {
                   {x, x<=1},
                   {x/2, True}
                  }}
                })2>=1/2},
             {3, True}
            }}
          }
           {
           {\[Piecewise], {
             {1, x>=√2||1/√2<=x<=1},
             {2, x>=√2||1/√2<=x<=1||x<=-(1/√2)},
             {3, True}
            }}
          }

Input     (* Min, Max, UnitStep, and Clip are piecewise functions of real arguments: *)
          PiecewiseExpand/@{
            Min[x,y],
            Max[x,y,z],
            UnitStep[x],
            Clip[x,{a,b}]
            }
Output    {{
            {\[Piecewise], {
              {x, x-y<=0},
              {y, True}
             }}
           },{
            {\[Piecewise], {
              {x, x-y>=0&&x-z>=0},
              {y, x-y<0&&y-z>=0},
              {z, True}
             }}
           },{
            {\[Piecewise], {
              {1, x>=0},
              {0, True}
             }}
           },{
```

```
              {\[Piecewise], {
                {a, a-x>0},
                {b, b-x<0&&a-x<=0},
                {x, True}
               }}
             }}

Input      (* Abs, Sign, and Arg are piecewise functions when their arguments are assumed to
           be real: *)
           Assuming[
            Element[x,Reals],
            PiecewiseExpand/@{
              Abs[x],
              Sign[x],
              Arg[x]}
            ]
Output            -x, x<0   ,      -1, x<0   ,      π, x<0
            {  {  x, True       {  1, x>0        {  0, True   }
                                   0, True

Input      (* KroneckerDelta and DiscreteDelta are piecewise functions of complex arguments:
           *)
           PiecewiseExpand/@{
             KroneckerDelta[x,y],
             DiscreteDelta[x,y]
             }
Output            1, x-y==0   ,      1, x==0&&y==0
            {  {  0, True        {  0, True           }

Input      (*Derivatives are computed piece-by-piece,unless the function is univariate in a
           real variable:*)
           D[
            Piecewise[
             {
              {(x^2-1)/(x-1),x!=1}
              },
             2
             ],
            x
            ]
Output
            {  (2 x)/(-1+x)-(-1+x2)/(-1+x)2,   x!=1
               0,                              True

Input      x=-4;
           y=Which[
             x>0,1/x,
             x<-3,x^2,
             True,0
             ]
Output     16

Input      a=-4;
           b=4;
           y=Switch[
             a^2,
             a b,1.0/a,
             b^2,1.0/b,
             _,01
             ]
Output     0.25
```

### Looping Statements

Mathematica has several looping functions [10], the most common of which is `Do[ ]`.

| | |
|---|---|
| `Do [expr,n]` | evaluates expr n times. |
| `Do [expr,{i,imax}]` | evaluates expr with the variable i successively taking on the values 1 through $i_{max}$ (in steps of 1). |
| `Do [expr,{i,imin,imax}]` | starts with $i = i_{min}$ |
| `Do [expr,{i,imin,imax},di}]` | uses steps di. |
| `Do [expr,{i,{i₁,i₂,…}}]` | uses the successive values $i_1, i_2$, …. |
| `Do [expr,{i,imin,imax},{j,jmin,jmax},…]` | evaluates expr looping over different values of j etc. for each i. |

***Mathematica Examples 1.37***

```
Input     Do[
           Print[5^k],
           {k,3,7}
           ]
Output    125
          625
          3125
          15625
          78125

Input     t=x;
          Do[
           Print[t=1/(1+k t)],
           {k,2,6,2}
           ]
Output    1/(1+2 x)
          1/(1 + 4/(1 + 2 x))
          1/(1 + 6/(1 + 4/(1 + 2x)))

Input     Do[
           Print[{i,j}],
           {i,4},
           {j,i}
           ]
Output    {1,1}
          {2,1}
          {2,2}
          {3,1}
          {3,2}
          {3,3}
          {4,1}
          {4,2}
          {4,3}
          {4,4}

Input     sum=0;
          Do[
            Print[sum=sum+i],
            {i,1,4}
            ];
          sum
Output    1
          3
          6
          10
```

```
            10

Input       fact=1;
            Do[
              Print[fact=fact*i],
              {i,1,4}
              ];
            fact
Output      1
            2
            6
            24
            24


Input       Do[
             Do[
              Do[
               Print["i= ",i," j= ",j," k= ",k],
               {i,1,2}
               ],
              {j,1,2}
              ],
             {k,1,2}
             ]
Output      i= 1 j= 1 k= 1
            i= 2 j= 1 k= 1
            i= 1 j= 2 k= 1
            i= 2 j= 2 k= 1
            i= 1 j= 1 k= 2
            i= 2 j= 1 k= 2
            i= 1 j= 2 k= 2
            i= 2 j= 2 k= 2

Input       sum=0;
            Do[
              Print[i,",",sum=sum+i^2],
              {i,1,6,2}
              ];
            sum
Output      1,1
            3,10
            5,35
            35

Input       Do[
             Do[
              If[
               Sqrt[i^2+j^2]∈Integers,
               Print[i," ",j]
               ],
              {j,i,10}],
             {i,1,10}
             ]
Output      3   4
            6   8

Input       Do[
             Print[k!],
             {k,3}
             ]
```

```
Output     1
           2
           6

Input      Do[
            Print[k," ",k^2," ",k^3],
            {k,3}
            ]
Output     1   1   1
           2   4   8
           3   9   27

Input      Do[
            Print[k," squared is ",k^2],
            {k,5}
            ]
Output     1 squared is 1
           2 squared is 4
           3 squared is 9
           4 squared is 16
           5 squared is 25

Input      Print["k k^2"]
           Print["-----"]
           Do[
            Print[k," ",k^2],
            {k,5}
            ]
Output     k k^2
           -----
           1 1
           2 4
           3 9
           4 16
           5 25

Input      Do[
            Print[k],
            {k,1.6,5.7,1.2}
            ]
Output     1.6
           2.8
           4.
           5.2

Input      Do[
            Print[k],
            {k,3(a+b),8(a+b),2(a+b)}
            ]
Output     3 (a+b)
           5 (a+b)
           7 (a+b)

Input      (* The step can be negative: *)
           Do[
            Print[i],
            {i,10,8,-1}
            ]
Output     10
           9
           8
```

```
Input       (* The values can be symbolic: *)
            Do[
            Print[n],
            {n,x,x+3 y,y}
            ]
Output      x
            x+y
            x+2 y
            x+3 y

Input       (* Loop over i and j, with j running up to i-1: *)
            Do[
            Print[{i,j}],
            {i,4},
            {j,i-1}
            ]
Output      {2,1}
            {3,1}
            {3,2}
            {4,1}
            {4,2}
            {4,3}


Input      (* The body can be a procedure: *)
           t=67;
           Do[
            Print[t];
            t=Floor[t/2],
            {3}
            ]
Output     67
           33
           16
```

| | |
|---|---|
| `Nest[f,expr,n]` | apply f to expr n times. |
| `FixedPoint[f,expr]` | start with expr, and apply f repeatedly until the result no longer changes. |
| `NestList[f,expr,n]` | gives a list of the results of applying f to expr 0 through n times. |
| `While[test,body]` | evaluate body repetitively, so long as test is. |
| `For[start,test,incr,body]` | executes start, then repeatedly evaluates body and incr until test fails to give True. |
| `Break[]` | exits the nearest enclosing Do, For, or While. |

***Mathematica Examples 1.38***

```
Input      Nest[
            f, x, 3
            ]
Output     f[f[f[x]]]

Input      Nest[
            Function[t,1/(1+t)], x, 3
            ]
Output     1/(1 + 1/(1 + 1/(1 + x)))

Input      NestList[
            f, x, 4
            ]
```

```
Output     {x,f[x],f[f[x]],f[f[f[x]]],f[f[f[f[x]]]]}

Input      NestList[
            Cos, 1.0, 10
            ]
Output     {1., 0.540302, 0.857553, 0.65429, 0.79348,
           0.701369, 0.76396, 0.722102, 0.750418, 0.731404, 0.744237}

Input      FixedPoint[
            Function[t,Print[t];Floor[t/2]],
            67
            ]
Output     67
           33
           16
           8
           4
           2
           1
           0
           0

Input      n=17;
           While[
            n=Floor[n/2]; n!=0,
            Print[n]
            ]
Output     8
           4
           2
           1

Input      n=1;
           While[
            n<4,
            Print[n];
            n=n+1
            ]
Output     1
           2
           3

Input      Do[
            Print[i];
            If[
             i>2,
             Break[]
             ],
            {i,10}
            ]
Output     1
           2
           3

Input      For[
            i=1; t=x,
            i^2<10,
            i=i+1,
            t=t^2+i;
            Print[t]
            ]
```

```
Output     1+x²
           2+(1+x²)²
           3+(2+(1+x²)²)²

Input      For[
            sum=0.0; x=1.0,
            (1/x)>0.15,
            x=x+1,
            sum=sum+1/x;
            Print[sum]
            ]
Output     1.
           1.5
           1.83333
           2.08333
           2.28333
           2.45
```

## 1.9 Modules, Blocks, and Local Variables

Global Variables are those variables declared in Main Program and can be used by Subprograms. Local Variables are those variables declared in Subprograms. The Wolfram Language normally assumes that all your variables are global [11]. This means that every time you use a name like x, the Wolfram Language normally assumes that you are referring to the same object. Particularly when you write subprograms, however, you may not want all your variables to be global. You may, for example, want to use the name x to refer to two quite different variables in two different subprograms. In this case, you need the x in each subprogram to be treated as a local variable. You can set up local variables in the Wolfram Language using modules. Within each module, you can give a list of variables which are to be treated as local to the module.

| | |
|---|---|
| `Module[{x,y,...},body]` | a module with local variables x, y, .... |
| `Module[{x=x0,y=y0,…},body]` | a module with initial values for local variables. |

***Mathematica Examples 1.39***

```
Input      k=25
Output     25

Input      Module[
            {k},
            Do[
             Print[k," ",2^k],
             {k,3}
             ]
            ]
Output     1  2
           2  4
           3  8

Input      k
Output     25
```

Thus, we can create programs as a series of modules, each performing a specific task. For subtasks, we can imbed modules within other modules to form a hierarchy of operations. The most common method for setting up modules is through function definitions,

***Mathematica Examples 1.40***

```
Input        k=25;
             integerPowers[x_Integer]:=Module[
               {k},
               Do[
                Print[k," ",x^k],
                {k,3}
                ]
               ]
             integerPowers[k]
Output       1 25
             2 625
             3 15625

Input        k
Output       25

Input[1]     t = 17
Output[1]    17

Input[2]     Module[
              {t},
              t=8;Print[t]
              ]
Output[2]    8

Input[3]     t
Output[3]    17

Input[4]     g[u_]:=Module[
               {t=u},
               t+=t/(1+u)
               ]
Input[5]     g[a]
Output[5]    a + a/(1 + a)

Input[6]     h[x_]:=Module[
               {t},
               t^2-1/;(t=x-4)>1
               ]
Input[7]     h[10]
Output[7]    35
```

The format of `Module` is `Module[{var1, var2, ...}, body]`, where `var1`, `var2`, ... are the variables we localize, and body is the body of the function. The value returned by `Module` is the value returned by the last operator in the body (unless an explicit `Return[]` statement is used within the body of `Module`. In this case, the argument of `Return[arg]` is returned). In particular, if one places the semicolon after this last operator, nothing (`Null`) is returned. As a variant, it is acceptable to initialize the local variables in the place of the declaration, with some global values: `Module[{var1 = value1, var2, ...}, body]`. However, one local variable (say, the one "just initialized" cannot be used in the initialization of another local variable inside the declaration list. The following would be a mistake: `Module[{var1 = value1, var2 = var1, ...}, body]`. Moreover, this will not result in an error, but just the global value for the symbol `var1` would be used in this example for the `var2` initialization (this is even more dangerous since no error message is generated and thus we don't see the problem.) In this case, it would be better to do initialization in steps: `Module[{var1=value1,var2,...}, var2=var1;body]`, that is, include the initialization of part of the variables in the body of `Module`. One can use `Return[value]` statement to return a value from anywhere within the `Module`. In this case, the rest of the code (if any) inside `Module` is slipped, and the result value is returned.

To show how this is done, the following code is an example of module which will simulate a single gambler playing the game until the goal is achieved or the money is gone.

***Mathematica Examples 1.41***

```
Input       GamblersRuin[a_,c_,p_]:=Module[
              {ranval,var1,var2,var3},
              var1=a;
              var2=c;
              var3=p;
              While[
               0<var1<var2,
               ranval=Random[];
               If[
                ranval<var3,
                var1=var1+1,
                var1=var1-1
                ]
               ];
              Return[
               var1==var2
               ]
              ]
```

There are several things to notice in this example. First, this is the same thing we have done in the past to define a function. That is, we have a function name `GamblersRuin` with three input variables, a, c, and p. The operator `:=` is used to start the definition. Secondly the function involves the Mathematica command `Module`. This just tells Mathematica to perform all the commands in the module (like a subroutine in Fortran or method in C++). There are some special features we need to understand in the `Module` command. The `Module` command has two arguments. The first argument is a list of all the local variables that will only be used inside the module. In the above example the local variable list is `{ranval,var1,var2,var3}`. These variables are only used in the module and are cleared once the module has executed. The second argument is all the commands that will be executed each time the module is called. There are some assignment commands at the beginning that are used to make things cleaner. The module uses temporary variables so that the values of the input variables are not overwritten when the module executes.

The last command is added to our list of input lines to return a result from the work done by the `Module`. Without this we would never get any results from out calculation. Any recognized variable type or structure within Mathematica can be returned by a `Module`. In this example, the returned value is the result of testing two variables in the code for equality. The code fragment `var1==var2` tests to determine if the variables, `var1` and `var2`, are equal. If the two variables are equal, then the line outputs `True` and if they are not equal, the line outputs `False`.

Again, all but the last command must be ended by a semicolon. This is to make sure that the commands are separated in the execution. Commands separated by blank space will be considered as terms to be multiplied together. Leaving out the semicolon will give rise to lots of error messages, wrong results or both.

Modules in Mathematica allow one to treat the names as local. When one uses `Block` then the names are global, but the values are local.

| | |
|---|---|
| `Block[{x,y,...},body]` | evaluate expr using local values for x,y, .... |
| `Block[{x=x0,y=y0,...},body]` | assign initial values to x,y; and evaluate ... as above. |

`Block[]` is automatically used to localize values of iterators in iteration constructs such as `Do`, `Sum` and `Table`.
`Block[]` may be used to pack several expressions into one unity.

***Mathematica Examples 1.42***

```
Input       Clear[x,t,a]
Input       x^2 + 3
Output      3+x2

Input       Block[ {x = a + 1}, %]
Output      3+(1+a)2

Input       x
Output      x

Input       t = 17
Output      17

Input       Module[ {t}, Print[t]]
Output      t$505

Input       t
Output      17

Input       Block[ {t}, Print[t] ]
Output      t
```

## 1.10 Requesting Input and Choice Dialog

The Wolfram Language usually works by taking whatever input you give, and then processing it. Sometimes, however, you may want to have a program you write explicitly request more input. You can do this using `Input` and `InputString`. Dialogs that display a choice can be created using `ChoiceDialog`, which blocks the kernel and returns `False` if the `CancelButton` is selected or `True` if the `DefaultButton` is selected.

| | |
|---|---|
| `Input[]` | interactively reads in one Wolfram Language expression. |
| `Input[prompt]` | requests input, displaying prompt as a "prompt". |
| `InputString[]` | interactively reads in a character string. |
| `InputString[prompt]` | requests input, displaying prompt as a "prompt". |
| `ChoiceDialog[expr]` | puts up a standard choice dialog that displays expr together with OK and Cancel buttons, and returns True if OK is clicked and False if Cancel is clicked. |
| `ChoiceDialog[expr,{lbl1->val1,lbl2->val2,…}]` | includes buttons with labels lbli, and returns the corresponding vali for the button clicked. |
| `Beep[]` | generates an audible beep when evaluated. |

***Mathematica Examples 1.43***

```
Input       ChoiceDialog["Make a choice..."]
Output      True

Input       ChoiceDialog["Pick one",{1->"a",2->"b",3->"c",4->"d"}]
Output      c

Input       Input[]^2
Output      4

Input       InputString[]
Output      Mohamed

Input       Plot[
             Evaluate[
```

```
      Input["Please input a function of x to plot on the unit interval"]
      ],
     {x,0,1}
    ]
```

Output

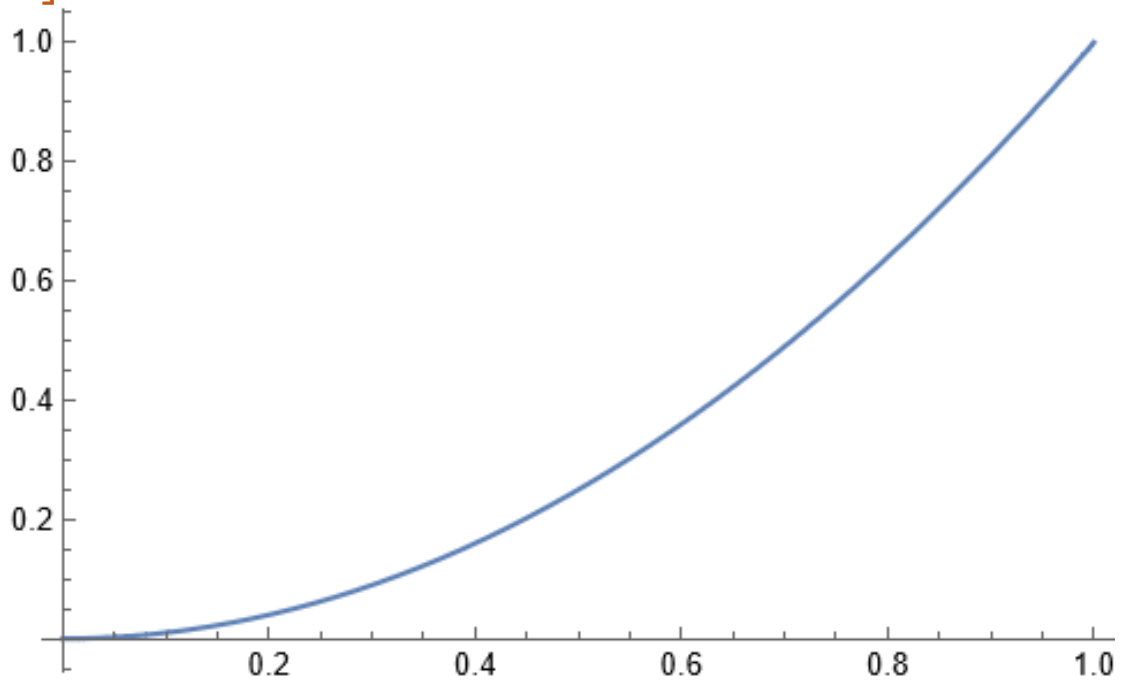

# CHAPTER 2

# MATRIX NOTATIONS, DECOMPOSITION, AND CALCULUS

In this chapter we collect five introductory topics needed for later use about matrices. we feel these topics are essential.

## 2.1 Vectors and Matrices with Dirac Notations

Throughout these lectures, we will distinguish scalars, vectors, and matrices by their typeface. We will let $M_{n\times m}$ denote the space of all real $n\times m$ matrices with $n$ rows and $m$ columns. Such matrices will be denoted using bold capital letters: **A**, **X**, **Y**, etc. An element of $M_{n\times 1}$ or $M_{1\times n}$, that is, a column or row vector, respectively, is denoted with a boldface lowercase letter: **a**, **x**, **y**, etc. An element of matrix **M** is a scalar, denoted with italic letter: $A$, $X$, $Y$, $a$, $x$, $y$, etc.

Moreover, we have chosen to use a Dirac matrix notation [1]. This choice was not made lightly. We are a strong advocate of index notation when appropriate. For example, index notation greatly simplifies the presentation and manipulation of differential geometry. As a rule-of-thumb, if your work primarily involves differentiation with respect to the spatial coordinates, then index notation is almost the appropriate choice. In the present case, however, we will be manipulating systems of equations in which matrix calculus is relatively simple, while matrix algebra and matrix arithmetic are messy and more involved. Thus, we have chosen to use the Dirac matrix notation.

The entities that Dirac called "kets" and "bras" are simply column vectors and row vectors, respectively. Of course, the elements of these vectors are generally complex numbers. In these lectures, for convenience, we will express ourselves in terms of vectors and matrices of real numbers. Hence, in the language of matrices, these two vectors are related by simply taking the transpose. In summary, Dirac refers to a "bra," which he denoted as $\langle\mathbf{a}|$, a "ket," which he denoted as $|\mathbf{b}\rangle$, and a square matrix **M**, we can associate these with vectors and matrices (in 3 dimensions) as follows

$$\langle\mathbf{a}| = \mathbf{a}^T = (a_1, a_2, a_3), \quad |\mathbf{b}\rangle = \mathbf{b} = \begin{pmatrix} b_1 \\ b_2 \\ b_3 \end{pmatrix}, \quad \mathbf{M} = \begin{pmatrix} M_{11} & M_{12} & M_{13} \\ M_{21} & M_{22} & M_{23} \\ M_{31} & M_{32} & M_{33} \end{pmatrix}. \tag{2.1}$$

The product of a bra and a ket, denoted by Dirac as $\langle\mathbf{a}||\mathbf{b}\rangle$ or, more commonly, by omitting one of the middle lines, as $\langle\mathbf{a}|\mathbf{b}\rangle$, is simply a number given by inner products of a row vector and a column vector in the usual way, i.e.,

$$\begin{aligned}\langle\mathbf{a}|\mathbf{b}\rangle &= \mathbf{a}^T\mathbf{b} \\ &= (a_1, a_2, a_3)\begin{pmatrix} b_1 \\ b_2 \\ b_3 \end{pmatrix} \\ &= a_1b_1 + a_2b_2 + a_3b_3 \\ &= \sum_{i=1}^{3} a_ib_i.\end{aligned} \tag{2.2}$$

We can also form the product of ket times a bra, which gives a square matrix, as shown below,

$$\begin{aligned}|\mathbf{b}\rangle\langle\mathbf{a}| &= \mathbf{b}\mathbf{a}^T \\ &= \begin{pmatrix} b_1 \\ b_2 \\ b_3 \end{pmatrix}(a_1, a_2, a_3) = \begin{pmatrix} b_1a_1 & b_1a_2 & b_1a_3 \\ b_2a_1 & b_2a_2 & b_2a_3 \\ b_3a_1 & b_3a_2 & b_3a_3 \end{pmatrix}.\end{aligned} \tag{2.3}$$

The product of square matrix times a ket corresponds to the product of square matrix times a column vector, yielding another column vector (i.e., a ket) as follows (using row-column multiplication or inner products of the rows with column)

$$\mathbf{M}|\mathbf{b}\rangle = \mathbf{Mb} = \begin{pmatrix} (M_{11} & M_{12} & M_{13}) \\ (M_{21} & M_{22} & M_{23}) \\ (M_{31} & M_{32} & M_{33}) \end{pmatrix} \begin{pmatrix} b_1 \\ b_2 \\ b_3 \end{pmatrix}$$
$$= \begin{pmatrix} M_{11}b_1 + M_{12}b_2 + M_{13}b_3 \\ M_{21}b_1 + M_{22}b_2 + M_{23}b_3 \\ M_{31}b_1 + M_{32}b_2 + M_{33}b_3 \end{pmatrix}$$
$$= \begin{pmatrix} \sum_{i=1}^{3} M_{1i}b_i \\ \sum_{i=1}^{3} M_{2i}b_i \\ \sum_{i=1}^{3} M_{3i}b_i \end{pmatrix}$$
$$= \begin{pmatrix} \langle \mathbf{m}_{(1)} | \mathbf{b}^{(1)} \rangle \\ \langle \mathbf{m}_{(2)} | \mathbf{b}^{(1)} \rangle \\ \langle \mathbf{m}_{(3)} | \mathbf{b}^{(1)} \rangle \end{pmatrix}, \tag{2.4}$$

where,

$$\langle \mathbf{m}_{(1)} | = (M_{11} \quad M_{12} \quad M_{13}), \qquad \langle \mathbf{m}_{(2)} | = (M_{21} \quad M_{22} \quad M_{23}), \qquad \langle \mathbf{m}_{(3)} | = (M_{31} \quad M_{32} \quad M_{33}), \tag{2.5}$$

$$|\mathbf{b}^{(1)}\rangle = \begin{pmatrix} b_1 \\ b_2 \\ b_3 \end{pmatrix}. \tag{2.6}$$

However, we also have column-row multiplication or (i.e., a linear combination of the columns using $b_i$)

$$\mathbf{M}|\mathbf{b}\rangle = \mathbf{Mb} = \left( \begin{pmatrix} M_{11} \\ M_{21} \\ M_{31} \end{pmatrix} \begin{pmatrix} M_{12} \\ M_{22} \\ M_{32} \end{pmatrix} \begin{pmatrix} M_{13} \\ M_{23} \\ M_{33} \end{pmatrix} \right) \begin{pmatrix} b_1 \\ b_2 \\ b_3 \end{pmatrix}$$
$$= \begin{pmatrix} M_{11} \\ M_{21} \\ M_{31} \end{pmatrix} b_1 + \begin{pmatrix} M_{12} \\ M_{22} \\ M_{32} \end{pmatrix} b_2 + \begin{pmatrix} M_{13} \\ M_{23} \\ M_{33} \end{pmatrix} b_3$$
$$= \begin{pmatrix} M_{11}b_1 \\ M_{21}b_1 \\ M_{31}b_1 \end{pmatrix} + \begin{pmatrix} M_{12}b_2 \\ M_{22}b_2 \\ M_{32}b_2 \end{pmatrix} + \begin{pmatrix} M_{13}b_3 \\ M_{23}b_3 \\ M_{33}b_3 \end{pmatrix}$$
$$= \begin{pmatrix} M_{11}b_1 + M_{12}b_2 + M_{13}b_3 \\ M_{21}b_1 + M_{22}b_2 + M_{23}b_3 \\ M_{31}b_1 + M_{32}b_2 + M_{33}b_3 \end{pmatrix}$$
$$= |\mathbf{m}^{(1)}\rangle\langle \mathbf{b}_{(1)}| + |\mathbf{m}^{(2)}\rangle\langle \mathbf{b}_{(2)}| + |\mathbf{m}^{(3)}\rangle\langle \mathbf{b}_{(3)}|$$
$$= \sum_{i=1}^{3} |\mathbf{m}^{(i)}\rangle\langle \mathbf{b}_{(i)}|, \tag{2.7}$$

where,

$$|\mathbf{m}^{(1)}\rangle = \begin{pmatrix} M_{11} \\ M_{21} \\ M_{31} \end{pmatrix}, \qquad |\mathbf{m}^{(2)}\rangle = \begin{pmatrix} M_{12} \\ M_{22} \\ M_{32} \end{pmatrix}, \qquad |\mathbf{m}^{(3)}\rangle = \begin{pmatrix} M_{13} \\ M_{23} \\ M_{33} \end{pmatrix}, \tag{2.8}$$

$$\langle \mathbf{b}_{(1)} | = (b_1), \qquad \langle \mathbf{b}_{(2)} | = (b_2), \qquad \langle \mathbf{b}_{(3)} | = (b_3). \tag{2.9}$$

Thus, $\mathbf{M}|\mathbf{b}\rangle$ is a linear combination of the columns of $\mathbf{M}$ [2]. This is fundamental. This thinking leads us to the column space of the matrix $\mathbf{M}$ and the idea of the rank of the matrix (as we will see in the next chapter). The key idea is to take all combinations of the columns. All real numbers $b_i$ are allowed - the space includes $\mathbf{M}|\mathbf{b}\rangle$ for all vectors $|\mathbf{b}\rangle$. In this way, we get infinitely many output vectors $\mathbf{M}|\mathbf{b}\rangle$.

The product of bra times a square matrix corresponds to the product of a row vector times a square matrix, which is again a row vector (i.e., a "bra"), written as

$$\begin{aligned}\langle \mathbf{a}|\mathbf{M} &= \mathbf{a}^T\mathbf{M}\\ &= (a_1, a_2, a_3)\begin{pmatrix} M_{11} & M_{12} & M_{13}\\ M_{21} & M_{22} & M_{23}\\ M_{31} & M_{32} & M_{33}\end{pmatrix}\\ &= \left(\sum_{1=1}^{3} a_i M_{i1}, \sum_{1=1}^{3} a_i M_{i2}, \sum_{1=1}^{3} a_i M_{i3}\right),\end{aligned} \tag{2.10}$$

or

$$\begin{aligned}\langle \mathbf{a}|\mathbf{M} = \mathbf{a}^T\mathbf{M} = (a_1, a_2, a_3)&\begin{pmatrix} (M_{11} & M_{12} & M_{13})\\ (M_{21} & M_{22} & M_{23})\\ (M_{31} & M_{32} & M_{33})\end{pmatrix}\\ &= a_1(M_{11} \quad M_{12} \quad M_{13}) + a_2(M_{21} \quad M_{22} \quad M_{23}) + a_3(M_{31} \quad M_{32} \quad M_{33})\\ &= (a_1M_{11} \quad a_1M_{12} \quad a_1M_{13}) + (a_2M_{21} \quad a_2M_{22} \quad a_2M_{23}) + (a_3M_{31} \quad a_3M_{32} \quad a_3M_{33})\\ &= \left(\sum_{1=1}^{3} a_i M_{i1}, \sum_{1=1}^{3} a_i M_{i2}, \sum_{1=1}^{3} a_i M_{i3}\right)\\ &= \sum_{i=1}^{3} \left|\mathbf{a}^{(i)}\right\rangle\left\langle \mathbf{m}_{(i)}\right|,\end{aligned} \tag{2.11}$$

where

$$\left|\mathbf{a}^{(1)}\right\rangle = (a_1), \qquad \left|\mathbf{a}^{(2)}\right\rangle = (a_2), \quad \left|\mathbf{a}^{(3)}\right\rangle = (a_3). \tag{2.12}$$

Thus, $\langle \mathbf{a}|\mathbf{M}$ is a linear combination of the rows of $\mathbf{M}$. This thinking leads us to the row space of the matrix $\mathbf{M}$ and the second definition of the rank of a matrix. The key idea is to take all combinations of the rows. All real numbers $a_i$ are allowed - the space includes $\langle \mathbf{a}|\mathbf{M}$ for all vectors $\langle \mathbf{a}|$. In this way, we get infinitely many output vectors $\langle \mathbf{a}|\mathbf{M}$.

Actually, we are seeing the clearest possible example of the first great theorem in linear algebra:

**Theorem 2.1:** The number of independent columns of a matrix equals the number of independent rows.

**Definition (Rank of a Matrix):** The rank of a matrix is the dimension of its column or row space.

This rank theorem is true for every matrix. Always columns and rows in linear algebra! The $m$ rows contain the same numbers $M_{ij}$ as the $n$ columns. But different vectors. This idea will explain the "big picture" of linear algebra. That picture shows how every $m$ by $n$ matrix $\mathbf{M}$ leads to four subspaces (two subspaces of $\mathbb{R}^m$ and two more of $\mathbb{R}^n$.)

Obviously, we can form an ordinary (real) number by taking the compound product of a bra, a square matrix, and a ket, which corresponds to forming the product of a row vector times a square matrix times a column vector

$$\begin{aligned}\langle \mathbf{a}|\mathbf{M}|\mathbf{b}\rangle &= \mathbf{a}^T\mathbf{M}\mathbf{b}\\ &= (a_1, a_2, a_3)\begin{pmatrix} M_{11} & M_{12} & M_{13}\\ M_{21} & M_{22} & M_{23}\\ M_{31} & M_{32} & M_{33}\end{pmatrix}\begin{pmatrix} b_1\\ b_2\\ b_3\end{pmatrix}\\ &= \left(\sum_{i=1}^{3} a_i M_{i1}, \sum_{i=1}^{3} a_i M_{i2}, \sum_{i=1}^{3} a_i M_{i3}\right)\begin{pmatrix} b_1\\ b_2\\ b_3\end{pmatrix}\\ &= \sum_{i=1}^{3} a_i M_{i1} b_1 + \sum_{i=1}^{3} a_i M_{i2} b_2 + \sum_{i=1}^{3} a_i M_{i3} b_3\\ &= \sum_{j=1}^{3}\sum_{i=1}^{3} a_i M_{ij} b_j,\end{aligned} \tag{2.13}$$

and

$$\mathbf{AB} = \begin{pmatrix} A_{11} & A_{12} & A_{13} \\ A_{21} & A_{22} & A_{23} \\ A_{31} & A_{32} & A_{33} \end{pmatrix} \begin{pmatrix} B_{11} & B_{12} & B_{13} \\ B_{21} & B_{22} & B_{23} \\ B_{31} & B_{32} & B_{33} \end{pmatrix}$$

$$= \begin{pmatrix} (A_{11} \quad A_{12} \quad A_{13}) \begin{pmatrix} B_{11} \\ B_{21} \\ B_{31} \end{pmatrix} & (A_{11} \quad A_{12} \quad A_{13}) \begin{pmatrix} B_{12} \\ B_{22} \\ B_{32} \end{pmatrix} & (A_{11} \quad A_{12} \quad A_{13}) \begin{pmatrix} B_{13} \\ B_{23} \\ B_{33} \end{pmatrix} \\ (A_{21} \quad A_{22} \quad A_{23}) \begin{pmatrix} B_{11} \\ B_{21} \\ B_{31} \end{pmatrix} & (A_{21} \quad A_{22} \quad A_{23}) \begin{pmatrix} B_{12} \\ B_{22} \\ B_{32} \end{pmatrix} & (A_{21} \quad A_{22} \quad A_{23}) \begin{pmatrix} B_{13} \\ B_{23} \\ B_{33} \end{pmatrix} \\ (A_{31} \quad A_{32} \quad A_{33}) \begin{pmatrix} B_{11} \\ B_{21} \\ B_{31} \end{pmatrix} & (A_{31} \quad A_{32} \quad A_{33}) \begin{pmatrix} B_{12} \\ B_{22} \\ B_{32} \end{pmatrix} & (A_{31} \quad A_{32} \quad A_{33}) \begin{pmatrix} B_{13} \\ B_{23} \\ B_{33} \end{pmatrix} \end{pmatrix}$$

$$= \begin{pmatrix} \sum_{p=1}^{3} A_{1p}B_{p1} & \sum_{p=1}^{3} A_{1p}B_{p2} & \sum_{p=1}^{3} A_{1p}B_{p3} \\ \sum_{p=1}^{3} A_{2p}B_{p1} & \sum_{p=1}^{3} A_{2p}B_{p2} & \sum_{p=1}^{3} A_{2p}B_{p3} \\ \sum_{p=1}^{3} A_{3p}B_{p1} & \sum_{p=1}^{3} A_{3p}B_{p2} & \sum_{p=1}^{3} A_{3p}B_{p3} \end{pmatrix}$$

$$= \begin{pmatrix} \langle \mathbf{a}_{(1)} | \mathbf{b}^{(1)} \rangle & \langle \mathbf{a}_{(1)} | \mathbf{b}^{(2)} \rangle & \langle \mathbf{a}_{(1)} | \mathbf{b}^{(3)} \rangle \\ \langle \mathbf{a}_{(2)} | \mathbf{b}^{(1)} \rangle & \langle \mathbf{a}_{(2)} | \mathbf{b}^{(2)} \rangle & \langle \mathbf{a}_{(2)} | \mathbf{b}^{(3)} \rangle \\ \langle \mathbf{a}_{(3)} | \mathbf{b}^{(1)} \rangle & \langle \mathbf{a}_{(3)} | \mathbf{b}^{(2)} \rangle & \langle \mathbf{a}_{(3)} | \mathbf{b}^{(3)} \rangle \end{pmatrix}, \tag{2.14}$$

$$\langle \mathbf{a}_{(1)} | = (A_{11} \quad A_{12} \quad A_{13}), \qquad \langle \mathbf{a}_{(2)} | = (A_{21} \quad A_{22} \quad A_{23}), \qquad \langle \mathbf{a}_{(3)} | = (A_{31} \quad A_{32} \quad A_{33}), \tag{2.15}$$

$$|\mathbf{b}^{(1)}\rangle = \begin{pmatrix} B_{11} \\ B_{21} \\ B_{31} \end{pmatrix}, \qquad |\mathbf{b}^{(2)}\rangle = \begin{pmatrix} B_{12} \\ B_{22} \\ B_{32} \end{pmatrix}, \qquad |\mathbf{b}^{(3)}\rangle = \begin{pmatrix} B_{13} \\ B_{23} \\ B_{33} \end{pmatrix}, \tag{2.16}$$

or

$$\mathbf{AB} = \begin{pmatrix} A_{11} \\ A_{21} \\ A_{31} \end{pmatrix} (B_{11} \quad B_{12} \quad B_{13}) + \begin{pmatrix} A_{12} \\ A_{22} \\ A_{32} \end{pmatrix} (B_{21} \quad B_{22} \quad B_{23}) + \begin{pmatrix} A_{13} \\ A_{23} \\ A_{33} \end{pmatrix} (B_{31} \quad B_{32} \quad B_{33})$$

$$= \begin{pmatrix} A_{11}B_{11} & A_{11}B_{12} & A_{11}B_{13} \\ A_{21}B_{11} & A_{21}B_{12} & A_{21}B_{13} \\ A_{31}B_{11} & A_{31}B_{12} & A_{31}B_{13} \end{pmatrix} + \begin{pmatrix} A_{12}B_{21} & A_{12}B_{22} & A_{12}B_{23} \\ A_{22}B_{21} & A_{22}B_{22} & A_{22}B_{23} \\ A_{32}B_{21} & A_{32}B_{22} & A_{32}B_{23} \end{pmatrix} + \begin{pmatrix} A_{13}B_{31} & A_{13}B_{32} & A_{13}B_{33} \\ A_{23}B_{31} & A_{23}B_{32} & A_{23}B_{33} \\ A_{33}B_{31} & A_{33}B_{32} & A_{33}B_{33} \end{pmatrix}$$

$$= \begin{pmatrix} \sum_{p=1}^{3} A_{1p}B_{p1} & \sum_{p=1}^{3} A_{1p}B_{p2} & \sum_{p=1}^{3} A_{1p}B_{p3} \\ \sum_{p=1}^{3} A_{2p}B_{p1} & \sum_{p=1}^{3} A_{2p}B_{p2} & \sum_{p=1}^{3} A_{2p}B_{p3} \\ \sum_{p=1}^{3} A_{3p}B_{p1} & \sum_{p=1}^{3} A_{3p}B_{p2} & \sum_{p=1}^{3} A_{3p}B_{p3} \end{pmatrix} = \sum_{p=1}^{3} \begin{pmatrix} A_{1p}B_{p1} & A_{1p}B_{p2} & A_{1p}B_{p3} \\ A_{2p}B_{p1} & A_{2p}B_{p2} & A_{2p}B_{p3} \\ A_{3p}B_{p1} & A_{3p}B_{p2} & A_{3p}B_{p3} \end{pmatrix} = \sum_{i=1}^{3} |\mathbf{a}^{(i)}\rangle\langle \mathbf{b}_{(i)}| \tag{2.17}$$

$$|\mathbf{a}^{(1)}\rangle = \begin{pmatrix} A_{11} \\ A_{21} \\ A_{31} \end{pmatrix}, \qquad |\mathbf{a}^{(2)}\rangle = \begin{pmatrix} A_{12} \\ A_{22} \\ A_{32} \end{pmatrix}, \qquad |\mathbf{a}^{(3)}\rangle = \begin{pmatrix} A_{13} \\ A_{23} \\ A_{33} \end{pmatrix}, \tag{2.18}$$

$$\langle \mathbf{b}_{(1)} | = (B_{11} \quad B_{12} \quad B_{13}), \qquad \langle \mathbf{b}_{(2)} | = (B_{21} \quad B_{22} \quad B_{23}), \qquad \langle \mathbf{b}_{(3)} | = (B_{31} \quad B_{32} \quad B_{33}). \tag{2.19}$$

## 2.2 Vectors and Matrices with Index Notations

Let $A_{ij} \in \mathbb{R}$, $i = 1,2,\dots,m$, $j = 1,2,\dots,n$. The matrix is the ordered rectangular array

$$\mathbf{A} = \begin{pmatrix} A_{11} & A_{12} & \cdots & A_{1n} \\ A_{21} & A_{22} & \cdots & A_{2n} \\ \vdots & \vdots & \ddots & \vdots \\ A_{m1} & A_{m2} & \dots & A_{mn} \end{pmatrix} = A_{ij}, \;\; i = 1,2,\dots,m; \; j = 1,2,\dots,n. \tag{2.20}$$

Note the first subscript locates the row in which the typical element lies while the second subscript locates the column. For example, $A_{jk}$ denotes the element lying in the $j^{\text{th}}$ row and $k^{\text{th}}$ column of the matrix **A**.

**Definition:** Let **A** and **B** are $m \times n$ matrices and let the sum $\mathbf{A} + \mathbf{B}$ be $\mathbf{C} = \mathbf{A} + \mathbf{B}$, then **C** is an $m \times n$ matrix, with the element $(i, j)$ given by

$$C_{ij} = A_{ij} + B_{ij}, \tag{2.21}$$

for all $i = 1,2,\dots,m$, $j = 1,2,\dots,n$.

**Definition:** Let **A** be an $m \times n$ matrix, and let $a$ be a scalar, then $\mathbf{C} = a\mathbf{A}$ is an $m \times n$ matrix, with the element $(i, j)$ given by

$$C_{ij} = aA_{ij}, \tag{2.22}$$

for all $i = 1,2,\dots,m$, $j = 1,2,\dots,n$.

**Definition:** Let **A** be an $m \times n$ matrix, then $\mathbf{A}^T$ is an $n \times m$ matrix (the transpose of **A**) with the element $(i, j)$ given by

$$(A^T)_{ij} = A_{ji}, \tag{2.23}$$

for all $i = 1,2,\dots,m$, $j = 1,2,\dots,n$.

***Example 2.1***

Let us consider,

$$\mathbf{A} = \begin{pmatrix} A_{11} & A_{12} & A_{13} \\ A_{21} & A_{22} & A_{23} \\ A_{31} & A_{32} & A_{33} \end{pmatrix} = A_{ij}, \quad \mathbf{B} = \begin{pmatrix} B_{11} & B_{12} & B_{13} \\ B_{21} & B_{22} & B_{23} \\ B_{31} & B_{32} & B_{33} \end{pmatrix} = B_{ij}.$$

Hence, we have

$$(A + B)_{ij} = A_{ij} + B_{ij} \Rightarrow \mathbf{A} + \mathbf{B} = \begin{pmatrix} A_{11} + B_{11} & A_{12} + B_{12} & A_{13} + B_{13} \\ A_{21} + B_{21} & A_{22} + B_{22} & A_{23} + B_{23} \\ A_{31} + B_{31} & A_{32} + B_{32} & A_{33} + B_{33} \end{pmatrix},$$

$$(cA)_{ij} = cA_{ij} \Rightarrow c\mathbf{A} = \begin{pmatrix} cA_{11} & cA_{12} & cA_{13} \\ cA_{21} & cA_{22} & cA_{23} \\ cA_{31} & cA_{32} & cA_{33} \end{pmatrix},$$

$$(A^T)_{ij} = A_{ji} \Rightarrow \mathbf{A}^T = \begin{pmatrix} A_{11} & A_{21} & A_{31} \\ A_{12} & A_{22} & A_{32} \\ A_{13} & A_{23} & A_{33} \end{pmatrix}.$$

**Definition:** Let **A** and **B** are $m \times n$ and $n \times p$ matrices, respectively, and let the product **AB** be

$$\mathbf{C} = \mathbf{AB}, \tag{2.24}$$

then **C** is an $m \times p$ matrix, with the element $(i, j)$ given by

$$C_{ij} = \sum_{k=1}^{n} A_{ik} B_{kj}, \tag{2.25}$$

for all $i = 1,2,\dots,m$, $j = 1,2,\dots,p$.

**Definition:** Let **A** be an $m \times n$ matrix, and **x** be $n \times 1$ vector, then the typical element of the product

$$\mathbf{z} = \mathbf{Ax}, \tag{2.26}$$

is given by

$$z_i = \sum_{k=1}^{n} A_{ik} x_k, \tag{2.27}$$

for all $i = 1,2, \dots, m$.

**Definition:** Let $\mathbf{A}$ be an $m \times n$ matrix, and $\mathbf{y}$ be an $m \times 1$ vector, then the typical element of the product

$$\mathbf{z}^T = \mathbf{y}^T\mathbf{A}, \tag{2.28}$$

is given by

$$z_i = \sum_{k=1}^{n} y_k A_{ki}, \tag{2.29}$$

for all $i = 1,2,\dots,n$.

**Definition:** Let $\mathbf{A}$ be an $m \times n$ matrix, $\mathbf{y}$ be an $m \times 1$ vector, and $\mathbf{x}$ be $n \times 1$ vector, then the scalar resulting from the product

$$\alpha = \mathbf{y}^T\mathbf{A}\mathbf{x}, \tag{2.30}$$

is given by

$$\alpha = \sum_{j=1}^{m}\sum_{k=1}^{n} y_j A_{jk} x_k. \tag{2.31}$$

**Definition:** The trace of an $n \times n$ square matrix $\mathbf{A}$ is defined as

$$\mathrm{tr}(\mathbf{A}) = \sum_{i=1}^{n} A_{ii} = A_{11} + A_{22} + \cdots + A_{nn}. \tag{2.32}$$

For example, let $\mathbf{A} = \begin{pmatrix} 1 & 2 & 3 \\ 4 & 5 & 6 \\ 7 & 8 & 9 \end{pmatrix}$, then $\mathrm{tr}(\mathbf{A}) = 1 + 5 + 9 = 15$.

**Theorem 2.2:**

1- The trace is a linear mapping. That is,

$$\mathrm{tr}(\mathbf{A} + \mathbf{B}) = \mathrm{tr}(\mathbf{A}) + \mathrm{tr}(\mathbf{B}), \tag{2.33}$$

$$\mathrm{tr}(c\mathbf{A}) = c\,\mathrm{tr}(\mathbf{A}), \tag{2.34}$$

$$\mathrm{tr}(\mathbf{A}) = \mathrm{tr}(\mathbf{A}^T), \tag{2.35}$$

for all square matrices $\mathbf{A}$ and $\mathbf{B}$, and all scalars $c$.

2- If $\mathbf{A}$ and $\mathbf{B}$ are $m \times n$ and $n \times m$ real matrices, respectively, then

$$\mathrm{tr}(\mathbf{AB}) = \mathrm{tr}(\mathbf{BA}). \tag{2.36}$$

More generally, the trace is invariant under cyclic permutations, that is,

$$\mathrm{tr}(\mathbf{ABC}) = \mathrm{tr}(\mathbf{BCA}) = \mathrm{tr}(\mathbf{CAB}). \tag{2.37}$$

**Proof:**

$$\mathrm{tr}(\mathbf{A} + \mathbf{B}) = \sum_{i=1}^{n} A_{ii} + B_{ii} = \sum_{i=1}^{n} A_{ii} + \sum_{i=1}^{n} B_{ii} = \mathrm{tr}(\mathbf{A}) + \mathrm{tr}(\mathbf{B}),$$

$$\mathrm{tr}(c\mathbf{A}) = \sum_{i=1}^{n} cA_{ii} = c\sum_{i=1}^{n} A_{ii} = c\mathrm{tr}(\mathbf{A}),$$

$$\mathrm{tr}(\mathbf{A}) = \sum_{i=1}^{n} A_{ii} = \mathrm{tr}(\mathbf{A}^T).$$

The trace of a matrix is the sum of its diagonal elements, but transposition leaves the diagonal elements unchanged.

$$\mathrm{tr}(\mathbf{AB}) = \sum_{i=1}^{n}\sum_{j=1}^{m} A_{ij}B_{ji} = \sum_{j=1}^{m}\sum_{i=1}^{n} B_{ji}A_{ij} = \mathrm{tr}(\mathbf{BA}),$$

$$\mathrm{tr}(\mathbf{ABC}) = \sum_{i=1}^{n}\sum_{j=1}^{m}\sum_{k=1}^{l} A_{ij}B_{jk}C_{ki} = \sum_{i=1}^{n}\sum_{j=1}^{m}\sum_{k=1}^{l} B_{jk}C_{ki}A_{ij} = \mathrm{tr}(\mathbf{BCA}).$$

■

## Some Types of Matrices

**Definition (Identity Matrix):** The identity matrix of size $n$ is the $n \times n$ square matrix with ones on the main diagonal and zeros elsewhere.

$$\mathbf{I} = \begin{pmatrix} 1 & 0 & \dots & 0 \\ 0 & 1 & \dots & \vdots \\ \vdots & \vdots & \ddots & 0 \\ 0 & 0 & \dots & 1 \end{pmatrix}. \tag{2.38}$$

**Definition (Real Matrix):** The conjugate complex of a matrix **M** is written as $\mathbf{M}^*$ and the elements of $\mathbf{M}^*$ are the conjugate complexes of the elements of **M** i.e.

$$(M^*)_{ij} = \left(M_{ij}\right)^*. \tag{2.39}$$

For a real matrix, all the elements are real and, therefore

$$\mathbf{M} = \mathbf{M}^*, \qquad (M)_{ij} = M_{ij}^*. \tag{2.40}$$

**Definition (Symmetric Matrix):** The transpose of a matrix **M** is obtained by changing rows into columns (or vice versa). For a symmetric matrix

$$\mathbf{M}^T = \mathbf{M},\ M_{ij} = M_{ji}. \tag{2.41}$$

**Definition (Skew-Symmetric Matrix):** The skew-symmetric matrix satisfies

$$\mathbf{M}^T = -\mathbf{M},\ M_{ij} = -M_{ji}. \tag{2.42}$$

**Definition (Orthogonal Matrix):** An orthogonal matrix satisfies

$$\mathbf{M}\mathbf{M}^T = \mathbf{M}^T\mathbf{M} = \mathbf{I}. \tag{2.43}$$

This leads to the equivalent characterization: a matrix **M** is orthogonal if its transpose is equal to its inverse:

$$\mathbf{M}^T = \mathbf{M}^{-1}. \tag{2.44}$$

**Definition (Conjugate Transpose or Hermitian Transpose):** Hermitian transpose of an $m \times n$ complex matrix **M** is an $n \times m$ matrix obtained by transposing **M** and applying complex conjugate on each entry.

$$\mathbf{M}^H = \mathbf{M}^\dagger = (\mathbf{M}^T)^*. \tag{2.45}$$

For real matrices, the conjugate transpose is just the transpose, $\mathbf{M}^H = \mathbf{M}^\dagger = \mathbf{M}^T$.

**Definition (Unitary Matrix):** A complex square matrix **M** is unitary if its conjugate transpose $\mathbf{M}^\dagger$ is also its inverse, that is, if

$$\mathbf{M}^\dagger\mathbf{M} = \mathbf{M}\mathbf{M}^\dagger = \mathbf{M}\mathbf{M}^{-1} = \mathbf{I}. \tag{2.46}$$

**Definition (Lower and Upper Triangular Matrix):** A matrix of the form

$$\begin{pmatrix} l_{11} & 0 & \dots & 0 \\ l_{21} & l_{22} & \ddots & \vdots \\ \vdots & \vdots & \ddots & 0 \\ l_{n1} & l_{n2} & \dots & l_{nn} \end{pmatrix}, \tag{2.47}$$

is called a lower triangular matrix or left triangular matrix, and analogously a matrix of the form

$$\begin{pmatrix} u_{11} & u_{12} & \dots & u_{1n} \\ 0 & u_{22} & \dots & u_{2n} \\ \vdots & \ddots & \ddots & \vdots \\ 0 & \dots & 0 & u_{nn} \end{pmatrix}, \tag{2.48}$$

is called an upper triangular matrix or right triangular matrix. In the lower triangular matrix, all elements above the diagonal are zeros, in the upper triangular matrix, all the elements below the diagonal are zeros.

For example, the matrix **A** is symmetric, but the matrix **B** is skew-symmetric.

$$\mathbf{A} = \begin{pmatrix} 0 & 1 & -2 \\ 1 & 3 & 0 \\ -2 & 0 & 5 \end{pmatrix}, \quad \mathbf{B} = \begin{pmatrix} 0 & 2 & -45 \\ -2 & 0 & -4 \\ 45 & 4 & 0 \end{pmatrix}.$$

## 2.3 Matrix Decomposition

A matrix decomposition or matrix factorization [2-5] is a factorization of a matrix into a product of matrices. There are many different matrix decompositions; each finds use among a particular class of problems. Table 2.1 represents some types of matrix factorization.

**Table 2.1.** Some types of matrix factorization.

| Method | Decomposition |
|---|---|
| LU decomposition | Applicable to: square matrix **A**, although rectangular matrices can be applicable.<br>Decomposition form:<br>$\mathbf{A} = \mathbf{LU}$,<br>where **L** is lower triangular, and **U** is upper triangular.<br><br>Related: the LDU decomposition is<br>$\mathbf{A} = \mathbf{LDU}$,<br>where **L** is lower triangular with ones on the diagonal, **U** is upper triangular with ones on the diagonal, and **D** is a diagonal matrix.<br><br>Related: the LUP decomposition is<br>$\mathbf{PA} = \mathbf{LU}$,<br>where **L** is lower triangular, **U** is upper triangular, and **P** is a permutation matrix. |
| Rank factorization | Applicable to: $m \times n$ matrix **A** of rank $r$.<br>Decomposition form:<br>$\mathbf{A} = \mathbf{CF}$,<br>where **C** is an $m \times r$ full column rank matrix, and **F** is an $r \times n$ full row rank matrix. |
| Cholesky decomposition | Applicable to: square, Hermitian, and positive definite matrix **A**.<br>Decomposition form:<br>$\mathbf{A} = \mathbf{U}^*\mathbf{U}$,<br>where **U** is a lower triangular matrix with real and positive diagonal entries, and $\mathbf{U}^*$ denotes the conjugate transpose of **U**. |
| QR decomposition | Applicable to: $m \times n$ matrix **A** with linearly independent columns.<br>Decomposition form:<br>$\mathbf{A} = \mathbf{QR}$,<br>where **Q** is a unitary matrix of size $m \times m$, and **R** is an upper triangular matrix of size $m \times n$. |
| Eigen decomposition | Applicable to: square matrix **A** with linearly independent eigenvectors (not necessarily distinct eigenvalues).<br>Decomposition form:<br>$\mathbf{A} = \mathbf{VDV}^{-1}$,<br>where **D** is a diagonal matrix formed from the eigenvalues of **A**, and the columns of **V** are the corresponding eigenvectors of **A**. |
| Jordan decomposition | Applicable to: square matrix **A**.<br><br>The Jordan normal form generalizes the eigen decomposition to cases where there are repeated eigenvalues and cannot be diagonalized. |
| Schur decomposition | Applicable to: square matrix $\mathbf{A}$.<br>Decomposition form:<br>$$\mathbf{A} = \mathbf{UTU}^*,$$<br>where $\mathbf{U}$ is a unitary matrix, $\mathbf{U}^*$ is the conjugate transpose of $\mathbf{U}$, and $\mathbf{T}$ is an upper triangular matrix called the complex Schur form, which has the eigenvalues of $\mathbf{A}$ along its diagonal. |
| Real Schur decomposition | Applicable to: square matrix $\mathbf{A}$.<br>Decomposition form:<br>$$\mathbf{A} = \mathbf{VSV}^T,$$<br>where $\mathbf{V}$ is a real orthogonal matrix, $\mathbf{V}^T$ is the transpose of $\mathbf{V}$, and $\mathbf{S}$ is a block upper triangular matrix called the real Schur form. |
| Takagi's factorization | Applicable to: square, complex, and symmetric matrix $\mathbf{A}$.<br>Decomposition form:<br>$$\mathbf{A} = \mathbf{VDV}^T,$$<br>where $\mathbf{D}$ is a real nonnegative diagonal matrix, and $\mathbf{V}$ is unitary. $\mathbf{V}^T$ denotes the matrix transpose of $\mathbf{V}$. |
| Singular value decomposition (SVD) | Applicable to: $m \times n$ matrix $\mathbf{A}$.<br>Decomposition form:<br>$$\mathbf{A} = \mathbf{UDV}^*,$$<br>where $\mathbf{D}$ is a nonnegative diagonal matrix, and $\mathbf{U}$ and $\mathbf{V}$ satisfy $\mathbf{U}^*\mathbf{U} = \mathbf{I}$, $\mathbf{V}^*\mathbf{V} = \mathbf{I}$. Here $\mathbf{V}^*$ and $\mathbf{U}^*$ are the conjugate transpose of $\mathbf{V}$ and $\mathbf{U}$, respectively (or simply the transpose, if $\mathbf{V}$ and $\mathbf{U}$ contain real numbers only), and $\mathbf{I}$ denotes the identity matrix (of some dimension). |

The first and most fundamental decomposition problem in linear algebra is the LU decomposition. In the remaining of the present section, we will briefly consider this problem. While in the next chapter, we will consider the SVD in detail.

### LU Factorization

Let $\mathbf{A}$ be a square matrix. The LU factorization refers to the factorization of $\mathbf{A}$, with proper row and/or column orderings or permutations, into two factors – a lower triangular matrix $\mathbf{L}$ and an upper triangular matrix $\mathbf{U}$:

$$\mathbf{A} = \mathbf{LU}. \tag{2.49}$$

For example, consider a $3 \times 3$ matrix $\mathbf{A}$; its LU decomposition looks like this:

$$\begin{pmatrix} a_{11} & a_{12} & a_{13} \\ a_{21} & a_{22} & a_{23} \\ a_{31} & a_{32} & a_{33} \end{pmatrix} = \begin{pmatrix} l_{11} & 0 & 0 \\ l_{21} & l_{22} & 0 \\ l_{31} & l_{32} & l_{33} \end{pmatrix} \begin{pmatrix} u_{11} & u_{12} & u_{13} \\ 0 & u_{22} & u_{23} \\ 0 & 0 & u_{33} \end{pmatrix}. \tag{2.50}$$

Without a proper ordering or permutations in the matrix, the factorization may fail to materialize. For example, it is easy to verify (by expanding the matrix multiplication) that $a_{11} = l_{11}u_{11}$. If $a_{11} = 0$, then at least one of $l_{11}$ and $u_{11}$ has to be zero, which implies that either $\mathbf{L}$ or $\mathbf{U}$ is singular. This is impossible if $\mathbf{A}$ is nonsingular (invertible). This is a procedural problem. It can be removed by simply reordering the rows of $\mathbf{A}$ so that the first element of the permuted matrix is non-zero.

Gaussian elimination or row reduction is an algorithm for giving $\mathbf{U}$ and also leads to a procedure for finding $\mathbf{L}$. It consists of a sequence of operations performed on the corresponding matrix of coefficients. To perform row reduction on a matrix, one uses a sequence of elementary row operations to modify the matrix until the lower left-hand corner of the matrix is filled with zeros, as much as possible. There are three types of elementary row operations:

- Swapping two rows.
- Multiplying a row by a non-zero number.
- Adding a multiple of one row to another row.

Using these operations, a matrix can always be transformed into an upper triangular matrix and in fact, one that is in row echelon form. For each row in a matrix, if the row does not consist of only zeros, then the leftmost non-zero entry is called the leading coefficient (or pivot) of that row. So, if two leading coefficients are in the same column, then a row operation of type 3 could be used to make one of those coefficients zero. Then by using the row-swapping operation, one can always order the rows so that for every non-zero row, the leading coefficient is to the right of the leading coefficient of the row above. If this is the case, then the matrix is said to be in row echelon form. So, the lower left part of the matrix contains only zeros, and all of the zero rows are below the non-zero rows.

A matrix is in reduced row echelon form if it satisfies the following conditions:

- It is in row echelon form.
- The leading entry in each non-zero row is a 1 (called a leading 1).
- Each column containing a leading 1 has zeros in all its other entries.

The reduced row echelon form is unique; in other words, it is independent of the sequence of row operations used. For example, in the following sequence of row operations (where two elementary operations on different rows are done at the first and third steps), the third and fourth matrices are the ones in row echelon form, and the final matrix is the unique reduced row echelon form.

$$\begin{pmatrix} 1 & 3 & 1 & 9 \\ 1 & 1 & -1 & 1 \\ 3 & 11 & 5 & 35 \end{pmatrix} \to \begin{pmatrix} 1 & 3 & 1 & 9 \\ 0 & -2 & -2 & -8 \\ 0 & 2 & 2 & 8 \end{pmatrix} \to \begin{pmatrix} 1 & 3 & 1 & 9 \\ 0 & -2 & -2 & -8 \\ 0 & 0 & 0 & 0 \end{pmatrix} \to \begin{pmatrix} 1 & 0 & -2 & -3 \\ 0 & 1 & 1 & 4 \\ 0 & 0 & 0 & 0 \end{pmatrix}.$$

Using row operations to convert a matrix into reduced row echelon form is sometimes called Gauss–Jordan elimination. In this case, Gaussian elimination refers to the process until it reaches its upper triangular or (unreduced) row echelon form.

If an LU factorization $\mathbf{A} = \mathbf{LU}$ does exist, then the gaussian algorithm gives $\mathbf{U}$ and also leads to a procedure for finding $\mathbf{L}$. Example 2.2 provides an illustration.

***Example 2.2***

Find an LU factorization for

$$\mathbf{A} = \begin{pmatrix} 5 & -5 & 10 & 0 & 5 \\ -3 & 3 & 2 & 2 & 1 \\ -2 & 2 & 0 & -1 & 0 \\ 1 & -1 & 10 & 2 & 5 \end{pmatrix}.$$

***Solution***

The reduction to row-echelon form is

$$\mathbf{A} = \begin{pmatrix} 5 & -5 & 10 & 0 & 5 \\ -3 & 3 & 2 & 2 & 1 \\ -2 & 2 & 0 & -1 & 0 \\ 1 & -1 & 10 & 2 & 5 \end{pmatrix} \to \begin{pmatrix} 1 & -1 & 2 & 0 & 1 \\ 0 & 0 & 8 & 2 & 4 \\ 0 & 0 & 4 & -1 & 2 \\ 0 & 0 & 8 & 2 & 4 \end{pmatrix} \to \begin{pmatrix} 1 & -1 & 2 & 0 & 1 \\ 0 & 0 & 1 & 1/4 & 1/2 \\ 0 & 0 & 0 & -2 & 0 \\ 0 & 0 & 0 & 0 & 0 \end{pmatrix}$$

$$\to \begin{pmatrix} 1 & -1 & 2 & 0 & 1 \\ 0 & 0 & 1 & 1/4 & 1/2 \\ 0 & 0 & 0 & 1 & 0 \\ 0 & 0 & 0 & 0 & 0 \end{pmatrix} = \mathbf{U}.$$

The blue columns are determined as follows: The first is the leading column of $\mathbf{A}$, and is used (by lower reduction) to create the first leading 1 and create zeros below it. This completes the work on the first row, and we repeat the procedure on the matrix consisting of the remaining rows. Thus, the second blue column is the leading column of this smaller matrix, which we use to create the second leading 1 and the zeros below it. If $\mathbf{U}$ denotes the final row-echelon matrix, then $\mathbf{A} = \mathbf{LU}$, where

$$\mathbf{L} = \begin{pmatrix} -5 & 0 & 0 & 0 \\ -3 & 8 & 0 & 0 \\ -2 & 4 & -2 & 0 \\ 1 & 8 & 0 & 1 \end{pmatrix}.$$

**Theorem 2.3:** Suppose an $m \times n$ matrix $\mathbf{A}$ is carried to a row-echelon matrix $\mathbf{U}$ via the gaussian algorithm. Let $\mathbf{P}_1$, $\mathbf{P}_2$, ..., $\mathbf{P}_s$ be the elementary matrices corresponding (in order) to the row interchanges used and write

$$\mathbf{P} = \mathbf{P}_s \cdots \mathbf{P}_2\mathbf{P}_1.$$

(If no interchanges are used, take $\mathbf{P} = \mathbf{I}_m$.) Then:
1. $\mathbf{PA}$ is the matrix obtained from $\mathbf{A}$ by doing these interchanges (in order) to $\mathbf{A}$.
2. $\mathbf{PA}$ has an LU factorization.

A matrix $\mathbf{P}$ that is the product of elementary matrices corresponding to row interchanges is called a permutation matrix. Such a matrix is obtained from the identity matrix by arranging the rows in a different order, so it has only one in each row and each column and zeros elsewhere. We regard the identity matrix as a permutation matrix. The elementary permutation matrices are those obtained from $\mathbf{I}$ by a single-row interchange, and every permutation matrix is a product of elementary ones.

***Example 2.3***

Given a square $4 \times 4$ matrix $\mathbf{A}$ defined as

$$\mathbf{A} = \begin{pmatrix} 0 & 0 & -1 & 2 \\ -1 & -1 & 1 & 2 \\ 2 & 1 & -3 & 6 \\ 0 & 1 & -1 & 4 \end{pmatrix}.$$

Find a permutation matrix $\mathbf{P}$ such that $\mathbf{PA}$ has an LU factorization.

***Solution***

Apply the gaussian algorithm to $\mathbf{A}$:

$$\begin{pmatrix} 0 & 0 & -1 & 2 \\ -1 & -1 & 1 & 2 \\ 2 & 1 & -3 & 6 \\ 0 & 1 & -1 & 4 \end{pmatrix} \overset{*}{\rightarrow} \begin{pmatrix} -1 & -1 & 1 & 2 \\ 0 & 0 & -1 & 2 \\ 2 & 1 & -3 & 6 \\ 0 & 1 & -1 & 4 \end{pmatrix} \rightarrow \begin{pmatrix} 1 & 1 & -1 & -2 \\ 0 & 0 & -1 & 2 \\ 0 & -1 & -1 & 10 \\ 0 & 1 & -1 & 4 \end{pmatrix} \overset{*}{\rightarrow} \begin{pmatrix} 1 & 1 & -1 & -2 \\ 0 & -1 & -1 & 10 \\ 0 & 0 & -1 & 2 \\ 0 & 1 & -1 & 4 \end{pmatrix} \rightarrow$$
$$\begin{pmatrix} 1 & 1 & -1 & -2 \\ 0 & 1 & 1 & -10 \\ 0 & 0 & -1 & 2 \\ 0 & 0 & -2 & 14 \end{pmatrix} \rightarrow \begin{pmatrix} 1 & 1 & -1 & -2 \\ 0 & 1 & 1 & -10 \\ 0 & 0 & 1 & -2 \\ 0 & 0 & 0 & 10 \end{pmatrix}.$$

Two-row interchanges were needed (marked with $*$), first and second rows, and then second and third. Hence, as in Theorem 2.3,

$$\mathbf{P} = \begin{pmatrix} 1 & 0 & 0 & 0 \\ 0 & 0 & 1 & 0 \\ 0 & 1 & 0 & 0 \\ 0 & 0 & 0 & 1 \end{pmatrix} \begin{pmatrix} 0 & 1 & 0 & 0 \\ 1 & 0 & 0 & 0 \\ 0 & 0 & 1 & 0 \\ 0 & 0 & 0 & 1 \end{pmatrix} = \begin{pmatrix} 0 & 1 & 0 & 0 \\ 0 & 0 & 1 & 0 \\ 1 & 0 & 0 & 0 \\ 0 & 0 & 0 & 1 \end{pmatrix}.$$

If we apply these interchanges (in order) to $\mathbf{A}$, the result is $\mathbf{PA}$. Now apply the LU algorithm to $\mathbf{PA}$:

$$\mathbf{PA} = \begin{pmatrix} -1 & -1 & 1 & 2 \\ 2 & 1 & -3 & 6 \\ 0 & 0 & -1 & 2 \\ 0 & 1 & -1 & 4 \end{pmatrix} \rightarrow \begin{pmatrix} 1 & 1 & -1 & -2 \\ 0 & -1 & -1 & 10 \\ 0 & 0 & -1 & 2 \\ 0 & 1 & -1 & 4 \end{pmatrix} \rightarrow \begin{pmatrix} 1 & 1 & -1 & -2 \\ 0 & 1 & 1 & -10 \\ 0 & 0 & -1 & 2 \\ 0 & 1 & -2 & 14 \end{pmatrix} \rightarrow \begin{pmatrix} 1 & 1 & -1 & -2 \\ 0 & 1 & 1 & -10 \\ 0 & 0 & 1 & -2 \\ 0 & 0 & 0 & 10 \end{pmatrix} \rightarrow$$
$$\begin{pmatrix} 1 & 1 & -1 & -2 \\ 0 & 1 & 1 & -10 \\ 0 & 0 & 1 & -2 \\ 0 & 1 & 0 & 1 \end{pmatrix} = \mathbf{U}.$$

Hence, $\mathbf{PA} = \mathbf{LU}$, where

$$\mathbf{L} = \begin{pmatrix} -1 & 0 & 0 & 0 \\ 2 & -1 & 0 & 0 \\ 0 & 0 & -1 & 0 \\ 0 & 1 & -2 & 10 \end{pmatrix} \text{ and } \mathbf{U} = \begin{pmatrix} 1 & 1 & -1 & -2 \\ 0 & 1 & 1 & -10 \\ 0 & 0 & 1 & -2 \\ 0 & 1 & 0 & 1 \end{pmatrix}.$$

Theorem 2.3 provides an essential general factorization theorem for matrices. If $\mathbf{A}$ is any $m \times n$ matrix, it asserts that there exists a permutation matrix $\mathbf{P}$ and an LU factorization $\mathbf{PA} = \mathbf{LU}$. Moreover, it shows that either $\mathbf{P} = \mathbf{I}$ or $\mathbf{P} = \mathbf{P}_s \cdots \mathbf{P}_2\mathbf{P}_1$, where $\mathbf{P}_1$, $\mathbf{P}_2$, ..., $\mathbf{P}_s$ are the elementary permutation matrices arising in the reduction of $\mathbf{A}$ to row-echelon form. Now observe that $\mathbf{P}_i^{-1} = \mathbf{P}_i$ for each $i$ (they are elementary row interchanges). Thus, $\mathbf{P}^{-1} = \mathbf{P}_1\mathbf{P}_2 \cdots \mathbf{P}_s$, so the matrix $\mathbf{A}$ can be factorized as

$$\mathbf{A} = \mathbf{P}^{-1}\mathbf{LU}, \tag{2.51}$$

where $\mathbf{P}^{-1}$ is a permutation matrix, $\mathbf{L}$ is lower triangular and invertible, and $\mathbf{U}$ is a row-echelon matrix. This is called a PLU factorization of $\mathbf{A}$.

**Theorem 2.4:** Let $\mathbf{A}$ be an $m \times n$ matrix that has an LU factorization

$$\mathbf{A} = \mathbf{LU}. \tag{2.52}$$

If $\mathbf{A}$ has rank $m$ (that is, $\mathbf{U}$ has no row of zeros), then $\mathbf{L}$ and $\mathbf{U}$ are uniquely determined by $\mathbf{A}$.

Finally, we can classify the LU decompositions according to the following cases:

- **LU factorization with partial pivoting**
  It turns out that a proper permutation in rows (or columns) is sufficient for LU factorization. LU factorization with partial pivoting (LUP) often refers to LU factorization with row permutations only:
  $$\mathbf{PA} = \mathbf{LU}, \tag{4.53}$$
  where $\mathbf{L}$ and $\mathbf{U}$ are again lower and upper triangular matrices, and $\mathbf{P}$ is a permutation matrix, which, when left-multiplied to $\mathbf{A}$, reorders the rows of $\mathbf{U}$. It turns out that all square matrices can be factorized in this form, and the factorization is numerically stable in practice.
- **LU factorization with full pivoting**
  An LU factorization with full pivoting involves both row and column permutations:
  $$\mathbf{PAQ} = \mathbf{LU}, \tag{4.54}$$
  where $\mathbf{L}$, $\mathbf{U}$, and $\mathbf{P}$ are defined as before, and $\mathbf{Q}$ is a permutation matrix that reorders the columns of $\mathbf{A}$.
- **Lower-diagonal-upper (LDU) decomposition**
  An LDU decomposition is a decomposition of the form
  $$\mathbf{A} = \mathbf{LDU}, \tag{4.55}$$
  where $\mathbf{D}$ is a diagonal matrix, and $\mathbf{L}$ and $\mathbf{U}$ are unit-triangular matrices, meaning that all the entries on the diagonals of $\mathbf{L}$ and $\mathbf{U}$ are one.

## 2.4 Basics of Matrix Calculus

In the world of single-variable functions, the options are limited for taking the derivative; for $f: \mathbb{R} \to \mathbb{R}$, $x \to f(x)$, the only derivative of our interest is $\frac{df}{dx}$. But with functions such as $\mathbf{g}(\mathbf{x}) = \mathbf{Ax}$ and $h(\mathbf{x}, \mathbf{A}) = \langle \mathbf{x}|\mathbf{A}|\mathbf{x}\rangle$, we can also consider derivatives such as $\frac{d\mathbf{g}}{d\mathbf{x}}$, $\frac{d\mathbf{g}}{dx_i}$, $\frac{dh}{d\mathbf{A}}$, $\frac{dh}{dA_{ij}}$, $\frac{dh}{d\mathbf{x}^T}$, etc. In particular, we have the following cases:

| | Scalar | Vector | Matrix |
|---|---|---|---|
| Scalar | $\frac{dy}{dx}$ | $\frac{dy}{d\mathbf{x}}$ | $\frac{dy}{d\mathbf{X}}$ |
| Vector | $\frac{d\mathbf{y}}{dx}$ | $\frac{d\mathbf{y}}{d\mathbf{x}}$ | $\frac{d\mathbf{y}}{d\mathbf{X}}$ |
| Matrix | $\frac{d\mathbf{Y}}{dx}$ | $\frac{d\mathbf{Y}}{d\mathbf{x}}$ | $\frac{d\mathbf{Y}}{d\mathbf{X}}$ |

There are many different versions of definitions [6,7], but here we use the denominator-layout notation. Also note that we use $d$ and $\partial$ interchangeably.

### Derivatives of Scalar

We first consider when we take the derivative of a scalar.

1. With respect to a scalar ($\frac{dy}{dx}$): We already know this case. This is simply the single-variable function case.

2. With respect to a vector ($\frac{dy}{d\mathbf{x}}$): An example of this case is when $y = \|\mathbf{x}\|$. This is the gradient we defined. That is, for $\mathbf{x} \in \mathbb{R}^{n\times 1}$,

$$\frac{dy}{d\mathbf{x}} = \begin{pmatrix} \frac{dy}{dx_1} \\ \vdots \\ \frac{dy}{dx_n} \end{pmatrix} \in \mathbb{R}^{n\times 1}. \tag{2.56}$$

We also define what happens when we take the derivative of a scalar with respect to a row vector $\mathbf{x}^T$:

$$\frac{dy}{d\mathbf{x}^T} = \left(\frac{dy}{dx_1}, \dots, \frac{dy}{dx_n}\right) \in \mathbb{R}^{1\times n}. \tag{2.57}$$

3. With respect to a matrix ($\frac{dy}{d\mathbf{X}}$): An example of this case is $y = \sqrt{\sum_{i=1}^{m}\sum_{j=1}^{n}\left|X_{ij}\right|^2}$, where $X_{ij}$ are the entries of the matrix. Expanding on the vector case, for $\mathbf{X} \in \mathbb{R}^{m\times n}$:

$$\frac{dy}{d\mathbf{X}} = \begin{pmatrix} \frac{dy}{dX_{11}} & \cdots & \frac{dy}{dX_{1n}} \\ \vdots & \ddots & \vdots \\ \frac{dy}{dX_{m1}} & \cdots & \frac{dy}{dX_{mn}} \end{pmatrix} \in \mathbb{R}^{m\times n}. \tag{2.58}$$

One thing to notice here is that when we take the derivative of a scalar, we end up with the same shape as the variable we took the derivative with respect to. For example, the shape of $\frac{dy}{d\mathbf{x}}$ is the same as the shape of $\mathbf{x}$. This is a nice property of the denominator-layout notation.

### Derivatives of Vector

Now we expand the scalar case to vectors, i.e., $\frac{d\mathbf{y}}{dx}$, $\frac{d\mathbf{y}}{d\mathbf{x}}$, and $\frac{d\mathbf{y}}{d\mathbf{X}}$. Note that $\mathbf{y}$ here does not necessarily have to be a column vector. The same definitions also apply to row vectors, including the resulting shapes.

1. With respect to a scalar ($\frac{d\mathbf{y}}{dx}$): An example of this case is $d(x\mathbf{v})/dx$ for a scalar $x$ and constant vector $\mathbf{v} \in \mathbb{R}^n$. For $\mathbf{y} \in \mathbb{R}^n$, this is defined as:

$$\frac{d\mathbf{y}}{dx} = \left(\frac{dy_1}{dx}, \dots, \frac{dy_n}{dx}\right) \in \mathbb{R}^{1\times n}. \tag{2.59}$$

2. With respect to a vector ($\frac{d\mathbf{y}}{d\mathbf{x}}$): An example of this case is $\mathbf{y} = \mathbf{A}\mathbf{x}$ for a constant matrix $\mathbf{A}$, and we evaluate $\frac{d\mathbf{y}}{d\mathbf{x}}$. For $\mathbf{y} \in \mathbb{R}^n$ and $\mathbf{x} \in \mathbb{R}^p$, this is defined as

$$\frac{d\mathbf{y}}{d\mathbf{x}} = \left(\nabla y_1(\mathbf{x}), \dots, \nabla y_n(\mathbf{x})\right) = \begin{pmatrix} \frac{dy_1}{dx_1} & \frac{dy_2}{dx_1} & \cdots & \frac{dy_n}{dx_1} \\ \frac{dy_1}{dx_2} & \frac{dy_2}{dx_2} & \cdots & \frac{dy_n}{dx_2} \\ \vdots & \vdots & \ddots & \vdots \\ \frac{dy_1}{dx_p} & \frac{dy_2}{dx_p} & \cdots & \frac{dy_n}{dx_p} \end{pmatrix} \in \mathbb{R}^{p\times n}. \tag{2.60}$$

Consider when $\mathbf{y} = \mathbf{A}\mathbf{x}$ for a constant matrix $\mathbf{A} \in \mathbb{R}^{n\times p}$. Explicit multiplication yields

$$\begin{aligned}
\mathbf{y} &= \mathbf{A}\mathbf{x} \\
&= \begin{pmatrix} A_{11} & \cdots & A_{1p} \\ \vdots & \ddots & \vdots \\ A_{n1} & \cdots & A_{np} \end{pmatrix} \begin{pmatrix} x_1 \\ \vdots \\ x_p \end{pmatrix} \\
&= \begin{pmatrix} A_{11}x_1 + A_{12}x_2 + \cdots + A_{1p}x_p \\ \vdots \\ A_{n1}x_1 + A_{n2}x_2 + \cdots + A_{np}x_p \end{pmatrix} \\
&= \begin{pmatrix} \sum_{k=1}^{p} A_{1k}x_k \\ \vdots \\ \sum_{k=1}^{p} A_{nk}x_k \end{pmatrix}.
\end{aligned} \tag{2.61}$$

This gives $y_i = \sum_{k=1}^{p} A_{ik}x_k$, and therefore $dy_i/dx_j = A_{ij}$. Hence, we have

$$\begin{aligned}
\frac{d\mathbf{y}}{d\mathbf{x}} &= \begin{pmatrix} \frac{dy_1}{dx_1} & \frac{dy_2}{dx_1} & \cdots & \frac{dy_n}{dx_1} \\ \frac{dy_1}{dx_2} & \frac{dy_2}{dx_2} & \cdots & \frac{dy_n}{dx_2} \\ \vdots & \vdots & \ddots & \vdots \\ \frac{dy_1}{dx_p} & \frac{dy_2}{dx_p} & \cdots & \frac{dy_n}{dx_p} \end{pmatrix} \\
&= \begin{pmatrix} A_{11} & A_{21} & \dots & A_{n1} \\ A_{12} & A_{22} & \dots & A_{n2} \\ \vdots & \vdots & \ddots & \vdots \\ A_{1p} & A_{2p} & \dots & A_{np} \end{pmatrix} \\
&= \mathbf{A}^T.
\end{aligned} \tag{2.62}$$

Hence, we have derived one helpful result:

$$\frac{d(\mathbf{A}\mathbf{x})}{d\mathbf{x}} = \mathbf{A}^T. \tag{2.63}$$

3. With respect to a matrix ($\frac{d\mathbf{y}}{d\mathbf{X}}$): An example of this case is $\mathbf{y} = \mathbf{X}\mathbf{v}$ for a constant vector $\mathbf{v}$, and we evaluate $\frac{d\mathbf{y}}{d\mathbf{X}}$. In general, this encodes three-dimensional information ($dy_i/dX_{jk}$) and is beyond the scope of these lectures. However, we define the following two specific cases:

$$\frac{d(\mathbf{X}\mathbf{v})}{d\mathbf{X}} = \mathbf{v}^T, \frac{d(\mathbf{v}^T\mathbf{X})}{d\mathbf{X}} = \mathbf{v}, \tag{2.64}$$

for a matrix $\mathbf{X}$ and a constant vector $\mathbf{v}$. Note that the second case is the derivative of a row vector with respect to a matrix.

The fundamental issue is that the derivative of a vector with respect to a vector, i.e., $\frac{d\mathbf{y}}{d\mathbf{x}}$, is often written in two competing ways. If the numerator $\mathbf{y}$ is of size $m$ and the denominator $\mathbf{x}$ of size $n$, then the result can be laid out as either an $m \times n$ matrix or $n \times m$ matrix, i.e., the elements of $\mathbf{y}$ laid out in columns and the elements of $\mathbf{x}$ laid out in rows, or vice versa. This leads to the following possibilities (see Table 2.2):

1- Numerator layout, i.e., lay out corresponding to $\mathbf{y}$ and $\mathbf{x}^T$. This is sometimes known as the Jacobian formulation. This corresponds to the $m \times n$ layout in the previous example.
2- Denominator layout, i.e., lay out corresponding to $\mathbf{y}^T$ and $\mathbf{x}$. This is sometimes known as the Hessian formulation. This corresponds to the $n \times m$ layout in the previous example.

**Table 2.2.** Derivatives of scalar, vector, and matrix based on numerator- and denominator-layout notations.

| Numerator-layout notation | Denominator-layout notation |
|---|---|
| $\frac{dy}{d\mathbf{x}} = \left(\frac{\partial y}{\partial x_1}, \frac{\partial y}{\partial x_2}, \dots, \frac{\partial y}{\partial x_n}\right) \in \mathbb{R}^{1\times n}$ | $\frac{dy}{d\mathbf{x}} = \begin{pmatrix} \frac{\partial y}{\partial x_1} \\ \frac{\partial y}{\partial x_2} \\ \vdots \\ \frac{\partial y}{\partial x_n} \end{pmatrix} \in \mathbb{R}^{n\times 1}$ |
| $\frac{\partial \mathbf{y}}{\partial x} = \begin{pmatrix} \frac{\partial y_1}{\partial x} \\ \frac{\partial y_2}{\partial x} \\ \vdots \\ \frac{\partial y_m}{\partial x} \end{pmatrix} \in \mathbb{R}^{m\times 1}$ | $\frac{\partial \mathbf{y}}{\partial x} = (\frac{\partial y_1}{\partial x}, \frac{\partial y_2}{\partial x}, \dots, \frac{\partial y_m}{\partial x}) \in \mathbb{R}^{1\times m}$ |
| $\frac{d\mathbf{y}}{d\mathbf{x}} = \begin{pmatrix} \frac{\partial y_1}{\partial x_1} & \frac{\partial y_1}{\partial x_2} & \cdots & \frac{\partial y_1}{\partial x_n} \\ \frac{\partial y_2}{\partial x_1} & \frac{\partial y_2}{\partial x_2} & \cdots & \frac{\partial y_2}{\partial x_n} \\ \vdots & \vdots & \ddots & \vdots \\ \frac{\partial y_m}{\partial x_1} & \frac{\partial y_m}{\partial x_2} & \cdots & \frac{\partial y_m}{\partial x_n} \end{pmatrix} \in \mathbb{R}^{m\times n}$ | $\frac{d\mathbf{y}}{d\mathbf{x}} = \begin{pmatrix} \frac{\partial y_1}{\partial x_1} & \frac{\partial y_2}{\partial x_1} & \cdots & \frac{\partial y_m}{\partial x_1} \\ \frac{\partial y_1}{\partial x_2} & \frac{\partial y_2}{\partial x_2} & \cdots & \frac{\partial y_m}{\partial x_2} \\ \vdots & \vdots & \ddots & \vdots \\ \frac{\partial y_1}{\partial x_n} & \frac{\partial y_2}{\partial x_n} & \cdots & \frac{\partial y_m}{\partial x_n} \end{pmatrix} \in \mathbb{R}^{n\times m}$ |
| $\frac{dy}{d\mathbf{X}} = \begin{pmatrix} \frac{\partial y}{\partial x_{11}} & \frac{\partial y}{\partial x_{21}} & \cdots & \frac{\partial y}{\partial x_{p1}} \\ \frac{\partial y}{\partial x_{12}} & \frac{\partial y}{\partial x_{22}} & \cdots & \frac{\partial y}{\partial x_{p2}} \\ \vdots & \vdots & \ddots & \vdots \\ \frac{\partial y}{\partial x_{1q}} & \frac{\partial y}{\partial x_{2q}} & \cdots & \frac{\partial y}{\partial x_{pq}} \end{pmatrix} \in \mathbb{R}^{q\times p}$ | $\frac{dy}{d\mathbf{X}} = \begin{pmatrix} \frac{\partial y}{\partial x_{11}} & \frac{\partial y}{\partial x_{12}} & \cdots & \frac{\partial y}{\partial x_{1q}} \\ \frac{\partial y}{\partial x_{21}} & \frac{\partial y}{\partial x_{22}} & \cdots & \frac{\partial y}{\partial x_{2q}} \\ \vdots & \vdots & \ddots & \vdots \\ \frac{\partial y}{\partial x_{p1}} & \frac{\partial y}{\partial x_{p2}} & \cdots & \frac{\partial y}{\partial x_{pq}} \end{pmatrix} \in \mathbb{R}^{p\times q}$ |
| $\frac{d\mathbf{Y}}{dx} = \begin{pmatrix} \frac{\partial y_{11}}{\partial x} & \frac{\partial y_{12}}{\partial x} & \cdots & \frac{\partial y_{1n}}{\partial x} \\ \frac{\partial y_{21}}{\partial x} & \frac{\partial y_{22}}{\partial x} & \cdots & \frac{\partial y_{2n}}{\partial x} \\ \vdots & \vdots & \ddots & \vdots \\ \frac{\partial y_{m1}}{\partial x} & \frac{\partial y_{m2}}{\partial x} & \cdots & \frac{\partial y_{mn}}{\partial x} \end{pmatrix} \in \mathbb{R}^{m\times n}$ | |

### Chain Rule

For a single-valued functions $f: \mathbb{R} \to \mathbb{R}$, $g: \mathbb{R} \to \mathbb{R}$, and $h: \mathbb{R} \to \mathbb{R}$, the derivative of $h(x) = f\big(g(x)\big)$ with respect to $x$, is obtained using a chain rule:

$$\frac{dh}{dx} = \frac{df}{dg}\frac{dg}{dx} = \frac{dg}{dx}\frac{df}{dg}. \tag{2.65}$$

For the multivariable case $h(x) = f(g_1(x), g_2(x))$, the chain rule is extended as

$$\frac{dh}{dx} = \frac{df}{dg_1}\frac{dg_1}{dx} + \frac{df}{dg_2}\frac{dg_2}{dx} = \frac{dg_1}{dx}\frac{df}{dg_1} + \frac{dg_2}{dx}\frac{df}{dg_2}. \tag{2.66}$$

Visually, we can represent the two chain rules in Figure 2.1:

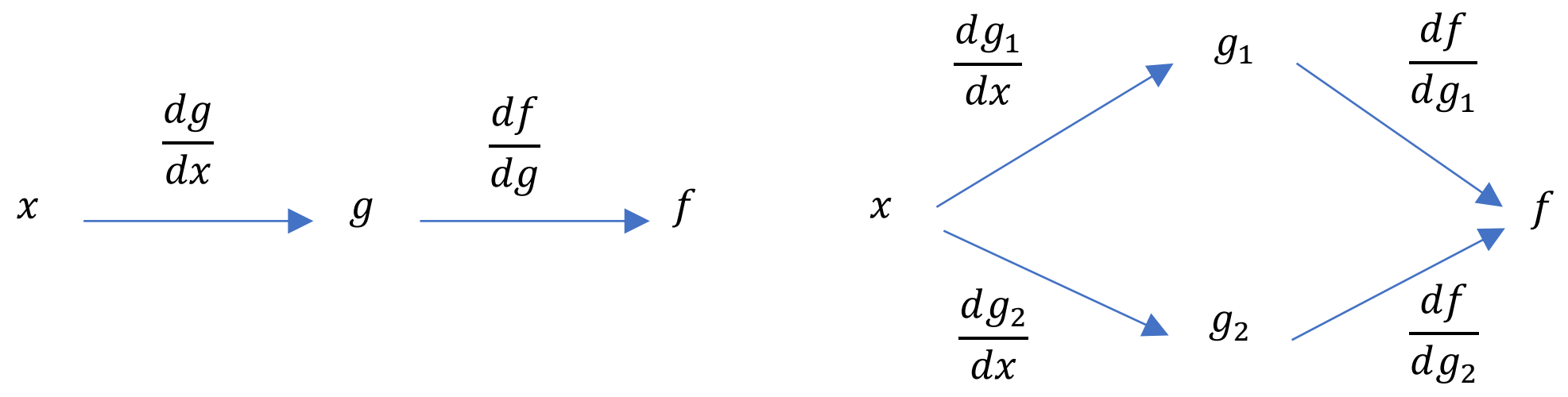

**Figure 2.1.** Chain rules are visualized.

This can be thought of as adding all components that contribute to the change of $h$. Building on this, we can extend the chain rule to also work in matrix calculus.

Consider $\mathbf{x} \in \mathbb{R}^p$, $\mathbf{y} \in \mathbb{R}^r$, $\mathbf{z} \in \mathbb{R}^n$ where $\mathbf{z}$ is a function of $\mathbf{y}$, and $\mathbf{y}$ is a function of $\mathbf{x}$; that is, $\mathbf{z} = f(\mathbf{y})$, $\mathbf{y} = g(\mathbf{x})$, and therefore $\mathbf{z} = f(g(\mathbf{x}))$. We can visualize this in Figure 2.2. Note how this figure considers the most general possible case.

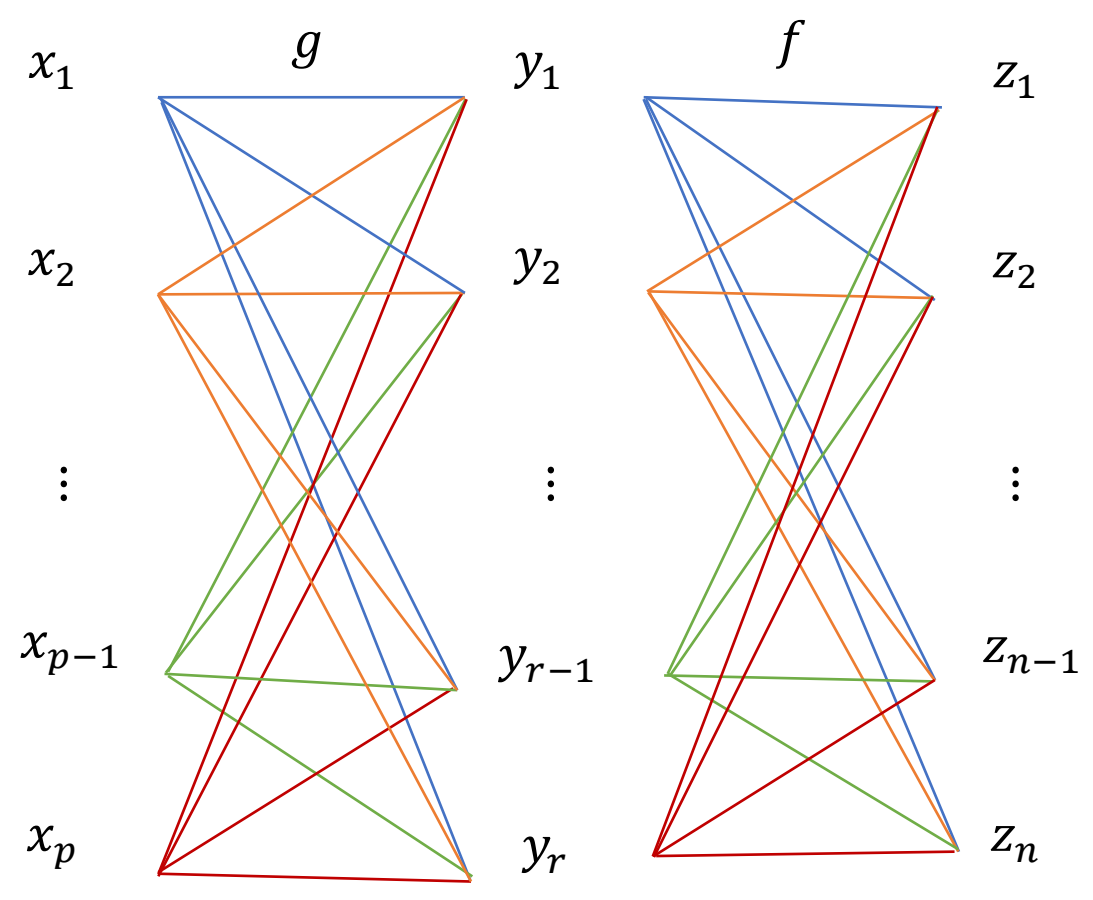

**Figure 2.2.** Visualization $\mathbf{z} = f(g(\mathbf{x}))$, where $\mathbf{z} = f(\mathbf{y})$, $\mathbf{y} = g(\mathbf{x})$.

Now we derive the chain rule for vectors in matrix calculus. Recall that we have previously defined $d\mathbf{z}/d\mathbf{x}$ as

$$\frac{d\mathbf{z}}{d\mathbf{x}} = \begin{pmatrix} \frac{dz_1}{dx_1} & \frac{dz_2}{dx_1} & \cdots & \frac{dz_n}{dx_1} \\ \frac{dz_1}{dx_2} & \frac{dz_2}{dx_2} & \cdots & \frac{dz_n}{dx_2} \\ \vdots & \vdots & \ddots & \vdots \\ \frac{dz_1}{dx_p} & \frac{dz_2}{dx_p} & \cdots & \frac{dz_n}{dx_p} \end{pmatrix} \in \mathbb{R}^{p \times n}. \tag{2.67}$$

By the chain rule,

$$\frac{dz_i}{dx_j} = \sum_{k=1}^{r} \frac{dz_i}{dy_k}\frac{dy_k}{dx_j} = \sum_{k=1}^{r} \frac{dy_k}{dx_j}\frac{dz_i}{dy_k}. \tag{2.68}$$

This directly follows from Figure 2.3, which can be obtained by isolating only $x_j$ and $z_i$ from Figure 2.2:

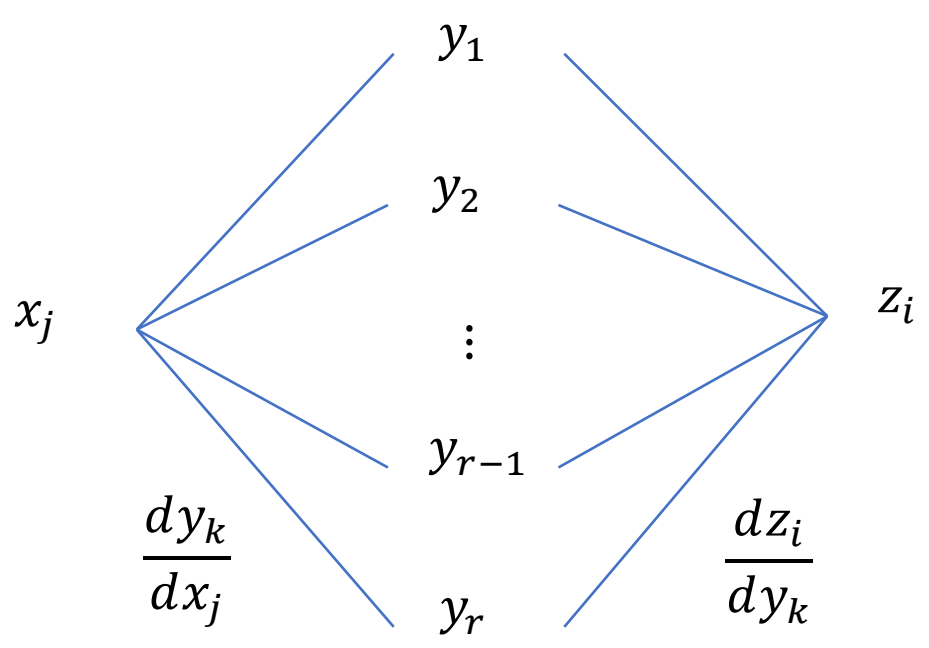

**Figure 2.3.** Chain rule visualized only considering $z_i$ and $x_j$.

Apply the scalar chain rule to each element of $d\mathbf{z}/d\mathbf{x}$. By the definition of matrix multiplication, observe that

$$\begin{aligned}
\left(\frac{d\mathbf{z}}{d\mathbf{x}}\right)^T &= \begin{pmatrix} \frac{dz_1}{dx_1} & \frac{dz_1}{dx_2} & \cdots & \frac{dz_1}{dx_p} \\ \frac{dz_2}{dx_1} & \frac{dz_2}{dx_2} & \cdots & \frac{dz_2}{dx_p} \\ \vdots & \vdots & \ddots & \vdots \\ \frac{dz_n}{dx_1} & \frac{dz_n}{dx_2} & \cdots & \frac{dz_n}{dx_p} \end{pmatrix} \in \mathbb{R}^{n\times p} \\
&= \begin{pmatrix} \sum_{k=1}^{r} \frac{dz_1}{dy_k}\frac{dy_k}{dx_1} & \sum_{k=1}^{r} \frac{dz_1}{dy_k}\frac{dy_k}{dx_2} & \cdots & \sum_{k=1}^{r} \frac{dz_1}{dy_k}\frac{dy_k}{dx_n} \\ \sum_{k=1}^{r} \frac{dz_2}{dy_k}\frac{dy_k}{dx_1} & \sum_{k=1}^{r} \frac{dz_2}{dy_k}\frac{dy_k}{dx_2} & \cdots & \sum_{k=1}^{r} \frac{dz_2}{dy_k}\frac{dy_k}{dx_n} \\ \vdots & \vdots & \ddots & \vdots \\ \sum_{k=1}^{r} \frac{dz_p}{dy_k}\frac{dy_k}{dx_1} & \sum_{k=1}^{r} \frac{dz_p}{dy_k}\frac{dy_k}{dx_2} & \cdots & \sum_{k=1}^{r} \frac{dz_p}{dy_k}\frac{dy_k}{dx_n} \end{pmatrix} \\
&= \begin{pmatrix} \frac{dz_1}{dy_1} & \frac{dz_1}{dy_2} & \cdots & \frac{dz_1}{dy_r} \\ \frac{dz_2}{dy_1} & \frac{dz_2}{dy_2} & \cdots & \frac{dz_2}{dy_r} \\ \vdots & \vdots & \ddots & \vdots \\ \frac{dz_n}{dy_1} & \frac{dz_n}{dy_2} & \cdots & \frac{dz_n}{dy_r} \end{pmatrix} \begin{pmatrix} \frac{dy_1}{dx_1} & \frac{dy_1}{dx_2} & \cdots & \frac{dy_1}{dx_p} \\ \frac{dy_2}{dx_1} & \frac{dy_2}{dx_2} & \cdots & \frac{dy_2}{dx_p} \\ \vdots & \vdots & \ddots & \vdots \\ \frac{dy_r}{dx_1} & \frac{dy_r}{dx_2} & \cdots & \frac{dy_r}{dx_p} \end{pmatrix} \\
&= \left(\frac{d\mathbf{z}}{d\mathbf{y}}\right)^T \left(\frac{d\mathbf{y}}{d\mathbf{x}}\right)^T.
\end{aligned} \tag{2.69}$$

Taking the transpose of both sides, we have that the chain rule extends to

$$\frac{d\mathbf{z}}{d\mathbf{x}} = \frac{d\mathbf{y}}{d\mathbf{x}}\frac{d\mathbf{z}}{d\mathbf{y}}. \tag{2.70}$$

Note the matrix multiplication order; $d\mathbf{y}/d\mathbf{x}$ comes first. The order did not matter for the scalar case, but we need to be mindful of the order for the matrix case.

The key idea for this derivation was to manipulate the matrices cleverly and use the scalar chain rule. When other types of derivatives are involved, this chain rule may change; some derivatives may be transposed, and the multiplication order may change. The chain rules also vary depending on how the derivatives are defined. However, the scalar chain rule must hold no matter what.

In the following, we represent some important derivative formulas in matrix calculus.

1- Suppose that $\mathbf{x}$ is an $n \times 1$ vector and $f(\mathbf{x})$ is a scalar function of the elements of $\mathbf{x}$. Then

$$\frac{\partial f}{\partial \mathbf{x}} = \begin{pmatrix} \frac{\partial f}{\partial x_1} \\ \vdots \\ \frac{\partial f}{\partial x_n} \end{pmatrix} \in \mathbb{R}^{n \times 1}. \tag{2.71}$$

2- Suppose $\mathbf{x}$ and $\mathbf{y}$ are $n$-element column vectors. Then

$$\mathbf{x}^T\mathbf{y} = \langle \mathbf{x}|\mathbf{y} \rangle = x_1y_1 + x_2y_2 + \cdots + x_ny_n. \tag{2.72}$$

Hence, we have

$$\begin{aligned} \frac{\partial (\mathbf{x}^T\mathbf{y})}{\partial \mathbf{x}} = \frac{\partial \langle \mathbf{x}|\mathbf{y} \rangle}{\partial \mathbf{x}} &= \begin{pmatrix} \frac{\partial \langle \mathbf{x}|\mathbf{y} \rangle}{\partial x_1} \\ \vdots \\ \frac{\partial \langle \mathbf{x}|\mathbf{y} \rangle}{\partial x_n} \end{pmatrix} \\ &= \begin{pmatrix} \frac{\partial (x_1y_1 + x_2y_2 + \cdots + x_ny_n)}{\partial x_1} \\ \vdots \\ \frac{\partial (x_1y_1 + x_2y_2 + \cdots + x_ny_n)}{\partial x_n} \end{pmatrix} \\ &= \begin{pmatrix} y_1 \\ \vdots \\ y_n \end{pmatrix} \\ &= \mathbf{y}. \end{aligned} \tag{2.73}$$

3- Also, we can obtain

$$\begin{aligned} \frac{\partial (\mathbf{x}^T\mathbf{y})}{\partial \mathbf{y}} = \frac{\partial \langle \mathbf{x}|\mathbf{y} \rangle}{\partial \mathbf{y}} &= \begin{pmatrix} \frac{\partial \langle \mathbf{x}|\mathbf{y} \rangle}{\partial y_1} \\ \vdots \\ \frac{\partial \langle \mathbf{x}|\mathbf{y} \rangle}{\partial y_n} \end{pmatrix} \\ &= \begin{pmatrix} \frac{\partial (x_1y_1 + x_2y_2 + \cdots + x_ny_n)}{\partial y_1} \\ \vdots \\ \frac{\partial (x_1y_1 + x_2y_2 + \cdots + x_ny_n)}{\partial y_n} \end{pmatrix} \\ &= \begin{pmatrix} x_1 \\ \vdots \\ x_n \end{pmatrix} \\ &= \mathbf{x}. \end{aligned} \tag{2.74}$$

4- Now we will compute the partial derivative of a quadratic form $\mathbf{x}^T\mathbf{A}\mathbf{x}$ with respect to a vector. First, write the quadratic form as follows

$$\begin{aligned} \mathbf{x}^T\mathbf{A}\mathbf{x} = \langle \mathbf{x}|\mathbf{A}|\mathbf{x} \rangle &= (x_1, x_2, \dots, x_n) \begin{pmatrix} A_{11} & \cdots & A_{1n} \\ \vdots & \ddots & \vdots \\ A_{n1} & \cdots & A_{nn} \end{pmatrix} \begin{pmatrix} x_1 \\ \vdots \\ x_n \end{pmatrix} \\ &= \left( \sum_{i=1}^{n} x_iA_{i1}, \sum_{i=1}^{n} x_iA_{i2}, \dots, \sum_{i=1}^{n} x_iA_{in} \right) \begin{pmatrix} x_1 \\ \vdots \\ x_n \end{pmatrix} \\ &= \sum_{j=1}^{n} \sum_{i=1}^{n} x_ix_jA_{ij}. \end{aligned} \tag{2.75}$$

Now take the partial derivative of the quadratic as follows:

$$\frac{\partial \mathbf{x}^T\mathbf{A}\mathbf{x}}{\partial \mathbf{x}} = \frac{\partial\langle\mathbf{x}|\mathbf{A}|\mathbf{x}\rangle}{\partial \mathbf{x}} = \begin{pmatrix} \frac{\partial\langle\mathbf{x}|\mathbf{A}|\mathbf{x}\rangle}{\partial x_1} \\ \vdots \\ \frac{\partial\langle\mathbf{x}|\mathbf{A}|\mathbf{x}\rangle}{\partial x_n} \end{pmatrix}$$
$$= \begin{pmatrix} \frac{\partial\left(\sum_{j=1}^{n}\sum_{i=1}^{n} x_i x_j A_{ij}\right)}{\partial x_1} \\ \vdots \\ \frac{\partial\left(\sum_{j=1}^{n}\sum_{i=1}^{n} x_i x_j A_{ij}\right)}{\partial x_n} \end{pmatrix}$$
$$= \begin{pmatrix} \frac{\partial \sum_{j=1}^{n} x_1 x_j A_{1j}}{\partial x_1} + \frac{\partial \sum_{i=1}^{n} x_i x_1 A_{i1}}{\partial x_1} \\ \vdots \\ \frac{\partial \sum_{j=1}^{n} x_n x_j A_{nj}}{\partial x_n} + \frac{\partial \sum_{i=1}^{n} x_i x_n A_{in}}{\partial x_n} \end{pmatrix}$$
$$= \begin{pmatrix} \sum_{j=1}^{n} x_j A_{1j} + \sum_{i=1}^{n} x_i A_{i1} \\ \vdots \\ \sum_{j=1}^{n} x_j A_{nj} + \sum_{i=1}^{n} x_i A_{in} \end{pmatrix}$$
$$= \begin{pmatrix} \sum_{j=1}^{n} x_j A_{1j} \\ \vdots \\ \sum_{j=1}^{n} x_j A_{nj} \end{pmatrix} + \begin{pmatrix} \sum_{i=1}^{n} x_i A_{i1} \\ \vdots \\ \sum_{i=1}^{n} x_i A_{in} \end{pmatrix} = \mathbf{A}\mathbf{x} + \mathbf{A}^T\mathbf{x} = (\mathbf{A} + \mathbf{A}^T)\mathbf{x}. \tag{2.76}$$

5- If $\mathbf{A}$ is symmetric, then $\mathbf{A} = \mathbf{A}^T$ and the above expression simplifies to

$$\frac{\partial\langle\mathbf{x}|\mathbf{A}|\mathbf{x}\rangle}{\partial \mathbf{x}} = 2\mathbf{A}\mathbf{x}. \tag{2.77}$$

6- Suppose $\mathbf{A}$ is an $m \times n$ matrix, $\mathbf{x}$ is an $n \times 1$ vector, and $\mathbf{y} = \mathbf{A}\mathbf{x}$;

$$\frac{\partial \mathbf{y}}{\partial \mathbf{x}} = \frac{\partial \mathbf{A}\mathbf{x}}{\partial \mathbf{x}} = \begin{pmatrix} \frac{dy_1}{dx_1} & \frac{dy_2}{dx_1} & \cdots & \frac{dy_m}{dx_1} \\ \frac{dy_1}{dx_2} & \frac{dy_2}{dx_2} & \cdots & \frac{dy_m}{dx_2} \\ \vdots & \vdots & \ddots & \vdots \\ \frac{dy_1}{dx_n} & \frac{dy_2}{dx_n} & \cdots & \frac{dy_m}{dx_n} \end{pmatrix} = \begin{pmatrix} A_{11} & A_{21} & \cdots & A_{m1} \\ A_{12} & A_{22} & \cdots & A_{m2} \\ \vdots & \vdots & \ddots & \vdots \\ A_{1n} & A_{2n} & \cdots & A_{mn} \end{pmatrix} = \mathbf{A}^T. \tag{2.78}$$

7- Now we suppose that $\mathbf{A}$ is an $m \times n$ matrix, $\mathbf{B}$ is an $n \times n$ matrix, and we want to compute the partial derivative of $\mathrm{Tr}(\mathbf{A}\mathbf{B}\mathbf{A}^T)$ with respect to $\mathbf{A}$. First, compute $\mathbf{A}\mathbf{B}\mathbf{A}^T$ as follows:

$$\mathbf{A}\mathbf{B}\mathbf{A}^T = \begin{pmatrix} A_{11} & \cdots & A_{1n} \\ \vdots & \ddots & \vdots \\ A_{m1} & \cdots & A_{mn} \end{pmatrix} \begin{pmatrix} B_{11} & \cdots & B_{1n} \\ \vdots & \ddots & \vdots \\ B_{n1} & \cdots & B_{nn} \end{pmatrix} \begin{pmatrix} A_{11} & \cdots & A_{m1} \\ \vdots & \ddots & \vdots \\ A_{1n} & \cdots & A_{mn} \end{pmatrix}$$
$$= \begin{pmatrix} \sum_{k,j} A_{1k}B_{kj}A_{1j} & \cdots & \sum_{k,j} A_{1k}B_{kj}A_{mj} \\ \vdots & \ddots & \vdots \\ \sum_{k,j} A_{mk}B_{kj}A_{1j} & \cdots & \sum_{k,j} A_{mk}B_{kj}A_{mj} \end{pmatrix}. \tag{2.79}$$

From this we see that the trace of $\mathbf{A}\mathbf{B}\mathbf{A}^T$ is given as

$$\mathrm{Tr}(\mathbf{ABA}^T) = \sum_{i,j,k} A_{ik} B_{kj} A_{ij}. \tag{2.80}$$

Its partial derivative with respect to **A** can be computed as

$$\frac{\partial \mathrm{Tr}(\mathbf{ABA}^T)}{\partial \mathbf{A}} = \begin{pmatrix} \frac{\partial \sum_{i,j,k} A_{ik}B_{kj}A_{ij}}{\partial A_{11}} & \cdots & \frac{\partial \sum_{i,j,k} A_{ik}B_{kj}A_{ij}}{\partial A_{1n}} \\ \vdots & \ddots & \vdots \\ \frac{\partial \sum_{i,j,k} A_{ik}B_{kj}A_{ij}}{\partial A_{m1}} & \cdots & \frac{\partial \sum_{i,j,k} A_{ik}B_{kj}A_{ij}}{\partial A_{mn}} \end{pmatrix}$$
$$= \begin{pmatrix} \frac{\partial \sum_j A_{11}B_{1j}A_{1j}}{\partial A_{11}} + \frac{\partial \sum_k A_{1k}B_{k1}A_{11}}{\partial A_{11}} & \cdots & \frac{\partial \sum_j A_{1n}B_{nj}A_{1j}}{\partial A_{1n}} + \frac{\partial \sum_k A_{1k}B_{kn}A_{1n}}{\partial A_{1n}} \\ \vdots & \ddots & \vdots \\ \frac{\partial \sum_j A_{m1}B_{1j}A_{mj}}{\partial A_{m1}} + \frac{\partial \sum_k A_{mk}B_{k1}A_{m1}}{\partial A_{m1}} & \cdots & \frac{\partial \sum_j A_{mn}B_{nj}A_{mj}}{\partial A_{mn}} + \frac{\partial \sum_k A_{mk}B_{kn}A_{mn}}{\partial A_{mn}} \end{pmatrix}$$
$$= \begin{pmatrix} \sum_j A_{1j}B_{1j} + \sum_k A_{1k}B_{k1} & \cdots & \sum_j A_{1j}B_{nj} + \sum_k A_{1k}B_{kn} \\ \vdots & \ddots & \vdots \\ \sum_j A_{mj}B_{1j} + \sum_k A_{mk}B_{k1} & \cdots & \sum_j A_{mj}B_{nj} + \sum_k A_{mk}B_{kn} \end{pmatrix}$$
$$= \begin{pmatrix} \sum_j A_{1j}B_{1j} & \cdots & \sum_j A_{1j}B_{nj} \\ \vdots & \ddots & \vdots \\ \sum_j A_{mj}B_{1j} & \cdots & \sum_j A_{mj}B_{nj} \end{pmatrix} + \begin{pmatrix} \sum_k A_{1k}B_{k1} & \cdots & \sum_k A_{1k}B_{kn} \\ \vdots & \ddots & \vdots \\ \sum_k A_{mk}B_{k1} & \cdots & \sum_k A_{mk}B_{kn} \end{pmatrix}$$
$$= \mathbf{AB}^T + \mathbf{AB}$$
$$= \mathbf{A}(\mathbf{B}^T + \mathbf{B}). \tag{2.81}$$

8- If **B** is symmetric, then this can be simplified to

$$\frac{\partial \mathrm{Tr}(\mathbf{ABA}^T)}{\partial \mathbf{A}} = 2\mathbf{AB}. \tag{2.82}$$

## 2.5 Numerical Differentiation and Finite Difference Methods

The process of estimating derivatives numerically is referred to as numerical differentiation [8]. Estimates can be derived in different ways from function evaluations. This section discusses finite difference methods.

The derivative $f^{(1)}(x)$ of a function $f(x)$ at the point $x = a$ is defined by:

$$\left.\frac{df(x)}{dx}\right|_{x=a} = f^{(1)}(a) = \lim_{x\to a} \frac{f(x) - f(a)}{x - a}. \tag{2.83}$$

Graphically, the definition is illustrated in Figure 2.4. The derivative is the value of the slope of the tangent line to the function at $x = a$. In finite difference approximations of the derivative, values of the function at different points in the neighborhood of the point $x = a$ are used for estimating the slope. The forward, backward, and central finite difference formulas are the simplest finite difference approximations of the derivative. In these approximations, illustrated in Figure 2.5, the derivative at the point $x_i$ is calculated from the values of two points. The derivative is estimated as the value of the slope of the line that connects the two points.

1- The forward difference is the slope of the line that connects points $(x_i, f(x_i))$ and $(x_{i+1}, f(x_{i+1}))$:

$$\left.\frac{df}{dx}\right|_{x=x_i} = \frac{f(x_{i+1}) - f(x_i)}{x_{i+1} - x_i}. \tag{2.84}$$

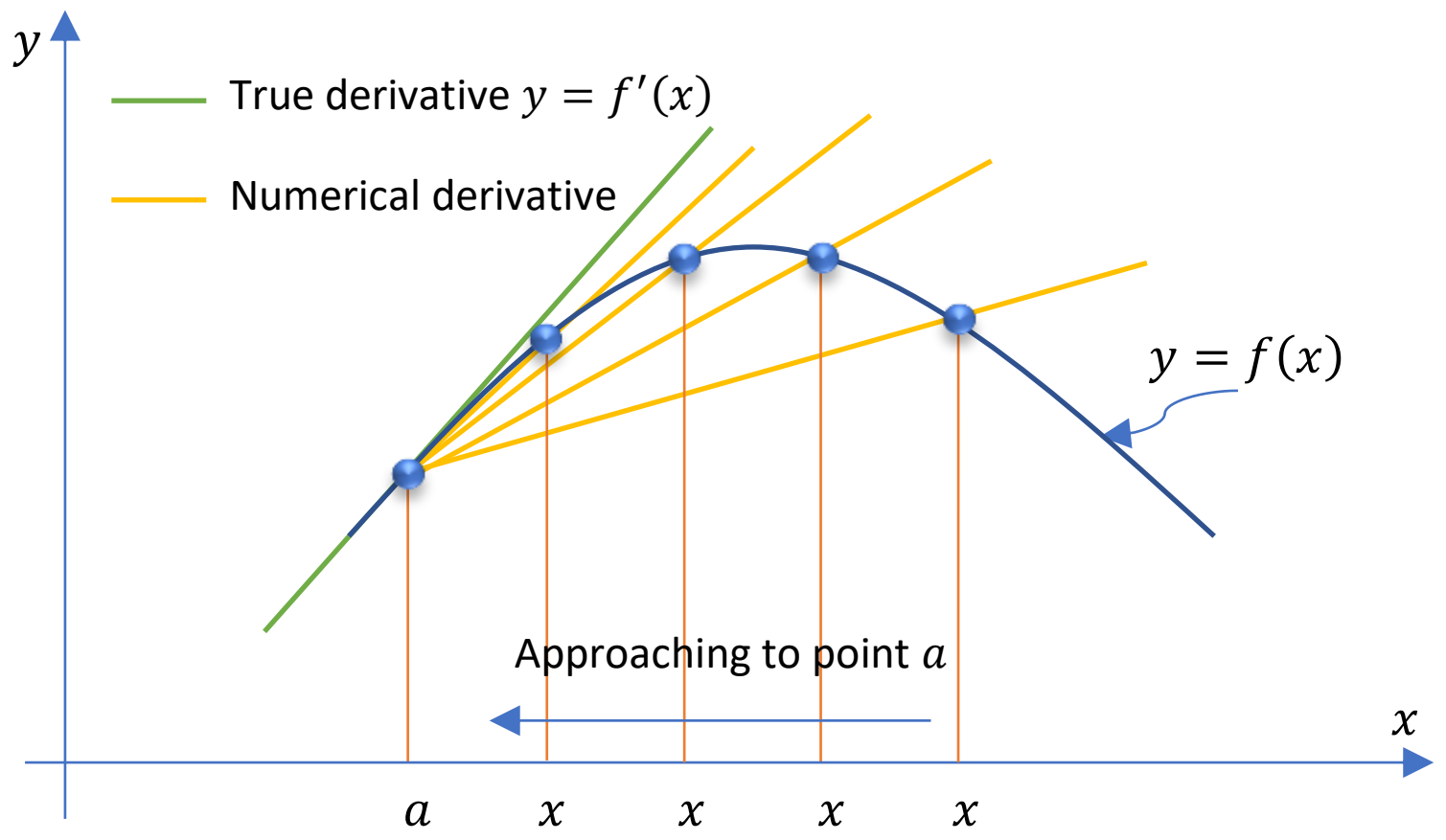

**Figure 2.4.** Definition of the derivative.

2- The backward difference is the slope of the line that connects points $(x_{i-1}, f(x_{i-1}))$ and $(x_i, f(x_i))$:

$$\left.\frac{df}{dx}\right|_{x=x_i} = \frac{f(x_i) - f(x_{i-1})}{x_i - x_{i-1}}. \tag{2.85}$$

3- The central difference is the slope of the line that connects points $(x_{i-1}, f(x_{i-1}))$ and $(x_{i+1}, f(x_{i+1}))$:

$$\left.\frac{df}{dx}\right|_{x=x_i} = \frac{f(x_{i+1}) - f(x_{i-1})}{x_{i+1} - x_{i-1}}. \tag{2.86}$$

The forward, backward, and central difference formulas, as well as many other finite difference formulas for approximating derivatives, can be derived by using Taylor series expansion. One advantage of using Taylor series expansion for deriving the formulas is that it also provides an estimate for the truncation error in the approximation. All the formulas derived in this section are for the case where the points are equally spaced.

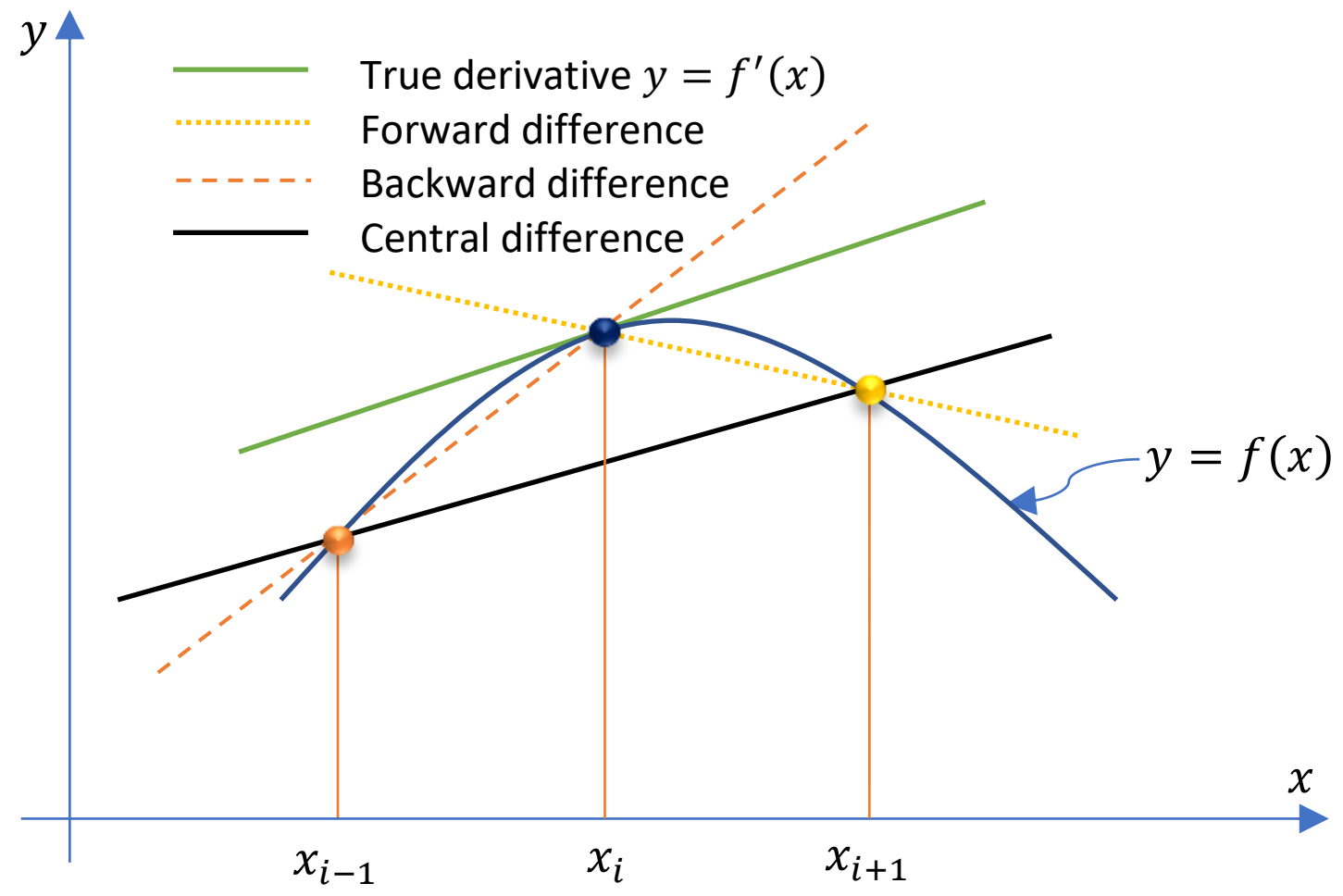

**Figure 2.5.** Finite difference approximation of derivative.

### Two-Point Forward Difference Formula for First Derivative

The value of a function at a point $x_{i+1}$ can be approximated by a Taylor series in terms of the value of the function and its derivatives at a point $x_i$:

$$f(x_{i+1}) = f(x_i) + f^{(1)}(x_i)h + \frac{f^{(2)}(x_i)}{2!}h^2 + \frac{f^{(3)}(x_i)}{3!}h^3 + \frac{f^{(4)}(x_i)}{4!}h^4 + \cdots, \tag{2.87}$$

where $h = x_{i+1} - x_i$ is the spacing between the points. By using a two-term Taylor series expansion with a remainder, (2.87) can be rewritten as:

$$f(x_{i+1}) = f(x_i) + f^{(1)}(x_i)h + \frac{f^{(2)}(\xi)}{2!}h^2, \tag{2.88}$$

where $\xi$ is a value of $x$ between $x_i$ and $x_{i+1}$. Solving (2.88) for $f^{(1)}(x_i)$ yields:

$$f^{(1)}(x_i) = \frac{f(x_{i+1}) - f(x_i)}{h} - \frac{f^{(2)}(\xi)}{2!}h. \tag{2.89}$$

An approximate value of the derivative $f^{(1)}(x_i)$ can now be calculated if the second term on the right-hand side of (2.89) is ignored. Ignoring this second term introduces a truncation (discretization) error. Since this term is proportional to $h$, the truncation error is said to be on the order of $h$ (written as $O(h)$):

$$\text{truncation error } = \frac{f^{(2)}(\xi)}{2!}h = O(h). \tag{2.90}$$

It should be pointed out here that the magnitude of the truncation error is not really known since the value of $f^{(2)}(\xi)$ is not known. Nevertheless, (2.90) is valuable since it implies that a smaller $h$ gives a smaller error. Using the notation of (2.90), the approximated value of the first derivative is:

$$f^{(1)}(x_i) = \frac{f(x_{i+1}) - f(x_i)}{h} + O(h). \tag{2.91}$$

### Two-Point Backward Difference Formula for First Derivative

The backward difference formula can also be derived by the application of Taylor series expansion. The value of the function at a point $x_{i-1}$ is approximated by a Taylor series:

$$f(x_{i-1}) = f(x_i) - f^{(1)}(x_i)h + \frac{f^{(2)}(x_i)}{2!}h^2 - \frac{f^{(3)}(x_i)}{3!}h^3 + \frac{f^{(4)}(x_i)}{4!}h^4 - \cdots, \tag{2.92}$$

where $h = x_i - x_{i-1}$. By using a two-term Taylor series expansion with a remainder, (2.92) can be rewritten as:

$$f(x_{i-1}) = f(x_i) - f^{(1)}(x_i)h + \frac{f^{(2)}(\xi)}{2!}h^2, \tag{2.93}$$

where $\xi$ is a value of $x$ between $x_{i-1}$ and $x_i$. Solving (2.93) for $f'(x_i)$ yields:

$$f^{(1)}(x_i) = \frac{f(x_i) - f(x_{i-1})}{h} - \frac{f^{(2)}(\xi)}{2!}h. \tag{2.94}$$

An approximate value of the derivative, $f^{(1)}(x_i)$, can be calculated if the second term on the right-hand side of (2.94) is ignored. This yield:

$$f^{(1)}(x_i) = \frac{f(x_i) - f(x_{i-1})}{h} + O(h). \tag{2.95}$$

### Two-Point Central Difference Formula for First Derivative

The central difference formula can be derived by using three terms in the Taylor series expansion and a remainder. The value of the function at a point $x_{i+1}$ in terms of the value of the function and its derivatives at a point $x_i$ is given by:

$$f(x_{i+1}) = f(x_i) + f^{(1)}(x_i)h + \frac{f^{(2)}(x_i)}{2!}h^2 + \frac{f^{(3)}(\xi_1)}{3!}h^3, \tag{2.96}$$

where $\xi_1$ is a value of $x$ between $x_i$ and $x_{i+1}$. The value of the function at a point $x_{i-1}$ is given by:

$$f(x_{i-1}) = f(x_i) - f^{(1)}(x_i)h + \frac{f^{(2)}(x_i)}{2!}h^2 - \frac{f^{(3)}(\xi_2)}{3!}h^3, \tag{2.97}$$

where $\xi_2$ is a value of $x$ between $x_{i-1}$ and $x_i$. In the last two equations, the spacing of the intervals is taken to be equal so that $h = x_{i+1} - x_i = x_i - x_{i-1}$. Subtracting (2.97) from (2.96) gives:

$$f(x_{i+1}) - f(x_{i-1}) = 2f^{(1)}(x_i)h + \frac{f^{(3)}(\xi_1)}{3!}h^3 + \frac{f^{(3)}(\xi_2)}{3!}h^3. \tag{2.98}$$

An estimate for the first derivative is obtained by solving (2.98) for $f^{(1)}(x_i)$ while neglecting the remainder terms, which introduces a truncation error, which is of the order of $h^2$:

$$f^{(1)}(x_i) = \frac{f(x_{i+1}) - f(x_{i-1})}{2h} + O(h^2). \tag{2.99}$$

**Table 2.3.** List of difference formulas that can be used for numerical evaluation of first derivatives.

| Method | First Derivative Formula | Truncation Error |
|---|---|---|
| Two-point forward difference | $f'(x_i) = \frac{f(x_{i+1}) - f(x_i)}{h}$ | $O(h)$ |
| Three-point forward difference | $f'(x_i) = \frac{-3f(x_i) + 4f(x_{i+1}) - f(x_{i+2})}{2h}$ | $O(h^2)$ |
| Two-point backward difference | $f'(x_i) = \frac{f(x_i) - f(x_{i-1})}{h}$ | $O(h)$ |
| Three-point backward difference | $f'(x_i) = \frac{f(x_{i-2}) - 4f(x_{i-1}) + 3f(x_i)}{2h}$ | $O(h^2)$ |
| Two-point central difference | $f'(x_i) = \frac{f(x_{i+1}) - f(x_{i-1})}{2h}$ | $O(h^2)$ |
| Four-point central difference | $f'(x_i) = \frac{f(x_{i-2}) - 8f(x_{i-1}) + 8f(x_{i+1}) - f(x_{i+2})}{12h}$ | $O(h^4)$ |

A comparison of (2.91), (2.95), and (2.99) shows that in the forward and backward difference approximation, the truncation error is of the order $h$, while in the central difference approximation, the truncation error is of the order $h^2$. This indicates that the central difference approximation gives a more accurate approximation of the derivative.

### Finite Difference Formulas for the Second Derivative

The same approach used to develop finite difference formulas for the first derivative can be used to develop expressions for higher-order derivatives. For example, for points $x_{i+1}$, and $x_{i-1}$, the four-term Taylor series expansion with a remainder is

$$f(x_{i+1}) = f(x_i) + f^{(1)}(x_i)h + \frac{f^{(2)}(x_i)}{2!}h^2 + \frac{f^{(3)}(x_i)}{3!}h^3 + \frac{f^{(4)}(\xi_1)}{4!}h^4, \tag{2.100}$$

$$f(x_{i-1}) = f(x_i) - f^{(1)}(x_i)h + \frac{f^{(2)}(x_i)}{2!}h^2 - \frac{f^{(2)}(x_i)}{3!}h^3 + \frac{f^{(4)}(\xi_2)}{4!}h^4, \tag{2.101}$$

where $\xi_1$ is a value of $x$ between $x_i$ and $x_{i+1}$, and $\xi_2$ is a value of $x$ between $x_i$ and $x_{i-1}$. Adding (2.100) and (2.101) gives

$$f(x_{i+1}) - f(x_{i-1}) = 2f(x_i) + 2\frac{f^{(2)}(x_i)}{2!}h^2 + \frac{f^{(4)}(\xi_1)}{4!}h^4 + \frac{f^{(4)}(\xi_2)}{4!}h^4. \tag{2.102}$$

An estimate for the second derivative can be obtained by solving (2.102) for $f^{(2)}(x_i)$ while neglecting the remainder terms. This introduces a truncation error of the order $h^2$.

$$f^{(2)}(x_i) = \frac{f(x_{i-1}) - 2f(x_i) + f(x_{i+1})}{h^2} + O(h^2). \tag{2.103}$$

Equation (2.103) is the three-point central difference formula with a truncation error of $O(h^2)$. The same procedure can be used to develop the three-point forward and backward difference formulas for the second derivative. Tables 2.3 and 2.4 list difference formulas of various accuracy that can be used for numerical evaluation of first and second derivatives.

### Numerical Partial Differentiation

For a function of several independent variables, the partial derivative of the function with respect to one of the variables represents the rate of change of the value of the function with respect to this variable, while all the other variables are kept constant. For a function $f(x, y)$ with two independent variables, the partial derivatives with respect to $x$ and $y$ at the point $(a, b)$ are defined as

$$\left.\frac{\partial f(x,y)}{\partial x}\right|_{\substack{x=a \\ y=b}} = \lim_{x \to a} \frac{f(x,b) - f(a,b)}{x - a}, \tag{2.104}$$

$$\left.\frac{\partial f(x,y)}{\partial y}\right|_{\substack{x=a \\ y=b}} = \lim_{y \to b} \frac{f(a,y) - f(a,b)}{y - b}. \tag{2.105}$$

**Table 2.4.** List of difference formulas that can be used for numerical evaluation of second derivatives.

| Method | Second Derivative Formula | Truncation Error |
|---|---|---|
| Three-point forward difference | $f''(x_i) = \frac{f(x_i) - 2f(x_{i+1}) + f(x_{i+2})}{h^2}$ | $O(h)$ |
| Four-point forward difference | $f''(x_i) = \frac{2f(x_i) - 5f(x_{i+1}) + 4f(x_{i+2}) - f(x_{i+3})}{h^2}$ | $O(h^2)$ |
| Three-point backward difference | $f''(x_i) = \frac{f(x_{i-2}) - 2f(x_{i-1}) + f(x_i)}{h^2}$ | $O(h)$ |
| Four-point backward difference | $f''(x_i) = \frac{-f(x_{i-3}) + 4f(x_{i-2}) - 5f(x_{i-1}) + 2f(x_i)}{h^2}$ | $O(h^2)$ |
| Three-point central difference | $f''(x_i) = \frac{f(x_{i-1}) - 2f(x_i) + f(x_{i+1})}{h^2}$ | $O(h^2)$ |
| Five-point central difference | $f''(x_i) = \frac{-f(x_{i-2}) + 16f(x_{i-1}) - 30f(x_i) + 16f(x_{i+1}) - f(x_{i+2})}{12h^2}$ | $O(h^4)$ |

This means that the finite difference formulas that are used for approximating the derivatives of functions with one independent variable can be adopted for calculating partial derivatives. The formulas are applied for one of the variables, while the other variables are kept constant. For example, consider a function of two independent variables $f(x, y)$ specified as a set of discrete $m \cdot n$ points $(x_1, y_1)$, $(x_1, y_2)$, ..., $(x_n, y_m)$. The spacing between the points in each direction is constant such that $h_x = x_{i+1} - x_i$ and $h_y = y_{i+1} - y_i$. Figure 2.6 shows a case where $n = 5$ and $m = 4$. An approximation for the partial derivative at a point $(x_i, y_i)$ with the two-point forward difference formula is

$$\left.\frac{\partial f}{\partial x}\right|_{\substack{x=x_i \\ y=y_i}} = \frac{f(x_{i+1}, y_i) - f(x_i, y_i)}{h_x}, \tag{2.106}$$

$$\left.\frac{\partial f}{\partial y}\right|_{\substack{x=x_i \\ y=y_i}} = \frac{f(x_i, y_{i+1}) - f(x_i, y_i)}{h_y}. \tag{2.107}$$

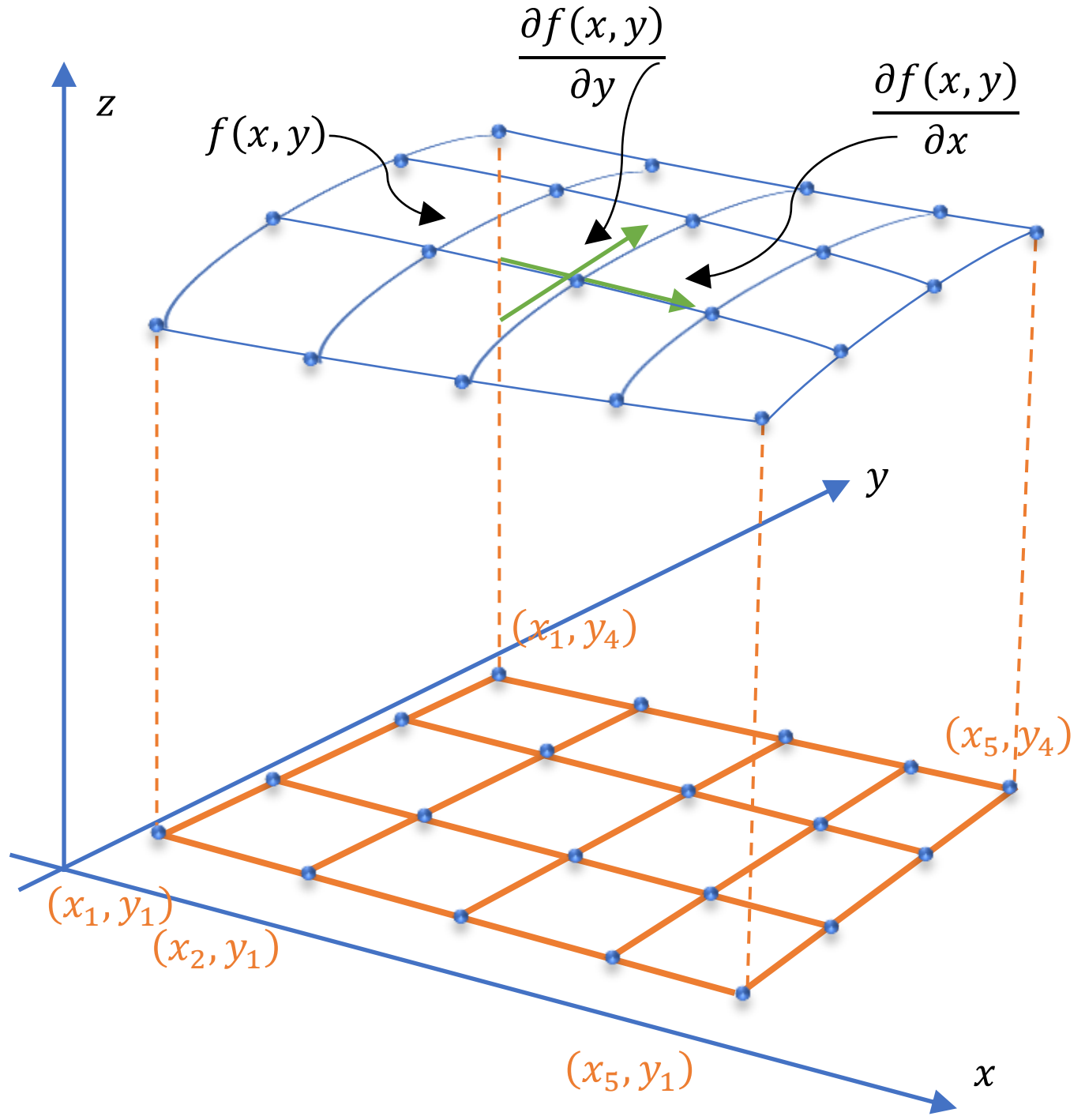

**Figure 2.6.** A function with two independent variables.

In the same way, the two-point backward and central difference formulas are

$$\left.\frac{\partial f}{\partial x}\right|_{\substack{x=x_i\\y=y_i}} = \frac{f(x_i, y_i) - f(x_{i-1}, y_i)}{h_x}, \tag{2.108}$$

$$\left.\frac{\partial f}{\partial y}\right|_{\substack{x=x_i\\y=y_i}} = \frac{f(x_i, y_i) - f(x_i, y_{i-1})}{h_y}, \tag{2.109}$$

$$\left.\frac{\partial f}{\partial x}\right|_{\substack{x=x_i\\y=y_i}} = \frac{f(x_{i+1}, y_i) - f(x_{i-1}, y_i)}{2h_x}, \tag{2.110}$$

$$\left.\frac{\partial f}{\partial x}\right|_{\substack{x=x_i\\y=y_i}} = \frac{f(x_i, y_{i+1}) - f(x_i, y_{i-1})}{2h_y}. \tag{2.111}$$

Moreover, the second partial derivatives with the three-point central difference formula are

$$\left.\frac{\partial^2 f}{\partial x^2}\right|_{\substack{x=x_i\\y=y_i}} = \frac{f(x_{i-1}, y_i) - 2f(x_i, y_i) + f(x_{i+1}, y_i)}{h_x^2}, \tag{2.112}$$

$$\left.\frac{\partial^2 f}{\partial y^2}\right|_{\substack{x=x_i\\y=y_i}} = \frac{f(x_i, y_{i-1}) - 2f(x_i, y_i) + f(x_i, y_{i+1})}{h_y^2}. \tag{2.113}$$

The second-order mixed four-point central finite difference formula is

$$\left.\frac{\partial^2 f}{\partial x \partial y}\right|_{\substack{x=x_i\\y=y_i}} = \frac{[f(x_{i+1}, y_{i+1}) - f(x_{i-1}, y_{i+1})] - [f(x_{i+1}, y_{i-1}) - f(x_{i-1}, y_{i-1})]}{2h_x \cdot 2h_y}. \tag{2.114}$$

***Mathematica Code 2.1***

```
Input [1]    (* Study the relation between approximation and exact solutions as a function
             of step size *)

             (*
             g(x):function
             x0:value at which the solution is desired
             h:step size value
             n:number of times step size is halved
             AV:Approximate value of the first derivative using Three point Forward
             Difference Approximation TPFDA
             Ev:Exact value of the first derivative
             Et:True error
             et:Absolute relative true percentage error
             *)

             (* Taking initial inputs from user *)
             x0=Input["Enter the intial point: For example 4 "] ;
             h=Input["Enter the step size : For example 0.4 "];
             n=Input["number of times step size is halved: For example 11 "];

             If[
               n<=0||h<=0,
               Beep[];
               MessageDialog[" n and h have to be postive number: "];
               Exit[];
               ];

             (* Taking the function from user *)
             g[x_] = Evaluate[Input["Please input a function of x and y to find the minimum:
             For example x*Exp[2*x] "]];

             (* Defination of the Three point Forward Difference Approximation function *)
             TPFDA[g_,x0_,h_]:=Module[
                {deriv},
                deriv=(-3*g[x0]+4*g[x0+h]-g[x0+2*h])/(2*h);
                deriv
                ];

             Ev=N[g'[x0]];

             (* Starting the calculations *)
             Do[
               Nn[i]=2^i;
               newh[i]=h/Nn[i];
               AV[i]=TPFDA[g,x0,newh[i]];
               Et[i]=Ev-AV[i];
               et[i]=Abs[(Et[i]/Ev)]*100.0;,
               {i,0,n-1,1}
               ];

             (* Results of each iteration *)
             data1=Table[
                {newh[i],AV[i],Et[i],et[i]},
                {i,0,n-1}
                ];

             Grid[
              Prepend[
               data1,
               {"step (h)","AV","Et","et"}
```

```
  ],
 Background->{None,{LightGray,{White,LightBlue}}}
 ]

(* Data visualization *)
data2=Table[
   {newh[i],AV[i]},
   {i,0,n-1}
   ];

plot2=ListLinePlot[
   data2,
   PlotStyle->Blue,
   Mesh->All,
   PlotRange->Full,
   PlotLabels->{"Approximate"}
   ];

data3=Table[
   {newh[i],Ev},
   {i,0,n-1}
   ];

plot3=ListLinePlot[
   data3,
   PlotStyle->Red,
   Mesh->All,
   PlotRange->Full,
   PlotLabels->{"Exact"}
   ];

Show[
 {plot2,plot3},
 PlotLabel->"Approximate Solution of the First Derivative using \n
Three point Forward Difference Approximation",
 AxesLabel->{"Step Size"}
 ]

data4=Table[
   {newh[i],et[i]},
   {i,0,n-1}
   ];

plot4=ListLinePlot[
  data4,
  PlotLabel->"Absolute Relative True Percentage Error \n
as a Function of Step Size",
  AxesLabel->{"Step Size","Absolute percentage error"},
  PlotStyle->Blue,
  Mesh->All,
  PlotRange->Full
  ]

(* Data Manipulation *)
Manipulate[
 ListPlot[
  data3,
  LabelStyle->Directive[Black,14],
  PlotStyle->Black,
  PlotLabels->Placed["Exact",Above],
  Joined->True,
```

```
      AxesLabel->{"Step","Apprx Sol."},
      Mesh->All,
      Epilog->{
        PointSize[0.013],
        Green,
        Arrow[{data2[[i]],data2[[i+1]]}],
        Red,
        Point[Table[data2[[j]],{j,1,i}]],
        Blue,
        Line[Table[data2[[j]],{j,1,i}]]
        }
      ],
    {i,1,n-2,1}
    ]
```

Output[1]

| step (h) | AV | Et | et |
|---|---|---|---|
| 0.4 | 12650.3 | 14178.3 | 52.8477 |
| 0.2 | 24371.4 | 2457.24 | 9.15903 |
| 0.1 | 26312.0 | 516.578 | 1.92547 |
| 0.05 | 26709.9 | 118.715 | 0.442493 |
| 0.025 | 26800.1 | 28.4721 | 0.106126 |
| 0.0125 | 26821.6 | 6.97288 | 0.0259904 |
| 0.00625 | 26826.9 | 1.72542 | 0.00643126 |
| 0.003125 | 26828.2 | 0.42915 | 0.0015996 |
| 0.0015625 | 26828.5 | 0.107013 | 0.000398878 |
| 0.00078125 | 26828.6 | 0.0267191 | 0.0000995919 |
| 0.000390625 | 26828.6 | 0.00667558 | 0.0000248823 |

Output[2]

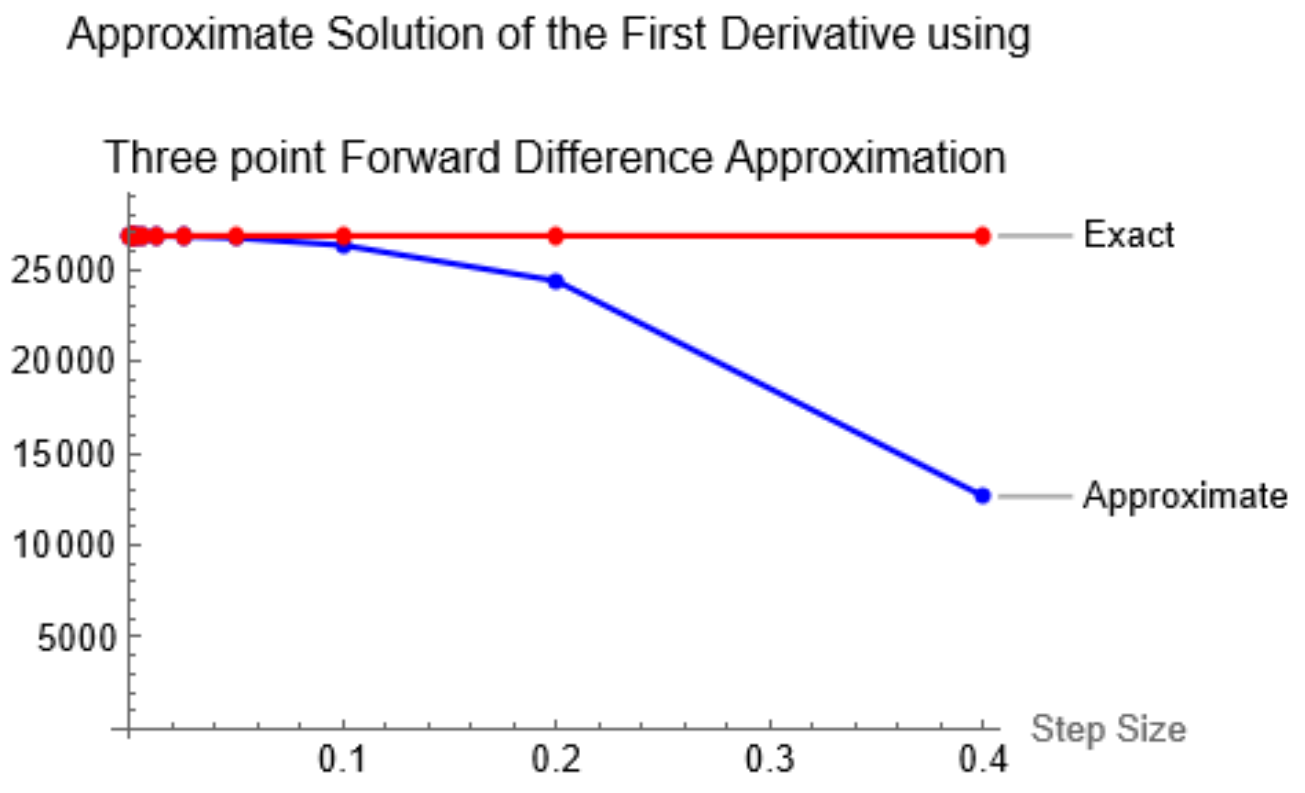

Output[3]

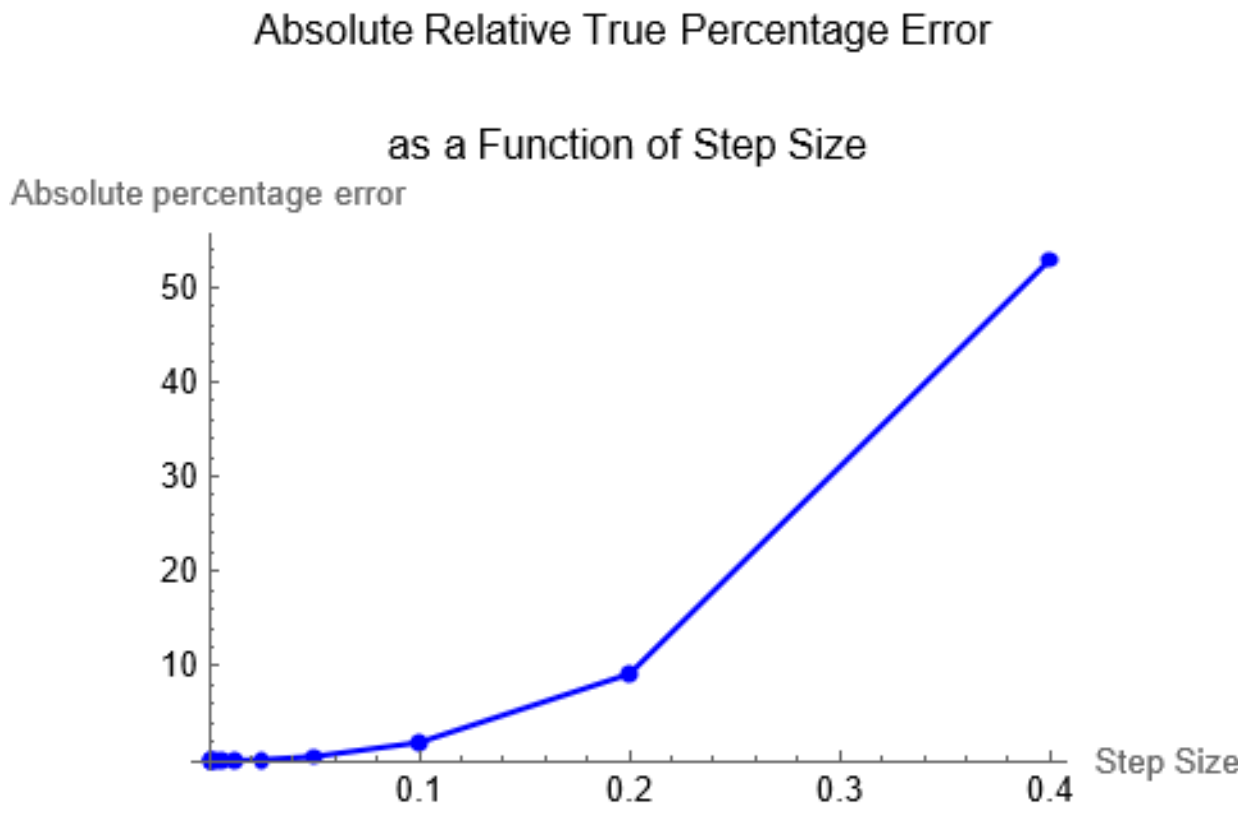

Output[3]

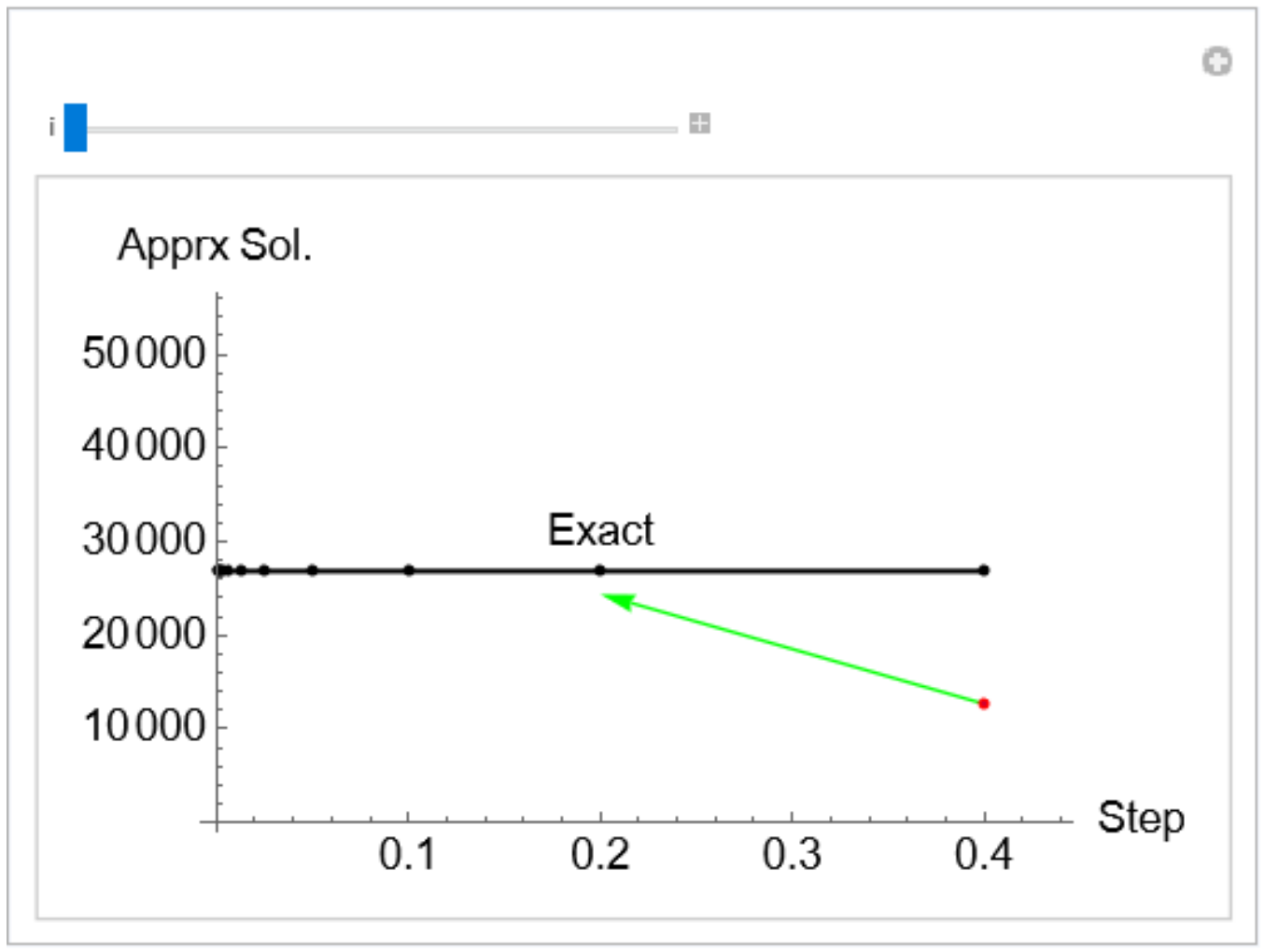

## 2.6 Mathematica Built-in Functions

Remember, in Chapter 1, we gave many functions for generating vectors and matrices.

### Vectors

Once a vector is defined, standard vector operations can be performed on it. The following list summarizes them.

| | |
|---|---|
| `{e₁,e₂,...}` | is a list of elements. |
| `Range[n]` | create the list $\{1,2,3,\ldots,n\}$. |
| `Range[n₁,n₂]` | create the list $\{n_1, n_1+1,\ldots,n_2\}$. |
| `Range[n₁,n₂,dn]` | create the list $\{n_1, n_1+dn,\ldots,n_2\}$. |
| `Table[f,{i,n}]` | build a length-n vector by evaluating f with $i = 1,2,\ldots,n$. |
| `Length[list]` | give the number of elements in list. |
| `List[[i]] or Part[list,i]` | give the $i^{th}$ element in the vector list. |
| `c v` | multiply a vector v by a scalar. |
| `a.b` | dot product of two vectors a. b. |
| `Cross[a,b]` | cross product of two vectors (also input as $a \times b$). |
| `Norm[v]` | Euclidean norm of a vector v. |
| `Normalize[v]` | gives the normalized form of a vector v. |
| `Orthogonalize [{v₁,v₂,…}]` | gives an orthonormal basis found by orthogonalizing the vectors $v_i$. |

### Matrix

The following list summarizes the functions for generating matrices.

| | |
|---|---|
| `{{a,b},{c,d}}` | matrix $\begin{pmatrix} a & b \\ c & d \end{pmatrix}$. |
| `Table[f,{i,m},{j,n}]` | build an $m \times n$ matrix by evaluating f with i ranging from 1 to m and j ranging. |
| `List[[i,j]] or Part[list,i,j]` | give the i, j th element in the matrix list. |
| `DiagonalMatrix[list]` | generate a square matrix with the elements in list on the main. |
| `Dimensions[list]` | give the dimensions of a matrix represented by list. |
| `Column[list]` | display the elements of list in a column. |
| `c m` | multiply a matrix m by a scalar. |
| `a.b` | dot product of two matrices a. b. |
| `Inverse[m]` | matrix inverse of a matrix m. |
| `MatrixPower[m,n]` | gives the $n^{th}$ power of a matrix m. |
| `Transpose[m]` | Transpose m. |
| `Det[m]` | Determinant of a matrix m. |
| `Tr[m]` | Trace m. |
| `Transpose[m]` | Transpose m. |

### Layout & Tables

| | |
|---|---|
| `MatrixForm[list]` | prints with the elements of list arranged in a regular array. |
| `TableForm[list]` | prints with the elements of list arranged in an array of rectangular cells. |

***Mathematica Examples 2.1***

```
Input     List[a,b,c,d]
Output    {a,b,c,d}

Input     Table[
           i^2,
           {i,10}
           ]
Output    {1,4,9,16,25,36,49,64,81,100}

Input     m = {{a, b},{c, d}}
Output    {{a, b},{c, d}}

Input     m[[1]]
Output    {a, b}

Input     m[[1, 2]]
Output    b
```

### Notes:

1- When using `MatrixForm`, it is very important to note that the `MatrixForm` is used for display purposes only [9]. If a matrix is defined with the `MatrixForm` in it, that matrix definition cannot be used in any subsequent calculation. For example, consider the following definition of matrix `m`:

`Input: m = MatrixForm[{{1, 2, 3,4}, {5, 6, 7,8}, {9,10,11,12}}]`

Output: $\begin{pmatrix} 1 & 2 & 3 & 4 \\ 5 & 6 & 7 & 8 \\ 9 & 10 & 11 & 12 \end{pmatrix}$

We cannot perform any operations on matrix m in this form. For example, using the `Transpose` function on `m` simply returns the initial matrix `m` wrapped in `Transpose`.

`Input: Transpose[m]`

Output: Transpose[$\begin{pmatrix} 1 & 2 & 3 & 4 \\ 5 & 6 & 7 & 8 \\ 9 & 10 & 11 & 12 \end{pmatrix}$]

The solution obviously is not to use the `MatrixForm` in the definition of matrices. After the definition, we can use the `MatrixForm` to get a nice-looking display.

`Input: m= {{1,2,3,4}, {5,6,7,8}, {9, 10, 11, 12}}; MatrixForm [m]`

Output: $\begin{pmatrix} 1 & 2 & 3 & 4 \\ 5 & 6 & 7 & 8 \\ 9 & 10 & 11 & 12 \end{pmatrix}$

`Input: Transpose[m]`

Output: $\begin{pmatrix} 1 & 5 & 9 \\ 2 & 6 & 10 \\ 3 & 7 & 11 \\ 4 & 8 & 12 \end{pmatrix}$

2- To explicitly define a $3 \times 1$ column vector is to enter it as a two-dimensional list. Here, each entry defines a row of a matrix with one column.

`Input: a = {{1}, {2}, {3}}; MatrixForm[a]`

Output: $\begin{pmatrix} 1 \\ 2 \\ 3 \end{pmatrix}$

`Input: b = {4, 5, 6}`

`Output: {4, 5, 6}`

In this case `Transpose[a] . b` makes sense. However, the result is a $1 \times 1$ matrix and not a scalar.

`Input: Transpose [a] . b`

`Output: {32}`

Obviously, now `a . b` will produce an error because the dimensions do not match.

`Input: a . b`

`Output: Dot : : dotsh : Tensors { {1} , {2} , {3} }`

`and {4, 5, 6} have incompatible shapes.`

`{{1}, {2}, {3}} . {4,5,6}`

To get a $3 \times 3$ matrix, we need to define `b` as a two-dimensional matrix with only one row, as follows:

`Input: b = {{4, 5, 6}}`

`Output: {{4,5,6}}`

Now a $(3 \times 1)$ **.** $(1 \times 3)$ can be evaluated directly, as one would expect.

`Input: MatrixForm [a . b]`

Output: $\begin{pmatrix} 4 & 5 & 6 \\ 8 & 10 & 12 \\ 12 & 15 & 18 \end{pmatrix}$

***Mathematica Examples 2.2***

```
Input       u={{x},{y}};
            MatrixForm[u]
            MatrixForm[Transpose[u]]
Output      ({
              {x},
              {y}
             })
Output      ({
              {x, y}
             })

Input       v={{x,y}};
            MatrixForm[v]
```

```
        MatrixForm[Transpose[v]]
Output  ({
          {x, y}
         })
Output  ({
          {x},
          {y}
         })

Input   (* Scalar product of vectors in three dimensions (row product column): *)
        {{a, b, c}} . {{x}, {y}, {z}}
Output  {{a x + b y + c z}}

Input   (* Scalar product of vectors in three dimensions (column product row): *)
        {{x},{y},{z}}.{{a,b,c}}//MatrixForm
Output  ({
          {a x, b x, c x},
          {a y, b y, c y},
          {a z, b z, c z}
         })

Input   (* Scalar product of vectors in two dimensions (row product column): *)
        u1 = {{1, 1}};
        v1 = {{-1}, {1}};
        u1 . v1
Output  {{0}}

Input   (* The product of a matrix and a vector: *)
        {{a, b}, {c, d}} . {{x}, {y}}
Output  {{a x + b y}, {c x + d y}}

Input   (* The product v^T M of a vector and a matrix: *)
        {{x, y}} . {{a, b}, {c, d}}
Output  {{a x + c y, b x + d y}}

Input   (* The product v^T Mw of a matrix and two vectors: *)
        {{x, y}} . {{a, b}, {c, d}} . {{r}, {s}}
Output  {{r (a x + c y) + s (b x + d y)}}

Input   (* The product of two matrices: *)
        {{a, b}, {c, d}} . {{1, 2}, {3, 4}} // MatrixForm
Output  ({
          {a+3 b, 2 a+4 b},
          {c+3 d, 2 c+4 d}
         })

Input   (* Multiply in the other order: *)
        {{1, 2}, {3, 4}} . {{a, b}, {c, d}} // MatrixForm
Output  ({
          {a+2 c, b+2 d},
          {3 a+4 c, 3 b+4 d}
         })

Input   {{1, 2, 3}, {4, 5, 6}} . {{a, b}, {c, d}, {e, f}} // MatrixForm
Output  ({
          {a+2 c+3 e, b+2 d+3 f},
          {4 a+5 c+6 e, 4 b+5 d+6 f}
         })

Input   Dot[{{3.2, 4.2, 5.2}}, {{0.75}, {1.1}, {0.0625}}]
Output  {{7.345}}
```

```
Input      (* Define a 2D-vector and a 3D-vector: *)
           v2 = {{0.618678, 0.213605}};
           v3 = {{0.804978}, {0.587651}, {0.2951}};
           {v2 // MatrixForm, v3 // MatrixForm}
Output     {({
              {0.618678, 0.213605}
             }),({
              {0.804978},
              {0.587651},
              {0.2951}
             })}

Input      (* Define a rectangular matrix of dimensions: *)
           r = {{0.187902, 0.498054, 0.767621}, {0.226789, 0.852257, 0.819982}};
           r // MatrixForm
Output     ({
             {0.187902, 0.498054, 0.767621},
             {0.226789, 0.852257, 0.819982}
            })

Input      v2 . r
Output     {{0.164694, 0.490181, 0.650062}}

Input      r . v2
Output     Dot::dotsh: Tensors {{0.187902,0.498054,0.767621},{0.226789,0.852257,0.819982}}
           and {{0.618678,0.213605}} have incompatible shapes.

           {{0.187902, 0.498054, 0.767621}, {0.226789, 0.852257, 0.819982}} . {{0.618678,
           0.213605}}

Input      r . v3
Output     {{0.670464}, {0.925367}}

Input      v2 . r . v3
Output     {{0.612464}}

Input      m = {{1, 2}, {3, 4}};
           v = {{5}, {6}};
Input      m . v
Output     {{17}, {39}}

Input      Transpose[v] . m
Output     {{23, 34}}

Input      Transpose[v] . m . v
Output     {{319}}

Input      (* Product of exact matrices: *)
           m = {{1, 2}, {3, 4}, {5, 6}};
           n = {{6, 5, 4}, {3, 2, 1}};
           m . n // MatrixForm
Output     ({
             {12, 9, 6},
             {30, 23, 16},
             {48, 37, 26}
            })

Input      n . m // MatrixForm
```

Output
```
({
  {41, 56},
  {14, 20}
 })
```

Input
```
(* Visualize the input and output matrices: *)
{ MatrixPlot[n], MatrixPlot[m . n] }
```

Output

{ 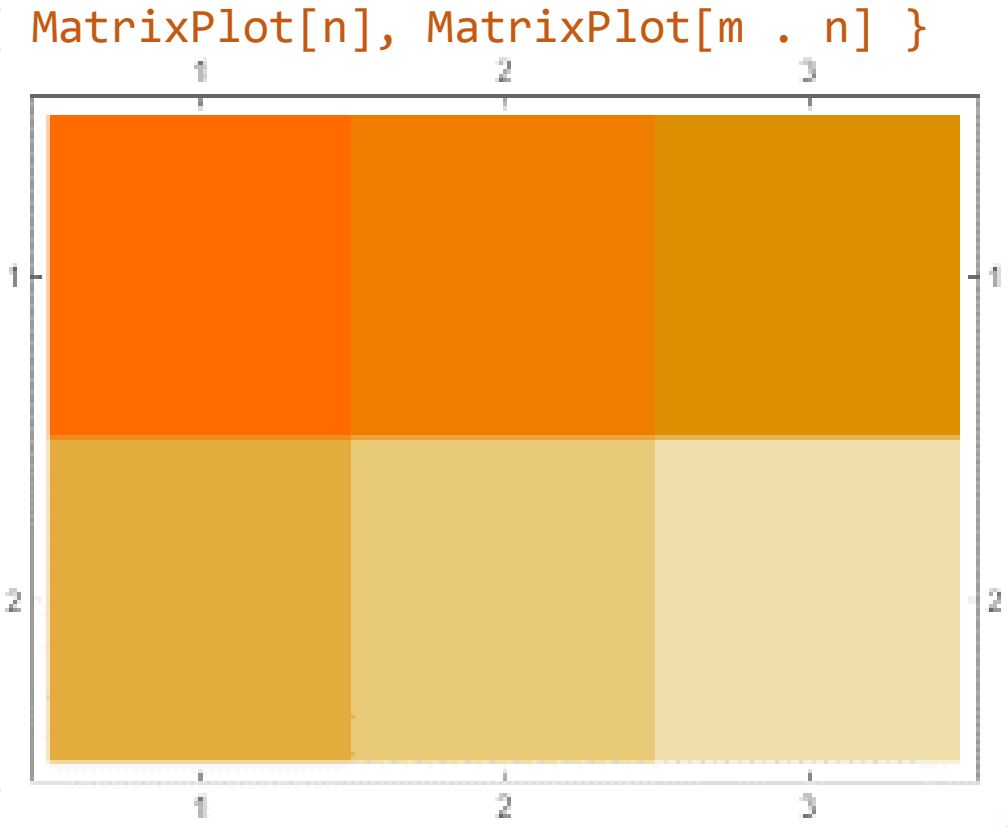 , 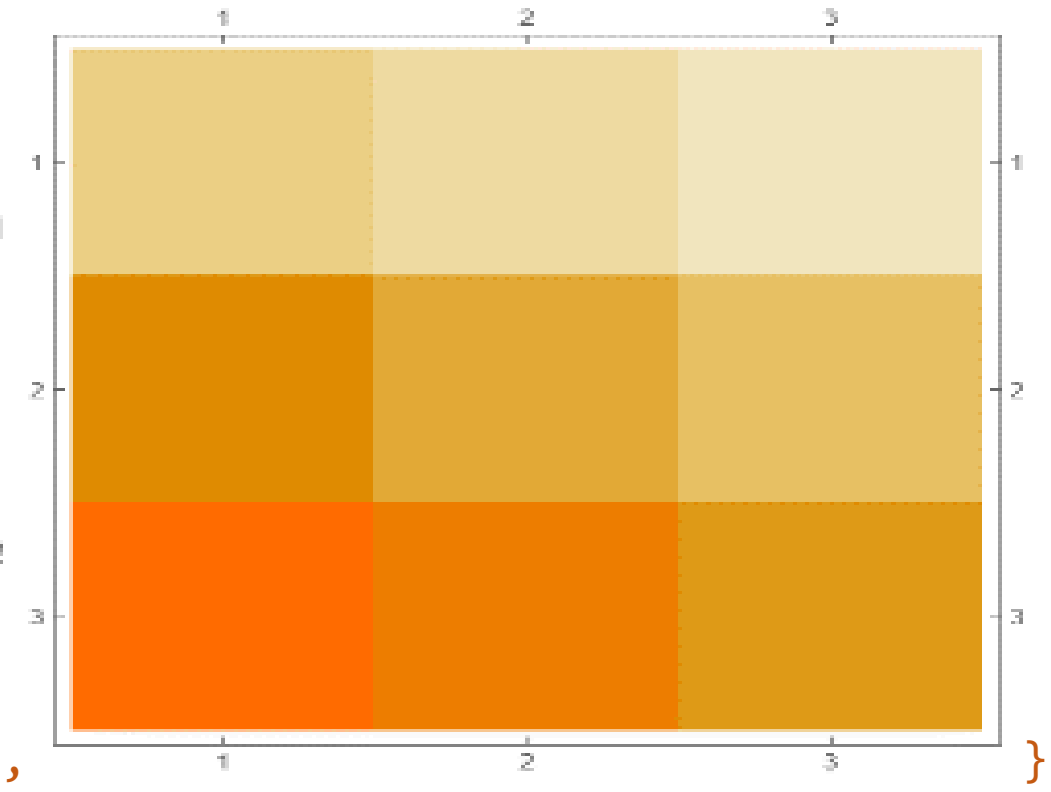 }

Input
```
MatrixPlot[
 Table[
  x^2-y^2,
  {x,-10,10},
  {y,-10,10}],
 ColorFunction->"TemperatureMap"
 ]
```

Output

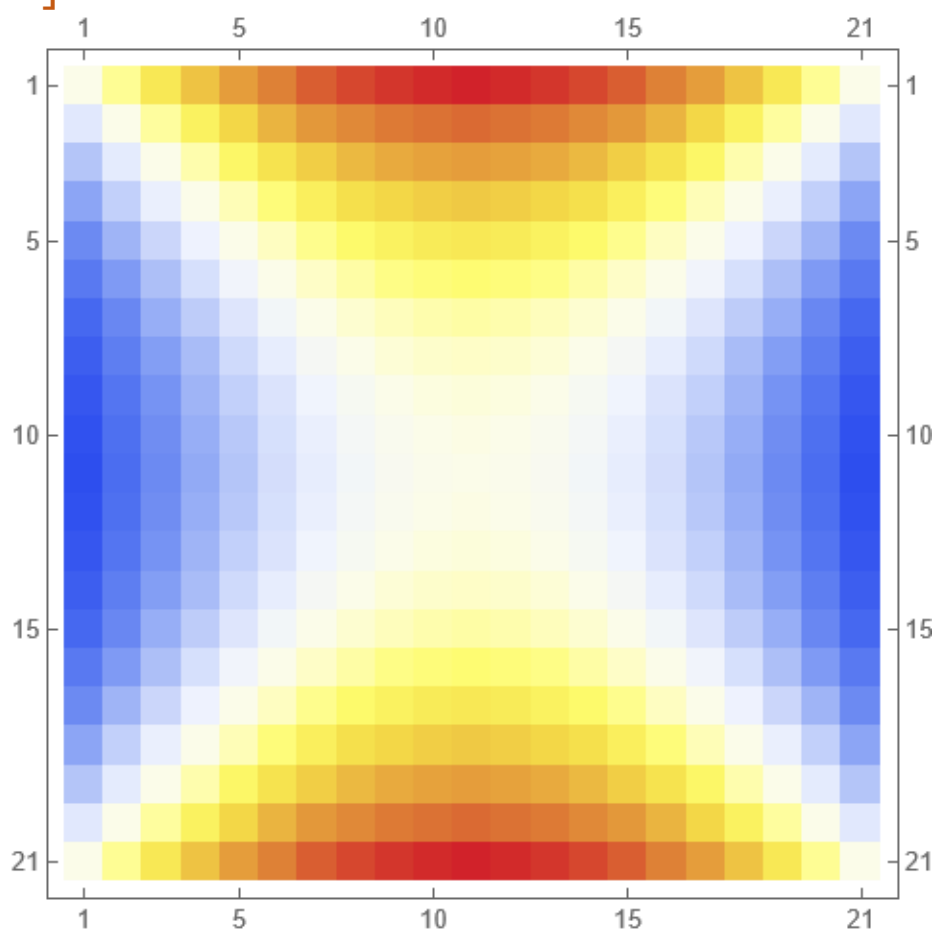

Input

(* A real symmetric matrix s gives a quadratic form $q{:}\mathbb{R}^n\text{—>}R$ by the formula $q(\mathrm{v}) = \mathrm{v}^T\,\mathrm{s}\,\mathrm{v}$: *)
```
s = {{12, 7, 13}, {7, 10, -8}, {13, -8, 4}};
s // MatrixForm
```

Output
```
({
  {12, 7, 13},
  {7, 10, -8},
  {13, -8, 4}
 })
```

Input
```
q[v_] := Transpose[v] . s . v
```
(* Quadratic forms have the property that $q(\alpha\,\mathrm{v}) = \alpha^2 q(\mathrm{v})$: *)

q[$\alpha$ {{x}, {y}, {z}}] == $\alpha^2$ q[{{x}, {y}, {z}}] // Simplify

Output
```
True
```

```
Input     (* Equivalently, they define a homogeneous quadratic polynomial in the variables
          of ℝ^n: *)
          q[{{x}, {y}, {z}}] // Expand
Output    {{12 x^2 + 14 x y + 10 y^2 + 26 x z - 16 y z + 4 z^2}}
```

3- Since Mathematica treats matrices as essentially lists [10], it does not distinguish between a column or a row vector. The actual form is determined from syntax in which it is used. For example, define two vectors `a` and `b` as follows:

```
Input : a={1, 2, 3}; b={4, 5, 6};
```

The inner product is evaluated simply as `a . b`, resulting in a scalar. Explicitly evaluating `Transpose[a]. b` will produce an error.

```
Input : a . b

Output: 32

Input : Transpose[a] . b

Output: Transpose :: imtx : The first two levels of the

one-dimensional list {1, 2, 3} cannot be transposed.

Transpose[{1, 2,3}] . {4, 5, 6}
```

If we want to treat `a` as a column vector ($3 \times 1$) and `b` as a row vector ($1 \times 3$) to get a $3 \times 3$ matrix from the product, we need to use the `Outer` function of Mathematica, as follows:

```
Input : ab = Outer[Times, a, b]; MatrixForm[ab]
```

Output: $\begin{pmatrix} 4 & 5 & 6 \\ 8 & 10 & 12 \\ 12 & 15 & 18 \end{pmatrix}$

***Mathematica Examples 2.3***

```
Input     (* Dot effectively treats vectors multiplied from the right as column vectors:
          *)
          a = {{1, 2}, {3, 4}, {5, 6}};
          a . {1, 1}
Output    {3, 7, 11}

Input     a . {{1}, {1}}
Output    {{3}, {7}, {11}}

Input     (* Dot effectively treats vectors multiplied from the left as row vectors: *)
          {1, 1, 1} . a
Output    {9, 12}

Input     {{1, 1, 1}} . a
Output    {{9, 12}}
```

### Matrix Decompositions

| | |
|---|---|
| `UpperTriangularMatrixQ[m]` | gives True if m is upper triangular, and False otherwise. |
| `UpperTriangularMatrixQ[m,k]` | gives True if m is upper triangular starting up from the $k^{th}$ diagonal, and False otherwise. |
| `LowerTriangularMatrixQ[m]` | gives True if m is lower triangular, and False otherwise. |
| `LowerTriangularMatrixQ[m,k]` | gives True if m is lower triangular starting down from the $k^{th}$ diagonal, and False otherwise. |
| `DiagonalMatrixQ[m]` | gives True if m is diagonal, and False otherwise. |
| `DiagonalMatrixQ[m,k]` | gives True if m has nonzero elements only on the $k^{th}$ diagonal, and False otherwise. |
| `IdentityMatrix[n]` | gives the $n \times n$ identity matrix. |
| `PermutationMatrix[permv]` | represents the permutation matrix given by permutation vector permv as a structured array. |
| `PermutationMatrix[pmat]` | converts a permutation matrix pmat to a structured array. |

***Mathematica Examples 2.4***

```
Input         (* Test if a matrix is upper triangular: *)
              UpperTriangularMatrixQ[({
                 {a, b, c},
                 {0, e, f},
                 {0, 0, g}
                })]
Output        True

Input         UpperTriangularMatrixQ[({
                 {1, 2, 3},
                 {4, 5, 6},
                 {7, 8, 9}
                })]
Output        False

Input         (*Test if a matrix is upper triangular starting from the first superdiagonal:*)
              MatrixForm[{{0,1,2},{0,0,3},{0,0,0}}]
              UpperTriangularMatrixQ[{{0,1,2},{0,0,3},{0,0,0}},1]
Output        ({
                {0, 1, 2},
                {0, 0, 3},
                {0, 0, 0}
               })
Output        True

Input         (* Test if a matrix is lower triangular: *)
              LowerTriangularMatrixQ[({
                 {a, 0, 0},
                 {b, c, 0},
                 {d, e, f}
                })]
Output        True

Input         LowerTriangularMatrixQ[({
                 {1, 2, 3},
                 {4, 5, 6},
                 {7, 8, 9}
                })]
Output        False

Input         (*Test if a matrix is lower triangular starting from the first superdiagonal:*)
              MatrixForm[{{1,2,0},{3,4,5},{6,7,8}}]
              LowerTriangularMatrixQ[{{1,2,0},{3,4,5},{6,7,8}},1]
Output        ({
                {1, 2, 0},
                {3, 4, 5},
                {6, 7, 8}
               })
```

```
Output      True

Input       MatrixForm[{{a,0,0},{0,b,0},{0,0,c}}]
            DiagonalMatrixQ[{{a,0,0},{0,b,0},{0,0,c}}]

Output      ({
              {a, 0, 0},
              {0, b, 0},
              {0, 0, c}
             })
Output      True

Input       MatrixForm[{{1,0,0},{0,0,2},{3,0,0}}]
            DiagonalMatrixQ[{{1,0,0},{0,0,2},{3,0,0}}]
Output      ({
              {1, 0, 0},
              {0, 0, 2},
              {3, 0, 0}
             })
Output      False

Input       MatrixForm[{{0,1,0},{0,0,2},{0,0,0}}]
            DiagonalMatrixQ[{{0,1,0},{0,0,2},{0,0,0}},1]
Output      ({
              {0, 1, 0},
              {0, 0, 2},
              {0, 0, 0}
             })
Output      True

Input       IdentityMatrix[3]//MatrixForm
Output      ({
              {1, 0, 0},
              {0, 1, 0},
              {0, 0, 1}
             })

Input       MatrixPlot[IdentityMatrix[3]]
Output
```

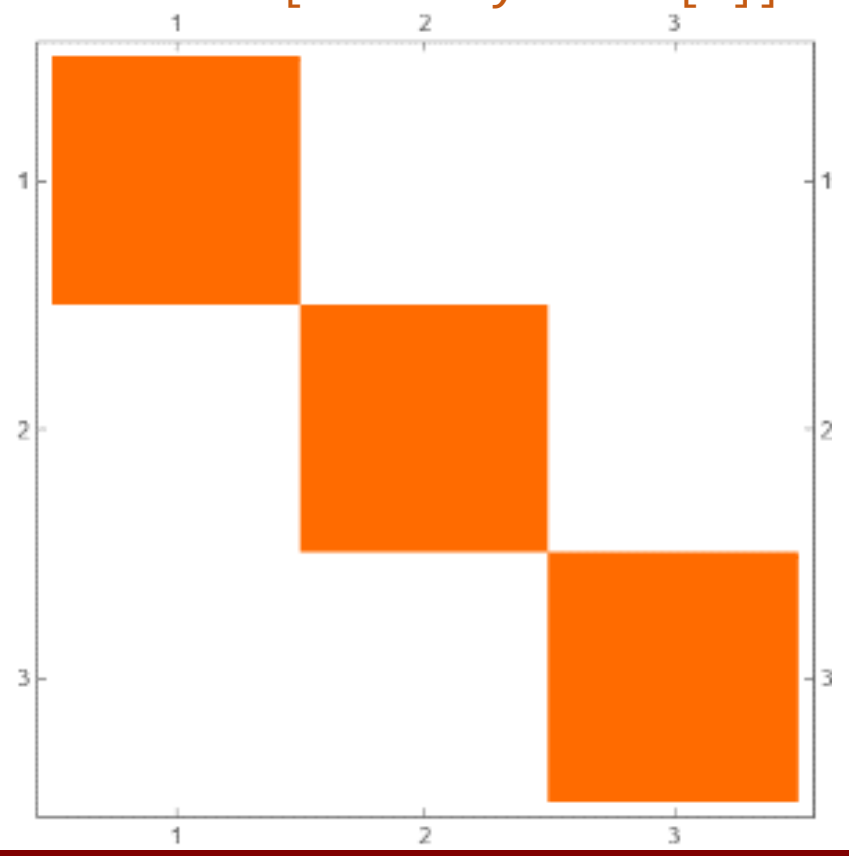

| | |
|---|---|
| RowReduce[m] | gives the row-reduced form of the matrix m. |
| LUDecomposition[m] | generates a representation of the LU decomposition of a square matrix m. |
| UpperTriangularize[m] | gives a matrix in which all but the upper triangular elements of m are replaced with zeros. |
| UpperTriangularize[m,k] | replaces with zeros only the elements below the $k^{th}$ subdiagonal of m. |
| LowerTriangularize[m] | gives a matrix in which all but the lower triangular elements of m are replaced with zeros. |
| LowerTriangularize[m,k] | replaces with zeros only the elements above the $k^{th}$ subdiagonal of m. |

**Notes:**

- RowReduce performs a version of Gaussian elimination, adding multiples of rows together so as to produce zero elements when possible. The final matrix is in reduced row echelon form.
- LUDecomposition returns a list of three elements. The first element is a combination of upper- and lower-triangular matrices, the second element is a vector specifying rows used for pivoting, and for approximate numerical matrices m the third element is an estimate of the $L^{\infty}$ condition number of m.
- UpperTriangularize[m] / LowerTriangularize[m] works even if m is not a square matrix.
- In UpperTriangularize[m,k]/ LowerTriangularize[m,k], positive k refers to subdiagonals above the main diagonal and negative k refers to subdiagonals below the main diagonal.

***Mathematica Examples 2.5***

```
Input       (* Get the upper triangular part of a matrix: *)
            MatrixForm[UpperTriangularize[{{1,2,3},{4,5,6},{7,8,9}}]]
Output      ({
              {1, 2, 3},
              {0, 5, 6},
              {0, 0, 9}
             })

Input       (* Get the strictly upper triangular part of a matrix: *)
            MatrixForm[UpperTriangularize[{{1,2,3},{4,5,6},{7,8,9}},1]]
Output      ({
              {0, 2, 3},
              {0, 0, 6},
              {0, 0, 0}
             })

Input       (* Get the lower triangular part of a matrix: *)
            MatrixForm[LowerTriangularize[{{1,2,3},{4,5,6},{7,8,9}}]]
Output      ({
              {1, 0, 0},
              {4, 5, 0},
              {7, 8, 9}
             })

Input       (* Get the strictly lower triangular part of a matrix: *)
            MatrixForm[LowerTriangularize[{{1,2,3},{4,5,6},{7,8,9}},-1]]
Output      ({
              {0, 0, 0},
              {4, 0, 0},
              {7, 8, 0}
             })

Input       (* Do row reduction on a square matrix: *)
            RowReduce[{{1,2,3},{4,5,6},{7,8,9}}];
            MatrixForm[%]
Output      ({
              {1, 0, -1},
              {0, 1, 2},
              {0, 0, 0}
             })

Input       (* Do row reduction on a rectangular matrix: *)
            RowReduce[{{1,2,3,1,0,0},{4,5,6,0,1,0},{7,8,9,0,0,1}}];
            MatrixForm[%]
Output      ({
              {1, 0, -1, 0, -(8/3), 5/3},
              {0, 1, 2, 0, 7/3, -(4/3)},
```

```
   {0, 0, 0, 1, -2, 1}
  })

Input    (* Compute the LU decomposition of a matrix: *)
         {lu, p, c} = LUDecomposition[({
              {1, 1, 1},
              {2, 4, 8},
              {3, 9, 27}
             })]
Output   {{{1,1,1},{2,2,6},{3,3,6}},{1,2,3},0}

Input    (* l is the strictly lower triangular part of lu with ones assumed along the
         diagonal: *)
         MatrixForm[l=LowerTriangularize[lu,-1]+IdentityMatrix[3]]
Output   ({
           {1, 0, 0},
           {2, 1, 0},
           {3, 3, 1}
          })

Input    (* u is the upper triangular part of lu: *)
         MatrixForm[u=UpperTriangularize[lu]]
Output   ({
           {1, 1, 1},
           {0, 2, 6},
           {0, 0, 6}
          })

Input    (* l.u reconstructs the original matrix: *)
         l.u//MatrixForm
Output   ({
           {1, 1, 1},
           {2, 4, 8},
           {3, 9, 27}
          })

Input    (* Find the LU decomposition of a symbolic matrix: *)
         lu = First[LUDecomposition[({
               {a, b, c},
               {d, e, f},
               {g, h, i}
              })]] // Simplify
Output   {{a,b,c},{d/a,-((b d)/a)+e,-((c d)/a)+f},{g/a,(b g-a h)/(b d-a e),(c e g-b f g-
         c d h+a f h+b d i-a e i)/(b d-a e)}}

Input    (* Extract the l and u matrices: *)
         {l,u}={LowerTriangularize[lu,-1]+IdentityMatrix[3],UpperTriangularize[lu]};
         {l//MatrixForm,u//MatrixForm}
Output   {({
            {1, 0, 0},
            {d/a, 1, 0},
            {g/a, (b g-a h)/(b d-a e), 1}
           }),({
            {a, b, c},
            {0, -((b d)/a)+e, -((c d)/a)+f},
            {0, 0, (c e g-b f g-c d h+a f h+b d i-a e i)/(b d-a e)}
           })}

Input    (* Verify that l.u equals the original matrix: *)
         l.u//Simplify//MatrixForm
```

```
Output      ({
              {a, b, c},
              {d, e, f},
              {g, h, i}
             })
```

| | |
|---|---|
| `CholeskyDecomposition[m]` | gives the Cholesky decomposition of a matrix m. |
| `QRDecomposition[m]` | yields the QR decomposition for a numerical matrix m. The result is a list {q,r}, where q is a unitary matrix and r is an upper-triangular matrix. |
| `SchurDecomposition[m]` | yields the Schur decomposition for a numerical matrix m, given as a list {q,t} where q is an orthonormal matrix and t is a block upper-triangular matrix. |
| `SchurDecomposition[{m,a}]` | gives the generalized Schur decomposition of m with respect to a. |
| `JordanDecomposition[m]` | yields the Jordan decomposition of a square matrix m. The result is a list {s,j} where s is a similarity matrix and j is the Jordan canonical form of m. |
| `HermiteDecomposition[m]` | gives the Hermite normal form decomposition of an integer matrix m. |
| `SmithDecomposition[m]` | gives the Smith normal form decomposition of an integer matrix m. |

***Mathematica Examples 2.6***

```
Input       (* Compute the Cholesky decomposition of a 2×2 real matrix: *)
            CholeskyDecomposition[({
               {2, 1},
               {1, 2}
              })]

            (* Verify the decomposition: *)
            ConjugateTranspose[%] . % // MatrixForm

            (* The original matrix is positive definite: *)
            PositiveDefiniteMatrixQ[({
               {2, 1},
               {1, 2}
              })]
Output      {{√2,1/√2},{0,√3/2}}
Output      ({
              {2, 1},
              {1, 2}
             })
Output      True

Input       (* The decomposition of a 2×2 matrix into a unitary (orthogonal) matrix and
            upper triangular matrix: *)
            m=({{1,2},{3,4}});
            {q,r}=QRDecomposition[m];
            {q//MatrixForm,r//MatrixForm}

            (* Verify that m=q†.r: *)
            m==ConjugateTranspose[q].r
Output      {({
               {1/√10, 3/√10},
               {3/√10, -(1/√10)}
              }),({
               {√10, 7√(2/5)},
               {0, √(2/5)}
              })}
Output      True

Input       (* Compute the QR decomposition for a 3×2 matrix with exact values: *)
            {q, r} = QRDecomposition[({
                 {1, 2},
```

```
                  {3, 4},
                  {5, 6}
                 })];
             {q // MatrixForm, r // MatrixForm}

             (* qT.r is the original matrix: *)
             Transpose[q].r//MatrixForm
Output       {({
                {1/ √35, 3/√35, √(5/7) },
                {13/√210, 2√(2/105), -√(5/42)}
               }),({
                {√35, 44/√35},
                {0, 2 √(6/35)}
               })}
             ({
               {1, 2},
               {3, 4},
               {5, 6}
              })

Input        (* Find the Schur decomposition of a real matrix: *)
             m={{2.7,4.8,8.1},{-.6,0,0},{.1,0,.3}};
             {q,t}=SchurDecomposition[m]

             (* Confirm the decomposition up to numerical rounding: *)
             m-q.t.ConjugateTranspose[q]//Chop

             (* Format q and t: *)
             {q//MatrixForm,t//MatrixForm}
Output       {{{0.890043,-0.454494,0.0354887},{-0.445021,-
             0.849324,0.283909},{0.0988936,0.268485,0.958194}},{{1.2,-
             2.90297,8.24441},{0.,1.2,-4.0942},{0.,0.,0.6}}}
Output       {{0,0,0},{0,0,0},{0,0,0}}
Output       {({
                {0.890043, -0.454494, 0.0354887},
                {-0.445021, -0.849324, 0.283909},
                {0.0988936, 0.268485, 0.958194}
               }),({
                {1.2, -2.90297, 8.24441},
                {0., 1.2, -4.0942},
                {0., 0., 0.6}
               })}

Input        (* Find the Jordan decomposition of a 3×3 matrix: *)
             JordanDecomposition[{{27,48,81},{-6,0,0},{1,0,3}}]

             (* Format the results: *)
             Map[MatrixForm,JordanDecomposition[{{27,48,81},{-6,0,0},{1,0,3}}]]
Output       {{{3,18,2},{-3,-9,-(1/4)},{1,2,0}},{{6,0,0},{0,12,1},{0,0,12}}}
Output       {({
                {3, 18, 2},
                {-3, -9, -(1/4)},
                {1, 2, 0}
               }),({
                {6, 0, 0},
                {0, 12, 1},
                {0, 0, 12}
               })}

Input        (* Decompose m into a unimodular matrix u and an upper-triangular matrix r: *)
             m=({{1,2,3},{5,4,3},{8,7,9}});
```

```
              {u,r}=HermiteDecomposition[m];
              {u//MatrixForm,r//MatrixForm}

              (* The determinant of u has absolute value 1: *)
              Abs[Det[u]]

              (* Confirm the decomposition: *)
              u.m==r
Output        {({
                 {1, 0, 0},
                 {3, 1, -1},
                 {-1, -3, 2}
                }),({
                 {1, 2, 3},
                 {0, 3, 3},
                 {0, 0, 6}
                })}
Output        1
Output        True

Input         (* Decompose m into unimodular matrices u and v and a diagonal matrix r: *)
              m=({{1,2,3},{5,4,3},{8,7,9}});
              {u,r,v}=SmithDecomposition[m];
              {u//MatrixForm,r//MatrixForm,v//MatrixForm}

              (* The determinants of u and v have absolute value 1: *)
              {Abs[Det[u]],Abs[Det[v]]}

              (* Each entry on the diagonal of r divides the successor: *)
              {Divisible[r[[2,2]],r[[1,1]]],Divisible[r[[3,3]],r[[2,2]]]}

              (* Confirm the decomposition: *)
              u.m.v==r
Output        {({
                 {1, 0, 0},
                 {3, 1, -1},
                 {-1, -3, 2}
                }),({
                 {1, 0, 0},
                 {0, 3, 0},
                 {0, 0, 6}
                }),({
                 {1, -2, -1},
                 {0, 1, -1},
                 {0, 0, 1}
                })}
Output        {1,1}
Output        {True,True}
Output        True
```

### Gradient function

| | |
|---|---|
| `Grad[f,{x1,…,xn}]` | gives the gradient $(\frac{\partial f}{\partial x_1}, \dots, \frac{\partial f}{\partial x_n})$. |

***Mathematica Examples 2.7***

```
Input         (* The gradient in three-dimensional Cartesian coordinates: *)
              Grad[f[x, y, z], {x, y, z}]
Output        {f^(1,0,0)[x,y,z], f^(0,1,0)[x,y,z], f^(0,0,1)[x,y,z]}
```

```
Input     (* The gradient in two dimensions: *)
          Grad[Sin[x^2 + y^2], {x, y}]
Output    {2 x Cos[x^2 + y^2], 2 y Cos[x^2 + y^2]}

Input     (*The gradient of a vector field in Cartesian coordinates, the Jacobian matrix*)
          Grad[{f[x, y, z], g[x, y, z], h[x, y, z]}, {x, y, z}] // MatrixForm
Output    ({
            {f^(1,0,0)[x,y,z], f^(0,1,0)[x,y,z], f^(0,0,1)[x,y,z]},
            {g^(1,0,0)[x,y,z], g^(0,1,0)[x,y,z], g^(0,0,1)[x,y,z]},
            {h^(1,0,0)[x,y,z], h^(0,1,0)[x,y,z], h^(0,0,1)[x,y,z]}
           })

Input     (* Compute the Hessian of a scalar function: *)
          grad1 = Grad[x*y*z, {x, y, z}]
Output    {y z, x z, x y}

Input     Grad[grad1, {x, y, z}] // MatrixForm
Output    ({
            {0, z, y},
            {z, 0, x},
            {y, x, 0}
           })

Input     (* Find the critical points of a function of two variables: *)
          f[x_, y_] := x^4 + y^4 - 20 x^2 - 10 x y – 25;
          grad2 = Grad[f[x, y], {x, y}]
Output    {-40 x + 4 x^3 - 10 y, -10 x + 4 y^3}

Input     sol = NSolve[grad2 == {0, 0}, {x, y}, Reals]
Output    {{x -> -3.39162, y -> -2.03915}, {x -> 3.39162, y -> 2.03915}, {x -> 0., y ->
          0.}}

Input     (* Compute the signs of ∂²f(x,y)/∂x² and the Hessian determinant: *)
          hessian2 = Grad[grad2, {x, y}]
Output    {{-40 + 12 x^2, -10}, {-10, 12 y^2}}

Input     Sign[hessian2[[1, 1]] /. sol]
Output    {1, 1, -1}

Input     Sign[Det[hessian2] /. sol]
Output    {1, 1, -1}


Input     (* Find the critical points of a function of three variables: *)
          f[x_, y_, z_] := 25 - x^2 + 2 x^4 - x y - y^2 + 4 x^2 y^2 + 2 y^4 - 2 x z - 2 y
          z -  z^2 + 4 x^2 z^2 + 4 y^2 z^2 + 2 z^4;
          grad3 = Grad[f[x, y, z], {x, y, z}]
Output    {-2 x + 8 x^3 - y + 8 x y^2 - 2 z + 8 x z^2, -x - 2 y + 8 x^2 y +
            8 y^3 - 2 z + 8 y z^2, -2 x - 2 y - 2 z + 8 x^2 z + 8 y^2 z + 8 z^3}

Input     sol = NSolve[grad3 == {0, 0, 0}, {x, y, z}, Reals]
Output    {{x -> -0.443969, y -> -0.443969, z -> -0.52661}, {x -> 0.443969, y -> 0.443969,
          z -> 0.52661}, {x -> -0.25, y -> 0.25, z -> 0.}, {x -> 0.25, y -> -0.25, z ->
          0.}, {x -> 0., y -> 0., z -> 0.}}

Input     (* Compute the Hessian matrix of f: *)
          hessian3 = Grad[grad3, {x, y, z}]
Output    {{-2 + 24 x^2 + 8 y^2 + 8 z^2, -1 + 16 x y, -2 + 16 x z}, {-1 +
             16 x y, -2 + 8 x^2 + 24 y^2 + 8 z^2, -2 + 16 y z}, {-2 +
```

```
              16 x z, -2 + 16 y z, -2 + 8 x^2 + 8 y^2 + 24 z^2}}

Input      Sign[Eigenvalues[hessian3] /. sol]
Output     {{1, 1, 1}, {1, 1, 1}, {-1, 1, 1}, {-1, 1, 1}, {-1, -1, 1}}

Input      Grad[x*y*Sin[z], {x, y, z}]
Output     {y Sin[z], x Sin[z], x y Cos[z]}

Input      D[x*y*Sin[z], {{x, y, z}}] // MatrixForm
Output     ({
             {y Sin[z]},
             {x Sin[z]},
             {x y Cos[z]}
            })

Input      Grad[{x*y, y*z, z*x}, {x, y, z}] // MatrixForm
Output     ({
             {y, x, 0},
             {0, z, y},
             {z, 0, x}
            })

Input      D[{x*y, y*z, z*x}, {{x, y, z}}] // MatrixForm
Output     ({
             {y, x, 0},
             {0, z, y},
             {z, 0, x}
            })
```

# CHAPTER 3

# LINEAR ALGEBRA – THE INCREDIBLE BEAUTY OF MATH

## 3.1 Vector Spaces

**Definition (Vector Space):** Let $K$ be a given field and let $V$ be a non-empty set with rules of addition and scalar multiplication, which assigns to any $|\mathbf{u}\rangle, |\mathbf{v}\rangle \in V$ a sum $|\mathbf{u}\rangle + |\mathbf{v}\rangle \in V$ and to any $|\mathbf{u}\rangle \in V, k \in K$ a product $k|\mathbf{u}\rangle \in V$. Then $V$ is called a vector space over $K$ (and the elements of $V$ are called vectors) if the following axioms hold.

1- There is a vector in $V$, denoted by $|\mathbf{0}\rangle$ and called the zero vector, for which
$$|\mathbf{u}\rangle + |\mathbf{0}\rangle = |\mathbf{u}\rangle, \quad |\mathbf{u}\rangle \in V. \tag{3.1}$$
2- For each vector $|\mathbf{u}\rangle \in V$ there is a vector in $V$, denoted by $|-\mathbf{u}\rangle$, for which
$$|\mathbf{u}\rangle + |-\mathbf{u}\rangle = |\mathbf{0}\rangle. \tag{3.2}$$
3- For any vectors $|\mathbf{u}\rangle, |\mathbf{v}\rangle \in V$,
$$|\mathbf{u}\rangle + |\mathbf{v}\rangle = |\mathbf{v}\rangle + |\mathbf{u}\rangle. \tag{3.3}$$
4- For any vectors $|\mathbf{u}\rangle, |\mathbf{v}\rangle, |\mathbf{w}\rangle \in V$,
$$(|\mathbf{u}\rangle + |\mathbf{v}\rangle) + |\mathbf{w}\rangle = |\mathbf{u}\rangle + (|\mathbf{v}\rangle + |\mathbf{w}\rangle). \tag{3.4}$$
5- For any scalar $k \in K$ and vectors $|\mathbf{u}\rangle, |\mathbf{v}\rangle \in V$,
$$k(|\mathbf{u}\rangle + |\mathbf{v}\rangle) = k|\mathbf{u}\rangle + k|\mathbf{v}\rangle. \tag{3.5}$$
6- For any scalars $a, b \in K$ and any vector $|\mathbf{u}\rangle \in V$,
$$(a+b)|\mathbf{u}\rangle = a|\mathbf{u}\rangle + b|\mathbf{u}\rangle. \tag{3.6}$$
7- For any scalars $a, b \in K$ and any vector $|\mathbf{u}\rangle \in V$,
$$(ab)|\mathbf{u}\rangle = a(b|\mathbf{u}\rangle). \tag{3.7}$$
8- For the unit scalar $1 \in K$, $1|\mathbf{u}\rangle = |\mathbf{u}\rangle$, for any vector $|\mathbf{u}\rangle \in V$.

Let us consider the following two important linear vector spaces:

1- The notation $\mathbb{R}^n$ is frequently used to denote the set of all $n$-tuples of elements in $\mathbb{R}$. Here $\mathbb{R}^n$ is viewed as a vector space over $\mathbb{R}$ where vector addition and scalar multiplication are defined by
$$(a_1, a_2, \dots, a_n) + (b_1, b_2, \dots, b_n) = (a_1 + b_1, a_2 + b_2, \dots, a_n + b_n), \tag{3.8}$$
and
$$k(a_1, a_2, \dots, a_n) = (ka_1, ka_2, \dots, ka_n). \tag{3.9}$$
The zero vector in $\mathbb{R}^n$ is the $n$-tuple of zeros,
$$(0,0, \dots ,0), \tag{3.10}$$
and the negative of a vector is defined by
$$-(a_1, a_2, \dots, a_n) = (-a_1, -a_2, \dots, -a_n). \tag{3.11}$$

2- The notation $M_{m\times n}$ or simply $M$, will be used to denote the set of all $m \times n$ matrices over $\mathbb{R}$. Here $M_{m\times n}$ is viewed as a vector space over $\mathbb{R}$ where matrix addition and scalar multiplication are defined by
$$\begin{pmatrix} a_{11} & a_{11} & \cdots & a_{1n} \\ a_{21} & a_{22} & \cdots & a_{2n} \\ \cdots & \cdots & \cdots & \cdots \\ a_{n1} & a_{n2} & \cdots & a_{nn} \end{pmatrix} + \begin{pmatrix} b_{11} & b_{11} & \cdots & b_{1n} \\ b_{21} & b_{22} & \cdots & b_{2n} \\ \cdots & \cdots & \cdots & \cdots \\ b_{n1} & b_{n2} & \cdots & b_{nn} \end{pmatrix} = \begin{pmatrix} a_{11} + b_{11} & \cdots \cdots & a_{1n} + b_{1n} \\ a_{21} + a_{21} & \cdots \cdots & a_{2n} + b_{2n} \\ \cdots & \cdots \cdots & \cdots \\ a_{n1} + b_{n1} & \cdots \cdots & a_{nn} + b_{nn} \end{pmatrix}, \tag{3.12}$$
and

$$k\begin{pmatrix} a_{11} & a_{11} & \cdots & a_{1n} \\ a_{21} & a_{22} & \cdots & a_{2n} \\ \cdots & \cdots & \cdots & \cdots \\ a_{n1} & a_{n2} & \cdots & a_{nn} \end{pmatrix} = \begin{pmatrix} ka_{11} & \cdots & \cdots & ka_{1n} \\ ka_{21} & \cdots & \cdots & ka_{2n} \\ \cdots & \cdots & \cdots & \cdots \\ ka_{n1} & \cdots & \cdots & ka_{nn} \end{pmatrix}. \tag{3.13}$$

The zero vector in $M_{m,n}$ is the matrix of zeros,

$$\begin{pmatrix} 0 & \cdots & \cdots & 0 \\ 0 & \cdots & \cdots & 0 \\ \cdots & \cdots & \cdots & \cdots \\ 0 & \cdots & \cdots & 0 \end{pmatrix}, \tag{3.14}$$

and the negative of a vector is defined by

$$-\begin{pmatrix} a_{11} & a_{11} & \cdots & a_{1n} \\ a_{21} & a_{22} & \cdots & a_{2n} \\ \cdots & \cdots & \cdots & \cdots \\ a_{n1} & a_{n2} & \cdots & a_{nn} \end{pmatrix} = \begin{pmatrix} -a_{11} & \cdots & \cdots & -a_{1n} \\ -a_{21} & \cdots & \cdots & -a_{2n} \\ \cdots & \cdots & \cdots & \cdots \\ -a_{n1} & \cdots & \cdots & -a_{nn} \end{pmatrix}. \tag{3.15}$$

**Definition (Linear Combination):** Let $V$ be a vector space over $\mathbb{R}$ and let $|\mathbf{u}_1\rangle, |\mathbf{u}_2\rangle, \ldots, |\mathbf{u}_n\rangle \in V$. Any vector in $V$ of the form

$$a_1|\mathbf{u}_1\rangle + a_2|\mathbf{u}_2\rangle + \cdots + a_n|\mathbf{u}_n\rangle, \tag{3.16}$$

where the $a_i \in \mathbb{R}$ is called a linear combination of $|\mathbf{u}_1\rangle, |\mathbf{u}_2\rangle, \ldots, |\mathbf{u}_n\rangle$.

For example: consider the following set of vectors in $\mathbb{R}^3$

$$S = \{|\mathbf{v}_1\rangle = \begin{pmatrix} 1 \\ 3 \\ 1 \end{pmatrix}, |\mathbf{v}_2\rangle = \begin{pmatrix} 0 \\ 1 \\ 2 \end{pmatrix}, |\mathbf{v}_3\rangle = \begin{pmatrix} 1 \\ 0 \\ -5 \end{pmatrix}\}.$$

The vector $|\mathbf{v}_1\rangle$ is a linear combination of $|\mathbf{v}_2\rangle$ and $|\mathbf{v}_3\rangle$ because $|\mathbf{v}_1\rangle = 3|\mathbf{v}_2\rangle + |\mathbf{v}_3\rangle$.

**Definition (Linear Span):** The vectors $|\mathbf{u}_1\rangle, |\mathbf{u}_2\rangle, \ldots, |\mathbf{u}_n\rangle$ span $V$ if, for every $|\mathbf{u}\rangle \in V$, there exist scalars $a_1, a_2, \ldots a_n$ such that

$$|\mathbf{u}\rangle = a_1|\mathbf{u}_1\rangle + a_2|\mathbf{u}_2\rangle + \cdots + a_n|\mathbf{u}_n\rangle, \tag{3.17}$$

that is, if $|\mathbf{u}\rangle$ is a linear combination of $|\mathbf{u}_1\rangle, |\mathbf{u}_2\rangle, \ldots, |\mathbf{u}_n\rangle$.

***Example 3.1***

The set of vectors

$$S = \{|\mathbf{v}_1\rangle = \begin{pmatrix} 1 \\ 0 \\ 0 \end{pmatrix}, |\mathbf{v}_2\rangle = \begin{pmatrix} 0 \\ 1 \\ 0 \end{pmatrix}, |\mathbf{v}_3\rangle = \begin{pmatrix} 0 \\ 0 \\ 1 \end{pmatrix}\},$$

spans $\mathbb{R}^3$ because any vector $|\mathbf{u}\rangle = \begin{pmatrix} u_1 \\ u_2 \\ u_3 \end{pmatrix}$ in $\mathbb{R}^3$ can be written as

$$|\mathbf{u}\rangle = u_1\begin{pmatrix} 1 \\ 0 \\ 0 \end{pmatrix} + u_2\begin{pmatrix} 0 \\ 1 \\ 0 \end{pmatrix} + u_3\begin{pmatrix} 0 \\ 0 \\ 1 \end{pmatrix} = \begin{pmatrix} u_1 \\ u_2 \\ u_3 \end{pmatrix}.$$

Also, let the field $K$ be the set $\mathbb{R}$ of real numbers and let the vector space $V$ be the Euclidean space $\mathbb{R}^n$. Consider the vectors

$$|\mathbf{e}_1\rangle = \begin{pmatrix} 1 \\ 0 \\ \vdots \\ 0 \end{pmatrix}, |\mathbf{e}_2\rangle = \begin{pmatrix} 0 \\ 1 \\ \vdots \\ 0 \end{pmatrix}, \ldots |\mathbf{e}_n\rangle = \begin{pmatrix} 0 \\ 0 \\ \vdots \\ 1 \end{pmatrix}.$$

Actually, any vector $|\mathbf{a}\rangle$ in $\mathbb{R}^n$ is a linear combination of $|\mathbf{e}_1\rangle, \ldots, |\mathbf{e}_n\rangle$, i.e.,

$$|\mathbf{a}\rangle = \begin{pmatrix} a_1 \\ a_2 \\ \vdots \\ n \end{pmatrix} = a_1|\mathbf{e}_1\rangle + a_2|\mathbf{e}_2\rangle + \cdots + a_n|\mathbf{e}_n\rangle.$$

**Definition (Linear Dependent and Linear Independent):** Let $V$ be a vector space over $\mathbb{R}$. The vectors $|\mathbf{u}_1\rangle, |\mathbf{u}_2\rangle, \ldots, |\mathbf{u}_n\rangle \in V$ are said to be linearly dependent over $\mathbb{R}$, or simply dependent if there exist scalars $a_1$, $a_2$, ..., $a_m \in \mathbb{R}$ not all of them 0, such that

$$a_1|\mathbf{u}_1\rangle + a_2|\mathbf{u}_2\rangle + \cdots + a_n|\mathbf{u}_n\rangle = |\mathbf{0}\rangle. \quad (3.18)$$

Otherwise, the vectors are said to be linearly independent over $\mathbb{R}$ or simply independent.

For example:

- The set $S = \{|\mathbf{v}_1\rangle = \begin{pmatrix}1\\2\end{pmatrix}, |\mathbf{v}_2\rangle = \begin{pmatrix}2\\4\end{pmatrix}\}$ in $\mathbb{R}^2$ is linearly dependent because $-2\begin{pmatrix}1\\2\end{pmatrix} + \begin{pmatrix}2\\4\end{pmatrix} = \begin{pmatrix}0\\0\end{pmatrix}$.
- However, the unit vectors $|\mathbf{e}_1\rangle, |\mathbf{e}_2\rangle, \ldots, |\mathbf{e}_n\rangle$ in $\mathbb{R}^n$ are linearly independent.

***Example 3.2***

Determine whether the following set of vectors in $\mathbb{R}^3$ is linearly independent or linearly dependent.

$$S = \{|\mathbf{v}_1\rangle = \begin{pmatrix}1\\2\\3\end{pmatrix}, |\mathbf{v}_2\rangle = \begin{pmatrix}0\\1\\2\end{pmatrix}, |\mathbf{v}_3\rangle = \begin{pmatrix}-2\\0\\1\end{pmatrix}\}.$$

***Solution***

To test for linear independence or linear dependence, form the vector equation

$$c_1|\mathbf{v}_1\rangle + c_2|\mathbf{v}_2\rangle + c_3|\mathbf{v}_3\rangle = |\mathbf{0}\rangle.$$

If the only solution of this equation is

$$c_1 = c_2 = c_3 = 0,$$

then the set $S$ is linearly independent. Otherwise, $S$ is linearly dependent. We have

$$c_1\begin{pmatrix}1\\2\\3\end{pmatrix} + c_2\begin{pmatrix}0\\1\\2\end{pmatrix} + c_3\begin{pmatrix}-2\\0\\1\end{pmatrix} = \begin{pmatrix}0\\0\\0\end{pmatrix},$$

$$\begin{pmatrix}c_1 - 2c_3\\2c_1 + c_2\\3c_1 + 2c_2 + c_3\end{pmatrix} = \begin{pmatrix}0\\0\\0\end{pmatrix}.$$

Which yields a homogeneous system of linear equations

$$c_1 - 2c_3 = 0,$$
$$2c_1 + c_2 = 0,$$
$$3c_1 + 2c_2 + c_3 = 0.$$

This implies that the only solution is the trivial solution

$$c_1 = c_2 = c_3 = 0.$$

So, $S$ is linearly independent.

**Theorem 3.1:** The vectors $|\mathbf{u}_1\rangle, |\mathbf{u}_2\rangle, \ldots, |\mathbf{u}_n\rangle$ are linearly dependent if and only if one of them is a linear combination of the others.

**Definition (Basis):** A set $S = \{|\mathbf{u}_1\rangle, |\mathbf{u}_2\rangle, \ldots, |\mathbf{u}_n\rangle\}$ of vectors is a basis of $V$ if the following two conditions hold:

1- $|\mathbf{u}_1\rangle, |\mathbf{u}_2\rangle, \ldots, |\mathbf{u}_n\rangle$ are linearly independent.

2- $|\mathbf{u}_1\rangle, |\mathbf{u}_2\rangle, \ldots, |\mathbf{u}_n\rangle$ span $V$.

In other words, a set $S = \{|\mathbf{u}_1\rangle, |\mathbf{u}_2\rangle, \ldots, |\mathbf{u}_n\rangle\}$ of vectors is a basis of $V$ if every vector $|\mathbf{u}\rangle \in V$ can be written uniquely as a linear combination of the basis vector.

For example, $S = \{|\mathbf{e}_1\rangle, |\mathbf{e}_2\rangle, \ldots, |\mathbf{e}_n\rangle\}$ is a basis of $\mathbb{R}^n$ and is called the standard basis.

***Example 3.3***

Show that $S = \left\{|\mathbf{u}_1\rangle = \begin{pmatrix}1\\-1\end{pmatrix}, |\mathbf{u}_2\rangle = \begin{pmatrix}3\\2\end{pmatrix}\right\}$ is a basis of $\mathbb{R}^2$.

***Solution***

Suppose that

$$a_1|\mathbf{u}_1\rangle + a_2|\mathbf{u}_2\rangle = \begin{pmatrix}0\\0\end{pmatrix} \Rightarrow \begin{pmatrix}a_1 + 3a_2\\ -a_1 + 2a_2\end{pmatrix} = \begin{pmatrix}0\\0\end{pmatrix}.$$

This leads to the equations

$$a_1 + 3a_2 = 0, \quad -a_1 + 2a_2 = 0.$$

Solving the system of equations, we have $a_1 = a_2 = 0$. So, the vectors $|\mathbf{u}_1\rangle$, $|\mathbf{u}_2\rangle$ are linearly independent. Next, we show that every vector in $\mathbb{R}^2$ can be written as a linear combination of $|\mathbf{u}_1\rangle$ and $|\mathbf{u}_2\rangle$. Let $|\mathbf{u}\rangle = \begin{pmatrix}x\\y\end{pmatrix} \in \mathbb{R}^2$, we must show that there exist scalars $b_1$ and $b_2$ such that

$$|\mathbf{u}\rangle = b_1|\mathbf{u}_1\rangle + b_2|\mathbf{u}_2\rangle \Rightarrow \begin{pmatrix}x\\y\end{pmatrix} = \begin{pmatrix}b_1 + 3b_2\\ -b_1 + 2b_2\end{pmatrix}.$$

This leads to the equations

$$b_1 + 3b_2 = x, \quad -b_1 + 2b_2 = y.$$

Therefore, for any values of $x$ and $y$, we have

$$b_1 = \frac{x+y}{5}, \quad b_2 = \frac{2x-3y}{5}.$$

**Definition (Dimension):** A vector space $V$ is said to be of finite dimension $n$ or to be $n$ dimensional, written

$$\dim V = n, \tag{3.19}$$

if $V$ has such a basis with $n$ elements.

**Definition (Coordinates):** Let $V$ be an $n$-dimensional vector space over $\mathbb{R}$, and suppose $S = \{|\mathbf{u}_1\rangle, |\mathbf{u}_2\rangle, \dots, |\mathbf{u}_n\rangle\}$ is a basis of $V$. Then any vector $|\mathbf{u}\rangle \in V$ can be expressed uniquely as a linear combination of the basis vectors in $S$; say

$$|\mathbf{u}\rangle = a_1|\mathbf{u}_1\rangle + a_2|\mathbf{u}_2\rangle + \cdots + a_n|\mathbf{u}_n\rangle. \tag{3.20}$$

These $n$ scalars $a_1, a_2, \dots, a_n$ are called the coordinates of $|\mathbf{u}\rangle$ relative to the basis $S$; and they form the $n$-tuple $(a_1, a_2, \dots a_n)^T$ in $\mathbb{R}^n$, called the coordinate vector of $|\mathbf{u}\rangle$ relative to $S$. We denote this vector by $[\mathbf{u}]_S$. Thus

$$[\mathbf{u}]_S = (a_1, a_2, \dots a_n)^T. \tag{3.21}$$

## 3.2 Inner Product

The inner product for a vector space $V$ is a map from $V \times V$ to the real numbers [1-6]. We can express this more clearly by saying that the inner product is a function on two vectors that produces a real number which we represent by $\langle\mathbf{u}|\mathbf{v}\rangle$. A vector space that also has an inner product is referred to as an inner product space.

**Definition (Inner Product):** Let $V$ be a vector space. Suppose each pair of vectors $|\mathbf{u}\rangle$, $|\mathbf{v}\rangle \in V$ there is assigned a number, denoted by $\langle\mathbf{u}|\mathbf{v}\rangle$. This function is called an inner product on $V$ if it satisfies the following axioms:

1- $\langle\mathbf{u}|a\mathbf{v}_1 + b\mathbf{v}_2\rangle = a\langle\mathbf{u}|\mathbf{v}_1\rangle + b\langle\mathbf{u}|\mathbf{v}_2\rangle$,

2- $\langle a\mathbf{u}_1 + b\mathbf{u}_2|\mathbf{v}\rangle = a\langle\mathbf{u}_1|\mathbf{v}\rangle + b\langle\mathbf{u}_2|\mathbf{v}\rangle$,

3- $\langle\mathbf{u}|\mathbf{v}\rangle = \langle\mathbf{v}|\mathbf{u}\rangle$,

4- $\langle\mathbf{u}|\mathbf{u}\rangle \geq 0$, and $\langle\mathbf{u}|\mathbf{u}\rangle = 0$ if and only if $|\mathbf{u}\rangle = 0$.

For a finite-dimensional vector space, the inner product can be written as a matrix multiplication of a row vector (bra) with a column vector (ket):

$$\langle\mathbf{u}|\mathbf{v}\rangle = u_1\mathrm{v}_1 + u_2\mathrm{v}_2 + \cdots + u_n\mathrm{v}_n = (u_1, u_2, \dots u_n)\begin{pmatrix}\mathrm{v}_1\\ \mathrm{v}_2\\ \dots\\ \mathrm{v}_n\end{pmatrix}. \tag{3.22}$$

**Example 3.4**

Show that the following function defines an inner product on $\mathbb{R}^2$, where $|\mathbf{u}\rangle = \begin{pmatrix} u_1 \\ u_2 \end{pmatrix}$ and $|\mathbf{v}\rangle = \begin{pmatrix} \mathrm{v}_1 \\ \mathrm{v}_2 \end{pmatrix}$.

$$\langle \mathbf{u}|\mathbf{v}\rangle = \mathbf{u}^T \cdot \mathbf{v} = u_1\mathrm{v}_1 + 2u_2\mathrm{v}_2.$$

***Solution***

Because the product of real numbers is commutative,

$$\begin{aligned}\langle \mathbf{u}|\mathbf{v}\rangle &= u_1\mathrm{v}_1 + 2u_2\mathrm{v}_2 \\ &= \mathrm{v}_1u_1 + 2\mathrm{v}_2u_2 \\ &= \langle \mathbf{v}|\mathbf{u}\rangle,\end{aligned}$$

Let $|\mathbf{w}\rangle = \begin{pmatrix} \mathrm{w}_1 \\ \mathrm{w}_2 \end{pmatrix}$. Then

$$\begin{aligned}\langle \mathbf{u}|\mathbf{v}+\mathbf{w}\rangle &= u_1(\mathrm{v}_1+\mathrm{w}_1) + 2u_2(\mathrm{v}_2+\mathrm{w}_2) \\ &= u_1\mathrm{v}_1 + u_1\mathrm{w}_1 + 2u_2\mathrm{v}_2 + 2u_2\mathrm{w}_2 \\ &= (u_1\mathrm{v}_1 + 2u_2\mathrm{v}_2) + (u_1\mathrm{w}_1 + 2u_2\mathrm{w}_2) \\ &= \langle \mathbf{u}|\mathbf{v}\rangle + \langle \mathbf{u}|\mathbf{w}\rangle.\end{aligned}$$

If $c$ is any scalar, then

$$\begin{aligned}c\langle \mathbf{u}|\mathbf{v}\rangle &= c(u_1\mathrm{v}_1 + 2u_2\mathrm{v}_2) \\ &= (cu_1)\mathrm{v}_1 + 2(cu_2)\mathrm{v}_2 \\ &= \langle c\mathbf{u}|\mathbf{v}\rangle.\end{aligned}$$

Because the square of a real number is non-negative,

$$\langle \mathbf{v}|\mathbf{v}\rangle = \mathrm{v}_1^2 + 2\mathrm{v}_2^2 \geq 0.$$

Moreover, this expression is equal to zero if and only if $|\mathbf{v}\rangle = |\mathbf{0}\rangle$ (that is, if and only if $\mathrm{v}_1 = \mathrm{v}_2 = 0$).

**Definition (Norm):** The square root of the inner product of a vector with itself is called the Euclidean norm, or length of $|\mathbf{u}\rangle$ and is designated by:

$$\|\mathbf{u}\| = \sqrt{\langle \mathbf{u}|\mathbf{u}\rangle}. \tag{3.23}$$

Also, a vector $|\mathbf{u}\rangle$ is said to be normalized if

$$\|\mathbf{u}\| = 1. \tag{3.24}$$

**Definition (Orthogonal Vectors):** Let $V$ be an inner product space. The vectors $|\mathbf{u}\rangle$, $|\mathbf{v}\rangle \in V$ are said to be orthogonal and $|\mathbf{u}\rangle$ is said to be orthogonal to $|\mathbf{v}\rangle$ if

$$\langle \mathbf{u}|\mathbf{v}\rangle = 0. \tag{3.25}$$

**Definition (Orthonormal Vectors):** A set $S = \{|\mathbf{u}_1\rangle, |\mathbf{u}_2\rangle, \dots, |\mathbf{u}_n\rangle\}$ of vectors in $V$ is called orthogonal if each pair of vectors in $S$ are orthogonal, and $S$ is called orthonormal if $S$ is orthogonal and each vector in $S$ has unit length. In other words, $S$ is orthogonal if

$$\langle \mathbf{u}_i|\mathbf{u}_j\rangle = 0, \text{ for } i \neq j, \tag{3.26}$$

and $S$ is orthonormal if

$$\langle \mathbf{u}_i|\mathbf{u}_j\rangle = \delta_{ij} = \begin{cases} 0 & i \neq j, \\ 1 & i = j. \end{cases} \tag{3.27}$$

The completeness, or closure, relation for this basis is given by

$$\textstyle\sum_i |\mathbf{u}_i\rangle\langle \mathbf{u}_i| = \mathbf{I}. \tag{3.28}$$

**Definition (Orthogonal and Orthonormal Basis):** A basis $S = \{|\mathbf{u}_1\rangle, |\mathbf{u}_2\rangle, \dots, |\mathbf{u}_n\rangle\}$ of a vector space $V$ is called an orthogonal basis or an orthonormal basis according as $S$ is an orthogonal set or an orthonormal set of vectors.

**Example 3.5**

Show that the following set is an orthonormal basis for $\mathbb{R}^3$

$$S = \left\{|\mathbf{v}_1\rangle = \frac{1}{\sqrt{2}}\begin{pmatrix} 1 \\ 1 \\ 0 \end{pmatrix}, |\mathbf{v}_2\rangle = \frac{\sqrt{2}}{6}\begin{pmatrix} -1 \\ 1 \\ 4 \end{pmatrix}, |\mathbf{v}_3\rangle = \frac{1}{3}\begin{pmatrix} 2 \\ -2 \\ 1 \end{pmatrix}\right\}.$$

***Solution***

First, we show that the three vectors are mutually orthogonal

$$\langle \mathbf{v}_1|\mathbf{v}_2\rangle = \mathbf{v}_1^T \cdot \mathbf{v}_2 = -\frac{1}{6}+\frac{1}{6}+0 = 0,$$
$$\langle \mathbf{v}_1|\mathbf{v}_3\rangle = \mathbf{v}_1^T \cdot \mathbf{v}_3 = -\frac{2}{3\sqrt{2}}+\frac{2}{3\sqrt{2}}+0 = 0,$$
$$\langle \mathbf{v}_2|\mathbf{v}_3\rangle = \mathbf{v}_2^T \cdot \mathbf{v}_3 = -\frac{\sqrt{2}}{9}-\frac{\sqrt{2}}{9}+2\frac{\sqrt{2}}{9} = 0.$$

Now, each vector is of length 1 because

$$\|\mathbf{v}_1\| = \sqrt{\langle \mathbf{v}_1|\mathbf{v}_1\rangle} = \sqrt{\frac{1}{2}+\frac{1}{2}+0} = 1,$$
$$\|\mathbf{v}_2\| = \sqrt{\langle \mathbf{v}_2|\mathbf{v}_2\rangle} = \sqrt{\frac{2}{36}+\frac{2}{36}+\frac{8}{9}} = 1,$$
$$\|\mathbf{v}_3\| = \sqrt{\langle \mathbf{v}_3|\mathbf{v}_3\rangle} = \sqrt{\frac{4}{9}+\frac{4}{9}+\frac{1}{9}} = 1.$$

**Theorem 3.2:** Suppose $S = \{|\mathbf{u}_1\rangle, |\mathbf{u}_2\rangle, \dots, |\mathbf{u}_n\rangle\}$ is an orthogonal set of nonzero vectors.

1- Then $S$ is linearly independent.

2- For any $|\mathbf{v}\rangle \in V$, we have

$$|\mathbf{v}\rangle = \frac{\langle \mathbf{u}_1|\mathbf{v}\rangle}{\langle \mathbf{u}_1|\mathbf{u}_1\rangle}|\mathbf{u}_1\rangle + \frac{\langle \mathbf{u}_2|\mathbf{v}\rangle}{\langle \mathbf{u}_2|\mathbf{u}_2\rangle}|\mathbf{u}_2\rangle + \cdots + \frac{\langle \mathbf{u}_n|\mathbf{v}\rangle}{\langle \mathbf{u}_n|\mathbf{u}_n\rangle}|\mathbf{u}_n\rangle. \tag{3.29}$$

- It is important to note that if $S = \{|\mathbf{u}_1\rangle, |\mathbf{u}_2\rangle, \dots, |\mathbf{u}_n\rangle\}$ is an orthonormal basis for $V$. Then, for any $|\mathbf{v}\rangle \in V$, we have

$$|\mathbf{v}\rangle = \langle \mathbf{u}_1|\mathbf{v}\rangle|\mathbf{u}_1\rangle + \langle \mathbf{u}_2|\mathbf{v}\rangle|\mathbf{u}_2\rangle + \cdots + \langle \mathbf{u}_n|\mathbf{v}\rangle|\mathbf{u}_n\rangle. \tag{3.30}$$

**Proof:**

We consider here the case of an orthonormal basis. Using the completeness relation, we have

$$\begin{aligned} |\mathbf{v}\rangle &= \mathbf{I}|\mathbf{v}\rangle \\ &= \left(\sum_i |\mathbf{u}_i\rangle\langle \mathbf{u}_i|\right)|\mathbf{v}\rangle \\ &= \sum_i \langle \mathbf{u}_i|\mathbf{v}\rangle|\mathbf{u}_i\rangle \\ &= \sum_i a_i|\mathbf{u}_i\rangle, \end{aligned}$$

where the coefficient $a_i$, which is equal to $\langle \mathbf{u}_i|\mathbf{v}\rangle$, represents the projection of $|\mathbf{v}\rangle$ onto $|\mathbf{u}_i\rangle$; $a_i$ is the component of $|\mathbf{v}\rangle$ along the vector $|\mathbf{u}_i\rangle$. So, within the orthonormal basis $\{|\mathbf{u}_i\rangle\}$, the ket $|\mathbf{v}\rangle$ is represented by the set of components, $a_1, a_2, \dots, a_n$ along $|\mathbf{u}_1\rangle, |\mathbf{u}_2\rangle, \dots, |\mathbf{u}_n\rangle$, respectively.

■

***Example 3.6***

Find the coordinate of $|\mathbf{w}\rangle = \begin{pmatrix} 5 \\ -5 \\ 2 \end{pmatrix}$ relative to the orthonormal basis for $\mathbb{R}^3$ shown below,

$$B = \{|\mathbf{v}_1\rangle = \frac{1}{5}\begin{pmatrix} 3 \\ 4 \\ 0 \end{pmatrix}, |\mathbf{v}_2\rangle = \frac{1}{5}\begin{pmatrix} -4 \\ 3 \\ 0 \end{pmatrix}, |\mathbf{v}_3\rangle = \begin{pmatrix} 0 \\ 0 \\ 1 \end{pmatrix}\}.$$

***Solution***

Because $B$ is orthonormal, we can use Theorem 3.2 to find the required coordinates

$$\langle \mathbf{w}|\mathbf{v}_1\rangle = (5,-5,2)\frac{1}{5}\begin{pmatrix}3\\4\\0\end{pmatrix} = -1, \qquad \langle \mathbf{w}|\mathbf{v}_2\rangle = (5,-5,2)\frac{1}{5}\begin{pmatrix}-4\\3\\0\end{pmatrix} = -7, \qquad \langle \mathbf{w}|\mathbf{v}_3\rangle = (5,-5,2)\begin{pmatrix}0\\0\\1\end{pmatrix} = 2.$$

So, the coordinate relative to $B$ is $[\mathbf{w}]_B = \begin{pmatrix}-1\\-7\\2\end{pmatrix}$.

## Change of basis matrix

Suppose $S = \{|\mathbf{v}_1\rangle, |\mathbf{v}_2\rangle, \dots, |\mathbf{v}_n\rangle\}$ is a basis of a vector space $V$ and suppose $\acute{S} = \{|\mathbf{u}_1\rangle, |\mathbf{u}_2\rangle, \dots, |\mathbf{u}_n\rangle\}$ is another basis. Since $\acute{S}$ is a basis, each vector in $S$ can be written uniquely as a linear combination of the elements in $\acute{S}$ [1-6]. Say,

$$\begin{aligned}|\mathbf{v}_1\rangle &= c_{11}|\mathbf{u}_1\rangle + c_{21}|\mathbf{u}_2\rangle + \cdots + c_{n1}|\mathbf{u}_n\rangle,\\ |\mathbf{v}_2\rangle &= c_{12}|\mathbf{u}_1\rangle + c_{22}|\mathbf{u}_2\rangle + \cdots + c_{n2}|\mathbf{u}_n\rangle,\\ &\dots \quad \dots \quad \dots\\ |\mathbf{v}_n\rangle &= c_{1n}|\mathbf{u}_1\rangle + c_{2n}|\mathbf{u}_2\rangle + \cdots + c_{nn}|\mathbf{u}_n\rangle.\end{aligned} \tag{3.31}$$

Let $\mathbf{Q}$ denote the above matrix of coefficients;

$$\mathbf{Q} = \begin{pmatrix} c_{11} & c_{12} & \cdots & c_{1n}\\ c_{21} & c_{22} & \cdots & c_{2n}\\ \dots & \dots & \dots & \dots\\ c_{n1} & c_{n2} & \cdots & c_{nn}\end{pmatrix}. \tag{3.32}$$

That is, $\mathbf{Q} = (c_{ij})$. Then $\mathbf{Q}$ is called the change of basis matrix from the old basis $S$ to the new basis $\acute{S}$. Using the completeness relation, we have

$$\begin{aligned}|\mathbf{v}_j\rangle &= \mathbf{I}|\mathbf{v}_j\rangle\\ &= \left(\sum_i |\mathbf{u}_i\rangle\langle \mathbf{u}_i|\right)|\mathbf{v}_j\rangle\\ &= \sum_i |\mathbf{u}_i\rangle\langle \mathbf{u}_i|\mathbf{v}_j\rangle = \sum_i c_{ij}\,|\mathbf{u}_i\rangle,\end{aligned} \tag{3.33}$$

where $c_{ij} = \langle \mathbf{u}_i|\mathbf{v}_j\rangle$.

**Theorem 3.3:** Let $\mathbf{Q}$ be the change of basis matrix from a basis $S$ to a basis $\acute{S}$ in a vector space $V$. Then, for any vector $|\mathbf{u}\rangle \in V$, we have

$$\mathbf{Q}[\mathbf{u}]_s = [\mathbf{u}]_{\acute{s}}\,. \tag{3.34}$$

Let $\mathbf{P} = \mathbf{Q}^{-1}$ be the change of basis matrix from a basis $\acute{S}$ to a basis $S$ in a vector space $V$. Then, for any vector $|\mathbf{u}\rangle \in V$, we have

$$\mathbf{Q}^{-1}\mathbf{Q}[\mathbf{u}]_s = [\mathbf{u}]_s = \mathbf{P}[\mathbf{u}]_{\acute{s}}\,. \tag{3.35}$$

***Example 3.7***

Let $S = \left\{|\mathbf{v}_1\rangle = \begin{pmatrix}1\\1\end{pmatrix}, |\mathbf{v}_2\rangle = \begin{pmatrix}2\\1\end{pmatrix}\right\}$ and $\acute{S} = \left\{|\mathbf{u}_1\rangle = \begin{pmatrix}1\\0\end{pmatrix}, |\mathbf{u}_2\rangle = \begin{pmatrix}0\\1\end{pmatrix}\right\}$ are bases of $\mathbb{R}^2$.

1- Find the transformation matrix $\mathbf{Q}$.

2- Find $[\mathbf{u}]_{\acute{s}}$, given that $[\mathbf{u}]_s = \begin{pmatrix}-3\\5\end{pmatrix}$.

***Solution***

It is easy to check that

$$|\mathbf{v}_1\rangle = |\mathbf{u}_1\rangle + |\mathbf{u}_2\rangle, \qquad |\mathbf{v}_2\rangle = 2|\mathbf{u}_1\rangle + |\mathbf{u}_2\rangle.$$

Hence, the transformation matrix $\mathbf{Q}$ from $S$ to $\acute{S}$ is the matrix

$$\mathbf{Q} = \begin{pmatrix}1 & 2\\1 & 1\end{pmatrix},$$

and,

$$\begin{pmatrix}1 & 2\\1 & 1\end{pmatrix}\begin{pmatrix}-3\\5\end{pmatrix} = \begin{pmatrix}x\\y\end{pmatrix} \;\Rightarrow\; [\mathbf{u}]_{\acute{s}} = \begin{pmatrix}x\\y\end{pmatrix} = \begin{pmatrix}7\\2\end{pmatrix}.$$

## 3.3 Subspace

**Definition (Subspace):** Let $W$ be a subset of a vector space $V$ over a field $K$. $W$ is called a subspace of $V$ if $W$ is itself a vector space over $K$ with respect to the operations of vector addition and scalar multiplication on $V$.

**Theorem 3.4:** Suppose $W$ is a subset of a vector space $V$. Then $W$ is a subspace of $V$ if and only if the following hold:

1- $|\mathbf{0}\rangle \in W$,
2- $W$ is closed under vector addition, that is: for every $|\mathbf{u}\rangle, |\mathbf{v}\rangle \in W$, the sum $|\mathbf{u}\rangle + |\mathbf{v}\rangle \in W$,
3- $W$ is closed under scalar multiplication, that is: for every $|\mathbf{u}\rangle \in W$, $k \in K$, the multiple $k|\mathbf{u}\rangle \in W$.

Because a subspace of a vector space is a vector space, it must contain the zero vector. In fact, the simplest subspace of a vector space is the one consisting of only the zero vector, $W = \{|\mathbf{0}\rangle\}$. This subspace is called the zero subspace. Another obvious subspace of $V$ is $V$ itself. Every vector space contains these two trivial subspaces, and subspaces other than these two are called proper (or nontrivial) subspaces.

***Example 3.8***

Let $U$ consist of all vectors in $\mathbb{R}^3$ whose entries are equal; that is,

$$U = \{\begin{pmatrix} a \\ b \\ c \end{pmatrix}, a = b = c\}.$$

For example, $\begin{pmatrix} 1 \\ 1 \\ 1 \end{pmatrix}, \begin{pmatrix} -3 \\ -3 \\ -3 \end{pmatrix}, \begin{pmatrix} 7 \\ 7 \\ 7 \end{pmatrix}$ are vectors in $U$. Geometrically, $U$ is the line through the origin $O$ and the point $\begin{pmatrix} 1 \\ 1 \\ 1 \end{pmatrix}$. Clearly $|\mathbf{0}\rangle = \begin{pmatrix} 0 \\ 0 \\ 0 \end{pmatrix}$ belongs to $U$, because all entries in 0 are equal. Further, suppose $|\mathbf{u}\rangle$ and $|\mathbf{v}\rangle$ are arbitrary vectors in $U$, say,

$$|\mathbf{u}\rangle = \begin{pmatrix} a \\ a \\ a \end{pmatrix} \text{ and } |\mathbf{v}\rangle = \begin{pmatrix} b \\ b \\ b \end{pmatrix}.$$

Then, for any scalar $k \in \mathbb{R}$, the following are also vectors in $U$:

$$|\mathbf{u}\rangle + |\mathbf{v}\rangle = \begin{pmatrix} a+b \\ a+b \\ a+b \end{pmatrix} \text{ and } k|\mathbf{u}\rangle = \begin{pmatrix} ka \\ ka \\ ka \end{pmatrix}.$$

Thus, $U$ is a subspace of $\mathbb{R}^3$.

**Theorem 3.5:** If $W_1$ and $W_2$ are both subspaces of a vector space $U$ then the intersection of $W_1$ and $W_2$ (denoted by $W_1 \cap W_2$ is also a subspace of $U$.

**Proof:**

Let $W_1$ and $W_2$ be subspaces of a vector space $V$. Clearly, $|\mathbf{0}\rangle \in W_1$ and $|\mathbf{0}\rangle \in W_2$, because $W_1$ and $W_2$ are subspaces; whence $|\mathbf{0}\rangle \in W_1 \cap W_2$. Now suppose $|\mathbf{v}_1\rangle$ and $|\mathbf{v}_2\rangle$ belong to the intersection $W_1 \cap W_2$. Then $|\mathbf{v}_1\rangle, |\mathbf{v}_2\rangle \in W_1$ and $|\mathbf{v}_1\rangle, |\mathbf{v}_2\rangle \in W_2$. Further, because $W_1$ and $W_2$ are subspaces, for any scalars $k_1, k_2 \in \mathbb{R}$, $k_1|\mathbf{v}_1\rangle + k_2|\mathbf{v}_2\rangle \in W_1$and $k_1|\mathbf{v}_1\rangle + k_2|\mathbf{v}_2\rangle \in W_2$. Thus, $k_1|\mathbf{v}_1\rangle + k_2|\mathbf{v}_2\rangle \in W_1 \cap W_2$. Therefore, $W_1 \cap W_2$ is a subspace of $V$.

■

**Definition (Sum of Subsets of Vector Space):** Let $V$ be a vector space. Suppose that $S_1$ and $S_2$ are non-empty subsets of $V$. The sum of $S_1$ and $S_2$, denoted by $S_1 + S_2$, is

$$S_1 + S_2 = \{|\mathbf{u}\rangle + |\mathbf{v}\rangle; |\mathbf{u}\rangle \in S_1, |\mathbf{v}\rangle \in S_2\}. \quad (3.36)$$

That is, $S_1 + S_2$ is the set of vectors of $V$ that can be obtained by adding a vector in $S_1$ to a vector in $S_2$.

**Theorem 3.6:** Let $V$ be a vector space and suppose $W_1$ and $W_2$ are subspaces of $V$. Then $W_1 + W_2$ is a subspace of $V$ that contains $W_1$ and $W_2$.

**Proof:**

Because $W_1$ and $W_2$ are subspaces, $|\mathbf{0}\rangle \in W_1$ and $|\mathbf{0}\rangle \in W_2$. Hence, $|\mathbf{0}\rangle = |\mathbf{0}\rangle + |\mathbf{0}\rangle$ belongs to $W_1 + W_2$. Now suppose $|\mathbf{v}_1\rangle, |\mathbf{v}_2\rangle \in W_1 + W_2$. Then $|\mathbf{v}_1\rangle = |\mathbf{u}\rangle + |\mathbf{v}\rangle$ and $|\mathbf{v}_2\rangle = |\acute{\mathbf{u}}\rangle + |\acute{\mathbf{v}}\rangle$, where $|\mathbf{u}\rangle, |\acute{\mathbf{u}}\rangle \in W_1$ and $|\mathbf{v}\rangle, |\acute{\mathbf{v}}\rangle \in W_2$. Then

$$k_1|\mathbf{v}_1\rangle + k_2|\mathbf{v}_2\rangle = k_1|\mathbf{u}\rangle + k_1|\mathbf{v}\rangle + k_2|\acute{\mathbf{u}}\rangle + k_2|\acute{\mathbf{v}}\rangle \in W_1 + W_2.$$

Thus, $W_1 + W_2$ is a subspace of $V$.

■

**Definition (Direct Sum of Subspaces):** A vector space $V$ is called the direct sum of $W_1$ and $W_2$ if $W_1$ and $W_2$ are subspaces of $V$ such that

$$W_1 \cap W_2 = \{|\mathbf{0}\rangle\}, \quad (3.37)$$

and

$$W_1 + W_2 = V. \quad (3.38)$$

We denote that $V$ is the direct sum of $W_1$ and $W_2$ by writing

$$V = W_1 \oplus W_2. \quad (3.39)$$

**Theorem 3.7:** Let $V$ be a vector space and suppose $W_1$ and $W_2$ are subspaces of $V$. Then $V = W_1 \oplus W_2$ if and only if each vector $|\mathbf{u}\rangle \in V$ can be written uniquely as

$$|\mathbf{u}\rangle = |\mathbf{u}_1\rangle + |\mathbf{u}_2\rangle, \quad (3.40)$$

where $|\mathbf{u}_1\rangle \in W_1$ and $|\mathbf{u}_2\rangle \in W_2$.

**Proof:**

To prove this theorem, we prove each direction separately:

1. Suppose $V = W_1 \oplus W_2$. Let $|\mathbf{u}\rangle \in V$ be an arbitrary vector. To see that $|\mathbf{u}\rangle$ can be written as $|\mathbf{u}\rangle = |\mathbf{u}_1\rangle + |\mathbf{u}_2\rangle$ where $|\mathbf{u}_1\rangle \in W_1$ and $|\mathbf{u}_2\rangle \in W_2$, we recall that since

$$V = W_1 \oplus W_2 = W_1 + W_2.$$

So,

$$|\mathbf{u}\rangle \in W_1 + W_2.$$

By the definition of $W_1 + W_2$, $|\mathbf{u}\rangle$ can be written as $|\mathbf{u}\rangle = |\mathbf{u}_1\rangle + |\mathbf{u}_2\rangle$ where $|\mathbf{u}_1\rangle \in W_1$ and $|\mathbf{u}_2\rangle \in W_2$.

To see that there is only one way that $|\mathbf{u}\rangle$ can be written as $|\mathbf{u}\rangle = |\mathbf{u}_1\rangle + |\mathbf{u}_2\rangle$, suppose there exists $|\mathbf{v}_1\rangle, |\acute{\mathbf{v}}_1\rangle \in W_1$ and $|\mathbf{v}_2\rangle, |\acute{\mathbf{v}}_2\rangle \in W_2$ such that

$$|\mathbf{u}\rangle = |\mathbf{v}_1\rangle + |\mathbf{v}_2\rangle, \qquad |\mathbf{u}\rangle = |\acute{\mathbf{v}}_1\rangle + |\acute{\mathbf{v}}_2\rangle.$$

Hence,

$$|\mathbf{v}_1\rangle - |\acute{\mathbf{v}}_1\rangle = |\acute{\mathbf{v}}_2\rangle - |\mathbf{v}_2\rangle.$$

Since, $|\mathbf{v}_1\rangle, |\acute{\mathbf{v}}_1\rangle \in W_1$ and $W_1$ is a subspace of $V$, we have $|\mathbf{v}_1\rangle - |\acute{\mathbf{v}}_1\rangle \in W_1$. Similarly, $|\acute{\mathbf{v}}_2\rangle - |\mathbf{v}_2\rangle \in W_2$. As a result,

$$|\mathbf{v}_1\rangle - |\acute{\mathbf{v}}_1\rangle = |\mathbf{v}_2\rangle - |\acute{\mathbf{v}}_2\rangle \in W_1 \cap W_2 = \{|\mathbf{0}\rangle\},$$

so $|\mathbf{v}_1\rangle = |\acute{\mathbf{v}}_1\rangle$ and $|\mathbf{v}_2\rangle = |\acute{\mathbf{v}}_2\rangle$. i.e., each vector $|\mathbf{u}\rangle \in V$ can be written uniquely as $|\mathbf{u}\rangle = |\mathbf{u}_1\rangle + |\mathbf{u}_2\rangle$ where $|\mathbf{u}_1\rangle \in W_1$ and $|\mathbf{u}_2\rangle \in W_2$.

2. We leave the proof of part (2) as an exercise.

■

**Theorem 3.8:** Let $V$ be a vector space and suppose $W_1$ and $W_2$ are finite-dimensional subspaces of $V$ such that $V = W_1 \oplus W_2$. Then $V$ is a finite-dimensional vector space. Moreover, if $\beta = \{|\mathbf{u}_1\rangle, |\mathbf{u}_2\rangle, \ldots, |\mathbf{u}_k\rangle\}$ is a basis for $W_1$ and $\gamma = \{|\mathbf{v}_1\rangle, |\mathbf{v}_2\rangle, \ldots, |\mathbf{v}_m\rangle\}$ is a basis for $W_2$, then

$$\alpha = \{|\mathbf{u}_1\rangle, |\mathbf{u}_2\rangle, \ldots, |\mathbf{u}_k\rangle, |\mathbf{v}_1\rangle, |\mathbf{v}_2\rangle, \ldots, |\mathbf{v}_m\rangle\}, \tag{3.41}$$

is a basis for $V$. Thus,

$$\dim V = \dim W_1 + \dim W_2. \tag{3.42}$$

**Theorem 3.9:** Suppose $W_1, W_2, \ldots, W_r$ are subspaces of $V$, and suppose

$$\begin{aligned} B_1 &= \{|\mathbf{u}_{11}\rangle, |\mathbf{u}_{21}\rangle, \ldots, |\mathbf{u}_{n_1 1}\rangle\}, \\ B_2 &= \{|\mathbf{u}_{12}\rangle, |\mathbf{u}_{22}\rangle, \ldots, |\mathbf{u}_{n_2 2}\rangle\}, \\ &\ldots \quad \ldots \quad \ldots \\ B_r &= \{|\mathbf{u}_{1r}\rangle, |\mathbf{u}_{2r}\rangle, \ldots, |\mathbf{u}_{n_r r}\rangle\}. \end{aligned} \tag{3.43}$$

are bases of $W_1, W_2, \ldots, W_r$, respectively. Then $V$ is the direct sum of $W_i$ if and only if the union

$$B = B_1 \cup B_2 \cup \ldots \cup B_r, \tag{3.44}$$

is a basis of $V$.

## 3.4 Linear Operator, Linear Transformation, and Matrix Representation

**Definition (Linear Transformation):** Let $V$ and $U$ be vector spaces over the same field $K$. A mapping $\mathbf{T}: V \to U$ is called a linear mapping (or linear transformation) if it satisfies the following two conditions [1-6]:

1- For any $|\mathbf{u}\rangle, |\mathbf{v}\rangle \in V$, $\mathbf{T}(|\mathbf{u}\rangle + |\mathbf{v}\rangle) = \mathbf{T}|\mathbf{u}\rangle + \mathbf{T}|\mathbf{v}\rangle$.

2- For any $k \in K$ and any $|\mathbf{u}\rangle \in V$, $\mathbf{T}(k|\mathbf{u}\rangle) = k\mathbf{T}|\mathbf{u}\rangle$.

In other words, let $V$ and $U$ be vector spaces over the same field $K$. A mapping $\mathbf{T}: V \to U$ is called a linear mapping (or linear transformation) if it satisfies the following condition, for any vectors $|\mathbf{u}\rangle, |\mathbf{v}\rangle \in V$ and any scalars $a, b \in K$

$$\mathbf{T}(a|\mathbf{u}\rangle + b|\mathbf{v}\rangle) = a\mathbf{T}|\mathbf{u}\rangle + b\mathbf{T}|\mathbf{v}\rangle. \tag{3.45}$$

***Example 3.9***

Show that $\mathbf{T}: \mathbb{R}^2 \to \mathbb{R}^3$ defined by $\mathbf{T}\begin{pmatrix} x \\ y \end{pmatrix} = \begin{pmatrix} x \\ x+y \\ x-y \end{pmatrix}$ is a linear transformation.

***Solution:***

We need to verify the two conditions of the definition. Given $\begin{pmatrix} x_1 \\ y_1 \end{pmatrix}$ and $\begin{pmatrix} x_2 \\ y_2 \end{pmatrix}$ in $\mathbb{R}^2$, compute

$$\mathbf{T}\left(\begin{pmatrix} x_1 \\ y_1 \end{pmatrix} + \begin{pmatrix} x_2 \\ y_2 \end{pmatrix}\right) = \mathbf{T}\begin{pmatrix} x_1 + x_2 \\ y_1 + y_2 \end{pmatrix} = \begin{pmatrix} x_1 + x_2 \\ x_1 + x_2 + y_1 + y_2 \\ x_1 + x_2 - y_1 - y_2 \end{pmatrix} = \begin{pmatrix} x_1 \\ x_1 + y_1 \\ x_1 - y_1 \end{pmatrix} + \begin{pmatrix} x_2 \\ x_2 + y_2 \\ x_2 - y_2 \end{pmatrix} = \mathbf{T}\begin{pmatrix} x_1 \\ y_1 \end{pmatrix} + \mathbf{T}\begin{pmatrix} x_2 \\ y_2 \end{pmatrix}.$$

This proves the first condition. For the second condition, we let $\alpha \in \mathbb{R}$ and compute

$$\mathbf{T}\left(\alpha\begin{pmatrix} x \\ y \end{pmatrix}\right) = \mathbf{T}\begin{pmatrix} \alpha x \\ \alpha y \end{pmatrix} = \begin{pmatrix} \alpha x \\ \alpha x + \alpha y \\ \alpha x - \alpha y \end{pmatrix} = \alpha\begin{pmatrix} x \\ x+y \\ x-y \end{pmatrix} = \alpha\mathbf{T}\begin{pmatrix} x \\ y \end{pmatrix}.$$

Hence $\mathbf{T}$ is a linear transformation.

**Definition (Null Space or Kernel of Linear Map):** Let $\mathbf{T}: V \to U$ be a linear map. Then the null space or kernel of **T**, denoted by null $\mathbf{T}$, is the set of all vectors in $V$ that map to zero:

$$\text{null}\, \mathbf{T} = \{|\mathbf{v}\rangle \in V : \mathbf{T}|\mathbf{v}\rangle = |\mathbf{0}\rangle\}. \tag{3.46}$$

**Theorem 3.10:** Let $\mathbf{T}: V \to U$ be a linear map. Then null $\mathbf{T}$ is a subspace of $V$.

**Definition (Injective Map):** The linear map $\mathbf{T}: V \to U$ is called injective if for all $|\mathbf{u}\rangle, |\mathbf{v}\rangle \in V$, the condition $\mathbf{T}|\mathbf{v}\rangle = \mathbf{T}|\mathbf{u}\rangle$ implies that $|\mathbf{v}\rangle = |\mathbf{u}\rangle$. In other words, different vectors in $V$ are mapped to different vectors in $U$.

**Definition (Range of Map):** Let $\mathbf{T}: V \to U$ be a linear map. The range of $\mathbf{T}$, denoted by range $\mathbf{T}$ , is the subset of vectors of $U$ that are in the image of $\mathbf{T}$

$$\text{range } \mathbf{T} = \{|\mathbf{u}\rangle \in U: \text{there exists } |\mathbf{v}\rangle \in V \text{ such that } \mathbf{T}|\mathbf{v}\rangle = |\mathbf{u}\rangle\}. \tag{3.47}$$

**Theorem 3.11:** Let $\mathbf{T}: V \to U$ be a linear map. Then range $\mathbf{T}$ is a subspace of $U$.

**Definition (Linear Operator):** Let $V$ be a vector space over a field $K$. We now consider the special case of linear mappings $\mathbf{F}: V \to V$, i.e., from $V$ into itself. They are also called linear operators on $V$.

Let $\mathbf{T}$ be a linear operator on a vector space $V$ over a field $K$ and suppose $S = \{|\mathbf{u}_1\rangle, |\mathbf{u}_2\rangle, \dots, |\mathbf{u}_n\rangle\}$ is a basis of $V$. Now $\mathbf{T}|\mathbf{u}_1\rangle,\dots, \mathbf{T}|\mathbf{u}_n\rangle$ are vectors in $V$ and so each is a linear combination of the vectors in the basis $S$; say,

$$\begin{aligned}
\mathbf{T}|\mathbf{u}_1\rangle &= a_{11}|\mathbf{u}_1\rangle + a_{21}|\mathbf{u}_2\rangle + \cdots + a_{n1}|\mathbf{u}_n\rangle,\\
\mathbf{T}|\mathbf{u}_2\rangle &= a_{12}|\mathbf{u}_1\rangle + a_{22}|\mathbf{u}_2\rangle + \cdots + a_{n2}|\mathbf{u}_n\rangle,\\
&\dots \dots \dots\\
\mathbf{T}|\mathbf{u}_n\rangle &= a_{1n}|\mathbf{u}_1\rangle + a_{2n}|\mathbf{u}_2\rangle + \cdots + a_{nn}|\mathbf{u}_n\rangle.
\end{aligned} \tag{3.48}$$

The following definition applies

**Definition (Matrix Representation of an Operator):** The above matrix of coefficients, denoted by $[\mathbf{T}]_S$, is called the matrix representation of $\mathbf{T}$ relative to the basis $S$ or simply the matrix of $\mathbf{T}$ in the basis $S$; that is,

$$[\mathbf{T}]_S = \begin{pmatrix} a_{11} & a_{12} & \cdots & a_{1n} \\ a_{21} & a_{22} & \cdots & a_{2n} \\ \cdots & \cdots & \cdots & \cdots \\ a_{n1} & a_{n2} & \cdots & a_{nn} \end{pmatrix}. \tag{3.49}$$

Using the completeness relation, we have

$$\begin{aligned}
\mathbf{T}|\mathbf{u}_j\rangle &= \mathbf{I}\,\mathbf{T}|\mathbf{u}_j\rangle\\
&= \left(\sum_i |\mathbf{u}_i\rangle\langle\mathbf{u}_i|\right)\mathbf{T}|\mathbf{u}_j\rangle\\
&= \sum_i |\mathbf{u}_i\rangle\langle\mathbf{u}_i|\mathbf{T}|\mathbf{u}_j\rangle\\
&= \sum_i a_{ij}|\mathbf{u}_i\rangle,
\end{aligned} \tag{3.50}$$

where $a_{ij} = \langle\mathbf{u}_i|\mathbf{T}|\mathbf{u}_j\rangle$.

***Example 3.10***

Let $\mathbf{T}$ be defined by $\mathbf{T}\begin{pmatrix} x \\ y \end{pmatrix} = \begin{pmatrix} x + 2y \\ 2x - y \end{pmatrix}$. Let $S = \left\{|\mathbf{e}_1\rangle = \begin{pmatrix} 1 \\ 0 \end{pmatrix}, |\mathbf{e}_2\rangle = \begin{pmatrix} 0 \\ 1 \end{pmatrix}\right\}$ be the standard basis of $\mathbb{R}^2$ . Find the matrix representation of $\mathbf{T}$ with respect to $S$.

***Solution***

We have the following computation

$$\mathbf{T}|\mathbf{e}_1\rangle = \mathbf{T}\begin{pmatrix} 1 \\ 0 \end{pmatrix} = \begin{pmatrix} 1 \\ 2 \end{pmatrix} = \begin{pmatrix} 1 \\ 0 \end{pmatrix} + 2\begin{pmatrix} 0 \\ 1 \end{pmatrix} = |\mathbf{e}_1\rangle + 2|\mathbf{e}_2\rangle,$$

$$\mathbf{T}|\mathbf{e}_2\rangle = \mathbf{T}\begin{pmatrix} 0 \\ 1 \end{pmatrix} = \begin{pmatrix} 2 \\ -1 \end{pmatrix} = 2\begin{pmatrix} 1 \\ 0 \end{pmatrix} - \begin{pmatrix} 0 \\ 1 \end{pmatrix} = 2|\mathbf{e}_1\rangle - |\mathbf{e}_2\rangle.$$

Therefore, the matrix representation of $\mathbf{T}$ is

$$[\mathbf{T}]_S = \begin{pmatrix} 1 & 2 \\ 2 & -1 \end{pmatrix}.$$

**Theorem 3.12:** Let $S = \{|\mathbf{u}_1\rangle, |\mathbf{u}_2\rangle, \dots, |\mathbf{u}_n\rangle\}$ be a basis for $V$ and let $\mathbf{T}$ be any linear operator on $V$. Then, for any vector $|\mathbf{u}\rangle \in V$,

$$[\mathbf{T}]_s[\mathbf{u}]_s = [\mathbf{T}|\mathbf{u}\rangle]_s. \tag{3.51}$$

That is if we multiply the matrix representation of $\mathbf{T}$ by the coordinate vector of $|\mathbf{u}\rangle$, then we obtain the coordinate vector of $\mathbf{T}|\mathbf{u}\rangle$.

The previous discussion shows that we can represent a linear operator by a matrix once we have chosen a basis. We ask the following natural question: How does our representation change if we select another basis?

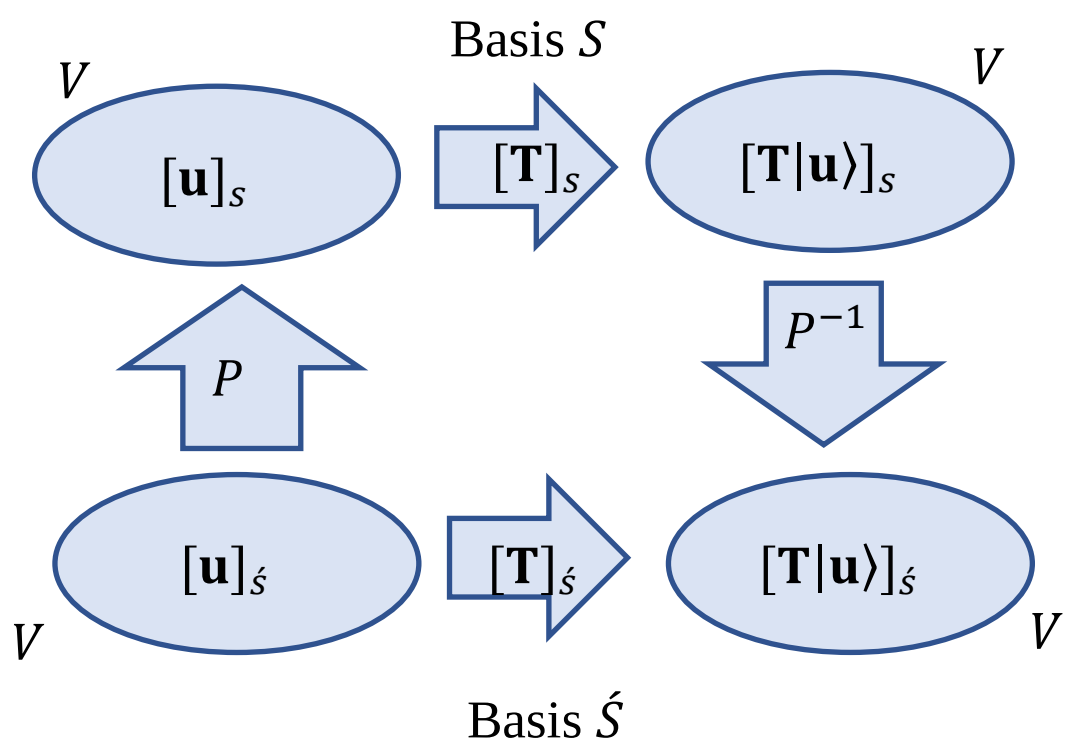

**Figure 3.1.** The relationships between coordinates of a vector in two different bases.

From Figure 3.1, there are two ways to get from the coordinate vector $[\mathbf{u}]_{\acute{s}}$ to the coordinate vector $[\mathbf{T}|\mathbf{u}\rangle]_{\acute{s}}$. One way is direct, using the matrix $[\mathbf{T}]_{\acute{s}}$ to obtain

$$[\mathbf{T}]_{\acute{s}}[\mathbf{u}]_{\acute{s}} = [\mathbf{T}|\mathbf{u}\rangle]_{\acute{s}}. \tag{3.52}$$

The other way is indirect, using the matrices $\mathbf{P}$, $[\mathbf{T}]_s$ and $\mathbf{P}^{-1}$ to obtain

$$\mathbf{P}^{-1}[\mathbf{T}]_s\mathbf{P}[\mathbf{u}]_{\acute{s}} = [\mathbf{T}|\mathbf{u}\rangle]_{\acute{s}}. \tag{3.53}$$

So that we get

$$[\mathbf{T}]_{\acute{s}} = \mathbf{P}^{-1}[\mathbf{T}]_s\mathbf{P}. \tag{3.54}$$

**Theorem 3.13:** Let $\mathbf{P}$ be the change of basis matrix from a basis $\acute{S}$ to a basis $S$ in a vector space $V$. Then, for any linear operator $\mathbf{T}$ on $V$,

$$[\mathbf{T}]_{\acute{s}} = \mathbf{P}^{-1}[\mathbf{T}]_s\mathbf{P}. \tag{3.55}$$

***Example 3.11***

Find the matrix representation of $\mathbf{T}\begin{pmatrix} x_1 \\ x_2 \end{pmatrix} = \begin{pmatrix} 2x_1 - 2x_2 \\ -x_1 + 3x_2 \end{pmatrix}$ with respect to $S$ where $S = \left\{|\mathbf{e}_1\rangle = \begin{pmatrix} 1 \\ 0 \end{pmatrix}, |\mathbf{e}_2\rangle = \begin{pmatrix} 0 \\ 1 \end{pmatrix}\right\}$.

Then, by using the result of the standard basis of $\mathbb{R}^2$, find the matrix $[\mathbf{T}]_{\acute{s}}$, where $\acute{S} = \left\{|\acute{\mathbf{e}}_1\rangle = \begin{pmatrix} 1 \\ 0 \end{pmatrix}, |\acute{\mathbf{e}}_2\rangle = \begin{pmatrix} 1 \\ 1 \end{pmatrix}\right\}$,

***Solution***

We have the following computation

$$\mathbf{T}|\mathbf{e}_1\rangle = \mathbf{T}\begin{pmatrix} 1 \\ 0 \end{pmatrix} = \begin{pmatrix} 2 \\ -1 \end{pmatrix} = 2\begin{pmatrix} 1 \\ 0 \end{pmatrix} - 1\begin{pmatrix} 0 \\ 1 \end{pmatrix} = 2|\mathbf{e}_1\rangle - 1|\mathbf{e}_2\rangle,$$
$$\mathbf{T}|\mathbf{e}_2\rangle = \mathbf{T}\begin{pmatrix} 0 \\ 1 \end{pmatrix} = \begin{pmatrix} -2 \\ 3 \end{pmatrix} = -2\begin{pmatrix} 1 \\ 0 \end{pmatrix} + 3\begin{pmatrix} 0 \\ 1 \end{pmatrix} = -2|\mathbf{e}_1\rangle + 3|\mathbf{e}_2\rangle.$$

Therefore, the matrix representation of $\mathbf{T}$ is

$$[\mathbf{T}]_s = \begin{pmatrix} 2 & -2 \\ -1 & 3 \end{pmatrix}.$$

It is easy to check that

$$|\acute{\mathbf{e}}_1\rangle = \begin{pmatrix}1\\0\end{pmatrix} = |\mathbf{e}_1\rangle + 0 \times |\mathbf{e}_2\rangle,$$
$$|\acute{\mathbf{e}}_2\rangle = \begin{pmatrix}1\\1\end{pmatrix} = |\mathbf{e}_1\rangle + |\mathbf{e}_2\rangle.$$

The transformation matrix **P** is

$$\mathbf{P} = \begin{pmatrix}1 & 1\\0 & 1\end{pmatrix},$$

and the matrix $\mathbf{P}^{-1}$ is

$$\mathbf{P}^{-1} = \begin{pmatrix}1 & -1\\0 & 1\end{pmatrix}.$$

Hence, we have

$$[\mathbf{T}]_{\acute{s}} = \mathbf{P}^{-1}[\mathbf{T}]_s\mathbf{P} = \begin{pmatrix}1 & -1\\0 & 1\end{pmatrix}\begin{pmatrix}2 & -2\\-1 & 3\end{pmatrix}\begin{pmatrix}1 & 1\\0 & 1\end{pmatrix} = \begin{pmatrix}3 & -2\\-1 & 2\end{pmatrix}.$$

## 3.5 Eigenvalues, Eigenvectors and Eigenspaces

**Definition (Similar Matrices):** If **A** and **B** are $n \times n$ matrices such that there is an invertible $n \times n$ matrix **P** with

$$\mathbf{B} = \mathbf{P}^{-1}\mathbf{A}\mathbf{P}, \tag{3.56}$$

then **A** and **B** are called similar.

**Remarks:**

1- With **P** being the change of basis matrix, similar matrices represent the same linear operator under two different bases.
2- If $\mathbf{B} = \mathbf{P}^{-1}\mathbf{A}\mathbf{P}$ then $\mathbf{A} = \mathbf{P}\mathbf{B}\mathbf{P}^{-1}$, so if **A** is similar to **B**, then **B** is similar to **A**.
3- If $\mathbf{A} = \mathbf{P}_1^{-1}\mathbf{B}\mathbf{P}_1$ and $\mathbf{B} = \mathbf{P}_2^{-1}\mathbf{C}\mathbf{P}_2$ then

$$\mathbf{A} = \mathbf{P}_1^{-1}(\mathbf{P}_2^{-1}\mathbf{C}\mathbf{P}_2)\mathbf{P}_1 = (\mathbf{P}_2\mathbf{P}_1)^{-1}\mathbf{C}(\mathbf{P}_2\mathbf{P}_1). \tag{3.57}$$

So, if **A** is similar to **B**, and **B** is similar to **C**, then **A** is similar to **C**. This allows us to put matrices into families in which all the matrices in a family are similar to each other. Then each family can be represented by a diagonal (or nearly diagonal) matrix.
4- If $\mathbf{B} = \mathbf{P}^{-1}\mathbf{A}\mathbf{P}$, then, $\det(\mathbf{A}) = \det(\mathbf{B})$.
5- If $\mathbf{B} = \mathbf{P}^{-1}\mathbf{A}\mathbf{P}$, then, $\mathrm{Tr}(\mathbf{A}) = \mathrm{Tr}(\mathbf{B})$.

Consider a polynomial $f(\lambda)$ over a field $K$; say

$$f(\lambda) = a_n\lambda^n + \cdots + a_1\lambda + a_0. \tag{3.58}$$

If **A** is the square matrix over $K$, then we define

$$f(\mathbf{A}) = a_n\mathbf{A}^n + \cdots + a_1\mathbf{A} + a_0\mathbf{I}, \tag{3.59}$$

where **I** is the identity matrix. In particular, we say that **A** is a root or zero of the polynomial $f(\lambda)$ if $f(\mathbf{A}) = \mathbf{0}$.

**Definition (Characteristic Matrix):** The matrix $\lambda\mathbf{I}_n - \mathbf{A}$, where $\mathbf{I}_n$ is the $n$-square identity matrix, and $\lambda$ is indeterminate, is called the characteristic matrix of **A**:

$$\lambda\mathbf{I}_n - \mathbf{A} = \begin{pmatrix}\lambda - A_{11} & -A_{21} & \cdots & -A_{n1}\\ -A_{12} & \lambda - A_{22} & \cdots & -A_{n2}\\ \cdots & \cdots & \cdots & \cdots\\ -A_{1n} & -A_{2n} & \cdots & \lambda - A_{nn}\end{pmatrix}. \tag{3.60}$$

Its determinant $\Delta_{\mathbf{A}}(\lambda) = \det(\lambda\mathbf{I}_n - \mathbf{A})$, which is a polynomial in $\lambda$, is called the characteristic polynomial of **A**. We also call

$$\Delta_{\mathbf{A}}(\lambda) = \det(\lambda\mathbf{I}_n - \mathbf{A}) = 0, \tag{3.61}$$

the characteristic equation of **A**.

**Theorem 4.14 (Cayley-Hamilton Theorem):** Every matrix is a zero of its characteristic polynomial.

For example, if $\mathbf{A} = \begin{pmatrix} 4 & -1 \\ 2 & 1 \end{pmatrix}$, the characteristic polynomial of $A$ is $\Delta_{\mathbf{A}}(\lambda) = \begin{vmatrix} \lambda - 4 & 1 \\ -2 & \lambda - 1 \end{vmatrix} = \lambda^2 - 5\lambda + 6$, and $\mathbf{A}$ is a zero of its characteristic polynomial $f(\mathbf{A}) = \mathbf{A}^2 - 5\mathbf{A} + 6\mathbf{I} = \begin{pmatrix} 0 & 0 \\ 0 & 0 \end{pmatrix}$.

**Theorem 4.15:** Similar matrices have the same characteristic polynomial.

**Definition (Eigenvalue and Eigenvector):** Let $\mathbf{A}$ be an $n$-square matrix over a field $K$. A scalar $\lambda \in K$ is called an eigenvalue of $\mathbf{A}$ if there exists a nonzero vector $|\mathbf{u}\rangle \in K^n$ for which

$$\mathbf{A}|\mathbf{u}\rangle = \lambda|\mathbf{u}\rangle. \tag{3.62}$$

Every vector satisfying this relation is then called an eigenvector of $\mathbf{A}$ belonging to the eigenvalue $\lambda$.

***Example 3.12***

For the matrix $\mathbf{A} = \begin{pmatrix} 2 & 0 \\ 0 & -1 \end{pmatrix}$, verify that $|\mathbf{x}_1\rangle = \begin{pmatrix} 1 \\ 0 \end{pmatrix}$, $|\mathbf{x}_2\rangle = \begin{pmatrix} 0 \\ 1 \end{pmatrix}$ are the eigenvector of $\mathbf{A}$ corresponding to eigenvalues $\lambda_1 = 2, \lambda_2 = -1$, respectively.

***Solution***

$$\mathbf{A}|\mathbf{x}_1\rangle = \begin{pmatrix} 2 & 0 \\ 0 & -1 \end{pmatrix}\begin{pmatrix} 1 \\ 0 \end{pmatrix} = \begin{pmatrix} 2 \\ 0 \end{pmatrix} = 2\begin{pmatrix} 1 \\ 0 \end{pmatrix} = 2|\mathbf{x}_1\rangle,$$

$$\mathbf{A}|\mathbf{x}_2\rangle = \begin{pmatrix} 2 & 0 \\ 0 & -1 \end{pmatrix}\begin{pmatrix} 0 \\ 1 \end{pmatrix} = \begin{pmatrix} 0 \\ -1 \end{pmatrix} = -\begin{pmatrix} 0 \\ 1 \end{pmatrix} = -|\mathbf{x}_2\rangle.$$

**Definition (Eigenspace):** The set $E_\lambda$ of all eigenvectors belonging to $\lambda$ is a subspace of $K^n$, called the eigenspace of $\lambda$.

**Definition (Singular Matrix):** A square matrix that is not invertible is called singular or degenerate. A square matrix is singular if and only if its determinant is zero.

**Theorem 3.16:** Let $\mathbf{A}$ be an $n$-square matrix over a field $K$. Then the following are equivalent.

1- A scalar $\lambda$ is an eigenvalue of $\mathbf{A}$.
2- The matrix $\mathbf{M} = \lambda\mathbf{I}_n - \mathbf{A}$ is singular.
3- The scalar $\lambda$ is a root of the characteristic polynomial $\Delta_{\mathbf{A}}(\lambda)$ of $\mathbf{A}$.

**Proof:**

The scalar $\lambda$ is an eigenvalue of $\mathbf{A}$ if and only if there exists a nonzero vector $|\mathbf{v}\rangle$ such that

$$\mathbf{A}|\mathbf{v}\rangle = \lambda|\mathbf{v}\rangle$$
$$\lambda\mathbf{I}|\mathbf{v}\rangle - \mathbf{A}|\mathbf{v}\rangle = |\mathbf{0}\rangle$$
$$(\lambda\mathbf{I} - \mathbf{A})|\mathbf{v}\rangle = |\mathbf{0}\rangle$$

or $\lambda\mathbf{I} - \mathbf{A} = \mathbf{M}$ is singular. In such a case $\lambda$ is a root of the characteristic polynomial $\Delta_{\mathbf{A}}(\lambda) = |\lambda\mathbf{I} - \mathbf{A}|$.

■

**Remarks:**

1- An eigenvalue of $\mathbf{A}$ is a scalar $\lambda$ such that

$$|\lambda\mathbf{I}_n - \mathbf{A}| = 0. \tag{3.63}$$

2- The eigenvectors of $\mathbf{A}$ corresponding to $\lambda$ are the nonzero solutions of

$$(\lambda\mathbf{I}_n - \mathbf{A})|\mathbf{x}\rangle = \mathbf{0}. \tag{3.64}$$

***Example 3.13***

Find the eigenvalues and corresponding eigenvectors of $\mathbf{A} = \begin{pmatrix} 2 & -12 \\ 1 & -5 \end{pmatrix}$.

***Solution***

The characteristic polynomial of $\mathbf{A}$ is

$$|\lambda\mathbf{I} - \mathbf{A}| = \begin{vmatrix} \lambda - 2 & 12 \\ -1 & \lambda + 5 \end{vmatrix} = (\lambda - 2)(\lambda + 5) + 12 = (\lambda + 1)(\lambda + 2).$$

So, the characteristic equation is $(\lambda + 1)(\lambda + 2) = 0$, which gives $\lambda_1 = -1, \lambda_2 = -2$ as eigenvalues of $\mathbf{A}$.

For $\lambda_1 = -1$

$$-\mathbf{I} - \mathbf{A} = \begin{pmatrix} -1 - 2 & 12 \\ -1 & -1 + 5 \end{pmatrix} = \begin{pmatrix} -3 & 12 \\ -1 & 4 \end{pmatrix}.$$

Let $\begin{pmatrix} -3 & 12 \\ -1 & 4 \end{pmatrix}\begin{pmatrix} x_1 \\ x_2 \end{pmatrix} = \begin{pmatrix} 0 \\ 0 \end{pmatrix}$. We get

$$\begin{pmatrix} x_1 \\ x_2 \end{pmatrix} = t\begin{pmatrix} 4 \\ 1 \end{pmatrix}, t \neq 0.$$

For $\lambda_1 = -2$

$$-2\mathbf{I} - \mathbf{A} = \begin{pmatrix} -2 - 2 & 12 \\ -1 & -2 + 5 \end{pmatrix} = \begin{pmatrix} -4 & 12 \\ -1 & 3 \end{pmatrix}.$$

Let $\begin{pmatrix} -4 & 12 \\ -1 & 3 \end{pmatrix}\begin{pmatrix} x_1 \\ x_2 \end{pmatrix} = \begin{pmatrix} 0 \\ 0 \end{pmatrix}$. We get

$$\begin{pmatrix} x_1 \\ x_2 \end{pmatrix} = t\begin{pmatrix} 3 \\ 1 \end{pmatrix}, t \neq 0.$$

**Definition (Diagonalizable Matrix):** A matrix $\mathbf{A}$ is said to be diagonalizable if there exists a non-singular matrix $\mathbf{P}$ such that $\mathbf{D} = \mathbf{P}^{-1}\mathbf{AP}$ is a diagonal matrix, i.e., if $\mathbf{A}$ is similar to a diagonal matrix $\mathbf{D}$.

For example, the matrix $\mathbf{A} = \begin{pmatrix} 1 & 3 & 0 \\ 3 & 1 & 0 \\ 0 & 0 & -2 \end{pmatrix}$ is diagonalizable because $\mathbf{P} = \begin{pmatrix} 1 & 1 & 0 \\ 1 & -1 & 0 \\ 0 & 0 & 1 \end{pmatrix}$ has the property

$$\mathbf{P}^{-1}\mathbf{AP} = \begin{pmatrix} 4 & 0 & 0 \\ 0 & -2 & 0 \\ 0 & 0 & -2 \end{pmatrix}.$$

**Theorem 3.17:** An $n$-square matrix $\mathbf{A}$ is similar to a diagonal matrix $\mathbf{D}$ if and only if $\mathbf{A}$ has $n$ linearly independent eigenvectors. In this case, the diagonal elements of $\mathbf{D}$ are the corresponding eigenvalues and $\mathbf{D} = \mathbf{P}^{-1}\mathbf{AP}$ where $\mathbf{P}$ is the matrix whose columns are the eigenvectors.

***Example 3.14***

The matrix $\mathbf{A} = \begin{pmatrix} 1 & 3 & 0 \\ 3 & 1 & 0 \\ 0 & 0 & -2 \end{pmatrix}$ has the eigenvalues and corresponding eigenvectors listed below

$$\lambda_1 = 4, |\mathbf{p}_1\rangle = \begin{pmatrix} 1 \\ 1 \\ 0 \end{pmatrix}, \qquad \lambda_2 = -2, |\mathbf{p}_2\rangle = \begin{pmatrix} 1 \\ -1 \\ 0 \end{pmatrix}, \qquad \lambda_3 = -2, |\mathbf{p}_3\rangle = \begin{pmatrix} 0 \\ 0 \\ 1 \end{pmatrix}.$$

The eigenvectors $|\mathbf{p}_1\rangle$, $|\mathbf{p}_2\rangle$ and $|\mathbf{p}_3\rangle$ are linearly independent. The matrix $P$ whose columns correspond to these eigenvectors, is

$$\mathbf{P} = \begin{pmatrix} 1 & 1 & 0 \\ 1 & -1 & 0 \\ 0 & 0 & 1 \end{pmatrix}.$$

Hence, we have

$$\mathbf{P}^{-1}\mathbf{AP} = \begin{pmatrix} 4 & 0 & 0 \\ 0 & -2 & 0 \\ 0 & 0 & -2 \end{pmatrix}.$$

**Theorem 3.18:** Let $|\mathbf{u}_1\rangle, |\mathbf{u}_2\rangle, \ldots, |\mathbf{u}_n\rangle$ be nonzero eigenvectors of a matrix **A** belonging to distinct eigenvalues $\lambda_1, \lambda_2, \ldots \lambda_n$. Then $|\mathbf{u}_1\rangle, |\mathbf{u}_2\rangle, \ldots, |\mathbf{u}_n\rangle$ are linearly independent.

**Theorem 3.19:** Suppose the characteristic polynomial $\Delta_{\mathbf{A}}(\lambda)$ of an $n$-square matrix **A** is a product of $n$ distinct factors, say, $\Delta_{\mathbf{A}}(\lambda) = (\lambda - \lambda_1)(\lambda - \lambda_2) \ldots (\lambda - \lambda_n)$. Then **A** is similar to a diagonal matrix whose diagonal elements are the $\lambda_i$.

## 3.6 The Four Fundamental Subspaces of a Matrix

### Row and Column Spaces of a Matrix

In this section, we will study the vector space spanned by the row vectors (or column vectors) of a matrix [7]. To begin, let us consider some terminology. For an $m \times n$ matrix **A**, the $n$-tuples corresponding to the rows of **A** are called the row vectors of **A**. Similarly, the columns of **A** are called the column vectors of **A**.

$$\mathbf{A} = \begin{pmatrix} a_{11} & a_{12} & \cdot\cdot & a_{1n} \\ a_{21} & a_{22} & \cdot\cdot & a_{2n} \\ \vdots & \vdots & & \vdots \\ a_{m1} & a_{m2} & \cdots & a_{mn} \end{pmatrix} \qquad \begin{matrix} \text{Row Vectors of } \mathbf{A} \\ (a_{11}, a_{12}, \ldots, a_{1n}) \\ (a_{21}, a_{22}, \ldots, a_{2n}) \\ \vdots \\ (a_{m1}, a_{m2}, \ldots, a_{mn}) \end{matrix} \qquad \begin{matrix} \text{Column Vectors of } \mathbf{A} \\ \begin{pmatrix} a_{11} \\ a_{21} \\ \vdots \\ a_{m1} \end{pmatrix} \begin{pmatrix} a_{12} \\ a_{22} \\ \vdots \\ a_{m2} \end{pmatrix} \cdots \begin{pmatrix} a_{1n} \\ a_{2n} \\ \vdots \\ a_{mn} \end{pmatrix}. \end{matrix} \tag{3.65}$$

**Definition (Row Space and Column Space of a Matrix):** Let **A** be an $m \times n$ matrix.

1. The row space of **A** is the subspace of $\mathbb{R}^n$ spanned by the row vectors of **A**.
2. The column space of **A** is the subspace of $\mathbb{R}^m$ spanned by the column vectors of **A**.

Recall that two matrices are row-equivalent if one can be obtained from the other by elementary row operations. The next theorem shows that the row-equivalent matrices have the same row space.

**Theorem 3.20:** If an $m \times n$ matrix **A** is row-equivalent to an $m \times n$ matrix **B**, then the row space of **A** is equal to the row space of **B**.

**Proof:**

Because the rows of **B** can be obtained from the rows of **A** by elementary row operations (scalar multiplication and addition), it follows that the row vectors of **B** can be written as linear combinations of the row vectors of **A**. The row vectors of **B** lie in the row space of **A**, and the subspace spanned by the row vectors of **B** is contained in the row space of **A**. But it is also true that the rows of **A** can be obtained from the rows of **B** by elementary row operations. So, you can conclude that the two-row spaces are subspaces of each other, making them equal.

■

**Remark:**

Theorem 3.20 says that the row space of a matrix is not changed by elementary row operations. Elementary row operations can, however, change the column space.

**Theorem 3.21:** If a matrix **A** is row-equivalent to a matrix **B** in row-echelon form, then the nonzero row vectors of **B** form a basis for the row space of **A**.

**Proof:**

If a matrix **B** is in row-echelon form, then its nonzero row vectors form a linearly independent set. Consequently, they form a basis for the row space of **B**, and by Theorem 3.20, they also form a basis for the row space of **A**.

■

***Example 3.15***

Find a basis for the row space of the matrix

$$\mathbf{A} = \begin{pmatrix} 1 & 3 & 1 & 3 \\ 0 & 1 & 1 & 0 \\ -3 & 0 & 6 & -1 \\ 3 & 4 & -2 & 1 \\ 2 & 0 & -4 & -2 \end{pmatrix}.$$

***Solution***

Using elementary row operations, rewrite $\mathbf{A}$ in row-echelon form as follows.

$$\mathbf{B} = \begin{pmatrix} 1 & 3 & 1 & 3 \\ 0 & 1 & 1 & 0 \\ 0 & 0 & 0 & 1 \\ 0 & 0 & 0 & 0 \\ 0 & 0 & 0 & 0 \end{pmatrix} \begin{matrix} \mathbf{w}_1 \\ \mathbf{w}_2 \\ \mathbf{w}_3. \\ \\ \\ \end{matrix}$$

By Theorem 3.21, you can conclude that the nonzero row vectors of $\mathbf{B}$,

$$\langle \mathbf{w}_1 | = (1,3,1,3), \qquad \langle \mathbf{w}_2 | = (0,1,1,0), \qquad \langle \mathbf{w}_3 | = (0,0,0,1),$$

form a basis for the row space of $\mathbf{A}$.

To find a basis for the column space of a matrix $\mathbf{A}$, we can use the fact that the column space of $\mathbf{A}$ is equal to the row space of $\mathbf{A}^T$ and apply the technique of Example 3.15 to the matrix $\mathbf{A}^T$.

***Example 3.16***

Find a basis for the column space of matrix $\mathbf{A}$ from Example 3.15.

$$\mathbf{A} = \begin{pmatrix} 1 & 3 & 1 & 3 \\ 0 & 1 & 1 & 0 \\ -3 & 0 & 6 & -1 \\ 3 & 4 & -2 & 1 \\ 2 & 0 & -4 & -2 \end{pmatrix}.$$

***Solution 1***

Take the transpose of $\mathbf{A}$ and use elementary row operations to write $\mathbf{A}^T$ in row-echelon form.

$$\mathbf{A}^T = \begin{pmatrix} 1 & 0 & -3 & 3 & 2 \\ 3 & 1 & 0 & 4 & 0 \\ 1 & 1 & 6 & -2 & -4 \\ 3 & 0 & -1 & 1 & -2 \end{pmatrix} \Rightarrow \begin{pmatrix} 1 & 0 & -3 & 3 & 2 \\ 0 & 1 & 9 & -5 & -6 \\ 0 & 0 & 1 & -1 & -1 \\ 0 & 0 & 0 & 0 & 0 \end{pmatrix} \begin{matrix} \mathbf{w}_1 \\ \mathbf{w}_2 \\ \mathbf{w}_3. \\ \\ \end{matrix}$$

So, $\langle \mathbf{w}_1 | = (1,0,-3,3,2)$, $\langle \mathbf{w}_2 | = (0,1,9,-5,-6)$, and $\langle \mathbf{w}_3 | = (0,0,1,-1,-1)$ form a basis for the row space of $\mathbf{A}^T$. This is equivalent to saying that the column vectors

$$\begin{pmatrix} 1 \\ 0 \\ -3 \\ 3 \\ 2 \end{pmatrix}, \begin{pmatrix} 0 \\ 1 \\ 9 \\ -5 \\ -6 \end{pmatrix}, \text{ and } \begin{pmatrix} 0 \\ 0 \\ 1 \\ -1 \\ -1 \end{pmatrix},$$

form a basis for the column space of $\mathbf{A}$.

**Theorem 3.22:** If $\mathbf{A}$ is an $m \times n$ matrix, then the row space and column space of $\mathbf{A}$ have the same dimension.

**Proof:**

Let $\langle \mathbf{v}_1 |, \langle \mathbf{v}_2 |, \ldots,$ and $\langle \mathbf{v}_m |$ be the row vectors and $|\mathbf{u}_1\rangle, |\mathbf{u}_2\rangle, \ldots,$ and $|\mathbf{u}_n\rangle$ be the column vectors of the matrix

$$\mathbf{A} = \begin{pmatrix} a_{11} & a_{12} & \cdots & a_{1n} \\ a_{21} & a_{22} & \cdots & a_{2n} \\ \vdots & \vdots & & \vdots \\ a_{m1} & a_{m2} & \cdots & a_{mn} \end{pmatrix}.$$

Suppose the row space of **A** has dimension $r$ and basis $S = \{\langle \mathbf{b}_1|, \langle \mathbf{b}_2|, \dots, \langle \mathbf{b}_r|\}$, where $\langle \mathbf{b}_\mathrm{i}| = (b_{i1}, b_{i2}, \dots, b_{in})$. Using this basis, you can write the row vectors of **A** as

$$\begin{aligned}\langle \mathbf{v}_1| &= c_{11}\langle \mathbf{b}_1| + c_{12}\langle \mathbf{b}_2| + \cdots + c_{1r}\langle \mathbf{b}_r|,\\ \langle \mathbf{v}_2| &= c_{21}\langle \mathbf{b}_1| + c_{22}\langle \mathbf{b}_2| + \cdots + c_{2r}\langle \mathbf{b}_r|,\\ &\vdots\\ \langle \mathbf{v}_\mathrm{m}| &= c_{m1}\langle \mathbf{b}_1| + c_{m2}\langle \mathbf{b}_2| + \cdots + c_{mr}\langle \mathbf{b}_r|.\end{aligned}$$

Rewrite this system of vector equations as follows

$$\begin{aligned}(a_{11}, a_{12}, \cdots, a_{1n}) &= c_{11}(b_{11}, b_{12}, \cdots, b_{1n}) + c_{12}(b_{21}, b_{22}, \cdots, b_{2n}) + \cdots + c_{1r}(b_{r1}, b_{r2}, \cdots, b_{rn}),\\ (a_{21}, a_{22}, \cdots, a_{2n}) &= c_{21}(b_{11}, b_{12}, \cdots, b_{1n}) + c_{22}(b_{21}, b_{22}, \cdots, b_{2n}) + \cdots + c_{2r}(b_{r1}, b_{r2}, \cdots, b_{rn}),\\ &\vdots\\ (a_{m1}, a_{m2}, \cdots, a_{mn}) &= c_{m1}(b_{11}, b_{12}, \cdots, b_{1n}) + c_{m2}(b_{21}, b_{22}, \cdots, b_{2n}) + \cdots + c_{mr}(b_{r1}, b_{r2}, \cdots, b_{rn}).\end{aligned}$$

Now, take only entries corresponding to the first column of matrix **A** to obtain the system of scalar equations shown below

$$\begin{aligned}a_{11} &= c_{11}b_{11} + c_{12}b_{21} + c_{13}b_{31} + \cdots + c_{1r}b_{r1},\\ a_{21} &= c_{21}b_{11} + c_{22}b_{21} + c_{23}b_{31} + \cdots + c_{2r}b_{r1},\\ a_{31} &= c_{31}b_{11} + c_{32}b_{21} + c_{33}b_{31} + \cdots + c_{3r}b_{r1},\\ &\vdots\\ a_{m1} &= c_{m1}b_{11} + c_{m2}b_{21} + c_{m3}b_{31} + \cdots + c_{mr}b_{r1}.\end{aligned}$$

Similarly, for the entries of the $j$th column, you can obtain the system below

$$\begin{aligned}a_{1j} &= c_{11}b_{1j} + c_{12}b_{2j} + c_{13}b_{3j} + \cdots + c_{1r}b_{rj},\\ a_{2j} &= c_{21}b_{1j} + c_{22}b_{2j} + c_{23}b_{3j} + \cdots + c_{2r}b_{rj},\\ a_{3j} &= c_{31}b_{1j} + c_{32}b_{2j} + c_{33}b_{3j} + \cdots + c_{3r}b_{rj},\\ &\vdots\\ a_{mj} &= c_{m1}b_{1j} + c_{m2}b_{2j} + c_{m3}b_{3j} + \cdots + c_{mr}b_{rj}.\end{aligned}$$

Now, let the vectors $|\mathbf{c}_i\rangle = (c_{1i}, c_{2i}, \cdots, c_{mi})^T$. Then the system for the $j^{\text{th}}$ column can be rewritten in a vector form as

$$|\mathbf{u}_j\rangle = b_{1j}|\mathbf{c}_1\rangle + b_{2j}|\mathbf{c}_2\rangle + \cdots + b_{rj}|\mathbf{c}_r\rangle.$$

Put all column vectors together to obtain

$$\begin{aligned}|\mathbf{u}_1\rangle &= (a_{11} \quad a_{21} \cdots a_{m1})^T = b_{11}|\mathbf{c}_1\rangle + b_{21}|\mathbf{c}_2\rangle + \cdots + b_{r\mathbf{c}}|\mathbf{c}_r\rangle,\\ |\mathbf{u}_2\rangle &= (a_{12} \quad a_{22} \cdots a_{m2})^T = b_{12}|\mathbf{c}_1\rangle + b_{22}|\mathbf{c}_2\rangle + \cdots + b_{r2}|\mathbf{c}_r\rangle,\\ &\vdots\\ |\mathbf{u}_n\rangle &= (a_{1n} \quad a_{2n} \cdots a_{mn})^T = b_{1n}|\mathbf{c}_1\rangle + b_{2n}|\mathbf{c}_2\rangle + \cdots + b_{rn}|\mathbf{c}_r\rangle.\end{aligned}$$

Because each column vector of **A** is a linear combination of $r$ vectors, the dimension of the column space of **A** is less than or equal to $r$ (the dimension of the row space of **A**). That is,

$$\dim(\text{column space of } \mathbf{A}) \le \dim(\text{row space of } \mathbf{A}).$$

Repeating this procedure for $\mathbf{A}^T$, we can conclude that the dimension of the column space of $\mathbf{A}^T$ is less than or equal to the dimension of the row space of $\mathbf{A}^T$. But this implies that the dimension of the row space of **A** is less than or equal to the dimension of the column space of **A**. That is,

$$\dim(\text{row space of } \mathbf{A}) \le \dim(\text{column space of } \mathbf{A}).$$

So, the two dimensions must be equal.

■

**Definition (Rank of a Matrix):** Let $\mathbf{A}$ be an $m \times n$ matrix. The dimension of the row (or column) space of a matrix $\mathbf{A}$ is called the rank of $\mathbf{A}$ and is denoted by $\mathrm{rank}(\mathbf{A})$.

***Example 3.17***

Find the rank of the matrix

$$\mathbf{A} = \begin{pmatrix} 1 & -2 & 0 & 1 \\ 2 & 1 & 5 & -3 \\ 0 & 1 & 3 & 5 \end{pmatrix}.$$

***Solution***

Convert to row-echelon form as follows

$$\mathbf{A} = \begin{pmatrix} 1 & -2 & 0 & 1 \\ 2 & 1 & 5 & -3 \\ 0 & 1 & 3 & 5 \end{pmatrix} \Rightarrow \mathbf{B} = \begin{pmatrix} 1 & -2 & 0 & 1 \\ 0 & 1 & 1 & -1 \\ 0 & 0 & 1 & 3 \end{pmatrix}.$$

Because $\mathbf{B}$ has three nonzero rows, the rank of $\mathbf{A}$ is 3.

### The Null-space of a Matrix

Remember, given a linear map $\mathbf{T}: V \to U$ between two vector spaces $V$ and $U$, the kernel of $\mathbf{T}$ (null space) is the vector subspace of all elements $|\mathbf{v}\rangle$ of $V$ such that $\mathbf{T}|\mathbf{v}\rangle = |\mathbf{0}\rangle$. null $\mathbf{T} = \{|\mathbf{v}\rangle \in V: \mathbf{T}|\mathbf{v}\rangle = |\mathbf{0}\rangle\}$.

**Theorem 3.23:** Consider a linear map represented as a $m \times n$ matrix $\mathbf{A}$ with coefficients in $\mathbb{R}$, that is operating on column vectors $|\mathbf{x}\rangle$ with $n$ components over $\mathbb{R}$. Then, the set of all solutions of the homogeneous system of linear equations

$$\mathbf{A}|\mathbf{x}\rangle = |\mathbf{0}\rangle, \tag{3.66}$$

is a subspace of $\mathbb{R}^n$ called the null-space of $\mathbf{A}$ and is denoted by $N(\mathbf{A})$. So,

$$N(\mathbf{A}) = \{|\mathbf{x}\rangle \in \mathbb{R}^n: \mathbf{A}|\mathbf{x}\rangle = |\mathbf{0}\rangle\}. \tag{3.67}$$

The dimension of the null-space of $\mathbf{A}$ is called the nullity of $\mathbf{A}$.

**Proof:**

Because $\mathbf{A}$ is an $m \times n$ matrix, the vector $|\mathbf{x}\rangle$ has size $n \times 1$. So, the set of all solutions of the system is a subset of $\mathbb{R}^n$. This set is clearly non-empty, because $\mathbf{A}|\mathbf{0}\rangle = |\mathbf{0}\rangle$. You can verify that it is a subspace by showing that it is closed under the operations of addition and scalar multiplication. Let $|\mathbf{x}_1\rangle$ and $|\mathbf{x}_2\rangle$ be two solution vectors of the system $\mathbf{A}|\mathbf{x}\rangle = |\mathbf{0}\rangle$, and let $c$ be a scalar. Because $\mathbf{A}|\mathbf{x}_1\rangle = |\mathbf{0}\rangle$ and $\mathbf{A}|\mathbf{x}_2\rangle = |\mathbf{0}\rangle$, we have

$$\mathbf{A}|\mathbf{x}_1 + \mathbf{x}_2\rangle = \mathbf{A}|\mathbf{x}_1\rangle + \mathbf{A}|\mathbf{x}_2\rangle = |\mathbf{0}\rangle + |\mathbf{0}\rangle = |\mathbf{0}\rangle,$$

and

$$\mathbf{A}|c\mathbf{x}_1\rangle = c(\mathbf{A}|\mathbf{x}_1\rangle) = c|\mathbf{0}\rangle = |\mathbf{0}\rangle.$$

So, both $|\mathbf{x}_1\rangle + |\mathbf{x}_2\rangle$ and $|c\mathbf{x}_1\rangle$ are solutions of $\mathbf{A}|\mathbf{x}\rangle = |\mathbf{0}\rangle$. Hence, the set of all solutions forms a subspace of $\mathbb{R}^n$.

■

**Theorem 3.24:** If $\mathbf{A}$ is an $m \times n$ matrix of rank $r$, then the dimension of the solution space of $\mathbf{A}|\mathbf{x}\rangle = |\mathbf{0}\rangle$ is $n - r$. That is,

$$n = \mathrm{rank}(\mathbf{A}) + \mathrm{nullity}(\mathbf{A}). \tag{3.68}$$

**Proof:**

Because $\mathbf{A}$ has rank $r$, it is row-equivalent to a reduced row-echelon matrix $\mathbf{B}$ with $r$ nonzero rows. No generality is lost by assuming that the upper left corner of $\mathbf{B}$ has the form of the $r \times r$ identity matrix $\mathbf{I}_r$. Moreover, because the zero rows of $\mathbf{B}$ contribute nothing to the solution, you can discard them to form the $r \times n$ matrix $\mathbf{B}'$, where

$$\mathbf{B}' = [\mathbf{I}_r \mid \mathbf{C}].$$

The matrix $\mathbf{C}$ has $n - r$ columns corresponding to the variables $x_{r+1}, x_{r+2}, \dots, x_n$. So, the solution space of $\mathbf{A}|\mathbf{x}\rangle = |\mathbf{0}\rangle$ can be represented by the system

$$
\begin{aligned}
x_1 + \quad & c_{11}x_{r+1} + c_{12}x_{r+2} + \cdots + c_{1,n-r}x_n = 0 \\
x_2 + \quad & c_{21}x_{r+1} + c_{22}x_{r+2} + \cdots + c_{2,n-r}x_n = 0 \\
\vdots \quad & \quad \vdots \qquad\qquad \vdots \qquad\qquad \vdots \qquad \vdots \\
x_r + & c_{r1}x_{r+1} + c_{r2}x_{r+2} + \cdots + c_{r,n-r}x_n = 0.
\end{aligned}
$$

Solving for the first $r$ variables in terms of the last $n-r$ variables produces $n-r$ vectors in the basis of the solution space. Consequently, the solution space has dimension $n-r$.

■

**Notes:**

- The product $\mathbf{A}|\mathbf{x}\rangle$ can be written in terms of the inner product of vectors as follows:
$$\mathbf{A}|\mathbf{x}\rangle = \begin{pmatrix} \langle \mathbf{a}_1|\mathbf{x}\rangle \\ \vdots \\ \langle \mathbf{a}_m|\mathbf{x}\rangle \end{pmatrix}. \tag{3.69}$$
Here, $\langle \mathbf{a}_1|, \ldots, \langle \mathbf{a}_m|$ denote the rows of the matrix $\mathbf{A}$. It follows that $|\mathbf{x}\rangle$ is in the kernel of $\mathbf{A}$, if and only if $|\mathbf{x}\rangle$ is orthogonal to each of the row vectors of $\mathbf{A}$.
- The row space of a matrix $\mathbf{A}$ is the span of the row vectors of $\mathbf{A}$. By the above reasoning, the kernel of $\mathbf{A}$ is the orthogonal complement to the row space. That is, a vector $|\mathbf{x}\rangle$ lies in the kernel of $\mathbf{A}$, if and only if it is perpendicular to every vector in the row space of $\mathbf{A}$.
- The left null space, or cokernel, of a matrix $\mathbf{A}$ consists of all column vectors $|\mathbf{x}\rangle$ such that $\langle \mathbf{x}|\mathbf{A} = \langle \mathbf{0}|$ or $\mathbf{A}^T|\mathbf{x}\rangle = |\mathbf{0}\rangle$. The left null space of $\mathbf{A}$ is the same as the kernel of $\mathbf{A}^T$. The left null space of $\mathbf{A}$ is the orthogonal complement to the column space of $\mathbf{A}$, and is dual to the cokernel of the associated linear transformation.
- The kernel (null space), the row space, the column space, and the left null space of $\mathbf{A}$ are the four fundamental subspaces associated to the matrix $\mathbf{A}$. (See Figure 3.2)

| | |
|---|---|
| The column space | contains all combinations of the columns of $\mathbf{A}$. |
| The row space | contains all combinations of the columns of $\mathbf{A}^T$. |
| The null space | contains all solutions $|\mathbf{x}\rangle$ to $\mathbf{A}|\mathbf{x}\rangle = |\mathbf{0}\rangle$. |
| The left null space | contains all solutions $|\mathbf{y}\rangle$ to $\mathbf{A}^T|\mathbf{y}\rangle = |\mathbf{0}\rangle$. |

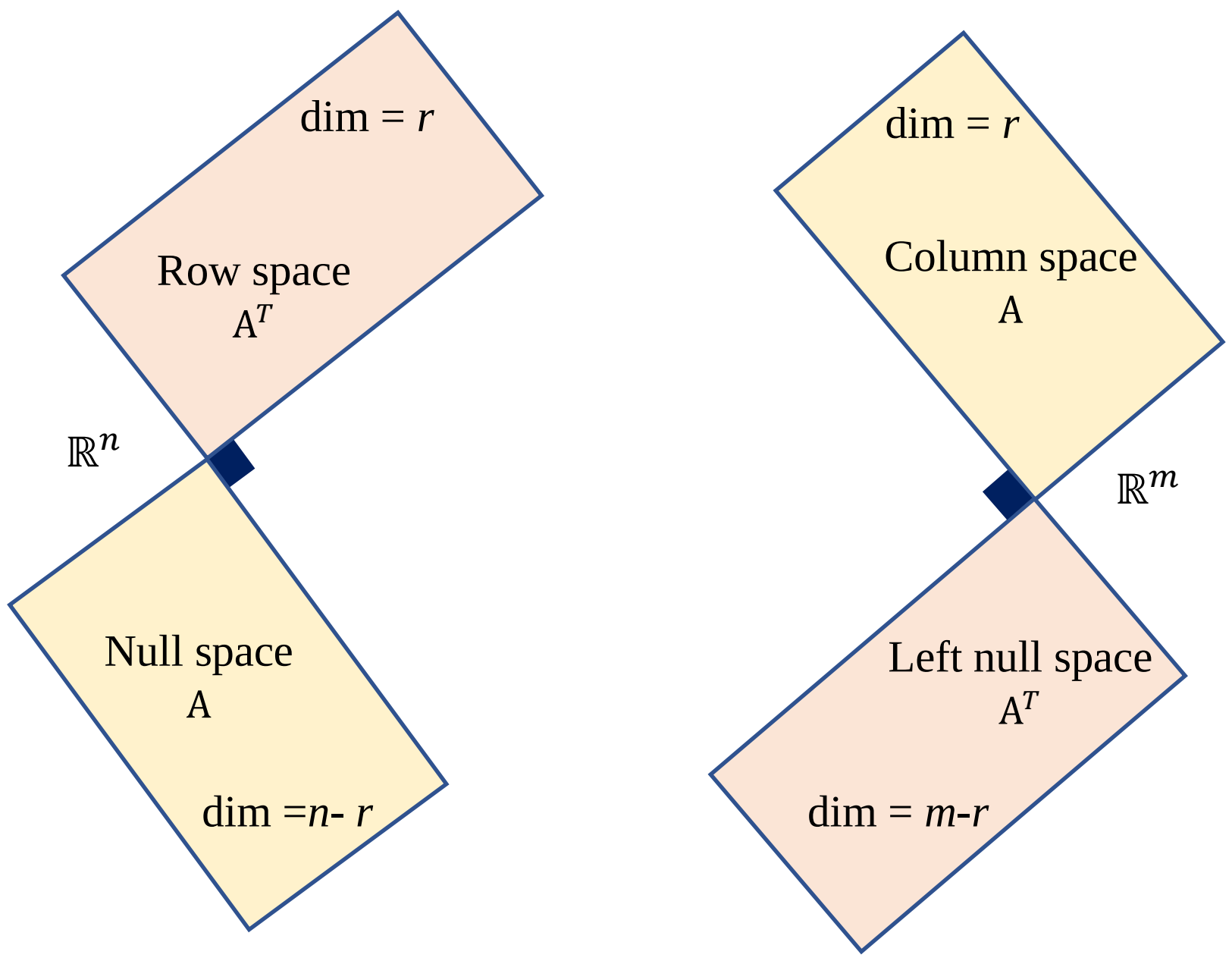

**Figure 3.2.** The four fundamental subspaces of a matrix.

## 3.7 Symmetric Matrices and Orthogonal Diagonalization

**Theorem 3.25:** If **A** is an $n \times n$ symmetric matrix, then the following properties are true [1-3].
1. **A** is diagonalizable.
2. All eigenvalues of **A** are real.
3. If $\lambda$ is an eigenvalue of **A** with multiplicity $k$, then $\lambda$ has $k$ linearly independent eigenvectors. That is, the eigenspace of $\lambda$ has dimension $k$.

Theorem 3.25 is called the real spectral theorem, and the set of eigenvalues of **A** is called the spectrum of **A**.

***Example 3.18***

Prove that a symmetric matrix $\mathbf{A} = \begin{pmatrix} a & c \\ c & b \end{pmatrix}$ is diagonalizable.

***Solution***

The characteristic polynomial of **A** is

$$|\lambda\mathbf{I} - \mathbf{A}| = \begin{vmatrix} \lambda - a & -c \\ -c & \lambda - b \end{vmatrix} = \lambda^2 - (a+b)\lambda + ab - c^2.$$

As a quadratic in $\lambda$, this polynomial has a discriminant of

$$\begin{aligned}(a+b)^2 - 4(ab - c^2) &= a^2 + 2ab + b^2 - 4ab + 4c^2 \\ &= a^2 - 2ab + b^2 + 4c^2 \\ &= (a-b)^2 + 4c^2.\end{aligned}$$

Because this discriminant is the sum of two squares, it must be either zero or positive. If $(a-b)^2 + 4c^2 = 0$, then $a = b$ and $c = 0$, which implies that is already diagonal.
That is,

$$\mathbf{A} = \begin{pmatrix} a & 0 \\ 0 & a \end{pmatrix}.$$

On the other hand, if $(a-b)^2 + 4c^2 > 0$, then by the quadratic formula, the characteristic polynomial of **A** has two distinct real roots, which implies that **A** has two distinct real eigenvalues. So, **A** is diagonalizable in this case.

**Theorem 3.26:** An $n \times n$ matrix **P** is orthogonal if and only if its column vectors form an orthonormal set.

**Proof:**

Suppose the column vectors of **P** form an orthonormal set:

$$\mathbf{P} = \begin{pmatrix} P_{11} & P_{12} & \dots & P_{1n} \\ P_{21} & P_{22} & \dots & P_{2n} \\ \vdots & \vdots & \ddots & \vdots \\ P_{n1} & P_{n2} & \dots & P_{nn} \end{pmatrix} = (|\mathbf{p}_1\rangle \vdots \ |\mathbf{p}_2\rangle \vdots \ \dots \ \vdots |\mathbf{p}_n\rangle).$$

Then the product $\mathbf{P}^T\mathbf{P}$ has the form

$$\mathbf{P}^T\mathbf{P} = \begin{pmatrix} P_{11} & P_{21} & \dots & P_{n1} \\ P_{12} & P_{22} & \dots & P_{n2} \\ \vdots & \vdots & \ddots & \vdots \\ P_{1n} & P_{2n} & \dots & P_{nn} \end{pmatrix} \begin{pmatrix} P_{11} & P_{12} & \dots & P_{1n} \\ P_{21} & P_{22} & \dots & P_{2n} \\ \vdots & \vdots & \ddots & \vdots \\ P_{n1} & P_{n2} & \dots & P_{nn} \end{pmatrix} = \begin{pmatrix} \langle\mathbf{p}_1|\mathbf{p}_1\rangle & \langle\mathbf{p}_1|\mathbf{p}_2\rangle & \dots & \langle\mathbf{p}_1|\mathbf{p}_n\rangle \\ \langle\mathbf{p}_2|\mathbf{p}_1\rangle & \langle\mathbf{p}_2|\mathbf{p}_2\rangle & \dots & \langle\mathbf{p}_2|\mathbf{p}_n\rangle \\ \vdots & \vdots & \ddots & \vdots \\ \langle\mathbf{p}_n|\mathbf{p}_1\rangle & \langle\mathbf{p}_n|\mathbf{p}_2\rangle & \dots & \langle\mathbf{p}_n|\mathbf{p}_n\rangle \end{pmatrix}.$$

Because the set $\{|\mathbf{p}_1\rangle, |\mathbf{p}_2\rangle, \dots, |\mathbf{p}_n\rangle\}$ is orthonormal; we have

$$\langle\mathbf{p}_i|\mathbf{p}_j\rangle = \delta_{ij}.$$

So, the matrix composed of dot products has the form

$$\mathbf{P}^T\mathbf{P} = \begin{pmatrix} 1 & 0 & \dots & 0 \\ 0 & 1 & \dots & 0 \\ \vdots & \vdots & \ddots & \vdots \\ 0 & 0 & \dots & 1 \end{pmatrix} = \mathbf{I}.$$

This implies that $\mathbf{P}^T = \mathbf{P}^{-1}$ and you can conclude that $\mathbf{P}$ is orthogonal.

Conversely, if $\mathbf{P}$ is orthogonal, you can reverse the steps above to verify that the column vectors of $\mathbf{P}$ form an orthonormal set.

■

**Theorem 3.27:** Let $\mathbf{A}$ be $n \times n$ matrix. Then $\mathbf{A}$ is orthogonally diagonalizable and has real eigenvalues if and only if $\mathbf{A}$ is symmetric.

**Proof:**

The proof of the theorem in one direction is fairly straightforward. That is, if you assume $\mathbf{A}$ is orthogonally diagonalizable, then there exists an orthogonal matrix $\mathbf{P}$ such that $\mathbf{D} = \mathbf{P}^{-1}\mathbf{AP}$ is diagonal. Moreover, because $\mathbf{P}^T = \mathbf{P}^{-1}$, you have $\mathbf{A} = \mathbf{PDP}^{-1} = \mathbf{PDP}^T$ which implies that

$$\begin{aligned}\mathbf{A}^T &= (\mathbf{PDP}^T)^T \\ &= (\mathbf{P}^T)^T\mathbf{D}^T\mathbf{P}^T \\ &= \mathbf{PDP}^T \\ &= \mathbf{A}.\end{aligned}$$

So, $\mathbf{A}$ is symmetric. We leave the remaining of the proof as an exercise.

■

### Orthogonal Diagonalization of a Symmetric Matrix

- Let $\mathbf{A}$ be an $n \times n$ symmetric matrix.
- Find all eigenvalues of $\mathbf{A}$ and determine the multiplicity of each.
- For each eigenvalue of multiplicity 1, choose a unit eigenvector. (Choose any eigenvector and then normalize it.)
- For each eigenvalue of multiplicity $k \geq 2$, find a set of $k$ linearly independent eigenvectors.
- The composite of steps 2 and 3 produces an orthonormal set of eigenvectors. Use these eigenvectors to form the columns of $\mathbf{P}$. The matrix $\mathbf{P}^{-1}\mathbf{AP} = \mathbf{P}^T\mathbf{AP} = \mathbf{D}$ will be diagonal. (The main diagonal entries of $\mathbf{D}$ are the eigenvalues of $\mathbf{A}$.)

***Example 3.19***

Find an orthogonal matrix $\mathbf{P}$ that orthogonally diagonalizes

$$\mathbf{A} = \begin{pmatrix} -2 & 2 \\ 2 & 1 \end{pmatrix}.$$

***Solution***

1. The characteristic polynomial of $\mathbf{A}$ is

$$|\lambda\mathbf{I} - \mathbf{A}| = \begin{vmatrix} \lambda + 2 & -2 \\ -2 & \lambda - 1 \end{vmatrix} = (\lambda + 3)(\lambda - 2).$$

So, the eigenvalues are $\lambda_1 = -3$ and $\lambda_2 = 2$.

2. For each eigenvalue, find an eigenvector. The eigenvectors $(-2,1)^T$ and $(1,2)^T$ form an orthogonal basis for $\mathbb{R}^2$. Each of these eigenvectors is normalized to produce an orthonormal basis.

$$|\mathbf{p}_1\rangle = \frac{(-2,1)^T}{\sqrt{5}}, |\mathbf{p}_2\rangle = \frac{(1,2)^T}{\sqrt{5}}.$$

3. Because each eigenvalue has a multiplicity of 1, go directly to step 4.

4. Using $|\mathbf{p}_1\rangle$ and $|\mathbf{p}_2\rangle$ as column vectors, construct the matrix $\mathbf{P}$

$$\mathbf{P} = \frac{1}{\sqrt{5}}\begin{pmatrix} -2 & 1 \\ 1 & 2 \end{pmatrix}.$$

Verify that is correct by computing $\mathbf{P}^{-1}\mathbf{AP} = \mathbf{P}^T\mathbf{AP}$

$$\mathbf{P}^T\mathbf{AP} = \frac{1}{\sqrt{5}}\begin{pmatrix} -2 & 1 \\ 1 & 2 \end{pmatrix}\begin{pmatrix} -2 & 2 \\ 2 & 1 \end{pmatrix}\frac{1}{\sqrt{5}}\begin{pmatrix} -2 & 1 \\ 1 & 2 \end{pmatrix} = \begin{pmatrix} -3 & 0 \\ 0 & 2 \end{pmatrix}.$$

## 3.8 Singular Values and Singular Vectors in the Singular Value Decomposition (SVD)

The best matrices (real symmetric matrices) have real eigenvalues and orthogonal eigenvectors [7]. But for other matrices, the eigenvalues are complex, or the eigenvectors are not orthogonal. If $\mathbf{A}$ is not square then $\mathbf{A}|\mathbf{x}\rangle = \lambda|\mathbf{x}\rangle$ is impossible and eigenvectors fail (left side in $\mathbb{R}^m$, right side in $\mathbb{R}^n$ ). We need an idea that succeeds for every matrix. The singular value decomposition fills this gap in a perfect way. The key point is that we need two sets of singular vectors, the $|\mathbf{u}\rangle$'s and the $|\mathbf{v}\rangle$'s. For a real $m \times n$ matrix, the $n$ right singular vectors $|\mathbf{v}_1\rangle, \dots, |\mathbf{v}_n\rangle$ are orthogonal in $\mathbb{R}^n$. The $m$ left singular vectors $|\mathbf{u}_1\rangle, \dots, |\mathbf{u}_m\rangle$ are perpendicular to each other in $\mathbb{R}^m$. The connection between $n$ $|\mathbf{v}\rangle$'s and $m$ $|\mathbf{u}\rangle$ 's is not $\mathbf{A}|\mathbf{x}\rangle = \lambda|\mathbf{x}\rangle$. That is for eigenvectors. For singular vectors, each $\mathbf{A}|\mathbf{v}\rangle$ equals $\sigma|\mathbf{u}\rangle$:

$$\mathbf{A}|\mathbf{v}_1\rangle = \sigma_1|\mathbf{u}_1\rangle, \cdots, \mathbf{A}|\mathbf{v}_r\rangle = \sigma_r|\mathbf{u}_r\rangle, \qquad \mathbf{A}|\mathbf{v}_{r+1}\rangle = \mathbf{0} \cdots \mathbf{A}|\mathbf{v}_n\rangle = \mathbf{0}. \tag{3.70}$$

We have separated the first $r$ $|\mathbf{v}\rangle$'s and $|\mathbf{u}\rangle$'s from the rest. That number $r$ is the rank of $\mathbf{A}$. We will have $r$ positive singular values in descending order $\sigma_1 \geq \sigma_2 \geq \cdots \geq \sigma_r > 0$. The last $n - r$ $|\mathbf{v}\rangle$'s are in the null-space of $\mathbf{A}$, and the last $m - r$ $|\mathbf{u}\rangle$'s are in the null-space of $\mathbf{A}^T$.

The first step is to write (3.70) in matrix form. All of the right singular vectors $|\mathbf{v}_1\rangle$ to $|\mathbf{v}_n\rangle$ go in the columns of $\mathbf{V}$. The left singular vectors $|\mathbf{u}_1\rangle$ to $|\mathbf{u}_m\rangle$ go in the columns of $\mathbf{U}$. Those are square orthogonal matrices ($\mathbf{V}^\mathrm{T} = \mathbf{V}^{-1}$ and $\mathbf{U}^\mathrm{T} = \mathbf{U}^{-1}$) because their columns are orthogonal unit vectors. Then (3.70) becomes the full SVD, with square matrices $\mathbf{V}$ and $\mathbf{U}$:

$$\mathbf{AV} = \mathbf{U\Sigma} \Rightarrow \mathbf{A}\begin{pmatrix} |\mathbf{v}_1\rangle & \dots & |\mathbf{v}_r\rangle & \dots & |\mathbf{v}_n\rangle \end{pmatrix} = \begin{pmatrix} |\mathbf{u}_1\rangle & \dots & |\mathbf{u}_r\rangle & \dots & |\mathbf{u}_m\rangle \end{pmatrix}\begin{pmatrix} \sigma_1 & 0 & \cdots & 0 \\ 0 & \ddots & 0 & 0 \\ \vdots & 0 & \sigma_r & \vdots \\ 0 & 0 & \cdots & 0 \end{pmatrix}. \tag{3.71}$$

You see $\mathbf{A}|\mathbf{v}_k\rangle = \sigma_k|\mathbf{u}_k\rangle$ in the first $r$ columns above. That is an important part of the SVD. It shows the basis of $|\mathbf{v}\rangle$'s for the row space of $\mathbf{A}$ and then $|\mathbf{u}\rangle$'s for the column space. After the positive numbers $\sigma_1, \dots, \sigma_r$ on the main diagonal of $\mathbf{\Sigma}$, the rest of that matrix is all zero from the null-spaces of $\mathbf{A}$ and $\mathbf{A}^\mathrm{T}$.

For example:

$$\mathbf{AV} = \mathbf{U\Sigma} \Rightarrow \begin{pmatrix} 3 & 0 \\ 4 & 5 \end{pmatrix} \frac{1}{\sqrt{2}} \begin{pmatrix} 1 & -1 \\ 1 & 1 \end{pmatrix} = \frac{1}{\sqrt{10}} \begin{pmatrix} 1 & -3 \\ 3 & 1 \end{pmatrix} \begin{pmatrix} 3\sqrt{5} & 0 \\ 0 & \sqrt{5} \end{pmatrix}. \tag{3.72}$$

The matrix $\mathbf{A}$ is not symmetric, so $\mathbf{V}$ is different from $\mathbf{U}$. The rank is 2, so there are two singular values $\sigma_1 = 3\sqrt{5}$ and $\sigma_2 = \sqrt{5}$. Their product $3 \cdot 5 = 15$ is the determinant of $\mathbf{A}$ (in this respect, singular values are like eigenvalues). The columns of $\mathbf{V}$ are orthogonal, and the columns of $\mathbf{U}$ are orthogonal. Those columns are unit vectors after the divisions by $\sqrt{2}$ and $\sqrt{10}$, so $\mathbf{V}$ and $\mathbf{U}$ are orthogonal matrices: $\mathbf{V}^\mathrm{T} = \mathbf{V}^{-1}$ and $\mathbf{U}^\mathrm{T} = \mathbf{U}^{-1}$.

That orthogonality allows us to go from $\mathbf{AV} = \mathbf{U\Sigma}$ to the famous expression of the SVD: Multiply both sides of $\mathbf{AV} = \mathbf{U\Sigma}$ by $\mathbf{V}^{-1} = \mathbf{V}^\mathrm{T}$. The singular value decomposition of $\mathbf{A}$ is

$$\mathbf{A} = \mathbf{U\Sigma V}^T. \tag{3.73}$$

Then column-row multiplication of $\mathbf{U\Sigma}$ times $\mathbf{V}^\mathrm{T}$ separates $\mathbf{A}$ into $r$ pieces of rank 1, i.e., pieces of the SVD is

$$\mathbf{A} = \mathbf{U\Sigma V}^T = \sigma_1|\mathbf{u}_1\rangle\langle\mathbf{v}_1| + \cdots + \sigma_r|\mathbf{u}_r\rangle\langle\mathbf{v}_r|. \tag{3.74}$$

In the above example, to recover $\mathbf{A}$, add the pieces $\sigma_1|\mathbf{u}_1\rangle\langle\mathbf{v}_1| + \sigma_2|\mathbf{u}_2\rangle\langle\mathbf{v}_2|$:

$$3\sqrt{5}\frac{1}{\sqrt{10}}\begin{pmatrix} 1 \\ 3 \end{pmatrix}\frac{1}{\sqrt{2}}\begin{pmatrix} 1 & 1 \end{pmatrix} + \sqrt{5}\frac{1}{\sqrt{10}}\begin{pmatrix} -3 \\ 1 \end{pmatrix}\frac{1}{\sqrt{2}}\begin{pmatrix} -1 & 1 \end{pmatrix} = \frac{3}{2}\begin{pmatrix} 1 & 1 \\ 3 & 3 \end{pmatrix} + \frac{1}{2}\begin{pmatrix} 3 & -3 \\ -1 & 1 \end{pmatrix} = \begin{pmatrix} 3 & 0 \\ 4 & 5 \end{pmatrix}. \tag{3.75}$$

Notice that the right singular vectors $\begin{pmatrix} 1 \\ 1 \end{pmatrix}$ and $\begin{pmatrix} -1 \\ 1 \end{pmatrix}$ in $\mathbf{V}$ are transposed to rows $\langle\mathbf{v}_1|, \langle\mathbf{v}_2|$ of $\mathbf{V}^T$.

**Proof:**

Our goal is $\mathbf{A} = \mathbf{U\Sigma V}^T$. We want to identify the two sets of singular vectors, the $|\mathbf{u}\rangle$'s and the $|\mathbf{v}\rangle$'s. One way to find those vectors is to form the symmetric matrices $\mathbf{A}^T\mathbf{A}$ and $\mathbf{AA}^T$:

$$\mathbf{A}^T\mathbf{A} = (\mathbf{V\Sigma}^T\mathbf{U}^T)(\mathbf{U\Sigma V}^T) = \mathbf{V}(\mathbf{\Sigma}^T\mathbf{\Sigma})\mathbf{V}^T, \tag{3.76}$$
$$\mathbf{AA}^T = (\mathbf{U\Sigma V}^T)(\mathbf{V\Sigma}^T\mathbf{U}^T) = \mathbf{U}(\mathbf{\Sigma\Sigma}^T)\mathbf{U}^T. \tag{3.77}$$

Both (3.76) and (3.77) produced symmetric matrices. Usually, $\mathbf{A}^T\mathbf{A}$ and $\mathbf{AA}^T$ are different. Both right hand sides have the special form $\mathbf{Q\Lambda Q}^T$. Eigenvalues are in $\mathbf{\Lambda} = \mathbf{\Sigma}^T\mathbf{\Sigma}$ or $\mathbf{\Sigma\Sigma}^T$. Eigenvectors are in $\mathbf{Q} = \mathbf{V}$ or $\mathbf{Q} = \mathbf{U}$. So, we know from (3.76) and (3.77) how $\mathbf{V}$, $\mathbf{U}$, and $\mathbf{\Sigma}$ connect to the symmetric matrices $\mathbf{A}^T\mathbf{A}$ and $\mathbf{AA}^T$.

$\mathbf{V}$ contains orthonormal eigenvectors of $\mathbf{A}^T\mathbf{A}$,
$\mathbf{U}$ contains orthonormal eigenvectors of $\mathbf{AA}^T$,
$\sigma_1^2$ to $\sigma_r^2$ are the nonzero eigenvalues of both $\mathbf{A}^T\mathbf{A}$ and $\mathbf{AA}^T$.

The plan is to start with the $|\mathbf{v}\rangle$'s. Choose orthonormal eigenvectors $|\mathbf{v}_1\rangle, \dots, |\mathbf{v}_r\rangle$ of $\mathbf{A}^T\mathbf{A}$. Then choose $\sigma_k = \sqrt{\lambda_k}$. To determine the $|\mathbf{u}\rangle$'s we require $\mathbf{A}|\mathbf{v}\rangle = \sigma|\mathbf{u}\rangle$. Hence, we have

$$\mathbf{A}^T\mathbf{A}|\mathbf{v}_k\rangle = \sigma_k^2|\mathbf{v}_k\rangle, \qquad |\mathbf{u}_k\rangle = \frac{\mathbf{A}|\mathbf{v}_k\rangle}{\sigma_k} \text{ for } k = 1, \dots, r. \tag{3.78}$$

We can prove that the vectors$|\mathbf{u}\rangle$ are eigenvectors of $\mathbf{AA}^T$ as follows:

$$\mathbf{AA}^T|\mathbf{u}_k\rangle = \mathbf{AA}^T\left(\frac{\mathbf{A}|\mathbf{v}_k\rangle}{\sigma_k}\right) = \mathbf{A}\left(\frac{\mathbf{A}^T\mathbf{A}|\mathbf{v}_k\rangle}{\sigma_k}\right) = \mathbf{A}\frac{\sigma_k^2|\mathbf{v}_k\rangle}{\sigma_k} = \sigma_k^2|\mathbf{u}_k\rangle. \tag{3.79}$$

Notice that $(\mathbf{AA}^T)\mathbf{A} = \mathbf{A}(\mathbf{A}^T\mathbf{A})$ was the key to (3.79). Moreover, the vectors$|\mathbf{u}\rangle$ are also orthonormal:

$$|\mathbf{u}_j\rangle\langle\mathbf{u}_k| = \left(\frac{\mathbf{A}|\mathbf{v}_j\rangle}{\sigma_j}\right)^T\left(\frac{\mathbf{A}|\mathbf{v}_k\rangle}{\sigma_k}\right) = \frac{|\mathbf{v}_j\rangle\mathbf{A}^T\mathbf{A}|\mathbf{v}_k\rangle}{\sigma_j\sigma_k} = \frac{\sigma_k}{\sigma_j}|\mathbf{v}_j\rangle\langle\mathbf{v}_k| = \begin{cases}1 & \text{if } j = k\\ 0 & \text{if } j \neq k\end{cases}. \tag{3.80}$$

Finally, we have to choose the last $n - r$ vectors $|\mathbf{v}_{r+1}\rangle$ to $|\mathbf{v}_n\rangle$ and the last $m - r$ vectors $|\mathbf{u}_{r+1}\rangle$ to $|\mathbf{u}_m\rangle$. These $|\mathbf{v}\rangle$'s and $|\mathbf{u}\rangle$'s are in the null-spaces of $\mathbf{A}$ and $\mathbf{A}^T$. We can choose any orthonormal bases for those null-spaces. They will automatically be orthogonal to the first $|\mathbf{v}\rangle$'s in the row space of $\mathbf{A}$ and the first $|\mathbf{u}\rangle$'s in the column space. This is because the spaces are orthogonal (see Figure 3.2).

■

***Example 3.20***

Find the matrices $\mathbf{U}, \mathbf{\Sigma}, \mathbf{V}$ for $\mathbf{A} = \begin{pmatrix}3 & 0\\4 & 5\end{pmatrix}$.

***Solution***

With rank 2, this $\mathbf{A}$ has two positive singular values $\sigma_1$ and $\sigma_2$. We will see that $\sigma_1$ is larger than $\lambda_{\max} = 5$, and $\sigma_2$ is smaller than $\lambda_{\min} = 3$. Begin with $\mathbf{A}^T\mathbf{A}$ and $\mathbf{AA}^T$:

$$\mathbf{A}^T\mathbf{A} = \begin{pmatrix}25 & 20\\20 & 25\end{pmatrix} \qquad \mathbf{AA}^T = \begin{pmatrix}9 & 12\\12 & 41\end{pmatrix}.$$

Those have the same trace (50) and the same eigenvalues $\sigma_1^2 = 45$ and $\sigma_2^2 = 5$. The square roots are $\sigma_1 = \sqrt{45}$ and $\sigma_2 = \sqrt{5}$. Then $\sigma_1\sigma_2 = 15$, and this is the determinant of $\mathbf{A}$.

A key step is to find the eigenvectors of $\mathbf{A}^T\mathbf{A}$ (with eigenvalues 45 and 5):

$$\begin{pmatrix}25 & 20\\20 & 25\end{pmatrix}\begin{pmatrix}1\\1\end{pmatrix} = 45\begin{pmatrix}1\\1\end{pmatrix} \qquad \begin{pmatrix}25 & 20\\20 & 25\end{pmatrix}\begin{pmatrix}-1\\1\end{pmatrix} = 5\begin{pmatrix}-1\\1\end{pmatrix}.$$

Then $|\mathbf{v}_1\rangle$ and $|\mathbf{v}_2\rangle$ are those orthogonal eigenvectors rescaled to length 1. Divide by $\sqrt{2}$, the right singular vectors are $|\mathbf{v}_1\rangle = \frac{1}{\sqrt{2}}\begin{pmatrix}1\\1\end{pmatrix}$ and $|\mathbf{v}_2\rangle = \frac{1}{\sqrt{2}}\begin{pmatrix}-1\\1\end{pmatrix}$. Also, the left singular vectors are $|\mathbf{u}_i\rangle = \frac{\mathbf{A}|\mathbf{v}_i\rangle}{\sigma_i}$.

Now compute $\mathbf{A}|\mathbf{v}_1\rangle$ and $\mathbf{A}|\mathbf{v}_2\rangle$ which will be $\sigma_1|\mathbf{u}_1\rangle = \sqrt{45}|\mathbf{u}_1\rangle$ and $\sigma_2|\mathbf{u}_2\rangle = \sqrt{5}|\mathbf{u}_2\rangle$:

$$\mathbf{A}|\mathbf{v}_1\rangle = \frac{3}{\sqrt{2}}\begin{pmatrix}1\\3\end{pmatrix} = \sqrt{\frac{45}{10}}\begin{pmatrix}1\\3\end{pmatrix} = \sigma_1|\mathbf{u}_1\rangle,$$

$$\mathbf{A}|\mathbf{v}_2\rangle = \frac{1}{\sqrt{2}}\begin{pmatrix}-3\\1\end{pmatrix} = \sqrt{\frac{5}{10}}\begin{pmatrix}-3\\1\end{pmatrix} = \sigma_2|\mathbf{u}_2\rangle.$$

The division by $\sqrt{10}$ makes $|\mathbf{u}_1\rangle$ and $|\mathbf{u}_2\rangle$ orthonormal. Then $\sigma_1 = \sqrt{45}$ and $\sigma_2 = \sqrt{5}$ as expected. The singular value decomposition of $\mathbf{A}$ is $\mathbf{U}$ times $\mathbf{\Sigma}$ times $\mathbf{V}^T$.

$$\mathbf{U} = \frac{1}{\sqrt{10}}\begin{pmatrix}1 & -3\\3 & 1\end{pmatrix}, \qquad \mathbf{\Sigma} = \begin{pmatrix}\sqrt{45} & 0\\0 & \sqrt{5}\end{pmatrix}, \qquad \mathbf{V} = \frac{1}{\sqrt{2}}\begin{pmatrix}1 & -1\\1 & 1\end{pmatrix}.$$

$\mathbf{U}$ and $\mathbf{V}$ contain orthonormal bases for the column space and the row space of $\mathbf{A}$ (both spaces are just $\mathbb{R}^2$). The real achievement is that those two bases diagonalize $\mathbf{A}$: $\mathbf{AV}$ equals $\mathbf{U\Sigma}$. The matrix $\mathbf{A} = \mathbf{U\Sigma V}^T$ splits into two rank-one matrices, columns times rows, with $\sqrt{2}\sqrt{10} = \sqrt{20}$.

$$\sigma_1|\mathbf{u}_1\rangle\langle\mathbf{v}_1| + \sigma_2|\mathbf{u}_2\rangle\langle\mathbf{v}_2| = \sqrt{\frac{45}{20}}\begin{pmatrix}1 & 1\\3 & 3\end{pmatrix} + \sqrt{\frac{5}{20}}\begin{pmatrix}3 & -3\\-1 & 1\end{pmatrix} = \begin{pmatrix}3 & 0\\4 & 5\end{pmatrix} = \mathbf{A}.$$

Every matrix is a sum of rank one matrices with orthogonal $|\mathbf{u}\rangle$'s and orthogonal $|\mathbf{v}\rangle$'s.

More generally, the singular value decomposition (SVD) is a factorization of a real or complex matrix. Specifically, the singular value decomposition of an $m \times n$ complex matrix $\mathbf{A}$ is a factorization of the form $\mathbf{A} = \mathbf{U\Sigma V}^*$ where $\mathbf{U}$ is an $m \times m$ complex unitary matrix, $\mathbf{\Sigma}$ is an $m \times n$ rectangular diagonal matrix with non-negative real numbers on the diagonal, $\mathbf{V}$ is an $n \times n$ complex unitary matrix, and $\mathbf{V}^*$ is the conjugate transpose of $\mathbf{V}$. Such decomposition always exists for any complex matrix (see Figure 3.3). If $\mathbf{A}$ is real, then $\mathbf{U}$ and $\mathbf{V}$ can be guaranteed to be real orthogonal matrices; in such contexts, the SVD is often denoted $\mathbf{A} = \mathbf{U\Sigma V}^T$. The diagonal entries $\sigma_i = \Sigma_{ii}$ of $\mathbf{\Sigma}$ are uniquely determined by $\mathbf{A}$ and are known as the singular values of $\mathbf{A}$. The number of non-zero singular values is equal to the rank of $\mathbf{A}$. The columns of $\mathbf{U}$ and the columns of $\mathbf{V}$ are called left-singular vectors and right-singular vectors of $\mathbf{A}$, respectively. They form two sets of orthonormal bases $|\mathbf{u}_1\rangle, \ldots, |\mathbf{u}_m\rangle$ and $|\mathbf{v}_1\rangle, \ldots, |\mathbf{v}_n\rangle$.

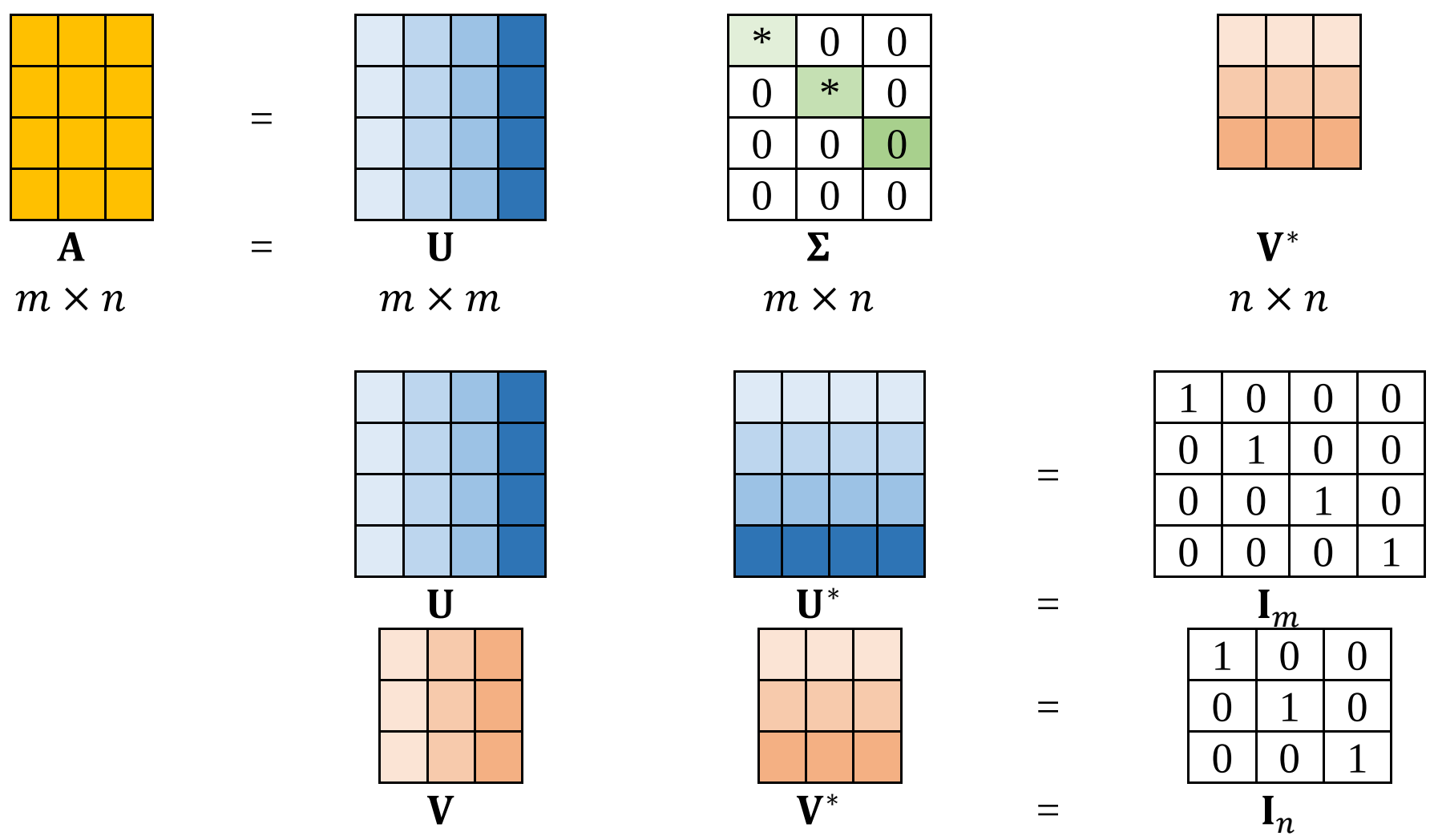

**Figure 3.3.** Visualization of the matrix multiplications in singular value decomposition

***Example 3.21***

$$\mathbf{A} = \begin{pmatrix} 1 & 0 & 0 & 0 & 2 \\ 0 & 0 & 3 & 0 & 0 \\ 0 & 0 & 0 & 0 & 0 \\ 0 & 2 & 0 & 0 & 0 \end{pmatrix}$$

A singular value decomposition of this matrix is given by $\mathbf{U\Sigma V}^*$

$$\mathbf{U} = \begin{pmatrix} 0 & -1 & 0 & 0 \\ -1 & 0 & 0 & 0 \\ 0 & 0 & 0 & -1 \\ 0 & 0 & -1 & 0 \end{pmatrix}, \quad \mathbf{\Sigma} = \begin{pmatrix} 3 & 0 & 0 & 0 & 0 \\ 0 & \sqrt{5} & 0 & 0 & 0 \\ 0 & 0 & 2 & 0 & 0 \\ 0 & 0 & 0 & 0 & 0 \end{pmatrix}, \quad \mathbf{V}^* = \begin{pmatrix} 0 & 0 & -1 & 0 & 0 \\ -\sqrt{0.2} & 0 & 0 & 0 & -\sqrt{0.8} \\ 0 & -1 & 0 & 0 & 0 \\ 0 & 0 & 0 & 1 & 0 \\ -\sqrt{0.8} & 0 & 0 & 0 & \sqrt{0.2} \end{pmatrix}$$

with

$$\mathbf{UU}^* = \begin{pmatrix} 1 & 0 & 0 & 0 \\ 0 & 1 & 0 & 0 \\ 0 & 0 & 1 & 0 \\ 0 & 0 & 0 & 1 \end{pmatrix} = \mathbf{I}, \quad \mathbf{VV}^* = \begin{pmatrix} 1 & 0 & 0 & 0 & 0 \\ 0 & 1 & 0 & 0 & 0 \\ 0 & 0 & 1 & 0 & 0 \\ 0 & 0 & 0 & 1 & 0 \\ 0 & 0 & 0 & 0 & 1 \end{pmatrix} = \mathbf{I}$$

The singular value decomposition is very general in the sense that it can be applied to any $m \times n$ matrix, whereas eigenvalue decomposition can only be applied to diagonalizable matrices. Nevertheless, the two decompositions are related. Given an SVD of **A**, as described above, the following two relations hold:

$$\mathbf{A}^*\mathbf{A} = \mathbf{V\Sigma}^*\mathbf{U}^*\mathbf{U\Sigma V}^* = \mathbf{V}(\mathbf{\Sigma}^*\mathbf{\Sigma})\mathbf{V}^* \tag{3.81}$$
$$\mathbf{AA}^* = \mathbf{U\Sigma V}^*\mathbf{V\Sigma}^*\mathbf{U}^* = \mathbf{U}(\mathbf{\Sigma\Sigma}^*)\mathbf{U}^* \tag{3.82}$$

The right-hand sides of these relations describe the eigenvalue decompositions of the left-hand sides. Consequently:

- The columns of **V** (right-singular vectors) are eigenvectors of $\mathbf{A}^*\mathbf{A}$.
- The columns of **U** (left-singular vectors) are eigenvectors of $\mathbf{AA}^*$.
- The non-zero elements of **Σ** (non-zero singular values) are the square roots of the non-zero eigenvalues of $\mathbf{A}^*\mathbf{A}$ or $\mathbf{AA}^*$.

### The Geometry of the SVD

The SVD separates a matrix into $\mathbf{A} = \mathbf{U\Sigma V}^T$: (orthogonal)×(diagonal)×(orthogonal). In two dimensions, we can draw those steps. The orthogonal matrices **U** and **V** rotate the plane. The diagonal matrix **Σ** stretches it along the axes. Figure 3.4 shows rotation times, stretching times rotation. Vectors $|\mathbf{x}\rangle$ on the unit circle go to $\mathbf{A}|\mathbf{x}\rangle$ on an ellipse. This picture applies to a two-by-two invertible matrix (because $\sigma_1 > 0$ and $\sigma_2 > 0$). First is a rotation of any $|\mathbf{x}\rangle$ to $\mathbf{V}^\mathbf{T}|\mathbf{x}\rangle$. Then **Σ** stretches that vector to $\mathbf{\Sigma V}^T|\mathbf{x}\rangle$. Then **U** rotates to $\mathbf{A}|\mathbf{x}\rangle = \mathbf{U\Sigma V}^T|\mathbf{x}\rangle$. We kept all determinants positive to avoid reflections. The four numbers $a, b, c, d$ in the matrix connect to two angles $\theta$ and $\phi$ and two numbers $\sigma_1$ and $\sigma_2$.

$$\begin{pmatrix} a & b \\ c & d \end{pmatrix} = \begin{pmatrix} \cos\theta & -\sin\theta \\ \sin\theta & \cos\theta \end{pmatrix} \begin{pmatrix} \sigma_1 & 0 \\ 0 & \sigma_2 \end{pmatrix} \begin{pmatrix} \cos\phi & \sin\phi \\ -\sin\phi & \cos\phi \end{pmatrix}. \tag{3.83}$$

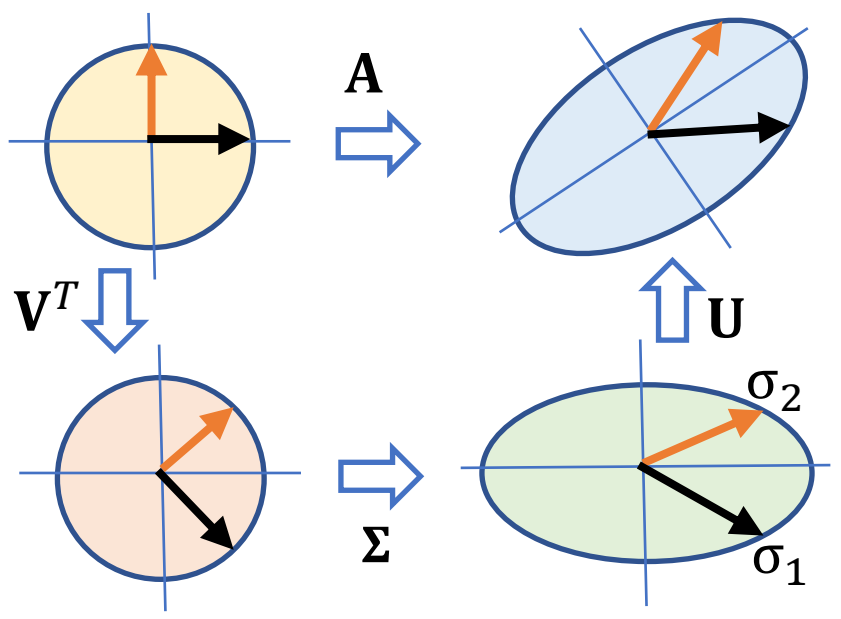

**Figure 3.4.** **U** and **V** are rotations and possible reflections. **Σ** stretches the circle to an ellipse.

### The Reduced Form of the SVD

The full form $\mathbf{AV} = \mathbf{U\Sigma}$ in (3.73) can have a lot of zeros in $\mathbf{\Sigma}$ when the rank of $\mathbf{A}$ is small and its null-space is large. Those zeros contribute nothing to matrix multiplication. The heart of the SVD is in the first $r$ $|\mathbf{v}\rangle$'s and $m$ $|\mathbf{u}\rangle$ 's and $\sigma$'s. We can reduce $\mathbf{AV} = \mathbf{U\Sigma}$ to $\mathbf{AV}_r = \mathbf{U}_r\mathbf{\Sigma}_r$ by removing the parts that are sure to produce zeros. This leaves the reduced SVD where $\mathbf{\Sigma}_r$ is now square:

$$\mathbf{AV}_r = \mathbf{U}_r\mathbf{\Sigma}_r \Rightarrow \mathbf{A}\begin{pmatrix} |\mathbf{v}_1\rangle & \underset{\text{row space}}{\ldots} & |\mathbf{v}_r\rangle \end{pmatrix} = \begin{pmatrix} |\mathbf{u}_1\rangle & \underset{\text{column space}}{\ldots} & |\mathbf{u}_r\rangle \end{pmatrix}\begin{pmatrix} \sigma_1 & & \\ & \ddots & \\ & & \sigma_r \end{pmatrix}. \tag{3.84}$$

We still have $\mathbf{V}_r^T\mathbf{V}_r = \mathbf{I}_r$ and $\mathbf{U}_r^T\mathbf{U}_r = \mathbf{I}_r$ from those orthogonal unit vectors $|\mathbf{v}\rangle$'s and $|\mathbf{u}\rangle$'s. But when $\mathbf{V}_r$ and $\mathbf{U}_r$ are not square, we can no longer have two-sided inverses: $\mathbf{V}_r\mathbf{V}_r^T \neq \mathbf{I}$ and $\mathbf{U}_r\mathbf{U}_r^T \neq \mathbf{I}$. (See Figure 3.5.)

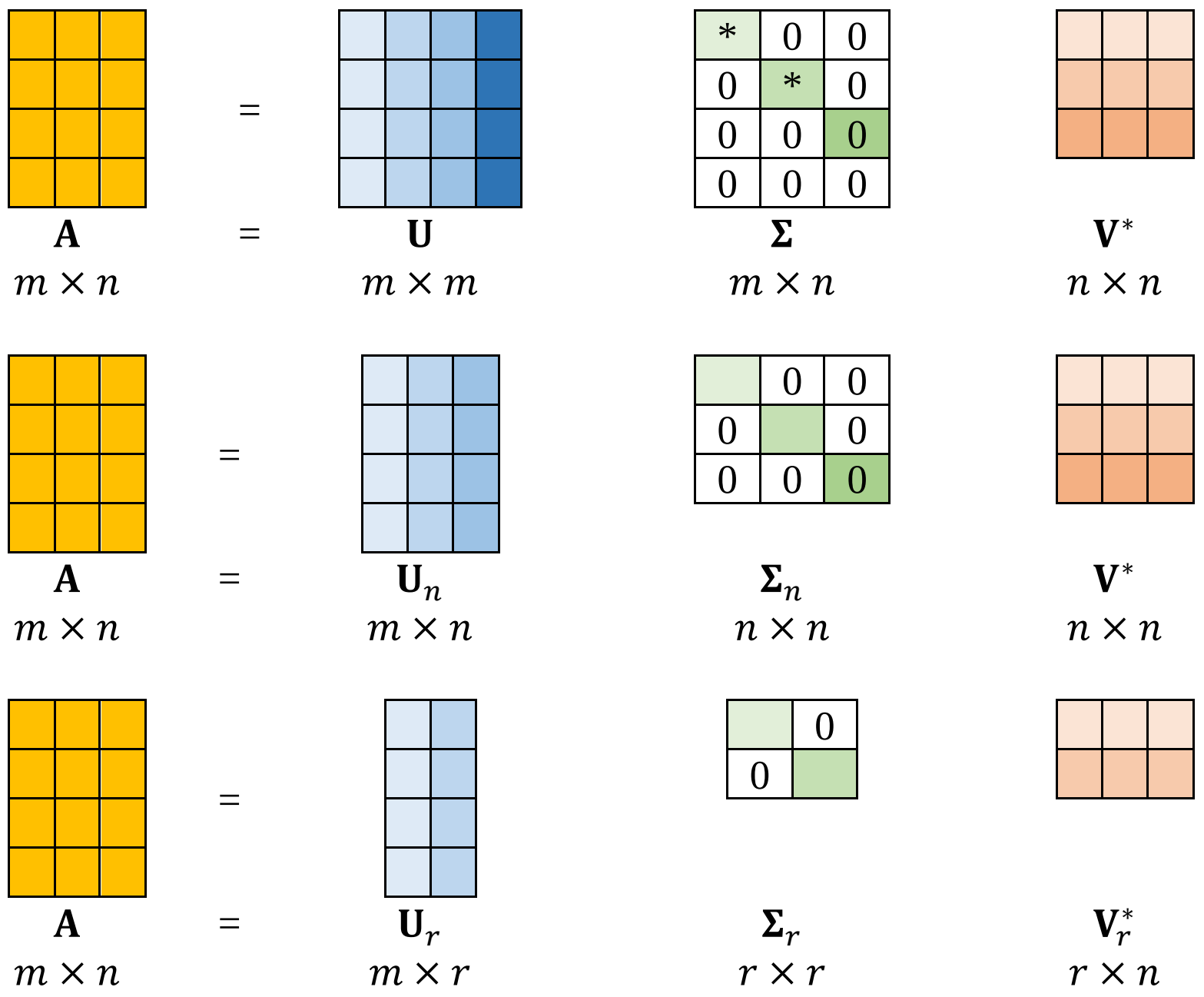

**Figure 3.5.** Visualization of the matrix multiplications in the reduced SVD.

## 3.9 Quadratic Form and Definite Matrix in $\mathbb{R}^n$

The reader should be well aware of the important roles of matrices in the study of linear equations, which can be expressed in the form

$$a_1x_1 + a_2x_2 + \ldots + a_nx_n = b. \tag{3.85}$$

The left side $a_1x_1 + a_2x_2 + \ldots + a_nx_n = \langle\mathbf{a}|\mathbf{x}\rangle$ of the equation is a (homogeneous) polynomial of degree 1 in $n$ variables, called a linear form. In this section, we study a (homogeneous) polynomial of degree 2 in several variables, called a quadratic form, and show that matrices also play an important role in the study of a quadratic form.

A quadratic equation in two variables $x$ and $y$ is an equation of the form

$$ax^2 + 2bxy + cy^2 + dx + ey + f = 0, \tag{3.86}$$

in which the left side consists of a constant term $f$, a linear form $dx + ey$, and a quadratic form $ax^2 + 2bxy + cy^2$. Note that this quadratic form may be written in matrix notation as

$$ax^2 + 2bxy + cy^2 = (x \quad y)\begin{pmatrix} a & b \\ b & c \end{pmatrix}\begin{pmatrix} x \\ y \end{pmatrix} = \langle \mathbf{x}|\mathbf{A}|\mathbf{x}\rangle, \tag{3.87}$$

where

$$|\mathbf{x}\rangle = \begin{pmatrix} x \\ y \end{pmatrix} \text{ and } \mathbf{A} = \begin{pmatrix} a & b \\ b & c \end{pmatrix}. \tag{3.88}$$

Note also that the matrix **A** is taken to be a (real) symmetric matrix. Geometrically, the solution set of a quadratic equation in $x$ and $y$ usually represents a conic section, such as an ellipse, a parabola, or a hyperbola in the $xy$-plane.

**Definitions:** An equation of the form

$$f|\mathbf{x}\rangle = \sum_{i=1}^{n}\sum_{j=1}^{n} A_{ij}x_i x_j + \sum_{i=1}^{n} b_i x_i + c = 0, \tag{3.89}$$

where $A_{ij}$, $b_i$ and $c$ are real constants, is called a quadratic equation in $n$ variables $x_1, x_2, ..., x_n$. In matrix form, it can be written as

$$f|\mathbf{x}\rangle = \langle \mathbf{x}|\mathbf{A}|\mathbf{x}\rangle + \langle \mathbf{b}|\mathbf{x}\rangle + c = 0, \tag{3.90}$$

where $\mathbf{A} = (A_{ij})$, $|\mathbf{x}\rangle = (x_1 \dots x_n)^T$ and $|\mathbf{b}\rangle = (b_1 \dots b_n)^T$ in $\mathbb{R}^n$.

(1) A linear form is a polynomial of degree 1 in $n$ variables $x_1, x_2, \dots, x_n$ of the form

$$\langle \mathbf{b}|\mathbf{x}\rangle = \sum_{i=1}^{n} b_i x_i, \tag{3.91}$$

where $|\mathbf{x}\rangle = (x_1 \dots x_n)^T$ and $|\mathbf{b}\rangle = (b_1 \dots b_n)^T$ in $\mathbb{R}^n$.

(2) A quadratic form is a (homogeneous) polynomial of degree 2 in $n$ variables $x_1, x_2, \dots, x_n$ of the form

$$q|\mathbf{x}\rangle = \langle \mathbf{x}|\mathbf{A}|\mathbf{x}\rangle = (x_1 \dots x_n)(A_{ij})\begin{pmatrix} x_1 \\ \vdots \\ x_n \end{pmatrix} = \sum_{i=1}^{n}\sum_{j=1}^{n} A_{ij}x_i x_j, \tag{3.92}$$

where $|\mathbf{x}\rangle = (x_1 \dots x_n)^T \in \mathbb{R}^n$ and $\mathbf{A} = (A_{ij})$ is a real $n \times n$ matrix.

***Example 3.22***

Consider the matrix $\mathbf{A} = \begin{pmatrix} 1 & 2 \\ 2 & 1 \end{pmatrix}$ and the vector $|\mathbf{x}\rangle = (x_1 \quad x_2)^T$.

***Solution***

$q$ is given by

$$\begin{aligned} q = \langle \mathbf{x}|\mathbf{A}|\mathbf{x}\rangle &= (x_1 \quad x_2)\begin{pmatrix} 1 & 2 \\ 2 & 1 \end{pmatrix}\begin{pmatrix} x_1 \\ x_2 \end{pmatrix} \\ &= (x_1 + 2x_2 \quad 2x_1 + x_2)\begin{pmatrix} x_1 \\ x_2 \end{pmatrix} \\ &= x_1^2 + 2x_1x_2 + 2x_1x_2 + x_2^2 = x_1^2 + 4x_1x_2 + x_2^2. \end{aligned}$$

***Example 3.23***

Consider the $3 \times 3$ diagonal matrix $\mathbf{A} = \begin{pmatrix} 1 & 0 & 0 \\ 0 & 2 & 0 \\ 0 & 0 & 4 \end{pmatrix}$ and a general 3-element vector $|\mathbf{x}\rangle = (x_1, \quad x_2, \quad x_3)^T$.

***Solution***

The general quadratic form is given by

$$\begin{aligned} q = \langle \mathbf{x}|\mathbf{A}|\mathbf{x}\rangle &= (x_1, \quad x_2, \quad x_3)\begin{pmatrix} 1 & 0 & 0 \\ 0 & 2 & 0 \\ 0 & 0 & 4 \end{pmatrix}\begin{pmatrix} x_1 \\ x_2 \\ x_3 \end{pmatrix} \\ &= (x_1, \quad 2x_2, \quad 4x_3)\begin{pmatrix} x_1 \\ x_2 \\ x_3 \end{pmatrix} \\ &= x_1^2 + 2x_2^2 + 4x_3^2. \end{aligned}$$

***Example 3.24***

Also, consider the following matrix $\mathbf{A} = \begin{pmatrix} -2 & 1 & 0 \\ 1 & -2 & 0 \\ 0 & 0 & -2 \end{pmatrix}$.

***Solution***

The general quadratic form is given by

$$\begin{aligned} q = \langle \mathbf{x}|\mathbf{A}|\mathbf{x}\rangle &= (x_1, \quad x_2, \quad x_3)\begin{pmatrix} -2 & 1 & 0 \\ 1 & -2 & 0 \\ 0 & 0 & -2 \end{pmatrix}\begin{pmatrix} x_1 \\ x_2 \\ x_3 \end{pmatrix} \\ &= (-2x_1 + x_2, \quad x_1 - 2x_2, \quad -2x_3)\begin{pmatrix} x_1 \\ x_2 \\ x_3 \end{pmatrix} \\ &= -2x_1^2 + x_1x_2 + x_1x_2 - 2x_2^2 - 2x_3^2 \\ &= -2x_1^2 + 2x_1x_2 - 2x_2^2 - 2x_3^2 \\ &= -2(x_1^2 - x_1x_2) - 2x_2^2 - 2x_3^2 \\ &= -2x_1^2 - 2(x_2^2 - x_1x_2) - 2x_3^2. \end{aligned}$$

**Theorem 3.28:** The matrix **A** in the definition of a quadratic form is any square matrix, but it can be restricted to be a symmetric matrix.

**Proof:**

In fact, any square matrix **A** is the sum of a symmetric part **B** ($\mathbf{B}^T = \mathbf{B}$) and a skew-symmetric part **C** ($\mathbf{C}^T = -\mathbf{C}$), say

$$\mathbf{A} = \mathbf{B} + \mathbf{C}, \qquad \mathbf{B} = \frac{1}{2}(\mathbf{A} + \mathbf{A}^T), \qquad \mathbf{C} = \frac{1}{2}(\mathbf{A} - \mathbf{A}^T),$$

where,

$$\mathbf{B}^T = \frac{1}{2}(\mathbf{A} + \mathbf{A}^T)^T = \frac{1}{2}(\mathbf{A}^T + (\mathbf{A}^T)^T) = \frac{1}{2}(\mathbf{A} + \mathbf{A}^T) = \mathbf{B},$$

$$\mathbf{C}^T = \frac{1}{2}(\mathbf{A} - \mathbf{A}^T)^T = \frac{1}{2}(\mathbf{A}^T - (\mathbf{A}^T)^T) = \frac{1}{2}(\mathbf{A}^T - \mathbf{A}) = -\frac{1}{2}(\mathbf{A} - \mathbf{A}^T) = -\mathbf{C}.$$

For the skew-symmetric matrix **C**, we have

$$\langle \mathbf{x}|\mathbf{C}|\mathbf{x}\rangle = \langle \mathbf{x}|\mathbf{C}|\mathbf{x}\rangle^T = \langle \mathbf{x}|\mathbf{C}^T|\mathbf{x}\rangle = -\langle \mathbf{x}|\mathbf{C}|\mathbf{x}\rangle.$$

Since, $\langle \mathbf{x}|\mathbf{C}|\mathbf{x}\rangle$ is a real number. Hence, $\langle \mathbf{x}|\mathbf{C}|\mathbf{x}\rangle = 0$. Therefore,

$$q(\mathbf{x}) = \langle \mathbf{x}|\mathbf{A}|\mathbf{x}\rangle = \langle \mathbf{x}|(\mathbf{B} + \mathbf{C})|\mathbf{x}\rangle = \langle \mathbf{x}|\mathbf{B}|\mathbf{x}\rangle.$$

This means that, without loss of generality, one may assume that the matrix **A** in the definition of a quadratic form is a symmetric matrix.

∎

**Remark:**

(1) A quadratic equation is said to be consistent if it has a solution, i.e., there is a vector $|\mathbf{x}\rangle \in \mathbb{R}^n$ such that $f|\mathbf{x}\rangle = 0$. Otherwise, it is said to be inconsistent. For instance, equation $2x^2 + 3y^2 = -1$ is inconsistent. In the following discussion, we will consider only consistent equations.

(2) From the definition of a quadratic form, one can see that fixing a basis like the standard basis for $\mathbb{R}^n$, a quadratic form is associated with a unique symmetric matrix, which is called the matrix representation of the quadratic form $q$ with respect to the basis chosen (the standard basis for $\mathbb{R}^n$). On the other hand, any (real) symmetric matrix **A** gives rise to a quadratic form $\langle\mathbf{x}|\mathbf{A}|\mathbf{x}\rangle$. For example, for a symmetric matrix $\begin{pmatrix}8 & 2\\ 2 & -1\end{pmatrix}$, the equation $(x_1 \quad x_2)\begin{pmatrix}8 & 2\\ 2 & -1\end{pmatrix}\begin{pmatrix}x_1\\ x_2\end{pmatrix}$ defines a quadratic form $8x_1^2 + 4x_1x_2 - x_2^2$.

### Diagonalization of a quadratic form

To study the solution of a quadratic equation $f|\mathbf{x}\rangle = 0$, we first consider an equation $\langle\mathbf{x}|\mathbf{A}|\mathbf{x}\rangle = c$ without a linear form. This quadratic form may be rewritten as the sum of two parts:

$$q|\mathbf{x}\rangle = \langle\mathbf{x}|\mathbf{A}|\mathbf{x}\rangle = \sum_{i=1}^{n} A_{ii}x_i^2 + 2\sum_{i<j} A_{ij}x_ix_j, \tag{3.93}$$

in which the first part $\sum_{i=1}^{n} A_{ii}x_i^2$ is called the (perfect) square terms, and the second part $\sum_{i\neq j} A_{ij}x_ix_j$ is called the cross-product terms. However, the symmetric matrix $\mathbf{A}$ can be orthogonally diagonalized, i. e., there exists an orthogonal matrix $\mathbf{P}$ ($\mathbf{P}^T = \mathbf{P}^{-1}$) such that

$$\mathbf{P}^T\mathbf{A}\mathbf{P} = \mathbf{P}^{-1}\mathbf{A}\mathbf{P} = \mathbf{D} = \begin{pmatrix}\lambda_1 & 0 & 0\\ 0 & \ddots & 0\\ 0 & 0 & \lambda_n\end{pmatrix}. \tag{3.94}$$

Here, the diagonal entries $\lambda_i$'s are the eigenvalues of $\mathbf{A}$, and the column vectors of $\mathbf{P}$ are their associated eigenvectors of $\mathbf{A}$. Then we get, by setting $|\mathbf{x}\rangle = \mathbf{P}|\mathbf{y}\rangle$,

$$\langle\mathbf{x}|\mathbf{A}|\mathbf{x}\rangle = \langle\mathbf{y}|\mathbf{P}^T\mathbf{A}\mathbf{P}|\mathbf{y}\rangle = \langle\mathbf{y}|\mathbf{D}|\mathbf{y}\rangle = \lambda_1y_1^2 + \lambda_2y_2^2 + \cdots + \lambda_ny_n^2, \tag{3.95}$$

which is a quadratic form without the cross-product terms. Consequently, we have proven the following theorem.

**Theorem 3.29:** Let $\langle\mathbf{x}|\mathbf{A}|\mathbf{x}\rangle$ be a quadratic form in $|\mathbf{x}\rangle = (x_1 \ldots x_n)^T \in \mathbb{R}^n$ for a symmetric matrix $\mathbf{A}$. Then there is a change of the coordinates of $|\mathbf{x}\rangle$ into $|\mathbf{y}\rangle = \mathbf{P}^T|\mathbf{y}\rangle = (y_1 \ldots y_n)^T$ such that

$$\langle\mathbf{x}|\mathbf{A}|\mathbf{x}\rangle = \langle\mathbf{y}|\mathbf{D}|\mathbf{y}\rangle = \lambda_1y_1^2 + \lambda_2y_2^2 + \cdots + \lambda_ny_n^2, \tag{3.96}$$

where $\mathbf{P}$ is an orthogonal matrix and $\mathbf{P}^T\mathbf{A}\mathbf{P} = \mathbf{D}$.

Remark:

(1) Recall that in the above theorem, the columns of the matrix $\mathbf{P}$ are the orthonormal eigenvectors of $\mathbf{A}$ and $|\mathbf{y}\rangle$ is just the coordinate expression of $|\mathbf{x}\rangle$ with respect to the orthonormal eigenvectors of $\mathbf{A}$.

(2) The solution set of a quadratic equation of the form $\langle\mathbf{x}|\mathbf{A}|\mathbf{x}\rangle = c$ is a hypersurface in $\mathbb{R}^n$, that is, a curved surface that can be parameterized in $n-1$ variables. These are called $n-1$ -dimensional quadratic surfaces, with axes in the directions of eigenvectors. In particular, if $n = 2$, the solution set of a quadratic equation is called a quadratic curve, or more commonly, a conic section. When $n = 3$, the quadratic surfaces are ellipsoids or hyperboloids depending on the signs of the eigenvalues of $\mathbf{A}$. Of course, a paraboloid is also a quadratic surface, but it appears when a linear form is present in the quadratic equation. The determination of the quadratic hypersurface depends on the signs of the eigenvalues of $\mathbf{A}$.

***Example 3.25***

Determine the conic section $3x^2 + 2xy + 3y^2 - 8 = 0$.

***Solution***

This equation can be written in the form

$$(x, y)\begin{pmatrix}3 & 1\\ 1 & 3\end{pmatrix}\begin{pmatrix}x\\ y\end{pmatrix} = 8.$$

The matrix $\mathbf{A} = \begin{pmatrix}3 & 1\\ 1 & 3\end{pmatrix}$ has eigenvalues $\lambda_1 = 2$ and $\lambda_2 = 4$ with associated unit eigenvectors

$$|\mathbf{v}_1\rangle = \frac{1}{\sqrt{2}}\begin{pmatrix}1\\ -1\end{pmatrix}, \qquad |\mathbf{v}_2\rangle = \frac{1}{\sqrt{2}}\begin{pmatrix}1\\ 1\end{pmatrix},$$

respectively, which form an orthonormal basis $\beta$. If $\alpha$ denotes the standard basis, then the transition matrix

$$\mathbf{P} = (|\mathbf{v}_1\rangle \quad |\mathbf{v}_2\rangle) = \frac{1}{\sqrt{2}}\begin{pmatrix}1 & 1\\ -1 & 1\end{pmatrix} = \begin{pmatrix}\cos 45 & \sin 45\\ -\sin 45 & \cos 45\end{pmatrix},$$

which is a rotation through 45° in a clockwise direction such that $\mathbf{P}^T = \mathbf{P}^{-1}$. It gives a change of coordinates,

$|\mathbf{x}\rangle = \mathbf{P}|\mathbf{y}\rangle$, i. e.,

$$\begin{pmatrix} x \\ y \end{pmatrix} = \frac{1}{\sqrt{2}}\begin{pmatrix} 1 & 1 \\ -1 & 1 \end{pmatrix}\begin{pmatrix} \acute{x} \\ \acute{y} \end{pmatrix} = \frac{1}{\sqrt{2}}\begin{pmatrix} \acute{x}+\acute{y} \\ -\acute{x}+\acute{y} \end{pmatrix}.$$

Thus, we get

$$\begin{aligned} 3x^2 + 2xy + 3y^2 &= \langle \mathbf{x}|\mathbf{A}|\mathbf{x}\rangle \\ &= \langle \mathbf{y}|\mathbf{P}^T\mathbf{A}\mathbf{P}|\mathbf{y}\rangle \\ &= \langle \mathbf{y}|\mathbf{D}|\mathbf{y}\rangle \\ &= \langle \mathbf{y}|\begin{pmatrix} 2 & 0 \\ 0 & 4 \end{pmatrix}|\mathbf{y}\rangle = 2(\acute{x})^2 + 4(\acute{y})^2 = 8, \end{aligned}$$

or

$$\frac{(\acute{x})^2}{4} + \frac{(\acute{y})^2}{2} = 1.$$

Its solution set is just an ellipse with axes $|\mathbf{v}_1\rangle = \mathbf{P}|\mathbf{e}_1\rangle$ and $|\mathbf{v}_2\rangle = \mathbf{P}|\mathbf{e}_2\rangle$.

**Definition:** Let $\mathbf{A} = (A_{ij}) \in M_{n\times n}(\mathbb{R})$ be a symmetric matrix, and let $|\mathbf{x}\rangle = (x_1, x_2, \ldots, x_n)^T \in \mathbb{R}^n$. Then, $\mathbf{A}$ is said to be

(1) positive definite if $\langle \mathbf{x}|\mathbf{A}|\mathbf{x}\rangle = \sum_{ij} A_{ij}x_ix_j > 0$ for all nonzero $|\mathbf{x}\rangle$,
(2) positive semidefinite if $\langle \mathbf{x}|\mathbf{A}|\mathbf{x}\rangle = \sum_{ij} A_{ij}x_ix_j \geq 0$ for all $|\mathbf{x}\rangle$,
(3) negative definite if $\langle \mathbf{x}|\mathbf{A}|\mathbf{x}\rangle = \sum_{ij} A_{ij}x_ix_j < 0$ for all nonzero $|\mathbf{x}\rangle$,
(4) negative semidefinite if $\langle \mathbf{x}|\mathbf{A}|\mathbf{x}\rangle = \sum_{ij} A_{ij}x_ix_j \leq 0$ for all $|\mathbf{x}\rangle$,
(5) indefinite if $\langle \mathbf{x}|\mathbf{A}|\mathbf{x}\rangle$ takes both positive and negative values.

***Example 3.26***

The real symmetric matrix

$$\begin{pmatrix} 2 & -1 & 0 \\ -1 & 2 & -1 \\ 0 & -1 & 2 \end{pmatrix},$$

is positive definite because the quadratic form satisfies

$$\begin{aligned} \langle \mathbf{x}|\mathbf{A}|\mathbf{x}\rangle &= (x_1 \quad x_2 \quad x_3)\begin{pmatrix} 2 & -1 & 0 \\ -1 & 2 & -1 \\ 0 & -1 & 2 \end{pmatrix}\begin{pmatrix} x_1 \\ x_2 \\ x_3 \end{pmatrix} \\ &= (x_1 \quad x_2 \quad x_3)\begin{pmatrix} 2x_1 - x_2 \\ -x_1 + 2x_2 - x_3 \\ -x_2 + 2x_3 \end{pmatrix} \\ &= x_1(2x_1 - x_2) + x_2(-x_1 + 2x_2 - x_3) + x_3(-x_2 + 2x_3) \\ &= 2x_1^2 - 2x_1x_2 + 2x_2^2 - 2x_2x_3 + 2x_3^2 \\ &= x_1^2 + (x_1 - x_2)^2 + (x_2 - x_3)^2 + x_3^2 > 0, \end{aligned}$$

unless $x_1 = x_2 = x_3 = 0$.

### Graphical analysis

When $|\mathbf{x}\rangle$ has only two elements, we can graphically represent $q$ in 3 dimensions. A positive definite quadratic form will always be positive except at the point where $|\mathbf{x}\rangle = |\mathbf{0}\rangle$. This gives a nice graphical representation where the plane at $|\mathbf{x}\rangle = |\mathbf{0}\rangle$ bounds the function from below. Figure 3.6 shows a positive definite quadratic form. Similarly, a negative definite quadratic form is bounded above by the plane $|\mathbf{x}\rangle = |\mathbf{0}\rangle$. Figure 3.7 shows a negative definite quadratic form. A positive semi-definite quadratic form is bounded below by the plane $|\mathbf{x}\rangle = |\mathbf{0}\rangle$ but will touch the plane at more than the single point $(0,0)^T$, it will touch the plane along a line. Figure 3.8 shows a positive semi-definite quadratic form. A negative semi-definite quadratic form is bounded above by the plane $|\mathbf{x}\rangle = |\mathbf{0}\rangle$ but will touch the plane at more than the single point $(0,0)^T$. It will touch the plane along a line. Figure 3.9 shows a negative-definite quadratic form. An indefinite quadratic form will not lie completely above or below the plane but will lie above for some values of $|\mathbf{x}\rangle$ and below for other values of $|\mathbf{x}\rangle$. Figure 3.10 shows an indefinite quadratic form.

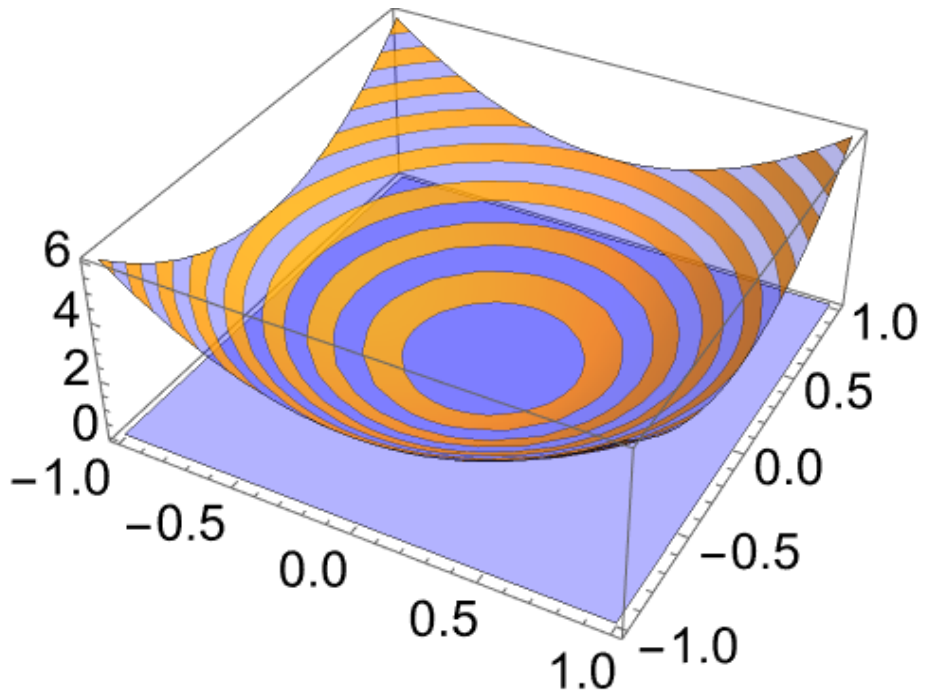

**Figure 3.6.** Positive definite quadratic form $3x_1^2 + 3x_2^2$.

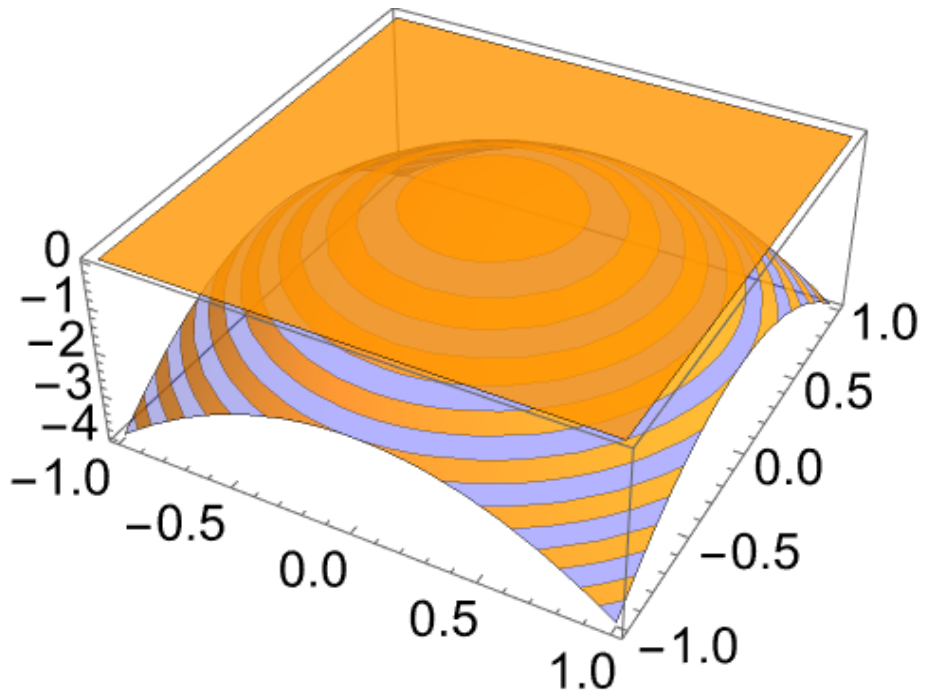

**Figure 3.7.** Negative definite quadratic form $-2x_1^2 - 2x_2^2$.

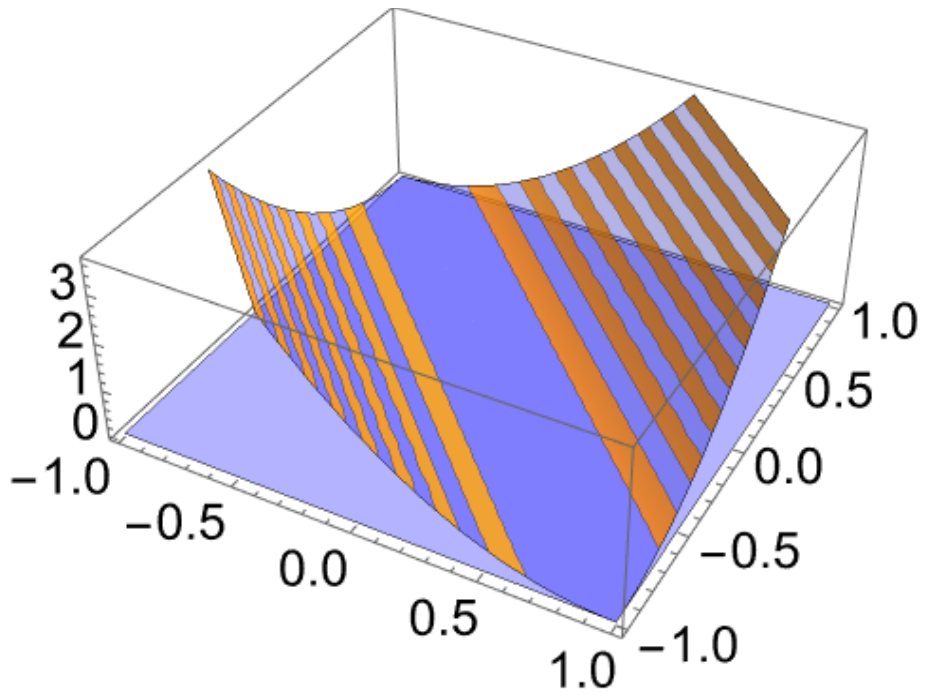

**Figure 3.8.** Positive semi-definite quadratic form $2x_1^2 + 4x_1x_2 + 2x_2^2$.

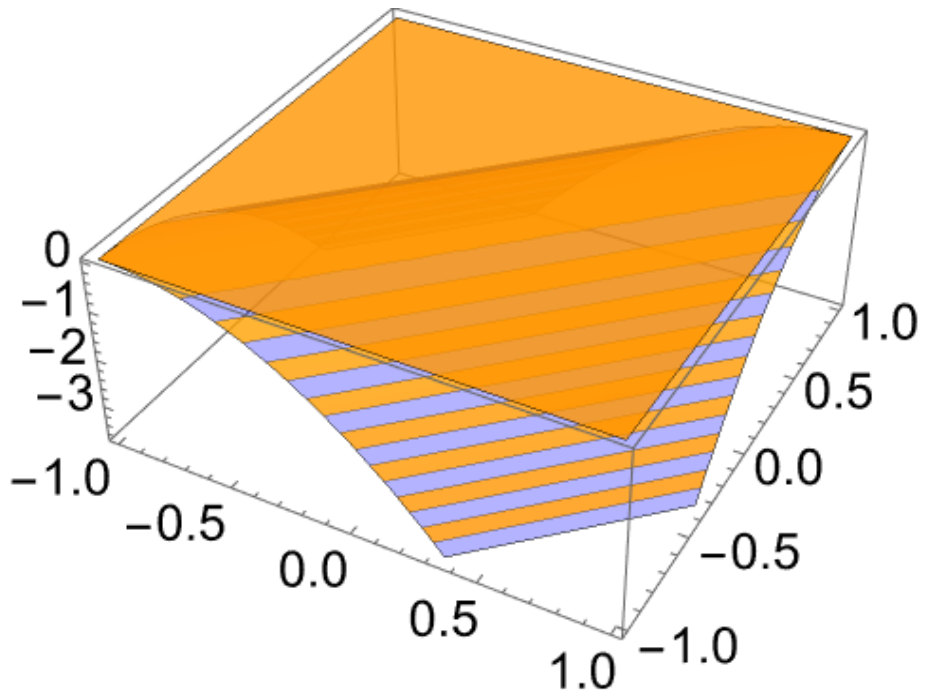

**Figure 3.9.** Negative semi-definite quadratic form $-2x_1^2 + 4x_1x_2 - 2x_2^2$.

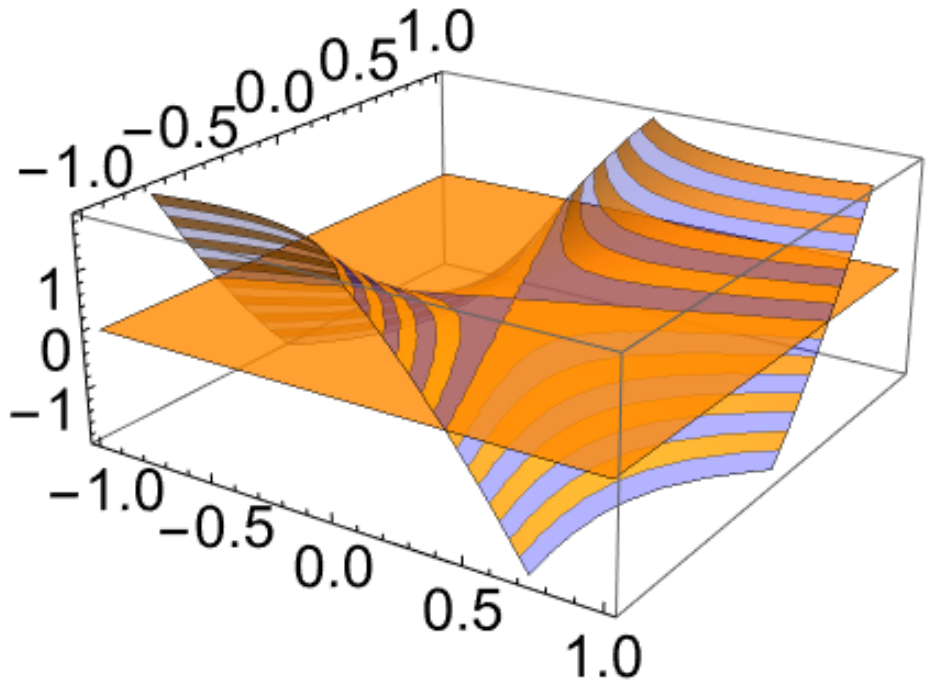

**Figure 3.10.** Indefinite quadratic form $-2x_1^2 + 4x_1x_2 + 2x_2^2$.

To determine whether or not a matrix **A** is positive definite, one can diagonalize **A** so that

$$\langle \mathbf{x}|\mathbf{A}|\mathbf{x}\rangle = \sum_{i,j=1}^{n} A_{ij}x_i x_j = \langle \mathbf{y}|\mathbf{D}|\mathbf{y}\rangle = \sum_{i=1}^{n} \lambda_i y_i^2, \tag{3.97}$$

where $|\mathbf{y}\rangle = \mathbf{P}^T|\mathbf{x}\rangle$ for an orthogonal matrix **P** and the $\lambda_i$ 's are eigenvalues of **A**. Therefore, $\langle \mathbf{x}|\mathbf{A}|\mathbf{x}\rangle > 0$ for all nonzero $|\mathbf{x}\rangle \in \mathbb{R}^n$ if and only if all the $\lambda_i$ 's are positive. Consequently, we have the following characterization of positive definite matrices:

**Theorem 3.30:** *A* real symmetric n × n matrix **A** is positive definite if and only if all the eigenvalues of **A** are positive.

So far, we have seen that it is important to determine whether or not a symmetric matrix **A** is positive definite. In most cases, the definition does not help much. But we have seen that Theorem 3.30 gives us a practical characterization of positive definite matrices: **A** is positive definite if and only if all eigenvalues of **A** are positive. We will find some other practical criteria in terms of the determinant of the matrix. For this, we again look at the quadratic form in two variables, $q|\mathbf{x}\rangle = ax^2 + 2bxy + cy^2$, $|\mathbf{x}\rangle = (x, y)^T$ which may be rewritten in a complete square form as

$$q|\mathbf{x}\rangle = ax^2 + 2bxy + cy^2 = a\left(x + \frac{b}{a}y\right)^2 + \left(c - \frac{b^2}{a}\right)y^2. \tag{3.98}$$

We see that $q$ is a positive definite, i.e., $q|\mathbf{x}\rangle = \langle \mathbf{x}|\mathbf{A}|\mathbf{x}\rangle > 0$ for any nonzero vector $|\mathbf{x}\rangle = (x, y)^T \in \mathbb{R}^2$, if and only if $a > 0$ and $ac > b^2$, or equivalently, the determinants of

$$(a) \text{ and } \begin{pmatrix} a & b \\ b & c \end{pmatrix}, \tag{3.99}$$

are positive. The natural generalization of the above conditions will involve all $n$-submatrices of **A**, called the principal submatrices of **A**, which are defined as the upper left square submatrices

$$\mathbf{A}_1 = (a_{11}), \qquad \mathbf{A}_2 = \begin{pmatrix} a_{11} & a_{12} \\ a_{21} & a_{22} \end{pmatrix}, \qquad \mathbf{A}_3 = \begin{pmatrix} a_{11} & a_{12} & a_{13} \\ a_{21} & a_{22} & a_{23} \\ a_{31} & a_{32} & a_{33} \end{pmatrix}, \dots, \mathbf{A}_n = \mathbf{A}, \tag{3.100}$$

or schematically represented in the following equation.

$$\begin{vmatrix} a_{11} & a_{12} & a_{13} & \dots & a_{1n} \\ a_{21} & a_{22} & a_{23} & \dots & a_{2n} \\ a_{31} & a_{32} & a_{33} & \dots & a_{3n} \\ & & & & \\ a_{n1} & a_{n2} & a_{n3} & \dots & a_{nn} \end{vmatrix} \tag{3.101}$$

With this construction, we have the following characterization of positive definite matrices.

**Theorem 3.31:** The following are equivalent for a real symmetric matrix **A**:
(1) **A** is positive definite, i. e., $\langle \mathbf{x}|\mathbf{A}|\mathbf{x}\rangle > 0$ for all nonzero vector $|\mathbf{x}\rangle$;
(2) all the eigenvalues of **A** are positive;

(3) all the principal submatrices $\mathbf{A}_k$'s have positive determinants;
(4) all the pivots (without row interchanges) are positive;
(5) there exists a non-singular matrix $\mathbf{W}$ such that $\mathbf{A} = \mathbf{W}^T\mathbf{W}$.

**Proof:**

(3) If $\mathbf{A}$ has positive eigenvalues $\lambda_1, \lambda_2, ..., \lambda_n$, then $\det \mathbf{A} = \lambda_1\lambda_2 \dots \lambda_n > 0$. To prove the same result for all the submatrices $\mathbf{A}_k$, we show that if $\mathbf{A}$ is positive definite, so is every $\mathbf{A}_k$. For each $k = 1, \dots, n$, consider all the vectors whose last $n - k$ components are zero, say $|\mathbf{x}\rangle = (x_1 \dots x_k, 0, \dots, 0)^T = \begin{pmatrix} |\mathbf{x}_k\rangle \\ |\mathbf{0}\rangle \end{pmatrix}$, where $|\mathbf{x}_k\rangle$ is any vector in $\mathbb{R}^k$. Then

$$\langle \mathbf{x}|\mathbf{A}|\mathbf{x}\rangle = (\langle \mathbf{x}_k|, \langle \mathbf{0}|) \begin{pmatrix} \mathbf{A}_k & * \\ & * \end{pmatrix} \begin{pmatrix} |\mathbf{x}_k\rangle \\ |\mathbf{0}\rangle \end{pmatrix} = \langle \mathbf{x}_k|\mathbf{A}_k|\mathbf{x}_k\rangle.$$

Thus $\langle \mathbf{x}|\mathbf{A}|\mathbf{x}\rangle > 0$ for all such nonzero $|\mathbf{x}\rangle$ if and only if $\langle \mathbf{x}_k|\mathbf{A}_k|\mathbf{x}_k\rangle > 0$ for all nonzero $|\mathbf{x}_k\rangle \in \mathbb{R}^k$ that is, $\mathbf{A}_k$'s are positive definite, all eigenvalues of $\mathbf{A}_k$ are positive, and its determinant is positive.

(4) Recall that the symmetric matrix $\mathbf{A}$ can be factorized uniquely into the form

$$\mathbf{A} = \mathbf{LDL}^T,$$

where $\mathbf{L}$ is a lower triangular matrix with 1's on its diagonal, and $\mathbf{D}$ is the diagonal matrix with the pivots $d_k$ of $\mathbf{A}$ on the diagonal. But the $k$-th pivot $d_k$ is exactly the ratio of $\det \mathbf{A}_k$ to $\det \mathbf{A}_{k-1}$: $d_k = \frac{\det \mathbf{A}_k}{\det \mathbf{A}_{k-1}}$. Hence, all $d_k$'s are positive.

(5) Let $\mathbf{A} = \mathbf{LDL}^T$ as above with

$$\mathbf{D} = \begin{pmatrix} d_1 & 0 & 0 \\ 0 & \ddots & 0 \\ 0 & 0 & d_n \end{pmatrix}, d_i > 0.$$

Define

$$\sqrt{\mathbf{D}} = \begin{pmatrix} \sqrt{d_1} & 0 & 0 \\ 0 & \ddots & 0 \\ 0 & 0 & \sqrt{d_n} \end{pmatrix}.$$

Then, clearly $\det \sqrt{\mathbf{D}} > 0$, $\mathbf{D} = \sqrt{\mathbf{D}}\sqrt{\mathbf{D}}$ and $\left(\sqrt{\mathbf{D}}\right)^T = \sqrt{\mathbf{D}}$. Hence,

$$\mathbf{A} = \mathbf{LDL}^T = \left(\mathbf{L}\sqrt{\mathbf{D}}\right)\left(\sqrt{\mathbf{D}}\mathbf{L}^T\right) = \left(\mathbf{L}\sqrt{\mathbf{D}}\right)\left(\mathbf{L}\sqrt{\mathbf{D}}\right)^T = \mathbf{W}^T\mathbf{W},$$

where $\mathbf{W} = \left(\mathbf{L}\sqrt{\mathbf{D}}\right)^T$, which is non-singular since $\mathbf{L}$ and $\sqrt{\mathbf{D}}$ are. If $\mathbf{A}$ is real symmetric and $\mathbf{A} = \mathbf{W}^T\mathbf{W}$, where $\mathbf{W}$ is non-singular, then for $|\mathbf{x}\rangle \neq |\mathbf{0}\rangle$ we have

$$\langle \mathbf{x}|\mathbf{A}|\mathbf{x}\rangle = \langle \mathbf{x}|\mathbf{W}^T\mathbf{W}|\mathbf{x}\rangle = \langle \mathbf{Wx}|\mathbf{Wx}\rangle = ||\mathbf{Wx}||^2 > 0,$$

because $|\mathbf{Wx}\rangle \neq |\mathbf{0}\rangle$.

■

***Example 3.27***

As an example consider the following matrix

$$\mathbf{G} = \begin{pmatrix} 4 & 2 & 0 \\ 2 & 9 & 0 \\ 0 & 0 & 2 \end{pmatrix}.$$

The element $a_{11} = 4 > 0$. Now consider the determinant of the first naturally occurring principal $2 \times 2$ submatrix $\begin{vmatrix} 4 & 2 \\ 2 & 9 \end{vmatrix} = 36 - 4 = 32 > 0$. Finally, consider the determinant of the entire matrix $\begin{vmatrix} 4 & 2 & 0 \\ 2 & 9 & 0 \\ 0 & 0 & 2 \end{vmatrix} = 64 > 0$. This matrix is then positive definite.

We now consider semidefinite matrices. One can easily establish the following analogous theorem.

**Theorem 3.32:** The following are equivalent for a real symmetric matrix $\mathbf{A}$:
(1) $\mathbf{A}$ is positive semidefinite, i.e., $\langle \mathbf{x}|\mathbf{A}|\mathbf{x}\rangle \geq 0$ for all vectors $|\mathbf{x}\rangle$;
(2) all the eigenvalues of $\mathbf{A}$ are non-negative.
(3) all the principal submatrices $\mathbf{A}_k$'s have nonnegative determinants.
(4) all the pivots (without row exchanges) are non-negative.
(5) there exists a matrix $\mathbf{W}$, possibly singular, such that $\mathbf{A} = \mathbf{W}^T\mathbf{W}$.

**Theorem 3.33:** Let $\mathbf{A}$ be a symmetric matrix of order $m$. Then $\mathbf{A}$ is negative definite iff its naturally ordered (leading) principal minors alternate in sign starting with a negative number.

***Example 3.28***

Consider the following matrix

$$\mathbf{E} = \begin{pmatrix} -2 & 1 & 0 \\ 1 & -2 & 0 \\ 0 & 0 & -2 \end{pmatrix}.$$

The element $a_{11} = -2 < 0$. Now consider the first naturally occurring principal $2 \times 2$ submatrix $\begin{vmatrix} -2 & 1 \\ 1 & -2 \end{vmatrix} = 4 - 1 = 3 > 0$. Finally, consider the determinant of the entire matrix $\begin{vmatrix} -2 & 1 & 0 \\ 1 & -2 & 0 \\ 0 & 0 & -2 \end{vmatrix} = -6 < 0$. This matrix is then negative definite.

## 3.10 Mathematica Built-in Functions

| | |
|---|---|
| `SymmetricMatrixQ[m]` | gives True if m is explicitly symmetric, and False otherwise. |
| `AntisymmetricMatrixQ[m]` | gives True if m is explicitly antisymmetric, and False otherwise. |
| `HermitianMatrixQ[m]` | gives True if m is explicitly Hermitian, and False otherwise. |
| `AntihermitianMatrixQ[m]` | gives True if m is explicitly antihermitian, and False otherwise. |
| `NormalMatrixQ[m]` | gives True if m is an explicitly normal matrix, and False otherwise. |
| `DiagonalizableMatrixQ[m]` | gives True if m is diagonalizable, and False otherwise. |
| `OrthogonalMatrixQ[m]` | gives True if m is an explicitly orthogonal matrix, and False otherwise. |
| `UnitaryMatrixQ[m]` | gives True if m is a unitary matrix, and False otherwise. |

***Mathematica Examples 3.1***

```
Input      SymmetricMatrixQ[{{1,2},{2,3}}]
Output     True

Input      AntisymmetricMatrixQ[{{0,a},{-a,0}}]
Output     True

Input      HermitianMatrixQ[{{1,3+4 I},{3-4 I,2}}]
Output     True

Input      AntihermitianMatrixQ[{{I,3+4 I},{-3+4 I,0}}]
Output     True

Input      (*Test if a matrix is normal:*)
           NormalMatrixQ[{{1,2,-1},{-1,1,2},{2,-1,1}}]
Output     True

Input      (*Test if a matrix is diagonalizable:*)
```

```
              DiagonalizableMatrixQ[{{2,1},{1,3}}]
Output        True

Input         (*Test if a matrix is orthogonal:*)
              OrthogonalMatrixQ[1/Sqrt[2] {{1,-1},{1,1}}]
Output        True

Input         (*Test if a matrix is unitary:*)
              UnitaryMatrixQ[1/Sqrt[2] {{1,I},{I,1}}]
Output        True
```

| | |
|---|---|
| MatrixRank[m] | gives the rank of the matrix m. |
| NullSpace[m] | gives a list of vectors that forms a basis for the null space of the matrix m. |
| Inverse[m] | gives the inverse of a square matrix m. |

***Mathematica Examples 3.2***

```
Input         (* Find the number of linearly independent rows of a numerical matrix: *)
              MatrixRank[{{1,2,3},{4,5,6},{7,8,9}}]
Output        2

Input         (* Find the number of linearly independent rows of a symbolic matrix: *)
              MatrixRank[{{a,b,c},{d,e,f},{g,h,i}}]
Output        3

Input         (* Find the null space of a 3×3 matrix: *)
              m={{1,2,3},{4,5,6},{7,8,9}};
              NullSpace[m]

              (* The action of m on the vector is the zero vector: *)
              m.{1,-2,1}
Output        {{1,-2,1}}
Output        {0,0,0}

Input         (* The following three vectors are not linearly independent: *)
              v1={1,2,3};
              v2={4,5,6};
              v3={7,8,9};
              Solve[a v1+b v2==v3,{a,b}]

              (* Therefore the matrix rank of the matrix whose rows are the vectors is 2: *)
              MatrixRank[{v1,v2,v3}]

              (* Therefore the null space of the matrix whose rows are the vectors is
              nonempty: *)
              NullSpace[{v1,v2,v3}]
Output        {{a->-1,b->2}}
Output        2
Output        {{1,-2,1}}

Input         (* The following three vectors are linearly independent: *)
              v1={1,2,3};
              v2={4,5,6};
              v3={7,8,10};
              Solve[a v1+b v2==v3,{a,b}]

              (* Therefore the matrix rank of the matrix whose rows are the vectors is 3: *)
              MatrixRank[{v1,v2,v3}]
```

```
        (*Therefore the null space of the matrix whose rows are the vectors is empty:*)
        NullSpace[{v1,v2,v3}]
Output  {}
Output  3
Output  {}

Input   (* Determine if the following vectors are linearly independent or not: *)
        Subscript[v,1]={108,90,252,186};
        Subscript[v,2]={240,260,520,420};
        Subscript[v,3]={264,340,536,468};
        Subscript[v,4]={522,705,1038,929};

        (* The rank of the matrix formed from the vectors is less than four,so they are
        not linearly independent: *)
        MatrixRank[{Subscript[v,1],Subscript[v,2],Subscript[v,3],Subscript[v,4]}]

        (* The matrix formed from the vectors has a nonempty null space,so they are not
        linearly independent: *)
        NullSpace[{Subscript[v,1],Subscript[v,2],Subscript[v,3],Subscript[v,4]}]
Output  2
Output  {{-44,-3,0,27},{-26,6,9,0}}

Input   (* Find the dimension of the column space of the following matrix: *)
        a=({{1,4,2,-9},{4,12,2,5},{6,7,-11,9},{5,15,10,12}});

        (* The dimension of the space of all linear combinations of the columns equals
        the matrix rank: *)
        MatrixRank[a]

        (* Since the null space is empty,the dimension of the column space equals the
        number of columns: *)
        NullSpace[a]
Output  4
Output  {}

Input   (* Find the dimension of the subspace spanned by the following vectors: *)
        Subscript[v,1]={8,8,0,-4,6};
        Subscript[v,2]={2,-2,-11,1,0};
        Subscript[v,3]={-4,-8,-12,4,-4};
        Subscript[v,4]={8,12,14,-6,6};
        Subscript[v,5]={4,4,6,-2,0};

        (* Since the matrix rank of the matrix formed by the vectors is three,that is
        the dimension of the subspace: *)
        MatrixRank[{Subscript[v,1],Subscript[v,2],Subscript[v,3],Subscript[v,4],Subscrip
        t[v,5]}]

        (* Since the matrix rank of the matrix formed by the vectors is three,that is
        the dimension of the subspace: *)
        With[{a={Subscript[v,1],Subscript[v,2],Subscript[v,3],Subscript[v,4],Subscript[v
        ,5]}},Length[a]-Length[NullSpace[a]]]
Output  3
Output  3

Input   (* Express a general vector in R^3 as a linear combination of the vectors
        {1,0,1}, {2,2,3} and {-1,-1,1}. First, verify the vectors are linearly
        independent by checking that their null space is empty: *)
        Subscript[b, 1] = {1, 0, 1};
        Subscript[b, 2] = {2, 2, 3};
        Subscript[b, 3] = {-1, -1, 1};
        NullSpace[{Subscript[b, 1], Subscript[b, 2], Subscript[b, 3]}]
```

```
             (* Form the matrix  whose columns are the basis vectors: *)
             p = Transpose[{Subscript[b, 1], Subscript[b, 2], Subscript[b, 3]}];

             (* The coefficients of a general vector  are given by TemplateBox[{p},
             Inverse].v: *)
             {Subscript[c, 1], Subscript[c, 2], Subscript[c, 3]} =  Inverse[p] . {x, y, z}

             (* Verify that  does indeed equal the linear combination : *)
             Simplify[
              Sum[
               Subscript[c, i] Subscript[b, i],
               {i, 1, 3}
               ]
              ]
Output       {}
Output       {x-y,-(x/5)+(2 y)/5+z/5,-((2 x)/5)-y/5+(2 z)/5}
Output       {x,y,z}

Input        (* Find the change-of-basis matrix that transforms coordinates with respect to
             the basis {bi} to coordinates with respect to the basis {di}: *)
             {Subscript[b,1],Subscript[b,2],Subscript[b,3]}={{2,5,-4},{1,0,3},{-3,3,-2}};
             {Subscript[d,1],Subscript[d,2],Subscript[d,3]}={{-2,4,1},{3,-4,-1},{3,3,-4}};

             (* The matrix p whose columns are the {bi}transforms from b-coordinates to
             standard coordinates: *)
             (p=Transpose[{Subscript[b,1],Subscript[b,2],Subscript[b,3]}])//MatrixForm

             (* The matrix q whose columns are the {di}transforms from d-coordinates to
             standard coordinates: *)
             (q=Transpose[{Subscript[d,1],Subscript[d,2],Subscript[d,3]}])//MatrixForm

             (* Its inverse converts from standard coordinates back to d-coordinates: *)
             Inverse[q]

             (* Therefore,q^-1.p converts from b-coordinates to d-coordinates: *)
             Inverse[q].p
Output       ({
               {2, 1, -3},
               {5, 0, 3},
               {-4, 3, -2}
              })
              ({
               {-2, 3, 3},
               {4, -4, 3},
               {1, -1, -4}
              })
Output       {{1,9/19,21/19},{1,5/19,18/19},{0,1/19,-(4/19)}}
Output       {{-(1/19),82/19,-(72/19)},{-(9/19),73/19,-(78/19)},{21/19,-(12/19),11/19}}
```

| | |
|---|---|
| `Eigenvalues[m]` | Eigenvalues m. |
| `Eigenvalues[m,k]` | gives the first k eigenvalues of m. |
| `Eigenvectors[m]` | Eigenvectors m. |
| `Eigenvectors[m,k]` | gives the first k eigenvectors of m. |
| `Eigensystem[m]` | gives a list {values,vectors} of the eigenvalues and eigenvectors of the square matrix m. |
| `Eigensystem[m,k]` | gives the eigenvalues and eigenvectors for the first k eigenvalues of m. |

***Mathematica Examples 3.3***

```
Input       r = Table[i+j+1, {i, 3}, {j, 3}]
Output      {{3, 4, 5}, {4, 5, 6}, {5, 6, 7}}

Input       Eigenvalues[r]
Output      {1/2 (15+√249),1/2 (15-√249), 0}

Input       Eigenvalues[{{a, b}, {c, d}}]
Output      {1/2 (a + d- √(a^2 + 4 b c- 2 a d + d^2 )), 1/2 (a+ d+√(a^2 + 4 b c - 2 a d +
            d^2 ))}

Input       (* Largest 5 eigenvalues: *)
            Eigenvalues[Table[N[1/(i + j + 1)], {i, 50}, {j, 50}], 5]
Output      {1.56835, 0.320271, 0.0506625, 0.00728517, 0.000970997}

Input       (* Three smallest eigenvalues *)
            Eigenvalues[Table[N[1/(i + j + 1), 20], {i, 10}, {j, 10}], -3]
Output      {1.14381991477649725×10^-10, 1.04594763489401510×10^-12, 4.43502120490729759×10^-15}

Input       Eigenvectors[{{a, b}, {c, d}}]
Output      {{-(-a+ d+ √(a^2 + 4 b c- 2 a d + d^2))/2 c, 1},{-(-a + d - √(a^2 + 4 b c - 2 a
            d + d^2))/ 2 c, 1}}

Input       Eigenvectors[{{1, 2, 3}, {4, 5, 6}, {7, 8, 9}}]
Output      {{-(-7-√33)/(2 (11 + 2 √33)), -(-29 - 5 √33)/(4 (11 + 2 √33)), 1}, {-(7 -
            √33)/(2 (-11+2 √33)), -(29 - 5 √33)/(4 (-11 + 2 √33)), 1}, {1, -2, 1}}

Input       N[%]
Output      {{0.283349, 0.641675, 1.}, {-1.28335, -0.141675, 1.}, {1., -2., 1.}}

Input       (* Eigenvectors computed using numerical methods *)
            Eigenvectors[N[{{1, 2, 3}, {4, 5, 6}, {7, 8, 9}}]]
Output      {{-0.231971, -0.525322, -0.818673}, {-0.78583, -0.0867513, 0.612328},
            {0.408248, -0.816497, 0.408248}}

Input[1]    m = Table[N[1/(i + j + 1)], {i, 10}, {j, 10}];
            (* Eigenvectors corresponding to the eigenvalues with smallest magnitude *)
Input[2]    vecs = Eigenvectors[m, -2]
Output      {{0.000199278, -0.0052208, 0.0482359, -0.213117, 0.487952, -0.525743,
            0.0417668, 0.478158, -0.437717, 0.12549}, {-0.0000188938, 0.000663109, -
            0.00844045, 0.054171, -0.200069, 0.451308, -0.63202, 0.536247, -0.252425,
            0.050584}}

Input       (* The corresponding eigenvalues *)
            Eigenvalues[m, -2]
Output      {1.045927766943346×10^-12, 4.450005850477068×10^-15}

Input       (* Find the eigenvectors corresponding to the 4 largest eigenvalues, or as many
            as there are if fewer *)
            Eigenvectors[Table[N[1/(i + j + 1)], {i, 3}, {j, 3}], UpTo[4]]
Output      {{-0.703153, -0.549268, -0.451532}, {-0.668535, 0.29444, 0.68291}, {0.242151, -
            0.782055, 0.574241}}

Input       Eigensystem[N[{{1, 2, 3}, {4, 5, 6}, {7, 8, 9}}]]
Output      {{16.1168, -1.11684, -1.30368×10-15}, {{-0.231971, -0.525322, -0.818673}, {-
            0.78583, -0.0867513, 0.612328}, {0.408248, -0.816497, 0.408248}}}

Input       Eigensystem[{{a, b}, {c, d}}]
```

```
Output        {{1/2 (a + d - √(a^2 + 4 b c - 2 a d + d^2)), 1/2 (a + d + √(a^2 + 4 b c - 2 a
              d + d^2))}, {{-(-a + d + √(a^2 + 4 b c - 2 a d + d^2))/2 c, 1}, {-(-a + d -
              √(a^2 + 4 b c - 2 a d + d^2))/2 c, 1}}}

Input         (* Find the 4 largest eigenvalues and their corresponding eigenvectors, or as
              many as there are if fewer *)
              Eigensystem[Table[N[1/(i + j + 1)], {i, 3}, {j, 3}], UpTo[4]]
Output        {{0.657051, 0.0189263, 0.000212737}, {{-0.703153, -0.549268, -0.451532}, {-
              0.668535, 0.29444, 0.68291}, {0.242151, -0.782055, 0.574241}}}

Input         (* Zero vectors are used when there are more eigenvalues than independent
              eigenvectors *)
              Eigenvectors[{{2, 1, 0}, {0, 2, 0}, {0, 0, 1}}]
Output        {1, 0, 0}, {0, 0, 0}, {0, 0, 1}}

Input         (*Eigenvectors with positive eigenvalues point in the same direction when acted
              on by the matrix:*)

              m={{1,2},{2,1}};
              {λ1,λ2}=Eigenvalues[m]
              {v1,v2}=Eigenvectors[m]
              Graphics[{{Thick,Arrow[{{0,0},v1}]},{Red,Arrow[{{0,0},m.v1}]}},Axes->True]

              (*Eigenvectors with negative eigenvalues point in the opposite direction when
              acted on by the matrix:*)
              Graphics[{{Thick,Arrow[{{0,0},v2}]},{Red,Arrow[{{0,0},m.v2}]}},Axes->True]
Output        {3,-1}
Output        {{1,1},{-1,1}}
Output
```

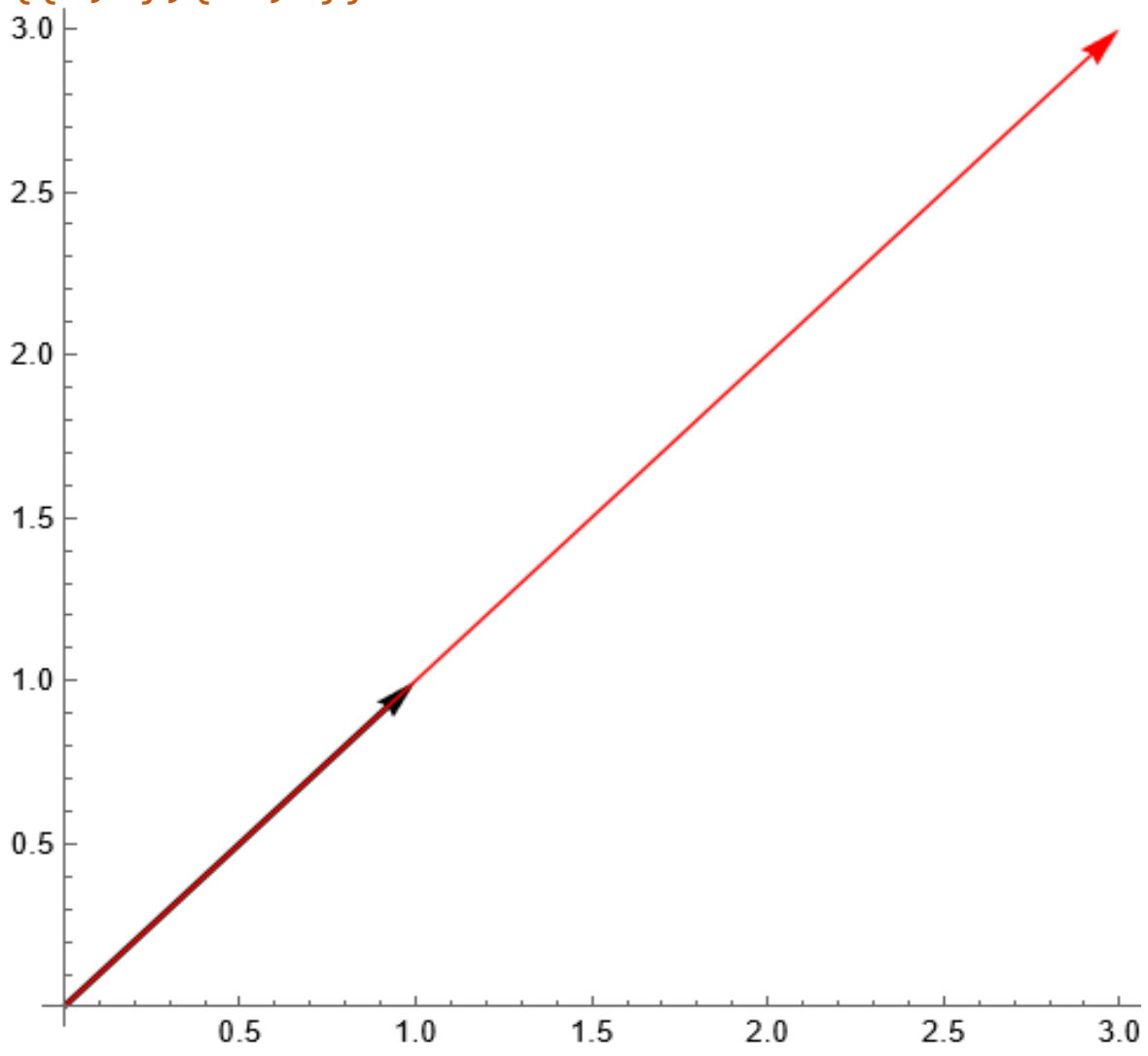

Output

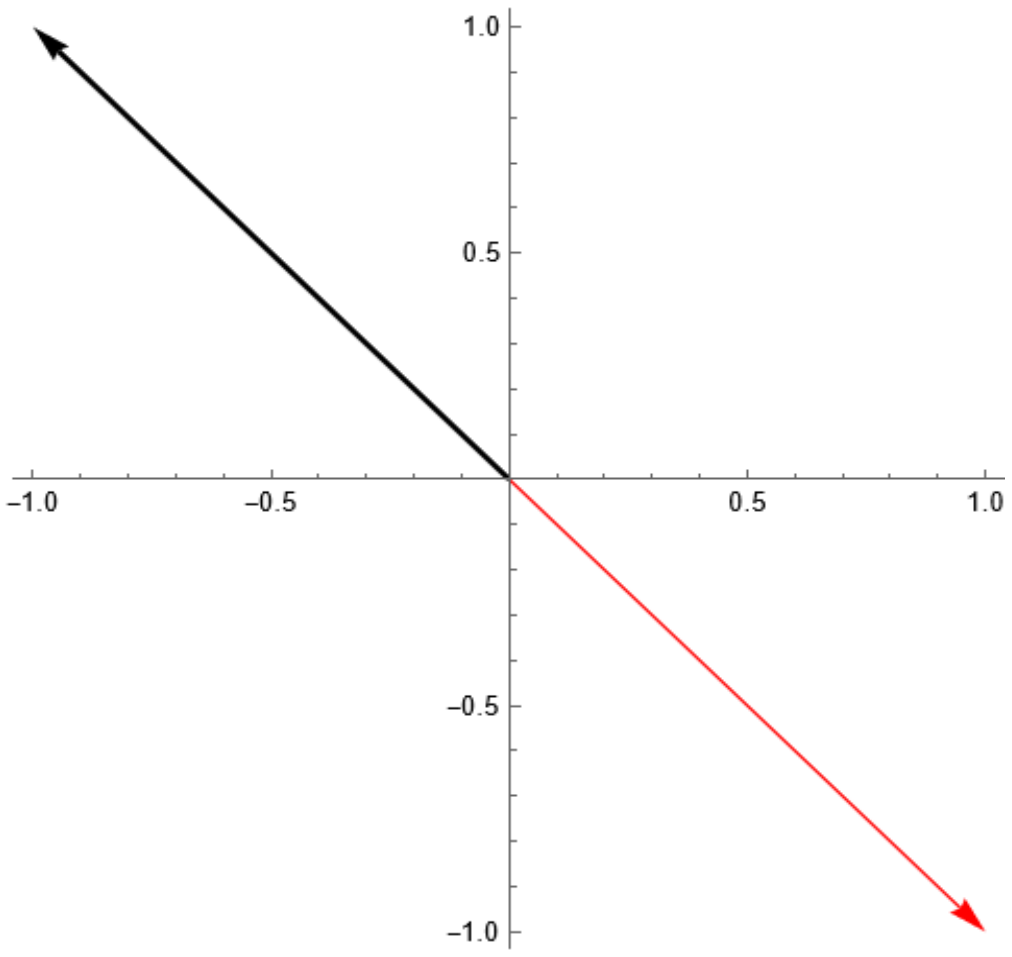

Input

```
(*Consider the following matrix and its associated quadratic form q=xT.a.x*)
a={{-2,2},{2,1}};
q={x,y}.a.{x,y}//Expand

(*The eigenvectors are the axes of the hyperbolas defined by q:*)
v=Eigenvectors[a]

(*The sign of the eigenvalue corresponds to the sign of the right-hand side of
the hyperbola equation:*)
λ=Eigenvalues[a]

ContourPlot[
 Table[
   q==n,
   {n,{-9,-4,-1,1,4,9}}
   ]//Evaluate,
 {x,-3,3},
 {y,-3,3},
 Epilog->(Arrow[{{0,0},#}]&/@v),
 PlotLegends->"Expressions"
 ]
```

Output -2 x2+4 x y+y2

Output {{-2,1},{1,2}}

Output {-3,2}

Output

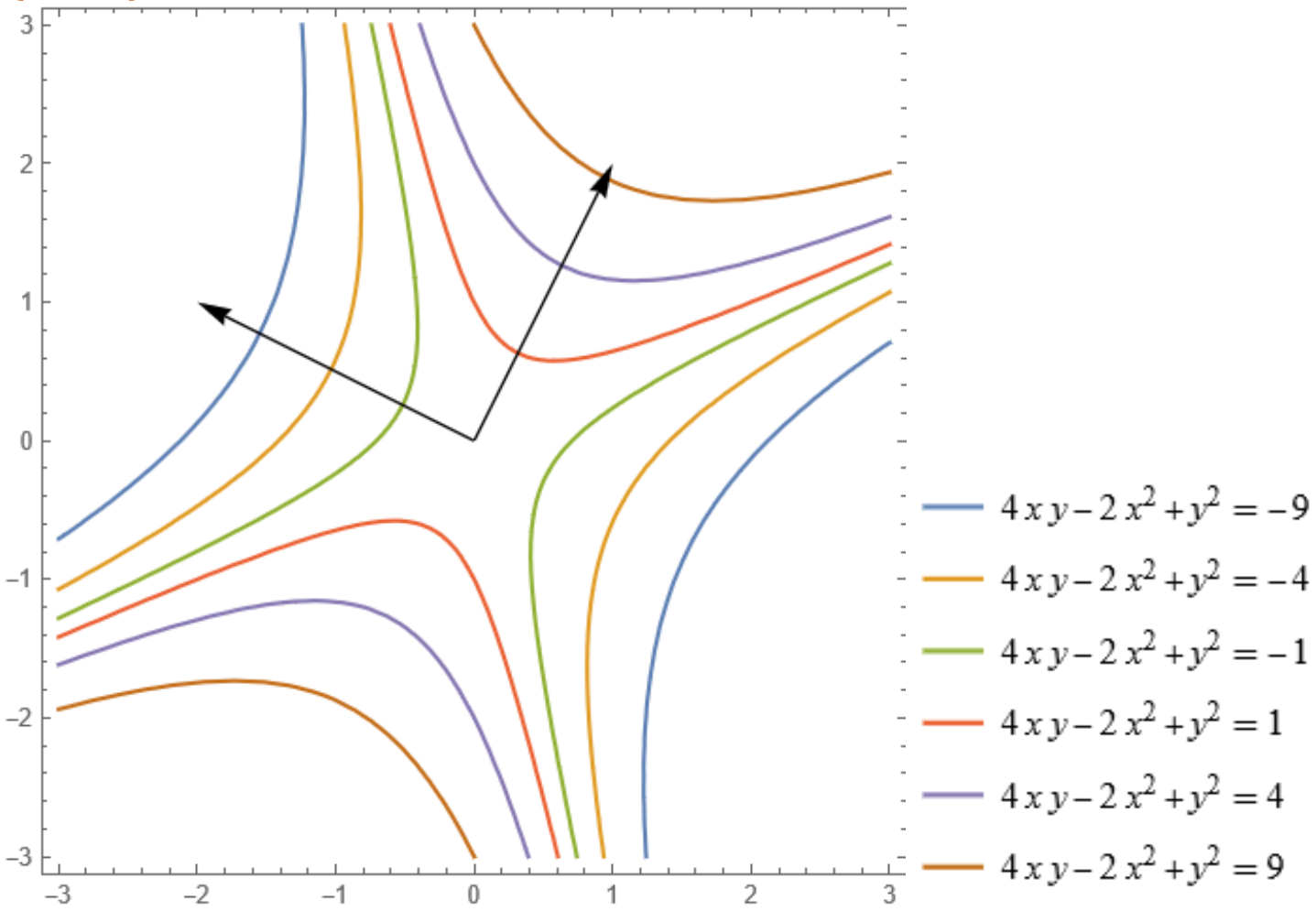

| | |
|---|---|
| `SingularValueList[m]` | gives a list of the non-zero singular values of a matrix m. |
| `SingularValueList[{m,a}]` | gives the generalized singular values of m with respect to a. |
| `SingularValueList[m,k]` | gives the k largest singular values of m. |
| `SingularValueList[{m,a},k]` | gives the k largest generalized singular values of m. |
| `SingularValueDecomposition[m]` | gives the singular value decomposition for a numerical matrix m as a list of matrices {u,σ,v}, where σ is a diagonal matrix and m can be written as u.σ.ConjugateTranspose[v]. |
| `SingularValueDecomposition[{m,a}]` | gives the generalized singular value decomposition of m with respect to a. |
| `SingularValueDecomposition[m,k]` | gives the singular value decomposition associated with the k largest singular values of m. |
| `SingularValueDecomposition[m,UpTo[k]]` | gives the decomposition for the k largest singular values, or as many as are available. |

The singular values of a matrix m are the square roots of the eigenvalues of m . m*, where * denotes Hermitian transpose. The number of such singular values is the smaller dimension of the matrix. `SingularValueList` sorts the singular values from largest to smallest. Very small singular values are usually numerically meaningless. With the option setting `Tolerance->t`, `SingularValueList` drops singular values that are less than a fraction `t` of the largest singular value [8,9].

- Repeated singular values appear with their appropriate multiplicity.
- The matrix m can be rectangular; the total number of singular values is always Min[Dimensions[m]].
- The singular values can be obtained from `Sqrt[Eigenvalues[ConjugateTranspose[m].m]]`.
- For `SingularValueDecomposition,` the diagonal elements of σ are the singular values of m. u and v are column orthonormal matrices, whose transposes can be considered as lists of orthonormal vectors.
- `SingularValueDecomposition` sets to zero any singular values that would be dropped by `SingularValueList`.

***Mathematica Examples 3.4***

```
Input     (* Compute the singular values of an invertible matrix: *)
          SingularValueList[({{1,2},{3,4}})]
Output    {√(15+√221),√(16-√221)}

Input     (* Compute the nonzero singular values of a singular matrix: *)
          SingularValueList[({
             {1, 2, 3},
             {4, 5, 6},
             {7, 8, 9}
            })]
Output    {√(3(95+√8881)/2)), √(3(95-√8881)/2))}

Input     (* Compute a singular value decomposition: *)
          {u, σ, v} = SingularValueDecomposition[({
              {1, 2},
              {1, 2}
             })]
Output    {{{1/√2,-(1/√2)},{1/√2,1/√}},{{√10,0},{0,0}},{{1/√5,-(2/√5)},{2/√5,1/√5}}}

Input     (* Reconstruct the input matrix: *)
          u.σ.Transpose[v]//MatrixForm
Output    ({
            {1, 2},
            {1, 2}
           })

Input     (* Compute a singular value decomposition for an invertible matrix: *)
```

```
              {u, σ, v} = SingularValueDecomposition[({
                  {0.5, 1},
                  {2, 2.5}
                 })]
              (* The matrix of singular values is also invertible: *)
              σ//MatrixForm
              (* Reconstruct the input matrix: *)
              u.σ.Transpose[v]
Output        {{{-0.324536,0.945873},{-0.945873,-0.324536}},{{3.38391,0.},{0.,0.221637}},{{-
              0.606994,-0.794707},{-0.794707,0.606994}}}
Output        ({
                {3.38391, 0.},
                {0., 0.221637}
               })
Output         {{0.5,1.},{2.,2.5}}
```

| | |
|---|---|
| `PositiveDefiniteMatrixQ[m]` | gives True if m is explicitly positive definite, and False otherwise. |
| `NegativeDefiniteMatrixQ[m]` | gives True if m is explicitly negative definite, and False otherwise. |
| `PositiveSemidefiniteMatrixQ[m]` | gives True if m is explicitly positive semidefinite, and False otherwise. |
| `NegativeSemidefiniteMatrixQ[m]` | gives True if m is explicitly negative semidefinite, and False otherwise. |
| `IndefiniteMatrixQ[m]` | gives True if m is explicitly indefinite, and False otherwise. |

***Mathematica Examples 3.5***

```
Input         (*Test if a matrix is explicitly positive definite:*)
              m={{2,3},{4,8}};
              PositiveDefiniteMatrixQ[m]
              (*This means that the quadratic form x.m.x>0 for all vectors x!=0:*)
              Plot3D[
               {x,y}.m.{x,y},
               {x,-1,1},
               {y,-1,1},
               LabelStyle->Directive[Black,20]
               ]
Output        True
Output
```

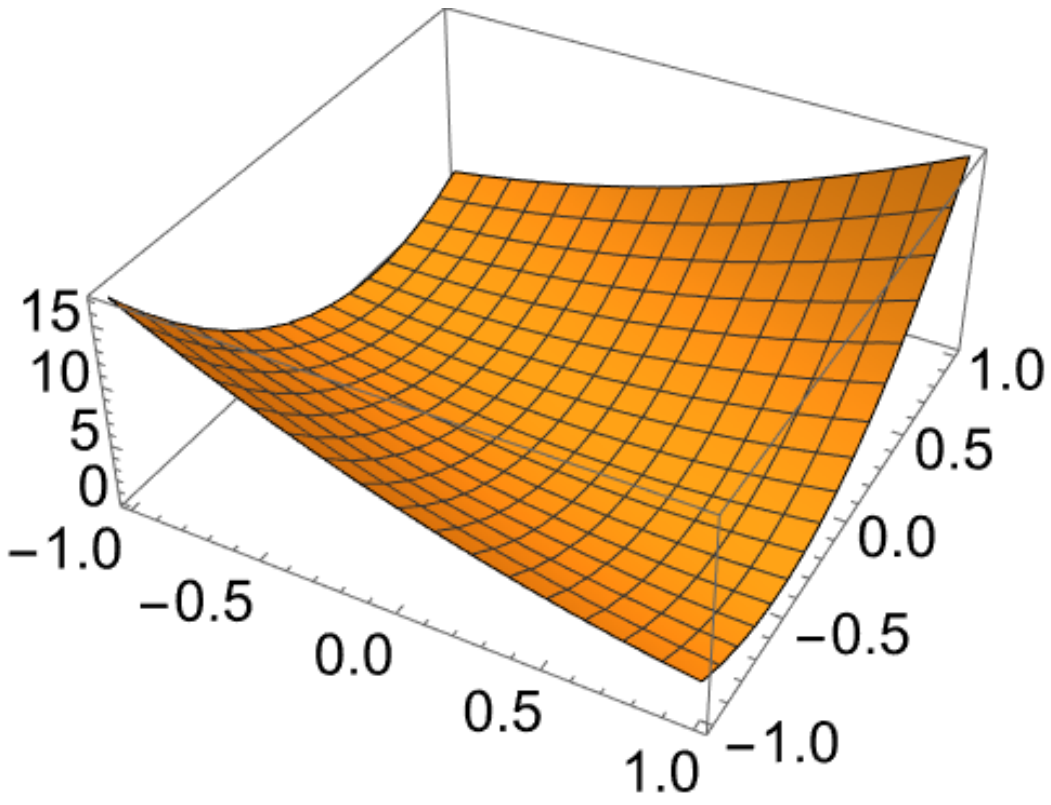

```
Input         (*Test if a matrix is explicitly negative definite:*)
              m={{-1,1},{1,-2}};
              NegativeDefiniteMatrixQ[m]

              (*This means that the quadratic form x.m.x<0 for all vectorsx!=0:*)
              Plot3D[
               {x,y}.m.{x,y},
               {x,-1,1},
               {y,-1,1},
```

```
               LabelStyle->Directive[Black,20]
               ]
Output        True
Output
```

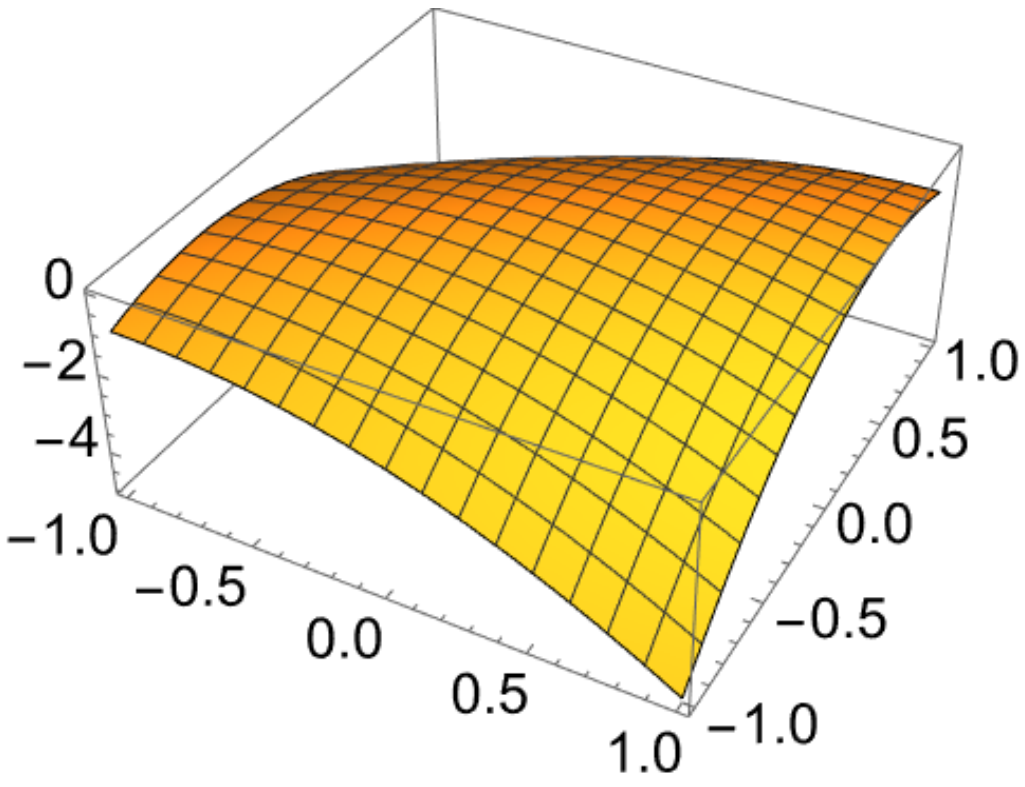

```
Input         (*Test if a matrix is explicitly positive semidefinite:*)
              m={{2,0},{0,0}};
              PositiveSemidefiniteMatrixQ[m]

              (*This means that the quadratic form for all vectors:*)
              Plot3D[
               {x,y}.m.{x,y},
               {x,-1,1},
               {y,-1,1},
               LabelStyle->Directive[Black,20]
               ]
Output        True
Output
```

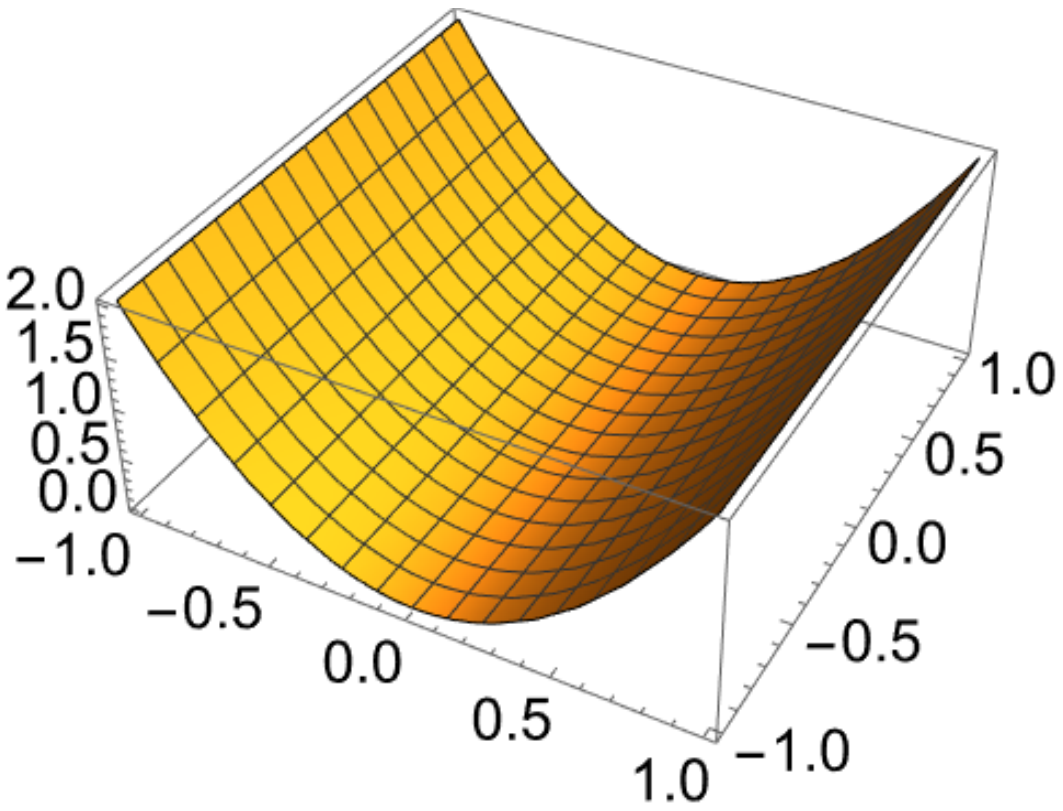

```
Input         (*Test if a matrix is explicitly negative semidefinite:*)
              m={{-2,0},{0,0}};
              NegativeSemidefiniteMatrixQ[m]

              (*This means that the quadratic form x.m.x<=0 for all vectorsx!=0:*)
              Plot3D[
               {x,y}.m.{x,y},
               {x,-1,1},
               {y,-1,1},
               LabelStyle->Directive[Black,20]
               ]
Output        True
```

Output

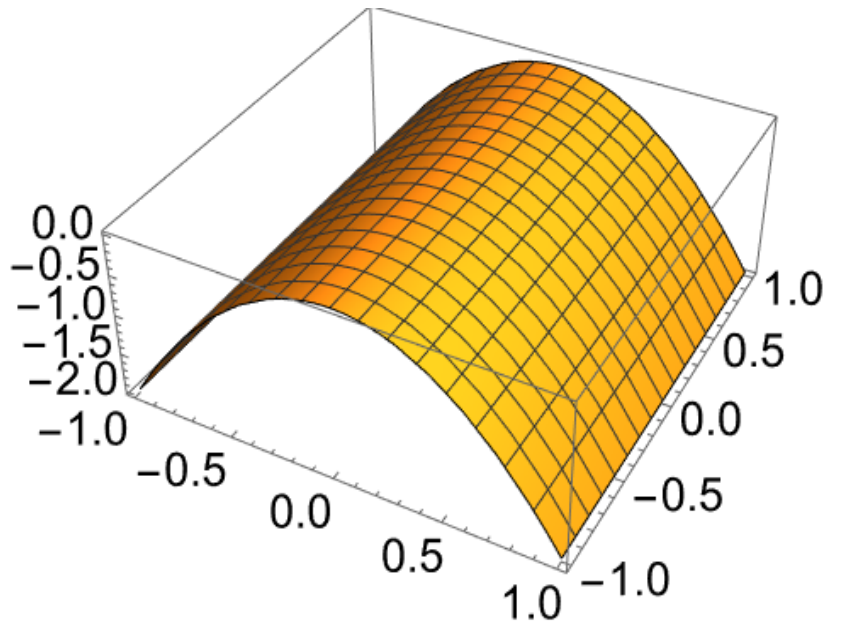

Input

```
(*m matrix is explicitly positive definite:*)
m={{2,3},{4,8}};
(*This means that the quadratic form x.m.x>0 for all vectors x!=0:*)
(*3d + contour plot*)

plot1=Plot3D[
   {x,y}.m.{x,y},
   {x,-1,1},
   {y,-1,1},
   ClippingStyle->None,
   MeshFunctions->{#3&},
   Mesh->15,
   MeshStyle->Opacity[.5],
   MeshShading->{{Opacity[.3],Blue},{Opacity[.8],Orange}},
   Lighting->"Neutral",
   LabelStyle->Directive[Black,20]
   ];
slice=SliceContourPlot3D[
   {x,y}.m.{x,y},
   z==0,
   {x,-1,1},
   {y,-1,1},
   {z,-1,1},
   Contours->15,
   Axes->False,
   PlotPoints->50,
   PlotRangePadding->0,
   ColorFunction->"Rainbow"
   ];
Show[
 plot1,
 slice,
 PlotRange->All,
 BoxRatios->{1,1,.6},
 FaceGrids->{Back,Left}
 ]
```

Output

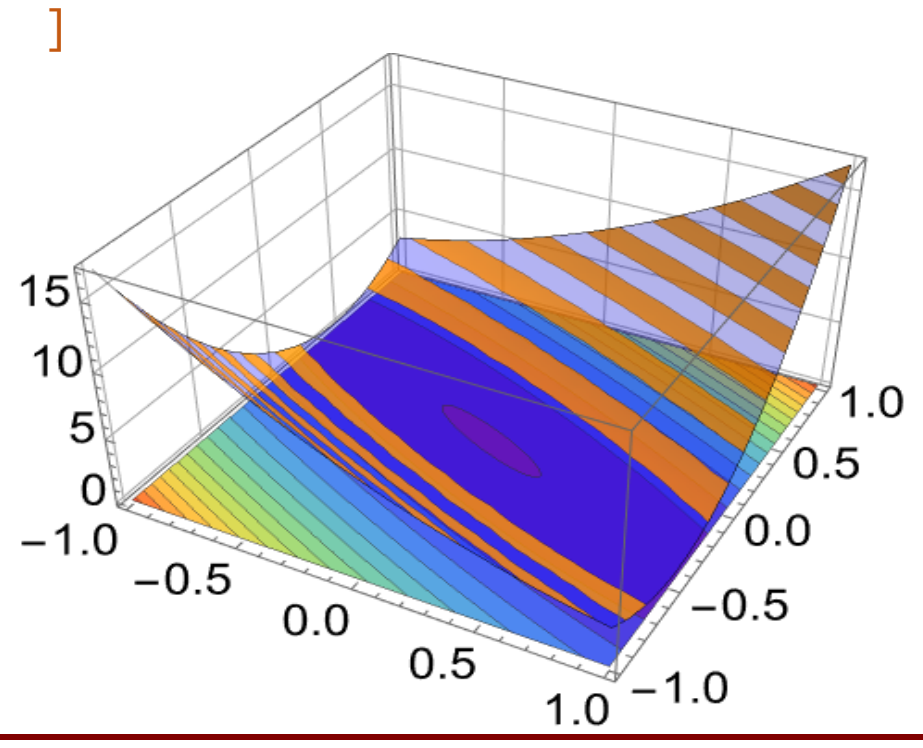

# CHAPTER 4

# CONVEX STES AND CONVEX FUNCTIONS

## 4.1 $p$-norm in finite dimensions

The length of a vector $|\mathbf{x}\rangle = (x_1, x_2, \dots, x_n)^T$ in the $n$-dimensional real vector space $\mathbb{R}^n$ is usually given by the Euclidean norm:

$$\|\mathbf{x}\| = (x_1^2 + x_2^2 + \cdots + x_n^2)^{\frac{1}{2}}. \tag{4.1}$$

It follows directly from the definition that $\|\mathbf{x}\| = \sqrt{\langle \mathbf{x}|\mathbf{x}\rangle}$.

The formula for the Euclidean distance between two points $|\mathbf{p}_1\rangle = (x_1, y_1, z_1)^T$ and $|\mathbf{p}_2\rangle = (x_2, y_2, z_2)^T$ in 3-dimensional space is

$$d(\mathbf{p}_1, \mathbf{p}_2) = ((x_1 - x_2)^2 + (y_1 - y_2)^2 + (z_1 - z_2)^2)^{\frac{1}{2}}. \tag{4.2}$$

The formula above, however, is not always the way we measure the distance between points laying in a 3-dimensional space. For example, we can think of the Earth as being a subset of three-dimensional space and the surface of the Earth as being the surface of a sphere. When measuring the distance between points on the surface of the Earth, say from New York City to Cairo. We do not draw a line through the two points and measure the distance along those points. This distance would correspond to the distance traveled if a tunnel was dug in a straight line through the Earth, starting at New York and ending at Cairo. This is not very practical; one usually does not travel between cities by tunneling from one to the other. Instead, we measure the distance between two points $|\mathbf{p}_1\rangle$ and $|\mathbf{p}_2\rangle$ on the Earth by measuring the short path starting at $|\mathbf{p}_1\rangle$ and ending at $|\mathbf{p}_2\rangle$ that stays on the surface of the Earth. In many situations, the Euclidean distance is insufficient for capturing the actual distances in a given space. An analogy to this is suggested by taxi drivers in a grid street plan who should measure distance not in terms of the length of the straight line to their destination but in terms of the rectilinear distance, which takes into account that streets are either orthogonal or parallel to each other (see Figure 4.1). The class of $p$-norms generalizes these two examples and has an abundance of applications in many parts of mathematics, physics, and computer science.

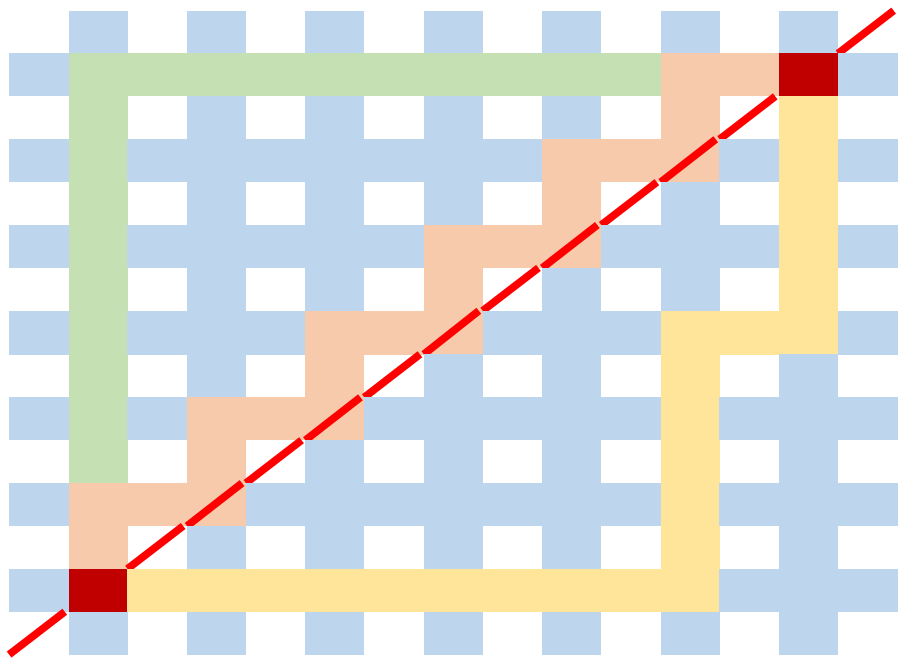

**Figure 4.1.** Taxicab or Manhattan distance versus Euclidean distance.

There are many useful measures of length (many different norms).

**Definition ($p$-Norm):** For a real number $p \geq 1$, the $p$-norm or $\ell^p$-norm of $|\mathbf{x}\rangle = (x_1, x_2, \dots, x_n)^T$ is defined by

$$\|\mathbf{x}\|_p = (|x_1|^p + |x_2|^p + \cdots + |x_n|^p)^{\frac{1}{p}}. \tag{4.3}$$

Some norms are

- $\ell^1$-norm, 1-norm, or the Manhattan norm: $\|\mathbf{x}\|_1 = |x_1| + |x_2| + \cdots + |x_n|$.
- $\ell^2$-norm, 2-norm, or the Euclidean norm: $\|\mathbf{x}\|_2 = \sqrt{|x_1|^2 + |x_2|^2 + \cdots + |x_n|^2}$,
- $\ell^\infty$-norm, $\infty$-norm or the Chebyshev norm: $\|\mathbf{x}\|_\infty = \max\{|x_1|, |x_2|, \dots |x_n|\}$.

Figure 4.2 gives a visual representation of these norms.

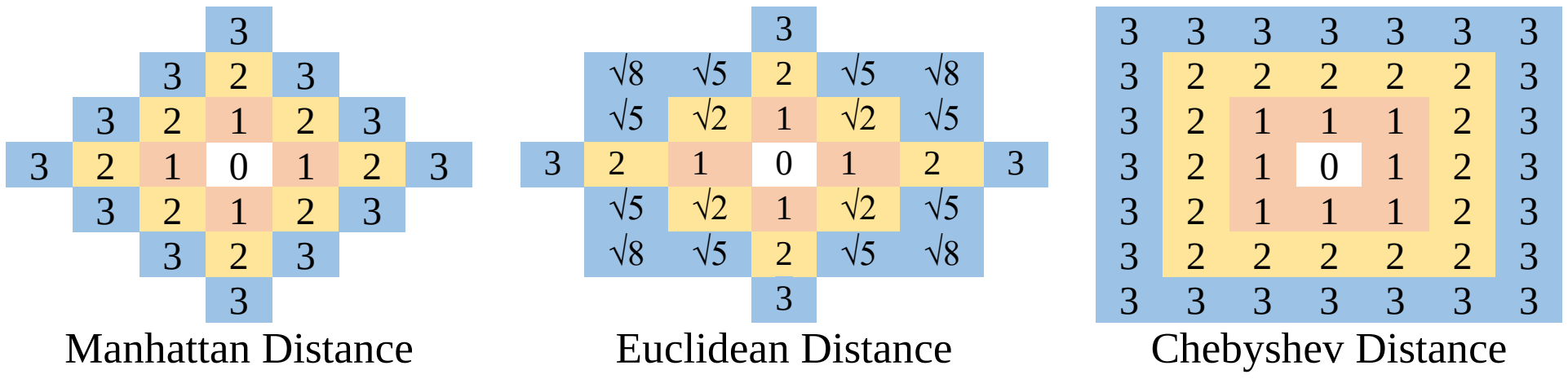

**Figure 4.2.** Manhattan, Euclidean, and Chebyshev Distance.

**Theorem 4.1:** $\|\mathbf{x}\|_\infty \le \|\mathbf{x}\|_2 \le \|\mathbf{x}\|_1$.

**Proof:**

$$\begin{aligned}\|\mathbf{x}\|_\infty &= \max\{|x_1|, |x_2|, \dots |x_n|\} \\ &= \max\left\{\sqrt{x_1^2}, \sqrt{x_2^2}, \dots, \sqrt{x_n^2}\right\} \\ &= \sqrt{x_k^2} \text{ for some } k \\ &\le \sqrt{x_1^2 + x_2^2 + \cdots + x_n^2} \\ &= \|\mathbf{x}\|_2 \\ &= \sqrt{|x_1|^2 + |x_2|^2 + \cdots + |x_n|^2} \\ &\le \sqrt{(|x_1| + |x_2| + \cdots + |x_n|)^2} \\ &= \|\mathbf{x}\|_1.\end{aligned}$$

■

## 4.2 Hyperplanes, Hyperspheres and Hyperellipsoid

### 4.2.1 Lines, Planes, Hyperplanes, and Half-Spaces in $\mathbb{R}^n$

The equation of a line means an equation in the variables $x$ and $y$ whose solution set is a line in the $xy$-plane. The most popular form in algebra is the (slope-intercept form)

$$y = mx + c. \tag{4.4}$$

This in effect uses $x$ as a parameter and writes $y$ as a function of $x$: $y = f(x) = mx + c$. When $x = 0$, and $y = c$, the point $(0, c)^T$ is the intersection of the line with the $y$-axis. Thinking of a line as a geometrical object and not the graph of a function, it makes sense to treat $x$ and $y$ more evenhandedly. The general equation for a line is

$$ax + by + d = 0, \tag{4.5}$$

with the condition that at least one of $a$ or $b$ is nonzero. This can easily be converted to slope-intercept form by solving for $y$:

$$y = -\frac{a}{b}x - \frac{d}{b}, \tag{4.6}$$

except for the special case $b = 0$, when the line is parallel to the $y$-axis.

If we set $x = t$ , $-\infty < t < \infty$ , then the solutions to (4.5) are

$$\begin{aligned} |\mathbf{y}\rangle &= (x, y)^T \\ &= \left(t, -\frac{a}{b}t - \frac{d}{b}\right)^T \\ &= t\left(1, -\frac{a}{b}\right)^T + \left(0, -\frac{d}{b}\right)^T \\ &= t|\mathbf{v}\rangle + |\mathbf{u}\rangle. \end{aligned} \tag{4.7}$$

Equation (4.7) is the line $L$ through $|\mathbf{u}\rangle$ in the direction of $|\mathbf{v}\rangle$. Equation (4.7) is called a vector equation for the line. Now let $|\mathbf{n}\rangle = (a, b)^T$, (4.5) becomes

$$\langle \mathbf{n}|\mathbf{y}\rangle + d = 0. \tag{4.8}$$

Moreover, if $|\mathbf{p}\rangle = (p_1, p_2)^T$ is a point on $L$, then

$$\langle \mathbf{n}|\mathbf{p}\rangle + d = 0, \tag{4.9}$$

which implies that $d = -\langle \mathbf{n}|\mathbf{p}\rangle$. Thus, we may write (4.8) as

$$\langle \mathbf{n}|\mathbf{y}\rangle - \langle \mathbf{n}|\mathbf{p}\rangle = 0, \tag{4.10}$$

and so we see that (4.10) is equivalent to the equation

$$\langle \mathbf{n}|\mathbf{y} - \mathbf{p}\rangle = 0. \tag{4.11}$$

Equation (4.11) is a normal equation for the line $L$ and $|\mathbf{n}\rangle$ is a normal vector for $L$. In words, (4.11) says that line $L$ consists of all points in $\mathbb{R}^2$ whose difference with $|\mathbf{p}\rangle$ is orthogonal to $|\mathbf{n}\rangle$. See Figure 4.3.

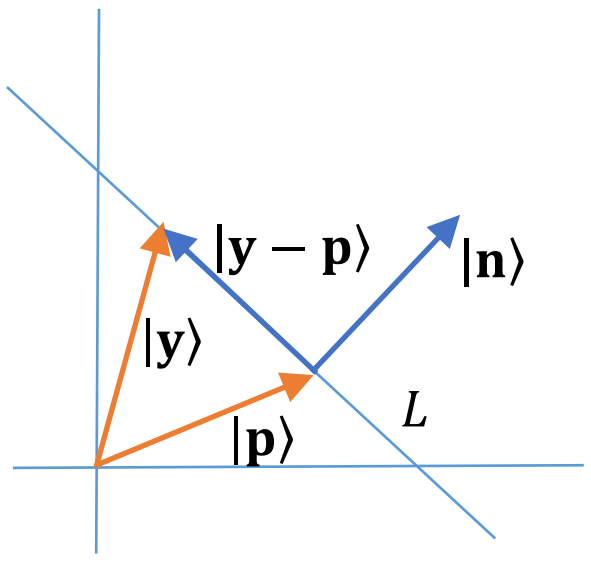

**Figure 4.3.** $L$ is the set of points $\mathbf{y}$ for which $|\mathbf{y} - \mathbf{p}\rangle$ is orthogonal to $|\mathbf{n}\rangle$.

Now, let us expand (4.11) to 3-dimensional space. A plane in $\mathbb{R}^3$ has the equation

$$ax + by + cz + d = 0, \tag{4.12}$$

where at least one of the numbers $a$, $b$, $c$ must be nonzero. If $c$ is not zero, it is often useful to think of the plane as the graph of a function $z$ of the variables $x$ and $y$. The equation can be rearranged like this:

$$z = -\left(\frac{a}{c}\right)x + \left(-\frac{b}{c}\right)y + \left(-\frac{d}{c}\right). \tag{4.13}$$

The equation of a plane (see Figure 4.4) with nonzero normal vector $|\mathbf{n}\rangle = (a, b, c)^T$ through the point $|\mathbf{x_0}\rangle = (x_0, y_0, z_0)^T$ is

$$\langle \mathbf{n}|\mathbf{x} - \mathbf{x_0}\rangle = 0, \tag{4.14}$$

where $|\mathbf{x}\rangle = (x, y, z)^T$. Plugging in gives the general equation of a plane,

$$ax + by + cz = -d, \tag{4.15}$$

where,

$$d = -ax_0 - by_0 - cz_0. \tag{4.16}$$

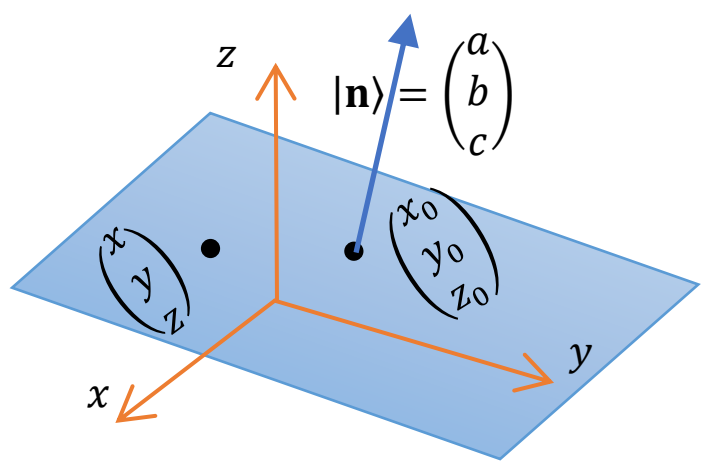

**Figure 4.4.** A plane with nonzero normal vector $|\mathbf{n}\rangle$ through the point $|\mathbf{x_0}\rangle = (x_0, y_0, z_0)^T$

A plane specified in this form, therefore, has $x$-, $y$-, and $z$-intercepts at

$$x = -\frac{d}{a}, y = -\frac{d}{b}, z = -\frac{d}{c}. \tag{4.17}$$

Now consider the case where $P$ is the set of all points $|\mathbf{y}\rangle = (x, y, z)^T$ in $\mathbb{R}^3$ that satisfies the equation

$$ax + by + cz + d = 0, \tag{4.18}$$

where $a$, $b$, $c$, and $d$ are scalars with at least one of $a$, $b$, and $c$ not being 0. If, for example, $a \neq 0$, then we may solve for $x$ to obtain

$$x = -\frac{b}{a}y - \frac{c}{a}z - \frac{d}{a}. \tag{4.19}$$

If we set $y = t$, $-\infty < t < \infty$, and $z = s$, $-\infty < s < \infty$, the solutions to (4.19) are

$$\begin{aligned} |\mathbf{y}\rangle &= (x, y, z)^T \\ &= \left(-\frac{b}{a}t - \frac{c}{a}s - \frac{d}{a}, t, s\right)^T \\ &= t(-\frac{b}{a}, 1, 0)^T + s(-\frac{c}{a}, 0, 1)^T + (-\frac{d}{a}, 0, 0)^T. \end{aligned} \tag{4.20}$$

Thus, we see that $P$ is a plane in $\mathbb{R}^3$. In analogy with the case of lines in $\mathbb{R}^2$, if we let $|\mathbf{n}\rangle = (a, b, c)^T$ and let $|\mathbf{p}\rangle = (p_1, p_2, p_3)^T$ be a point on $P$, then we have

$$\langle \mathbf{n}|\mathbf{y}\rangle + d = ax + by + cz + d = 0, \tag{4.21}$$

from which we see that $\langle \mathbf{n}|\mathbf{p}\rangle = -d$, and so we may write (4.18) as

$$\langle \mathbf{n}|\mathbf{y} - \mathbf{p}\rangle = 0. \tag{4.22}$$

We call (4.21) a normal equation for $P$ and we call $|\mathbf{n}\rangle$ a normal vector for $P$. In words, (4.22) says that the plane $P$ consists of all points in $\mathbb{R}^3$ whose difference with $|\mathbf{p}\rangle$ is orthogonal to $|\mathbf{n}\rangle$. See Figure 4.5.

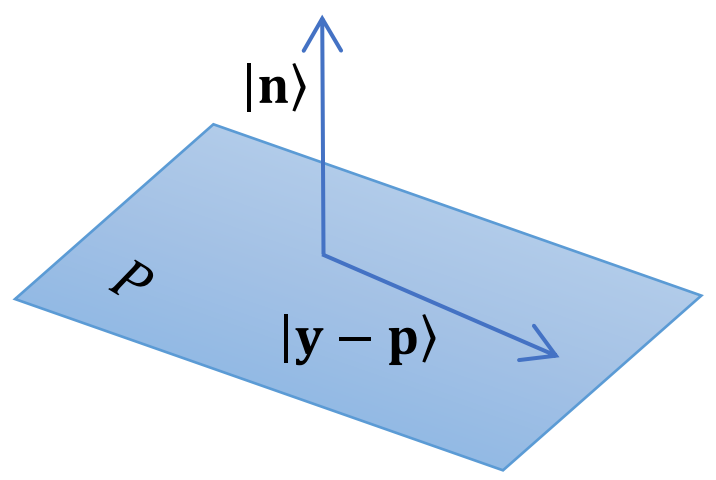

**Figure 4.5.** The plane $P$ is the set of points $\mathbf{y}$ for which $|\mathbf{y} - \mathbf{p}\rangle$ is orthogonal to $|\mathbf{n}\rangle$.

**Definition (Hyperplane):** Suppose $|\mathbf{n}\rangle$ and $|\mathbf{p}\rangle$ are vectors in $\mathbb{R}^n$ with $|\mathbf{n}\rangle \neq |\mathbf{0}\rangle$. The set of all vectors $|\mathbf{y}\rangle$ in $\mathbb{R}^n$ which satisfies the equation

$$\langle \mathbf{n}|\mathbf{y} - \mathbf{p}\rangle = 0, \tag{4.23}$$

or

$$\langle \mathbf{n}|\mathbf{y}\rangle = b, \quad \text{where } b = \langle \mathbf{n}|\mathbf{p}\rangle. \tag{4.24}$$

is called a hyperplane through the point $|\mathbf{p}\rangle$. We call $|\mathbf{n}\rangle$ a normal vector for the hyperplane, and we call (4.23) a normal equation for the hyperplane.
In this terminology, a line in $\mathbb{R}^2$ is a hyperplane and a plane in $\mathbb{R}^3$ is a hyperplane.

**Remarks**

1- If we let $|\mathbf{n}\rangle = (a_1, a_2, \dots, a_n)^T$, $|\mathbf{p}\rangle = (p_1, p_2, \dots, p_n)^T$, and $|\mathbf{y}\rangle = (y_1, y_2, \dots, y_n)^T$, then we may write (4.23) as

$$a_1(y_1 - p_1) + a_2(y_2 - p_2) + \dots + a_n(y_n - p_n) = 0, \tag{4.25}$$

or

$$a_1 y_1 + a_2 y_2 + \dots + a_n y_n + d = 0, \tag{4.26}$$

where $d = -\langle \mathbf{n}|\mathbf{p}\rangle$.

2- The normal equation description of a hyperplane simplifies a number of geometric calculations. For example, given a hyperplane $H$ through $|\mathbf{p}\rangle$ with normal vector $|\mathbf{n}\rangle$ and a point $|\mathbf{q}\rangle$ in $\mathbb{R}^n$, the distance from $|\mathbf{q}\rangle$ to $H$ is simply the length of the projection of $|\mathbf{q} - \mathbf{p}\rangle$ onto $|\mathbf{n}\rangle$. Thus if $|\mathbf{u}\rangle$ is the unit vector in the direction of $|\mathbf{n}\rangle$, then the distance from $|\mathbf{q}\rangle$ to $H$ is $|\langle \mathbf{q} - \mathbf{p}|\mathbf{u}\rangle|$, (see Figure 4.6). Moreover, if we let $d = -\langle \mathbf{p}|\mathbf{n}\rangle$ as in (4.26), then we have

$$|\langle \mathbf{q} - \mathbf{p}|\mathbf{u}\rangle| = |\langle \mathbf{q}|\mathbf{u}\rangle - \langle \mathbf{p}|\mathbf{u}\rangle| = \frac{\langle \mathbf{q}|\mathbf{n}\rangle - \langle \mathbf{p}|\mathbf{n}\rangle}{\| \mathbf{n} \|} = \frac{|\langle \mathbf{q}|\mathbf{n}\rangle + d|}{\| \mathbf{n} \|}. \tag{4.27}$$

In particular, (4.27) may be used to find the distance from a point to a line in $\mathbb{R}^2$ and from a point to a plane in $\mathbb{R}^3$.

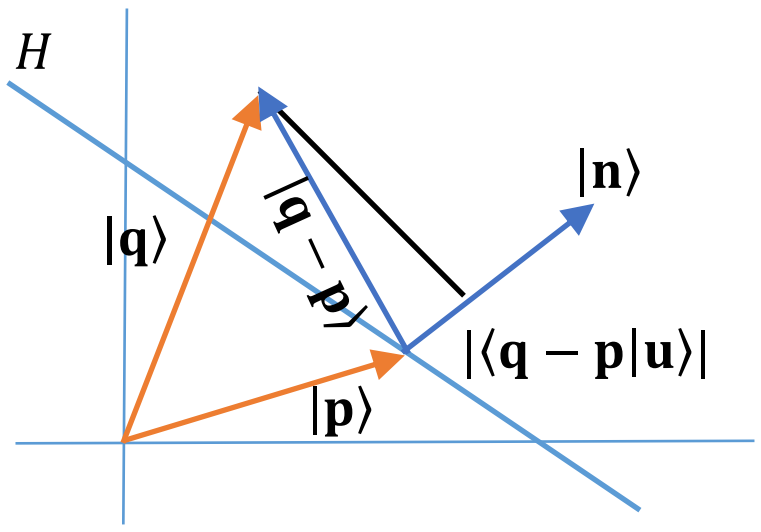

**Figure 4.6.** Distance from a point $|\mathbf{q}\rangle$ to a hyperplane $H$.

3- Geometrically, a hyperplane $H = \{|\mathbf{x}\rangle : \langle \mathbf{a}|\mathbf{x}\rangle = b\}$, with $\|\mathbf{a}\|_2=1$, is a translation of the set of vectors orthogonal to $|\mathbf{a}\rangle$. The direction of the translation is determined by $|\mathbf{a}\rangle$, and the amount by $b$. Precisely, $|b|$ is the length of the closest point $|\mathbf{x_0}\rangle$ on $H$ from the origin, and the sign of $b$ determines if $H$ is away from the origin along the direction $|\mathbf{a}\rangle$ or $-|\mathbf{a}\rangle$. As we increase the magnitude of $b$, the hyperplane is shifting further away along $\pm|\mathbf{a}\rangle$, depending on the sign of $b$. In Figure 4.7, the scalar $b$ is positive, as $|\mathbf{x_0}\rangle$ and $|\mathbf{a}\rangle$ point to the same direction.

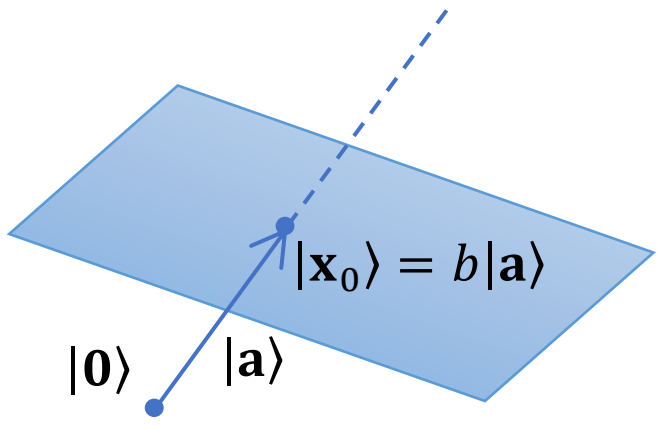

**Figure 4.7.** The geometry of hyperplanes.

A half-space is a subset of $\mathbb{R}^n$ defined by a single inequality involving a scalar product.

**Definition (Half-Space):** A half-space in $\mathbb{R}^n$ is a set of the form

$$H = \{|\mathbf{x}\rangle : \langle \mathbf{a}|\mathbf{x}\rangle \geq b\}, \tag{4.28}$$

where $|\mathbf{a}\rangle \in \mathbb{R}^n$, $|\mathbf{a}\rangle \neq |\mathbf{0}\rangle$, and $b \in \mathbb{R}$ are given.

Geometrically, in Figure 4.8, the half-space above is the set of points such that $\langle \mathbf{a}|\mathbf{x}-\mathbf{x_0}\rangle \geq 0$, that is, the angle between $|\mathbf{x}-\mathbf{x_0}\rangle$ and $|\mathbf{a}\rangle$ is acute (in $[-90, +90]$). Here $|\mathbf{x_0}\rangle$ is the point closest to the origin on the hyperplane defined by the equality $\langle \mathbf{a}|\mathbf{x}\rangle = b$. (When $|\mathbf{a}\rangle$ is normalized, as in the picture, $|\mathbf{x_0}\rangle = b|\mathbf{a}\rangle$.)

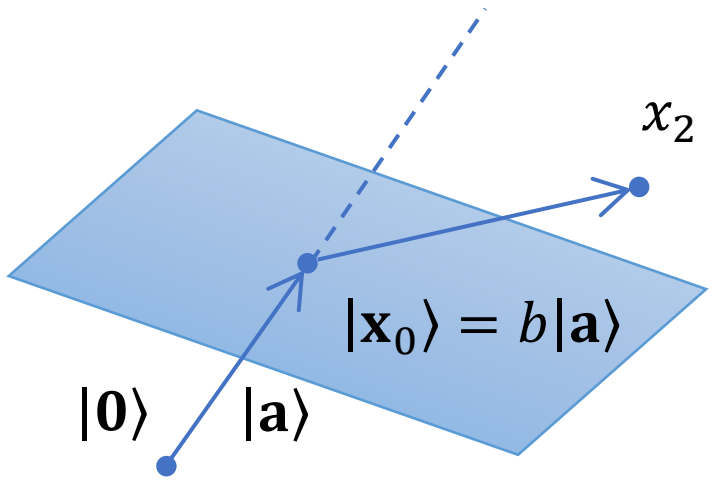

**Figure 4.8.** The half-space $\{|\mathbf{x}\rangle : \langle \mathbf{a}|\mathbf{x}\rangle \geq b\}$ is the set of points such that $|\mathbf{x}-\mathbf{x_0}\rangle$ forms an acute angle with $|\mathbf{a}\rangle$, where $|\mathbf{x_0}\rangle$ is the projection of the origin on the boundary of the half-space.

### 4.2.2 Circle, Spheres, Hyperspheres, and Balls in $\mathbb{R}^n$

In two dimensions, the equation for a circle of radius $r$ is

$$x^2 + y^2 = r^2. \tag{4.29}$$

In three dimensions, the equation for a 3-dimensional sphere (2-sphere) is

$$x^2 + y^2 + z^2 = r^2. \tag{4.30}$$

In hyperspace (space with more than three dimensions), the equation for an $n$-dimensional sphere ($(n-1)$-sphere) becomes

$$x_1^2 + x_2^2 + x_3^2 + \ldots + x_n^2 = r^2. \tag{4.31}$$

An $n$-sphere (a hypersphere) is the set of points in $(n+1)$-dimensional Euclidean space that are situated at a constant distance $r$ from a fixed point, called the center. It is the generalization of an ordinary sphere in ordinary three-dimensional space. The "radius" of a sphere is the constant distance of its points to the center. When the sphere has a unit radius, it is usual to call it the unit $n$-sphere or (the $n$-sphere).

**Definition 1 ($n$-sphere):** In terms of the standard norm, the $n$-sphere is defined as

$$S^n = \{|\mathbf{x}\rangle \in \mathbb{R}^{n+1} : \|\mathbf{x}\|_2 = 1\}, \tag{4.32}$$

and an $n$-sphere of radius $r$ can be defined as

$$S^n = \{|\mathbf{x}\rangle \in \mathbb{R}^{n+1} : \|\mathbf{x}\|_2 = r\}. \tag{4.33}$$

The dimension of $n$-sphere is $n$, and must not be confused with the dimension $(n+1)$ of the Euclidean space in which it is naturally embedded.

**Definition 2 ($n$-sphere):** The set of points in $(n+1)$-space, $(x_1, x_2, \ldots, x_{n+1})^T$, that define an $n$-sphere, $S^n(r)$, is represented by the equation:

$$r^2 = \sum_{i=1}^{n+1} (x_i - c_i)^2, \tag{4.34}$$

where $|c\rangle = (c_1, c_2, \ldots, c_{n+1})^T$ is a center point, and $r$ is the radius.

In particular:

- The pair of points at the ends of a (one-dimensional) line segment is a 0-sphere.
- A circle, which is the one-dimensional circumference of a (two-dimensional) disk, is a 1-sphere.
- The two-dimensional surface of a three-dimensional ball is a 2-sphere, often simply called a sphere.

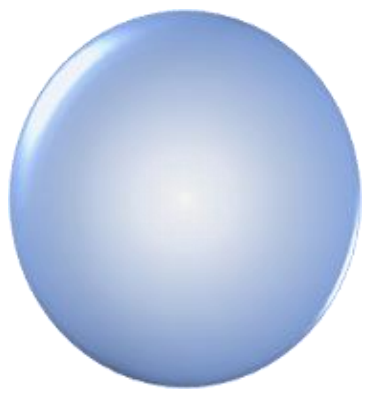

**Figure 4.9.** In Euclidean space, a ball is a volume bounded by a sphere.

An $n$-sphere is the surface or boundary of an $(n+1)$-dimensional ball. In other words, a ball is a solid figure bounded by a sphere; it is also called a solid sphere, see Figure 4.9. It may be a closed ball (including the boundary points that constitute the sphere) or an open ball (excluding them). These concepts are defined not only in three-dimensional Euclidean space but also for lower and higher dimensions spaces in general. A ball in $n$ dimensions is called a hyperball or $n$-ball and is bounded by a hypersphere or $(n-1)$- sphere. Thus, for example,

- A 1-ball, a line segment, is the interior of a 0-sphere.
- A 2-ball, a disk, is the interior of a circle (1-sphere).
- A 3-ball, an ordinary ball, is the interior of a sphere (2-sphere).

Any normed vector space $V$ with the norm $\|\ \ \|$ is a metric space with the metric (distance function) $d(\mathbf{x},\mathbf{y}) = \|\mathbf{x}-\mathbf{y}\|$. In such spaces, an arbitrary ball $B_r(\mathbf{y})$ of points $|\mathbf{x}\rangle$ around a point $|\mathbf{y}\rangle$ with a distance of less than $r$ may be viewed as scaled (by $r$) and translated (by $|\mathbf{y}\rangle$) copy of a unit ball $B_1(\mathbf{0})$. Such "centered" balls with $|\mathbf{y}\rangle = |\mathbf{0}\rangle$ are denoted with $B(r)$.

**Definition (Ball):** An open ball around the origin with radius $r$ is given by the set (with the $p$-norm)

$$B(r) = \left\{|\mathbf{x}\rangle \in \mathbb{R}^n : \|\mathbf{x}\|_p = (|x_1|^p + |x_2|^p + \cdots + |x_n|^p)^{\frac{1}{p}} < r\right\}. \tag{4.35}$$

The Euclidean ball (or just ball) in $\mathbb{R}^n$ has the form

$$B(\mathbf{x}_c, r) = \{|\mathbf{x}\rangle : \|\mathbf{x}-\mathbf{x}_c\|_2 \le r\} = \{|\mathbf{x}\rangle : \langle\mathbf{x}-\mathbf{x}_c|\mathbf{x}-\mathbf{x}_c\rangle \le r^2\}. \tag{4.36}$$

For $n = 2$, in a 2-dimensional space $\mathbb{R}^2$, "balls" according to the $\ell^1$-norm are bounded by squares with their diagonals parallel to the coordinate axes; those according to the $\ell^\infty$-norm have squares with their sides parallel to the coordinate axes as their boundaries. The $\ell^2$-norm generates the well-known discs within circles; see Figure 4.10.

For $n = 3$, the $\ell^1$-balls are within octahedra with axes-aligned body diagonals, the $\ell^\infty$-balls are within cubes with axes-aligned edges. Obviously, $\ell^2$ generates the inner of usual spheres.

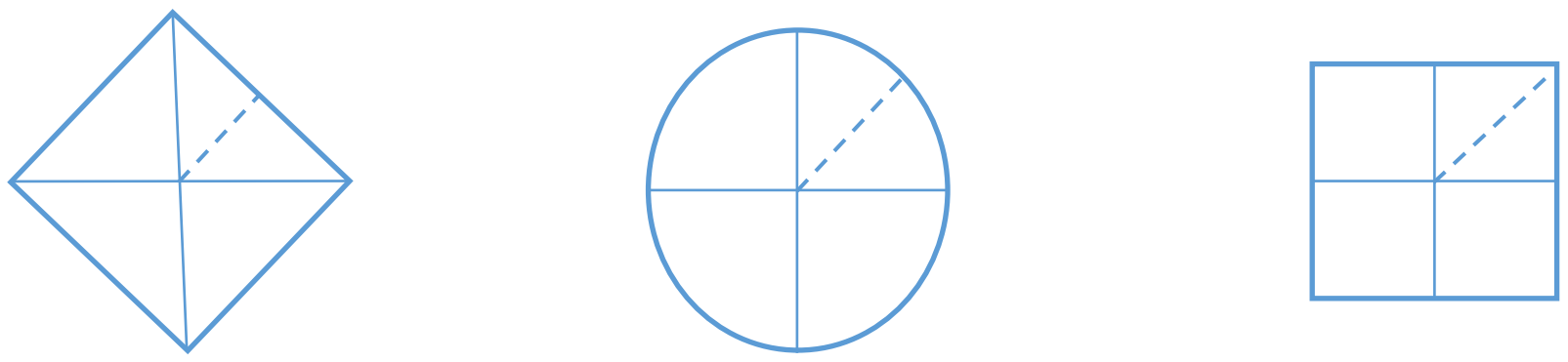

Manhattan Distance. (square) | Euclidean Distance. (circle) | Chebyshev Distance. (square)

**Figure 4.10.** In a 2-dimensional plane $\mathbb{R}^2$, "balls" according to the $\ell^1$-, $\ell^2$-, $\ell^\infty$- norms.

### 4.2.3 Ellipse, Ellipsoid, and Hyper-Ellipsoid in $\mathbb{R}^n$

Using a Cartesian coordinate system in which the origin is the center of the ellipse, and the coordinate axes are axes of the ellipse, the implicit equation of the ellipse has the standard form

$$\frac{x_1^2}{a^2}+\frac{x_2^2}{b^2}=1, \tag{4.37}$$

where $a$, $b$ are positive real numbers. The points $(a,0)^T$ and $(0,b)^T$ lie on the curve. The line segments from the origin to these points are called the principal semi-axes of the ellipse, because $a$, $b$ are half the length of the principal axes. They correspond to the semi-major axis and semi-minor axis of an ellipse. Also, the equation of the ellipse can be written as

$$\begin{aligned}\langle \mathbf{x}|\mathbf{A}^{-1}|\mathbf{x}\rangle &= (x_1 \quad x_2)\begin{pmatrix}\frac{1}{a^2} & 0\\ 0 & \frac{1}{b^2}\end{pmatrix}\begin{pmatrix}x_1\\ x_2\end{pmatrix}\\ &= (x_1 \quad x_2)\begin{pmatrix}\frac{1}{a^2}x_1\\ \frac{1}{b^2}x_2\end{pmatrix}\\ &= \frac{x_1^2}{a^2}+\frac{x_2^2}{b^2},\end{aligned} \tag{4.38}$$

where $|\mathbf{x}\rangle = \begin{pmatrix}x_1\\ x_2\end{pmatrix}$, $\mathbf{A} = \begin{pmatrix}a^2 & 0\\ 0 & b^2\end{pmatrix}$, and $\mathbf{A}^{-1} = \begin{pmatrix}\frac{1}{a^2} & 0\\ 0 & \frac{1}{b^2}\end{pmatrix}$. The standard equation in $\mathbb{R}^n$ is

$$\frac{x_1^2}{a_1^2}+\frac{x_2^2}{a_2^2}+\cdots+\frac{x_n^2}{a_n^2}=1. \tag{4.39}$$

**Definition (Ellipsoids):** The ellipsoids is the set

$$\mathcal{E} = \{|\mathbf{x}\rangle : \langle \mathbf{x} - \mathbf{x}_c|\mathbf{P}^{-1}|\mathbf{x}-\mathbf{x}_c\rangle \le 1\}, \tag{4.40}$$

where $\mathbf{P}$ is symmetric and positive definite (i.e., $\mathbf{P}=\mathbf{P}^T$ and $\langle \mathbf{u}|\mathbf{P}|\mathbf{u}\rangle > 0$ whenever $|\mathbf{u}\rangle \neq |\mathbf{0}\rangle$). The vector $|\mathbf{x}_c\rangle \in \mathbb{R}^n$ is the center of the ellipsoid.

**Remarks**

- The matrix $\mathbf{P}$ determines how far the ellipsoid extends in every direction from $|\mathbf{x}_c\rangle$; the lengths of the semi-axes of $\mathcal{E}$ are given by $\sqrt{\lambda_i}$, where $\lambda_i$ are the eigenvalues of $\mathbf{P}$.
- A ball is an ellipsoid with $\mathbf{P}=r^2\mathbf{I}$.
- Figure 4.11 shows an ellipsoid in $\mathbb{R}^2$.
- A hyper-ellipsoid, or ellipsoid of dimension $n-1$ in a Euclidean space of dimension $n$, is a quadric hypersurface defined by a polynomial of degree two that has a homogeneous part of degree two which is a positive definite quadratic form.

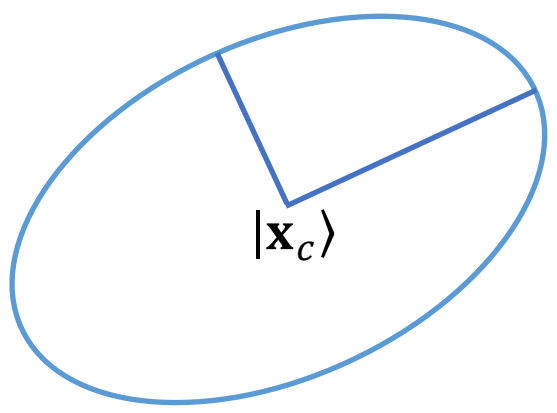

**Figure 4.11.** The geometry of an ellipsoid in $\mathbb{R}^2$.

## 4.3 Convex sets

Convexity, or convex analysis [1-6], is an area of mathematics where one studies questions related to two basic objects, namely convex sets and convex functions. Triangles, rectangles, and "certain" polygons are examples of convex sets in the plane, and the quadratic function $f(x) = ax^2 + bx + c$ is convex provided that $a \geq 0$. Actually, the points in the plane on or above the graph of this quadratic function is another example of a convex set. But one may also consider convex sets in $\mathbb{R}^n$, for any $n$, and convex functions of several variables. Convexity is the mathematical core of optimization, and it plays an important role in many other mathematical areas such as statistics, approximation theory, differential equations, and mathematical economics.

A basic optimization problem is to minimize a real-valued function $f$ of $n$ variables, say $f|\mathbf{x}\rangle$ where $|\mathbf{x}\rangle = (x_1, \ldots, x_n)^T \in A$ and $A$ is the domain of $f$. A global minimum of $f$ is a point $|\mathbf{x}^*\rangle$ with $f|\mathbf{x}^*\rangle \leq f|\mathbf{x}\rangle$ for all $|\mathbf{x}\rangle \in A$. Often it is hard to find a global minimum, so one settles with a local minimum point which satisfies $f|\mathbf{x}^*\rangle \leq f|\mathbf{x}\rangle$ for all $|\mathbf{x}\rangle$ in $A$ that are sufficiently close to $|\mathbf{x}^*\rangle$. There are several optimization algorithms that can locate a local minimum of $f$. Unfortunately, the function value in a local minimum may be much larger than the global minimum value. This raises the question: Are there functions where a local minimum point is also a global minimum? The main answer to this question is:

If $f$ is a convex function, and the domain $A$ is a convex set, then a local minimum point is also a global minimum point! Thus, one can find the global minimum of convex functions, whereas this may be hard (or even impossible) in other situations.

We start the study of convexity with sets. Geometric ideas play an underlying role in convex analysis, its extensions, and applications.

**Definition (Line Segment):** Given two elements $|\mathbf{a}\rangle$ and $|\mathbf{b}\rangle$ in $\mathbb{R}^n$, define the interval/line segment

$$[\mathbf{a}, \mathbf{b}] = \{\lambda|\mathbf{a}\rangle + (1-\lambda)|\mathbf{b}\rangle : \lambda \in [0,1]\}. \tag{4.41}$$

The parameter value $\lambda = 0$ corresponds to $|\mathbf{b}\rangle$, and the parameter value $\lambda = 1$ corresponds to $|\mathbf{a}\rangle$. Values of the parameter $\lambda$ between 0 and 1 correspond to the (closed) line segment between $|\mathbf{a}\rangle$ and $|\mathbf{b}\rangle$. Note that if $|\mathbf{a}\rangle = |\mathbf{b}\rangle$, then this interval reduces to a singleton $[\mathbf{a}, \mathbf{b}] = \{\mathbf{a}\}$. Equivalently, we have

$$[\mathbf{a}, \mathbf{b}] = \{|\mathbf{b}\rangle + \lambda(|\mathbf{a}\rangle - |\mathbf{b}\rangle) : \lambda \in [0,1]\} = \{|\mathbf{b}\rangle + \lambda|\mathbf{a} - \mathbf{b}\rangle : \lambda \in [0,1]\}. \tag{4.42}$$

This expression gives another interpretation: the line segment is the sum of the base point $|\mathbf{b}\rangle$ (corresponding to $\lambda = 0$) and the direction $(|\mathbf{a}\rangle - |\mathbf{b}\rangle = |\mathbf{a} - \mathbf{b}\rangle)$ (which points from $|\mathbf{b}\rangle$ to $|\mathbf{a}\rangle$) scaled by the parameter $\lambda$. Thus, $\lambda$ gives the fraction of the way from $|\mathbf{b}\rangle$ to $|\mathbf{a}\rangle$. As $\lambda$ increases from 0 to 1, the point moves from $|\mathbf{b}\rangle$ to $|\mathbf{a}\rangle$. See Figure 4.12.

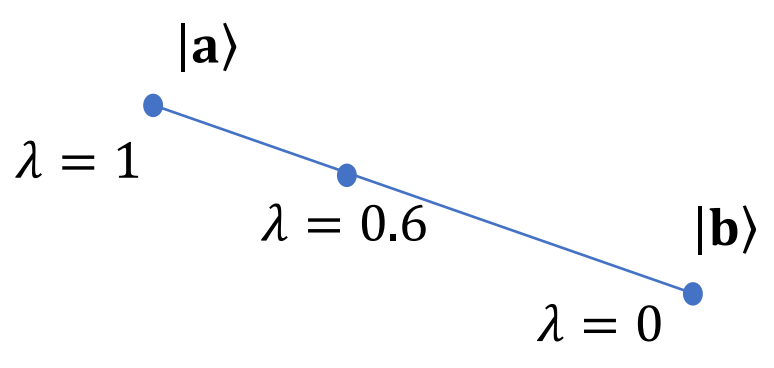

**Figure 4.12.** The line segment between $|\mathbf{a}\rangle$ and $|\mathbf{b}\rangle$ in $\mathbb{R}^n$.

**Definition (Convex Set):** A set $C$ is convex if the line segment between any two points in $C$ lies in $C$. Equivalently, $C$ is convex if

$$\lambda|\mathbf{a}\rangle + (1-\lambda)|\mathbf{b}\rangle \in C, \tag{4.43}$$

for all $|\mathbf{a}\rangle, |\mathbf{b}\rangle \in C$ and $\lambda \in [0,1]$.

Loosely speaking, a convex set in $\mathbb{R}^2$ (or $\mathbb{R}^n$) is a set "with no holes." See Figure 4.13.

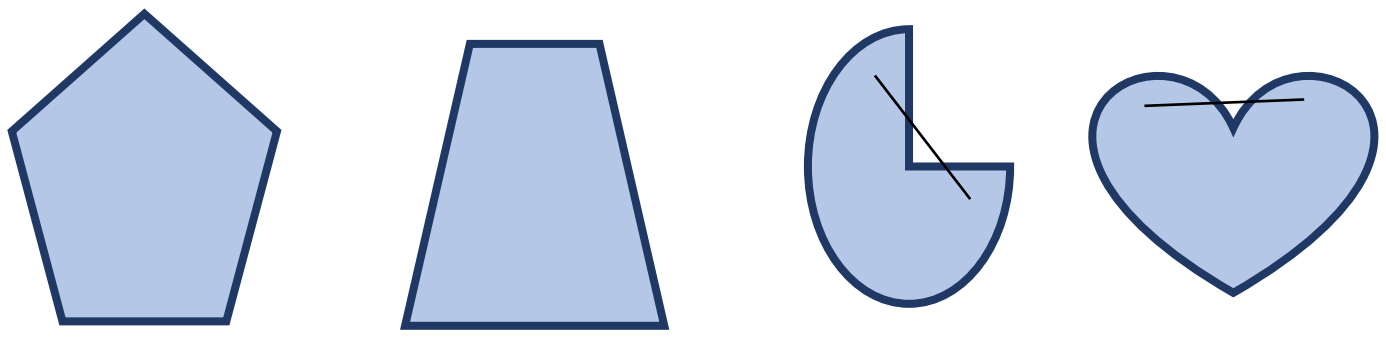

**Figure 4.13.** Some simple convex and non-convex sets.

Some simple examples

- The empty set $\emptyset$ is convex.
- Any single point $\{x_0\}$ is convex.
- The whole space $\mathbb{R}^n$ is convex.
- A line segment is convex.

### Hyperplanes, Half-Spaces, and Support Line/Hyperplane

Remember, a hyperplane is a set of the form $H = \{|\mathbf{x}\rangle : \langle \mathbf{a}|\mathbf{x}\rangle = b\}$, where $|\mathbf{a}\rangle \in \mathbb{R}^n$, $|\mathbf{a}\rangle \neq |\mathbf{0}\rangle$, and $b \in \mathbb{R}$. Geometrically, the hyperplane $H$ can be interpreted as the set of points with a constant inner product to a given vector $|\mathbf{a}\rangle$, or as a hyperplane with a normal vector $|\mathbf{a}\rangle$; the constant $b \in \mathbb{R}$ determines the offset of the hyperplane from the origin. This geometric interpretation can be understood by expressing the hyperplane in the form $\{|\mathbf{x}\rangle : \langle \mathbf{a}|\mathbf{x} - \mathbf{x}_0\rangle = 0\}$, where $|\mathbf{x}_0\rangle$ is any point in the hyperplane (i.e., any point that satisfies $\langle \mathbf{a}|\mathbf{x}_0\rangle = b$). A hyperplane divides $\mathbb{R}^n$ into two half-spaces. A (closed) half-space is a set of the form $\{|\mathbf{x}\rangle : \langle \mathbf{a}|\mathbf{x}\rangle \leq b\}$, where $|\mathbf{a}\rangle \neq |\mathbf{0}\rangle$. Half-spaces are convex. This is illustrated in Figure 4.14 and Figure 4.15.

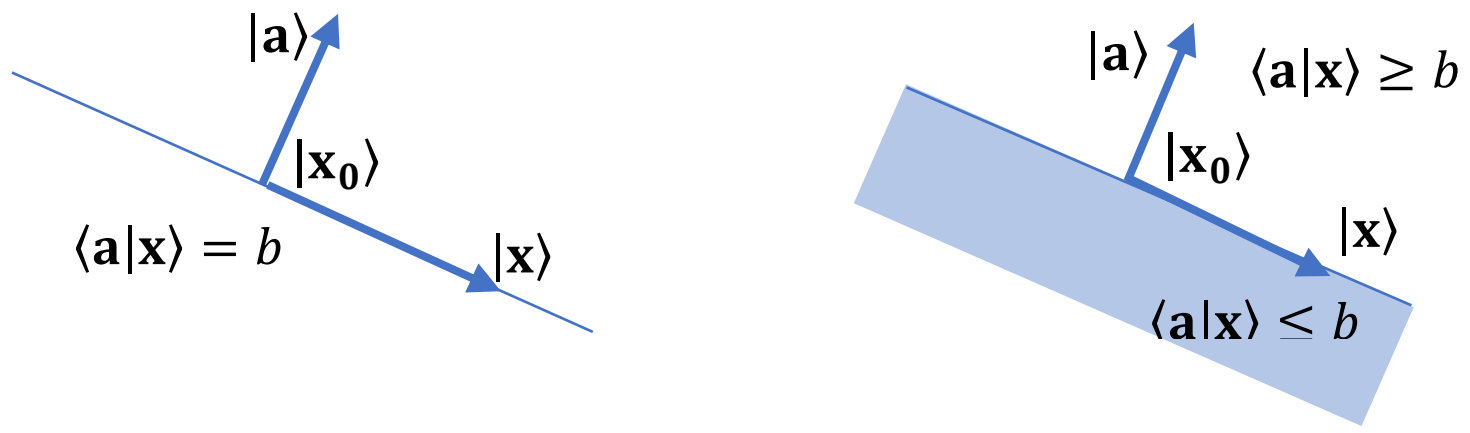

**Figure 4.14.** A hyperplane defined by $\langle \mathbf{a}|\mathbf{x}\rangle = b$ in $\mathbb{R}^2$ determines two half-spaces.

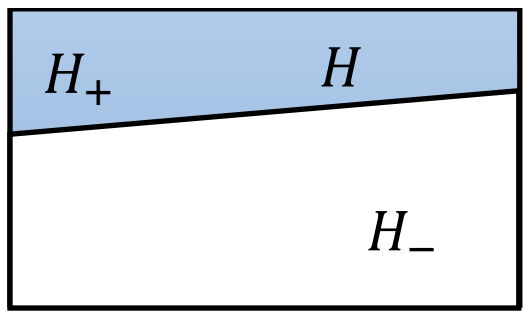

**Figure 4.15.** The two half-spaces are determined by a hyperplane, $H$.

**Theorem 4.2:** The hyperplane is a convex set.

**Proof:**

Let $\langle \mathbf{a}|\mathbf{x}\rangle = b$ be a hyperplane. $|\mathbf{x}_1\rangle$, and $|\mathbf{x}_2\rangle$ are any points on the hyperplane. Then $\langle \mathbf{a}|\mathbf{x}_1\rangle = b$, $\langle \mathbf{a}|\mathbf{x}_2\rangle = b$. Therefore, for $0 \leq \lambda \leq 1$,

$$\begin{aligned}\langle \mathbf{a}|\lambda\mathbf{x}_1 + (1-\lambda)\mathbf{x}_2\rangle &= \langle \mathbf{a}|\lambda\mathbf{x}_1\rangle + \langle \mathbf{a}|(1-\lambda)\mathbf{x}_2\rangle \\ &= \lambda\langle \mathbf{a}|\mathbf{x}_1\rangle + (1-\lambda)\langle \mathbf{a}|\mathbf{x}_2\rangle \\ &= \lambda b + (1-\lambda)b \\ &= b.\end{aligned}$$

Hence, $\lambda|\mathbf{x}_1\rangle + (1-\lambda)|\mathbf{x}_2\rangle$, for $0 \leq \lambda \leq 1$ lies in the hyperplane. So, the hyperplane is convex.

∎

**Definition (Support Line/Hyperplane):** An important property of any convex set $K$ in the plane is that at every point $|\mathbf{p}\rangle$ on the boundary of $K$, there exists at least one line $l$ (or generally a $(n-1)$-dimensional hyperplane in higher dimensions) that passes through $|\mathbf{p}\rangle$ such that $K$ lies entirely in one of the closed half-planes (half-spaces) defined by $l$ (see Figure 4.16). Such a line is called a support line/hyperplane for $K$.

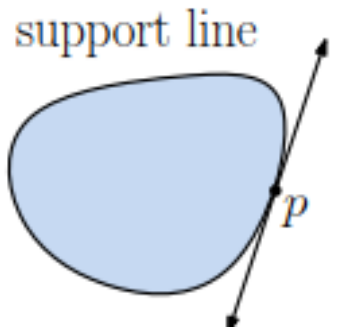

**Figure 4.16.** Support line/hyperplane.

### Euclidean ball

**Theorem 4.3:** The Euclidean ball is a convex set.

**Proof:**

The Euclidean ball in $\mathbb{R}^n$ has the form

$$B(\mathbf{x}_c, r) = \{|\mathbf{x}\rangle : \|\mathbf{x}-\mathbf{x}_c\|_2 \le r\} = \{|\mathbf{x}\rangle : \langle \mathbf{x}-\mathbf{x}_c|\mathbf{x}-\mathbf{x}_c\rangle \le r^2\}.$$

The point $|\mathbf{x}_c\rangle$ is the center of the ball, and the scalar $r$ is its radius; $B(\mathbf{x}_c, r)$ consists of all points within a distance $r$ of the center $|\mathbf{x}_c\rangle$. If $\|\mathbf{x}_1-\mathbf{x}_c\|_2 \le r$, $\|\mathbf{x}_2-\mathbf{x}_c\|_2 \le r$, and $0 \le \theta \le 1$, then

$$\begin{aligned}\|\theta\mathbf{x}_1 + (1-\theta)\mathbf{x}_2 - \mathbf{x}_c\|_2 &= \|\theta\mathbf{x}_1 - \theta\mathbf{x}_c + (1-\theta)\mathbf{x}_2 - \mathbf{x}_c + \theta\mathbf{x}_c\|_2 \\ &= \|\theta(\mathbf{x}_1-\mathbf{x}_c) + (1-\theta)(\mathbf{x}_2-\mathbf{x}_c)\|_2 \\ &\le \theta\|(\mathbf{x}_1-\mathbf{x}_c)\|_2 + (1-\theta)\|(\mathbf{x}_2-\mathbf{x}_c)\|_2 \\ &\le \theta r + (1-\theta) r \\ &\le r.\end{aligned}$$

■

### Positive Semi-definite Matrices

**Theorem 4.4:** The set of positive semi-definite matrices is convex.

**Proof:**

Let $\mathbf{M}_1$ and $\mathbf{M}_2$ are two positive semidefinite matrices. So that, $\langle\mathbf{x}|\mathbf{M}_1|\mathbf{x}\rangle \ge 0$ and $\langle\mathbf{x}|\mathbf{M}_2|\mathbf{x}\rangle \ge 0$. If

$$\mathbf{M} = \theta\mathbf{M}_1 + (1-\theta)\mathbf{M}_2,$$

we have

$$\begin{aligned}\langle\mathbf{x}|\mathbf{M}|\mathbf{x}\rangle &= \langle\mathbf{x}|(\theta\mathbf{M}_1 + (1-\theta)\mathbf{M}_2)|\mathbf{x}\rangle \\ &= \theta\langle\mathbf{x}|\mathbf{M}_1|\mathbf{x}\rangle + (1-\theta)\langle\mathbf{x}|\mathbf{M}_2|\mathbf{x}\rangle \\ &\ge 0.\end{aligned}$$

Hence, the set of positive semi-definite matrices is convex.

■

***Example 4.1***

Show that

$$C = \{(x_1, x_2)^T : 2x_1 + 3x_2 = 7\} \subset \mathbb{R}^2,$$

is a convex set.

***Solution***

Let $|\mathbf{x}\rangle, |\mathbf{y}\rangle \in C$, where $|\mathbf{x}\rangle = (x_1, x_2)^T$, $|\mathbf{y}\rangle = (y_1, y_2)^T$. The line segment joining $|\mathbf{x}\rangle$ and $|\mathbf{y}\rangle$ is the set

$$W = \{|\mathbf{w}\rangle : |\mathbf{w}\rangle = \lambda|\mathbf{x}\rangle + (1-\lambda)|\mathbf{y}\rangle, 0 \le \lambda \le 1\}.$$

For some $0 \le \lambda \le 1$, let $|\mathbf{w}\rangle = (w_1, w_2)^T$ be a point of set $W$, so that

$$w_1 = \lambda x_1 + (1-\lambda)y_1, \qquad w_2 = \lambda x_2 + (1-\lambda)y_2.$$

Since $|\mathbf{x}\rangle, |\mathbf{y}\rangle \in C$, $2x_1 + 3x_2 = 7$ and $2y_1 + 3y_2 = 7$. But,

$$\begin{aligned} 2w_1 + 3w_2 &= 2(\lambda x_1 + (1-\lambda)y_1) + 3(\lambda x_2 + (1-\lambda)y_2) \\ &= \lambda(2x_1 + 3x_2) + (1-\lambda)(2y_1 + 3y_2) \\ &= 7\lambda + 7(1-\lambda) \\ &= 7. \end{aligned}$$

Therefore, $|\mathbf{w}\rangle = (w_1, w_2)^T \in C$. Since $|\mathbf{w}\rangle$ is any point of $C$, $|\mathbf{x}\rangle, |\mathbf{y}\rangle \in C \Rightarrow [\mathbf{x}:\mathbf{y}] \subset C$. Hence $C$ is convex.

***Example 4.2***

Show that in $\mathbb{R}^3$, the closed ball $x_1^2 + x_2^2 + x_3^2 \le 1$ is a convex set.

***Solution***

Let $S = \{(x_1, x_2, x_3)^T : x_1^2 + x_2^2 + x_3^2 \le 1\}$. Also, let $|\mathbf{x}\rangle, |\mathbf{y}\rangle \in S$, where,

$$|\mathbf{x}\rangle = (x_1, x_2, x_3)^T, \qquad |\mathbf{y}\rangle = (y_1, y_2, y_3)^T.$$

Then, by giving condition, we have

$$x_1^2 + x_2^2 + x_3^2 \le 1 \text{ and } y_1^2 + y_2^2 + y_3^2 \le 1.$$

Now, for some scalar $\lambda$, $0 \le \lambda \le 1$, we have

$$\begin{aligned} \|\lambda\mathbf{x} + (1-\lambda)\mathbf{y}\|^2 &= (\lambda x_1 + (1-\lambda)y_1)^2 + (\lambda x_2 + (1-\lambda)y_2)^2 + (\lambda x_3 + (1-\lambda)y_3)^2 \\ &= \lambda^2(x_1^2 + x_2^2 + x_3^2) + (1-\lambda)^2(y_1^2 + y_2^2 + y_3^2) \\ &\quad + 2\lambda(1-\lambda)(x_1y_1 + x_2y_2 + x_3y_3). \end{aligned}$$

By Schwartz's inequality,

$$x_1y_1 + x_2y_2 + x_3y_3 \le \sqrt{x_1^2 + x_2^2 + x_3^2}\sqrt{y_1^2 + y_2^2 + y_3^2},$$

we have

$$\|\lambda\mathbf{x} + (1-\lambda)\mathbf{y}\|^2 \le \lambda^2 + (1-\lambda)^2 + 2\lambda(1-\lambda) = \big(\lambda + (1-\lambda)\big)^2 = 1.$$

Therefore, $\lambda|\mathbf{x}\rangle + (1-\lambda)|\mathbf{y}\rangle$ is a point in $S$. Thus, $|\mathbf{x}\rangle, |\mathbf{y}\rangle \in S \Rightarrow [\mathbf{x}:\mathbf{y}] \subset S$. Hence $S$ is convex.

**Definition (Convex Combination):** Given $|\boldsymbol{\omega}_1\rangle, \dots |\boldsymbol{\omega}_m\rangle \in \mathbb{R}^n$, the element

$$|\mathbf{x}\rangle = \lambda_1|\boldsymbol{\omega}_1\rangle + \cdots + \lambda_m|\boldsymbol{\omega}_m\rangle = \sum_{i=1}^{m} \lambda_i|\boldsymbol{\omega}_i\rangle, \tag{4.44}$$

where,

$$\sum_{i=1}^{m} \lambda_i = 1, \tag{4.45}$$

and $\lambda_i \ge 0$ for some $m \in \mathbb{N}$, is called a convex combination of $|\boldsymbol{\omega}_1\rangle, \dots |\boldsymbol{\omega}_m\rangle$.

***Example 4.3***

Determine whether the vector $(0,7)^T$ is a convex combination of the set $\{(3,6)^T, (-6,9)^T, (2,1)^T, (-1,1)^T\}$.

***Solution***

For these vectors, we have

$$\binom{0}{7} = \binom{3}{6}\lambda_1 + \binom{-6}{9}\lambda_2 + \binom{2}{1}\lambda_3 + \binom{-1}{1}\lambda_4,$$

or

$$3\lambda_1 - 6\lambda_2 + 2\lambda_3 - \lambda_4 = 0, \qquad 6\lambda_1 + 9\lambda_2 + \lambda_3 + \lambda_4 = 7.$$

To these equations, we add a third condition,

$$\lambda_1 + \lambda_2 + \lambda_3 + \lambda_4 = 1.$$

We must determine whether there exist non-negative values of $\lambda_1$, $\lambda_2$, $\lambda_3$, and $\lambda_4$ that simultaneously satisfy the above equations. Solving these equations, we obtain

$$\lambda_1 = \frac{2}{3} + \frac{1}{2}\lambda_4, \quad \lambda_2 = \frac{1}{3} - \frac{5}{16}\lambda_4, \quad \lambda_3 = -\frac{19}{16}\lambda_4,$$

with $\lambda_4$ arbitrary. The choice $\lambda_4 = 0$ is forced, giving

$$\lambda_1 = \frac{2}{3}, \quad \lambda_2 = \frac{1}{3}, \quad \lambda_3 = 0, \quad \lambda_4 = 0,$$

as an acceptable set of constants. Thus, $(0,7)^T$ is a convex combination of the given set of four vectors.

Next, we proceed with intersections of convex sets [1-6].

**Theorem 4.5:** The intersection of an arbitrary collection of convex sets is convex. Equivalently, let $\{\Omega_\alpha\}_{\alpha\in I}$ be a collection of convex subsets of $\mathbb{R}^n$, where $I$ is an arbitrary index set. Then $\bigcap_{\alpha\in I}\Omega_\alpha$ is also a convex subset of $\mathbb{R}^n$.

**Proof:**

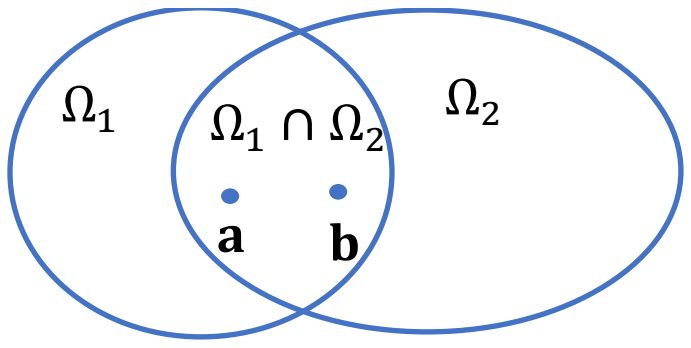

**Figure 4.17.** The intersection of two convex sets is convex

Taking any $|\mathbf{a}\rangle, |\mathbf{b}\rangle \in \bigcap_{\alpha\in I}\Omega_\alpha$. We get that $|\mathbf{a}\rangle, |\mathbf{b}\rangle \in \Omega_\alpha$ for all $\alpha \in I$, see Figure 4.17. The convexity of each $\Omega_\alpha$ ensures that

$$\lambda|\mathbf{a}\rangle + (1-\lambda)|\mathbf{b}\rangle \in \Omega_\alpha,$$

for any $\lambda \in (0,1)$. Thus

$$\lambda|\mathbf{a}\rangle + (1-\lambda)|\mathbf{b}\rangle \in \textstyle\bigcap_{\alpha\in I}\Omega_\alpha,$$

and the intersection $\bigcap_{\alpha\in I}\Omega_\alpha$ is convex.

■

**Theorem 4.6:** Let $|\mathbf{b}_i\rangle \in \mathbb{R}^n$ and $\beta_i \in \mathbb{R}$ for $i \in I$, where $I$ is an arbitrary index set. Then the set

$$C = \{|\mathbf{x}\rangle \in \mathbb{R}^n : \langle\mathbf{x}|\mathbf{b}_i\rangle \le \beta_i, i \in I\}, \tag{4.46}$$

is convex.

Note: Theorem 4.6 would still be valid, of course, if some of the inequalities $\le$ were replaced by $\ge, >, <$ or $=$.

**Proof:**

Let $C_i = \{|\mathbf{x}\rangle : \langle\mathbf{x}|\mathbf{b}_i\rangle \le \beta_i\}$. Then $C_i$ is a closed half-space and $C = \bigcap_{i\in I} C_i$.

■

**Definition (Polyhedral Convex Set):** A set which can be expressed as the intersection of finitely many closed half-spaces and hyperplanes of $\mathbb{R}^n$ is called a polyhedral convex set.

$$\mathcal{P} = \{|\mathbf{x}\rangle : \langle\mathbf{a}_j|\mathbf{x}\rangle \le b_j, j = 1, \dots, m, \langle\mathbf{c}_j|\mathbf{x}\rangle = d_j, j = 1, \dots, p\}, \tag{4.47}$$

or a polyhedron is defined as the solution set of a finite number of linear equalities and inequalities,
or, a polyhedron $\mathcal{P} \subseteq \mathbb{R}^n$ is defined as the solution set of a system of linear inequalities. Thus, $\mathcal{P}$ has the form

$$\mathcal{P} = \{|\mathbf{x}\rangle \in \mathbb{R}^n : \mathbf{A}|\mathbf{x}\rangle \le |\mathbf{b}\rangle\}, \tag{4.48}$$

where $\mathbf{A}$ is a real $m \times n$ matrix, $|\mathbf{b}\rangle \in \mathbb{R}^m$ and where the vector inequality means it holds for every component.

It is easily shown that polyhedron are convex sets. See Figure 4.18. Such sets are considerably better behaved than general convex sets, mostly because of their lack of "curvature."

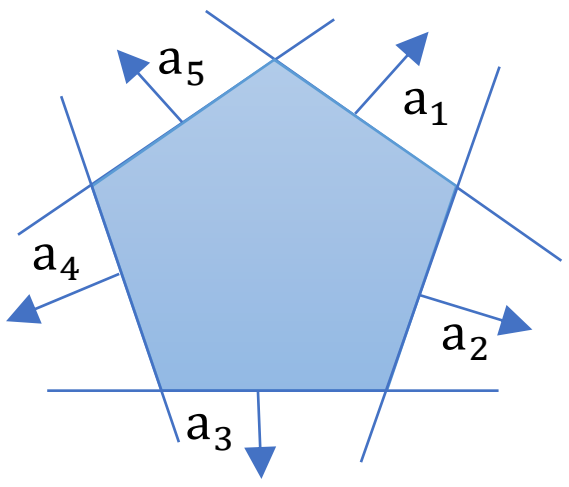

**Figure 4.18.** The polyhedron is the intersection of five half-spaces.

**Theorem 4.7:** Every polyhedron is a convex set.

**Proof:**

Consider a polyhedron $\mathcal{P} = \{|\mathbf{x}\rangle \in \mathbb{R}^n : \mathbf{A}|\mathbf{x}\rangle \leq |\mathbf{b}\rangle\}$ and let $|\mathbf{x}_1\rangle, |\mathbf{x}_2\rangle \in \mathcal{P}$ and $0 \leq \lambda \leq 1$. Then

$$\mathbf{A}\big((1-\lambda)|\mathbf{x}_1\rangle + \lambda|\mathbf{x}_2\rangle\big) = (1-\lambda)\mathbf{A}|\mathbf{x}_1\rangle + \lambda\mathbf{A}|\mathbf{x}_2\rangle \leq (1-\lambda)|\mathbf{b}\rangle + \lambda|\mathbf{b}\rangle = |\mathbf{b}\rangle,$$

which shows that $(1-\lambda)|\mathbf{x}_1\rangle + \lambda|\mathbf{x}_2\rangle \in \mathcal{P}$ and the convexity of $\mathcal{P}$ follows.

■

**Theorem 4.8:** A subset $\Omega$ of $\mathbb{R}^n$ is convex if and only if it contains all convex combinations of its elements.

**Proof:**

Actually, by definition, a set $\Omega$ is convex if and only if $\lambda_1|\boldsymbol{\omega}_1\rangle + \lambda_2|\boldsymbol{\omega}_2\rangle \in \Omega$ whenever $|\boldsymbol{\omega}_1\rangle, |\boldsymbol{\omega}_2\rangle \in \Omega, \lambda_1, \lambda_2 \geq 0$ and $\lambda_1 + \lambda_2 = 1$. In other words, the convexity of $\Omega$ means that $\Omega$ is closed undertaking convex combinations with $m = 2$. We must show that this implies $\Omega$ is also closed undertaking convex combinations with $m > 2$. Fix now a positive integer $m \geq 2$ and suppose that every convex combination of $k \in \mathbb{N}$ elements from $\Omega$, where $k \leq m$, belongs to $\Omega$. Form the convex combination

$$|\mathbf{y}\rangle = \sum_{i=1}^{m+1} \lambda_i |\boldsymbol{\omega}_i\rangle, \qquad \sum_{i=1}^{m+1} \lambda_i = 1, \qquad \lambda_i \geq 0,$$

and observe that if $\lambda_{m+1} = 1$, then $\lambda_1 = \lambda_2 = \cdots = \lambda_m = 0$, so $|\mathbf{y}\rangle = |\boldsymbol{\omega}_{m+1}\rangle \in \Omega$. In the case where $\lambda_{m+1} < 1$ we get the representations

$$\sum_{i=1}^{m} \lambda_i + \lambda_{m+1} = 1,$$

and

$$\sum_{i=1}^{m} \lambda_i = 1 - \lambda_{m+1}, \quad \text{i.e.,} \qquad \sum_{i=1}^{m} \frac{\lambda_i}{1-\lambda_{m+1}} = 1,$$

which implies in turn the inclusion

$$|\mathbf{z}\rangle = \sum_{i=1}^{m} \frac{\lambda_i}{1-\lambda_{m+1}} |\boldsymbol{\omega}_i\rangle \in \Omega.$$

It yields, therefore, the relationships

$$\begin{aligned}|\mathbf{y}\rangle &= (1-\lambda_{m+1}) \sum_{i=1}^{m} \left(\frac{\lambda_i}{1-\lambda_{m+1}}\right) |\boldsymbol{\omega}_i\rangle + \lambda_{m+1}|\boldsymbol{\omega}_{m+1}\rangle \\ &= (1-\lambda_{m+1})|\mathbf{z}\rangle + \lambda_{m+1}|\boldsymbol{\omega}_{m+1}\rangle \in \Omega,\end{aligned}$$

and thus completes the proof of the proposition.

■

**Theorem 4.9:** Let $\Omega_1$ be a convex subset of $\mathbb{R}^n$ and let $\Omega_2$ be a convex subset of $\mathbb{R}^p$. Then the Cartesian product $\Omega_1 \times \Omega_2$ is a convex subset of $\mathbb{R}^n \times \mathbb{R}^p$.

**Proof:**

Fix $a = (|\mathbf{a}_1\rangle, |\mathbf{a}_2\rangle), b = (|\mathbf{b}_1\rangle, |\mathbf{b}_2\rangle) \in \Omega_1 \times \Omega_2, \lambda \in (0,1)$, $|\mathbf{a}_1\rangle, |\mathbf{b}_1\rangle \in \Omega_1$ and $|\mathbf{a}_2\rangle, |\mathbf{b}_2\rangle \in \Omega_2$. The convexity of $\Omega_1$ and $\Omega_2$ gives us

$$\lambda|\mathbf{a}_i\rangle + (1-\lambda)|\mathbf{b}_i\rangle \in \Omega_i \text{ for } i = 1,2,$$

which implies, therefore that

$$\begin{aligned}\lambda a + (1-\lambda)b &= \lambda(|\mathbf{a}_1\rangle, |\mathbf{a}_2\rangle) + (1-\lambda)(|\mathbf{b}_1\rangle, |\mathbf{b}_2\rangle) \\ &= (\lambda|\mathbf{a}_1\rangle + (1-\lambda)|\mathbf{b}_1\rangle, \lambda|\mathbf{a}_2\rangle + (1-\lambda)|\mathbf{b}_2\rangle) \\ &\in \Omega_1 \times \Omega_2.\end{aligned}$$

Thus, the Cartesian product $\Omega_1 \times \Omega_2$ is convex.

■

**Theorem 4.10:** Let $\Omega_1, \Omega_2 \subset \mathbb{R}^n$ be convex and let $\lambda \in \mathbb{R}$. Then both sets $\lambda\Omega_1$ and $\Omega_1 + \Omega_2$ are also convex in $\mathbb{R}^n$.

**Proof:**

Let $\lambda|\mathbf{v}_1\rangle, \lambda|\mathbf{v}_2\rangle \in \lambda\Omega_1$, where $|\mathbf{v}_1\rangle, |\mathbf{v}_2\rangle \in \Omega_1$. Because $\Omega_1$ is convex, we have

$$\alpha|\mathbf{v}_1\rangle + (1-\alpha)|\mathbf{v}_2\rangle \in \Omega_1,$$

for any $\alpha \in (0,1)$. Hence,

$$\alpha\lambda|\mathbf{v}_1\rangle + (1-\alpha)\lambda|\mathbf{v}_2\rangle = \lambda(\alpha|\mathbf{v}_1\rangle + (1-\alpha)|\mathbf{v}_2\rangle) \in \lambda\Omega_1,$$

and thus $\lambda\Omega_1$ is convex.

Let $|\mathbf{v}_1\rangle, |\mathbf{v}_2\rangle \in \Omega_1 + \Omega_2$. Then $|\mathbf{v}_1\rangle = |\dot{\mathbf{v}}_1\rangle + |\ddot{\mathbf{v}}_1\rangle$ and $|\mathbf{v}_2\rangle = |\dot{\mathbf{v}}_2\rangle + |\ddot{\mathbf{v}}_2\rangle$, where $|\dot{\mathbf{v}}_1\rangle$, $|\dot{\mathbf{v}}_2\rangle \in \Omega_1$ and $|\ddot{\mathbf{v}}_1\rangle$, $|\ddot{\mathbf{v}}_2\rangle \in \Omega_2$. Because $\Omega_1$ and $\Omega_2$ are convex, for all $\alpha \in (0,1)$, we have

$$|\mathbf{x}_1\rangle = \alpha|\dot{\mathbf{v}}_1\rangle + (1-\alpha)|\dot{\mathbf{v}}_2\rangle \in \Omega_1,$$

and

$$|\mathbf{x}_2\rangle = \alpha|\ddot{\mathbf{v}}_1\rangle + (1-\alpha)|\ddot{\mathbf{v}}_2\rangle \in \Omega_2,$$

By definition of $\Omega_1 + \Omega_2$, $|\mathbf{x}_1\rangle + |\mathbf{x}_2\rangle \in \Omega_1 + \Omega_2$. Now,

$$\alpha|\mathbf{v}_1\rangle + (1-\alpha)|\mathbf{v}_2\rangle = \alpha(|\dot{\mathbf{v}}_1\rangle + |\ddot{\mathbf{v}}_1\rangle) + (1-\alpha)(|\dot{\mathbf{v}}_2\rangle + |\ddot{\mathbf{v}}_2\rangle) = |\mathbf{x}_1\rangle + |\mathbf{x}_2\rangle \in \Omega_1 + \Omega_2.$$

Hence, $\Omega_1 + \Omega_2$ is convex.

■

Let us continue with the definition of affine mappings.

**Definition (Affine Transformation):** $\mathbf{B}: \mathbb{R}^n \to \mathbb{R}^p$ is affine mapping (affine transformation) if and only if

$$\mathbf{B}(\lambda|\mathbf{x}\rangle + (1-\lambda)|\mathbf{y}\rangle) = \lambda\mathbf{B}|\mathbf{x}\rangle + (1-\lambda)\mathbf{B}|\mathbf{y}\rangle \text{ for all } |\mathbf{x}\rangle, |\mathbf{y}\rangle \in \mathbb{R}^n \text{ and } \lambda \in \mathbb{R}. \tag{4.49}$$

**Theorem 4.11:** A mapping $\mathbf{B}: \mathbb{R}^n \to \mathbb{R}^p$ is affine if there exists a linear mapping $\mathbf{A}: \mathbb{R}^n \to \mathbb{R}^p$ and an element $|\mathbf{b}\rangle \in \mathbb{R}^p$ such that

$$\mathbf{B}|\mathbf{x}\rangle = \mathbf{A}|\mathbf{x}\rangle + |\mathbf{b}\rangle, \tag{4.50}$$

for all $|\mathbf{x}\rangle \in \mathbb{R}^n$.

**Proof:**

If $\mathbf{B}|\mathbf{x}\rangle = \mathbf{A}|\mathbf{x}\rangle + \mathbf{b}$ where $\mathbf{A}$ is linear, one has

$$\begin{aligned}\mathbf{B}\big((1-\lambda)|\mathbf{x}\rangle + \lambda|\mathbf{y}\rangle\big) &= \mathbf{A}\big((1-\lambda)|\mathbf{x}\rangle + \lambda|\mathbf{y}\rangle\big) + |\mathbf{b}\rangle \\ &= (1-\lambda)\mathbf{A}|\mathbf{x}\rangle + \lambda\mathbf{A}|\mathbf{y}\rangle + |\mathbf{b}\rangle \\ &= (1-\lambda)\mathbf{A}|\mathbf{x}\rangle + \lambda\mathbf{A}|\mathbf{y}\rangle + |\mathbf{b}\rangle + \lambda|\mathbf{b}\rangle - \lambda|\mathbf{b}\rangle \\ &= (1-\lambda)\mathbf{A}|\mathbf{x}\rangle + \lambda\mathbf{A}|\mathbf{y}\rangle + \lambda|\mathbf{b}\rangle + (1-\lambda)|\mathbf{b}\rangle \\ &= \big((1-\lambda)\mathbf{A}|\mathbf{x}\rangle + (1-\lambda)|\mathbf{b}\rangle\big) + (\lambda\mathbf{A}|\mathbf{y}\rangle + \lambda|\mathbf{b}\rangle) \\ &= (1-\lambda)(\mathbf{A}|\mathbf{x}\rangle + |\mathbf{b}\rangle) + \lambda(\mathbf{A}|\mathbf{y}\rangle + |\mathbf{b}\rangle) \\ &= (1-\lambda)\mathbf{B}|\mathbf{x}\rangle + \lambda\mathbf{B}|\mathbf{y}\rangle.\end{aligned}$$

Thus, $\mathbf{B}$ is affine. Conversely, if $\mathbf{B}$ is affine, let $|\mathbf{b}\rangle = \mathbf{B}|\mathbf{0}\rangle$ and $\mathbf{A}|\mathbf{x}\rangle = \mathbf{B}|\mathbf{x}\rangle - |\mathbf{b}\rangle$. Then $\mathbf{A}$ is an affine transformation with $\mathbf{A}|\mathbf{0}\rangle = |\mathbf{0}\rangle$. Then $\mathbf{A}$ is actually linear.

■

Now we show that set convexity is preserved under affine operations.

**Theorem 4.12:** Let $\mathbf{B}: \mathbb{R}^n \to \mathbb{R}^p$ be an affine mapping. Suppose that $\Omega$ is a convex subset of $\mathbb{R}^n$ and $\Theta$ is a convex subset of $\mathbb{R}^p$. Then $\mathbf{B}(\Omega)$ is a convex subset of $\mathbb{R}^p$ and $\mathbf{B}^{-1}(\Theta)$ is a convex subset of $\mathbb{R}^n$.

**Proof:**

Fix any $|\mathbf{a}\rangle, |\mathbf{b}\rangle \in \mathbf{B}(\Omega)$ and $\lambda \in (0,1)$. Then $|\mathbf{a}\rangle = \mathbf{B}|\mathbf{x}\rangle$ and $|\mathbf{b}\rangle = \mathbf{B}|\mathbf{y}\rangle$ for $|\mathbf{x}\rangle, |\mathbf{y}\rangle \in \Omega$. Since $\Omega$ is convex, we have $\lambda|\mathbf{x}\rangle + (1-\lambda)|\mathbf{y}\rangle \in \Omega$. Then

$$\begin{aligned}\lambda|\mathbf{a}\rangle + (1-\lambda)|\mathbf{b}\rangle &= \lambda\mathbf{B}|\mathbf{x}\rangle + (1-\lambda)\mathbf{B}|\mathbf{y}\rangle \\ &= \mathbf{B}(\lambda|\mathbf{x}\rangle + (1-\lambda)|\mathbf{y}\rangle) \in \mathbf{B}(\Omega),\end{aligned}$$

which justifies the convexity of the image $\mathbf{B}(\Omega)$.

Taking now any $|\mathbf{x}\rangle, |\mathbf{y}\rangle \in \mathbf{B}^{-1}(\Theta)$ and $\lambda \in (0,1)$, we get $\mathbf{B}|\mathbf{x}\rangle$ and $\mathbf{B}|\mathbf{y}\rangle \in \Theta$. This gives us

$$\lambda\mathbf{B}|\mathbf{x}\rangle + (1-\lambda)\mathbf{B}|\mathbf{y}\rangle = \mathbf{B}(\lambda|\mathbf{x}\rangle + (1-\lambda)|\mathbf{y}\rangle) \in \Theta,$$

by the convexity of $\Theta$. Thus, we have $\lambda|\mathbf{x}\rangle + (1-\lambda)|\mathbf{y}\rangle \in \mathbf{B}^{-1}(\Theta)$, which verifies the convexity of the inverse image $\mathbf{B}^{-1}(\Theta)$.

■

## 4.4 Convex Hull

Indeed, there are many problems that are comparatively easy to solve for convex sets but very hard in general. However, not all sets are convex, and a discrete set of points is never convex, unless it consists of at most one point only. In such a case, it is useful to make a given set $P$ convex, that is, approximate $P$ with or, rather, encompass $P$ within a convex set $H \supseteq P$. Ideally, $H$ differs from $P$ as little as possible; that is, we want $H$ to be the smallest convex set enclosing $P$.

At this point, let us step back for a second and ask ourselves whether this wish makes sense at all: Does such a set $H$ (always) exist? Fortunately, we are on the safe side because the whole space $\mathbb{R}^n$ is certainly convex. It is less obvious, but we will see below that $H$ is actually unique. Therefore, it is legitimate to refer to $H$ as the smallest convex set enclosing $P$ or, shortly, the convex hull of $P$.

The convex hull [1-6] of a shape is the smallest convex set that contains it. The convex hull may be defined either as the intersection of all convex sets containing a given subset of Euclidean space or, equivalently, as the set of all convex combinations of points in the subset. For a bounded subset of the plane, the convex hull may be visualized as the shape enclosed by a rubber band stretched around the subset. (See Figure 4.19)

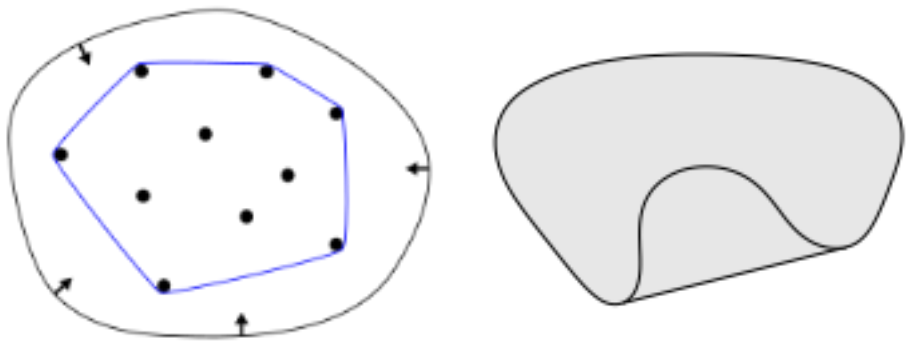

**Figure 4.19.** The convex hulls of two sets in $\mathbb{R}^2$. Left. The convex hull of a set of 10 points. Right. The convex hull of the non-convex shape.

**Definition (Convex Hull):** Let $\Omega$ be a subset of $\mathbb{R}^n$. The convex hull of $\Omega$ is defined by

$$\text{co}\,\Omega = \bigcap\{C : C \text{ is convex and } \Omega \subset C\}. \tag{4.51}$$

From the definition of the convex hull of a given set $\Omega$, the convex hull is

- the set of all convex combinations of points in $\Omega$, or
- the (unique) minimal convex set containing $\Omega$, or
- the intersection of all convex sets containing $\Omega$.

**Theorem 4.13:** For any subset $\Omega$ of $\mathbb{R}^n$, its convex hull admits the representation

$$\text{co}\,\Omega = \left\{\sum_{i=1}^{m} \lambda_i|\mathbf{a}_i\rangle : \sum_{i=1}^{m} \lambda_i = 1\,, \lambda_i \geq 0, |\mathbf{a}_i\rangle \in \Omega, m \in \mathbb{N}\right\}. \tag{4.52}$$

**Proof:**

Denoting by $C$ the right-hand side of the representation to prove, we obviously have $\Omega \subset C$. Let us check that the set $C$ is convex. Take any $|\mathbf{a}\rangle, |\mathbf{b}\rangle \in C$ and get

$$|\mathbf{a}\rangle = \sum_{i=1}^{p} \alpha_i|\mathbf{a}_i\rangle\,, \qquad |\mathbf{b}\rangle = \sum_{i=1}^{q} \beta_i|\mathbf{b}_i\rangle,$$

where $|\mathbf{a}_i\rangle, |\mathbf{b}_i\rangle \in \Omega$, $\beta \geq 0$ with $\sum_{i=1}^{p} \alpha_i = 1$, $\sum_{i=1}^{q} \beta_i = 1$, and $p, q \in \mathbb{N}$. It follows easily that for every number $\lambda \in (0,1)$, we have

$$\lambda|\mathbf{a}\rangle + (1-\lambda)|\mathbf{b}\rangle = \sum_{i=1}^{p} \lambda\alpha_i|\mathbf{a}_i\rangle + \sum_{i=1}^{q} (1-\lambda)\beta_i|\mathbf{b}_i\rangle$$
$$= \lambda\alpha_1|\mathbf{a}_1\rangle + \lambda\alpha_2|\mathbf{a}_2\rangle + \cdots + \lambda\alpha_p|\mathbf{a}_p\rangle + (1-\lambda)\beta_1|\mathbf{b}_1\rangle + (1-\lambda)\beta_2|\mathbf{b}_2\rangle + (1-\lambda)\beta_q|\mathbf{b}_q\rangle.$$

Then the resulting equality

$$\sum_{i=1}^{p} \lambda\alpha_i + \sum_{i=1}^{q} (1-\lambda)\beta_i = \lambda\sum_{i=1}^{p} \alpha_i + (1-\lambda)\sum_{i=1}^{q} \beta_i = 1,$$

ensures that $\lambda|\mathbf{a}\rangle + (1-\lambda)|\mathbf{b}\rangle \in C$, and thus co $\Omega \subset C$ by the definition of co $\Omega$.

Fix now any $|\mathbf{a}\rangle = \sum_{i=1}^{m} \lambda_i|\mathbf{a}_i\rangle \in C$ with $|\mathbf{a}_i\rangle \in \Omega \subset \text{co}\,\Omega$ for $i = 1, \dots, m$. Since the set co $\Omega$ is convex, we conclude by the fact “a subset of $\mathbb{R}^n$ is convex if and only if it contains all convex combinations of its elements” that $|\mathbf{a}\rangle \in \text{co}\,\Omega$ and thus co $\Omega = C$.

■

**Theorem 4.14:** The convex hull co $\Omega$ is the smallest convex set containing $\Omega$.

**Proofs:**

**Proof 1:** First of all, co($\Omega$) contains $\Omega$: for every $|\mathbf{x}\rangle \in \Omega$, $1|\mathbf{x}\rangle$ is a convex combination of size 1, so $|\mathbf{x}\rangle \in \text{co}(\Omega)$.

Second, co($\Omega$) is a convex set: if we take $|\mathbf{x}\rangle, |\mathbf{y}\rangle \in \text{co}(\Omega)$, which are the convex combinations of points in $\Omega$, then $t|\mathbf{x}\rangle + (1-t)|\mathbf{y}\rangle$ can be expanded to get another convex combination of points in $\Omega$.

All convex sets containing $\Omega$ must contain co($\Omega$), and co($\Omega$) is itself a convex set containing $\Omega$; therefore, it is the smallest such set.

■

**Proof 2:** The convexity of the set co $\Omega \supset \Omega$ follows from the fact “the intersection of an arbitrary collection of convex sets is convex.” On the other hand, for any convex set $C$ that contains $\Omega$, we clearly have co $\Omega \subset C$, which verifies the theorem.

■

**Proof 3:** The intersection of an arbitrary collection of convex sets is convex. Since any set is contained in at least one convex set (the whole vector space in which it sits), it follows that any set, $A$, is contained in the smallest convex set, namely the intersection of all the convex sets that contain $A$.

■

## 4.5 Convex Cone

**Definition (Cone):** A subset $K$ of $\mathbb{R}^n$ is called a cone if it is closed under positive scalar multiplication, i.e., $\lambda|\mathbf{x}\rangle \in K$ when $|\mathbf{x}\rangle \in K$ and $\lambda > 0$. Such a set is a union of half-lines emanating from the origin. The origin itself may or may not be included.

**Definition (Convex Cone):** A convex cone is a cone which is a convex set which means that for any $|\mathbf{x}_1\rangle, |\mathbf{x}_2\rangle \in C$ and $\theta_1, \theta_2 \geq 0$, we have

$$\theta_1|\mathbf{x}_1\rangle + \theta_2|\mathbf{x}_2\rangle \in C. \tag{4.53}$$

Points of this form can be described geometrically as forming the two-dimensional pie slice with an apex $|\mathbf{0}\rangle$ and edges passing through $|\mathbf{x}_1\rangle$ and $|\mathbf{x}_2\rangle$. (See Figure 4.20)

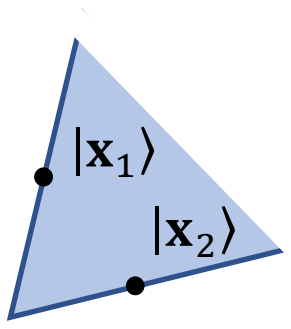

**Figure 4.20.** The pie slice shows all points of the form $\theta_1|\mathbf{x}_1\rangle + \theta_2|\mathbf{x}_2\rangle$, where $\theta_1, \theta_2 \geq 0$.

**Examples:**

1- Subspaces of $\mathbb{R}^n$ are, in particular, convex cones.
2- So are the open and closed half-spaces corresponding to a hyperplane through the origin.
3- Two of the most important convex cones are the non-negative orthant of $\mathbb{R}^n$,

$$\{|\mathbf{x}\rangle = (\xi_1, \xi_2, \dots, \xi_n)^T : \xi_1 \geq 0, \dots, \xi_n \geq 0\}, \tag{4.54}$$

and the positive orthant

$$\{|\mathbf{x}\rangle = (\xi_1, \xi_2, \dots, \xi_n)^T : \xi_1 > 0, \dots, \xi_n > 0\}. \tag{4.55}$$

**Theorem 4.15:** The intersection of an arbitrary collection of convex cones is a convex cone.

**Theorem 4.16:** Let $|\mathbf{b}_i\rangle \in \mathbb{R}^n$ for $i \in I$, where $I$ is an arbitrary index set. Then

$$K = \{|\mathbf{x}\rangle \in \mathbb{R}^n : \langle\mathbf{x}|\mathbf{b}_i\rangle \leq 0, i \in I\}, \tag{4.56}$$

is a convex cone.

**Theorem 4.17:** A subset of $\mathbb{R}^n$ is a convex cone if and only if it is closed under addition and positive scalar multiplication.

**Proof:**

Let $K$ be a cone. Let $|\mathbf{x}\rangle \in K$ and $|\mathbf{y}\rangle \in K$. If $K$ is convex, the vector $|\mathbf{z}\rangle = \frac{1}{2}(|\mathbf{x}\rangle + |\mathbf{y}\rangle)$ belongs to $K$, and hence $|\mathbf{x}\rangle + |\mathbf{y}\rangle = 2|\mathbf{z}\rangle \in K$. On the other hand, if $K$ is closed under addition, and if $0 < \lambda < 1$, the vectors $(1-\lambda)|\mathbf{x}\rangle$ and $\lambda|\mathbf{y}\rangle$ belong to $K$, and hence $(1-\lambda)|\mathbf{x}\rangle + \lambda|\mathbf{y}\rangle \in K$. Thus, $K$ is convex if and only if it is closed under addition.

∎

**Theorem 4.18:** A subset of $\mathbb{R}^n$ is a convex cone if and only if it contains all the positive linear combinations of its elements (i.e., linear combinations $\lambda_1|\mathbf{x}_1\rangle + \lambda_2|\mathbf{x}_2\rangle + \cdots + \lambda_m|\mathbf{x}_m\rangle$ in which the coefficients are all positive).

**Theorem 4.19:** Let $S$ be an arbitrary subset of $\mathbb{R}^n$, and let $K$ be the set of all positive linear combinations of $S$. Then $K$ is the smallest convex cone which includes $S$.

**Proof:**

Clearly, $K$ is closed under addition and positive scalar multiplication, and $K \supset S$. Every convex cone including $S$ must, on the other hand, include $K$.

■

A simpler description is possible when $S$ is convex, as follows.

**Theorem 4.20:** Let $C$ be a convex set, and let

$$K = \{\lambda|\mathbf{x}\rangle : \lambda > 0, |\mathbf{x}\rangle \in C\}. \tag{4.57}$$

Then $K$ is the smallest convex cone which includes $C$.

**Proof:**

This follows from the preceding theorem. Namely, every positive linear combination of elements of $C$ is a positive scalar multiple of a convex combination of elements of $C$ and hence is an element of $K$.

■

A vector $|\mathbf{x}^*\rangle$ is said to be normal to a convex set $C$ at a point $|\mathbf{a}\rangle$, where $|\mathbf{a}\rangle \in C$, if $|\mathbf{x}^*\rangle$ does not make an acute angle with any line segment in $C$ with $|\mathbf{a}\rangle$ as an endpoint, i.e., if $\langle \mathbf{v} - \mathbf{a}|\mathbf{x}^*\rangle \le 0$ for every $|\mathbf{x}\rangle \in C$. For instance, if $C$ is a half-space $\{|\mathbf{x}\rangle : \langle \mathbf{x}|\mathbf{b}\rangle \le \beta\}$ and $|\mathbf{a}\rangle$ satisfies $\langle \mathbf{a}|\mathbf{b}\rangle = \beta$, then $|\mathbf{b}\rangle$ is normal to $C$ at $|\mathbf{a}\rangle$. In general, the set of all vectors $|\mathbf{x}^*\rangle$ normal to $C$ at $|\mathbf{a}\rangle$ is called the normal cone to $C$ at $|\mathbf{a}\rangle$. The reader can verify easily that this cone is always convex.

Another easily verified example of a convex cone is the barrier cone of a convex set $C$. This is defined as the set of all vectors $|\mathbf{x}^*\rangle$ such that, for some $\beta \in \mathbb{R}$, $\langle \mathbf{x}|\mathbf{x}^*\rangle \le \beta$ for every $|\mathbf{x}\rangle \in C$.

## 4.6 Convex Functions

This section collects basic facts about general (extended-real-valued) convex functions, including their analytic and geometric characterizations [1-6].

**Definition (Extended Real Valued Function):** Let $f: \Omega \to \overline{\mathbb{R}}$, ($\overline{\mathbb{R}} = (-\infty, \infty]$), be an extended real-valued function defined on a convex set $\Omega \subset \mathbb{R}^n$. Then the function $f$ is convex on $\Omega$ if

$$f(\lambda|\mathbf{x}\rangle + (1-\lambda)|\mathbf{y}\rangle) \le \lambda f|\mathbf{x}\rangle + (1-\lambda) f|\mathbf{y}\rangle, \tag{4.58}$$

for all $|\mathbf{x}\rangle, |\mathbf{y}\rangle \in \Omega$ and $\lambda \in (0,1)$.

Geometrically, this inequality means that the line segment between $(|\mathbf{x}\rangle, f|\mathbf{x}\rangle)$ and $(|\mathbf{y}\rangle, f|\mathbf{y}\rangle)$, which is the chord from $|\mathbf{x}\rangle$ to $|\mathbf{y}\rangle$, lies above the graph of $f$ (see Figures 4.21 and 4.22).

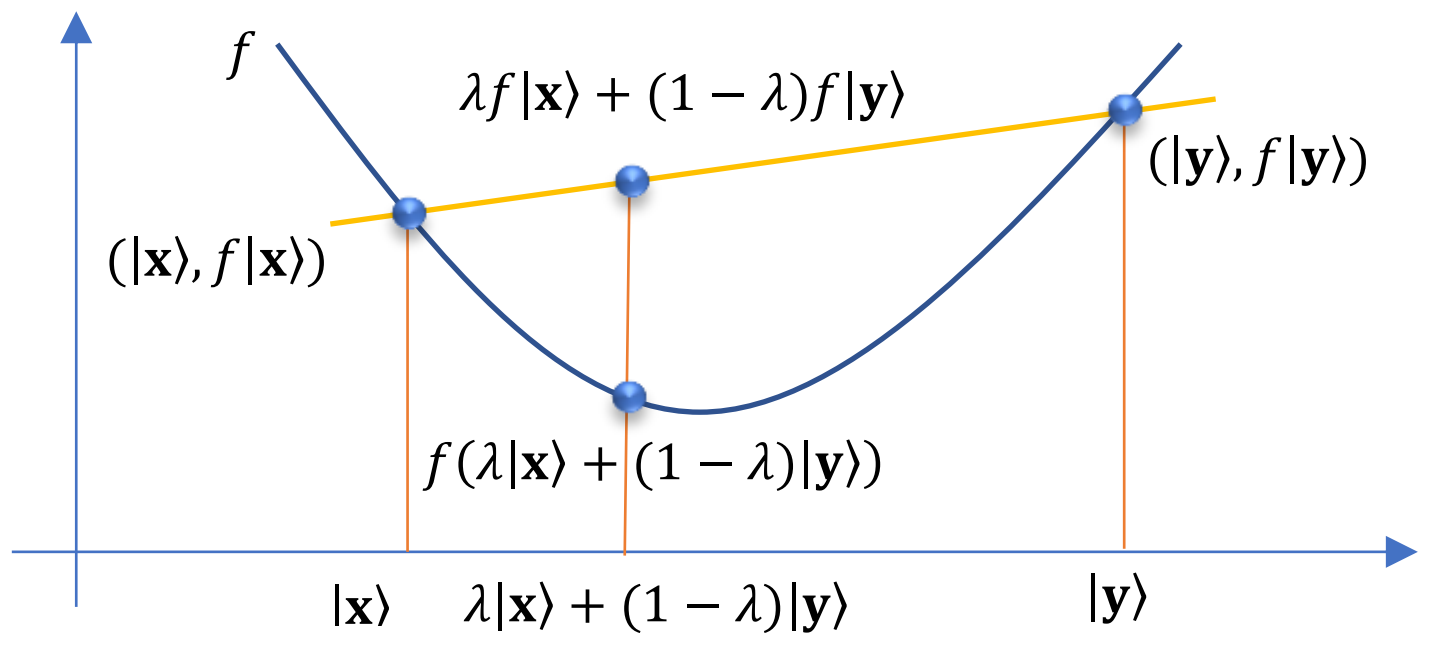

**Figure 4.21.** Convex function. The line segment lies above the graph of $f$.

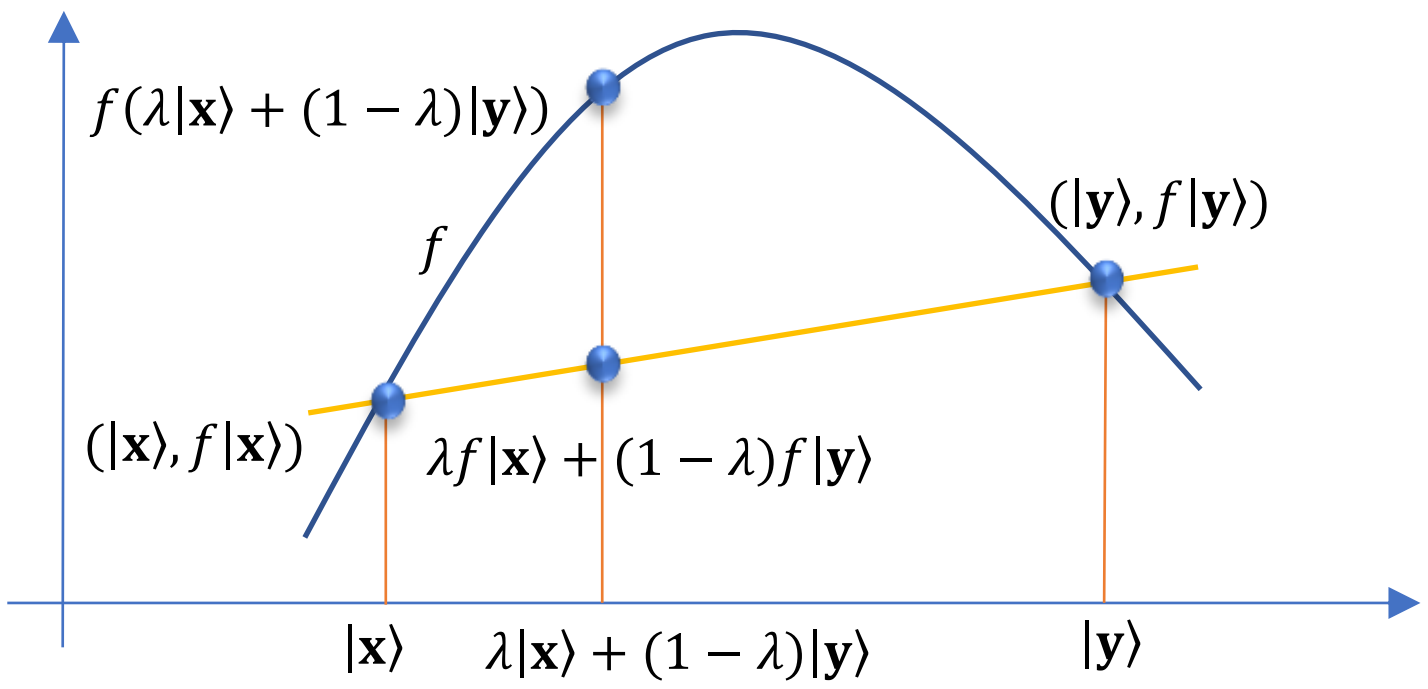

**Figure 4.22.** Concave function. The line segment lies below the graph of $f$.

**Definition (Strictly Convex Function):** A function $f$ is strictly convex if strict inequality holds

$$f(\lambda|\mathbf{x}\rangle + (1-\lambda)|\mathbf{y}\rangle) < \lambda f|\mathbf{x}\rangle + (1-\lambda) f|\mathbf{y}\rangle, \tag{4.59}$$

for all $|\mathbf{x}\rangle \neq |\mathbf{y}\rangle \in \Omega$ and $\lambda \in (0,1)$.

A strictly convex function $f$ is a function that the straight line between any pair of points on the curve $f$ is above the curve $f$ except for the intersection points between the straight line and the curve. We say $f$ is concave if $-f$ is convex and strictly concave if $-f$ is strictly convex.

It is often convenient to extend a convex function to all of $\mathbb{R}^n$ by defining its value to be $\infty$ outside its domain. If $f$ is convex, we define its extended-value extension $\tilde{f}: \mathbb{R}^n \to \mathbb{R} \cup \{\infty\}$ by

$$\tilde{f}|\mathbf{x}\rangle = \begin{cases} f|\mathbf{x}\rangle & |\mathbf{x}\rangle \in \text{dom } f, \\ \infty & |\mathbf{x}\rangle \notin \text{dom } f. \end{cases} \tag{4.60}$$

The extension $\tilde{f}$ is defined on all $\mathbb{R}^n$, and takes values in $\mathbb{R} \cup \{\infty\}$. We can recover the domain of the original function $f$ from the extension $\tilde{f}$ as $\text{dom } f = \{|\mathbf{x}\rangle : \tilde{f}|\mathbf{x}\rangle < \infty\}$. The extension can simplify notation since we do not need to explicitly describe the domain or add the qualifier for all $|\mathbf{x}\rangle \in \text{dom } \tilde{f}$ every time we refer to $f|\mathbf{x}\rangle$. Consider, for example, the basic defining inequality (4.58). In terms of the extension $\tilde{f}$, we can express it as: for $0 < \lambda < 1$,

$$\tilde{f}(\lambda|\mathbf{x}\rangle + (1-\lambda)|\mathbf{y}\rangle) \leq \lambda \tilde{f}|\mathbf{x}\rangle + (1-\lambda)\tilde{f}|\mathbf{y}\rangle, \tag{4.61}$$

for any $|\mathbf{x}\rangle$ and $|\mathbf{y}\rangle$. (For $\lambda = 0$ or $\lambda = 1$ the inequality always holds.)

Let us illustrate the convexity of functions by examples. In the following, we give a few more examples of convex and concave functions. We start with some functions on $\mathbb{R}$, with variable $x$.

***Example 4.4***

Is the function $f(x) = |x|$, $x \in \mathbb{R}$ convex? Is it strictly convex?

***Solution***

To see if $f$ is convex, we need to see if (4.58) is true; to see if it is strictly convex, we need to check (4.59). Furthermore, these inequalities have to be true for every $x_1, x_2 \in \mathbb{R}$, and every $\lambda \in [0,1]$.

$$\begin{aligned} f\big((1-\lambda)x_1 + \lambda x_2\big) &= |(1-\lambda)x_1 + \lambda x_2| \\ &\leq |(1-\lambda)x_1| + |\lambda x_2| \\ &= (1-\lambda)|x_1| + \lambda|x_2| \\ &= (1-\lambda)f(x_1) + \lambda f(x_2). \end{aligned}$$

Therefore (4.58) is satisfied, so $f$ is convex. To show that it is strictly convex, we would have to show that the inequality

$$|(1-\lambda)x_1 + \lambda x_2| \leq |(1-\lambda)x_1| + |\lambda x_2|,$$

can be replaced by strict inequality $<$. However, we can't do this: for example, if $x_1 = 1$, $x_2 = 2$, $\lambda = 0.5$, the left side of the inequality ($|\frac{1}{2} + \frac{2}{2}| = \frac{3}{2}$) is exactly equal to the right side ($|\frac{1}{2}| + |\frac{2}{2}| = \frac{3}{2}$). So $f$ is not strictly convex.

Note that proving that $f(x)$ is convex requires a general argument, whereas proving that $f(x)$ was not strictly convex only required a single counterexample. This is because the definition of convexity is a "for all" or "for every" type of argument. To prove convexity, you need an argument that allows for all possible values of $x_1$, $x_2$, and $\lambda$, whereas to disprove it, you only need to give one set of values where the necessary condition does not hold.

***Example 4.5***

Show that every affine function $f(x) = ax + b$, $x \in \mathbb{R}$ is convex, but not strictly convex.

***Solution***

$$\begin{aligned} f\big((1-\lambda)x_1 + \lambda x_2\big) &= a\big((1-\lambda)x_1 + \lambda x_2\big) + b \\ &= a\big((1-\lambda)x_1 + \lambda x_2\big) + \big((1-\lambda) + \lambda\big)b \\ &= (1-\lambda)(ax_1 + b) + \lambda(ax_2 + b) \\ &= (1-\lambda)f(x_1) + \lambda f(x_2). \end{aligned}$$

So, we see that inequality (4.58) is, in fact, satisfied as equality. That's fine, so every affine function in one variable is convex. However, this means we can't replace the inequality $\leq$ with the strict inequality $<$, so affine functions are not strictly convex.

***Example 4.6***

Show that $f(x) = x^2$, $x \in \mathbb{R}$ is strictly convex.

***Solution***

Pick $x_1$, $x_2$ so that $x_1 \neq x_2$, and pick $\lambda \in (0,1)$.

$$\begin{aligned} f\big((1-\lambda)x_1 + \lambda x_2\big) &= \big((1-\lambda)x_1 + \lambda x_2\big)^2 \\ &= (1-\lambda)^2 x_1^2 + \lambda^2 x_2^2 + 2(1-\lambda)\lambda x_1 x_2 \end{aligned}$$

Since $x_1 \neq x_2$, $(x_1 - x_2)^2 > 0$. Expanding, this means that $x_1^2 + x_2^2 > 2x_1x_2$. This means that

$$\begin{aligned} (1-\lambda)^2 x_1^2 + \lambda^2 x_2^2 + 2(1-\lambda)\lambda x_1 x_2 &< (1-\lambda)^2 x_1^2 + \lambda^2 x_2^2 + (1-\lambda)(\lambda)(x_1^2 + x_2^2) \\ &= (1 - 2\lambda + \lambda^2 + \lambda - \lambda^2)x_1^2 + (\lambda - \lambda^2 + \lambda^2)x_2^2 \\ &= (1-\lambda)x_1^2 + \lambda x_2^2 \\ &= (1-\lambda)f(x_1) + \lambda f(x_2), \end{aligned}$$

which proves strict convexity.

More examples of convex functions on $\mathbb{R}$ are

- Exponential: $e^{ax}$ is convex on $\mathbb{R}$, for any $a \in \mathbb{R}$.
- Powers: $x^a$ is convex on $\mathbb{R}_{++}$ when $a \geq 1$ or $a \leq 0$.
- Powers of absolute value: $|x|^p$, for $p \geq 1$, is convex on $\mathbb{R}$.
- Logarithm: $\log x$ is concave on $\mathbb{R}_{++}$ ($\mathbb{R}_{++}$: The set of positive real numbers).

***Example 4.7***

The following functions are convex:

$$f|\mathbf{x}\rangle = \langle \mathbf{a}|\mathbf{x}\rangle + b,$$

for $|\mathbf{x}\rangle \in \mathbb{R}^n$, where $|\mathbf{a}\rangle \in \mathbb{R}^n$ and $b \in \mathbb{R}$.

***Solution***

Indeed, the function $f$ is convex since

$$\begin{aligned} f(\lambda|\mathbf{x}\rangle + (1-\lambda)|\mathbf{y}\rangle) &= \langle \mathbf{a}|\lambda\mathbf{x} + (1-\lambda)\mathbf{y}\rangle + b \\ &= \lambda\langle \mathbf{a}|\mathbf{x}\rangle + (1-\lambda)\langle \mathbf{a}|\mathbf{y}\rangle + b \\ &= \lambda(\langle \mathbf{a}|\mathbf{x}\rangle + b) + (1-\lambda)(\langle \mathbf{a}|\mathbf{y}\rangle + b) \\ &= \lambda f|\mathbf{x}\rangle + (1-\lambda)f|\mathbf{y}\rangle, \end{aligned}$$

for all $|\mathbf{x}\rangle, |\mathbf{y}\rangle \in \mathbb{R}^n$ and $\lambda \in (0,1)$.

***Example 4.8***

The following function is convex:

$$g|\mathbf{x}\rangle = \|\mathbf{x}\| \text{ for } |\mathbf{x}\rangle \in \mathbb{R}^n.$$

***Solution***

The function $g$ is convex since for $|\mathbf{x}\rangle, |\mathbf{y}\rangle \in \mathbb{R}^n$ and $\lambda \in (0,1)$, we have

$$\begin{aligned} g(\lambda|\mathbf{x}\rangle + (1-\lambda)|\mathbf{y}\rangle) &= \|\lambda\mathbf{x} + (1-\lambda)\mathbf{y}\| \\ &\leq \lambda\|\mathbf{x}\| + (1-\lambda)\|\mathbf{y}\| \\ &= \lambda g|\mathbf{x}\rangle + (1-\lambda)g|\mathbf{y}\rangle, \end{aligned}$$

which follows from the triangle inequality and the fact that $\|\alpha\,\mathbf{u}\| = |\alpha|\,\|\mathbf{u}\|$, whenever $\alpha \in \mathbb{R}$ and $|\mathbf{u}\rangle \in \mathbb{R}^n$.

***Example 4.9***

Let $\mathbf{A}$ be an $n \times n$ symmetric matrix. It is called positive semidefinite if $\langle \mathbf{Au}|\mathbf{u}\rangle \geq 0$ for all $|\mathbf{u}\rangle \in \mathbb{R}^n$. Let us check that $\mathbf{A}$ is positive semidefinite if and only if the function $f: \mathbb{R}^n \to \mathbb{R}$ defined by

$$f|\mathbf{x}\rangle = \frac{1}{2}\langle \mathbf{Ax}|\mathbf{x}\rangle, \qquad |\mathbf{x}\rangle \in \mathbb{R}^n,$$

is convex.

***Solution***

Indeed, a direct calculation shows that for any $|\mathbf{x}\rangle, |\mathbf{y}\rangle \in \mathbb{R}^n$ and $\lambda \in (0,1)$ we have

$$\begin{aligned} &\lambda f|\mathbf{x}\rangle + (1-\lambda)f|\mathbf{y}\rangle - f(\lambda|\mathbf{x}\rangle + (1-\lambda)|\mathbf{y}\rangle) \\ &= \lambda\frac{1}{2}\langle \mathbf{Ax}|\mathbf{x}\rangle + (1-\lambda)\frac{1}{2}\langle \mathbf{Ay}|\mathbf{y}\rangle - \frac{1}{2}\langle \mathbf{A}(\lambda\mathbf{x} + (1-\lambda)\mathbf{y})|\lambda\mathbf{x} + (1-\lambda)\mathbf{y}\rangle \\ &= \frac{1}{2}\lambda\langle \mathbf{Ax}|\mathbf{x}\rangle + \frac{1}{2}(1-\lambda)\langle \mathbf{Ay}|\mathbf{y}\rangle - \frac{1}{2}\langle \mathbf{A}\lambda\mathbf{x} + \mathbf{A}(1-\lambda)\mathbf{y}|\lambda\mathbf{x} + (1-\lambda)\mathbf{y}\rangle \\ &= \frac{1}{2}\{\lambda\langle \mathbf{Ax}|\mathbf{x}\rangle + (1-\lambda)\langle \mathbf{Ay}|\mathbf{y}\rangle - \langle \mathbf{A}\lambda\mathbf{x}|\lambda\mathbf{x}\rangle - \langle \mathbf{A}\lambda\mathbf{x}|(1-\lambda)\mathbf{y}\rangle - \langle \mathbf{A}(1-\lambda)\mathbf{y}|\lambda\mathbf{x}\rangle \\ &\qquad - \langle \mathbf{A}(1-\lambda)\mathbf{y}|(1-\lambda)\mathbf{y}\rangle\} \\ &= \frac{1}{2}\{\lambda\langle \mathbf{Ax}|\mathbf{x}\rangle - \lambda^2\langle \mathbf{Ax}|\mathbf{x}\rangle + (1-\lambda)\langle \mathbf{Ay}|\mathbf{y}\rangle - (1-\lambda)^2\langle \mathbf{Ay}|\mathbf{y}\rangle - \lambda(1-\lambda)\langle \mathbf{Ax}|\mathbf{y}\rangle - \lambda(1-\lambda)\langle \mathbf{Ay}|\mathbf{x}\rangle\} \\ &= \frac{1}{2}\{\lambda(1-\lambda)\langle \mathbf{Ax}|\mathbf{x}\rangle + \lambda(1-\lambda)\langle \mathbf{Ay}|\mathbf{y}\rangle - \lambda(1-\lambda)\langle \mathbf{Ax}|\mathbf{y}\rangle - \lambda(1-\lambda)\langle \mathbf{Ay}|\mathbf{x}\rangle\} \\ &= \frac{1}{2}\lambda(1-\lambda)\{\langle \mathbf{Ax}|\mathbf{x}\rangle - \langle \mathbf{Ax}|\mathbf{y}\rangle + \langle \mathbf{Ay}|\mathbf{y}\rangle - \langle \mathbf{Ay}|\mathbf{x}\rangle\} \\ &= \frac{1}{2}\lambda(1-\lambda)\{\langle \mathbf{Ax}|\mathbf{x}-\mathbf{y}\rangle - \langle \mathbf{Ay}|\mathbf{x}-\mathbf{y}\rangle\} \\ &= \frac{1}{2}\lambda(1-\lambda)\langle \mathbf{A}(\mathbf{x}-\mathbf{y})|\mathbf{x}-\mathbf{y}\rangle = \frac{1}{2}\lambda(1-\lambda)\langle \mathbf{A}(\mathbf{x}-\mathbf{y})|\mathbf{x}-\mathbf{y}\rangle. \end{aligned}$$

If the matrix $\mathbf{A}$ is positive semidefinite, then

$$\langle \mathbf{A}(\mathbf{x}-\mathbf{y})|\mathbf{x}-\mathbf{y}\rangle \geq 0,$$

so the function $f$ is convex by (4.58). Conversely, assuming the convexity of $f$ and it is easy to verify that $\mathbf{A}$ is positive semidefinite.

For example, let $h(x) := x^2$ for $x \in \mathbb{R}$. The convexity of the simplest quadratic function $h$ follows from a more general result for the quadratic function on $\mathbb{R}^n$.

The following characterization of convexity is known as the Jensen inequality.

**Theorem 4.21 (Jensen inequality):** A function $f: \mathbb{R}^n \to \overline{\mathbb{R}}$ is convex if and only if for any numbers $\lambda_i \geq 0$ as $i = 1, \ldots, m$ with $\sum_{i=1}^m \lambda_i = 1$ and for any elements $|\mathbf{x}_i\rangle \in \mathbb{R}^n$, $i = 1, \ldots, m$, it holds that

$$f\left(\sum_{i=1}^m \lambda_i|\mathbf{x}_i\rangle\right) \leq \sum_{i=1}^m \lambda_i f|\mathbf{x}_i\rangle. \tag{4.62}$$

**Proof:**

Since (4.62) immediately implies the convexity of $f$, we only need to prove that any convex function $f$ satisfies the Jensen inequality (4.62). Arguing by induction and taking into account that for $m = 1$, inequality (4.62) holds trivially and for $m = 2$, inequality (4.62) holds by the definition of convexity, we suppose that it holds for an integer $m = k$ with $k \geq 2$. Fix numbers $\lambda_i \geq 0$, $i = 1, \dots, k+1$, with $\sum_{i=1}^{k+1} \lambda_i = 1$ and elements $|\mathbf{x}_i\rangle \in \mathbb{R}^n$, $i = 1, \dots, k+1$. We obviously have the relationship

$$\sum_{i=1}^{k} \lambda_i = 1 - \lambda_{k+1}.$$

If $\lambda_{k+1} = 1$, then $\lambda_i = 0$ for all $i = 1, \dots, k$ and (4.62) obviously holds for $m = k+1$ in this case. Supposing now that $0 \leq \lambda_{k+1} < 1$, we get

$$\sum_{i=1}^{k} \frac{\lambda_i}{1 - \lambda_{k+1}} = 1,$$

and by direct calculations based on convexity, arrive at

$$\begin{aligned} f\left(\sum_{i=1}^{k+1} \lambda_i |\mathbf{x}_i\rangle\right) &= f\left((1-\lambda_{k+1})\frac{\sum_{i=1}^{k} \lambda_i |\mathbf{x}_i\rangle}{1-\lambda_{k+1}} + \lambda_{k+1}|\mathbf{x}_{k+1}\rangle\right) \\ &\leq (1-\lambda_{k+1}) f\left(\frac{\sum_{i=1}^{k} \lambda_i |\mathbf{x}_i\rangle}{1-\lambda_{k+1}}\right) + \lambda_{k+1} f|\mathbf{x}_{k+1}\rangle \\ &= (1-\lambda_{k+1}) f\left(\sum_{i=1}^{k} \frac{\lambda_i}{1-\lambda_{k+1}} |\mathbf{x}_i\rangle\right) + \lambda_{k+1} f|\mathbf{x}_{k+1}\rangle \\ &\leq (1-\lambda_{k+1}) \sum_{i=1}^{k} \frac{\lambda_i}{1-\lambda_{k+1}} f|\mathbf{x}_i\rangle + \lambda_{k+1} f|\mathbf{x}_{k+1}\rangle \\ &= \sum_{i=1}^{k+1} \lambda_i f|\mathbf{x}_i\rangle. \end{aligned}$$

This justifies inequality (4.62) and completes the proof of the theorem.

■

**Definition (Domain and Epigraph of a Function):** The domain and epigraph of $f: \mathbb{R}^n \to \overline{\mathbb{R}}$ are defined, respectively, by

$$\text{dom } f = \{|\mathbf{x}\rangle \in \mathbb{R}^n : f|\mathbf{x}\rangle < \infty\}. \tag{4.63}$$

The graph of a function $f: \mathbb{R}^n \to \mathbb{R}$ is defined as

$$\text{graph } f = \{(|\mathbf{x}\rangle, f|\mathbf{x}\rangle) : |\mathbf{x}\rangle \in \text{dom } f\}. \tag{4.64}$$

The epigraph of a function $f: \mathbb{R}^n \to \mathbb{R}$ is defined as

$$\begin{aligned} \text{epi } f &= \{(|\mathbf{x}\rangle, t) \in \mathbb{R}^n \times \mathbb{R} : |\mathbf{x}\rangle \in \mathbb{R}^n, t \geq f|\mathbf{x}\rangle\} \\ &= \{(|\mathbf{x}\rangle, t) \in \mathbb{R}^n \times \mathbb{R} : |\mathbf{x}\rangle \in \text{dom } f, t \geq f|\mathbf{x}\rangle\}, \end{aligned} \tag{4.65}$$

(*i.e.*, all the points in the Cartesian product $\mathbb{R}^n \times \mathbb{R}$ lying on or above its graph.)

('Epi' means 'above,' so epigraph means 'above the graph.') The definition is illustrated in Figure 4.23.

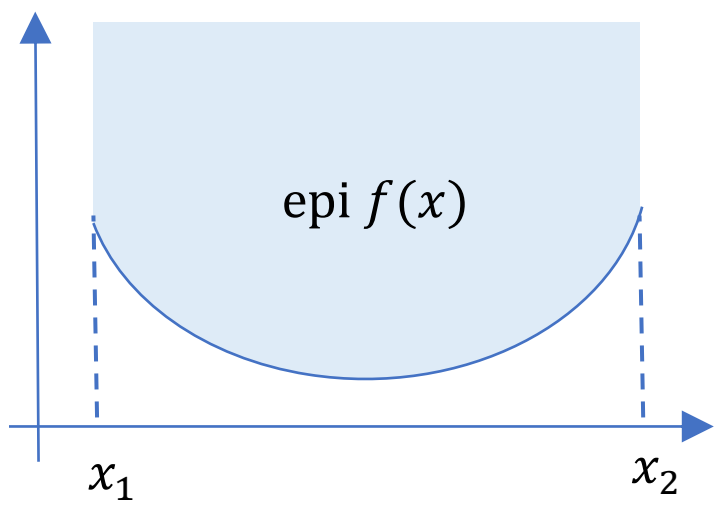

**Figure 4.23.** Epigraphs of a convex function.

The next theorem gives a geometric characterization of the function convexity via the convexity of the associated epigraphical set.

**Theorem 4.22:** A function $f: \mathbb{R}^n \to \overline{\mathbb{R}}$ is convex if and only if its epigraph, epi $f$, is a convex subset of the product space $\mathbb{R}^n \times \mathbb{R}$.

**Proof:**

Assuming that $f$ is convex, fix pairs $(|\mathbf{x}_1\rangle, t_1), (|\mathbf{x}_2\rangle, t_2) \in \text{epi } f$ and a number $\lambda \in (0,1)$. Then we have $f|\mathbf{x}_i\rangle \le t_i$ for $i = 1,2$. Thus, the convexity of $f$ ensures that

$$f(\lambda|\mathbf{x}_1\rangle + (1-\lambda)|\mathbf{x}_2\rangle) \le \lambda f|\mathbf{x}_1\rangle + (1-\lambda) f|\mathbf{x}_2\rangle \le \lambda t_1 + (1-\lambda)t_2.$$

This implies that

$$\begin{aligned}\lambda(|\mathbf{x}_1\rangle, t_1) + (1-\lambda)(|\mathbf{x}_2\rangle, t_2) &= (\lambda|\mathbf{x}_1\rangle, \lambda t_1) + \big((1-\lambda)|\mathbf{x}_2\rangle, (1-\lambda)t_2\big)\\ &= (\lambda|\mathbf{x}_1\rangle + (1-\lambda)|\mathbf{x}_2\rangle, \lambda t_1 + (1-\lambda)t_2)\\ &\in \text{epi } f,\end{aligned}$$

and thus, the epigraph, epi $f$, is a convex subset of $\mathbb{R}^n \times \mathbb{R}$.

Conversely, suppose that the set epi $f$ is convex and fix $|\mathbf{x}_1\rangle, |\mathbf{x}_2\rangle \in \text{dom } f$ and a number $\lambda \in (0,1)$. Then $(|\mathbf{x}_1\rangle, f|\mathbf{x}_1\rangle)$ , $(|\mathbf{x}_2\rangle, f|\mathbf{x}_2\rangle) \in \text{epi } f$. This tells us that

$$\lambda(|\mathbf{x}_1\rangle, f|\mathbf{x}_1\rangle) + (1-\lambda)(|\mathbf{x}_2\rangle, f|\mathbf{x}_2\rangle) = (\lambda|\mathbf{x}_1\rangle + (1-\lambda)|\mathbf{x}_2\rangle, \lambda f|\mathbf{x}_1\rangle + (1-\lambda) f|\mathbf{x}_2\rangle) \in \text{epi } f,$$

and thus implies the inequality

$$f(\lambda|\mathbf{x}_1\rangle + (1-\lambda)|\mathbf{x}_2\rangle) \le \lambda f|\mathbf{x}_1\rangle + (1-\lambda) f|\mathbf{x}_2\rangle,$$

which justifies the convexity of the function $f$.

■

Now we show that convexity is preserved under some important operations.

**Theorem 4.23:** Let $f: \mathbb{R}^n \to \overline{\mathbb{R}}$ be convex functions for all $i = 1, \dots, m$. Then the following functions are convex as well:

(i) The multiplication by scalars $\lambda f$ for any $\lambda > 0$.

(ii) The sum function $\sum_{i=1}^{m} f_i$.

(iii) The maximum function $\max_{1 \le i \le m} f_i$.

**Proof:**

The convexity of $\lambda f$ in (i) follows directly from the definition. It is sufficient to prove (ii) and (iii) for $m = 2$, since the general cases immediately follow by induction.

(ii) Fix any $|\mathbf{x}\rangle, |\mathbf{y}\rangle \in \mathbb{R}^n$ and $\lambda \in (0,1)$. Then we have

$$\begin{aligned}(f_1 + f_2)(\lambda|\mathbf{x}\rangle + (1-\lambda)|\mathbf{y}\rangle) &= f_1(\lambda|\mathbf{x}\rangle + (1-\lambda)|\mathbf{y}\rangle) + f_2(\lambda|\mathbf{x}\rangle + (1-\lambda)|\mathbf{y}\rangle)\\ &\le \lambda f_1|\mathbf{x}\rangle + (1-\lambda) f_1|\mathbf{y}\rangle + \lambda f_2|\mathbf{x}\rangle + (1-\lambda) f_2|\mathbf{y}\rangle\\ &= \lambda(f_1 + f_2)|\mathbf{x}\rangle + (1-\lambda)(f_1 + f_2)|\mathbf{y}\rangle,\end{aligned}$$

which verifies the convexity of the sum function $f_1 + f_2$.

(iii) Denote $g := \max\{f_1, f_2\}$ and get for any $|\mathbf{x}\rangle, |\mathbf{y}\rangle \in \mathbb{R}^n$ and $\lambda \in (0,1)$ that

$$f_i(\lambda\mathbf{x} + (1-\lambda)\mathbf{y}) \le \lambda f_i(\mathbf{x}) + (1-\lambda) f_i(\mathbf{y}) \le \lambda g(\mathbf{x}) + (1-\lambda) g(\mathbf{y}),$$

for $i = 1,2$. This directly implies that

$$\begin{aligned}g(\lambda|\mathbf{x}\rangle + (1-\lambda)|\mathbf{y}\rangle) &= \max\{f_1(\lambda|\mathbf{x}\rangle + (1-\lambda)|\mathbf{y}\rangle), f_2(\lambda|\mathbf{x}\rangle + (1-\lambda)|\mathbf{y}\rangle)\}\\ &\le \lambda g|\mathbf{x}\rangle + (1-\lambda) g|\mathbf{y}\rangle,\end{aligned}$$

which shows that the maximum function $g|\mathbf{x}\rangle = \max\{f_1|\mathbf{x}\rangle, f_2|\mathbf{x}\rangle\}$ is convex.

■

### Convexity conditions, case 1, $f$ is a function of one variable only

**Theorem 4.24:** Given a convex function $f:\overline{\mathbb{R}} \to \overline{\mathbb{R}}$, assume that its domain is an open interval $I$. For any $a, b \in I$ and $a < x < b$, we have the inequalities

$$\frac{f(x)-f(a)}{x-a} \leq \frac{f(b)-f(a)}{b-a} \leq \frac{f(b)-f(x)}{b-x}. \tag{4.66}$$

**Proof:**

Fix $a, b, x$ as above and form the numbers $t = \frac{x-a}{b-a} \in (0,1)$. Then

$$\begin{aligned} f(x) &= f\big(a + (x-a)\big) \\ &= f\left(a + \frac{x-a}{b-a}(b-a)\right) \\ &= f\big(a + t(b-a)\big) \\ &= f(tb + (1-t)a). \end{aligned}$$

This gives us the inequalities

$$f(x) \leq tf(b) + (1-t)f(a),$$

and

$$\begin{aligned} f(x) - f(a) &\leq tf(b) + (1-t)f(a) - f(a) \\ &= t[f(b) - f(a)] \\ &= \frac{x-a}{b-a}\big(f(b) - f(a)\big), \end{aligned}$$

which can be equivalently written as

$$\frac{f(x)-f(a)}{x-a} \leq \frac{f(b)-f(a)}{b-a}.$$

Similarly, we have the estimate

$$\begin{aligned} f(x) - f(b) &\leq tf(b) + (1-t)f(a) - f(b) \\ &= (t-1)[f(b) - f(a)] \\ &= \frac{x-b}{b-a}\big(f(b) - f(a)\big), \end{aligned}$$

which finally implies that

$$\frac{f(b)-f(a)}{b-a} \leq \frac{f(b)-f(x)}{b-x}.$$

and thus completes the proof of the theorem.

■

**Theorem 4.25:** Suppose that $f:\mathbb{R} \to \overline{\mathbb{R}}$ is differentiable on its domain, which is an open interval $I$. Then $f$ is convex if and only if the derivative $f'$ is nondecreasing on $I$.

**Proof:**

Suppose that $f$ is convex and fix $a < b$ with $a, b \in I$. By Theorem 4.24, we have

$$\frac{f(x)-f(a)}{x-a} \leq \frac{f(b)-f(a)}{b-a},$$

for any $x \in (a, b)$. This implies by the derivative definition that

$$f'(a) \leq \frac{f(b)-f(a)}{b-a}.$$

Similarly, we arrive at the estimate

$$\frac{f(b)-f(a)}{b-a} \le f'(b),$$

and conclude that $f'(a) \le f'(b)$, i.e., $f'$ is a nondecreasing function.

To prove the converse, suppose that $f'$ is nondecreasing on $I$ and fix $x_1 < x_2$ with $x_1, x_2 \in I$ and $t \in (0,1)$. Then

$$x_1 < x_t < x_2 \text{ for } x_t := tx_1 + (1-t)x_2.$$

By the mean value theorem, we find $c_1$, $c_2$ such that $x_1 < c_1 < x_t < c_2 < x_2$ and

$$f(x_t) - f(x_1) = f'(c_1)(x_t - x_1) = f'(c_1)(1-t)(x_2 - x_1),$$
$$f(x_t) - f(x_2) = f'(c_2)(x_t - x_2) = f'(c_2)t(x_1 - x_2).$$

This can be equivalently rewritten as

$$tf(x_t) - tf(x_1) = f'(c_1)t(1-t)(x_2 - x_1) \quad \Rightarrow \frac{tf(x_t) - tf(x_1)}{t(1-t)(x_2 - x_1)} = f'(c_1),$$

$$(1-t)f(x_t) - (1-t)f(x_2) = f'(c_2)t(1-t)(x_1 - x_2) \quad \Rightarrow \frac{(1-t)f(x_t) - (1-t)f(x_2)}{t(1-t)(x_1 - x_2)} = f'(c_2).$$

Since $f'(c_1) \le f'(c_2)$, we have

$$\Rightarrow \frac{tf(x_t) - tf(x_1)}{t(1-t)(x_2 - x_1)} \le \frac{(1-t)f(x_t) - (1-t)f(x_2)}{t(1-t)(x_1 - x_2)},$$
$$\Rightarrow \frac{f(x_t) - f(x_1)}{(1-t)} \le \frac{f(x_2) - f(x_t)}{t},$$
$$\Rightarrow tf(x_t) + (1-t)f(x_t) \le tf(x_1) + (1-t)f(x_2),$$
$$\Rightarrow f(x_t) \le tf(x_1) + (1-t)f(x_2),$$
$$\Rightarrow f(tx_1 + (1-t)x_2) \le tf(x_1) + (1-t)f(x_2),$$

which justifies the convexity of the function $f$.

∎

As a result, we have

**Theorem 4.26:** Let $f: \mathbb{R} \to \overline{\mathbb{R}}$ be a function and let $f$ be differentiable on its domain. Then $f$ is convex if and only if

$$f(x_2) \ge f(x_1) + f'(x_1)(x_2 - x_1), \tag{4.67}$$

for all $x_1, x_2 \in \mathbb{R}$.

**Proof:**

First, we will assume $f$ is convex and try to prove the inequality. Take any $x_1, x_2 \in \mathbb{R}$, and assume $x_1 \neq x_2$ because otherwise, the inequality is already satisfied: it just says $f(x) \ge f(x)$. We have

$$f\big((1-t)x_1 + tx_2\big) \le (1-t)f(x_1) + tf(x_2),$$

whenever $t \in (0,1)$, which we can rewrite as

$$f\big(x_1 + t(x_2 - x_1)\big) \le (1-t)f(x_1) + tf(x_2),$$
$$\Rightarrow f\big(x_1 + t(x_2 - x_1)\big) - f(x_1) \le tf(x_2) - tf(x_1),$$
$$\Rightarrow \frac{f\big(x_1 + t(x_2 - x_1)\big) - f(x_1)}{t(x_2 - x_1)}(x_2 - x_1) \le f(x_2) - f(x_1).$$

If we take the limit as $t \to 0$, then $t(x_2 - x_1) \to 0$ as well, which means the left-hand side of this inequality approaches $f'(x_1)(x_2 - x_1)$. The right-hand side does not depend on $t$, so it remains the same, and we get

$$f'(x_1)(x_2 - x_1) \le f(x_2) - f(x_1),$$
$$\Rightarrow f(x_2) \ge f(x_1) + f'(x_1)(x_2 - x_1).$$

Next, we will assume that the inequality holds and try to prove that $f$ is convex. Let $u, v \in \mathbb{R}$ and let $w = tu + (1 - t)v$ with $t \in (0,1)$. Then we have

$$f(u) \geq f(w) + f'(w)(u - w),$$
$$f(v) \geq f(w) + f'(w)(v - w),$$

So, if we add $t$ times the first inequality and $(1 - t)$ times the second inequality, we get

$$\begin{aligned} tf(u) + (1-t)f(v) &\geq tf(w) + (1-t)f(w) + f'(w)(tu - tw + (1-t)v - (1-t)w) \\ &= f(w) + f'(w)(tu + (1-t)v - w) \\ &= f(w) + f'(w)(w - w) = f(w), \end{aligned}$$

and since $w = tu + (1 - t)v$, this is exactly the inequality

$$tf(u) + (1-t)f(v) \geq f(tu + (1-t)v),$$

that proves that $f$ is convex.

■

**Theorem 4.27:** Let $f: \mathbb{R} \to \overline{\mathbb{R}}$ be twice differentiable on its domain, which is an open interval $I$. Then $f$ is convex if and only if $f''(x) \geq 0$ for all $x \in I$.

**Proof:**

Since $f''(x) \geq 0$ for all $x \in I$ if and only if the derivative function $f'$ is nondecreasing on this interval. Then the conclusion follows directly from Theorem 4.25.

■

Much simpler! If $f$ is differentiable (or, better yet, twice differentiable), checking these conditions is almost always easier.

Equivalent conditions for strict convexity can be obtained in a natural way, changing $\geq$ to $>$ and requiring that $x_1$, and $x_2$ be distinct in Theorem 4.26. Essentially, Theorem 4.26 says that $f$ lies above its tangent lines, see Figure 4.24, while Theorem 4.27 says that $f$ is always "curving upward." (A convex function lies above its tangents, but below its secants). These conditions are usually easier to verify than that of the definition of convexity.

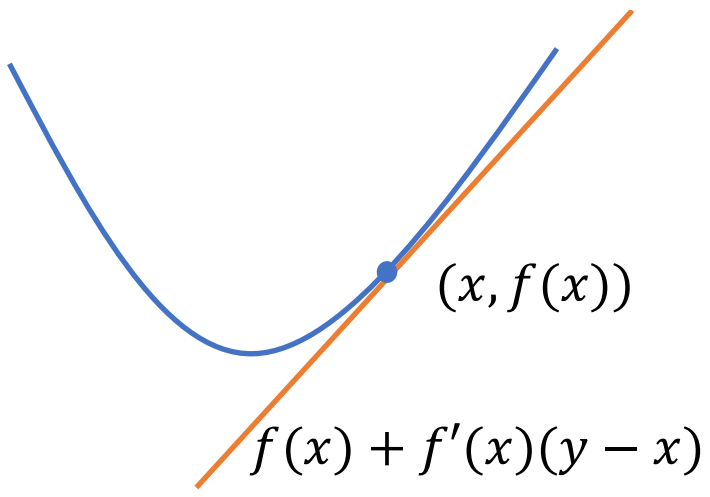

**Figure 4.24.** A convex function lies above its tangents

***Example 4.10***

Show that the following function is convex on $\mathbb{R}$

$$f(x) = \begin{cases} 1/x & \text{if } x > 0, \\ \infty & \text{otherwise.} \end{cases}$$

***Solution***

To verify its convexity, we get that $f''(x) = \frac{2}{x^3} > 0$ for all $x$ belonging to the domain of $f$, which is $I = (0, \infty)$. Thus, this function is convex on $\mathbb{R}$ by Theorem 4.27.

***Example 4.11***

Show that $f(x) = x^2$ is strictly convex using Theorem 4.26.

***Solution***

Pick any $x_1, x_2 \in \mathbb{R}$ with $x_1 \neq x_2$. We have $f'(x_1) = 2x_1$, so we need to show that

$$x_2^2 > x_1^2 + 2x_1(x_2 - x_1).$$

Expanding the right-hand side and rearranging terms, we see this is equivalent to

$$x_1^2 - 2x_1x_2 + x_2^2 > 0,$$

or

$$(x_1 - x_2)^2 > 0,$$

which is clearly true since $x_1 \neq x_2$. Thus, $f$ is strictly convex.

***Example 4.12***

Show that $f(x) = x^2$ is strictly convex using Theorem 4.27.

***Solution***

$f''(x) = 2 > 0$ for all $x \in \mathbb{R}$, so $f$ is strictly convex.

## 4.7 General Case (Multi-Variable Function) and Sub-gradient of Convex Functions

If $f$ is a differentiable function, then $|\mathbf{y}\rangle = f|\mathbf{a}\rangle + \langle \nabla f(\mathbf{a})|\mathbf{x} - \mathbf{a}\rangle$, is the equation of a hyperplane that is tangent to the surface $|\mathbf{y}\rangle = f|\mathbf{x}\rangle$ at the point $(|\mathbf{a}\rangle, f|\mathbf{a}\rangle)$. And if $f$ is also convex, then $f|\mathbf{x}\rangle \geq f|\mathbf{a}\rangle + \langle \nabla f(\mathbf{a})|\mathbf{x} - \mathbf{a}\rangle$ for all $|\mathbf{x}\rangle$ in the domain of the function, so the tangent plane lies below the graph of the function and is a supporting hyperplane of the epigraph.

The epigraph of an arbitrary convex function is a convex set by definition. Hence, through each boundary point belonging to the epigraph of a convex function there passes a supporting hyperplane. The supporting hyperplanes of a convex one-variable function $f$, defined on an open interval, which says that the line $y = f(x_0) + a(x - x_0)$ supports the epigraph at the point $(x_0, f(x_0))$ if (and only if) $f'_-(x_0) \leq a \leq f'_+(x_0)$. The existence of supporting hyperplanes characterizes convexity.

What if $f$ is not differentiable at $|\mathbf{x}_t\rangle$ (i.e., $\nabla f|\mathbf{x}_t\rangle$ does not exist)? Can we still apply this algorithm? In this section, we are going to answer this question. We will introduce the sub-gradient, which is a concept closely related to the gradient. And for convex functions, we will show that even if the gradient does not exist, the sub-gradient always exists. Then, in convex optimization problems, when the function is non-differential at a certain point, we can use its sub-gradient as an alternative to the gradient. First, we give the formal definition of sub-gradient [1-4].

**Definition (Sub-gradient):** A vector $|\mathbf{g}\rangle \in \mathbb{R}^d$ is a sub-gradient of $f : \mathbb{R}^d \longrightarrow \mathbb{R}$ at $|\mathbf{x}\rangle \in \text{dom}\, f$, if for any $|\mathbf{y}\rangle \in \text{dom} f$, see Figure 4.25, we have

$$f|\mathbf{y}\rangle \geq f|\mathbf{x}\rangle + \langle \mathbf{g}|\mathbf{y} - \mathbf{x}\rangle. \tag{4.68}$$

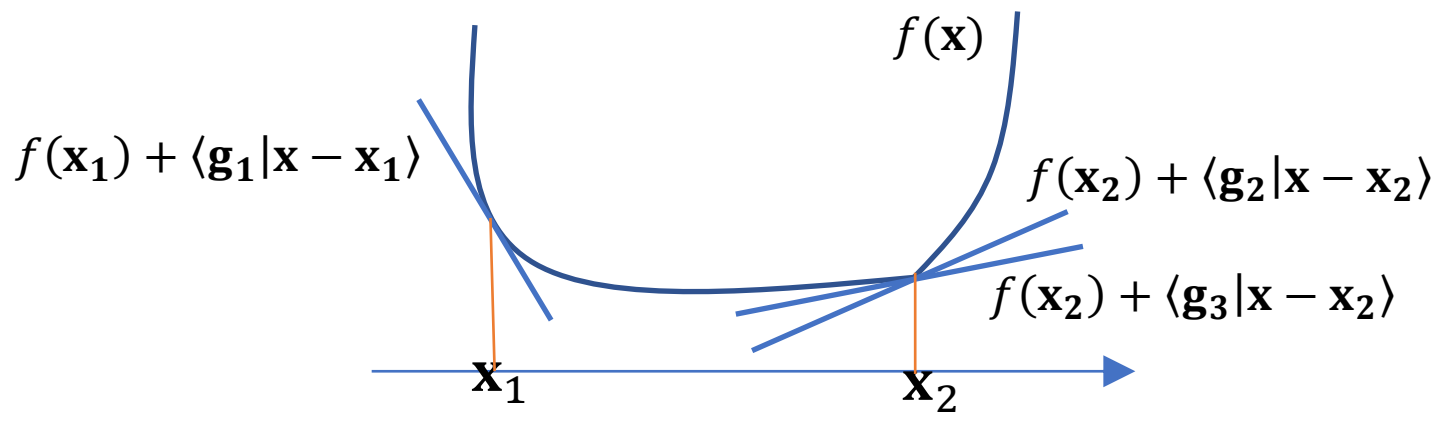

**Figure 4.25.** The sub-gradient of a nondifferentiable convex function. Where $|\mathbf{g}_1\rangle$ is a sub-gradient at point $|\mathbf{x}_1\rangle$ and $|\mathbf{g}_2\rangle$ and $|\mathbf{g}_3\rangle$ are sub-gradients at point $|\mathbf{x}_2\rangle$. Notice that when the function is differentiable, such as at point $|\mathbf{x}_1\rangle$, the sub-gradient, $|\mathbf{g}_1\rangle$, just becomes the gradient to the function.

**Remark:**

Clearly, if $f$ is differentiable at $|\mathbf{x}\rangle \in \operatorname{dom} f$, then its sub-gradient at $|\mathbf{x}\rangle$, $|\mathbf{g}\rangle$ is equal to $\nabla f|\mathbf{x}\rangle$. A function $f$ is called sub-differentiable at $|\mathbf{x}\rangle \in \operatorname{dom} f$ if there exists at least one sub-gradient at $|\mathbf{x}\rangle$.

**Definition (Sub-differential):** The set of sub-gradients of $f$ at $|\mathbf{x}\rangle \in \operatorname{dom} f$ is called the sub-differential of $f$ at $|\mathbf{x}\rangle$, denoted by $\partial f|\mathbf{x}\rangle$.

**Theorem 4.28:** The sub-differential of $f: \mathbb{R}^d \longrightarrow \mathbb{R}$ at $|\mathbf{x}\rangle \in \operatorname{dom} f$ is a closed and convex set.

**Proof:**

For any $|\mathbf{g}_1\rangle, |\mathbf{g}_2\rangle \in \partial f|\mathbf{x}\rangle$ and any $\alpha \in [0,1]$, we want to show $\alpha|\mathbf{g}_1\rangle + (1-\alpha)|\mathbf{g}_2\rangle \in \partial f|\mathbf{x}\rangle$. By definition of sub-gradients, we have

$$f|\mathbf{y}\rangle \geq f|\mathbf{x}\rangle + \langle \mathbf{g}_1|\mathbf{y}-\mathbf{x}\rangle,$$
$$f|\mathbf{y}\rangle \geq f|\mathbf{x}\rangle + \langle \mathbf{g}_2|\mathbf{y}-\mathbf{x}\rangle,$$

Multiplying the first inequality by $\alpha$ and the second one by $(1-\alpha)$, and then adding them together leads to

$$\alpha f|\mathbf{y}\rangle + (1-\alpha) f|\mathbf{y}\rangle \geq \alpha f|\mathbf{x}\rangle + (1-\alpha) f|\mathbf{x}\rangle + \langle \alpha \mathbf{g}_1 + (1-\alpha)\mathbf{g}_2|\mathbf{y}-\mathbf{x}\rangle$$
$$= f|\mathbf{x}\rangle + \langle \alpha \mathbf{g}_1 + (1-\alpha)\mathbf{g}_2|\mathbf{y}-\mathbf{x}\rangle.$$

Thus, by definition, we have $\alpha|\mathbf{g}_1\rangle + (1-\alpha)|\mathbf{g}_2\rangle \in \partial f(\mathbf{x})$

■

The following theorem is very useful and is simple to proof.

**Theorem 4.29:** A point $|\mathbf{x}^*\rangle$ is a minimizer of a convex function $f: \mathbb{R}^d \longrightarrow \mathbb{R}$, if and only if $|\mathbf{0}\rangle \in \partial f(\mathbf{x}^*)$.

**Proof:**

If $|\mathbf{x}^*\rangle$ is the minimizer of $f$, then for any $|\mathbf{y}\rangle \in \operatorname{dom} f$, we have $f|\mathbf{y}\rangle \geq f|\mathbf{x}^*\rangle$, which can be rewritten as

$$f|\mathbf{y}\rangle \geq f|\mathbf{x}^*\rangle + \langle \mathbf{0}|\mathbf{y}-\mathbf{x}^*\rangle.$$

Therefore, by definition, $|\mathbf{0}\rangle$ is a sub-gradient of $f$ at $|\mathbf{x}^*\rangle$ and $|\mathbf{0}\rangle \in \partial f(\mathbf{x}^*)$.

Conversely, if $|\mathbf{0}\rangle \in \partial f(\mathbf{x}^*)$, then for any $|\mathbf{y}\rangle \in \operatorname{dom} f$, we have

$$f|\mathbf{y}\rangle \geq f|\mathbf{x}^*\rangle + \langle \mathbf{0}|\mathbf{y}-\mathbf{x}^*\rangle = f|\mathbf{x}^*\rangle.$$

Clearly, by definition, $|\mathbf{x}^*\rangle$ is a minimizer of $f$.

■

**Theorem 4.30:** Let $f: \mathbb{R}^n \longrightarrow \overline{\mathbb{R}}$ convex and let $|\tilde{\mathbf{x}}\rangle \in \operatorname{dom} f$ be a local minimizer of $f$. Then $f$ attains its global minimum at this point.

**Proof:**

Since $|\tilde{\mathbf{x}}\rangle$ is a local minimizer of $f$, there is $\delta > 0$ such that

$$f|\mathbf{u}\rangle \geq f|\tilde{\mathbf{x}}\rangle,$$

for all $|\mathbf{u}\rangle \in IB(\tilde{\mathbf{x}}, \delta)$, where $IB$ means a closed unit ball. Fix $|\mathbf{x}\rangle \in \mathbb{R}^n$ and construct a sequence of $|\mathbf{x}_k\rangle = (1-k^{-1})|\tilde{\mathbf{x}}\rangle + k^{-1}|\tilde{\mathbf{x}}\rangle$ as $k \in \mathbb{N}$. Thus, we have $|\mathbf{x}_k\rangle \in IB(\tilde{\mathbf{x}}, \delta)$ when $k$ is sufficiently large. It follows from the convexity of $f$ that

$$f|\tilde{\mathbf{x}}\rangle \leq f|\mathbf{x}_k\rangle \leq (1-k^{-1}) f|\tilde{\mathbf{x}}\rangle + k^{-1} f|\mathbf{x}\rangle,$$

which readily implies that $k^{-1} f|\tilde{\mathbf{x}}\rangle \leq k^{-1} f|\mathbf{x}\rangle$, and hence $f|\tilde{\mathbf{x}}\rangle \leq f|\mathbf{x}\rangle$ for all $|\mathbf{x}\rangle \in \mathbb{R}^n$.

■

**Note that:** $|\tilde{\mathbf{x}}\rangle$ is an intrtior point of $\operatorname{dom} f$, $|\tilde{\mathbf{x}}\rangle \in \operatorname{int}(\operatorname{dom} f)$, if there is $\delta > 0$ such that $IB(|\tilde{\mathbf{x}}\rangle, \delta) \subset \operatorname{dom} f$.

**Theorem 4.31:** Let $f: \mathbb{R}^n \longrightarrow \overline{\mathbb{R}}$ be convex and differentiable at $|\tilde{\mathbf{x}}\rangle \in \text{int}\,(\text{dom}\, f)$. Then we have

$$\partial f|\tilde{\mathbf{x}}\rangle = \{\nabla f|\tilde{\mathbf{x}}\rangle\}, \tag{4.69}$$

and

$$f|\mathbf{x}\rangle - f|\tilde{\mathbf{x}}\rangle \geq \langle \nabla f(\tilde{\mathbf{x}})|\mathbf{x} - \tilde{\mathbf{x}}\rangle \quad \text{for all } |\mathbf{x}\rangle \in \mathbb{R}^n. \tag{4.70}$$

**Proof:**

It follows from the differentiability of $f$ at $|\tilde{\mathbf{x}}\rangle$ that for any $\epsilon > 0$ there is $\delta > 0$ with

$$-\epsilon\|\mathbf{x} - \tilde{\mathbf{x}}\| \leq f|\mathbf{x}\rangle - f|\tilde{\mathbf{x}}\rangle - \langle \nabla f(\tilde{\mathbf{x}})|\mathbf{x} - \tilde{\mathbf{x}}\rangle \leq \epsilon\|\mathbf{x} - \tilde{\mathbf{x}}\| \text{ whenever } \|\mathbf{x} - \tilde{\mathbf{x}}\| < \delta. \tag{*}$$

Consider further the convex function

$$\varphi|\mathbf{x}\rangle = f|\mathbf{x}\rangle - f|\tilde{\mathbf{x}}\rangle - \langle \nabla f(\tilde{\mathbf{x}})|\mathbf{x} - \tilde{\mathbf{x}}\rangle + \epsilon\|\mathbf{x} - \tilde{\mathbf{x}}\|, \qquad |\mathbf{x}\rangle \in \mathbb{R}^n,$$

and observe that $\varphi|\mathbf{x}\rangle \geq \varphi|\tilde{\mathbf{x}}\rangle = 0$ for all $|\mathbf{x}\rangle \in IB(\tilde{\mathbf{x}}, \delta)$. The convexity of $\varphi$ ensures that $\varphi|\mathbf{x}\rangle \geq \varphi|\tilde{\mathbf{x}}\rangle$ for all $|\mathbf{x}\rangle \in \mathbb{R}^n$. Thus

$$\langle \nabla f(\tilde{\mathbf{x}})|\mathbf{x} - \tilde{\mathbf{x}}\rangle \leq f|\mathbf{x}\rangle - f|\tilde{\mathbf{x}}\rangle + \epsilon\|\mathbf{x} - \tilde{\mathbf{x}}\|, \quad \text{whenever } |\mathbf{x}\rangle \in \mathbb{R}^n,$$

which yields (4.70) by letting $\epsilon \downarrow 0$.

It follows from (4.70) that $\nabla f|\tilde{\mathbf{x}}\rangle \in \partial f|\tilde{\mathbf{x}}\rangle$. Picking now $|\mathbf{v}\rangle \in \partial f|\tilde{\mathbf{x}}\rangle$, we get

$$\begin{aligned}\langle \mathbf{v}|\mathbf{x} - \tilde{\mathbf{x}}\rangle \leq f|\mathbf{x}\rangle - f|\tilde{\mathbf{x}}\rangle \quad &\Longrightarrow \quad \langle \mathbf{v}|\mathbf{x} - \tilde{\mathbf{x}}\rangle - \langle \nabla f(\tilde{\mathbf{x}})|\mathbf{x} - \tilde{\mathbf{x}}\rangle \leq f|\mathbf{x}\rangle - f|\tilde{\mathbf{x}}\rangle - \langle \nabla f(\tilde{\mathbf{x}})|\mathbf{x} - \tilde{\mathbf{x}}\rangle \\ &\Longrightarrow \quad \langle \mathbf{v} - \nabla f(\tilde{\mathbf{x}})|\mathbf{x} - \tilde{\mathbf{x}}\rangle \leq f|\mathbf{x}\rangle - f|\tilde{\mathbf{x}}\rangle - \langle \nabla f(\tilde{\mathbf{x}})|\mathbf{x} - \tilde{\mathbf{x}}\rangle\end{aligned}$$

Then the second part of (*) gives us that

$$\langle \mathbf{v} - \nabla f(\tilde{\mathbf{x}})|\mathbf{x} - \tilde{\mathbf{x}}\rangle \leq \epsilon\|\mathbf{x} - \tilde{\mathbf{x}}\|, \quad \text{whenever } \|\mathbf{x} - \tilde{\mathbf{x}}\| < \delta.$$

Finally, we observe that $\|\mathbf{v} - \nabla f(\tilde{\mathbf{x}})\| \leq \epsilon$, which yields $|\mathbf{v}\rangle = \nabla f|\tilde{\mathbf{x}}\rangle$, since $\epsilon > 0$ was chosen arbitrarily. Thus $\partial f|\tilde{\mathbf{x}}\rangle = \{\nabla f|\tilde{\mathbf{x}}\rangle\}$.

■

**Theorem 4.32:** Let $f: \mathbb{R}^n \longrightarrow \overline{\mathbb{R}}$ be a strictly convex function. Then the differentiability of $f$ at $|\tilde{\mathbf{x}}\rangle \in \text{int}\,(\text{dom}\, f)$ implies the strict inequality

$$f|\mathbf{x}\rangle - f|\tilde{\mathbf{x}}\rangle > \langle \nabla f(\tilde{\mathbf{x}})|\mathbf{x} - \tilde{\mathbf{x}}\rangle \quad \text{whenever } |\mathbf{x}\rangle \neq |\tilde{\mathbf{x}}\rangle. \tag{4.71}$$

**Proof:**

Since $f$ is convex, we get from Theorem 4.31 that

$$\langle \nabla f(\tilde{\mathbf{x}})|\mathbf{u} - \tilde{\mathbf{x}}\rangle \leq f|\mathbf{u}\rangle - f|\tilde{\mathbf{x}}\rangle \quad \text{for all } |\mathbf{u}\rangle \in \mathbb{R}^n.$$

Fix $|\mathbf{x}\rangle \neq |\tilde{\mathbf{x}}\rangle$ and let $|\mathbf{u}\rangle = \frac{1}{2}(|\mathbf{x}\rangle + |\tilde{\mathbf{x}}\rangle)$. It follows from the above that

$$\left\langle \nabla f(\tilde{\mathbf{x}}) \middle| \frac{1}{2}(\mathbf{x} + \tilde{\mathbf{x}}) - \tilde{\mathbf{x}} \right\rangle \leq f\left|\frac{1}{2}(\mathbf{x} + \tilde{\mathbf{x}})\right\rangle - f|\tilde{\mathbf{x}}\rangle < \frac{1}{2} f|\mathbf{x}\rangle + \frac{1}{2} f|\tilde{\mathbf{x}}\rangle - f|\tilde{\mathbf{x}}\rangle.$$

Thus, we arrive at the strict inequality (4.71) and complete the proof.

■

**Remarks:**

1. The function $f|\mathbf{y}\rangle$ given by $f|\mathbf{x}\rangle + \langle \nabla f(\mathbf{x})|\mathbf{y} - \mathbf{x}\rangle$ is, of course, the first-order Taylor approximation of $f$ near $|\mathbf{x}\rangle$. The inequality (4.70) states that for a convex function, the first-order Taylor approximation is, in fact, a global under-estimator of the function, see Figure 4.24. Conversely, if the first-order Taylor approximation of a function is always a global under-estimator of the function, then the function is convex.
2. The inequality (4.70) shows that from local information about a convex function (i.e., its value and derivative at a point) we can derive global information (i.e., a global under-estimator of it). This is perhaps the most important property of convex functions and explains some of the remarkable properties of convex functions and convex optimization problems. As one simple example, the inequality (4.70) shows that if $\nabla f|\mathbf{x}\rangle = 0$, then for all $|\mathbf{y}\rangle \in \text{dom}\, f$, $f|\mathbf{y}\rangle \geq f|\mathbf{x}\rangle$, *i.e.*, $|\mathbf{x}\rangle$ is a global minimizer of the function $f$.

***Example 4.13***

Let $p|\mathbf{x}\rangle = \|\mathbf{x}\|$ be the Euclidean norm function on $\mathbb{R}^n$. Then we have

$$\partial p|\mathbf{x}\rangle = \begin{cases} IB & \text{if } |\mathbf{x}\rangle = |\mathbf{0}\rangle, \\ |\mathbf{x}\rangle/\|\mathbf{x}\| & \text{otherwise.} \end{cases}$$

***Solution***

To verify this, observe first that the Euclidean norm function $p$ is differentiable at any nonzero point with

$$\begin{aligned} \nabla p|\mathbf{x}\rangle &= \nabla(x_1^2 + \cdots + x_n^2)^{\frac{1}{2}} \\ &= \frac{1}{2}(x_1^2 + \cdots + x_n^2)^{-\frac{1}{2}} \times 2 \times \begin{pmatrix} x_1 \\ \vdots \\ x_n \end{pmatrix} \\ &= |\mathbf{x}\rangle/\|\mathbf{x}\|, \end{aligned}$$

as $|\mathbf{x}\rangle \neq |\mathbf{0}\rangle$.

It remains to calculate its subdifferential at $|\mathbf{x}\rangle = |\mathbf{0}\rangle$. To proceed by definition (4.68), we have that $|\mathbf{v}\rangle \in \partial p|\mathbf{0}\rangle$ if and only if

$$\langle\mathbf{v}|\mathbf{x}\rangle = \langle\mathbf{v}|\mathbf{x} - \mathbf{0}\rangle \le p|\mathbf{x}\rangle - p|\mathbf{0}\rangle = \|\mathbf{x}\| \text{ for all } |\mathbf{v}\rangle \in \mathbb{R}^n.$$

Letting $|\mathbf{x}\rangle = |\mathbf{v}\rangle$ gives us

$$\langle\mathbf{v}|\mathbf{v}\rangle = \|\mathbf{v}\|^2 \le \|\mathbf{v}\|,$$

which implies that $\|\mathbf{v}\| \le 1$, i.e., $|\mathbf{v}\rangle \in IB$. Now take $|\mathbf{v}\rangle \in IB$ and deduce from the Cauchy-Schwarz inequality that

$$\begin{aligned} \langle\mathbf{v}|\mathbf{x} - \mathbf{0}\rangle &= \langle\mathbf{v}|\mathbf{x}\rangle \\ &\le \|\mathbf{v}\|\,\|\mathbf{x}\| \\ &\le \|\mathbf{x}\| \\ &= p|\mathbf{x}\rangle - p|\mathbf{0}\rangle \ \text{ for all } |\mathbf{x}\rangle \in \mathbb{R}^n, \end{aligned}$$

and thus $|\mathbf{v}\rangle \in \partial p|\mathbf{0}\rangle$, which shows that $\partial p|\mathbf{0}\rangle = IB$.

**Theorem 4.33:** Let $f: \mathbb{R}^n \longrightarrow \overline{\mathbb{R}}$ be a differentiable function on its domain $D$, which is an open convex set. Then $f$ is convex if and only if

$$\langle\nabla f(\mathbf{u})|\mathbf{x} - \mathbf{u}\rangle \le f|\mathbf{x}\rangle - f|\mathbf{u}\rangle \quad \text{for all } |\mathbf{u}\rangle, |\mathbf{x}\rangle \in D. \tag{4.72}$$

**Proof:**

The "only if" part follows from Theorem 4.31. To justify the converse, suppose that (4.72) holds and then fix any $|\mathbf{x}_1\rangle, |\mathbf{x}_2\rangle \in D$ and $t \in (0,1)$. Denoting $|\mathbf{x}_t\rangle = t|\mathbf{x}_1\rangle + (1-t)|\mathbf{x}_2\rangle$, we have $|\mathbf{x}_t\rangle \in D$ by the convexity of $D$. Then

$$\begin{aligned} \langle\nabla f(\mathbf{x}_t)|\mathbf{x}_1 - \mathbf{x}_t\rangle &\le f|\mathbf{x}_1\rangle - f|\mathbf{x}_t\rangle, \\ \langle\nabla f(\mathbf{x}_t)|\mathbf{x}_2 - \mathbf{x}_t\rangle &\le f|\mathbf{x}_2\rangle - f|\mathbf{x}_t\rangle. \end{aligned}$$

It follows furthermore that

$$\begin{aligned} t\langle\nabla f(\mathbf{x}_t)|\mathbf{x}_1 - \mathbf{x}_t\rangle &\le tf|\mathbf{x}_1\rangle - tf|\mathbf{x}_t\rangle, \\ (1-t)\langle\nabla f(\mathbf{x}_t)|\mathbf{x}_2 - \mathbf{x}_t\rangle &\le (1-t)f|\mathbf{x}_2\rangle - (1-t)f|\mathbf{x}_t\rangle. \end{aligned}$$

Summing up these inequalities, we arrive at

$$0 \le tf|\mathbf{x}_1\rangle + (1-t)f|\mathbf{x}_2\rangle - f|\mathbf{x}_t\rangle,$$

which ensures that $f|\mathbf{x}_t\rangle \le tf|\mathbf{x}_1\rangle + (1-t)f|\mathbf{x}_2\rangle$, and so verifies the convexity of $f$.

■

**Theorem 4.34:** Let $f: \mathbb{R}^n \longrightarrow \overline{\mathbb{R}}$ be convex and twice continuously differentiable on an open subset $V$ of its domain containing $\tilde{\mathbf{x}}$. Then we have

$$\langle\nabla^2 f(\tilde{\mathbf{x}})\mathbf{u}|\mathbf{u}\rangle \ge 0 \quad \text{for all } |\mathbf{u}\rangle \in \mathbb{R}^n, \tag{4.73}$$

where $\nabla^2 f|\tilde{\mathbf{x}}\rangle$ is the Hessian matrix of $f$ at $|\tilde{\mathbf{x}}\rangle$.

**Proof:**

Let $\mathbf{A} = \nabla^2 f|\tilde{\mathbf{x}}\rangle$, which is a symmetric matrix. Then

$$\lim_{|\mathbf{h}\rangle\to|\mathbf{0}\rangle} \frac{f|\tilde{\mathbf{x}}+\mathbf{h}\rangle - f|\tilde{\mathbf{x}}\rangle - \langle\nabla f(\tilde{\mathbf{x}})|\mathbf{h}\rangle - \frac{1}{2}\langle\mathbf{Ah}|\mathbf{h}\rangle}{\|\mathbf{h}\|^2} = 0. \quad (*)$$

It follows from (*) that for any $\epsilon > 0$ there is $\delta > 0$ such that

$$-\epsilon\|\mathbf{h}\|^2 \le f|\tilde{\mathbf{x}}+\mathbf{h}\rangle - f|\tilde{\mathbf{x}}\rangle - \langle\nabla f(\tilde{\mathbf{x}})|\mathbf{h}\rangle - \frac{1}{2}\langle\mathbf{Ah}|\mathbf{h}\rangle \le \epsilon\|\mathbf{h}\|^2 \text{ for all } \|\mathbf{h}\| \le \delta.$$

By Theorem 4.31, we readily have

$$f|\tilde{\mathbf{x}}+\mathbf{h}\rangle - f|\tilde{\mathbf{x}}\rangle - \langle\nabla f(\tilde{\mathbf{x}})|\mathbf{h}\rangle \ge 0.$$

Combining the above inequalities ensures that

$$-\epsilon\|\mathbf{h}\|^2 \le \frac{1}{2}\langle\mathbf{Ah}|\mathbf{h}\rangle \text{ whenever } \|\mathbf{h}\| \le \delta. \quad (**)$$

Observe further that for any $|\mathbf{0}\rangle \ne |\mathbf{u}\rangle \in \mathbb{R}^n$ the element $|\mathbf{h}\rangle = \frac{|\mathbf{u}\rangle}{\|\mathbf{u}\|}$ satisfies $\|\mathbf{h}\| \le \delta$ and, being substituted into (**), gives us the estimate

$$-\epsilon\delta^2 \le \frac{1}{2}\delta^2\left\langle \mathbf{A}\frac{\mathbf{u}}{\|\mathbf{u}\|}\middle|\frac{\mathbf{u}}{\|\mathbf{u}\|}\right\rangle.$$

It shows therefore (since the case where $|\mathbf{u}\rangle = |\mathbf{0}\rangle$ is trivial) that

$$-2\epsilon\|\mathbf{u}\|^2 \le \langle\mathbf{Au}|\mathbf{u}\rangle \text{ whenever } |\mathbf{u}\rangle \in \mathbb{R}^n,$$

which implies by letting $\epsilon \downarrow 0$ that $0 \le \langle\mathbf{Au}|\mathbf{u}\rangle$ for all $|\mathbf{u}\rangle \in \mathbb{R}^n$.

■

**Theorem 4.35:** Let $f\colon \mathbb{R}^n \longrightarrow \overline{\mathbb{R}}$ be a function twice continuously differentiable on its domain $D$, which is an open convex subset of $\mathbb{R}^n$. Then $f$ is convex if and only if $\nabla^2 f|\tilde{\mathbf{x}}\rangle$ is positive semidefinite for every $|\tilde{\mathbf{x}}\rangle \in D$.

**Proof:**

Taking Theorem 4.34 into account, we only need to verify that if $\nabla^2 f|\tilde{\mathbf{x}}\rangle$ is positive semidefinite for every $|\tilde{\mathbf{x}}\rangle \in D$, then $f$ is convex. To proceed, for any $|\mathbf{x}_1\rangle, |\mathbf{x}_2\rangle \in D$ define $|\mathbf{x}_t\rangle = t|\mathbf{x}_1\rangle + (1-t)|\mathbf{x}_2\rangle$ and consider the function

$$\varphi(t) := f(t|\mathbf{x}_1\rangle + (1-t)|\mathbf{x}_2\rangle) - tf|\mathbf{x}_1\rangle - (1-t)f|\mathbf{x}_2\rangle, \qquad t \in \mathbb{R}.$$

It is clear that $\varphi$ is well defined on an open interval $I$ containing $(0,1)$. Then

$$\varphi'(t) = \langle\nabla f(\mathbf{x}_t)|\mathbf{x}_1 - \mathbf{x}_2\rangle - f|\mathbf{x}_1\rangle + f|\mathbf{x}_2\rangle,$$

and

$$\varphi''(t) = \langle\nabla^2 f(\mathbf{x}_t)(\mathbf{x}_1 - \mathbf{x}_2)|\mathbf{x}_1 - \mathbf{x}_2\rangle \ge 0,$$

for every $t \in I$ since $\nabla^2 f|\mathbf{x}_t\rangle$ is positive semidefinite. Hence, $\varphi$ is convex on $I$. Since $\varphi(0) = \varphi(1) = 0$, for any $t \in (0,1)$ we have

$$\varphi(t) = \varphi(t(1) + (1-t)0) \le t\varphi(1) + (1-t)\varphi(0) = 0,$$

which implies, in turn, the inequality

$$f(t|\mathbf{x}_1\rangle + (1-t)|\mathbf{x}_2\rangle) \le tf|\mathbf{x}_1\rangle + (1-t)f|\mathbf{x}_2\rangle.$$

This justifies that the function $f$ is convex on its domain and, thus, on $\mathbb{R}^n$.

■

For a function on $\mathbb{R}$, Theorem 4.35 reduces to the simple condition $\mathrm{f}''(x) \ge 0$ (and dom $f$ convex,), which means that the derivative is nondecreasing. The condition $\nabla^2 f \succcurlyeq 0$ can be interpreted geometrically as the requirement that the graph of the function has positive (upward) curvature at $x$.

***Example 4.14***

Determine whether the following functions are convex, concave or neither:

1. $f: \mathbb{R} \to \mathbb{R}, f(x) = -8x^2$;
2. $f: \mathbb{R}^3 \to \mathbb{R}, f|\mathbf{x}\rangle = 4x_1^2 + 3x_2^2 + 5x_3^2 + 6x_1x_2 + x_1x_3 - 3x_1 - 2x_2 + 15$;
3. $f: \mathbb{R}^3 \to \mathbb{R}, f|\mathbf{x}\rangle = 2x_1x_2 - x_1^2 - x_2^2$.

***Solution***

1. We will use Theorem 4.35. We first compute the Hessian, which in this case is just the second derivative:

$$\frac{d^2}{dx^2} f(x) = -16 < 0$$

for all $x \in \mathbb{R}$. Hence, $f$ is concave over $\mathbb{R}$.

2. The Hessian matrix of $f$ is

$$\mathbf{H} = \begin{pmatrix} 8 & 6 & 1 \\ 6 & 6 & 0 \\ 1 & 0 & 10 \end{pmatrix}.$$

The leading principal minors of **H** are

$$\Delta_1 = 8 > 0,$$
$$\Delta_2 = \det \begin{pmatrix} 8 & 6 \\ 6 & 6 \end{pmatrix} = 12 > 0,$$
$$\Delta_2 = \det \begin{pmatrix} 8 & 6 & 1 \\ 6 & 6 & 0 \\ 1 & 0 & 10 \end{pmatrix} = 114 > 0,$$

Hence, **H** is positive definite for all $|\mathbf{x}\rangle \in \mathbb{R}^3$. Therefore, $f$ is a convex function over $\mathbb{R}^3$.

3. The Hessian of $f$ is

$$\mathbf{H} = \begin{pmatrix} -2 & 2 \\ 2 & -2 \end{pmatrix}.$$

which is negative semi-definite for all $|\mathbf{x}\rangle \in \mathbb{R}^2$. Hence, $f$ is concave on $\mathbb{R}^2$.

We now turn to an important theorem, which shows that for a convex function, its sub-differential is always nonempty set.

**Theorem 4.36 (Existence of sub-gradient for convex functions):** If $f: \mathbb{R}^d \longrightarrow \mathbb{R}$ is convex, and $|\mathbf{x}\rangle \in$ int dom $f$, then $\partial f|\mathbf{x}\rangle$ is nonempty.

Now let us see some examples of sub-gradients.

***Example 4.15***

For $f(x) = |x|$,

$$\partial f(x) = \begin{cases} \{1\} & \text{if } x > 0, \\ \{-1\} & \text{if } x < 0, \\ [-1.1] & \text{if } x = 0. \end{cases}$$

***Solution***

Note that $f(x)$ is differentiable when $x > 0$ or $x < 0$, so the subgradient of $f$ is equal to its gradient (1 and $-1$, respectively). At $x = 0$, for any $y \in \mathbb{R}$, its sub-gradient $g$ should satisfy

$$|y| \geq |0| + g(y - 0). \qquad \text{i.e.,} \qquad g \cdot y \leq |y|.$$

Thus, $g \in [-1.1]$.

***Example 4.16***

Let $f|\mathbf{x}\rangle = \|\mathbf{x}\|_1 = \sum_{i=1}^d |x_i|$, and $|\mathbf{g}\rangle = (g_1. \cdots . g_d)^T$ be a sub-gradient of $f$ at $|\mathbf{x}\rangle$.

***Solution***

Then we have

$$g_i = \begin{cases} 1 & \text{if } x_i > 0, \\ -1 & \text{if } x_i < 0. i = 1. \cdots . d, \\ \in [-1.1] & \text{if } x_i = 0. \end{cases}$$

Finally, consider the subdifferential of the point-wise maximum function.

***Example 4.17***

Let $f|\mathbf{x}\rangle = \max(f_1|\mathbf{x}\rangle, f_2|\mathbf{x}\rangle)$, where $f_1$ and $f_2$ are convex and differentiable.

***Solution***

Then we have

$$\partial f|\mathbf{x}\rangle = \begin{cases} \{\nabla f_1|\mathbf{x}\rangle\} & \text{if } f_1|\mathbf{x}\rangle > f_2|\mathbf{x}\rangle, \\ \{\nabla f_2|\mathbf{x}\rangle\} & \text{if } f_2|\mathbf{x}\rangle > f_1|\mathbf{x}\rangle, \\ \mathrm{conv}(\langle\nabla f_1(\mathbf{x})|\nabla f_2(\mathbf{x})\rangle) & \text{if } f_1|\mathbf{x}\rangle = f_2|\mathbf{x}\rangle. \end{cases}$$

Clearly, when $f_1|\mathbf{x}\rangle \neq f_2|\mathbf{x}\rangle$, $f|\mathbf{x}\rangle$ is equal to either $f_1|\mathbf{x}\rangle$ or $f_2|\mathbf{x}\rangle$ and is therefore differentiable.
When $f_1|\mathbf{x}\rangle = f_2|\mathbf{x}\rangle$, both $\nabla f_1|\mathbf{x}\rangle$ and $\nabla f_2|\mathbf{x}\rangle$ are sub-gradients of $f|\mathbf{x}\rangle$. Since subdifferential is a convex set, $\partial f|\mathbf{x}\rangle$ should be a convex hull of $\nabla f_1|\mathbf{x}\rangle$ and $\nabla f_2|\mathbf{x}\rangle$.

It is essential to calculate the sub-gradient for non-differential functions. We have the following rules in calculating the sub-gradient.

### 1. Non-negative Scaling

Suppose $f|\mathbf{x}\rangle$ is convex, its domain $\mathrm{dom}f$ and that $\alpha > 0$, then

$$\partial(\alpha f|\mathbf{x}\rangle) = \alpha\partial f|\mathbf{x}\rangle. \tag{4.74}$$

Note that here " = " indicates set equality, i.e., for two sets $A$ and $B$, $A = B$ if and only if $A \subseteq B$ and $B \subseteq A$.

### 2. Summation

Suppose $f|\mathbf{x}\rangle = \sum_{i=1}^{n} f_i|\mathbf{x}\rangle$ where all of the $f_i's$ are convex. Then

$$\partial f|\mathbf{x}\rangle = \sum_{i=1}^{n} \partial f_i|\mathbf{x}\rangle = \partial f_1|\mathbf{x}\rangle + \partial f_2|\mathbf{x}\rangle + \cdots + \partial f_n|\mathbf{x}\rangle. \tag{7.75}$$

Note that here "+" indicates set sum, i.e., for two sets $A$ and $B$, $\mathbf{A} + \mathbf{B} = \{(a+b): a \in \mathbf{A}, \mathrm{b} \in \boldsymbol{B}\}$.

## 4.8 Mathematica Built-in Functions

The Wolfram Language provides built-in functions for many standard distance measures, [7,8].

| | |
|---|---|
| `Norm[expr]` | gives the norm of a number, vector, or matrix. |
| `Norm[expr,p]` | gives the p-norm. |
| `Normalize[v]` | gives the normalized form of a vector v. |
| `UnitVector[k]` | gives the two-dimensional unit vector in the $k^{th}$ direction. |
| `UnitVector[n,k]` | gives the n-dimensional unit vector in the $k^{th}$ direction. |
| `Abs[z]` | gives the absolute value of the real or complex number z. |
| `EuclideanDistance[u,v]` | gives the Euclidean distance between vectors u and v. |
| `SquaredEuclideanDistance[u,v]` | gives the squared Euclidean distance between vectors u and v. |
| `NormalizedSquaredEuclideanDistance[u,v]` | gives the normalized squared Euclidean distance between vectors u and v. |
| `RootMeanSquare[list]` | gives the root mean square of values in list. |
| `ManhattanDistance[u,v]` | gives the Manhattan or "city block" distance between vectors u and v. |
| `ChessboardDistance[u,v]` | gives the chessboard, Chebyshev, or sup norm distance between vectors u and v. |

***Mathematica Examples 3.1***

```
Input     (*Norm of a vector:*)
          Norm[{x,y,z}]
Output    √(Abs[x]^2+√Abs[y]^2+√Abs[z]^2)
```

```
Input     (*The p-norm:*)
          Norm[{x,y,z},p]
Output    (Abs[x]^p+√Abs[y]^p+√Abs[z]^p)^(1/p)

Input     (*The Infinity-norm:*)
          Norm[{x,y,z},Infinity]
           Max[Abs[x],Abs[y],Abs[z]]
           Normalize[{1,5,1}]
Output    {1/(3 √3),5/(3 √3),1/(3 √3)}

Input     (*Symbolic vectors:*)
          Normalize[{x,y}]
Output    {x/√(Abs[x]^2+√Abs[y]^2),y/√(Abs[x]^2+√Abs[y]^2)}

Input     (*The unit vector in the direction in three dimensions:*)
          UnitVector[3,2]
Output    {0,1,0}

Input     (*Euclidean distance between two vectors:*)
          EuclideanDistance[{a,b,c},{x,y,z}]
Output    √(Abs[a-x]^2+√Abs[b-y]^2+√Abs[c-z]^2)

Input     (*Euclidean distance between numeric vectors:*)
          EuclideanDistance[{1,2,3},{2,4,6}]
Output    √14

Input     (*Squared Euclidean distance between numeric vectors:*)
          SquaredEuclideanDistance[{a,b,c},{x,y,z}]
Output    Abs[a-x]2+Abs[b-y]2+Abs[c-z]2

Input     (*Manhattan distance between two vectors:*)
          ManhattanDistance[{a,b,c},{x,y,z}]
Output    Abs[a-x]+Abs[b-y]+Abs[c-z]

Input     (*Manhattan distance between numeric vectors:*)
          ManhattanDistance[{1,2,3},{2,4,6}]
Output    6

Input     (*The chessboard distance between two vectors:*)
          ChessboardDistance[{a,b,c},{x,y,z}]
Output    Max[Abs[a-x],Abs[b-y],Abs[c-z]]

Input     (*The normalized squared Euclidean distance between two vectors:*)
          NormalizedSquaredEuclideanDistance[{a,b},{x,y}]
Output    (Abs[a+1/2 (-a-b)-x+(x+y)/2]2+Abs[1/2 (-a-b)+b-y+(x+y)/2]2)/(2 (Abs[a+1/2 (-a-
          b)]2+Abs[1/2 (-a-b)+b]2+Abs[x+1/2 (-x-y)]2+Abs[1/2 (-x-y)+y]2))

Input     (*RootMeanSquare of a list:*)
          RootMeanSquare[{a,b,c,d}]
Output    1/2 √(a^2+b^2+c^2)
```

| | |
|---|---|
| `Hyperplane[n,p]` | represents the hyperplane with normal n passing through the point p. |
| `Hyperplane[n,c]` | represents the hyperplane with normal n given by the points that satisfy n.x=c. |
| `HalfSpace[n,p]` | represents the half-space of points x such that n.(x-p)<=0. |
| `HalfSpace[n,c]` | represents the half-space of points x such that n.x<=c. |

***Mathematica Examples 3.1***

```
Input       (*A Hyperplane in 2D:*)
            Graphics[
             Hyperplane[{2,1}]
             ]

            (*And in 3D:*)
            Graphics3D[
             Hyperplane[{-1,-1,1},{1,2,3}]
             ]
```

Output

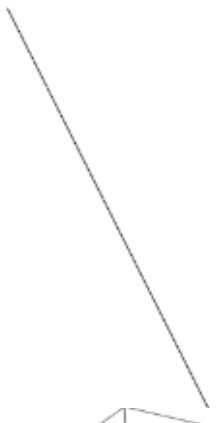

Output

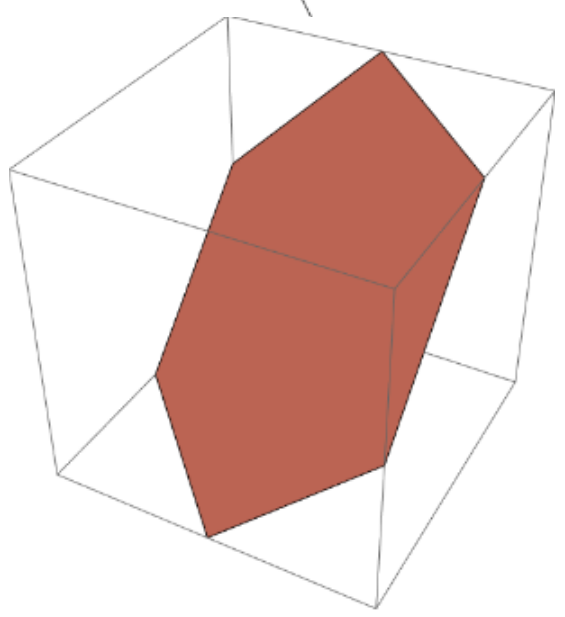

```
Input       (*A hyperplane in 2D defined by a normal vector and a point:*)
            ill=Graphics[
               {PointSize[Medium],Point[{{0,0}}],Arrowheads[Medium],Thick,Arrow[{{0,0},{-
            1,1}}]},
               PlotRange->2,
               Axes->True
               ];
            Show[
             Graphics[
              {Blue,Hyperplane[{-1,1},{0,0}]}
              ],
             ill,
             PlotRange->2
             ]
```

Output

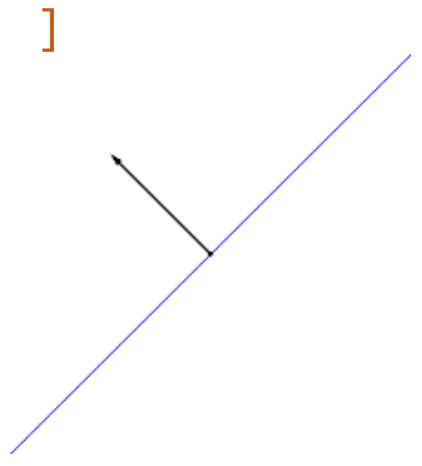

```
Input       (*Define a hyperplane in 3D using a normal vector and a constant:*)
            ill=Graphics3D[
               {Arrowheads[Medium],Thick,Arrow[{{1,0,0},{2,0,0}}]},
               PlotRange->2,
               Axes->True
               ];
            Show[
             ill,
             Graphics3D[
              Hyperplane[{1,0,0},1]
```

```
          ]
         ]
```

Output

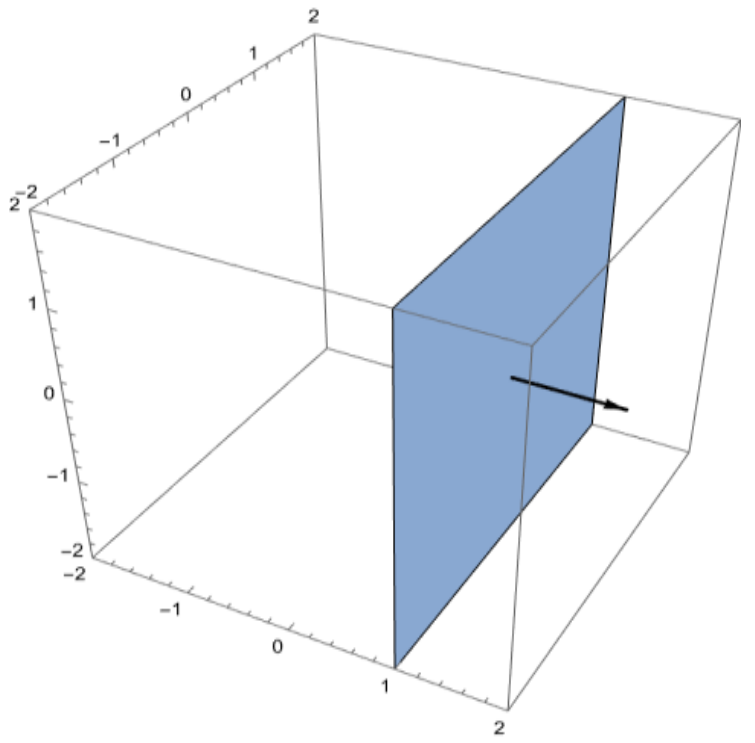

Input

```
(*Orthogonal arrangements of hyperplanes:*)
Graphics3D[
{Opacity[0.5],{Red,Hyperplane[{1,0,0},0]},{Green,Hyperplane[{0,1,0},0]},{Blue,Hyperplane[{0,0,1},0]}},
 Lighting->"Neutral"
 ]
```

Output

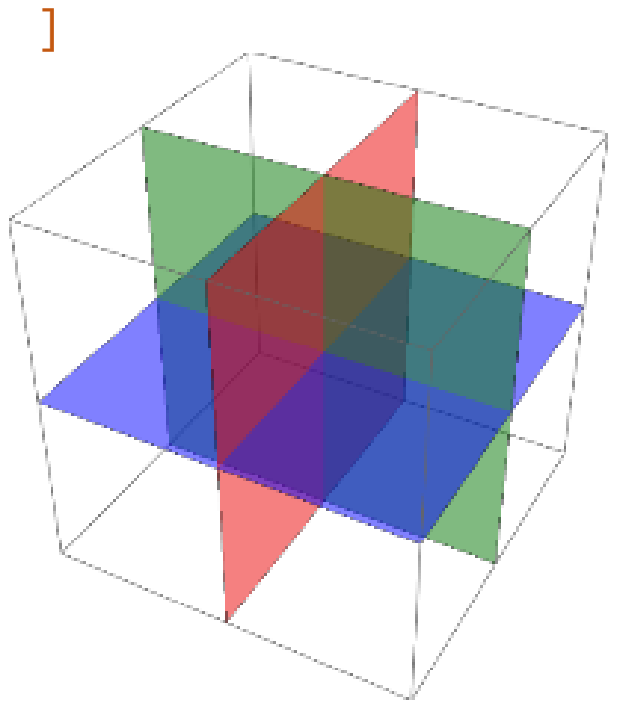

Input

```
(*A tangent plane to an implicitly defined curve in 2D is given by its normal Gradf(x,y) at a point on the curve.Start by finding points on the curve:*)
f=x^2/9+y^2/4-1;
pts=FindInstance[f==0,{x,y},Reals,2]
(*Find tangent lines at each of the points:*)
tangents=Hyperplane[Grad[f,{x,y}],{x,y}]/. pts
(*Visualize the solution:*)
Show[
 {ContourPlot[
   f==0,
   {x,-4,4},
   {y,-3,3},
   LabelStyle->Directive[Black,20]
   ],
  Graphics[
   {tangents,{Red,Point[{x,y}/. pts]}}
   ]
  },
 AspectRatio->Automatic,
 PlotRange->All
 ]
```

Output

```
{{x->5/34,y->-(√10379/51)},{x->-(101/34),y->√203/51}}
```

Output

```
{Hyperplane[{5/153,-( √10379/102)},{5/34,-( √10379/51)}],Hyperplane[{-(101/153),√203/102},{-(101/34), √203/51}]}
```

Output

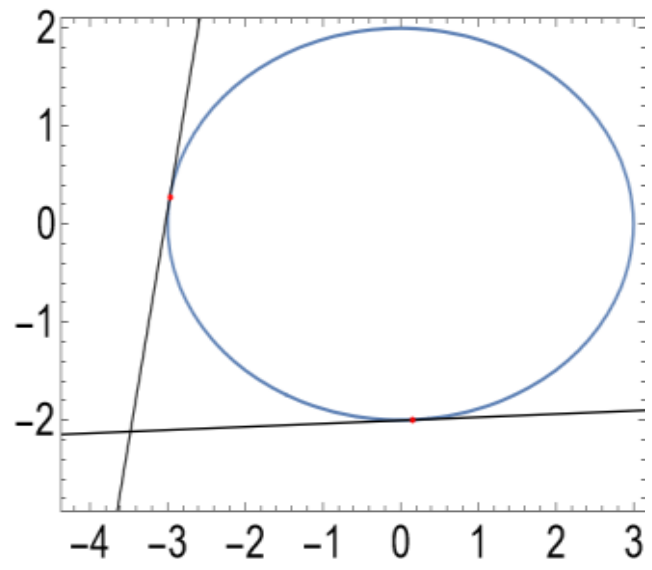

Input

```
(*A HalfSpace in 2D:*)
Graphics[
 HalfSpace[{2,1},{0,0}]
 ]
(*And in 3D:*)
Graphics3D[
 HalfSpace[{-1,-1,1},{0,0,0}]
 ]
```

Output

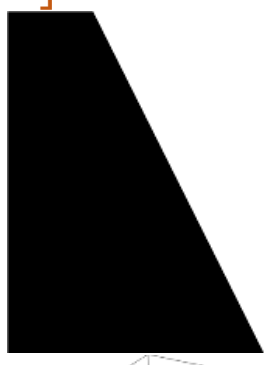

Output

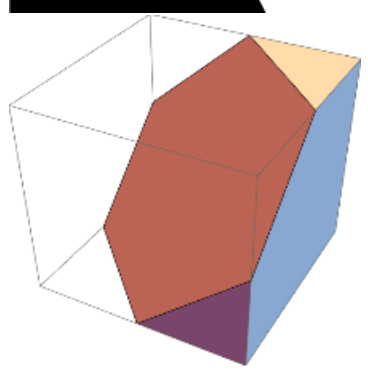

Input

```
(*A half-space in 2D defined by a normal vector and a point:*)
ill=Graphics[
   {PointSize[Medium],Point[{{0,0}}],Arrowheads[Medium],Thick,Arrow[{{0,0},{-
1,1}}]]},
   PlotRange->2,
   Axes->True
   ];
Show[
 Graphics[
  {Blue,HalfSpace[{-1,1},{0,0}]}
  ],
 ill,
 PlotRange->2
 ]
```

Output

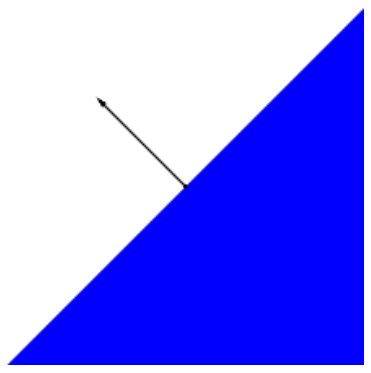

Input

```
(*Define a half-space in 3D using a normal vector and a constant:*)
ill=Graphics3D[
   {Arrowheads[Medium],Thick,Arrow[{{1,0,2},{-1,0,2}}]]},
   PlotRange->2,
   Axes->True
   ];
```

```
Show[
 ill,
 Graphics3D[
  HalfSpace[{-1,0,0},-1]
  ]
 ]
```

Output

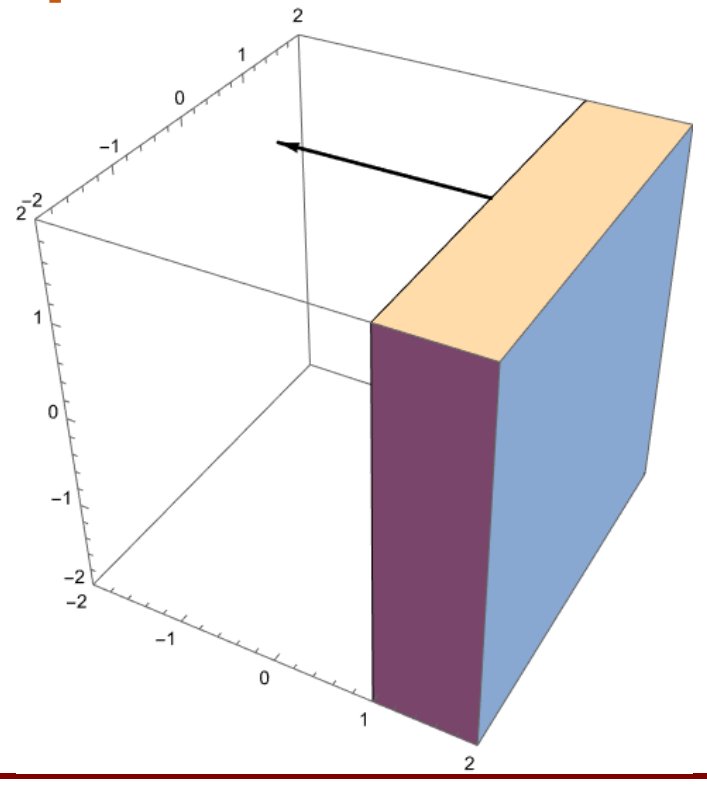

| | |
|---|---|
| `Sphere[p]` | represents a unit sphere centered at the point p. |
| `Sphere[p,r]` | represents a sphere of radius r centered at the point p. |
| `Sphere[{p1,p2,…},r]` | represents a collection of spheres of radius r. |
| `Ball[p]` | represents the unit ball centered at the point p. |
| `Ball[p,r]` | represents the ball of radius r centered at the point p. |
| `Ball[{p1,p2,…},r]` | represents a collection of balls of radius r. |
| `SphericalShell[c,{rinner,router}]` | represents a filled spherical shell centered at c with inner radius rinner and outer radius router. |
| `Ellipsoid[p,{r1,…}]` | represents an axis-aligned ellipsoid centered at the point p and with semiaxes lengths ri. |
| `Ellipsoid[p,Σ]` | represents an ellipsoid centered at p and weight matrix Σ. |
| `Simplex[{p1,…,pk}]` | represents the simplex spanned by points pi. |

***Mathematica Examples 3.1***

Input

```
(*A unit ball in 3D:*)
Graphics3D[
 Ball[]
 ]
```

Output

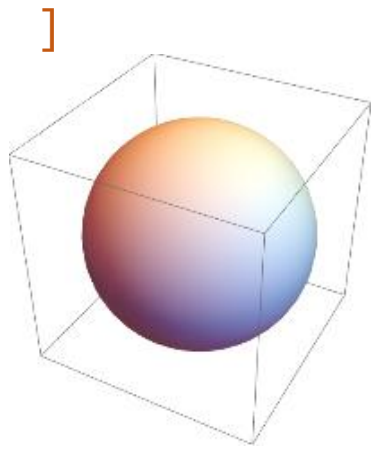

Input

```
(*The standard spherical shell at the origin:*)
SphericalShell[];
Graphics3D[
 {Opacity[0.5],%}
 ]
```

Output

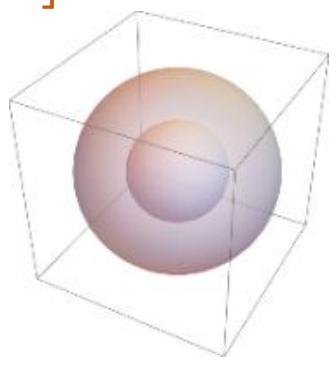

```
Input     (*A Simplex in 3D:*)
          Graphics3D[
           Simplex[{{0,0,1},{1,0,0},{1,0,1},{1,1,1}}]
           ]
Output
```

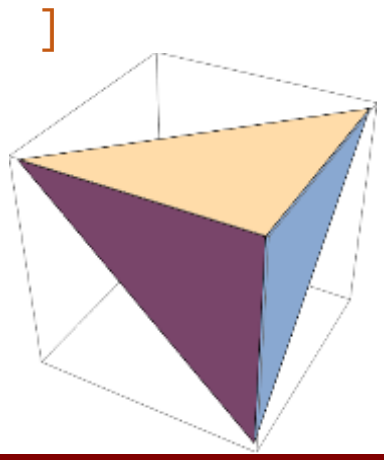

| | |
|---|---|
| `BoundedRegionQ[reg]` | gives True if reg is a bounded region and False otherwise. |
| `ConvexRegionQ[reg]` | gives True if reg is a convex region and False otherwise. |
| `ConvexHullRegion[{p1,p2,…}]` | gives the convex hull from the points p1, p2, …. |
| `ConvexHullRegion[reg]` | gives the convex hull of the region reg. |
| `ConvexHullMesh[{p1,p2,…}]` | gives a BoundaryMeshRegion representing the convex hull from the points p1, p2, …. |
| `ConvexHullMesh[mreg]` | gives the convex hull of the mesh region mreg. |

***Mathematica Examples 3.1***

```
Input     (*A bounded region:*)
          BoundedRegionQ[Disk[]]
          Region[
           Disk[],
           LabelStyle->Directive[Black,20]
           ]
Output    True
Output
```

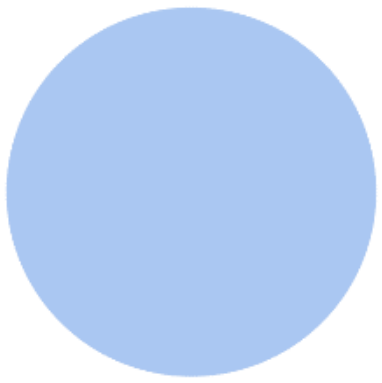

```
Input     (*An unbounded region:*)
          BoundedRegionQ[HalfPlane[{{0,0,0},{-1,1,1}},{1,1,1}]]
          Region[
           HalfPlane[{{0,0,0},{-1,1,1}},{1,1,1}],
           Boxed->True,
           LabelStyle->Directive[Black,20]
           ]
Output    False
Output
```

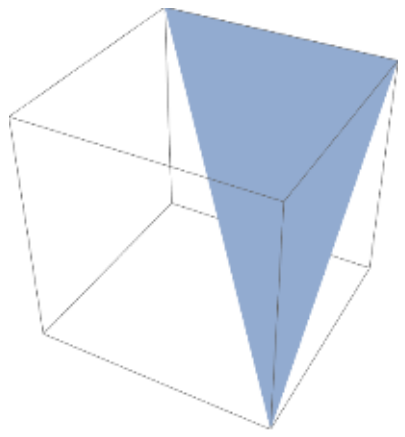

```
Input     (*Test whether a rectangle is convex:*)
          ConvexRegionQ[Rectangle[]]
Output    True
```

```
Input       (*A circle is not a convex region:*)
            Region[
              Circle[]
              ];
            ConvexRegionQ[%]
Output      False

Input       (*A 2D convex hull mesh from points:*)
            pts=RandomReal[{-1,1},{50,2}];
            ℛ=ConvexHullRegion[pts]
            (*The region is the smallest convex region that includes the points:*)
            Region[
             ℛ,
             Epilog->{Black,Point[pts]},
             LabelStyle->Directive[Black,20]
             ]
Output      Polygon[          ]
```

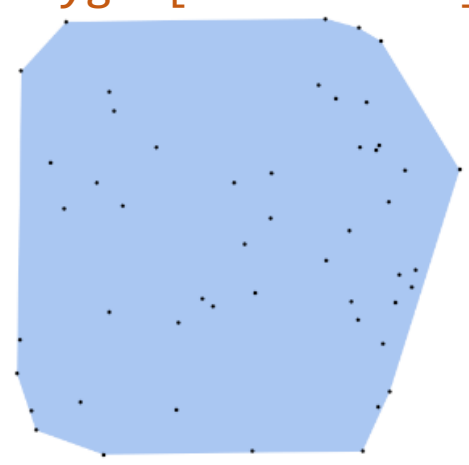

```
Input       (*A 2D convex hull mesh:*)
            pts=RandomReal[{-1,1},{50,2}];
            ConvexHullMesh[pts]
            (*The region is the smallest convex region that includes the points:*)
            Show[
             %,
             Graphics[Point[pts]]
             ]
Output
```

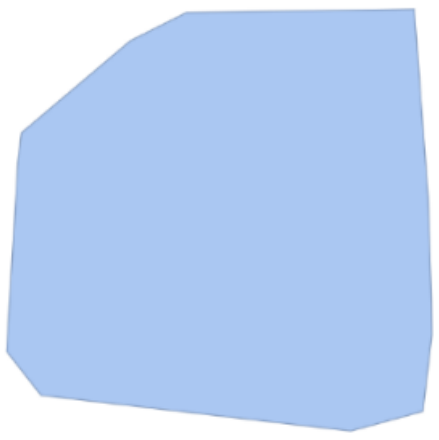

```
Output
```

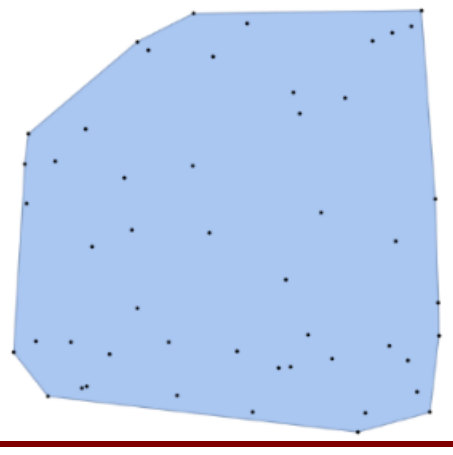

| | |
|---|---|
| `RegionPlot[pred,{x,xmin,xmax},{y,ymin,ymax}]` | makes a plot showing the region in which pred is True. |
| `RegionPlot[{pred1,pred2,…},…]` | plots several regions corresponding to the predi. |
| `RegionPlot[{…,w[predi,…],…},…]` | plots predi with features defined by the symbolic wrapper w. |
| `RegionPlot3D[pred,{x,xmin,xmax},{y,ymin,ymax},{z,zmin,zmax}]` | makes a plot showing the three-dimensional region in which pred is True. |
| `RegionPlot3D[{pred1,pred2,…},…]` | plots several regions corresponding to the predi. |

***Mathematica Examples 3.1***

Input
```
(*Plot a region defined by logical combinations of inequalities:*)
RegionPlot[
 x^2+y^3<2&&x+y<1,
 {x,-2,2},
 {y,-2,2},
 LabelStyle->Directive[Black,20]
 ]
```

Output

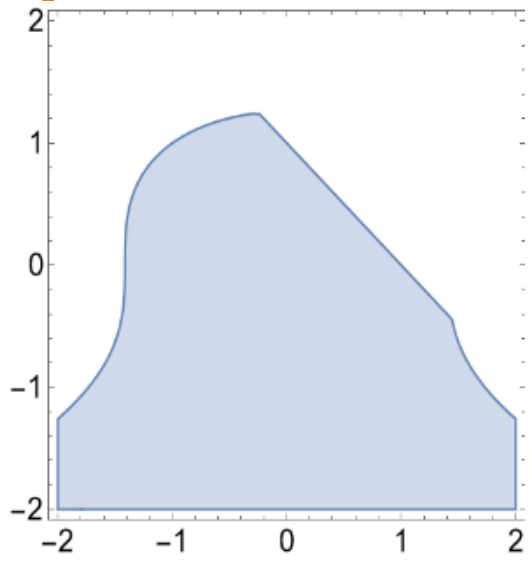

Input
```
(*Simple regions including a disk:*)
RegionPlot[
 x^2+y^2<=1,
 {x,-1,1},
 {y,-1,1},
 LabelStyle->Directive[Black,20]
 ]
(*Disk annulus:*)
RegionPlot[
 1/4<=x^2+y^2<=1,
 {x,-1,1},
 {y,-1,1},
 LabelStyle->Directive[Black,20]
 ]
```

Output

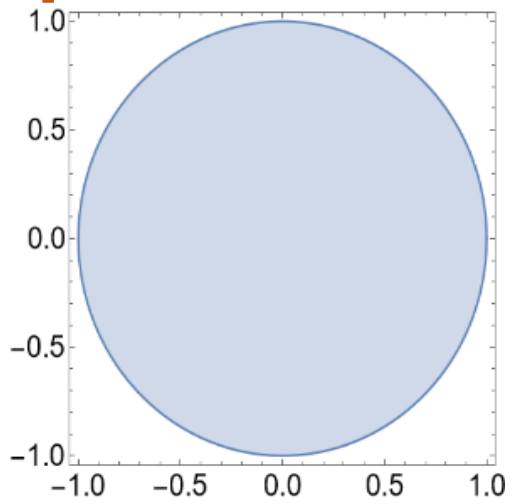

Output

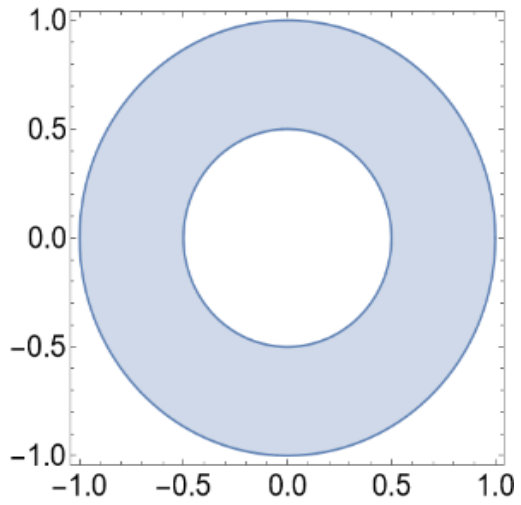

Input
```
(*Plot a 3D region:*)
RegionPlot3D[
 x^2+y^3-z^2>0,
 {x,-2,2},
 {y,-2,2},
 {z,-2,2},
 LabelStyle->Directive[Black,20]
 ]
```

Output

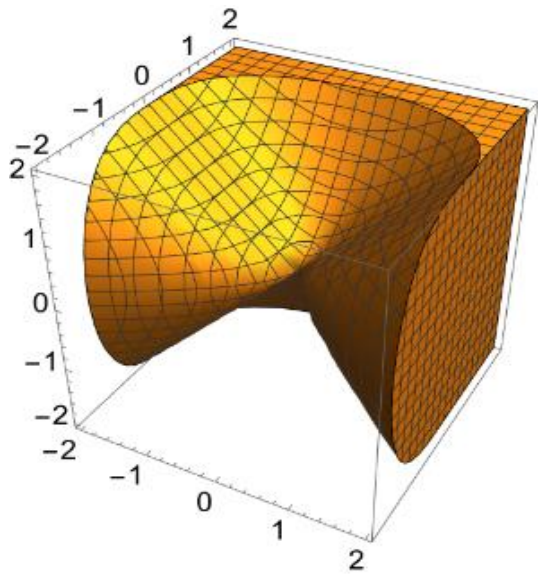

Input

```
(*Plot multiple regions:*)
RegionPlot3D[
 {x+y+z<-2,x+y+z>2},
 {x,-2,2},
 {y,-2,2},
 {z,-2,2},
 LabelStyle->Directive[Black,20]
 ]
```

Output

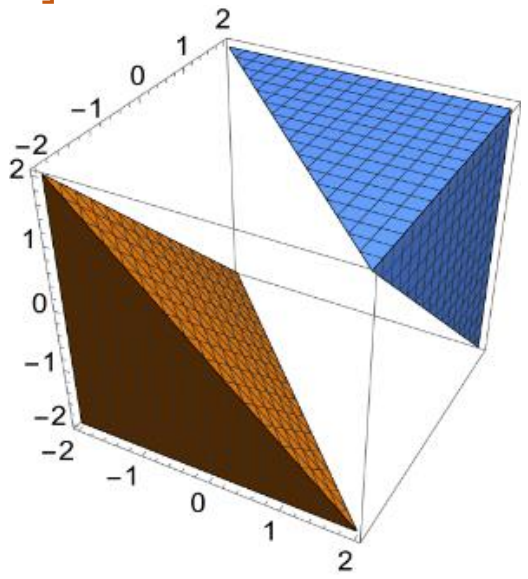

Input

```
(*Plot 3D regions defined by logical combinations of inequalities:*)
RegionPlot3D[
 x^2+y^2+z^2<1&&x^2+y^2<z^2,
 {x,-1,1},
 {y,-1,1},
 {z,-1,1},
 PlotPoints->35,
 PlotRange->All,
 LabelStyle->Directive[Black,20]
 ]
```

Output

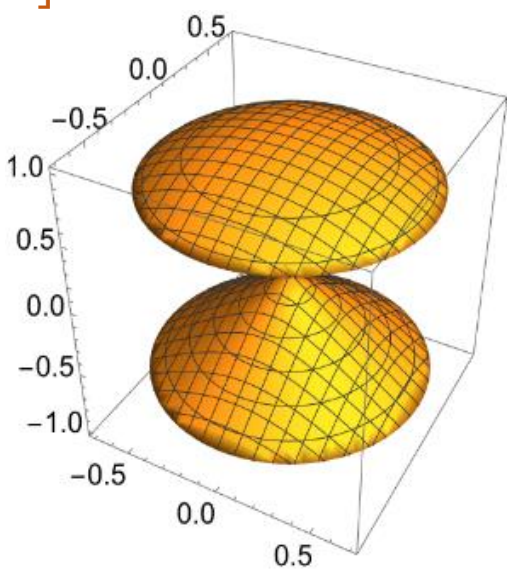

Input

```
(*Simple regions including a Ball:*)
RegionPlot3D[
 x^2+y^2+z^2<=3,
 {x,-2,2},
 {y,-2,2},
 {z,-2,2},
 Mesh->None,
 LabelStyle->Directive[Black,20]
 ]
```

Output

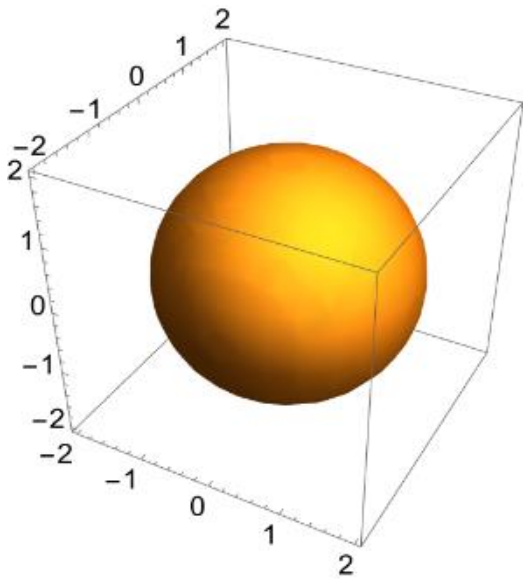

Input
```
(*Half of a spherical shell:*)
RegionPlot3D[
 1<=x^2+y^2+z^2<=3&&y>=0,
 {x,-2,2},
 {y,-2,2},
 {z,-2,2},
 Mesh->None,
 PlotPoints->50,
 LabelStyle->Directive[Black,20]
 ]
```
Output

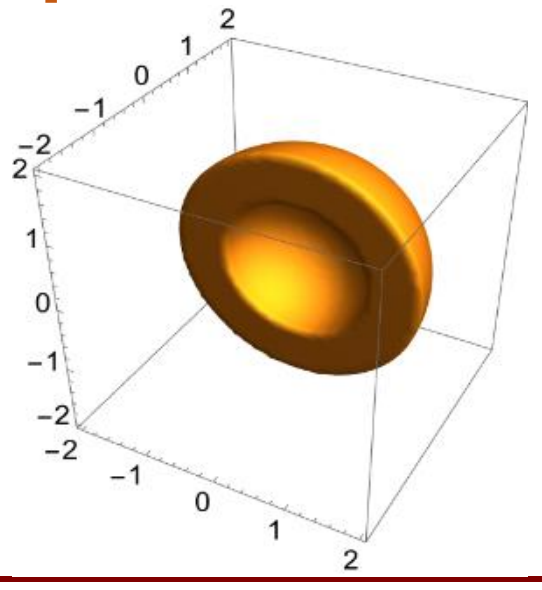

| | |
|---|---|
| `FunctionConvexity[f,{x1,x2,…}]` | finds the convexity of the function f with variables x1,x2,… over the reals. |
| `FunctionConvexity[{f,cons},{x1,x2,…}]` | finds the convexity when variables are restricted by the constraints cons representing a convex region. |

Convexity is also known as convex, concave, strictly convex and strictly concave [9]. By default, the following definitions are used: (+1: convex), (0: affine), (-1: concave), (Indeterminate : neither convex nor concave).

***Mathematica Examples 3.1***

Input
```
(*Find the convexity of a univariate function:*)
FunctionConvexity[Abs[x],x]
1
```

Input
```
(*Find the convexity of a multivariate function:*)
FunctionConvexity[-Exp[x+y],{x,y}]
```
Output
```
-1
```

Input
```
(*Find the convexity of a function with variables restricted by constraints:*)
FunctionConvexity[{x^3+y^3,x>=0&&y>=0},{x,y}]
```
Output
```
1
```

Input
```
(*Univariate functions:*)
FunctionConvexity[Sin[x],x]
FunctionConvexity[2^x,x]
FunctionConvexity[x-Sqrt[x^2],x]
Plot[
 {Sin[x],2^x,x-Sqrt[x^2]},
 {x,-3,3},
```

```
               LabelStyle->Directive[Black,20]
               ]
Output        Indeterminate
Output        1
Output        -1
Output
```

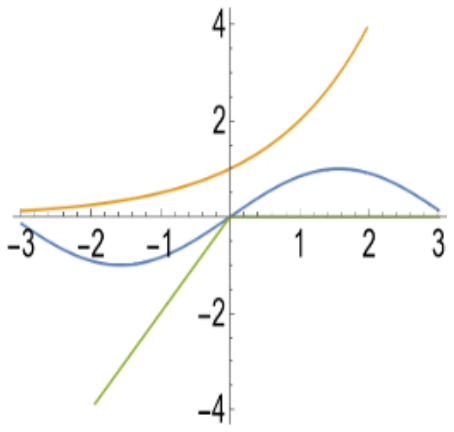

```
Input         (*Multivariate functions:*)
              FunctionConvexity[x^2+y^2,{x,y}]
              Plot3D[
               x^2+y^2,
               {x,-2,2},
               {y,-2,2},
               LabelStyle->Directive[Black,20]
               ]
Output        1
Output
```

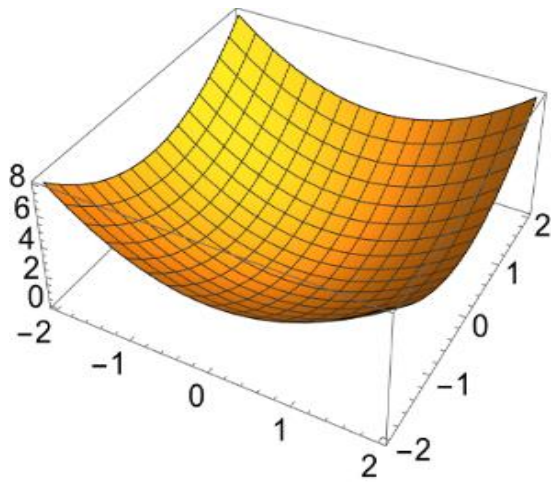

```
Input         (*FunctionConvexity gives a conditional answer here:*)
              FunctionConvexity[x^4+a x^2,x]
Output        1  if a∈ Reals && a>=0

Input         (*Check the convexity of:*)
              FunctionConvexity[x^2-1,x]
              (*The segment connecting any two points on the graph lies above the graph:*)
              Plot[
               x^2-1,
               {x,-5,5},
               Epilog->{Green,Line[{{-2,3},{3,8}}]},
               LabelStyle->Directive[Black,20]
               ]
Output        1
Output
```

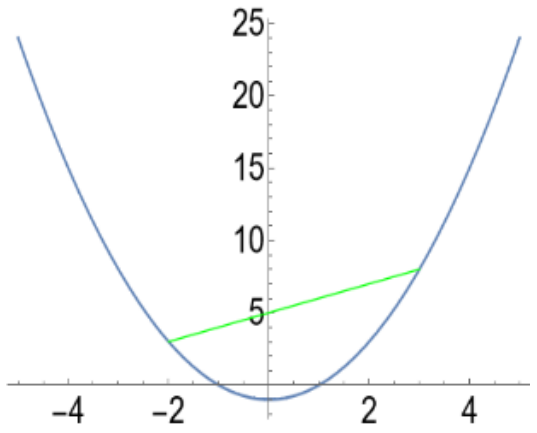

```
Input         (*Norm[v,p]} is convex for 1<=p<=infinity,but not strictly convex:*)
              FunctionConvexity[Norm[{x,y},p],{x,y},Assumptions->p>=1]
              FunctionConvexity[Norm[{x,y,z},Infinity],{x,y,z}]
              FunctionConvexity[Norm[{x,y},p],{x,y},Assumptions->p>=1,StrictInequalities-
              >True]
              Plot3D[
```

```
             {Norm[{x,y},1],Norm[{x,y},2],Norm[{x,y},Infinity]},
             Sequence[{x,-1,1},{y,-1,1},
              Mesh->None,
              PlotStyle->Opacity[0.7],
              PlotPoints->50,
              LabelStyle->Directive[Black,20]]
             ]
Output      1
Output      1
Output      Indeterminate
Output
```

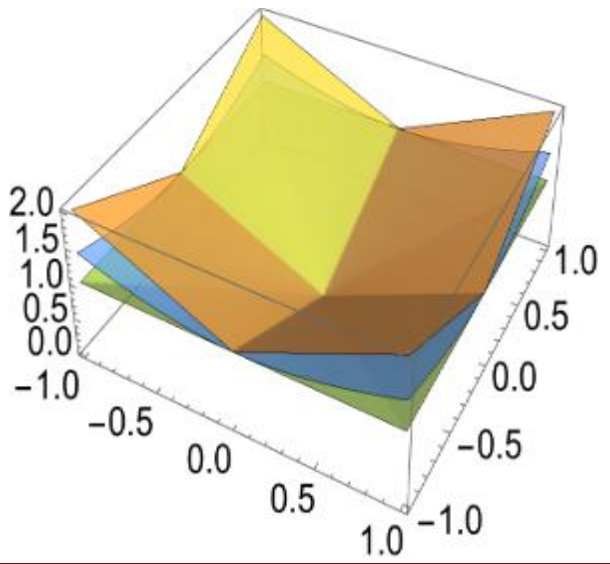

# CHAPTER 5

# SINGLE-VARIABLE OPTIMIZATION WITHOUT CONSTRAINTS

## 5.1 Single-variable Functions

Three general classes of nonlinear optimization problems can be identified as follows:

- One-dimensional unconstrained problems.
- Multidimensional unconstrained problems.
- Multidimensional constrained problems.

The optimization problem in which the performance measure is a function of one variable is the most elementary type of optimization problem. Yet, it is of central importance to optimization theory and practice because single-variable optimization often arises as a subproblem within the iterative procedures for solving multivariable optimization problems. Before we present conditions for a point to be an optimal point, we define three different types of optimal points.

**Definitions (Global and Local Minima):**

- A function $f$ has a global maximum at $c$ if $f(c) \geq f(x)$ for all $x$ in $D$, where $D$ is the domain of $f$. Similarly, $f$ has a global minimum at $c$ if $f(c) \leq f(x)$ for all $x$ in $D$.
- A function $f$ has a local maximum at $c$ if $f(c) \geq f(x)$ for all $x$ in some open interval containing $c$. Similarly, $f$ has a local minimum at $c$ if $f(c) \leq f(x)$ when $x$ is near $c$.
- An inflection point or saddle point is a point that does not correspond to a local optimum (minimum or maximum).

The maximum and minimum values of $f$ are called the extreme values of $f$.

Figure 5.1 shows the graph of a function $f$ with global maximum at $d$ and global minimum at $a$. Note that, the point $(d, f(d))$ is the highest point on the graph, and $(a, f(a))$ is the lowest point. If we consider only values of $x$ near $b$ [for instance, if we restrict our attention to the interval $(a, c)$], then $f(b)$ is the largest of those values of $f(x)$ and is the local maximum value of $f$. In Figure 5.1 multiple local optima are possible, and the global minimum/ maximum can be found only by locating all local optima and selecting the best one.

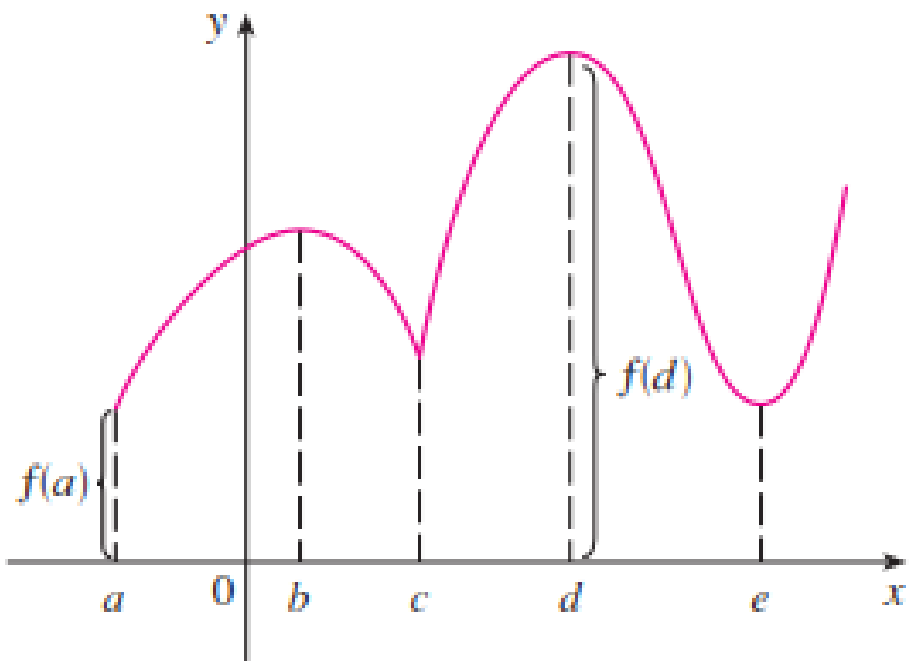

**Figure 5.1**. Minimum value $f(a)$, maximum value $f(d)$.

**Definition (Increasing and Decreasing Functions):** A function $f$ is called increasing on an interval $I$ if

$$f(x_1) < f(x_2) \text{ whenever } x_1 < x_2 \text{ in } I.$$

It is called decreasing on $I$ if

$$f(x_1) > f(x_2) \text{ whenever } x_1 < x_2 \text{ in } I.$$

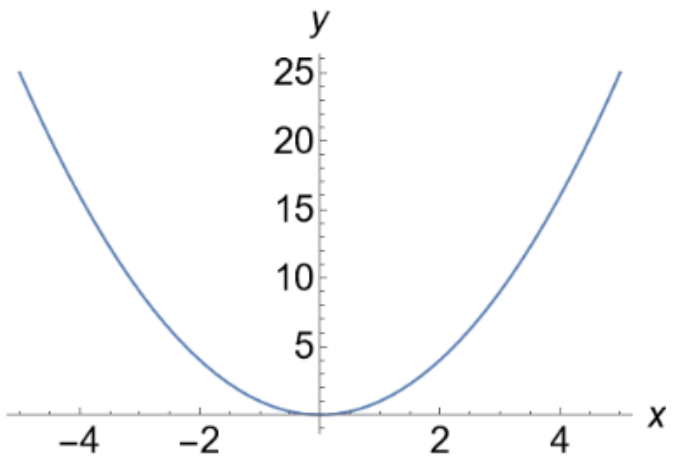

**Figure 5.2.** The function is monotonically decreasing on the interval $[-\infty, 0]$ and monotonically increasing on the interval $[0, \infty]$. It is a unimodal function.

**Definition (Monotonic Functions):** A function $f(x)$ is monotonic (either increasing or decreasing) if for any two points $x_1$ and $x_2$, with $x_1 \leq x_2$, it follows that

1- $f(x_1) \leq f(x_2)$ (Monotonically increasing)
2- $f(x_1) \geq f(x_2)$ (Monotonically decreasing)

**Definition (Unimodal Functions):** A function $f(x)$ is unimodal on the interval $a \leq x \leq b$ if and only if it is monotonic on either side of the single optimal point $x^*$ in the interval. See Figure 5.2.

**Theorem 5.1 (Increasing/Decreasing Test):**
(a) If $f'(x) > 0$ on an interval, then $f$ is increasing on that interval.
(b) If $f'(x) < 0$ on an interval, then $f$ is decreasing on that interval.

**Proof:**

Let $x_1$ and $x_2$ be any two numbers in the interval with $x_1 < x_2$. According to the definition of an increasing function, we have to show that $f(x_1) < f(x_2)$. Because we are given that $f'(x) > 0$, we know that $f$ is differentiable on $[x_1, x_2]$. So, by the mean value theorem, there is a number $c$ between $x_1$ and $x_2$ such that

$$f'(c) = \frac{f(x_2) - f(x_1)}{x_2 - x_1}.$$

Now, $f'(c) > 0$ by assumption and because $x_2 - x_1 > 0$. Thus, $f(x_1) < f(x_2)$. This shows that $f$ is increasing.

Part (b) is proved similarly.

■

**Theorem 5.2 (Necessary Condition):** If a function $f(x)$ is defined in the interval $a \leq x \leq b$ and has a local minimum at $x = x^*$, where $a < x^* < b$, and the derivative $\frac{d}{dx} f(x) = f'(x)$ exists as a finite number at $x = x^*$, then

$$f'(x^*) = 0. \quad (5.1)$$

**Proof:**

The definition of derivative provides

$$f'(x^*) = \lim_{h \to 0} \frac{f(x^* + h) - f(x^*)}{h}.$$

Since $x^*$ is a relative minimum, we have $f(x^*) \leq f(x^* + h)$ for all values of $h > 0$ (sufficiently close to zero). Hence,

$$\frac{f(x^*+h)-f(x^*)}{h} \geq 0, \qquad \text{if } h>0,$$

$$\frac{f(x^*+h)-f(x^*)}{h} \leq 0, \qquad \text{if } h<0.$$

The limit as $h$ tends to zero through positive values as follows:

$$f'(x^*) \geq 0.$$

However, it gives the limit as $h$ tends to zero through negative values as follows:

$$f'(x^*) \leq 0.$$

The only way to satisfy both conditions is to have

$$f'(x^*) = 0.$$

■

**Definition (Critical Point):** A critical point of a function $f$ is a point $c$ in the domain of $f$ such that $f'(c)=0$ or $f'(c)$ does not exist. See Figure 5.3.

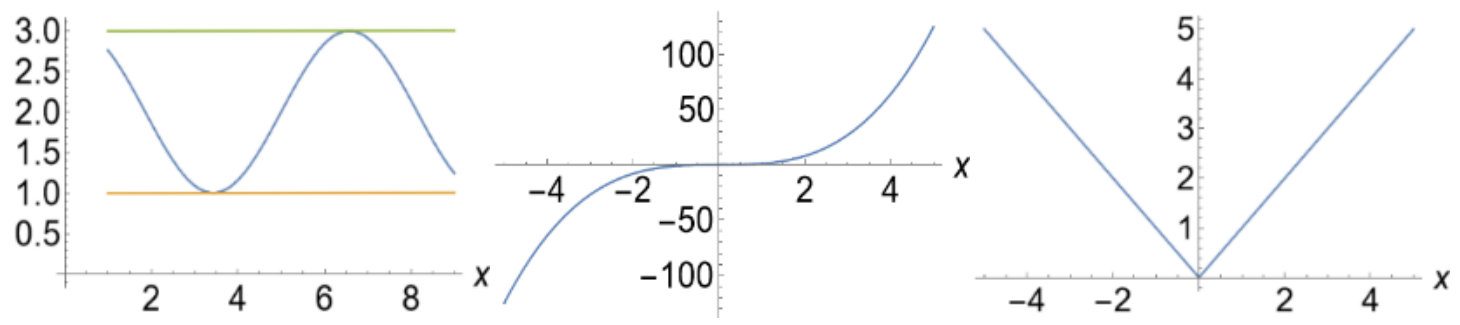

**Figure 5.3.** Critical points.

**The First Derivative Test**
Suppose that $c$ is a critical point of a continuous function $f$.
a) If $f'$ changes from positive to negative at $c$, then $f$ has a local maximum at $c$.
(b) If $f'$ changes from negative to positive at $c$, then $f$ has a local minimum at $c$.
(c) If $f'$ does not change sign at $c$ (that is, $f$ is positive on both sides of $c$ or negative on both sides), then $f$ has no local maximum or minimum at $c$.
See Figure 5.4.

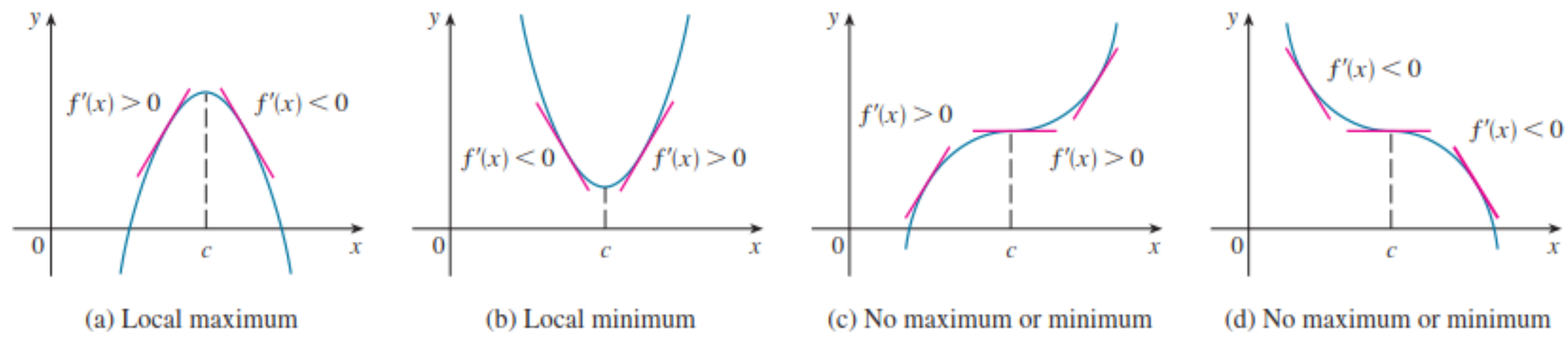

**Figure 5.4.** Local maximum and minimum.

Remember, a function (or its graph) is convex (concave upward) on an interval $I$ if $f'$ is an increasing function on $I$. It is concave (concave downward) on $I$ if $f'$ is decreasing on $I$. In Figure 5.5, the slopes of the tangent lines increase from left to right on the interval $(a,b)$, so $f'$ is increasing and $f$ is convex on $(a,b)$. [It can be proved that this is equivalent to saying that the graph of $f$ lies above all of its tangent lines on $(a,b)$.] Similarly, the slopes of the tangent lines decrease from left to right on $(b,c)$, so $f'$ is decreasing and $f$ is concave on $(b,c)$. [Review Theorems 4.25 and 4.26]

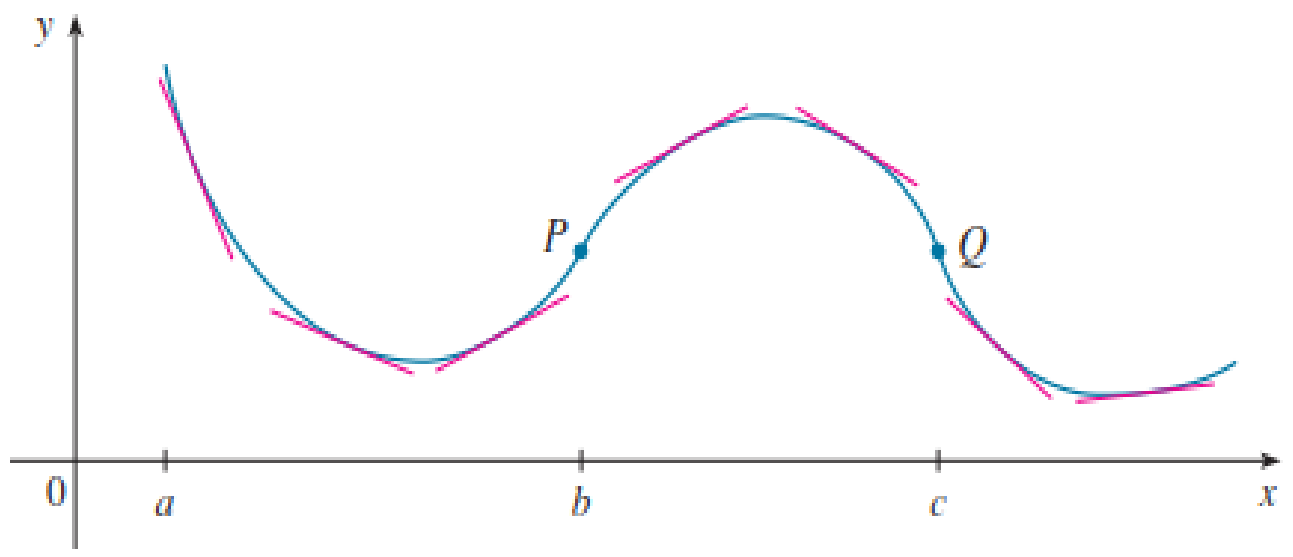

**Figure 5.5.** Concave upward and downward of a function.

Because $f'' = (f')'$, we know that if $f''(x)$ is positive, then $f'$ is an increasing function, and so $f$ is convex. Similarly, if $f''$ is negative, then $f'$ is decreasing, and $f$ is concave. Thus, we have the following test for concavity. [Review Theorem 4.27]

**Concavity Test**
(a) If $f'' > 0$ for all $x$ in $I$, then the graph of $f$ is convex on $I$.
(b) If $f'' < 0$ for all $x$ in $I$, then the graph of $f$ is concave on $I$.

**Theorem 5.3.1 (Sufficient Condition), in the case $f$ is twice differentiable:**
Suppose $f$ is continuous near $c$.
(a) If $f'(c) = 0$ and $f''(c) > 0$, then $f$ has a local minimum at $c$.
(b) If $f'(c) = 0$ and $f''(c) < 0$, then $f$ has a local maximum at $c$.

**Theorem 5.3.2 (Sufficient Condition):**
Suppose at a point $x^*$ the first derivative is zero and the first nonzero higher-order derivative is denoted by $n$.

1- If $n$ is odd, then $x^*$ is a point of inflection.
2- If $n$ is even, then $x^*$ is a local optimum. Moreover:
   a- If that derivative is positive, then the point $x^*$ is a local minimum.
   b- If that derivative is negative, then the point $x^*$ is a local maximum.

**Proof:**

This result is easily verified by using the Taylor series expansion of the function $f(x)$. we have

$$f(x^* + h) = f(x^*) + hf'(x^*) + \frac{h^2}{2} f''(x^*) + \cdots + \frac{h^n}{n!} f^{(n)}(x^* + \theta h), \qquad \text{for } 0 < \theta < 1.$$

Since the first nonvanishing higher-order derivative is $n$, the Taylor expansion becomes

$$f(x^* + h) - f(x^*) = \frac{h^n}{n!} f^{(n)}(x^* + \theta h).$$

If $n$ is odd, then the value of $f(x^* + h) - f(x^*)$ can be made positive or negative by choosing $h$ positive or negative. This implies that depending on the sign of $h$, $f(x^* + h) - f(x^*)$ could be positive or negative. Hence the function does not attain a minimum or a maximum at $x^*$, and $x^*$ is an inflection point.

Now, consider the case when $n$ is even. Then the term $h^n$ is always positive, and for all $h$ sufficiently small, the sign of $f(x^* + h) - f(x^*)$ will be dominated by the $\frac{d^{(n)}}{dx^{(n)}} f$. Hence, if $\frac{d^{(n)}}{dx^{(n)}} f$ is positive, $f(x^* + h) - f(x^*) > 0$, and $x^*$ corresponds to a local minimum. A similar argument can be applied in the case of a local maximum.

∎

***Example 5.1***

Maximize: $z = xe^{-x^2}$.

***Solution***

Here $f'(x) = e^{-x^2} - 2x^2e^{-x^2} = e^{-x^2}(1-2x^2)$ which is defined for all $x$ and which vanishes only at $x = \pm\frac{1}{\sqrt{2}}$.

Since $x$ is unrestricted, the values of the objective function at the stationary points (critical points),

$$f(\pm 1/\sqrt{2}) = \pm(1/\sqrt{2})e^{-\frac{1}{2}} = \pm 0.429$$

must be compared to the limiting values of $f(x)$ as $x \to \pm\infty$, which are both 0 in this case. Recording these results,

| $x$ | $x \to -\infty$ | $-1/\sqrt{2}$ | $1/\sqrt{2}$ | $x \to +\infty$ |
|---|---|---|---|---|
| $f(x)$ | 0 | $-0.429$ | 0.429 | 0 |

We see that a global maximum exists at $x^* = 1/\sqrt{2}$ and is $f(x^*) = 0.429$.

## 5.2 Optimization Problem and Algorithms

In this section, we consider the optimization algorithms for single-variable and unconstrained functions. These algorithms are repeatedly used as a subtask of many multi-variable optimization methods. Therefore, a clear understanding of these algorithms will help readers learn complex algorithms discussed in subsequent chapters.

The algorithms described in this chapter can be used to solve minimization problems of the following type:

$$\text{Minimize } f(x), \tag{5.2}$$

where $f(x)$ is the objective function, and $x$ is a real variable. The purpose of an optimization algorithm is to find a solution $x$, for which the function $f(x)$ is minimum. Two distinct types of algorithms are presented in this chapter: Direct search methods use only objective function values to locate the minimum point, and gradient-based methods use the first and/or the second-order derivatives of the objective function to locate the minimum point.

Even though the optimization methods described here are for minimization problems, they can also be used to solve a maximization problem by adopting the following procedure. In this case, an equivalent dual problem $(-f(x))$ is formulated and minimized. Hence, the algorithms described in this chapter can be directly used to solve a maximization problem. The approach is to use an iterative search algorithm, also called a line-search method. There are many of these algorithms [1-10], for examples

- **Bracketing Methods**
  - Exhaustive Search Method
  - Bounding Phase Method
- **Region-Elimination Methods**
  - Interval Halving Method
  - Fibonacci Search Method
  - Golden Section Search Method
- **Point-Estimation Method**
  - Successive Quadratic Estimation Method
- **Gradient-based Methods**
  - Newton-Raphson Method
  - Bisection Method
  - Secant Method
  - Cubic Search Method

In an iterative algorithm, we start with an initial candidate solution $x^{(0)}$ and generate a sequence of iterates $x^{(1)}, x^{(2)}, \ldots$. For each iteration $k = 0,1,2,\ldots$, the next point $x^{(k+1)}$ depends on $x^{(k)}$ and the objective function $f$. The algorithm may use only the value of $f$ at specific points, or perhaps its first derivative $f'$, or its second derivative $f''$.

## 5.3 Bracketing Methods

### 5.3.1 Exhaustive Search Method

It is the most basic of all the search strategies [1-10]. The optimum of a function is bracketed by computing the values of the function at a number of evenly spaced locations, as shown in Figure 5.6. Typically, the search starts with a lower limit on the variable and compares three successive function values at a time, based on the premise that the function is unimodal. Depending on the outcome of the comparison, the search is either terminated or continued by replacing one of the three points with a new point. The search continues until the minimum is bracketed.

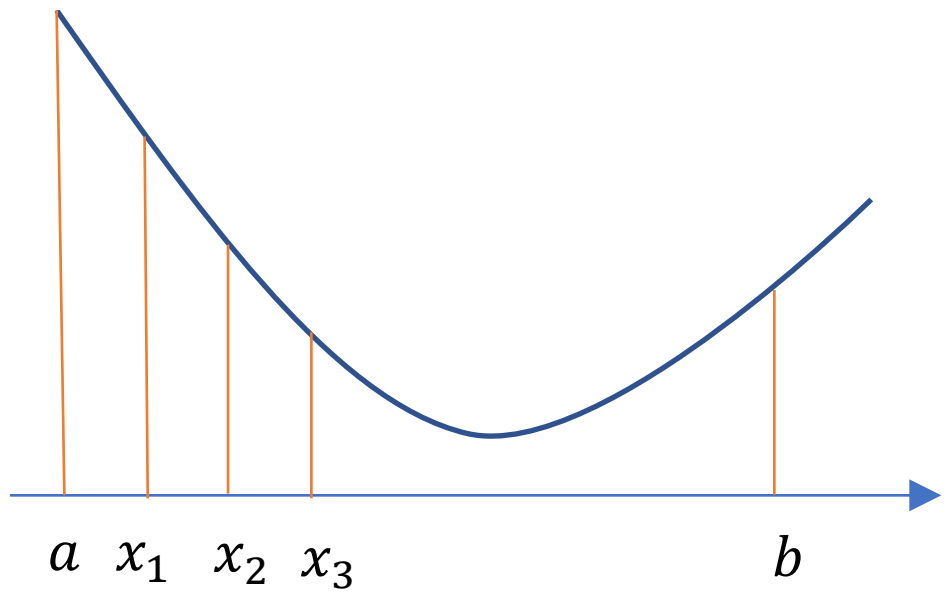

**Figure 5.6.** The exhaustive search method uses equally spaced points.

**Algorithm**

**Step 1:** Set $x_1 = a$, $\Delta x = (b - a)/n$ ($n$ is the number of intermediate points), $x_2 = x_1 + \Delta x$, and $x_3 = x_2 + \Delta x$.

**Step 2:** If $f(x_1) \geq f(x_2) \leq f(x_3)$, the minimum point lies in $(x_1, x_3)$, Terminate;
Else $x_1 = x_2$, $x_2 = x_3$, $x_3 = x_2 + \Delta x$, and go to Step 3.

**Step 3:** Is $x_3 \leq b$? If yes, go to Step 2;
Else no minimum exists in $(a, b)$ or a boundary point ($a$ or $b$) is the minimum point.

***Example 5.2***

Consider the problem:

$$\text{Minimize } f(x) = (x - 5)^2 + 6,$$

in the interval [1,9] and $n = 20$.

***Solution***

A plot of the function is shown in Figure 5.7. The plot shows that the minimum lies at $x^* = 5$, $f(x^*) = 6$.

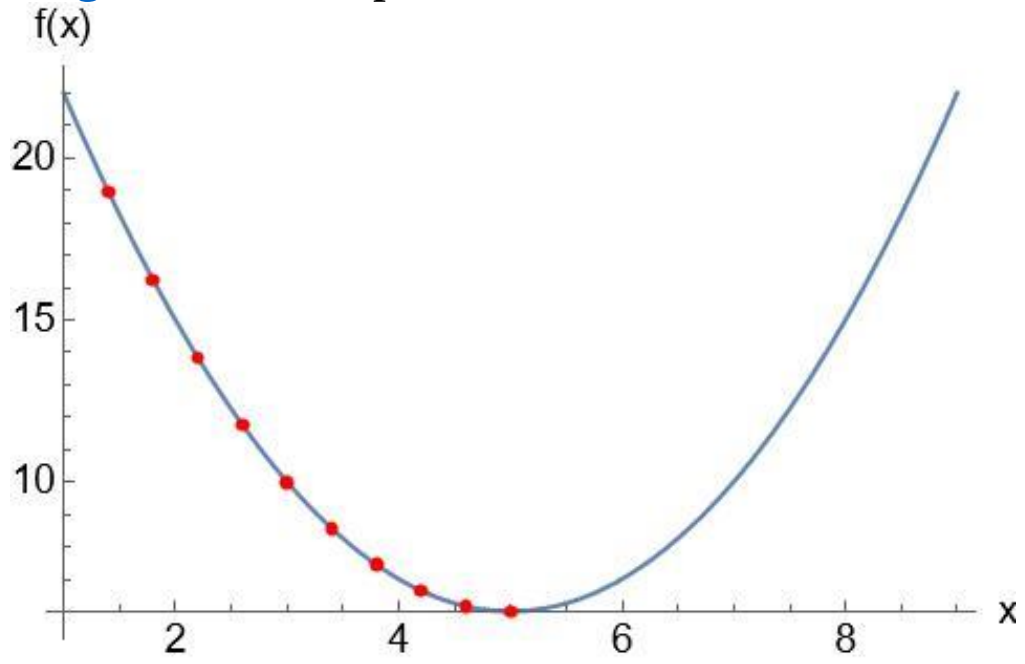

**Figure 5.7.** The results of 10 iterations of the exhaustive search method for $f(x) = (x - 5)^2 + 6$.

From Table 5.1, it is obtained that at iteration 10, the optimal condition is attained. Hence, the minimum lies in $(x_1, x_3)$, which is 4.6 and 5.4, and the objective function will possess a minimum at $x_2 = 5$.

**Table 5.1.** The results were produced by Mathematica code 5.1.

| No. of iters. | $x_1$ | $x_2$ | $x_3$ | $f(x_1)$ | $f(x_2)$ | $f(x_3)$ |
|---|---|---|---|---|---|---|
| 1 | 1 | 1.4 | 1.8 | 22 | 18.96 | 16.24 |
| 2 | 1.4 | 1.8 | 2.2 | 18.96 | 16.24 | 13.84 |
| 3 | 1.8 | 2.2 | 2.6 | 16.24 | 13.84 | 11.76 |
| 4 | 2.2 | 2.6 | 3 | 13.84 | 11.76 | 10 |
| 5 | 2.6 | 3 | 3.4 | 11.76 | 10 | 8.56 |
| 6 | 3 | 3.4 | 3.8 | 10 | 8.56 | 7.44 |
| 7 | 3.4 | 3.8 | 4.2 | 8.56 | 7.44 | 6.64 |
| 8 | 3.8 | 4.2 | 4.6 | 7.44 | 6.64 | 6.16 |
| 9 | 4.2 | 4.6 | 5 | 6.64 | 6.16 | 6 |
| 10 | 4.6 | 5 | 5.4 | 6.16 | 6 | 6.16 |

***Mathematica Code 5.1*** Exhaustive Search Method

```
(* Exhaustive Search Method *)

(*
Notations
a,b        :Lower and upper limits
n          :The number of intermediate points
f          :Objective function in variable x
delta      :Increment in x
lii        :The last iteration index
result[k] :The results of iteration k
*)

(* Taking Input from User *)
a = Input["Please input the lower limit of the domain where the function is defined"];
b = Input["Please input the upper limit of the domain where the function is defined"];

If[
    a > b,
    Beep[];
    MessageDialog["The value of lower limit a must be less than the upper limit b"];
    Exit[];
    ];

n = Input["Please input the the number of intermediate points"];

If[
    ! IntegerQ[n] || n < 0,
    Beep[];
    MessageDialog["The value of n must be postive integer"];
    Exit[];
    ];

(* Taking the Function from User *)
f[x_] = Evaluate[Input["Please input a unimodal function of x to find the minimum "]];
(* The user types in, for instance x^2 *)

(* Initiating Required Variables *)
delta = (b - a)/n;
x1 = a;
x2 = x1 + delta;
x3 = x2 + delta;
```

```
(* Starting the Algorithm *)
Do[
    If[
      x3 > b,
      Break[]
      ];

    y1[k_] = f[x1];
    y2[k_] = f[x2];
    y3[k_] = f[x3];

    lii = k;
    result[k] = N[{k, x1, x2, x3, y1[k], y2[k], y3[k]}];

    If[
      y1[k] > y2[k] && y2[k] < y3[k],
      Break[],
      x1 = x2;
      x2 = x3;
      x3 = x2 + delta;
      ];,
    {k, 1, ∞}
    ];

(* Final Result *)
Which[
  x3 > b && f[a] > f[b],
  Print["The boundary point b=", N[b], " is the mimimum point" , "\nThe solution is
(approximately) ", N[f[b]]];,

  x3 > b && f[b] > f[a],
  Print["The boundary point a=", N[a], " is the mimimum point" , "\nThe solution is
(approximately) ", N[f[a]]];,

  x3 <= b,
  Print["The solution lies between " , N[x1], " and ", N[x3] , "\nThe solution is
(approximately) ", N[y2[lii]], " at x= ", N[x2]];
  ]

(* Results of Each Iteration *)
table = TableForm[
    Table[
      result[i],
      {i, 1, lii}
      ],
    TableHeadings -> {None, {"No. of iters.", "x1", "x2", "x3", "f[x1]", "f[x2]", "f[x3]"}}
    ]

Export["example52.xls", table, "XLS"];

(* Data Visualization *)
Which[
  x3 > b && f[a] > f[b],
  Plot[
    f[x],
    {x, a, b},
    AxesLabel -> {"x", "f(x)"},
    LabelStyle -> Directive[Black, 14],
    Epilog -> {PointSize[0.015], Red, Point[Prepend[Table[{result[i][[3]], result[i][[6]]},
{i, 1, lii}], {b, f[b]}]]}
```

```
    ],

  x3 > b && f[b] > f[a],
  Plot[
    f[x],
    {x, a, b},
    AxesLabel -> {"x", "f(x)"},
    LabelStyle -> Directive[Black, 14],
    Epilog -> {PointSize[0.015], Red, Point[Prepend[Table[{result[i][[3]], result[i][[6]]},
{i, 1, lii}], {a, f[a]}]]}
    ],

  x3 <= b,
  Plot[
  f[x],
  {x, a, b},
  AxesLabel -> {"x", "f(x)"},
  LabelStyle -> Directive[Black, 14],Epilog ->
{PointSize[0.015],Red,Point[Table[{result[i][[3]], result[i][[6]]},{i, 1, lii}]]}
  ]
 ]

(* Data Manipulation *)
Manipulate[
  Plot[
   f[x],
  {x, a, b},
  AxesLabel -> {"x", "f(x)"},
  LabelStyle -> Directive[Black, 14],
  Epilog -> { PointSize[0.02],Red,Point[{result[i][[3]], result[i][[6]]}] }
  ],
  {i, 1, lii, 1}
  ]
```

### 5.3.2 Bounding Phase Method

**Definition (strictly unimodal):** A function is strictly unimodal if it is unimodal and has no intervals of finite length in which the function is of constant value.

**Theorem (Elimination Property):** Suppose $f$ is strictly unimodal on the interval $a \le x \le b$ with a minimum at $x^*$. Let $x_1$ and $x_2$ be two points in the interval such that $a < x_1 < x_2 < b$. Comparing the functional values at $x_1$ and $x_2$, we can conclude:
(i) If $f(x_1) > f(x_2)$, then the minimum of $f(x)$ does not lie in the interval $(a, x_1)$. In other words, $x^* \in (x_1, b)$ (see Figure 5.8).
(ii) If $f(x_1) < f(x_2)$, then the minimum does not lie in the interval $(x_2, b)$ or $x^* \in (a, x_2)$ (see Figure 5.8).

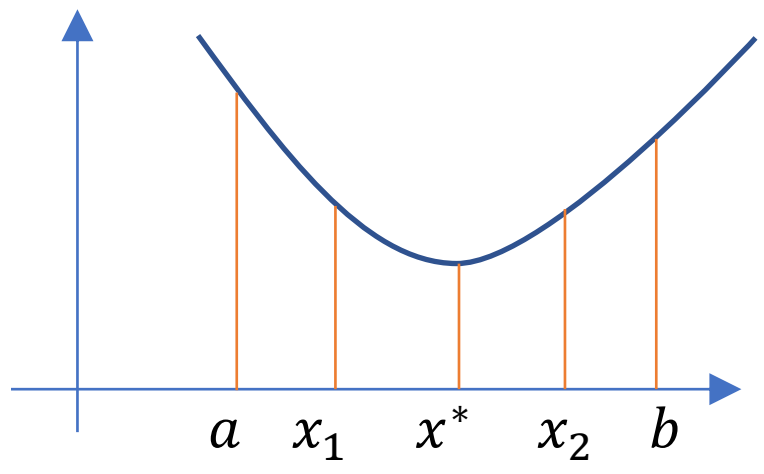

**Figure 5.8.** Case (i) and case (ii) of elimination property.

In the initial phase, starting at some selected trial point, the optimum is roughly bracketed within a finite interval by using the elimination property. Typically, this bounding search is conducted using some heuristic expanding pattern, although extrapolation methods have also been devised. An example of an expanding pattern is Swann's method, in which the $(k+1)$ test point is generated using the recursion

$$x_{k+1} = x_k + 2^k\Delta \text{ for } k = 0, 1, 2, \ldots, \tag{5.3}$$

where $x_0$ is an arbitrarily selected starting point and $\Delta$ is a step-size parameter of suitably chosen magnitude. The sign of $\Delta$ is determined by comparing $f(x_0)$, $f(x_0 + |\Delta|)$, and $f(x_0 - |\Delta|)$. If

$$f(x_0 - |\Delta|) \geq f(x_0) \geq f(x_0 + |\Delta|), \tag{5.4}$$

then, because of the unimodality assumption, the minimum must lie to the right of $x_0$, and $\Delta$ is chosen to be positive. If the inequalities are reversed, $\Delta$ is chosen to be negative; if

$$f(x_0 - |\Delta|) \geq f(x_0) \leq f(x_0 + |\Delta|), \tag{5.5}$$

the minimum has been bracketed between $x_0 - |\Delta|$ and $x_0 + |\Delta|$ and the bounding search can be terminated. The remaining case,

$$f(x_0 - |\Delta|) \leq f(x_0) \geq f(x_0 + |\Delta|), \tag{5.6}$$

is ruled out by the unimodality assumption. However, the occurrence of the above condition indicates that the given function is not unimodal [1-10].

**Algorithm**

**Step 1:** Choose an initial guess $x^{(0)}$ and an increment $\Delta$. Set $k = 0$.

**Step 2:** If $f(x^{(0)} - |\Delta|) \geq f(x^{(0)}) \geq f(x^{(0)} + |\Delta|)$, then $\Delta$ is positive;
Else if $f(x^{(0)} - |\Delta|) \leq f(x^{(0)}) \leq f(x^{(0)} + |\Delta|)$, then $\Delta$ is negative;
Else go to Step 1.

**Step 3:** Set $x^{(k+1)} = x^{(k)} + 2^k\Delta$.

**Step 4:** If $f(x^{(k+1)}) < f(x^{(k)})$, set $k = k + 1$ and go to Step 3;
Else the minimum lies in the interval $(x^{(k-1)}, x^{(k+1)})$ and Terminate.

If the chosen $\Delta$ is large, the bracketing accuracy of the minimum point is poor, but the bracketing of the minimum is faster. On the other hand, if the chosen $\Delta$ is small, the bracketing accuracy is better, but more function evaluations may be necessary to bracket the minimum. This method of bracketing the optimum is usually faster than the exhaustive search method. The algorithm approaches the optimum exponentially, but the accuracy in the obtained interval may not be very good, whereas in the exhaustive search method, the iterations required to attain near the optimum may be large, but the obtained accuracy is good. An algorithm with a mixed strategy may be more desirable.

***Example 5.3***

Consider the problem:

$$\text{Minimize } f(x) = (x-5)^2 + 6,$$

given the starting point $x_0 = 9$ and a step size $|\Delta| = 0.1$.

***Solution***

From Table 5.2, it is obtained that at iteration 5, the optimal condition is attained. Consequently, in five evaluations $x^*$ has been bracketed within the interval $2.7 \leq x \leq 7.5$ .

**Table 5.2.** The results were produced by Mathematica code 5.2.

| No. of iters. $k$ | $x[k]$ | $x[k+1]$ | $\lvert x[k] - x[k+1] \rvert$ | $f(x[k])$ | $f(x[k+1])$ |
|---|---|---|---|---|---|
| 0 | 9 | 8.9 | 0.1 | 22 | 21.21 |
| 1 | 8.9 | 8.7 | 0.2 | 21.21 | 19.69 |
| 2 | 8.7 | 8.3 | 0.4 | 19.69 | 16.89 |
| 3 | 8.3 | 7.5 | 0.8 | 16.89 | 12.25 |
| 4 | 7.5 | 5.9 | 1.6 | 12.25 | 6.81 |
| 5 | 5.9 | 2.7 | 3.2 | 6.81 | 11.29 |

A plot of the function is shown in Figure 5.9. The plot shows that the minimum lies at $x^* = 5$, $f(x^*) = 6$.

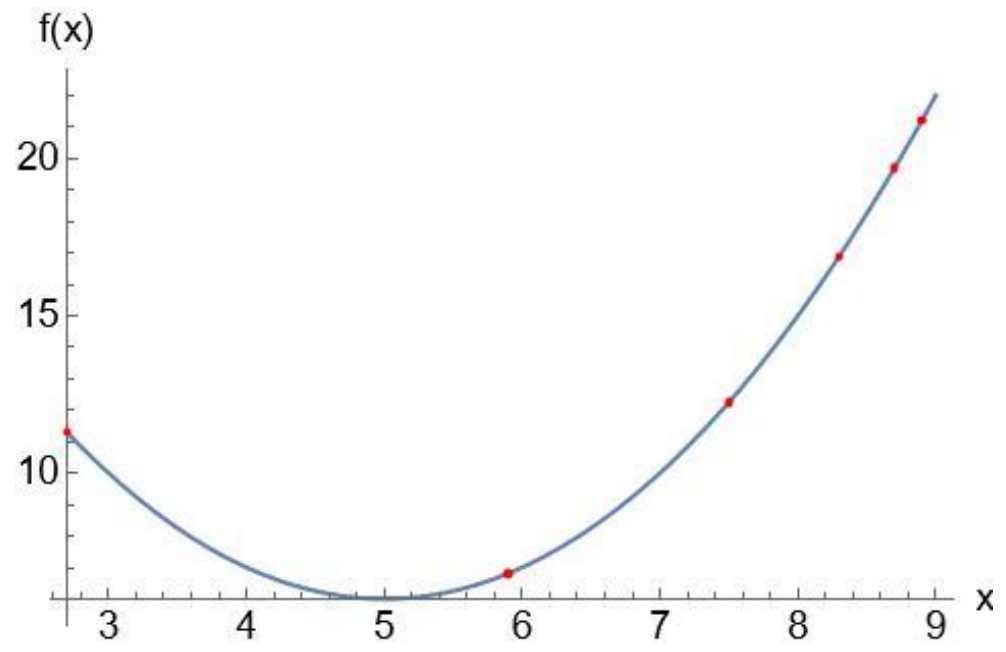

**Figure 5.9.** The results of 5 iterations of the bounding phase method for $f(x) = (x-5)^2 + 6$.

***Mathematica Code 5.2*** Bounding Phase Method

```
(* Bounding Phase Method *)

(*
Notations
x0        :Initial guess of the min
f         :Objective function in varible x
delta     :Increment in x
lii       :The last iteration index
result[k] :The results of iteration k
*)

(* Taking Input from User *)
x0=Input["Enter the initial guess:"];
delta=Input["Enter the increment:"];

(* Taking the Function from User *)
f[x_] = Evaluate[Input["Please input a function of x to find the minimum "]];
(* the user types in, for instance x^2 *)

(* Initiating Required Variables *)
y1 = f[x0-Abs[delta]];
y2 = f[x0];
y3 = f[x0+Abs[delta]];

(* Determining Whether the Inicrement Is Positive or Negative*)
Which[
  y1>y2&&y2>y3,
  increment=Abs[delta],
  y1<y2&&y2<y3,
  increment=-Abs[delta],
  y1>y2&&y2<y3,
  Print["The minimum has been bracketed between ",x0-Abs[delta]," and ",x0+Abs[delta]];
  Exit[];
  ];

(* Starting the Algorithm *)
Do[
  x[0]=x0;
  x[k+1]=x[k]+2^k*increment;

  b2=x[k+1];
```

```
  lii=k;
  result[k]=N[{k,x[k],x[k+1],Abs[x[k]-x[k+1]],f[x[k]],f[x[k+1]]}];

  Which[
   f[x[k]]<f[x[k+1]],

   (* Final Result *)
   Print["The solution is between ", x[k-1], " and  ", x[k+1],"\nThe solution is
(approximately) ", f[x[k]]];
   Break [],

   k>50,
   Print["After 50 iteration the function decrease"];
   Exit[]
   ];,
  {k,0,∞}
  ];

(* Results of Each Iteration *)
table=TableForm[
  Table[
   result[i],
   {i,0,lii}
   ],
  TableHeadings->{None,{"No. of iters.","x[k]","x[k+1]","Abs[x[k]-
x[k+1]]","f[x[k]]","f[x[k+1]]"}}
  ]
Export["example53.xls",table,"XLS"];

(* Data Visualization *)
Plot[
 f[x],
 {x, x0, b2},
 AxesLabel->{"x","f(x)"},
 LabelStyle->Directive[Black,14],
 Epilog->{PointSize[0.01],Red,Point[Table[{result[i][[3]],result[i][[6]]},{i,0,lii}]]}
 ]

(* Data Manipulation *)
Manipulate[
 Plot[
  f[x],
  {x,x0, b2},
  AxesLabel->{"x","f(x)"},
  LabelStyle->Directive[Black,14],
  Epilog->{PointSize[0.02],Red,Point[{result[i][[3]],result[i][[6]]}]}
  ],
 {i,0,lii,1}
 ]
```

## 5.4 Region-Elimination Methods

### 5.4.1 Interval Halving Method

This method [6-9] deletes exactly one-half the interval at each stage. This is also called a three-point equal-interval search since it works with three equally spaced trial points in the search interval. Figure 5.10 shows these three points in the interval. (Three points chosen in the interval $(a, b)$ are all equidistant from each other and equidistant from the boundaries by the same amount.) The three points divide the search space into four regions. The fundamental region elimination rule is used to eliminate a portion of the search space based on function values at the three chosen points.

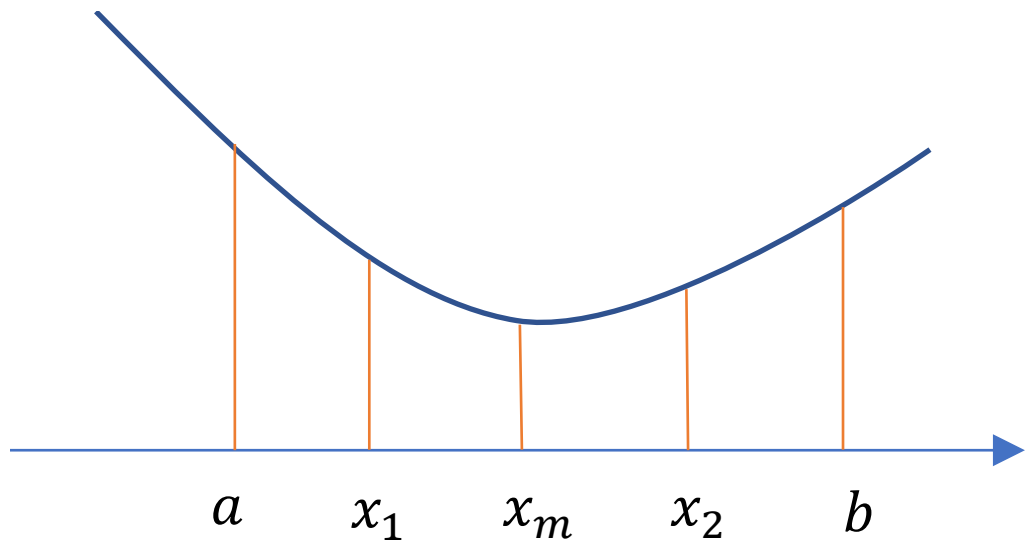

**Figure 5.10.** The three equidistant trial points of the interval halving method.

Two of the function values are compared at a time, and some region is eliminated. There are three scenarios that may occur.

- If $f(x_1) < f(x_m)$, then the minimum cannot lie beyond $x_m$. Therefore, we reduce the interval from $(a, b)$ to $(a, x_m)$. The point $x_m$ being the middle of the search space, this elimination reduces the search space to 50 percent of the original search space.
- On the other hand, if $f(x_1) > f(x_m)$, the minimum cannot lie in the interval $(a, x_1)$. The point $x_1$ being at one-fourth point in the search space, this reduction is only 25 percent. Thereafter, we compare function values at $x_m$ and $x_2$ to eliminate a further 25 percent of the search space. This process continues until a small enough interval is found.
- Finally, if $f(x_1) = f(x_m)$, we can conclude that regions $(a, x_1)$ and $(b, x_2)$ can be eliminated with the assumption that there exists only one local minimum in the search space $(a, b)$. We must have $x_1 < x^* < x_m$.

The basic steps of the search procedure for finding the minimum of a function $f(x)$ over the interval $(a, b)$ are as follows:

**Algorithm**

**Step 1:** Choose a lower bound $a$ and an upper bound $b$. Choose also a small number $\epsilon$. Let $x_m = (a + b)/2, L_0 = L = b - a$. Compute $f(x_m)$.

**Step 2:** Set $x_1 = a + L/4, x_2 = b - L/4$. Compute $f(x_1)$ and $f(x_2)$.

**Step 3:** If $f(x_1) < f(x_m)$ set $b = x_m;\ x_m = x_1$; go to Step 5;
Else go to Step 4.

**Step 4:** If $f(x_2) < f(x_m)$ set $a = x_m; x_m = x_2$; go to Step 5;
Else set $a = x_1, b = x_2$; go to Step 5.

**Step 5:** Calculate $L = b - a$. If $|L| < \epsilon$, Terminate;
Else go to Step 2.

***Example 5.4***

Consider the problem:

$$\text{Minimize } f(x) = (x - 5)^2 + 6,$$

in the interval [3,9] and $\epsilon = 0.01$.

***Solution***

A plot of the function is shown in Figure 5.11. The plot shows that the minimum lies at $x^* = 5, f(x^*) = 6$.

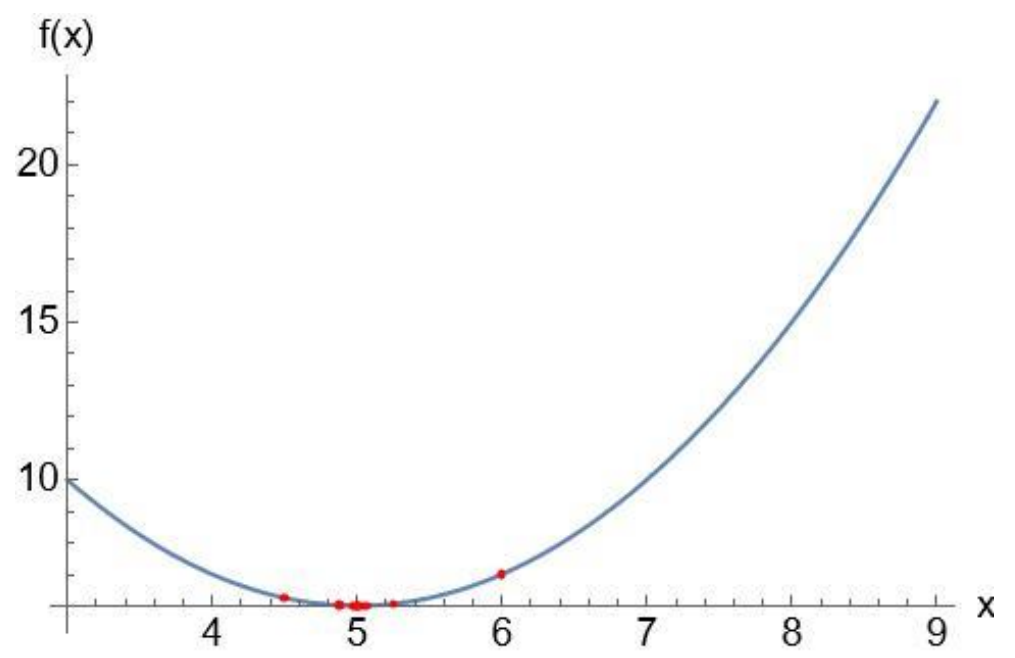

**Figure 5.11.** The results of 10 iterations of the interval halving method for $f(x) = (x-5)^2 + 6$.

From Table 5.3, it is obtained that at iteration 10, the optimal condition is attained. Consequently, in ten evaluations $x^*$ has been bracketed within the interval $4.995 \leq x \leq 5.0009$ .

**Table 5.3.** The results were produced by Mathematica code 5.3.

| No. it. | $a$ | $x_1$ | $x_m$ | $x_2$ | $b$ | $f(x_1)$ | $f(x_m)$ | $f(x_2)$ | $l$ |
|---|---|---|---|---|---|---|---|---|---|
| 1 | 3 | 4.5 | 6 | 7.5 | 9 | 6.25 | 7 | 12.25 | 6 |
| 2 | 3 | 3.75 | 4.5 | 5.25 | 6 | 7.5625 | 6.25 | 6.0625 | 3 |
| 3 | 4.5 | 4.875 | 5.25 | 5.625 | 6 | 6.015625 | 6.0625 | 6.390625 | 1.5 |
| 4 | 4.5 | 4.6875 | 4.875 | 5.0625 | 5.25 | 6.097656 | 6.015625 | 6.003906 | 0.75 |
| 5 | 4.875 | 4.96875 | 5.0625 | 5.15625 | 5.25 | 6.000977 | 6.003906 | 6.024414 | 0.375 |
| 6 | 4.875 | 4.921875 | 4.96875 | 5.015625 | 5.0625 | 6.006104 | 6.000977 | 6.000244 | 0.1875 |
| 7 | 4.96875 | 4.992188 | 5.015625 | 5.039063 | 5.0625 | 6.000061 | 6.000244 | 6.001526 | 0.09375 |
| 8 | 4.96875 | 4.980469 | 4.992188 | 5.003906 | 5.015625 | 6.000381 | 6.000061 | 6.000015 | 0.046875 |
| 9 | 4.992188 | 4.998047 | 5.003906 | 5.009766 | 5.015625 | 6.000004 | 6.000015 | 6.000095 | 0.023438 |
| 10 | 4.992188 | 4.995117 | 4.998047 | 5.000977 | 5.003906 | 6.000024 | 6.000004 | 6.000001 | 0.011719 |

Thus, in 10 stages, the initial search interval of length 6 has been reduced exactly to 0.0117.

***Mathematica Code 5.3*** Interval Halving Method

```
(* Interval Halving Method *)

(*
Notations
a,b        :Lower and upper limits
f          :Objective function in varible x
epsilon    :Tolerance
lii        :The last iteration index
result[k] :The results of iteration k
*)

(* Taking Input from User *)
a = Input["Please input the lower limit of the domain where the function is defined"];
b = Input["Please input the upper limit of the domain where the function is defined"];

If[
  a>b,
  Beep[];
  MessageDialog["The value of lower limit a must be less than b"];
  Exit[];
  ];

a0=a;
```

```
b0=b;

epsilon=Input["Please enter a tolerance (small postive number):"];
(* A very small number to check whether the loop would be executed or not*)

If[
  epsilon<=0,
  Beep[];
  MessageDialog["Tolerance value has to be small postive number: "];
  Exit[];
  ];

(* Taking the Function from User *)
f[x_] = Evaluate[Input["Please input a function of x to find the minimum "]]; (* the user
types in, for instance x^2 *)

(* Initiating Required Variables *)
l=b-a;

(* Starting the Algorithm *)
Do[
 xm=(a+b)/2;
 ym=f[xm];
 x1=a+(l/4);
 x2=b-(l/4);
 y1=f[x1];
 y2=f[x2];

 lii=k;
 result[k]=N[{k,a,x1,xm,x2,b,f[x1],f[xm],f[x2],l,Abs[x1-x2]}];

 If[
  f[x1]<f[xm],
  b=xm;
  xm=x1;,

  If[f[x2]<f[xm],
   a=xm;
   xm=x2;,
   a=x1;
   b=x2;
   ]
  ];

 l=b-a;
 If[
  Abs[l]<epsilon,
  Break[]
  ];,
 {k,1,∞}
 ]

(* Final Result *)
Print["The solution lies between " ,N[a], " and ", N[b] , "\nThe solution is (approximately)
", N[f[xm]]];

(* Results of Each Iteration *)
table=TableForm[
  Table[
   result[i],
   {i,1,lii}
```

```
    ],
   TableHeadings->{None,{"No. of
iters.","a","x1","xm","x2","b","f[x1]","f[xm]","f[x2]","l","Abs[x1-x2]"}}
   ]

Export["example54.xls",table,"XLS"];

(* Data Visualization *)
Plot[
 f[x],
 {x, a0, b0},
 AxesLabel->{"x","f(x)"},
 LabelStyle->Directive[Black,14],
 Epilog->{PointSize[0.01],Red,Point[Table[{result[i][[4]],result[i][[8]]},{i,1,lii}]]}
 ]

(* Data Manipulation *)
Manipulate[
 Plot[
  f[x],
  {x,a0,b0},
  AxesLabel->{"x","f(x)"},
  LabelStyle->Directive[Black,14],
  Epilog->{PointSize[0.02],Red,Point[{result[i][[4]],result[i][[8]]}]}
  ],
 {i,1,lii,1}
 ]
```

### 5.4.2 Fibonacci Search Method

Consider an interval of uncertainty

$$I^k = [x_L^k, x_U^k], \tag{5.7}$$

and assume that two points $x_a^k$ and $x_b^k$ are located in $I^k$, as depicted in Figure 5.12.

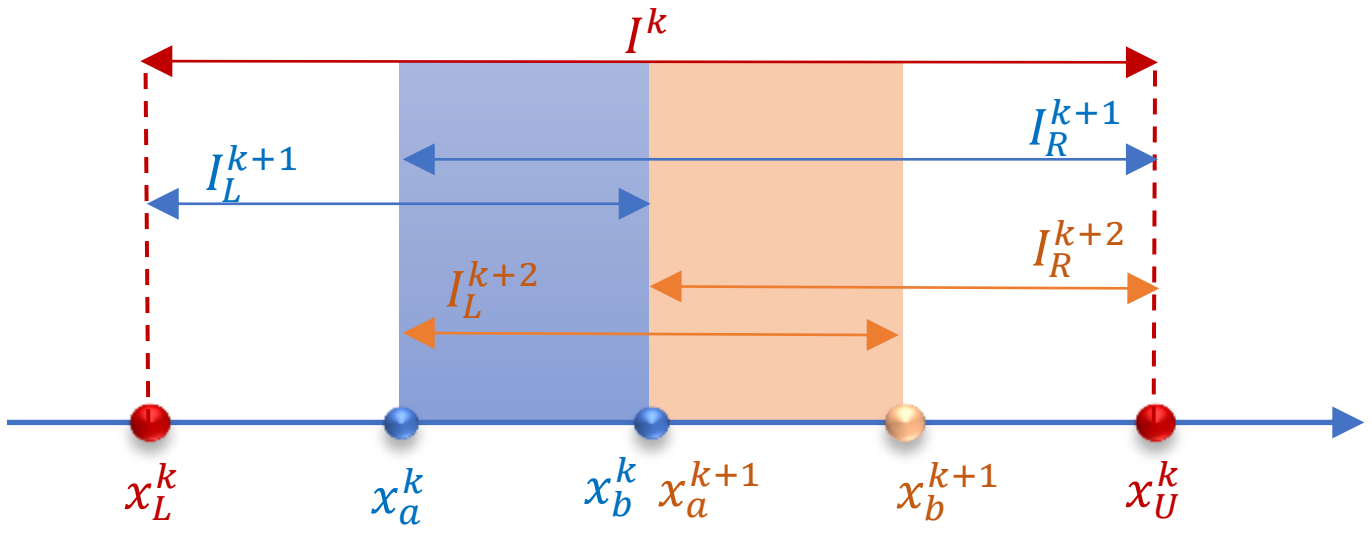

**Figure 5.12.** Reduction of range of uncertainty**.**

The values of $f(x)$ at $x_a^k$ and $x_b^k$, namely, $f(x_a^k)$ and $f(x_b^k)$, can be used to select the left interval [2-4]

$$I_L^{k+1} = [x_L^k, x_b^k], \tag{5.8}$$

if $f(x_a^k) < f(x_b^k)$, the right interval

$$I_R^{k+1} = [x_a^k, x_U^k], \tag{5.9}$$

if $f(x_a^k) > f(x_b^k)$, or either of $I_R^{k+1}$ and $I_L^{k+1}$ if $f(x_a^k) = f(x_b^k)$.

If the right interval $I_R^{k+1}$ is selected, it contains the minimizer and, in addition, the value of $f(x)$ is known at one interior point of $I_R^{k+1}$ namely, at the point $x_b^k$. If $f(x)$ is evaluated at one more interior point of $I_R^{k+1}$, say, at point$x_b^{k+1}$, sufficient information is available to allow a further reduction in the region of uncertainty, and the above cycle of events can be repeated. One of the two new sub-intervals $I_L^{k+2}$ and $I_R^{k+2}$, shown in Figure 5.12, can be selected as before, and so on. From Figure 5.12

$$I^k = I_L^{k+1} + I_R^{k+2}, \tag{5.10}$$

and if we assume equal intervals, then

$$\begin{aligned} I_L^{k+1} &= I_R^{k+1} = I^{k+1}, \\ I_L^{k+2} &= I_R^{k+2} = I^{k+2}. \end{aligned} \tag{5.11}$$

Equation (5.10) gives the recursive relation

$$I^k = I^{k+1} + I^{k+2}. \tag{5.12}$$

If the above procedure is repeated a number of times, a sequence of intervals $\{I_1, I_2, \dots, I_n\}$ will be generated as follows:

$$\begin{aligned} I^1 &= I^2 + I^3, \\ I^2 &= I^3 + I^4, \\ \dots &= \dots \\ I^n &= I^{n+1} + I^{n+2}. \end{aligned} \tag{5.13}$$

The Fibonacci sequence is generated by assuming that the interval for iteration $n+2$ vanishes, that is, $I^{n+2} = 0$. From (5.12), we can write

$$\begin{aligned} I^{n+1} &= I^n - I^{n+2} = I^n \equiv F_0 I^n, \\ I^n &= I^{n+1} + I^{n+2} = I^n \equiv F_1 I^n, \\ I^{n-1} &= I^n + I^{n+1} = 2I^n \equiv F_2 I^n, \\ I^{n-2} &= I^{n-1} + I^n = 3I^n \equiv F_3 I^n, \\ I^{n-3} &= I^{n-2} + I^{n-1} = 5I^n \equiv F_4 I^n, \\ I^{n-4} &= I^{n-3} + I^{n-2} = 8I^n \equiv F_5 I^n, \\ \dots &= \dots \\ I^k &= I^{k+1} + I^{k+2} = F_{n-k+1} I^n, \\ \dots &= \dots \\ I^1 &= I^2 + I^3 = F_n I^n. \end{aligned} \tag{5.14}$$

The sequence generated, namely,

$$\{1,1,2,3,5,8,13,\dots\} = \{F_0, F_1, F_2, F_3, F_4, F_5, F_6 \dots\}, \tag{5.15}$$

is the well-known Fibonacci sequence. It can be generated by using the recursive relation

$$F_k = F_{k-1} + F_{k-2} \quad \text{for } k \geq 2, \tag{5.16}$$

where $F_0 = F_1 = 1$.

If the number of iterations is assumed to be $n$, then from (5.14), the Fibonacci search reduces the interval of uncertainty to

$$I^n = \frac{I^1}{F_n}. \tag{5.17}$$

From (5.14), it follows that $I^n = I^{n-1}/2$. The Fibonacci sequence of intervals can be generated only if $n$ is known. If the objective of the optimization is to find $x^*$ to within a prescribed tolerance, the required $n$ can be readily deduced by using (5.17). Let us assume that the initial bounds of the minimizer, namely, $x_L^1$ and $x_U^1$, and the value of $n$ are given, and a mathematical description of $f(x)$ is available. The implementation consists of computing the successive intervals, evaluating $f(x)$, and selecting the appropriate intervals. At the $k$th iteration, the quantities $x_L^k$, $x_a^k$, $x_b^k$ , $x_U^k$, $I^{k+1}$ and

$$f_a^k = f(x_a^k), \qquad f_b^k = f(x_b^k), \tag{5.18}$$

are known, and the quantities $x_L^{k+1}$, $x_a^{k+1}$, $x_b^{k+1}$, $x_U^{k+1}$, $I^{k+2}$, $f_a^{k+1}$, and $f_b^{k+1}$ are required. Interval $I^{k+2}$ can be obtained from (5.14) as

$$\frac{I^{k+2}}{I^{k+1}} = \frac{F_{n-k-1}I^n}{F_{n-k}I^n} \quad \Rightarrow \quad I^{k+2} = \frac{F_{n-k-1}}{F_{n-k}} I^{k+1}. \tag{5.19}$$

The remaining quantities can be computed as follows.

If $f_a^k > f_b^k$, then $x^*$ is in interval $[x_a^k, x_U^k]$. In this case, as illustrated in Figure 5.13, we assign

$$x_L^{k+1} = x_a^k, \tag{5.20.a}$$
$$x_U^{k+1} = x_U^k, \tag{5.20.b}$$
$$x_a^{k+1} = x_b^k, \tag{5.20.c}$$
$$x_b^{k+1} = x_L^{k+1} + I^{k+2}, \tag{5.20.d}$$
$$f_a^{k+1} = f_b^k, \tag{5.20.e}$$
$$f_b^{k+1} = f\left(x_b^{k+1}\right). \tag{5.20.f}$$

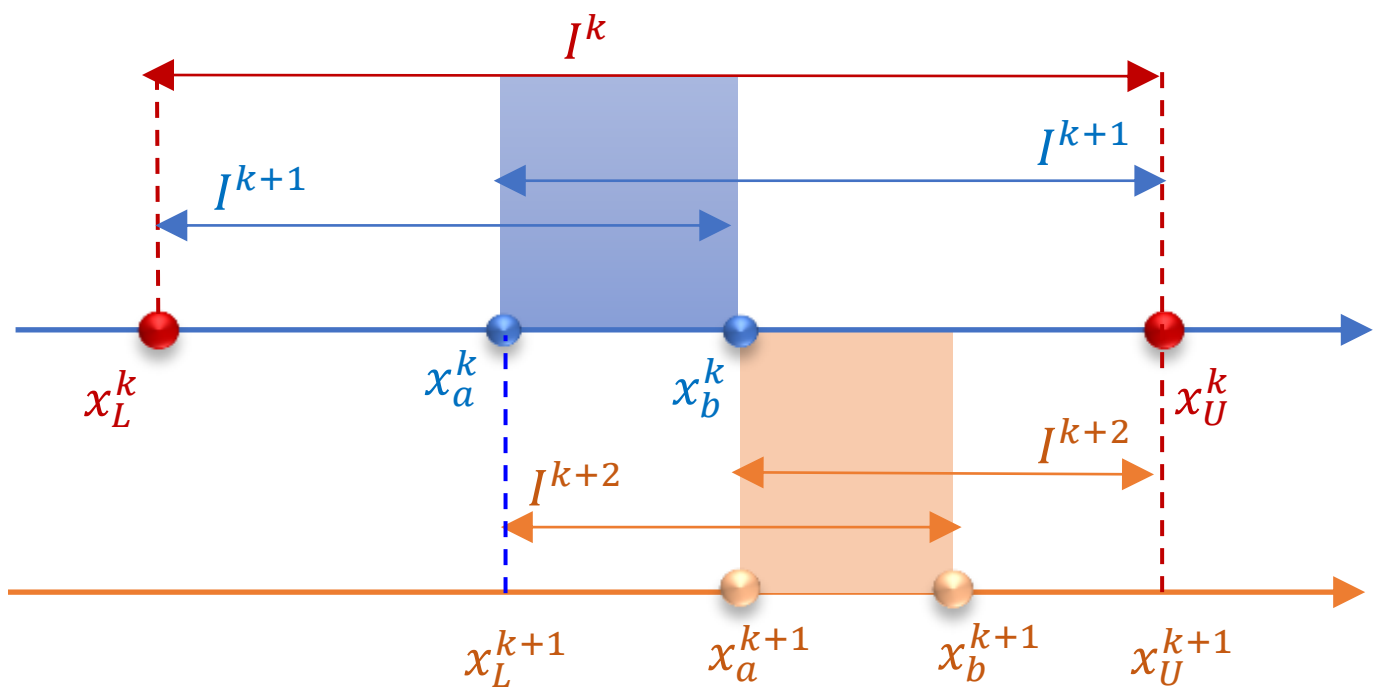

**Figure 5.13.** Assignments in $k^{\text{th}}$ iteration of the Fibonacci method if $f_a^k > f_b^k$.

On the other hand, if $f_a^k < f_b^k$, then $x^*$ is in interval $[x_L^k, x_b^k]$. In this case, as depicted in Figure 5.14, we assign

$$x_L^{k+1} = x_L^k, \tag{5.21.a}$$
$$x_U^{k+1} = x_b^k, \tag{5.21.b}$$
$$x_a^{k+1} = x_U^{k+1} - I^{k+2}, \tag{5.21.c}$$
$$x_b^{k+1} = x_a^k, \tag{5.21.d}$$
$$f_b^{k+1} = f_a^k, \tag{5.21.e}$$
$$f_a^{k+1} = f(x_a^{k+1}). \tag{5.21.f}$$

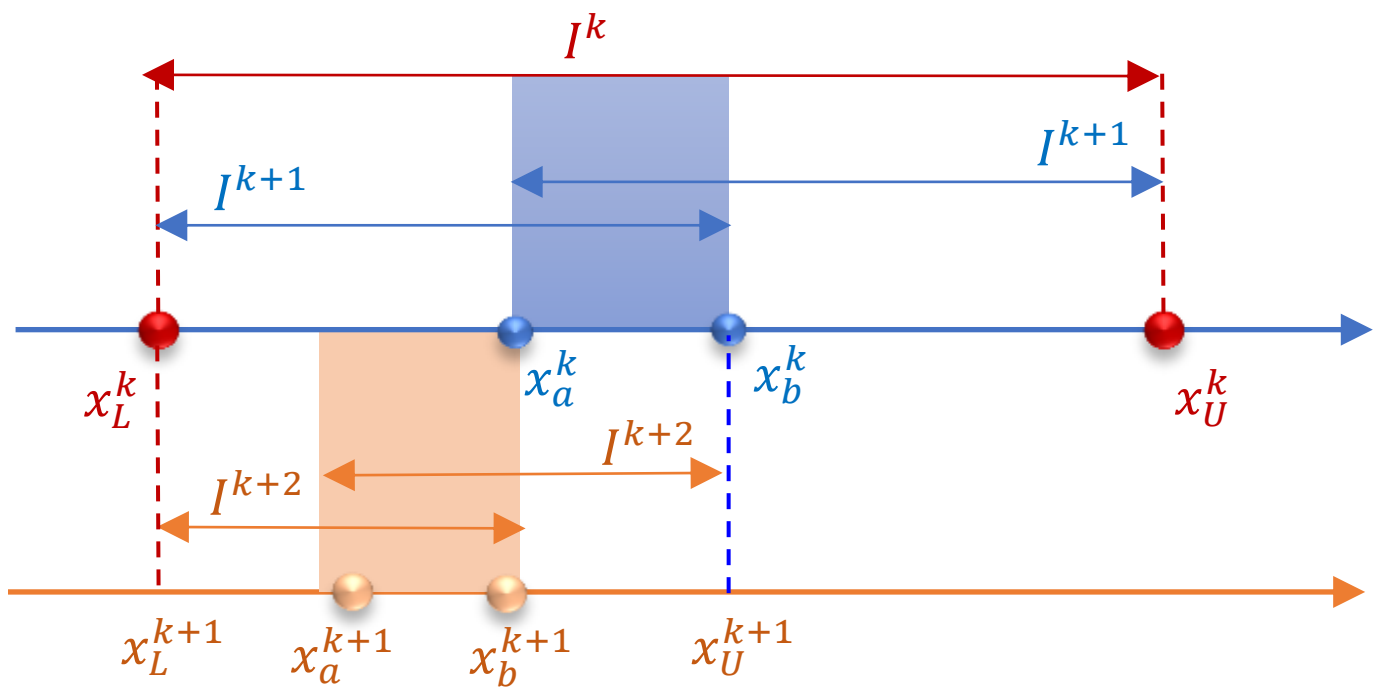

**Figure 5.14.** Assignments in $k^{\text{th}}$ iteration of the Fibonacci method if $f_a^k < f_b^k$.

For the case $f_a^k = f_b^k$, either of the above sets of assignments can be used since $x^*$ is contained by both intervals $[x_L^k, x_b^k]$ and $[x_a^k, x_U^k]$. The above procedure is repeated until $k = n - 2$ in which case $I^{k+2} = I^n$ and

$$x^* = x_a^{k+1} = x_b^{k+1}, \tag{5.22}$$

as depicted in Figure 5.15.

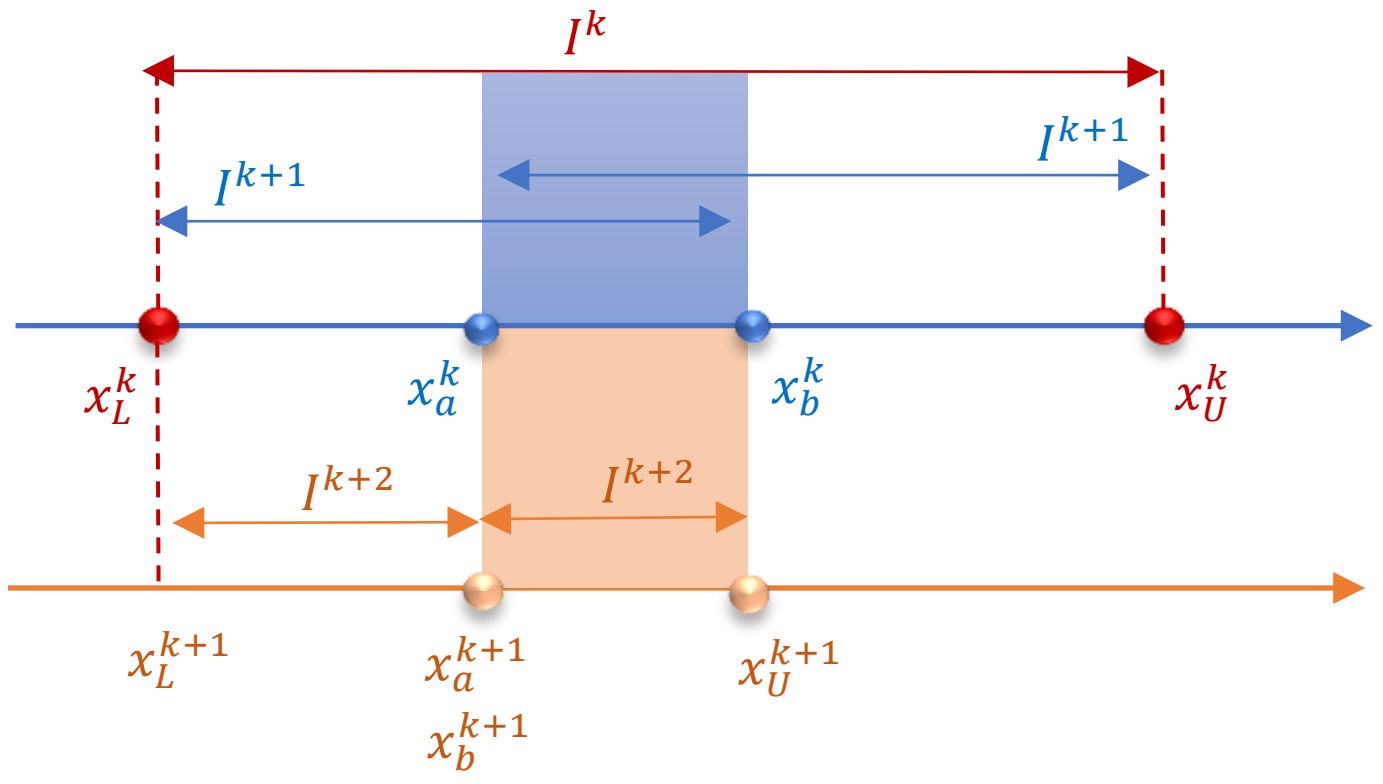

**Figure 5.15.** Assignments in iteration $n - 2$ of the Fibonacci method if $f_a^k < f_b^k$.

Evidently, the minimizer is determined to be within a tolerance $\pm 1/F_n$.

**Algorithm**

**Step 1:** Choose a lower bound $a$ and an upper bound $b$. Set $L = b - a$.
Assume the desired number of function evaluations to be $n$. Set $k = 2$.

**Step 2:** Compute $L_k^* = (F_{n-k+1}/F_{n+1})L$. Set $x_1 = a + L_k^*$ and $x_2 = b - L_k^*$.

**Step 3:** Compute one of $f(x_1)$ or $f(x_2)$, which was not evaluated earlier.
Use the fundamental region-elimination rule to eliminate a region. Set new $a$ and $b$.

**Step 4:** Is $k = n$? If no, set $k = k + 1$ and go to Step 2;
Else Terminate.

***Example 5.5***

Consider the problem:

$$\text{Minimize } f(x) = (x-5)^2 + 6,$$

in the interval $[1,9]$ and $n = 10$.

***Solution***

A plot of the function is shown in Figure 5.16. The plot shows that the minimum lies at $x^* = 5$, $f(x^*) = 6$.

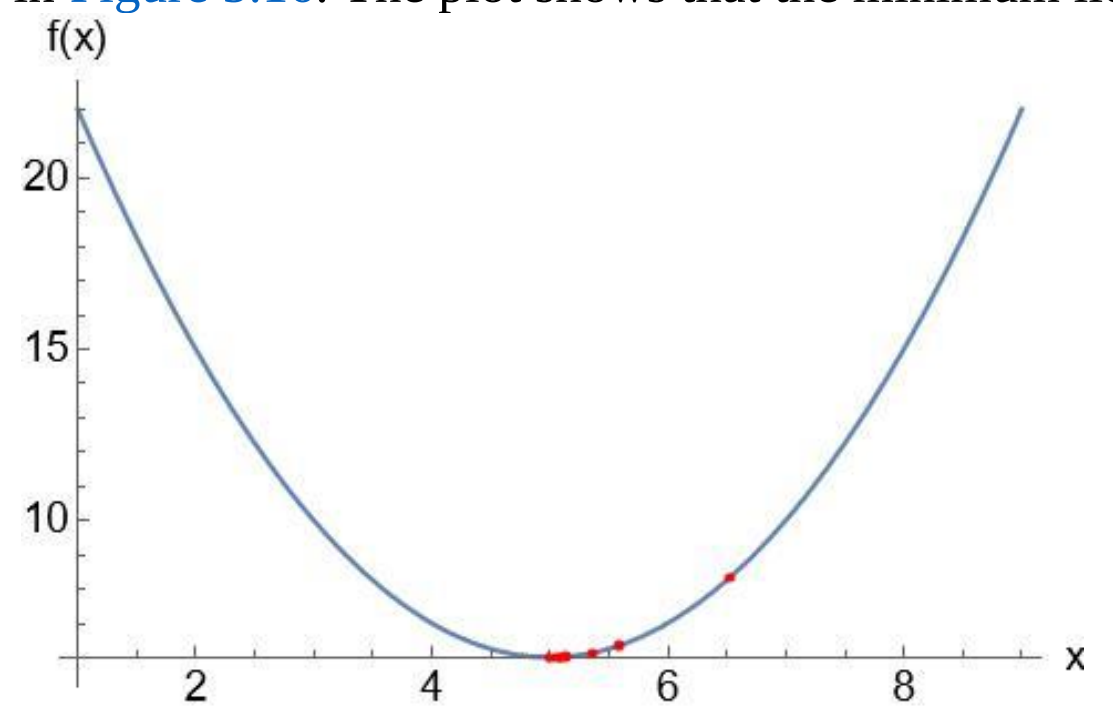

**Figure 5.16.** The results of 10 iterations of the Fibonacci method for $f(x) = (x-5)^2 + 6$.

From Table 5.4, it is obtained that at iteration 10, the optimal condition is attained. Consequently, in ten evaluations $x^*$ has been bracketed within the interval $5.04494 \le x \le 5.13483$ .

**Table 5.4.** The results were produced by Mathematica code 5.4.

| No. it. | $l^*[k]$ | $x_1$ | $x_2$ | $f(x_1)$ | $f(x_2)$ | $a$ | $x^*[k]$ | $b$ | $f(x^*[k])$ |
|---|---|---|---|---|---|---|---|---|---|
| 1 | 3.05618 | 4.05618 | 5.94382 | 6.890797 | 6.890797 | 4.05618 | 6.52809 | 9 | 8.335059 |
| 2 | 1.88764 | 5.94382 | 7.11236 | 6.890797 | 10.46206 | 4.05618 | 5.58427 | 7.11236 | 6.341371 |
| 3 | 1.168539 | 5.224719 | 5.94382 | 6.050499 | 6.890797 | 4.05618 | 5 | 5.94382 | 6 |
| 4 | 0.719101 | 4.775281 | 5.224719 | 6.050499 | 6.050499 | 4.775281 | 5.359551 | 5.94382 | 6.129277 |
| 5 | 0.449438 | 5.224719 | 5.494382 | 6.050499 | 6.244414 | 4.775281 | 5.134831 | 5.494382 | 6.01818 |
| 6 | 0.269663 | 5.044944 | 5.224719 | 6.00202 | 6.050499 | 4.775281 | 5 | 5.224719 | 6 |
| 7 | 0.179775 | 4.955056 | 5.044944 | 6.00202 | 6.00202 | 4.955056 | 5.089888 | 5.224719 | 6.00808 |
| 8 | 0.089888 | 5.044944 | 5.134831 | 6.00202 | 6.01818 | 4.955056 | 5.044944 | 5.134831 | 6.00202 |
| 9 | 0.089888 | 5.044944 | 5.044944 | 6.00202 | 6.00202 | 5.044944 | 5.089888 | 5.134831 | 6.00808 |
| 10 | 0 | 5.044944 | 5.134831 | 6.00202 | 6.01818 | 5.044944 | 5.089888 | 5.134831 | 6.00808 |

***Mathematica Code 5.4*** Fibonacci Search Method

```
(* Fibonacci Search Method *)

(*
Notations
a,b          :Lower and upper limits
n            :The number of function evaluation
f            :Objective function
xstar[k]     :Final design solution of iteration k
f[star[k]]   :Final objective function value of iteration k
lii          :The last iteration index
result[k]    :The results of iteration k
*)

(* Taking Input from User *)
a = Input["Please input the lower limit of the domain where the function is defined"];
b = Input["Please input the upper limit of the domain where the function is defined"];

If[
  a>b,
  Beep[];
  MessageDialog["The value of lower limit a must be less than the upper limit b"];
  Exit[];
  ];

n = Input["Please input the desired number of function evaluation (greater than 2):"];

If[
  !IntegerQ[n]||n<=2,
  Beep[];
  MessageDialog["The value of n must be postive integer greater than 2"];
  Exit[];
  ];

epslion=(2*(b-a))/Fibonacci[n+1];

(* Taking the Function from User *)
f[x_] = Evaluate[Input["Please input a unimodal function of x to find the minimum "]];
(* the user types in, for instance x^2 *)
```

```
(* Initiating Required Variables *)
a0=a;
b0=b;
l=b-a;

(* Starting the Algorithm *)
Do[
  lstar[k_]=(Fibonacci[n-k+1]/Fibonacci[n+1])*l;
  x1=a+lstar[k];
  x2=b-lstar[k];
  f1=f[x1];
  f2=f[x2];

  If[
   f1>=f2,
   a=x1;,
   b=x2;
   ];

  xstar[k]=0.5(a+b);

  lii=k;
  result[k]=N[{k-1,lstar[k],x1,x2,f1,f2,a,xstar[k],b,f[xstar[k]]}];

  If[
   k>n,
   Break[]
   ];,
  {k,2,∞}
  ];

(* Final Result *)
Print["The desired accuracy is ϵ=2(b-a)/F(n+1)= ",N[epslion],"\nThe solution lies between
",N[a]," and ",N[b], " at ",N[xstar[lii]],"\nThe solution is (approximately)",
N[f[xstar[lii]]]];

(* Results of Each Iteration *)
table=TableForm[
  Table[
   result[i], {i,2,lii}
   ],
  TableHeadings->{None,{"No. of
iters.","lstar[k]","x1","x2","f[x1]","f[x2]","a","xstar[k]","b","f[xstar[k]]"}}
  ]

Export["example55.xls",table,"XLS"];

(* Data Visualization *)
Plot[
 f[x],
  {x, a0,b0},
 AxesLabel->{"x","f(x)"},
 LabelStyle->Directive[Black,14],
 Epilog->{PointSize[0.01],Red,Point[Table[{result[i][[8]],result[i][[10]]},{i,2,lii}]]}
 ]

(* Data Manipulation *)
Manipulate[
 Plot[
  f[x],
  {x,a0,b0},
```

```
  AxesLabel->{"x","f(x)"},
  LabelStyle->Directive[Black,14],
  Epilog->{PointSize[0.02],Red,Point[{result[i][[8]],result[i][[10]]}]}
  ],
 {i,2,lii,1}
 ]
```

### 5.4.3 Golden Section Search Method

One difficulty of the Fibonacci search method is that the Fibonacci numbers have to be calculated and stored. Another problem is that at every iteration, the proportion of the eliminated region is not the same. In order to overcome these two problems and yet calculate one new function evaluation per iteration, the golden section search method [1-10] is used. Consider a unimodal function $f$ of one variable and the interval $[a_0, b_0]$. We have to evaluate $f$ at two intermediate points, as illustrated in Figure 5.17. We choose the intermediate points in such a way that the reduction in the range is symmetric, in a sense that

$$a_1 - a_0 = b_0 - b_1 = \rho(b_0 - a_0), \tag{5.23}$$

where,

$$\rho < \frac{1}{2}. \tag{5.24}$$

We then evaluate $f$ at the intermediate points. If $f(a_1) < f(b_1)$, then the minimizer must lie in the range $[a_0, b_1]$. If, on the other hand, $f(a_1) \geq f(b_1)$, then the minimizer is located in the range $[a_1, b_0]$ (see Figure 5.18).

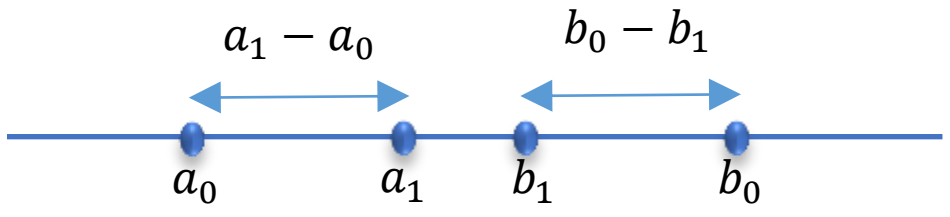

**Figure 5.17.** Evaluating the objective function at two intermediate points.

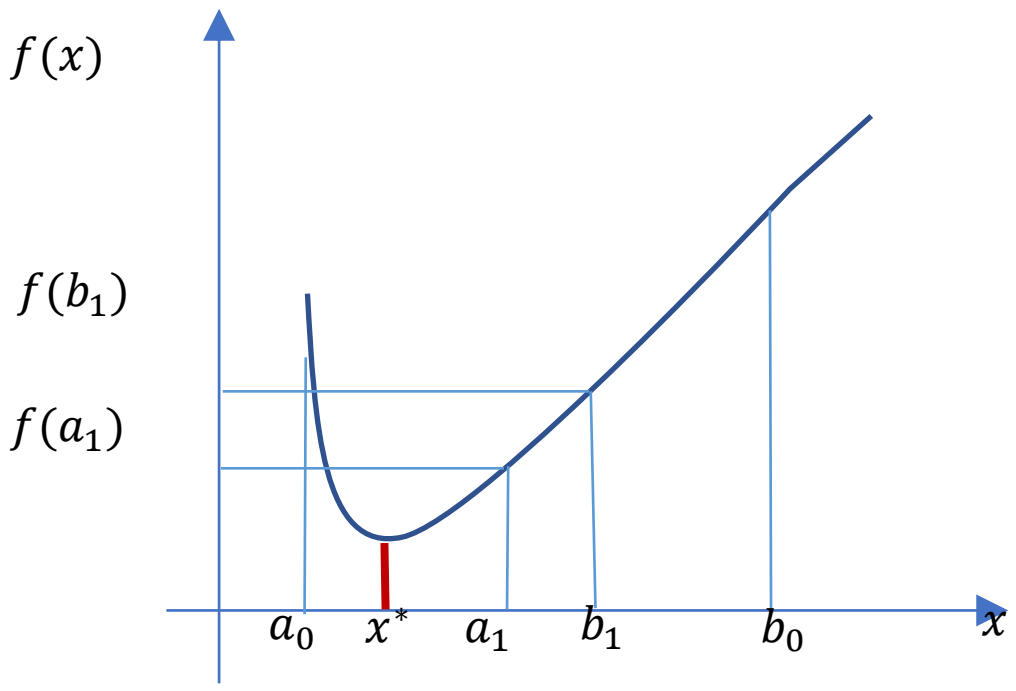

**Figure 5.18.** The case where $f(a_1) < f(b_1)$; the minimizer $x^* \in [a_0, b_1]$.

Starting with the reduced range of uncertainty, we can repeat the process and similarly find two new points, say $a_2$ and $b_2$, using the same value of $\rho < \frac{1}{2}$ as before. However, we would like to minimize the number of objective function evaluations while reducing the width of the uncertainty interval. Suppose, for example, that $f(a_1) < f(b_1)$, as in Figure 5.18. Then, we know that $x^* \in [a_0, b_1]$. Because $a_1$ is already in the uncertainty interval and $f(a_1)$ is already known, we can make $a_1$ coincide with $b_2$. Thus, only one new evaluation of $f$ at $a_2$ would be necessary. To find the value of $\rho$ that results in only one new evaluation of $f$, see Figure 5.19. Without loss of generality, imagine that the original range $[a_o, b_o]$ is of unit length. Then, to have only one new evaluation of $f$ it is enough to choose $\rho$ so that

$$\rho(b_1 - a_0) = b_1 - b_2. \tag{5.25}$$

Because $b_1 - a_0 = 1 - \rho$ and $b_1 - b_2 = 1 - 2\rho$, we have

$$\rho(1 - \rho) = 1 - 2\rho. \tag{5.26}$$

We write the quadratic equation above as

$$\rho^2 - 3\rho + 1 = 0. \tag{5.27}$$

The solutions are

$$\rho_1 = \frac{3 + \sqrt{5}}{2}, \qquad \rho_2 = \frac{3 - \sqrt{5}}{2}. \tag{5.28}$$

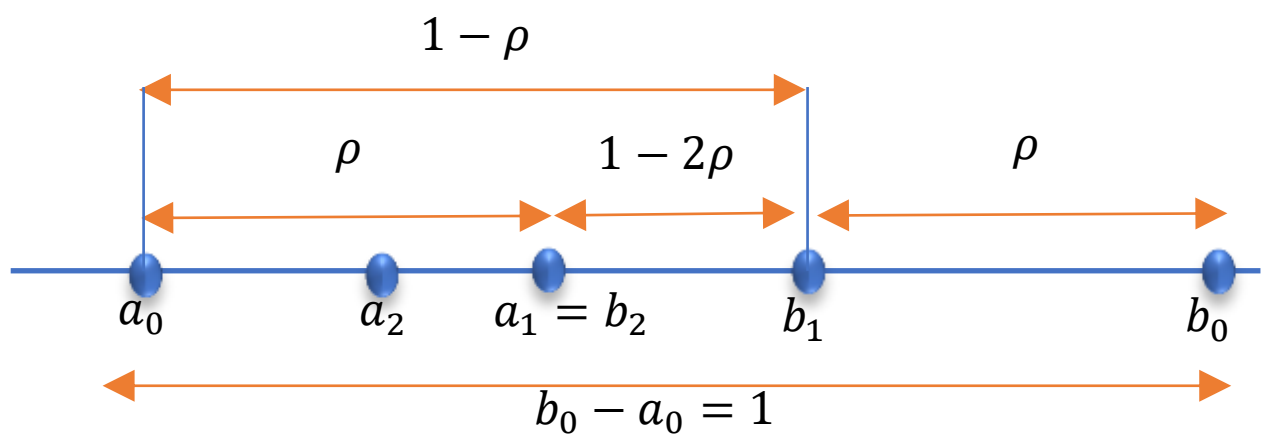

**Figure 2.19.** Finding the value of $\rho$ resulted in only one new evaluation of $f$.

Because we require that $\rho < \frac{1}{2}$, we take

$$\rho = \frac{3 - \sqrt{5}}{2} \approx 0.382. \tag{5.29}$$

Observe that

$$1 - \rho = \frac{\sqrt{5} - 1}{2}, \tag{5.30}$$

and

$$\frac{\rho}{1 - \rho} = \frac{3 - \sqrt{5}}{\sqrt{5} - 1} = \frac{\sqrt{5} - 1}{2} = \frac{1 - \rho}{1}, \tag{5.31}$$

that is,

$$\frac{\rho}{1 - \rho} = \frac{1 - \rho}{1}. \tag{5.32}$$

Thus, dividing a range in the ratio of $\rho$ to $1 - \rho$ has the effect that the ratio of the shorter segment to the longer equals the ratio of the longer to the sum of the two. This rule was referred to as the golden section.

Using the golden section rule means that at every stage of the uncertainty range reduction (except the first), the objective function $f$ needs only be evaluated at one new point. The uncertainty range is reduced by the ratio $1 - \rho \approx 0.61803$ at every stage. Hence, $n$ steps of reduction using the golden section method reduce the range by the factor

$$(1 - \rho)^n \approx (0.61803)^n. \tag{5.33}$$

In this algorithm, the search space $(a, b)$ is first linearly mapped to a unit interval search space $(0,1)$. Thereafter, two points at $\rho$ from either end of the search space are chosen so that at every iteration, the eliminated region is $(1 - \rho)$ to that in the previous iteration

**Algorithm**

**Step 1:** Choose a lower bound $a$ and an upper bound $b$. Also, choose a small number $\epsilon$. Normalize the variable $x$ by using the equation $x_{norm} = (x - a)/(b - a)$. Thus, $a_{norm} = 0$, $b_{norm} = 1$, and $L_{norm} = 1$. Set $k = 1$.

**Step 2:** Set $x_{1\,norm} = a_{norm} + (0.382)L_{norm}$ and $x_{2\,norm} = b_{norm} - (0.382)L_{norm}$. Compute $f(x_{1\,norm})$ or $f(x_{2\,norm})$, depending on whichever of the two was not evaluated earlier. Use the fundamental region-elimination rule to eliminate a region. Set new $a_{norm}$ and $b_{norm}$.

**Step 3:** Is $|L_{norm}| < \epsilon$ small? If no, set $k = k + 1$, go to Step 2;
Else Terminate.

***Example 5.6***

Consider the problem:

$$\text{Minimize } f(x) = (x - 5)^2 + 6,$$

in the interval [3,9] and $\epsilon = 0.01$.

***Solution***

A plot of the function is shown in Figure 5.20. The plot shows that the minimum lies at $x^* = 5$, $f(x^*) = 6$.

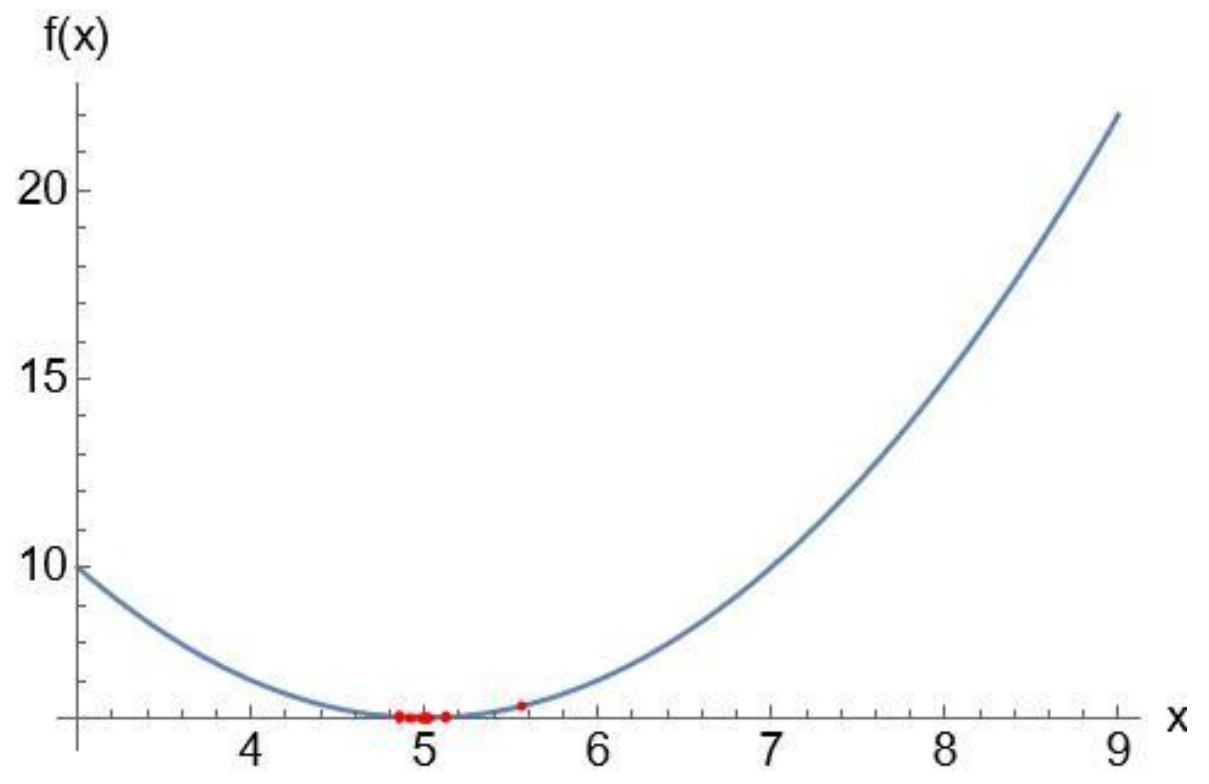

**Figure 5.20.** The results of 11 iterations of the golden section search method for $f(x) = (x - 5)^2 + 6$.

From Table 5.5, it is obtained that at iteration 11, the optimal condition is attained. Consequently, the $x^*$ has been bracketed within the interval $4.98752 \leq x \leq 5.01765$ .

**Table 5.5.a.** The results were produced by Mathematica code 5.5.

| No. it. | $a_{norm}$ | $b_{norm}$ | $l_{norm}$ | $x_{1\,norm}$ | $x_{2\,norm}$ |
|---|---|---|---|---|---|
| 1 | 0 | 1 | 1 | 0.382 | 0.618 |
| 2 | 0 | 0.618 | 0.618 | 0.236076 | 0.381924 |
| 3 | 0.236076 | 0.618 | 0.381924 | 0.381971 | 0.472105 |
| 4 | 0.236076 | 0.472105 | 0.236029 | 0.326239 | 0.381942 |
| 5 | 0.236076 | 0.381942 | 0.145866 | 0.291797 | 0.326221 |
| 6 | 0.291797 | 0.381942 | 0.090145 | 0.326232 | 0.347506 |
| 7 | 0.291797 | 0.347506 | 0.05571 | 0.313078 | 0.326225 |
| 8 | 0.313078 | 0.347506 | 0.034429 | 0.32623 | 0.334355 |
| 9 | 0.32623 | 0.347506 | 0.021277 | 0.334357 | 0.339379 |
| 10 | 0.32623 | 0.339379 | 0.013149 | 0.331253 | 0.334356 |
| 11 | 0.331253 | 0.339379 | 0.008126 | 0.334357 | 0.336275 |

**Table 5.5.b.**

| No. it. | $x_1$ | $x_2$ | $f_1$ | $f_2$ | $a_{new}$ | $b_{new}$ | $x_{star}$ | $f(x_{star})$ |
|---|---|---|---|---|---|---|---|---|
| 1 | 5.292 | 6.708 | 6.085264 | 8.917264 | 3 | 6.708 | 4.854 | 6.021316 |
| 2 | 4.416456 | 5.291544 | 6.340524 | 6.084998 | 4.416456 | 6.708 | 5.562228 | 6.3161 |
| 3 | 5.291826 | 5.83263 | 6.085162 | 6.693273 | 4.416456 | 5.83263 | 5.124543 | 6.015511 |
| 4 | 4.957435 | 5.291652 | 6.001812 | 6.085061 | 4.416456 | 5.291652 | 4.854054 | 6.0213 |
| 5 | 4.750781 | 4.957327 | 6.06211 | 6.001821 | 4.750781 | 5.291652 | 5.021216 | 6.00045 |
| 6 | 4.957393 | 5.085039 | 6.001815 | 6.007232 | 4.750781 | 5.085039 | 4.91791 | 6.006739 |
| 7 | 4.878467 | 4.957352 | 6.01477 | 6.001819 | 4.878467 | 5.085039 | 4.981753 | 6.000333 |
| 8 | 4.957378 | 5.006129 | 6.001817 | 6.000038 | 4.957378 | 5.085039 | 5.021208 | 6.00045 |
| 9 | 5.006144 | 5.036272 | 6.000038 | 6.001316 | 4.957378 | 5.036272 | 4.996825 | 6.00001 |
| 10 | 4.987515 | 5.006135 | 6.000156 | 6.000038 | 4.987515 | 5.036272 | 5.011894 | 6.000141 |
| 11 | 5.006141 | 5.017647 | 6.000038 | 6.000311 | 4.987515 | 5.017647 | 5.002581 | 6.000007 |

***Mathematica Code 5.5*** Golden Section Search Method

```
(* Golden Section Search Method *)

(*
Notations
a,b         :Lower and upper limits
epsilon     :Tolerance
f           :Objective function
lii         :The last iteration index
result[k]   :The results of iteration k
*)

(* Taking Input from User *)
a = Input["Please input the lower limit of the domain where the function is defined"];
b = Input["Please input the upper limit of the domain where the function is defined"];

If[
  a>b,
  Beep[];
  MessageDialog["The value of lower limit a must be less than the upper limit b"];
  Exit[];
  ];

epsilon=Input["Please enter a tolerance (small postive number):"];(* A very small number to
check whether the loop would be executed or not*)

If[
  epsilon<=0,
  Beep[];
  MessageDialog["Tolerance value has to be small postive number: "];
  Exit[];
  ];

(* Taking the Function from User *)
f[x_] = Evaluate[Input["Please input a unimodal function of x to find the minimum "]];
(* The user types in, for instance x^2 *)

(*Normalize the Variable x and Initiating Required Variables*)
a0=a;
b0=b;

anew=a;
```

```
bnew=b;

(* Starting the Algorithm *)
Do[
  anorm=(anew-a)/(b-a);
  bnorm=(bnew-a)/(b-a);

  lnorm=bnorm-anorm;

  x1norm=anorm+0.382*lnorm;
  x2norm=bnorm-0.382*lnorm;

  x1=x1norm(b0-a0)+a0;
  x2=x2norm(b0-a0)+a0;

  f1=f[x1];
  f2=f[x2];

  Which [f1>f2,
   anew=x1(*move lower bound to x1*);,
   f1<f2,
   bnew=x2(*move upper bound to x2*);,
   f1==f2,
   anew=x1(*move lower bound to x1*);
   bnew=x2(*move upper bound to x2*);];

  xstar=0.5*(anew+bnew);

  lii=k;
  result[k]=N[{k,anorm,bnorm,lnorm,x1norm,x2norm,x1,x2,f1,f2,anew,bnew,xstar,f[xstar]}];

  If[
   Abs[lnorm]<epsilon,
   Break[]
   ];,
  {k,1,∞}
  ];

(* Final Result *)
Print["The solution lies between ",N[anew]," and ",N[bnew], " at ",N[xstar],"\nThe solution 
is (approximately)", N[f[xstar]]];

(* Results of Each Iteration *)
table=TableForm[
  Table[
   result[i], {i,1,lii}
   ],
  TableHeadings->{None,{"No. of 
iters.","anorm","bnorm","lnorm","x1norm","x2norm","x1","x2","f1","f2","anew","bnew","xstar",
"f[xstar]"}}
  ]

Export["example56.xls",table,"XLS"];

(* Data Visualization *)
Plot[
 f[x],
  {x, a0,b0},
 AxesLabel->{"x","f(x)"},
 LabelStyle->Directive[Black,14],
 Epilog->{PointSize[0.01],Red,Point[Table[{result[i][[13]],result[i][[14]]},{i,1,lii}]]}
```

```
 ]

(* Data Manipulation *)
Manipulate[
 Plot[
  f[x],
  {x,a0,b0},
  AxesLabel->{"x","f(x)"},
  LabelStyle->Directive[Black,14],
  Epilog->{PointSize[0.02],Red,Point[{result[i][[13]],result[i][[14]]}]}
  ],
 {i,1,lii,1}
 ]
```

## 5.5 Successive Quadratic Estimation Method (Powell Method)

In the approximation approach to one-dimensional optimization, an approximate expression for the objective function is assumed, usually in the form of a low-order polynomial. If a second-order polynomial of the form, [11]

$$q(x) = a_0 + a_1 x + a_2 x^2, \tag{5.34}$$

is assumed where $a_0$, $a_1$, and $a_2$ are constants, a quadratic interpolation method is obtained. In other words, a quadratic estimation scheme assumes that within the bounded interval the function can be approximated by a quadratic function, and this approximation will improve as the points used to construct the approximation approach the actual minimum.

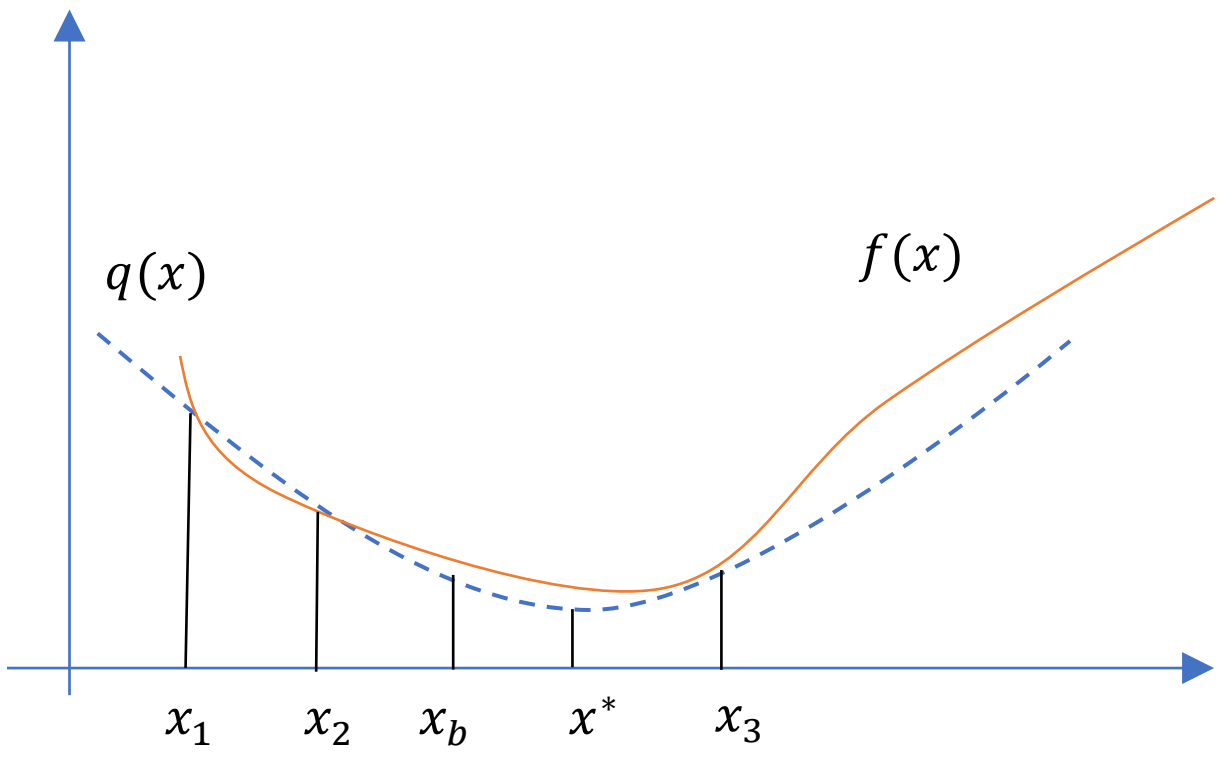

**Figure 5.21.** The function $f(x)$ and the interpolated quadratic function.

Since any quadratic function can be defined with three points, the algorithm begins with three initial points. Figure 5.21 shows the original function and three initial points $x_1$, $x_2$, and $x_3$. The fitted quadratic curve through these three points is also plotted with a dashed line. The minimum $\bar{x}$ of this curve is used as one of the candidate points for the next iteration. For non-quadratic functions, a number of iterations of this algorithm is necessary, whereas for quadratic objective functions the exact minimum can be found in one iteration only.

Given three consecutive points $x_1$, $x_2$, $x_3$ and their corresponding function values $f_1 = f(x_1)$, $f_2 = f(x_2)$, $f_3 = f(x_3)$, we seek to determine three constants $a_1$, $a_2$, and $a_3$ such that the quadratic function

$$q(x) = a_0 + a_1(x - x_1) + a_2(x - x_1)(x - x_2), \tag{5.35}$$

agrees with $f(x)$ at these three points. We proceed by evaluating $q(x)$ at each of the three given points. We have

$$\begin{aligned} q(x_1) &= f_1 = a_0 + 0 \times a_1 + 0 \times a_2, \\ q(x_2) &= f_2 = a_0 + a_1(x_2 - x_1) + 0 \times a_2, \\ q(x_3) &= f_3 = a_0 + a_1(x_3 - x_1) + a_2(x_3 - x_1)(x_3 - x_2). \end{aligned} \tag{5.36}$$

In the matrix form

$$\begin{pmatrix} f_1 \\ f_2 \\ f_3 \end{pmatrix} = \begin{pmatrix} 1 & 0 & 0 \\ 1 & x_2 - x_1 & 0 \\ 1 & x_3 - x_1 & (x_3 - x_1)(x_3 - x_2) \end{pmatrix} \begin{pmatrix} a_0 \\ a_1 \\ a_2 \end{pmatrix}. \tag{5.37}$$

We have

$$a_0 = f_1, \qquad a_1 = \frac{f_2 - f_1}{x_2 - x_1}, \qquad a_2 = \frac{1}{(x_3 - x_2)} \left( \frac{f_3 - f_1}{x_3 - x_1} - \frac{f_2 - f_1}{x_2 - x_1} \right). \tag{5.38}$$

Hence, given three points and their function values, the quadratic estimate can be constructed by simply evaluating the expressions for $a_0$, $a_1$, $a_2$ given above.

Now following the proposed strategy, if the approximating quadratic is a good approximation to the function to be optimized over the interval $x_1$ to $x_3$, then it can be used to predict the location of the optimum. Recall that the stationary points of a single-variable function can be determined by setting its first derivative to zero and solving for the roots of the resulting equation. In the case of our quadratic approximating function,

$$\frac{dq}{dx} = a_1 + a_2(x - x_2) + a_2(x - x_1) = 0, \tag{5.39}$$

can be solved to yield the estimate

$$\bar{x} = \frac{x_2 + x_1}{2} - \frac{a_1}{2a_2}. \tag{5.40}$$

Since the function $f(x)$ is unimodal over the interval in question and since the approximating quadratic is also a unimodal function, it is reasonable to expect that $\bar{x}$ will be a good estimate of the desired exact optimum $x^*$.

The above approach constitutes one iteration of the quadratic interpolation method. If $f(x)$ cannot be represented accurately by a second-order polynomial, a number of such iterations can be performed. The appropriate strategy is to attempt to reduce the interval of uncertainty in each iteration as was done in the search methods. This can be achieved by rejecting either $x_1$ or $x_3$ and then using the two remaining points along with the point $\bar{x}$ for a new interpolation. This procedure continues until two consecutive estimates are close to each other. After a number of iterations, the three points will be in the neighborhood of $x^*$. Based on these results, we present the following Powell's algorithm.

**Algorithm**

**Step 1:** Let $x_1$ be an initial point and $\Delta$ be the step size. Compute $x_2 = x_1 + \Delta$.

**Step 2:** Evaluate $f(x_1)$ and $f(x_2)$.

**Step 3:** If $f(x_1) > f(x_2)$, let $x_3 = x_1 + 2\Delta$;
Else let $x_3 = x_1 - \Delta$. Evaluate $f(x_3)$.

**Step 4:** Determine $F_{\min} = \min(f_1, f_2, f_3)$ and $X_{\min}$ is the point $x_i$ that corresponds to $F_{\min}$.

**Step 5:** Use points $x_1$, $x_2$, and $x_3$ to calculate $\bar{x}$ using (5.40).

**Step 6:** Are $e_f = |(F_{\min} - f(\bar{x}))/f(\bar{x})|$ and $e_x = |(X_{\min} - \bar{x})/\bar{x}|$ small? If not, go to Step 7;
Else the optimum is the best of the current four points, and Terminate.

**Step 7:** Save the best point and two bracketing it, if possible; otherwise, save the best three points.
Relabel them according to $x_1 < x_2 < x_3$ and go to Step 4.

***Example 5.7.1***

Consider the problem:

$$\text{Minimize } f(x) = (x-5)^2 + 6$$

given the starting point $x_0 = 9$, a step size $|\Delta| = 0.1$, $\epsilon x = 0.001$ and $\epsilon f = 0.001$.

***Solution***

A plot of the function is shown in Figure 5.22. The plot shows that the minimum lies at $x^* = 5$, $f(x^*) = 6$.

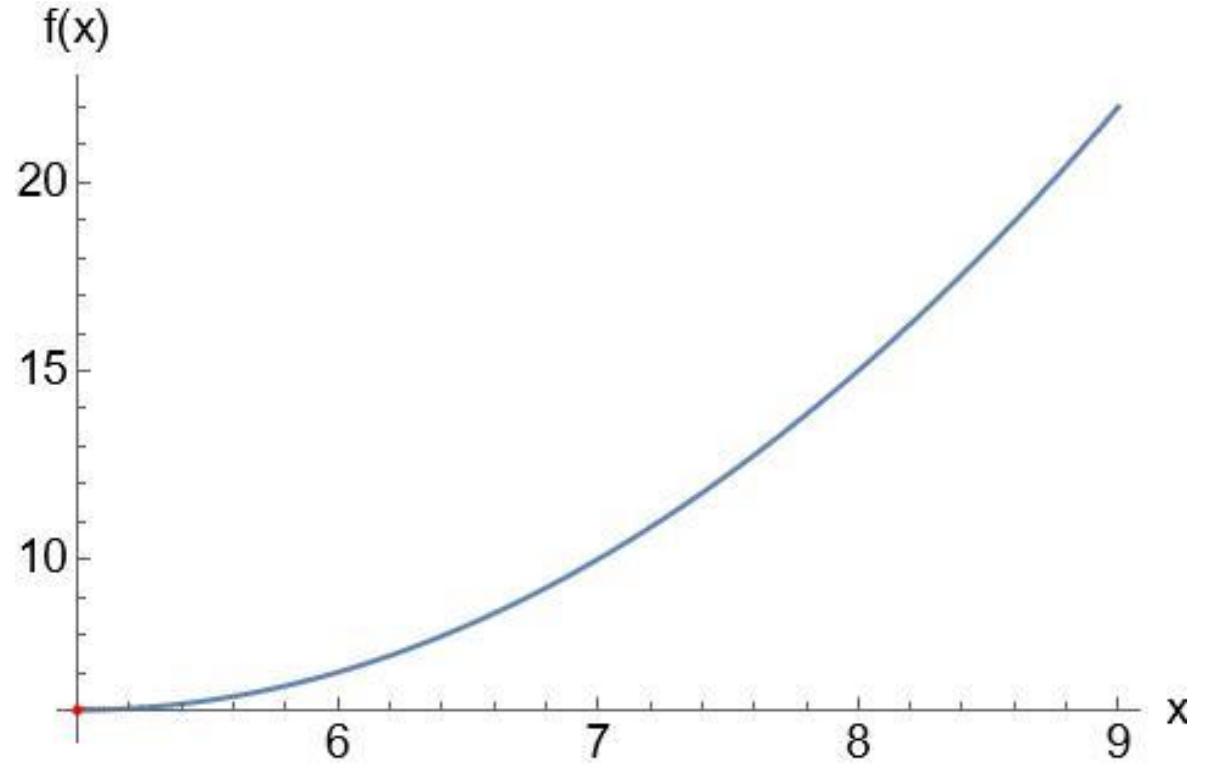

**Figure 5.22.** The results of 2 iterations of the successive quadratic method for $f(x) = (x-5)^2 + 6$.

From Table 5.6, it is obtained that at iteration 2, the optimal condition is attained. Consequently, $x^* = 5$.

**Table 5.6.**

| No. it. | $x_1$ | $x_2$ | $x_3$ | $\bar{x}$ | $x_{min}$ | $f(x_1)$ | $f(x_2)$ | $f(x_3)$ | $f(\bar{x})$ | $f(x_{min})$ | $e_f$ | $e_x$ |
|---|---|---|---|---|---|---|---|---|---|---|---|---|
| 1 | 8.9 | 9 | 9.1 | 5 | 8.9 | 21.21 | 22 | 22.81 | 6 | 21.21 | 2.535 | 0.78 |
| 2 | 5 | 8.9 | 9 | 5 | 5 | 6 | 21.21 | 22 | 6 | 6 | 0 | 1.97E-14 |

***Example 5.7.2***

Consider the problem:

$$\text{Minimize } f(x) = 0.5x^2 + 50/x$$

given the starting point $x_0 = 10$, a step size $|\Delta| = 0.1$, $\epsilon x = 0.001$ and $\epsilon f = 0.001$.

***Solution***

A plot of the function is shown in Figure 5.23. The plot shows that the minimum lies at $x^* = 3.6903$.

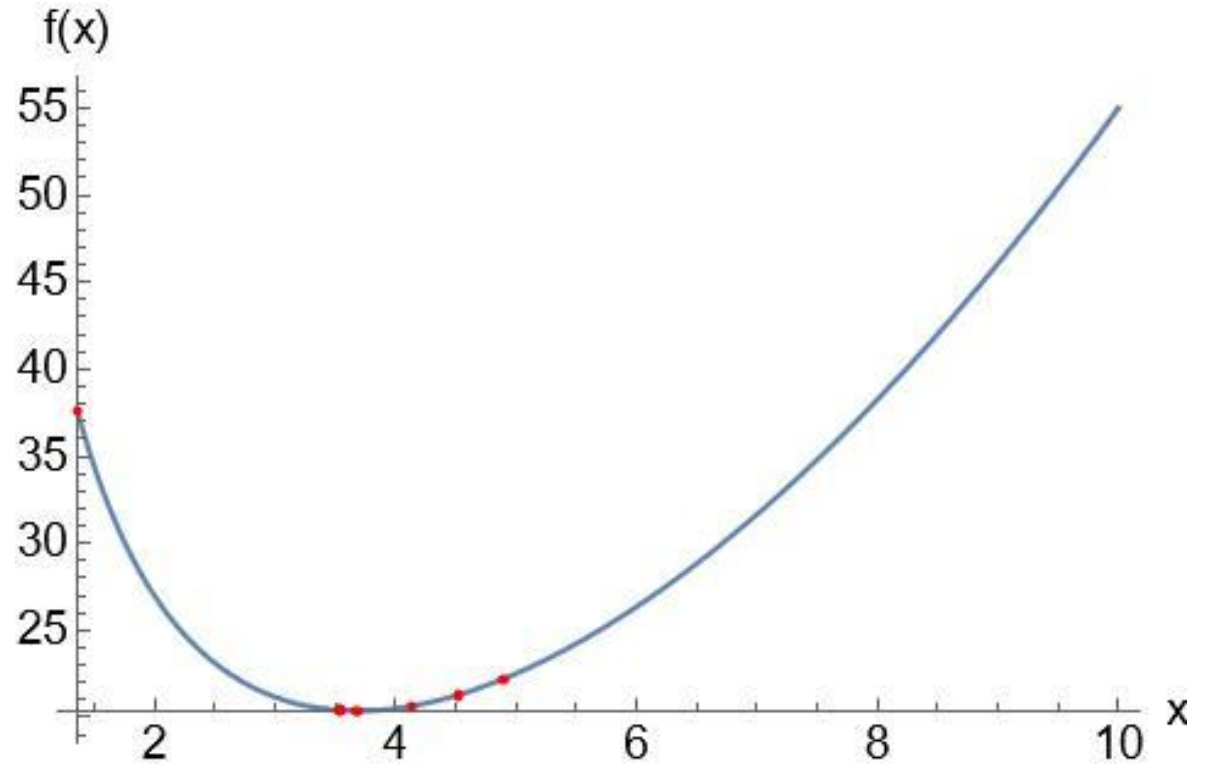

**Figure 5.23.** The results of 7 iterations of the successive quadratic method for $f(x) = 0.5x^2 + 50/x$.

From Table 5.7, it is obtained that at iteration 7, the optimal condition is attained. Consequently, $x^* = 3.6903$.

**Table 5.7.1.** The results were produced by Mathematica code 5.6.

| No. it. | $x_1$ | $x_2$ | $x_3$ | $\bar{x}$ | $x_{\min}$ |
|---|---|---|---|---|---|
| 1 | 9.9 | 10 | 10.1 | 1.36376 | 9.9 |
| 2 | 1.36376 | 9.9 | 10 | 4.523968 | 1.36376 |
| 3 | 1.36376 | 4.523968 | 9.9 | 4.900613 | 4.523968 |
| 4 | 1.36376 | 4.523968 | 4.900613 | 4.141879 | 4.523968 |
| 5 | 4.141879 | 4.523968 | 4.900613 | 3.53613 | 4.141879 |
| 6 | 3.53613 | 4.141879 | 4.523968 | 3.669568 | 3.53613 |
| 7 | 3.53613 | 3.669568 | 4.141879 | 3.690382 | 3.669568 |

**Table 5.7.2.**

| No. it. | $f(x_1)$ | $f(x_2)$ | $f(x_3)$ | $f(\bar{x})$ | $f(x_{\min})$ | $e_f$ | $e_x$ |
|---|---|---|---|---|---|---|---|
| 1 | 54.05551 | 55 | 55.9555 | 37.59325 | 54.05551 | 0.437904 | 6.25934 |
| 2 | 37.59325 | 54.05551 | 55 | 21.28539 | 37.59325 | 0.766153 | 0.698548 |
| 3 | 37.59325 | 21.28539 | 54.05551 | 22.21081 | 21.28539 | 0.041665 | 0.076857 |
| 4 | 37.59325 | 21.28539 | 22.21081 | 20.6494 | 21.28539 | 0.030799 | 0.09225 |
| 5 | 20.6494 | 21.28539 | 22.21081 | 20.39186 | 20.6494 | 0.012629 | 0.171303 |
| 6 | 20.39186 | 20.6494 | 21.28539 | 20.35845 | 20.39186 | 0.001641 | 0.036363 |
| 7 | 20.39186 | 20.35845 | 20.6494 | 20.35819 | 20.35845 | 1.25E-05 | 0.00564 |

***Mathematica Code 5.6*** Successive Quadratic Method

```
(* Successive Quadratic Method *)

(*
Notations
x1         :Initial guess of the min
delta      :The step size
f          :Objective function of the varible x
epsilonf   :Small number to check whether the loop would be executed or not
epsilonx   :Small number to check whether the loop would be executed or not
lii        :The last iteration index
result[k] :The results of iteration k
*)

(* Taking Input from User *)
x1=Input["Enter the initial point:"];
delta=Input["Enter the step size:"];

epsilonf=Input["Please enter a small number epsilonf: function tolerance "];
epsilonx=Input["Please enter a small number epsilonx: x variable tolerance "];
(* A very small number to check whether the loop would be executed or not*)

If[
  epsilonf<=0||epsilonx<=0,
  Beep[];
  MessageDialog["Tolerance value has to be small postive number: "];
  Exit[];
  ];

(* Taking the Function from User *)
f[x_] = Evaluate[Input["Please input a function of x to find the minimum "]]; (* the user
types in, for instance x^2 *)

(* Initiating Required Variables *)
a0=x1;
```

```
x2=x1+delta;
If[
  f[x1]>f[x2],
  x3=x1+2*delta,
  x3=x1-delta
  ];

{x1,x2,x3}=Sort[{x1,x2,x3}];

(* Starting the Algorithm *)
Do[
 y1=f[x1];
 y2=f[x2];
 y3=f[x3];

 fmin=Min[y1,y2,y3];

 Which[
  fmin==y1,
  xmin=x1,
  fmin==y2,
  xmin=x2,
  fmin==y3,
  xmin=x3
  ];

 a1=(y2-y1)/(x2-x1);
 a2=1/(x3-x2)*((y3-y1)/(x3-x1)-(y2-y1)/(x2-x1));

 xbar=(x2+x1)/2-a1/(2*a2);
 ybar=f[xbar];

 ef=Abs[(fmin-ybar)/ybar];
 ex=Abs[(xmin-xbar)/xbar];

 lii=k;
 result[k]=N[{k,x1,x2,x3,xbar,xmin,f[x1],f[x2],f[x3],f[xbar],fmin,ef,ex}];

 If[
  ef<=epsilonf||ex<=epsilonx,
  Break[];,
  {x1,x2,x3}=Sort[
     Take[
      SortBy[
       {x1,x2,x3,xbar},
       f](* Sort the four points (x1, x2, x3, xbar) by the values of the function f[x] *),
      3]
     ](* Take the first three points *);
  ];,
 {k,1,∞}
 ]

(* Final Result *)
Print["The solution at ",N[xbar],"\nThe solution is (approximately)", N[f[xbar]]];

(* Results of Each Iteration *)
table=TableForm[
  Table[
   result[i],
   {i,1,lii}
```

```
    ],
  TableHeadings->{None,{"No. of
iters.","x1","x2","x3","xbar","N[xmin]","f[x1]","f[x2]","f[x3]","f[xbar]","N[fmin]","N[ef]",
"N[ex]"}}
  ]

Export["example57.xls",table,"XLS"];

(* Data Visualization *)
llimt=Min[Table[result[i][[5]],{i,1,lii}]];

If[
  a0<llimt,
  {a0,llimit}={llimit,a0}
  ];

Plot[
 f[x],
 {x, llimt, a0},
 AxesLabel->{"x","f(x)"},
 LabelStyle->Directive[Black,14],
 Epilog->{PointSize[0.01],Red,Point[Table[{result[i][[5]],result[i][[10]]},{i,1,lii}]]}
 ]

(* Data Manipulation *)
Manipulate[
 Plot[
  f[x],
  {x,llimt, a0},
  AxesLabel->{"x","f(x)"},
  LabelStyle->Directive[Black,14],
  Epilog->{PointSize[0.02],Red,Point[{result[i][[5]],result[i][[10]]}]}
  ],
 {i,1,lii,1}
 ]
```

## 5.6 Gradient-based Methods

### 5.6.1 Newton-Raphson Method

The Newton-Raphson method [1-10] requires that the function $f$ be twice differentiable, i.e., at each measurement point $x^{(k)}$ we can determine $f\left(x^{(k)}\right)$, $f'\left(x^{(k)}\right)$, and $f''\left(x^{(k)}\right)$. In this method, a linear approximation to the first derivative of the function is made at a point using Taylor's series expansion.

Using Taylor's series expansion, we have

$$f(x+h) = f(x) + f'(x)h + \frac{1}{2}f''(x)h^2 + \cdots. \tag{5.41a}$$

$$x + h = x^{(k+1)}, \quad x = x^{(k)} \quad \text{and} \quad h = x^{(k+1)} - x^{(k)}, \tag{5.41b}$$

$$f\left(x^{(k+1)}\right) = q(x^{(k+1)}) = f\left(x^{(k)}\right) + f'\left(x^{(k)}\right)\left(x^{(k+1)} - x^{(k)}\right) + \frac{1}{2}f''\left(x^{(k)}\right)\left(x^{(k+1)} - x^{(k)}\right)^2 + \cdots. \tag{5.41c}$$

Instead of minimizing $f$, we minimize its approximation. The first-order necessary condition for a minimizer yields

$$f'\left(x^{(k+1)}\right) \approx f'\left(x^{(k)}\right) + f''\left(x^{(k)}\right)\left(x^{(k+1)} - x^{(k)}\right) = 0, \tag{5.42}$$

we obtain

$$x^{(k+1)} = x^{(k)} - \frac{f'\left(x^{(k)}\right)}{f''\left(x^{(k)}\right)}. \tag{5.43}$$

Newton's method works well if $f''(x) > 0$ everywhere (see Figure 5.24). However, if $f''(x) < 0$ for some $x$, Newton's method may fail to converge to the minimizer (see Figure 5.25).

Indeed, if we set $g(x) = f'(x)$, then we obtain a formula for the iterative solution of the equation $g(x) = 0$:

$$x^{(k+1)} = x^{(k)} - \frac{g\left(x^{(k)}\right)}{g'\left(x^{(k)}\right)}. \tag{5.44}$$

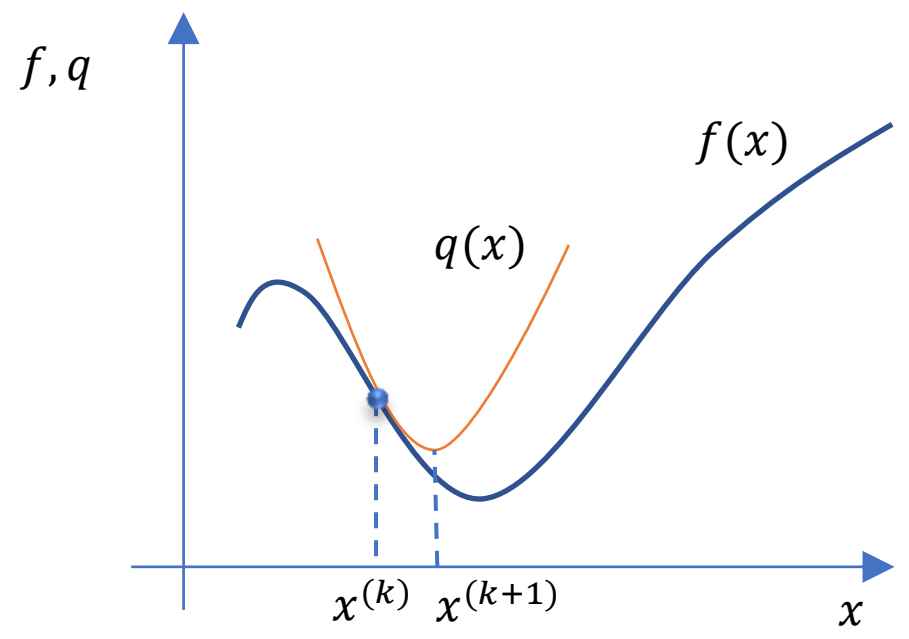

**Figure 5.24.** Newton's algorithm with $f''(x) > 0$.

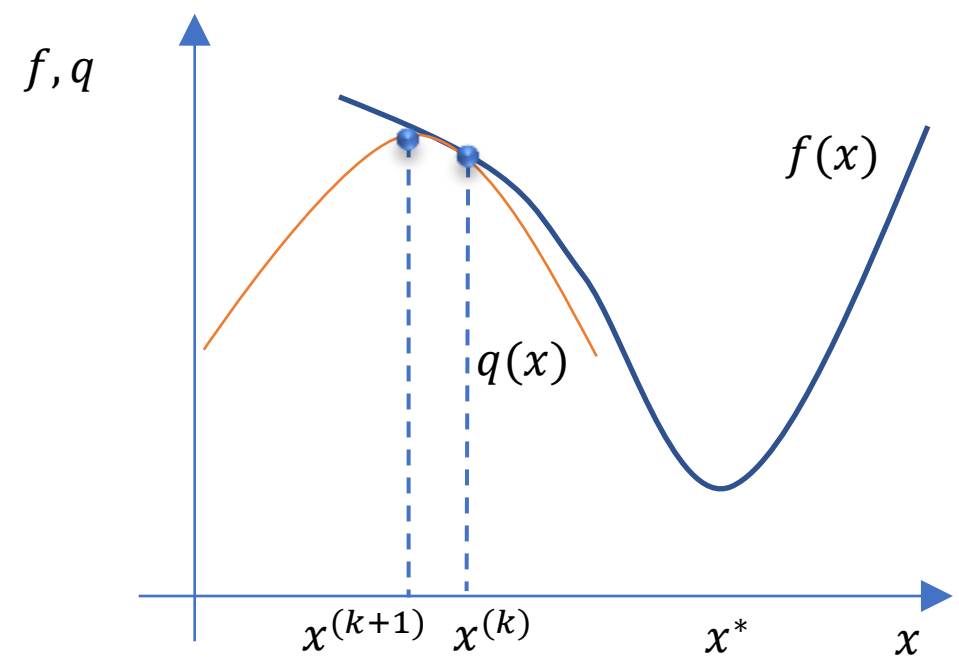

**Figure 5.25.** Newton's algorithm with $f''(x) < 0$.

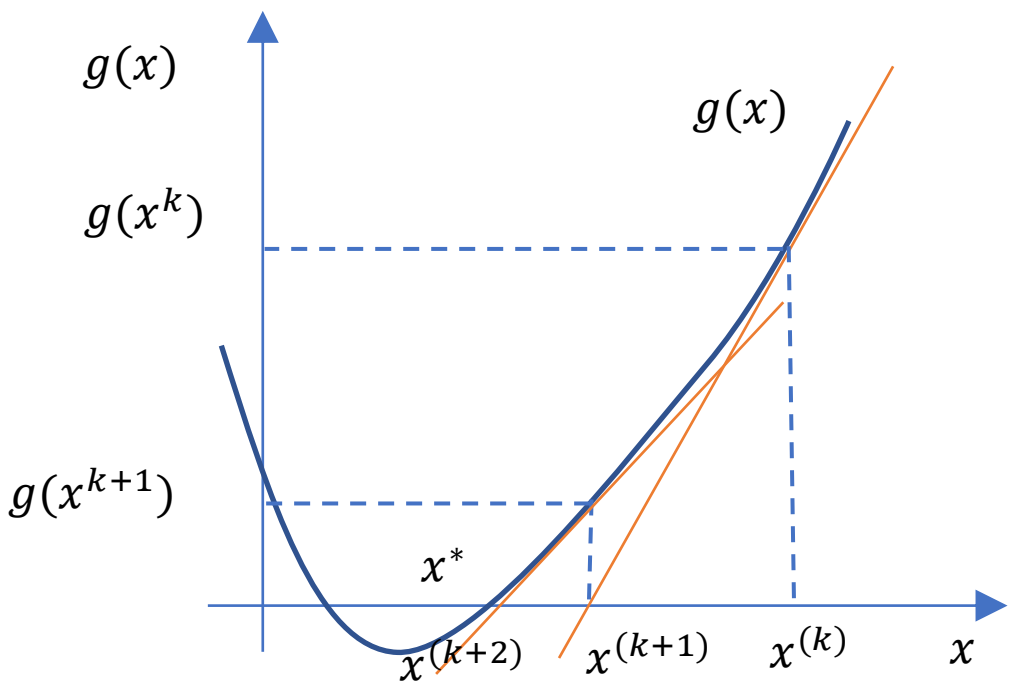

**Figure 5.26.** Newton's method of tangents.

Newton's method for solving equations of the form $g(x) = 0$ is also referred to as Newton's method of tangents. This name is easily justified if we look at a geometric interpretation of the method when applied to the solution of the equation $g(x) = 0$ (see Figure 5.26). If we draw a tangent to $g(x)$ at the given point $x^{(k)}$, then the tangent line intersects the $x$-axis at the point $x^{(k+1)}$ which we expect to be closer to the root $x^*$ of $g(x) = 0$. Note that the slope of $g(x)$ at $x^{(k)}$ is

$$g'\left(x^{(k)}\right) = \frac{g\left(x^{(k)}\right)}{x^{(k)} - x^{(k+1)}}, \tag{5.45.a}$$

$$x^{(k+1)} = x^{(k)} - \frac{g\left(x^{(k)}\right)}{g'\left(x^{(k)}\right)}. \tag{5.45.b}$$

The convergence of the algorithm depends on the initial point and the nature of the objective function. For mathematical functions, the derivative may be easy to compute, but in practice, the gradients have to be computed numerically.

**Algorithm**

**Step 1:** Choose initial guess $x^{(1)}$ and a small number $\epsilon$. Set $k = 1$. Compute $f'(x^{(1)})$.

**Step 2:** Compute $f''(x^{(k)})$.

**Step 3:** Calculate $x^{(k+1)} = x^{(k)} - f'(x^{(k)})/f''(x^{(k)})$. Compute $f'(x^{(k+1)})$.

**Step 4:** If $|f'(x^{(k+1)})| < \epsilon$, Terminate;
Else set $k = k + 1$ and go to Step 2.

***Example 5.8.1***

Consider the problem:

$$\text{Minimize } f(x) = (x - 5)^2 + 6,$$

given the starting point $x_0 = 9$, $\epsilon x = 0.001$ and $\epsilon f = 0.001$.

***Solution***

A plot of the function is shown in Figure 5.27. The plot shows that the minimum lies at $x^* = 5$, $f(x^*) = 6$.

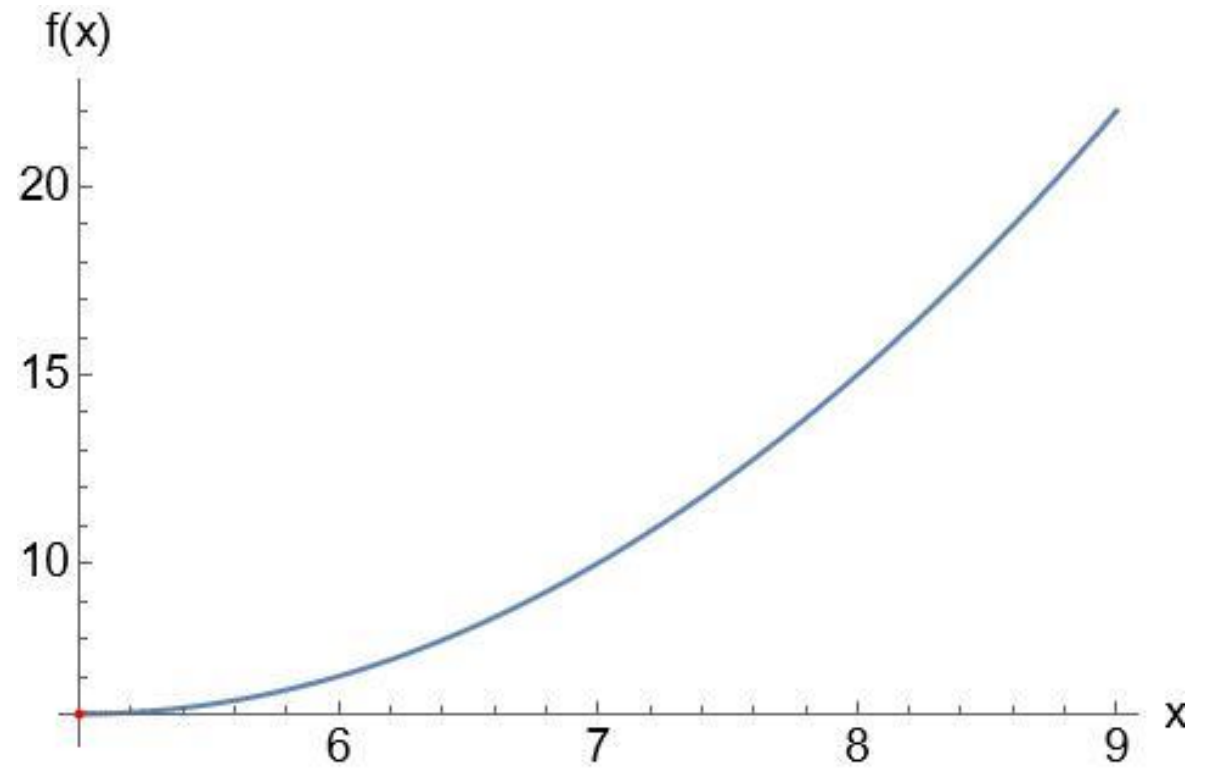

**Figure 5.27.** The results of 1 iteration of the Newton-Raphson method for $f(x) = (x - 5)^2 + 6$.

From Table 5.8, it is obtained that at iteration 1, the optimal condition is attained. Consequently, the $x^* = 5$.

**Table 5.8.**

| No. it. | $x_1$ | $x_2$ | $f(x_1)$ | $f(x_2)$ | $f'(x_1)$ | $f'(x_2)$ | $\lvert x_2 - x_1\rvert$ |
|---|---|---|---|---|---|---|---|
| 1 | 9 | 5 | 22 | 6 | 8 | 0 | 4 |

***Example 5.8.2***

Using Newton's method, find the minimizer of $f(x) = 2x^2 + \frac{100}{x}$, given the starting point $x_0 = 10$, $\epsilon x = 0.001$ and $\epsilon f = 0.001$.

***Solution***

A plot of the function is shown in Figure 5.28. The plot shows that the minimum lies at $x^* = 2.92402$, $f(x^*) = 51.2992$.

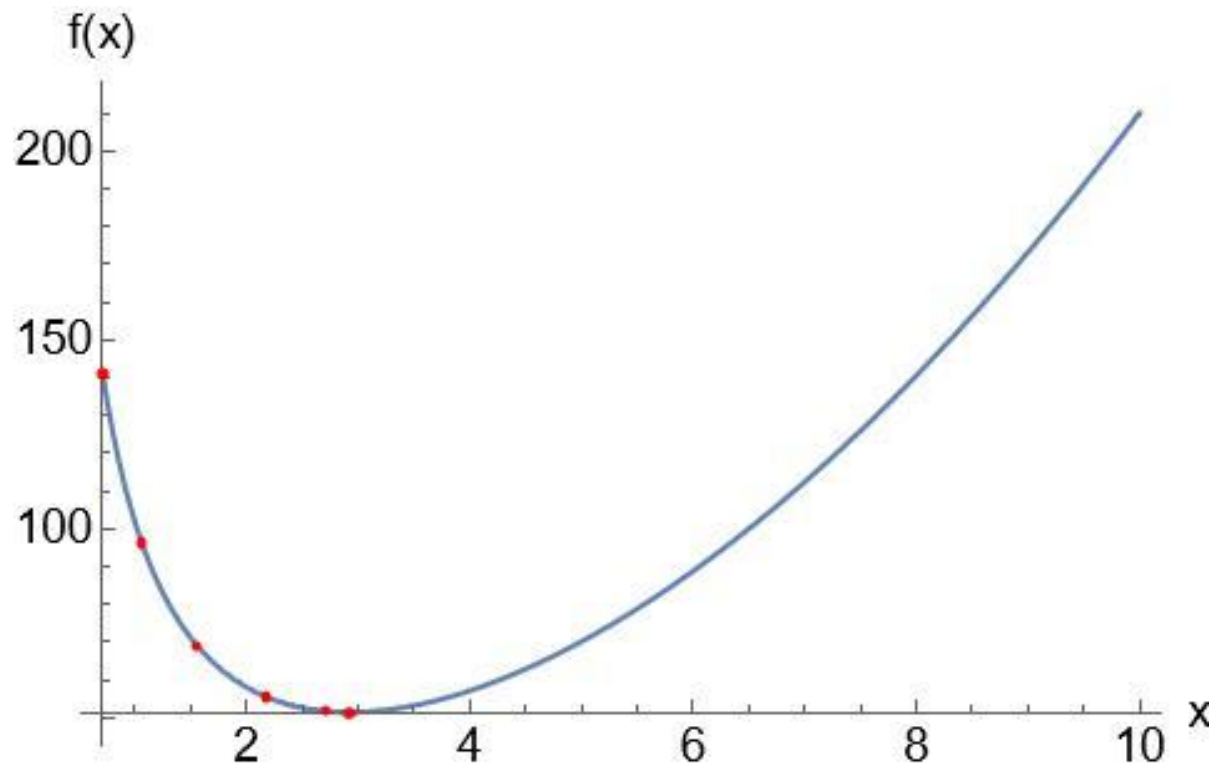

**Figure 5.28.** The results of 8 iterations of the Newton-Raphson method for $f(x) = 2x^2 + \frac{100}{x}$.

From Table 5.9, it is obtained that at iteration 8, the optimal condition is attained. Consequently, $x^* = 2.92402$.

**Table 5.9.** The results were produced by Mathematica code 5.7.

| No. it. | $x_1$ | $x_2$ | $f(x_1)$ | $f(x_2)$ | $f'(x_1)$ | $f'(x_2)$ | $\|x_2 - x_1\|$ |
|---|---|---|---|---|---|---|---|
| 1 | 10 | 0.714286 | 210 | 141.0204 | 39 | -193.143 | 9.285714 |
| 2 | 0.714286 | 1.063676 | 141.0204 | 96.27642 | -193.143 | -84.1309 | 0.34939 |
| 3 | 1.063676 | 1.558014 | 96.27642 | 69.03909 | -84.1309 | -34.9642 | 0.494338 |
| 4 | 1.558014 | 2.172682 | 69.03909 | 55.46716 | -34.9642 | -12.4933 | 0.614668 |
| 5 | 2.172682 | 2.704303 | 55.46716 | 51.60461 | -12.4933 | -2.85659 | 0.531621 |
| 6 | 2.704303 | 2.906717 | 51.60461 | 51.30108 | -2.85659 | -0.20885 | 0.202414 |
| 7 | 2.906717 | 2.923915 | 51.30108 | 51.29928 | -0.20885 | -0.00123 | 0.017198 |
| 8 | 2.923915 | 2.924018 | 51.29928 | 51.29928 | -0.00123 | -4.3E-08 | 0.000103 |

***Mathematica Code 5.7*** Newton Raphson Method

```
(* Newton Raphson Method *)

(*
Notations
x1         :Initial guess of the min
f          :Objective function of the varible x
epsilonf   :Small number to check whether the loop would be executed or not
epsilonx   :Small number to check whether the loop would be executed or not
lii        :The last iteration index
result[k]  :The results of iteration k
*)

(* Taking Input from User *)
x1=Input["Enter the initial guess:"];

epsilonf=Input["Please enter a small number epsilonf: function tolerance"];
(* A very small number to check whether the loop would be executed or not *)
epsilonx=Input["Please enter a small number epsilonx: x variable tolerance"];
(* A very small number to check whether the loop would be executed or not *)
```

```
If[
  epsilonf<=0||epsilonx<=0,
  Beep[];
  MessageDialog["Tolerance value has to be small postive number: "];
  Exit[];
  ];

(* Taking the Function from User *)
f[x_] = Evaluate[Input["Please input a function of x to find the minimum "]]; (* The user
types in, for instance x^2 *)

(* Initiating Required Variables *)
a0=x1;

(* Starting the Algorithm *)
Do[
  dfx1=f'[x1];

  If[
   N[dfx1]==0,
   Break[]
   ];

  ddfx1=f''[x1];

  If[
   ddfx1==0,
   Print["The Newton-Raphson scheme requires that the function f be twice differentiable and
f''[x]!=0"];(*Ending Program*)
   Exit[];
   ];
  x2=x1-dfx1/ddfx1;

  lii=k;
  result[k]=N[{k,x1,x2,f[x1],f[x2],f'[x1],f'[x2],Abs[x2-x1]}];

  If[
   Abs[x2-x1]<epsilonx&&Abs[f'[x2]]<epsilonf,
   Break[],
   x1=x2;
   ];,
  {k,1,∞}
  ];

(* Final Result *)
Print["The solution is x= ",  N[x2],"\nThe solution is (approximately)= ", N[f[x2]]];

(* Results of Each Iteration *)
table=TableForm[
  Table[
   result[i], {i,1,lii}
   ],
  TableHeadings->{None,{"No. of iters.","x1","x2","f[x1]","f[x2]","f'[x1]","f'[x2]","Abs[x2-
x1]"}}
  ]

Export["example58.xls",table,"XLS"];

(* Data Visualization *)
llimt=Min[Table[result[i][[3]],{i,1,lii}]];
```

```
If[
  a0<llimt,
  {a0,llimt}={llimt,a0}
  ];

Plot[
 f[x],
  {x, a0,llimt},
 AxesLabel->{"x","f(x)"},
 LabelStyle->Directive[Black,14],
 Epilog->{PointSize[0.01],Red,Point[Table[{result[i][[3]],result[i][[5]]},{i,1,lii}]]}
 ]

(* Data Manipulation *)
Manipulate[
 Plot[
  f[x],
  {x,a0,llimt},
  AxesLabel->{"x","f(x)"},
  LabelStyle->Directive[Black,14],
  Epilog->{PointSize[0.02],Red,Point[{result[i][[3]],result[i][[5]]}]}
  ],
 {i,1,lii,1}
 ]
```

### 5.6.2 Bisection Method (Bolzano search)

Again, we consider finding the minimizer of an objective function $f: \mathbb{R} \rightarrow \mathbb{R}$ over an interval $[a_0, b_0]\cdot$ As before, we assume that the objective function $f$ is unimodal. Further, suppose that $f$ is continuously differentiated and that we can use the sign of the derivative $f'$ as a basis for reducing the uncertainty interval.

The Newton-Raphson method involves the computation of the second derivative, a numerical computation of which requires three function evaluations. In the bisection method [1-10], the computation of the second derivative is avoided; instead, only the first derivative is used. The method is similar to the region elimination method, but here derivatives are used to make the decision about the region to be eliminated. If both the function value and the derivative of the function are available, then an efficient region elimination search can be conducted using just a single point rather than a pair of points to identify a point at which $f'(x) = 0$.

To begin, let $x^{(0)} = (a_0 + b_0)/2$ is the midpoint of the initial uncertainty interval. Next, evaluate $f'\left(x^{(0)}\right)$.

(1) If $f'\left(x^{(0)}\right) > 0$, then we deduce that the minimizer lies to the left of $x^{(0)}$. In other words, we reduce the uncertainty interval to $\left[a_0, x^{(0)}\right]$.
(2) On the other hand, if $f'\left(x^{(0)}\right) < 0$, then we deduce that the minimizer lies to the right of $x^{(0)}$. In this case, we reduce the uncertainty interval to $\left[x^{(0)}, b_0\right]\cdot$
(3) Finally, if $f'\left(x^{(0)}\right) = 0$, then we declare $x^{(0)}$ to be the minimizer and terminate our search.

With the new uncertainty interval computed, we repeat the process iteratively. At each iteration $k$, we compute the midpoint of the uncertainty interval. Call this point $x^{(k)}$. Depending on the sign of $f'\left(x^{(k)}\right)$ (assuming that it is nonzero), we reduce the uncertainty interval to the left or right of $x^{(k)}$. If, at any iteration $k$ we find that $f'\left(x^{(k)}\right) = 0$, then we declare $x^{(k)}$ to be the minimizer and terminate our search.

Two salient features distinguish the bisection method from the golden section.

- First, instead of using values of $f$, the bisection method uses the sign of $f'$.
- Second, at each iteration, the length of the uncertainty interval is reduced by a factor of 1/2. Hence, after $N$ steps, the range is reduced by a factor of $(1/2)^N$. This factor is smaller than in the golden section.

**Remark:**

Optimization means finding a maximum or minimum. In mathematical terms, optimization means finding where the derivative is zero. One can compare the goal of optimization in the equation $\frac{df(x)}{dx} = 0$ with the goal of root-finding in equation $f(x) = 0$. The two equations are essentially the same, both setting a function equal to zero. This motivates the idea that we can use our existing single-equation root-finding tools to optimize nonlinear equations. If the function is simple, we can perform the differentiation by hand and obtain the functional form of the derivative. At that point, we can apply any single-equation root-finding technique, such as the bisection method or the Newton-Raphson method, without modification, using as an input $f'(x)$ rather than $f(x)$. If we either cannot or will not differentiate the function analytically, we can still use the framework of the bisection method or the Newton-Raphson method, where we use the finite difference formula to provide the first derivative of the function. If we are using a technique like the bisection method, that is all we require as input. If we are using a technique like the Newton-Raphson method, which requires derivatives of the function, then we shall require the second derivative as well. The finite difference formulae for the second derivative are simple; see Chapter 2.

**Algorithm**

**Step 1:** Choose two points $a$ and $b$ such that $f'(a) < 0$ and $f'(b) > 0$. Also, choose a small number $\epsilon$. Set $x_1 = a$ and $x_2 = b$.

**Step 2:** Calculate $z = (x_2 + x_1)/2$ and evaluate $f'(z)$.

**Step 3:** $\epsilon f = |f'(z)|$, If $|f'(z)| \le \epsilon$, Terminate;
Else if $f'(z) < 0$ set $x_1 = z$ and go to Step 2;
Else if $f'(z) > 0$ set $x_2 = z$ and go to Step 2.

The only change in the bisection code necessary to convert it from a root-finding routine to an optimization routine is that, where we previously evaluated the function at the brackets, we now evaluate the first derivative of the function at the brackets using a finite difference formula. The only changes in the Newton-Raphson with numerical derivatives method to convert it from a root-finding routine to an optimization routine are that,

(i) where we previously evaluated the function, we now evaluate the first derivative of the function using a finite difference formula and

(ii) where we previously evaluated the first derivative of the function, we now evaluate the second derivative of the function using a finite difference formula.

Note that the search logic of this region elimination method is based purely on the sign of the derivative and does not use its magnitude.

***Example 5.9***

Use the bisection method to find the minimum of

$$f(x) = 0.5(x-2)^4 + 2e^x,$$

on the interval $x \in [-1,5]$, $\epsilon f = 0.01$.

***Solution***

A plot of the function is shown in Figure 5.29. The plot shows that the minimum lies at $x^* = 0.7256$, $f(x^*) = 5.4508$.

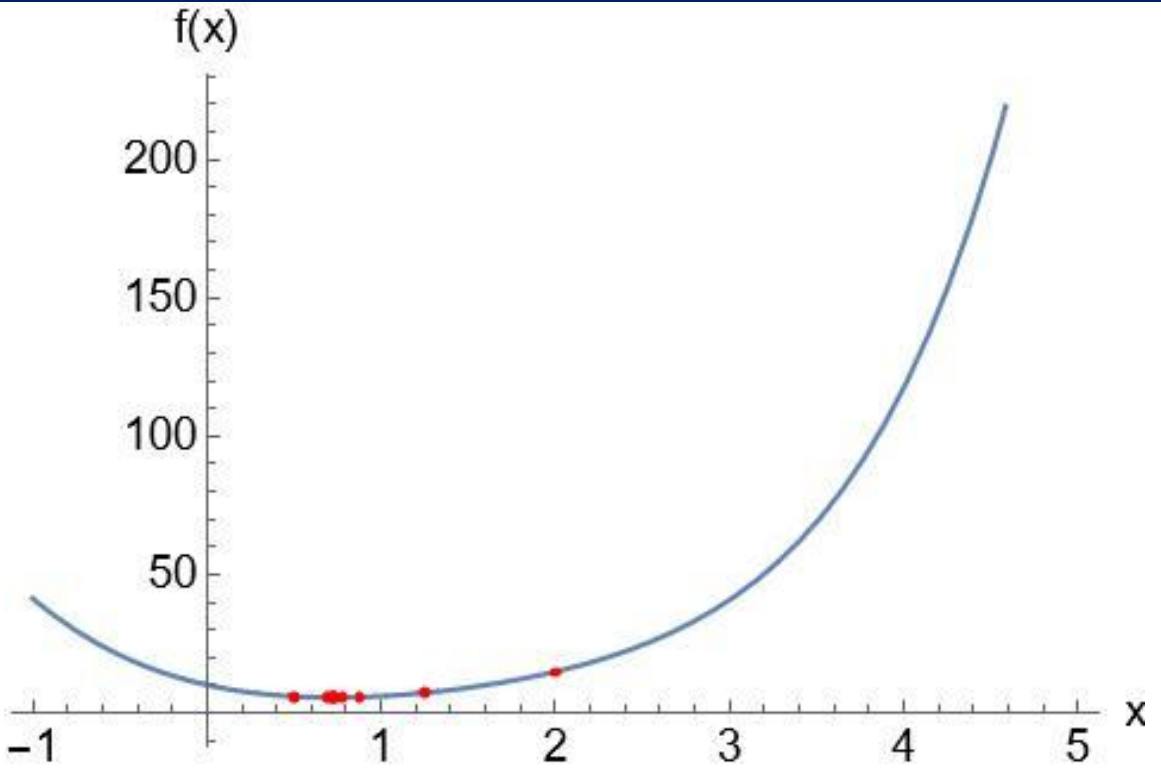

**Figure 5.29.** The results of 11 iterations of the bisection method for $f(x) = 0.5(x-2)^4 + 2e^x$.

You may find it useful to follow the data in Table 5.10.

**Table 5.10.** The results were produced by Mathematica code 5.8.

| No. it. | $x_1$ | $x_3$ | $x_2$ | $\lvert x_2 - x_1 \rvert$ | $f'(x_1)$ | $f'(x_3)$ | $f'(x_2)$ | $f(x_3)$ |
|---|---|---|---|---|---|---|---|---|
| 1 | -1 | 2 | 5 | 6 | -53.2642 | 14.77811 | 350.8263 | 14.77811 |
| 2 | -1 | 0.5 | 2 | 3 | -53.2642 | -3.45256 | 14.77811 | 5.828693 |
| 3 | 0.5 | 1.25 | 2 | 1.5 | -3.45256 | 6.136936 | 14.77811 | 7.138889 |
| 4 | 0.5 | 0.875 | 1.25 | 0.75 | -3.45256 | 1.950094 | 6.136936 | 5.598654 |
| 5 | 0.5 | 0.6875 | 0.875 | 0.375 | -3.45256 | -0.5445 | 1.950094 | 5.461247 |
| 6 | 0.6875 | 0.78125 | 0.875 | 0.1875 | -0.5445 | 0.747857 | 1.950094 | 5.471536 |
| 7 | 0.6875 | 0.734375 | 0.78125 | 0.09375 | -0.5445 | 0.113785 | 0.747857 | 5.45125 |
| 8 | 0.6875 | 0.710938 | 0.734375 | 0.046875 | -0.5445 | -0.21223 | 0.113785 | 5.452392 |
| 9 | 0.710938 | 0.722656 | 0.734375 | 0.023438 | -0.21223 | -0.04845 | 0.113785 | 5.450866 |
| 10 | 0.722656 | 0.728516 | 0.734375 | 0.011719 | -0.04845 | 0.032858 | 0.113785 | 5.45082 |
| 11 | 0.722656 | 0.725586 | 0.728516 | 0.005859 | -0.04845 | -0.00775 | 0.032858 | 5.450784 |

***Mathematica Code 5.8*** Bisection Method

```
(* Bisection Method *)

(*
Notations
a,b        :Two initial guesses
f          :Objective function of the varible x
epsilonf   :Small number to check whether the loop would be executed or not
lii        :The last iteration index
result[k]  :The results of iteration k
*)

(* Taking Input from User *)
a = Input["Please input the lower initial guess"];
b = Input["Please input the upper initial guess"];

If[
  a>b,
  Beep[];
  MessageDialog["The value of lower guess a must be less than b"];
  Exit[];
  ];

epsilonf=Input["Please enter a small number epsilonf: function tolerance"];
```

```
(* A very small number to check whether the loop would be executed or not*)

If[
  epsilonf<=0,
  Beep[];
  MessageDialog["Tolerance value has to be small postive number: "];
  Exit[];
  ];

(* Taking the function as input from user *)
f[x_] = Evaluate[Input["Please input a function of x to find the minimum "]];
(* The user types in, for instance x^2 *)

Which[
  f'[a]==0,
  Print["The solution is x= ",  N[a],"\nThe solution is (approximately)= ", N[f[a]]];
  Exit[];,
  f'[b]==0,
  Print["The solution is x= ",  N[b],"\nThe solution is (approximately)= ", N[f[b]]];
  Exit[];
  ];

Do[
  If[
   f'[a]*f'[b]<0 ,
   Break[];
   ];

  c=(a+b)/2;(* Center of the interval [a,b] *)
  w=(b-a)/2;(* Half width of the interval [a,b] *)
  nw=w*2;(* New half width of the interval *)
  a=c-nw;(* New a of the interval *)
  b=c+nw;(* New b of the interval *)

  If[
   u>10 ,
   Beep[];
   MessageDialog["The values of f'[a] and f'[b] must be satisfy f'[a]*f'[b]<0, \nAfter 10
iterations, we can not find values of a and b satisfy f'[a]*f'[b]<0"];
   Exit[];
   ];
  ,{u,1,∞}
  ];

(* Initiating Required Variables *)
a0=a;
b0=b;

x1=a;
x2=b;

(* Starting the Algorithm *)
Do[
  x3=(x1+x2)/2;
  dfx1=f'[x1];
  dfx2=f'[x2];
  dfx3=f'[x3];

  lii=k;
  result[k]=N[{k,x1,x3,x2,Abs[x2-x1],dfx1,dfx3,dfx2,f[x3]}];
```

```
   Which[
    Abs[dfx3]<epsilonf,
    Break[];,
    dfx3<0,
    x1=x3;,
    dfx3>0,
    x2=x3;
    ],
   {k,1,∞}
   ];

(* Final Result *)
Print["The solution is x= ",  N[x3],"\nThe solution is (approximately)= ", N[f[x3]]];

(* Results of Each Iteration *)
table=TableForm[
   Table[
    result[i], {i,1,lii}
    ],
   TableHeadings->{None,{"No. of iters.","x1","x3","x2","Abs[x2-
x1]","f'[x1]","f'[x3]","f'[x2]","f[x3]"}}
   ]

Export["example59.xls",table,"XLS"];

(* Data Visualization *)
Plot[
 f[x],
  {x, a0,b0},
 AxesLabel->{"x","f(x)"},
 LabelStyle->Directive[Black,14],
 Epilog->{PointSize[0.01],Red,Point[Table[{result[i][[3]],result[i][[9]]},{i,1,lii}]]}
 ]

(* Data Manipulation *)
Manipulate[
 Plot[
  f[x],
  {x,a0,b0},
  AxesLabel->{"x","f(x)"},
  LabelStyle->Directive[Black,14],
  Epilog->{PointSize[0.02],Red,Point[{result[i][[3]],result[i][[9]]}]}
  ],
 {i,1,lii,1}
 ]
```

### 5.6.3 Secant Method

This is also a gradient-based method, which follows the same condition $f'(a)f'(b) < 0$ defined in the bisection method. The secant method [1-10] combines Newton's method with a region elimination scheme for finding a root of the equation $f'(x) = 0$ in the interval $(a, b)$. In this method, both the magnitude and sign of derivatives are used to create a new point. The derivative of the function is assumed to vary linearly between the two chosen boundary points. Since boundary points have derivatives with opposite signs and the derivatives vary linearly between the boundary points, there exists a point between these two points with a zero derivative. Knowing the derivatives at the boundary points, the point with zero derivative can be easily found.

Newton's method for minimizing $f$ uses second derivatives of $f$. If the second derivative is not available, we may attempt to approximate it using the first derivative information. In particular, we may approximate $f''\left(x^{(k)}\right)$ with

$$f''\left(x^{(k)}\right) = \frac{f'\left(x^{(k)}\right) - f'\left(x^{(k-1)}\right)}{x^{(k)} - x^{(k-1)}}. \tag{5.46}$$

Using the foregoing approximation of the second derivative, we obtain the algorithm

$$x^{(k+1)} = x^{(k)} - \frac{x^{(k)} - x^{(k-1)}}{f'(x^{(k)}) - f'(x^{(k-1)})} f'\left(x^{(k)}\right), \tag{5.47}$$

called the secant method. Note that the algorithm requires two initial points to start it, which we denote $x^{(-1)}$ and $x^{(0)}$. The secant algorithm can be represented in the following equivalent form:

$$x^{(k+1)} = \frac{f'\left(x^{(k)}\right)x^{(k-1)} - f'\left(x^{(k-1)}\right)x^{(k)}}{f'(x^{(k)}) - f'(x^{(k-1)})}. \tag{5.48}$$

Observe that, like Newton's method, the secant method does not directly involve values of $f\left(x^{(k)}\right)$. Instead, it tries to drive the derivative $f'$ to zero. In fact, as we did for Newton's method, we can interpret the secant method as an algorithm for solving equations of the form $g(x) = 0$. Specifically, the secant algorithm for finding a root of the equation $g(x) = 0$ takes the form

$$x^{(k+1)} = x^{(k)} - \frac{x^{(k)} - x^{(k-1)}}{g(x^{(k)}) - g(x^{(k-1)})} g\left(x^{(k)}\right), \tag{5.49}$$

or, equivalently,

$$x^{(k+1)} = \frac{g\left(x^{(k)}\right)x^{(k-1)} - g\left(x^{(k-1)}\right)x^{(k)}}{g(x^{(k)}) - g(x^{(k-1)})}. \tag{5.50}$$

The secant method for root finding is illustrated in Figure 5.30. Unlike Newton's method, which uses the slope of $g$ to determine the next point, the secant method uses the "secant" between the $(k-1)$th and $k$th points to determine the $(k+1)$th point.

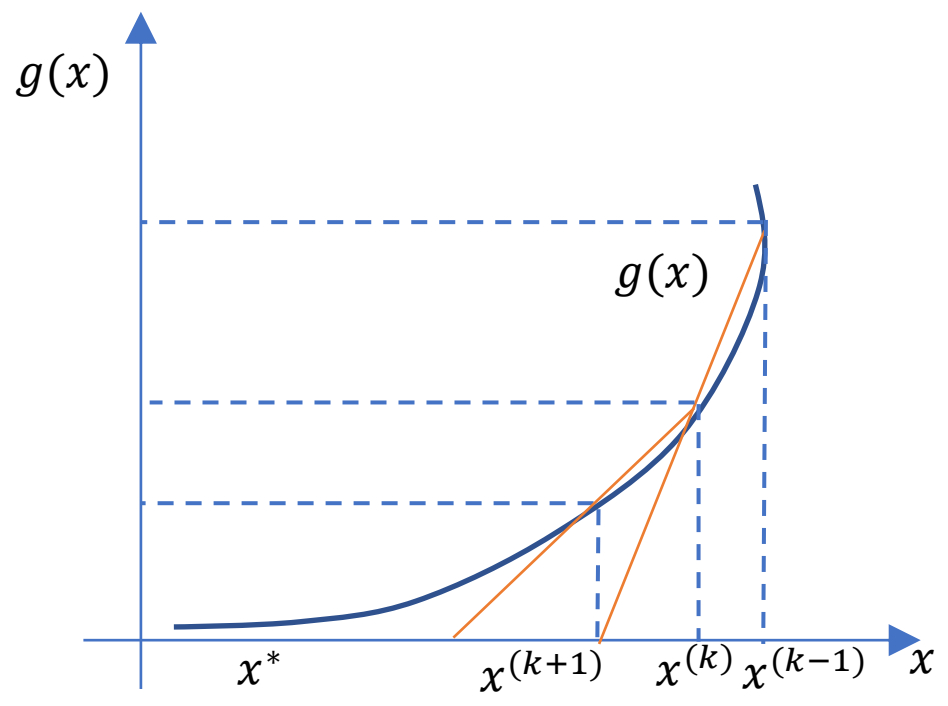

**Figure 5.30.** Secant method for root finding.

**Algorithm**

| | |
|---|---|
| **Step 1:** | Choose two points $a$ and $b$ such that $f'(a) < 0$ and $f'(b) > 0$. Also, choose a small number $\epsilon$. Set $x_1 = a$ and $x_2 = b$. |
| **Step 2:** | Calculate the new point z using (5.47) and evaluate $f'(z)$. |
| **Step 3:** | If $\lvert f'(z)\rvert \le \epsilon$, Terminate;<br>Else if $f'(z) < 0$ set $x_1 = z$ and go to Step 2;<br>Else if $f'(z) > 0$ set $x_2 = z$ and go to Step 2. |

***Example 5.10***

Use the bisection method to find the minimum of $f(x) = 0.5(x-2)^4 + 2e^x$ on the interval $x \in [-1,2]$, $\epsilon f = 0.01$.

***Solution***

A plot of the function is shown in Figure 5.31. The plot shows that the minimum lies at $x^* = 0.7266$, $f(x^*) = 5.4508$.

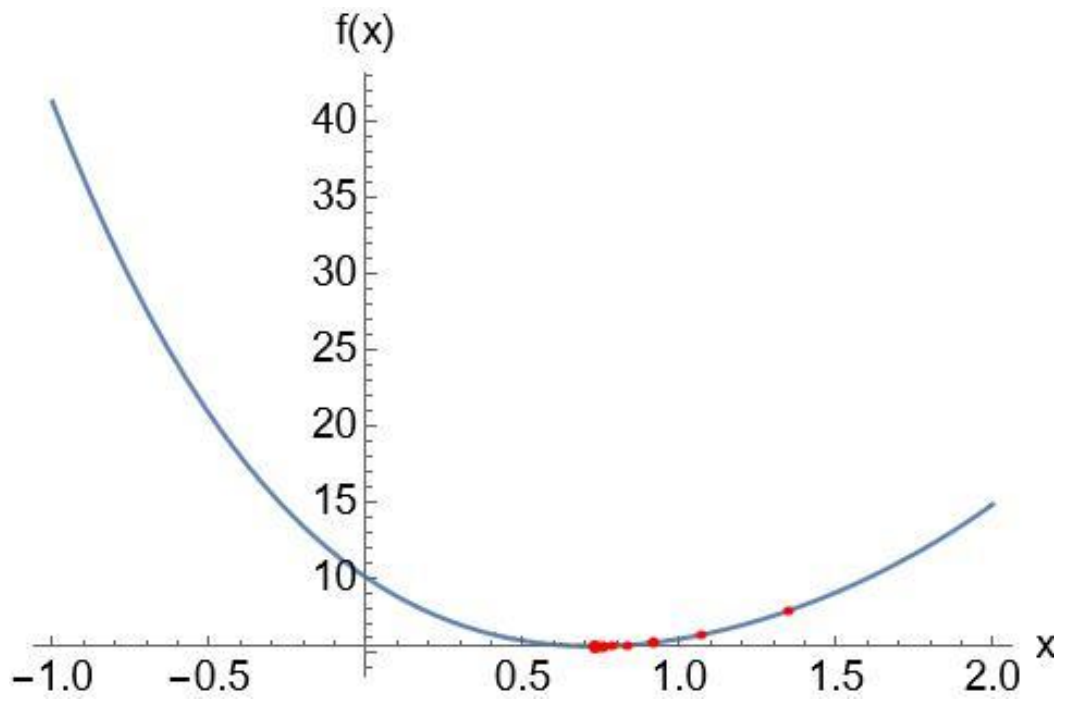

**Figure 5.31.** The results of 13 iterations of secant method for $f(x) = 0.5(x-2)^4 + 2e^x$.

You may find it useful to follow the data in Table 5.11.

**Table 5.11.** The results were produced by Mathematica code 5.9.

| No. it. | $x_1$ | $x_3$ | $x_2$ | $\lvert x_2 - x_1 \rvert$ | $f'(x_1)$ | $f'(x_3)$ | $f'(x_2)$ | $f(x_3)$ |
|---|---|---|---|---|---|---|---|---|
| 1 | -1 | 1.34843 | 2 | 3 | -53.2642 | 7.149512 | 14.77811 | 7.792869 |
| 2 | -1 | 1.070511 | 1.34843 | 2.34843 | -53.2642 | 4.227679 | 7.149512 | 6.206945 |
| 3 | -1 | 0.918256 | 1.070511 | 2.070511 | -53.2642 | 2.478185 | 4.227679 | 5.694485 |
| 4 | -1 | 0.832974 | 0.918256 | 1.918256 | -53.2642 | 1.421443 | 2.478185 | 5.527753 |
| 5 | -1 | 0.78533 | 0.832974 | 1.832974 | -53.2642 | 0.801957 | 1.421443 | 5.474698 |
| 6 | -1 | 0.758848 | 0.78533 | 1.78533 | -53.2642 | 0.447749 | 0.801957 | 5.458135 |
| 7 | -1 | 0.744186 | 0.758848 | 1.758848 | -53.2642 | 0.248452 | 0.447749 | 5.453028 |
| 8 | -1 | 0.736088 | 0.744186 | 1.744186 | -53.2642 | 0.137379 | 0.248452 | 5.451465 |
| 9 | -1 | 0.731622 | 0.736088 | 1.736088 | -53.2642 | 0.075813 | 0.137379 | 5.450989 |
| 10 | -1 | 0.729161 | 0.731622 | 1.731622 | -53.2642 | 0.041792 | 0.075813 | 5.450844 |
| 11 | -1 | 0.727805 | 0.729161 | 1.729161 | -53.2642 | 0.023023 | 0.041792 | 5.450801 |
| 12 | -1 | 0.727059 | 0.727805 | 1.727805 | -53.2642 | 0.01268 | 0.023023 | 5.450787 |
| 13 | -1 | 0.726648 | 0.727059 | 1.727059 | -53.2642 | 0.006982 | 0.01268 | 5.450783 |

***Mathematica Code 5.9*** Secant Method

```
(* Secant Method *)

(*
Notations
a,b         :Two initial guesses
f           :Objective function of the varible x
epsilonf    :Small number to check whether the loop would be executed or not
lii         :The last iteration index
result[k]   :The results of iteration k
*)

(* Taking Input from User *)
a = Input["Please input the lower initial guess"];
b = Input["Please input the upper initial guess"];

If[
  a>b,
```

```
  Beep[];
  MessageDialog["The value of lower guess a must be less than b"];
  Exit[];
  ];

epsilonf=Input["Please enter a small number epsilonf: function tolerance"];
(* A very small number to check whether the loop would be executed or not*)

If[
  epsilonf<=0,
  Beep[];
  MessageDialog["Tolerance value has to be small postive number: "];
  Exit[];
  ];

(* Taking the function as input from user *)
f[x_] = Evaluate[Input["Please input a function of x to find the minimum "]];
(* The user types in, for instance x^2 *)

Which[
  f'[a]==0,
  Print["The solution is x= ",  N[a],"\nThe solution is (approximately)= ", N[f[a]]];
  Exit[];,
  f'[b]==0,
  Print["The solution is x= ",  N[b],"\nThe solution is (approximately)= ", N[f[b]]];
  Exit[];
  ];

Do[
  If[
   f'[a]*f'[b]<0 ,
   Break[];
   ];

  c=(a+b)/2;(* Center of the interval [a,b] *)
  w=(b-a)/2;(* Half width of the interval [a,b] *)
  nw=w*2;(* New half width of the interval *)
  a=c-nw;(* New a of the interval *)
  b=c+nw;(* New b of the interval *)

  If[
   u>10,
   Beep[];
   MessageDialog["The values of f'[a] and f'[b] must be satisfy f'[a]*f'[b]<0, \nAfter 10
iterations, we can not find values of a and b satisfy f'[a]*f'[b]<0"];
   Exit[];
   ];
  ,{u,1,∞}
  ];

(* Initiating Required Variables *)
a0=a;
b0=b;

x1=a;
x2=b;

(* Starting the Algorithm *)
Do[
  If[
   f'[x2]-f'[x1]==0,
```

```
   Beep[];
   MessageDialog["The derivative of the function is assumed to vary linearly between the two chosen boundary points. \nwe have (f'[x2]=f'[x1]). "];
   Exit[];
   ];
  x3=x2-((x2-x1)f'[x2])/(f'[x2]-f'[x1]);

  dfx1=f'[x1];
  dfx2=f'[x2];
  dfx3=f'[x3];

  lii=k;
  result[k]=N[{k,x1,x3,x2,Abs[x2-x1],dfx1,dfx3,dfx2,f[x3]}];

  Which[
   Abs[dfx3]<epsilonf,
   Break[];,
   dfx3<0,
   x1=x3;,
   dfx3>0,
   x2=x3;
   ],
  {k,1,∞}
  ];

(* Final Result *)
Print["The solution is x= ",  N[x3],"\nThe solution is (approximately)= ", N[f[x3]]];

(* Results of Each Iteration *)
table=TableForm[
  Table[
   result[i], {i,1,lii}
   ],
  TableHeadings->{None,{"No. of iters.","x1","x3","x2","Abs[x2-x1]","f'[x1]","f'[x3]","f'[x2]","f[x3]"}}
  ]

Export["example510.xls",table,"XLS"];

(* Data Visualization *)
Plot[
 f[x],
  {x, a0,b0},
 AxesLabel->{"x","f(x)"},
 LabelStyle->Directive[Black,14],
 Epilog->{PointSize[0.01],Red,Point[Table[{result[i][[3]],result[i][[9]]},{i,1,lii}]]}
 ]
(* Data Manipulation *)
Manipulate[
 Plot[
  f[x],
  {x,a0,b0},
  AxesLabel->{"x","f(x)"},
  LabelStyle->Directive[Black,14],
  Epilog->{PointSize[0.02],Red,Point[{result[i][[3]],result[i][[9]]}]}
  ],
 {i,1,lii,1}
 ]
```

## 5.7 Cubic Search Method

Another one-dimensional optimization method that is sometimes quite useful is the cubic interpolation method [7]. This is based on the third-order polynomial

$$\bar{f}(x) = a_0 + a_1 x + a_2 x^2 + a_3 x^3. \tag{5.51}$$

The basic logic is similar to the quadratic approximation scheme. However, in this instance, because both the function value and the derivative value are available at each point, the approximating polynomial can be constructed using fewer points. The coefficients $a_i$ can be determined such that $\bar{f}(x)$ and/or its derivatives at certain points are equal to $f(x)$ and/or its derivatives. Since there are four coefficients in (5.51), four equations are needed for the complete characterization of $\bar{f}(x)$. These equations can be chosen in a number of ways, and several cubic interpolation formulas can be generated.

The cubic search starts with an arbitrary point $x_1$ and finds another point $x_2$ by a bounding search such that the derivatives $f'(x_1)$ and $f'(x_2)$ are of opposite signs. In other words, the stationary point $\bar{x}$ where $f'(x) = 0$ is bracketed between $x_1$ and $x_2$. A cubic approximation function of the form

$$\bar{f}(x) = a_0 + a_1(x - x_1) + a_2(x - x_1)(x - x_2) + a_3(x - x_1)^2(x - x_2), \tag{5.52}$$

is fitted such that (5.52) agrees with $f(x)$ at the two points $x_1$ and $x_2$. The first derivative of $\bar{f}(x)$ is given by

$$\frac{d}{dx}\bar{f}(x) = a_1 + a_2(x - x_1) + a_2(x - x_2) + a_3(x - x_1)^2 + 2a_3(x - x_1)(x - x_2) \tag{5.53}$$

The coefficients $a_0$, $a_1$, $a_2$, and $a_3$ of (5.52) can now be determined using the values of $f(x_1)$, $f(x_2)$, $f'(x_1)$, and $f'(x_2)$ by solving the following linear equations:

$$\begin{aligned} f_1 &= f(x_1) = a_0, \\ f_2 &= f(x_2) = a_0 + a_1(x_2 - x_1), \\ f_1' &= f'(x_1) = a_1 + a_2(x_1 - x_2), \\ f_2' &= f'(x_2) = a_1 + a_2(x_2 - x_1) + a_3(x_2 - x_1)^2. \end{aligned} \tag{5.54}$$

As in the quadratic case discussed earlier, given these coefficients, an estimate of the stationary point of $f$ can be obtained from the approximating cubic of (5.52). In this case, when we set the derivative of $\bar{f}(x)$ given by (5.53) to zero, we get a quadratic equation. By applying the formula for the root of the quadratic equation, a closed-form solution to the stationary point $\bar{x}$ of the approximating cubic is obtained as follows:

$$\bar{x} = \begin{cases} x_2 & \text{if } \mu < 0, \\ x_2 - \mu(x_2 - x_1) & \text{if } 0 \le \mu \le 1, \\ x_1 & \text{if } \mu > 1, \end{cases} \tag{5.55}$$

where

$$\mu = \frac{f_2' + w - z}{f_2' - f_1' + 2w}, \tag{5.56}$$

$$z = \frac{3(f_1 - f_2)}{x_2 - x_1} + f_1' + f_2', \tag{5.57}$$

$$w = \begin{cases} (z^2 - f_1' f_2')^{1/2} & \text{if } x_1 < x_2, \\ -(z^2 - f_1' f_2')^{1/2} & \text{if } x_1 > x_2. \end{cases} \tag{5.58}$$

Hence, similar to Powell's successive quadratic estimation method, the minimum of the approximation function $\bar{f}(x)$ can be used as an estimate of the true minimum of the objective function. This estimate and the earlier two points ($x_1$ and $x_2$) may be used to find the next estimate of the true minimum point. Two points ($x_1$ and $x_2$) are so chosen that the product of their first derivative is negative. This procedure may be continued until the desired accuracy is obtained as follows:

**Algorithm**

**Step 1:** Choose an initial point $x^{(0)}$, a step size $\Delta$, and two termination parameters $\epsilon_1$ and $\epsilon_2$. Compute $f'(x^{(0)})$. If $f'(x^{(0)}) > 0$, set $\Delta = -\Delta$. Set $k = 0$.

**Step 2:** Compute $x^{(k+1)} = x^{(k)} + 2^k\Delta$.

**Step 3:** Evaluate $f'(x^{(k+1)})$.
If $f'(x^{(k+1)})f'(x^{(k)}) \leq 0$, set $x_1 = x^{(k)}$, $x_2 = x^{(k+1)}$, and go to Step 4;
Else set $k = k + 1$ and go to Step 2.

**Step 4:** Calculate the point $\bar{x}$ using (5.55).

**Step 5:** If $f(\bar{x}) < f(x_1)$, go to Step 6;
Else set $\bar{x} = \bar{x} - \frac{1}{2}(\bar{x} - x_1)$ until $f(\bar{x}) \leq f(x_1)$ is achieved.

**Step 6:** Compute $f'(\bar{x})$. If $|f'(\bar{x})| \leq \epsilon_1$ and $|(\bar{x} - x_1)/\bar{x}| \leq \epsilon_2$, Terminate;
Else if $f'(\bar{x})f'(x_1) < 0$, set $x_2 = \bar{x}$;
Else if $f'(\bar{x})f'(x_2) < 0$ set $x_1 = \bar{x}$.
Go to Step 4.

***Example 5.11***

Use the cubic method to find the minimum of

$$f(x) = 0.5(x-2)^4 + 2e^x$$

with $x_0 = 10$, $\Delta = 0.1$, $\epsilon_f = 0.01$, and $\epsilon_x = 0.01$.

***Solution***

A plot of the function is shown in Figure 5.32. The plot shows that the minimum lies at $x^* = 0.7261$, $f(x^*) = 5.45078$.

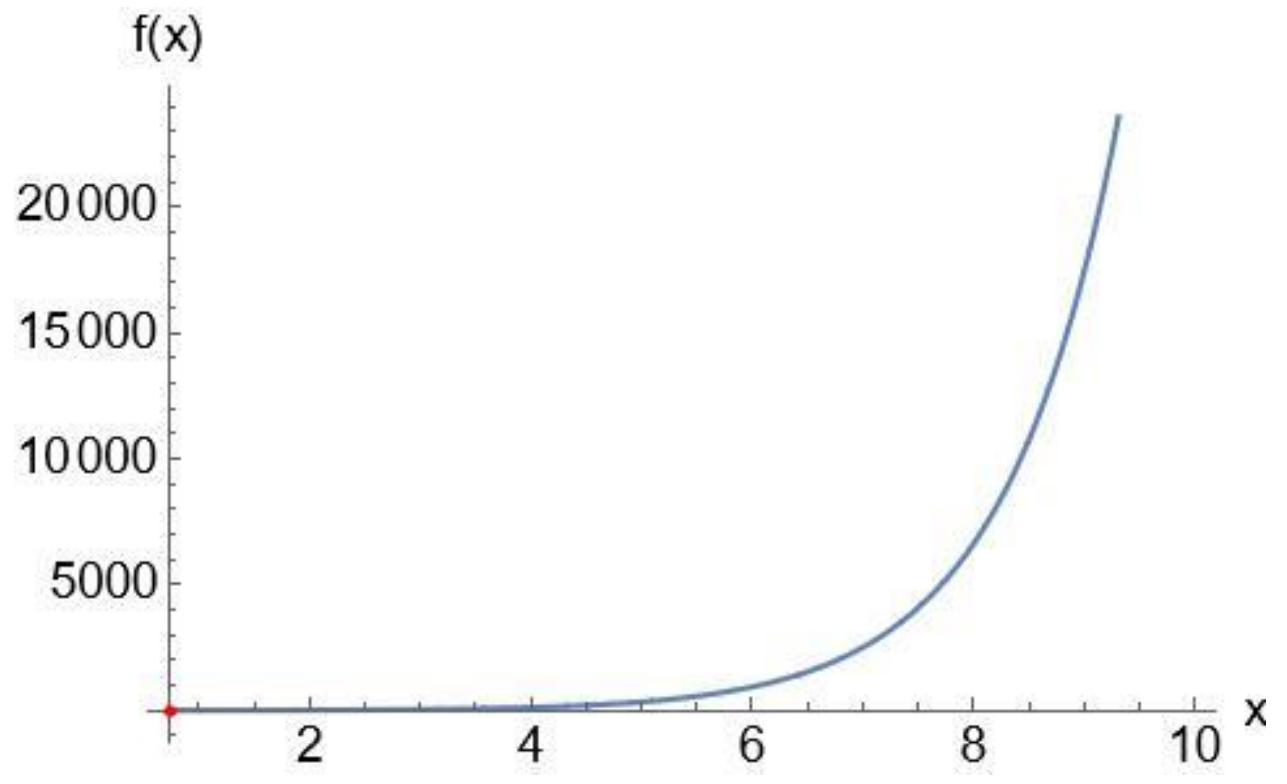

**Figure 5.32.** The results of 3 iterations of the cubic search method for $f(x) = 0.5(x-2)^4 + 2e^x$.

You may find it useful to follow the data in Table 5.12.

**Table 5.12.** The results were produced by Mathematica code 5.10.

| No. it. | $x_1$ | $x_2$ | $\bar{x}$ | $f(x_1)$ | $f(x_2)$ | $f(\bar{x})$ | $\lvert f'(\bar{x})\rvert$ | $\lvert(\bar{x}-x_1)/\bar{x}\rvert$ |
|---|---|---|---|---|---|---|---|---|
| 1 | 3.8 | -2.6 | 0.752585 | 94.65117 | 224.0213 | 5.455596 | 0.362898 | 4.049261 |
| 2 | 0.752585 | -2.6 | 0.723063 | 5.455596 | 224.0213 | 5.450847 | 1.320085 | 0.186114 |
| 3 | 0.752585 | 0.723063 | 0.726145 | 5.455596 | 5.450847 | 5.450781 | 2.57E-06 | 0.036412 |

**Mathematica Code 5.10** Cubic Search Method

```
(* Cubic Search Method *)

(*
Notations
x0         :Initial point
delta      :Postive step size
epsilonf   :Small number to check whether the loop would be executed or not
epsilonx   :Small number to check whether the loop would be executed or not
f          :Objective function
xbar[k]    :Final design solution of iteration k
f[xbar[k]]:Final objective function value of iteration k
lii        :The last iteration index
result[k] :The results of iteration k
*)

(* Taking Input from User *)
x0=Input["Enter the initial point:"];
delta=Input["Enter the postive step size:"];

epsilonf=Input["Please enter a small number epsilonf: function tolerance"];
epsilonx=Input["Please enter a small number epsilonx: x variable tolerance"];

If[
  epsilonf<=0||epsilonx<=0,
  Beep[];
  MessageDialog["Tolerance value has to be small postive number: "];
  Exit[];
  ];

(* Taking the function as input from user *)
f[x_] = Evaluate[Input["Please input a function of x to find the minimum "]];
(* The user types in, for instance x^2*)

(* Initiating Required Variables *)
dfx0=f'[x0];

Which[
  dfx0==0,
  Print["The solution is x= ",  N[x0],"\nThe solution is (approximately)= ", N[f[x0]]];
  Exit[];,
  dfx0>0,
  delta=-delta;
  ];

a=x0;

Do[
  b=a+(2^i)*delta;
  dfa=f'[a];
  dfb=f'[b];

  If[
   dfb*dfa<0,
   x1=a;
   x2=b;
   Break[];
   ];

  a=b;
```

```
  If[
   i>20 ,
   Beep[];
   MessageDialog["The values of f'[b] and f'[a] must be satisfy f'[b]f'[a]<0, \nAfter 20 iterations, we can not find values of b and a satisfy f'[b]f'[a]<0"];
   Exit[];
   ];,
  {i,1,∞}
  ];

(* Starting the Algorithm *)
Do[
  dfx2=f'[x2];
  dfx1=f'[x1];

  z=(3(f[x1]-f[x2]))/(x2-x1)+dfx1+dfx2;

  Which[
   x1<x2,
   w= Sqrt[z^2 - dfx1*dfx2],
   x1>x2,
   w=-Sqrt[z^2 - dfx1*dfx2]
   ];

  μ=(dfx2+w-z)/(dfx2-dfx1+2w);

  Which[
   μ<0,
   xbar=x2,
   0<=μ&&μ<=1,
   xbar=x2-μ*(x2-x1),
   μ>1,
   xbar=x1
   ];

  fxbar=f[xbar];
  dfxbar=f'[xbar];

  errordf=Abs[dfxbar];
  errorx=Abs[(xbar-x1)/xbar];

  Do[
   If[
    fxbar<=f[x1],
    Break[];
    ];
   xbar=xbar-(xbar-x1)/2;
   fxbar=f[xbar];,
   {i,1,∞}
   ];

  lii=k;
  result[k]=N[{k,x1,x2,xbar,f[x1],f[x2],f[xbar],errordf,errorx}];


  Which[
   errordf<=epsilonf||errorx<=epsilonx,
   Break[];,
   dfxbar*dfx1<0,
   x2=xbar;,
   dfxbar*dfx2<0,
```

```
    x1=xbar;,
    k>30,
    Print["After 30 iteration the error greater than the epsilon"];
    Exit[];
    ];,
   {k,1,∞}
   ];

(* Final Result *)
Print["The solution lies between ",N[x1]," and ",N[x2], " at ",N[xbar],"\nThe solution is
(approximately)", N[f[xbar]]];

(* Results of Each Iteration *)
table=TableForm[
  Table[
   result[i], {i,1,lii}
   ],
  TableHeadings->{None,{"No
it.","x1","x2","xbar","f[x1]","f[x2]","f[xbar]","errordf","errorx"}}
  ]

Export["example511.xls",table,"XLS"];

(* Data Visualization *)
Plot[
 f[x],
  {x,x0,xbar},
 AxesLabel->{"x","f(x)"},
 LabelStyle->Directive[Black,14],
 Epilog->{PointSize[0.01],Red,Point[Table[{result[i][[4]],result[i][[7]]},{i,1,lii}]]}
 ]

(* Data Manipulation *)
Manipulate[
 Plot[
  f[x],
  {x,x0,xbar},
  AxesLabel->{"x","f(x)"},
  LabelStyle->Directive[Black,14],
  Epilog->{PointSize[0.02],Red,Point[{result[i][[4]],result[i][[7]]}]}
  ],
 {i,1,lii,1}
 ]
```

## 5.8 Mathematica Built-in Functions

The Wolfram Language has a collection of commands that do unconstrained optimization [12,13]

| | |
|---|---|
| `FindMinimum[f,x]` | searches for a local minimum in f, starting from an automatically selected point. |
| `FindMinimum[f,{x,x0}]` | searches for a local minimum in f, starting from the point x=x0. |
| `FindMinimum[f,{{x,x0},{y,y0},…}]` | searches for a local minimum in a function of several variables. |
| `FindMinimum[{f,cons},{{x,x0},{y,y0},…}]` | searches for a local minimum subject to the constraints cons. |
| `FindMinimum[{f,cons},{x,y,…}]` | starts from a point within the region defined by the constraints. |

**Details and Options**

- `FindMinimum` returns a list of the form `{fmin ,{x->xmin }}`, where `fmin` is the minimum value of `f` found, and `xmin` is the value of `x` for which it is found.
- If the starting point for a variable is given as a list, the values of the variable are taken to be lists with the same dimensions.
- The constraints `cons` can contain equations, inequalities or logical combinations of these.
- The constraints `cons` can be any logical combination of:

| | |
|---|---|
| `lhs==rhs` | equations |
| `lhs>rhs or lhs>=rhs` | Inequalities |
| `{x,y,…}∈reg` | region specification |

- `FindMinimum[f,{x,x0,x1}]` searches for a local minimum in `f` using `x0` and `x1` as the first two values of `x`, avoiding the use of derivatives.
- `FindMinimum[f,{x,x0,xmin,xmax}]` searches for a local minimum, stopping the search if `x` ever gets outside the range `xmin` to `xmax`.
- Except when `f` and `cons` are both linear, the results found by `FindMinimum` may correspond only to local, but not global, minima.
- By default, all variables are assumed to be real.
- For linear `f` and `cons, x∈Integers` can be used to specify that a variable can take on only integer values.
- The following options can be given:

| | | |
|---|---|---|
| `AccuracyGoal` | `Automatic` | the accuracy sought |
| `EvaluationMonitor` | `None` | expression to evaluate whenever $f$ is evaluated |
| `Gradient` | `Automatic` | the list of gradient components for $f$ |
| `MaxIterations` | `Automatic` | maximum number of iterations to use |
| `Method` | `Automatic` | method to use |
| `PrecisionGoal` | `Automatic` | the precision sought |
| `StepMonitor` | `None` | expression to evaluate whenever a step is taken |
| `WorkingPrecision` | `MachinePrecision` | the precision used in internal computations |

- The settings for `AccuracyGoal` and `PrecisionGoal` specify the number of digits to seek in both the value of the position of the minimum, and the value of the function at the minimum.
- `FindMinimum` continues until either of the goals specified by `AccuracyGoal` or `PrecisionGoal` is achieved.
- Possible settings for Method include `"ConjugateGradient"`, `"PrincipalAxis"`, `"LevenbergMarquardt"`, `"Newton"`, `"QuasiNewton"`, `"InteriorPoint"`, and `"LinearProgramming"`, with the default being Automatic.
- The `FindMinimum` function in the Wolfram Language has five essentially different ways of choosing this model, controlled by the method option.

| | |
|---|---|
| `"Newton"` | use the exact Hessian or a finite difference approximation if the symbolic derivative cannot be computed |
| `"QuasiNewton"` | use the quasi-Newton BFGS approximation to the Hessian built up by updates based on past steps |
| `"LevenbergMarquardt"` | a Gauss–Newton method for least-squares problems; the Hessian is approximated by , where  is the Jacobian of the residual function |
| `"ConjugateGradient"` | a nonlinear version of the conjugate gradient method for solving linear systems; a model Hessian is never formed explicitly |
| `"PrincipalAxis"` | works without using any derivatives, not even the gradient, by keeping values from past steps; it requires two starting conditions in each variable |

- With `Method->Automatic`, the Wolfram Language uses the quasi-Newton method unless the problem is structurally a sum of squares, in which case the Levenberg–Marquardt variant of the Gauss–Newton method is used. When given two starting conditions in each variable, the principal axis method is used.

***Mathematica Examples 5.1***

Input
```
(*With different starting points,you may get different local minima*)
FindMinimum[x Cos[x],{x,2}]
Plot[
 x Cos[x],
 {x,0,20},
 LabelStyle->Directive[Black,20]
 ]
```

Output
```
{-3.28837,{x->3.42562}}
```

Output

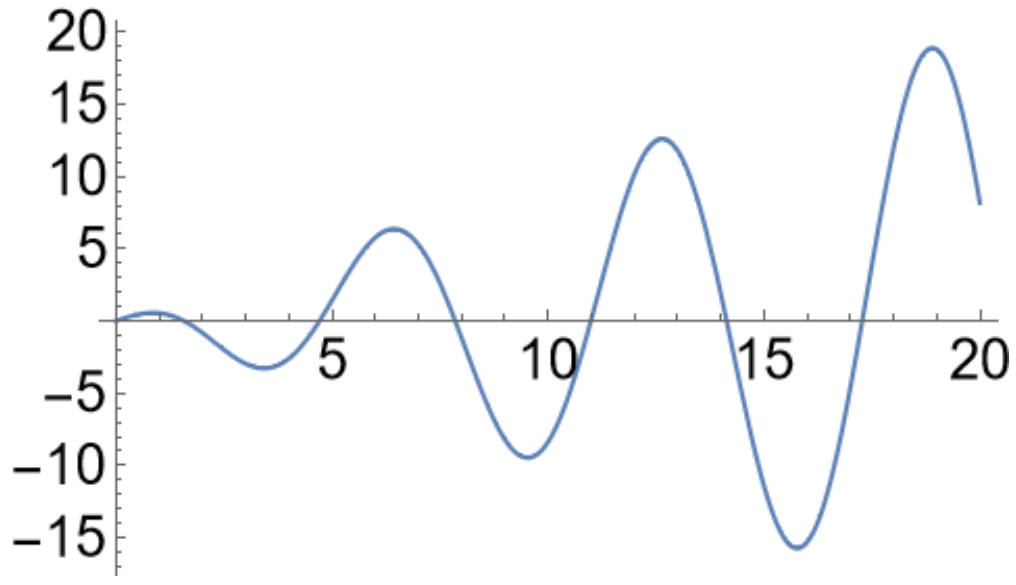

Input
```
FindMinimum[x Cos[x],{x,5}]
```
Output
```
{-3.28837,{x->3.42562}}
```

Input
```
FindMinimum[x Cos[x],{x,10}]
```
Output
```
{-9.47729,{x->9.52933}}
```

Input
```
f=2 x^2-3 x+5;
Plot[f,
 {x,-10,10},
 LabelStyle->Directive[Black,20]
 ]
FindMinimum[f,x,Method->"ConjugateGradient"]
FindMinimum[f,x,Method->"PrincipalAxis"]
FindMinimum[f,x,Method->"Newton"]
FindMinimum[f,x,Method->"QuasiNewton"]
FindMinimum[f,x,Method->"LevenbergMarquardt"]
```

Output

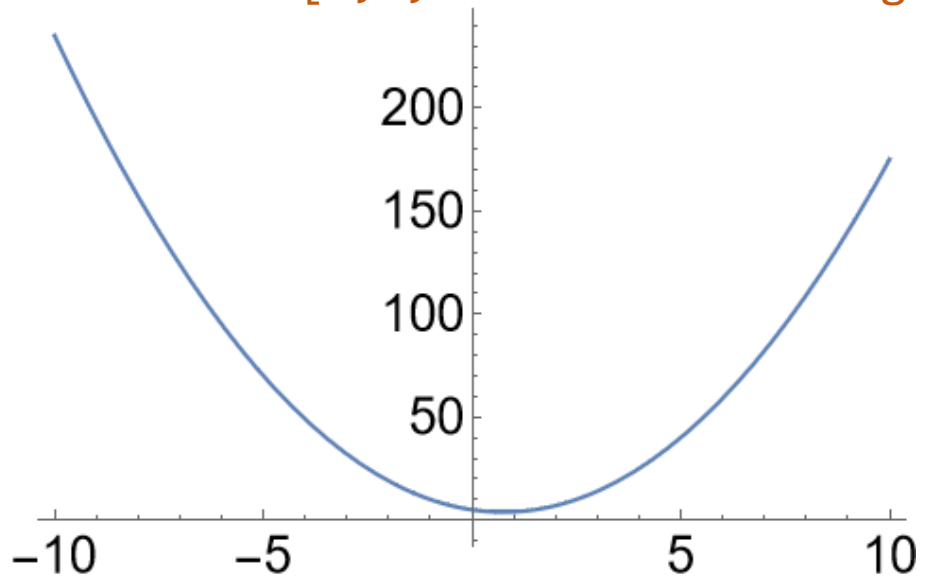

Output
```
{3.875,{x->0.75}}
```
Output
```
{3.875,{x->0.75}}
```
Output
```
{3.875,{x->0.75}}
```
Output
```
{3.875,{x->0.75}}
```
Output

FindMinimum::notlm: The objective function for the method LevenbergMarquardt must be in a least-squares form: Sum[f[i][x]^2, {i, 1, n}] or Sum[w[i] f[i][x]^2, {i, 1, n}] with positive w[i].

```
FindMinimum[f,x,Method->LevenbergMarquardt]
```

| | |
|---|---|
| `FindMinimumPlot[f,{x,xst}]` | plots the steps and the points at which the function f and any of its derivatives are evaluated in FindMinimum[f,{x,xst}], superimposed on a plot of f versus x. |
| `FindMinimumPlot[f,{{x,xst},{y,yst}}]` | plots the steps and the points at which the bivariate function f and any of its derivatives are evaluated, superimposed on a contour plot of f as a function of x and y. |
| `FindMinimumPlot[f,range,property]` | returns the specified property. |

The `FindMinimumPlot` command is defined in the ``Optimization`UnconstrainedProblems` `` package loaded automatically by this notebook. It runs `FindMinimum`, keeps track of the function and gradient evaluations and steps taken during the search (using the `EvaluationMonitor` and `StepMonitor` options), and shows them superimposed on a plot of the function.

- Steps and evaluation points are color coded for easy detection as follows:
  - Steps are shown with blue lines and blue points.
  - Function evaluations are shown with green points.
  - Gradient evaluations are shown with red points.
  - Hessian evaluations are shown with cyan points.
  - Residual function evaluations are shown with yellow points.
  - Jacobian evaluations are shown with purple points.
  - The search termination is shown with a large black point.
- `FindMinimumPlot` and `FindRootPlot` return a list containing `{result,summary,plot}`, where:
  - result is the result of `FindMinimum` or `FindRoot`.
  - summary is a list of rules showing the number of steps and evaluations of the function and its derivatives.
  - plot is the graphics object shown.
- From the plot, it is clear that `FindMinimum` has found a local minimum point.
- With the setting `PlotLegends→Automatic, FindMinimumPlot` shows a legend for the evaluation points.

***Mathematica Examples 5.2***

Input
```
(*This loads a package that contains some utility functions.*)
<<Optimization`UnconstrainedProblems`
```
Input
```
(*This shows a plot of the function x*Sin[x+1].*)
Plot[
 x*Sin[x+1],
 {x,-10,10},
 LabelStyle->Directive[Black,20]
 ]
(*This shows the steps taken by FindMinimum for the function x*Sin[x+1] starting at x=0.*)
FindMinimumPlot[
 x*Sin[x+1],
 {x,0},
 LabelStyle->Directive[Black,20]
 ]
```
Output

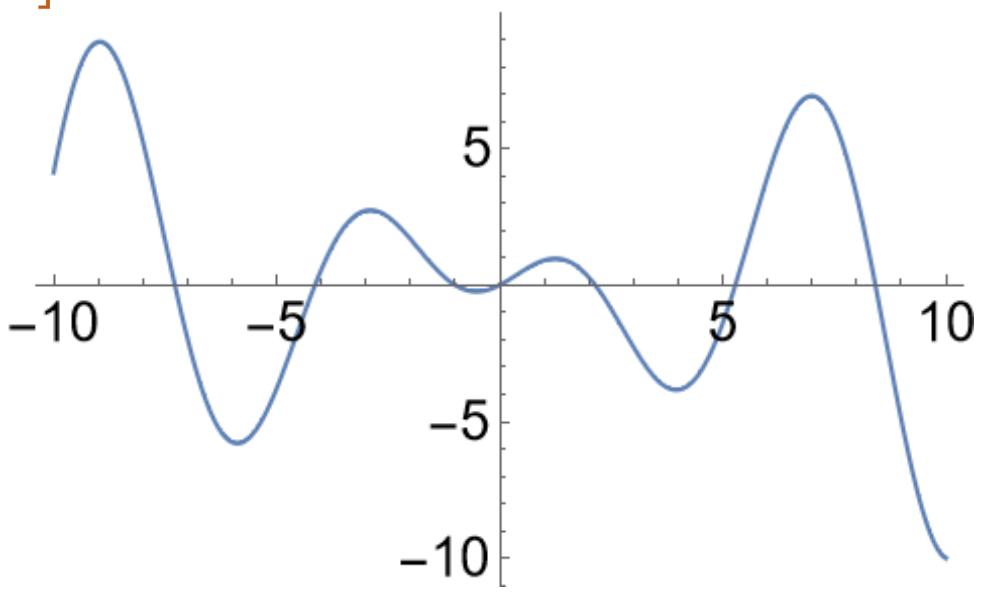

```
{{-0.240125,{x->-0.520269}},{Steps->5,Function->6,Gradient->6},
```

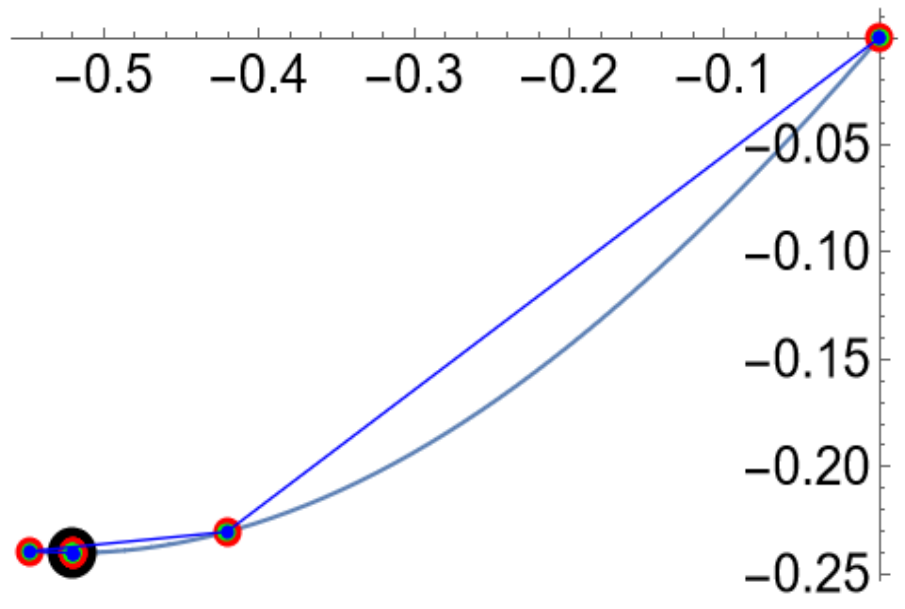

}

Input
```
(*This shows the steps taken by FindMinimum for the function x*Sin[x+1] starting at x=2.*)
FindMinimumPlot[
 x*Sin[x+1],
 {x,2},
 LabelStyle->Directive[Black,20]
 ]
```

Output
```
{{-3.83922,{x->3.95976}},{Steps->4,Function->9,Gradient->9},
```

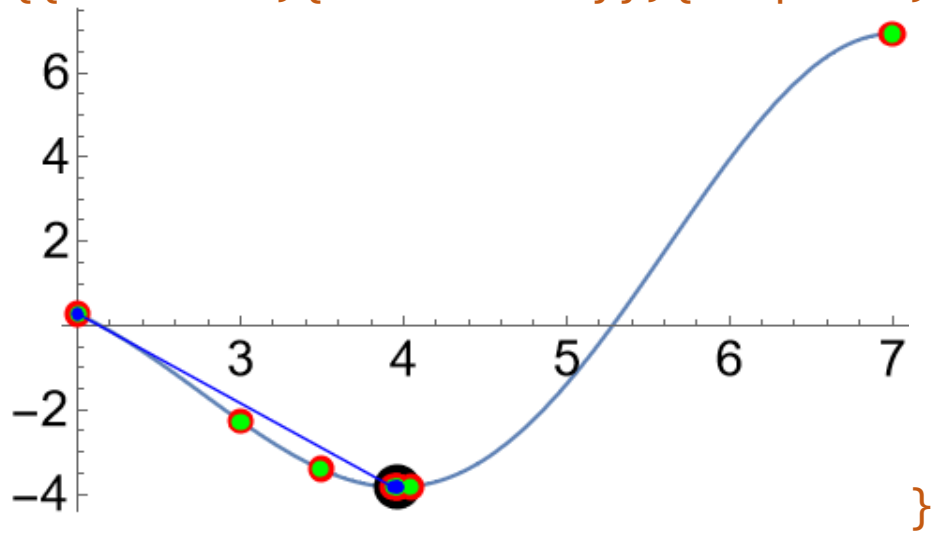

}

Input
```
(*Show the steps and function evaluations used in finding a local minimum of the function exp(x)+1/x*)
FindMinimumPlot[
 Exp[x]+1/x,
 {x,1,1.1},
 LabelStyle->Directive[Black,20]
 ]
```

Output
```
{{3.44228,{x->0.703467}},{Steps->6,Function->14},
```

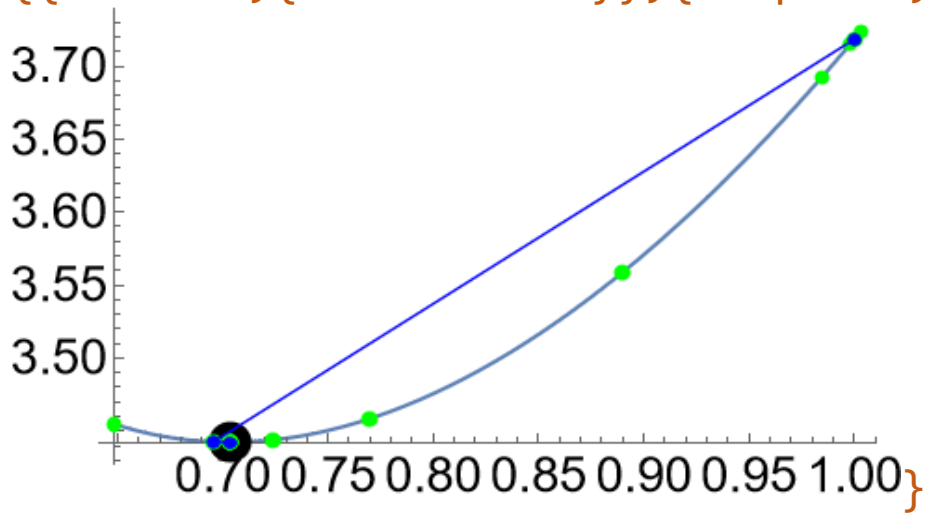

}

| | |
|---|---|
| `Minimize[f,x]` | minimizes f exactly with respect to x. |
| `Minimize[f,{x,y,…}]` | minimizes f exactly with respect to x, y, …. |
| `Minimize[{f,cons},{x,y,…}]` | minimizes f exactly subject to the constraints cons. |
| `Minimize[…,x∈rdom]` | constrains x to be in the region or domain rdom. |
| `Minimize[…,…,dom]` | constrains variables to the domain dom, typically Reals or Integers. |

## Details and Options

- Minimize finds the global minimum of `f` subject to the constraints given.
- Minimize returns a list of the form `{fmin,{x->xmin,y->ymin,…}}`.

- If `f` and `cons` are linear or polynomial, `Minimize` will always find a global minimum.
- The constraints `cons` can be any logical combination of:

| | |
|---|---|
| `lhs==rhs` | equations |
| `lhs>rhs, lhs≥rhs, lhs<rhs, lhs≤rhs` | inequalities (LessEqual,…) |
| `lhs≻rhs, lhs≽rhs, lhs≺rhs, lhs≼rhs` | vector inequalities (VectorLessEqual,…) |
| `Exists[…], ForAll[…]` | quantified conditions |
| `{x,y,…}∈rdom` | region or domain specification |

- By default, all variables are assumed to be real.
- `Minimize` will return exact results if given exact input. With approximate input, it automatically calls `NMinimize`.
- `Minimize` will return the following forms:

| | |
|---|---|
| `{fmin,{x→xmin,…}}` | finite minimum |
| `{∞,{x→Indeterminate,…}}` | infeasible, i.e. the constraint set is empty |
| `{-∞,{x→xmin,…}}` | unbounded, i.e. the values of f can be arbitrarily small |

- Even if the same minimum is achieved at several points, only one is returned.
- `Minimize[f,x,WorkingPrecision->n]` uses `n` digits of precision while computing a result.

***Mathematica Examples 5.3***

Input
```
Plot[
 2 x^2-3 x+5,
 {x,-10,10},
 LabelStyle->Directive[Black,20]
 ]
Minimize[2 x^2-3 x+5,x]
```
Output

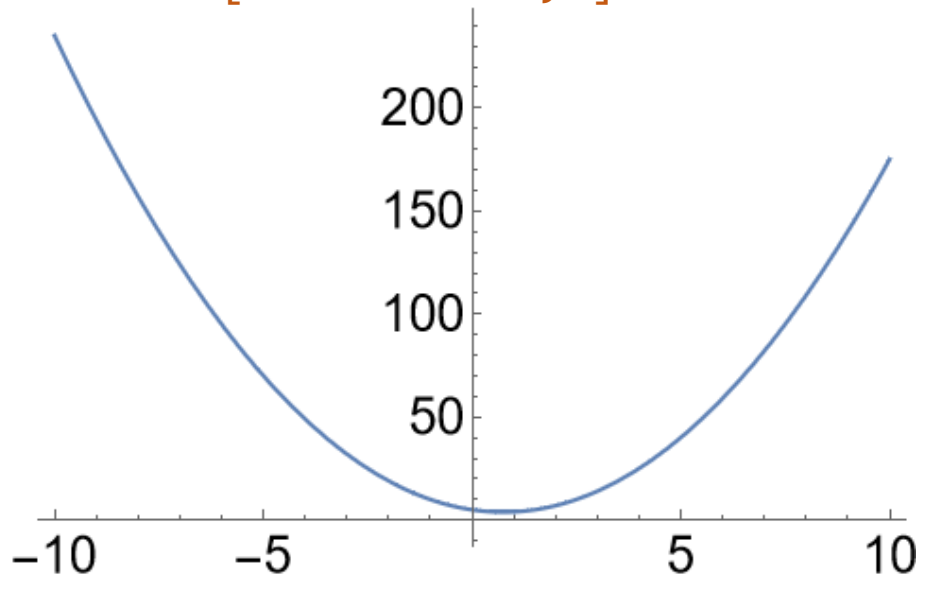

```
{31/8,{x->3/4}}
```

Input
```
Minimize[a x^2+b x+c,x]
```
Output
```
{{
  {\[Piecewise], {
    {c, (b==0&&a==0)||(b==0&&a>0)},
    {(-b2+4 a c)/(4 a), (b>0&&a>0)||(b<0&&a>0)},
    {-∞, True}
   }}
 },{x->{
    {\[Piecewise], {
      {-(b/(2 a)), (b>0&&a>0)||(b<0&&a>0)},
      {0, (b==0&&a==0)||(b==0&&a>0)},
      {Indeterminate, True}
     }}
   }}}
```

Input
```
Plot[
 Exp[x],
 {x,-10,10},
 LabelStyle->Directive[Black,20]
```

```
              ]
             Minimize[Exp[x],x]
Output
```

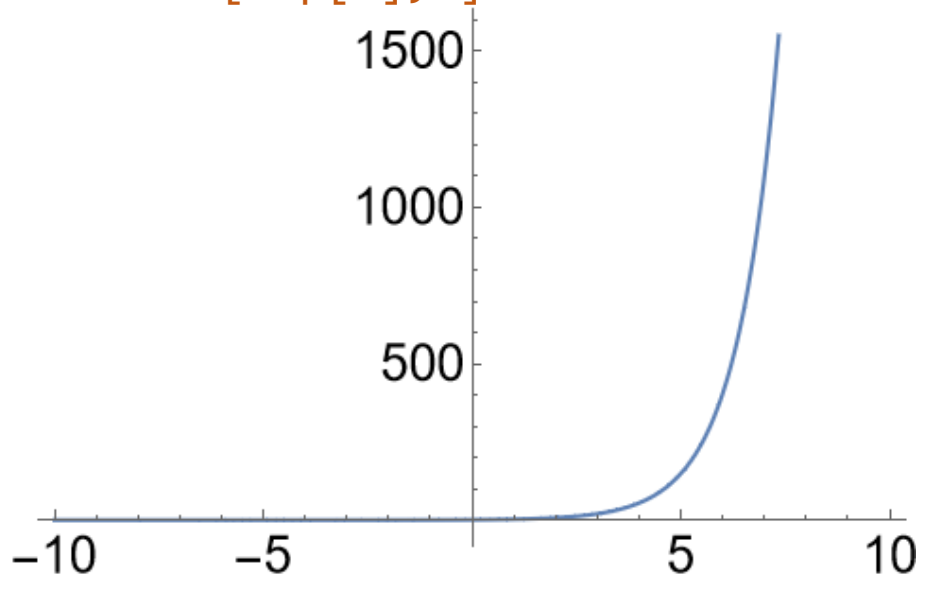

```
Output       Minimize::natt: The minimum is not attained at any point satisfying the given constraints.
             {0,{x->-∞}}

Input        f=Expand[(x-1) (x-2) (x-4) (x-7) (x-9) (x-10)];
             Minimize[f,x]
             Plot[
              f,
              {x,0.5,10.5},Epilog->{PointSize[Medium],Red,Point[{x/. %[[2]],%[[1]]}]},
              LabelStyle->Directive[Black,20]
              ]
Output       {-(35721/64),{x->11/2}}
Output
```

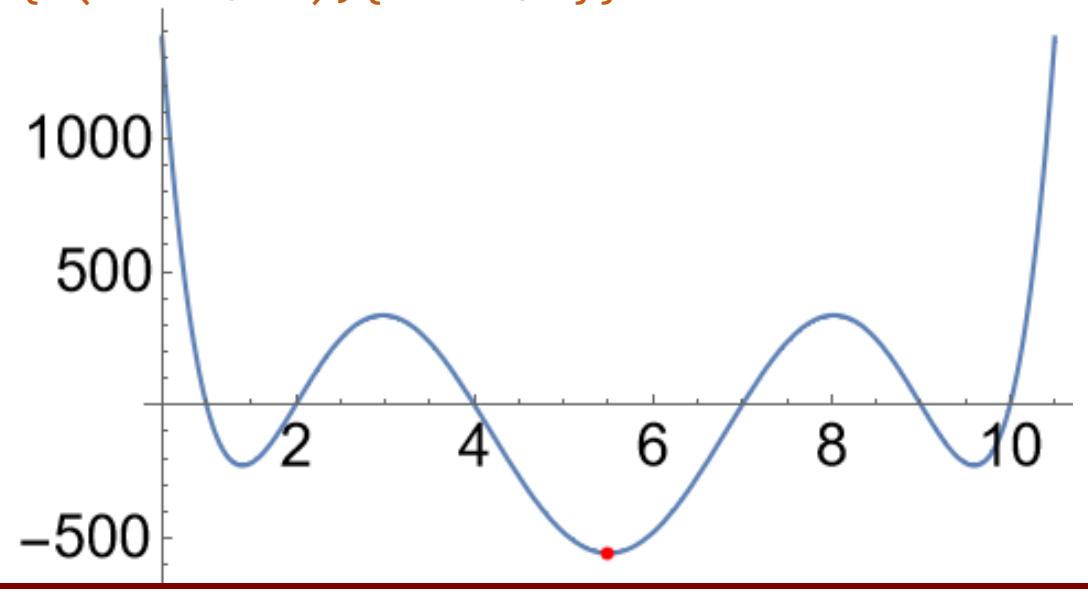

| | |
|---|---|
| `NMinimize[f,x]` | minimizes f numerically with respect to x. |
| `NMinimize[f,{x,y,…}]` | minimizes f numerically with respect to x, y, …. |
| `NMinimize[{f,cons},{x,y,…}]` | minimizes f numerically subject to the constraints cons. |
| `NMinimize[…,x∈rdom]` | constrains x to be in the region or domain rdom. |

- `NMinimize` always attempts to find a global minimum of `f` subject to the constraints given.
- `NMinimize` returns a list of the form `{fmin,{x->xmin,y->ymin,…}}`.
- If `f` and `cons` are linear or convex, the result given by `NMinimize` will be the global minimum, over both real and integer values; otherwise, the result may sometimes only be a local minimum.
- If `NMinimize` determines that the constraints cannot be satisfied, it returns `{Infinity,{x->Indeterminate,…}}`.

***Mathematica Examples 5.4***

```
Input        Plot[
              x^4-3 x^2-x,
              {x,-2,2},
              LabelStyle->Directive[Black,20]
              ]
             NMinimize[x^4-3 x^2-x,x]
```

Output

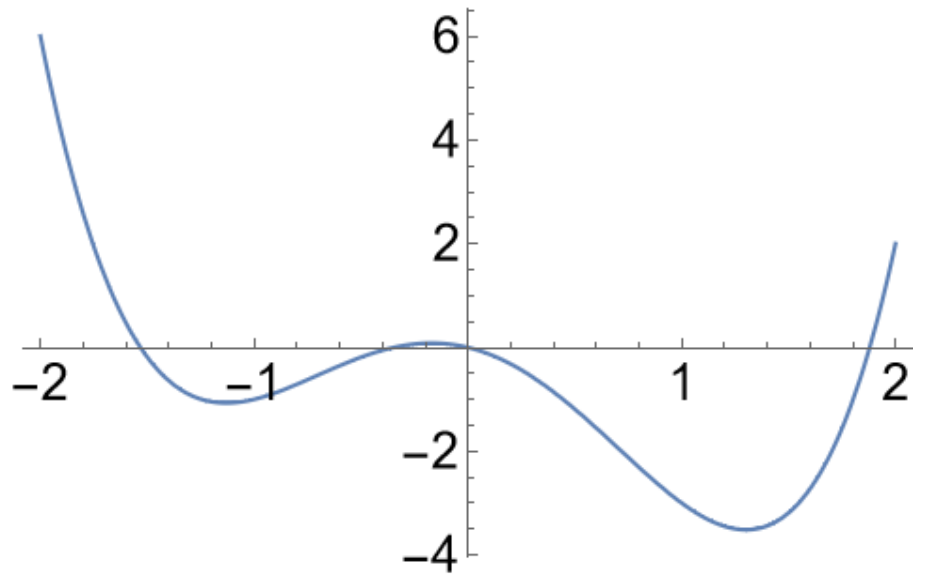

Output {-3.51391,{x->1.30084}}

| | |
|---|---|
| `FindRoot[f,{x,x0}]` | searches for a numerical root of f, starting from the point x=x0. |
| `FindRoot[lhs==rhs,{x,x0}]` | searches for a numerical solution to the equation lhs==rhs. |
| `FindRoot[{f1,f2,…},{{x,x0},{y,y0},…}]` | searches for a simultaneous numerical root of all the fi. |
| `FindRoot[{eqn1,eqn2,…},{{x,x0},{y,y0},…}]` | searches for a numerical solution to the simultaneous equations eqni. |

### Details and Options

- If the starting point for a variable is given as a list, the values of the variable are taken to be lists with the same dimensions.
- `FindRoot` returns a list of replacements for `x, y, …`, in the same form as obtained from `Solve`.
- `FindRoot[lhs==rhs,{x,x0,x1}]` searches for a solution using `x0` and `x1` as the first two values of `x`, avoiding the use of derivatives.
- If you specify only one starting value of `x, FindRoot` searches for a solution using Newton methods. If you specify two starting values, `FindRoot` uses a variant of the secant method.
- If all equations and starting values are real, then `FindRoot` will search only for real roots. If any are complex, it will also search for complex roots.
- `FindRoot` continues until either of the goals specified by `AccuracyGoal` or `PrecisionGoal` is achieved.

***Mathematica Examples 5.5***

```
Input       FindRoot[Sin[x]+Exp[x],{x,0}]
Output      {x->-0.588533}

Input       FindRoot[Cos[x]==x,{x,0}]
Output      {x->0.739085}
```

### Resource Functions

**BisectionMethodFindRoot:** Determine the root of an equation using the bisection method [14].

| | |
|---|---|
| `ResourceFunction["BisectionMethodFindRoot"][f,{x,xa,xb},tol,n]` | searches for a numerical root of f between the points xa and xb using tol digits and up to n steps. |
| `ResourceFunction["BisectionMethodFindRoot"][lhs==rhs,{x,xa,xb},tol,n]` | searches for a numerical solution to the equation lhs==rhs. |
| `ResourceFunction["BisectionMethodFindRoot"][f,{x,xa,xb},tol,n,property]` | returns a property of the search for the root of f. |

Details and Options

- BisectionMethodFindRoot supports two options for property:

| | |
|---|---|
| `"Solution"` | return the root of f |
| `"Steps"` | return a table of steps taken to reach the root |

- `"PropertyAssociation"` can be used to return an Association of the properties.
- `BisectionMethodFindRoot` terminates when the result is correct to the requested tolerance or the maximum number of steps has been taken, whichever comes first.

***Mathematica Examples 5.6***

```
Input     ResourceFunction["BisectionMethodFindRoot"][
           x-Sqrt[30],
           {x,5,6},
           5,
           100
           ]
Output    {x->5.4772}

Input     ResourceFunction["BisectionMethodFindRoot"][
           Cos[x]==x,
           {x,0,1},
           3,
           100,
           "Steps"
           ]
Output
```

| "steps" | "a" | "f[a]" | "b" | "f[b]" |
|---|---|---|---|---|
| 1 | 1.00 | -0.459698 | 0 | 1. |
| 2 | 1.00 | -0.459698 | 0.500 | 0.377583 |
| 3 | 0.750 | -0.0183111 | 0.500 | 0.377583 |
| 4 | 0.750 | -0.0183111 | 0.625 | 0.185963 |
| 5 | 0.750 | -0.0183111 | 0.688 | 0.0853349 |
| 6 | 0.750 | -0.0183111 | 0.719 | 0.0338794 |
| 7 | 0.750 | -0.0183111 | 0.734 | 0.00787473 |
| 8 | 0.742 | -0.00519571 | 0.734 | 0.00787473 |
| 9 | 0.742 | -0.00519571 | 0.738 | 0.00134515 |
| 10 | 0.740 | -0.00192387 | 0.738 | 0.00134515 |
| 11 | 0.739 | 0. | 0.739 | 0. |

**NewtonsMethodFindRoot:** Determine the root of an equation using Newton's method [15].

| | |
|---|---|
| `ResourceFunction["NewtonsMethodFindRoot"][f,{x,x0},tol]` | searches for a numerical root of f starting at x0 with digits equal to tol. |
| `ResourceFunction["NewtonsMethodFindRoot"][lhs==rhs,{x,x0},tol]` | searches for a numerical solution to the equation lhs==rhs. |
| `ResourceFunction["NewtonsMethodFindRoot"][f,{x,x0},tol,property]` | returns a property of the search for the root of f. |

***Mathematica Examples 5.7***

```
Input     ResourceFunction["NewtonsMethodFindRoot"][
           x-Sqrt[30],
           {x,5},
           10
           ]
Output    {x->5.477225575}
```

Input
```
ResourceFunction["NewtonsMethodFindRoot"][
 Cos[x]==x,
 {x,1},
 8,
 "Steps"
 ]
```
Output

| step | x | residual | derivative |
|---|---|---|---|
| 0 | 1. | 0.4596977 | -1.84147098 |
| 1 | 0.750364 | 0.0189231 | -1.68190495 |
| 2 | 0.739113 | 0.0000465 | -1.67363254 |
| 3 | 0.739085 | 0.*10^-8 | |

**NewtonMethod:** Approximate the root of a function using Newton's method [16].

| | |
|---|---|
| `ResourceFunction["NewtonMethod"][f,{x,x0},n]` | returns the root approximation obtained by applying Newton's method at most n times to a differentiable function f with starting value x0. |

***Mathematica Examples 5.8***

Input
```
Plot[
 x Sin[4 x],
 {x,0,3},
 LabelStyle->Directive[Black,20]
 ]
```
Output

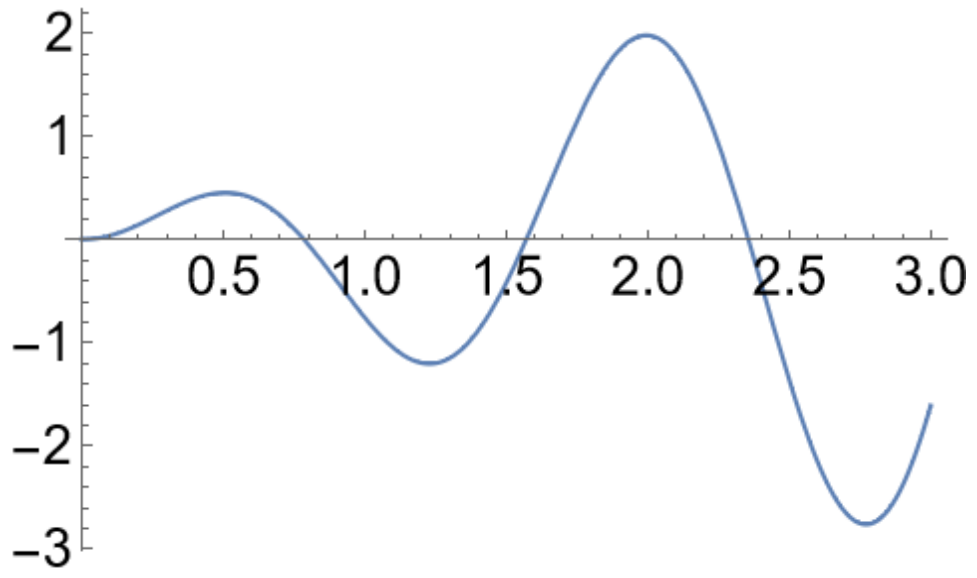

Input
```
ResourceFunction["NewtonMethod"][
 x Sin[4 x],
 {x,1.5}
 ]
```
Output
```
{x -> 1.5708}
```

**NewtonMethodPlot:** Plot the function together with a graphical display of the Newton iterations approximating its root [17].

| | |
|---|---|
| `ResourceFunction["NewtonMethodPlot"][f,{x,xmin,xmax},pt]` | returns a plot of f from x=xmin to x=xmax, together with illustrations representing the iterations of Newton's root-finding method, starting at x=pt. |

***Mathematica Examples 5.9***

Input
```
ResourceFunction["NewtonMethodPlot"][
 x^3-2+3 x,
 {x,-3,5},
 4,
 LabelStyle->Directive[Black,20]
 ]
```

```
Output     {
            {4.},
            {2.54902},
            {1.56161},
            {0.932191},
            {0.645648},
            {0.597163}
           }
Output
```

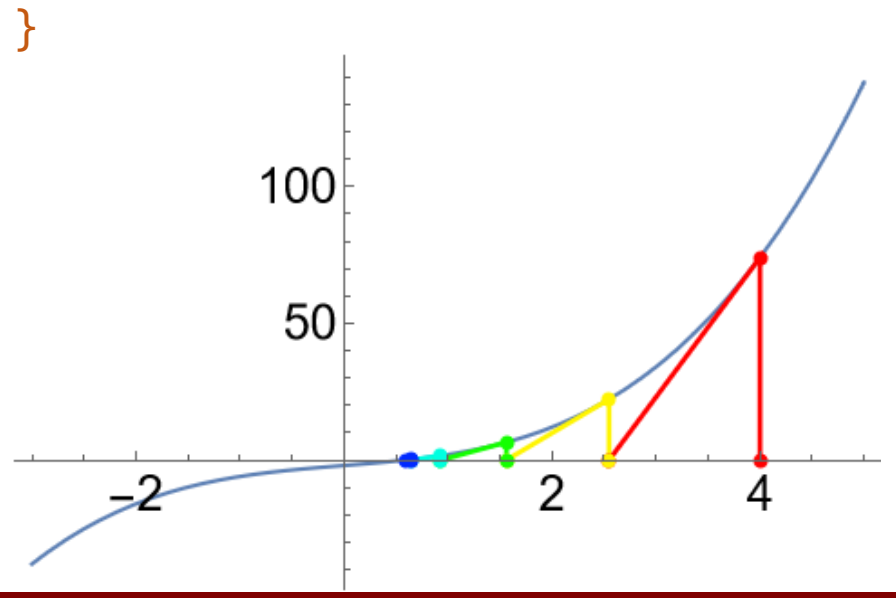

**SecantMethodFindRoot:** Determine the root of an equation using the secant method [18].

| | |
|---|---|
| `ResourceFunction["SecantMethodFindRoot"][f,{x,x0,x1},prec]` | searches for a numerical root of f starting at x0 with digits equal to prec. |
| `ResourceFunction["SecantMethodFindRoot"][lhs==rhs,{x,x0,x1},prec]` | searches for a numerical solution to the equation lhs==rhs. |
| `ResourceFunction["SecantMethodFindRoot"][f,{x,x0,x1},prec,property]` | returns a property of the search for the root of f. |

***Mathematica Examples 5.10***

```
Input      ResourceFunction["SecantMethodFindRoot"][
            x^2-30,
            {x,5,6},
            10
            ]
Output     {x->5.477225575}
```

**NumericalMethodFindRoot:** Find the root of an equation or number using a specified numerical method [19].

| | |
|---|---|
| `ResourceFunction["NumericalMethodFindRoot"][f,x,method]` | searches for a numerical root of f as a function of x, using the specified method. |
| `ResourceFunction["NumericalMethodFindRoot"][f,{x,x0}, method]` | searches for a numerical root of f, starting from the point x= x0. |
| `ResourceFunction["NumericalMethodFindRoot"][f,{x, x0},method,property]` | returns the specified property for the numerical search. |

- `ResourceFunction["NumericalMethodFindRoot"]` supports "Bisection", "Newton" and "Secant" methods.

***Mathematica Examples 5.11***

```
Input      ResourceFunction["NumericalMethodFindRoot"][
            x Cos[x],
            x,
            "Newton"
            ]
Output     {x->-7.85398}
```

Input
```
ResourceFunction["NumericalMethodFindRoot"][
 x Cos[x],
 {x,2},
 "Newton"
 ]
```
Output
```
{x->1.5708}
```

Input
```
grid=ResourceFunction["NumericalMethodFindRoot"][
  x Cos[x],
  {x,2},
  "Newton",
  "Steps"
  ]
```
Output

| "steps" | "x" | "f[x]" |
|---|---|---|
| 1 | 2. | -0.832294 |
| 2 | 1.62757 | -0.0923469 |
| 3 | 1.57265 | -0.00291944 |
| 4 | 1.5708 | -3.43469*10^-6 |
| 5 | 1.5708 | -4.78107*10^-12 |
| 6 | 1.5708 | 9.61835*10^-17 |
| 7 | 1.5708 | 9.61835*10^-17 |

**FindRootPlot:** Visualize the function evaluations done by `FindRoot`, [20].

| | |
|---|---|
| `ResourceFunction["FindRootPlot"][f,{x,xst}]` | plots the steps and the points at which the function f and any of its derivatives are evaluated in FindRoot[f,{x,xst}], superimposed on a plot of f versus x. |
| `ResourceFunction["FindRootPlot"][{f1,f2},{{x,xst},{y,yst}}]` | plots the steps and the points at which the pair of functions and their derivatives are evaluated, superimposed on a contour plot of the merit function. |
| `ResourceFunction["FindRootPlot"][f,range,property]` | returns the specified property. |

***Mathematica Examples 5.12***

Input
```
ResourceFunction["FindRootPlot"][
 Cos[x],
 {x,3,6},
 PlotRange->All,
 LabelStyle->Directive[Black,20]
 ]
```
Output

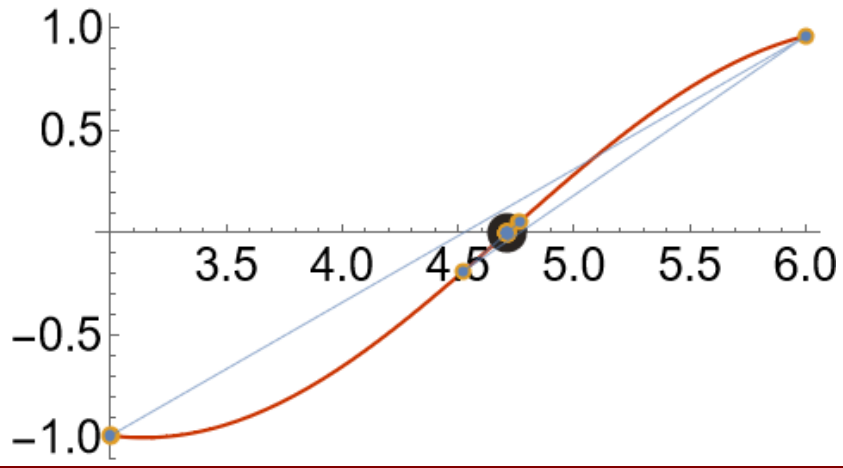

# CHAPTER 6

# MULTIVARIABLE UNCONSTRAINED OPTIMIZATION

This chapter includes the necessary and sufficient conditions for the minimum of an unconstrained function of several variables. We examine fundamental concepts (for instance, tangent planes, directional derivative, and Taylor series expansion of multivariable functions). First, we focus on the two-variables functions. The visualization of these functions in three dimensions is essential. Then, we expand these ideas to $n$-variables case.

## 6.1 Level Curves and Local Linearization of the Two-Variable Functions

One way to visualize a function of two variables is through its graph. The graph of $f$ is the surface with equation $z = f|\mathbf{x}\rangle$, $|\mathbf{x}\rangle = (x, y)^T$ [1]. See, for example, Figure 6.1.

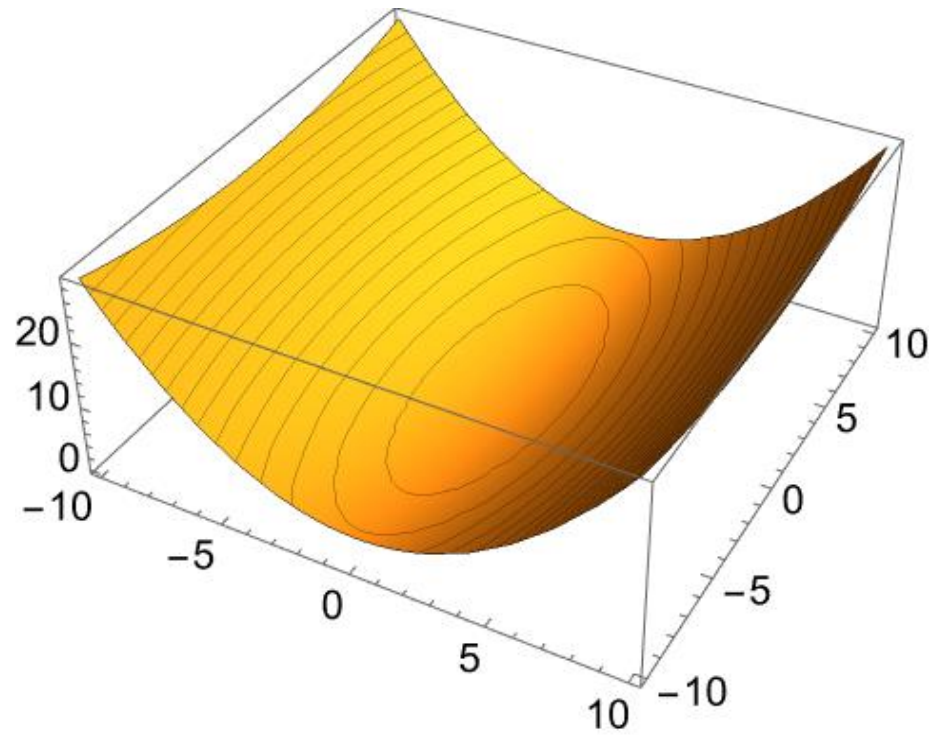

**Figure 6.1.** The graph of $f|\mathbf{x}\rangle = \frac{x^2}{4} + \frac{y^2}{25}$.

Another method for visualizing functions is a contour map on which points of constant elevation are joined to form contour lines or level curves.

**Definition (Level Curves):** The level curves of a function of two variables are the curves in the $xy$-plane with equations $f|\mathbf{x}\rangle = k$, where $k$ is a constant in the range of $f$.

More generally, a contour line for a function of two variables is a curve connecting points where the function has the same particular value (a constant value). It is a plane section of the three-dimensional graph of the function $f|\mathbf{x}\rangle$ parallel to the $xy$-plane. The surface is steep where the level curves are close together. It is somewhat flatter where they are farther apart, Figure 6.2.

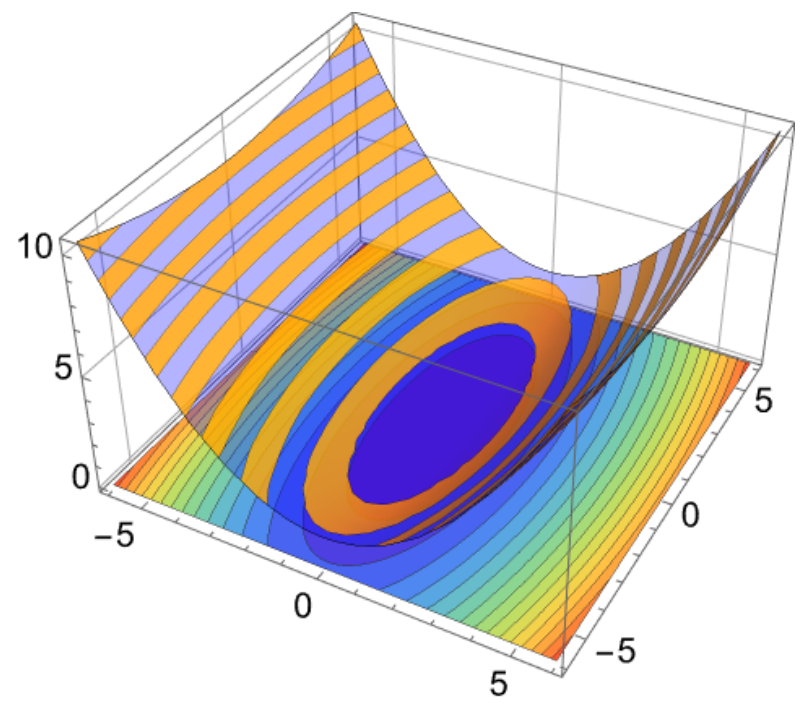

**Figure 6.2.** The level curve C with equation $f|\mathbf{x}\rangle = k$ is the projection of the trace of $f$ in the plane $z = k$ onto the $xy$-plane.

***Example 6.1***

Sketch a contour map for the surface described by $f|\mathbf{x}\rangle = x^2 + y^2$, using the level curves corresponding to $k = 0$, 1, 4, 9, and 16.

***Solution***

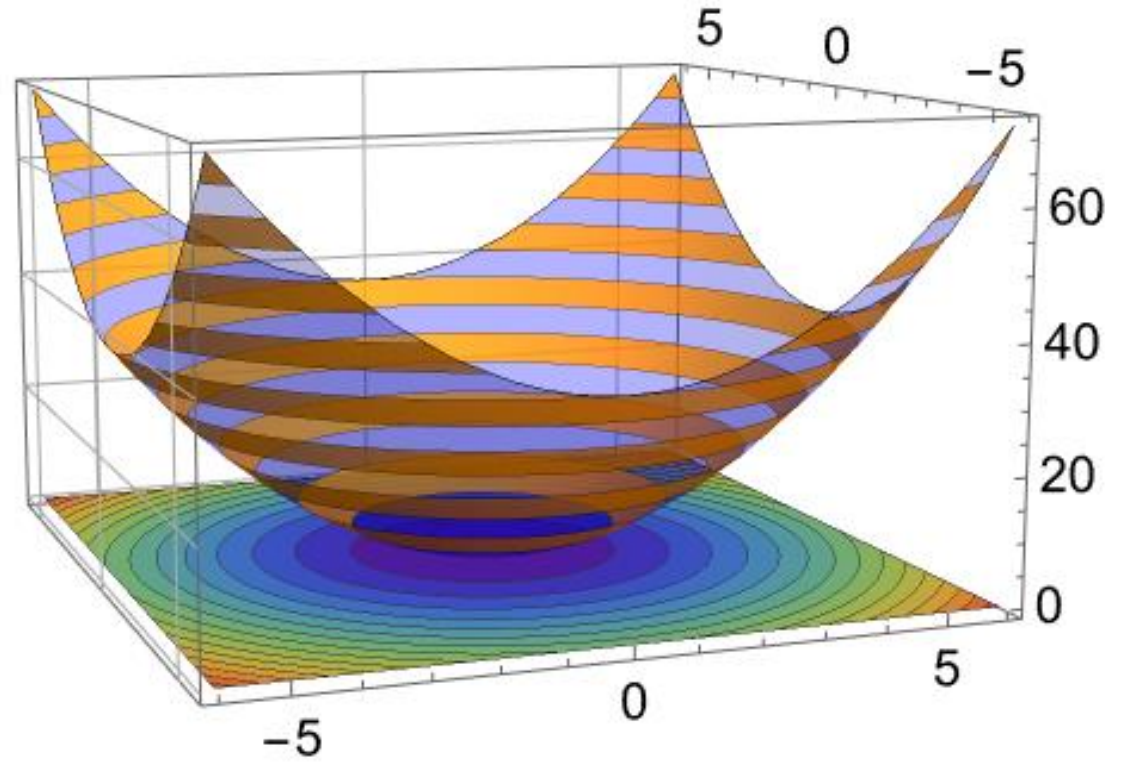

**Figure 6.3.** The level curve of the function $f|\mathbf{x}\rangle = x^2 + y^2$.

The level curve of $f$ corresponding to each value of $k$ is a circle $x^2 + y^2 = k$ of radius $\sqrt{k}$, centered at the origin, Figure 6.3. For example, if $k = 4$, the level curve is the circle with the equation $x^2 + y^2 = 4$, centered at the origin and having a radius 2.

## Taylor series for the two-variable functions and linearization

It will be helpful to review the Taylor series for a function of one variable and see how it extends to functions of more than one variable. Recall that the Taylor series for a function $f(x)$, based at a point $x = a$, is given by the following, where we assume that $f$ is analytic:

$$f(x) = f(a) + f'(a)(x - a) + \frac{1}{2!}f''(a)(x - a)^2 + \cdots. \tag{6.1}$$

Therefore, we can approximate $f$ using a constant:

$$f(x) \approx f(a), \tag{6.2}$$

or using a linear approximation (which is the tangent line to $f$ at $a$):

$$f(x) \approx f(a) + f'(a)(x - a), \tag{6.3}$$

or using a quadratic approximation:

$$f(x) \approx f(a) + f'(a)(x - a) + \frac{1}{2!}f''(a)(x - a)^2. \tag{6.4}$$

We can do something similar if $f$ depends on more than one variable as follows [1].

**Definition (Tangent Plane):** Provided $f$ is differentiable at $(a, b)$, the approximation $f(x, y)$:

$$f(x, y) \approx f(a, b) + f_x(a, b)(x - a) + f_y(a, b)(y - b), \tag{6.5.1}$$

is called the linear approximation or the tangent plane approximation of $f(x, y)$. In ket (column) notation

$$f|\mathbf{x}\rangle \approx f|\hat{\mathbf{x}}\rangle + f_x|\hat{\mathbf{x}}\rangle(x - a) + f_y|\hat{\mathbf{x}}\rangle(y - b), \tag{6.5.2}$$

where $|\mathbf{x}\rangle = (x, y)^T$, and $|\hat{\mathbf{x}}\rangle = (a, b)^T$.

**Remarks:**

- Since the tangent plane lies close to the surface near the point at which they meet, $z$-values on the tangent plane are close to values of $f|\mathbf{x}\rangle$ for points near $|\hat{\mathbf{x}}\rangle$. See Figure 6.4.
- We are thinking of $a$ and $b$ as fixed, so the expression on the right side of (6.5) is linear in $x$ and $y$. The right side of this approximation gives the local linearization of $f$ near $x = a$, $y = b$.
- For a function of one variable, local linearity means that as we zoom in on the graph, it looks like a straight line. As we zoom in on the graph of a two-variable function, the graph usually looks like a plane.

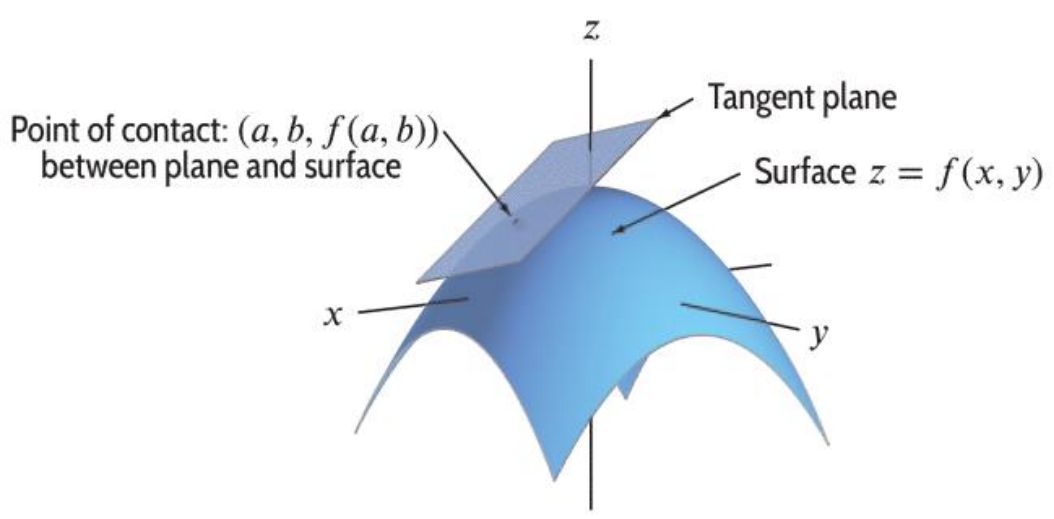

**Figure 6.4.** The tangent plane to the surface $z = f(x, y)$ at the point $(a, b)$

If we want to go further with a second-order (quadratic) approximation, it looks very similar. First, if $z = f(x, y)$ at $(a, b)$, the quadratic approximation looks like this:

$$f(x,y) \approx f(a,b) + f_x(a,b)(x-a) + f_y(a,b)(y-b) + \frac{1}{2}\left(f_{xx}(a,b)(x-a)^2 + 2f_{xy}(a,b)(x-a)(y-b) + f_{yy}(a,b)(y-b)^2\right), \tag{6.6.1}$$

where we assume that $f_{xy}(a,b) = f_{yx}(a,b)$. In ket notation, we have

$$f|\mathbf{x}\rangle \approx f|\hat{\mathbf{x}}\rangle + f_x|\hat{\mathbf{x}}\rangle(x-a) + f_y|\hat{\mathbf{x}}\rangle(y-b) + \frac{1}{2}\left[f_{xx}|\hat{\mathbf{x}}\rangle(x-a)^2 + 2f_{xy}|\hat{\mathbf{x}}\rangle(x-a)(y-b) + f_{yy}|\hat{\mathbf{x}}\rangle(y-b)^2\right], \tag{6.6.2}$$

where $|\mathbf{x}\rangle = (x,y)^T$, and $|\hat{\mathbf{x}}\rangle = (a,b)^T$. Moreover, the gradient of $f$, in 2D, is a vector of first partial derivatives:

$$\nabla f = \begin{pmatrix} f_x \\ f_y \end{pmatrix}, \tag{6.7}$$

and the $2 \times 2$ matrix of second partial derivatives is the Hessian matrix

$$\mathbf{H}_f = \begin{pmatrix} f_{xx} & f_{yx} \\ f_{xy} & f_{yy} \end{pmatrix}. \tag{6.8}$$

Using these notations, the linear approximation to $f$ at $|\hat{\mathbf{x}}\rangle = (a,b)^T$ is

$$f|\mathbf{x}\rangle \approx f|\hat{\mathbf{x}}\rangle + \langle \nabla f(\hat{\mathbf{x}})|\mathbf{x} - \hat{\mathbf{x}}\rangle, \tag{6.9}$$

and the quadratic approximation to $f$ is:

$$f(\mathbf{x}) \approx f|\hat{\mathbf{x}}\rangle + \langle \nabla f(\hat{\mathbf{x}})|\mathbf{x} - \hat{\mathbf{x}}\rangle + \frac{1}{2}\langle \mathbf{x} - \hat{\mathbf{x}}|\mathbf{H}_f(\hat{\mathbf{x}})|\mathbf{x} - \hat{\mathbf{x}}\rangle. \tag{6.10}$$

## 6.2 The Directional Derivative of Two-Variable Function

Suppose that $f$ is a function defined by the equation $z = f|\mathbf{x}\rangle$, $|\mathbf{x}\rangle = (x,y)^T$ and let $|\mathbf{p}\rangle \equiv (a,b)^T$ be a point in the domain $D$ of $f$. Furthermore, let $|\mathbf{u}\rangle$ be a unit vector having a specified direction. Then the vertical plane containing the line $L$ passing through $|\mathbf{p}\rangle$ and having the same direction as $|\mathbf{u}\rangle$ will intersect the surface $z = f|\mathbf{x}\rangle$ along a curve $C$ (see Figure 6.5). Intuitively, we see that the rate of change of $z$ at the point $|\mathbf{p}\rangle$ with respect to the distance measured along $L$ is given by the slope of the tangent line $T$ to the curve $C$ at the point $|\acute{\mathbf{p}}\rangle \equiv (a, b, f(a,b))^T$. To find the slope of $T$, first, observe that $|\mathbf{u}\rangle$ may be specified by writing $|\mathbf{u}\rangle = (u_1, u_2)^T$ for appropriate components $u_1$ and $u_2$. Equivalently, we may specify $|\mathbf{u}\rangle$ by giving the angle $\theta$ that it makes with the positive $x$-axis, in which case $u_1 = \cos\theta$ and $u_2 = \sin\theta$ (see Figure 6.6).

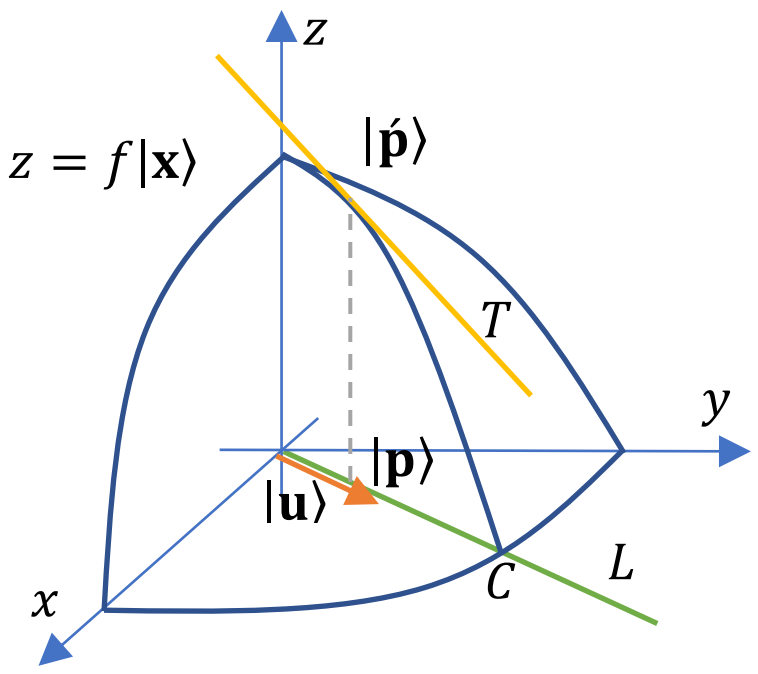

**Figure 6.5.** The rate of change of $z$ at $|\mathbf{p}\rangle$ with respect to the distance measured along $L$ is given by the slope of $T$.

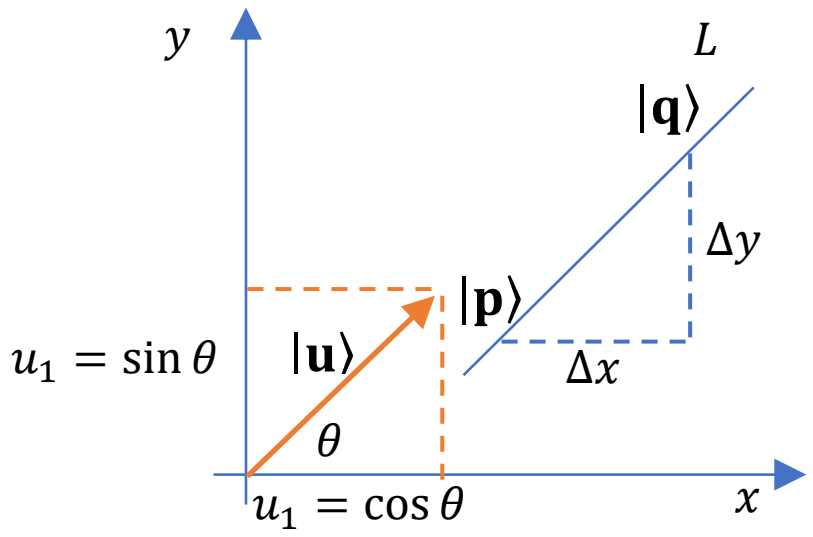

**Figure 6.6.** Any direction in the plane can be specified in terms of a unit vector $|\mathbf{u}\rangle$.

Next, let $|\mathbf{q}\rangle \equiv (a + \Delta x, b + \Delta y)^T$ be any point distinct from $|\mathbf{p}\rangle$ lying on the line $L$ passing through $|\mathbf{p}\rangle$ and having the same direction as $|\mathbf{u}\rangle$ (Figure 6.6). Since the vector $|\mathbf{pq}\rangle$ is parallel to $|\mathbf{u}\rangle$, it must be a scalar multiple of $|\mathbf{u}\rangle$. In other words, there exists a nonzero number $h$ such that

$$|\mathbf{pq}\rangle = h|\mathbf{u}\rangle = (hu_1, hu_2)^T. \tag{6.11}$$

But $|\mathbf{pq}\rangle$ is also given by $(\Delta x, \Delta y)^T$, and therefore,

$$\Delta x = hu_1, \qquad \Delta y = hu_2, \qquad h = \sqrt{(\Delta x)^2 + (\Delta y)^2}. \tag{6.12}$$

So, the point $|\mathbf{q}\rangle$ can be expressed as $|\mathbf{q}\rangle \equiv (a + hu_1, b + hu_2)^T$. Therefore, the slope of the secant line $S$ passing through the points $|\acute{\mathbf{p}}\rangle$ and $|\acute{\mathbf{q}}\rangle$ (see Figure 6.7) is given by

$$\frac{\Delta z}{h} = \frac{f|\mathbf{q}\rangle - f|\mathbf{p}\rangle}{h} = \frac{f(a + hu_1, b + hu_2)^T - f(a,b)^T}{h}. \tag{6.13}$$

Observe that (6.13) also gives the average rate of change of $z = f|\mathbf{x}\rangle$ from $|\mathbf{p}\rangle \equiv (a,b)^T$ to $|\mathbf{q}\rangle \equiv (a + \Delta x, b + \Delta y)^T = (a + hu_1, b + hu_2)^T$ in the direction of $|\mathbf{u}\rangle$.

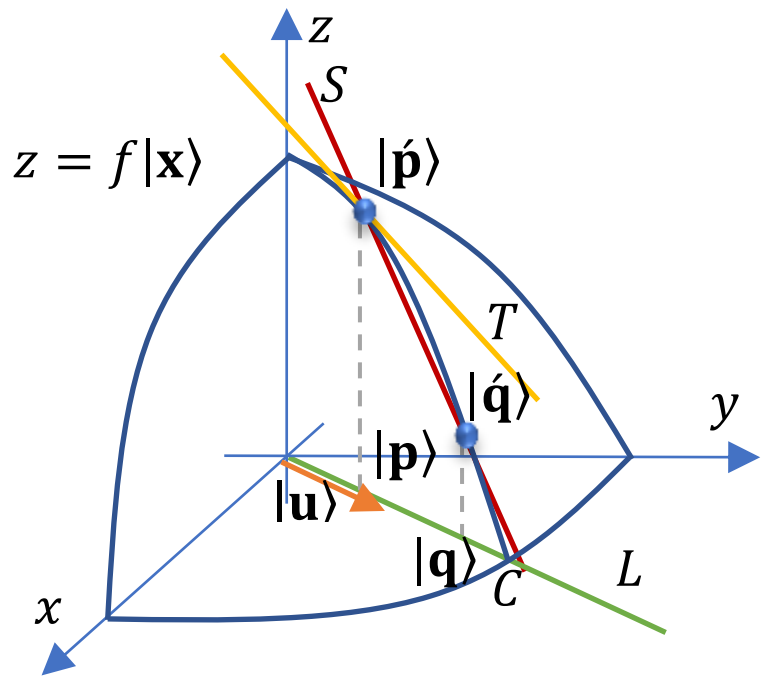

**Figure 6.7.** The secant line passes through the points $|\acute{\mathbf{p}}\rangle$ and $|\acute{\mathbf{q}}\rangle$ on the curve $C$.

If we let $h$ approach zero in (6.13), we see that the slope of the secant line $S$ approaches the slope of the tangent line at $|\acute{\mathbf{p}}\rangle$. Also, the average rate of change of $z$ approaches the (instantaneous) rate of change of $z$ at $|\mathbf{p}\rangle$ in the direction of $|\mathbf{u}\rangle$. This limit, whenever it exists, is called the directional derivative of $f$ at $|\mathbf{p}\rangle \equiv (a,b)^T$ in the direction of $|\mathbf{u}\rangle$. Since the point $|\mathbf{p}\rangle$ is arbitrary, we can replace it by $|\mathbf{p}\rangle \equiv (x,y)^T$ and define the directional derivative of $f$ at any point as follows.

**Definition (Directional Derivative):** Let $f$ be a function of $|\mathbf{x}\rangle = (x, y)^T$, and let $|\mathbf{u}\rangle = (u_1, u_2)^T$ be a unit vector. Then the directional derivative of $f$ at $|\mathbf{x}\rangle$ in the direction of $|\mathbf{u}\rangle$ is

$$D_{|\mathbf{u}\rangle} f|\mathbf{x}\rangle = \lim_{h\to 0} \frac{f|\mathbf{x} + h\mathbf{u}\rangle - f|\mathbf{x}\rangle}{h}, \qquad (6.14.a)$$

$$D_{|\mathbf{u}\rangle} f|\mathbf{x}\rangle = \lim_{h\to 0} \frac{f(x + hu_1, y + hu_2)^T - f(x, y)^T}{h}, \qquad (6.14.b)$$

if this limit exists.

**Remark:** If $|\mathbf{u}\rangle = |\mathbf{i}\rangle$ ($u_1 = 1$ and $u_2 = 0$), then (6.14) gives

$$D_{|\mathbf{i}\rangle} f|\mathbf{x}\rangle = \lim_{h\to 0} \frac{f(x + h, y)^T - f(x, y)^T}{h} = f_x|\mathbf{x}\rangle. \qquad (6.15)$$

That is, the directional derivative of $f$ in the $x$-direction is the partial derivative of $f$ in the $x$-direction. Similarly, we can show that $D_{|\mathbf{j}\rangle} f|\mathbf{x}\rangle = f_y|\mathbf{x}\rangle$.

**Theorem 6.1:** If $f$ is a differentiable function of $|\mathbf{x}\rangle = (x, y)^T$, then $f$ has a directional derivative in the direction of any unit vector $|\mathbf{u}\rangle = (u_1, u_2)^T$ and

$$D_{|\mathbf{u}\rangle} f|\mathbf{x}\rangle = f_x|\mathbf{x}\rangle u_1 + f_y|\mathbf{x}\rangle u_2. \qquad (6.16)$$

**Proof:**

Fix the point $|\hat{\mathbf{x}}\rangle = (a, b)^T$. Using local linearity, we have

$$\Delta f = f|\mathbf{x}\rangle - f|\hat{\mathbf{x}}\rangle \approx f_x|\hat{\mathbf{x}}\rangle \Delta x + f_y|\hat{\mathbf{x}}\rangle \Delta y = f_x|\hat{\mathbf{x}}\rangle h u_1 + f_y|\hat{\mathbf{x}}\rangle h u_2.$$

Thus, dividing by $h$ gives

$$\frac{\Delta f}{h} \approx f_x|\hat{\mathbf{x}}\rangle u_1 + f_y|\hat{\mathbf{x}}\rangle u_2.$$

This approximation becomes exact as $h \to 0$, so we have the following formula:

$$D_{|\mathbf{u}\rangle} f|\hat{\mathbf{x}}\rangle = f_x|\hat{\mathbf{x}}\rangle\, u_1 + f_y|\hat{\mathbf{x}}\rangle\, u_2.$$

Finally, since $|\hat{\mathbf{x}}\rangle = (a, b)^T$ is arbitrary, we may replace it by $|\mathbf{x}\rangle = (x, y)^T$ and the result follows.

■

## 6.3 The Gradient and Tangent Planes of Functions of Two- and Three-Variables

The directional derivative $D_{|\mathbf{u}\rangle} f|\mathbf{x}\rangle$ can be written as the dot product of the unit vector $|\mathbf{u}\rangle = (u_1, u_2)^T$ and the gradient vector $\nabla f = (f_x, f_y)^T$. We have,

$$D_{|\mathbf{u}\rangle} f|\mathbf{x}\rangle = f_x|\mathbf{x}\rangle u_1 + f_y|\mathbf{x}\rangle u_2 = \langle \nabla f(\mathbf{x})|\mathbf{u}\rangle. \qquad (6.17)$$

**Theorem 6.2:** If $f$ is a differentiable function of $|\mathbf{x}\rangle = (x, y)^T$, then $f$ has a directional derivative in the direction of any unit vector $|\mathbf{u}\rangle$, and

$$D_{|\mathbf{u}\rangle} f|\mathbf{x}\rangle = \langle \nabla f(\mathbf{x})|\mathbf{u}\rangle. \qquad (6.18)$$

Suppose $\theta$ is the angle between the vectors $\nabla f|\mathbf{x}\rangle$ and $|\mathbf{u}\rangle$. At the point $(a, b)^T$, we have

$$\langle \nabla f(\mathbf{x})|\mathbf{u}\rangle = \|\nabla f\|\|\mathbf{u}\| \cos\theta = \|\nabla f\| \cos\theta. \qquad (6.19)$$

Imagine that $\nabla f$ is fixed and that $|\mathbf{u}\rangle$ can rotate. (see Figure 6.8). The maximum value of $D_{|\mathbf{u}\rangle} f$ occurs when $\cos\theta = 1$, so $\theta = 0$ and $|\mathbf{u}\rangle$ is pointing in the direction of $\nabla f|\mathbf{x}\rangle$ .

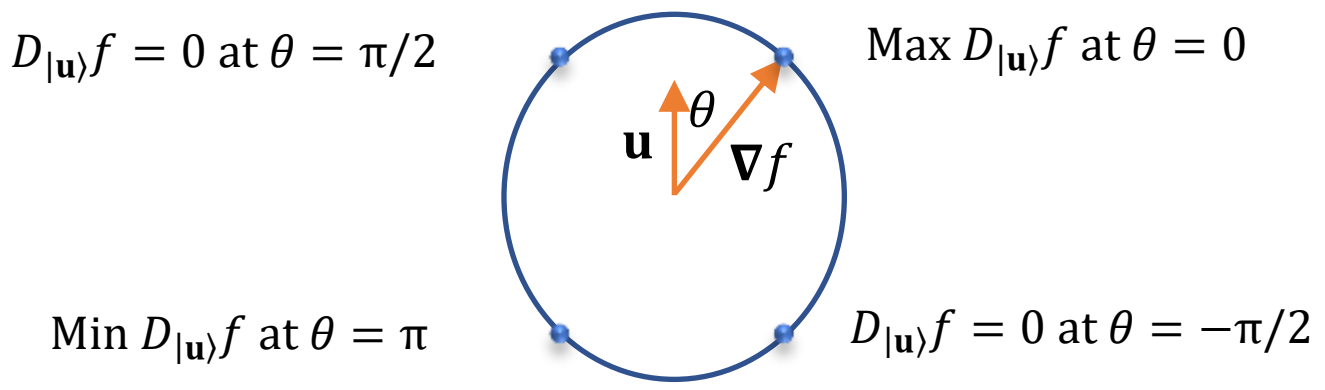

**Figure 6.8.** Values of the directional derivative at different angles to the gradient.

Then

$$\text{Max } D_{|\mathbf{u}\rangle} f = \|\nabla f\|. \tag{6.20}$$

The minimum value of $D_{|\mathbf{u}\rangle} f$ occurs when $\cos\theta = -1$, so $\theta = \pi$ and $|\mathbf{u}\rangle$ is pointing in the direction opposite to $\nabla f$. Then

$$\text{Min } D_{|\mathbf{u}\rangle} f = -\|\nabla f\|. \tag{6.21}$$

Hence, we have

**Theorem 6.3:** Suppose $f$ is differentiable at the point $|\mathbf{x}\rangle = (x, y)^T$.
1. If $\nabla f|\mathbf{x}\rangle = |\mathbf{0}\rangle$, then $D_{|\mathbf{u}\rangle} f|\mathbf{x}\rangle = 0$ for every $|\mathbf{u}\rangle$.
2. The maximum value of $D_{|\mathbf{u}\rangle} f|\mathbf{x}\rangle$ is $|\nabla f(\mathbf{x})|$, and this occurs when $|\mathbf{u}\rangle$ has the same direction as $\nabla f$.
3. The minimum value of $D_{|\mathbf{u}\rangle} f|\mathbf{x}\rangle$ is $-|\nabla f(\mathbf{x})|$, and this occurs when $|\mathbf{u}\rangle$ has the direction of $-\nabla f$.

**Remarks:**

1. Property (2) of Theorem 6.3 tells us that $f$ increases most rapidly in the direction of $\nabla f$. This direction is called the direction of the steepest ascent.
2. Property (3) of Theorem 6.3 says that $f$ decreases most rapidly in the direction of $-\nabla f$. This direction is called the direction of the steepest descent.

We can now give the geometric interpretation of the gradient in 2D. Suppose that the curve $C$ is represented by the vector function

$$|\mathbf{r}(t)\rangle = (g(t), h(t))^T, \tag{6.22}$$

where $g$ and $h$ are differentiable functions, $a = g(t_0)$ and $b = h(t_0)$, and $t_0$ lies in the parameter interval (Figure 6.9). Since the point $(x, y)^T = (g(t), h(t))^T$ lies on $C$, we have

$$f|\mathbf{x}\rangle = f|\mathbf{r}(t)\rangle = f(g(t), h(t))^T = c, \tag{6.23}$$

for all $t$ in the parameter interval.

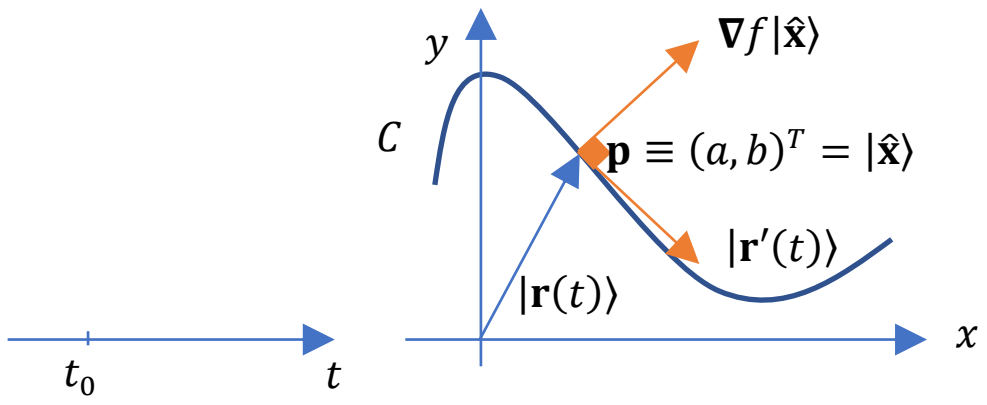

**Figure 6.9.** The curve may be represented by $|\mathbf{r}(t)\rangle = (x, y)^T = (g(t), h(t))^T$.

Differentiating both sides of this equation with respect to $t$ and using the chain rule for a function of two variables, we obtain

$$\frac{\partial f}{\partial x}\frac{dx}{dt} + \frac{\partial f}{\partial y}\frac{dy}{dt} = 0. \tag{6.24}$$

Using $\nabla f|\mathbf{x}\rangle = \left(\frac{\partial f}{\partial x}, \frac{\partial f}{\partial y}\right)^T$ and $|\mathbf{r}'(t)\rangle = \left(\frac{dx}{dt}, \frac{dy}{dt}\right)^T$, we can write (6.24) in the form

$$\langle \nabla f(\mathbf{x}) | \mathbf{r}'(t)\rangle = 0. \tag{6.25}$$

In particular, when $t = t_0$, i.e., $|\hat{\mathbf{x}}\rangle = (a, b)^T$, we have

$$\langle \nabla f(\hat{\mathbf{x}}) | \mathbf{r}'(t_0)\rangle = 0. \tag{6.26}$$

Thus, if $|\mathbf{r}'(t_0)\rangle \neq |\mathbf{0}\rangle$, the vector $\nabla f|\hat{\mathbf{x}}\rangle$ is orthogonal to the tangent vector $|\mathbf{r}'(t_0)\rangle$ at $(a, b)^T$. Loosely speaking, we have demonstrated the following:

**Theorem 6.4.1:** $\nabla f$ is orthogonal to the level curve $f|\mathbf{x}\rangle = c$ at point $|\hat{\mathbf{x}}\rangle$.

***Example 6.2***

Let $f|\mathbf{x}\rangle = x^2 - y^2$. Find the level curve of $f$ passing through the point $|\hat{\mathbf{x}}\rangle = (5,3)^T$. Also, find the gradient of $f$ at that point, and make a sketch of both the level curve and the gradient vector.

***Solution***

Since $f|\hat{\mathbf{x}}\rangle = 25 - 9 = 16$, the required level curve is the hyperbola $x^2 - y^2 = 16$. The gradient of $f$ at any point $|\mathbf{x}\rangle = (x, y)^T$ is

$$\nabla f|\mathbf{x}\rangle = (2x, -2y)^T,$$

and, in particular, the gradient of $f$ at $|\hat{\mathbf{x}}\rangle$ is

$$\nabla f|\hat{\mathbf{x}}\rangle = (10, -6)^T.$$

The level curve and $\nabla f|\mathbf{x}\rangle$ are shown in Figure 6.10.

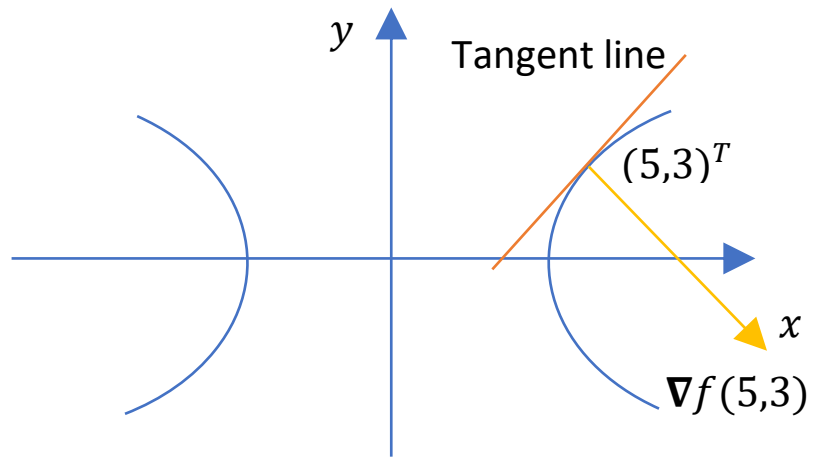

**Figure 6.10.** The gradient $\nabla f|\mathbf{x}\rangle$ is orthogonal to the level curve $x^2 - y^2 = 16$ at $|\hat{\mathbf{x}}\rangle = (5,3)^T$.

### Functions of Three Variables and Level Surfaces

A function $F$ of three variables is a rule that assigns to each ordered triple $|\mathbf{x}\rangle = (x, y, z)^T$ in a domain $D = \{|\mathbf{x}\rangle = (x, y, z)^T : x, y, z \in \mathbb{R}\}$ a unique real number $w$ denoted by $F|\mathbf{x}\rangle = w$. For example, $F|\mathbf{x}\rangle = xyz$. Since the graph of a function of three variables is composed of the points $(x, y, z, w)^T$, where $w = F|\mathbf{x}\rangle$, lying in four-dimensional space, we cannot draw the graphs of such functions. But by examining the level surfaces, which are the surfaces with equations $F|\mathbf{x}\rangle = k$, $k$ a constant, we are often able to gain some insight into the nature of $F$.

***Example 6.3***

Find the level surfaces of the function defined by $F|\mathbf{x}\rangle = x^2 + y^2 + z^2$.

***Solution***

The required level surfaces of $F$ are the graphs of the equations $x^2 + y^2 + z^2 = k$, where $k \geq 0$. These surfaces are concentric spheres of radius $\sqrt{k}$ centered at the origin (see Figure 6.11). Observe that $F$ has the same value for all points $(x, y, z)^T$ lying on any such sphere.

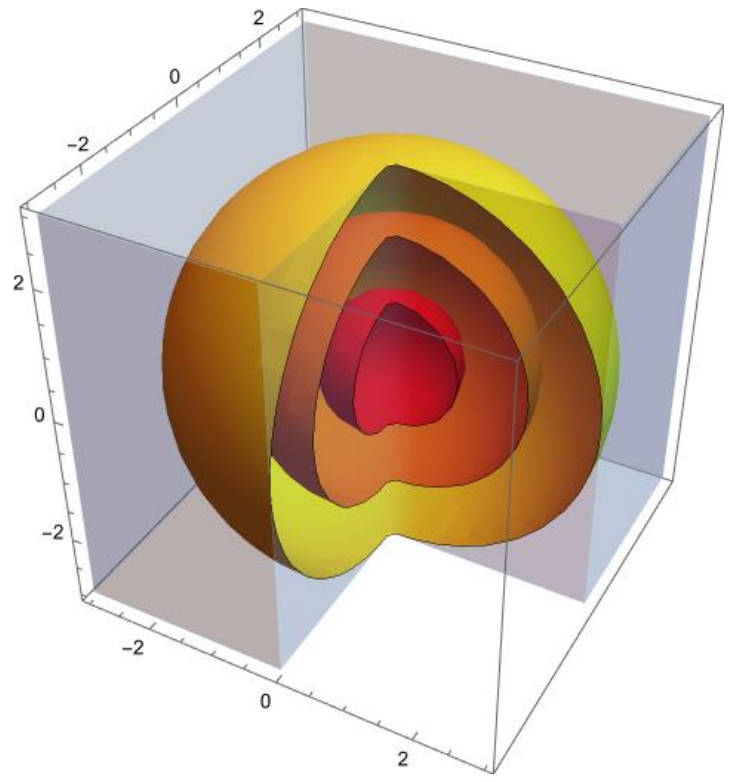

**Figure 6.11.** The level surfaces of $F|\mathbf{x}\rangle = x^2 + y^2 + z^2$ corresponding to $k = 1, 4, 9$.

Suppose that $F|\mathbf{x}\rangle = k$ is the level surface $S$ of a differentiable function $F$ defined by $T = F|\mathbf{x}\rangle$, $|\mathbf{x}\rangle = (x, y, z)^T$. Suppose that $|\hat{\mathbf{x}}\rangle \equiv (a, b, c)^T$ is a point on $S$ and let $C$ be a smooth curve on $S$ passing through $|\hat{\mathbf{x}}\rangle$. Then $C$ can be described by the vector function $|\mathbf{r}(t)\rangle = (f(t), g(t), h(t))^T$ where $f(t_0) = a$, $g(t_0) = b$ , $h(t_0) = c$, and $t_0$ is a point in the parameter interval (see Figure 6.12).

Since the point $|\mathbf{x}\rangle = (x, y, z)^T = (f(t), g(t), h(t))^T$ lies on $S$, we have, $F|\mathbf{x}\rangle = F(f(t), g(t), h(t))^T = k$, for all $t$ in the parameter interval. If $|\mathbf{r}\rangle$ is differentiable, then we can use the chain rule to differentiate both sides of this equation to obtain,

$$\frac{\partial F}{\partial x}\frac{dx}{dt} + \frac{\partial F}{\partial y}\frac{dy}{dt} + \frac{\partial F}{\partial z}\frac{dz}{dt} = 0. \tag{6.27}$$

This is the same as

$$\left(F_x, F_y, F_z\right) \cdot \begin{pmatrix} \frac{dx}{dt} \\ \frac{dy}{dt} \\ \frac{dz}{dt} \end{pmatrix} = 0, \tag{6.28}$$

or, in an even more abbreviated form, $\langle \nabla F | \mathbf{r}'(t)\rangle = 0$. At $t = t_0$, i.e., $|\hat{\mathbf{x}}\rangle = (a, b, c)^T$, we have $\langle \nabla F(\hat{\mathbf{x}}) | \mathbf{r}'(t_0)\rangle = 0$.

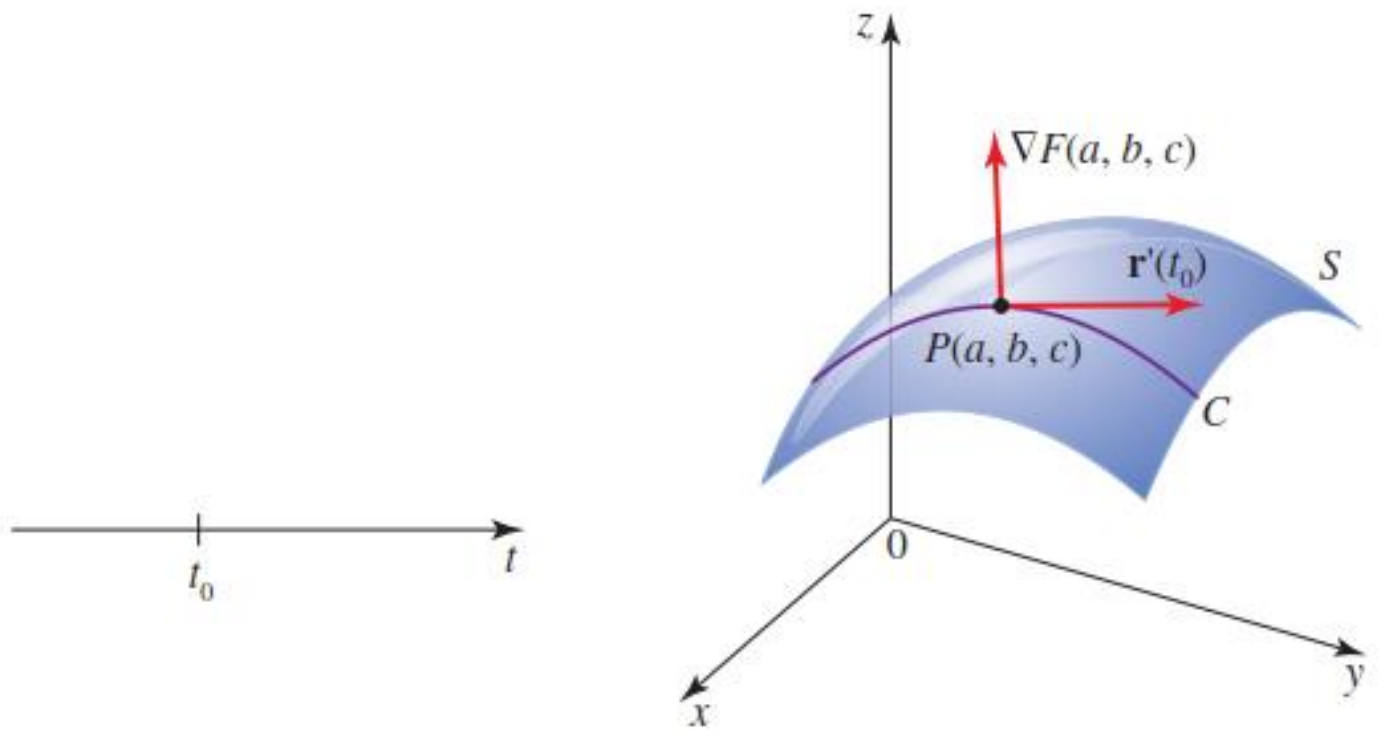

**Figure 6.12.** The curve $C$ is described by $|\mathbf{r}(t)\rangle = (f(t), g(t), h(t))^T$ with $P \equiv |\hat{\mathbf{x}}\rangle = (a, b, c)^T$ corresponding to $t_0$.

This shows that if $|\mathbf{r}'(t_0)\rangle \neq |\mathbf{0}\rangle$, then the gradient vector $\nabla F|\hat{\mathbf{x}}\rangle$ is orthogonal to the tangent vector $|\mathbf{r}'(t_0)\rangle$ to $C$ at $|\hat{\mathbf{x}}\rangle$ (see Figure 6.12). Since this argument holds for any differentiable curve passing through $|\hat{\mathbf{x}}\rangle$ on $C$, we have shown that $\nabla F|\hat{\mathbf{x}}\rangle$ is orthogonal to the tangent vector of every curve on $S$ passing through $|\hat{\mathbf{x}}\rangle$. Thus, we have demonstrated the following result.

**Theorem 6.4.2:** $\nabla F$ is orthogonal to the level surface $F|\mathbf{x}\rangle = 0$ at $|\hat{\mathbf{x}}\rangle$.

The gradient $\nabla F|\hat{\mathbf{x}}\rangle$ is orthogonal to the tangent vector of every curve on $S$ passing through $\mathbf{p} \equiv |\hat{\mathbf{x}}\rangle = (a, b, c)^T$ (Figure 6.12). This suggests that we define the tangent plane to $S$ at $\mathbf{p}$ to be the plane passing through $\mathbf{p}$ and containing all these tangent vectors. Equivalently, the tangent plane should have $\nabla F|\hat{\mathbf{x}}\rangle$ as a normal vector.

**Definition (Normal Line):** Let $\mathbf{p} \equiv |\hat{\mathbf{x}}\rangle = (a, b, c)^T$ be a point on the surface $S$ described by $F|\mathbf{x}\rangle = 0$, where $F$ is differentiable at $\mathbf{p}$, and suppose that $\nabla F|\hat{\mathbf{x}}\rangle \neq |\mathbf{0}\rangle$. Then the tangent plane to $S$ at $\mathbf{p}$ is the plane that passes through $\mathbf{p}$ and has normal vector $\nabla F|\hat{\mathbf{x}}\rangle$. The normal line to $S$ at $\mathbf{p}$ is the line that passes through $\mathbf{p}$ and has the same direction as $\nabla F|\hat{\mathbf{x}}\rangle$.

The equation of the tangent plane is

$$(x - a)F_x|\hat{\mathbf{x}}\rangle + (y - b)F_y|\hat{\mathbf{x}}\rangle + (z - c)F_z|\hat{\mathbf{x}}\rangle = 0. \tag{6.29}$$

In the following section, we can easily extend all of the above results and definitions for $n$-variables case and study the optimality criteria.

## 6.4 Optimality Criteria (General Case $n$-Variables)

The present section will largely consist of a generalization of the results of Chapter 5 to the case of $n$-variables [2-8]. We have,

$$\text{optimize: } z = f|\mathbf{x}\rangle, \qquad \text{where } |\mathbf{x}\rangle = (x_1, x_2, \dots, x_n)^T. \tag{6.30}$$

The optimality criteria can be derived using the definition of a locally optimal point and the Taylor series expansion of a multivariable function.

**Definition ($\epsilon$-Neighborhood):** An $\epsilon$-neighborhood ($\epsilon > 0$) around $|\hat{\mathbf{x}}\rangle$ is the set of all vectors $|\mathbf{x}\rangle$ such that

$$\begin{aligned} \|\mathbf{x} - \hat{\mathbf{x}}\| &= \langle \mathbf{x} - \hat{\mathbf{x}} | \mathbf{x} - \hat{\mathbf{x}} \rangle \\ &= (x_1 - \hat{x}_1)^2 + (x_2 - \hat{x}_2)^2 + \cdots + (x_n - \hat{x}_n)^2 \le \epsilon^2. \end{aligned} \tag{6.31}$$

In geometrical terms, an $\epsilon$-neighborhood around $|\hat{\mathbf{x}}\rangle$ is the interior and boundary of an $n$-dimensional sphere of radius $\epsilon$ centered at $|\hat{\mathbf{x}}\rangle$.

**Definition (Local and Global Minimizer):** An objective function $f|\mathbf{x}\rangle$ has a local minimizer at $|\hat{\mathbf{x}}\rangle$ if there exists an $\epsilon$-neighborhood around $|\hat{\mathbf{x}}\rangle$ such that $f|\mathbf{x}\rangle \ge f|\hat{\mathbf{x}}\rangle$ for all $|\mathbf{x}\rangle$ in this $\epsilon$-neighborhood at which the function is defined. If the condition is met for every positive $\epsilon$ (no matter how large), then $f|\mathbf{x}\rangle$ has a global minimizer at $|\hat{\mathbf{x}}\rangle$.

**Definition (Saddle Point):** A saddle point or minimax point is a point on the surface of the graph of a function where the slopes (derivatives) in orthogonal directions are all zero (a critical point) but which is not a local extremum of the function.

In the most general terms, a saddle point for a smooth function (whose graph is a curve, surface, or hypersurface) is a stationary point such that the curve/surface/etc. in the neighborhood of that point is not entirely on any side of the tangent space at that point. (see Figure 6.13)

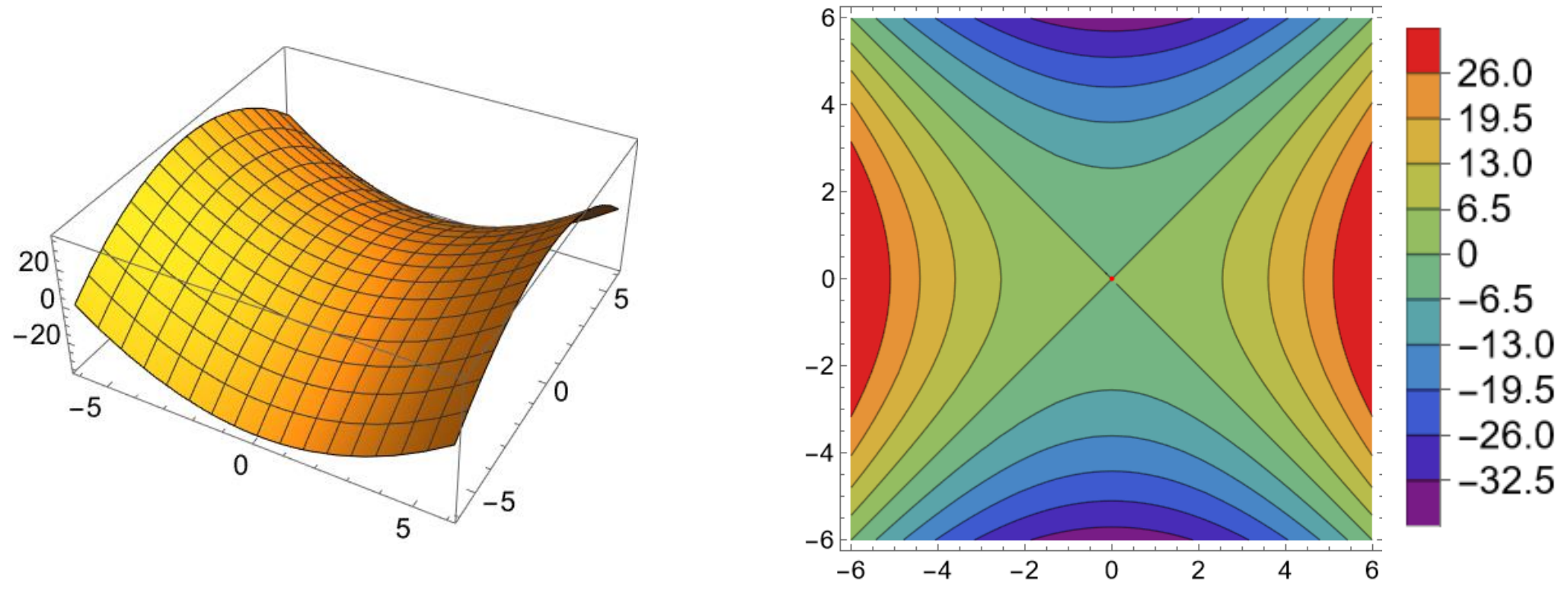

**Figure 6.13.** A saddle point (in red) on the graph of $z = x^2 - y^2$.

Remember, the gradient vector at any point $|\hat{\mathbf{x}}\rangle$ is represented by $\nabla f|\hat{\mathbf{x}}\rangle$ which is an $n$-dimensional vector given as follows:

$$\nabla f = \begin{pmatrix} \dfrac{\partial f}{\partial x_1} \\ \dfrac{\partial f}{\partial x_2} \\ \vdots \\ \dfrac{\partial f}{\partial x_n} \end{pmatrix}. \tag{6.32}$$

The first-order partial derivatives can be calculated numerically using finite difference. The second-order derivatives in multivariable functions form a matrix, $\nabla^2 f|\hat{\mathbf{x}}\rangle$ (the Hessian matrix) given as follows:

$$\mathbf{H}_f|\hat{\mathbf{x}}\rangle \equiv \nabla^2 f|\hat{\mathbf{x}}\rangle = \begin{pmatrix} \frac{\partial^2 f}{\partial x_1^2} & \frac{\partial^2 f}{\partial x_1 \partial x_2} & \cdots & \frac{\partial^2 f}{\partial x_1 \partial x_n} \\ \frac{\partial^2 f}{\partial x_2 \partial x_1} & \frac{\partial^2 f}{\partial x_2^2} & \cdots & \\ \cdots & \cdots & \ddots & \\ \frac{\partial^2 f}{\partial x_n \partial x_1} & \frac{\partial^2 f}{\partial x_n \partial x_2} & \cdots & \frac{\partial^2 f}{\partial x_n^2} \end{pmatrix} = \left(\frac{\partial^2 f}{\partial x_i \partial x_j}\right), \quad (i,j = 1,2,\ldots,n). \tag{6.33}$$

The second-order partial derivatives can also be calculated numerically using finite difference.

**Theorem 6.5 (Generalized of Theorem 6.3):** For small displacements from $|\hat{\mathbf{x}}\rangle$ in various directions, the direction of maximum increase in $f|\hat{\mathbf{x}}\rangle$ is the direction of the vector $\nabla f|\hat{\mathbf{x}}\rangle$.

**Proof:**

For any fixed vector $|\hat{\mathbf{x}}\rangle$ and any unit vector $|\mathbf{u}\rangle$, the directional derivative,

$$D_{|\mathbf{u}\rangle} f|\hat{\mathbf{x}}\rangle = \langle \nabla f(\hat{\mathbf{x}})|\mathbf{u}\rangle,$$

gives the rate of change of $f|\mathbf{x}\rangle$ at $|\hat{\mathbf{x}}\rangle$ in the direction of $|\mathbf{u}\rangle$. Since

$$\langle \nabla f|\mathbf{u}\rangle = \|\nabla f\|\|\mathbf{u}\| \cos\theta = \|\nabla f\| \cos\theta,$$

the greatest increase in $f|\mathbf{x}\rangle$ occurs when $\theta = 0$, i.e., when $|\mathbf{u}\rangle$ is in the same direction as $\nabla f$. Therefore, any (small) movement from $|\hat{\mathbf{x}}\rangle$ in the direction of $\nabla f|\hat{\mathbf{x}}\rangle$ will, initially, increase the function over $f|\hat{\mathbf{x}}\rangle$ as rapidly as possible.

■

**Definition (Feasible Direction):** A vector $|\mathbf{d}\rangle \in \mathbb{R}^n$, $|\mathbf{d}\rangle \neq |\mathbf{0}\rangle$, is a feasible direction at $|\mathbf{x}\rangle \in \Omega$ if there exists $\alpha_0 > 0$ such that $|\mathbf{x}\rangle + \alpha|\mathbf{d}\rangle \in \Omega$ for all $\alpha \in [0, \alpha_0]$. See Figure 6.14.

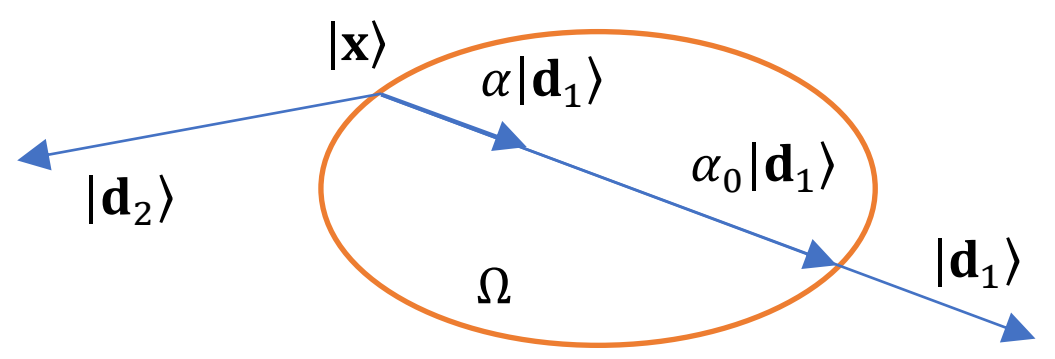

**Figure 6.14.** $|\mathbf{d}_1\rangle$ is a feasible direction, however $|\mathbf{d}_2\rangle$ is not a feasible direction.

**Theorem 6.6 (First-Order Necessary Condition):** Let $\Omega$ be a subset of $\mathbb{R}^n$ and $f \in C^1$ a real-valued function on $\Omega$. If $|\mathbf{x}^*\rangle$ is a local minimizer of $f$ over $\Omega$, then for any feasible direction $|\mathbf{d}\rangle$ at $|\mathbf{x}^*\rangle$, we have

$$\langle \mathbf{d}|\nabla f(\mathbf{x}^*)\rangle \geq 0. \tag{6.34}$$

If $|\mathbf{x}^*\rangle$ is located in the interior of $\Omega$, then

$$\nabla f(\mathbf{x}^*) = |\mathbf{0}\rangle. \tag{6.35}$$

**Proof:**

If $|\mathbf{d}\rangle$ is a feasible direction at $|\mathbf{x}^*\rangle$, then we have

$$|\mathbf{x}\rangle = |\mathbf{x}^*\rangle + \alpha|\mathbf{d}\rangle \in \Omega.$$

From the Taylor theorem,

$$f|\mathbf{x}\rangle = f|\mathbf{x}^*\rangle + \alpha\langle \mathbf{d}|\nabla f(\mathbf{x}^*)\rangle + o(\alpha).$$

If

$$\langle \mathbf{d}|\nabla f(\mathbf{x}^*)\rangle < 0,$$

then as $\alpha \to 0$

$$\alpha\langle\mathbf{d}|\nabla f(\mathbf{x}^*)\rangle + o(\alpha) < 0,$$

and so

$$f|\mathbf{x}\rangle < f|\mathbf{x}^*\rangle.$$

This contradicts the assumption that $|\mathbf{x}^*\rangle$ is a local minimizer. So, a necessary condition for $|\mathbf{x}^*\rangle$ to be a minimizer is $\langle\mathbf{d}|\nabla f(\mathbf{x}^*)\rangle \geq 0$.

If $|\mathbf{x}^*\rangle$ is in the interior of $\Omega$, vectors exist in all directions which are feasible. Thus, a direction $|\mathbf{d}\rangle = |\mathbf{d}_1\rangle$ yields

$$\langle\mathbf{d}_1|\nabla f(\mathbf{x}^*)\rangle \geq 0.$$

Similarly, for a direction $|\mathbf{d}\rangle = -|\mathbf{d}_1\rangle$

$$-\langle\mathbf{d}_1|\nabla f(\mathbf{x}^*)\rangle \geq 0.$$

In this case, a necessary condition for $|\mathbf{x}^*\rangle$ to be a local minimizer is

$$\nabla f(\mathbf{x}^*) = |\mathbf{0}\rangle.$$

■

**Theorem 6.7 (Second-Order Necessary Condition):** Let $\Omega$ be a subset of $\mathbb{R}^n$ and $f \in C^2$ a real-valued function on $\Omega$, $|\mathbf{x}^*\rangle$ is a local minimizer of $f$ over $\Omega$, and $|\mathbf{d}\rangle$ a feasible direction at $|\mathbf{x}^*\rangle$. If $\langle\mathbf{d}|\nabla f(\mathbf{x}^*)\rangle = 0$, then

$$\langle\mathbf{d}|\mathbf{H}(\mathbf{x}^*)|\mathbf{d}\rangle \geq 0, \tag{6.36}$$

where $\mathbf{H}$ is the Hessian of $f$.
If $|\mathbf{x}^*\rangle$ is a local minimizer in the interior of $\Omega$, then $\nabla f(\mathbf{x}^*) = |\mathbf{0}\rangle$ and $\langle\mathbf{d}|\mathbf{H}(\mathbf{x}^*)|\mathbf{d}\rangle \geq 0$ for all $|\mathbf{d}\rangle \neq |\mathbf{0}\rangle$.
(This condition is equivalent to stating that $\mathbf{H}(\mathbf{x}^*)$ is positive semidefinite.)

**Proof:**

If $|\mathbf{d}\rangle$ is a feasible direction at $|\mathbf{x}^*\rangle$, then we have

$$|\mathbf{x}\rangle = |\mathbf{x}^*\rangle + \alpha|\mathbf{d}\rangle \in \Omega.$$

From the Taylor theorem,

$$f|\mathbf{x}\rangle = f|\mathbf{x}^*\rangle + \alpha\langle\mathbf{d}|\nabla f(\mathbf{x}^*)\rangle + \frac{1}{2}\alpha^2\langle\mathbf{d}|\mathbf{H}(\mathbf{x}^*)|\mathbf{d}\rangle + o(\alpha^2).$$

If

$$\langle\mathbf{d}|\mathbf{H}(\mathbf{x}^*)|\mathbf{d}\rangle < 0,$$

then as $\alpha \to 0$

$$\frac{1}{2}\alpha^2\langle\mathbf{d}|\mathbf{H}(\mathbf{x}^*)|\mathbf{d}\rangle + o(\alpha^2) < 0,$$

and so

$$f|\mathbf{x}\rangle < f|\mathbf{x}^*\rangle.$$

This contradicts the assumption that $|\mathbf{x}^*\rangle$ is a local minimizer. So that, if $\langle\mathbf{d}|\nabla f(\mathbf{x}^*)\rangle = 0$ then $\langle\mathbf{d}|\mathbf{H}(\mathbf{x}^*)|\mathbf{d}\rangle \geq 0$.

If $|\mathbf{x}^*\rangle$ is a local minimizer in the interior of $\Omega$, then all vectors $|\mathbf{d}\rangle$ are feasible directions and, therefore, the second part holds.

■

***Example 6.4***

Let $f|\mathbf{x}\rangle = x_1^2 - x_2^2$. The Theorem 6.6 requires that $\nabla f|\mathbf{x}\rangle = (2x_1, -2x_2)^T = |\mathbf{0}\rangle$. Thus, $|\mathbf{x}\rangle = (0,0)^T$ satisfies the Theorem 6.6. The Hessian matrix of $f$ is

$$\mathbf{H} = \begin{pmatrix} 2 & 0 \\ 0 & -2 \end{pmatrix}.$$

The Hessian matrix is indefinite; that is, for some $|\mathbf{d}_1\rangle \in \mathbb{R}^2$ we have $\langle \mathbf{d}_1|\mathbf{H}|\mathbf{d}_1\rangle > 0$, e.g., $|\mathbf{d}_1\rangle = (1,0)^T$ and, for some $|\mathbf{d}_2\rangle$, we have $\langle \mathbf{d}_2|\mathbf{H}|\mathbf{d}_2\rangle < 0$, e.g., $|\mathbf{d}_2\rangle = (0,1)^T$. Thus, $|\mathbf{x}\rangle = (0,0)^T$ does not satisfy Theorem 6.7, and hence it is not a minimizer. Figure 6.15 shows the graph of $f|\mathbf{x}\rangle = x_1^2 - x_2^2$.

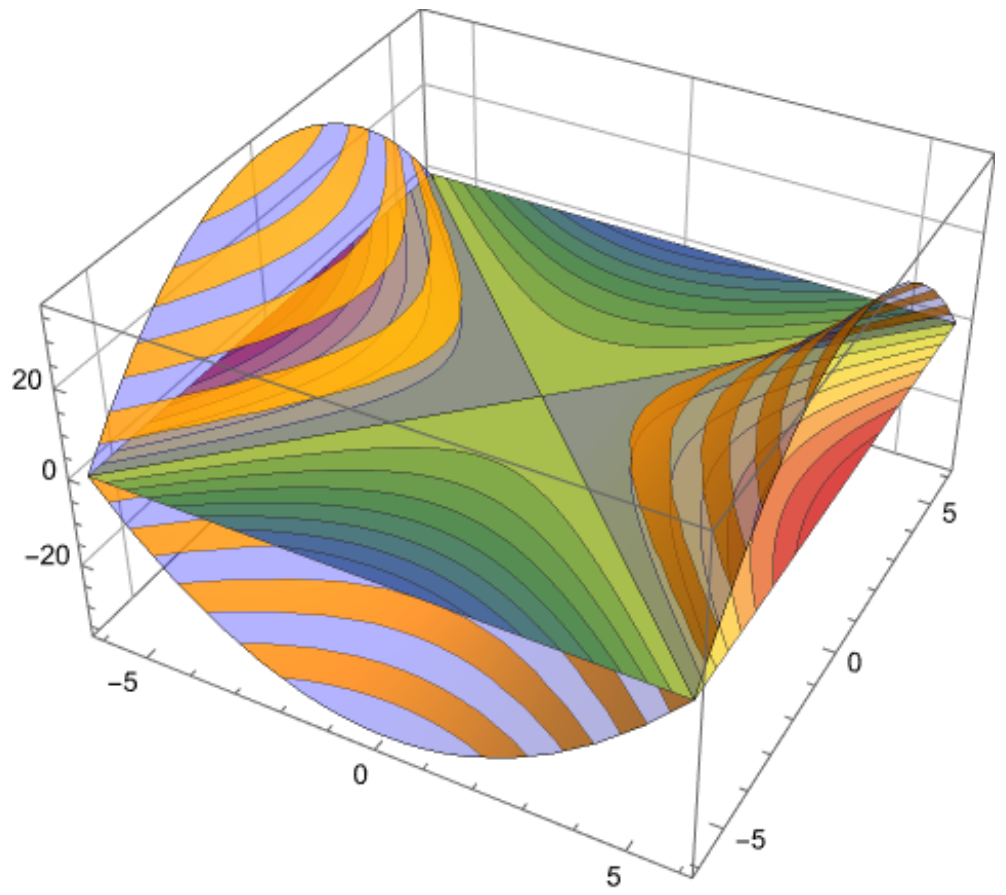

**Figure 6.15.** 3D and level curves of the function $f|\mathbf{x}\rangle = x_1^2 - x_2^2$.

***Example 6.5***

Let $f|\mathbf{x}\rangle = x_1^2 - x_1 + x_2 + x_1x_2, x_1 \geq 0, x_2 \geq 0$. We will show that the necessary conditions for $|\mathbf{x}^*\rangle = (0.5,0)^T$ to be a local minimizer are satisfied.

The gradient of the function $f|\mathbf{x}\rangle$ is

$$\nabla f|\mathbf{x}\rangle = (2x_1 - 1 + x_2, 1 + x_1)^T.$$

Let $|\mathbf{d}\rangle = (d_1, d_2)^T$ is a feasible direction, we obtain

$$\langle \mathbf{d}|\nabla f(\mathbf{x})\rangle = (2x_1 - 1 + x_2)d_1 + (1 + x_1)d_2,$$

and

$$\langle \mathbf{d}|\nabla f(\mathbf{x}^*)\rangle = \frac{3}{2}d_2.$$

Since $d_2 \geq 0$ for $|\mathbf{d}\rangle$ to be a feasible direction, we get

$$\langle \mathbf{d}|\nabla f(\mathbf{x})\rangle \geq 0.$$

So, Theorem 6.6, for a minimum, is satisfied.

if $d_2 = 0$, we obtain

$$\langle \mathbf{d}|\nabla f(\mathbf{x}^*)\rangle = 0.$$

The Hessian is

$$\mathbf{H}(\mathbf{x}^*) = \begin{pmatrix} 2 & 1 \\ 1 & 0 \end{pmatrix}.$$

so

$$\langle \mathbf{d}|\mathbf{H}(\mathbf{x}^*)|\mathbf{d}\rangle = 2d_1^2 + 2d_1d_2,$$

for $d_2 = 0$, we get

$$\langle \mathbf{d}|\mathbf{H}(\mathbf{x}^*)|\mathbf{d}\rangle = 2d_1^2 \geq 0,$$

for every feasible value of $d_1$. Therefore, Theorem 6.7 are satisfied.

**Theorem 6.8 (Second-Order Necessary Condition):** Let $\Omega$ be a subset of $\mathbb{R}^n$ and $f \in C^2$ a real-valued function on $\Omega$, $|\mathbf{x}^*\rangle$ is located in the interior of $\Omega$, then the conditions

1. $\nabla \boldsymbol{f}(\mathbf{x}^*) = |\mathbf{0}\rangle$

2. $\mathbf{H}(\mathbf{x}^*)$ is positive definite
are sufficient for $|\mathbf{x}^*\rangle$ to be a strict local minimizer.

**Proof:**

If $|\mathbf{d}\rangle$ is a feasible direction at $|\mathbf{x}^*\rangle$, then we have

$$|\mathbf{x}\rangle = |\mathbf{x}^*\rangle + \alpha|\mathbf{d}\rangle \in \Omega.$$

From the Taylor theorem,

$$f|\mathbf{x}\rangle = f|\mathbf{x}^*\rangle + \alpha\langle\mathbf{d}|\nabla f(\mathbf{x}^*)\rangle + \frac{1}{2}\alpha^2\langle\mathbf{d}|\mathbf{H}(\mathbf{x}^*)|\mathbf{d}\rangle + o(\alpha^2),$$

and if condition (1) is satisfied, we have,

$$f|\mathbf{x}\rangle = f|\mathbf{x}^*\rangle + \frac{1}{2}\alpha^2\langle\mathbf{d}|\mathbf{H}(\mathbf{x}^*)|\mathbf{d}\rangle + o(\alpha^2).$$

Now, if condition (2) is satisfied, then

$$\frac{1}{2}\alpha^2\langle\mathbf{d}|\mathbf{H}(\mathbf{x}^*)|\mathbf{d}\rangle + o(\alpha^2) > 0, \qquad \text{as } \|\mathbf{d}\| \to 0.$$

Therefore,

$$f|\mathbf{x}\rangle > f|\mathbf{x}^*\rangle.$$

that is, $|\mathbf{x}^*\rangle$ is a strict local minimizer.

■

***Example 6.6***

Let $f|\mathbf{x}\rangle = x_1^2 + x_2^2$. We have $\nabla f|\mathbf{x}^*\rangle = 2x_1 + 2x_2 = |\mathbf{0}\rangle$, if and only if $|\mathbf{x}\rangle = (0,0)^T$. For all $|\mathbf{x}\rangle \in \mathbb{R}^2$, we have

$$\mathbf{H} = \begin{pmatrix} 2 & 0 \\ 0 & 2 \end{pmatrix}.$$

The point $|\mathbf{x}\rangle = (0,0)^T$ satisfies Theorems 6.6, 6.7, and 6.8. It is a strict local minimizer. Actually $|\mathbf{x}\rangle = (0,0)^T$ is a strict global minimizer. Figure 6.16 shows the graph of $f|\mathbf{x}\rangle = x_1^2 + x_2^2$.

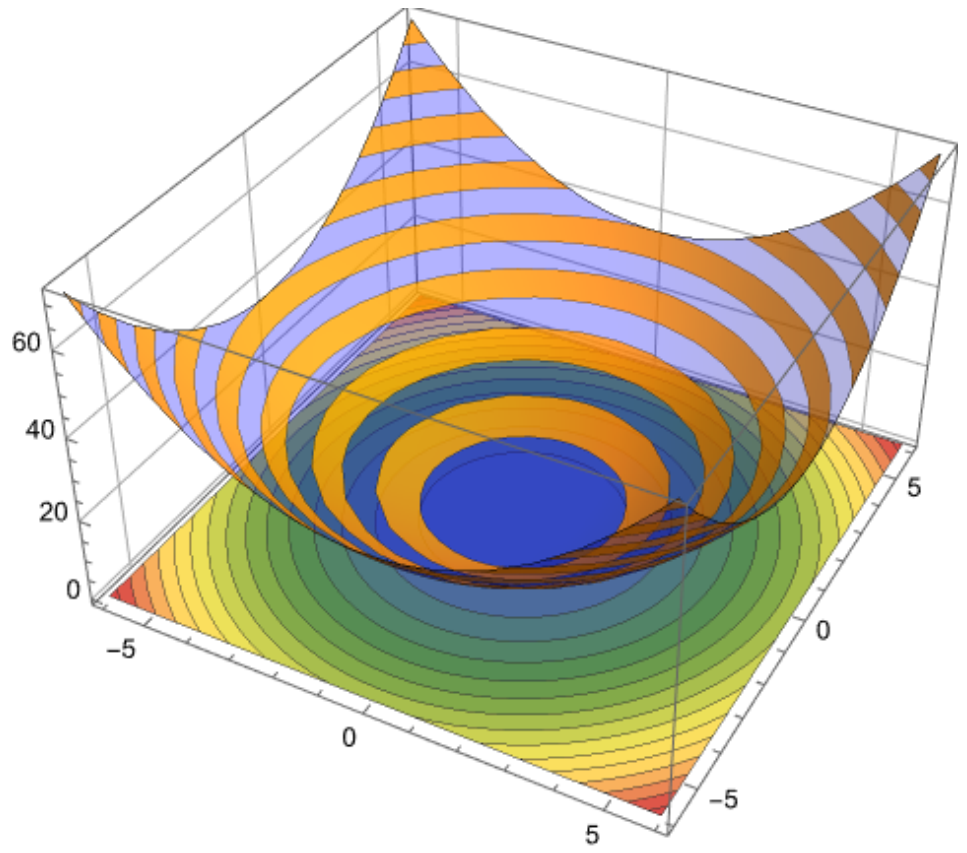

**Figure 6.16.** 3D and level curves of the function $f|\mathbf{x}\rangle = x_1^2 + x_2^2$.

Analytical solutions based on calculus are even harder to obtain for multivariable programs than for single-variable programs, and so, once again, numerical methods are used to approximate (local) maxima to within prescribed tolerances. In subsequent chapters, we will discuss numerical algorithms in detail.

## 6.5 Mathematica Built-in Functions

The commands `FindMinimum` and `NMinimize,` and `FindMinimumPlot` can do unconstrained optimization for both single-variable and multivariable functions. This section represents some examples in two dimensions [9,10].

| | |
|---|---|
| `ContourPlot[f,{x,xmin,xmax},{y,ymin,ymax}]` | generates a contour plot of f as a function of x and y. |
| `ContourPlot3D[f,{x,xmin,xmax},{y,ymin,ymax},{z,zmin,zmax}]` | produces a three-dimensional contour plot of f as a function of x, y, and z. |

***Mathematica Examples 6.1***

Input
```
ContourPlot3D[
 x^3+y^2-z^2==0,
 {x,-2,2},
 {y,-2,2},
 {z,-2,2},
 LabelStyle->Directive[Black,16]
 ]
```
Output

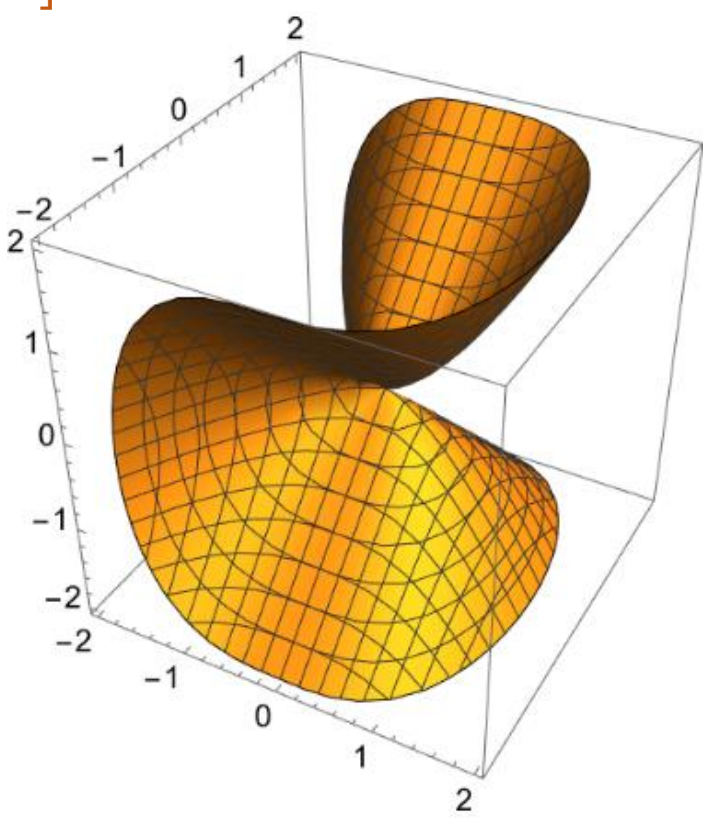

Input
```
ContourPlot3D[
 x^2+y^2+z^2,
 {x,-2,2},
 {y,-2,2},
 {z,-2,2},
 Contours->5,
 RegionFunction->Function[{x,y,z},x<0||y>0],
 ContourStyle->{Red,Orange,Yellow,Green,Blue},
 Mesh->None,
 LabelStyle->Directive[Black,16]
 ]
```
Output

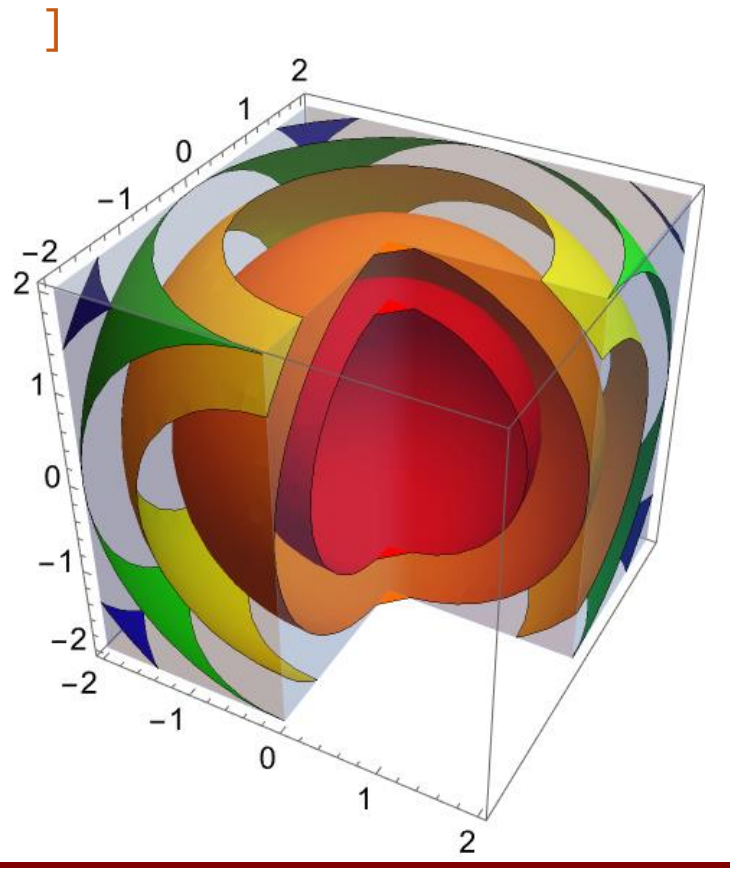

| | |
|---|---|
| `ArgMin[f,x]` | gives a position xmin at which f is minimized. |
| `ArgMin[f,{x,y,…}]` | gives a position {xmin,ymin,…} at which f is minimized. |
| `MinValue[f,x]` | gives the minimum value of f with respect to x. |
| `MinValue[f,{x,y,…}]` | gives the exact minimum value of f with respect to x, y, …. |
| `NMinimize[f,x]` | minimizes f numerically with respect to x. |
| `NMinimize[f,{x,y,…}]` | minimizes f numerically with respect to x, y, …. |
| `NMinValue[f,x]` | gives the minimum value of f with respect to x. |
| `NMinValue[f,{x,y,…}]` | gives the minimum value of f with respect to x, y, …. |
| `Sow[e]` | specifies that e should be collected by the nearest enclosing Reap. |
| `Reap[expr]` | gives the value of expr together with all expressions to which Sow has been applied during its evaluation. Expressions sown using Sow[e] or Sow[e,tagi] with different tags are given in different lists. |

***Mathematica Examples 6.2***

```
Input     (* Find a minimizer point for a univariate function: *)
          ArgMin[
           2 x^2-3 x+5,
           x
           ]
Output    ¾

Input     (* Find a minimizer point for a multivariate function: *)
          ArgMin[
           (x y-3)^2+1,
           {x,y}
           ]
Output    {-1,-3}

Input     (* Find the minimum value of a univariate function: *)
          MinValue[
           2 x^2-3 x+5,
           x
           ]
Output    31/8

Input     (* Find the minimum value of a multivariate function: *)
          MinValue[
           (x y-3)^2+1,
           {x,y}
           ]
Output    1

Input     (* Find the global minimum of an unconstrained problem: *)
          NMinimize[
           x^4-3 x^2-x,
           x
           ]
Output    {-3.51391,{x->1.30084}}

Input     (* Find the global minimum value of a univariate function: *)
          NMinValue[
           2 x^2-3 x+5,
           x
           ]
Output    3.875
```

```
Input     (* Find the global minimum value of a multivariate function: *)
          NMinValue[
           (x y-3)^2+1,
           {x,y}
           ]
Output    1.
```

| | |
|---|---|
| `FindMinimum[f,{{x,x0},{y,y0},…}]` | searches for a local minimum in a function of several variables. |
| `FindMinimumPlot[f,{{x,xst},{y,yst}}]` | plots the steps and the points at which the bivariate function f and any of its derivatives are evaluated, superimposed on a contour plot of f as a function of x and y. |
| `Minimize[f,{x,y,…}]` | minimizes f exactly with respect to x, y, …. |
| `NMinimize[f,{x,y,…}]` | minimizes f numerically with respect to x, y, …. |

***Mathematica Examples 6.3***

```
Input     <<Optimization`UnconstrainedProblems`

Input     Plot3D[
           Sin[x] Sin[2 y],
           {x,-3,3},
           {y,-3,3},
           ColorFunction->"Rainbow",
           PlotLegends->BarLegend[Automatic],
           LabelStyle->Directive[Black,16]
           ]

          ContourPlot[
           Sin[x] Sin[2 y],
           {x,-3,3},
           {y,-3,3},
           PlotLegends->Automatic,
           Contours->10,
           ColorFunction->"Rainbow",
           LabelStyle->Directive[Black,16]
           ]

          FindMinimum[
           Sin[x] Sin[2 y],
           {{x,2},{y,2}}
           ]

          FindMinimumPlot[
           Sin[x] Sin[2 y],
           {{x,2},{y,2}},
           Method->"Newton"
           ]

          FindMinimumPlot[
           Sin[x] Sin[2 y],
           {{x,2},{y,2}},
           Method->"QuasiNewton"
           ]

          FindMinimumPlot[
           Sin[x] Sin[2 y],
           {{x,2},{y,2}},
           Method->"PrincipalAxis"
           ]
```

Output

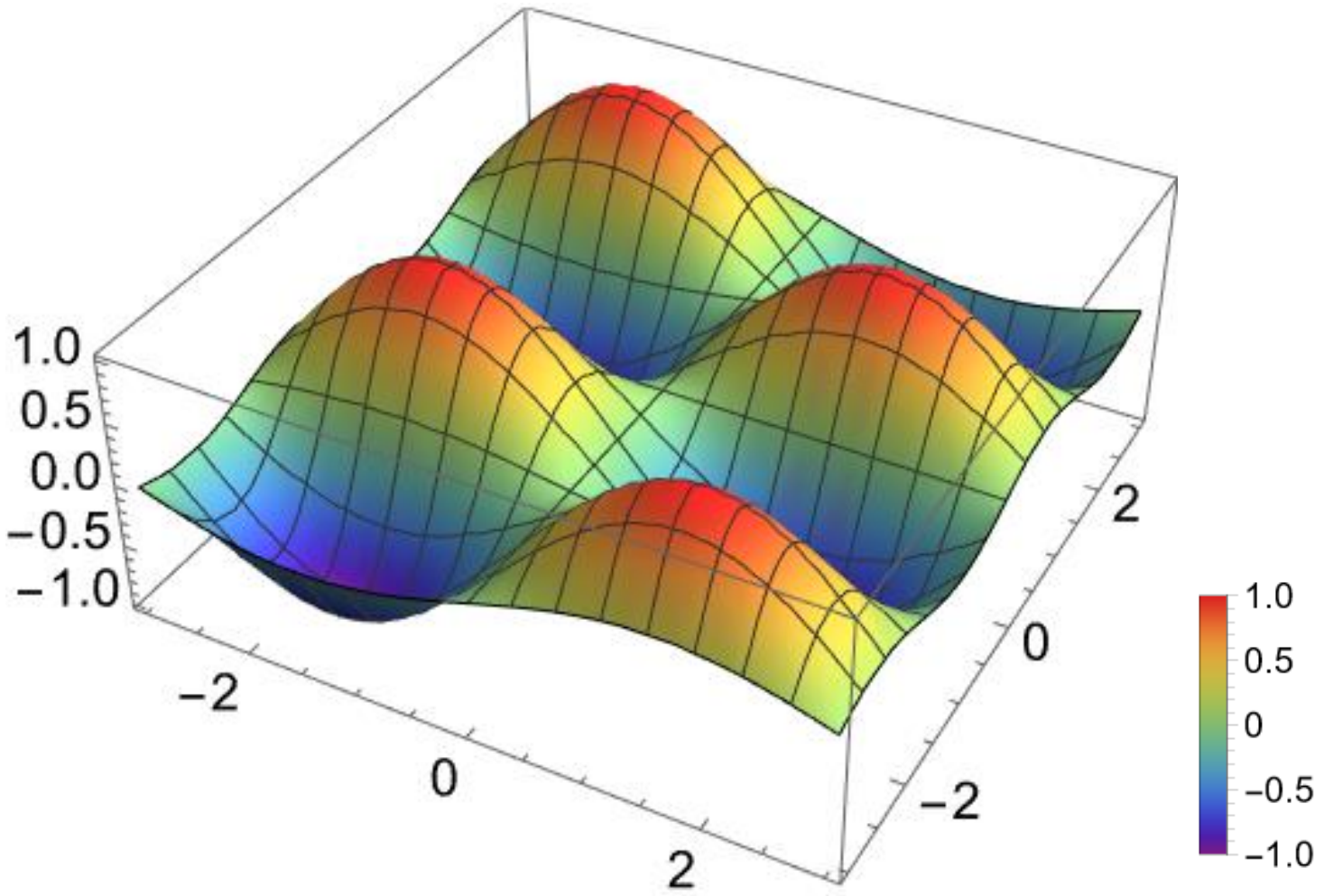

Output

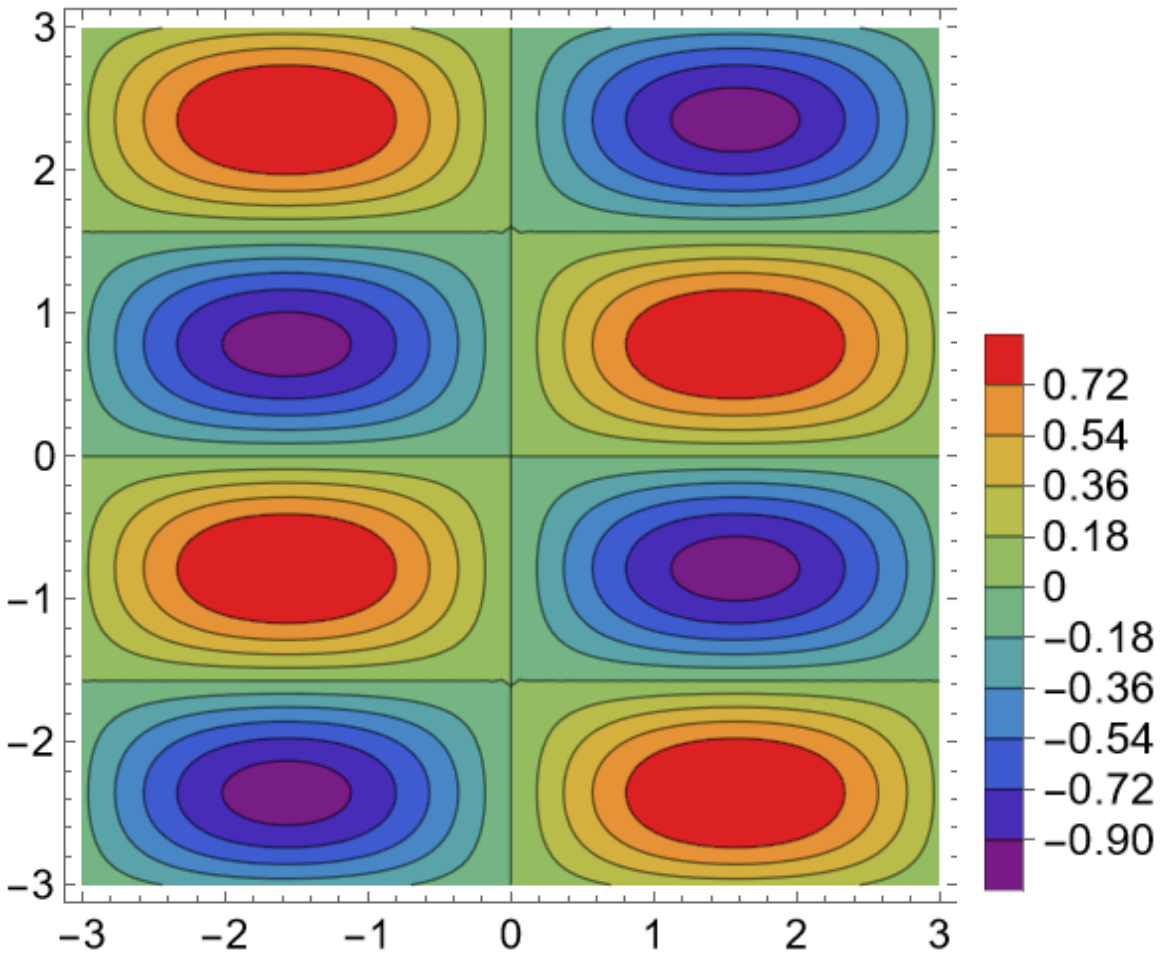

Output {-1.,{x->1.5708,y->2.35619}}

Output {{-1.,{x->1.5708,y->2.35619}},{Steps->3,Function->5,Gradient->5},

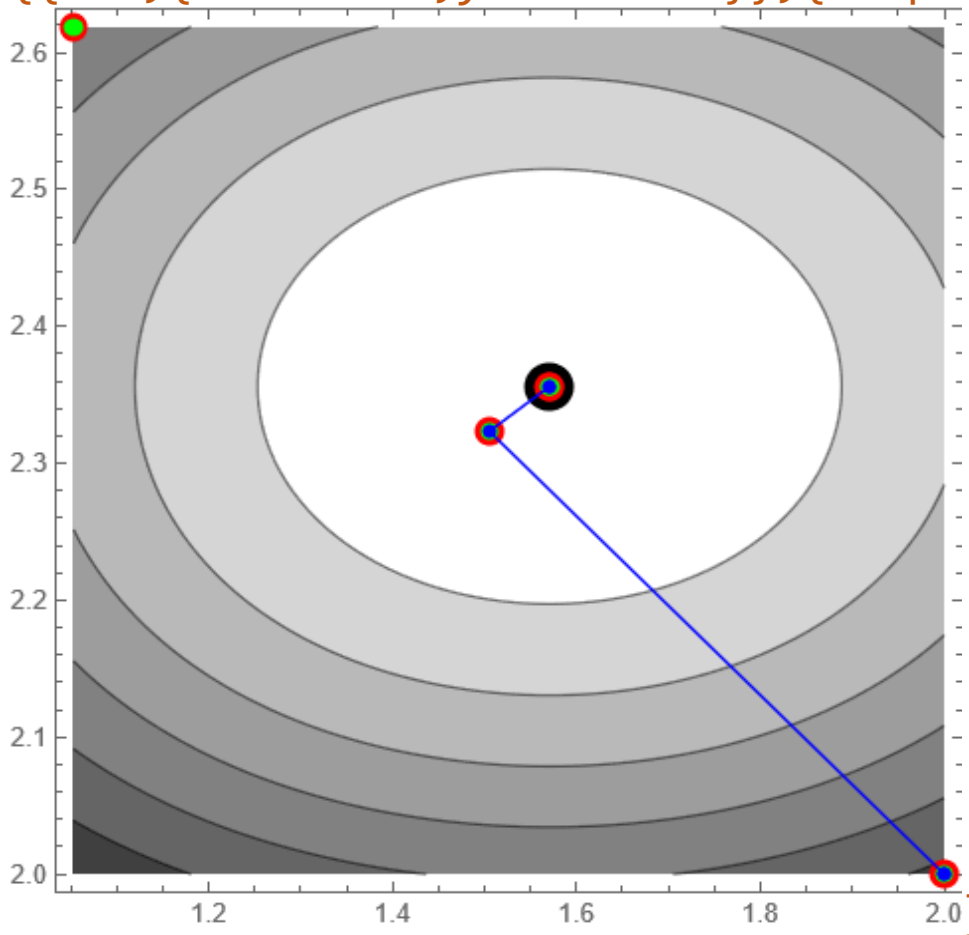

}

Output {{-1.,{x->1.5708,y->2.35619}},{Steps->7,Function->9,Gradient->9},

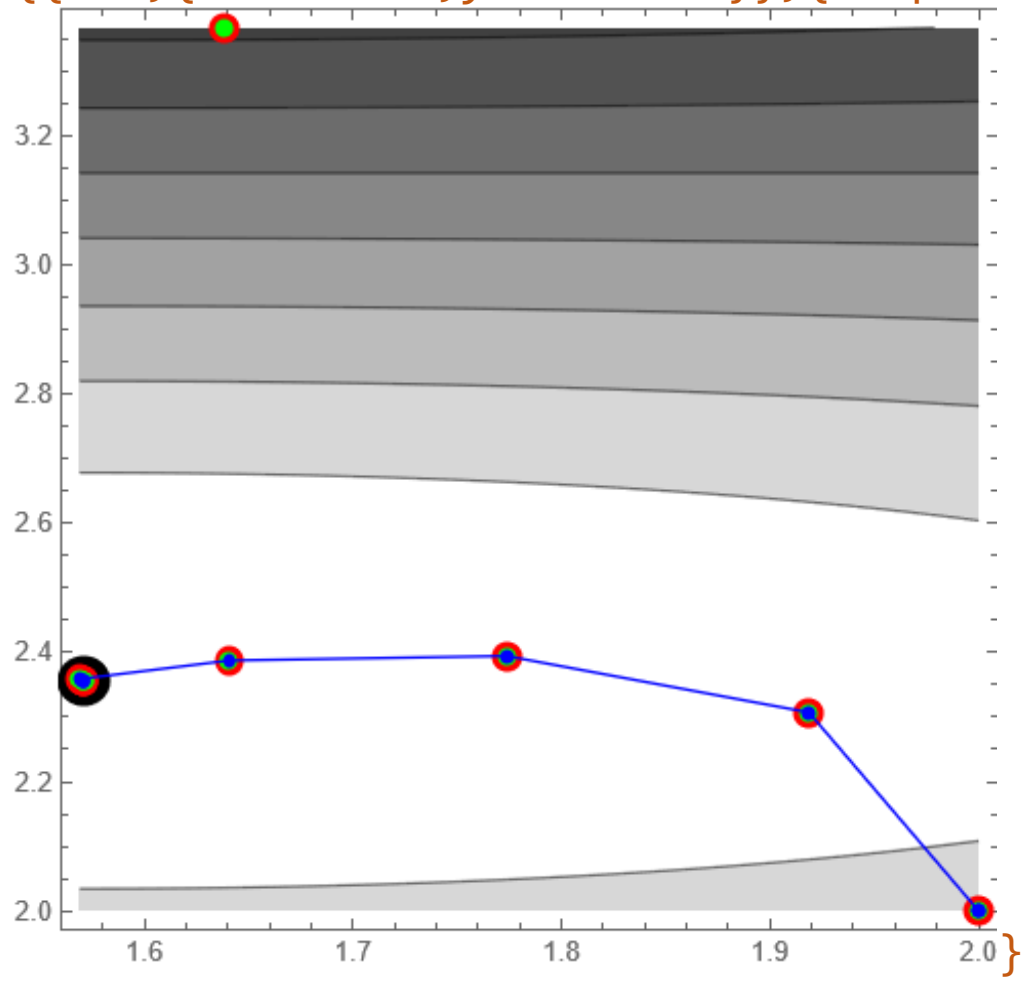

}

Output {{-1.,{x->1.5708,y->2.35619}},{Steps->2,Function->72},

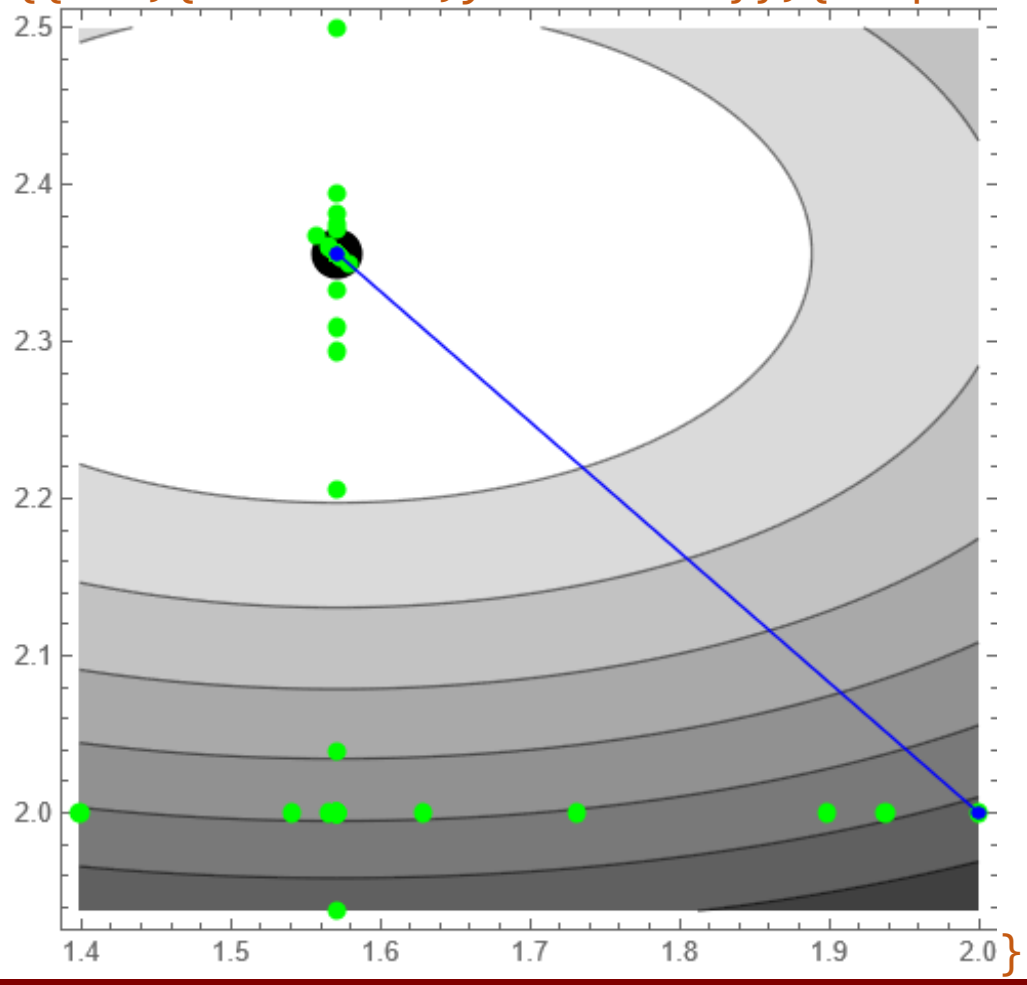

}

***Mathematica Examples 6.4***

Input

```
Plot3D[
 (x y-3)^2+1,
 {x,-4,4},
 {y,-4,4},
 ColorFunction->"Rainbow",
 PlotLegends->BarLegend[Automatic],
 LabelStyle->Directive[Black,16]
 ]

ContourPlot[
 (x y-3)^2+1,
 {x,-4,4},
 {y,-4,4},
 PlotLegends->Automatic,
 Contours->10,
 ColorFunction->"Rainbow",
 LabelStyle->Directive[Black,16]
 ]

Minimize[
 (x y-3)^2+1,
```

```
 {x,y}
 ]

FindMinimumPlot[
 (x y-3)^2+1,
 {{x,-1},{y,-1}},
 Method->"Newton"
 ]
```

Output

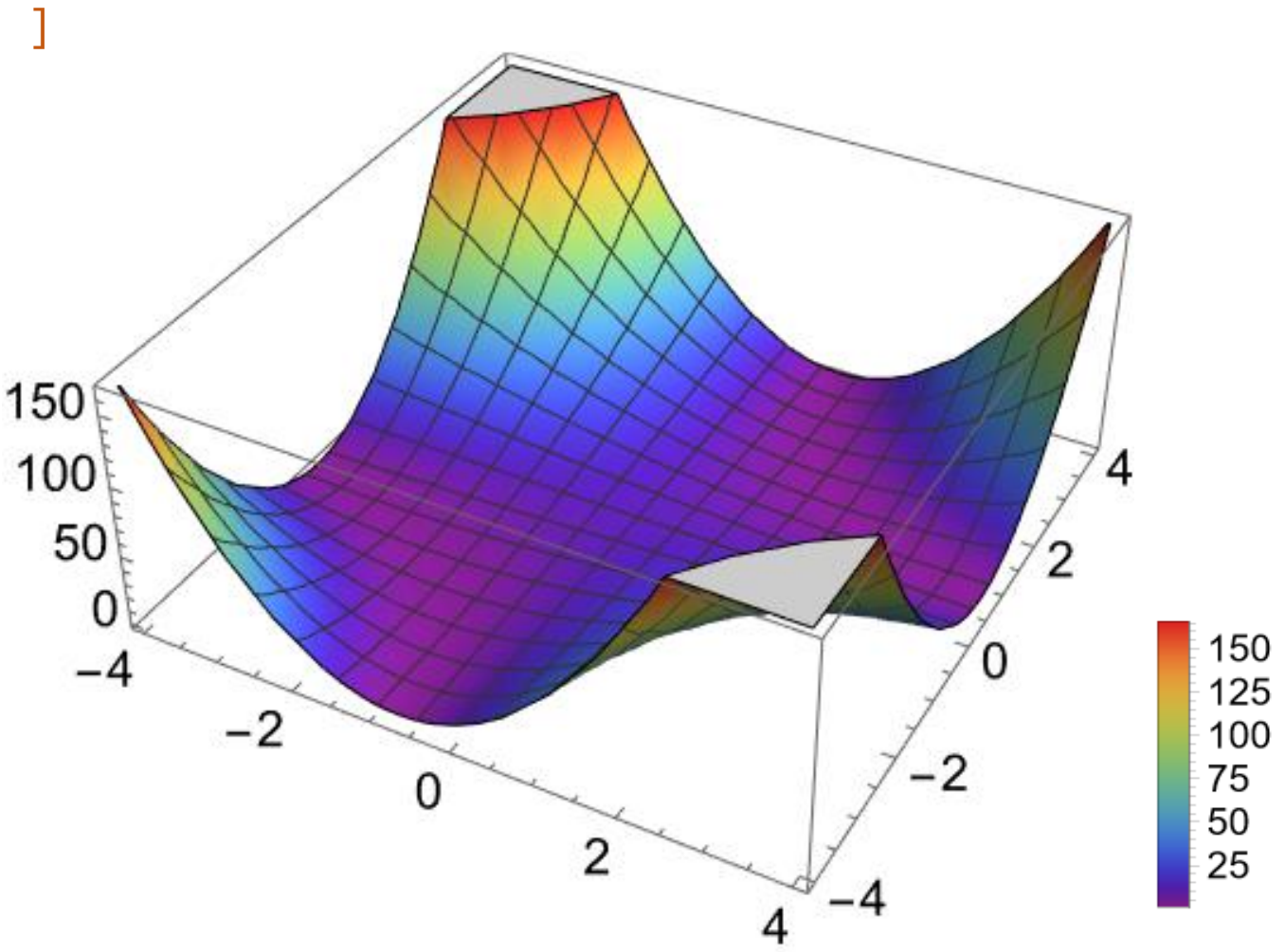

Output

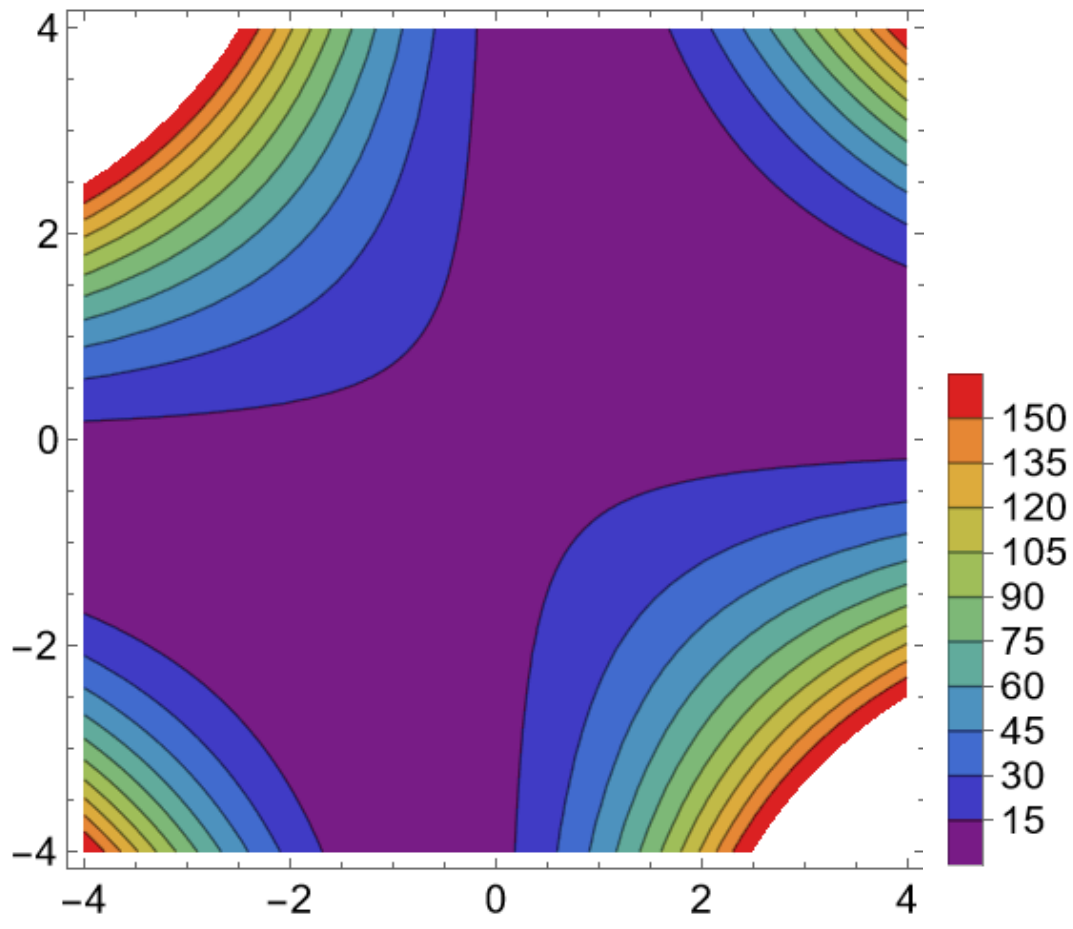

```
Output    {1,{x->-1,y->-3}}
Output    {{1.,{x->-1.39926,y->-2.14398}},{Steps->5,Function->12,Gradient->12},
```

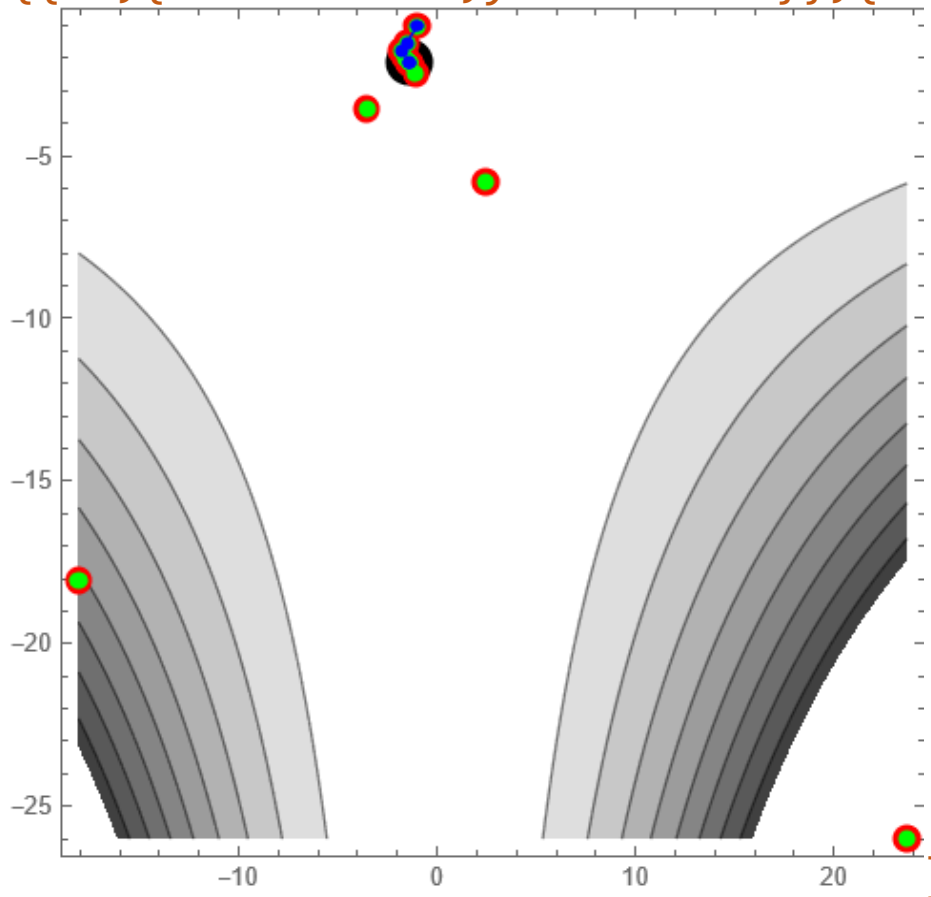

}

***Mathematica Examples 6.5***

Input

```
h=Cos[x^2-3 y]+Sin[x^2+y^2]

Plot3D[
 h,
 {x,-3,3},
 {y,-3,3},
 ColorFunction->"Rainbow",
 PlotLegends->BarLegend[Automatic],
 LabelStyle->Directive[Black,16]
 ]

ContourPlot[
 h,
 {x,-3,3},
 {y,-3,3},
 PlotLegends->Automatic,
 Contours->10,
 ColorFunction->"Rainbow",
 LabelStyle->Directive[Black,16]
 ]

Minimize[
 h,
 {x,y}
 ]

FindMinimum[
 h,
 {{x,1},{y,1}},
 Method->"Newton"
 ]

FindMinimumPlot[
 h,
 {{x,1},{y,1}},
 Method->"Newton"
 ]

FindMinimumPlot[
 h,
 {{x,1},{y,1}},
 Method->"QuasiNewton"
 ]

FindMinimumPlot[
 h,
 {{x,1},{y,1}},
 Method->"ConjugateGradient"
 ]

FindMinimumPlot[
 h,
 {{x,1},{y,1}},
 Method->"PrincipalAxis"
 ]
```

Output

```
Cos[x2-3 y]+Sin[x2+y2]
```

Output 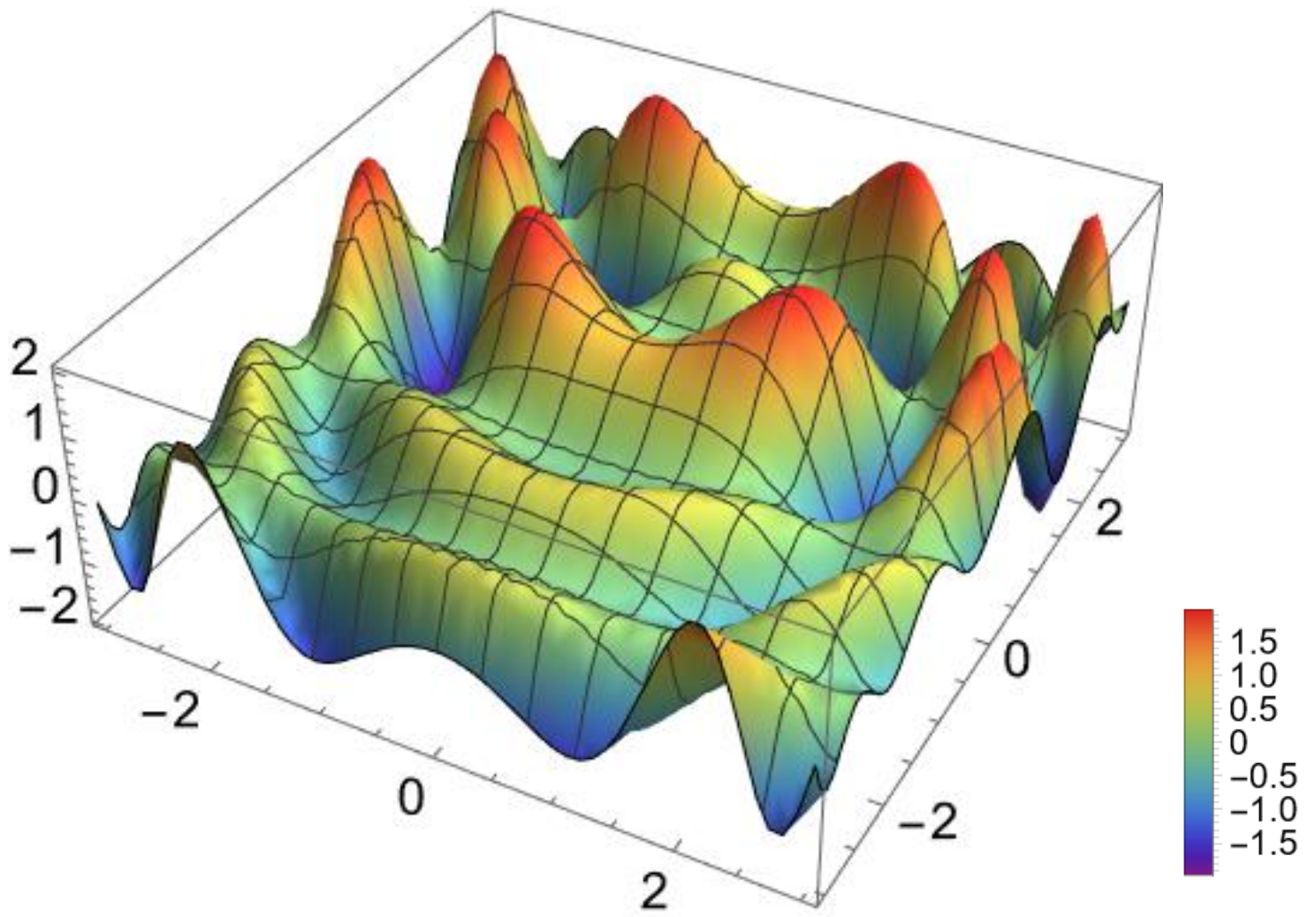

Output 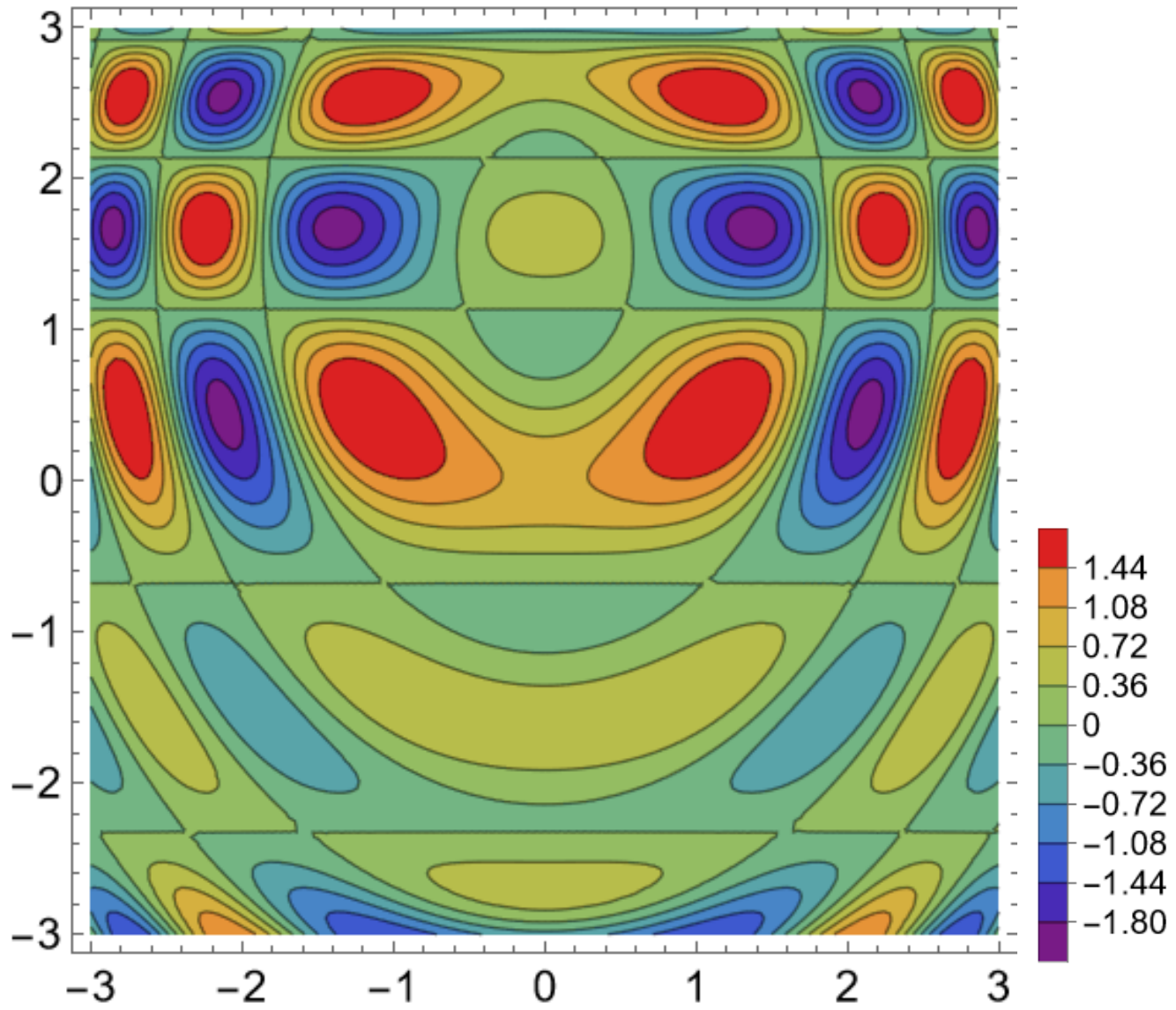

```
Output     {-2,{x->-√(-9-2π+9√(1+82π))/2,y->3/2 (-1+√(1+82π))}}
Output     {-2.,{x->1.37638,y->1.67868}}
Output     {{-2.,{x->1.37638,y->1.67868}},{Steps->5,Function->6,Gradient->6},
```

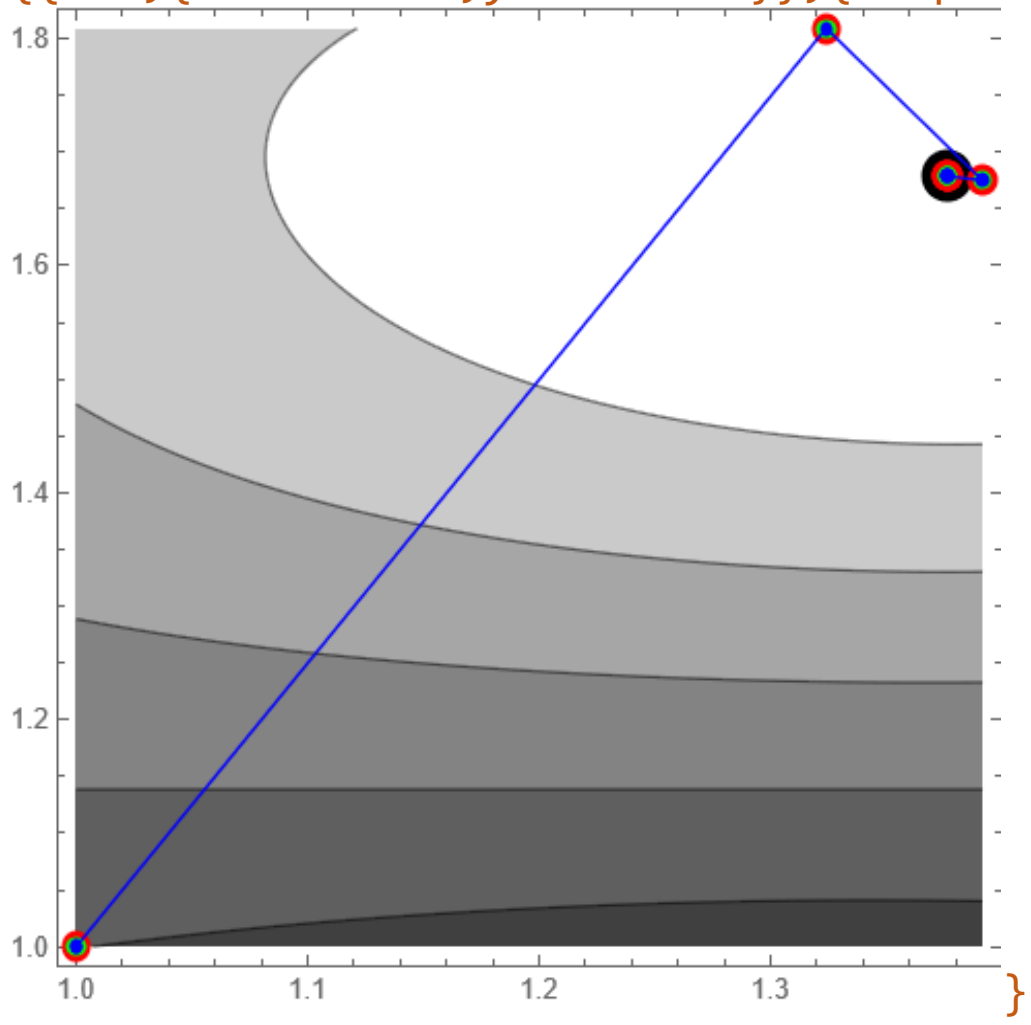

```
}
Output     {{-2.,{x->1.37638,y->1.67868}},{Steps->9,Function->13,Gradient->13},
```

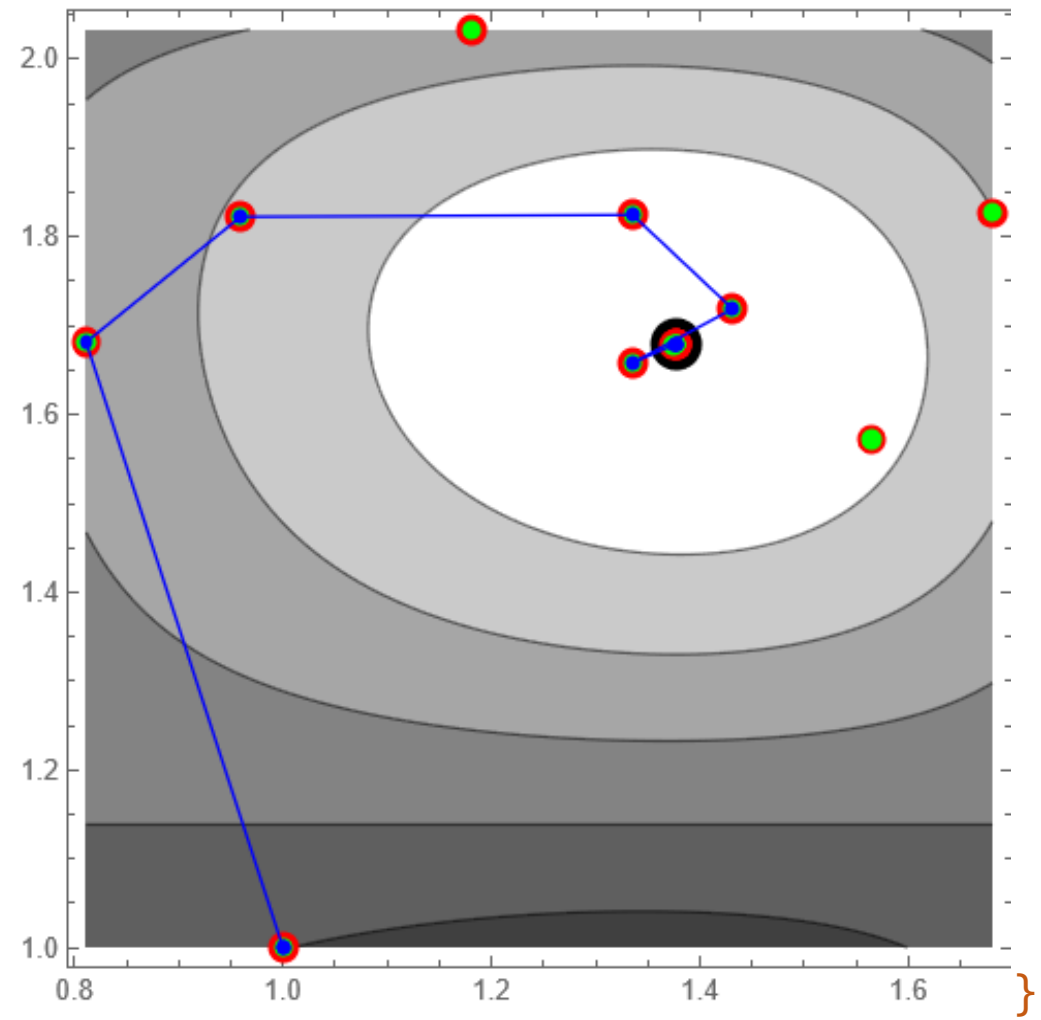
}

Output {{-2.,{x->1.37638,y->1.67868}},{Steps->10,Function->32,Gradient->32},

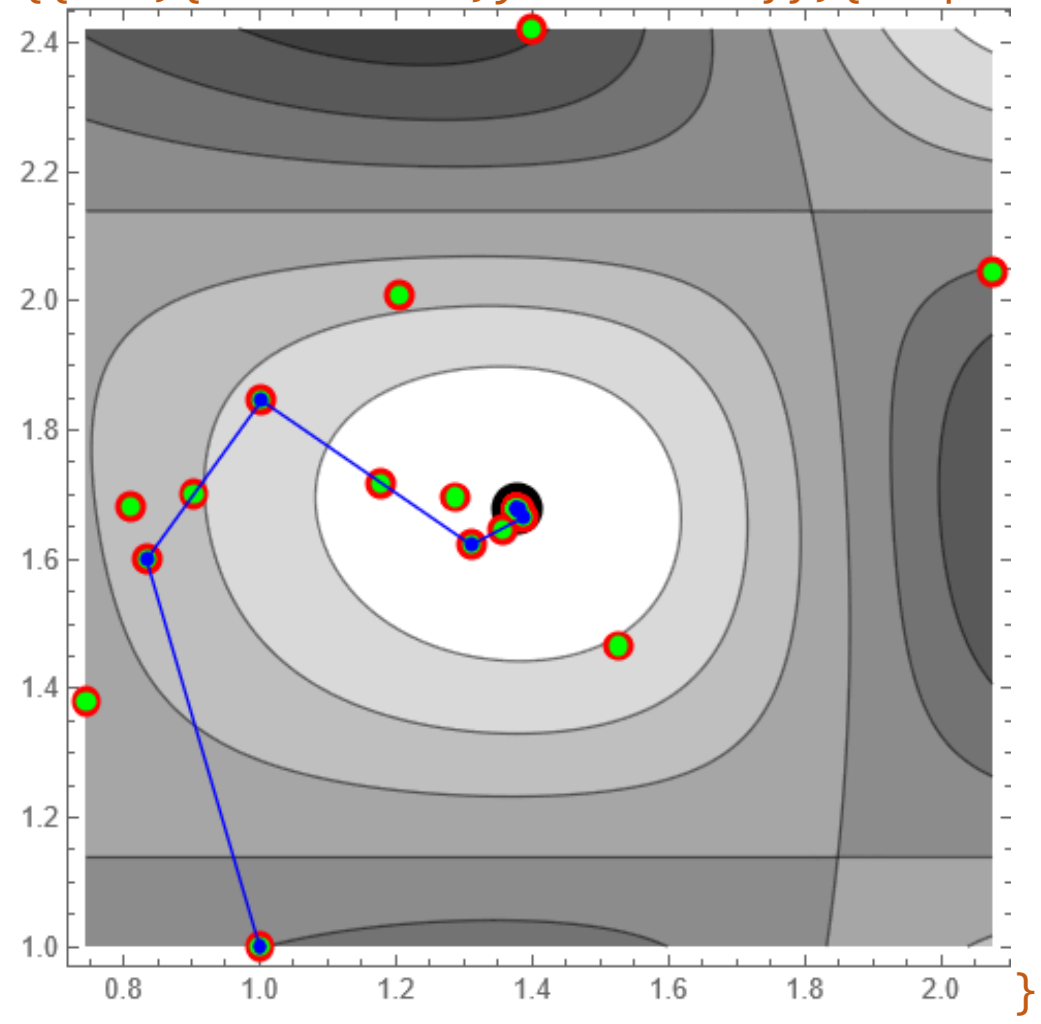
}

Output {{-0.179902,{x->4.94463*10-9,y->0.905726}},{Steps->1,Function->56},

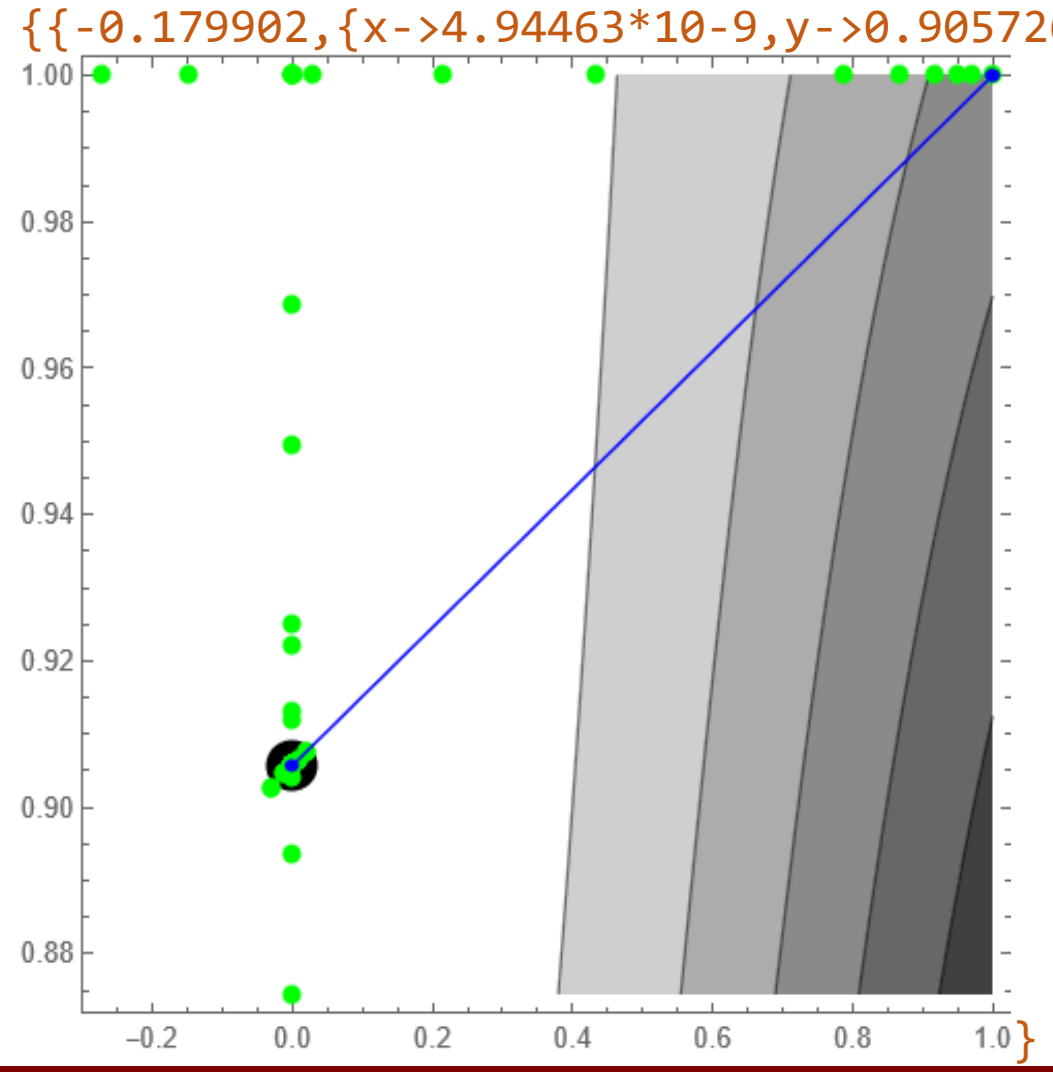
}

*Mathematica Examples 6.6*

Input
```
(* Evaluate a sequence of expressions,"sowing" some to be collected by Reap: *)
Reap[Sow[a];b;Sow[c];Sow[d];e]
```

```
Output    {e,{{a,c,d}}}

Input     (* Compute a sum,"sowing" i^2 at each step: *)
          Reap[
           Sum[
            Sow[i^2]+1,
            {i,10}
            ]
           ]
Output    {395,{{1,4,9,16,25,36,49,64,81,100}}}

Input     (* Use Reap and Sow to collect step data: *)
          f[x_]=Exp[x]+1/x;
          {res,steps}=Reap[
            FindMinimum[
             f[x],
             {x,3},
             StepMonitor:>Sow[{x,f[x]}]
             ]
            ]
          (* Show steps on a plot of the function: *)
          Plot[
           f[x],
           {x,.1,1.5},
           Epilog->{Red,Map[Point,steps[[1]]]},
           LabelStyle->Directive[Black,16]
           ]
Output    {{3.44228,{x-
          >0.703467}},{{{1.5,5.14836},{1.12002,3.95775},{0.633045,3.463},{0.736537,3.44639},
          {0.706923,3.44232},{0.703306,3.44228},{0.703468,3.44228},{0.703467,3.44228}}}}

Output
```

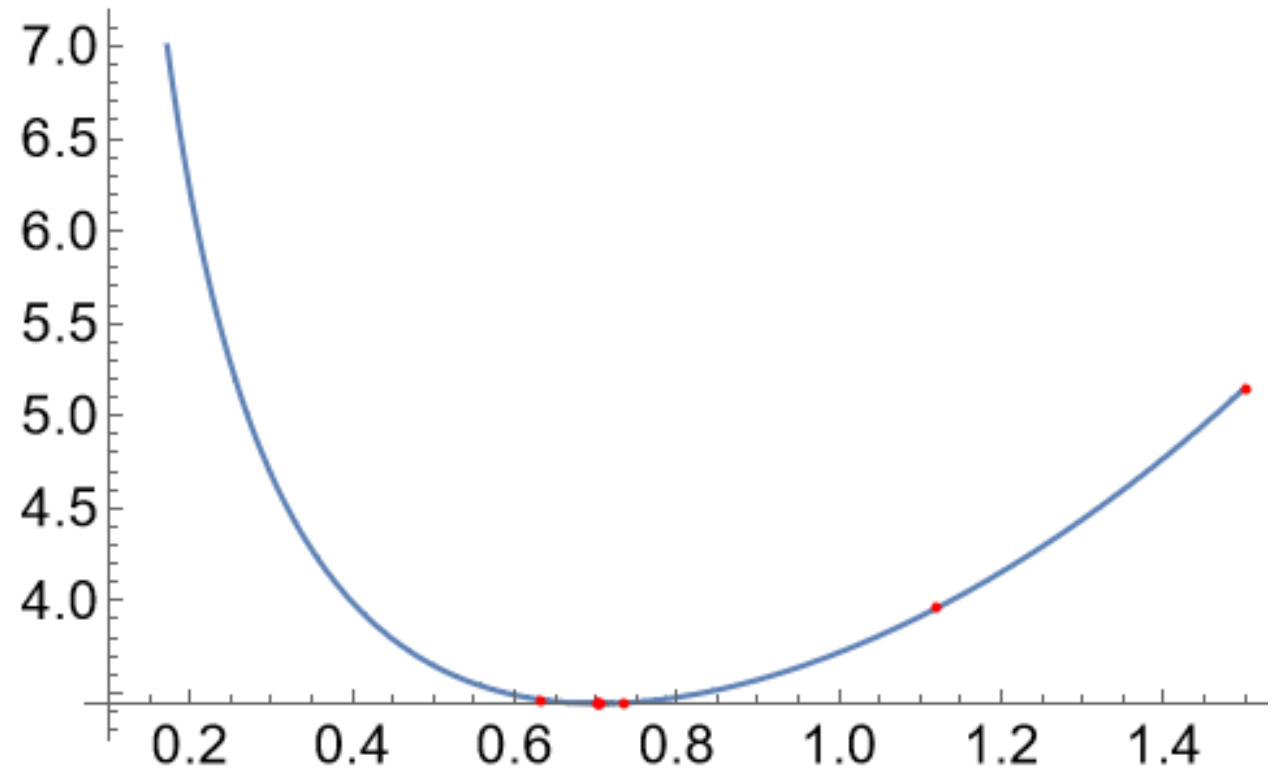

```
Input     (* Investigate steps and evaluations for a numerical minimization: *)
          {sol,data}=Reap[
            FindMinimum[
             (x-1)^2+100 (y-x^2)^2,
             {{x,-1},{y,1}},
             EvaluationMonitor:>Sow[{x,y},"Evaluations"],
             StepMonitor:>Sow[{x,y},"Steps"]
             ],
            _,
            Rule
            ];sol
Output    {0.,{x->1.,y->1.}}

Input     (* Show evaluations in red,steps in yellow,and the final point in green: *)
          ContourPlot[
```

```
 (x-1)^2+100 (y-x^2)^2,
 {x,-1,1},
 {y,-1,1},
 Epilog->{{Green,PointSize[0.05],Point[{x,y}/.
sol[[2]]]},{Yellow,PointSize[0.03],Point["Steps"/.
data]},{Red,PointSize[0.015],Point["Evaluations"/. data]}},
 LabelStyle->Directive[Black,16]
 ]
```

Output

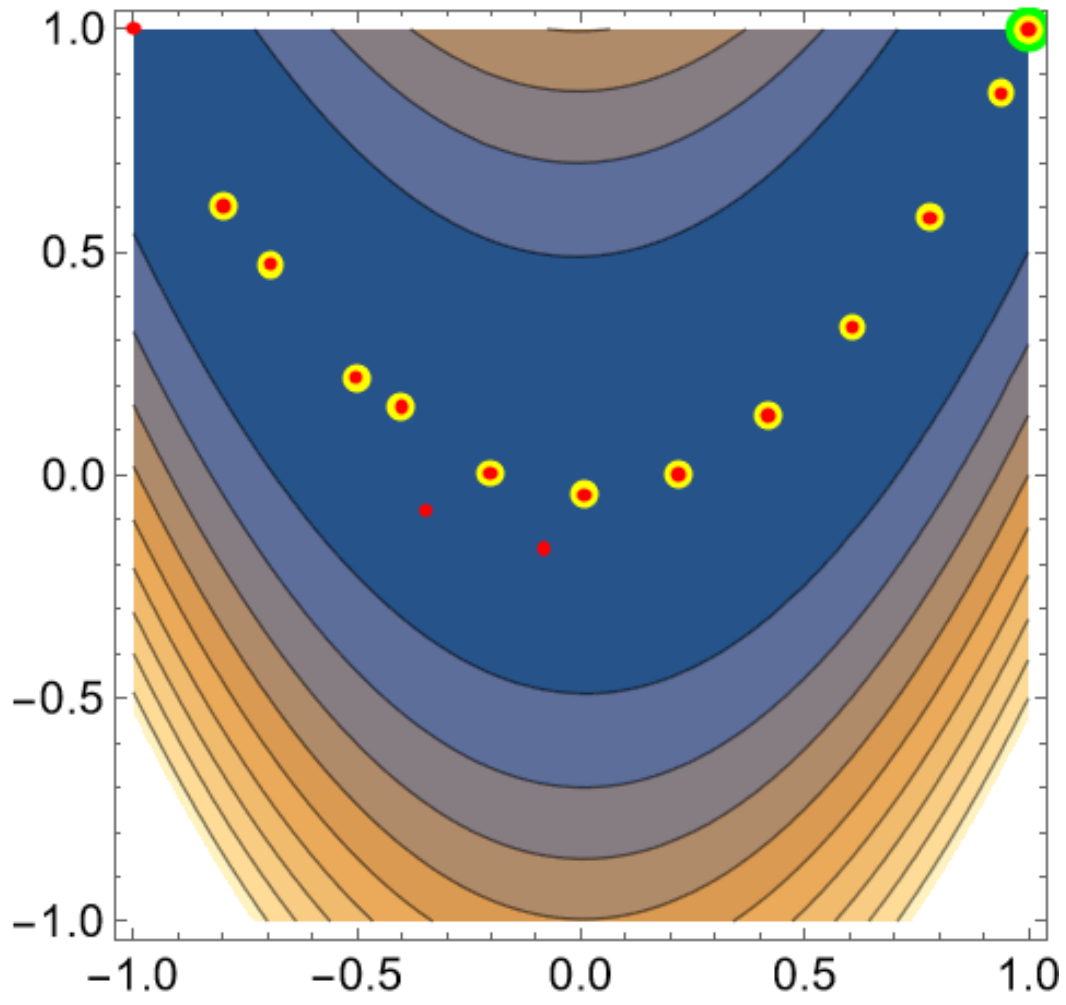

Input

```
(* Steps taken by FindMinimum in finding the minimum of a function: *)
pts=Reap[
    FindMinimum[
     (1-x)^2+100 (-x^2-y)^2+1,
     {{x,-1.2},{y,1}},
     StepMonitor:>Sow[{x,y}]
     ]
    ][[2,1]];
pts=Join[{{-1.2,1}},pts];
ContourPlot[
 (1-x)^2+100 (-x^2-y)^2+1//Log,
 {x,-1.3,1.5},
 {y,-1.5,1.4},
 Epilog->{Red,Line[pts],Point[pts]},
 LabelStyle->Directive[Black,16]
 ]
```

Output

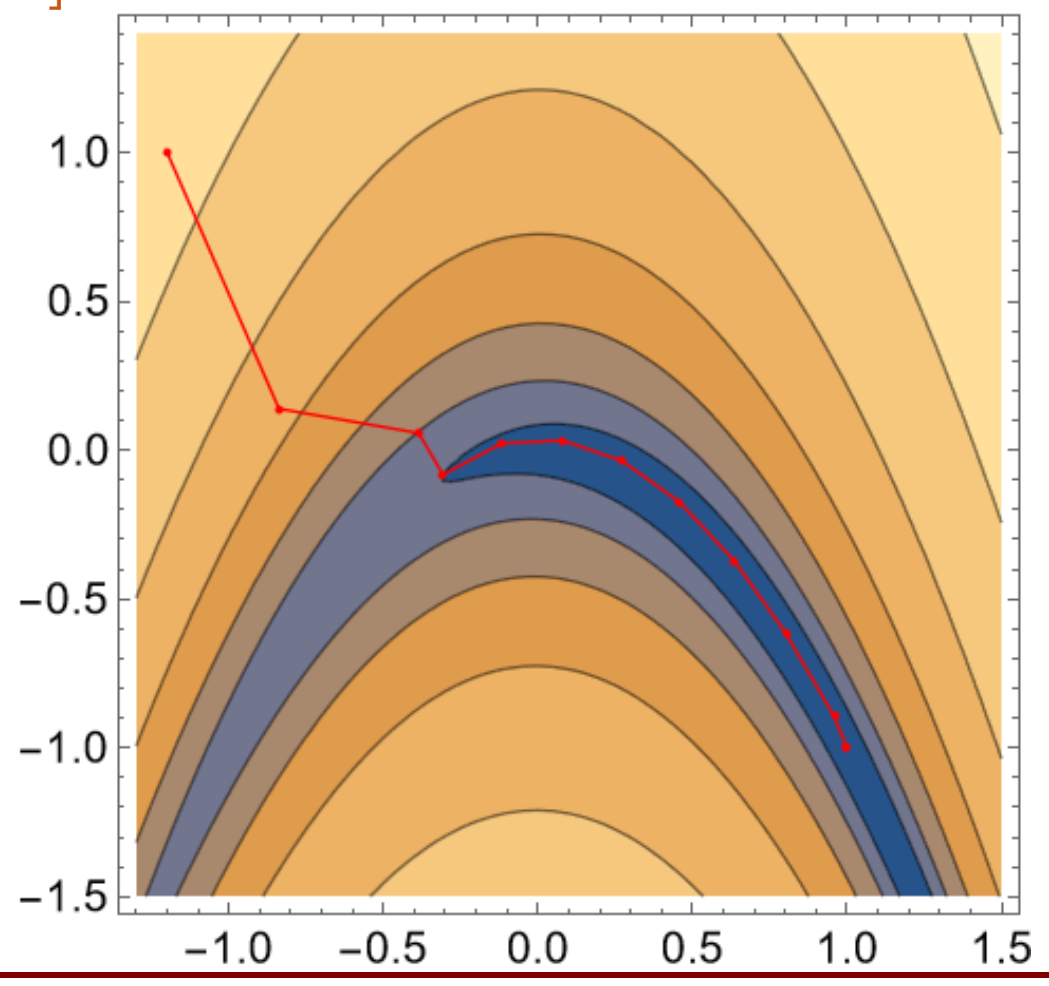

# CHAPTER 7

# MULTI-VARIABLE DIRECT-SEARCH ALGORITHMS

## 7.1 Introduction

In this section, we examine methods or algorithms that iteratively produce estimates of $|\mathbf{x}^*\rangle$, that set of design variables that causes $f|\mathbf{x}\rangle$ to take on its minimum value. All methods lie on several common principles:

- Selection of a basis point: a set of values of the variables must be found that belong to the domain.
- Calculation of the objective function at this point.
- Choice of a second feasible point according to a given method.
- Calculation of the objective function at this point.
- Comparison of the value of the objective function at the second point to that at the basis point.
- If the second point is better, it constitutes the new basis. If the initial point is better, modify the search direction, strategy, or stop.

The methods can be classified into three categories:

1. Direct-search methods, which use only function values.
2. Gradient methods, which require estimates of the first derivative of $f|\mathbf{x}\rangle$).
3. Second-order methods, which require estimates of the first and second derivatives of $f|\mathbf{x}\rangle$.

Direct-search methods are not very efficient. As a result, their application is restricted to problems where gradient information is difficult to obtain, for example, if the objective function is not continuous. Gradient methods are based on gradient information. First-order methods are based on the linear approximation of the Taylor series, and hence they entail the gradient **g**. On the other hand, second-order methods are based on the quadratic approximation of the Taylor series. They entail the gradient **g** as well as the Hessian **H**. We will consider examples from each class since no one class of methods can be expected to solve all problems with equal efficiency uniformly. In this chapter, we consider the following algorithms

**Direct-Search Methods**

1. Box evolutionary method.
2. Simplex search method.
3. Hooke-Jeeves pattern search method.
4. Powell conjugate direction method.
5. Cyclic coordinate search

**First-order Approximations Methods (Gradient-Based Methods)**

1. Steepest descent method.
2. Steepest descent without line search (formula 1 and 2).
3. Steepest descent without line search (Barzilai–Borwein two-point formula).

**Second-order Approximations Methods (Hessian-Based Methods)**

1. Newton method.
2. Modified Newton method.
3. Gauss-Newton method.

4. Marquardt method.

**Conjugate-Direction Methods**

1. Conjugate-Gradient method.
2. Hestenes-Stiefel, Polak-Ribiere, and Fletcher-Reeves methods.

**Quasi-Newton Methods**

1. Rank one correction formula.
2. The Davidon-Fletcher-Powell (DFP) and Broyden-Fletcher-Goldfarb-Shanno (BFGS) methods.

## 7.2 Artificial Landscapes

The test functions, known as artificial landscapes, are useful for evaluating characteristics of optimization algorithms, such as convergence rate, precision, robustness, and general performance. In Table 7.1, the formula in 2D, global minimum, and search domain of 8 well-known functions are presented. The contour plots of the functions are shown in Figure 7.1.

**Table 7.1.**

| | | |
|---|---|---|
| Sphere function | Formula in 2D | $f\,(x,y) = x^2 + y^2$ |
| | Global minimum | $f(0,0) = 0$ |
| | Search domain | $x, y \in [-\infty, \infty]$ |
| Rosenbrock function | Formula in 2D | $f\,(x,y) = 100(y - x^2)^2 + (1 - x)^2$ |
| | Global minimum | $f(1,1) = 0$ |
| | Search domain | $x, y \in [-\infty, \infty]$ |
| Beale function | Formula in 2D | $f\,(x,y) = (1.5 - x + xy)^2 + (2.25 - x + xy^2)^2 + (2.625 - x + xy^3)^2$ |
| | Global minimum | $f(3,0.5) = 0$ |
| | Search domain | $x, y \in [-4.5, 4.5]$ |
| Booth function | Formula in 2D | $f(x,y) = (x + 2y - 7)^2 + (2x + y - 5)^2$ |
| | Global minimum | $f(1,3) = 0$ |
| | Search domain | $x, y \in [-10, 10]$ |
| Matyas function | Formula in 2D | $f(x,y) = 0.26(x^2 + y^2) - 0.48xy$ |
| | Global minimum | $f(0,0) = 0$ |
| | Search domain | $x, y \in [-10, 10]$ |
| Himmelblau function | Formula in 2D | $f(x,y) = (x^2 + y - 11)^2 + (x + y^2 - 7)^2$ |
| | Global minimum | $f(3,2) = f(-2.805, 3.131) = f(-3.779, -3.283) = f(3.584, -1.848) = 0$ |
| | Search domain | $x, y \in [-5, 5]$ |
| Three-hump camel function | Formula in 2D | $f(x,y) = 2x^2 - 1.05x^4 + x^6/6 + xy + y^2$ |
| | Global minimum | $f(0,0) = 0$ |
| | Search domain | $x, y \in [-5, 5]$ |
| Styblinski-Tang function | Formula in 2D | $f(x,y) = (x^4 - 16x^2 + 5x)/2 + (y^4 - 16y^2 + 5y)/2$ |
| | Global minimum | $-39.16617 < f(-2.9035, -2.9035) < -39.16616$ |
| | Search domain | $x, y \in [-5, 5]$ |

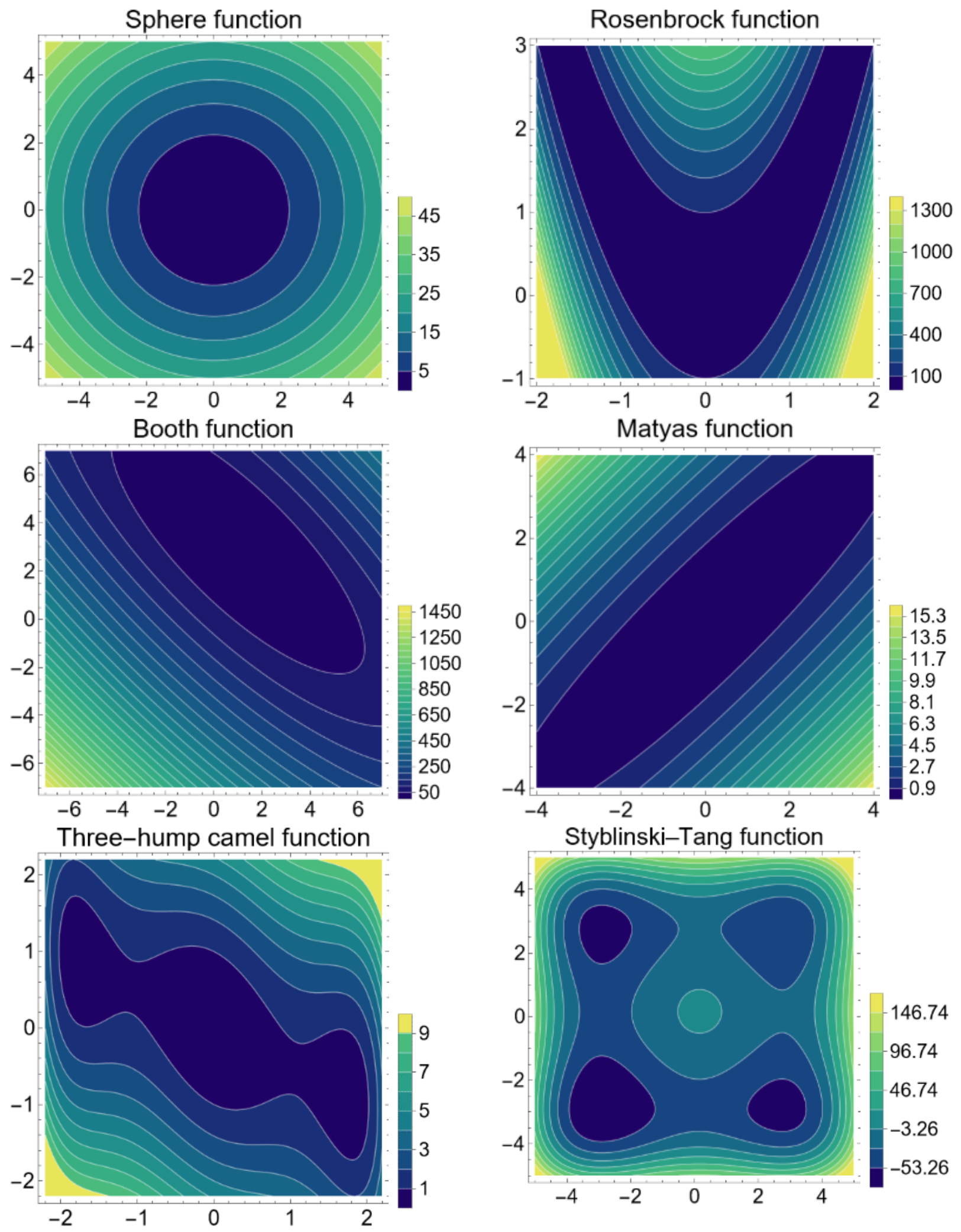

**Figure 7.1.** The contour plots some test functions

The results were produced by Mathematica code 7.1.

***Mathematica Code 7.1*** Contour plot of the Styblinski–Tang function

```
(* Styblinski-Tang function *)

ContourPlot[
 (x^4-16x^2+5x)/2+(y^4-16y^2+5y)/2,
 {x,-5,5},
 {y,-5,5},
 Contours->Function[{min,max},Range[min,max,25]],
 ContourStyle->{White},
 ClippingStyle->Automatic,
 ColorFunction->"BlueGreenYellow",
 PlotLegends->Automatic,
 LabelStyle->Directive[Black,20],
 PlotLabel->"Styblinski-Tang function"
 ]
```

## 7.3 Unidirectional Search

One of the very basic search techniques for multi-variable optimization is a unidirectional search method. The idea of this method is to convert the multi-variable function to a single variable and then investigate the optimum. Many multi-variable optimization techniques use successive unidirectional search techniques to find the minimum point along a particular search direction. Since unidirectional searches will be mentioned in the different algorithms, we illustrate here how a unidirectional search can be performed on a multi-variable function. A unidirectional search is a one-dimensional search performed by comparing function values only along a specific direction. Usually, a unidirectional search is performed from a point $|\mathbf{x}_t\rangle$ and in a specified direction $|\mathbf{s}_t\rangle$. That is, only points that lie on a line (in an N-dimensional space) pass through the point $|\mathbf{x}_t\rangle$ and oriented along the search direction $|\mathbf{s}_t\rangle$ are allowed to be considered in the search process. Any arbitrary point on that line can be expressed as follows:

$$|\mathbf{x}(\alpha)\rangle = |\mathbf{x}_t\rangle + \alpha\big|\mathbf{s}^{(t)}\big\rangle. \tag{7.1}$$

The parameter $\alpha$ is a scalar quantity, specifying a relative measure of the distance of the point $|\mathbf{x}(\alpha)\rangle$ from $|\mathbf{x}_t\rangle$. Note, however, that the above equation is a vector equation specifying all design variables $x_i(\alpha)$. In order to find the minimum point on the specified line, we can rewrite the multi-variable objective function in terms of a single variable $\alpha$ by substituting each variable $x_i$ by the expression $x_i(\alpha)$ given in (7.1) and by using a suitable single-variable search method described in Chapter 5. Once the optimal value $\alpha^*$ is found, the corresponding point can also be found using (7.1). In the following code, Mathematica function, Mathematica Code 7.2, starts by bracketing the minimum of a single variable function $f|\mathbf{x}(\alpha)\rangle$, using bounding phase method, then isolating the minimum using the golden section search method.

***Mathematica Code 7.2*** Unidirectional Search

```
unidirectionalsearch[α_,delt_,eps_]:=Module[

{α0=α,delta=delt,epsilon=eps,y1,y2,y3,αα,a,b,increment,a0,b0,anew,bnew,anorm,bnorm,lnorm,α1n
orm,α2norm,α1,α2,f1,f2,αstar},

   (* Bounding Phase Method *)

   (* Initiating required variables *)
   y1 = f[α0-Abs[delta]];
   y2 = f[α0];
   y3 = f[α0+Abs[delta]];

   (*Determining whether the inicrement is positive or negative*)
   Which[
    y1==y2,
    a=α0-Abs[delta];
    b=α0;
    Goto[end];,
    y2==y3,
    a=α0;
    b=α0+Abs[delta];
    Goto[end];,
    y1==y3||(y1>y2&&y2<y3),
    a=α0-Abs[delta];
    b=α0+Abs[delta];
    Goto[end];
    ];

   Which[
     y1>y2&&y2>y3,
     increment=Abs[delta];,
```

```
      y1<y2&&y2<y3,
      increment=-Abs[delta];
      ]

     (* Starting the algorithm *)
     Do[
      αα[0]=α0;
      αα[k+1]=αα[k]+2^k*increment;

      Which[
       f[αα[k]]<f[αα[k+1]],(* Evidently, it is impossible the condition to hold for k=0 *)
       a=αα[k-1];
       b= αα[k+1];
       Break [],

       k>50,
       Print["After 50 iterations the bounding phase method can not braketing the min of
alpha"];
       Exit[]
       ];,
      {k,0,∞}
      ];

    Label[end];

    If[
     a>b,
     {a,b}={b,a}
     ];

    (* Golden Section Search Method *)

    (* Initiating required variables*)
    a0=a;
    b0=b;
    anew=a;
    bnew=b;

    If[
     a0==b0,
     αstar=a;
     Goto[final]
     ];

    (* Starting the algorithm *)
    Do[
     (* Normalize the variable α *)
     anorm=(anew-a)/(b-a);
     bnorm=(bnew-a)/(b-a);

     lnorm=bnorm-anorm;

     α1norm=anorm+0.382*lnorm;
     α2norm=bnorm-0.382*lnorm;

     α1=α1norm(b0-a0)+a0;
     α2=α2norm(b0-a0)+a0;

     f1=f[α1];
     f2=f[α2];
```

```
     Which [
      f1>f2,
      anew=α1(*move lower bound to α1*);,
      f1<f2,
      bnew=α2(*move upper bound to α2*);,
      f1==f2,
      anew=α1(*move lower bound to α1*);
      bnew=α2(*move upper bound to α2*);
      ];

     αstar=0.5*(anew+bnew);

     If[
      Abs[lnorm]<epsilon,
      Break[]
      ];,
     {k,1,∞}
     ];

    Label[final];

    (* Final result *)
    N[αstar]
    ]

(*Taking the function as input from user*)
f[α_] = Evaluate[Input["Please input a function of x to find the minimum "]];(* the user
types in for instance x^2 *)

 unidirectionalsearch[-6,1,0.00001]
 -6.15674*10^-6
```

## 7.4 Direct Search Methods

In this section, we present a few minimization algorithms that use function values only. Direct methods are also called zero-order, black box, pattern search, or derivative-free methods [1-8]. We assume here that $f|\mathbf{x}\rangle$ is continuous, and $\nabla f(\mathbf{x})$ may or may not exist but certainly is not available. Also, we assume that $f|\mathbf{x}\rangle$ has a single minimum in the domain of interest.

In a single-variable function optimization, there are only two search directions in a point that can be modified—either in the positive $x$-direction or the negative $x$-direction. The extent of increment or decrement in each direction depends on the current point and the objective function. In multi-variable optimization, each variable can be modified either in the positive or in the negative direction, thereby totaling $2^N$ different ways, where $N$ is the number of variables. Typically, the domain of the performance index (variables) will be subdivided into a grid of points, and then various strategies for shrinking the area in which the solution lay will be applied. Often this led to nearly exhaustive enumeration and hence proved unsatisfactory for all problems. The more useful idea was to pick a base point and, rather than attempting to cover the entire range of the variables, evaluate the performance index in some pattern about the base point.

### 7.4.1 Box Evolutionary Method

The basic idea of this method is to construct a $N$-dimensional hypercube (called a factorial design pattern) where the center point is the current point. The Box algorithm requires $(2^N + 1)$ points, of which $2^N$ are corner points of an $N$-dimensional hypercube centered on the other point. For all $2^N + 1$ points, function values are obtained and compared. Then, the best point is identified. In the next iteration, another hypercube is formed around this best point. If at any iteration, an improved point is not found, the size of the hypercube is reduced. This process continues until the hypercube becomes very small [1-5].

For example, in two dimensions, a square pattern, such as in Figure 7.2, is located. Then the best of the five points is selected as the next base point around which to locate the next pattern of points. If none of the corner points is better than the base point, the scale of the grid is reduced, and the search continues.

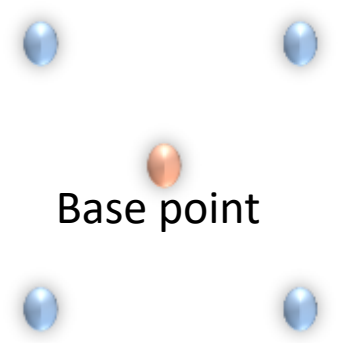

**Figure 7.2.** Factorial design pattern in 2D.

**Algorithm**

**Step 1:** Choose an initial point $\mathbf{x}^{(0)}$ and size reduction parameters $\Delta_i$ for all design variables, $i = 1,2,\ldots,N$. Choose a termination parameter $\epsilon$. Set $|\bar{\mathbf{x}}\rangle = |\mathbf{x}_0\rangle$.

**Step 2:** If $||\Delta|| < \epsilon$, Terminate;
Else create $2^N$ points by adding and subtracting $\Delta_i/2$ from each variable at the point $\bar{\mathbf{x}}$.

**Step 3:** Compute function values at all $(2^N + 1)$ points.
Find the point having the minimum function value.
Designate the minimum point to be $|\bar{\mathbf{x}}\rangle$.

**Step 4:** If $|\bar{\mathbf{x}}\rangle = |\mathbf{x}_0\rangle$, reduce size parameters $\Delta_i = \Delta_i/2$ and go to Step 2;
Else set $|\mathbf{x}_0\rangle = |\bar{\mathbf{x}}\rangle$ and go to Step 2.

## Remarks

- It is evident from the algorithm that at most $2^N$ functions are evaluated at each iteration. Thus, the required number of function evaluations increases exponentially with $N$.
- The algorithm convergence depends on the initial hypercube size, location, and the chosen size reduction parameter $\Delta_i$.
- Starting with a large $\Delta_i$ is good, but the convergence to the minimum may require more iterations and hence more function evaluations.
- On the other hand, starting with a small hypercube may lead to premature convergence on a suboptimal point, especially in the case of highly nonlinear functions.
- It is important to note that the reduction of the size parameter $(\Delta_i)$ by a factor of two in Step 4 is not always necessary. A smaller or larger reduction can be used. However, a smaller reduction (a factor smaller than two and greater than one) is usually recommended for better convergence.

***Example 7.1***

Consider the problem:

$$\text{Minimize } f|\mathbf{x}\rangle = x^2 + y^2,$$

using, $|\mathbf{x}_0\rangle = (5,5)$, $\Delta_i = (1,1)$, and $\epsilon = 0.1$.

***Solution***

The 3D and contour plots of the function are shown in Figure 7.3. The plots show that the minimum lies at $|\mathbf{x}^*\rangle = (0,0)^T$, $f|\mathbf{x}^*\rangle = 0$.

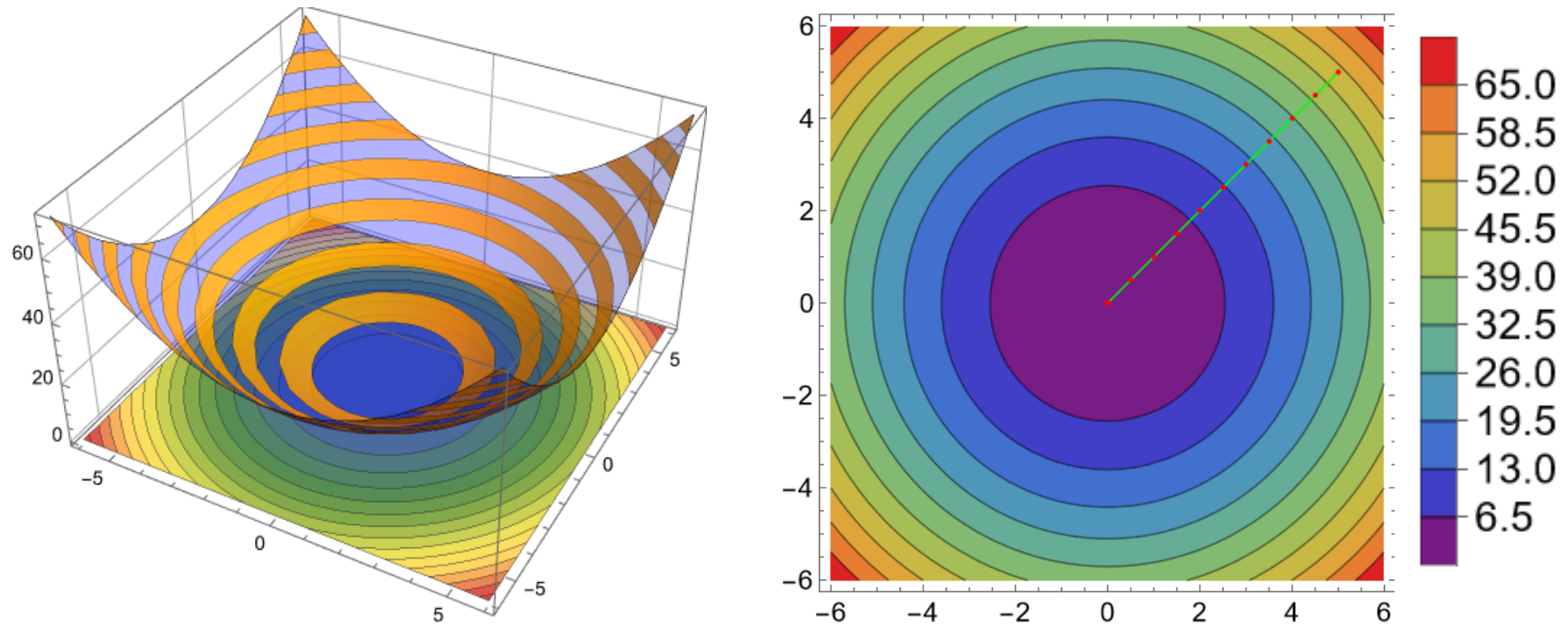

**Figure 7.3.** The results of 14 iterations of the Box evolutionary method for $f|\mathbf{x}\rangle = x^2 + y^2$.

From Table 7.2, after 14 iterations, the optimal condition is attained. The results were produced by Mathematica code 7.3.

**Table 7.2.a.**

| No. of iters. | corner4 | Corner1 | x0 | Corner2 | Corner3 |
|---|---|---|---|---|---|
| 1 | {4.5, 5.5} | {4.5, 4.5} | {5., 5.} | {5.5, 4.5} | {5.5, 5.5} |
| 2 | {4., 5.} | {4., 4.} | {4.5, 4.5} | {5., 4.} | {5., 5.} |
| 3 | {3.5, 4.5} | {3.5, 3.5} | {4., 4.} | {4.5, 3.5} | {4.5, 4.5} |
| 4 | {3., 4.} | {3., 3.} | {3.5, 3.5} | {4., 3.} | {4., 4.} |
| 5 | {2.5, 3.5} | {2.5, 2.5} | {3., 3.} | {3.5, 2.5} | {3.5, 3.5} |
| 6 | {2., 3.} | {2., 2.} | {2.5, 2.5} | {3., 2.} | {3., 3.} |
| 7 | {1.5, 2.5} | {1.5, 1.5} | {2., 2.} | {2.5, 1.5} | {2.5, 2.5} |
| 8 | {1., 2.} | {1., 1.} | {1.5, 1.5} | {2., 1.} | {2., 2.} |
| 9 | {0.5, 1.5} | {0.5, 0.5} | {1., 1.} | {1.5, 0.5} | {1.5, 1.5} |
| 10 | {0., 1.} | {0., 0.} | {0.5, 0.5} | {1., 0.} | {1., 1.} |
| 11 | {-0.5, 0.5} | {-0.5, -0.5} | {0., 0.} | {0.5, -0.5} | {0.5, 0.5} |
| 12 | {-0.25, 0.25} | {-0.25, -0.25} | {0., 0.} | {0.25, -0.25} | {0.25, 0.25} |
| 13 | {-0.125, 0.125} | {-0.125, -0.125} | {0., 0.} | {0.125, -0.125} | {0.125, 0.125} |
| 14 | {-0.0625, 0.0625} | {-0.0625, -0.0625} | {0., 0.} | {0.0625, -0.0625} | {0.0625, 0.0625} |

**Table 7.2.b.**

| No. of iters. | f(corner4) | f(corner1) | f(x0) | f(corner2) | f(corner3) | xstar |
|---|---|---|---|---|---|---|
| 1 | 50.5 | 40.5 | 50 | 50.5 | 60.5 | {4.5, 4.5} |
| 2 | 41 | 32 | 40.5 | 41 | 50 | {4., 4.} |
| 3 | 32.5 | 24.5 | 32 | 32.5 | 40.5 | {3.5, 3.5} |
| 4 | 25 | 18 | 24.5 | 25 | 32 | {3., 3.} |
| 5 | 18.5 | 12.5 | 18 | 18.5 | 24.5 | {2.5, 2.5} |
| 6 | 13 | 8 | 12.5 | 13 | 18 | {2., 2.} |
| 7 | 8.5 | 4.5 | 8 | 8.5 | 12.5 | {1.5, 1.5} |
| 8 | 5 | 2 | 4.5 | 5 | 8 | {1., 1.} |
| 9 | 2.5 | 0.5 | 2 | 2.5 | 4.5 | {0.5, 0.5} |
| 10 | 1 | 0 | 0.5 | 1 | 2 | {0., 0.} |
| 11 | 0.5 | 0.5 | 0 | 0.5 | 0.5 | {0., 0.} |
| 12 | 0.125 | 0.125 | 0 | 0.125 | 0.125 | {0., 0.} |
| 13 | 0.03125 | 0.03125 | 0 | 0.03125 | 0.03125 | {0., 0.} |
| 14 | 0.007813 | 0.007813 | 0 | 0.007813 | 0.007813 | {0., 0.} |

***Mathematica Code 7.3*** Box Evolutionary Search Method

```
(* Box's Evolutionary Search Method *)

(*
Notations
x0          :Intial vector
deltai       :Size reduction parameter delta
epsilon      :Small number to check the accuracy of the box method
f[x,y]       :Objective function
xbar[k]      :Final design solution of iteration k
f[xbar[k]] :Final objective function value of iteration k
lii          :The last iteration index
result[k]   :The results of iteration k
*)

(* Taking Initial Inputs from User *)
x0=Input["Enter the intial point in the format {x, y}; for example {5,5} "] ;
deltai=Input["Enter the size reduction parameter delta in the format {Δx,Δy}; for example
{1,1}"];
epsilon=Input["Please enter accuracy of the box method: for example 0.1 "];

domainx=Input["Please enter domain of x variable for 3D and contour plots; for example {-
6,6}"];
domainy=Input["Please enter domain of y variable for 3D and contour plots; for example {-
6,6}"];

If[
  epsilon<=0||deltai[[1]]<=0||deltai[[2]]<=0,
  Beep[];
  MessageDialog["Tolerance value has to be small postive number: "];
  Exit[];
  ];

dd=deltai;
xx=x0;

(* Taking the Function from User *)
f[{x_,y_}] = Evaluate[Input["Please input a function of x and y to find the minimum "]];
(* The user types in, for instance x^2+y^2 *)

(* Starting the algorithm *)
Do[
 delta =Norm[deltai];

 If[
  delta<epsilon||k>50,
  Break[]
  ];

 corner1=x0+{-deltai[[1]],-deltai[[2]]}/2;
 corner2=x0+{deltai[[1]],-deltai[[2]]}/2;
 corner3=x0+{deltai[[1]],deltai[[2]]}/2;
 corner4=x0+{-deltai[[1]],deltai[[2]]}/2;

 xstar=Flatten[
   Take[
    SortBy[
     {x0,corner1,corner2,corner3,corner4},
     f
     ](* Sort the five points by the values of the function f[x] *),1
    ],1
```

```
    ](* Take the first points *);

  lii=k;

 result[k]=N[{k,Row[corner4,","],Row[corner1,","],Row[x0,","],Row[corner2,","],Row[corner3,",
 "],f[corner4],f[corner1],f[x0],f[corner2],f[corner3],Row[xstar,","]}];
  plotresult[k]=N[{x0,corner1,corner2,corner3,corner4}];

  If[
   xstar==x0,
   deltai=deltai/2,
   x0=xstar
   ],
  {k,1,∞}
  ]

 (* Final Result *)
 If[
   lii==50,
   Print[" After 50 iterations the mimimum point around the intial point is " ,N[xstar],
 "\nThe solution is (approximately) ", N[f[xstar]]];,
  Print["The solution is x= ",  N[x0],"\nThe solution is (approximately)= ", N[f[x0]]];
   ]

 (* Results of Each Iteration *)
 table=TableForm[
   Table[
    result[i],
    {i,1,lii}
    ],
   TableHeadings->{None,{"No. of
 iters.","corner4","corner1","x0","corner2","corner3","f[corner4]","f[corner1]","f[x0]","f[co
 rner2]","f[corner3]","xstar"}}
   ]

 Export["example73.xls",table,"XLS"];

 (* Data visualization *)

 (* Domain of Varibles*)
 xleft=domainx[[1]];
 xright=domainx[[2]];
 ydown=domainy[[1]];
 yup=domainy[[2]];

 (* 3D+ Contour Plot *)
 plot1=Plot3D[
    f[{x,y}],
    {x,xleft,xright},
    {y,ydown,yup},
    ClippingStyle->None,
    MeshFunctions->{#3&},
    Mesh->15,
    MeshStyle->Opacity[.5],
    MeshShading->{{Opacity[.3],Blue},{Opacity[.8],Orange}},
    Lighting->"Neutral"
    ];

 slice=SliceContourPlot3D[
    f[{x,y}],
    z==0,
```

```
   {x,xleft,xright},
   {y,ydown,yup},
   {z,-1,1},
   Contours->15,
   Axes->False,
   PlotPoints->50,
   PlotRangePadding->0,
   ColorFunction->"DarkRainbow"
   ];
Show[
 plot1,
 slice,
 PlotRange->All,
 BoxRatios->{1,1,.6},
 FaceGrids->{Back,Left}
 ]

(* Contour Plot with Step Iterations *)
ContourPlot[
 f[{x,y}],
 {x,xleft,xright},
 {y,ydown,yup},
 LabelStyle->Directive[Black,16],
 ColorFunction->"Rainbow",
 PlotLegends->Automatic,
 Contours->10,
 Epilog-
>{Green,Line[Table[plotresult[i][[1]],{i,1,lii}]],Red,Point[Table[plotresult[i][[1]],{i,1,li
i}]]}
 ]

(* Data Manipulation *)
Manipulate[
 ContourPlot[
  f[{x,y}],
  {x,xleft,xright},
  {y,ydown,yup},
  LabelStyle->Directive[Black,14],
  ColorFunction->"Rainbow",
  PlotLegends->Automatic,
  Contours->10,
  Epilog->{
    PointSize[0.015],
    Green,
    Arrow[{plotresult[i][[1]],plotresult[i+1][[1]]}],
    Red,
    Point[plotresult[i]],
    Blue,
    Line[Table[plotresult[j][[1]],{j,1,i}]]
    }
  ],
 {i,1,lii-1,1}
 ]
```

### 7.4.2 Simplex Search Method (Nelder–Mead Simplex Method)

In the box search method, the number of evaluations, even for problems of modest dimension, rapidly becomes overwhelming. Hence, more efficient methods must be investigated if we are to use these direct-search methods on problems of reasonable size. In the simplex search method [1-5,9], the number of points in the initial simplex is much less compared to that in Box evolutionary optimization method. This reduces the number of function evaluations required in each iteration. With $N$ variables, only $(N+1)$ points are used in the initial simplex.

In $N$ dimensions, a regular simplex is a polyhedron composed of $N+1$ equidistant points, which form its vertices. For example, an equilateral triangle is a simplex in two dimensions; a tetrahedron is a simplex in three dimensions. The main property of a simplex employed by their algorithm is that a new simplex can be generated on any face of the old one by projecting any chosen vertex a suitable distance through the centroid of the remaining vertices of the old simplex. The new simplex is then formed by replacing the old vertex with the newly generated projected point. In this way, each new simplex is generated with a single evaluation of the objective.

For simplicity, let us consider the two-variable problem. We make three initial estimates of the position of the minimum. These will define a starting simplex. Suppose we label the vertices of the simplex $a$, $b$, $c$ and call the corresponding function values $f_a$, $f_b$, $f_c$. The vertex with the highest function value is said to be the worst, and this point must be replaced with a better one. The basic move in the simplex method is reflection. A new trial point is obtained by reflecting the worst point in the centroid of the remaining vertices. This is a heuristic way of placing a new solution estimate in a region where lower function values are likely to occur. Suppose, for instance, that $f_a > f_b > f_c$. Then the vertex $a$ would be reflected in the centroid of vertices $b$ and $c$. Let the new point be labeled as $d$ and let $f_d$ be the associated function value. If $f_d < f_a$ the new point is an improvement on the vertex $a$, and a new simplex is defined by deleting the old worst point and renaming vertex $d$ as $a$. This process is demonstrated in Figure 7.4.

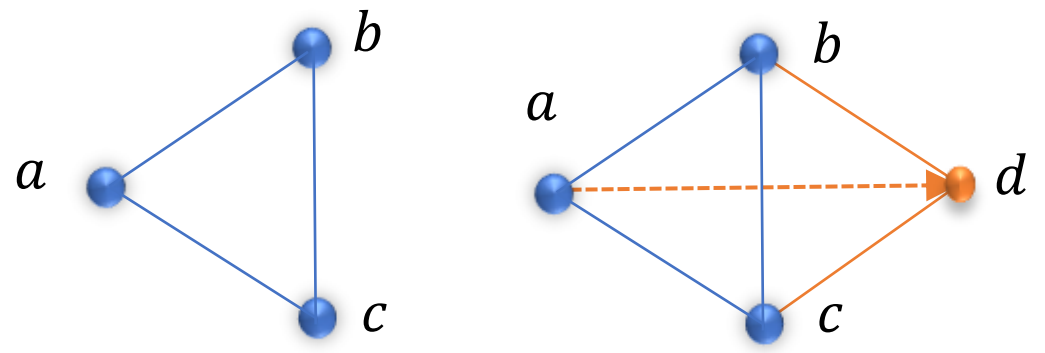

**Figure 7.4.** Construction of new simplex.

Hence, the method begins by setting up a regular simplex in the space of the independent variables and evaluating the function at each vertex. At each iteration, the worst point in the simplex is found first. Then, a new simplex is formed from the old simplex by some fixed rules that steer the search away from the worst point in the simplex. The extent of steering depends on the relative function values of the simplex.

Four different situations may arise depending on the function values. At first, the centroid $\mathbf{x}_c$ of all but the worst point $\mathbf{x}_h$ is determined. Thereafter, the worst point in the simplex is reflected about the centroid, and a new point $\mathbf{x}_r$ is found. If the function value at this point is better than the best point in the simplex, the reflection is considered to have taken the simplex to a good region in the search space. The amount of expansion is controlled by the factor $\gamma$. On the other hand, if the function value at the reflected point is worse than the worst point in the simplex, the reflection is considered to have taken the simplex to a bad region in the search space. Thus, a contraction in the direction from the centroid to the reflected point is made. The amount of contraction is controlled by a factor $\beta$ (a negative value of $\beta$ is used). Finally, if the function value at the reflected point is better than the worst and worse than the next-to-worst point in the simplex, a contraction is made with a positive $\beta$ value. The default scenario is the reflected point itself. The obtained new point replaces the worst point in the simplex, and the algorithm continues with the new simplex.

So, a single iteration then evaluates four simplex operations (the situations are depicted in Figure 7.5):

**Reflection.** $|\mathbf{x}_{\text{new ref}}\rangle = |\mathbf{x}_c\rangle + \alpha|\mathbf{x}_c - \mathbf{x}_h\rangle = 2|\mathbf{x}_c - \mathbf{x}_h\rangle$, reflects the highest-valued point over the centroid. This typically moves the simplex from high regions toward lower regions. Here, $\alpha > 0$ and is typically set to 1.

**Expansion.** $|\mathbf{x}_{\text{new exp}}\rangle = |\mathbf{x}_c\rangle + \gamma|\mathbf{x}_r - \mathbf{x}_c\rangle = (1+\gamma)|\mathbf{x}_c\rangle - \gamma|\mathbf{x}_h\rangle$, like reflection, but the reflected point is sent even further. This is done when the reflected point has an objective function value less than all points in the simplex. Here, $\gamma > max(1, \alpha)$ and is typically set to 2.

**Contraction.** $|\mathbf{x}_{\text{new Con}}\rangle = |\mathbf{x}_c\rangle \pm \beta|\mathbf{x}_h - \mathbf{x}_c\rangle$, the simplex is shrunk down by moving away from the worst point. See Figure 7.5.

**Shrinkage.** All points are moved toward the best point, typically halving the separation distance.

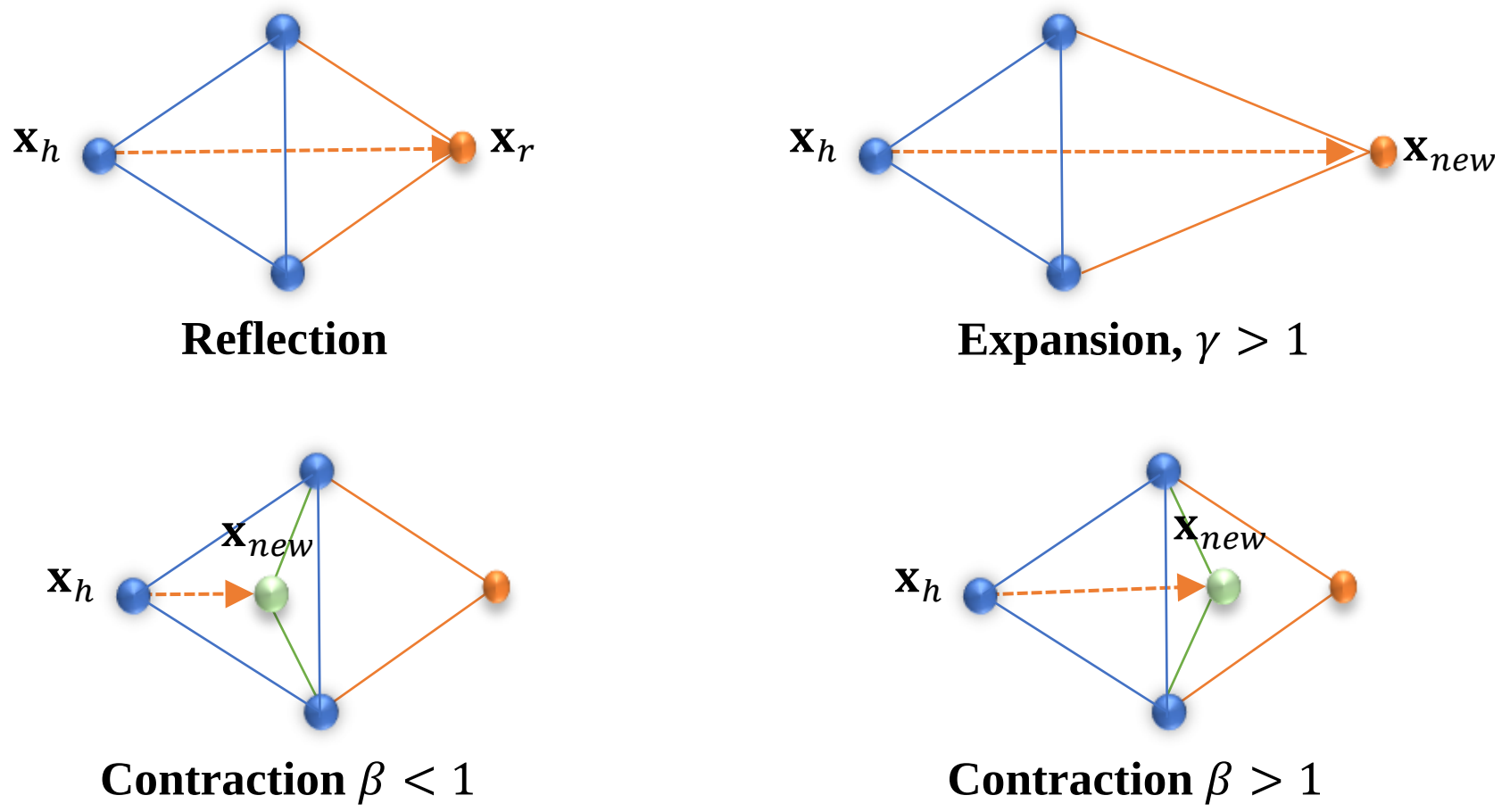

**Figure 7.5.** Four simplex operations, reflection, expansion, contraction, and shrinkage

**Algorithm**

**Step 1:** Choose $\gamma > 1$, $\beta \in (0,1)$, and a termination parameter $\epsilon$. Create an initial simplex.

**Step 2:** Find $|\mathbf{x}_h\rangle$ (the worst point), $|\mathbf{x}_l\rangle$ (the best point), and $|\mathbf{x}_g\rangle$ (next to the worst point). Calculate

$$|\mathbf{x}_c\rangle = \frac{1}{N}\sum_{i=1,i\neq h}^{N+1} |\mathbf{x}_i\rangle$$

(the mean of all vertices except the highest point $\mathbf{x}_h$.)

**Step 3:** Calculate the reflected point $|\mathbf{x}_r\rangle = 2|\mathbf{x}_c\rangle - |\mathbf{x}_h\rangle$.
Set $|\mathbf{x}_{\text{new}}\rangle = |\mathbf{x}_{\text{new}}\rangle$.
If $f|\mathbf{x}_r\rangle < f|\mathbf{x}_l\rangle$, set $|\mathbf{x}_{\text{new}}\rangle = (1+\gamma)|\mathbf{x}_c\rangle - \gamma|\mathbf{x}_h\rangle$ (Expansion);
Else if $f|\mathbf{x}_r\rangle \geq f|\mathbf{x}_l\rangle$, set $|\mathbf{x}_{\text{new}}\rangle = (1-\beta)|\mathbf{x}_c\rangle + \beta|\mathbf{x}_h\rangle$ (contraction);
Else if $f|\mathbf{x}_g\rangle < f|\mathbf{x}_r\rangle < f|\mathbf{x}_h\rangle$, set $|\mathbf{x}_{\text{new}}\rangle = (1+\beta)|\mathbf{x}_c\rangle - \beta|\mathbf{x}_h\rangle$ (contraction).
Calculate $f|\mathbf{x}_{\text{new}}\rangle$ and replace $|\mathbf{x}_h\rangle$ by $|\mathbf{x}_{\text{new}}\rangle$.

**Step 4:** If $\left(\sum_{i=1}^{N+1}\frac{(f|\mathbf{x}_i\rangle - f|\mathbf{x}_c\rangle)^2}{N+1}\right)^{\frac{1}{2}} \leq \epsilon$, Terminate; Else, go to Step 2.

**Remarks:**

- Any other termination criteria may also be used.
- The performance of the above algorithm depends on the values of $\beta$ and $\gamma$. If a large value of $\gamma$ or $1/\beta$ is used, the approach to the optimum point may be faster, but the convergence to the optimum point may be difficult. On the other hand, smaller values of $\gamma$ or $1/\beta$ may require more function evaluations to converge near the optimum point. The recommended values for parameters are $\gamma \approx 2.0$ and $|\beta| \approx 0.5$.
- Even though some guidelines are suggested for choosing the initial simplex, it should be kept in mind that the points chosen for the initial simplex should not form a zero-volume $N$-dimensional hypercube. Thus, in a function with two variables, the chosen three points in the simplex should not lie along a line. Similarly, in a function with three variables, four points in the initial simplex should not lie on a plane.

- It can be shown from elementary geometry that given an $N$-dimensional starting or base point $|\mathbf{x}^0\rangle$ and a scale factor $\alpha$, the other $N$ vertices of the simplex in $N$ dimensions are given by

$$|\mathbf{x}_j^i\rangle = \begin{cases} |\mathbf{x}_j^i\rangle + \delta_1 \text{ if } j = i, \\ |\mathbf{x}_j^i\rangle + \delta_2 \text{ if } j \neq i. \end{cases} \tag{7.2}$$

For $i$ and $j = 1, 2, 3, \dots, N$. The increments $\delta_1$ and $\delta_2$, which depend only on $N$ and the selected scale factor $\sigma$, are calculated from

$$\delta_1 = \left(\frac{\sqrt{N+1} + N - 1}{N\sqrt{2}}\right)\sigma, \quad \delta_2 = \left(\frac{\sqrt{N+1} - 1}{N\sqrt{2}}\right)\sigma. \tag{7.3}$$

Note that the scale factor $\sigma$ is chosen by the user to suit the problem at hand. The choice $\sigma = 1$ leads to a regular simplex with sides of unit length.

- For $N = 2$, we have

$$\delta_1 = \left(\frac{\sqrt{3}+1}{2\sqrt{2}}\right)\sigma, \quad \delta_2 = \left(\frac{\sqrt{3}-1}{2\sqrt{2}}\right)\sigma. \tag{7.4}$$

The vertices are $(x_1^0, x_2^0)^T$, $(x_1^1, x_2^1)^T = (x_1^0 + \delta_1, x_2^0 + \delta_2)^T$ and $(x_1^2, x_2^2)^T = (x_1^0 + \delta_2, x_2^0 + \delta_1)^T$ (see Figure 7.6.)

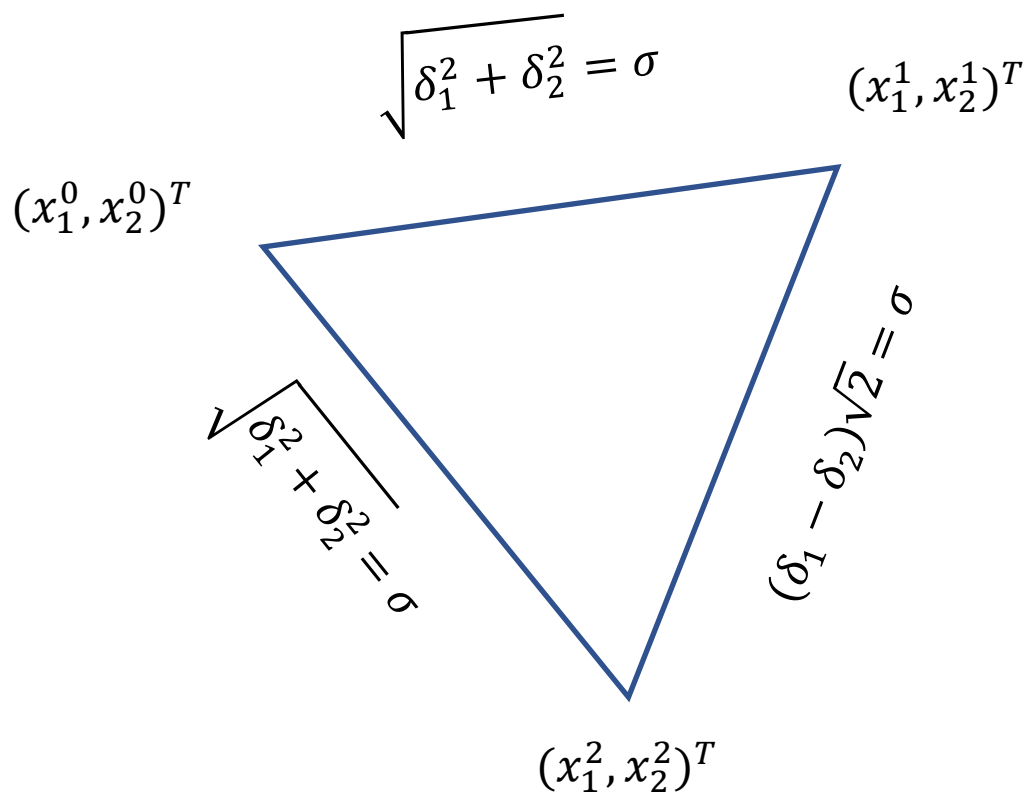

**Figure 7.6.** The vertices of regular simplex in 2D.

- One of the ways to create a simplex is to choose a base point $|\mathbf{x}^0\rangle$ and a scale factor $C$. Then $(N+1)$ points are $|\mathbf{x}^0\rangle$ and for $i, j = 1, 2, \dots, N$,

$$|\mathbf{x}_j^i\rangle = \begin{cases} |\mathbf{x}_j^0\rangle + C & \text{if } j = i, \\ |\mathbf{x}_j^0\rangle + C\Delta & \text{otherwise.} \end{cases} \tag{7.5}$$

where

$$\Delta = \begin{cases} 0.25 & \text{if } N = 3, \\ \dfrac{\sqrt{N+1} - 2}{N - 3} & \text{otherwise.} \end{cases} \tag{7.6}$$

Any other initial simplex may be used.

- This method has no relationship to the simplex method of linear programming. The similarity in name is indeed unfortunate.
- The simplex derives its name from the fact that it is the simplest possible polytope in any given space.

***Example 7.2***

Consider the problem:

$$\text{Minimize } f|\mathbf{x}\rangle = (x^2 + y - 11)^2 + (x + y^2 - 7)^2,$$

using, $|\mathbf{x}_0\rangle = (4,3)$, and $\epsilon = 0.1$.

***Solution***

The 3D and contour plots of the function are shown in Figure 7.7. The plots show that the minimum lies at $|\mathbf{x}^*\rangle = (2.989, 2.075)^T$, $f|\mathbf{x}^*\rangle = 0.0860$.

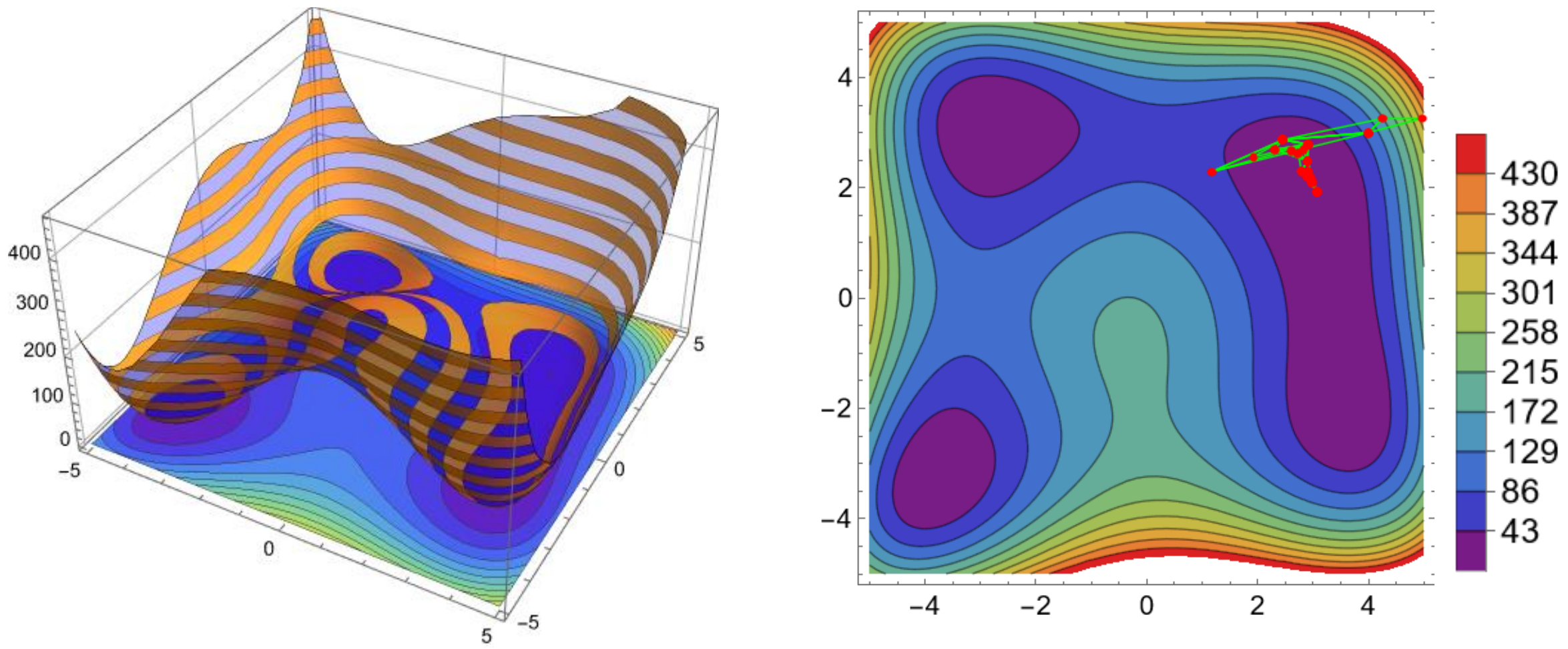

**Figure 7.7.** The results of 18 iterations of the Nelder–Mead simplex method for $f|\mathbf{x}\rangle = (x^2 + y - 11)^2 + (x + y^2 - 7)^2$.

From Table 7.3, after 18 iterations, the optimal condition is attained. The results were produced by Mathematica code 7.4.

**Table 7.3.a.**

| No. of iters. | $x_l$ | $x_g$ | $x_h$ | $x_{new}$ | $x_r$ |
|---|---|---|---|---|---|
| 1 | {4., 3.} | {4.259, 3.259} | {4.966, 3.259} | {2.456, 2.871} | {3.293, 3.} |
| 2 | {2.456, 2.871} | {4., 3.} | {4.259, 3.259} | {1.167, 2.288} | {2.198, 2.612} |
| 3 | {2.456, 2.871} | {1.167, 2.288} | {4., 3.} | {2.906, 2.789} | {-0.377, 2.159} |
| 4 | {2.906, 2.789} | {2.456, 2.871} | {1.167, 2.288} | {1.924, 2.559} | {4.195, 3.372} |
| 5 | {2.906, 2.789} | {2.456, 2.871} | {1.924, 2.559} | {2.303, 2.695} | {3.438, 3.101} |
| 6 | {2.906, 2.789} | {2.303, 2.695} | {2.456, 2.871} | {2.899, 2.485} | {2.752, 2.614} |
| 7 | {2.899, 2.485} | {2.906, 2.789} | {2.303, 2.695} | {2.603, 2.666} | {3.503, 2.580} |
| 8 | {2.899, 2.485} | {2.603, 2.666} | {2.906, 2.789} | {2.829, 2.683} | {2.597, 2.362} |
| 9 | {2.899, 2.485} | {2.829, 2.683} | {2.603, 2.666} | {2.733, 2.625} | {3.126, 2.502} |
| 10 | {2.899, 2.485} | {2.733, 2.625} | {2.829, 2.683} | {2.793, 2.300} | {2.805, 2.428} |
| 11 | {2.793, 2.300} | {2.899, 2.485} | {2.733, 2.625} | {3.072, 1.928} | {2.959, 2.161} |
| 12 | {3.072, 1.928} | {2.793, 2.300} | {2.899, 2.485} | {2.916, 2.299} | {2.965, 1.743} |
| 13 | {3.072, 1.928} | {2.916, 2.299} | {2.793, 2.300} | {2.894, 2.207} | {3.195, 1.928} |
| 14 | {3.072, 1.928} | {2.894, 2.207} | {2.916, 2.299} | {2.949, 2.184} | {3.049, 1.836} |
| 15 | {3.072, 1.928} | {2.949, 2.184} | {2.894, 2.207} | {2.952, 2.132} | {3.128, 1.905} |
| 16 | {3.072, 1.928} | {2.952, 2.132} | {2.949, 2.184} | {2.981, 2.107} | {3.075, 1.876} |
| 17 | {2.981, 2.107} | {3.072, 1.928} | {2.952, 2.132} | {2.989, 2.075} | {3.101, 1.904} |
| 18 | {2.989, 2.075} | {2.981, 2.107} | {3.072, 1.928} | {3.029, 2.009} | {2.898, 2.253} |

**Table 7.3.b.**

| No. of iters. | $f(x_r)$ | $f(x_l)$ | $f(x_g)$ | $f(x_h)$ | $f(x_{new})$ |
|---|---|---|---|---|---|
| 1 | 36.0982 | 100.00 | 170.159 | 359.977 | 18.057 |
| 2 | 16.7422 | 18.057 | 100.00 | 170.159 | 54.3794 |
| 3 | 83.0544 | 18.057 | 54.3794 | 100.00 | 13.6581 |
| 4 | 172.825 | 13.6581 | 18.057 | 54.3794 | 24.6291 |
| 5 | 52.0472 | 13.6581 | 18.057 | 24.6291 | 15.5945 |
| 6 | 7.33708 | 13.6581 | 15.5945 | 18.057 | 4.32513 |
| 7 | 24.8349 | 4.32513 | 13.6581 | 15.5945 | 9.78196 |
| 8 | 4.97272 | 4.32513 | 9.78196 | 13.6581 | 9.25445 |
| 9 | 7.30935 | 4.32513 | 9.25445 | 9.78196 | 7.70397 |
| 10 | 3.38373 | 4.32513 | 7.70397 | 9.25445 | 1.98503 |
| 11 | 0.40026 | 1.98503 | 4.32513 | 7.70397 | 0.177807 |
| 12 | 1.20910 | 0.17781 | 1.98503 | 4.32513 | 1.49152 |
| 13 | 1.30372 | 0.17781 | 1.49152 | 1.98503 | 0.762333 |
| 14 | 0.35569 | 0.17781 | 0.762333 | 1.49152 | 0.529752 |
| 15 | 0.53482 | 0.17781 | 0.529752 | 0.762333 | 0.269427 |
| 16 | 0.27332 | 0.17781 | 0.269427 | 0.529752 | 0.176198 |
| 17 | 0.34442 | 0.17619 | 0.177807 | 0.269427 | 0.0860848 |
| 18 | 1.07163 | 0.08609 | 0.176198 | 0.177807 | 0.0374227 |

***Mathematica Code 7.4*** Nelder–Mead Simplex Method

```
(* Simplex (Nelder and Mead) Search Method *)

(*
Notations :
x0        :Intial vector
delta     :Size reduction parameter delta
epsilon   :Small number to check the accuracy of the simplex (Nelder and Mead) search method
f[x,y]    :Objective function
lii       :The last iteration index
result[k] :The results of iteration k
*)

(* Taking Initial Inputs from User *)
x0=Input["Enter the intial point in the format {x, y}; for example {4,3} "] ;
epsilon=Input["Please enter accuracy of the simplex (Nelder and Mead ) search method; for 
example 0.1 "];

domainx=Input["Please enter domain of x variable for 3D and contour plots; for example {-
5,5}"];
domainy=Input["Please enter domain of y variable for 3D and contour plots; for example {-
5,5}"];

If[
  epsilon<=0,
  Beep[];
  MessageDialog["epsilon, n and σ have to be postive number: "];
  Exit[];
  ];

xx=x0;
σ=1;(* Scale factor to generate simplex:The choice σ=1 leads to a regular simplex with sides 
of unit length *)
γ=2;(* The value of the parameter γ: γ>1 *)
```

```
β=0.5;(* The value of the parameter β: 0<β<1 *)

(* Taking the Function from User *)
f[{x_,y_}] = Evaluate[Input["Please input a function of x and y to find the minimum "]];
(* For example (x^2+y-11)^2+(x+y^2-7)^2 *)
δ1=N[(1+Sqrt[3])/(2*Sqrt[2])*σ];
δ2=N[(-1+Sqrt[3])/(2*Sqrt[2])*σ];
corner1=x0;
corner2=x0+{δ1,δ2};
corner3=x0+{δ2,δ2};
{xl,xg,xh}=SortBy[
  {corner1,corner2,corner3},
  f
  ];(* Sort the three points by the values of the function f[x] *)

(* Starting the Algorithm *)
Do[
 xc=1/2*(xl+xg);
 xr=2*xc-xh;

 Which[
  f[xr]<f[xl],
  xnew=(1+γ)*xc-γ*xh;,
  f[xr]>=f[xl],
  xnew=(1-β)*xc+β*xh;,
  f[xg]<f[xr]<f[xh],
  xnew=(1+β)*xc-β*xh;
  ];

 lii=k;

result[k]=N[{k,Row[xl,","],Row[xg,","],Row[xh,","],Row[xnew,","],Row[xr,","],f[xr],f[xl],f[x
g],f[xh],f[xnew]}];
 plotresult[k]=N[{xl,xg,xh,xl}];

 xh=xnew;
 error=Sqrt[1/3 ((f[xl]-f[xc])2+(f[xg]-f[xc])2+(f[xh]-f[xc])2)];

 If[
  error<epsilon||k>50,
  Break[]
  ];

 {xl,xg,xh}=SortBy[
   {xl,xg,xh},
   f
   ];(* Sort the three points by the values of the function f[x] *),
 {k,1,∞}
 ];

(* Final Result *)
If[
  lii==50,
  Print[" After 50 iterations the mimimum point around the intial point is " ,N[xl], "\nThe
solution is (approximately) ", N[f[xl]]];,
 Print["The solution is x= ",  N[xl],"\nThe solution is (approximately)= ", N[f[xl]]];
  ]

(* Results of Each Iteration *)
table=TableForm[
  Table[
```

```
   result[i],
   {i,1,lii}
   ],
  TableHeadings->{None,{"No. of
iters.","xl","xg","xh","xnew","xr","f[xr]","f[xl]","f[xg]","f[xh]","f[xnew]"}}
  ]

Export["example74.xls",table,"XLS"];
(* Data Visualization *)
(* Domain of Varibles*)
xleft=domainx[[1]];
xright=domainx[[2]];
ydown=domainy[[1]];
yup=domainy[[2]];

(* 3D+ Contour Plot *)
plot1=Plot3D[
   f[{x,y}],
   {x,xleft,xright},
   {y,ydown,yup},
   ClippingStyle->None,
   MeshFunctions->{#3&},
   Mesh->15,
   MeshStyle->Opacity[.5],
   MeshShading->{{Opacity[.3],Blue},{Opacity[.8],Orange}},
   Lighting->"Neutral"
   ];
slice=SliceContourPlot3D[
   f[{x,y}],
   z==0,
   {x,xleft,xright},
   {y,ydown,yup},
   {z,-1,1},
   Contours->15,
   Axes->False,
   PlotPoints->50,
   PlotRangePadding->0,
   ColorFunction->"Rainbow"
   ];
Show[
 plot1,
 slice,
 PlotRange->All,
 BoxRatios->{1,1,.6},
 FaceGrids->{Back,Left}
 ]

(* Contour Plot with Step Iterations *)
ContourPlot[
 f[{x,y}],
 {x,xleft,xright},
 {y,ydown,yup},
 LabelStyle->Directive[Black,16],
 ColorFunction->"Rainbow",
 PlotLegends->Automatic,
 Contours->10,
 Epilog-
>{PointSize[0.015],Green,Line[Table[plotresult[i][[1]],{i,1,lii}]],Red,Point[Table[plotresul
t[i][[1]],{i,1,lii}]]}
 ]
```

```
ContourPlot[
 f[{x,y}],
 {x,xleft,xright},
 {y,ydown,yup},
 LabelStyle->Directive[Black,16],
 ColorFunction->"Rainbow",
 PlotLegends->Automatic,
 Contours->10,
 Epilog-
>{PointSize[0.015],Green,Line[Flatten[Table[plotresult[i],{i,1,lii}],1]],Red,Point[Flatten[T
able[plotresult[i],{i,1,lii}],1]]}
 ]

(* Data Manipulation *)
Manipulate[
 ContourPlot[
  f[{x,y}],
  {x,xleft,xright},
  {y,ydown,yup},
  LabelStyle->Directive[Black,14],
  ColorFunction->"Rainbow",
  PlotLegends->Automatic,
  Contours->10,
  Epilog->{
    PointSize[0.015],
    Yellow,
    Arrow[{plotresult[i][[1]],plotresult[i+1][[1]]}],
    Red,
    Point[plotresult[i][[1]]],
    Green,
    Line[Table[plotresult[j][[1]],{j,1,i}]]
    }
  ],
 {i,1,lii-1,1}
 ]

Manipulate[
 ContourPlot[
  f[{x,y}],
  {x,xleft,xright},
  {y,ydown,yup},
  LabelStyle->Directive[Black,14],
  ColorFunction->"Rainbow",
  PlotLegends->Automatic,
  Contours->10,
  Epilog->{
    PointSize[0.015],
    Yellow,
    Arrow[{plotresult[i][[3]],plotresult[i+1][[1]]}],
    Red,
    Point[Flatten[Table[plotresult[j],{j,1,i}],1]],
    Green,
    Line[Flatten[Table[plotresult[j],{j,1,i}],1]]
    }
  ],
 {i,1,lii-1,1}
 ]
```

### 7.4.3 Hooke-Jeeves Pattern Search Method

The pattern search method works by creating a set of search directions iteratively. In a $N$-dimensional problem, this requires at least $N$ linearly independent search directions. For example, in a two-variable function, at least two search directions are required to go from any one point to any other point. In this approach, the set of directions is chosen to be the coordinate directions in the space of the problem variables. Each of the coordinate directions is sequentially searched using a single-variable optimization scheme until no further improvement is possible.

The Hooke-Jeeves method [1-5] traverses the search space based on evaluations at small steps in each coordinate direction. At every iteration, the Hooke-Jeeves method evaluates $f|\mathbf{x}\rangle$ and $f|\mathbf{x} \pm \Delta\mathbf{e}^{(i)}\rangle$ for a given step size $\Delta$ in every coordinate direction from a base point $|\mathbf{x}\rangle$. It accepts any improvement it may find. If no improvements are found, it will decrease the step size. The process repeats until the step size is sufficiently small. Figure 7.8 shows a few iterations of the algorithm. One step of the Hooke-Jeeves method requires $2N$ function evaluations for an $N$-dimensional problem.

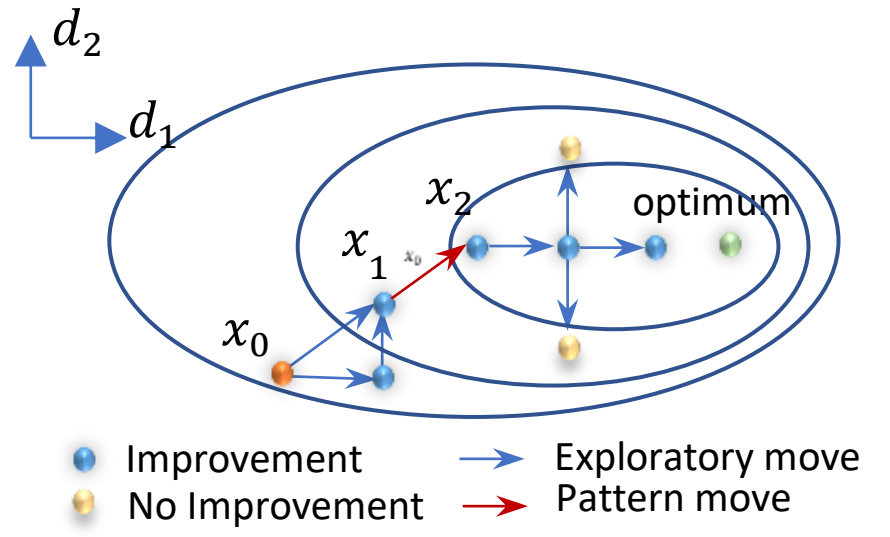

**Figure 7.8.** Illustration of exploratory and pattern moves in the Hooke-Jeeves method.

In the 2-dimensional Hooke and Jeeves method, start from an initial point denoted by $|\mathbf{x}_0^1\rangle$ (Figure 7.9). Start with a fixed variation $\pm\Delta$ parallel to $\mathrm{Ox}_1$. The first point giving a better result is called $|\mathbf{x}_1^1\rangle$. Start from this new point and make a variation $\pm\Delta$ parallel to $\mathrm{Ox}_2$. If there is some improvement, the new point is noted $|\mathbf{x}_2^1\rangle$. It is supposed that at least one of the two directions gave a better result, called $|\mathbf{x}_2^1\rangle$ distinct from $|\mathbf{x}_0^1\rangle$. This final point of the first stage will serve to determine the initial point for the second stage, noted $|\mathbf{x}_0^2\rangle$.

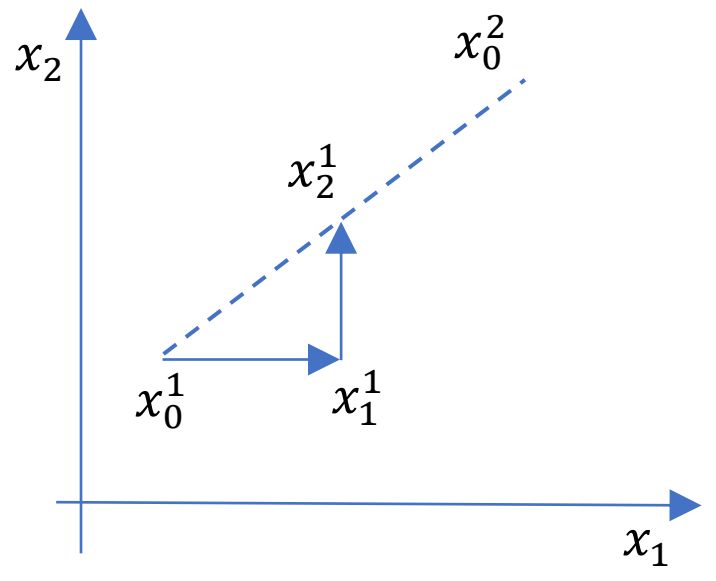

**Figure 7.9.** 2-dimensional Hooke and Jeeves method.

For the second stage, the direction of search is modified by joining the first point $|\mathbf{x}_0^1\rangle$ to the new point $|\mathbf{x}_2^1\rangle$. The acceleration method is introduced at this level by multiplying the distance from $|\mathbf{x}_0^1\rangle$ to $|\mathbf{x}_2^1\rangle$ by a factor $\alpha$ yielding the point $|\mathbf{x}_0^2\rangle$

$$|\mathbf{x}_0^2\rangle = |\mathbf{x}_0^1\rangle + \alpha|\mathbf{x}_2^1 - \mathbf{x}_0^1\rangle = \alpha|\mathbf{x}_2^1\rangle + (1-\alpha)|\mathbf{x}_0^1\rangle. \tag{7.7}$$

The recommended factor is $\alpha = 2$. From this point $|\mathbf{x}_0^2\rangle$, a search is made completely identical to the first stage (parallel to the axes with $\pm\Delta$). A new better point $\left|\mathbf{x}_2^{2)}\right\rangle$ will be obtained, which will serve by acceleration between $|\mathbf{x}_2^2\rangle$ and $|\mathbf{x}_0^2\rangle$ to find the initial point $|\mathbf{x}_0^3\rangle$ of the third stage, and so on. In the case of a failure, that is, when the final point of a stage is not better than the initial point, the search is started again at the very beginning, but in a cautious way, with a small step.

In other words, basically, the Hooke and Jeeves pattern search procedure consists of a sequence of exploratory moves about a base point which, if successful, are followed by pattern moves. The exploratory moves examine the local behavior of the function. An exploratory move is performed in the vicinity of the current point systematically to find the best point around the current point. The pattern moves utilize the information generated in the exploration to step rapidly along the valleys. The procedure is summarized as follows.

**Exploratory move**

The purpose of an exploratory move is to acquire information about the function $f|\mathbf{x}\rangle$ in the neighborhood of the current base point $|\mathbf{x}^c\rangle$. Each variable $x_i^c$, in turn, is given an increment $\Delta_i$ and a check is made of the new function value. If any move is a success (i.e., results in a reduced function value), the new value of that variable will be retained. After all the variables have been considered, a new base point $|\mathbf{x}\rangle$ will be reached. If $|\mathbf{x}\rangle = |\mathbf{x}^c\rangle$, no function reduction has been achieved. The step length $\Delta_i$ is reduced, and the procedure is repeated. If $|\mathbf{x}\rangle \neq |\mathbf{x}^c\rangle$, a pattern moves from $|\mathbf{x}^c\rangle$ is made. Set $i = 1$ and $|\mathbf{x}\rangle = |\mathbf{x}^c\rangle$.

**Algorithm**

**Step 1:** Calculate $f = f|\mathbf{x}\rangle$, $f^+ = f(x_i + \Delta_i)$ and $f^- = f(x_i - \Delta_i)$.

**Step 2:** Find $f_{\min} = \min(f, f^+, f^-)$. Set $|\mathbf{x}\rangle$ corresponds to $f_{\min}$.

**Step 3:** Is $i = N$? If no, set $i = i + 1$ and go to Step 1;
Else $|\mathbf{x}\rangle$ is the result, and go to Step 4.

**Step 4:** If $|\mathbf{x}\rangle \neq |\mathbf{x}^c\rangle$, success;
Else failure.

The current point is changed to the best point at the end of each variable perturbation. If the point found at the end of all variable perturbations is different from the original point, the exploratory move is a success; otherwise, the exploratory move is a failure. In any case, the best point is considered to be the outcome of the exploratory move.

**Pattern move**

A pattern move attempts to speed up the search by using the information already acquired about $f|\mathbf{x}\rangle$ so as to identify the best search direction. By intuition, a move is made from $|\mathbf{x}_k\rangle$ in the direction $|\mathbf{x}_k\rangle - |\mathbf{x}_{k-1}\rangle$, since a move in this direction has led to a decrease in the function value. In other words, a new point is found by jumping from the current best point $|\mathbf{x}^c\rangle$ along a direction connecting the previous best point $|\mathbf{x}^{(k-1)}\rangle$ and the current base point $|\mathbf{x}^{(k)}\rangle$ as follows:

$$|\mathbf{x}_{k+1}^p\rangle = |\mathbf{x}_k\rangle + (|\mathbf{x}_k\rangle - |\mathbf{x}_{k-1}\rangle). \tag{7.8}$$

The algorithm works as follows:

**Algorithm**

**Step 1:** Choose a starting point $|\mathbf{x}_0\rangle$, variable increments $\Delta_i$ $(i = 1, 2, \ldots, N)$, a step reduction factor $\alpha > 1$, and a termination parameter, $\epsilon$. Set $k = 0$.

**Step 2:** Perform an exploratory move with $|\mathbf{x}_k\rangle$ as the base point. Say $|\mathbf{x}\rangle$ is the outcome of the exploratory move. If the exploratory move is a success, set $|\mathbf{x}_{k+1}\rangle = |\mathbf{x}\rangle$ and go to Step 4;
Else go to Step 3.

**Step 3:** Is $\|\Delta\| < \epsilon$? If yes, Terminate;
Else set $\Delta_i = \Delta_i/\alpha$ for $i = 1, 2, \ldots, N$ and go to Step 2.

**Step 4:** Set $k = k + 1$ and perform the pattern move:

$$|\mathbf{x}_{k+1}^p\rangle = |\mathbf{x}_k\rangle + (|\mathbf{x}_k\rangle - |\mathbf{x}_{k-1}\rangle).$$

**Step 5:** Perform another exploratory move using $|\mathbf{x}_{k+1}^p\rangle$ as the base point. Let the result be $|\mathbf{x}_{k+1}\rangle$.

| | |
|---|---|
| **Step 6:** | Is $f|\mathbf{x}_{k+1}\rangle < f|\mathbf{x}_k\rangle$? If yes, go to Step 4;<br>Else go to Step 3. |

***Example 7.3***

Consider the problem:

$$\text{Minimize } f|\mathbf{x}\rangle = (x+2y-7)^2 + (2x+y-5)^2,$$

using, $|\mathbf{x}_0\rangle = (5,5)$, $\Delta_i = (1,1)$, $\alpha = 2$ and $\epsilon = 0.01$.

***Solution***

The 3D and contour plots of the function are shown in Figure 7.10. The plots show that the minimum lies at $|\mathbf{x}^*\rangle = (0.9961, 3.003)^T$, $f|\mathbf{x}^*\rangle = 0.00003$.

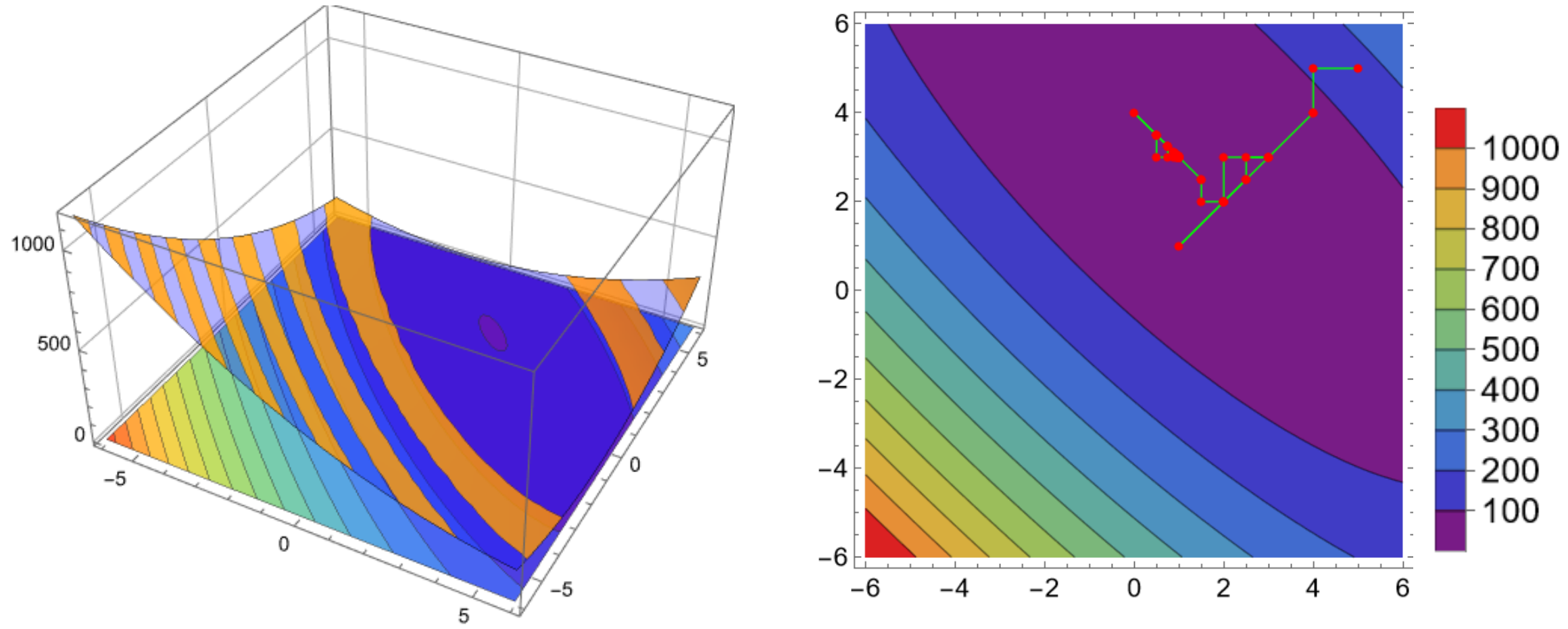

**Figure 7.10.** The results of 12 iterations of the Hooke-Jeeves pattern search method for $f|\mathbf{x}\rangle = (x+2y-7)^2 + (2x+y-5)^2$.

From Table 7.4, after 12 iterations, the optimal condition is attained. The results were produced by Mathematica code 7.5.

**Table 7.4.**

| No. iters. | $x_0$ | $p_{11}$ | $x_1$ | $x_2$ | $f(x_0)$ | $f(p_{11})$ | $f(x_1)$ | $f(x_2)$ |
|---|---|---|---|---|---|---|---|---|
| 1 | {5., 5.} | {4., 5.} | {4., 4.} | {3., 3.} | 164 | 113 | 74 | 20 |
| 2 | {3., 3.} | {2., 3.} | {2., 2.} | {1., 1.} | 20 | 5 | 2 | 20 |
| 3 | {3., 3.} | {2.5, 3.} | {2.5, 2.5} | {2., 2.} | 20 | 11.25 | 6.5 | 2 |
| 4 | {2., 2.} | {1.5, 2.} | {1.5, 2.5} | {1., 3.} | 2 | 2.250 | 0.5 | 0 |
| 5 | {1., 3.} | {0.5, 3.} | {0.5, 3.5} | {0., 4.} | 0 | 1.250 | 0.5 | 2 |
| 6 | {1., 3.} | {0.75, 3.} | {0.75, 3.25} | {0.5, 3.5} | 0 | 0.313 | 0.125 | 0.5 |
| 7 | {1., 3.} | {0.875, 3.} | {0.875, 3.125} | {0.75, 3.25} | 0 | 0.078 | 0.031 | 0.125 |
| 8 | {1., 3.} | {0.938, 3.} | {0.938, 3.063} | {0.875, 3.125} | 0 | 0.019 | 0.008 | 0.031 |
| 9 | {1., 3.} | {0.969, 3.} | {0.969, 3.031} | {0.938, 3.063} | 0 | 0.005 | 0.002 | 0.008 |
| 10 | {1., 3.} | {0.984, 3.} | {0.984, 3.016} | {0.969, 3.031} | 0 | 0.001 | 0.001 | 0.002 |
| 11 | {1., 3.} | {0.992, 3.} | {0.992, 3.008} | {0.984, 3.016} | 0 | 0.0003 | 0.0001 | 0.001 |
| 12 | {1., 3.} | {0.996, 3.} | {0.996, 3.004} | {0.992, 3.008} | 0 | 7.63E-05 | 3.05E-05 | 0.0001 |

***Mathematica Code 7.5*** Hooke-Jeeves Pattern Search Method

```
(* Hooke-Jeeves Pattern Search Method *)

(*
Notations :x0:Intial vector
```

```
deltai     :Size reduction parameter delta
epsilon    :Small number to check the accuracy of the Hooke-Jeeves Pattern Search Method
f[x,y]     :Objective function
α          :The factor α of pattern move:(α=2)
lii        :The last iteration index
result[k] :The results of iteration k
*)

(* Taking Initial Inputs from User *)
x0=Input["Enter the intial point in the format {x, y}; for example {5,5} "] ;

deltai=Input["Enter the size reduction parameter delta in the format {Δx,Δy}; for example
{1,1}"];

epsilon=Input["Please enter accuracy of the Hooke-Jeeves Pattern Search Method; for example
0.01 "];

α=Input["Please enter the factor α: The recommended factor is α=2 "];

domainx=Input["Please enter domain of x variable for 3D and contour plots; for example {-
6,6}"];

domainy=Input["Please enter domain of y variable for 3D and contour plots; for example {-
6,6}"];

If[
  epsilon<=0||deltai[[1]]<=0||deltai[[2]]<=0||α<=0,
  Beep[];
  MessageDialog["the values of epsilon, Δx, Δy and α have to be postive number: "];
  Exit[];
  ];

dd=deltai;
xx=x0;

(* Taking the Function from User *)
f[{x_,y_}] = Evaluate[Input["Please input a function of x and y to find the minimum "]];
(*for example: (x+2y-7)^2+(2x+y-5)^2 *)

(* Defination of Exploratory Search Function *)
exploratorySearch[x_,del_]:=Module[
   {x0=x,deltai=del,cornerPx,cornerNx,cornerPy,cornerNy},

   cornerPx=x0+{deltai[[1]],0};
   cornerNx=x0+{-deltai[[1]],0};

   p11=Flatten[
     Take[
      SortBy[
       {cornerPx,cornerNx},
       f
       ](* Sort the two points by the values of the function f[x] *),1
      ],1
     ](* Take the first points *);

   cornerPy=p11+{0,deltai[[2]]};
   cornerNy=p11+{0,-deltai[[2]]};

   p12=Flatten[
     Take[
      SortBy[
```

```
       {cornerPy,cornerNy},
       f
       ](* Sort the two points by the values of the function f[x] *),1
      ],1
     ](* Take the first points *);

   If[
    p12!=x0,
    xstar=p12;,
    xstar=x0;
    ];
   xstar
   ];

(* Starting the Main Algorithm *)

(* The Pattern Move *)
Do[
 x1=exploratorySearch[x0,deltai];
 x2=x0+α*(x1-x0);

 lii=k;
 result[k]=N[{k,Row[x0,","],Row[p11,","],Row[x1,","],Row[x2,","],f[x0],f[p11],f[x1],f[x2]}];
 plotresult[k]=N[{x0,p11,p12,x1,x2}];

 If[
  f[x2]<f[x1],
  x0=x2;,

  delta =Norm[deltai];
  If[
   delta<epsilon||k>50,
   Break[],
   deltai=deltai/2
   ];
  ],
 {k,1,∞}
 ]

(* Final Result *)

Which[
 lii==50,
 Print[" After 50 iterations the mimimum point around the intial point is " ,N[x1], "\nThe
solution is (approximately) ", N[f[x1]]],
 delta<epsilon,
 Print[" We reashed the condition delta<epsilon; the mimimum point around the intial point
is " ,N[x1], "\nThe solution is (approximately) ", N[f[x1]]],
 True,
 Print["The solution is x= ",  N[x1],"\nThe solution is (approximately)= ", N[f[x1]]]
 ]

(* Results of Each Iteration *)
table=TableForm[
  Table[
   result[i],
   {i,1,lii}
   ],
  TableHeadings->{None,{"No. of
iters.","x0","p11","xl","x2","f[x0]","f[p11]","f[x1]","f[x2]"}}
  ]
```

```
Export["example75.xls",table,"XLS"];

(* Data Visualization *)

(* Domain of Varibles*)
xleft=domainx[[1]];
xright=domainx[[2]];
ydown=domainy[[1]];
yup=domainy[[2]];

(* 3D+ Contour Plot *)
plot1=Plot3D[
   f[{x,y}],
   {x,xleft,xright},
   {y,ydown,yup},
   ClippingStyle->None,
   MeshFunctions->{#3&},
   Mesh->15,
   MeshStyle->Opacity[.5],
   MeshShading->{{Opacity[.3],Blue},{Opacity[.8],Orange}},
   Lighting->"Neutral"
   ];

slice=SliceContourPlot3D[
   f[{x,y}],
   z==0,
   {x,xleft,xright},
   {y,ydown,yup},
   {z,-1,1},
   Contours->15,
   Axes->False,
   PlotPoints->50,
   PlotRangePadding->0,
   ColorFunction->"Rainbow"
   ];

Show[
 plot1,
 slice,
 PlotRange->All,
 BoxRatios->{1,1,.6},
 FaceGrids->{Back,Left}
 ]

(* Contour Plot with Step Iterations *)
ContourPlot[
 f[{x,y}],
 {x,xleft,xright},
 {y,ydown,yup},
 LabelStyle->Directive[Black,16],
 ColorFunction->"Rainbow",
 PlotLegends->Automatic,
 Contours->10,
 Epilog-
>{PointSize[0.015],Green,Line[Flatten[Table[plotresult[i],{i,1,lii}],1]],Red,Point[Flatten[T
able[plotresult[i],{i,1,lii}],1]]}
 ]

(* Data Manipulation *)
Manipulate[
```

```
ContourPlot[
 f[{x,y}],
 {x,xleft,xright},
 {y,ydown,yup},
 LabelStyle->Directive[Black,14],
 ColorFunction->"Rainbow",
 PlotLegends->Automatic,
 Contours->10,
 Epilog->{
   PointSize[0.015],
   Yellow,
   Arrow[{plotresult[i][[3]],plotresult[i+1][[3]]}],
   Red,
   Point[Flatten[Table[plotresult[j],{j,1,i}],1]],
   Green,
   Line[Flatten[Table[plotresult[j],{j,1,i}],1]]
   }
 ],
{i,1,lii-1,1}
]
```

### 7.4.4 Powell Conjugate Direction Method

The Powell method [7,10] is developed for the convex quadratic problem, but it can be applied successfully to nonquadratic problems.

**Definition (Conjugate Directions):**
(a) Two distinct nonzero vectors $|\mathbf{a}_1\rangle$ and $|\mathbf{a}_2\rangle$ are said to be conjugate with respect to a real symmetric matrix $\mathbf{H}$, if

$$\langle \mathbf{a}_1|\mathbf{H}|\mathbf{a}_2\rangle = 0. \tag{7.9}$$

(b) A finite set of distinct nonzero vectors $\{|\mathbf{a}_0\rangle, |\mathbf{a}_1\rangle, \dots, |\mathbf{a}_k\rangle\}$ is said to be conjugate with respect to a real symmetric matrix $\mathbf{H}$, if

$$\langle \mathbf{a}_i|\mathbf{H}|\mathbf{a}_j\rangle = 0 \text{ for all } i \neq j. \tag{7.10}$$

If $\mathbf{H} = \mathbf{I}_n$, where $\mathbf{I}_n$ is the $n \times n$ identity matrix, then (7.9) can be expressed as

$$\langle \mathbf{a}_i|\mathbf{H}|\mathbf{a}_j\rangle = \langle \mathbf{a}_i|\mathbf{I}_n|\mathbf{a}_j\rangle = \langle \mathbf{a}_i|\mathbf{a}_j\rangle = 0 \text{ for } i \neq j. \tag{7.11}$$

This is the well-known condition for orthogonality between vectors $|\mathbf{a}_i\rangle$ and $|\mathbf{a}_j\rangle$. Hence, conjugacy is a generalization of orthogonality.

If $|\mathbf{a}_j\rangle$ for $j = 0,1,\dots,k$ are eigenvectors of $\mathbf{H}$ then

$$\mathbf{H}|\mathbf{a}_j\rangle = \lambda_j|\mathbf{a}_j\rangle. \tag{7.12}$$

where the $\lambda_j$ are the eigenvalues of $\mathbf{H}$. Hence, we have

$$\langle \mathbf{a}_i|\mathbf{H}|\mathbf{a}_j\rangle = \lambda_j\langle \mathbf{a}_i|\mathbf{a}_j\rangle = 0 \text{ for } i \neq j. \tag{7.13}$$

since $|\mathbf{a}_i\rangle$ and $|\mathbf{a}_j\rangle$ for $i \neq j$ are orthogonal. In effect, the set of eigenvectors $|\mathbf{a}_j\rangle$ constitutes a set of conjugate directions with respect to $\mathbf{H}$.

**Theorem 7.1 (Linear Independence of Conjugate Vectors):**
If nonzero vectors $\{|\mathbf{a}_0\rangle, |\mathbf{a}_1\rangle, \dots, |\mathbf{a}_k\rangle\}$ form a conjugate set with respect to a positive definite matrix $\mathbf{H}$, then they are linearly independent.

**Theorem 7.2 (Parallel Subspace Property and Generation of Conjugate Directions in Powell's Method):**
Given a quadratic function $q(\mathbf{x}) = a + \langle \mathbf{b}|\mathbf{x}\rangle + \frac{1}{2}\langle \mathbf{x}|\mathbf{C}|\mathbf{x}\rangle$ of two variables (where $a$ is a scalar quantity, $|\mathbf{b}\rangle$ is a vector, and $\mathbf{C}$ is a $2 \times 2$ matrix). $|\mathbf{x}^1\rangle$ and $|\mathbf{x}^2\rangle$ are two arbitrary but distinct points, $|\mathbf{d}\rangle$ is a direction.
If $|\bar{\mathbf{x}}^1\rangle$ is the solution to the problem

$$\text{minimize } q(|\mathbf{x}^1\rangle + \lambda|\mathbf{d}\rangle), \tag{7.14}$$

and $|\bar{\mathbf{x}}^2\rangle$ is the solution to the problem

$$\text{minimize } q(|\mathbf{x}^2\rangle + \lambda|\mathbf{d}\rangle), \tag{7.15}$$

then the direction $|\bar{\mathbf{x}}^2 - \bar{\mathbf{x}}^1\rangle$ is conjugate to $|\mathbf{d}\rangle$ or, in other words, the quantity $\langle \bar{\mathbf{x}}^2 - \bar{\mathbf{x}}^1|\mathbf{C}|\mathbf{d}\rangle$ is zero.

**Proof:**

The points along the direction $|\mathbf{d}\rangle$ from $|\mathbf{x}^1\rangle$ are

$$|\mathbf{x}\rangle = |\mathbf{x}^1\rangle + \lambda|\mathbf{d}\rangle,$$

and the minimum of $q|\mathbf{x}\rangle$ along $|\mathbf{d}\rangle$ is obtained by finding $\lambda^*$ such that $dq/d|\mathbf{x}\rangle = 0$. The derivative is calculated using the chain rule:

$$\frac{\partial q}{\partial \lambda} = \frac{\partial q}{\partial |\mathbf{x}\rangle}\frac{\partial |\mathbf{x}\rangle}{\partial \lambda} = \langle \mathbf{b}|\mathbf{d}\rangle + \langle \mathbf{x}|\mathbf{C}|\mathbf{d}\rangle.$$

By hypothesis, the minimum occurs at $|\bar{\mathbf{x}}^1\rangle$; hence,

$$\langle \mathbf{b}|\mathbf{d}\rangle + \langle \bar{\mathbf{x}}^1|\mathbf{C}|\mathbf{d}\rangle = 0.$$

Similarly, since the minimum of $q|\mathbf{x}\rangle$ along $|\mathbf{d}\rangle$ from $|\mathbf{x}^2\rangle$ is attained at $|\bar{\mathbf{x}}^2\rangle$ we have

$$\langle \mathbf{b}|\mathbf{d}\rangle + \langle \bar{\mathbf{x}}^2|\mathbf{C}|\mathbf{d}\rangle = 0.$$

Subtracting the last two equations, we have

$$\langle \bar{\mathbf{x}}^2 - \bar{\mathbf{x}}^1|\mathbf{C}|\mathbf{d}\rangle = 0.$$

Accordingly, by definition, the directions $|\mathbf{d}\rangle$ and $|\bar{\mathbf{x}}^2 - \bar{\mathbf{x}}^1\rangle$ are $\mathbf{C}$ conjugate, and the parallel subspace property of quadratic functions has been verified.

■

Thus, if two arbitrary points $|\mathbf{x}^1\rangle$ and $|\mathbf{x}^2\rangle$ and an arbitrary search direction $|\mathbf{d}\rangle$ are chosen, two unidirectional searches, one from each point, will create two points $|\bar{\mathbf{x}}^1\rangle$ and $|\bar{\mathbf{x}}^2\rangle$. For quadratic functions, we can say that the minimum of the function lies on the line joining the points $|\bar{\mathbf{x}}^1\rangle$ and $|\bar{\mathbf{x}}^2\rangle$, as depicted in Figure 7.11(a). The vector $|\bar{\mathbf{x}}^2 - \bar{\mathbf{x}}^1\rangle$ forms a conjugate direction with the original direction vector $|\mathbf{d}\rangle$.

**Extended parallel subspace property**

Instead of using two points $|\mathbf{x}^1\rangle$ and $|\mathbf{x}^2\rangle$ and a direction vector $|\mathbf{d}\rangle$ to create one pair of conjugate directions, one point $|\mathbf{x}^1\rangle$ and both coordinate directions ($(1,0)^T$ and $(0,1)^T$ ) can be used to create a pair of conjugate directions $|\mathbf{d}\rangle$ and $|\bar{\mathbf{x}}^2 - \bar{\mathbf{x}}^1\rangle$ (Figure 7.11(b)). The point $|\bar{\mathbf{x}}^1\rangle$ is obtained by performing a unidirectional search along $(1,0)^T$ from the point $|\mathbf{x}^1\rangle$. Then, the point $|\mathbf{x}^2\rangle$ is obtained by performing a unidirectional search along $(0,1)^T$ from $|\bar{\mathbf{x}}^1\rangle$ and finally the point $|\bar{\mathbf{x}}^2\rangle$ is found by a unidirectional search along the direction $(1,0)^T$ from the point $|\mathbf{x}^2\rangle$.

By comparing Figures 7.11(a) and 7.11(b), we notice that both figures follow the parallel subspace property. The former approach requires two unidirectional searches to find a pair of conjugate directions, whereas the latter approach requires three unidirectional searches. For a quadratic function, the minimum lies in the direction $|\bar{\mathbf{x}}^2 - \bar{\mathbf{x}}^1\rangle$, but for higher-order polynomials the true minimum may not lie in the above direction. Thus, in the case of the quadratic function, four unidirectional searches will find the minimum point, and in higher-order polynomials, more than four unidirectional searches may be necessary. In the latter case, a few iterations of this procedure are required to find the true minimum point. This concept of parallel subspace property can also be extended to higher dimensions.

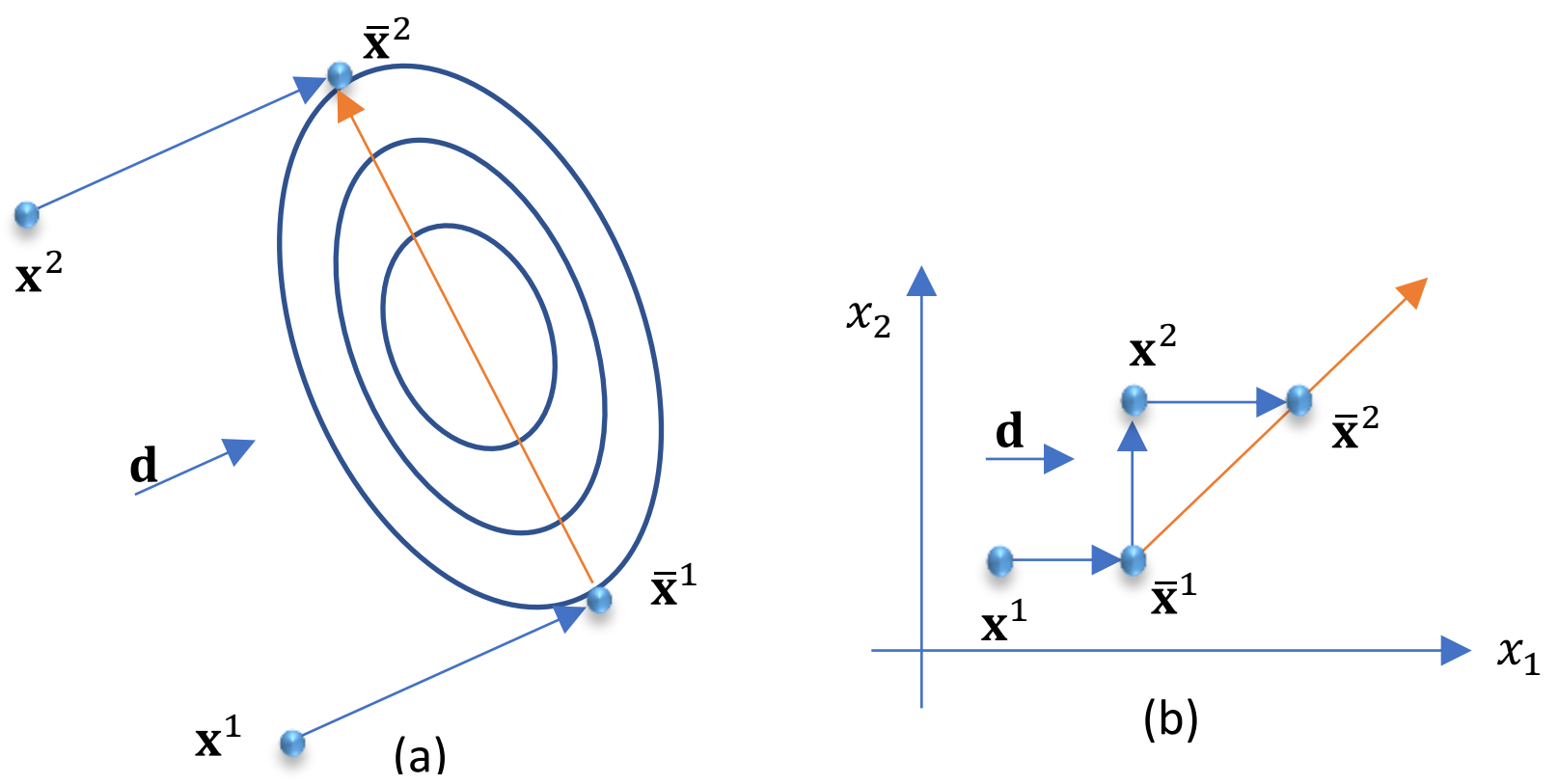

**Figure 7.11.** An illustration of parallel subspace.

Hence, the basic idea is to create a set of $N$ linearly independent search directions and perform a series of unidirectional searches along each of these search directions, starting each time from the previous best point.

**Algorithm**

**Step 1:** Choose a starting point $|\mathbf{x}_0\rangle$ and a set of $N$ linearly independent directions; possibly $|\mathbf{s}_i\rangle = |\mathbf{e}_i\rangle$ for $i = 1,2,\ldots,N$.

**Step 2:** Minimize along $N$ unidirectional search directions using the previous minimum point to begin the next search. Begin with the search direction $|\mathbf{s}_1\rangle$ and end with $|\mathbf{s}_N\rangle$. Thereafter, perform another unidirectional search along $|\mathbf{s}_1\rangle$.

**Step 3:** Form a new conjugate direction $|\mathbf{d}\rangle$ using the extended parallel subspace property.

**Step 4:** If $\|\mathbf{d}\|$ is small or search directions are linearly dependent, Terminate;
Else replace $|\mathbf{s}_j\rangle = |\mathbf{s}_{j-1}\rangle$ for all $j = N, N-1, \ldots, 2$. Set $|\mathbf{s}_1\rangle = |\mathbf{d}\rangle / \|\mathbf{d}\|$ and go to Step 2.

If the function is quadratic, exactly $(N-1)$ loops through Steps 2 to 4 are required. Since in every iteration of the above algorithm exactly $(N+1)$ unidirectional searches are necessary, a total of $(N-1)\times(N+1)$ or $(N^2-1)$ unidirectional searches are necessary to find $N$ conjugate directions. Thereafter, one final unidirectional search is necessary to obtain the minimum point. Thus, in order to find the minimum of a quadratic objective function, the conjugate direction method requires a total of $N^2$ unidirectional searches. For other functions, more loops of the above algorithm may be required. One difficulty with this algorithm is that since unidirectional searches are carried out numerically by using a single-variable search method, the computation of the minimum for unidirectional searches may not be exact. Thus, the resulting directions may not be exactly conjugate to each other. To calculate the extent of deviation, linear independence of the conjugate directions is usually checked. If the search directions are not found to be linearly independent, a completely new set of search directions (possibly conjugate to each other) may be created at the current point. Here, to make the implementation simpler, the coordinate directions can be used again as search directions at the current point.

***Example 7.4***

Consider the problem:

$$\text{Minimize } f|\mathbf{x}\rangle = x^4 - 10xy + y^3,$$

using, $|\mathbf{x}_0\rangle = (5,5)$, and $\epsilon = 0.01$.

***Solution***

The 3D and contour plots of the function are shown in Figure 7.12. The plots show that the minimum lies at $|\mathbf{x}^*\rangle = (1.836, 2.469)^T$, $f|\mathbf{x}^*\rangle = -18.917$.

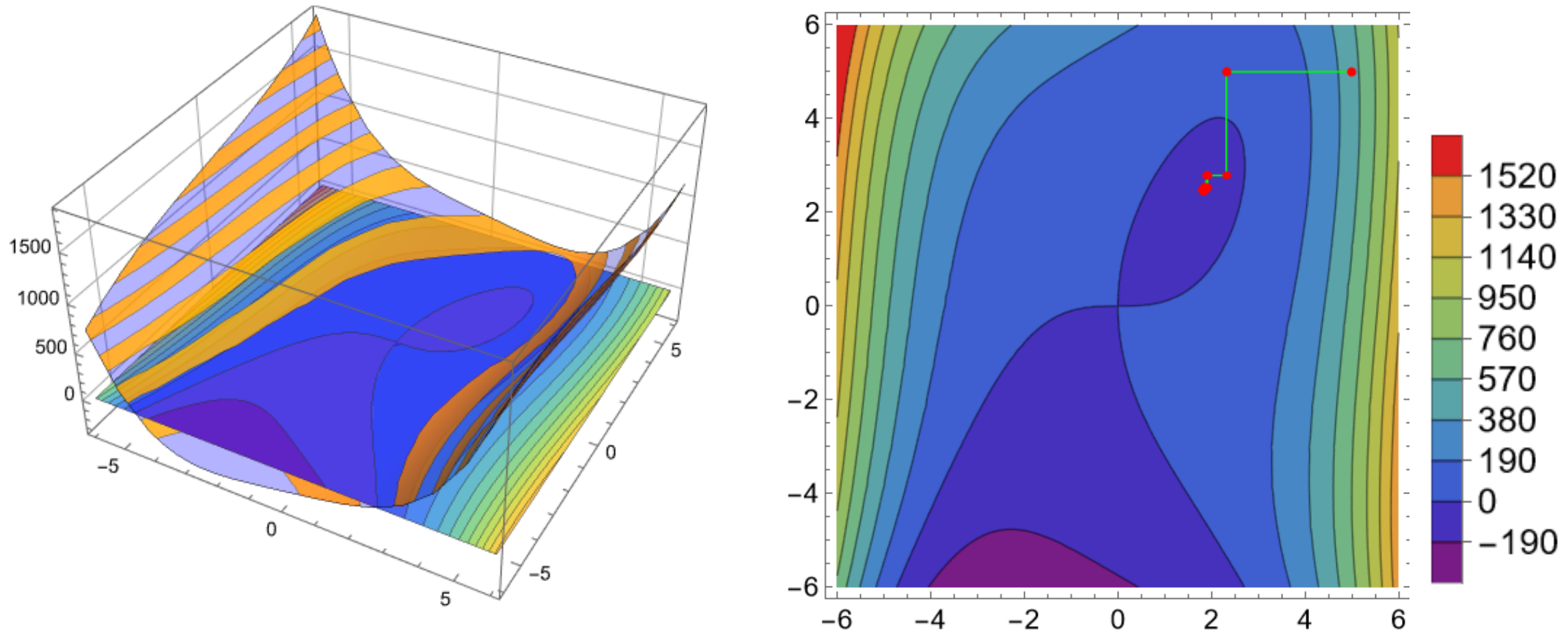

**Figure 7.12.** The results of 5 iterations of the Powell conjugate direction method for $f|\mathbf{x}\rangle = x^4 - 10xy + y^3$.

From Table 7.5, after 12 iterations, the optimal condition is attained. The results were produced by Mathematica code 7.6.

**Table 7.5.a.**

| No. of iters. | $\alpha_1$ | $\alpha_2$ | $\alpha_3$ | $\alpha_4$ |
|---|---|---|---|---|
| 1 | -2.68312 | -2.21248 | -0.40121 | 1.141534 |
| 2 | 0.000583 | -0.25487 | -0.06682 | 1.229789 |
| 3 | 0.000583 | -0.04819 | -0.01804 | 1.165926 |
| 4 | 0.000583 | -0.01804 | 0.000583 | 0.709552 |
| 5 | 0.000583 | 0.000583 | 0.000583 | 4.292127 |

**Table 7.5.b.**

| | $x_{11}$ | $x_2$ | $x_{22}$ | $x_{final}$ | $f(x_{11})$ | $f(x_2)$ | $f(x_{22})$ | $f(x_{final})$ |
|---|---|---|---|---|---|---|---|---|
| 1 | {2.317, 5.} | {2.317, 2.788} | {1.916, 2.788} | {1.859, 2.474} | 37.97 | -14.11 | -18.27 | -18.91 |
| 2 | {1.916, 2.788} | {1.916, 2.533} | {1.849, 2.533} | {1.834, 2.474} | -18.27 | -18.80 | -18.90 | -18.92 |
| 3 | {1.850, 2.533} | {1.850, 2.484} | {1.832, 2.484} | {1.829, 2.477} | -18.90 | -18.91 | -18.92 | -18.92 |
| 4 | {1.833, 2.484} | {1.833, 2.466} | {1.833, 2.466} | {1.833, 2.472} | -18.92 | -18.92 | -18.92 | -18.92 |
| 5 | {1.834, 2.466} | {1.834, 2.466} | {1.834, 2.466} | {1.836, 2.469} | -18.92 | -18.92 | -18.92 | -18.92 |

***Mathematica Code 7.6*** Powell’s Conjugate Direction Method

```
(* Powell's Conjugate Direction Method *)

(*
Notations:
x0           :Intial vector
epsilonPcd0  :Small number to check the accuracy of the Powell's Conjugate Direction Method
f[x,y]       :Objective function
lii          :The last iteration index
result[k]    :The results of iteration k
*)

(* Taking Initial Inputs from User *)
```

```
x0=Input["Enter the intial point in the format {x, y}; for example {5,5} "] ;
epsilonPcd0=Input["Please enter accuracy of the Powell's Conjugate Direction Method; for
example 0.01 "];

domainx=Input["Please enter domain of x variable for 3D and contour plots; for example {-
6,6}"];
domainy=Input["Please enter domain of y variable for 3D and contour plots; for example {-
6,6}"];

If[
  epsilonPcd0<=0,
  Beep[];
  MessageDialog["The value of epsilonPcd has to be postive number: "];
  Exit[];
  ];

(* Taking the Function from User *)
f[{x_,y_}] = Evaluate[Input["Please input a function of x and y to find the minimum "]];
(* For example: x^4-10*x*y+y^3 *)

(*Defination of the Unidirectionalsearch Function*)
(*This function start by bracketing the minimum (using Bounding Phase Method) then isolating
the minimum (using Golden Section Search Method) *)

unidirectionalsearch[α_,delt_,eps_]:=Module[


{α0=α,delta=delt,epsilon=eps,y1,y2,y3,αα,a,b,increment,a0,b0,anew,bnew,anorm,bnorm,lnorm,α1n
orm,α2norm,α1,α2,φ1,φ2,αstar},

   (* Bounding Phase Method *)

   (* Initiating Required Variables *)
   y1 =φ[α0-Abs[delta]];
   y2 =φ[α0];
   y3 =φ[α0+Abs[delta]];

   (*Determining Whether the Inicrement Is Positive or Negative*)
   Which[
    y1==y2,
    a=α0-Abs[delta];
    b=α0;
    Goto[end];,
    y2==y3,
    a=α0;
    b=α0+Abs[delta];
    Goto[end];,
    y1==y3||(y1>y2&&y2<y3),
    a=α0-Abs[delta];
    b=α0+Abs[delta];
    Goto[end];
    ];

   Which[
     y1>y2&&y2>y3,
     increment=Abs[delta];,
     y1<y2&&y2<y3,
     increment=-Abs[delta];
     ]

    (* Starting the Algorithm *)
```

```
     Do[
      αα[0]=α0;
      αα[k+1]=αα[k]+2^k*increment;

      Which[
       φ[αα[k]]<φ[αα[k+1]],(* Evidently, it is impossible the condition to hold for k=0 *)
       a=αα[k-1];
       b= αα[k+1];
       Break [],

       k>50,
       Print["After 50 iterations the bounding phase method can not braketing the min of
alpha"];
       Exit[]
       ];,
      {k,0,∞}
      ];

    Label[end];

    If[
     a>b,
     {a,b}={b,a}
     ];

    (* Golden Section Search Method *)

    (* Initiating Required Variables*)
    a0=a;
    b0=b;
    anew=a;
    bnew=b;

    If[
     a0==b0,
     αstar=a;
     Goto[final]
     ];

    (* Starting the Algorithm *)
    Do[
     (* Normalize the Variable α *)
     anorm=(anew-a)/(b-a);
     bnorm=(bnew-a)/(b-a);

     lnorm=bnorm-anorm;

     α1norm=anorm+0.382*lnorm;
     α2norm=bnorm-0.382*lnorm;

     α1=α1norm(b0-a0)+a0;
     α2=α2norm(b0-a0)+a0;

     φ1=φ[α1];
     φ2=φ[α2];

     Which [
      φ1>φ2,
      anew=α1(*move lower bound to α1*);,
      φ1<φ2,
      bnew=α2(*move upper bound to α2*);,
```

```
      φ1==φ2,
      anew=α1(*move lower bound to α1*);
      bnew=α2(*move upper bound to α2*);
      ];

     αstar=0.5*(anew+bnew);

     If[
      Abs[lnorm]<epsilon,
      Break[]
      ];,
     {k,1,∞}
     ];

    Label[final];

    (* Final result *)
    N[αstar]
    ];

α00=2;(* The intial point of α; for example 2*)
delt0=1;(* The parameter delta of Bounding Phase Method; for example 1 *)
eps0=0.01;(* The accuracy of the Golden Section Search Method; for example 0.01 *)

(* Linearly Independent Directions *)
d1={1,0};
d2={0,1};

(* Main Loop *)
Do[
  φ[α_]=f[x0+α*d1];
  α1=unidirectionalsearch[α00,delt0,eps0];

  x11=x0+α1*d1;
  f11=f[x11];

  φ[α_]=f[x11+α*d2];
  α2=unidirectionalsearch[α00,delt0,eps0];

  x2=x11+α2*d2;
  f2=f[x2];

  φ[α_]=f[x2+α*d1];
  α3=unidirectionalsearch[α00,delt0,eps0];

  x22=x2+α3*d1;
  f22=f[x22];

  d=x22-x11;

  φ[α_]=f[x11+α*d];
  α4=unidirectionalsearch[α00,delt0,eps0];

  xfinal=x11+α4*d;
  ffinal=f[xfinal];

  lii=k;

result[k]=N[{k,α1,α2,α3,α4,Row[x11,","],Row[x2,","],Row[x22,","],Row[xfinal,","],f[x11],f[x2
],f[x22],f[xfinal]}];
  plotresult[k]=N[{x0,x11,x2,x22,xfinal}];
```

```
  If[
   Norm[d]<epsilonPcd0||k>50,
   Break[];,
   x0=x22
   ];,
  {k,1,∞}
  ];

(* Final Result *)

If[
  lii==50,
  Print[" After 50 iterations the mimimum point around the intial point is " ,N[xfinal],
"\nThe solution is (approximately) ", N[ffinal]];,
 Print["The solution is x= ",  N[xfinal],"\nThe solution is (approximately)= ", N[ffinal]];
  ]

(* Results of Each Iteration *)
table=TableForm[
  Table[
   result[i],
   {i,1,lii}
   ],
  TableHeadings->{None,{"No. of
iters.","α1","α2","α3","α4","x11","x2","x22","xfinal","f[x11]","f[x2]","f[x22]","f[xfinal]"}
}
  ]

Export["example76.xls",table,"XLS"];

(* Data Visualization *)
(* Domain of Varibles*)
xleft=domainx[[1]];
xright=domainx[[2]];
ydown=domainy[[1]];
yup=domainy[[2]];

(* 3D+ Contour Plot *)
plot1=Plot3D[
   f[{x,y}],
   {x,xleft,xright},
   {y,ydown,yup},
   ClippingStyle->None,
   MeshFunctions->{#3&},
   Mesh->15,
   MeshStyle->Opacity[.5],
   MeshShading->{{Opacity[.3],Blue},{Opacity[.8],Orange}},
   Lighting->"Neutral"
   ];

slice=SliceContourPlot3D[
   f[{x,y}],
   z==0,
   {x,xleft,xright},
   {y,ydown,yup},
   {z,-1,1},
   Contours->15,
   Axes->False,
   PlotPoints->50,
   PlotRangePadding->0,
```

```
   ColorFunction->"Rainbow"
   ];

Show[
 plot1,
 slice,
 PlotRange->All,
 BoxRatios->{1,1,.6},
 FaceGrids->{Back,Left}
 ]

ContourPlot[
 f[{x,y}],
 {x,xleft,xright},
 {y,ydown,yup},
 LabelStyle->Directive[Black,14],
 ColorFunction->"Rainbow",
 PlotLegends->Automatic,
 Contours->10,
 Epilog-
>{PointSize[0.015],Green,Line[Flatten[Table[plotresult[i],{i,1,lii}],1]],Red,Point[Flatten[T
able[plotresult[i],{i,1,lii}],1]]}
 ]

(* Data Manipulation *)
Manipulate[
 ContourPlot[
  f[{x,y}],
  {x,xleft,xright},
  {y,ydown,yup},
  LabelStyle->Directive[Black,14],
  ColorFunction->"Rainbow",
  PlotLegends->Automatic,
  Contours->10,
  Epilog->{
    PointSize[0.015],
    Yellow,
    Arrow[{plotresult[i][[2]],plotresult[i][[4]]}],
    Red,
    Point[Flatten[Table[plotresult[j],{j,1,i}],1]],
    Green,
    Line[Flatten[Table[plotresult[j],{j,1,i}],1]]
    }
  ],
 {i,1,lii,1}
 ]
```

### 7.4.5 Cyclic Coordinate Search

Coordinate-based algorithms solve optimization problems by advancing along coordinate directions toward a solution. Namely, they iteratively (and approximately) minimize the objective along one or a handful of coordinates at a time. These methods have a long history, being one of the first algorithms proposed for solving optimization problems computationally. Before the rise of modern-day huge-scale applications, these methods were regarded as too simple. They continued to be used in derivative-free optimization, though, where access to derivatives is impossible/too expensive.

Cyclic coordinate search [7], also known as coordinate-descent or taxicab search, simply alternates between coordinate directions for its line search. In this algorithm, an initial point $|\mathbf{x}_1\rangle$ is assumed, and $f|\mathbf{x}\rangle$ is minimized in direction $|\mathbf{d}_1\rangle$ to obtain a new point $|\mathbf{x}_2\rangle$. The procedure is repeated for points $|\mathbf{x}_2\rangle$, $|\mathbf{x}_3\rangle$, ... and when $k = n$, the algorithm is reinitialized, and the procedure is repeated until convergence is achieved. In particular, the search directions $|\mathbf{d}_1\rangle$, ..., $|\mathbf{d}_n\rangle$, where the $|\mathbf{d}_i\rangle$ are vectors of zeros, except for a 1 in the $i$th position. Therefore, along each search direction $|\mathbf{d}_i\rangle$ the corresponding variable $\mathbf{x}_i$ is changed only, with all remaining variables being kept constant to their previous values. For example, the third basis function denoted $|\mathbf{d}_3\rangle$, in four-dimensional space is: $|\mathbf{d}_3\rangle = (0,0,1,0)^T$. The algorithm of the cyclic coordinate method can be summarized as follows:

**Algorithm**

**Step 1:** Input $|\mathbf{x}_1\rangle$ and initialize the tolerance ε.
Set $k = 1$.

**Step 2:** Set $|\mathbf{d}_k\rangle = (0\ 0\ \cdots\ 0\ 1\ 0\ \cdots\ 0)^T$.

**Step 3:** Find $\alpha_k$ , the value of $\alpha$ that minimizes $f(|\mathbf{x}_k\rangle + \alpha|\mathbf{d}_k\rangle)$, using a line search.
Set $|\mathbf{x}_{k+1}\rangle = |\mathbf{x}_k\rangle + \alpha_k|\mathbf{d}_k\rangle$
Calculate $f_{k+1} = f|\mathbf{x}_{k+1}\rangle$.

**Step 4:** If $\|\alpha_k \mathbf{d}_k\|_2 < \varepsilon$ then output $|\mathbf{x}^*\rangle = |\mathbf{x}_{k+1}\rangle$ and $f|\mathbf{x}^*\rangle = f_{k+1}$, and stop.

**Step 5:** If $k = n$, set $|\mathbf{x}_1\rangle = |\mathbf{x}_{k+1}\rangle$, $k = 1$ and repeat from Step 2;
otherwise, set $k = k + 1$ and repeat from Step 2.

**Remark:**

The coordinate-descent algorithm is not very effective or reliable in practice since an oscillatory behavior can sometimes occur.

***Example 7.5***

Consider the problem:

$$\text{Minimize } f|\mathbf{x}\rangle = 5x^2 - 6xy + 5y^2,$$

using, $|\mathbf{x}^0\rangle = (5,5)$, and $\epsilon = 0.01$.

***Solution***

The 3D and contour plots of the function are shown in Figure 7.13. The plots show that the minimum lies at $|\mathbf{x}^*\rangle = (\{0.003, 0.0072)^T$, $f|\mathbf{x}^*\rangle = 0.0002$.

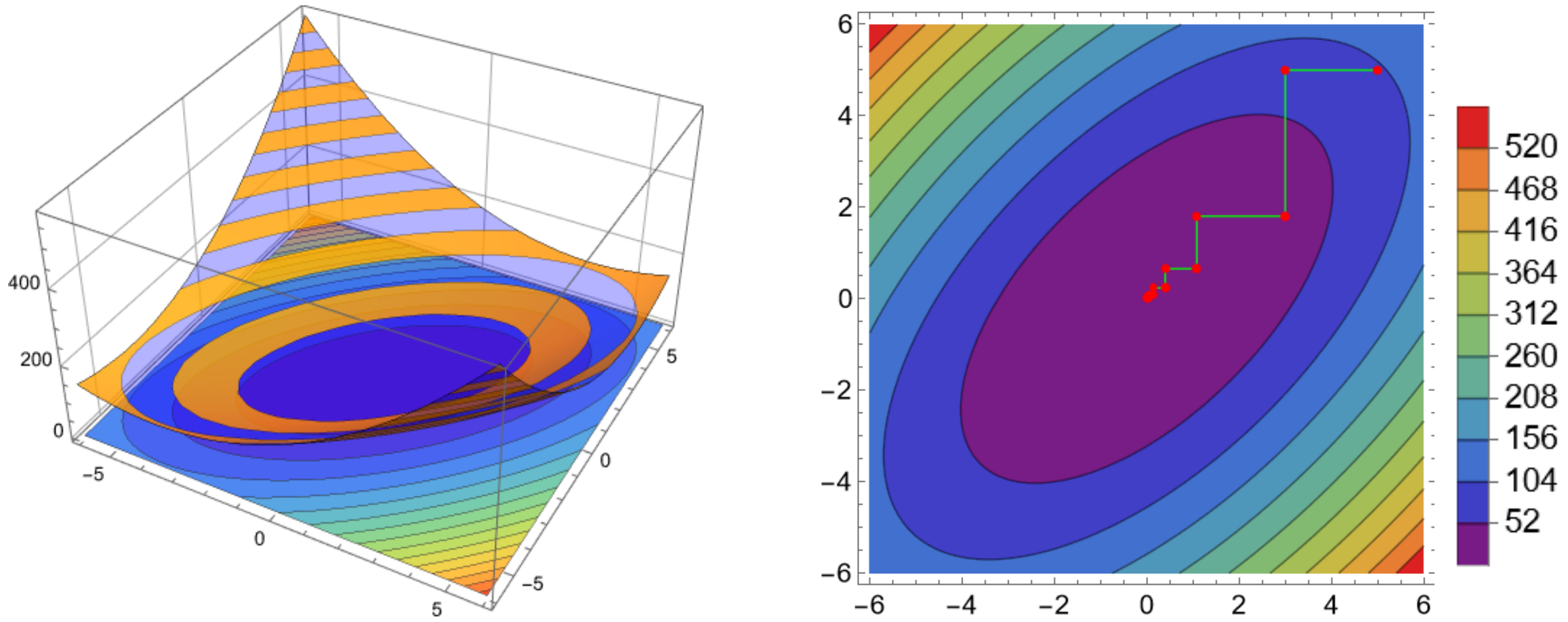

**Figure 7.13.** The results of 7 iterations of the cyclic coordinate search method for $f|\mathbf{x}\rangle = 5x^2 - 6xy + 5y^2$.

From Table 7.6, after 7 iterations, the optimal condition is attained. The results were produced by Mathematica code 7.7.

**Table 7.6.**

| No. iters. | $\alpha_1$ | $\alpha_2$ | $x_1$ | $x_2$ | $x_3$ | $f[x_2]$ | $f[x_3]$ |
|---|---|---|---|---|---|---|---|
| 1 | -2 | -3.19392 | {5., 5.} | {2.999, 5.} | {2.999, 1.806} | 80 | 28.80018 |
| 2 | -1.91522 | -1.14895 | {2.999, 1.806} | {1.085, 1.806} | {1.085, 0.657} | 10.43816 | 3.765782 |
| 3 | -0.68679 | -0.41984 | {1.085, 0.657} | {0.398, 0.657} | {0.398, 0.237} | 1.381876 | 0.506868 |
| 4 | -0.25487 | -0.14581 | {0.398, 0.237} | {0.143, 0.237} | {0.143, 0.092} | 0.18018 | 0.065697 |
| 5 | -0.08545 | -0.04819 | {0.143, 0.092} | {0.058, 0.092} | {0.058, 0.043} | 0.026814 | 0.011018 |
| 6 | -0.03667 | -0.03667 | {0.058, 0.043} | {0.021, 0.043} | {0.021, 0.007} | 0.006117 | 0.00159 |
| 7 | -0.01804 | 0.000583 | {0.021, 0.007} | {0.003, 0.007} | {0.003, 0.007} | 0.000145 | 0.000175 |

***Mathematica Code 7.7*** Cyclic Coordinate Search

```
(* Cyclic Coordinate Search *)

(*
Notations:
x1          :Intial vector
epsilon0cda :Small number to check the accuracy of the Powell's Conjugate Direction Method
f[x,y]      :Objective function
lii         :The last iteration index
result[k]   :The results of iteration k
*)

(* Taking Initial Inputs from User *)
x1=Input["Enter the intial point in the format {x, y}; for example {5,5} "] ;
epsilon0cda=Input["Please enter accuracy of the Cyclic Coordinate Search Method; for example
0.01 "];

If[
  epsilon0cda<=0,
  Beep[];
  MessageDialog["The values of eps, epsilonPcd, Δx, and Δy have to be postive number: "];
  Exit[];
  ];

domainx=Input["Please enter domain of x variable for 3D and contour plots; for example {-
6,6}"];
domainy=Input["Please enter domain of y variable for 3D and contour plots; for example {-
6,6}"];

(* Taking the Function from User *)
f[{x_,y_}]=Evaluate[Input["Please input a function of x and y to find the minimum "]];
(* For example: 5x^2-6x*y+5y^2 *)

(*Defination of the Unidirectionalsearch Function*)
(*This function start by bracketing the minimum (using Bounding Phase Method) then isolating
the minimum (using Golden Section Search Method) *)

unidirectionalsearch[α_,delt_,eps_]:=Module[


{α0=α,delta=delt,epsilon=eps,y1,y2,y3,αα,a,b,increment,a0,b0,anew,bnew,anorm,bnorm,lnorm,α1n
orm,α2norm,α1,α2,φ1,φ2,αstar},

   (* Bounding Phase Method *)

   (* Initiating Required Variables *)
   y1 =φ[α0-Abs[delta]];
```

```
    y2 =φ[α0];
    y3 =φ[α0+Abs[delta]];

    (*Determining Whether the Inicrement Is Positive or Negative*)
    Which[
     y1==y2,
     a=α0-Abs[delta];
     b=α0;
     Goto[end];,
     y2==y3,
     a=α0;
     b=α0+Abs[delta];
     Goto[end];,
     y1==y3||(y1>y2&&y2<y3),
     a=α0-Abs[delta];
     b=α0+Abs[delta];
     Goto[end];
     ];

    Which[
      y1>y2&&y2>y3,
      increment=Abs[delta];,
      y1<y2&&y2<y3,
      increment=-Abs[delta];
      ]

     (* Starting the Algorithm *)
     Do[
      αα[0]=α0;
      αα[k+1]=αα[k]+2^k*increment;

      Which[
       φ[αα[k]]<φ[αα[k+1]],(* Evidently, it is impossible the condition to hold for k=0 *)
       a=αα[k-1];
       b= αα[k+1];
       Break [],

       k>50,
       Print["After 50 iterations the bounding phase method can not braketing the min of
alpha"];
       Exit[]
       ];,
      {k,0,∞}
      ];

    Label[end];

    If[
     a>b,
     {a,b}={b,a}
     ];

    (* Golden Section Search Method *)
    (* Initiating Required Variables*)
    a0=a;
    b0=b;
    anew=a;
    bnew=b;

    If[
     a0==b0,
```

```
     αstar=a;
     Goto[final]
     ];

    (* Starting the Algorithm *)
    Do[
     (* Normalize the Variable α *)
     anorm=(anew-a)/(b-a);
     bnorm=(bnew-a)/(b-a);

     lnorm=bnorm-anorm;

     α1norm=anorm+0.382*lnorm;
     α2norm=bnorm-0.382*lnorm;

     α1=α1norm(b0-a0)+a0;
     α2=α2norm(b0-a0)+a0;

     φ1=φ[α1];
     φ2=φ[α2];

     Which [
      φ1>φ2,
      anew=α1(*move lower bound to α1*);,
      φ1<φ2,
      bnew=α2(*move upper bound to α2*);,
      φ1==φ2,
      anew=α1(*move lower bound to α1*);
      bnew=α2(*move upper bound to α2*);
      ];

     αstar=0.5*(anew+bnew);

     If[
      Abs[lnorm]<epsilon,
      Break[]
      ];,
     {k,1,∞}
     ];

    Label[final];

    (* Final result *)
    N[αstar]
    ];

α00=2;(* The intial point of α; for example 2*)
delt0=1;(* The parameter delta of Bounding Phase Method; for example 1 *)
eps0=0.01;(* The accuracy of the Golden Section Search Method; for example 0.01 *)

(* Linearly Independent Directions *)
d1={1,0};
d2={0,1};

(* Main Loop *)
Do[
  φ[α_]=f[x1+α*d1];
  α1=unidirectionalsearch[α00,delt0,eps0];

  x2=x1+α1*d1;
```

```
  If[
   Norm[α1*d1]<epsilon0cda,
   Break[];
   ];

  φ[α_]=f[x2+α*d2];
  α2=unidirectionalsearch[α00,delt0,eps0];

  x3=x2+α2*d2;

  lii=k;
  result[k]=N[{k,α1,α2,Row[x1,","],Row[x2,","],Row[x3,","],f[x2],f[x3]}];
  plotresult[k]=N[{x1,x2,x3}];

  If[
   Norm[α2*d2]<epsilon0cda||k>50,
   Break[];,
   x1=x3;
   ];,

  {k,1,∞}
  ];

(* Final Result *)

If[
  lii==50,
  Print[" After 50 iterations the mimimum point around the intial point is " ,N[x3], "\nThe
solution is (approximately) ", N[f[x3]]];,
 Print["The solution is x= ",  N[x3],"\nThe solution is (approximately)= ", N[f[x3]]];
  ]

(* Results of Each Iteration *)
table=TableForm[
  Table[
   result[i],
   {i,1,lii}
   ],
  TableHeadings->{None,{"No. of iters.","α1","α2","x1","x2","x3","f[x2]","f[x3]"}}
  ]

Export["example77.xls",table,"XLS"];

(* Data Visualization *)
(* Domain of Varibles*)
xleft=domainx[[1]];
xright=domainx[[2]];
ydown=domainy[[1]];
yup=domainy[[2]];

(* 3D+ Contour Plot *)
plot1=Plot3D[
   f[{x,y}],
   {x,xleft,xright},
   {y,ydown,yup},
   ClippingStyle->None,
   MeshFunctions->{#3&},
   Mesh->15,
   MeshStyle->Opacity[.5],
   MeshShading->{{Opacity[.3],Blue},{Opacity[.8],Orange}},
   Lighting->"Neutral"
```

```
   ];
slice=SliceContourPlot3D[
   f[{x,y}],
   z==0,
   {x,xleft,xright},
   {y,ydown,yup},
   {z,-1,1},
   Contours->15,
   Axes->False,
   PlotPoints->50,
   PlotRangePadding->0,
   ColorFunction->"Rainbow"
   ];
Show[
 plot1,
 slice,
 PlotRange->All,
 BoxRatios->{1,1,.6},
 FaceGrids->{Back,Left}
 ]

ContourPlot[
 f[{x,y}],
 {x,xleft,xright},
 {y,ydown,yup},
 LabelStyle->Directive[Black,16],
 ColorFunction->"Rainbow",
 PlotLegends->Automatic,
 Contours->10,
 Epilog-
>{PointSize[0.015],Green,Line[Flatten[Table[plotresult[i],{i,1,lii}],1]],Red,Point[Flatten[T
able[plotresult[i],{i,1,lii}],1]]}
 ]

(* Data Manipulation *)
Manipulate[
 ContourPlot[
  f[{x,y}],
  {x,xleft,xright},
  {y,ydown,yup},
  LabelStyle->Directive[Black,14],
  ColorFunction->"Rainbow",
  PlotLegends->Automatic,
  Contours->10,
  Epilog->{
    PointSize[0.015],
    Yellow,
    Arrow[{plotresult[i][[2]],plotresult[i][[3]]}],
    Red,
    Point[Flatten[Table[plotresult[j],{j,1,i}],1]],
    Green,
    Line[Flatten[Table[plotresult[j],{j,1,i}],1]]
    }
  ],
 {i,1,lii,1}
 ]
```

# CHAPTER 8

# MULTI-VARIABLE GRADIENT- and HESSIAN- BASED ALGORITHMS

## 8.1 First-order Approximations Methods (Gradient-Based Methods)

This section discusses a variety of algorithms that use first-order methods to select the appropriate descent direction.

### 8.1.1. Cauchy Method (Steepest Descent Method)

The search direction used in the simple gradient method is the negative of the gradient at any particular point $|\mathbf{x}_k\rangle$.

$$|\mathbf{x}_{k+1}\rangle = |\mathbf{x}_k\rangle - \alpha \boldsymbol{\nabla} f|\mathbf{x}_k\rangle, \tag{8.1}$$

where $\alpha$ is a fixed positive parameter. Since this direction gives maximum descent in function values. The method has two disadvantages: the need to make an appropriate choice for $\alpha$ and inherent sluggishness near the minimum due to the corrections vanishing as $\boldsymbol{\nabla} f$ goes to zero.

Accordingly, we are led to adjust $\alpha$ at each iteration:

$$|\mathbf{x}_{k+1}\rangle = |\mathbf{x}_k\rangle - \alpha_k \boldsymbol{\nabla} f|\mathbf{x}_k\rangle. \tag{8.2}$$

The method of steepest descent [1-3] is a gradient algorithm where the step size $\alpha_k$ is chosen to achieve the maximum amount of decrease of the objective function at each individual step. Hence, at every iteration, the derivative is computed at the current point, and a unidirectional single-variable searching method is performed in the negative to this derivative direction to find the minimum point along that direction. The minimum point becomes the current point, and the search is continued from this point. The algorithm continues until a point having a small enough gradient vector is found. This algorithm guarantees improvement in the function value at every iteration. Specifically, $\alpha_k$ is chosen to minimize $\phi_k(\alpha) \equiv f(|\mathbf{x}_k\rangle - \alpha \boldsymbol{\nabla} f|\mathbf{x}_k\rangle)$. In other words,

$$\alpha_k = \underset{\alpha \geq 0}{\arg\ \min} f(|\mathbf{x}_k\rangle - \alpha \boldsymbol{\nabla} f|\mathbf{x}_k\rangle). \tag{8.3}$$

To summarize, the steepest descent algorithm proceeds as follows: at each step, starting from the point $|\mathbf{x}_k\rangle$, we conduct a line search in the direction $-\boldsymbol{\nabla} f|\mathbf{x}_k\rangle$ until a minimizer, $|\mathbf{x}_{k+1}\rangle$, is found. A typical sequence resulting from the method of steepest descent is depicted in Figure 8.1. Observe that the method of steepest descent moves in orthogonal steps.

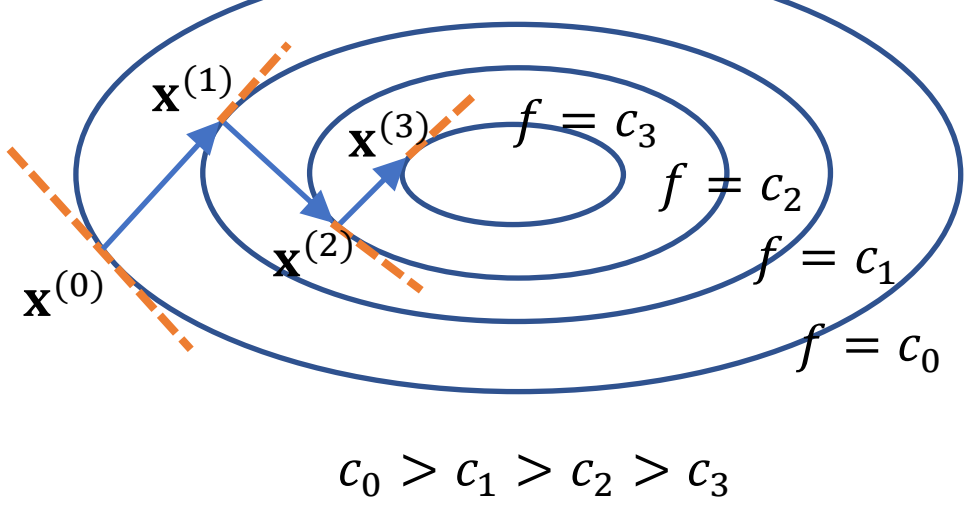

**Figure 8.1.** Typical sequence resulting from the method of steepest descent.

***Example 8.1***

The initial point is $|\mathbf{x}_0\rangle = (1,1,1)^T$. We use the method of steepest descent to find the minimizer of

$$f(x_1,x_2,x_3) = (-5+x_1)^4 + (1+x_2)^2 + 4(-2+x_3)^4.$$

***Solution***

$$\nabla f|\mathbf{x}\rangle = (4(-5+x_1)^3, 2(1+x_2), 16(-2+x_3)^3)^T.$$

We perform six iterations. We find

**Iteration 1:**

| | |
|---|---|
| $\nabla f|\mathbf{x}_0\rangle$ | $:(-256,4,-16)$ |
| $|\mathbf{x}_0\rangle - \alpha\nabla f|\mathbf{x}_0\rangle$ | $:(1+256\,\alpha, 1-4\,\alpha, 1+16\,\alpha)$ |
| $\phi_0(\alpha)$ | $:(2-4\,\alpha)^2 + 4\,(-1+16\,\alpha)^4 + (-4+256\,\alpha)^4$ |
| $\alpha_0$ | $:0.0174879$ |
| $|\mathbf{x}_1\rangle$ | $:(5.4769, 0.930048, 1.27981)$ |

**Iteration 2:**

| | |
|---|---|
| $\nabla f|\mathbf{x}_1\rangle$ | $:(0.433863, 3.8601, -5.97678)$ |
| $|\mathbf{x}_1\rangle - \alpha\nabla f|\mathbf{x}_1\rangle$ | $:(5.4769\ -0.433863\,\alpha, 0.930048\ -3.8601\,\alpha, 1.27981\ +5.97678\,\alpha)$ |
| $\phi_1(\alpha)$ | $:(1.93005\ -3.8601\,\alpha)^2 + (0.476904\ -0.433863\,\alpha)^4 + 4\,(-0.720194+5.97678\,\alpha)^4$ |
| $\alpha_1$ | $:0.196981$ |
| $|\mathbf{x}_2\rangle$ | $:(5.39144, 0.169683, 2.45712)$ |

**Iteration 3:**

| | |
|---|---|
| $\nabla f|\mathbf{x}_2\rangle$ | $:(0.239916, 2.33937, 1.52829\}$ |
| $|\mathbf{x}_2\rangle - \alpha\nabla f|\mathbf{x}_2\rangle$ | $:(5.39144\ -0.239916\,\alpha, 0.169683\ -2.33937\,\alpha, 2.45712\ -1.52829\,\alpha)$ |
| $\phi_2(\alpha)$ | $:(1.16968\ -2.33937\,\alpha)^2 + 4\,(0.457118\ -1.52829\,\alpha)^4 + (0.391441\ -0.239916\,\alpha)^4$ |
| $\alpha_2$ | $:0.465307$ |
| $|\mathbf{x}_3\rangle$ | $:(5.27981, -0.918841, 1.74599)$ |

**Iteration 4:**

| | |
|---|---|
| $\nabla f|\mathbf{x}_3\rangle$ | $:(0.087626, 0.162317, -0.262215)$ |
| $|\mathbf{x}_3\rangle - \alpha\nabla f|\mathbf{x}_3\rangle$ | $:(5.27981\ -0.087626\,\alpha, -0.918841-0.162317\,\alpha, 1.74599\ +0.262215\,\alpha)$ |
| $\phi_3(\alpha)$ | $:(0.0811585\ -0.162317\,\alpha)^2 + (0.279806\ -0.087626\,\alpha)^4 + 4\,(-0.254007+0.262215\,\alpha)^4$ |
| $\alpha_3$ | $:0.630721$ |
| $|\mathbf{x}_4\rangle$ | $:(5.22454, -1.02122, 1.91138)$ |

**Iteration 5:**

| | |
|---|---|
| $\nabla f|\mathbf{x}_4\rangle$ | $:(0.0452829, -0.0424364, -0.0111366)$ |
| $|\mathbf{x}_4\rangle - \alpha\nabla f|\mathbf{x}_4\rangle$ | $:(5.22454\ -0.0452829\,\alpha, -1.02122+0.0424364\,\alpha, 1.91138\ +0.0111366\,\alpha)$ |
| $\phi_4(\alpha)$ | $:(0.224539\ -0.0452829\,\alpha)^4 + 4\,(-0.0886225+0.0111366\,\alpha)^4 + (-0.0212182+0.0424364\,\alpha)^2$ |
| $\alpha_4$ | $:0.848719$ |
| $|\mathbf{x}_5\rangle$ | $:(5.18611, -0.985202, 1.92083)$ |

**Iteration 6:**

| | |
|---|---|
| $\nabla f|\mathbf{x}_5\rangle$ | $:(0.0257836, 0.0295968, -0.00793986)$ |
| $|\mathbf{x}_5\rangle - \alpha\nabla f|\mathbf{x}_5\rangle$ | $:(5.18611\ -0.0257836\,\alpha, -0.985202-0.0295968\,\alpha, 1.92083\ +0.00793986\,\alpha)$ |
| $\phi_5(\alpha)$ | $:(0.0147984\ -0.0295968\,\alpha)^2 + (0.186106\ -0.0257836\,\alpha)^4 + 4\,(-0.0791707+0.00793986\,\alpha)^4$ |
| $\alpha_5$ | $:0.795376$ |
| $|\mathbf{x}_6\rangle$ | $:(5.1656, -1.00874, 1.92714)$ |

We obtain $\alpha_0, \dots, \alpha_5$, by using the line search method. In Figure 8.2, we show a plot of $\phi_i(\alpha)$ versus $\alpha$. Note that the minimizer of $f$ is $(5,-1,2)^T$, and it appears that we have arrived at the minimizer in six iterations.

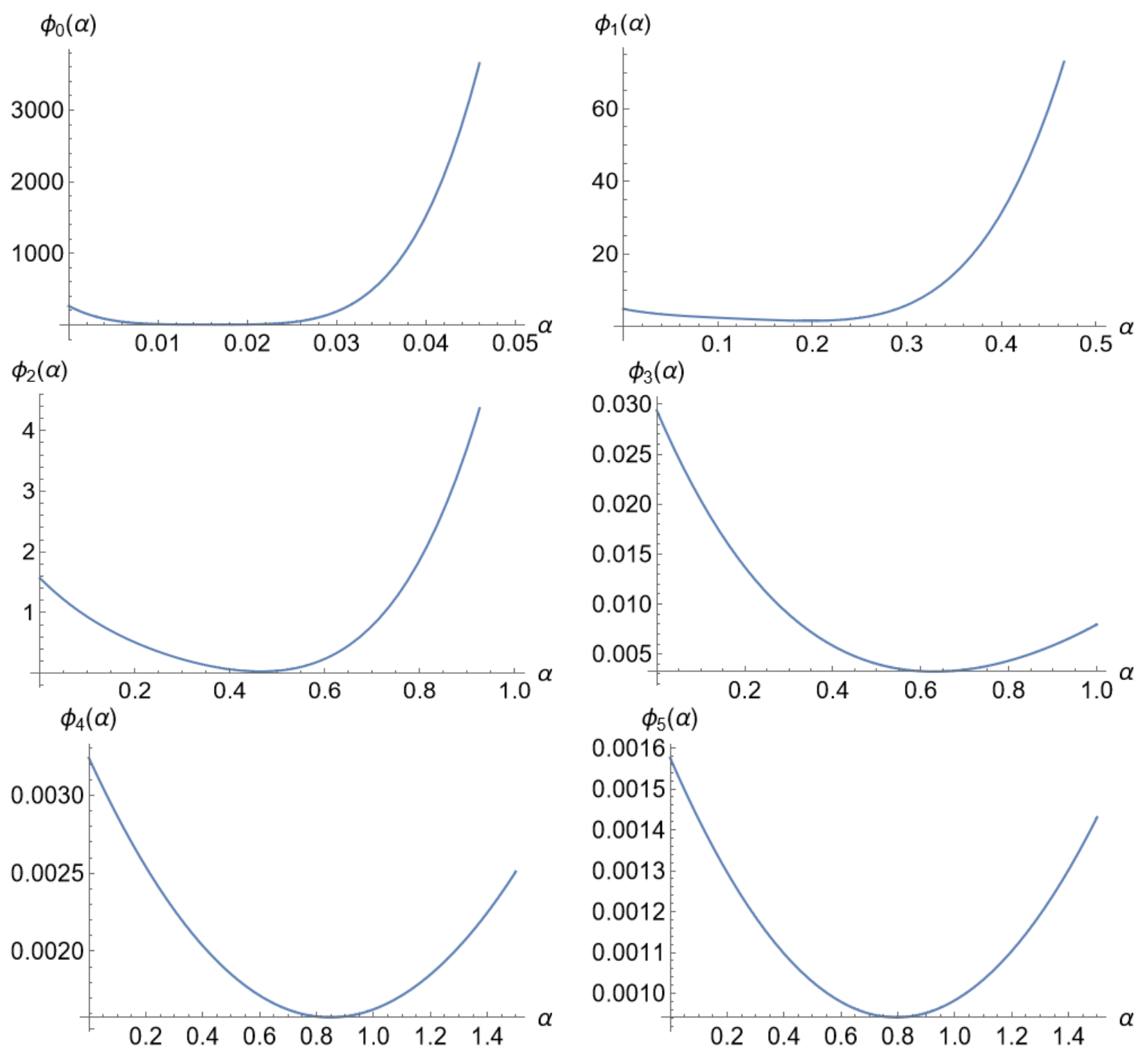

**Figure 8.2.** $\phi_i(\alpha)$ versus $\alpha$.

**Theorem 8.1:** If $\{|\mathbf{x}_k\rangle\}_{k=0}^{\infty}$ is a steepest descent sequence for a given function $f: \mathbb{R}^n \to \mathbb{R}$, then for each $k$ the vector $|\mathbf{x}_{k+1}\rangle - |\mathbf{x}_k\rangle$ is orthogonal to the vector $|\mathbf{x}_{k+2}\rangle - |\mathbf{x}_{k+1}\rangle$.

**Proof:**

**First, note that:**

Let $f$ be a function of $x$ and $y$ and $|\mathbf{u}\rangle = (u_1, u_2)^T$ be a unit vector. Then the directional derivative of $f$ at $|\mathbf{x}\rangle = (x, y)^T$ in the direction of $|\mathbf{u}\rangle$ is

$$D_{|\mathbf{u}\rangle} f|\mathbf{x}\rangle = \lim_{h \to 0} \frac{f|\mathbf{x} + h\mathbf{u}\rangle - f|\mathbf{x}\rangle}{h} = \langle \nabla f|\mathbf{u}\rangle. \tag{8.4}$$

where $|\mathbf{x} + h\mathbf{u}\rangle = (x + hu_1, y + hu_2)^T$,

$$D_{|\mathbf{u}\rangle} f|\mathbf{x}\rangle = \lim_{h \to 0} \frac{f|\mathbf{x} + h\mathbf{u}\rangle - f|\mathbf{x}\rangle}{h} = \frac{\partial f|\mathbf{x}\rangle}{\partial|\mathbf{u}\rangle} = \langle \mathbf{u}|\nabla f\rangle. \tag{8.5}$$

But

$$\frac{d}{d\alpha} f|\mathbf{x} + \alpha\mathbf{u}\rangle = \lim_{h \to 0} \frac{f|\mathbf{x} + (\alpha + h)\mathbf{u}\rangle - f|\mathbf{x} + \alpha\mathbf{u}\rangle}{h}. \tag{8.6}$$

So that

$$D_{|\mathbf{u}\rangle} f|\mathbf{x}\rangle = \frac{\partial f|\mathbf{x}\rangle}{\partial|\mathbf{u}\rangle} = \frac{d}{d\alpha} f|\mathbf{x} + \alpha\mathbf{u}\rangle\Big|_{\alpha=0} = \langle \nabla f|\mathbf{u}\rangle. \tag{8.7}$$

From the iterative formula of the method of steepest descent, it follows that

$$\begin{aligned}\langle \mathbf{x}_{k+1} - \mathbf{x}_k|\mathbf{x}_{k+2} - \mathbf{x}_{k+1}\rangle &= \langle \alpha_k \nabla f(\mathbf{x}_k)|\alpha_{k+1}\nabla f(\mathbf{x}_{k+1})\rangle \\ &= \alpha_k\alpha_{k+1}\langle \nabla f(\mathbf{x}_k)|\nabla f(\mathbf{x}_{k+1})\rangle.\end{aligned}$$

To complete the proof, it is enough to show that

$$\langle \nabla f(\mathbf{x}_k) | \nabla f(\mathbf{x}_{k+1}) \rangle = 0.$$

To this end, observe that $\alpha_k$ is a nonnegative scalar that minimizes $\phi_k(\alpha) \equiv f(|\mathbf{x}_k\rangle - \alpha \nabla f|\mathbf{x}_k\rangle)$. Hence,

$$\begin{aligned} 0 &= \phi_k'(\alpha_k) \\ &= \frac{d}{d\alpha}\phi_k(\alpha_k) \\ &= \frac{d}{d\alpha} f(|\mathbf{x}_k\rangle - \alpha_k \nabla f|\mathbf{x}_k\rangle) \\ &= \big\langle \nabla f\big(\mathbf{x}_k - \alpha_k \nabla f(\mathbf{x}_k)\big)\big| -\nabla f(\mathbf{x}_k)\big\rangle \\ &= -\langle \nabla f(\mathbf{x}_{k+1}) | \nabla f(\mathbf{x}_k)\rangle, \end{aligned}$$

and the proof is completed.

■

The above proposition implies that $\nabla f|\mathbf{x}_k\rangle$ is parallel to the tangent plane to the level set $\{f|\mathbf{x}\rangle = f|\mathbf{x}_{k+1}\rangle\}$ at $|\mathbf{x}_{k+1}\rangle$.

**Theorem 8.2:** If $\{|\mathbf{x}_k\rangle\}_{k=0}^{\infty}$ is a steepest descent sequence for $f: \mathbb{R}^n \to \mathbb{R}$ and if $\nabla f|\mathbf{x}_k\rangle \neq |\mathbf{0}\rangle$, then

$$f|\mathbf{x}_{k+1}\rangle < f|\mathbf{x}_k\rangle. \tag{8.8}$$

**Proof:**

First, recall that

$$|\mathbf{x}_{k+1}\rangle = |\mathbf{x}_k\rangle - \alpha_k \nabla f|\mathbf{x}_k\rangle,$$

where $\alpha_k \geq 0$ is the minimizer of

$$\phi_k(\alpha) = f(|\mathbf{x}_k\rangle - \alpha \nabla f|\mathbf{x}_k\rangle),$$

over all $\alpha \geq 0$. Thus, for $\alpha \geq 0$, we have

$$\phi_k(\alpha_k) \leq \phi_k(\alpha).$$

By the chain rule,

$$\begin{aligned} \phi_k'(0) &= \frac{d}{d\alpha}\phi_k(0) \\ &= \frac{d}{d\alpha} f(|\mathbf{x}_k\rangle - 0 \times \nabla f|\mathbf{x}_k\rangle) \\ &= \langle \nabla f((\mathbf{x}_k) - 0 \times \nabla f(\mathbf{x}_k)) | -\nabla f(\mathbf{x}_k)\rangle \\ &= -\|\nabla f|\mathbf{x}_k\rangle\|^2 < 0, \end{aligned}$$

because $\nabla f|\mathbf{x}_k\rangle \neq |\mathbf{0}\rangle$ by assumption. Thus, $\phi_k'(0) < 0$ and this implies that there is an $\bar{\alpha} > 0$ such that $\phi_k(0) > \phi_k(\alpha)$ for all $\alpha \in (0, \bar{\alpha}]$. Hence,

$$f|\mathbf{x}_{k+1}\rangle = \phi_k(\alpha_k) \leq \phi_k(\bar{\alpha}) < \phi_k(0) = f|\mathbf{x}_k\rangle,$$

and the proof of the statement is completed.

■

In the above, we proved that the algorithm possesses the descent property: $f|\mathbf{x}_{k+1}\rangle < f|\mathbf{x}_k\rangle$ if $\nabla f|\mathbf{x}_k\rangle \neq |\mathbf{0}\rangle$.
If for some $k$, we have $\nabla f|\mathbf{x}_k\rangle = |\mathbf{0}\rangle$, in this case, $|\mathbf{x}_{k+1}\rangle = |\mathbf{x}_k\rangle$.

1- We can use the above as the basis for a stopping criterion for the algorithm. The condition $\nabla f|\mathbf{x}_{k+1}\rangle = |\mathbf{0}\rangle$, however, is not directly suitable as a practical stopping criterion because the numerical computation of the gradient will rarely be identically equal to zero.

2- A practical stopping criterion is to check if the norm $\|\nabla f|\mathbf{x}_k\rangle\|$ of the gradient is less than a prespecified threshold, in which case we stop.

3- Alternatively, we may compute the absolute difference $|f|\mathbf{x}_{k+1}\rangle - f|\mathbf{x}_k\rangle|$ between objective function values for every two successive iterations, and if the difference is less than some prespecified threshold, then we stop; that is, we stop when

$$|f|\mathbf{x}_{k+1}\rangle - f|\mathbf{x}_k\rangle| < \varepsilon, \tag{8.9}$$

where $\varepsilon > 0$ is a prespecified threshold.

4- Yet another alternative is to compute the norm $\|\mathbf{x}_{k+1} - \mathbf{x}_k\|$ of the difference between two successive iterates, and we stop if the norm is less than a prespecified threshold:

$$\|\mathbf{x}_{k+1} - \mathbf{x}_k\| < \varepsilon. \tag{8.10}$$

5- Alternatively, we may check "relative" values of the above quantities, for example,

$$\frac{|f|\mathbf{x}_{k+1}\rangle - f|\mathbf{x}_k\rangle|}{|f|\mathbf{x}_k\rangle|} < \varepsilon, \tag{8.11}$$

or

$$\frac{\|\mathbf{x}_{k+1} - \mathbf{x}_k\|}{\|\mathbf{x}_k\|} < \varepsilon. \tag{8.12}$$

6- The above relative stopping criteria are preferable because the relative criteria are "scale-independent." For example, scaling the objective function does not change the satisfaction of the criterion $|f|\mathbf{x}_{k+1}\rangle - f|\mathbf{x}_k\rangle|/|f|\mathbf{x}_k\rangle| < \varepsilon$. Similarly, scaling the decision variable does not change the satisfaction of the criterion $\|\mathbf{x}_{k+1} - \mathbf{x}_k\|/\|\mathbf{x}_k\| < \varepsilon$.

7- To avoid dividing by very small numbers, we can modify these stopping criteria as follows:

$$\frac{|f|\mathbf{x}_{k+1}\rangle - f|\mathbf{x}_k\rangle|}{\max(1, |f|\mathbf{x}_k\rangle|)} < \varepsilon, \tag{8.13}$$

$$\frac{\|\mathbf{x}_{k+1} - \mathbf{x}_k\|}{\max(1, \|\mathbf{x}_k\|)} < \varepsilon. \tag{8.14}$$

Note that the above stopping criteria are relevant to all the iterative algorithms we discuss in this part.

**Algorithm**

**Step 1:** Choose a maximum number of iterations $M$ to be performed, an initial point $|\mathbf{x}_0\rangle$, two termination parameters $\epsilon_1$, $\epsilon_2$, and set $k = 0$.

**Step 2:** Calculate $\nabla f|\mathbf{x}_k\rangle$, the first derivative at the point $x^{(k)}$.

**Step 3:** If $\|\nabla f(\mathbf{x}_k)\| \le \epsilon_1$, Terminate;
Else if $k \ge M$; Terminate;
Else go to Step 4.

**Step 4:** Perform a unidirectional search to find $\alpha_k$ using $\epsilon_2$ such that $f|\mathbf{x}_{k+1}\rangle = f(|\mathbf{x}_k\rangle - \alpha_k \nabla f|\mathbf{x}_k\rangle)$ is minimum. One criterion for termination is when $|\nabla f|\mathbf{x}_{k+1}\rangle \cdot \nabla f|\mathbf{x}_k\rangle| \le \epsilon_2$.

**Step 5:** Is $\|\mathbf{x}_{k+1} - \mathbf{x}_k\|/\|\mathbf{x}_k\| \le \epsilon_1$? If yes, Terminate;
Else set $k = k + 1$ and go to Step 2.

**Remark:**

Cauchy method works well when $|\mathbf{x}_0\rangle$ is far away from $|\mathbf{x}^*\rangle$. When the current point is very close to the minimum, the change in the gradient vector is small. Thus, the new point created by the unidirectional search is also close to the current point (the change in the variables is directly related to the magnitude of the gradient, which is itself going to zero). There is, then, no mechanism to produce an acceleration toward the minimum in the final iterations. This slows the convergence process near the true minimum. Convergence can be made faster by using second-order derivatives.

***Example 8.2***

Consider the problem:

$$\text{Minimize } f|\mathbf{x}\rangle = \frac{x^2}{4} + \frac{y^2}{25},$$

using, $|\mathbf{x}_0\rangle = (5,5)^T$, and $\epsilon = 0.01$.

***Solution***

The 3D and contour plots of the function are shown in Figure 8.3. The plots show that the minimum lies at $|\mathbf{x}^*\rangle = (0.0049, 0.0052)^T$, $f|\mathbf{x}^*\rangle = 0.000007$.

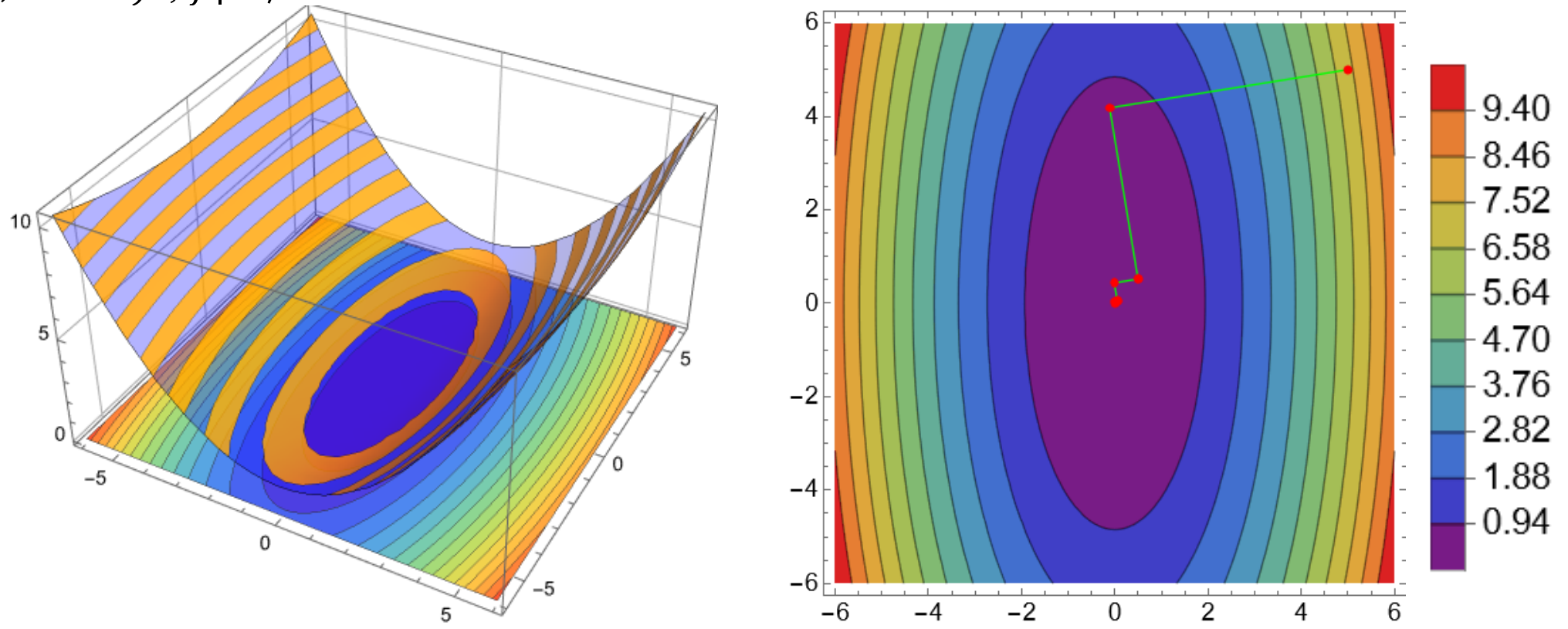

**Figure 8.3.** The results of 6 iterations of the steepest descent method for $f|\mathbf{x}\rangle = x^2/4 + y^2/25$.

From Table 8.1, after 6 iterations, the optimal condition is attained. The results were produced by Mathematica code 8.1.

**Table 8.1.**

| No. of iters. | $\alpha 1$ | $x0$ | $x1$ | $f(x0)$ | $f(x1)$ | error |
|---|---|---|---|---|---|---|
| 1 | 2.04452 | {5., 5.} | {-0.111, 4.182} | 7.25 | 0.702726 | 0.732041 |
| 2 | 10.95101 | {-0.111, 4.182} | {0.498, 0.5183} | 0.702726 | 0.072777 | 0.887803 |
| 3 | 2.04452 | {0.498, 0.518} | {-0.011, 0.434} | 0.072777 | 0.007547 | 0.516225 |
| 4 | 11.049 | {-0.011, 0.434} | {0.050, 0.050} | 0.007547 | 0.000731 | 0.388033 |
| 5 | 2.04452 | {0.050, 0.050} | {-0.001, 0.042} | 0.000731 | 7.12E-05 | 0.051943 |
| 6 | 10.95101 | {-0.001, 0.042} | {0.005, 0.005} | 7.12E-05 | 7.33E-06 | 0.037377 |

***Mathematica Code 8.1*** Steepest Descent Method

```
(*Cauchy Steepest Descent Method*)

(*
Notations:
x0        :Intial vector
epsilonSD :Small number to check the accuracy of the Steepest Descent Search Method
f[x,y]    :Objective function
lii       :The last iteration index
result[k] :The results of iteration k
*)

(* Taking Initial Inputs from User *)
x0=Input["Enter the intial point in the format {x, y}; for example {5,5} "] ;
epsilonSD=Input["Please enter accuracy of the Steepest Descent Method; for example 0.01 "];

domainx=Input["Please enter domain of x variable for 3D and contour plots; for example {-
6,6}"];
domainy=Input["Please enter domain of y variable for 3D and contour plots; for example {-
6,6}"];
```

```
If[
  epsilonSD<=0,
  Beep[];
  MessageDialog["The value of epsilonSD has to be postive number: "];
  Exit[];
  ];

(* Taking the Function from User *)
f[{x_,y_}] = Evaluate[Input["Please input a function of x and y to find the minimum "]];
(* For example : x^2/4+y^2/25 *)

(*Defination of the Unidirectionalsearch Function*)
(*This function start by bracketing the minimum (using Bounding Phase Method) then isolating
the minimum (using Golden Section Search Method) *)

unidirectionalsearch[α_,delt_,eps_]:=Module[


{α0=α,delta=delt,epsilon=eps,y1,y2,y3,αα,a,b,increment,a0,b0,anew,bnew,anorm,bnorm,lnorm,α1n
orm,α2norm,α1,α2,φ1,φ2,αstar},

   (* Bounding Phase Method *)
   (* Initiating Required Variables *)
   y1 =φ[α0-Abs[delta]];
   y2 =φ[α0];
   y3 =φ[α0+Abs[delta]];

   (*Determining Whether the Inicrement Is Positive or Negative*)
   Which[
    y1==y2,
    a=α0-Abs[delta];
    b=α0;
    Goto[end];,
    y2==y3,
    a=α0;
    b=α0+Abs[delta];
    Goto[end];,
    y1==y3||(y1>y2&&y2<y3),
    a=α0-Abs[delta];
    b=α0+Abs[delta];
    Goto[end];
    ];

   Which[
     y1>y2&&y2>y3,
     increment=Abs[delta];,
     y1<y2&&y2<y3,
     increment=-Abs[delta];
     ]

    (* Starting the Algorithm *)
    Do[
     αα[0]=α0;
     αα[k+1]=αα[k]+2^k*increment;

     Which[
      φ[αα[k]]<φ[αα[k+1]],(* Evidently, it is impossible the condition to hold for k=0 *)
      a=αα[k-1];
      b= αα[k+1];
      Break [],
```

```
      k>50,
      Print["After 50 iterations the bounding phase method can not braketing the min of
alpha"];
      Exit[]
      ];,
     {k,0,∞}
     ];

   Label[end];

   If[
    a>b,
    {a,b}={b,a}
    ];

   (* Golden Section Search Method *)
   (* Initiating Required Variables*)
   a0=a;
   b0=b;
   anew=a;
   bnew=b;

   If[
    a0==b0,
    αstar=a;
    Goto[final]
    ];

   (* Starting the Algorithm *)
   Do[
    (* Normalize the Variable α *)
    anorm=(anew-a)/(b-a);
    bnorm=(bnew-a)/(b-a);

    lnorm=bnorm-anorm;

    α1norm=anorm+0.382*lnorm;
    α2norm=bnorm-0.382*lnorm;

    α1=α1norm(b0-a0)+a0;
    α2=α2norm(b0-a0)+a0;

    φ1=φ[α1];
    φ2=φ[α2];

    Which [
     φ1>φ2,
     anew=α1(*move lower bound to α1*);,
     φ1<φ2,
     bnew=α2(*move upper bound to α2*);,
     φ1==φ2,
     anew=α1(*move lower bound to α1*);
     bnew=α2(*move upper bound to α2*);
     ];

    αstar=0.5*(anew+bnew);

    If[
     Abs[lnorm]<epsilon,
     Break[]
     ];,
```

```
     {k,1,∞}
     ];

    Label[final];

    (* Final Result *)
    N[αstar]
    ];

α00=2;(* The intial point of α; for example 2*)
delt0=1;(* The parameter delta of Bounding Phase Method; for example 1 *)
eps0=0.01;(* The accuracy of the Golden Section Search Method; for example 0.01 *)

gradfx[x_,y_]=Grad[f[{x,y}],{x,y}];

(* Main Loop *)
Do[
  gradfx0=gradfx[x0[[1]],x0[[2]]];

  φ[α_]=f[x0-α*gradfx0];
  α1=unidirectionalsearch[α00,delt0,eps0];
  x1=x0-α1*gradfx0;
  error=Norm[(x1-x0)]/Max[1,Norm[x0]];

  lii=k;
  result[k]=N[{k,α1,Row[x0,","],Row[x1,","],f[x0],f[x1],error}];
  plotresult[k]=N[{x0,x1}];

  If[
   error<epsilonSD||k>50||Norm[gradfx0]<epsilonSD,
   Break[],
   x0=x1;
   ];,
  {k,1,∞}
  ];

(* Final Result *)

If[
  lii==50,
  Print[" After 50 iterations the mimimum point around the intial point is " ,N[x1], "\nThe
solution is (approximately) ", N[f[x1]]];,
 Print["The solution is x= ",  N[x1],"\nThe solution is (approximately)= ", N[f[x1]]];
  ]

(* Results of Each Iteration *)
table=TableForm[
  Table[
   result[i],
   {i,1,lii}
   ],
  TableHeadings->{None,{"No. of iters.","α1","x0","x1","f[x0]","f[xl]","error"}}
  ]

Export["example81.xls",table,"XLS"];

(* Data Visualization *)
(* Domain of Varibles*)
xleft=domainx[[1]];
xright=domainx[[2]];
ydown=domainy[[1]];
```

```
yup=domainy[[2]];

(* 3D+ Contour Plot *)
plot1=Plot3D[
   f[{x,y}],
   {x,xleft,xright},
   {y,ydown,yup},
   ClippingStyle->None,
   MeshFunctions->{#3&},
   Mesh->15,
   MeshStyle->Opacity[.5],
   MeshShading->{{Opacity[.3],Blue},{Opacity[.8],Orange}},
   Lighting->"Neutral"
   ];
slice=SliceContourPlot3D[
   f[{x,y}],
   z==0,
   {x,xleft,xright},
   {y,ydown,yup},
   {z,-1,1},
   Contours->15,
   Axes->False,
   PlotPoints->50,
   PlotRangePadding->0,
   ColorFunction->"Rainbow"
   ];
Show[
 plot1,
 slice,
 PlotRange->All,
 BoxRatios->{1,1,.6},
 FaceGrids->{Back,Left}
 ]

(* Contour Plot with Step Iterations *)
ContourPlot[
 f[{x,y}],
 {x,xleft,xright},
 {y,ydown,yup},
 LabelStyle->Directive[Black,14],
 ColorFunction->"Rainbow",
 PlotLegends->Automatic,
 Contours->10,
 Epilog-
>{PointSize[0.015],Green,Line[Flatten[Table[plotresult[i],{i,1,lii}],1]],Red,Point[Flatten[T
able[plotresult[i],{i,1,lii}],1]]}
 ]

(* Data Manipulation *)
Manipulate[
 ContourPlot[
  f[{x,y}],
  {x,xleft,xright},
  {y,ydown,yup},
  LabelStyle->Directive[Black,16],
  ColorFunction->"Rainbow",
  PlotLegends->Automatic,
  Contours->10,
  Epilog->{
    PointSize[0.015],
    Yellow,
```

```
    Arrow[{plotresult[i][[1]],plotresult[i][[2]]}],
    Red,
    Point[Flatten[Table[plotresult[j],{j,1,i}],1]],
    Green,
    Line[Flatten[Table[plotresult[j],{j,1,i}],1]]
    }
  ],
 {i,1,lii,1}
 ]
```

### 8.1.2. Steepest Descent Without Line Search (Formula 1)

If the Hessian of $f|\mathbf{x}_k\rangle$ is available, the value of $\alpha$, $\alpha_k$, that minimizes $f|\mathbf{x}_k+\alpha\mathbf{d}\rangle$ can be determined by using an analytical method. If $\mathbf{H}_k=\mathbf{H}_f|\mathbf{x}_k\rangle=\nabla^2 f|\mathbf{x}_k\rangle$ and $|\mathbf{g}_k\rangle=\nabla f|\mathbf{x}_k\rangle$, $|\boldsymbol{\delta}_k\rangle=\alpha|\mathbf{d}_k\rangle$ , the Taylor series yields

$$f|\mathbf{x}_k+\boldsymbol{\delta}_k\rangle \approx f|\mathbf{x}_k\rangle+\langle\boldsymbol{\delta}_k|\mathbf{g}_k\rangle+\frac{1}{2}\langle\boldsymbol{\delta}_k|\mathbf{H}_k|\boldsymbol{\delta}_k\rangle, \tag{8.15}$$

and if $|\mathbf{d}_k\rangle$ is the steepest-descent direction, i.e.,

$$|\boldsymbol{\delta}_k\rangle=-\alpha|\mathbf{g}_k\rangle, \tag{8.16}$$

we obtain

$$f|\mathbf{x}_k-\alpha\mathbf{g}_k\rangle \approx f|\mathbf{x}_k\rangle-\alpha\langle\mathbf{g}_k|\mathbf{g}_k\rangle+\frac{1}{2}\alpha^2\langle\mathbf{g}_k|\mathbf{H}_k|\mathbf{g}_k\rangle. \tag{8.17}$$

By differentiating and setting the result to zero, we get

$$\frac{d}{d\alpha}f|\mathbf{x}_k-\alpha\mathbf{g}_k\rangle \approx -\langle\mathbf{g}_k|\mathbf{g}_k\rangle+\alpha\langle\mathbf{g}_k|\mathbf{H}_k|\mathbf{g}_k\rangle=0, \tag{8.18}$$

or

$$\alpha=\alpha_k\approx\frac{\langle\mathbf{g}_k|\mathbf{g}_k\rangle}{\langle\mathbf{g}_k|\mathbf{H}_k|\mathbf{g}_k\rangle}. \tag{8.19}$$

Note that the value of $\alpha_k$ is positive if $\langle\mathbf{g}_k|\mathbf{H}_k|\mathbf{g}_k\rangle>0$, which is guaranteed to hold if $\mathbf{H}_k$ is positive definite. Now if we assume that $\alpha=\alpha_k$ minimizes $f|\mathbf{x}_k-\alpha\mathbf{d}_k\rangle$, (8.2) can be expressed as

$$|\mathbf{x}_{k+1}\rangle=|\mathbf{x}_k\rangle-\frac{\langle\mathbf{g}_k|\mathbf{g}_k\rangle}{\langle\mathbf{g}_k|\mathbf{H}_k|\mathbf{g}_k\rangle}|\mathbf{g}_k\rangle. \tag{8.20}$$

The accuracy of $\alpha_k$ will depend heavily on the magnitude of $|\boldsymbol{\delta}_k\rangle$ since the quadratic approximation of the Taylor series is valid only in the neighborhood of the point $|\mathbf{x}_k\rangle$. At the start of the optimization, $\|\boldsymbol{\delta}_k\|$ will be relatively large and so $\alpha_k$ will be inaccurate. Nevertheless, the reduction will be achieved in $f|\mathbf{x}\rangle$ since $f|\mathbf{x}_k+\alpha\mathbf{d}_k\rangle$ is minimized in the steepest-descent direction. As the solution is approached, $\|\boldsymbol{\delta}_k\|$ is decreased, and consequently, the accuracy of $\alpha_k$ will progressively be improved, and the maximum reduction in $f|\mathbf{x}\rangle$ will eventually be achieved in each iteration. For quadratic functions, (8.15) is satisfied with the equal sign, and hence $\alpha=\alpha_k$ yields the maximum reduction in $f|\mathbf{x}\rangle$ in every iteration.

### 8.1.3. Steepest Descent Without Line Search (Formula 2)

If the Hessian is not available, the value of $\alpha_k$ can be determined by calculating $f|\mathbf{x}\rangle$ at points $|\mathbf{x}_k\rangle$ and $|\mathbf{x}_k\rangle-\hat{\alpha}|\mathbf{g}_k\rangle$ where $\hat{\alpha}$ is an estimate of $\alpha_k$. If

$$f_k=f|\mathbf{x}_k\rangle \text{ and } \hat{f}=f|\mathbf{x}_k-\hat{\alpha}\mathbf{g}_k\rangle. \tag{8.21}$$

Equation (8.17) gives

$$\hat{f}\approx f_k-\hat{\alpha}\langle\mathbf{g}_k|\mathbf{g}_k\rangle+\frac{1}{2}\hat{\alpha}^2\langle\mathbf{g}_k|\mathbf{H}_k|\mathbf{g}_k\rangle, \tag{8.22}$$

or

$$\langle \mathbf{g}_k|\mathbf{H}_k|\mathbf{g}_k\rangle \approx \frac{2\left(\hat{f} - f_k + \hat{\alpha}\langle \mathbf{g}_k|\mathbf{g}_k\rangle\right)}{\hat{\alpha}^2}. \tag{8.23}$$

Now from (8.19) and (8.23), we get

$$\alpha_k \approx \frac{\hat{\alpha}^2\langle \mathbf{g}_k|\mathbf{g}_k\rangle}{2\left(\hat{f} - f_k + \hat{\alpha}\langle \mathbf{g}_k|\mathbf{g}_k\rangle\right)}. \tag{8.24}$$

Note that (8.24) yields a positive $\alpha_k$ only if $\left(\hat{f} - f_k + \hat{\alpha}\langle \mathbf{g}_k|\mathbf{g}_k\rangle\right)$ is positive. From (8.23), this is guaranteed to be true if $\mathbf{H}_k$ is positive definite. A suitable value for $\hat{\alpha}$ is the optimum value of $\alpha$ in the previous iteration, namely, $\alpha_{k-1}$. For the first iteration, the value $\hat{\alpha} = 1$ can be used. An algorithm that uses an estimated step size instead of a line search is as follows:

**Algorithm**

**Step 1:** Input $|\mathbf{x}_1\rangle$ and initialize the tolerance $\varepsilon$.
Set $k = 1$ and $\alpha_0 = 1$.
Compute $f_1 = f|\mathbf{x}_1\rangle$.

**Step 2:** Compute $|\mathbf{g}_k\rangle$.

**Step 3:** Set $|\mathbf{d}_k\rangle = -|\mathbf{g}_k\rangle$ and $\hat{\alpha} = \alpha_{k-1}$.
Compute $\hat{f} = f(|\mathbf{x}_k\rangle - \hat{\alpha}|\mathbf{g}_k\rangle)$.
Compute $\alpha_k$ from $\alpha_k \approx \frac{\hat{\alpha}^2\langle \mathbf{g}_k|\mathbf{g}_k\rangle}{2(\hat{f} - f_k + \hat{\alpha}\langle \mathbf{g}_k|\mathbf{g}_k\rangle)}$.

**Step 4:** Set $|\mathbf{x}_{k+1}\rangle = |\mathbf{x}_k\rangle + \alpha_k|\mathbf{d}_k\rangle$ and calculate $f_{k+1} = f|\mathbf{x}_{k+1}\rangle$.

**Step 5:** If $\|\alpha_k \mathbf{d}_k\| < \varepsilon$, then do:
Output $|\mathbf{x}^*\rangle = |\mathbf{x}_{k+1}\rangle$ and $f|\mathbf{x}^*\rangle = f_{k+1}$, and stop.
Otherwise, set $k = k + 1$ and repeat from Step 2.

***Example 8.3***

Consider the problem:

$$\text{Minimize } f|\mathbf{x}\rangle = \frac{x^2}{4} + \frac{y^2}{25},$$

using, $|\mathbf{x}_0\rangle = (5,5)^T$, and $\epsilon = 0.01$.

***Solution***

The 3D and contour plots of the function are shown in Figure 8.4. The plots show that the minimum lies at $|\mathbf{x}^*\rangle = (-0.00009,0.003)^T$, $f|\mathbf{x}^*\rangle = 5.82 \times 10^{-7}$.

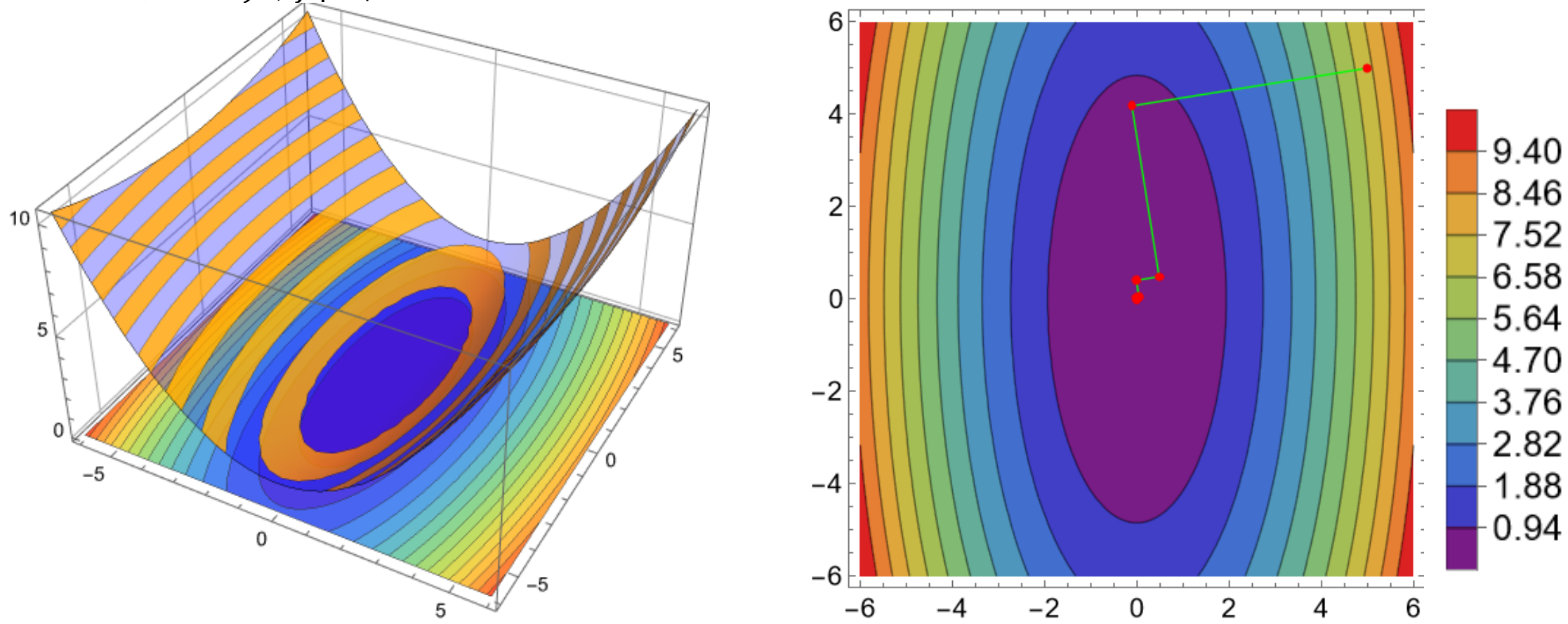

**Figure 8.4.** The results of 7 iterations of the steepest descent method (Formula 2) for $f|\mathbf{x}\rangle = x^2/4 + y^2/25$.

From Table 8.2, after 7 iterations, the optimal condition is attained. The results were produced by Mathematica code 8.2.

**Table 8.2.**

| iters. | $\alpha1$ | $\alpha2$ | $x1$ | $x2$ | $f(x1)$ | $f(x2)$ | error |
|---|---|---|---|---|---|---|---|
| 1 | 1 | 2.042833 | {5., 5.} | {-0.107, 4.183} | 7.25 | 0.702722 | 5.172039 |
| 2 | 2.042833 | 11.05172 | {-0.107, 4.183} | {0.485, 0.485} | 0.702722 | 0.068113 | 3.74527 |
| 3 | 11.05172 | 2.042833 | {0.485, 0.485} | {-0.010, 0.405} | 0.068113 | 0.006602 | 0.501311 |
| 4 | 2.042833 | 11.05172 | {-0.010, 0.405} | {0.047, 0.047} | 0.006602 | 0.00064 | 0.363018 |
| 5 | 11.05172 | 2.042833 | {0.047, 0.047} | {-0.001, 0.039} | 0.00064 | 6.2E-05 | 0.048591 |
| 6 | 2.042833 | 11.05172 | {-0.001, 0.039} | {0.005, 0.005} | 6.2E-05 | 6.01E-06 | 0.035186 |
| 7 | 11.05172 | 2.042833 | {0.005, 0.0046} | {-0.00009, 0.003} | 6.01E-06 | 5.83E-07 | 0.00471 |

***Mathematica Code 8.2*** Steepest Descent Without Line Search (Formula 2)

```
(* Steepest Descent Method Without Unidirectional Search(Formula 2) *)

(*
Notations:
x0         :Intial vector
epsilonSD   :Small number to check the accuracy of the Steepest Descent Search Method
f[x,y]      :Objective function
lii         :The last iteration index
result[k]   :The results of iteration k
*)

(* Taking Initial Inputs from User *)
x1=Input["Enter the intial point in the format {x, y}; for example {5,5} "] ;
epsilonSD=Input["Please enter accuracy of the steepest descent method (Formula 2); for
example 0.01 "];

domainx=Input["Please enter domain of x variable for 3D and contour plots; for example {-
6,6}"];
domainy=Input["Please enter domain of y variable for 3D and contour plots; for example {-
6,6}"];

If[
  epsilonSD<=0,
  Beep[];
  MessageDialog["The value epsilonSD has to be postive number: "];
  Exit[];
  ];

(* Taking the Function from User *)
f[{x_,y_}] = Evaluate[Input["Please input a function of x and y to find the minimum "]];
(* For example : x^2/4+y^2/25 *)

α1=1;
gradfx[x_,y_]=Grad[f[{x,y}],{x,y}];

(* Main Loop *)
Do[
  f1=f[x1];
  gradfx1=gradfx[x1[[1]],x1[[2]]];
  d=-1*gradfx1;
  fbar=f[x1+α1*d];
  α2=((α1)^2*(Dot[gradfx1,gradfx1]))/(2*(fbar-f1+α1*(Dot[gradfx1,gradfx1])));
```

```
  x2=x1+α2*d;

  error=Norm[α2*d];

  lii=k;
  result[k]=N[{k,α1,α2,Row[x1,","],Row[x2,","],f[x1],f[x2],error}];
  plotresult[k]=N[{x1,x2}];

  If[
   error<epsilonSD||k>50,
   Break[],
   α1=α2;
   x1=x2;
   ];,
  {k,1,∞}
  ];

(* Final Result *)
If[
  lii==50,
  Print[" After 50 iterations the mimimum point around the intial point is " ,N[x2], "\nThe 
solution is (approximately) ", N[f[x2]]];,
 Print["The solution is x= ",  N[x2],"\nThe solution is (approximately)= ", N[f[x2]]];
  ]

(* Results of Each Iteration *)
table=TableForm[
  Table[
   result[i],
   {i,1,lii}
   ],
  TableHeadings->{None,{"No. of iters.","α1","α2","x1","x2","f[x1]","f[x2]","error"}}
  ]

Export["example82.xls",table,"XLS"];

(* Data Visualization *)
(* Domain of Varibles*)
xleft=domainx[[1]];
xright=domainx[[2]];
ydown=domainy[[1]];
yup=domainy[[2]];

(* 3D+ Contour Plot *)
plot1=Plot3D[
   f[{x,y}],
   {x,xleft,xright},
   {y,ydown,yup},
   ClippingStyle->None,
   MeshFunctions->{#3&},
   Mesh->15,
   MeshStyle->Opacity[.5],
   MeshShading->{{Opacity[.3],Blue},{Opacity[.8],Orange}},
   Lighting->"Neutral"
   ];
slice=SliceContourPlot3D[
   f[{x,y}],
   z==0,
   {x,xleft,xright},
   {y,ydown,yup},
   {z,-1,1},
```

```
   Contours->15,
   Axes->False,
   PlotPoints->50,
   PlotRangePadding->0,
   ColorFunction->"Rainbow"
   ];
Show[
 plot1,
 slice,
 PlotRange->All,
 BoxRatios->{1,1,.6},
 FaceGrids->{Back,Left}
 ]

(* Contour Plot with Step Iterations *)
ContourPlot[
 f[{x,y}],
 {x,xleft,xright},
 {y,ydown,yup},
 LabelStyle->Directive[Black,16],
 ColorFunction->"Rainbow",
 PlotLegends->Automatic,
 Contours->10,
 Epilog-
>{PointSize[0.015],Green,Line[Flatten[Table[plotresult[i],{i,1,lii}],1]],Red,Point[Flatten[T
able[plotresult[i],{i,1,lii}],1]]}
 ]

(* Data Manipulation *)
Manipulate[
 ContourPlot[
  f[{x,y}],
  {x,xleft,xright},
  {y,ydown,yup},
  LabelStyle->Directive[Black,14],
  ColorFunction->"Rainbow",
  PlotLegends->Automatic,
  Contours->10,
  Epilog->{
    PointSize[0.015],
    Yellow,
    Arrow[{plotresult[i][[1]],plotresult[i][[2]]}],
    Red,
    Point[Flatten[Table[plotresult[j],{j,1,i}],1]],
    Green,
    Line[Flatten[Table[plotresult[j],{j,1,i}],1]]
    }
  ],
 {i,1,lii,1}
 ]
```

### 8.1.4. Steepest Descent Without Line Search (Formula 3: Barzilai–Borwein Two-Point Formulas)

If the Hessian in (8.15) is positive definite, then the convex quadratic function at the right-hand side of (8.15) has a unique minimizer which can be obtained as

$$|\boldsymbol{\delta}_k\rangle = -\mathbf{H}_k^{-1}|\mathbf{g}_k\rangle, \tag{8.25}$$

by setting the gradient of the quadratic function to zero. On comparing the $|\boldsymbol{\delta}_k\rangle$ in (8.25) with that obtained in the steepest-descent method, i.e.,

$$|\boldsymbol{\delta}_k\rangle = -\alpha_k|\mathbf{g}_k\rangle, \tag{8.26}$$

we note that $\alpha_k$ in (8.26) corresponds to $\mathbf{H}_k^{-1}$ in (8.25). Thus, we can approximate the inverse Hessian matrix $\mathbf{H}_k^{-1}$ by $\alpha_k\mathbf{I}$ or, equivalently, we can approximate $\mathbf{H}_k$ by $\beta_k\mathbf{I}$ I with $\beta_k = 1/\alpha_k$. Hence, we can write

$$|\mathbf{g}_k\rangle \approx |\mathbf{g}_{k-1}\rangle + \mathbf{H}_k|\boldsymbol{\delta}_{k-1}\rangle, \tag{8.27}$$

which implies that

$$\mathbf{H}_k^{-1}|\boldsymbol{\gamma}_{k-1}\rangle \approx |\boldsymbol{\delta}_{k-1}\rangle, \tag{8.28}$$

and

$$\mathbf{H}_k|\boldsymbol{\delta}_{k-1}\rangle \approx |\boldsymbol{\gamma}_{k-1}\rangle, \tag{8.29}$$

where

$$|\boldsymbol{\delta}_{k-1}\rangle = |\mathbf{x}_k\rangle - |\mathbf{x}_{k-1}\rangle, \qquad |\boldsymbol{\gamma}_{k-1}\rangle = |\mathbf{g}_k\rangle - |\mathbf{g}_{k-1}\rangle, \tag{8.30}$$

In other words, $|\boldsymbol{\delta}_{k-1}\rangle$ and $|\boldsymbol{\gamma}_{k-1}\rangle$ are the differences between the two points $|\mathbf{x}_k\rangle$ and $|\mathbf{x}_{k-1}\rangle$ and the gradients at these two points. From (8.28), it follows that $\alpha_k\mathbf{I}$ approximates $\mathbf{H}_k^{-1}$ if $\alpha = \alpha_k$ is a solution of the problem

$$\underset{\alpha}{\text{minimize}}\ \|\alpha\boldsymbol{\gamma}_{k-1} - \boldsymbol{\delta}_{k-1}\|, \tag{8.31}$$

which is a single-variable least-squares problem that involves the objective function

$$\|\boldsymbol{\gamma}_{k-1}\|^2\alpha^2 - 2\langle\boldsymbol{\delta}_{k-1}|\boldsymbol{\gamma}_{k-1}\rangle\alpha + \|\boldsymbol{\delta}_{k-1}\|^2. \tag{8.32}$$

By setting the derivative of the above function with respect to $\alpha$ to zero, the minimizer of (8.31) is obtained as

$$\alpha_k = \frac{\langle\boldsymbol{\delta}_{k-1}|\boldsymbol{\gamma}_{k-1}\rangle}{\langle\boldsymbol{\gamma}_{k-1}|\boldsymbol{\gamma}_{k-1}\rangle}. \tag{8.33}$$

Alternatively, from (8.29) it follows that $\beta_k\mathbf{I}$ approximates $\mathbf{H}_k$ if $\beta = \beta_k$ is a solution of the single-variable least-squares optimization problem

$$\underset{\beta}{\text{minimize}}\ \|\beta\boldsymbol{\delta}_{k-1} - \boldsymbol{\gamma}_{k-1}\|. \tag{8.34}$$

The minimizer of the problem in (8.34) is given by

$$\beta_k = \frac{\langle\boldsymbol{\delta}_{k-1}|\boldsymbol{\delta}_{k-1}\rangle}{\langle\boldsymbol{\delta}_{k-1}|\boldsymbol{\gamma}_{k-1}\rangle}, \tag{8.35}$$

and, in effect, an alternative formula for $\alpha_k$ is obtained as

$$\alpha_k = \frac{\langle\boldsymbol{\delta}_{k-1}|\boldsymbol{\delta}_{k-1}\rangle}{\langle\boldsymbol{\delta}_{k-1}|\boldsymbol{\gamma}_{k-1}\rangle}. \tag{8.36}$$

The two formulas for $\alpha_k$ (8.33) and (8.36) can be used to implement Step 3 of the above algorithm, which uses an estimated step size [4-6].

***Example 8.4***

Consider the problem:

$$\text{Minimize } f|\mathbf{x}\rangle = \frac{x^2}{4} + \frac{y^2}{25},$$

using, $|\mathbf{x}_0\rangle = (5,5)$, and $\epsilon = 0.01$.

***Solution***

The 3D and contour plots of the function are shown in Figure 8.5. The plots show that the minimum lies at $|\mathbf{x}^*\rangle = (2.00952 \times 10^{-10}, 0.00002)^T$, $f|\mathbf{x}^*\rangle = 1.91749 \times 10^{-11}$.

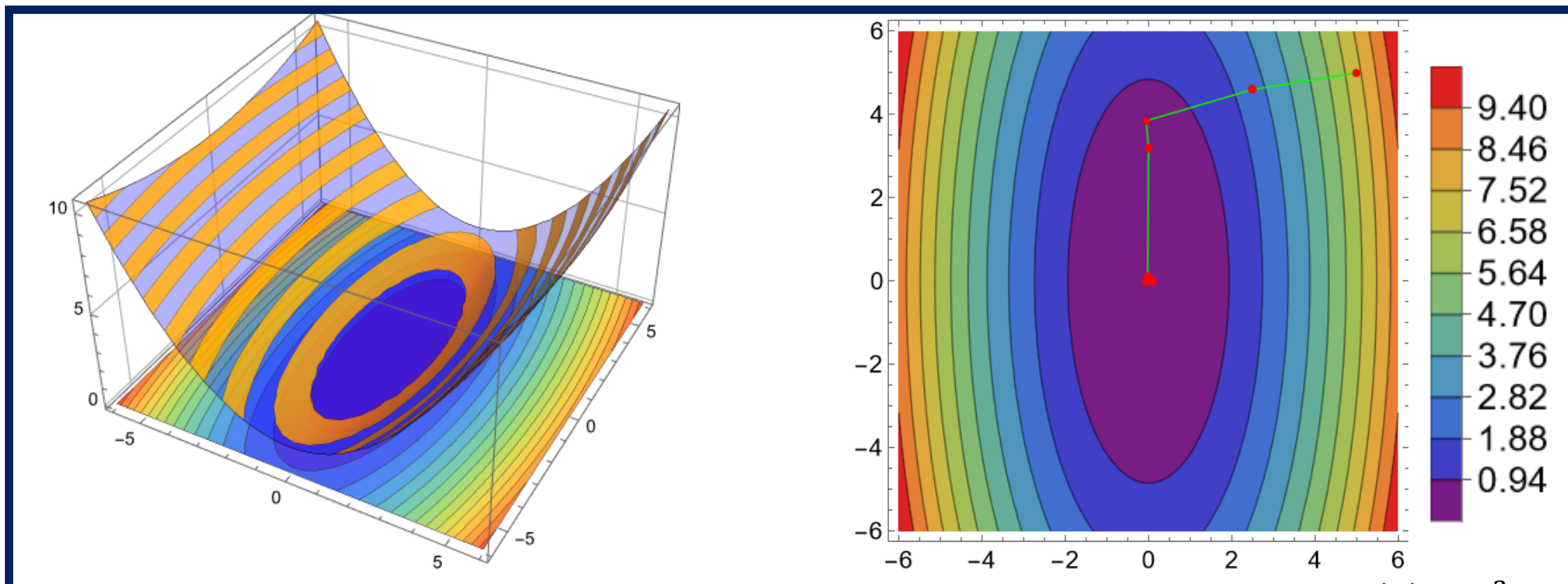

**Figure 8.5.** The results of 6 iterations of steepest descent, Barzilai–Borwein Formula, for $f|\mathbf{x}\rangle = x^2/4 + y^2/25$.

From Table 8.3, after 6 iterations, the optimal condition is attained. The results were produced by Mathematica code 8.3.

**Table 8.3.**

| No. iters. | α2 | α3 | x2 | x3 | $f(x2)$ | $f(x3)$ | error |
|---|---|---|---|---|---|---|---|
| 0 | 1 | 2.042833 | {5., 5.} | {2.5, 4.6} | 7.25 | 2.4089 | 2.66190 |
| 1 | 2.042833 | 2.143616 | {2.5, 4.6} | {-0.053, 3.848} | 2.4089 | 0.59307 | 0.66242 |
| 2 | 2.143616 | 12.02617 | {-0.054, 3.848} | {0.004, 3.188} | 0.59307 | 0.40662 | 3.06754 |
| 3 | 12.02617 | 12.49627 | {0.004, 3.188} | {-0.019, 0.121} | 0.40662 | 0.00068 | 0.17059 |
| 4 | 12.49627 | 3.456497 | {-0.019, 0.121} | {0.101, 0.00003} | 0.00068 | 0.00256 | 0.17481 |
| 5 | 3.456497 | 2 | {0.101, 0.00004} | {-0.074, 0.00003} | 0.00256 | 0.00136 | 0.07366 |
| 6 | 2 | 2 | {-0.074, 0.00003} | {2.009 E-10, 0.00002} | 0.00136 | 1.92E-11 | 3.5E-06 |

***Mathematica Code 8.3*** Steepest Descent Without Line Search (Barzilai–Borwein Formula)

```
(* Steepest Descent Method Without Unidirectional Search(Barzilai-Borwein Two-Point Formula)
*)

(*
Notations:
x0        :Intial vector
epsilonSD :Small number to check the accuracy of the Steepest Descent Search Method
f[x,y]    :Objective function
lii       :The last iteration index
result[k] :The results of iteration k
*)

(* Taking Initial Inputs from User *)
x1=Input["Enter the intial point in the format {x, y}; for example {5,5} "] ;
epsilonSD=Input["Please enter accuracy of the Steepest Descent Method (Barzilai-Borwein Two-
Point Formula); for example 0.01 "];

domainx=Input["Please enter domain of x variable for 3D and contour plots; for example {-
6,6}"];
domainy=Input["Please enter domain of y variable for 3D and contour plots; for example {-
6,6}"];

If[
  epsilonSD<=0,
  Beep[];
  MessageDialog["The value of epsilonSD has to be postive number: "];
```

```
  Exit[];
  ];

(* Taking the Function from User *)
f[{x_,y_}] = Evaluate[Input["Please input a function of x and y to find the minimum "]];
(* For example :x^2/4+y^2/25 *)

gradfx[x_,y_]=Grad[f[{x,y}],{x,y}];

(* Iteration Number Zero *)
α1=1;
gradfx1=gradfx[x1[[1]],x1[[2]]];
x2=x1-α1*gradfx1;
gradfx2=gradfx[x2[[1]],x2[[2]]];
δ=x2-x1;
γ=gradfx2-gradfx1;
α2=Dot[δ,δ]/Dot[δ,γ];
result[0]=N[{0,α1,α2,Row[x1,","],Row[x2,","],f[x1],f[x2],Norm[α2*gradfx2]}];
plotresult[0]=N[{x1,x2}];

(* Main Loop *)
Do[
  x3=x2-α2*gradfx2;
  gradfx3=gradfx[x3[[1]],x3[[2]]];

  δ=x3-x2;
  γ=gradfx3-gradfx2;
  α3=Dot[δ,δ]/Dot[δ,γ];
  error=Norm[α3*gradfx3];

  lii=k;
  result[k]=N[{k,α2,α3,Row[x2,","],Row[x3,","],f[x2],f[x3],error}];
  plotresult[k]=N[{x2,x3}];

  If[
   error<epsilonSD||k>50,
   Break[],
   α2=α3;
   x2=x3;
   gradfx2=gradfx3;
   ];,
  {k,1,∞}
  ];

(* Final Result *)
If[
  lii==50,
  Print[" After 50 iterations the mimimum point around the intial point is " ,N[x3], "\nThe solution is (approximately) ", N[f[x3]]];,
 Print["The solution is x= ",  N[x3],"\nThe solution is (approximately)= ", N[f[x3]]];
  ]

(* Results of Each Iteration *)
table=TableForm[
  Table[
   result[i],
   {i,0,lii}
   ],
  TableHeadings->{None,{"No. of iters.","α2","α3","x2","x3","f[x2]","f[x3]","error"}}
  ]
Export["example83.xls",table,"XLS"];
```

```
(* Data Visualization *)
(* Domain of Varibles*)
xleft=domainx[[1]];
xright=domainx[[2]];
ydown=domainy[[1]];
yup=domainy[[2]];

(* 3D+ Contour Plot *)
plot1=Plot3D[
   f[{x,y}],
   {x,xleft,xright},
   {y,ydown,yup},
   ClippingStyle->None,
   MeshFunctions->{#3&},
   Mesh->15,
   MeshStyle->Opacity[.5],
   MeshShading->{{Opacity[.3],Blue},{Opacity[.8],Orange}},
   Lighting->"Neutral"
   ];
slice=SliceContourPlot3D[
   f[{x,y}],
   z==0,
   {x,xleft,xright},
   {y,ydown,yup},
   {z,-1,1},
   Contours->15,
   Axes->False,
   PlotPoints->50,
   PlotRangePadding->0,
   ColorFunction->"Rainbow"
   ];
Show[
 plot1,
 slice,
 PlotRange->All,
 BoxRatios->{1,1,.6},
 FaceGrids->{Back,Left}
 ]

(* Contour Plot with Step Iterations *)
ContourPlot[
 f[{x,y}],
 {x,xleft,xright},
 {y,ydown,yup},
 LabelStyle->Directive[Black,16],
 ColorFunction->"Rainbow",
 PlotLegends->Automatic,
 Contours->10,
 Epilog-
>{PointSize[0.015],Green,Line[Flatten[Table[plotresult[i],{i,0,lii}],1]],Red,Point[Flatten[T
able[plotresult[i],{i,0,lii}],1]]}
 ]

(* Data Manipulation *)
Manipulate[
 ContourPlot[
  f[{x,y}],
  {x,xleft,xright},
  {y,ydown,yup},
  LabelStyle->Directive[Black,14],
  ColorFunction->"Rainbow",
```

```
  PlotLegends->Automatic,
  Contours->10,
  Epilog->{
    PointSize[0.015],
    Yellow,
    Arrow[{plotresult[i][[1]],plotresult[i][[2]]}],
    Red,
    Point[Flatten[Table[plotresult[j],{j,0,i}],1]],
    Green,
    Line[Flatten[Table[plotresult[j],{j,0,i}],1]]
    }
  ],
 {i,0,lii,1}
 ]
```

### 8.1.5. Momentum

Gradient descent will take a long time to traverse a nearly flat surface (Regions that are nearly flat have gradients with small magnitudes and can thus require many iterations of gradient descent to traverse.). Allowing momentum to accumulate is one way to speed progress. We can modify gradient descent to incorporate momentum. The momentum update equations are [7]:

$$|\mathbf{v}_{k+1}\rangle = \beta|\mathbf{v}_k\rangle - \alpha|\mathbf{g}_k\rangle, \tag{8.37}$$
$$|\mathbf{x}_{k+1}\rangle = |\mathbf{x}_k\rangle + |\mathbf{v}_{k+1}\rangle. \tag{8.38}$$

For $\beta = 0$, we recover gradient descent. Momentum can be interpreted as a ball rolling down a nearly horizontal incline. The ball naturally gathers momentum as gravity causes it to accelerate, just as the gradient causes momentum to accumulate in this descent method.

### 8.1.6. Nesterov Momentum

One issue of momentum is that the steps do not slow down enough at the bottom of a valley and tend to overshoot the valley floor. Nesterov momentum modifies the momentum algorithm to use the gradient at the projected future position [7]:

$$|\mathbf{v}_{k+1}\rangle = \beta|\mathbf{v}_k\rangle - \alpha\nabla f(|\mathbf{x}_k\rangle + \beta|\mathbf{v}_k\rangle), \tag{8.39}$$
$$|\mathbf{x}_{k+1}\rangle = |\mathbf{x}_k\rangle + |\mathbf{v}_{k+1}\rangle. \tag{8.40}$$

## 8.2 Second-order Approximations Methods (Hessian-Based Methods)

This section focuses on the second-order approximations that use the Hessian in multivariate optimization to direct the search [1-3,7].

### 8.2.1. Newton Method

Recall that the method of steepest descent uses only first derivatives (gradients) in selecting a suitable search direction. This strategy is not always the most effective. If higher derivatives are used, the resulting iterative algorithm may perform better than the steepest descent method. Newton method (sometimes called the Newton-Raphson method) uses first and second derivatives and, indeed, does perform better than the steepest descent method if the initial point is close to the minimizer. The idea behind this method is as follows. Given a starting point, we construct a quadratic approximation to the objective function that matches the first and second derivative values at that point. We then minimize the approximate (quadratic) function instead of the original objective function. We use the minimizer of the approximate function as the starting point in the next step and repeat the procedure iteratively. If the objective function is quadratic, then the approximation is exact, and the method yields the true minimizer in one step. If, on the other hand, the objective function is not quadratic, then the approximation will provide only an estimate of the position of the true minimizer. Figure 8.6 illustrates the above idea.

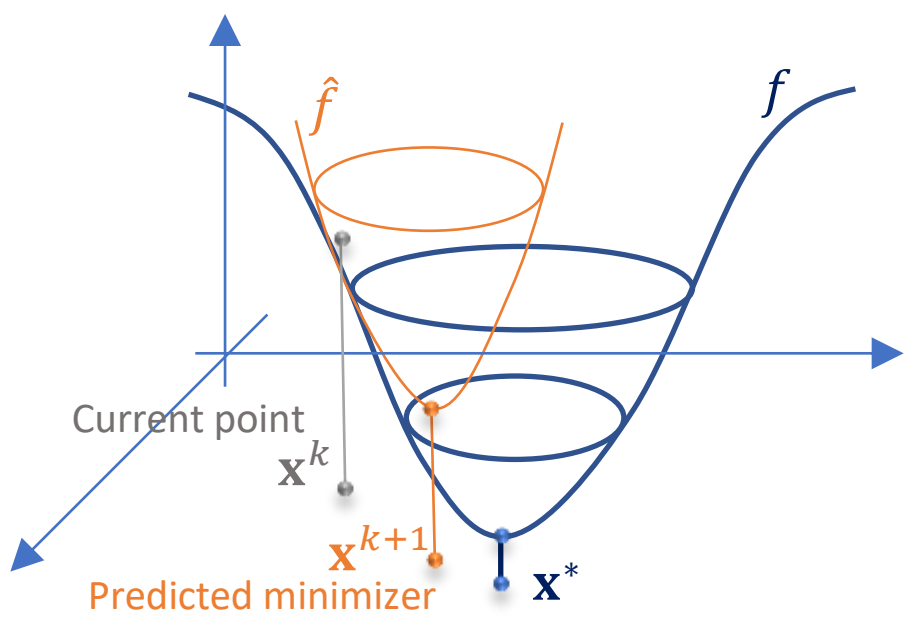

**Figure 8.6.** Quadratic approximation to the objective function using first and second derivatives.

Consider once again the Taylor expansion of the objective:

$$f|\mathbf{x}\rangle = f|\mathbf{x}_k\rangle + \langle \nabla f(\mathbf{x}_k)|\Delta\mathbf{x}\rangle + \frac{1}{2}\langle \Delta\mathbf{x}|\nabla^2 f(\mathbf{x}_k)|\Delta\mathbf{x}\rangle + O(\Delta\mathbf{x}^3). \tag{8.41}$$

We form a quadratic approximation to $f|\mathbf{x}\rangle$ by dropping terms of order 3 and above:

$$\hat{f}(\mathbf{x};\mathbf{x}_k) = f|\mathbf{x}_k\rangle + \langle \nabla f(\mathbf{x}_k)|\Delta\mathbf{x}\rangle + \frac{1}{2}\langle \Delta\mathbf{x}|\nabla^2 f(\mathbf{x}_k)|\Delta\mathbf{x}\rangle, \tag{8.42}$$

where $\hat{f}(\mathbf{x};\mathbf{x}_k)$ denotes an approximating function constructed at $|\mathbf{x}_k\rangle$ which is itself a function of $|\mathbf{x}\rangle$. Now let us use this quadratic approximation of $f|\mathbf{x}\rangle$ to form an iteration sequence by forcing $|\mathbf{x}_{k+1}\rangle$ the next point in the sequence, to be a point where the gradient of the approximation is zero. Therefore,

$$\nabla\hat{f}(\mathbf{x};\mathbf{x}_k) = \nabla f|\mathbf{x}_k\rangle + \langle \nabla^2 f(\mathbf{x}_k)|\Delta\mathbf{x}\rangle = 0, \tag{8.43}$$

so, the search direction becomes

$$|\Delta\mathbf{x}\rangle = |\mathbf{s}_k\rangle = |\mathbf{x}_{k+1}\rangle - |\mathbf{x}_k\rangle = -\frac{\nabla f|\mathbf{x}_k\rangle}{\nabla^2 f|\mathbf{x}_k\rangle}. \tag{8.44}$$

Accordingly, this successive quadratic approximation scheme produces Newton's optimization method:

$$|\mathbf{x}_{k+1}\rangle = |\mathbf{x}_k\rangle - \frac{\nabla f|\mathbf{x}_k\rangle}{\nabla^2 f|\mathbf{x}_k\rangle} = |\mathbf{x}_k\rangle - \mathbf{H}_k^{-1}|\mathbf{g}_k\rangle, \tag{8.45}$$

where, $\mathbf{H}_k = \mathbf{H}_f|\mathbf{x}_k\rangle = \nabla^2 f|\mathbf{x}_k\rangle$ and $|\mathbf{g}_k\rangle = \nabla f|\mathbf{x}_k\rangle$.

This solution exists if and only if the following conditions hold:

(a) The Hessian is nonsingular. (b) The approximation in (8.41) is valid.

It can also be shown that if the matrix $\mathbf{H}_k^{-1}$ is positive-semidefinite, the direction $|\mathbf{s}_k\rangle$ must be descent. But if the matrix $\mathbf{H}_k^{-1}$ is not positive-semidefinite, the direction $|\mathbf{s}_k\rangle$ may or may not be descent, depending on whether the quantity $\langle \mathbf{g}_k|\mathbf{H}_k^{-1}|\mathbf{g}_k\rangle$ is positive or not. Thus, the above search direction may not always guarantee a decrease in the function value in the vicinity of the current point. But the second-order optimality condition suggests that $\nabla^2 f|\mathbf{x}^*\rangle$ be positive-definite for the minimum point. Thus, it can be assumed that the matrix $\nabla^2 f|\mathbf{x}^*\rangle$ is positive-definite in the vicinity of the minimum point (i.e., for $\|\mathbf{x}-\mathbf{x}^*\| < \varepsilon$), and the method is suitable and efficient when the initial point is close to the optimum point.

### Modified Newton Method

Any quadratic function has a Hessian, which is constant for any $|\mathbf{x}\rangle$. If the function has a minimum and the second-order sufficiency conditions for a minimum hold, then $\mathbf{H}$ is positive definite and, therefore, nonsingular at any point $|\mathbf{x}\rangle$. Since any quadratic function is represented exactly by the quadratic approximation of the Taylor series, the solution in (8.44) exists. Furthermore, for any point $|\mathbf{x}\rangle$ one iteration will yield the solution. For nonquadratic functions, the Newton step will often be large when $|\mathbf{x}_0\rangle$ is far from $|\mathbf{x}^*\rangle$, and there is the real possibility of divergence.

It is possible to modify the method in a logical and simple way to ensure descent by adding a line search as in the Cauchy method. That is, we form the sequence of iterates

$$|\mathbf{x}_{k+1}\rangle = |\mathbf{x}_k\rangle + \alpha_k|\mathbf{d}_k\rangle = |\mathbf{x}_k\rangle - \alpha_k \frac{\nabla f|\mathbf{x}_k\rangle}{\nabla^2 f|\mathbf{x}_k\rangle}, \quad (8.46)$$

by choosing $\alpha_k$ such that

$$f|\mathbf{x}_{k+1}\rangle \to \min, \quad (8.47)$$

which ensures that

$$f|\mathbf{x}_{k+1}\rangle \le f|\mathbf{x}_k\rangle, \quad (8.48)$$

where,

$$|\mathbf{d}_k\rangle = -\frac{\nabla f|\mathbf{x}_k\rangle}{\nabla^2 f|\mathbf{x}_k\rangle}. \quad (8.49)$$

This is the modified Newton method, and we find it reliable and efficient when the first and second derivatives are accurately and inexpensively calculated.

| **Algorithm** | **Newton Algorithm** |
|---|---|
| **Step 1:** | Choose a maximum number of iterations $M$ to be performed, an initial point $|\mathbf{x}_0\rangle$, two termination parameters $\epsilon_1$, $\epsilon_2$, and set $k = 0$. |
| **Step 2:** | Compute $\mathbf{g}_k$ and $\mathbf{H}_k$.<br>If $\|\mathbf{g}_k\| \le \epsilon_1$, Terminate;<br>If $\mathbf{H}_k$ is not a positive definite, force it to become a positive definite.<br>Compute $\mathbf{H}_k^{-1}$ and $|\mathbf{d}_k\rangle = -\mathbf{H}_k^{-1}|\mathbf{g}_k\rangle$ . |
| **Step 3:** | Perform a unidirectional search to find $\alpha^{(k)}$ using $\epsilon_2$ such that<br>$f|\mathbf{x}_{k+1}\rangle = f(|\mathbf{x}_k\rangle + \alpha_k|\mathbf{d}_k\rangle)$<br>is minimum. |
| **Step 4:** | Set $|\mathbf{x}_{k+1}\rangle = |\mathbf{x}_k\rangle + \alpha_k|\mathbf{d}_k\rangle$<br>$f_{k+1} = f|\mathbf{x}_{k+1}\rangle$ |
| **Step 5:** | Is $\left\|\alpha^{(k)}\mathbf{d}_k\right\| \le \epsilon_1$? If yes, then output $|\mathbf{x}^*\rangle = |\mathbf{x}_{k+1}\rangle$ and $f|\mathbf{x}^*\rangle = f_{k+1}$, and Terminate;<br>Else set $k = k + 1$ and go to Step 2. |

***Example 8.5***

Use Newton method to minimize the function:

$$f|\mathbf{x}\rangle = (x_1 + 10x_2)^2 + 5(x_3 - x_4)^2 + (x_2 - 2x_3)^4 + 10(x_1 - x_4)^4.$$

Use as the starting point $|\mathbf{x}_0\rangle = (3, -1, 0, 1)^T$. Perform three iterations.

***Solution***

Note that $f|\mathbf{x}_0\rangle = 215$. We have

$$\nabla f|\mathbf{x}\rangle = \begin{pmatrix} 2(x_1 + 10x_2) + 40(x_1 - x_4)^3 \\ 20(x_1 + 10x_2) + 4(x_2 - 2x_3)^3 \\ 10(x_3 - x_4) - 8(x_2 - 2x_3)^3 \\ -10(x_3 - x_4) - 40(x_1 - x_4)^3 \end{pmatrix},$$

and $\mathbf{H}_f|\mathbf{x}\rangle$ is given by

$$\begin{pmatrix} 2 + 120(x_1 - x_4)^2 & 20 & 0 & -120(x_1 - x_4)^2 \\ 20 & 200 + 12(x_2 - 2x_3)^2 & -24(x_2 - 2x_3)^2 & 0 \\ 0 & -24(x_2 - 2x_3)^2 & 10 + 48(x_2 - 2x_3)^2 & -10 \\ -120(x_1 - x_4)^2 & 0 & -10 & 10 + 120(x_1 - x_4)^2 \end{pmatrix}.$$

**Iteration 1.**

$$|\mathbf{g}_0\rangle = (306, -144, -2, -310)^T,$$

$$\mathbf{H}_0 = \begin{pmatrix} 482 & 20 & 0 & -480 \\ 20 & 212 & -24 & 0 \\ 0 & -24 & 58 & -10 \\ -480 & 0 & -10 & 490 \end{pmatrix},$$

$${\mathbf{H}_0}^{-1} = \begin{pmatrix} .1126 & -.0089 & .0154 & .1106 \\ -.0089 & .0057 & .0008 & -.0087 \\ .0154 & .0008 & .0203 & .0155 \\ .1106 & -.0087 & .0155 & .1107 \end{pmatrix},$$

$${\mathbf{H}_0}^{-1}|\mathbf{g}_0\rangle = (1.4127, -0.8413, -0.2540, 0.7460)^T.$$

Hence,

$$|\mathbf{x}_1\rangle = |\mathbf{x}_0\rangle - {\mathbf{H}_0}^{-1}|\mathbf{g}_0\rangle = (1.5873, -0.1587, 0.2540, 0.2540)^T, \qquad f|\mathbf{x}_1\rangle = 31.8.$$

**Iteration 2.**

$$|\mathbf{g}_1\rangle = (94.81, -1.179, 2.371, -94.81)^T,$$

$$\mathbf{H}_1 = \begin{pmatrix} 215.3 & 20 & 0 & -213.3 \\ 20 & 205.3 & -10.67 & 0 \\ 0 & -10.67 & 31.34 & -10 \\ -213.3 & 0 & -10 & 223.3 \end{pmatrix},$$

$${\mathbf{H}_1}^{-1}|\mathbf{g}_1\rangle = (0.5291, -0.0529, -0.0846, 0.0846)^T,$$

Hence,

$$|\mathbf{x}_2\rangle = |\mathbf{x}_1\rangle - {\mathbf{H}_1}^{-1}|\mathbf{g}_1\rangle = (1.0582, -0.1058, 0.1694, 0.1694)^T, \qquad f|\mathbf{x}_2\rangle = 6.28.$$

**Iteration 3.**

$$|\mathbf{g}_2\rangle = (28.09, -0.3475, 0.7031, -28.08)^T,$$

$$\mathbf{H}_2 = \begin{pmatrix} 96.80 & 20 & 0 & -94.80 \\ 20 & 202.4 & -4.744 & 0 \\ 0 & -4.744 & 19.49 & -10 \\ -94.80 & 0 & -10 & 104.80 \end{pmatrix},$$

$$|\mathbf{x}_3\rangle = (0.7037, -0.0704, 0.1121, 0.1111)^T, \qquad f|\mathbf{x}_3\rangle = 1.24.$$

***Example 8.6***

Consider the problem:

$$\text{Minimize } f|\mathbf{x}\rangle = (x^2 + y - 11)^2 + (x + y^2 - 7),$$

using, $|\mathbf{x}_0\rangle = (5,5)^T$, and $\epsilon = 0.01$.

***Solution***

The 3D and contour plots of the function are shown in Figure 8.7. The plots show that the minimum lies at $|\mathbf{x}^*\rangle = (3.294, 0.076)^T$, $f|\mathbf{x}^*\rangle = -3.695$.

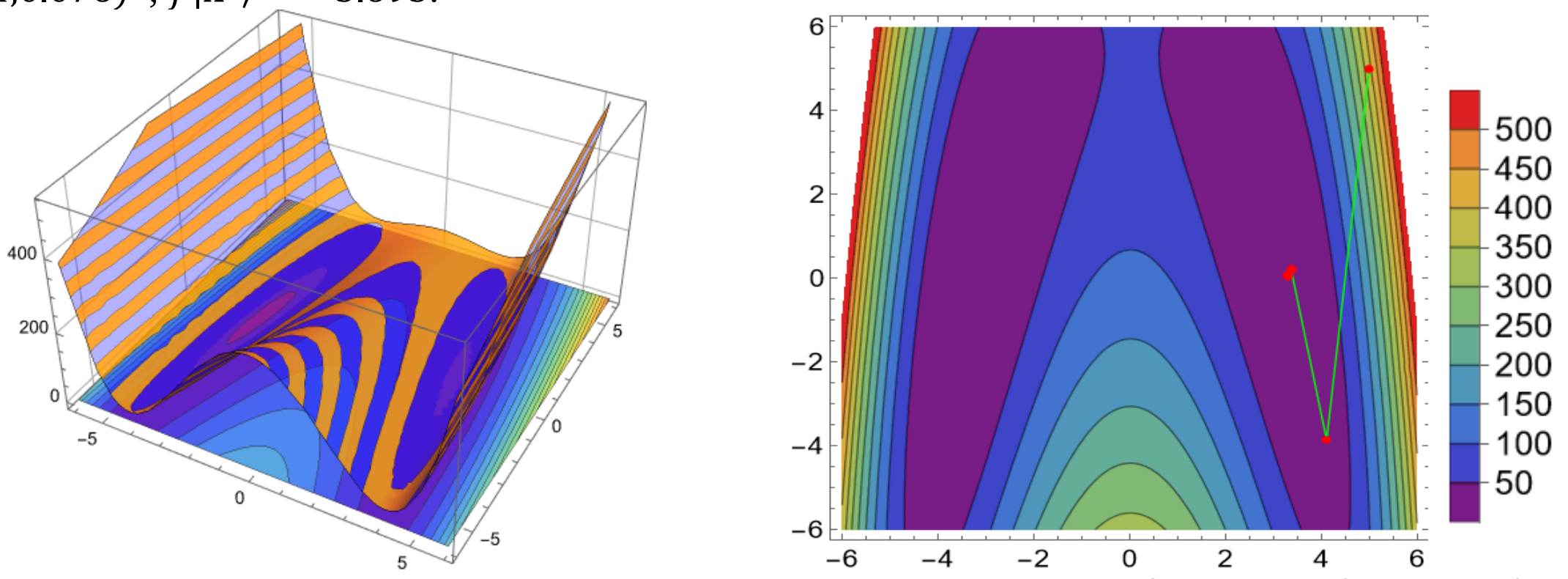

**Figure 8.7.** The results of 5 iterations of the Newton method, for $f|\mathbf{x}\rangle = (x^2 + y - 11)^2 + (x + y^2 - 7)$.

From Table 8.4, after 5 iterations, the optimal condition is attained. The results were produced by Mathematica code 8.4.

**Table 8.4.**

| No. of iters. | αstar | $x0$ | $x1$ | $f(x0)$ | $f(x1)$ | error |
|---|---|---|---|---|---|---|
| 1 | 1.107836 | {5., 5.} | {4.113, -3.856} | 384 | 16.21152 | 8.900749 |
| 2 | 1.132219 | {4.113, -3.856} | {3.372, 0.209} | 16.21152 | -3.25051 | 4.13176 |
| 3 | 1.02885 | {3.372, 0.209} | {3.297, 0.057} | -3.25051 | -3.69447 | 0.168754 |
| 4 | 0.998709 | {3.297, 0.057} | {3.294, 0.076} | -3.69447 | -3.69482 | 0.018747 |
| 5 | 1.001587 | {3.294, 0.076} | {3.294, 0.076} | -3.69482 | -3.69482 | 2.6E-05 |

***Mathematica Code 8.4*** Basic Newton Method

```
(* Basic Newton Method *)

(*
Notations:
x0        :Intial vector
epsilonNM :Small number to check the accuracy of the basic Newton method
f[x,y]    :Objective function
lii       :The last iteration index
result[k] :The results of iteration k
*)

(* Taking Initial Inputs from User *)
x0=Input["Enter the intial point in the format {x, y}; for example {5,5} "] ;
epsilonNM=Input["Please enter accuracy of the basic Newton method; for example 0.01 "];

If[
  epsilonNM<=0,
  Beep[];
  MessageDialog["The value of epsilonNM has to be postive number: "];
  Exit[];
  ];

domainx=Input["Please enter domain of x variable for 3D and contour plots; for example {-
6,6}"];
domainy=Input["Please enter domain of y variable for 3D and contour plots; for example {-
6,6}"];

(* Taking the Function from User *)
f[{x_,y_}] = Evaluate[Input["Please input a function of x and y to find the minimum "]];
(* For example: (x^2+y-11)^2+(x+y^2-7) *)

(*Defination of the Unidirectionalsearch Function*)
(*This function start by bracketing the minimum (using Bounding Phase Method) then isolating
the minimum (using Golden Section Search Method) *)

unidirectionalsearch[α_,delt_,eps_]:=Module[


{α0=α,delta=delt,epsilon=eps,y1,y2,y3,αα,a,b,increment,a0,b0,anew,bnew,anorm,bnorm,lnorm,α1n
orm,α2norm,α1,α2,φ1,φ2,αstar},

   (* Bounding Phase Method *)

   (* Initiating Required Variables *)
   y1 =φ[α0-Abs[delta]];
   y2 =φ[α0];
```

```
   y3 =φ[α0+Abs[delta]];

   (*Determining Whether the Inicrement Is Positive or Negative*)
   Which[
    y1==y2,
    a=α0-Abs[delta];
    b=α0;
    Goto[end];,
    y2==y3,
    a=α0;
    b=α0+Abs[delta];
    Goto[end];,
    y1==y3||(y1>y2&&y2<y3),
    a=α0-Abs[delta];
    b=α0+Abs[delta];
    Goto[end];
    ];

   Which[
     y1>y2&&y2>y3,
     increment=Abs[delta];,
     y1<y2&&y2<y3,
     increment=-Abs[delta];
     ]

    (* Starting the Algorithm *)
    Do[
     αα[0]=α0;
     αα[k+1]=αα[k]+2^k*increment;

     Which[
      φ[αα[k]]<φ[αα[k+1]],(* Evidently, it is impossible the condition to hold for k=0 *)
      a=αα[k-1];
      b= αα[k+1];
      Break [],

      k>50,
      Print["After 50 iterations the bounding phase method can not braketing the min of 
alpha"];
      Exit[]
      ];,
     {k,0,∞}
     ];

   Label[end];

   If[
    a>b,
    {a,b}={b,a}
    ];

   (* Golden Section Search Method *)

   (* Initiating Required Variables*)
   a0=a;
   b0=b;
   anew=a;
   bnew=b;

   If[
    a0==b0,
```

```
     αstar=a;
     Goto[final]
     ];

    (* Starting the Algorithm *)
    Do[
     (* Normalize the Variable α *)
     anorm=(anew-a)/(b-a);
     bnorm=(bnew-a)/(b-a);

     lnorm=bnorm-anorm;

     α1norm=anorm+0.382*lnorm;
     α2norm=bnorm-0.382*lnorm;

     α1=α1norm(b0-a0)+a0;
     α2=α2norm(b0-a0)+a0;

     φ1=φ[α1];
     φ2=φ[α2];

     Which [
      φ1>φ2,
      anew=α1(*move lower bound to α1*);,
      φ1<φ2,
      bnew=α2(*move upper bound to α2*);,
      φ1==φ2,
      anew=α1(*move lower bound to α1*);
      bnew=α2(*move upper bound to α2*);
      ];

     αstar=0.5*(anew+bnew);

     If[
      Abs[lnorm]<epsilon,
      Break[]
      ];,
     {k,1,∞}
     ];

    Label[final];

    (* Final Result *)
    N[αstar]
    ];

α00=2;(* The intial point of α; for example 2*)
delt0=1;(* The parameter delta of Bounding Phase Method; for example 1 *)
eps0=0.01;(* The accuracy of the Golden Section Search Method; for example 0.01 *)

gradfx[x_,y_]=Grad[f[{x,y}],{x,y}];
hessianfx[x_,y_]=Grad[gradfx[x,y],{x,y}];

(* Main Loop *)
Do[
  gradfx0=gradfx[x0[[1]],x0[[2]]];
  hessianfx0=hessianfx[x0[[1]],x0[[2]]];

  If[
   N[Norm[gradfx0]]==0,
```

```
   Print["The Newton-Raphson scheme requires that the gradiant of the function f
!=0"];(*Ending program*)
   Exit[];
   ];

  If[
   !PositiveDefiniteMatrixQ[hessianfx0],
   hessianfx0=IdentityMatrix[2]
   ];

  dx0=-(Inverse[hessianfx0].gradfx0);
  φ[α_]=f[x0+α*dx0];
  αstar=unidirectionalsearch[α00,delt0,eps0];
  x1=x0+αstar*dx0;

  error=Norm[αstar*dx0];

  lii=k;
  result[k]=N[{k,αstar,Row[x0,","],Row[x1,","],f[x0],f[x1],error}];
  plotresult[k]=N[{x0,x1}];

  If[
   error<epsilonNM||k>50,
   Break[],
   x0=x1;
   ];,
  {k,1,∞}
  ];

(* Final Result *)
If[
  lii==50,
  Print[" After 50 iterations the mimimum point around the intial point is " ,N[x1], "\nThe
solution is (approximately) ", N[f[x1]]];,
 Print["The solution is x= ",  N[x1],"\nThe solution is (approximately)= ", N[f[x1]]];
  ]

(* Results of Each Iteration *)
table=TableForm[
  Table[
   result[i],
   {i,1,lii}
   ],
  TableHeadings->{None,{"No. of iters.","αstar","x0","x1","f[x0]","f[x1]","error"}}
  ]

Export["example84.xls",table,"XLS"];

(* Data Visualization *)
(* Domain of Varibles*)
xleft=domainx[[1]];
xright=domainx[[2]];
ydown=domainy[[1]];
yup=domainy[[2]];

(* 3D+ Contour Plot *)
plot1=Plot3D[
   f[{x,y}],
   {x,xleft,xright},
   {y,ydown,yup},
   ClippingStyle->None,
```

```
   MeshFunctions->{#3&},
   Mesh->15,
   MeshStyle->Opacity[.5],
   MeshShading->{{Opacity[.3],Blue},{Opacity[.8],Orange}},
   Lighting->"Neutral"
   ];
slice=SliceContourPlot3D[
   f[{x,y}],
   z==0,
   {x,xleft,xright},
   {y,ydown,yup},
   {z,-1,1},
   Contours->15,
   Axes->False,
   PlotPoints->50,
   PlotRangePadding->0,
   ColorFunction->"Rainbow"
   ];
Show[
 plot1,
 slice,
 PlotRange->All,
 BoxRatios->{1,1,.6},
 FaceGrids->{Back,Left}
 ]
(* Contour Plot with Step Iterations *)
ContourPlot[
 f[{x,y}],
 {x,xleft,xright},
 {y,ydown,yup},
 LabelStyle->Directive[Black,16],
 ColorFunction->"Rainbow",
 PlotLegends->Automatic,
 Contours->10,
 Epilog-
>{PointSize[0.015],Green,Line[Flatten[Table[plotresult[i],{i,1,lii}],1]],Red,Point[Flatten[T
able[plotresult[i],{i,1,lii}],1]]}
 ]
(* Data Manipulation *)
Manipulate[
 ContourPlot[
  f[{x,y}],
  {x,xleft,xright},
  {y,ydown,yup},
  LabelStyle->Directive[Black,14],
  ColorFunction->"Rainbow",
  PlotLegends->Automatic,
  Contours->10,
  Epilog->{
    PointSize[0.015],
    Yellow,
    Arrow[{plotresult[i][[1]],plotresult[i][[2]]}],
    Red,
    Point[Flatten[Table[plotresult[j],{j,1,i}],1]],
    Green,
    Line[Flatten[Table[plotresult[j],{j,1,i}],1]]
    }
  ],
 {i,1,lii,1}
 ]
```

### 8.2.2. Marquardt Method and Modification of the Hessian

If the Hessian is not positive definite in any iteration of the Newton algorithm, it is forced to become positive definite in Step 2 of the algorithm. This modification of

$$\mathbf{H}_k = \nabla^2 f|\mathbf{x}_k\rangle, \tag{8.50}$$

can be accomplished in one of several ways. One approach is to replace the matrix $\mathbf{H}_k$ by the $n \times n$ identity matrix $\mathbf{I}_n$ whenever it becomes nonpositive definite [8]. Since $\mathbf{I}_n$ is positive definite, the problem of a nonsingular $\mathbf{H}_k$ is eliminated.

Another approach would be to let

$$\widehat{\mathbf{H}}_k = \frac{1}{1+\beta}(\mathbf{H}_k + \beta \mathbf{I}_n), \tag{8.51}$$

where $\beta$ is a positive scalar that is slightly larger than the absolute value of the most negative eigenvalue of $\mathbf{H}_k$ so as to assure the positive definiteness of $\mathbf{H}_k$ in (8.51). If $\beta$ is large, then

$$\widehat{\mathbf{H}}_k \approx \mathbf{I}_n, \tag{8.52}$$

and from (8.49)

$$|\mathbf{d}_k\rangle \approx -\nabla f|\mathbf{x}_k\rangle, \tag{8.53}$$

In effect, the modification in (8.51) converts the Newton method into the steepest-descent method. A nonpositive definite $\mathbf{H}_k$ is likely to arise at points far from the solution where the steepest-descent method is most effective in reducing the value of $f|\mathbf{x}\rangle$. Therefore, the modification in (8.51) leads to an algorithm that combines the complementary convergence properties of the Newton and steepest-descent methods.

#### Marquardt Method

Cauchy method works well when the initial point is far away from the minimum point, and the Newton method works well when the initial point is near the minimum point. In any given problem, it is usually not known whether the chosen initial point is away from the minimum or close to the minimum, but wherever be the minimum point, a method can be devised to take advantage of both these methods.

The Marquardt method [3,9-11] combines Cauchy and Newton methods in a convenient manner that exploits the strengths of both but does require second-order information. Marquardt specified the search direction to be

$$|\mathbf{s}_k\rangle = -\left[\mathbf{H}_f(\mathbf{x}_k) + \lambda_k \mathbf{I}\right]^{-1} \nabla f|\mathbf{x}_k\rangle, \tag{8.54}$$

and set $\alpha_k = +1$ in $|\mathbf{x}_{k+1}\rangle = |\mathbf{x}_k\rangle + \alpha_k|\mathbf{s}_k\rangle$, since $\lambda$ is used to control both the direction of the search and the length of the step. To begin the search, let $\lambda_0$ be a large constant, say $10^4$, such that

$$\left[\mathbf{H}_f(\mathbf{x}_0) + \lambda_0 \mathbf{I}\right]^{-1} \cong [\lambda_0 \mathbf{I}]^{-1} = \frac{\mathbf{I}}{\lambda_0}. \tag{8.55}$$

We notice from (8.54) that as $\lambda$ decreases from a large value to zero, $|\mathbf{s}\rangle$ goes from the gradient to the Newton direction. Hence, in the Marquardt method, the Cauchy method is initially followed. Thereafter, the Newton method is adopted.

**Algorithm**

**Step 1:** Choose a starting point, $|\mathbf{x}_0\rangle$, the maximum number of iterations, $M$, and a termination parameter, $\epsilon$. Set $k = 0$ and $\lambda^{(0)} = 10^4$ (a large number).

**Step 2:** Calculate $\nabla f|\mathbf{x}_k\rangle$ and $\mathbf{H}^{(k)}$.

**Step 3:** If $\|\nabla f(\mathbf{x}_k)\| \le \epsilon$ or $k \ge M$? Terminate;
Else go to Step 4.

| | |
|---|---|
| **Step 4:** | Calculate $\mathbf{s}^{(k)} = -\left[\mathbf{H}^{(k)} + \lambda\mathbf{I}\right]^{-1}\nabla f|\mathbf{x}_k\rangle$. Set $|\mathbf{x}_{k+1}\rangle = |\mathbf{x}_k\rangle + |\mathbf{s}_k\rangle$. |
| **Step 5:** | Is $f|\mathbf{x}_{k+1}\rangle < f|\mathbf{x}_k\rangle$? If yes, go to Step 6; Else go to Step 7. |
| **Step 6:** | Set $\lambda^{(k+1)} = \frac{1}{2}\lambda^{(k)}$, $k = k + 1$, and go to Step 2. |
| **Step 7:** | Set $\lambda^{(k)} = 2\lambda^{(k)}$ and go to Step 4. |

The algorithm can be made faster by performing a unidirectional search while finding the new point in Step 4: $|\mathbf{x}_{k+1}\rangle = |\mathbf{x}_k\rangle + \alpha_k|\mathbf{s}_k\rangle$. Since the computations of the Hessian matrix and its inverse are computationally expensive, the unidirectional search along $|\mathbf{s}_k\rangle$ is not usually performed. For simpler objective functions, however, a unidirectional search in Step 4 can be achieved to find the new point $|\mathbf{x}_{k+1}\rangle$.

***Example 8.7***

Consider the problem:

$$\text{Minimize } f|\mathbf{x}\rangle = (x^2 + y - 11)^2 + (x + y^2 - 7),$$

using, $|\mathbf{x}_0\rangle = (5,5)^T$, and $\epsilon = 0.01$.

***Solution***

The 3D and contour plots of the function are shown in Figure 8.8. The plots show that the minimum lies at $|\mathbf{x}^*\rangle = (\{3.294, 0.076)^T$, $f|\mathbf{x}^*\rangle = -3.695$.

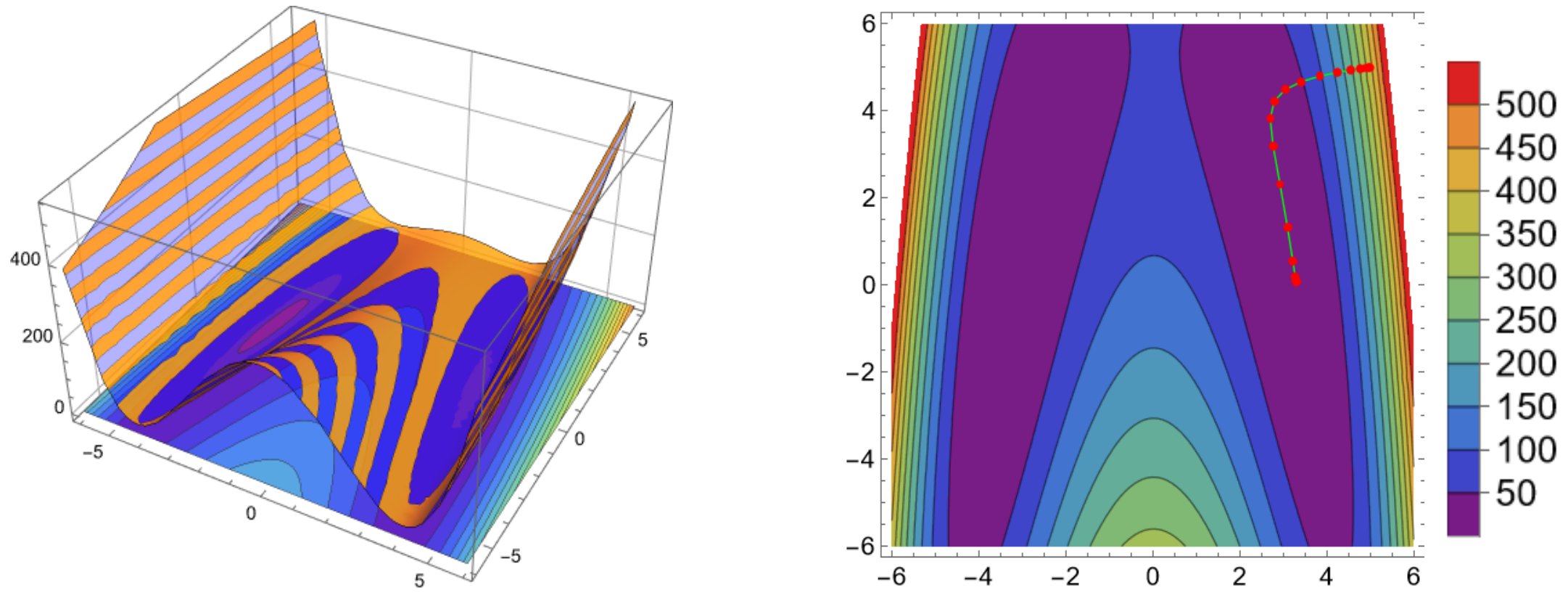

**Figure 8.8.** The results of 18 iterations of the Marquardt method, for $f|\mathbf{x}\rangle = (x^2 + y - 11)^2 + (x + y^2 - 7)$.

From Table 8.5, after 18 iterations, the optimal condition is attained. The results were produced by Mathematica code 8.5.

**Table 8.5.**

| No. of iters. | x0 | x1 | $f(x0)$ | $f(x1)$ | λ | grad $f(x0)$ |
|---|---|---|---|---|---|---|
| 1 | {5., 5.} | {4.962, 4.995} | 384 | 369.8427 | 10000 | 384.0117 |
| 2 | {4.962, 4.995} | {4.892, 4.986} | 369.8427 | 344.0235 | 5000 | 373.7556 |
| 3 | {4.892, 4.986} | {4.765, 4.968} | 344.0235 | 300.6203 | 2500 | 354.7538 |
| 4 | {4.765, 4.968} | {4.553, 4.937} | 300.6203 | 237.0638 | 1250 | 321.8436 |
| 5 | {4.553, 4.937} | {4.238, 4.884} | 237.0638 | 161.5026 | 625 | 270.9743 |
| 6 | {4.238, 4.884} | {3.840, 4.799} | 161.5026 | 92.96678 | 312.5 | 204.6507 |
| 7 | {3.840, 4.799} | {3.417, 4.673} | 92.96678 | 46.93404 | 156.25 | 134.9973 |
| 8 | {3.417, 4.673} | {3.046, 4.491} | 46.93404 | 23.8945 | 78.125 | 76.86545 |
| 9 | {3.046, 4.491} | {2.792, 4.225} | 23.8945 | 14.70148 | 39.0625 | 37.66962 |
| 10 | {2.792, 4.225} | {2.700, 3.823} | 14.70148 | 10.33225 | 19.53125 | 16.28916 |
| 11 | {2.700, 3.823} | {2.765, 3.199} | 10.33225 | 6.027281 | 9.765625 | 8.188342 |
| 12 | {2.765, 3.199} | {2.925, 2.314} | 6.027281 | 1.297385 | 4.882813 | 6.134885 |
| 13 | {2.925, 2.314} | {3.097, 1.323} | 1.297385 | -2.14442 | 2.441406 | 4.401651 |
| 14 | {3.097, 1.323} | {3.220, 0.557} | -2.14442 | -3.46404 | 1.220703 | 2.478397 |
| 15 | {3.220, 0.557} | {3.276, 0.190} | -3.46404 | -3.68166 | 0.610352 | 0.966446 |
| 16 | {3.276, 0.190} | {3.291, 0.091} | -3.68166 | -3.69458 | 0.305176 | 0.229988 |
| 17 | {3.291, 0.091} | {3.293, 0.077} | -3.69458 | -3.69482 | 0.152588 | 0.030814 |
| 18 | {3.293, 0.077} | {3.293, 0.075} | -3.69482 | -3.69482 | 0.076294 | 0.002229 |

***Mathematica Code 8.5*** Marquardt Method

```
(* Marquardt Method and Modification of the Hessian *)

(*
Notations:
x0          :Intial vector
epsilonMM    :Small number to check the accuracy of the Marquardt Method and Modification of
the Hessian
f[x,y]       :Objective function
lii          :The last iteration index
result[k]    :The results of iteration k
*)

(* Taking Initial Inputs from User *)
x0=Input["Enter the intial point in the format {x, y}; for example {5,5} "] ;
epsilonMM=Input["Please enter accuracy of the Marquardt Method and Modification of the
Hessian; for example 0.01 "];

domainx=Input["Please enter domain of x variable for 3D and contour plots; for example {-
6,6}"];
domainy=Input["Please enter domain of y variable for 3D and contour plots; for example {-
6,6}"];

If[
  epsilonMM<=0,
  Beep[];
  MessageDialog["The value of epsilonMM has to be postive number: "];
  Exit[];
  ];
(* Taking the Function from User *)
f[{x_,y_}] = Evaluate[Input["Please input a function of x and y to find the minimum "]];
(* For example: (x^2+y-11)^2+(x+y^2-7) *)

gradfx[x_,y_]=Grad[f[{x,y}],{x,y}];
hessianfx[x_,y_]=Grad[gradfx[x,y],{x,y}];
λ=10000;

(* Main Loop *)
Do[
  gradfx0=gradfx[x0[[1]],x0[[2]]];
  hessianfx0=hessianfx[x0[[1]],x0[[2]]];

  If[
   N[Norm[gradfx0]]==0,
   Print["The Marquardt's Method requires that the gradiant of the function f !=0"];
(* Ending program *)
   Exit[];
   ];
```

```
  dx0=-N[Inverse[hessianfx0+λ*IdentityMatrix[2]].gradfx0];
  x1=x0+dx0;

  lii=k;
  result[k]=N[{k,Row[x0,","],Row[x1,","],f[x0],f[x1],λ,N[Norm[gradfx0]]}];
  plotresult[k]=N[{x0,x1}];

  If[
   f[x1]<f[x0],
   λ=0.5*λ;
   x0=x1;,
   λ=2*λ;
   x0=x1;
   ];

  If[
   N[Norm[gradfx0]]<epsilonMM||k>50,
   Break[];
   ];,
  {k,1,∞}
  ];

(* Final Result *)
If[
  lii==50,
  Print[" After 50 iterations the mimimum point around the intial point is " ,N[x1], "\nThe solution is (approximately) ", N[f[x1]]];,
 Print["The solution is x= ",  N[x1],"\nThe solution is (approximately)= ", N[f[x1]]];
  ]

(* Results of Each Iteration *)
table=TableForm[
  Table[
   result[i],
   {i,1,lii}
   ],
  TableHeadings->{None,{"No. of iters.","x0","x1","f[x0]","f[x1]","λ","N[gradfx0]"}}
  ]

Export["example85.xls",table,"XLS"];

(* Data Visualization *)
(* Domain of Variables *)
xleft=domainx[[1]];
xright=domainx[[2]];
ydown=domainy[[1]];
yup=domainy[[2]];

(* 3D+ contour plot *)
plot1=Plot3D[
   f[{x,y}],
   {x,xleft,xright},
   {y,ydown,yup},
   ClippingStyle->None,
   MeshFunctions->{#3&},
   Mesh->15,
   MeshStyle->Opacity[.5],
   MeshShading->{{Opacity[.3],Blue},{Opacity[.8],Orange}},
   Lighting->"Neutral"
   ];
```

```
slice=SliceContourPlot3D[
   f[{x,y}],
   z==0,
   {x,xleft,xright},
   {y,ydown,yup},
   {z,-1,1},
   Contours->15,
   Axes->False,
   PlotPoints->50,
   PlotRangePadding->0,
   ColorFunction->"Rainbow"
   ];
Show[
 plot1,
 slice,
 PlotRange->All,
 BoxRatios->{1,1,.6},
 FaceGrids->{Back,Left}
 ]

(* Contour Plot with Step Iterations *)
ContourPlot[
 f[{x,y}],
 {x,xleft,xright},
 {y,ydown,yup},
 LabelStyle->Directive[Black,16],
 ColorFunction->"Rainbow",
 PlotLegends->Automatic,
 Contours->10,
 Epilog-
>{PointSize[0.015],Green,Line[Flatten[Table[plotresult[i],{i,1,lii}],1]],Red,Point[Flatten[T
able[plotresult[i],{i,1,lii}],1]]}
 ]

(* Data Manipulation *)
Manipulate[
 ContourPlot[
  f[{x,y}],
  {x,xleft,xright},
  {y,ydown,yup},
  LabelStyle->Directive[Black,14],
  ColorFunction->"Rainbow",
  PlotLegends->Automatic,
  Contours->10,
  Epilog->{
    PointSize[0.015],
    Yellow,
    Arrow[{plotresult[i][[1]],plotresult[i][[2]]}],
    Red,
    Point[Flatten[Table[plotresult[j],{j,1,i}],1]],
    Green,
    Line[Flatten[Table[plotresult[j],{j,1,i}],1]]
    }
  ],
 {i,1,lii,1}
 ]
```

### 8.2.3. Gauss-Newton Method

In many optimization problems, the objective function is in the form of a vector of functions given by

$$|\mathbf{f}\rangle = [\, f_1|\mathbf{x}\rangle \;\; f_2|\mathbf{x}\rangle \;\; \cdots \;\; f_m|\mathbf{x}\rangle]^T, \tag{8.56}$$

where $f_p|\mathbf{x}\rangle$ for $p = 1,2,\ldots,m$ are independent functions of $|\mathbf{x}\rangle$. The solution sought is a point $|\mathbf{x}\rangle$ such that all $f_p|\mathbf{x}\rangle$ are reduced to zero simultaneously. In problems of this type, a real-valued function can be formed as

$$F = \sum_{p=1}^{m} (f_p|\mathbf{x}\rangle)^2 = \langle \mathbf{f}|\mathbf{f}\rangle. \tag{8.57}$$

If $F$ is minimized by using a multidimensional unconstrained algorithm, then the individual functions $f_p|\mathbf{x}\rangle$ are minimized in the least-squares sense. A method for the solution of this class of problems, known as theGauss–Newton method [3], can be readily developed by applying the Newton method. Since there are a number of functions $f_p|\mathbf{x}\rangle$ and each one depends on $x_i$ for $i = 1,2,\ldots,n$ a gradient matrix (the Jacobian) can be formed as

$$\mathbf{J} = \begin{pmatrix} \frac{\partial f_1}{\partial x_1} & \frac{\partial f_1}{\partial x_2} & \cdots & \frac{\partial f_1}{\partial x_n} \\ \frac{\partial f_2}{\partial x_1} & \frac{\partial f_2}{\partial x_2} & \cdots & \frac{\partial f_2}{\partial x_n} \\ \vdots & \vdots & \ddots & \vdots \\ \frac{\partial f_m}{\partial x_1} & \frac{\partial f_m}{\partial x_2} & \cdots & \frac{\partial f_m}{\partial x_n} \end{pmatrix}. \tag{8.58}$$

The number of functions $m$ may exceed the number of variables $n$; that is, the Jacobian need not be a square matrix. By differentiating $F$ in (8.57) with respect to $x_i$, we obtain

$$\frac{\partial F}{\partial x_i} = \sum_{p=1}^{m} 2 f_p|\mathbf{x}\rangle \frac{\partial f_p}{\partial x_i}, \tag{8.59}$$

for $i = 1,2,\ldots,n$. Alternatively, in matrix form

$$\begin{pmatrix} \frac{\partial F}{\partial x_1} \\ \frac{\partial F}{\partial x_2} \\ \vdots \\ \frac{\partial F}{\partial x_n} \end{pmatrix} = 2 \begin{pmatrix} \frac{\partial f_1}{\partial x_1} & \frac{\partial f_1}{\partial x_2} & \cdots & \frac{\partial f_1}{\partial x_n} \\ \frac{\partial f_2}{\partial x_1} & \frac{\partial f_2}{\partial x_2} & \cdots & \frac{\partial f_2}{\partial x_n} \\ \vdots & \vdots & \ddots & \vdots \\ \frac{\partial f_m}{\partial x_1} & \frac{\partial f_m}{\partial x_2} & \cdots & \frac{\partial f_m}{\partial x_n} \end{pmatrix} \begin{pmatrix} f_1|\mathbf{x}\rangle \\ f_2|\mathbf{x}\rangle \\ \vdots \\ f_m|\mathbf{x}\rangle \end{pmatrix}. \tag{8.60}$$

Hence the gradient of $F$, designated by $|\mathbf{g}_F\rangle$, can be expressed as

$$|\mathbf{g}_F\rangle = \; 2\mathbf{J}^T|\mathbf{f}\rangle. \tag{8.61}$$

Assuming that $f_p|\mathbf{x}\rangle \in \mathrm{C}^2$, (8.59) yields

$$\frac{\partial^2 F}{\partial x_i \partial x_j} = 2\sum_{p=1}^{m} \frac{\partial f_p}{\partial x_i}\frac{\partial f_p}{\partial x_j} + 2\sum_{p=1}^{m} f_p|\mathbf{x}\rangle \frac{\partial^2 f_p}{\partial x_i \partial x_j}, \tag{8.62}$$

for $i, j = 1,2,\ldots,n$. If the second derivative of $f_p|\mathbf{x}\rangle$ are neglected, we have

$$\frac{\partial^2 F}{\partial x_i \partial x_j} \approx 2\sum_{p=1}^{m} \frac{\partial f_p}{\partial x_i}\frac{\partial f_p}{\partial x_j}. \tag{8.63}$$

Thus, the Hessian of $F$, designated by $\mathbf{H}_F$, can be deduced as

$$\mathbf{H}_F \approx 2\mathbf{J}^T\mathbf{J}. \tag{8.64}$$

Since the gradient and Hessian of $F$ are now known, the Newton method can be applied to the solution of the problem. The necessary recursive relation is given by (8.46), (8.49), (8.61), and (8.64) as

$$\begin{aligned} |\mathbf{x}_{k+1}\rangle &= |\mathbf{x}_k\rangle - \alpha_k (2\mathbf{J}^T\mathbf{J})^{-1}(2\mathbf{J}^T|\mathbf{f}\rangle) \\ &= |\mathbf{x}_k\rangle - \alpha_k (\mathbf{J}^T\mathbf{J})^{-1}(\mathbf{J}^T|\mathbf{f}\rangle), \end{aligned} \tag{8.65}$$

where $\alpha_k$ is the value of $\alpha$ that minimizes $F(|\mathbf{x}_k\rangle + \alpha|\mathbf{d}_k\rangle)$. When $|\mathbf{x}_k\rangle$ is in the neighborhood of $|\mathbf{x}^*\rangle$, the matrix in (8.64) becomes an accurate representation of the Hessian of $F_k$, and the method converges rapidly. If functions $f_p|\mathbf{x}\rangle$ are linear, $F$ is quadratic, the matrix in (8.64) is the Hessian, and the problem is solved in one iteration. The method breaks down if $\mathbf{H}_F$ becomes singular, as in the case of the Newton method.

**Algorithm** **Gauss-Newton algorithm**

**Step 1:** Input $\mathbf{x}^{(0)}$ and initialize the tolerance $\varepsilon$. Set $k = 0$.

**Step 2:** Compute $f_{pk} = f_p|\mathbf{x}_k\rangle$ for $p = 1, 2, \ldots, m$ and $F_k$.

**Step 3:** Compute $\mathbf{J}_k$, $|\mathbf{g}_k\rangle = 2\mathbf{J}_k^T|\mathbf{f}\rangle$, and $\mathbf{H}_k = 2\mathbf{J}_k^T\mathbf{J}_k$.

**Step 4:** Compute $|\mathbf{d}_k\rangle = -(\mathbf{J}^T\mathbf{J})^{-1}(\mathbf{J}^T|\mathbf{f}\rangle) = -\mathbf{H}_F|\mathbf{g}_k\rangle$.

**Step 5:** Find $\alpha_k$, the value of $\alpha$ that minimizes $F(|\mathbf{x}_k\rangle + \alpha|\mathbf{d}_k\rangle)$.

**Step 6:** Set $|\mathbf{x}_{k+1}\rangle = |\mathbf{x}_k\rangle + \alpha_k|\mathbf{d}_k\rangle$. Compute $f_{p\,(k+1)}$ for $p = 1, 2, \ldots, m$ and $F_{k+1}$.

**Step 7:** If $|F_{k+1} - F_k| < \varepsilon$, then output $|\mathbf{x}^*\rangle = |\mathbf{x}_{k+1}\rangle$, $f_{p\,(k+1)}|\mathbf{x}^*\rangle$ for $p = 1, 2, \ldots, m$, and $F_{k+1}$ and stop; otherwise, set $k = k + 1$ and repeat from Step 3.

***Example 8.8***

Consider the problem:

$$\text{Minimize } \mathbf{f}|\mathbf{x}\rangle = \{x^2 + y - 11, x + y^2 - 7\},$$

using, $|\mathbf{x}_0\rangle = (5,5)^T$, and $\epsilon = 0.01$, $\mathbf{f}|\mathbf{x}\rangle$ is a vector-valued function.

***Solution***

The 3D and contour plots of the function are shown in Figure 8.9. The plots show that the minimum lies at $|\mathbf{x}^*\rangle = (\{3,2)^T, f|\mathbf{x}^*\rangle = 5.475256215385999 \times 10^{-14}$.

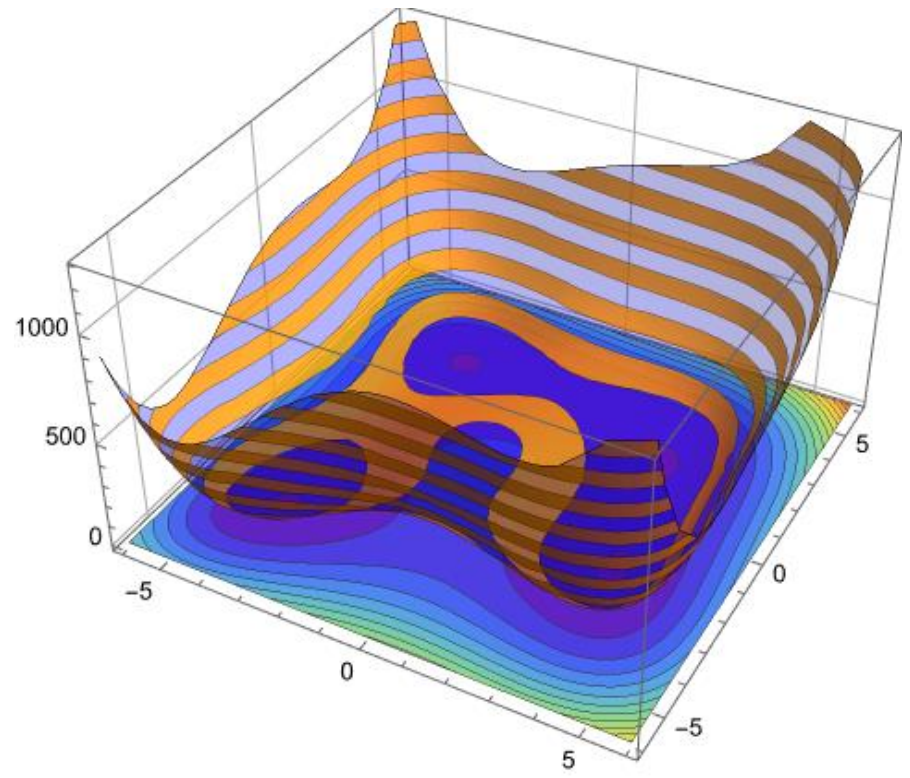

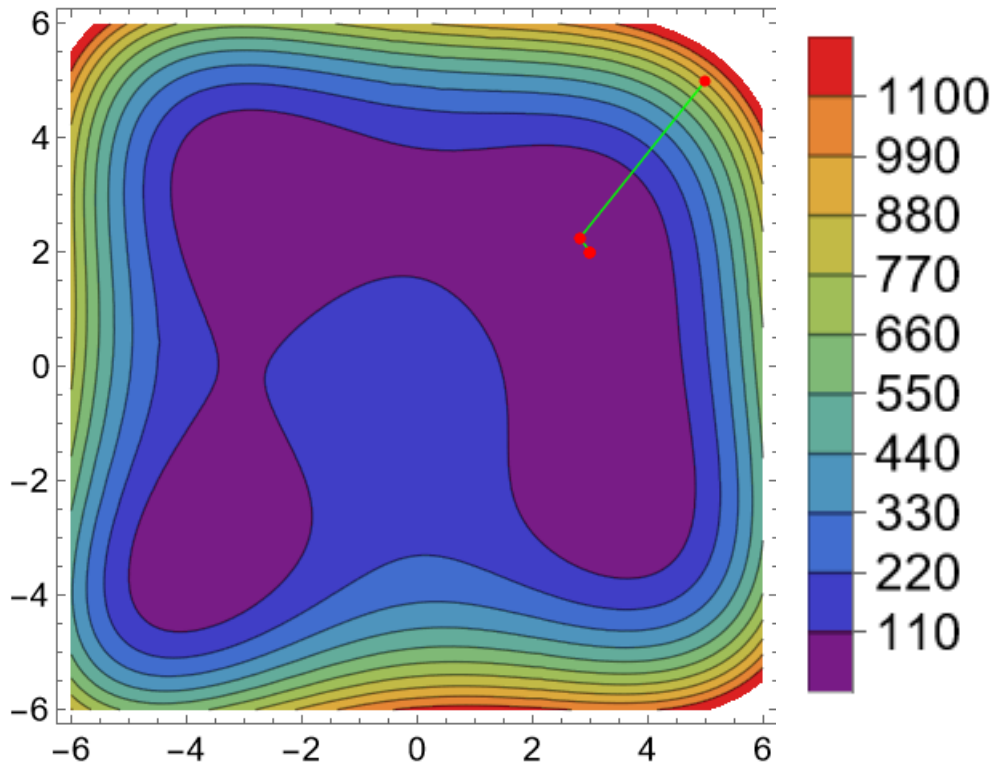

**Figure 8.9.** The results of 4 iterations of the Gauss–Newton method, for $\mathbf{f}|\mathbf{x}\rangle = \{x^2 + y - 11, x + y^2 - 7\}$.

From Table 8.6, after 4 iterations, the optimal condition is attained. The results were produced by Mathematica code 8.6.

**Table 8.6.**

| No. of iters. | αstar | $x0$ | $x1$ | $F(x0)$ | $F(x1)$ | error |
|---|---|---|---|---|---|---|
| 1 | 1.290081 | {5., 5.} | {2.824, 2.250} | 890 | 1.390712 | 3.506557 |
| 2 | 1.013781 | {2.824, 2.250} | {3.006, 2.009} | 1.390712 | 0.003992 | 0.301471 |
| 3 | 0.998709 | {3.006, 2.009} | {3.000, 2.000} | 0.003992 | 3.25E-08 | 0.011407 |
| 4 | 0.998709 | {3.000, 2.000} | {3, 2} | 3.25E-08 | 5.48E-14 | 3.81E-05 |

***Mathematica Code 8.6*** Gauss–Newton Algorithm

```
(* Gauss–Newton Method with unidirectional search (Approximate Hessian) (Vector Valued
Function) (Sum of Squares) *)

(*
Notations:
x0          :Intial vector
epsilonNM  :Small number to check the accuracy of the Gauss-Newton Method
f[x,y]      :Objective function
lii         :The last iteration index
result[k]  :The results of iteration k
*)

(* Taking Initial Inputs from User *)
x0=Input["Enter the intial point in the format {x, y}; for example {5,5} "] ;
epsilonNM=Input["Please enter accuracy of the Gauss-Newton Method; for example 0.01 "];

domainx=Input["Please enter domain of x variable for 3D and contour plots; for example {-
6,6}"];
domainy=Input["Please enter domain of y variable for 3D and contour plots; for example {-
6,6}"];

If[
  epsilonNM<=0,
  Beep[];
  MessageDialog["The value of epsilonNM has to be postive number: "];
  Exit[];
  ];

(* Taking the Function from User *)
f[{x_,y_}] = Evaluate[Input["Please input a function of x and y  as a vector valued
function; for example {x^2+y-11,x+y^2-7} "]];
(* VIP VIP Input the f as a vector valued function; for example: {x^2+y-11,x+y^2-7} *)

(*Defination of the Unidirectionalsearch Function*)
(*This function start by bracketing the minimum (using Bounding Phase Method) then isolating
the minimum (using Golden Section Search Method) *)

unidirectionalsearch[α_,delt_,eps_]:=Module[

{α0=α,delta=delt,epsilon=eps,y1,y2,y3,αα,a,b,increment,a0,b0,anew,bnew,anorm,bnorm,lnorm,α1n
orm,α2norm,α1,α2,φ1,φ2,αstar},

   (* Bounding Phase Method *)
   (* Initiating Required Variables *)
   y1 =φ[α0-Abs[delta]];
   y2 =φ[α0];
   y3 =φ[α0+Abs[delta]];

   (*Determining Whether the Inicrement Is Positive or Negative*)
   Which[
    y1==y2,
    a=α0-Abs[delta];
    b=α0;
    Goto[end];,
    y2==y3,
    a=α0;
    b=α0+Abs[delta];
    Goto[end];,
```

```
     y1==y3||(y1>y2&&y2<y3),
     a=α0-Abs[delta];
     b=α0+Abs[delta];
     Goto[end];
     ];

    Which[
      y1>y2&&y2>y3,
      increment=Abs[delta];,
      y1<y2&&y2<y3,
      increment=-Abs[delta];
      ]

     (* Starting the Algorithm *)
     Do[
      αα[0]=α0;
      αα[k+1]=αα[k]+2^k*increment;

      Which[
       φ[αα[k]]<φ[αα[k+1]],(* Evidently, it is impossible the condition to hold for k=0 *)
       a=αα[k-1];
       b= αα[k+1];
       Break [],

       k>50,
       Print["After 50 iterations the bounding phase method can not braketing the min of
alpha"];
       Exit[]
       ];,
      {k,0,∞}
      ];

    Label[end];

    If[
     a>b,
     {a,b}={b,a}
     ];

    (* Golden Section Search Method *)
    (* Initiating Required Variables*)
    a0=a;
    b0=b;
    anew=a;
    bnew=b;

    If[
     a0==b0,
     αstar=a;
     Goto[final]
     ];

    (* Starting the Algorithm *)
    Do[
     (* Normalize the Variable α *)
     anorm=(anew-a)/(b-a);
     bnorm=(bnew-a)/(b-a);

     lnorm=bnorm-anorm;

     α1norm=anorm+0.382*lnorm;
```

```
     α2norm=bnorm-0.382*lnorm;

     α1=α1norm(b0-a0)+a0;
     α2=α2norm(b0-a0)+a0;

     φ1=φ[α1];
     φ2=φ[α2];

     Which [
      φ1>φ2,
      anew=α1(*move lower bound to α1*);,
      φ1<φ2,
      bnew=α2(*move upper bound to α2*);,
      φ1==φ2,
      anew=α1(*move lower bound to α1*);
      bnew=α2(*move upper bound to α2*);
      ];

     αstar=0.5*(anew+bnew);

     If[
      Abs[lnorm]<epsilon,
      Break[]
      ];,
     {k,1,∞}
     ];

    Label[final];

    (* Final Result *)
    N[αstar]
    ];

α00=2;(* The intial point of α; for example 2*)
delt0=1;(* The parameter delta of Bounding Phase Method; for example 1 *)
eps0=0.01;(* The accuracy of the Golden Section Search Method; for example 0.01 *)

v={x,y};
jacobianf=Grad[f[{x,y}],v];
functionF[{x_,y_}]:=f[{x,y}].f[{x,y}];

gradF[x_,y_]=2*Transpose[jacobianf].f[{x,y}];
HessianF[x_,y_]=2*Transpose[jacobianf].jacobianf;

(* Main Loop *)
Do[
  gradFx0=gradF[x0[[1]],x0[[2]]];
  hessianFx0=HessianF[x0[[1]],x0[[2]]];

  If[
   N[Norm[gradFx0]]==0,
   Print["The Newton-Raphson scheme requires that the gradiant of the function f
!=0"];(*Ending program*)
   Exit[];
   ];

  If[
   !PositiveDefiniteMatrixQ[hessianFx0],
   hessianFx0=IdentityMatrix[2]
   ];
```

```
  inverseHessianF=Inverse[hessianFx0];
  dx0=-(inverseHessianF.gradFx0);
  φ[α_]=functionF[x0+α*dx0];
  αstar=unidirectionalsearch[α00,delt0,eps0];
  x1=x0+αstar*dx0;

  error=Norm[αstar*dx0];

  lii=k;
  result[k]=N[{k,αstar,Row[x0,","],Row[x1,","],functionF[x0],functionF[x1],error}];
  plotresult[k]=N[{x0,x1}];

  If[
   error<epsilonNM||k>50,
   Break[],
   x0=x1;
   ];,
  {k,1,∞}
  ];

(* Final Result *)
If[
  lii==50,
  Print[" After 50 iterations the mimimum point around the intial point is " ,N[x1], "\nThe
solution is (approximately) ", N[functionF[x1]]];,
 Print["The solution is x= ",  N[x1],"\nThe solution is (approximately)= ",
N[functionF[x1]]];
  ]

(* Results of Each Iteration *)
table=TableForm[
  Table[
   result[i],
   {i,1,lii}
   ],
  TableHeadings->{None,{"No. of
iters.","αstar","x0","x1","functionF[x0]","functionF[x1]","error"}}
  ]

Export["example86.xls",table,"XLS"];

(* Data Visualization *)
(* Domain of Varibles*)
xleft=domainx[[1]];
xright=domainx[[2]];
ydown=domainy[[1]];
yup=domainy[[2]];

(* 3D+ Contour Plot *)
plot1=Plot3D[
   functionF[{x,y}],
   {x,xleft,xright},
   {y,ydown,yup},
   ClippingStyle->None,
   MeshFunctions->{#3&},
   Mesh->15,
   MeshStyle->Opacity[.5],
   MeshShading->{{Opacity[.3],Blue},{Opacity[.8],Orange}},
   Lighting->"Neutral"
   ];
slice=SliceContourPlot3D[
```

```
    functionF[{x,y}],
    z==0,
    {x,xleft,xright},
    {y,ydown,yup},
    {z,-1,1},
    Contours->15,
    Axes->False,
    PlotPoints->50,
    PlotRangePadding->0,
    ColorFunction->"Rainbow"
    ];
Show[
 plot1,
 slice,
 PlotRange->All,
 BoxRatios->{1,1,.6},
 FaceGrids->{Back,Left}
 ]

(* Contour Plot with Step Iterations *)
ContourPlot[
 functionF[{x,y}],
 {x,xleft,xright},
 {y,ydown,yup},
 LabelStyle->Directive[Black,16],
 ColorFunction->"Rainbow",
 PlotLegends->Automatic,
 Contours->10,
 Epilog-
>{PointSize[0.015],Green,Line[Flatten[Table[plotresult[i],{i,1,lii}],1]],Red,Point[Flatten[T
able[plotresult[i],{i,1,lii}],1]]}
 ]

(* Data Manipulation *)
Manipulate[
 ContourPlot[
  functionF[{x,y}],
  {x,xleft,xright},
  {y,ydown,yup},
  LabelStyle->Directive[Black,14],
  ColorFunction->"Rainbow",
  PlotLegends->Automatic,
  Contours->10,
  Epilog->{
    PointSize[0.015],
    Yellow,
    Arrow[{plotresult[i][[1]],plotresult[i][[2]]}],
    Red,
    Point[Flatten[Table[plotresult[j],{j,1,i}],1]],
    Green,
    Line[Flatten[Table[plotresult[j],{j,1,i}],1]]
    }
  ],
 {i,1,lii,1}
 ]
```

***Example 8.9***

Consider the problem:

$$\text{Minimize } \mathbf{f}|\mathbf{x}\rangle = \{x^2 + y - 11, x + y^2 - 7\},$$

using, $|\mathbf{x}_0\rangle = (5,5)^T$, and $\epsilon = 0.01$, $\mathbf{f}|\mathbf{x}\rangle$ is a vector-valued function.

***Solution***

The 3D and contour plots of the function are shown in Figure 8.10. The plots show that the minimum lies at $|\mathbf{x}^*\rangle = (3,2)^T$, $f|\mathbf{x}^*\rangle = 2.471036864817747 \times 10^{-12}$.

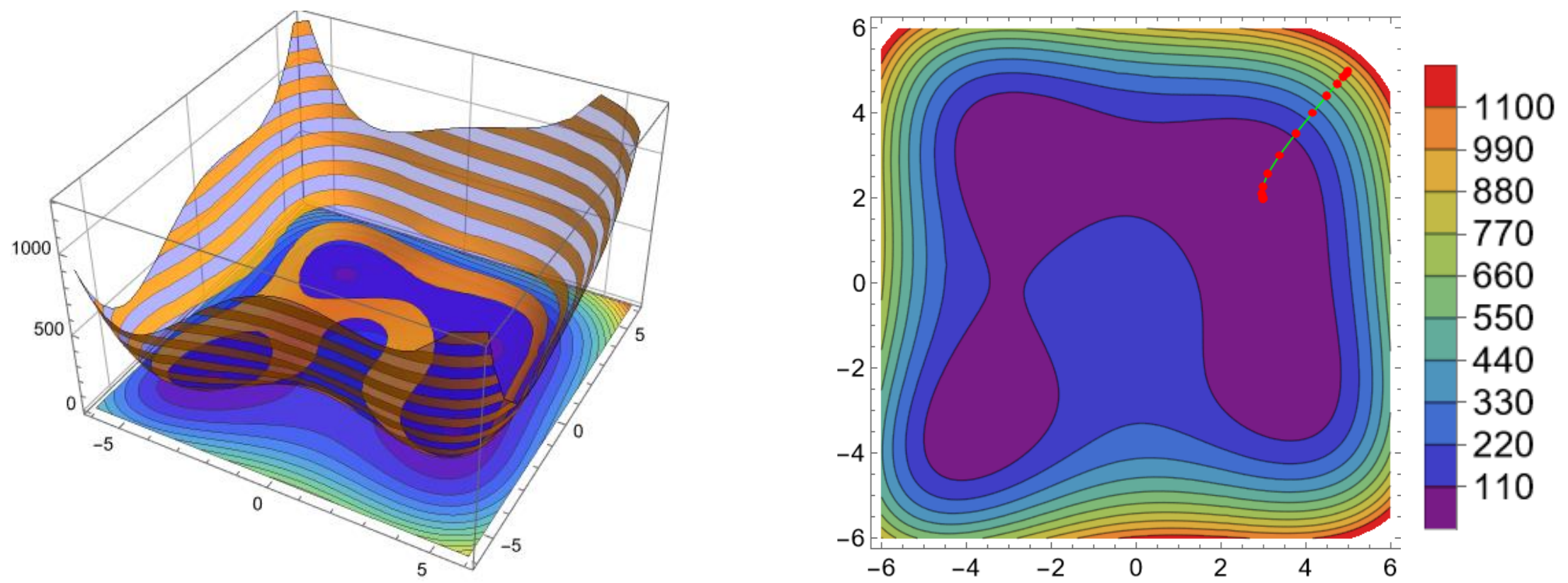

**Figure 8.10.** The results of 15 iterations of the Marquardt method of a vector-valued function, for $\mathbf{f}|\mathbf{x}\rangle = \{x^2 + y - 11, x + y^2 - 7\}$.

From Table 8.7, after 15 iterations, the optimal condition is attained. The results were produced by Mathematica code 8.7.

**Table 8.7.**

| No. of iters. | $x0$ | $x1$ | $f(x0)$ | $f(x1)$ | $\lambda$ | grad $f(x0)$ |
|---|---|---|---|---|---|---|
| 1 | {5., 5.} | {4.958, 4.951} | 890 | 848.7297 | 10000 | 655.3472 |
| 2 | {4.958, 4.951} | {4.879, 4.859} | 848.7297 | 774.1588 | 5000 | 634.6301 |
| 3 | {4.879, 4.859} | {4.738, 4.692} | 774.1588 | 651.0815 | 2500 | 596.4431 |
| 4 | {4.738, 4.692} | {4.504, 4.416} | 651.0815 | 477.1697 | 1250 | 530.9784 |
| 5 | {4.504, 4.416} | {4.166, 4.013} | 477.1697 | 283.7561 | 625 | 431.8054 |
| 6 | {4.166, 4.013} | {3.762, 3.515} | 283.7561 | 127.7414 | 312.5 | 307.3011 |
| 7 | {3.762, 3.515} | {3.379, 3.008} | 127.7414 | 41.25237 | 156.25 | 184.773 |
| 8 | {3.379, 3.008} | {3.109, 2.580} | 41.25237 | 9.220859 | 78.125 | 92.13854 |
| 9 | {3.109, 2.580} | {2.989, 2.279} | 9.220859 | 1.456565 | 39.0625 | 37.54547 |
| 10 | {2.989, 2.279} | {2.975, 2.106} | 1.456565 | 0.170816 | 19.53125 | 12.3052 |
| 11 | {2.975, 2.106} | {2.990, 2.028} | 0.170816 | 0.011524 | 9.765625 | 3.400272 |
| 12 | {2.990, 2.028} | {2.998, 2.004} | 0.011524 | 0.000294 | 4.882813 | 0.793858 |
| 13 | {2.998, 2.004} | {2.999, 2.000} | 0.000294 | 2.23E-06 | 2.441406 | 0.123629 |
| 14 | {2.999, 2.000} | {3, 2} | 2.23E-06 | 4.6E-09 | 1.220703 | 0.010723 |
| 15 | {3, 2} | {3, 2} | 4.6E-09 | 2.47E-12 | 0.610352 | 0.000486 |

***Mathematica Code 8.7*** Marquardt Method

```
(* Marquardt Method and Modification of the Hessian (Approximate) (Vector Valued Function)
(Sum of Squares) *)

(*
Notations  :x0:Intial vector
```

```
epsilonMM  :Small number to check the accuracy of the Marquardt Method and Vector Valued
Function
f[x,y]     :Objective function
lii        :The last iteration index
result[k]  :The results of iteration k
*)

(* Taking Initial Inputs from User *)
x0=Input["Enter the intial point in the format {x, y}; for example {5,5} "] ;
epsilonMM=Input["Please enter accuracy of the Marquardt Method and Vector Valued Function;
for example 0.01 "];

domainx=Input["Please enter domain of x variable for 3D and contour plots; for example {-
6,6}"];
domainy=Input["Please enter domain of y variable for 3D and contour plots; for example {-
6,6}"];

If[
  epsilonMM<=0,
  Beep[];
  MessageDialog["The value of epsilonMM has to be postive number: "];
  Exit[];
  ];

(* Taking the Function from User *)
f[{x_,y_}] = Evaluate[Input["Please input a function of x and y  as a vector valued
function; for example {x^2+y-11,x+y^2-7} "]];
(* VIP VIP Input the f as a vector valued function; for example: {x^2+y-11,x+y^2-7} *)

v={x,y};
jacobianf=Grad[f[{x,y}],v];
functionF[{x_,y_}]:=f[{x,y}].f[{x,y}];

gradF[x_,y_]=2*Transpose[jacobianf].f[{x,y}];
HessianF[x_,y_]=2*Transpose[jacobianf].jacobianf;
λ=10000;

(* Main Loop *)
Do[
  gradFx0=gradF[x0[[1]],x0[[2]]];
  hessianFx0=HessianF[x0[[1]],x0[[2]]];

  If[
   N[Norm[gradFx0]]==0,
   Print["The Marquardt's Method requires that the gradiant of the function f !=0"];(*Ending
program*)
   Exit[];
   ];
  inverseHessianF=Inverse[hessianFx0+λ*IdentityMatrix[2]];
  dx0=-(inverseHessianF.gradFx0);
  x1=x0+dx0;

  lii=k;
  result[k]=N[{k,Row[x0,","],Row[x1,","],functionF[x0],functionF[x1],λ,N[Norm[gradFx0]]}];
  plotresult[k]=N[{x0,x1}];

  If[
   functionF[x1]<functionF[x0],
   λ=0.5*λ;
   x0=x1;,
   λ=2*λ;
```

```
    x0=x1;
    ];

   If[
    N[Norm[gradFx0]]<epsilonMM||k>50,
    Break[];
    ];,
   {k,1,∞}
   ];

(* Final result *)
If[
  lii==50,
  Print[" After 50 iterations the mimimum point around the intial point is " ,N[x1], "\nThe
solution is (approximately) ", N[functionF[x1]]];,
 Print["The solution is x= ",  N[x1],"\nThe solution is (approximately)= ",
N[functionF[x1]]];
  ]

(* Results of Each Iteration *)
table=TableForm[
  Table[
   result[i],
   {i,1,lii}
   ],
  TableHeadings->{None,{"No. of iters.","x0","x1","f[x0]","f[x1]","λ","N[gradfx0]"}}
  ]

Export["example87.xls",table,"XLS"];

(* Data Visualization *)
(* Domain of Varibles*)
xleft=domainx[[1]];
xright=domainx[[2]];
ydown=domainy[[1]];
yup=domainy[[2]];

(* 3D+ Contour Plot *)
plot1=Plot3D[
   functionF[{x,y}],
   {x,xleft,xright},
   {y,ydown,yup},
   ClippingStyle->None,
   MeshFunctions->{#3&},
   Mesh->15,
   MeshStyle->Opacity[.5],
   MeshShading->{{Opacity[.3],Blue},{Opacity[.8],Orange}},
   Lighting->"Neutral"
   ];
slice=SliceContourPlot3D[
   functionF[{x,y}],
   z==0,
   {x,xleft,xright},
   {y,ydown,yup},
   {z,-1,1},
   Contours->15,
   Axes->False,
   PlotPoints->50,
   PlotRangePadding->0,
   ColorFunction->"Rainbow"
   ];
```

```
Show[
 plot1,
 slice,
 PlotRange->All,
 BoxRatios->{1,1,.6},
 FaceGrids->{Back,Left}
 ]

(* Contour Plot with Step Iterations *)
ContourPlot[
 functionF[{x,y}],
 {x,xleft,xright},
 {y,ydown,yup},
 LabelStyle->Directive[Black,16],
 ColorFunction->"Rainbow",
 PlotLegends->Automatic,
 Contours->10,
 Epilog-
>{PointSize[0.015],Green,Line[Flatten[Table[plotresult[i],{i,1,lii}],1]],Red,Point[Flatten[T
able[plotresult[i],{i,1,lii}],1]]}
 ]

(* Data Manipulation *)
Manipulate[
 ContourPlot[
  functionF[{x,y}],
  {x,xleft,xright},
  {y,ydown,yup},
  LabelStyle->Directive[Black,14],
  ColorFunction->"Rainbow",
  PlotLegends->Automatic,
  Contours->10,
  Epilog->{
    PointSize[0.015],
    Yellow,
    Arrow[{plotresult[i][[1]],plotresult[i][[2]]}],
    Red,
    Point[Flatten[Table[plotresult[j],{j,1,i}],1]],
    Green,
    Line[Flatten[Table[plotresult[j],{j,1,i}],1]]
    }
  ],
 {i,1,lii,1}
 ]
```

# CHAPTER 9

# MULTI-VARIABLE CONJUGATE-DIRECTION AND QUASI-NEWTON ALGORITHMS

We have seen before that the steepest descent method tends to be effective far from the minimum but becomes less so as $|\mathbf{x}^*\rangle$ is approached, whereas the Newton method can be unreliable far from $|\mathbf{x}^*\rangle$ but is very efficient as $|\mathbf{x}_k\rangle$ approaches the minimum. In this chapter, we discuss methods that tend to exhibit the positive characteristics of the steepest descent and Newton methods using only first-order derivative [1-5]. An important class of methods of this type is the conjugate-direction method. The conjugate-direction methods will be reliable far from $|\mathbf{x}^*\rangle$ and will accelerate as the sequence of iterates approaches the minimum. In the steepest descent and Newton methods, the direction of search in each iteration depends on the local properties of the objective function. However, a relationship may exist between successive search directions. In conjugate-direction methods, the optimization is performed by using a sequence of search directions that bear a relationship to one another. In this chapter, we also will discuss the quasi-Newton methods. The foundation of the quasi-Newton methods is the classical Newton method. The basic idea in quasi-Newton methods is that the direction of search is based on an $n \times n$ direction matrix $\mathbf{S}$ which is the approximate of the inverse Hessian matrix by some other matrix that should be positive definite. This saves the work of computation of second derivatives and also avoids the difficulties associated with the loss of positive definiteness.

## 9.1 Conjugate Direction Method

We have previously given the conjugacy conditions for a set of directions $|\mathbf{d}_i\rangle$, $i = 1, 2, 3, \ldots, r \leq n$, and an $n \times n$ symmetric matrix $\mathbf{H}$. A finite set of distinct nonzero vectors $\{|\mathbf{d}_0\rangle, |\mathbf{d}_1\rangle, \ldots, |\mathbf{d}_m\rangle\}$ is said to be conjugate with respect to a real symmetric matrix $\mathbf{H}$, if

$$\langle \mathbf{d}_i | \mathbf{H} | \mathbf{d}_j \rangle = 0 \quad \text{for all } i \neq j. \tag{9.1}$$

The $\mathbf{H}$-conjugate directions $|\mathbf{d}_0\rangle, |\mathbf{d}_1\rangle, \ldots, |\mathbf{d}_m\rangle$ form a set of linearly independent vectors. Moreover, if $\mathbf{H}$ is the symmetric positive definite matrix, then the $\mathbf{H}$-conjugate vectors $|\mathbf{d}_i\rangle$, $i = 0,1,2, \cdots, n-1$ form a basis for $\mathbb{R}^n$. To see this fact, let $\alpha_i$, $i = 0,1, \ldots, k$ be scalar, such that,

$$\alpha_0 |\mathbf{d}_0\rangle + \alpha_1 |\mathbf{d}_1\rangle + \cdots + \alpha_k |\mathbf{d}_k\rangle = \sum_{i=0}^{k} \alpha_i |\mathbf{d}_i\rangle = |\mathbf{0}\rangle.$$

Hence

$$\sum_{i=0}^{k} \alpha_i \langle \mathbf{d}_j | \mathbf{H} | \mathbf{d}_i \rangle = \alpha_j \langle \mathbf{d}_j | \mathbf{H} | \mathbf{d}_j \rangle = \langle \mathbf{d}_j | \mathbf{H} | \mathbf{0} \rangle = 0.$$

Since $\mathbf{H}$ is a symmetric positive definite matrix and $|\mathbf{d}_j\rangle \neq |\mathbf{0}\rangle$, we get $\alpha_j = 0$, that is,

$$\alpha_0 = 0, \alpha_1 = 0, \ldots, \alpha_k = 0.$$

As a result, the set of distinct nonzero vectors $\{|\mathbf{d}_0\rangle, |\mathbf{d}_1\rangle, \ldots, |\mathbf{d}_k\rangle\}$ are linearly independent.

***Example 9.1***

Construct a set of **H**-conjugate vectors $|\mathbf{d}_0\rangle, |\mathbf{d}_1\rangle, |\mathbf{d}_2\rangle$, where $|\mathbf{d}_0\rangle = (1,0,0)^T$ and the matrix $\mathbf{H} = \begin{pmatrix} 3 & 0 & 1 \\ 0 & 4 & 1 \\ 1 & 1 & 3 \end{pmatrix}$.

***Solution***

First, note that the matrix **H** is positive definite. Let $|\mathbf{d}_1\rangle = (d_{11}, d_{12}, d_{13})^T$, and $|\mathbf{d}_2\rangle = (d_{21}, d_{22}, d_{23})^T$. We require $\langle\mathbf{d}_0|\mathbf{H}|\mathbf{d}_1\rangle = 0$. We have

$$\langle\mathbf{d}_0|\mathbf{H}|\mathbf{d}_1\rangle = (1,0,0)\begin{pmatrix} 3 & 0 & 1 \\ 0 & 4 & 1 \\ 1 & 1 & 3 \end{pmatrix}\begin{pmatrix} d_{11} \\ d_{12} \\ d_{13} \end{pmatrix} = (3,0,1)\begin{pmatrix} d_{11} \\ d_{12} \\ d_{13} \end{pmatrix} = 3d_{11} + d_{13}.$$

Let $d_{11} = 1$, $d_{12} = 0$ and $d_{13} = -3$. Then $|\mathbf{d}_1\rangle = (1,0,-3)^T$, and thus $\langle\mathbf{d}_0|\mathbf{H}|\mathbf{d}_1\rangle = 0$.

To find the third vector $|\mathbf{d}_2\rangle = (d_{21}, d_{22}, d_{23})^T$, which would be **H**-conjugate with $|\mathbf{d}_0\rangle$ and $|\mathbf{d}_1\rangle$, we require $\langle\mathbf{d}_0|\mathbf{H}|\mathbf{d}_2\rangle = 0$ and $\langle\mathbf{d}_1|\mathbf{H}|\mathbf{d}_2\rangle = 0$. We have

$$\langle\mathbf{d}_0|\mathbf{H}|\mathbf{d}_2\rangle = (1,0,0)\begin{pmatrix} 3 & 0 & 1 \\ 0 & 4 & 1 \\ 1 & 1 & 3 \end{pmatrix}\begin{pmatrix} d_{21} \\ d_{22} \\ d_{23} \end{pmatrix} = (3,0,1)\begin{pmatrix} d_{21} \\ d_{22} \\ d_{23} \end{pmatrix} = 3d_{21} + d_{23} = 0,$$

and

$$\langle\mathbf{d}_1|\mathbf{H}|\mathbf{d}_2\rangle = (1,0,-3)\begin{pmatrix} 3 & 0 & 1 \\ 0 & 4 & 1 \\ 1 & 1 & 3 \end{pmatrix}\begin{pmatrix} d_{21} \\ d_{22} \\ d_{23} \end{pmatrix} = (0,-3,-8)\begin{pmatrix} d_{21} \\ d_{22} \\ d_{23} \end{pmatrix} = -3d_{22} - 8d_{23} = 0.$$

If we take $|\mathbf{d}_2\rangle = (1,8,-3)^T$, then the resulting set of vectors is mutually conjugate.

Consider the quadratic objective function $f\colon \mathbb{R}^n \to \mathbb{R}$ defined by

$$f|\mathbf{x}\rangle = \frac{1}{2}\langle\mathbf{x}|\mathbf{H}|\mathbf{x}\rangle + \langle\mathbf{b}|\mathbf{x}\rangle + c, \tag{9.2}$$

where **H** is an $n \times n$ symmetric positive definite matrix, $|\mathbf{b}\rangle$ is an $n \times 1$ vector, and $c$ is a real number. The function $f$ has a gradient vector $|\mathbf{g}\rangle$ given by: $\mathbf{g}|\mathbf{x}\rangle = \nabla f|\mathbf{x}\rangle = \mathbf{H}|\mathbf{x}\rangle + |\mathbf{b}\rangle$. If $f$ has the strict global minimizer $|\mathbf{x}^*\rangle$ in $\mathbb{R}^n$, then $\mathbf{g}|\mathbf{x}^*\rangle = |\mathbf{0}\rangle$. For a quadratic function of $n$ variables, the conjugate direction method reaches the solution after $n$ steps, as in the following theorem.

**Theorem 9.1 (Conjugate Direction Algorithm):** Let $\{|\mathbf{d}_i\rangle\}_{i=0}^{n-1}$ be a set of nonzero **H**-conjugate directions. For any $|\mathbf{x}_0\rangle \in \mathbb{R}^n$ the sequence $\{|\mathbf{x}_k\rangle\}$ generated according to

$$|\mathbf{x}_{k+1}\rangle = |\mathbf{x}_k\rangle + \alpha_k|\mathbf{d}_k\rangle, \qquad k \geq 0, \tag{9.3}$$

with

$$\alpha_k = -\frac{\langle\mathbf{g}_k|\mathbf{d}_k\rangle}{\langle\mathbf{d}_k|\mathbf{H}|\mathbf{d}_k\rangle}, \tag{9.4}$$

and

$$|\mathbf{g}_k\rangle = \mathbf{H}|\mathbf{x}_k\rangle + |\mathbf{b}\rangle, \tag{9.5}$$

converges to a unique solution, $|\mathbf{x}^*\rangle$ after $n$ steps, that is, $|\mathbf{x}_n\rangle = |\mathbf{x}^*\rangle$.

**Proof:**

Since the $|\mathbf{d}_k\rangle$ are linearly independent, then there exist constants $\alpha_i$, $i = 0,1,\dots,n-1$ such that,

$$|\mathbf{x}^*\rangle - |\mathbf{x}_0\rangle = \alpha_0|\mathbf{d}_0\rangle + \alpha_1|\mathbf{d}_1\rangle + \cdots + \alpha_{n-1}|\mathbf{d}_{n-1}\rangle = \sum_{i=0}^{n-1}\alpha_i|\mathbf{d}_i\rangle,$$

We multiply by **H** and take the scalar product with $|\mathbf{d}_k\rangle$ to find

$$\begin{aligned}\langle\mathbf{d}_k|\mathbf{H}|\mathbf{x}^* - \mathbf{x}_0\rangle &= \langle\mathbf{d}_k|\mathbf{H}\sum_{i=0}^{n-1}\alpha_i|\mathbf{d}_i\rangle \\ &= \sum_{i=0}^{n-1}\alpha_i\langle\mathbf{d}_k|\mathbf{H}|\mathbf{d}_i\rangle = \alpha_k\langle\mathbf{d}_k|\mathbf{H}|\mathbf{d}_k\rangle,\end{aligned}$$

or

$$\alpha_k = \frac{\langle \mathbf{d}_k|\mathbf{H}|\mathbf{x}^* - \mathbf{x}_0\rangle}{\langle \mathbf{d}_k|\mathbf{H}|\mathbf{d}_k\rangle}.$$

Using $|\mathbf{x}_{k+1}\rangle = |\mathbf{x}_k\rangle + \alpha_k|\mathbf{d}_k\rangle$, the iterative process from $|\mathbf{x}_0\rangle$ up to $|\mathbf{x}_k\rangle$ gives

$$\begin{aligned}
&|\mathbf{x}_1\rangle = |\mathbf{x}_0\rangle + \alpha_0|\mathbf{d}_0\rangle,\\
&|\mathbf{x}_2\rangle = |\mathbf{x}_1\rangle + \alpha_1|\mathbf{d}_1\rangle = |\mathbf{x}_0\rangle + \alpha_0|\mathbf{d}_0\rangle + \alpha_1|\mathbf{d}_1\rangle,\\
&|\mathbf{x}_3\rangle = |\mathbf{x}_2\rangle + \alpha_2|\mathbf{d}_2\rangle = |\mathbf{x}_0\rangle + \alpha_0|\mathbf{d}_0\rangle + \alpha_1|\mathbf{d}_1\rangle + \alpha_2|\mathbf{d}_2\rangle,\\
&\dots\ \dots\ \dots\\
&|\mathbf{x}_k\rangle = |\mathbf{x}_0\rangle + \sum_{i=0}^{k-1}\alpha_i|\mathbf{d}_i\rangle.
\end{aligned}$$

Moreover, by the $\mathbf{H}$ -orthogonality of the $|\mathbf{d}_k\rangle$ it follows that

$$\langle \mathbf{d}_k|\mathbf{H}|\mathbf{x}_k - \mathbf{x}_0\rangle = \langle \mathbf{d}_k|\mathbf{H}\sum_{i=0}^{k-1}\alpha_i|\mathbf{d}_i\rangle = \sum_{i=0}^{k-1}\alpha_i\langle \mathbf{d}_k|\mathbf{H}|\mathbf{d}_i\rangle = 0,$$

and

$$\langle \mathbf{d}_k|\mathbf{H}|\mathbf{x}_k\rangle = \langle \mathbf{d}_k|\mathbf{H}|\mathbf{x}_0\rangle.$$

Hence, we get

$$\begin{aligned}
\alpha_k &= \frac{\langle \mathbf{d}_k|\mathbf{H}|\mathbf{x}^* - \mathbf{x}_0\rangle}{\langle \mathbf{d}_k|\mathbf{H}|\mathbf{d}_k\rangle}\\
&= \frac{\langle \mathbf{d}_k|\mathbf{H}|\mathbf{x}^*\rangle - \langle \mathbf{d}_k|\mathbf{H}|\mathbf{x}_0\rangle}{\langle \mathbf{d}_k|\mathbf{H}|\mathbf{d}_k\rangle}\\
&= \frac{\langle \mathbf{d}_k|\mathbf{H}|\mathbf{x}^*\rangle - \langle \mathbf{d}_k|\mathbf{H}|\mathbf{x}_k\rangle}{\langle \mathbf{d}_k|\mathbf{H}|\mathbf{d}_k\rangle}\\
&= \frac{\langle \mathbf{d}_k|\mathbf{H}|\mathbf{x}^* - \mathbf{x}_k\rangle}{\langle \mathbf{d}_k|\mathbf{H}|\mathbf{d}_k\rangle}.
\end{aligned}$$

From $|\mathbf{g}_k\rangle = \mathbf{H}|\mathbf{x}_k\rangle + |\mathbf{b}\rangle$, we have $\mathbf{H}|\mathbf{x}_k\rangle = |\mathbf{g}_k\rangle - |\mathbf{b}\rangle$, and since $|\mathbf{g}_k\rangle = |\mathbf{0}\rangle$ at minimizer $|\mathbf{x}^*\rangle$ , we get $\mathbf{H}|\mathbf{x}^*\rangle = -|\mathbf{b}\rangle$. Therefore,

$$\begin{aligned}
\alpha_k &= \frac{\langle \mathbf{d}_k|\mathbf{H}|\mathbf{x}^* - \mathbf{x}_k\rangle}{\langle \mathbf{d}_k|\mathbf{H}|\mathbf{d}_k\rangle}\\
&= \frac{\langle \mathbf{d}_k|\mathbf{H}\mathbf{x}^* - \mathbf{H}\mathbf{x}_k\rangle}{\langle \mathbf{d}_k|\mathbf{H}|\mathbf{d}_k\rangle}\\
&= \frac{\langle \mathbf{d}_k|-\mathbf{b} - \mathbf{H}\mathbf{x}_k\rangle}{\langle \mathbf{d}_k|\mathbf{H}|\mathbf{d}_k\rangle}\\
&= -\frac{\langle \mathbf{d}_k|\mathbf{b} + \mathbf{H}\mathbf{x}_k\rangle}{\langle \mathbf{d}_k|\mathbf{H}|\mathbf{d}_k\rangle}\\
&= -\frac{\langle \mathbf{d}_k|\mathbf{g}_k\rangle}{\langle \mathbf{d}_k|\mathbf{H}|\mathbf{d}_k\rangle}.
\end{aligned}$$

Finally, we have

$$\begin{aligned}
|\mathbf{x}^*\rangle &= |\mathbf{x}_0\rangle + \alpha_0|\mathbf{d}_0\rangle + \alpha_1|\mathbf{d}_1\rangle + \cdots + \alpha_{n-1}|\mathbf{d}_{n-1}\rangle\\
&= |\mathbf{x}_0\rangle + \sum_{i=0}^{n-1}\alpha_i|\mathbf{d}_i\rangle\\
&= |\mathbf{x}_0\rangle + |\mathbf{x}_n\rangle - |\mathbf{x}_0\rangle\\
&= |\mathbf{x}_n\rangle.
\end{aligned}$$

■

**Algorithm** (Conjugate Direction Algorithm for the Quadratic Function)

**Step 1:** Given an initial point $|\mathbf{x}_0\rangle \in \mathbb{R}^n$ and a tolerance value $0 < \varepsilon < 1$.

**Step 2:** Choose $\mathbf{H}$-conjugate directions $|\mathbf{d}_0\rangle, |\mathbf{d}_1\rangle, \dots, |\mathbf{d}_{n-1}\rangle$.

**Step 3:** Set the iteration $k = 0$.

**Step 4:** Compute $|\mathbf{g}_k\rangle$ the gradient of the function at $|\mathbf{x}_k\rangle$.

**Step 5:** Compute $\alpha_k$ by the formula $\alpha_k = -\frac{\langle \mathbf{g}_k|\mathbf{d}_k\rangle}{\langle \mathbf{d}_k|\mathbf{H}|\mathbf{d}_k\rangle}$.

**Step 6:** Compute a new point $|\mathbf{x}_{k+1}\rangle$ by the formula $|\mathbf{x}_{k+1}\rangle = |\mathbf{x}_k\rangle + \alpha_k|\mathbf{d}_k\rangle$.

**Step 7:** If $\|\mathbf{g}_{k+1}\| < \varepsilon$, then $|\mathbf{x}_{k+1}\rangle = |\mathbf{x}^*\rangle$, which is the minimizer.

**Step 8:** Set $k = k + 1$ and go to step 4.

***Example 9.2***

Let $f\colon \mathbb{R}^2 \to \mathbb{R}$ be defined by

$$f|\mathbf{x}\rangle = \frac{1}{2}\langle \mathbf{x}|\mathbf{H}|\mathbf{x}\rangle + \langle \mathbf{b}|\mathbf{x}\rangle + c$$

where

$$\mathbf{H} = \begin{pmatrix} 2 & 0 \\ 0 & 4 \end{pmatrix}, |\mathbf{b}\rangle = \begin{pmatrix} -1 \\ 1 \end{pmatrix}, c = 1, |\mathbf{x}_0\rangle = \begin{pmatrix} 0 \\ 0 \end{pmatrix},$$
$$|\mathbf{d}_0\rangle = \begin{pmatrix} 1 \\ 0 \end{pmatrix}, |\mathbf{d}_1\rangle = \begin{pmatrix} 0 \\ 1 \end{pmatrix}.$$

Verify that the iteration $|\mathbf{x}_2\rangle$ generated from the general conjugate direction method is the minimizer of $f$. Take $\varepsilon = 0.0001$.

***Solution***

We can write $f|\mathbf{x}\rangle$ as follows

$$\begin{aligned} f|\mathbf{x}\rangle &= \frac{1}{2}(x_1, x_2)\begin{pmatrix} 2 & 0 \\ 0 & 4 \end{pmatrix}\begin{pmatrix} x_1 \\ x_2 \end{pmatrix} - (-1,1)\begin{pmatrix} x_1 \\ x_2 \end{pmatrix} + 1 \\ &= \frac{1}{2}(2x_1, 4x_2)\begin{pmatrix} x_1 \\ x_2 \end{pmatrix} + x_1 - x_2 + 1 \\ &= \frac{1}{2}(2x_1^2 + 4x_2^2) + x_1 - x_2 + 1 \\ &= x_1^2 + 2x_2^2 + x_1 - x_2 + 1. \end{aligned}$$

Now, the gradient vector $\mathbf{g}|\mathbf{x}\rangle$ is

$$\mathbf{g}|\mathbf{x}\rangle = \left(\frac{\partial f}{\partial x_1}, \frac{\partial f}{\partial x_2}\right)^T = (2x_1 + 1, 4x_2 - 1)^T,$$

so that

$$\mathbf{g}_0 = \mathbf{g}|\mathbf{x}_0\rangle = (1, -1)^T.$$

The values of $\alpha_0$ and $|\mathbf{x}_1\rangle$ are

$$\alpha_0 = -\frac{\langle \mathbf{g}_0|\mathbf{d}_0\rangle}{\langle \mathbf{d}_0|\mathbf{H}|\mathbf{d}_0\rangle} = -\frac{(1,-1)\begin{pmatrix} 1 \\ 0 \end{pmatrix}}{(1,0)\begin{pmatrix} 2 & 1 \\ 1 & 2 \end{pmatrix}\begin{pmatrix} 1 \\ 0 \end{pmatrix}} = -\frac{1}{2},$$

$$|\mathbf{x}_1\rangle = |\mathbf{x}_0\rangle + \alpha_0|\mathbf{d}_0\rangle = \begin{pmatrix} 0 \\ 0 \end{pmatrix} - \frac{1}{2}\begin{pmatrix} 1 \\ 0 \end{pmatrix} = \begin{pmatrix} -\frac{1}{2} \\ 0 \end{pmatrix},$$

and therefore

$$\mathbf{g}_1 = \mathbf{g}|\mathbf{x}_1\rangle = (0, -1)^T \Rightarrow \|\mathbf{g}_1\| = 1 > \varepsilon.$$

Compute $\alpha_1$ and $|\mathbf{x}_2\rangle$ by

$$\alpha_1 = -\frac{\langle \mathbf{g}_1|\mathbf{d}_1\rangle}{\langle \mathbf{d}_1|\mathbf{H}|\mathbf{d}_1\rangle} = -\frac{(0,-1)\begin{pmatrix}0\\1\end{pmatrix}}{(0,1)\begin{pmatrix}2 & 0\\0 & 4\end{pmatrix}\begin{pmatrix}0\\1\end{pmatrix}} = -\frac{-1}{4} = \frac{1}{4},$$

$$|\mathbf{x}_2\rangle = |\mathbf{x}_1\rangle + \alpha_1|\mathbf{d}_1\rangle = \begin{pmatrix}-\frac{1}{2}\\0\end{pmatrix} + \frac{1}{4}\begin{pmatrix}0\\1\end{pmatrix} = \begin{pmatrix}-\frac{1}{2}\\0\end{pmatrix} + \begin{pmatrix}0\\\frac{1}{4}\end{pmatrix} = \begin{pmatrix}-\frac{1}{2}\\\frac{1}{4}\end{pmatrix},$$

and

$$\mathbf{g}_2 = \mathbf{g}|\mathbf{x}_2\rangle = \begin{pmatrix}0\\0\end{pmatrix} \Rightarrow \|\mathbf{g}_2\| = 0 < \varepsilon.$$

Therefore, $|\mathbf{x}_2\rangle$ is the minimizer of $f$. Since $f$ is a quadratic function of two variables, therefore it converges to the minimum point$|\mathbf{x}_2\rangle$ and the minimum value of the objective function $f|\mathbf{x}^*\rangle = 5/8$ in two iterations.

As we shall see below, the method also possesses an important property in the intermediate steps [1]. To see this, suppose that we start at $|\mathbf{x}_0\rangle$ and search in the direction $|\mathbf{d}_0\rangle$ to obtain

$$|\mathbf{x}_1\rangle = |\mathbf{x}_0\rangle - \left(\frac{\langle \mathbf{g}_0|\mathbf{d}_0\rangle}{\langle \mathbf{d}_0|\mathbf{H}|\mathbf{d}_0\rangle}\right)|\mathbf{d}_0\rangle. \tag{9.6}$$

We claim that

$$\langle \mathbf{g}_1|\mathbf{d}_0\rangle = 0. \tag{9.7}$$

To see this,

$$\begin{aligned}\langle \mathbf{g}_1|\mathbf{d}_0\rangle &= \langle \mathbf{H}\mathbf{x}_1 - \mathbf{b}|\mathbf{d}_0\rangle \\ &= \left\langle \mathbf{H}\left(\mathbf{x}_0 - \left(\frac{\langle \mathbf{g}_0|\mathbf{d}_0\rangle}{\langle \mathbf{d}_0|\mathbf{H}|\mathbf{d}_0\rangle}\right)\mathbf{d}_0\right) - \mathbf{b}\middle|\mathbf{d}_0\right\rangle \\ &= \langle \mathbf{H}\mathbf{x}_0|\mathbf{d}_0\rangle - \left\langle\left(\frac{\langle \mathbf{g}_0|\mathbf{d}_0\rangle}{\langle \mathbf{d}_0|\mathbf{H}|\mathbf{d}_0\rangle}\right)\mathbf{H}\mathbf{d}_0\middle|\mathbf{d}_0\right\rangle - \langle \mathbf{b}|\mathbf{d}_0\rangle \\ &= \langle \mathbf{x}_0|\mathbf{H}|\mathbf{d}_0\rangle - \left(\frac{\langle \mathbf{g}_0|\mathbf{d}_0\rangle}{\langle \mathbf{d}_0|\mathbf{H}|\mathbf{d}_0\rangle}\right)\langle \mathbf{d}_0|\mathbf{H}|\mathbf{d}_0\rangle - \langle \mathbf{b}|\mathbf{d}_0\rangle \\ &= \langle \mathbf{H}\mathbf{x}_0 - \mathbf{b}|\mathbf{d}_0\rangle - \langle \mathbf{g}_0|\mathbf{d}_0\rangle \\ &= \langle \mathbf{g}_0|\mathbf{d}_0\rangle - \langle \mathbf{g}_0|\mathbf{d}_0\rangle = 0.\end{aligned} \tag{9.8}$$

The equation $\langle \mathbf{g}_1|\mathbf{d}_0\rangle = 0$ implies that $\alpha_0$ has the property that $\alpha_0 = \arg\min \phi_0(\alpha)$, where $\phi_0(\alpha) = f|\mathbf{x}_0 + \alpha\mathbf{d}_0\rangle$. To see this, apply the chain rule to get

$$\frac{d}{d\alpha}\phi_0(\alpha) = \langle \nabla f(\mathbf{x_0} + \alpha\mathbf{d}_0)|\mathbf{d}_0\rangle, \tag{9.9}$$

at $\alpha = \alpha_0$, we get

$$\frac{d}{d\alpha}\phi_0(\alpha_0) = \langle \mathbf{g}_1|\mathbf{d}_0\rangle = 0. \tag{9.10}$$

Because $\phi_0$ is a quadratic function of $\alpha$, and the coefficient of the $\alpha^2$ term in $\phi_0$ is $\langle \mathbf{d}_0|\mathbf{H}|\mathbf{d}_0\rangle > 0$, the above implies that $\alpha_0 = \arg\min_{\alpha\in\mathbb{R}} \phi_0(\alpha)$. Using a similar argument, we can show that for all $k$,

$$\langle \mathbf{g}_{k+1}|\mathbf{d}_k\rangle = 0, \tag{9.11}$$

and hence

$$\alpha_k = \arg\min f|\mathbf{x}_k + \alpha\mathbf{d}_k\rangle. \tag{9.12}$$

In fact, an even stronger condition holds, as given by the following theorem.

**Theorem 9.2 (Orthogonality of Gradient to a Set of Conjugate Directions):** In the conjugate direction algorithm,

$$\langle \mathbf{g}_{k+1}|\mathbf{d}_i\rangle = 0, \tag{9.13}$$

for all $k$, $0 \le k \le n-1$, and $0 \le i \le k$.

**Proof:**

First, note that,

$$\mathbf{H}|\mathbf{x}_{k+1}-\mathbf{x}_k\rangle = \mathbf{H}|\mathbf{x}_{k+1}\rangle - |\mathbf{b}\rangle - (\mathbf{H}|\mathbf{x}_k\rangle - |\mathbf{b}\rangle) = |\mathbf{g}_{k+1}\rangle - |\mathbf{g}_k\rangle,$$

because $|\mathbf{g}_k\rangle = \mathbf{H}|\mathbf{x}_k\rangle - |\mathbf{b}\rangle$. Thus,

$$|\mathbf{g}_{k+1}\rangle = |\mathbf{g}_k\rangle + \alpha_k \mathbf{H}|\mathbf{d}_k\rangle,$$

where $|\mathbf{x}_{k+1}-\mathbf{x}_k\rangle = \alpha_k|\mathbf{d}_k\rangle$. We prove this theorem by induction. The result is true for $k=0$ because $\langle \mathbf{g}_1|\mathbf{d}_0\rangle = 0$, as shown before. We now show that if the result is true for $k-1$ (i.e., $\langle \mathbf{g}_k|\mathbf{d}_i\rangle = 0, i \le k-1$ ), then it is true for $k$ (i.e., $\langle \mathbf{g}_{k+1}|\mathbf{d}_i\rangle = 0$, $i \le k$ ). Fix $k > 0$ and $0 \le i < k$. By the induction hypothesis, $\langle \mathbf{g}_k|\mathbf{d}_i\rangle = 0$. Because

$$|\mathbf{g}_{k+1}\rangle = |\mathbf{g}_k\rangle + \alpha_k \mathbf{H}|\mathbf{d}_k\rangle,$$

and $\langle \mathbf{d}_k|\mathbf{H}|\mathbf{d}_i\rangle = 0$ by $\mathbf{H}$-conjugacy, we have

$$\langle \mathbf{g}_{k+1}|\mathbf{d}_i\rangle = \langle \mathbf{g}_k|\mathbf{d}_i\rangle + \alpha_k\langle \mathbf{d}_k|\mathbf{H}|\mathbf{d}_i\rangle = 0.$$

It remains to be shown that $\langle \mathbf{g}_{k+1}|\mathbf{d}_k\rangle = 0$. Indeed,

$$\begin{aligned}
\langle \mathbf{g}_{k+1}|\mathbf{d}_k\rangle &= \langle \mathbf{H}\mathbf{x}_{k+1} - \mathbf{b}|\mathbf{d}_k\rangle \\
&= \langle \mathbf{x}_{k+1}|\mathbf{H}|\mathbf{d}_k\rangle - \langle \mathbf{b}|\mathbf{d}_k\rangle \\
&= \left\langle \mathbf{x}_k - \left(\frac{\langle \mathbf{g}_k|\mathbf{d}_k\rangle}{\langle \mathbf{d}_k|\mathbf{H}|\mathbf{d}_k\rangle}\right)\mathbf{d}_k \middle| \mathbf{H} \middle| \mathbf{d}_k \right\rangle - \langle \mathbf{b}|\mathbf{d}_k\rangle \\
&= \langle \mathbf{x}_k|\mathbf{H}|\mathbf{d}_k\rangle - \left(\frac{\langle \mathbf{g}_k|\mathbf{d}_k\rangle}{\langle \mathbf{d}_k|\mathbf{H}|\mathbf{d}_k\rangle}\right)\langle \mathbf{d}_k|\mathbf{H}|\mathbf{d}_k\rangle - \langle \mathbf{b}|\mathbf{d}_k\rangle \\
&= \langle \mathbf{H}\mathbf{x}_k|\mathbf{d}_k\rangle - \langle \mathbf{b}|\mathbf{d}_k\rangle - \langle \mathbf{g}_k|\mathbf{d}_k\rangle \\
&= \langle \mathbf{H}\mathbf{x}_k - \mathbf{b}|\mathbf{d}_k\rangle - \langle \mathbf{g}_k|\mathbf{d}_k\rangle = 0,
\end{aligned}$$

because $\mathbf{H}|\mathbf{x}_k\rangle - \mathbf{b} = |\mathbf{g}_k\rangle$. Therefore, by induction, for all $0 \le k \le n-1$ and $0 \le i \le k$,

$$\langle \mathbf{g}_{k+1}|\mathbf{d}_i\rangle = 0.$$

∎

Theorem 9.2 says that under given conditions, the gradient vector at $|\mathbf{x}_{k+1}\rangle$ is orthogonal to all of the preceding directions $|\mathbf{d}_i\rangle$, $i = 0,1,\cdots,k$. Figure 9.1 illustrates this statement.

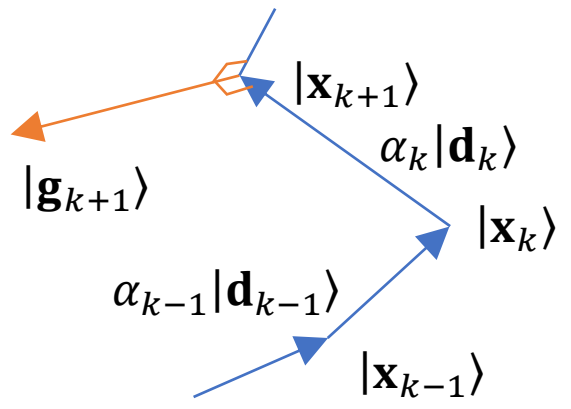

**Figure 9.1.** Illustration of Theorem 9.2.

## 9.2 Conjugate Gradient Method

The conjugate gradient method [6] is the modification of the steepest descent method. It uses conjugate directions to minimize a quadratic function, but the steepest descent method uses local gradients to minimize the quadratic function. The conjugate gradient method does not repeat any previous search directions and converges in $n$ iterations for $n$ variables of quadratic functions.

Consider the quadratic function $f(\mathbf{x}) = \frac{1}{2}\langle \mathbf{x}|\mathbf{H}|\mathbf{x}\rangle - \langle \mathbf{x}|\mathbf{b}\rangle$ where $\mathbf{H}$ is an $n \times n$ symmetric positive definite matrix, and $|\mathbf{b}\rangle$ is an $n \times 1$ vector. Our first search direction from an initial point $|\mathbf{x}_0\rangle$ is in the direction of the steepest descent; that is $|\mathbf{d}_0\rangle = -|\mathbf{g}_0\rangle$. Thus

$$|\mathbf{x}_1\rangle = |\mathbf{x}_0\rangle + \alpha_0|\mathbf{d}_0\rangle, \quad (9.14)$$

where $\alpha_0 = \min_{\alpha \ge 0} f(|\mathbf{x}_0\rangle + \alpha|\mathbf{d}_0\rangle) = -\frac{\langle \mathbf{g}_0|\mathbf{d}_0\rangle}{\langle \mathbf{d}_0|\mathbf{H}|\mathbf{d}_0\rangle}$. In the next step, we search for a direction $|\mathbf{d}_1\rangle$ that is $\mathbf{H}$-conjugate to $|\mathbf{d}_0\rangle$. We choose $|\mathbf{d}_1\rangle$ as a linear combination of $|\mathbf{g}_1\rangle$ and $|\mathbf{d}_0\rangle$ that can be expressed as $|\mathbf{d}_1\rangle = -|\mathbf{g}_1\rangle + \beta_0|\mathbf{d}_0\rangle$.

In general, at the $(k+1)$th step, we have $|\mathbf{x}_{k+1}\rangle = |\mathbf{x}_k\rangle + \alpha_k|\mathbf{d}_k\rangle$, where

$$\alpha_k = -\frac{\langle \mathbf{g}_k|\mathbf{d}_k\rangle}{\langle \mathbf{d}_k|\mathbf{H}|\mathbf{d}_k\rangle}, \tag{9.15}$$

and we choose $|\mathbf{d}_{k+1}\rangle$ to be a linear combination of $|\mathbf{g}_{k+1}\rangle$ and $|\mathbf{d}_k\rangle$, which can be presented

$$|\mathbf{d}_{k+1}\rangle = -|\mathbf{g}_{k+1}\rangle + \beta_k|\mathbf{d}_k\rangle, \quad k = 0,1,2,\cdots. \tag{9.16}$$

The coefficient $\beta_k$, $k = 0,1,\cdots$, are chosen in such a way that $|\mathbf{d}_{k+1}\rangle$ is $\mathbf{H}$-conjugate to $|\mathbf{d}_0\rangle, |\mathbf{d}_1\rangle, \cdots, |\mathbf{d}_k\rangle$. This is accomplished by choosing $\beta_k$ as follows:

$$\beta_k = \frac{\langle \mathbf{g}_{k+1}|\mathbf{H}|\mathbf{d}_k\rangle}{\langle \mathbf{d}_k|\mathbf{H}|\mathbf{d}_k\rangle}. \tag{9.17}$$

Now, we present the conjugate gradient algorithm.

**Algorithm** (Conjugate Gradient Algorithm for Quadratic Function)

**Step 1:** Choose a starting point $|\mathbf{x}_0\rangle \in \mathbb{R}^n$, tolerance value $0 < \epsilon < 1$, set $k = 0$.

**Step 2:** Compute $|\mathbf{g}_0\rangle$, which is the gradient vector at the point $|\mathbf{x}_0\rangle$. If $|\mathbf{g}_0\rangle = |\mathbf{0}\rangle$ stop else go to step 3.

**Step 3:** Set $|\mathbf{d}_0\rangle = -|\mathbf{g}_0\rangle$.

**Step 4:** Compute $\mathbf{H}$, which is the Hessian matrix for the function $f|\mathbf{x}\rangle$ at $|\mathbf{x}_k\rangle$.

**Step 5:** Compute $\alpha_k$ as: $\alpha_k = -\frac{\langle \mathbf{g}_k|\mathbf{d}_k\rangle}{\langle \mathbf{d}_k|\mathbf{H}|\mathbf{d}_k\rangle}$.

**Step 6:** Compute $|\mathbf{x}_{k+1}\rangle$ as: $|\mathbf{x}_{k+1}\rangle = |\mathbf{x}_k\rangle + \alpha_k|\mathbf{d}_k\rangle$.

**Step 7:** Compute $|\mathbf{g}_{k+1}\rangle$. If $\|\mathbf{g}_{k+1}\| < \varepsilon$ stop, else go to step 8.

**Step 8:** Compute $\beta_k$ by $\beta_k = \frac{\langle \mathbf{g}_{k+1}|\mathbf{H}|\mathbf{d}_k\rangle}{\langle \mathbf{d}_k|\mathbf{H}|\mathbf{d}_k\rangle}$ ($\beta_k$ is chosen such that $|\mathbf{d}_k\rangle$ is $\mathbf{H}$-conjugate with $|\mathbf{d}_{k+1}\rangle$).

**Step 9:** Compute $|\mathbf{d}_{k+1}\rangle = -|\mathbf{g}_{k+1}\rangle + \beta_k|\mathbf{d}_k\rangle$.

**Step 10:** Set $k = k + 1$, and go to step 4.

***Example 9.3***

Find the minimizer of the quadratic function

$$f|\mathbf{x}\rangle = \frac{1}{2}\langle \mathbf{x}|\mathbf{H}|\mathbf{x}\rangle - \langle \mathbf{x}|\mathbf{b}\rangle,$$

where

$$\mathbf{H} = \begin{pmatrix} 3 & 0 & \sqrt{3} \\ 0 & 4 & 2 \\ \sqrt{3} & 2 & 3 \end{pmatrix} \text{ and } \boldsymbol{b} = \begin{pmatrix} 2 \\ 0 \\ 1 \end{pmatrix},$$

by using the conjugate gradient algorithm, where $|\mathbf{x}_0\rangle = (0,0,0)^T$, and $\varepsilon = 0.0001$.

***Solution***

1. Find the gradient vector of the function $f|\mathbf{x}\rangle$:

$$\mathbf{g}|\mathbf{x}\rangle = \mathbf{H}|\mathbf{x}\rangle - \mathbf{b} = \begin{pmatrix} 3 & 0 & \sqrt{3} \\ 0 & 4 & 2 \\ \sqrt{3} & 2 & 3 \end{pmatrix}\begin{pmatrix} x_1 \\ x_2 \\ x_3 \end{pmatrix} - \begin{pmatrix} 2 \\ 0 \\ 1 \end{pmatrix} = \begin{pmatrix} 3x_1 + \sqrt{3}x_3 \\ 4x_2 + 2x_3 \\ \sqrt{3}x_1 + 2x_2 + 3x_3 \end{pmatrix} - \begin{pmatrix} 2 \\ 0 \\ 1 \end{pmatrix} = \begin{pmatrix} 3x_1 + \sqrt{3}x_3 - 2 \\ 4x_2 + 2x_3 \\ \sqrt{3}x_1 + 2x_2 + 3x_3 - 1 \end{pmatrix}.$$

2. Find $|\mathbf{g}_0\rangle$: $|\mathbf{g}_0\rangle = \mathbf{g}|\mathbf{x}_0\rangle = (-2,0,-1)^T$

3. Find $|\mathbf{d}_0\rangle$: $|\mathbf{d}_0\rangle = -|\mathbf{g}_0\rangle = (2,0,1)^T$

4. Find $\alpha_0$:

$$\alpha_0 = -\frac{\langle \mathbf{g_0}|\mathbf{d}_0\rangle}{\langle \mathbf{d}_0|\mathbf{H}|\mathbf{d}_0\rangle} = -\frac{(-2,0,-1)\begin{pmatrix}2\\0\\1\end{pmatrix}}{(2,0,1)\begin{pmatrix}3 & 0 & \sqrt{3}\\0 & 4 & 2\\\sqrt{3} & 2 & 3\end{pmatrix}\begin{pmatrix}2\\0\\1\end{pmatrix}} = -\frac{-5}{15+4\sqrt{3}} = 0.228017.$$

5. Find $|\mathbf{x}_1\rangle$:

$$|\mathbf{x}_1\rangle = |\mathbf{x}_0\rangle + \alpha_0|\mathbf{d}_0\rangle = \begin{pmatrix}0\\0\\0\end{pmatrix} + 0.228017\begin{pmatrix}2\\0\\1\end{pmatrix} = \begin{pmatrix}0.456034\\0\\0.228017\end{pmatrix}.$$

6. Find $|\mathbf{g_1}\rangle$**:** $|\mathbf{g_1}\rangle = \begin{pmatrix}-0.237053\\0.456034\\0.473872\end{pmatrix}$ where $\|\mathbf{g_1}\| = 0.699082 > \varepsilon$.

7. Find $\beta_0$:

$$\beta_0 = \frac{\langle \mathbf{g}_1|\mathbf{H}|\mathbf{d}_0\rangle}{\langle \mathbf{d}_0|\mathbf{H}|\mathbf{d}_0\rangle} = \frac{(-0.237053,0.456034,0.473872)\begin{pmatrix}3 & 0 & \sqrt{3}\\0 & 4 & 2\\\sqrt{3} & 2 & 3\end{pmatrix}\begin{pmatrix}2\\0\\1\end{pmatrix}}{(2,0,1)\begin{pmatrix}3 & 0 & \sqrt{3}\\0 & 4 & 2\\\sqrt{3} & 2 & 3\end{pmatrix}\begin{pmatrix}2\\0\\1\end{pmatrix}} = \frac{2.14232}{15+4\sqrt{3}} = 0.097697.$$

8. Find $\boldsymbol{d}^{(1)}$:

$$|\mathbf{d}_1\rangle = -|\mathbf{g}_1\rangle + \beta_0|\mathbf{d}_0\rangle = \begin{pmatrix}0.237053\\-0.456034\\-0.473872\end{pmatrix} + 0.097697\begin{pmatrix}2\\0\\1\end{pmatrix} = \begin{pmatrix}0.237053\\-0.456034\\-0.473872\end{pmatrix} + \begin{pmatrix}0.195394\\0\\0.097697\end{pmatrix} = \begin{pmatrix}0.432447\\-0.456034\\-0.376175\end{pmatrix}.$$

9. Find $\alpha_1$:

$$\alpha_1 = -\frac{\langle \mathbf{g}_1|\mathbf{d}_1\rangle}{\langle \mathbf{d}_1|\mathbf{H}|\mathbf{d}_1\rangle} = -\frac{(-0.2371,0.4560,0.4739)\begin{pmatrix}0.4325\\-0.4560\\-0.3762\end{pmatrix}}{(0.4325,-0.4560,-0.3762)\begin{pmatrix}3 & 0 & \sqrt{3}\\0 & 4 & 2\\\sqrt{3} & 2 & 3\end{pmatrix}\begin{pmatrix}0.4325\\-0.4560\\-0.3762\end{pmatrix}} = -\frac{-0.4887}{1.9401} = 0.2519.$$

10. Find $\boldsymbol{x}^{(2)}$:

$$|\mathbf{x}_2\rangle = |\mathbf{x}_1\rangle + \alpha_1|\mathbf{d}_1\rangle = \begin{pmatrix}0.4560\\0\\0.2280\end{pmatrix} + 0.2519\begin{pmatrix}0.4325\\-0.4560\\-0.3762\end{pmatrix} = \begin{pmatrix}0.4560\\0\\0.2280\end{pmatrix} + \begin{pmatrix}0.1089\\-0.1149\\-0.0948\end{pmatrix} = \begin{pmatrix}0.5649\\-0.1149\\0.1332\end{pmatrix}.$$

11. Find $|\mathbf{g}_2\rangle$:

$$|\mathbf{g}_2\rangle = \begin{pmatrix}-0.0743\\-0.1930\\0.1486\end{pmatrix}, \text{ with } \|\mathbf{g}_2\| = 0.2547 > \varepsilon.$$

12. Find $\beta_1$:

$$\beta_1 = \frac{\langle \mathbf{g}_2|\mathbf{H}|\mathbf{d}_1\rangle}{\langle \mathbf{d}_1|\mathbf{H}|\mathbf{d}_1\rangle} = \frac{(-0.0743,-0.1930,0.1486)\begin{pmatrix}3 & 0 & \sqrt{3}\\0 & 4 & 2\\\sqrt{3} & 2 & 3\end{pmatrix}\begin{pmatrix}0.4325\\-0.4560\\-0.3762\end{pmatrix}}{(0.4325,-0.4560,-0.3762)\begin{pmatrix}3 & 0 & \sqrt{3}\\0 & 4 & 2\\\sqrt{3} & 2 & 3\end{pmatrix}\begin{pmatrix}0.4325\\-0.4560\\-0.3762\end{pmatrix}} = \frac{0.2575}{1.9401} = 0.1327.$$

13. Find $\boldsymbol{d}^{(2)}$:

$$|\mathbf{d}_2\rangle = -|\mathbf{g}_2\rangle + \beta_1|\mathbf{d}_1\rangle = \begin{pmatrix}0.0743\\0.1930\\-0.1486\end{pmatrix} + 0.1327\begin{pmatrix}0.4325\\-0.4560\\-0.3762\end{pmatrix} = \begin{pmatrix}0.0743\\0.1930\\-0.1486\end{pmatrix} + \begin{pmatrix}0.0574\\-0.0605\\-0.0499\end{pmatrix} = \begin{pmatrix}0.1317\\0.1325\\-0.1985\end{pmatrix}.$$

14. Find $\alpha_2$:

$$\alpha_2 = -\frac{\langle \mathbf{g}_2|\mathbf{d}_2\rangle}{\langle \mathbf{d}_2|\mathbf{H}|\mathbf{d}_2\rangle} = -\frac{(-0.0743,-0.1930,0.1486)\begin{pmatrix}0.1317\\0.1325\\-0.1985\end{pmatrix}}{(0.1317,0.1325,-0.1985)\begin{pmatrix}3 & 0 & \sqrt{3}\\ 0 & 4 & 2\\ \sqrt{3} & 2 & 3\end{pmatrix}\begin{pmatrix}0.1317\\0.1325\\-0.1985\end{pmatrix}} = -\frac{-0.0648}{0.0447} = 1.4508.$$

15. Find $\boldsymbol{x}^{(3)}$:

$$|\mathbf{x}_3\rangle = |\mathbf{x}_2\rangle + \alpha_2|\mathbf{d}_2\rangle = \begin{pmatrix}0.5649\\-0.1149\\0.1333\end{pmatrix} + 1.4508\begin{pmatrix}0.1317\\0.1325\\-0.1985\end{pmatrix} = \begin{pmatrix}0.5649\\-0.1149\\0.1333\end{pmatrix} + \begin{pmatrix}0.1910\\0.1922\\-0.2880\end{pmatrix} = \begin{pmatrix}0.7560\\0.0774\\-0.1547\end{pmatrix}.$$

Note that

$$|\mathbf{g}_3\rangle = \begin{pmatrix}0.00004\\-0.00002\\0.00001\end{pmatrix}, \text{ with } \|\mathbf{g}_3\| = 0.000042 < \varepsilon.$$

***Example 9.4***

Consider the problem:

$$\text{Minimize } f|\mathbf{x}\rangle = 3x^2 + 12y^2,$$

using, $|\mathbf{x}_0\rangle = (5,5)^T$, and $\epsilon = 0.01$.

***Solution***

The 3D and contour plots of the function are shown in Figure 9.2. The plots show that the minimum lies at $|\mathbf{x}^*\rangle = (0.,0.)^T$, $f|\mathbf{x}^*\rangle = 0$.

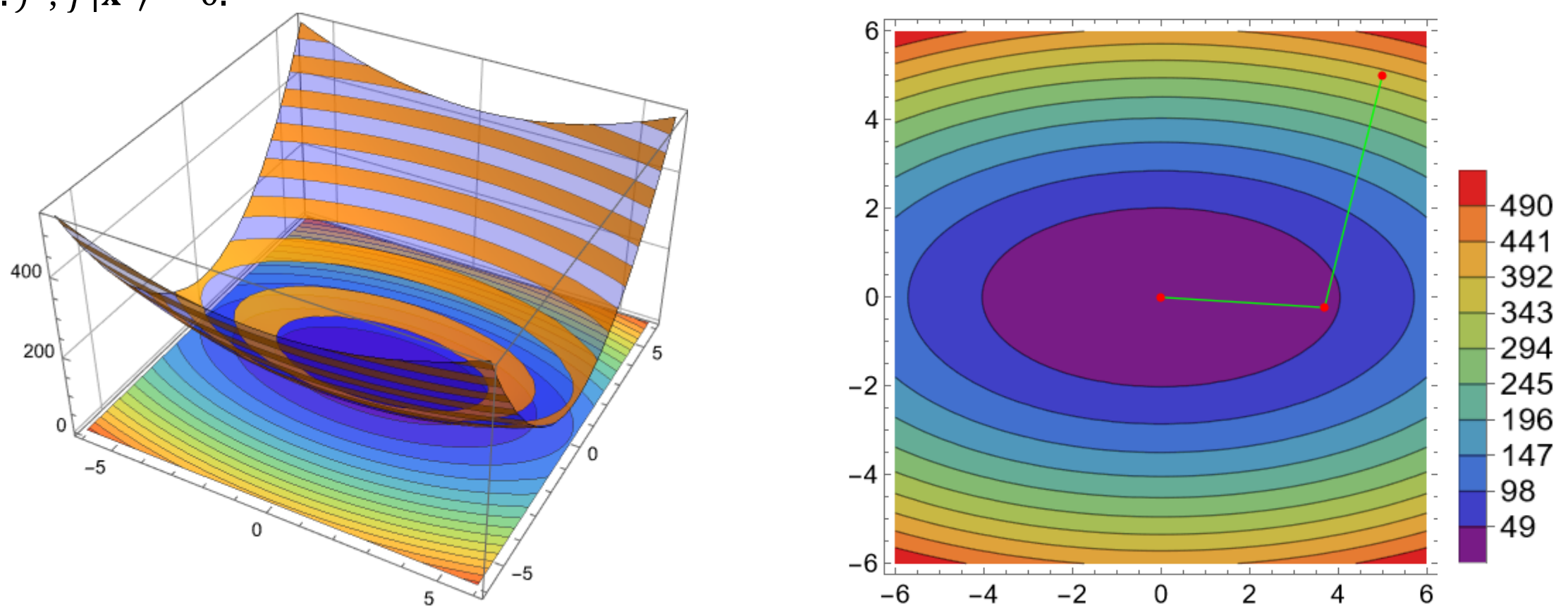

**Figure 9.2.** The results of 2 iterations of the conjugate gradient algorithm for $f|\mathbf{x}\rangle = 3x^2 + 12y^2$.

From Table 9.1, after 2 iterations, the optimal condition is attained. The results were produced by Mathematica code 9.1.

**Table 9.1.**

| No. of iters. | $\alpha x0$ | $\beta x0$ | $\alpha x1$ | $x0$ | $x1$ | $x2$ | $f(x0)$ | $f(x1)$ | $f(x2)$ |
|---|---|---|---|---|---|---|---|---|---|
| 1 | 0.044 | 0.034 | 0.159 | {5.,5.} | {3.692,-0.230} | {0.,0.} | 375. | 41.5385 | 0. |

***Mathematica Code 9.1*** Conjugate Gradient Algorithm (for Quadratic Function)

```
(* Conjugate Gradient Algorithm for Quadratic Functions Only *)

(*
Notations:
x0          :Intial vector
epsiloncga   :Small number to check the accuracy of the Conjugate Gradient Algorithm
f[x,y]       :Objective function
```

```
lii          :The last iteration index
result[k]    :The results of iteration k
*)

(* Taking Initial Inputs from User *)
x0=Input["Enter the intial point in the format {x, y}; for example {5,5} "] ;
epsiloncga=Input["Please enter accuracy of the conjugate gradient algorithm; for example
0.01 "];

domainx=Input["Please enter domain of x variable for 3D and contour plots; for example {-
6,6}"];
domainy=Input["Please enter domain of y variable for 3D and contour plots; for example {-
6,6}"];

If[
  epsiloncga<=0,
  Beep[];
  MessageDialog["The values of epsiloncga has to be postive number: "];
  Exit[];
  ];

(* Taking the Function from User *)
f[{x_,y_}] = Evaluate[Input["Please input a function of x and y to find the minimum "]];
(* For example: 3x^2+12y^2*)

gradfx[x_,y_]=Grad[f[{x,y}],{x,y}];
hessianfx=Grad[gradfx[x,y],{x,y}];

gradfx0=gradfx[x0[[1]],x0[[2]]];

If[
  N[Norm[gradfx0]]==0,
  Print["The conjugate gradient algorithm requires that the gradiant of the function f
!=0"];(*Ending program*)
  Exit[];
  ];

dx0=-gradfx0;

If[
  N[Transpose[dx0].(hessianfx.dx0)]==0,
  Print["The conjugate gradient algorithm requires that the Transpose[dx0].hessianfx.dx0 of
the function f !=0"];(*Ending program*)
  Exit[];
  ];

αx0=-((Transpose[gradfx0].dx0)/Transpose[dx0].(hessianfx.dx0));

x1=x0+αx0*dx0;
gradfx1=gradfx[x1[[1]],x1[[2]]];

If[
 Norm[gradfx1]<epsiloncga,

 (* Final Result *)
 lii=1;
 result=N[{αx0,Row[x0,","],Row[x1,","],f[x0],f[x1]}];
 plotresult=N[{x0,x1}];

 Print["The solution is x= ",  N[x1],"\nThe solution is (approximately)= ", N[f[x1]]];
```

```
 table=TableForm[
   {result},
   TableHeadings->{None,{"αx0","x0","x1","f[x0]","f[x1]"}}
   ],

 lii=2;
 βx0=(Transpose[gradfx1].hessianfx.dx0)/Transpose[dx0].(hessianfx.dx0);
 dx1=-gradfx1+βx0*dx0;

 If[
  N[Transpose[dx1].(hessianfx.dx1)]==0,
  Print["The conjugate gradient algorithm requires that the Transpose[dx1].hessianfx.dx1 of
the function f !=0"];(*Ending program*)
  Exit[];
  ];

 αx1=-((Transpose[gradfx1].dx1)/Transpose[dx1].(hessianfx.dx1));
 x2=x1+αx1*dx1;

 result=N[{αx0,βx0,αx1,Row[x0,","],Row[x1,","],Row[x2,","],f[x0],f[x1],f[x2]}];
 plotresult=N[{x0,x1,x2}];
 (* Final result *)
 Print["The solution is x= ",  N[x2],"\nThe solution is (approximately)= ", N[f[x2]]];
 table=TableForm[
   Table[result,{i,1,1}],
   TableHeadings->{None,{"αx0","βx0","αx1","x0","x1","x2","f[x0]","f[x1]","f[x2]"}}
   ]
 ]
Export["example91.xls",table,"XLS"];

(* Data Visualization *)
(* Domain of Varibles*)
xleft=domainx[[1]];
xright=domainx[[2]];
ydown=domainy[[1]];
yup=domainy[[2]];
(* 3D+ Contour Plot *)
plot1=Plot3D[
   f[{x,y}],
   {x,xleft,xright},
   {y,ydown,yup},
   ClippingStyle->None,
   MeshFunctions->{#3&},
   Mesh->15,
   MeshStyle->Opacity[.5],
   MeshShading->{{Opacity[.3],Blue},{Opacity[.8],Orange}},
   Lighting->"Neutral"
   ];
slice=SliceContourPlot3D[
   f[{x,y}],
   z==0,
   {x,xleft,xright},
   {y,ydown,yup},
   {z,-1,1},
   Contours->15,
   Axes->False,
   PlotPoints->50,
   PlotRangePadding->0,
   ColorFunction->"Rainbow"
   ];
Show[
```

```
 plot1,
 slice,
 PlotRange->All,
 BoxRatios->{1,1,.6},
 FaceGrids->{Back,Left}
 ]
(* Contour Plot with Step Iterations *)
ContourPlot[
 f[{x,y}],
 {x,xleft,xright},
 {y,ydown,yup},
 LabelStyle->Directive[Black,16],
 ColorFunction->"Rainbow",
 PlotLegends->Automatic,
 Contours->10,
 Epilog->{PointSize[0.015],Green,Line[plotresult],Red,Point[plotresult]}
 ]
(* Data Manipulation *)
Manipulate[
 ContourPlot[
  f[{x,y}],
  {x,xleft,xright},
  {y,ydown,yup},
  LabelStyle->Directive[Black,14],
  ColorFunction->"Rainbow",
  PlotLegends->Automatic,
  Contours->10,
  Epilog->{
    PointSize[0.015],
    Yellow,
    Arrow[{plotresult[[i]],plotresult[[i+1]]}],
    Red,
    Point[Table[plotresult[[j]],{j,1,i+1}]],
    Green,
    Line[Table[plotresult[[j]],{j,1,i+1}]]
    }
  ],
 {i,1,lii,1}
 ]
```

### 9.2.1. Convergence Analysis of Conjugate Gradient method

**Theorem 9.3 (Convergence of Conjugate Gradient Method):** (a) If $\mathbf{H}$ is a positive definite matrix, then for any initial point $|\mathbf{x}_0\rangle$ and an initial direction

$$|\mathbf{d}_0\rangle = -|\mathbf{g}_0\rangle = -|\mathbf{b} + \mathbf{H}\mathbf{x}_0\rangle, \tag{9.18}$$

the sequence generated by the recursive relation

$$|\mathbf{x}_{k+1}\rangle = |\mathbf{x}_k\rangle + \alpha_k|\mathbf{d}_k\rangle, \tag{9.19}$$

where

$$\alpha_k = -\frac{\langle \mathbf{g}_k|\mathbf{d}_k\rangle}{\langle \mathbf{d}_k|\mathbf{H}|\mathbf{d}_k\rangle}, \tag{9.20a}$$

$$|\mathbf{g}_k\rangle = |\mathbf{b} + \mathbf{H}\mathbf{x}_k\rangle, \tag{9.20b}$$

$$|\mathbf{d}_{k+1}\rangle = -|\mathbf{g}_{k+1}\rangle + \beta_k|\mathbf{d}_k\rangle, \tag{9.20c}$$

$$\beta_k = \frac{\langle \mathbf{g}_{k+1}|\mathbf{H}|\mathbf{d}_k\rangle}{\langle \mathbf{d}_k|\mathbf{H}|\mathbf{d}_k\rangle}, \tag{9.20d}$$

converges to the unique solution $|\mathbf{x}^*\rangle$ for the minimization of the quadratic function $f|\mathbf{x}\rangle = a + \langle \mathbf{b}|\mathbf{x}\rangle + \frac{1}{2}\langle \mathbf{x}|\mathbf{H}|\mathbf{x}\rangle$.
(b) The gradient $|\mathbf{g}_k\rangle$ is orthogonal to $\{|\mathbf{g}_0\rangle, |\mathbf{g}_1\rangle, \dots, |\mathbf{g}_{k-1}\rangle\}$, i.e.,

$$\langle \mathbf{g}_k|\mathbf{g}_i\rangle = 0, \quad 0 \le i < k. \tag{9.21}$$

**Proof:**

(a) The proof of convergence is the same as in Theorem 9.1. What remains to prove is that the directions $|\mathbf{d}_0\rangle, |\mathbf{d}_1\rangle, \dots, |\mathbf{d}_{n-1}\rangle$ form a conjugate set (i.e., are $\mathbf{H}$-conjugate), that is,

$$\langle \mathbf{d}_k|\mathbf{H}|\mathbf{d}_i\rangle = 0,\ \ 0 \le i < k \text{ and } 1 \le k \le n.$$

The proof is by induction. We assume that

$$\langle \mathbf{d}_k|\mathbf{H}|\mathbf{d}_i\rangle = 0 \ \text{ for } \ 0 \le i < k, \tag{9.22}$$

and show that

$$\langle \mathbf{d}_{k+1}|\mathbf{H}|\mathbf{d}_i\rangle = 0,\ \ \ 0 \le i < k+1.$$

Let $S(|\mathbf{v}_0\rangle, |\mathbf{v}_1\rangle, \dots, |\mathbf{v}_k\rangle)$ be the subspace spanned by vectors $|\mathbf{v}_0\rangle, |\mathbf{v}_1\rangle, \dots, |\mathbf{v}_k\rangle$. From

$$|\mathbf{g}_{k+1}\rangle = |\mathbf{g}_k\rangle + \alpha_k\mathbf{H}|\mathbf{d}_k\rangle,$$

and hence for $k = 0$, we have

$$|\mathbf{g}_1\rangle = |\mathbf{g}_0\rangle + \alpha_0\mathbf{H}|\mathbf{d}_0\rangle = |\mathbf{g}_0\rangle - \alpha_0\mathbf{H}|\mathbf{g}_0\rangle,$$

since $|\mathbf{d}_0\rangle = -|\mathbf{g}_0\rangle$. In addition, (9.20c) yields

$$|\mathbf{d}_1\rangle = -|\mathbf{g}_1\rangle + \beta_0|\mathbf{d}_0\rangle = -(1+\beta_0)|\mathbf{g}_0\rangle + \alpha_0\mathbf{H}|\mathbf{g}_0\rangle,$$

that is, $|\mathbf{g}_1\rangle$ and $|\mathbf{d}_1\rangle$ are linear combinations of $|\mathbf{g}_0\rangle$ and $\mathbf{H}|\mathbf{g}_0\rangle$, and so

$$S(|\mathbf{g}_0\rangle, |\mathbf{g}_1\rangle) = S(|\mathbf{d}_0\rangle, |\mathbf{d}_1\rangle) = S(|\mathbf{g}_0\rangle, \mathbf{H}|\mathbf{g}_0\rangle).$$

Similarly, for $k = 2$, we get

$$|\mathbf{g}_2\rangle = |\mathbf{g}_0\rangle - [\alpha_0 + \alpha_1(1+\beta_0)]\mathbf{H}|\mathbf{g}_0\rangle + \alpha_0\alpha_1\mathbf{H}^2|\mathbf{g}_0\rangle,$$
$$|\mathbf{d}_2\rangle = -[1 + (1+\beta_0)\beta_1]|\mathbf{g}_0\rangle + [\alpha_0 + \alpha_1(1+\beta_0) + \alpha_0\beta_1]\mathbf{H}|\mathbf{g}_0\rangle - \alpha_0\alpha_1\mathbf{H}^2|\mathbf{g}_0\rangle,$$

and hence

$$S(|\mathbf{g}_0\rangle, |\mathbf{g}_1\rangle, |\mathbf{g}_2\rangle) = S(|\mathbf{g}_0\rangle, \mathbf{H}|\mathbf{g}_0\rangle, \mathbf{H}^2|\mathbf{g}_0\rangle),$$
$$S(|\mathbf{d}_0\rangle, |\mathbf{d}_1\rangle, |\mathbf{d}_2\rangle) = S(|\mathbf{g}_0\rangle, \mathbf{H}|\mathbf{g}_0\rangle, \mathbf{H}^2|\mathbf{g}_0\rangle).$$

By continuing the induction, we can show that

$$S(|\mathbf{g}_0\rangle, |\mathbf{g}_1\rangle, \dots, |\mathbf{g}_k\rangle) = S(|\mathbf{g}_0\rangle, \mathbf{H}|\mathbf{g}_0\rangle, \dots, \mathbf{H}^k|\mathbf{g}_0\rangle), \tag{9.23.a}$$
$$S(|\mathbf{d}_0\rangle, |\mathbf{d}_1\rangle, \dots, |\mathbf{d}_k\rangle) = S(|\mathbf{g}_0\rangle, \mathbf{H}|\mathbf{g}_0\rangle, \dots, \mathbf{H}^k|\mathbf{g}_0\rangle). \tag{9.23.b}$$

Now from (9.20c)

$$\langle \mathbf{d}_{k+1}|\mathbf{H}|\mathbf{d}_i\rangle = -\langle \mathbf{g}_{k+1}|\mathbf{H}|\mathbf{d}_i\rangle + \beta_k\langle \mathbf{d}_k|\mathbf{H}|\mathbf{d}_i\rangle. \tag{9.24}$$

For $i = k$, the use of (9.20d) gives

$$\langle \mathbf{d}_{k+1}|\mathbf{H}|\mathbf{d}_k\rangle = -\langle \mathbf{g}_{k+1}|\mathbf{H}|\mathbf{d}_k\rangle + \beta_k\langle \mathbf{d}_k|\mathbf{H}|\mathbf{d}_k\rangle = 0. \tag{9.25}$$

For $i < k$, (9.23.b) shows that

$$\mathbf{H}|\mathbf{d}_i\rangle \in S(|\mathbf{d}_0\rangle, |\mathbf{d}_1\rangle, \dots, |\mathbf{d}_k\rangle),$$

and thus $\mathbf{H}|\mathbf{d}_i\rangle$ can be represented by the linear combination

$$\mathbf{H}|\mathbf{d}_i\rangle = \sum_{i=0}^{k} a_i|\mathbf{d}_i\rangle, \tag{9.26}$$

where $a_i$ for $i = 0,1,\dots,k$ are constants. Now from (9.24) and (9.26)

$$\langle \mathbf{d}_{k+1}|\mathbf{H}|\mathbf{d}_i\rangle = -\sum_{i=0}^{k} a_i\langle \mathbf{g}_{k+1}|\mathbf{d}_i\rangle + \beta_k\langle \mathbf{d}_k|\mathbf{H}|\mathbf{d}_i\rangle = 0,\ \ \ i < k. \tag{9.27}$$

The first term is zero from the orthogonality property of Theorem 9.2, whereas the second term is zero from the assumption in (9.22). By combining (9.25) and (9.27), we have

$$\langle \mathbf{d}_{k+1}|\mathbf{H}|\mathbf{d}_i\rangle = 0,\ \ \ 0 \le i < k+1. \tag{9.28}$$

For $k = 0$, (9.28) gives

$$\langle \mathbf{d}_1|\mathbf{H}|\mathbf{d}_i\rangle = 0,\ \ 0 \le i < 1,$$

and, therefore, from (9.22) and (9.28), we have

$$\langle \mathbf{d}_2|\mathbf{H}|\mathbf{d}_i\rangle = 0,\ \ 0 \le i < 2,$$
$$\langle \mathbf{d}_3|\mathbf{H}|\mathbf{d}_i\rangle = 0,\ \ 0 \le i < 3,$$
$$\dots\ \dots\ \dots$$
$$\langle \mathbf{d}_k|\mathbf{H}|\mathbf{d}_i\rangle = 0,\ \ 0 \le i < k.$$

(b) From (9.23a) and (9.23b), $|\mathbf{g}_0\rangle, |\mathbf{g}_1\rangle, \ldots, |\mathbf{g}_k\rangle$ span the same subspace as $|\mathbf{d}_0\rangle, |\mathbf{d}_1\rangle, \ldots, |\mathbf{d}_k\rangle$ and, consequently, they are linearly independent. We can write

$$|\mathbf{g}_i\rangle = \sum_{j=0}^{i} a_j |\mathbf{d}_j\rangle,$$

where $a_j$ for $j = 0,1, \ldots, i$ are constants. Therefore, from Theorem 9.2

$$\langle \mathbf{g}_k | \mathbf{g}_i \rangle = \sum_{j=0}^{i} a_j \langle \mathbf{g}_k | \mathbf{d}_j \rangle = 0, \quad 0 \le i < k.$$

■

From the proof of Theorem 9.3, we have the following theorem.

**Theorem 9.4:** In the conjugate gradient algorithm, the directions $|\mathbf{d}_0\rangle, |\mathbf{d}_1\rangle, \ldots, |\mathbf{d}_{n-1}\rangle$ are **H**-conjugate.

**Remarks:**

1- The expression for $\alpha_k$ in Theorem 9.3 can be simplified as follows. From (9.20.c)

$$-\langle \mathbf{g}_k | \mathbf{d}_k \rangle = \langle \mathbf{g}_k | \mathbf{g}_k \rangle - \beta_{k-1} \langle \mathbf{g}_k | \mathbf{d}_{k-1} \rangle, \tag{9.29}$$

where

$$\langle \mathbf{g}_k | \mathbf{d}_{k-1} \rangle = 0, \tag{9.30}$$

according to Theorem 9.2. Hence

$$-\langle \mathbf{g}_k | \mathbf{d}_k \rangle = \langle \mathbf{g}_k | \mathbf{g}_k \rangle, \tag{9.31}$$

and, therefore, the expression for $\alpha_k$ in (9.20.a) is modified as

$$\alpha_k = \frac{\langle \mathbf{g}_k | \mathbf{g}_k \rangle}{\langle \mathbf{d}_k | \mathbf{H} | \mathbf{d}_k \rangle}. \tag{9.32}$$

2- The expression for $\beta_k$ in Theorem 9.3 can also be simplified as follows. We have,

$$\mathbf{H}|\mathbf{d}_k\rangle = \frac{1}{\alpha_k} |\mathbf{g}_{k+1} - \mathbf{g}_k\rangle, \tag{9.33}$$

and so

$$\langle \mathbf{g}_{k+1} | \mathbf{H} | \mathbf{d}_k \rangle = \frac{1}{\alpha_k} (\langle \mathbf{g}_{k+1} | \mathbf{g}_{k+1} \rangle - \langle \mathbf{g}_{k+1} | \mathbf{g}_k \rangle). \tag{9.34}$$

Now from (9.23.a) and (9.23.b)

$$|\mathbf{g}_k\rangle \in S(|\mathbf{d}_1\rangle, |\mathbf{d}_2\rangle, \ldots, |\mathbf{d}_k\rangle), \tag{9.35}$$

or

$$|\mathbf{g}_k\rangle = \sum_{i=0}^{k} a_i |\mathbf{d}_i\rangle, \tag{9.36}$$

and as a result

$$\langle \mathbf{g}_{k+1} | \mathbf{g}_k \rangle = \sum_{i=0}^{k} a_i \langle \mathbf{g}_{k+1} | \mathbf{d}_i \rangle = 0, \tag{9.37}$$

by virtue of Theorem 9.2. Therefore,

$$\beta_k = \frac{\langle \mathbf{g}_{k+1} | \mathbf{H} | \mathbf{d}_k \rangle}{\langle \mathbf{d}_k | \mathbf{H} | \mathbf{d}_k \rangle} = \frac{1}{\alpha_k} \langle \mathbf{g}_{k+1} | \mathbf{g}_{k+1} \rangle \frac{\alpha_k}{\langle \mathbf{g}_k | \mathbf{g}_k \rangle} = \frac{\langle \mathbf{g}_{k+1} | \mathbf{g}_{k+1} \rangle}{\langle \mathbf{g}_k | \mathbf{g}_k \rangle}. \tag{9.38}$$

3- Hence, the iterative scheme of the conjugate gradient method is

$$|\mathbf{x}_{k+1}\rangle = |\mathbf{x}_k\rangle + \alpha_k |\mathbf{d}_k\rangle, \tag{9.39.a}$$

$$\alpha_k = \frac{\langle \mathbf{g}_k | \mathbf{g}_k \rangle}{\langle \mathbf{d}_k | \mathbf{H} | \mathbf{d}_k \rangle}, \tag{9.39.b}$$

$$|\mathbf{d}_{k+1}\rangle = -|\mathbf{g}_{k+1}\rangle + \beta_k |\mathbf{d}_k\rangle, \tag{9.39.c}$$

$$\beta_k = \frac{\langle \mathbf{g}_{k+1} | \mathbf{g}_{k+1} \rangle}{\langle \mathbf{g}_k | \mathbf{g}_k \rangle}, \tag{9.39.d}$$

## 9.3 Method of Conjugate for General Functions

The conjugate gradient algorithm can be extended to general nonlinear functions [7] by interpreting $f|\mathbf{x}\rangle = \frac{1}{2}\langle\mathbf{x}|\mathbf{H}|\mathbf{x}\rangle - \langle\mathbf{x}|\mathbf{b}\rangle$ as a second-order Taylor series approximation of the objective function. Near the solution, such functions behave approximately as quadratics. For a quadratic, the matrix $\mathbf{H}$, the Hessian of the quadratic, is constant. However, for a general nonlinear function, the Hessian is a matrix that has to be reevaluated at each iteration of the algorithm. This can be computationally very expensive. Thus, an efficient implementation of the conjugate gradient algorithm that eliminates the Hessian evaluation at each step is desirable.

Observe that $\mathbf{H}$ appears only in the computation of the scalars $\alpha_k$ and $\beta_k$. Because $\alpha_k = \min_{\alpha\geq 0} f(|\mathbf{x}_k\rangle + \alpha|\mathbf{d}_k\rangle)$ the closed-form formula for $\alpha_k$ in the algorithm can be replaced by a numerical line search procedure. Therefore, we only need to concern ourselves with the formula for $\beta_k$. Here, we will discuss three well-known formulas for $\beta_k$.

### 9.3.1. The Hestenes-Stiefel Formula

Recall that $\beta_k = \frac{\langle\mathbf{g}_{k+1}|\mathbf{H}|\mathbf{d}_k\rangle}{\langle\mathbf{d}_k|\mathbf{H}|\mathbf{d}_k\rangle}$. We have, $|\mathbf{x}_{k+1}\rangle = |\mathbf{x}_k\rangle + \alpha_k|\mathbf{d}_k\rangle$. Then, multiply both sides by $\mathbf{H}$ yields

$$\mathbf{H}|\mathbf{x}_{k+1}\rangle = \mathbf{H}|\mathbf{x}_k\rangle + \alpha_k\mathbf{H}|\mathbf{d}_k\rangle, \tag{9.40}$$

which equivalent to

$$\mathbf{H}|\mathbf{x}_{k+1}\rangle - |\mathbf{b}\rangle = \mathbf{H}|\mathbf{x}_k\rangle - |\mathbf{b}\rangle + \alpha_k\mathbf{H}|\mathbf{d}_k\rangle. \tag{9.41}$$

Since $|\mathbf{g}_k\rangle = \mathbf{H}|\mathbf{x}_k\rangle - |\mathbf{b}\rangle$ and $|\mathbf{g}_{k+1}\rangle = \mathbf{H}|\mathbf{x}_{k+1}\rangle - |\mathbf{b}\rangle$, we get $|\mathbf{g}_{k+1}\rangle = |\mathbf{g}_k\rangle + \alpha_k\mathbf{H}|\mathbf{d}_k\rangle$. Then $\mathbf{H}|\mathbf{d}_k\rangle = \frac{1}{\alpha_k}|\mathbf{g}_{k+1} - \mathbf{g}_k\rangle$. Substituting this into the original equation for $\beta_k$ gives

$$\beta_k = \frac{\langle\mathbf{g}_{k+1}|\mathbf{g}_{k+1} - \mathbf{g}_k\rangle}{\langle\mathbf{d}_k|\mathbf{g}_{k+1} - \mathbf{g}_k\rangle}, \tag{9.42}$$

which is called the Hestenes-Stiefel formula [7].

### 9.3.2. The Polak-Ribiere Formula

Starting from Hestenes-Stiefel formula

$$\beta_k = \frac{\langle\mathbf{g}_{k+1}|\mathbf{g}_{k+1} - \mathbf{g}_k\rangle}{\langle\mathbf{d}_k|\mathbf{g}_{k+1}\rangle - \langle\mathbf{d}_k|\mathbf{g}_k\rangle}. \tag{9.43}$$

Since $\langle\mathbf{d}_k|\mathbf{g}_{k+1}\rangle = 0$ and $|\mathbf{d}_k\rangle = -|\mathbf{g}_k\rangle + \beta_{k-1}|\mathbf{d}_{k-1}\rangle$. Then, multiplying both sides by $\langle\mathbf{g}_k|$ implies

$$\begin{aligned}\langle\mathbf{g}_k|\mathbf{d}_k\rangle &= -\langle\mathbf{g}_k|\mathbf{g}_k\rangle + \beta_{k-1}\langle\mathbf{g}_k|\mathbf{d}_{k-1}\rangle,\\ &= -\langle\mathbf{g}_k|\mathbf{g}_k\rangle,\end{aligned} \tag{9.44}$$

and thus, the expression for $\beta_k$ becomes

$$\beta_k = \frac{\langle\mathbf{g}_{k+1}|\mathbf{g}_{k+1} - \mathbf{g}_k\rangle}{\langle\mathbf{g}_k|\mathbf{g}_k\rangle}, \tag{9.45}$$

which is known as the Polak-Ribiere formula.

### 9.3.3. The Fletcher-Reeves Formula

Starting with Polak-Ribiere formula

$$\beta_k = \frac{\langle\mathbf{g}_{k+1}|\mathbf{g}_{k+1}\rangle - \langle\mathbf{g}_{k+1}|\mathbf{g}_k\rangle}{\langle\mathbf{g}_k|\mathbf{g}_k\rangle}. \tag{9.46}$$

Since $\langle\mathbf{g}_{k+1}|\mathbf{g}_k\rangle = 0$, so that

$$\beta_k = \frac{\langle\mathbf{g}_{k+1}|\mathbf{g}_{k+1}\rangle}{\langle\mathbf{g}_k|\mathbf{g}_k\rangle}, \tag{9.47}$$

which is called the Fletcher-Reeves formula.

**Algorithm** (Hestenes-Stiefel, Polak-Ribiere, or Fletcher-Reeves Algorithm for General Nonlinear Functions)

**Step 1:** Choose a starting point $|\mathbf{x}_0\rangle \in \mathbb{R}^n$ and a tolerance $0 < \varepsilon < 1$. Set $k = 0$.

**Step 2:** Compute the gradient $|\mathbf{g}_0\rangle$ of $f|\mathbf{x}\rangle$ at $|\mathbf{x}_0\rangle$.

**Step 3:** Set $|\mathbf{d}_0\rangle = -|\mathbf{g}_0\rangle$.

**Step 4:** Compute $\alpha_k$ as: $\alpha_k$ by using line search $\alpha_k = \min_{\alpha\geq 0} f|\mathbf{x}_k + \alpha\mathbf{d}_k\rangle$.

**Step 5:** Compute $|\mathbf{x}_{k+1}\rangle$ as: $|\mathbf{x}_{k+1}\rangle = |\mathbf{x}_k\rangle + \alpha_k|\mathbf{d}_k\rangle$.

**Step 6:** Compute $|\mathbf{g}_{k+1}\rangle$. If $\|\mathbf{g}_{k+1}\| < \varepsilon$ stop, else go to the next step.

**Step 7:** Compute $\beta_k$, (9.42), (9.45), or (9.47).

**Step 8:** Compute $|\mathbf{d}_{k+1}\rangle$ as: $|\mathbf{d}_{k+1}\rangle = -|\mathbf{g}_{k+1}\rangle + \beta_k|\mathbf{d}_k\rangle$.

**Step 9:** Set $k = k + 1$ and go to step 4.

***Example 9.5***

Consider the problem:

$$\text{Minimize } f|\mathbf{x}\rangle = (x^2 - x + 2)^2 + (y^2 - y + 1)^2,$$

using, $|\mathbf{x}_0\rangle = (5,5)^T$, and $\epsilon = 0.01$.

***Solution***

The 3D and contour plots of the function are shown in Figure 9.3. The plots show that the minimum lies at $|\mathbf{x}^*\rangle = (0.500599, 0.500632)^T$, $f|\mathbf{x}^*\rangle = 3.625$.

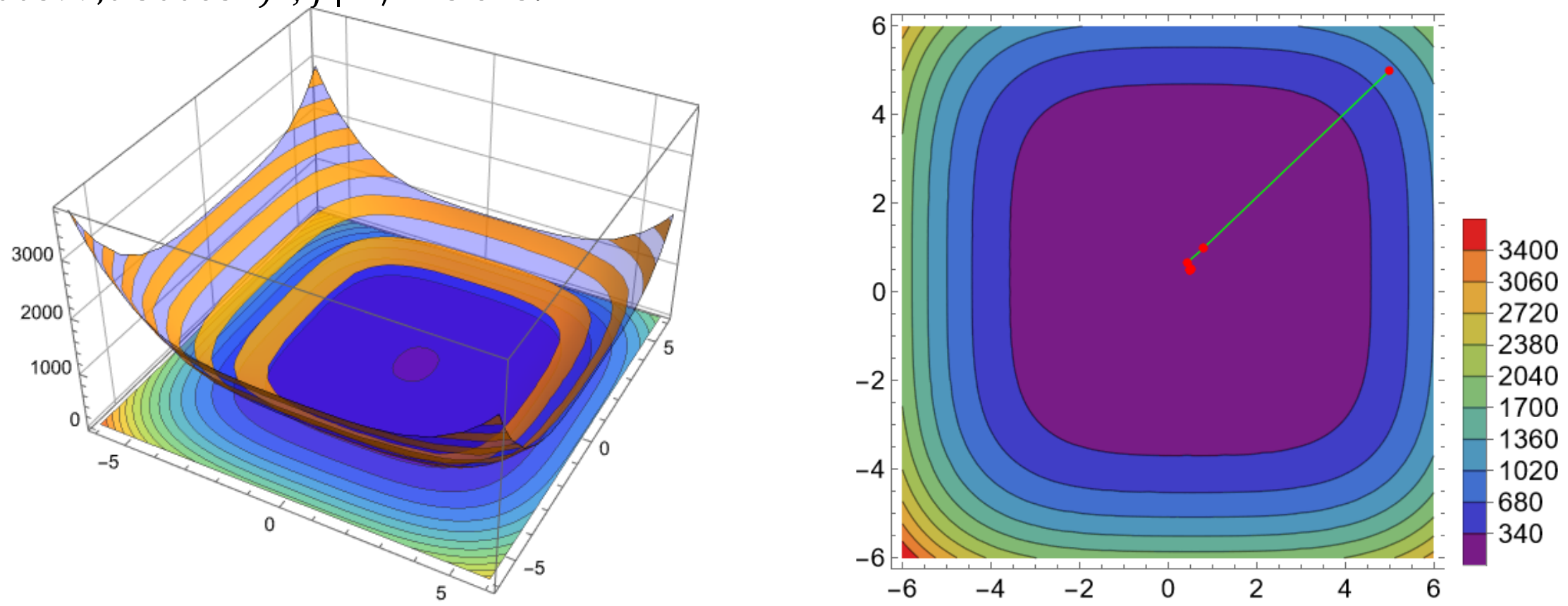

**Figure 9.3.** The results of 5 iterations of the conjugate gradient algorithm for $f|\mathbf{x}\rangle = (x^2 - x + 2)^2 + (y^2 - y + 1)^2$.

From Table 9.2, after 5 iterations, the optimal condition is attained. The results were produced by Mathematica code 9.2.

**Table 9.2.**

| No. of iters. | $\alpha star$ | $x0$ | $x1$ | $f(x0)$ | $f(x1)$ | error |
|---|---|---|---|---|---|---|
| 1 | 0.010613 | {5., 5.} | {0.797, 0.988} | 925 | 4.356 | 2.916 |
| 2 | 0.162708 | {0.797, 0.988} | {0.439, 0.672} | 4.356 | 3.683 | 0.683 |
| 3 | 0.241639 | {0.439, 0.672} | {0.512, 0.516} | 3.683 | 3.625 | 0.101 |
| 4 | 0.177781 | {0.512, 0.516} | {0.498, 0.505} | 3.625 | 3.625 | 0.020 |
| 5 | 0.250952 | {0.498, 0.505} | {0.500, 0.500} | 3.625 | 3.625 | 0.004 |

***Mathematica Code 9.2*** Conjugate Gradient Algorithm (Fletcher-Reeves Formula)

```
(* Conjugate Gradient Algorithm for General Functions Fletcher-Reeves Formula *)
```

```
(*
Notations:
x0       :Intial vector
epsilonCGM :Small number to check the accuracy of the conjugate gradient method for general
functions (Fletcher-Reeves Formula)
f[x,y]    :Objective function
lii       :The last iteration index
result[k] :The results of iteration k
*)

(* Taking Initial Inputs from User *)
x0=Input["Enter the intial point in the format {x, y}; for example {5,5} "] ;
epsilonCGM=Input["Please enter accuracy of the conjugate gradient method (Fletcher-Reeves
Formula); for example 0.01 "];

If[
  epsilonCGM<=0,
  Beep[];
  MessageDialog["The value of epsilonCGM has to be postive number: "];
  Exit[];
  ];

domainx=Input["Please enter domain of x variable for 3D and contour plots; for example {-
6,6}"];
domainy=Input["Please enter domain of y variable for 3D and contour plots; for example {-
6,6}"];

(* Taking the Function from User *)
f[{x_,y_}] = Evaluate[Input["Please input a function of x and y to find the minimum "]];
(* For example: (x^2-x+2)^2+(y^2-y+1)^2 *)

(*Defination of the Unidirectionalsearch Function*)
(*This function start by bracketing the minimum (using Bounding Phase Method) then isolating
the minimum (using Golden Section Search Method) *)

unidirectionalsearch[α_,delt_,eps_]:=Module[


{α0=α,delta=delt,epsilon=eps,y1,y2,y3,αα,a,b,increment,a0,b0,anew,bnew,anorm,bnorm,lnorm,α1n
orm,α2norm,α1,α2,φ1,φ2,αstar},

   (* Bounding Phase Method *)
   (* Initiating Required Variables *)
   y1 =φ[α0-Abs[delta]];
   y2 =φ[α0];
   y3 =φ[α0+Abs[delta]];

   (*Determining Whether the Inicrement Is Positive or Negative*)
   Which[
    y1==y2,
    a=α0-Abs[delta];
    b=α0;
    Goto[end];,
    y2==y3,
    a=α0;
    b=α0+Abs[delta];
    Goto[end];,
    y1==y3||(y1>y2&&y2<y3),
    a=α0-Abs[delta];
    b=α0+Abs[delta];
    Goto[end];
```

```
    ];

   Which[
     y1>y2&&y2>y3,
     increment=Abs[delta];,
     y1<y2&&y2<y3,
     increment=-Abs[delta];
     ]

    (* Starting the Algorithm *)
    Do[
     αα[0]=α0;
     αα[k+1]=αα[k]+2^k*increment;

     Which[
      φ[αα[k]]<φ[αα[k+1]],(* Evidently, it is impossible the condition to hold for k=0 *)
      a=αα[k-1];
      b= αα[k+1];
      Break [],

      k>50,
      Print["After 50 iterations the bounding phase method can not braketing the min of
alpha"];
      Exit[]
      ];,
     {k,0,∞}
     ];

   Label[end];

   If[
    a>b,
    {a,b}={b,a}
    ];

   (* Golden Section Search Method *)
   (* Initiating Required Variables*)
   a0=a;
   b0=b;
   anew=a;
   bnew=b;

   If[
    a0==b0,
    αstar=a;
    Goto[final]
    ];

   (* Starting the Algorithm *)
   Do[
    (* Normalize the Variable α *)
    anorm=(anew-a)/(b-a);
    bnorm=(bnew-a)/(b-a);

    lnorm=bnorm-anorm;

    α1norm=anorm+0.382*lnorm;
    α2norm=bnorm-0.382*lnorm;

    α1=α1norm(b0-a0)+a0;
    α2=α2norm(b0-a0)+a0;
```

```
    φ1=φ[α1];
    φ2=φ[α2];

    Which [
     φ1>φ2,
     anew=α1(*move lower bound to α1*);,
     φ1<φ2,
     bnew=α2(*move upper bound to α2*);,
     φ1==φ2,
     anew=α1(*move lower bound to α1*);
     bnew=α2(*move upper bound to α2*);
     ];

    αstar=0.5*(anew+bnew);

    If[
     Abs[lnorm]<epsilon,
     Break[]
     ];,
    {k,1,∞}
    ];

   Label[final];

   (* Final Result *)
   N[αstar]
   ];

α00=2;(* The intial point of α; for example 2*)
delt0=1;(* The parameter delta of Bounding Phase Method; for example 1 *)
eps0=0.01;(* The accuracy of the Golden Section Search Method; for example 0.01 *)

gradfx[x_,y_]=Grad[f[{x,y}],{x,y}];

gradfx0=gradfx[x0[[1]],x0[[2]]];

If[
  N[Norm[gradfx0]]==0,
  Print["The conjugate gradient method requires that the gradiant of the function f
!=0"];(*Ending program*)
  Exit[];
  ];

dx0=-gradfx0;

(* Main Loop *)
Do[
  φ[α_]=f[x0+α*dx0];
  αstar=unidirectionalsearch[α00,delt0,eps0];

  x1=x0+αstar*dx0;
  gradfx1=gradfx[x1[[1]],x1[[2]]];
  error=Norm[gradfx1];

  lii=k;
  result[k]=N[{k,αstar,Row[x0,","],Row[x1,","],f[x0],f[x1],error}];
  plotresult[k]=N[{x0,x1}];

  If[
   error<epsilonCGM||k>50,
```

```
   Break[],
   βx0=(Transpose[gradfx1].gradfx1)/(Transpose[gradfx0].gradfx0);
   dx1=-gradfx1+βx0*dx0;
   x0=x1;
   dx0=dx1;
   gradfx0=gradfx1;
   ],
  {k,1,∞}
  ];

(* Final Result *)
Print["The solution is x= ",  N[x1],"\nThe solution is (approximately)= ", N[f[x1]]]

(* Results of Each Iteration *)
table=TableForm[
  Table[
   result[i],
   {i,1,lii}
   ],
  TableHeadings->{None,{"No. of iters.","αstar","x0","x1","f[x0]","f[x1]","error"}}
  ]

Export["example92.xls",table,"XLS"];

(* Data Visualization *)
(* Domain of Varibles*)
xleft=domainx[[1]];
xright=domainx[[2]];
ydown=domainy[[1]];
yup=domainy[[2]];

(* 3D+ Contour Plot *)
plot1=Plot3D[
   f[{x,y}],
   {x,xleft,xright},
   {y,ydown,yup},
   ClippingStyle->None,
   MeshFunctions->{#3&},
   Mesh->15,
   MeshStyle->Opacity[.5],
   MeshShading->{{Opacity[.3],Blue},{Opacity[.8],Orange}},
   Lighting->"Neutral"
   ];

slice=SliceContourPlot3D[
   f[{x,y}],
   z==0,
   {x,xleft,xright},
   {y,ydown,yup},
   {z,-1,1},
   Contours->15,
   Axes->False,
   PlotPoints->50,
   PlotRangePadding->0,
   ColorFunction->"Rainbow"
   ];
Show[
 plot1,
 slice,
 PlotRange->All,
 BoxRatios->{1,1,.6},
```

```
 FaceGrids->{Back,Left}
 ]

(* Contour Plot with Step Iterations *)
ContourPlot[
 f[{x,y}],
 {x,xleft,xright},
 {y,ydown,yup},
 LabelStyle->Directive[Black,16],
 ColorFunction->"Rainbow",
 PlotLegends->Automatic,
 Contours->10,
 Epilog-
>{PointSize[0.015],Green,Line[Flatten[Table[plotresult[i],{i,1,lii}],1]],Red,Point[Flatten[T
able[plotresult[i],{i,1,lii}],1]]}
 ]

(* Data Manipulation *)
Manipulate[
 ContourPlot[
  f[{x,y}],
  {x,xleft,xright},
  {y,ydown,yup},
  LabelStyle->Directive[Black,14],
  ColorFunction->"Rainbow",
  PlotLegends->Automatic,
  Contours->10,
  Epilog->{
    PointSize[0.015],
    Yellow,
    Arrow[{plotresult[i][[1]],plotresult[i][[2]]}],
    Red,
    Point[Flatten[Table[plotresult[j],{j,1,i}],1]],
    Green,
    Line[Flatten[Table[plotresult[j],{j,1,i}],1]]
    }
  ],
 {i,1,lii,1}
 ]
```

## 9.4 The Basic Quasi-Newton Approach

In Section 9.1, multidimensional optimization methods were considered, in which the search for the minimizer is carried out by using a set of conjugate directions. An important feature of some of these methods is that explicit expressions for the second derivatives of the objective function $f|\mathbf{x}\rangle$ are not required. Another class of methods that do not require explicit expressions for the second derivatives is the class of quasi-Newton methods [1-5]. Quasi-Newton methods rank among the most efficient methods available and are used very extensively in numerous applications. Quasi-Newton methods, like most other methods, are developed for the convex quadratic problem and are then extended to the general problem. The foundation of these methods is the classical Newton method described in Section 8.2. In the Newton method, we require to compute the inverse of the Hessian at every iteration, which is a very expensive computation. Furthermore, Newton method fails to converge to the minimizer $|\mathbf{x}^*\rangle$ for the following cases.

1. The search direction $|\mathbf{d}_k\rangle$ is orthogonal to the gradient vector $|\mathbf{g}_k\rangle$.
2. The inverse of the Hessian matrix exists, and it is not positive definite.
3. The search direction $|\mathbf{d}_k\rangle$ is so large that $f|\mathbf{x}_{k+1}\rangle > f|\mathbf{x}_k\rangle$.

4. The inverse of the Hessian matrix does not exist.

These drawbacks of Newton method gave the motivation to develop the quasi-Newton methods. The basic principle in quasi-Newton methods is that the direction of search is based on an $n \times n$ direction matrix $\mathbf{S}$ which serves the same purpose as the inverse Hessian in the Newton method.

Recall that the idea behind Newton method is to locally approximate the function $f$ being minimized, at every iteration, by a quadratic function. The basic recursive formula representing Newton method is

$$|\mathbf{x}_{k+1}\rangle = |\mathbf{x}_k\rangle - \mathbf{H}_k^{-1}|\mathbf{g}_k\rangle. \tag{9.48}$$

However, if the initial point is not sufficiently close to the solution, then the algorithm may not possess the descent property (i.e., $f|\mathbf{x}_{k+1}\rangle \not< f|\mathbf{x}_k\rangle$ for some $k$). To avoid this and guarantee that the algorithm has descent property, we modify the original algorithm as follows:

$$|\mathbf{x}_{k+1}\rangle = |\mathbf{x}_k\rangle - \alpha_k \mathbf{H}_k^{-1}|\mathbf{g}_k\rangle, \tag{9.49}$$

where $\alpha_k$ is chosen to ensure that $f|\mathbf{x}_{k+1}\rangle < f|\mathbf{x}_k\rangle$. For example, we may choose $\alpha_k = \arg\min_{\alpha\geq 0} f|\mathbf{x}_k - \alpha\mathbf{H}_k^{-1}\mathbf{g}_k\rangle$. We can then determine an appropriate value of $\alpha_k$ by performing a line search in the direction $-\mathbf{H}_k^{-1}|\mathbf{g}_k\rangle$. Note that although the line search is simply the minimization of the real variable function $\phi_k(\alpha) = f(|\mathbf{x}_k\rangle - \alpha\mathbf{H}_k^{-1}|\mathbf{g}_k\rangle)$, it is not a trivial problem to solve because we need to evaluate $\mathbf{H}_k$ and solve the equation $\mathbf{H}_k|\mathbf{d}_k\rangle = -|\mathbf{g}_k\rangle$.

***Example 9.6***

Let $f: \mathbb{R}^2 \to \mathbb{R}$ be a function defined by

$$f|\mathbf{x}\rangle = x_1 - x_2 + 2x_1^2 + 2x_1x_2 + x_2^2, \qquad |\mathbf{x}\rangle = (x_1, x_2)^T.$$

Use Newton method with line search to find the minimizer of the objective function $f$. Take $|\mathbf{x}_0\rangle = (0,0)^T$ and $\varepsilon = 0.0001$.

***Solution***

First, we find the gradient vector $|\mathbf{g}\rangle$ as:

$$|\mathbf{g}\rangle = \nabla f|\mathbf{x}\rangle = \left(\frac{\partial f}{\partial x_1}, \frac{\partial f}{\partial x_2}\right)^T = (1 + 4x_1 + 2x_2, -1 + 2x_1 + 2x_2)^T.$$

Now, find the Hessian matrix $\mathbf{H}$ as follows

$$\mathbf{H} = \mathbf{H}_f|\mathbf{x}\rangle = \nabla^2 f|\mathbf{x}\rangle = \begin{pmatrix} \frac{\partial^2 f}{\partial x_1^2} & \frac{\partial^2 f}{\partial x_1 \partial x_2} \\ \frac{\partial^2 f}{\partial x_2 \partial x_1} & \frac{\partial^2 f}{\partial x_2^2} \end{pmatrix} = \begin{pmatrix} 4 & 2 \\ 2 & 2 \end{pmatrix}.$$

**Iteration 1.**

$$|\mathbf{g}_0\rangle = (1, -1)^T,$$
$$\mathbf{H}_0^{-1} = \begin{pmatrix} 1/2 & -1/2 \\ -1/2 & 1 \end{pmatrix},$$
$$|\mathbf{d}_0\rangle = -\mathbf{H}_0^{-1}|\mathbf{g}_0\rangle = (-1, 3/2)^T.$$

Hence

$$|\mathbf{x}_1\rangle = |\mathbf{x}_0\rangle + \alpha_0|\mathbf{d}_0\rangle = (-\alpha_0, 3\alpha_0/2)^T,$$
$$f|\mathbf{x}_1\rangle = -\frac{5}{2}\alpha_0 + \frac{5}{4}\alpha_0^2,$$

so that, $\frac{df|\mathbf{x}_1\rangle}{d\alpha_0} = -\frac{5}{2} + \frac{5}{2}\alpha_0$. Now, to find the minimizer of $f|\mathbf{x}_1\rangle$, set $\frac{df|\mathbf{x}_1\rangle}{d\alpha_0} = 0$, which implies $\alpha_0 = 1$. Note that $\frac{d^2f|\mathbf{x}_1\rangle}{d\alpha_0^2} = \frac{5}{2} > 0$, so that $\alpha_0 = 1$ is the minimizer of $f|\mathbf{x}_1\rangle$, which yield

$$|\mathbf{x}_1\rangle = (-1, 3/2)^T,$$
$$|\mathbf{g}_1\rangle = (0,0)^T.$$

Since $\|\mathbf{g}_1\| = \|\nabla f|\mathbf{x}_1\rangle\| = 0 < \varepsilon$. Thus, $|\mathbf{x}_1\rangle = (-1,3/2)^T$ is the minimizer of $f|\mathbf{x}\rangle$.

### 9.4.1. Approximating the Inverse Hessian (Computation of $\mathbf{S}_k$)

The idea of the quasi-Newton method is to approximate the inverse of Hessian by some other matrix which should be positive definite. This saves the work of computation of second derivatives and avoids the difficulties associated with the loss of positive definiteness. Positive definiteness is essential; otherwise, the search direction $|\mathbf{d}_k\rangle$ is not guaranteed to be a descent direction (i.e., $f|\mathbf{x}_{k+1}\rangle \nless f|\mathbf{x}_k\rangle$.) So, we approximate $\mathbf{H}_k^{-1}$ by another matrix $\mathbf{S}_k$, using only the first partial derivatives of the objective function $f$, [1-5].

Let $f|\mathbf{x}\rangle \in C^2$ be a function in $\mathbb{R}^n$ and assume that the gradients of $f|\mathbf{x}\rangle$ at points $|\mathbf{x}_k\rangle$ and $|\mathbf{x}_{k+1}\rangle$ are $|\mathbf{g}_k\rangle$ and $|\mathbf{g}_{k+1}\rangle$, respectively. Now, suppose that

$$|\boldsymbol{\gamma}_k\rangle = |\mathbf{g}_{k+1}\rangle - |\mathbf{g}_k\rangle, \tag{9.50}$$

and

$$|\boldsymbol{\delta}_k\rangle = |\mathbf{x}_{k+1}\rangle - |\mathbf{x}_k\rangle. \tag{9.51}$$

Consider the Taylor expansion of the objective:

$$f|\mathbf{x}_{k+1}\rangle = f|\mathbf{x}_k\rangle + \langle\nabla f(\mathbf{x}_k)|\boldsymbol{\delta}_k\rangle + \frac{1}{2}\langle\boldsymbol{\delta}_k|\nabla^2 f(\mathbf{x}_k)|\boldsymbol{\delta}_k\rangle + \cdots, \tag{9.52}$$

then the Taylor series gives the elements of $|\mathbf{g}_{k+1}\rangle$ as

$$|\mathbf{g}_{k+1}\rangle = \nabla f|\mathbf{x}_{k+1}\rangle = \nabla f|\mathbf{x}_k\rangle + \nabla^2 f(\mathbf{x}_k)|\boldsymbol{\delta}_k\rangle + \frac{1}{2}\langle\boldsymbol{\delta}_k|\nabla^3 f(\mathbf{x}_k)|\boldsymbol{\delta}_k\rangle + \cdots. \tag{9.53}$$

Now if the objective function $f|\mathbf{x}\rangle$ is quadratic, the second derivative of $f|\mathbf{x}\rangle$ are constant (independent of $|\mathbf{x}\rangle$) and, in turn, the second derivatives of $|\mathbf{g}_k\rangle$ are zero (i.e., $\nabla^3 f|\mathbf{x}_k\rangle = 0$). Thus

$$\nabla f|\mathbf{x}_{k+1}\rangle = \nabla f|\mathbf{x}_k\rangle + \nabla^2 f(\mathbf{x}_k)|\boldsymbol{\delta}_k\rangle, \tag{9.54}$$

or

$$|\mathbf{g}_{k+1}\rangle = |\mathbf{g}_k\rangle + \mathbf{H}|\boldsymbol{\delta}_k\rangle, \tag{9.55}$$

where $\mathbf{H}$ is the Hessian of $f|\mathbf{x}\rangle$. Hence, we have

$$|\boldsymbol{\gamma}_k\rangle = \mathbf{H}|\boldsymbol{\delta}_k\rangle. \tag{9.56}$$

If the gradient is evaluated sequentially at $n+1$ points, say, $|\mathbf{x}_0\rangle, |\mathbf{x}_1\rangle, \dots, |\mathbf{x}_n\rangle$ such that the changes in $|\mathbf{x}\rangle$, namely,

$$\begin{aligned} |\boldsymbol{\delta}_0\rangle &= |\mathbf{x}_1 - \mathbf{x}_0\rangle, \\ |\boldsymbol{\delta}_1\rangle &= |\mathbf{x}_2 - \mathbf{x}_1\rangle, \\ &\vdots \\ |\boldsymbol{\delta}_{n-1}\rangle &= |\mathbf{x}_n - \mathbf{x}_{n-1}\rangle, \end{aligned} \tag{9.57}$$

form a set of linearly independent vectors; then sufficient information is obtained to determine $\mathbf{H}$ uniquely. To illustrate this fact, $n$ equations of (9.56) can be re-arranged as

$$(|\boldsymbol{\gamma}_0\rangle, |\boldsymbol{\gamma}_1\rangle, \dots, |\boldsymbol{\gamma}_{n-1}\rangle) = \mathbf{H}(|\boldsymbol{\delta}_0\rangle, |\boldsymbol{\delta}_1\rangle, \dots, |\boldsymbol{\delta}_{n-1}\rangle), \tag{9.58}$$

and, therefore,

$$\mathbf{H} = (|\boldsymbol{\gamma}_0\rangle, |\boldsymbol{\gamma}_1\rangle, \dots, |\boldsymbol{\gamma}_{n-1}\rangle)(|\boldsymbol{\delta}_0\rangle, |\boldsymbol{\delta}_1\rangle, \dots, |\boldsymbol{\delta}_{n-1}\rangle)^{-1}. \tag{9.59}$$

The solution exists if $|\boldsymbol{\delta}_0\rangle, |\boldsymbol{\delta}_1\rangle, \dots, |\boldsymbol{\delta}_{n-1}\rangle$ form a set of linearly independent vectors.

The strategy of the quasi-Newton methods is as follows.

1- Assume that a positive definite real symmetric matrix $\mathbf{S}_k$ is available for the minimization of the quadratic problem $f|\mathbf{x}\rangle = a + \langle\mathbf{b}|\mathbf{x}\rangle + \frac{1}{2}\langle\mathbf{x}|\mathbf{H}|\mathbf{x}\rangle$, which is an approximation of $\mathbf{H}^{-1}$.
2- Compute a quasi-Newton direction

$$|\mathbf{d}_k\rangle = -\mathbf{S}_k|\mathbf{g}_k\rangle. \tag{9.60}$$

3- Find $\alpha_k$, the value of $\alpha$ that minimizes $f|\mathbf{x}_k + \alpha\mathbf{d}_k\rangle$, by differentiating $f|\mathbf{x}_k - \alpha\mathbf{S}_k\mathbf{g}_k\rangle$ with respect to $\alpha$ and then setting the result to zero, the value of $\alpha$ that minimizes $f|\mathbf{x}_k - \alpha\mathbf{S}_k\mathbf{g}_k\rangle$ can be deduced as

$$\alpha_k = \frac{\langle \mathbf{g}_k|\mathbf{S}_k|\mathbf{g}_k\rangle}{\langle \mathbf{g}_k|\mathbf{S}_k\mathbf{H}\mathbf{S}_k|\mathbf{g}_k\rangle} = \frac{\langle \mathbf{g}_k|\mathbf{S}_k|\mathbf{g}_k\rangle}{\langle \mathbf{S}_k\mathbf{g}_k|\mathbf{H}|\mathbf{S}_k\mathbf{g}_k\rangle}, \tag{9.61}$$

where $\mathbf{S}_k$ and $\mathbf{H}$ are positive definite and $|\mathbf{g}_k\rangle = |\mathbf{b}\rangle + \mathbf{H}|\mathbf{x}_k\rangle$ is the gradient of $f|\mathbf{x}\rangle$ at $|\mathbf{x}\rangle = |\mathbf{x}_k\rangle$.

4- Determine a change in $|\mathbf{x}\rangle$ as

$$|\boldsymbol{\delta}_k\rangle = \alpha_k|\mathbf{d}_k\rangle, \tag{9.62}$$

and deduce a new point $|\mathbf{x}_{k+1}\rangle$, using (9.51).

5- The change in the gradient, $|\boldsymbol{\gamma}_k\rangle$, can be computed using (9.50) (by computing the gradient at points $|\mathbf{x}_k\rangle$ and $|\mathbf{x}_{k+1}\rangle$.)

6- Apply a correction to $\mathbf{S}_k$ and generate

$$\mathbf{S}_{k+1} = \mathbf{S}_k + \mathbf{C}_k, \tag{9.63}$$

where $\mathbf{C}_k$ is an $n \times n$ correction matrix which can be computed from available data. Note that the equations that the matrices $\mathbf{S}_k$ are required to satisfy, do not determine those matrices uniquely. Thus, we have some freedom in the way we compute the $\mathbf{S}_k$ matrix.

### 9.4.2. Convergence Analysis of quasi-Newton algorithms for quadratic problem

Let $\mathbf{S}_0, \mathbf{S}_1, \mathbf{S}_2, \ldots$ be successive approximations of the inverse $\mathbf{H}_k^{-1}$ of the Hessian. On applying the above procedure iteratively starting with an initial point $|\mathbf{x}_0\rangle$ and an initial positive definite matrix $\mathbf{S}_0$, say, $\mathbf{S}_0 = \mathbf{I}_n$, the sequences $|\boldsymbol{\delta}_0\rangle, |\boldsymbol{\delta}_1\rangle, \ldots, |\boldsymbol{\delta}_k\rangle$, $|\boldsymbol{\gamma}_0\rangle, |\boldsymbol{\gamma}_1\rangle, \ldots, |\boldsymbol{\gamma}_k\rangle$, and $\mathbf{S}_1, \mathbf{S}_2, \ldots, \mathbf{S}_{k+1}$ can be generated. If the approximation $\mathbf{S}_{k+1}$ of the Hessian satisfy

$$\mathbf{S}_{k+1}|\boldsymbol{\gamma}_i\rangle = |\boldsymbol{\delta}_i\rangle, \quad 0 \le i \le k, \tag{9.64}$$

then for $n$ steps $(k = n - 1)$, we can write

$$\begin{gathered} \mathbf{S}_n|\boldsymbol{\gamma}_0\rangle = |\boldsymbol{\delta}_0\rangle, \\ \mathbf{S}_n|\boldsymbol{\gamma}_1\rangle = |\boldsymbol{\delta}_1\rangle, \\ \vdots \\ \mathbf{S}_n|\boldsymbol{\gamma}_{n-1}\rangle = |\boldsymbol{\delta}_{n-1}\rangle. \end{gathered} \tag{9.65}$$

This set of equations can be represented as

$$\mathbf{S}_n(|\boldsymbol{\gamma}_0\rangle, |\boldsymbol{\gamma}_1\rangle, \ldots, |\boldsymbol{\gamma}_{n-1}\rangle) = (|\boldsymbol{\delta}_0\rangle, |\boldsymbol{\delta}_1\rangle, \ldots, |\boldsymbol{\delta}_{n-1}\rangle), \tag{9.66}$$

or

$$\mathbf{S}_n = (|\boldsymbol{\delta}_0\rangle, |\boldsymbol{\delta}_1\rangle, \ldots, |\boldsymbol{\delta}_{n-1}\rangle)(|\boldsymbol{\gamma}_0\rangle, |\boldsymbol{\gamma}_1\rangle, \ldots, |\boldsymbol{\gamma}_{n-1}\rangle)^{-1}, \tag{9.67}$$

and from (9.58) and (9.67), we find that if $(|\boldsymbol{\gamma}_0\rangle, |\boldsymbol{\gamma}_1\rangle, \ldots, |\boldsymbol{\gamma}_{n-1}\rangle)$ is nonsingular, then $\mathbf{H}^{-1}$ is determined uniquely after $n$ steps, via

$$\mathbf{S}_n = \mathbf{H}^{-1} = (|\boldsymbol{\delta}_0\rangle, |\boldsymbol{\delta}_1\rangle, \ldots, |\boldsymbol{\delta}_{n-1}\rangle)(|\boldsymbol{\gamma}_0\rangle, |\boldsymbol{\gamma}_1\rangle, \ldots, |\boldsymbol{\gamma}_{n-1}\rangle)^{-1}. \tag{9.68}$$

As a consequence, we conclude that if $\mathbf{S}_n$ satisfies the equations $\mathbf{S}_n|\boldsymbol{\gamma}_i\rangle = |\boldsymbol{\delta}_i\rangle$, $0 \le i \le n-1$, then the algorithm $|\mathbf{x}_{k+1}\rangle = |\mathbf{x}_k\rangle - \alpha_k\mathbf{S}_k|\mathbf{g}_k\rangle$, $\alpha_k = \arg\min_{\alpha\ge 0} f|\mathbf{x}_k - \alpha\mathbf{S}_k\mathbf{g}_k\rangle$, is guaranteed to solve problems with quadratic objective functions in $n+1$ steps, because the update $|\mathbf{x}_{n+1}\rangle = |\mathbf{x}_n\rangle - \alpha_n\mathbf{S}_n|\mathbf{g}_n\rangle$ is equivalent to Newton algorithm.

Hence, quasi-Newton algorithms have the form.

$$|\mathbf{x}_{k+1}\rangle = |\mathbf{x}_k\rangle + \alpha_k|\mathbf{d}_k\rangle, \tag{9.69.1}$$

$$|\mathbf{d}_k\rangle = -\mathbf{S}_k|\mathbf{g}_k\rangle, \tag{9.69.2}$$

$$\alpha_k = \arg\min_{\alpha\ge 0} f|\mathbf{x}_k + \alpha\mathbf{d}_k\rangle, \tag{9.69.3}$$

where the matrices $\mathbf{S}_0, \mathbf{S}_1, \ldots$ are symmetric.

**Theorem 9.5:** Consider a quasi-Newton algorithm applied to a quadratic function with Hessian $\mathbf{H} = \mathbf{H}^T$ such that for $0 \le k < n-1$,

$$\mathbf{S}_{k+1}|\boldsymbol{\gamma}_i\rangle = |\boldsymbol{\delta}_i\rangle,\ 0 \le i \le k, \tag{9.70}$$

where $\mathbf{S}_{k+1} = \mathbf{S}_{k+1}^T$. If $\alpha_i \neq 0$, $0 \le i \le k$, then $|\mathbf{d}_0\rangle, \dots, |\mathbf{d}_{k+1}\rangle$ are $\mathbf{H}$-conjugate.

**Proof:**

We proceed by induction. We begin with the $k = 0$ case: that $|\mathbf{d}_0\rangle$ and $|\mathbf{d}_1\rangle$ are $\mathbf{H}$-conjugate. Because $\alpha_0 \neq 0$, we can write $|\mathbf{d}_0\rangle = |\boldsymbol{\delta}_0\rangle/\alpha_0$. Hence,

$$\begin{aligned}\langle\mathbf{d}_1|\mathbf{H}|\mathbf{d}_0\rangle &= -\langle\mathbf{g}_1|\mathbf{S}_1\mathbf{H}|\mathbf{d}_0\rangle\\ &= -\frac{1}{\alpha_0}\langle\mathbf{g}_1|\mathbf{S}_1\mathbf{H}|\boldsymbol{\delta}_0\rangle\\ &= -\frac{1}{\alpha_0}\langle\mathbf{g}_1|\mathbf{S}_1|\boldsymbol{\gamma}_0\rangle\\ &= -\frac{1}{\alpha_0}\langle\mathbf{g}_1|\boldsymbol{\delta}_0\rangle = -\langle\mathbf{g}_1|\mathbf{d}_0\rangle,\end{aligned}$$

where we used $|\mathbf{d}_1\rangle = -\mathbf{S}_1|\mathbf{g}_1\rangle \Rightarrow \langle\mathbf{d}_1| = -(\mathbf{S}_1|\mathbf{g}_1\rangle)^T = -\langle\mathbf{g}_1|\mathbf{S}_1^T = -\langle\mathbf{g}_1|\mathbf{S}_1$, $\mathbf{H}|\boldsymbol{\delta}_0\rangle = |\boldsymbol{\gamma}_0\rangle$, and $\mathbf{S}_1|\boldsymbol{\gamma}_0\rangle = |\boldsymbol{\delta}_0\rangle$. But

$$\langle\mathbf{d}_0|\mathbf{g}_1\rangle = \langle\mathbf{d}_0|\boldsymbol{\nabla} f(\mathbf{x}_0 + \alpha_0\mathbf{d}_0)\rangle = \phi_0'(\alpha_0).$$

Since $\alpha_0 = \arg\min_{\alpha\ge 0} f|\mathbf{x}_0 + \alpha\mathbf{d}_0\rangle > 0$, we have $\phi_0'(\alpha_0) = 0$. Hence, $\langle\mathbf{g}_1|\mathbf{d}_0\rangle = 0$. As a consequence of $\alpha_0 > 0$ being the minimizer of $\phi(\alpha) = f|\mathbf{x}_0 + \alpha\mathbf{d}_0\rangle$. As a result, $\langle\mathbf{d}_1|\mathbf{H}|\mathbf{d}_0\rangle = 0$. Assume that the result is true for $k-1$ (where $k < n-1$ ). It suffices to show that $\langle\mathbf{d}_{k+1}|\mathbf{H}|\mathbf{d}_i\rangle = 0$, $0 \le i \le k$. Given $i$, $0 \le i \le k$, using the same steps as in the $k = 0$ case, and using the assumption that $\alpha_i \neq 0$, we obtain

$$\begin{aligned}\langle\mathbf{d}_{k+1}|\mathbf{H}|\mathbf{d}_i\rangle &= -\langle\mathbf{g}_{k+1}|\mathbf{S}_{k+1}\mathbf{H}|\mathbf{d}_i\rangle\\ &= -\frac{1}{\alpha_i}\langle\mathbf{g}_{k+1}|\mathbf{S}_{k+1}\mathbf{H}|\boldsymbol{\delta}_i\rangle\\ &= -\frac{1}{\alpha_i}\langle\mathbf{g}_{k+1}|\mathbf{S}_{k+1}|\boldsymbol{\gamma}_i\rangle\\ &= -\frac{1}{\alpha_i}\langle\mathbf{g}_{k+1}|\boldsymbol{\delta}_i\rangle = -\langle\mathbf{g}_{k+1}|\mathbf{d}_i\rangle,\end{aligned}$$

where we used $|\mathbf{d}_{k+1}\rangle = -\mathbf{S}_{k+1}|\mathbf{g}_{k+1}\rangle \Rightarrow \langle\mathbf{d}_{k+1}| = -(\mathbf{S}_{k+1}|\mathbf{g}_{k+1}\rangle)^T = -\langle\mathbf{g}_{k+1}|\mathbf{S}_{k+1}$, $|\mathbf{d}_i\rangle = |\boldsymbol{\delta}_i\rangle/\alpha_i$, $\mathbf{H}|\boldsymbol{\delta}_i\rangle = |\boldsymbol{\gamma}_i\rangle$, and $\mathbf{S}_{k+1}|\boldsymbol{\gamma}_i\rangle = |\boldsymbol{\delta}_i\rangle$. Because $|\mathbf{d}_0\rangle, \dots, |\mathbf{d}_k\rangle$ are $\mathbf{H}$-conjugate by assumption, we conclude that $\langle\mathbf{g}_{k+1}|\mathbf{d}_i\rangle = 0$. Hence, $\langle\mathbf{d}_{k+1}|\mathbf{H}|\mathbf{d}_i\rangle = 0$.

■

Now if $k = n$, (9.60)-(9.62) yield

$$|\mathbf{d}_n\rangle = -\mathbf{H}^{-1}|\mathbf{g}_n\rangle, \tag{9.71.1}$$
$$\alpha_n = 1, \tag{9.71.2}$$
$$|\boldsymbol{\delta}_n\rangle = -\mathbf{H}^{-1}|\mathbf{g}_n\rangle, \tag{9.71.3}$$

respectively, and, therefore, from (9.51)

$$|\mathbf{x}_{n+1}\rangle = |\mathbf{x}_n\rangle - \mathbf{H}^{-1}|\mathbf{g}_n\rangle = |\mathbf{x}^*\rangle, \tag{9.72}$$

as in the Newton method.

In any derivation of $\mathbf{C}_n$, $\mathbf{S}_{k+1}$ must satisfy (9.65), and the following properties are highly desirable:

1. Vectors $|\boldsymbol{\delta}_0\rangle, |\boldsymbol{\delta}_1\rangle, \dots, |\boldsymbol{\delta}_k\rangle$ should form a set of conjugate directions.
2. A positive definite matrix $\mathbf{S}_k$ should give rise to a positive definite matrix $\mathbf{S}_{k+1}$.

The first property will ensure that the properties of conjugate-direction methods apply to the quasi-Newton method as well. The second property will ensure that $|\mathbf{d}_k\rangle$ is a descent direction in every iteration, i.e., for $k = 0,1, \ldots$.

Consider the point $|\mathbf{x}_k + \boldsymbol{\delta}_k\rangle$, and let

$$|\boldsymbol{\delta}_k\rangle = \alpha|\mathbf{d}_k\rangle,$$

where

$$|\mathbf{d}_k\rangle = -\mathbf{S}_k|\mathbf{g}_k\rangle.$$

For $\alpha > 0$, the Taylor series gives

$$f|\mathbf{x}_k + \boldsymbol{\delta}_k\rangle = f|\mathbf{x}_k\rangle + \langle\mathbf{g}_k|\boldsymbol{\delta}_k\rangle + \frac{1}{2}\langle\boldsymbol{\delta}_k|\mathbf{H}(\mathbf{x}_k + c\boldsymbol{\delta}_k)|\boldsymbol{\delta}_k\rangle,$$

where $c$ is a constant in the range $0 \le c < 1$. On eliminating $|\boldsymbol{\delta}_k\rangle$, we obtain

$$\begin{aligned} f|\mathbf{x}_k + \boldsymbol{\delta}_k\rangle &= f|\mathbf{x}_k\rangle - \alpha\langle\mathbf{g}_k|\mathbf{S}_k|\mathbf{g}_k\rangle + O\left(\alpha\|\mathbf{d}_k\|_2\right) \\ &= f|\mathbf{x}_k\rangle - \left[\alpha\langle\mathbf{g}_k|\mathbf{S}_k|\mathbf{g}_k\rangle - O\left(\alpha\|\mathbf{d}_k\|_2\right)\right], \end{aligned}$$

where $O\left(\alpha\|\mathbf{d}_k\|_2\right)$ is the remainder, which approaches zero faster than $\alpha\|\mathbf{d}_k\|_2$. Now if $\mathbf{S}_k$ is positive definite, then for a sufficiently small $\alpha > 0$, we have

$$\alpha\langle\mathbf{g}_k|\mathbf{S}_k|\mathbf{g}_k\rangle - O\left(\alpha\|\mathbf{d}_k\|_2\right) > 0,$$

since $\alpha > 0$, $\langle\mathbf{g}_k|\mathbf{S}_k|\mathbf{g}_k\rangle > 0$, and $O\left(\alpha\|\mathbf{d}_k\|_2\right) \to 0$. Therefore,

$$f|\mathbf{x}_k + \boldsymbol{\delta}_k\rangle < f|\mathbf{x}_k\rangle,$$

that is, if $\mathbf{S}_k$ is positive definite, then $|\mathbf{d}_k\rangle$ is a descent direction.

## 9.5 The Rank One Correction Formula

The general formula for correcting the matrix $\mathbf{S}_k$ has been defined in (9.63), where $\mathbf{C}_k$ is considered to be the correction matrix added to $\mathbf{S}_k$. To derive a rank one formula [8-12], we choose a scaled outer product of a vector $|\mathbf{z}\rangle$ for $\mathbf{C}_k$ as

$$\mathbf{C}_k = \beta_k|\mathbf{z}_k\rangle\langle\mathbf{z}_k|. \tag{9.73}$$

Hence, the update formula is

$$\mathbf{S}_{k+1} = \mathbf{S}_k + \beta_k|\mathbf{z}_k\rangle\langle\mathbf{z}_k|, \tag{9.74}$$

where $\beta_k \in \mathbb{R}$ and $|\mathbf{z}_k\rangle \in \mathbb{R}^n$. Note that, the rank of the outer product $|\mathbf{z}_k\rangle\langle\mathbf{z}_k|$ equal 1.

$$\begin{aligned} \operatorname{rank}|\mathbf{z}_k\rangle\langle\mathbf{z}_k| &= \operatorname{rank}\left[\begin{pmatrix} z_{1,k} \\ \vdots \\ z_{n,k} \end{pmatrix}\left(z_{1,k}, \ldots, z_{n,k}\right)\right] \\ &= \operatorname{rank}\begin{pmatrix} z_{1,k}z_{1,k} & z_{1,k}z_{2,k} & \cdots & z_{1,k}z_{n,k} \\ z_{2,k}z_{1,k} & z_{2,k}z_{2,k} & \cdots & z_{2,k}z_{n,k} \\ \vdots & \vdots & \ddots & \vdots \\ z_{n,k}z_{1,k} & z_{n,k}z_{2,k} & \cdots & z_{n,k}z_{n,k} \end{pmatrix} \\ &= \operatorname{rank}\begin{pmatrix} z_{1,k}\left(z_{1,k}, z_{2,k}, \ldots, z_{n,k}\right) \\ z_{2,k}\left(z_{1,k}, z_{2,k}, \ldots, z_{n,k}\right) \\ \vdots \\ z_{n,k}\left(z_{1,k}, z_{2,k}, \ldots, z_{n,k}\right) \end{pmatrix} = 1. \end{aligned} \tag{9.75}$$

Our goal now is to determine $\beta_k$ and $|\mathbf{z}_k\rangle$, given $\mathbf{S}_k$, $|\boldsymbol{\gamma}_k\rangle$, $|\boldsymbol{\delta}_k\rangle$, so that the required relationship $\mathbf{S}_{k+1}|\boldsymbol{\gamma}_i\rangle = |\boldsymbol{\delta}_i\rangle$, $i = 0, \ldots, k$, is satisfied. To begin, let us first consider the relation,

$$\mathbf{S}_{k+1}|\boldsymbol{\gamma}_k\rangle = |\boldsymbol{\delta}_k\rangle. \tag{9.76}$$

In other words, given $\mathbf{S}_k$, $|\boldsymbol{\gamma}_k\rangle$, and $|\boldsymbol{\delta}_k\rangle$, we wish to find $\beta_k$ and $|\mathbf{z}_k\rangle$ to ensure that

$$\mathbf{S}_{k+1}|\boldsymbol{\gamma}_k\rangle = \left(\mathbf{S}_k + \beta_k|\mathbf{z}_k\rangle\langle\mathbf{z}_k|\right)|\boldsymbol{\gamma}_k\rangle = |\boldsymbol{\delta}_k\rangle. \tag{9.77}$$

From (9.76)

$$|\boldsymbol{\delta}_k\rangle = \mathbf{S}_k|\boldsymbol{\gamma}_k\rangle + \beta_k|\mathbf{z}_k\rangle\langle\mathbf{z}_k|\boldsymbol{\gamma}_k\rangle, \tag{9.78}$$

and hence

$$\langle\boldsymbol{\gamma}_k|(|\boldsymbol{\delta}_k\rangle - \mathbf{S}_k|\boldsymbol{\gamma}_k\rangle) = \beta_k\langle\boldsymbol{\gamma}_k|\mathbf{z}_k\rangle\langle\mathbf{z}_k|\boldsymbol{\gamma}_k\rangle = \beta_k\langle\mathbf{z}_k|\boldsymbol{\gamma}_k\rangle^2. \tag{9.79}$$

From (9.77), we have

$$|\boldsymbol{\delta}_k\rangle - \mathbf{S}_k|\boldsymbol{\gamma}_k\rangle = \beta_k|\mathbf{z}_k\rangle\langle\mathbf{z}_k|\boldsymbol{\gamma}_k\rangle = (\beta_k\langle\mathbf{z}_k|\boldsymbol{\gamma}_k\rangle)|\mathbf{z}_k\rangle, \tag{9.80.a}$$

$$(|\boldsymbol{\delta}_k\rangle - \mathbf{S}_k|\boldsymbol{\gamma}_k\rangle)^T = (\beta_k|\mathbf{z}_k\rangle\langle\mathbf{z}_k|\boldsymbol{\gamma}_k\rangle)^T = \beta_k\langle\boldsymbol{\gamma}_k|\mathbf{z}_k\rangle\langle\mathbf{z}_k| = (\beta_k\langle\mathbf{z}_k|\boldsymbol{\gamma}_k\rangle)\langle\mathbf{z}_k|, \tag{9.80.b}$$

since $\langle\mathbf{z}_k|\boldsymbol{\gamma}_k\rangle$ is a scalar. Hence

$$(|\boldsymbol{\delta}_k\rangle - \mathbf{S}_k|\boldsymbol{\gamma}_k\rangle)(|\boldsymbol{\delta}_k\rangle - \mathbf{S}_k|\boldsymbol{\gamma}_k\rangle)^T = \beta_k\langle\mathbf{z}_k|\boldsymbol{\gamma}_k\rangle^2\beta_k|\mathbf{z}_k\rangle\langle\mathbf{z}_k|, \tag{9.81}$$

and from (9.78) and (9.80), we have

$$\begin{aligned}\mathbf{C}_k &= \beta_k|\mathbf{z}_k\rangle\langle\mathbf{z}_k| \\ &= \frac{(|\boldsymbol{\delta}_k\rangle - \mathbf{S}_k|\boldsymbol{\gamma}_k\rangle)(|\boldsymbol{\delta}_k\rangle - \mathbf{S}_k|\boldsymbol{\gamma}_k\rangle)^T}{\beta_k\langle\mathbf{z}_k|\boldsymbol{\gamma}_k\rangle^2} \\ &= \frac{(|\boldsymbol{\delta}_k\rangle - \mathbf{S}_k|\boldsymbol{\gamma}_k\rangle)(|\boldsymbol{\delta}_k\rangle - \mathbf{S}_k|\boldsymbol{\gamma}_k\rangle)^T}{\langle\boldsymbol{\gamma}_k|(|\boldsymbol{\delta}_k\rangle - \mathbf{S}_k|\boldsymbol{\gamma}_k\rangle)} = \frac{|\boldsymbol{\delta}_k - \mathbf{S}_k\boldsymbol{\gamma}_k\rangle\langle\boldsymbol{\delta}_k - \mathbf{S}_k\boldsymbol{\gamma}_k|}{\langle\boldsymbol{\gamma}_k|\boldsymbol{\delta}_k - \mathbf{S}_k\boldsymbol{\gamma}_k\rangle}.\end{aligned} \tag{9.82}$$

With the correction matrix known, $\mathbf{S}_{k+1}$ can be deduced from (9.75) as

$$\begin{aligned}\mathbf{S}_{k+1} &= \mathbf{S}_k + \mathbf{C}_k \\ &= \mathbf{S}_k + \frac{|\boldsymbol{\delta}_k - \mathbf{S}_k\boldsymbol{\gamma}_k\rangle\langle\boldsymbol{\delta}_k - \mathbf{S}_k\boldsymbol{\gamma}_k|}{\langle\boldsymbol{\gamma}_k|\boldsymbol{\delta}_k - \mathbf{S}_k\boldsymbol{\gamma}_k\rangle}.\end{aligned} \tag{9.83}$$

This formula is known as the unique rank one update formula for $\mathbf{S}_{k+1}$. To implement (9.82), an initial symmetric positive definite matrix is selected for $\mathbf{S}_0$ at the start of the algorithm, and the next point $|\mathbf{x}_1\rangle$ is computed using (9.69). Then the new matrix $\mathbf{S}_1$ is computed using (9.82) and the new point $|\mathbf{x}_2\rangle$ is determined from (9.69). This iterative process is continued until convergence is achieved. If $\mathbf{S}_k$ is symmetric, (9.82) ensures that $\mathbf{S}_{k+1}$ is also symmetric. However, there is no guarantee that $\mathbf{S}_{k+1}$ remains positive definite even if $\mathbf{S}_k$ is positive definite. This might lead to a breakdown of the procedure, especially when used for the optimization of nonquadratic functions.

There are two problems associated with the rank one method. First, a positive definite matrix $\mathbf{S}_k$ may not yield a positive definite matrix $\mathbf{S}_{k+1}$. Second, the denominator in formula (9.82) may become zero, and the method will break down since $\mathbf{S}_{k+1}$ will become undefined. From (9.82), we can write

$$\begin{aligned}\langle\boldsymbol{\gamma}_i|\mathbf{S}_{k+1}|\boldsymbol{\gamma}_i\rangle &= \langle\boldsymbol{\gamma}_i|\mathbf{S}_k|\boldsymbol{\gamma}_i\rangle + \langle\boldsymbol{\gamma}_i|\frac{|\boldsymbol{\delta}_k - \mathbf{S}_k\boldsymbol{\gamma}_k\rangle\langle\boldsymbol{\delta}_k - \mathbf{S}_k\boldsymbol{\gamma}_k|}{\langle\boldsymbol{\gamma}_k|\boldsymbol{\delta}_k - \mathbf{S}_k\boldsymbol{\gamma}_k\rangle}|\boldsymbol{\gamma}_i\rangle \\ &= \langle\boldsymbol{\gamma}_i|\mathbf{S}_k|\boldsymbol{\gamma}_i\rangle + \frac{\langle\boldsymbol{\gamma}_i|\boldsymbol{\delta}_k - \mathbf{S}_k\boldsymbol{\gamma}_k\rangle\langle\boldsymbol{\delta}_k - \mathbf{S}_k\boldsymbol{\gamma}_k|\boldsymbol{\gamma}_i\rangle}{\langle\boldsymbol{\gamma}_k|\boldsymbol{\delta}_k - \mathbf{S}_k\boldsymbol{\gamma}_k\rangle} \\ &= \langle\boldsymbol{\gamma}_i|\mathbf{S}_k|\boldsymbol{\gamma}_i\rangle + \frac{[\langle\boldsymbol{\gamma}_i|\boldsymbol{\delta}_k\rangle - \langle\boldsymbol{\gamma}_i|\mathbf{S}_k|\boldsymbol{\gamma}_k\rangle][\langle\boldsymbol{\delta}_k|\boldsymbol{\gamma}_i\rangle - \langle\boldsymbol{\gamma}_k|\mathbf{S}_k|\boldsymbol{\gamma}_i\rangle]}{\langle\boldsymbol{\gamma}_k|\boldsymbol{\delta}_k - \mathbf{S}_k\boldsymbol{\gamma}_k\rangle} \\ &= \langle\boldsymbol{\gamma}_i|\mathbf{S}_k|\boldsymbol{\gamma}_i\rangle + \frac{[\langle\boldsymbol{\gamma}_i|\boldsymbol{\delta}_k\rangle - \langle\boldsymbol{\gamma}_i|\mathbf{S}_k|\boldsymbol{\gamma}_k\rangle]^2}{\langle\boldsymbol{\gamma}_k|\boldsymbol{\delta}_k - \mathbf{S}_k\boldsymbol{\gamma}_k\rangle}.\end{aligned} \tag{9.84}$$

Hence, if $\mathbf{S}_k$ is positive definite, a sufficient condition for $\mathbf{S}_{k+1}$ to be positive and definite is

$$\langle\boldsymbol{\gamma}_k|\boldsymbol{\delta}_k - \mathbf{S}_k\boldsymbol{\gamma}_k\rangle > 0. \tag{9.85}$$

### 9.5.1. Convergence Analysis of Rank One Correction Formula

**Theorem 9.6:** For the rank one algorithm, if $\mathbf{H} = \mathbf{H}^T$ is the Hessian of a convex quadratic problem and

$$|\boldsymbol{\gamma}_i\rangle = \mathbf{H}|\boldsymbol{\delta}_i\rangle, \quad 0 \le i \le k, \tag{9.86}$$

where $|\boldsymbol{\delta}_i\rangle$ for $i = 0, \dots, k$ are given linearly independent vectors, then for any initial symmetric matrix $\mathbf{S}_0$

$$|\boldsymbol{\delta}_i\rangle = \mathbf{S}_{k+1}|\boldsymbol{\gamma}_i\rangle, \quad 0 \le i \le k, \tag{9.87}$$

where

$$\mathbf{S}_{i+1} = \mathbf{S}_i + \frac{|\boldsymbol{\delta}_i - \mathbf{S}_i\boldsymbol{\gamma}_i\rangle\langle\boldsymbol{\delta}_i - \mathbf{S}_i\boldsymbol{\gamma}_i|}{\langle\boldsymbol{\gamma}_i|\boldsymbol{\delta}_i - \mathbf{S}_i\boldsymbol{\gamma}_i\rangle}. \tag{9.88}$$

**Proof:**

We prove the theorem by induction. Suppose now that the theorem is true for $i < k$; that is

$$|\boldsymbol{\delta}_i\rangle = \mathbf{S}_k|\boldsymbol{\gamma}_i\rangle, \quad 0 \le i \le k-1,$$

and show that

$$|\boldsymbol{\delta}_i\rangle = \mathbf{S}_{k+1}|\boldsymbol{\gamma}_i\rangle, \quad 0 \le i \le k.$$

Our construction of the correction matrix ensures that $|\boldsymbol{\delta}_k\rangle = \mathbf{S}_{k+1}|\boldsymbol{\gamma}_k\rangle$. So we only have to show that

$$|\boldsymbol{\delta}_i\rangle = \mathbf{S}_{k+1}|\boldsymbol{\gamma}_i\rangle, \quad i < k.$$

If $0 \le i \le k-1$, since $\mathbf{S}_k$ is symmetric, (9.83) yields

$$\begin{aligned}\mathbf{S}_{k+1}|\boldsymbol{\gamma}_i\rangle &= \left(\mathbf{S}_k + \frac{|\boldsymbol{\delta}_k - \mathbf{S}_k\boldsymbol{\gamma}_k\rangle\langle\boldsymbol{\delta}_k - \mathbf{S}_k\boldsymbol{\gamma}_k|}{\langle\boldsymbol{\gamma}_k|\boldsymbol{\delta}_k - \mathbf{S}_k\boldsymbol{\gamma}_k\rangle}\right)|\boldsymbol{\gamma}_i\rangle \\ &= \mathbf{S}_k|\boldsymbol{\gamma}_i\rangle + |\boldsymbol{\zeta}_k\rangle\langle\boldsymbol{\delta}_k - \mathbf{S}_k\boldsymbol{\gamma}_k|\boldsymbol{\gamma}_i\rangle \\ &= \mathbf{S}_k|\boldsymbol{\gamma}_i\rangle + |\boldsymbol{\zeta}_k\rangle(\langle\boldsymbol{\delta}_k|\boldsymbol{\gamma}_i\rangle - \langle\boldsymbol{\gamma}_k|\mathbf{S}_k|\boldsymbol{\gamma}_i\rangle) \\ &= |\boldsymbol{\delta}_i\rangle + |\boldsymbol{\zeta}_k\rangle(\langle\boldsymbol{\delta}_k|\boldsymbol{\gamma}_i\rangle - \langle\boldsymbol{\gamma}_k|\boldsymbol{\delta}_i\rangle),\end{aligned}$$

where

$$|\boldsymbol{\zeta}_k\rangle = \frac{|\boldsymbol{\delta}_k - \mathbf{S}_k\boldsymbol{\gamma}_k\rangle}{\langle\boldsymbol{\gamma}_k|\boldsymbol{\delta}_k - \mathbf{S}_k\boldsymbol{\gamma}_k\rangle},$$

and we use $|\boldsymbol{\delta}_i\rangle = \mathbf{S}_k|\boldsymbol{\gamma}_i\rangle$, $0 \le i \le k-1$ (from induction assumption).
For $0 \le i \le k$

$$|\boldsymbol{\gamma}_i\rangle = \mathbf{H}|\boldsymbol{\delta}_i\rangle,$$

and

$$\langle\boldsymbol{\gamma}_k| = \langle\boldsymbol{\delta}_k|\mathbf{H}.$$

Hence for $0 \le i \le k-1$, we have

$$\langle\boldsymbol{\delta}_k|\boldsymbol{\gamma}_i\rangle - \langle\boldsymbol{\gamma}_k|\boldsymbol{\delta}_i\rangle = \langle\boldsymbol{\delta}_k|\mathbf{H}|\boldsymbol{\delta}_i\rangle - \langle\boldsymbol{\delta}_k|\mathbf{H}|\boldsymbol{\delta}_i\rangle = 0,$$

and consequently, we get

$$\mathbf{S}_{k+1}|\boldsymbol{\gamma}_i\rangle = |\boldsymbol{\delta}_i\rangle, \quad 0 \le i \le k-1.$$

By combining $|\boldsymbol{\delta}_k\rangle = \mathbf{S}_{k+1}|\boldsymbol{\gamma}_k\rangle$ with the above equation, we obtain

$$\mathbf{S}_{k+1}|\boldsymbol{\gamma}_i\rangle = |\boldsymbol{\delta}_i\rangle, \quad 0 \le i \le k.$$

To complete the induction, we have

$$\mathbf{S}_1|\boldsymbol{\gamma}_i\rangle = |\boldsymbol{\delta}_i\rangle, \quad 0 \le i \le 0,$$

and we can write

$$\begin{aligned}\mathbf{S}_2|\boldsymbol{\gamma}_i\rangle &= |\boldsymbol{\delta}_i\rangle, \quad 0 \le i \le 1, \\ \mathbf{S}_3|\boldsymbol{\gamma}_i\rangle &= |\boldsymbol{\delta}_i\rangle, \quad 0 \le i \le 2, \\ &\vdots \\ \mathbf{S}_{k+1}|\boldsymbol{\gamma}_i\rangle &= |\boldsymbol{\delta}_i\rangle, \quad 0 \le i \le k.\end{aligned}$$

■

Quasi-Newton Algorithms:
The Rank One Correction.
The Davidon Fletcher Powell (DFP) Method (see section 9.6.1)

**Algorithm** The Broyden Fletcher Goldfarb Shanno (BFGS) Method (see section 9.6.2)

**Step 1:** Input $|\mathbf{x}_0\rangle$ and initialize the tolerance $\varepsilon$.

**Step 2:** Set $k = 0$ and select a real symmetric positive definite matrix $\mathbf{S}_0$ (put $\mathbf{S}_0 = \mathbf{I}_n$). Compute $|\mathbf{g}_0\rangle$.

**Step 3:** If $|\mathbf{g}_k\rangle = |\mathbf{0}\rangle$, stop; else, $|\mathbf{d}_k\rangle = -\mathbf{S}_k|\mathbf{g}_k\rangle$.
Compute $\alpha_k = \arg\min_{\alpha \ge 0} f|\mathbf{x}_k + \alpha\mathbf{d}_k\rangle$, the value of $\alpha$ that minimizes $f|\mathbf{x}_k + \alpha\mathbf{d}_k\rangle$, using a line search.
Set $|\boldsymbol{\delta}_k\rangle = \alpha_k|\mathbf{d}_k\rangle$ and update $|\mathbf{x}_{k+1}\rangle = |\mathbf{x}_k\rangle + |\boldsymbol{\delta}_k\rangle$.

**Step 4:** If $\|\boldsymbol{\delta}_k\|_2 < \varepsilon$, output $|\mathbf{x}^*\rangle = |\mathbf{x}_{k+1}\rangle$ and $f|\mathbf{x}^*\rangle = f|\mathbf{x}_{k+1}\rangle$, and stop.

**Step 5:** Compute $|\mathbf{g}_{k+1}\rangle$ and set $|\boldsymbol{\gamma}_k\rangle = |\mathbf{g}_{k+1}\rangle - |\mathbf{g}_k\rangle$.

**The rank one correction method**

Compute $\mathbf{S}_{k+1}$ by the formula

$$\mathbf{S}_{k+1} = \mathbf{S}_k + \frac{|\boldsymbol{\delta}_k - \mathbf{S}_k\boldsymbol{\gamma}_k\rangle\langle\boldsymbol{\delta}_k - \mathbf{S}_k\boldsymbol{\gamma}_k|}{\langle\boldsymbol{\gamma}_k|\boldsymbol{\delta}_k - \mathbf{S}_k\boldsymbol{\gamma}_k\rangle}.$$

**The DFP method**

Compute $\mathbf{S}_{k+1}$ by the formula

$$\mathbf{S}_{k+1} = \mathbf{S}_k + \frac{|\boldsymbol{\delta}_k\rangle\langle\boldsymbol{\delta}_k|}{\langle\boldsymbol{\delta}_k|\boldsymbol{\gamma}_k\rangle} - \frac{|\mathbf{S}_k\boldsymbol{\gamma}_k\rangle\langle\mathbf{S}_k\boldsymbol{\gamma}_k|}{\langle\boldsymbol{\gamma}_k|\mathbf{S}_k|\boldsymbol{\gamma}_k\rangle}.$$

**The BFGS method**

Compute $\mathbf{S}_{k+1}$ by the formula

$$\mathbf{S}_{k+1} = \mathbf{S}_k + \left(1 + \frac{\langle\boldsymbol{\gamma}_k|\mathbf{S}_k|\boldsymbol{\gamma}_k\rangle}{\langle\boldsymbol{\gamma}_k|\boldsymbol{\delta}_k\rangle}\right)\frac{|\boldsymbol{\delta}_k\rangle\langle\boldsymbol{\delta}_k|}{\langle\boldsymbol{\delta}_k|\boldsymbol{\gamma}_k\rangle} - \frac{|\mathbf{S}_k\boldsymbol{\gamma}_k\rangle\langle\boldsymbol{\delta}_k| + |\boldsymbol{\delta}_k\rangle\langle\mathbf{S}_k\boldsymbol{\gamma}_k|}{\langle\boldsymbol{\gamma}_k|\boldsymbol{\delta}_k\rangle}.$$

**Step 6:** Set $k = k + 1$ and go to step 3.

In step 3, the value of $\alpha_k$ is obtained by using a line search to render the algorithm more amenable to nonquadratic problems. However, for convex quadratic problems, $\alpha_k$ should be calculated by using (9.61).

***Example 9.7***

Let

$$f|\mathbf{x}\rangle = a + \langle\mathbf{b}|\mathbf{x}\rangle + \frac{1}{2}\langle\mathbf{x}|\mathbf{H}|\mathbf{x}\rangle, \quad |\mathbf{x}\rangle \in \mathbb{R}^2,$$

with $a = 7$, $|\mathbf{b}\rangle = \begin{pmatrix}-1\\1\end{pmatrix}$, $\mathbf{H} = \begin{pmatrix}1 & 0\\0 & 2\end{pmatrix}$, and $|\mathbf{x}_0\rangle = \begin{pmatrix}0\\0\end{pmatrix}$. Use the rank one correction method to generate two $\mathbf{H}$-conjugate directions.

***Solution***

We first compute the gradient $|\mathbf{g}\rangle = \nabla f|\mathbf{x}\rangle$ and evaluate it at $|\mathbf{x}_0\rangle$,

$$|\mathbf{g}_0\rangle = |\mathbf{b}\rangle + \mathbf{H}|\mathbf{x}_0\rangle = \begin{pmatrix}-1\\1\end{pmatrix}.$$

It is a nonzero vector, so we proceed with the first iteration.

**Iteration 1.**

Let $\mathbf{S}_0 = \mathbf{I}_2$. Then,

$$|\mathbf{d}_0\rangle = -\mathbf{S}_0|\mathbf{g}_0\rangle = \begin{pmatrix}1\\-1\end{pmatrix}.$$

The step size $\alpha_0$ is

$$\alpha_0 = -\frac{\langle\mathbf{g}_0|\mathbf{d}_0\rangle}{\langle\mathbf{d}_0|\mathbf{H}|\mathbf{d}_0\rangle} = 2/3.$$

Hence,

$$|\mathbf{x}_1\rangle = |\mathbf{x}_0\rangle + \alpha_0|\mathbf{d}_0\rangle = \begin{pmatrix}2/3\\-2/3\end{pmatrix}.$$

We evaluate the gradient $|\mathbf{g}\rangle = \nabla f|\mathbf{x}\rangle$ at $|\mathbf{x}_1\rangle$ to obtain

$$|\mathbf{g}_1\rangle = |\mathbf{b}\rangle + \mathbf{H}|\mathbf{x}_1\rangle = \begin{pmatrix}-1/3\\-1/3\end{pmatrix}.$$

It is a nonzero vector, so we proceed with the second iteration.

**Iteration 2.**

Compute $\mathbf{S}_1$,

$$\mathbf{S}_1 = \mathbf{S}_0 + \frac{|\boldsymbol{\delta}_0 - \mathbf{S}_0\boldsymbol{\gamma}_0\rangle\langle\boldsymbol{\delta}_0 - \mathbf{S}_0\boldsymbol{\gamma}_0|}{\langle\boldsymbol{\gamma}_0|\boldsymbol{\delta}_0 - \mathbf{S}_0\boldsymbol{\gamma}_0\rangle}.$$

To find $\mathbf{S}_1$ we need to compute

$$|\boldsymbol{\gamma}_0\rangle = |\mathbf{g}_1\rangle - |\mathbf{g}_0\rangle = \begin{pmatrix} 2/3 \\ -4/3 \end{pmatrix}.$$

Using the above, we determine,

$$|\boldsymbol{\delta}_0\rangle - \mathbf{S}_0|\boldsymbol{\gamma}_0\rangle = \begin{pmatrix} 0 \\ 2/3 \end{pmatrix},$$
$$\langle\boldsymbol{\gamma}_0|(|\boldsymbol{\delta}_0\rangle - \mathbf{S}_0|\boldsymbol{\gamma}_0\rangle) = -8/9.$$

Then, we obtain

$$\mathbf{S}_1 = \begin{pmatrix} 1 & 0 \\ 0 & 1/2 \end{pmatrix},$$

and

$$|\mathbf{d}_1\rangle = -\mathbf{S}_1|\mathbf{g}_1\rangle = \begin{pmatrix} 1/3 \\ 1/6 \end{pmatrix}.$$

Note that $\langle\mathbf{d}_0|\mathbf{H}|\mathbf{d}_1\rangle = 0$, that is, $|\mathbf{d}_0\rangle$ and $|\mathbf{d}_1\rangle$ are **H**-conjugate.

***Example 9.8***

Consider the problem:

$$\text{Minimize } f|\mathbf{x}\rangle = (x^2 - 5x + 5)^2 + (y^2 - 2y + 5)^2,$$

using, $|\mathbf{x}_0\rangle = (5,5)^T$, and $\epsilon = 0.01$.

***Solution***

The 3D and contour plots of the function are shown in Figure 9.4. The plots show that the minimum lies at $|\mathbf{x}^*\rangle = (3.6176, 0.999846)^T$, $f|\mathbf{x}^*\rangle = 16$.

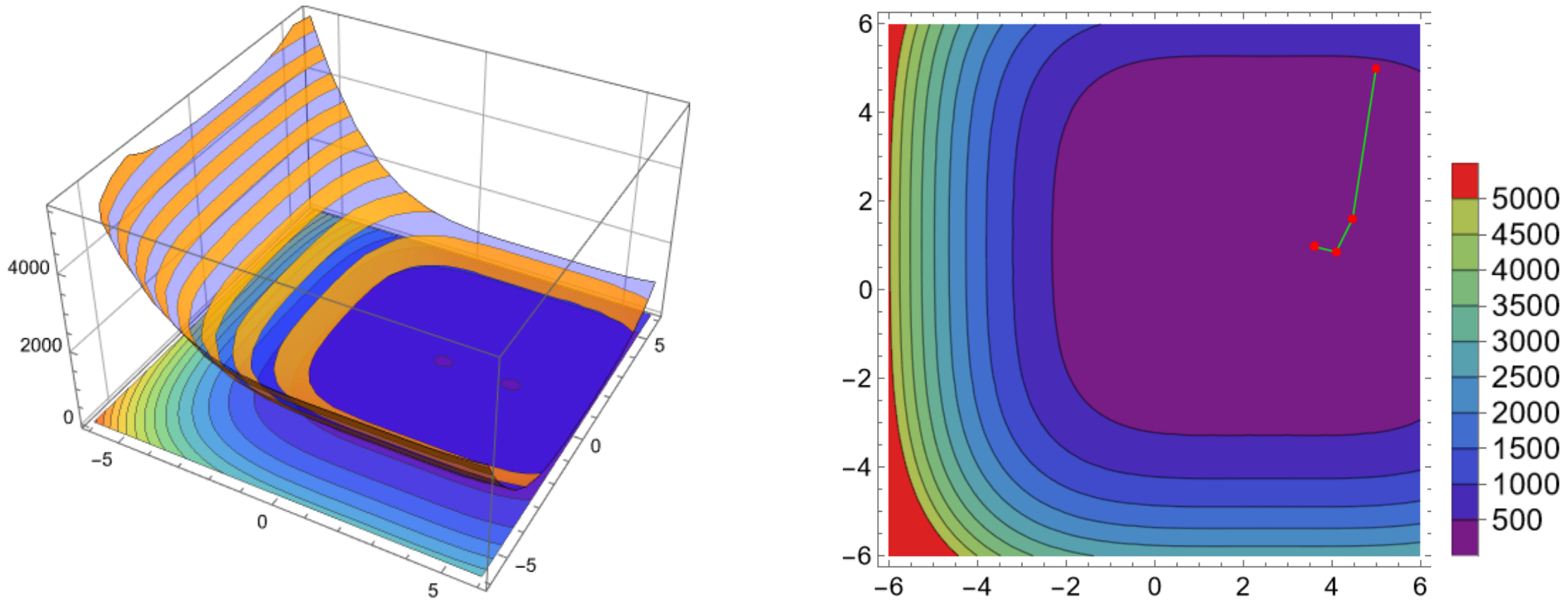

**Figure 9.4.** The results of 5 iterations of the quasi-Newton algorithms: the rank one correction method for $f|\mathbf{x}\rangle = (x^2 - 5x + 5)^2 + (y^2 - 2y + 5)^2$.

From Table 9.3, after 5 iterations, the optimal condition is attained. The results were produced by Mathematica code 9.3.

**Table 9.3.**

| No. of iters. | *astar* | *x*0 | *x*1 | $f(x0)$ | $f(x1)$ | error |
|---|---|---|---|---|---|---|
| 1 | 0.010613 | {5., 5.} | {4.469, 1.603} | 425 | 25.95832 | 23.2341 |
| 2 | -0.03667 | {4.469, 1.603} | {4.095, 0.857} | 25.95832 | 17.84122 | 8.578855 |
| 3 | 0.059401 | {4.095, 0.857} | {3.592, 0.981} | 17.84122 | 16.00582 | 0.381213 |
| 4 | -0.09696 | {3.592, 0.981} | {3.617, 1.001} | 16.00582 | 16.00003 | 0.031206 |
| 5 | 0.065162 | {3.617, 1.001} | {3.617, 0.999} | 16.00003 | 16 | 0.005005 |

***Mathematica Code 9.3*** Quasi-Newton Algorithm (The Rank One Correction)

```
(* Quasi Newton Algorithm (The Rank One Correction Method) *)
```

```
(*
Notations:
x0       :Intial vector
epsilonQNA :Small number to check the accuracy of the quasi-Newton algorithm (the rank one
correction method)
f[x,y]    :Objective function
lii       :The last iteration index
result[k] :The results of iteration k
*)

(* Taking Initial Inputs from User *)
x0=Input["Enter the intial point in the format {x, y}; for example {5,5} "] ;
epsilonQNA=Input["Please enter accuracy of the quasi-Newton algorithm (the rank one
correction method); for example 0.01 "];

If[
  epsilonQNA<=0,
  Beep[];
  MessageDialog["The value of epsilonQNA has to be postive number: "];
  Exit[];
  ];

domainx=Input["Please enter domain of x variable for 3D and contour plots; for example {-
6,6}"];
domainy=Input["Please enter domain of y variable for 3D and contour plots; for example {-
6,6}"];

(* Taking the Function from User *)
f[{x_,y_}] = Evaluate[Input["Please input a function of x and y to find the minimum "]];
(* For example: (x^2-5x+5)^2+(y^2-2y+5)^2 *)

(*Defination of the Unidirectionalsearch Function*)
(*This function start by bracketing the minimum (using Bounding Phase Method) then isolating
the minimum (using Golden Section Search Method) *)

unidirectionalsearch[α_,delt_,eps_]:=Module[


{α0=α,delta=delt,epsilon=eps,y1,y2,y3,αα,a,b,increment,a0,b0,anew,bnew,anorm,bnorm,lnorm,α1n
orm,α2norm,α1,α2,φ1,φ2,αstar},

   (* Bounding Phase Method *)
   (* Initiating Required Variables *)
   y1 =φ[α0-Abs[delta]];
   y2 =φ[α0];
   y3 =φ[α0+Abs[delta]];

   (*Determining Whether the Inicrement Is Positive or Negative*)
   Which[
    y1==y2,
    a=α0-Abs[delta];
    b=α0;
    Goto[end];,
    y2==y3,
    a=α0;
    b=α0+Abs[delta];
    Goto[end];,
    y1==y3||(y1>y2&&y2<y3),
    a=α0-Abs[delta];
    b=α0+Abs[delta];
    Goto[end];
```

```
    ];

   Which[
     y1>y2&&y2>y3,
     increment=Abs[delta];,
     y1<y2&&y2<y3,
     increment=-Abs[delta];
     ]

    (* Starting the Algorithm *)
    Do[
     αα[0]=α0;
     αα[k+1]=αα[k]+2^k*increment;

     Which[
      φ[αα[k]]<φ[αα[k+1]],(* Evidently, it is impossible the condition to hold for k=0 *)
      a=αα[k-1];
      b= αα[k+1];
      Break [],

      k>50,
      Print["After 50 iterations the bounding phase method can not braketing the min of
alpha"];
      Exit[]
      ];,
     {k,0,∞}
     ];

   Label[end];

   If[
    a>b,
    {a,b}={b,a}
    ];

   (* Golden Section Search Method *)
   (* Initiating Required Variables*)
   a0=a;
   b0=b;
   anew=a;
   bnew=b;

   If[
    a0==b0,
    αstar=a;
    Goto[final]
    ];

   (* Starting the Algorithm *)
   Do[
    (* Normalize the Variable α *)
    anorm=(anew-a)/(b-a);
    bnorm=(bnew-a)/(b-a);

    lnorm=bnorm-anorm;

    α1norm=anorm+0.382*lnorm;
    α2norm=bnorm-0.382*lnorm;

    α1=α1norm(b0-a0)+a0;
    α2=α2norm(b0-a0)+a0;
```

```
    φ1=φ[α1];
    φ2=φ[α2];

    Which [
     φ1>φ2,
     anew=α1(*move lower bound to α1*);,
     φ1<φ2,
     bnew=α2(*move upper bound to α2*);,
     φ1==φ2,
     anew=α1(*move lower bound to α1*);
     bnew=α2(*move upper bound to α2*);
     ];

    αstar=0.5*(anew+bnew);

    If[
     Abs[lnorm]<epsilon,
     Break[]
     ];,
    {k,1,∞}
    ];

   Label[final];

   (* Final Result *)
   N[αstar]
   ];

α00=2;(* The intial point of α; for example 2*)
delt0=1;(* The parameter delta of Bounding Phase Method; for example 1 *)
eps0=0.01;(* The accuracy of the Golden Section Search Method; for example 0.01 *)

gradfx[x_,y_]=Grad[f[{x,y}],{x,y}];
gradfx0=gradfx[x0[[1]],x0[[2]]];

If[
  N[Norm[gradfx0]]==0,
  Print["The quasi-Newton algorithm requires that the gradiant of the function f
!=0"];(*Ending program*)
  Exit[];
  ];

s0=IdentityMatrix[2];
dx0=-IdentityMatrix[2].gradfx0;
(* Main Loop *)
Do[
  φ[α_]=f[x0+α*dx0];
  αstar=unidirectionalsearch[α00,delt0,eps0];
  delta0=αstar*dx0;

  x1=x0+delta0;
  gradfx1=gradfx[x1[[1]],x1[[2]]];
  gamax0=gradfx1-gradfx0;


  s1=s0+((delta0-s0.gamax0).Transpose[delta0-s0.gamax0])/(Transpose[gamax0].(delta0-
s0.gamax0));
  dx1=-(s1.gradfx1);
  error=Norm[gradfx1];
```

```
  lii=k;
  result[k]=N[{k,αstar,Row[x0,","],Row[x1,","],f[x0],f[x1],error}];
  plotresult[k]=N[{x0,x1}];

  If[
   error<epsilonQNA||k>50,
   Break[],
   x0=x1;
   dx0=dx1;
   s0=s1;
   gradfx0=gradfx1;
   ],
  {k,1,∞}
  ];

(* Final Result *)
Print["The solution is x= ",  N[x1],"\nThe solution is (approximately)= ", N[f[x1]]]

(* Results of Each Iteration *)
table=TableForm[
  Table[
   result[i],
   {i,1,lii}
   ],
  TableHeadings->{None,{"No. of iters.","αstar","x0","x1","f[x0]","f[x1]","error"}}
  ]

Export["example93.xls",table,"XLS"];

(* Data Visualization *)

(* Domain of Varibles*)
xleft=domainx[[1]];
xright=domainx[[2]];
ydown=domainy[[1]];
yup=domainy[[2]];

(* 3D+ Contour Plot *)
plot1=Plot3D[
   f[{x,y}],
   {x,xleft,xright},
   {y,ydown,yup},
   ClippingStyle->None,
   MeshFunctions->{#3&},
   Mesh->15,
   MeshStyle->Opacity[.5],
   MeshShading->{{Opacity[.3],Blue},{Opacity[.8],Orange}},
   Lighting->"Neutral"
   ];

slice=SliceContourPlot3D[
   f[{x,y}],
   z==0,
   {x,xleft,xright},
   {y,ydown,yup},
   {z,-1,1},
   Contours->15,
   Axes->False,
   PlotPoints->50,
   PlotRangePadding->0,
   ColorFunction->"Rainbow"
```

```
  ];

Show[
 plot1,
 slice,
 PlotRange->All,
 BoxRatios->{1,1,.6},
 FaceGrids->{Back,Left}
 ]

(* Contour Plot with Step Iterations *)
ContourPlot[
 f[{x,y}],
 {x,xleft,xright},
 {y,ydown,yup},
 LabelStyle->Directive[Black,16],
 ColorFunction->"Rainbow",
 PlotLegends->Automatic,
 Contours->10,
 Epilog-
>{PointSize[0.015],Green,Line[Flatten[Table[plotresult[i],{i,1,lii}],1]],Red,Point[Flatten[T
able[plotresult[i],{i,1,lii}],1]]}
 ]

(* Data Manipulation *)
Manipulate[
 ContourPlot[
  f[{x,y}],
  {x,xleft,xright},
  {y,ydown,yup},
  LabelStyle->Directive[Black,14],
  ColorFunction->"Rainbow",
  PlotLegends->Automatic,
  Contours->10,
  Epilog->{
    PointSize[0.015],
    Yellow,
    Arrow[{plotresult[i][[1]],plotresult[i][[2]]}],
    Red,
    Point[Flatten[Table[plotresult[j],{j,1,i}],1]],
    Green,
    Line[Flatten[Table[plotresult[j],{j,1,i}],1]]
    }
  ],
 {i,1,lii,1}
 ]
```

## 9.6 The Rank Two Correction Formula

The rank two update formulas guarantee both symmetry and positive definiteness of the matrix $\mathbf{S}_{k+1}$ and are more robust in minimizing general nonlinear functions, hence are preferred in practical applications. The Davidon Fletcher Powell (DFP) and Broyden Fletcher Goldfarb Shanno (BFGS) iterative methods are described in detail in the following sections.

### 9.6.1. The DFP Method

The DFP algorithm [9,13] is superior to the rank one correction algorithm because it preserves the positive definiteness of the matrix $\mathbf{S}_k$. In rank two updates, we choose the update matrix $\mathbf{C}_k$ as the sum of two “rank one” updates as

$$\mathbf{C}_k = \beta_{1k}|\mathbf{z}_{1,k}\rangle\langle\mathbf{z}_{1,k}| + \beta_{2k}|\mathbf{z}_{2,k}\rangle\langle\mathbf{z}_{2,k}|, \tag{9.89}$$

where $\beta_{1k}, \beta_{2k} \in \mathbb{R}$ and the vectors $|\mathbf{z}_{1,k}\rangle, |\mathbf{z}_{2,k}\rangle \in \mathbb{R}^n$. Therefore, (9.63) and (9.89) lead to the new update formula

$$\mathbf{S}_{k+1} = \mathbf{S}_k + \beta_{1k}|\mathbf{z}_{1,k}\rangle\langle\mathbf{z}_{1,k}| + \beta_{2k}|\mathbf{z}_{2,k}\rangle\langle\mathbf{z}_{2,k}|. \tag{9.90}$$

By forcing (9.89) to satisfy the quasi-Newton condition, $|\boldsymbol{\delta}_k\rangle = \mathbf{S}_{k+1}|\boldsymbol{\gamma}_k\rangle$, we obtain

$$\begin{aligned} |\boldsymbol{\delta}_k\rangle &= \left(\mathbf{S}_k + \beta_{1k}|\mathbf{z}_{1,k}\rangle\langle\mathbf{z}_{1,k}| + \beta_{2k}|\mathbf{z}_{2,k}\rangle\langle\mathbf{z}_{2,k}|\right)|\boldsymbol{\gamma}_k\rangle \\ &= \mathbf{S}_k|\boldsymbol{\gamma}_k\rangle + \beta_{1k}\langle\mathbf{z}_{1,k}|\boldsymbol{\gamma}_k\rangle|\mathbf{z}_{1,k}\rangle + \beta_{2k}\langle\mathbf{z}_{2,k}|\boldsymbol{\gamma}_k\rangle|\mathbf{z}_{2,k}\rangle, \end{aligned} \tag{9.91}$$

where $\langle\mathbf{z}_{1,k}|\boldsymbol{\gamma}_k\rangle$ and $\langle\mathbf{z}_{2,k}|\boldsymbol{\gamma}_k\rangle$ can be identified as scalars. Although the vectors $|\mathbf{z}_{1,k}\rangle$ and $|\mathbf{z}_{2,k}\rangle$ in (9.90) are not unique, the following choices can be made to satisfy (9.90):

$$\begin{aligned} |\mathbf{z}_{1,k}\rangle &= |\boldsymbol{\delta}_k\rangle, \\ |\mathbf{z}_{2,k}\rangle &= \mathbf{S}_k|\boldsymbol{\gamma}_k\rangle, \\ \beta_{1k} &= \frac{1}{\langle\mathbf{z}_{1,k}|\boldsymbol{\gamma}_k\rangle}, \\ \beta_{2k} &= -\frac{1}{\langle\mathbf{z}_{2,k}|\boldsymbol{\gamma}_k\rangle}. \end{aligned} \tag{9.92}$$

So, the correction matrix (9.89) becomes

$$\begin{aligned} \mathbf{C}_k &= \frac{1}{\langle\mathbf{z}_{1,k}|\boldsymbol{\gamma}_k\rangle}|\boldsymbol{\delta}_k\rangle\langle\mathbf{z}_{1,k}| - \frac{1}{\langle\mathbf{z}_{2,k}|\boldsymbol{\gamma}_k\rangle}\mathbf{S}_k|\boldsymbol{\gamma}_k\rangle\langle\mathbf{z}_{2,k}| \\ &= \frac{|\boldsymbol{\delta}_k\rangle\langle\boldsymbol{\delta}_k|}{\langle\boldsymbol{\delta}_k|\boldsymbol{\gamma}_k\rangle} - \frac{(\mathbf{S}_k|\boldsymbol{\gamma}_k\rangle)(\mathbf{S}_k|\boldsymbol{\gamma}_k\rangle)^T}{(\mathbf{S}_k|\boldsymbol{\gamma}_k\rangle)^T|\boldsymbol{\gamma}_k\rangle} \\ &= \frac{|\boldsymbol{\delta}_k\rangle\langle\boldsymbol{\delta}_k|}{\langle\boldsymbol{\delta}_k|\boldsymbol{\gamma}_k\rangle} - \frac{\mathbf{S}_k|\boldsymbol{\gamma}_k\rangle\langle\boldsymbol{\gamma}_k|\mathbf{S}_k}{\langle\boldsymbol{\gamma}_k|\mathbf{S}_k|\boldsymbol{\gamma}_k\rangle} \\ &= \frac{|\boldsymbol{\delta}_k\rangle\langle\boldsymbol{\delta}_k|}{\langle\boldsymbol{\delta}_k|\boldsymbol{\gamma}_k\rangle} - \frac{|\mathbf{S}_k\boldsymbol{\gamma}_k\rangle\langle\mathbf{S}_k\boldsymbol{\gamma}_k|}{\langle\boldsymbol{\gamma}_k|\mathbf{S}_k|\boldsymbol{\gamma}_k\rangle}, \end{aligned} \tag{9.93}$$

and thus, the rank two update formula (9.90) can be expressed as

$$\begin{aligned} \mathbf{S}_{k+1} &= \mathbf{S}_k + \mathbf{C}_k \\ &= \mathbf{S}_k + \frac{|\boldsymbol{\delta}_k\rangle\langle\boldsymbol{\delta}_k|}{\langle\boldsymbol{\delta}_k|\boldsymbol{\gamma}_k\rangle} - \frac{|\mathbf{S}_k\boldsymbol{\gamma}_k\rangle\langle\mathbf{S}_k\boldsymbol{\gamma}_k|}{\langle\boldsymbol{\gamma}_k|\mathbf{S}_k|\boldsymbol{\gamma}_k\rangle}, \end{aligned} \tag{9.94}$$

where the correction is an $n \times n$ symmetric matrix of rank two. This equation is known as the DFP formula. The validity of (9.93) can be demonstrated by post-multiplying both sides by $|\boldsymbol{\gamma}_k\rangle$, that is,

$$\mathbf{S}_{k+1}|\boldsymbol{\gamma}_k\rangle = \mathbf{S}_k|\boldsymbol{\gamma}_k\rangle + \frac{|\boldsymbol{\delta}_k\rangle\langle\boldsymbol{\delta}_k|\boldsymbol{\gamma}_k\rangle}{\langle\boldsymbol{\delta}_k|\boldsymbol{\gamma}_k\rangle} - \frac{\mathbf{S}_k|\boldsymbol{\gamma}_k\rangle\langle\boldsymbol{\gamma}_k|\mathbf{S}_k|\boldsymbol{\gamma}_k\rangle}{\langle\boldsymbol{\gamma}_k|\mathbf{S}_k|\boldsymbol{\gamma}_k\rangle}.$$

Since $\langle\boldsymbol{\delta}_k|\boldsymbol{\gamma}_k\rangle$ and $\langle\boldsymbol{\gamma}_k|\mathbf{S}_k|\boldsymbol{\gamma}_k\rangle$ are scalars, we have

$$\mathbf{S}_{k+1}|\boldsymbol{\gamma}_k\rangle = |\boldsymbol{\delta}_k\rangle,$$

as required.

Using the iteration formula $|\mathbf{x}_{k+1}\rangle = |\mathbf{x}_k\rangle + \alpha_k|\mathbf{d}_k\rangle$, such that $|\mathbf{d}_k\rangle$ related to $|\boldsymbol{\delta}_k\rangle = |\mathbf{x}_{k+1}\rangle - |\mathbf{x}_k\rangle$ by the formula $|\boldsymbol{\delta}_k\rangle = \alpha_k|\mathbf{d}_k\rangle$. Thus (9.94) can be expressed as

$$\mathbf{S}_{k+1} = \mathbf{S}_k + \alpha_k\frac{|\mathbf{d}_k\rangle\langle\mathbf{d}_k|}{\langle\mathbf{d}_k|\boldsymbol{\gamma}_k\rangle} - \frac{|\mathbf{S}_k\boldsymbol{\gamma}_k\rangle\langle\mathbf{S}_k\boldsymbol{\gamma}_k|}{\langle\boldsymbol{\gamma}_k|\mathbf{S}_k|\boldsymbol{\gamma}_k\rangle}. \tag{9.95}$$

Note that, (9.83) and (9.93) are known as inverse update formulas since these equations approximate the inverse of the Hessian matrix of $f$.

***Example 9.9***

Locate the minimizer of

$$f|\mathbf{x}\rangle = \langle\mathbf{b}|\mathbf{x}\rangle + \frac{1}{2}\langle\mathbf{x}|\mathbf{H}|\mathbf{x}\rangle, \ \ |\mathbf{x}\rangle \in \mathbb{R}^2,$$

with $\mathbf{H} = \begin{pmatrix} 4 & 1 \\ 1 & 2 \end{pmatrix}$ and $|\mathbf{b}\rangle = \begin{pmatrix} 1 \\ -1 \end{pmatrix}$. Use the initial point $|\mathbf{x}_0\rangle = \begin{pmatrix} 0 \\ 0 \end{pmatrix}$ and $\mathbf{S}_0 = \mathbf{I}_2$.

***Solution***

Compute the gradient $|\mathbf{g}_0\rangle$,

$$|\mathbf{g}_0\rangle = |\mathbf{b}\rangle + \mathbf{H}|\mathbf{x}_0\rangle = \begin{pmatrix} 1 \\ -1 \end{pmatrix}.$$

It is a nonzero vector, so we proceed with the first iteration.

**Iteration 1.** Given that $\mathbf{S}_0 = \mathbf{I}_2$. Then,

$$|\mathbf{d}_0\rangle = -\mathbf{S}_0|\mathbf{g}_0\rangle = \begin{pmatrix} -1 \\ 1 \end{pmatrix}.$$

Because $f$ is a quadratic function,

$$\alpha_0 = \arg\min_{\alpha \geq 0} f|\mathbf{x}_0 + \alpha\mathbf{d}_0\rangle = -\frac{\langle\mathbf{g}_0|\mathbf{d}_0\rangle}{\langle\mathbf{d}_0|\mathbf{H}|\mathbf{d}_0\rangle} = \frac{1}{2}.$$

Therefore,

$$|\mathbf{x}_1\rangle = |\mathbf{x}_0\rangle + \alpha_0|\mathbf{d}_0\rangle = \begin{pmatrix} -1/2 \\ 1/2 \end{pmatrix}.$$

We then compute the gradient $|\mathbf{g}_1\rangle$

$$|\mathbf{g}_1\rangle = |\mathbf{b}\rangle + \mathbf{H}|\mathbf{x}_1\rangle = \begin{pmatrix} -1/2 \\ -1/2 \end{pmatrix}.$$

It is a nonzero vector, so we proceed with the second iteration.

**Iteration 2.** Compute

$$|\boldsymbol{\delta}_0\rangle = |\mathbf{x}_1\rangle - |\mathbf{x}_0\rangle = \begin{pmatrix} -1/2 \\ 1/2 \end{pmatrix},$$

$$|\boldsymbol{\gamma}_0\rangle = |\mathbf{g}_1\rangle - |\mathbf{g}_0\rangle = \begin{pmatrix} -3/2 \\ 1/2 \end{pmatrix}.$$

Using the above, we determine,

$$|\boldsymbol{\delta}_0\rangle\langle\boldsymbol{\delta}_0| = \begin{pmatrix} 1/4 & -1/4 \\ -1/4 & 1/4 \end{pmatrix},$$

$$\langle\boldsymbol{\delta}_0|\boldsymbol{\gamma}_0\rangle = 1,$$

$$\mathbf{S}_0|\boldsymbol{\gamma}_0\rangle = \begin{pmatrix} -3/2 \\ 1/2 \end{pmatrix}.$$

Thus,

$$\mathbf{S}_0|\boldsymbol{\gamma}_0\rangle\langle\boldsymbol{\gamma}_0|\mathbf{S}_0 = \begin{pmatrix} -3/2 \\ 1/2 \end{pmatrix}(-3/2, 1/2) = \begin{pmatrix} 9/4 & -3/4 \\ -3/4 & 1/4 \end{pmatrix},$$

$$\langle\boldsymbol{\gamma}_0|\mathbf{S}_0|\boldsymbol{\gamma}_0\rangle = (-3/2, 1/2)\begin{pmatrix} 1 & 0 \\ 0 & 1 \end{pmatrix}\begin{pmatrix} -3/2 \\ 1/2 \end{pmatrix} = 2.5.$$

Using the above, we now compute $\mathbf{S}_1$:

$$\mathbf{S}_1 = \mathbf{S}_0 + \frac{|\boldsymbol{\delta}_0\rangle\langle\boldsymbol{\delta}_0|}{\langle\boldsymbol{\delta}_0|\boldsymbol{\gamma}_0\rangle} - \frac{\mathbf{S}_0|\boldsymbol{\gamma}_0\rangle\langle\boldsymbol{\gamma}_0|\mathbf{S}_0}{\langle\boldsymbol{\gamma}_0|\mathbf{S}_0|\boldsymbol{\gamma}_0\rangle} = \begin{pmatrix} 0.35 & 0.05 \\ 0.05 & 1.15 \end{pmatrix}.$$

We now compute $|\mathbf{d}_1\rangle = -\mathbf{S}_1|\mathbf{g}_1\rangle = \begin{pmatrix} 0.2 \\ 0.6 \end{pmatrix}$ and

$$\alpha_1 = \arg\min_{\alpha \geq 0} f|\mathbf{x}_1 + \alpha\mathbf{d}_1\rangle = -\frac{\langle\mathbf{g}_1|\mathbf{d}_1\rangle}{\langle\mathbf{d}_1|\mathbf{H}|\mathbf{d}_1\rangle} = 0.357143.$$

Hence, the new update

$$|\mathbf{x}_2\rangle = |\mathbf{x}_1\rangle + \alpha_1|\mathbf{d}_1\rangle = \begin{pmatrix} -0.428571 \\ 0.714286 \end{pmatrix} = |\mathbf{x}^*\rangle,$$

because $f$ is a quadratic function of two variables.

Note that we have $\langle\mathbf{d}_0|\mathbf{H}|\mathbf{d}_1\rangle = \langle\mathbf{d}_1|\mathbf{H}|\mathbf{d}_0\rangle = 0$; that is, $|\mathbf{d}_0\rangle$ and $|\mathbf{d}_1\rangle$ are $\mathbf{H}$-conjugate directions.

#### 9.6.1.1. Convergence Analysis of the DFP Method

**Theorem 9.7:** If the DFP algorithm is applied to a quadratic function with the Hessian $\mathbf{H} = \mathbf{H}^T > 0$, then we have

$$\mathbf{S}_{k+1}|\boldsymbol{\gamma}_i\rangle = |\boldsymbol{\delta}_i\rangle,\ 0 \le i \le k. \tag{9.96}$$

**Proof:**

We use mathematical induction to prove this theorem. For $k = 0$, we have

$$\begin{aligned}\mathbf{S}_1|\boldsymbol{\gamma}_0\rangle &= \mathbf{S}_0|\boldsymbol{\gamma}_0\rangle + \frac{|\boldsymbol{\delta}_0\rangle\langle\boldsymbol{\delta}_0|}{\langle\boldsymbol{\delta}_0|\boldsymbol{\gamma}_0\rangle}|\boldsymbol{\gamma}_0\rangle - \frac{\mathbf{S}_0|\boldsymbol{\gamma}_0\rangle\langle\boldsymbol{\gamma}_0|\mathbf{S}_0}{\langle\boldsymbol{\gamma}_0|\mathbf{S}_0|\boldsymbol{\gamma}_0\rangle}|\boldsymbol{\gamma}_0\rangle\\ &= \mathbf{S}_0|\boldsymbol{\gamma}_0\rangle + |\boldsymbol{\delta}_0\rangle\frac{\langle\boldsymbol{\delta}_0|\boldsymbol{\gamma}_0\rangle}{\langle\boldsymbol{\delta}_0|\boldsymbol{\gamma}_0\rangle} - \mathbf{S}_0|\boldsymbol{\gamma}_0\rangle\frac{\langle\boldsymbol{\gamma}_0|\mathbf{S}_0|\boldsymbol{\gamma}_0\rangle}{\langle\boldsymbol{\gamma}_0|\mathbf{S}_0|\boldsymbol{\gamma}_0\rangle}\\ &= \mathbf{S}_0|\boldsymbol{\gamma}_0\rangle + |\boldsymbol{\delta}_0\rangle - \mathbf{S}_0|\boldsymbol{\gamma}_0\rangle.\end{aligned}$$

Therefore,

$$\mathbf{S}_1|\boldsymbol{\gamma}_0\rangle = |\boldsymbol{\delta}_0\rangle.$$

Assume that the result is true for $k-1$, that is, $\mathbf{S}_k|\boldsymbol{\gamma}_i\rangle = |\boldsymbol{\delta}_i\rangle$, where $0 \le i \le k-1$. Now we prove that the result is true for $k$, i.e., $\mathbf{S}_{k+1}|\boldsymbol{\gamma}_i\rangle = |\boldsymbol{\delta}_i\rangle$, where $0 \le i \le k$. First, we take $i = k$:

$$\begin{aligned}\mathbf{S}_{k+1}|\boldsymbol{\gamma}_k\rangle &= \mathbf{S}_k|\boldsymbol{\gamma}_k\rangle + \frac{|\boldsymbol{\delta}_k\rangle\langle\boldsymbol{\delta}_k|}{\langle\boldsymbol{\delta}_k|\boldsymbol{\gamma}_k\rangle}|\boldsymbol{\gamma}_k\rangle - \frac{\mathbf{S}_k|\boldsymbol{\gamma}_k\rangle\langle\boldsymbol{\gamma}_k|\mathbf{S}_k}{\langle\boldsymbol{\gamma}_k|\mathbf{S}_k|\boldsymbol{\gamma}_k\rangle}|\boldsymbol{\gamma}_k\rangle\\ &= \mathbf{S}_k|\boldsymbol{\gamma}_k\rangle + |\boldsymbol{\delta}_k\rangle\frac{\langle\boldsymbol{\delta}_k|\boldsymbol{\gamma}_k\rangle}{\langle\boldsymbol{\delta}_k|\boldsymbol{\gamma}_k\rangle} - \mathbf{S}_k|\boldsymbol{\gamma}_k\rangle\frac{\langle\boldsymbol{\gamma}_k|\mathbf{S}_k|\boldsymbol{\gamma}_k\rangle}{\langle\boldsymbol{\gamma}_k|\mathbf{S}_k|\boldsymbol{\gamma}_k\rangle}\\ &= \mathbf{S}_k|\boldsymbol{\gamma}_k\rangle + |\boldsymbol{\delta}_k\rangle - \mathbf{S}_k|\boldsymbol{\gamma}_k\rangle\\ &= |\boldsymbol{\delta}_k\rangle.\end{aligned}$$

For $i < k$:

$$\begin{aligned}\mathbf{S}_{k+1}|\boldsymbol{\gamma}_i\rangle &= \mathbf{S}_k|\boldsymbol{\gamma}_i\rangle + \frac{|\boldsymbol{\delta}_k\rangle\langle\boldsymbol{\delta}_k|}{\langle\boldsymbol{\delta}_k|\boldsymbol{\gamma}_k\rangle}|\boldsymbol{\gamma}_i\rangle - \frac{\mathbf{S}_k|\boldsymbol{\gamma}_k\rangle\langle\boldsymbol{\gamma}_k|\mathbf{S}_k}{\langle\boldsymbol{\gamma}_k|\mathbf{S}_k|\boldsymbol{\gamma}_k\rangle}|\boldsymbol{\gamma}_i\rangle\\ &= |\boldsymbol{\delta}_i\rangle + |\boldsymbol{\delta}_k\rangle\frac{\langle\boldsymbol{\delta}_k|\boldsymbol{\gamma}_i\rangle}{\langle\boldsymbol{\delta}_k|\boldsymbol{\gamma}_k\rangle} - \mathbf{S}_k|\boldsymbol{\gamma}_k\rangle\frac{\langle\boldsymbol{\gamma}_k|\mathbf{S}_k|\boldsymbol{\gamma}_i\rangle}{\langle\boldsymbol{\gamma}_k|\mathbf{S}_k|\boldsymbol{\gamma}_k\rangle}\\ &= |\boldsymbol{\delta}_i\rangle + |\boldsymbol{\delta}_k\rangle\frac{\langle\boldsymbol{\delta}_k|\boldsymbol{\gamma}_i\rangle}{\langle\boldsymbol{\delta}_k|\boldsymbol{\gamma}_k\rangle} - \mathbf{S}_k|\boldsymbol{\gamma}_k\rangle\frac{\langle\boldsymbol{\gamma}_k|\boldsymbol{\delta}_i\rangle}{\langle\boldsymbol{\gamma}_k|\mathbf{S}_k|\boldsymbol{\gamma}_k\rangle}.\end{aligned}$$

Since

$$\langle\boldsymbol{\delta}_k|\boldsymbol{\gamma}_i\rangle = \langle\boldsymbol{\delta}_k|\mathbf{H}|\boldsymbol{\delta}_i\rangle = \langle\alpha_k\mathbf{d}_k|\mathbf{H}|\alpha_i\mathbf{d}_i\rangle = \alpha_i\alpha_k\langle\mathbf{d}_k|\mathbf{H}|\mathbf{d}_i\rangle = 0,$$

and

$$\langle\boldsymbol{\gamma}_k|\boldsymbol{\delta}_i\rangle = \langle\mathbf{H}\boldsymbol{\delta}_k|\boldsymbol{\delta}_i\rangle = \langle\boldsymbol{\delta}_k|\mathbf{H}|\boldsymbol{\delta}_i\rangle = 0,$$

then we get

$$\mathbf{S}_{k+1}|\boldsymbol{\gamma}_i\rangle = |\boldsymbol{\delta}_i\rangle.$$

This completes the proof.

■

From Theorems 9.5 and 9.7, we conclude that the DFP algorithm is a conjugate direction algorithm, that is if the DFP algorithm is applied to a quadratic with Hessian $\mathbf{H} = \mathbf{H}^T > 0$ such that $\mathbf{S}_{k+1}|\boldsymbol{\gamma}_i\rangle = |\boldsymbol{\delta}_i\rangle$, $0 \le i \le k$ with $\mathbf{S}_{k+1} = \mathbf{S}_{k+1}^T$ and $\alpha_i \ne 0$, where $0 \le i \le k$, then $|\mathbf{d}_0\rangle, |\mathbf{d}_1\rangle, \dots, |\mathbf{d}_{n-1}\rangle$ are $\mathbf{H}$-conjugate.

**Theorem 9.8:** Suppose that $|\mathbf{g}_k\rangle \ne |\mathbf{0}\rangle$. In the DFP algorithm if $\mathbf{S}_k$ is positive definite, then $\mathbf{S}_{k+1}$ is also positive definite.

**Proof:**

We have to show that

$$\langle\mathbf{x}|\mathbf{S}_{k+1}|\mathbf{x}\rangle > 0,$$

for all $|\mathbf{x}\rangle \ne 0$. For this, we have to write the following quadratic form:

$$\begin{aligned}\langle\mathbf{x}|\mathbf{S}_{k+1}|\mathbf{x}\rangle &= \langle\mathbf{x}|\mathbf{S}_k|\mathbf{x}\rangle + \langle\mathbf{x}|\frac{|\boldsymbol{\delta}_k\rangle\langle\boldsymbol{\delta}_k|}{\langle\boldsymbol{\delta}_k|\boldsymbol{\gamma}_k\rangle}|\mathbf{x}\rangle - \langle\mathbf{x}|\frac{\mathbf{S}_k|\boldsymbol{\gamma}_k\rangle\langle\boldsymbol{\gamma}_k|\mathbf{S}_k}{\langle\boldsymbol{\gamma}_k|\mathbf{S}_k|\boldsymbol{\gamma}_k\rangle}|\mathbf{x}\rangle \\ &= \langle\mathbf{x}|\mathbf{S}_k|\mathbf{x}\rangle + \frac{\langle\mathbf{x}|\boldsymbol{\delta}_k\rangle^2}{\langle\boldsymbol{\delta}_k|\boldsymbol{\gamma}_k\rangle} - \frac{\langle\mathbf{x}|\mathbf{S}_k|\boldsymbol{\gamma}_k\rangle\langle\boldsymbol{\gamma}_k|\mathbf{S}_k|\mathbf{x}\rangle}{\langle\boldsymbol{\gamma}_k|\mathbf{S}_k|\boldsymbol{\gamma}_k\rangle}.\end{aligned}$$

Therefore,

$$\langle\mathbf{x}|\mathbf{S}_{k+1}|\mathbf{x}\rangle = \langle\mathbf{x}|\mathbf{S}_k|\mathbf{x}\rangle + \frac{\langle\mathbf{x}|\boldsymbol{\delta}_k\rangle^2}{\langle\boldsymbol{\delta}_k|\boldsymbol{\gamma}_k\rangle} - \frac{\langle\mathbf{x}|\mathbf{S}_k|\boldsymbol{\gamma}_k\rangle^2}{\langle\boldsymbol{\gamma}_k|\mathbf{S}_k|\boldsymbol{\gamma}_k\rangle}. \tag{9.97}$$

Set $|\mathbf{a}\rangle \equiv \mathbf{S}_k^{1/2}|\mathbf{x}\rangle$, and $|\mathbf{b}\rangle \equiv \mathbf{S}_k^{1/2}|\boldsymbol{\gamma}_k\rangle$, where $\mathbf{S}_k = \mathbf{S}_k^{1/2}\mathbf{S}_k^{1/2}$. Since $\mathbf{S}_k > 0$, and $\mathbf{S}_k^{1/2}$ is well defined, then

$$\begin{aligned}\langle\mathbf{x}|\mathbf{S}_k|\mathbf{x}\rangle &= \left\langle\mathbf{x}\middle|\mathbf{S}_k^{1/2}\mathbf{S}_k^{1/2}\middle|\mathbf{x}\right\rangle = \left\langle\mathbf{S}_k^{1/2}\mathbf{x}\middle|\mathbf{S}_k^{1/2}\mathbf{x}\right\rangle = \langle\mathbf{a}|\mathbf{a}\rangle, \\ \langle\mathbf{x}|\mathbf{S}_k|\boldsymbol{\gamma}_k\rangle &= \left\langle\mathbf{x}\middle|\mathbf{S}_k^{1/2}\mathbf{S}_k^{1/2}\middle|\boldsymbol{\gamma}_k\right\rangle = \left\langle\mathbf{S}_k^{1/2}\mathbf{x}\middle|\mathbf{S}_k^{1/2}\boldsymbol{\gamma}_k\right\rangle = \langle\mathbf{a}|\mathbf{b}\rangle,\end{aligned}$$

and

$$\langle\boldsymbol{\gamma}_k|\mathbf{S}_k|\boldsymbol{\gamma}_k\rangle = \left\langle\boldsymbol{\gamma}_k\middle|\mathbf{S}_k^{1/2}\mathbf{S}_k^{1/2}\middle|\boldsymbol{\gamma}_k\right\rangle = \left\langle\mathbf{S}_k^{1/2}\boldsymbol{\gamma}_k\middle|\mathbf{S}_k^{1/2}\boldsymbol{\gamma}_k\right\rangle = \langle\mathbf{b}|\mathbf{b}\rangle.$$

Therefore, (9.97) is reduced in the following form:

$$\begin{aligned}\langle\mathbf{x}|\mathbf{S}_{k+1}|\mathbf{x}\rangle &= \langle\mathbf{a}|\mathbf{a}\rangle + \frac{\langle\mathbf{x}|\boldsymbol{\delta}_k\rangle^2}{\langle\boldsymbol{\delta}_k|\boldsymbol{\gamma}_k\rangle} - \frac{\langle\mathbf{a}|\mathbf{b}\rangle^2}{\langle\mathbf{b}|\mathbf{b}\rangle} \\ &= \|\mathbf{a}\|^2 - \frac{\langle\mathbf{a}|\mathbf{b}\rangle^2}{\|\mathbf{b}\|^2} + \frac{\langle\mathbf{x}|\boldsymbol{\delta}_k\rangle^2}{\langle\boldsymbol{\delta}_k|\boldsymbol{\gamma}_k\rangle} \\ &= \frac{\|\mathbf{a}\|^2\|\mathbf{b}\|^2 - \langle\mathbf{a}|\mathbf{b}\rangle^2}{\|\mathbf{b}\|^2} + \frac{\langle\mathbf{x}|\boldsymbol{\delta}_k\rangle^2}{\langle\boldsymbol{\delta}_k|\boldsymbol{\gamma}_k\rangle}.\end{aligned} \tag{9.98}$$

Note that

$$\begin{aligned}\langle\boldsymbol{\delta}_k|\boldsymbol{\gamma}_k\rangle &= \langle\boldsymbol{\delta}_k|\mathbf{g}_{k+1} - \mathbf{g}_k\rangle \\ &= \langle\boldsymbol{\delta}_k|\mathbf{g}_{k+1}\rangle - \langle\boldsymbol{\delta}_k|\mathbf{g}_k\rangle = \alpha_k\langle\mathbf{d}_k|\mathbf{g}_{k+1}\rangle - \langle\boldsymbol{\delta}_k|\mathbf{g}_k\rangle.\end{aligned}$$

Since $\langle\mathbf{d}_k|\mathbf{g}_{k+1}\rangle = 0$, then

$$\langle\boldsymbol{\delta}_k|\boldsymbol{\gamma}_k\rangle = -\langle\boldsymbol{\delta}_k|\mathbf{g}_k\rangle.$$

Since $|\boldsymbol{\delta}_k\rangle = \alpha_k|\mathbf{d}_k\rangle = -\alpha_k\mathbf{S}_k|\mathbf{g}_k\rangle$, then

$$\langle\boldsymbol{\delta}_k|\boldsymbol{\gamma}_k\rangle = \alpha_k\langle\mathbf{g}_k|\mathbf{S}_k|\mathbf{g}_k\rangle. \tag{9.99}$$

Using (9.99) in (9.98), we obtain

$$\langle\mathbf{x}|\mathbf{S}_{k+1}|\mathbf{x}\rangle = \frac{\|\mathbf{a}\|^2\|\mathbf{b}\|^2 - \langle\mathbf{a}|\mathbf{b}\rangle^2}{\|\mathbf{b}\|^2} + \frac{\langle\mathbf{x}|\boldsymbol{\delta}_k\rangle^2}{\alpha_k\langle\mathbf{g}_k|\mathbf{S}_k|\mathbf{g}_k\rangle}. \tag{9.100}$$

In (9.100), both terms of RHS are nonnegative. The first term is nonnegative due to the Cauchy-Schwarz inequality

$$\langle\mathbf{x}|\mathbf{y}\rangle \le \|\mathbf{x}\|\|\mathbf{y}\|, \quad \forall|\mathbf{x}\rangle, |\mathbf{y}\rangle \in \mathbb{R}^n,$$

holds if $|\mathbf{x}\rangle = \alpha|\mathbf{y}\rangle$, where $\alpha \in \mathbb{R}$. The second term is nonnegative due to the fact that $\mathbf{S}_k$ is positive definite and $\alpha_k > 0$. Therefore, in order to show that

$$\langle\mathbf{x}|\mathbf{S}_{k+1}|\mathbf{x}\rangle > 0, \quad |\mathbf{x}\rangle \ne |\mathbf{0}\rangle,$$

we need to show that the first term is zero only if $|\mathbf{a}\rangle$ and $|\mathbf{b}\rangle$ are proportional, that is, if $|\mathbf{a}\rangle = \beta|\mathbf{b}\rangle$ for some $\beta \in \mathbb{R}$. Thus, to complete the proof, it is enough to show that if

$$|\mathbf{a}\rangle = \beta|\mathbf{b}\rangle,$$

then

$$\frac{\langle\mathbf{x}|\boldsymbol{\delta}_k\rangle^2}{\alpha_k\langle\mathbf{g}_k|\mathbf{S}_k|\mathbf{g}_k\rangle} > 0.$$

Note that

$$\mathbf{S}_k^{1/2}|\mathbf{x}\rangle = |\mathbf{a}\rangle = \beta|\mathbf{b}\rangle = \beta\mathbf{S}_k^{1/2}|\boldsymbol{\gamma}_k\rangle = \mathbf{S}_k^{1/2}(\beta|\boldsymbol{\gamma}_k\rangle),$$

which yields

$$|\mathbf{x}\rangle = \beta|\boldsymbol{\gamma}_k\rangle.$$

Using value of $|\mathbf{x}\rangle$ in the second term of the right-hand side of (9.100), we get

$$\begin{aligned}\frac{\langle \mathbf{x}|\boldsymbol{\delta}_k\rangle^2}{\alpha_k\langle \mathbf{g}_k|\mathbf{S}_k|\mathbf{g}_k\rangle} &= \frac{\langle \beta\boldsymbol{\gamma}_k|\boldsymbol{\delta}_k\rangle^2}{\alpha_k\langle \mathbf{g}_k|\mathbf{S}_k|\mathbf{g}_k\rangle} \\ &= \beta^2\frac{\langle \boldsymbol{\gamma}_k|\boldsymbol{\delta}_k\rangle^2}{\alpha_k\langle \mathbf{g}_k|\mathbf{S}_k|\mathbf{g}_k\rangle} \\ &= \beta^2\frac{(\alpha_k\langle \mathbf{g}_k|\mathbf{S}_k|\mathbf{g}_k\rangle)^2}{\alpha_k\langle \mathbf{g}_k|\mathbf{S}_k|\mathbf{g}_k\rangle} \\ &= \beta^2\alpha_k\langle \mathbf{g}_k|\mathbf{S}_k|\mathbf{g}_k\rangle > 0.\end{aligned}$$

Thus, from (9.100), for all $|\mathbf{x}\rangle \neq |\mathbf{0}\rangle$, we get

$$\langle \mathbf{x}|\mathbf{S}_{k+1}|\mathbf{x}\rangle > 0.$$

■

### 9.6.2. The BFGS Method

Recall the updating formulas for the approximation of the inverse $\mathbf{H}^{-1}$ of the Hessian matrix which were based on the following equations:

$$\mathbf{S}_{k+1}|\boldsymbol{\gamma}_i\rangle = |\boldsymbol{\delta}_i\rangle,\ \ 0 \le i \le k. \tag{9.101}$$

We then formulated update formulas for the approximations to the inverse $\mathbf{H}^{-1}$ of the Hessian. It is possible to derive a family of direct update formulas in which approximations to the Hessian matrix $\mathbf{H}$ are considered instead of approximating $\mathbf{H}^{-1}$. To do this, let $\mathbf{E}_k$ be our estimate of $\mathbf{H}_k$ at the $k$th step. We need $\mathbf{E}_{k+1}$ to satisfy the following set of equations:

$$|\boldsymbol{\gamma}_i\rangle = \mathbf{E}_{k+1}|\boldsymbol{\delta}_i\rangle,\ \ 0 \le i \le k. \tag{9.102}$$

Note that, (9.101) and (9.102) are similar; the only difference is that $|\boldsymbol{\delta}_i\rangle$ and $|\boldsymbol{\gamma}_i\rangle$ are interchanged. Thus, given any update formula for $\mathbf{S}_k$, a corresponding update formula for $\mathbf{E}_k$ can be found by interchanging the roles of $\mathbf{E}_k$ and $\mathbf{S}_k$ and of $|\boldsymbol{\gamma}_k\rangle$ and $|\boldsymbol{\delta}_k\rangle$. In particular, the BFGS update [8,14-18] for $\mathbf{E}_k$ corresponding to the DFP update for $\mathbf{S}_k$. Formulas related in this way are said to be dual or complementary. Recall the DFP update for the approximation $\mathbf{S}_k$ of the inverse $\mathbf{H}^{-1}$ Hessian from (9.94) as

$$\mathbf{S}_{k+1}^{\mathrm{DFP}} = \mathbf{S}_k + \frac{|\boldsymbol{\delta}_k\rangle\langle\boldsymbol{\delta}_k|}{\langle\boldsymbol{\delta}_k|\boldsymbol{\gamma}_k\rangle} - \frac{|\mathbf{S}_k\boldsymbol{\gamma}_k\rangle\langle\mathbf{S}_k\boldsymbol{\gamma}_k|}{\langle\boldsymbol{\gamma}_k|\mathbf{S}_k|\boldsymbol{\gamma}_k\rangle}. \tag{9.103}$$

We apply the complementarity concept in the above equation. In other words, the procedure used in deriving (9.83) and (9.94) can be followed by using $\mathbf{E}_k$, $|\boldsymbol{\delta}_k\rangle$, and $|\boldsymbol{\gamma}_k\rangle$ in place of $\mathbf{S}_k$, $|\boldsymbol{\gamma}_k\rangle$, and $|\boldsymbol{\delta}_k\rangle$, respectively. This leads to the rank two update formula, similar to (9.94), known as the BFGS update of $\mathbf{E}_k$:

$$\mathbf{E}_{k+1} = \mathbf{E}_k + \frac{|\boldsymbol{\gamma}_k\rangle\langle\boldsymbol{\gamma}_k|}{\langle\boldsymbol{\gamma}_k|\boldsymbol{\delta}_k\rangle} - \frac{|\mathbf{E}_k\boldsymbol{\delta}_k\rangle\langle\mathbf{E}_k\boldsymbol{\delta}_k|}{\langle\boldsymbol{\delta}_k|\mathbf{E}_k|\boldsymbol{\delta}_k\rangle}, \tag{9.104}$$

that represents an update equation for the approximation $\mathbf{E}_k$ of the Hessian itself. In practical computations, (9.104) is rewritten more conveniently in terms of $\mathbf{S}_k$, using the following formulas

$$\mathbf{S}_{k+1}^{\mathrm{BFGS}} = (\mathbf{E}_{k+1})^{-1} = \left(\mathbf{E}_k + \frac{|\boldsymbol{\gamma}_k\rangle\langle\boldsymbol{\gamma}_k|}{\langle\boldsymbol{\gamma}_k|\boldsymbol{\delta}_k\rangle} - \frac{|\mathbf{E}_k\boldsymbol{\delta}_k\rangle\langle\mathbf{E}_k\boldsymbol{\delta}_k|}{\langle\boldsymbol{\delta}_k|\mathbf{E}_k|\boldsymbol{\delta}_k\rangle}\right)^{-1}, \tag{9.105}$$

or

$$\mathbf{S}_{k+1}^{\mathrm{BFGS}} = \mathbf{S}_k + \left(1 + \frac{\langle\boldsymbol{\gamma}_k|\mathbf{S}_k|\boldsymbol{\gamma}_k\rangle}{\langle\boldsymbol{\gamma}_k|\boldsymbol{\delta}_k\rangle}\right)\frac{|\boldsymbol{\delta}_k\rangle\langle\boldsymbol{\delta}_k|}{\langle\boldsymbol{\delta}_k|\boldsymbol{\gamma}_k\rangle} - \frac{|\mathbf{S}_k\boldsymbol{\gamma}_k\rangle\langle\boldsymbol{\delta}_k| + |\boldsymbol{\delta}_k\rangle\langle\mathbf{S}_k\boldsymbol{\gamma}_k|}{\langle\boldsymbol{\gamma}_k|\boldsymbol{\delta}_k\rangle}. \tag{9.106}$$

This represents the BFGS formula for updating $\mathbf{S}_k$. Numerical experience indicates that the BFGS method is the best unconstrained quasi-Newton method and is less influenced by errors in finding $\alpha_k$ compared to the DFP method.

***Example 9.10***

Use the BFGS method to minimize

$$f|\mathbf{x}\rangle = a + \langle \mathbf{b}|\mathbf{x}\rangle + \frac{1}{2}\langle \mathbf{x}|\mathbf{H}|\mathbf{x}\rangle, \ \ |\mathbf{x}\rangle \in \mathbb{R}^2,$$

where

$$\mathbf{H} = \begin{pmatrix} 5 & -1 \\ -1 & 4 \end{pmatrix}, \qquad |\mathbf{b}\rangle = \begin{pmatrix} 0 \\ -1 \end{pmatrix}, \qquad a = -1.$$

Use the initial point $|\mathbf{x}_0\rangle = \begin{pmatrix} 0 \\ 0 \end{pmatrix}$ and take $\mathbf{S}_0 = \mathbf{I}_2$. Verify that $\mathbf{S}_2^{\mathrm{BFGS}} = \mathbf{H}^{-1}$.

***Solution***

We compute the gradient $|\mathbf{g}_0\rangle$

$$|\mathbf{g}_0\rangle = |\mathbf{b}\rangle + \mathbf{H}|\mathbf{x}_0\rangle = \begin{pmatrix} 0 \\ -1 \end{pmatrix}.$$

It is a nonzero vector, so we proceed with the first iteration.

**Iteration 1.** Given that $\mathbf{S}_0 = \mathbf{I}_2$. Then,

$$|\mathbf{d}_0\rangle = -\mathbf{S}_0|\mathbf{g}_0\rangle = \begin{pmatrix} 0 \\ 1 \end{pmatrix}.$$

The objective function is quadratic, and hence we can use the following formula to compute $\alpha_0$:

$$\alpha_0 = \arg\min_{\alpha \geq 0} f|\mathbf{x}_0 + \alpha\mathbf{d}_0\rangle = -\frac{\langle \mathbf{g}_0|\mathbf{d}_0\rangle}{\langle \mathbf{d}_0|\mathbf{H}|\mathbf{d}_0\rangle} = \frac{1}{4}.$$

Therefore,

$$|\mathbf{x}_1\rangle = |\mathbf{x}_0\rangle + \alpha_0|\mathbf{d}_0\rangle = \begin{pmatrix} 0 \\ 1/4 \end{pmatrix}.$$

We then compute the gradient $|\mathbf{g}_1\rangle$

$$|\mathbf{g}_1\rangle = |\mathbf{b}\rangle + \mathbf{H}|\mathbf{x}_1\rangle = \begin{pmatrix} -1/4 \\ 0 \end{pmatrix}.$$

It is a nonzero vector, so we proceed with the second iteration.

**Iteration 2.** To compute $\mathbf{S}_1^{\mathrm{BFGS}}$, we need to evaluate the following quantities:

$$|\boldsymbol{\delta}_0\rangle = |\mathbf{x}_1\rangle - |\mathbf{x}_0\rangle = \begin{pmatrix} 0 \\ 1/4 \end{pmatrix},$$

$$|\boldsymbol{\gamma}_0\rangle = |\mathbf{g}_1\rangle - |\mathbf{g}_0\rangle = \begin{pmatrix} -1/4 \\ 1 \end{pmatrix}.$$

Using the above, we now compute $\mathbf{S}_1^{\mathrm{BFGS}}$:

$$\mathbf{S}_1^{\mathrm{BFGS}} = \mathbf{S}_0 + \left(1 + \frac{\langle \boldsymbol{\gamma}_0|\mathbf{S}_0|\boldsymbol{\gamma}_0\rangle}{\langle \boldsymbol{\gamma}_0|\boldsymbol{\delta}_0\rangle}\right)\frac{|\boldsymbol{\delta}_0\rangle\langle\boldsymbol{\delta}_0|}{\langle\boldsymbol{\delta}_0|\boldsymbol{\gamma}_0\rangle} - \frac{\mathbf{S}_0|\boldsymbol{\gamma}_0\rangle\langle\boldsymbol{\delta}_0| + |\boldsymbol{\delta}_0\rangle\langle\boldsymbol{\gamma}_0|\mathbf{S}_0}{\langle\boldsymbol{\gamma}_0|\boldsymbol{\delta}_0\rangle} = \begin{pmatrix} 1 & 0.25 \\ 0.25 & 0.3125 \end{pmatrix}.$$

We now compute $|\mathbf{d}_1\rangle$

$$|\mathbf{d}_1\rangle = -\mathbf{S}_1^{\mathrm{BFGS}}|\mathbf{g}_1\rangle = \begin{pmatrix} 0.25 \\ 0.0625 \end{pmatrix},$$

and

$$\alpha_1 = \arg\min_{\alpha \geq 0} f|\mathbf{x}_1 + \alpha\mathbf{d}_1\rangle = -\frac{\langle \mathbf{g}_1|\mathbf{d}_1\rangle}{\langle \mathbf{d}_1|\mathbf{H}|\mathbf{d}_1\rangle} = 0.210526.$$

Therefore, the new update

$$|\mathbf{x}_2\rangle = |\mathbf{x}_1\rangle + \alpha_1|\mathbf{d}_1\rangle = \begin{pmatrix} 0.052632 \\ 0.263158 \end{pmatrix}.$$

Because our objective function $f$ is a quadratic on $\mathbb{R}^2$, $|\mathbf{x}_2\rangle$ is the minimizer. Note that the gradient at $|\mathbf{x}_2\rangle$ is $|\mathbf{g}_2\rangle = |\mathbf{0}\rangle$.

To verify that $\mathbf{S}_2^{\mathrm{BFGS}} = \mathbf{H}^{-1}$, we compute

$$|\boldsymbol{\delta}_1\rangle = |\mathbf{x}_2\rangle - |\mathbf{x}_1\rangle = \begin{pmatrix} 0.0526316 \\ 0.0131579 \end{pmatrix},$$

$$|\boldsymbol{\gamma}_1\rangle = |\mathbf{g}_2\rangle - |\mathbf{g}_1\rangle = \begin{pmatrix} 1/4 \\ 0 \end{pmatrix}.$$

Hence,

$$\mathbf{S}_2^{\mathrm{BFGS}} = \mathbf{S}_1 + \left(1 + \frac{\langle \boldsymbol{\gamma}_1|\mathbf{S}_1|\boldsymbol{\gamma}_1\rangle}{\langle \boldsymbol{\gamma}_1|\boldsymbol{\delta}_1\rangle}\right)\frac{|\boldsymbol{\delta}_1\rangle\langle\boldsymbol{\delta}_1|}{\langle\boldsymbol{\delta}_1|\boldsymbol{\gamma}_1\rangle} - \frac{\mathbf{S}_1|\boldsymbol{\gamma}_1\rangle\langle\boldsymbol{\delta}_1| + |\boldsymbol{\delta}_1\rangle\langle\boldsymbol{\gamma}_1|\mathbf{S}_1}{\langle\boldsymbol{\gamma}_1|\boldsymbol{\delta}_1\rangle} = \begin{pmatrix} 0.210526 & 0.052632 \\ 0.052632 & 0.263158 \end{pmatrix}.$$

Note that indeed $\mathbf{S}_2^{\mathrm{BFGS}}\mathbf{H} = \mathbf{H}\mathbf{S}_2^{\mathrm{BFGS}} = \mathbf{I}_2$, and hence $\mathbf{S}_2^{\mathrm{BFGS}} = \mathbf{H}^{-1}$.

## 9.7 Other Families of Quasi-Newton Algorithms

### 9.7.1. The Huang Family

This is a more quasi-Newton general family [19] which encompasses the rank one, DFP, BFGS, as well as some other formulas. It is of the form

$$\mathbf{S}_{k+1} = \mathbf{S}_k + \frac{\theta|\boldsymbol{\delta}_k\rangle\langle\boldsymbol{\delta}_k| + \phi|\boldsymbol{\delta}_k\rangle\langle\boldsymbol{\gamma}_k|\mathbf{S}_k}{\theta\langle\boldsymbol{\delta}_k|\boldsymbol{\gamma}_k\rangle + \phi\langle\boldsymbol{\gamma}_k|\mathbf{S}_k|\boldsymbol{\gamma}_k\rangle} - \frac{\psi\mathbf{S}_k|\boldsymbol{\gamma}_k\rangle\langle\boldsymbol{\delta}_k| + \omega\mathbf{S}_k|\boldsymbol{\gamma}_k\rangle\langle\boldsymbol{\gamma}_k|\mathbf{S}_k}{\psi\langle\boldsymbol{\delta}_k|\boldsymbol{\gamma}_k\rangle + \omega\langle\boldsymbol{\gamma}_k|\mathbf{S}_k|\boldsymbol{\gamma}_k\rangle}, \tag{9.107}$$

where $\theta$, $\phi$, $\psi$, and $\omega$ are independent parameters. Different choices of $\theta$, $\phi$, $\psi$, and $\omega$ in (9.107) lead to different algorithms. The formulas that can be generated from the Huang formula are given in Table 9.4.

**Table 9.4.** The Huang family

| Formula | Parameters |
|---|---|
| Rank one | $\theta = 1$, $\phi = -1$, $\psi = 1$, and $\omega = -1$, |
| DFP | $\theta = 1$, $\phi = 0$, $\psi = 0$, and $\omega = 1$, |
| BFGS | $\phi/\theta = -\langle\boldsymbol{\delta}_k|\boldsymbol{\gamma}_k\rangle/(\langle\boldsymbol{\delta}_k|\boldsymbol{\gamma}_k\rangle + \langle\boldsymbol{\gamma}_k|\mathbf{S}_k|\boldsymbol{\gamma}_k\rangle)$, $\psi = 1$, and $\omega = 0$, |
| McCormick | $\theta = 1$, $\phi = 0$, $\psi = 1$, and $\omega = 0$, |
| Pearson | $\theta = 0$, $\phi = 1$, $\psi = 0$, and $\omega = 1$. |

### 9.7.2. Hoshino Method

Like the DFP and BFGS formulas, the Hoshino formula [20] is of rank two. It is given by

$$\mathbf{S}_{k+1} = \mathbf{S}_k + \theta_k|\boldsymbol{\delta}_k\rangle\langle\boldsymbol{\delta}_k| - \psi_k(|\boldsymbol{\delta}_k\rangle\langle\boldsymbol{\gamma}_k|\mathbf{S}_k + \mathbf{S}_k|\boldsymbol{\gamma}_k\rangle\langle\boldsymbol{\delta}_k| + \mathbf{S}_k|\boldsymbol{\gamma}_k\rangle\langle\boldsymbol{\gamma}_k|\mathbf{S}_k), \tag{9.108}$$

where

$$\theta_k = \frac{\langle\boldsymbol{\gamma}_k|\boldsymbol{\delta}_k\rangle + 2\langle\boldsymbol{\gamma}_k|\mathbf{S}_k|\boldsymbol{\gamma}_k\rangle}{\langle\boldsymbol{\gamma}_k|\boldsymbol{\delta}_k\rangle^2 + \langle\boldsymbol{\gamma}_k|\boldsymbol{\delta}_k\rangle\langle\boldsymbol{\gamma}_k|\mathbf{S}_k|\boldsymbol{\gamma}_k\rangle} \text{ and } \psi_k = \frac{1}{\langle\boldsymbol{\gamma}_k|\boldsymbol{\delta}_k\rangle + \langle\boldsymbol{\gamma}_k|\mathbf{S}_k|\boldsymbol{\gamma}_k\rangle}. \tag{9.109}$$

The inverse of $\mathbf{S}_{k+1}$, designated as $\mathbf{E}_{k+1}$, can be obtained by applying the Sherman-Morrison formula.

### 9.7.3. The Broyden Family

This formula entails an independent parameter $\phi_k$ and is given by

$$\mathbf{S}_{k+1} = (1 - \phi_k)\mathbf{S}_{k+1}^{\text{DFP}} + \phi_k\mathbf{S}_{k+1}^{\text{BFGS}}, \tag{9.110}$$

which is known as the Broyden formula [8]. Evidently, if $\phi_k = 1$ the Broyden formula reduces to the BFGS formula, and when $\phi_k = 0$, the Broyden formula reduces to the DFP formula, and if

$$\phi_k = \frac{\langle\boldsymbol{\delta}_k|\boldsymbol{\gamma}_k\rangle}{\langle\boldsymbol{\delta}_k|\boldsymbol{\gamma}_k\rangle \pm \langle\boldsymbol{\gamma}_k|\mathbf{S}_k|\boldsymbol{\gamma}_k\rangle}, \tag{9.111}$$

the rank one or Hoshino formula is obtained.

# CHAPTER 10

# MULTI-VARIABLE CONSTRAINED OPTIMIZATION

## 10.1 Constrained Optimization Problems

In this chapter, equality-constrained optimization problems are discussed. Almost all decision-making problems aim to minimize or maximize a function and simultaneously have a requirement for satisfying some constraints. Applying constraints to a problem can affect the solution, but this need not be the case, as shown in Figure 10.1.

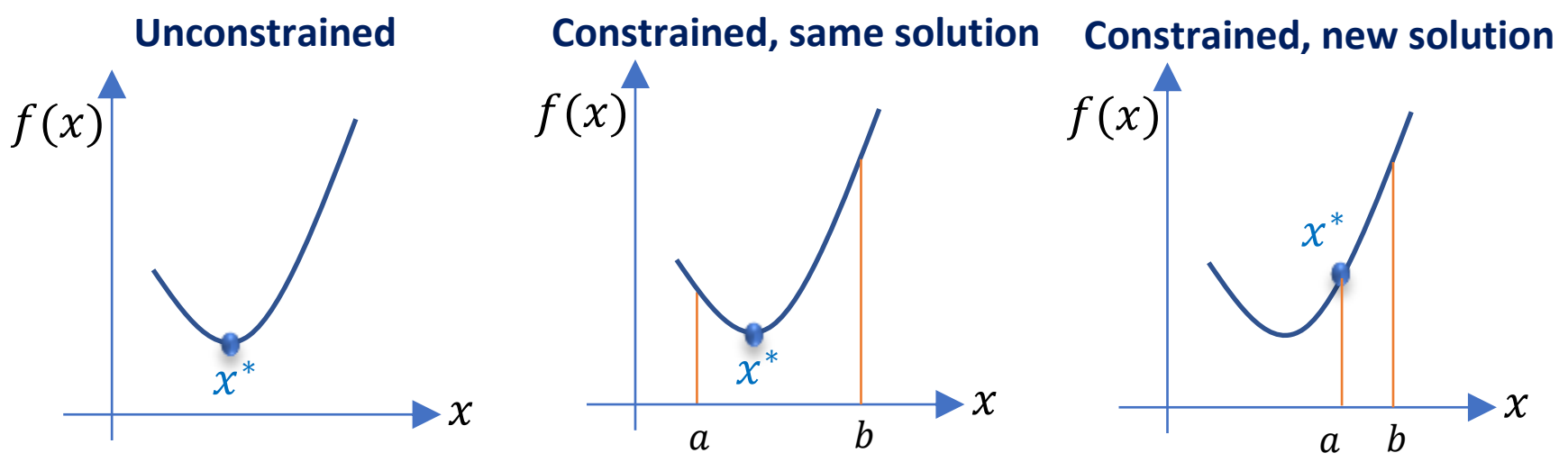

**Figure 10.1.** Constraints can change the solution to a problem but do not have to.

The general nonlinear constrained optimization problems can be formulated as [1]:

$$\begin{array}{lll} \text{Minimize} & f|\mathbf{x}\rangle & \\ \text{Subject to} & h_1|\mathbf{x}\rangle = 0, & \phi_1|\mathbf{x}\rangle \le 0, \\ & \vdots & \vdots \\ & h_m|\mathbf{x}\rangle = 0, & \phi_p|\mathbf{x}\rangle \le 0, \end{array} \tag{10.1.a}$$

where,

$$\mathbf{x} \equiv |\mathbf{x}\rangle \in \mathbb{R}^n, \tag{10.1.b}$$

$$f: \mathbb{R}^n \to \mathbb{R}, \tag{10.1.c}$$

$$h_i: \mathbb{R}^n \to \mathbb{R}, \quad i = 1, \dots m, \tag{10.1.d}$$

$$\phi_j: \mathbb{R}^n \to \mathbb{R}, \quad j = 1, \dots p. \tag{10.1.e}$$

In vector notation, the problem above can be represented in the following standard form:

$$\begin{array}{ll} \text{Minimize} & f|\mathbf{x}\rangle, \\ \text{Subject to} & \mathbf{h}|\mathbf{x}\rangle = |\mathbf{0}\rangle, \\ & \boldsymbol{\phi}|\mathbf{x}\rangle \le |\mathbf{0}\rangle, \end{array} \tag{10.2.a}$$

where,

$$\mathbf{h}: \mathbb{R}^n \to \mathbb{R}^m, \quad \mathbf{h}|\mathbf{x}\rangle = \begin{pmatrix} h_1|\mathbf{x}\rangle \\ \vdots \\ h_m|\mathbf{x}\rangle \end{pmatrix}, \tag{10.2.b}$$

$$\boldsymbol{\phi}: \mathbb{R}^n \to \mathbb{R}^p, \quad \boldsymbol{\phi}|\mathbf{x}\rangle = \begin{pmatrix} \phi_1|\mathbf{x}\rangle \\ \vdots \\ \phi_p|\mathbf{x}\rangle \end{pmatrix}. \tag{10.2.c}$$

In this chapter, we consider only the optimization problems with equality constraints such that

$$\begin{aligned} &\text{Minimize} && f|\mathbf{x}\rangle, \\ &\text{Subject to} && h_1|\mathbf{x}\rangle = 0, \\ & && \vdots \\ & && h_m|\mathbf{x}\rangle = 0. \end{aligned} \tag{10.3}$$

For solving this problem, there are several methods, e.g.

1- Direct substitution method

2- Lagrange multiplier method.

## 10.2 Direct Substitution Method

For an optimization problem with $n$ variables and $m$ equality constraints, it is possible (theoretically) to express any set of $m$ variables in terms of the remaining $(n-m)$ variables. When these expressions are substituted into the original objective function, then the reduced objective function involves only $n-m$ variables. This reduced objective function is not subjected to any constraint, so its optimum value can be found by using unconstrained optimization techniques. Theoretically, the method of direct substitution is very simple. However, from a practical point of view, it is not convenient. In most of the practical problems, the constraint equations are highly non-linear in nature. In those cases, it becomes impossible to solve them and express any $m$ variables in terms of the remaining $(n-m)$ variables from the given constraints [2].

***Example 10.1***

Solve

$$\begin{aligned} &\text{Minimize } z = 9 - 8x_1 - 6x_2 - 4x_3 + 2x_1^2 + 2x_2^2 + x_3^2 + 2x_1x_2 + 2x_1x_3 \\ &\text{subject to } x_1 + x_2 + 2x_3 = 3 \end{aligned}$$

***Solution***

From the given constraint, we have $x_2 = 3 - x_1 - 2x_3$, and substituting it in the objective function, we have

$$\begin{aligned} z &= 9 - 8x_1 - 6(3 - x_1 - 2x_3) - 4x_3 + 2x_1^2 + 2(3 - x_1 - 2x_3)^2 + x_3^2 + 2x_1(3 - x_1 - 2x_3) + 2x_1x_3 \\ &= 2x_1^2 + 9x_3^2 + 6x_1x_3 - 8x_1 - 16x_3 + 9. \end{aligned}$$

Now we have to optimize $z$ with respect to the variables $x_1$ and $x_3$.

The necessary conditions for the optimality of $z$ are given by

$$\frac{\partial z}{\partial x_1} = 0 \text{ and } \frac{\partial z}{\partial x_3} = 0,$$

or

$$\begin{aligned} 2x_1 + 3x_3 &= 4, \\ 3x_1 + 9x_3 &= 8. \end{aligned}$$

Solving these equations, we have $x_1 = \frac{4}{3}, x_3 = \frac{4}{9}$. Now,

$$x_2 = 3 - x_1 - 2x_3 = 7/9.$$

The Hessian matrix is given by

$$\mathbf{H} = \begin{pmatrix} \frac{\partial^2 z}{\partial x_1^2} & \frac{\partial^2 z}{\partial x_1 \partial x_3} \\ \frac{\partial^2 z}{\partial x_3 \partial x_1} & \frac{\partial^2 z}{\partial x_3^2} \end{pmatrix} = \begin{pmatrix} 4 & 6 \\ 6 & 18 \end{pmatrix}.$$

The leading principal minors are

$$H_1 = |4| = 4 \text{ and } H_2 = \begin{vmatrix} 4 & 6 \\ 6 & 18 \end{vmatrix} = 36.$$

Since $H_1 > 0$ and $H_2 > 0$, $z$ is the minimum for $x_1 = 4/3, x_2 = 7/9$ and $x_3 = 4/9$, and the minimum value of $z$ is $1/9$.

## 10.3 Lagrange Multipliers Illustration

Figure 10.2a shows the graph of a function $f$ defined by the equation $z = f|\mathbf{x}\rangle$. Observe that $f$ has an absolute minimum at $(0,0)^T$ and an absolute minimum value of 0. However, if the independent variables $x$ and $y$ are subjected to a constraint of the form $h|\mathbf{x}\rangle = k$, $|\mathbf{x}\rangle = (x,y)^T$, then the points $(x,y,z)^T$ that satisfy both $z = f|\mathbf{x}\rangle$ and $h|\mathbf{x}\rangle = k$ lie on the curve $C$, the intersection of the surface $z = f|\mathbf{x}\rangle$ and the $h|\mathbf{x}\rangle = k$ (Figure 10.2b). From the Figure 10.2b, we can see that the absolute minimum of $f$ subject to the constraint $h|\mathbf{x}\rangle = k$ occurs at the point $|\mathbf{x}^*\rangle = (a,b)^T$. Furthermore, $f$ has the constrained absolute minimum value $f|\mathbf{x}^*\rangle$ rather than the unconstrained absolute minimum value of 0 at $(0,0)^T$.

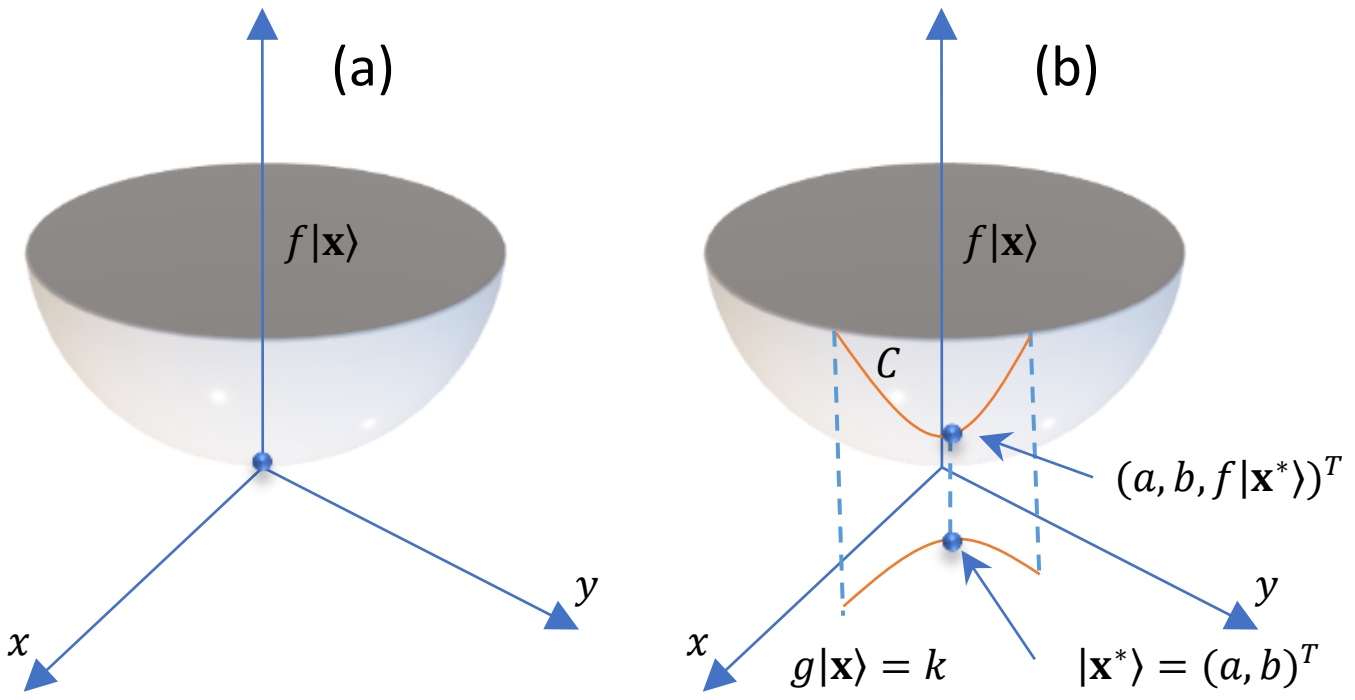

**Figure 10.2.** The function has an unconstrained minimum value of 0, but it has a constrained minimum value of $f|\mathbf{x}^*\rangle$ when subjected to the constraint $h|\mathbf{x}\rangle = k$.

We will now consider a method called the method of Lagrange multipliers. To see how this method works [3], let us start by studying the problem of finding the minimum of the objective function $f$ subject to the constraint $h|\mathbf{x}\rangle = k$. Figure 10.3 shows the level curves of $f$ drawn in the $xyz$-coordinate system. These level curves are reproduced in the $xy$-plane in Figure 10.3. Observe that the level curves of $f$ with equations $f|\mathbf{x}\rangle = c$, where $c < f|\mathbf{x}^*\rangle$, $|\mathbf{x}^*\rangle = (a,b)^T$, have no points in common with the graph of the constraint equation $h|\mathbf{x}\rangle = k$ (for example, the level curves $f|\mathbf{x}\rangle = c_1$ and $f|\mathbf{x}\rangle = c_2$ shown in Figure 10.3). Thus, points lying on these curves are not candidates for the constrained minimum of $f$.

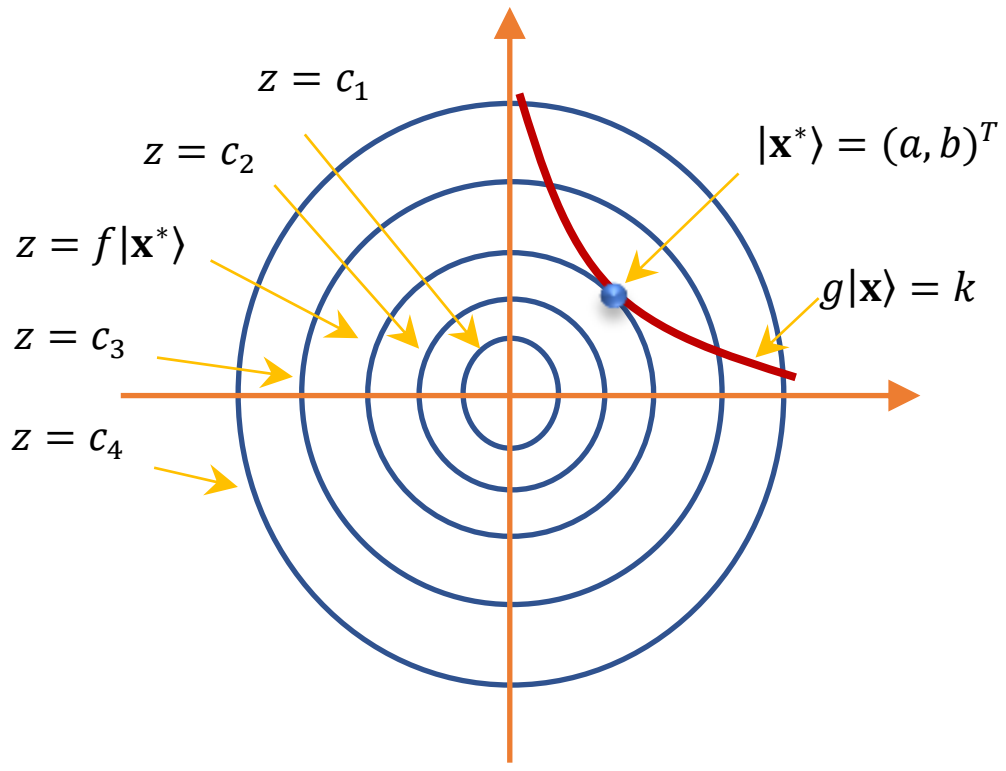

**Figure 10.3.** The level curves of $f$ in the $xyz$-plane.

On the other hand, the level curves of $f$ with equation $f|\mathbf{x}\rangle = c$, where $c \geq f|\mathbf{x}^*\rangle$, do intersect the graph of the constraint equation $h|\mathbf{x}\rangle = k$ (such as the level curves of $f|\mathbf{x}\rangle = c_3$ and $f|\mathbf{x}\rangle = c_4$). These points of intersection are candidates for the constrained minimum of $f$.

Finally, observe that the larger $c$ is for $c \geq f|\mathbf{x}^*\rangle$, the larger the value $f|\mathbf{x}\rangle$ is for $|\mathbf{x}\rangle$ lying on the level curve $h|\mathbf{x}\rangle = k$. This observation suggests that we can find the constrained minimum of $f$ by choosing the smallest value of $c$ so that the level curve $f|\mathbf{x}\rangle = c$ still intersects the curve $h|\mathbf{x}\rangle = k$. At such a point $|\mathbf{x}^*\rangle = (a, b)^T$ the level curve of $f$ just touches the graph of the constraint equation $h|\mathbf{x}\rangle = k$. That is, the two curves have a common tangent at $|\mathbf{x}^*\rangle$ (see Figure 10.3). Equivalently, their normal lines at this point coincide. Putting it yet another way, the gradient vectors $\nabla f|\mathbf{x}^*\rangle$ and $\nabla h|\mathbf{x}^*\rangle$ have the same direction, so $\nabla f|\mathbf{x}^*\rangle = \lambda \nabla h|\mathbf{x}^*\rangle$ for some scalar $\lambda$. These geometric arguments suggest the following theorem.

**Theorem 10.1:** Let $f$ and $h$ have continuous first partial derivatives in some region $D$ in the plane. If $f$ has an extremum at a point $|\mathbf{x}^*\rangle = (a, b)^T$ on the smooth constraint curve $h|\mathbf{x}\rangle = c$ lying in $D$ and

$$\nabla h|\mathbf{x}^*\rangle \neq |\mathbf{0}\rangle, \tag{10.4.a}$$

then there is a real number $\lambda$ such that

$$\nabla f|\mathbf{x}^*\rangle = \lambda \nabla h|\mathbf{x}^*\rangle. \tag{10.4.b}$$

The number $\lambda$ is called a Lagrange multiplier.

**Proof:**

Suppose that the smooth curve $C$ described by $h|\mathbf{x}\rangle = c$, $|\mathbf{x}\rangle = (x, y)^T$, is represented by the vector function

$$|\mathbf{r}(t)\rangle = (x(t), y(t))^T, \qquad |\mathbf{r}'(t)\rangle = (x'(t)\,, y'(t))^T \neq |\mathbf{0}\rangle,$$

where $x'$ and $y'$ are continuous on an open interval $I$ (Figure 10.4). Then the values assumed by $f$ on $C$ are given by

$$h(t) = f|\mathbf{r}(t)\rangle.$$

Suppose that $f$ has an extreme value at $|\mathbf{x}^*\rangle = (a, b)^T$. If $t_0$ is the point in $I$ corresponding to the point $|\mathbf{x}^*\rangle$, then $h$ has an extreme value at $t_0$. Therefore, $h'(t_0) = 0$. Using the Chain Rule, we have

$$\begin{aligned}
\frac{d}{dt}h(t_0) &= h'(t_0) \\
&= \frac{\partial}{\partial x} f|\mathbf{r}(t_0)\rangle \frac{d}{dt} x(t_0) + \frac{\partial}{\partial y} f|\mathbf{r}(t_0)\rangle \frac{d}{dt} y(t_0) \\
&= f_x|\mathbf{r}(t_0)\rangle x'(t_0) + f_y|\mathbf{r}(t_0)\rangle y'(t_0) \\
&= f_x|\mathbf{x}^*\rangle x'(t_0) + f_y|\mathbf{x}^*\rangle y'(t_0) \\
&= \left(f_x|\mathbf{x}^*\rangle, f_y|\mathbf{x}^*\rangle\right) \begin{pmatrix} x'(t_0) \\ y'(t_0) \end{pmatrix} \\
&= \langle \nabla f(\mathbf{x}^*) | \mathbf{r}'(t_0)\rangle \\
&= 0.
\end{aligned}$$

This shows that $\nabla f|\mathbf{x}^*\rangle$ is orthogonal to $\mathbf{r}'(t_0)$. But, $\nabla h|\mathbf{x}^*\rangle$ is orthogonal to $\mathbf{r}'(t_0)$. Therefore, the gradient vectors $\nabla f|\mathbf{x}^*\rangle$ and $\nabla h|\mathbf{x}^*\rangle$ are parallel, so there is a number $\lambda$ such that $\nabla f|\mathbf{x}^*\rangle = \lambda \nabla h|\mathbf{x}^*\rangle$.

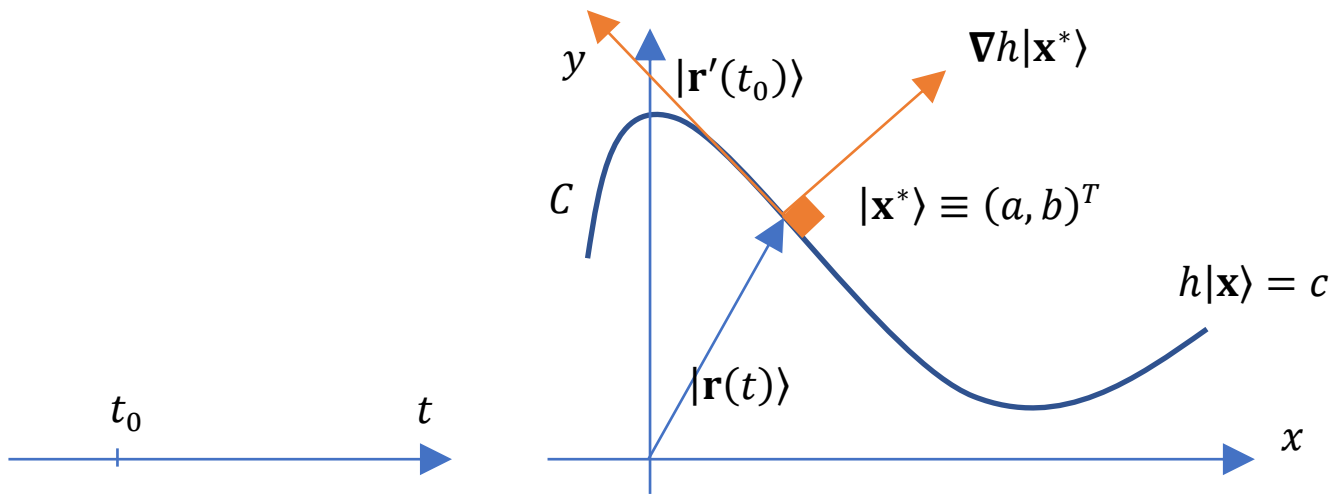

**Figure 10.4.** The curve $C$ is represented by the vector $|\mathbf{r}(t)\rangle$.

■

**The Method of Lagrange Multipliers**

Suppose $f$ and $h$ have continuous first partial derivatives. To find the maximum and minimum values of $f$ subject to the constraint $h$ (assuming that these extreme values exist and that $\nabla h \neq |\mathbf{0}\rangle$ on $h|\mathbf{x}\rangle = k$):

1. Solve the equations $\nabla f|\mathbf{x}\rangle = \lambda \nabla h|\mathbf{x}\rangle$ and $h|\mathbf{x}\rangle = k$ for $x$, $y$, and $\lambda$.
2. Evaluate $f$ at each solution point found in Step 1. The largest value yields the constrained maximum of $f$, and the smallest value yields the constrained minimum of $f$.

Note that. Since

$$\nabla f|\mathbf{x}\rangle = \begin{pmatrix} f_x|\mathbf{x}\rangle \\ f_y|\mathbf{x}\rangle \end{pmatrix}, \qquad \nabla h|\mathbf{x}\rangle = \begin{pmatrix} h_x|\mathbf{x}\rangle \\ h_y|\mathbf{x}\rangle \end{pmatrix}, \tag{10.5}$$

we see, by equating like components, that the vector equation

$$\nabla f|\mathbf{x}\rangle = \lambda \nabla h|\mathbf{x}\rangle, \tag{10.6}$$

is equivalent to the two scalar equations

$$f_x|\mathbf{x}\rangle = \lambda h_x|\mathbf{x}\rangle, \qquad f_y|\mathbf{x}\rangle = \lambda h_y|\mathbf{x}\rangle. \tag{10.7}$$

These scalar equations together with the constraint equation $h|\mathbf{x}\rangle = k$ give a system of three equations to be solved for the three unknowns $x$, $y$, and $\lambda$.

***Example 10.2***

Find the maximum and minimum values of the function $f|\mathbf{x}\rangle = x^2 - 2y$ subject to $x^2 + y^2 = 9$.

***Solution***

The constraint equation is $h|\mathbf{x}\rangle = x^2 + y^2 - 9 = 0$. Since

$$\nabla f|\mathbf{x}\rangle = \begin{pmatrix} 2x \\ -2 \end{pmatrix}, \qquad \nabla h|\mathbf{x}\rangle = \begin{pmatrix} 2x \\ 2y \end{pmatrix}.$$

The equation $\nabla f|\mathbf{x}\rangle = \lambda \nabla h|\mathbf{x}\rangle$ becomes

$$\begin{pmatrix} 2x \\ -2 \end{pmatrix} = \lambda \begin{pmatrix} 2x \\ 2y \end{pmatrix} = \begin{pmatrix} 2\lambda x \\ 2\lambda y \end{pmatrix}.$$

Hence, we have the following system of three equations in the three variables $x$, $y$, and $\lambda$:

$$2x = 2\lambda x, \quad -2 = 2\lambda y, \qquad x^2 + y^2 = 9.$$

From the first equation, we have

$$2x(1 - \lambda) = 0.$$

so $x = 0$, or $\lambda = 1$. If $x = 0$, then the third equation gives $y = \pm 3$. If $\lambda = 1$, then the second equation gives $y = -1$, which upon substitution into the third equation yields $x = \pm 2\sqrt{2}$. Therefore, $f$ has possible extreme values at the points $(0, -3)^T$, $(0,3)^T$, $(-2\sqrt{2}, -1)^T$, and $(2\sqrt{2}, -1)^T$. Evaluating $f$ at each of these points gives

$$f(0, -3)^T = 6, \qquad f(0,3)^T = -6, \qquad f\left(-2\sqrt{2}, -1\right)^T = 10 \quad \text{and} \quad f\left(2\sqrt{2}, -1\right)^T = 10.$$

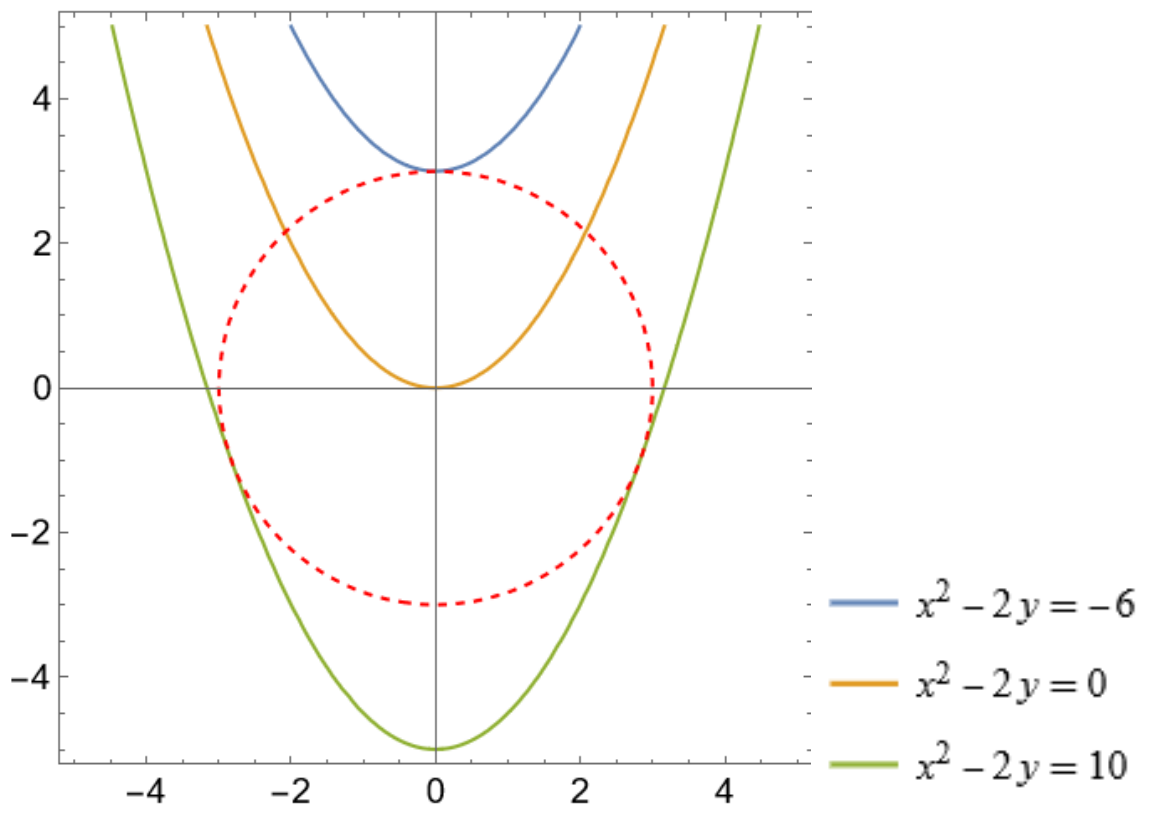

**Figure 10.5.** The extreme values of $f$ occur at the points where the level curves of $f$ are tangent to the graph of the constraint equation (the circle).

We conclude that the maximum value of $f$ on the circle $x^2 + y^2 = 9$ is 10, attained at the points $(-2\sqrt{2}, -1)^T$ and $(2\sqrt{2}, -1)^T$, and that the minimum value of $f$ on the circle is $-6$, attained at the point $(0,3)^T$. Figure 10.5 shows the graph of the constraint equation $x^2 + y^2 = 9$ and some level curves of the objective function $f$. Observe that the extreme values of $f$ are attained at the points where the level curves of $f$ are tangent to the graph of the constraint equation.

### 10.3.1. Optimizing a Function of Three Variables Subject to One Constraint

The proof of the Lagrange Theorem for functions of three variables [3] is similar to that for functions of two variables.

Since

$$\nabla f|\mathbf{x}\rangle = \begin{pmatrix} f_x|\mathbf{x}\rangle \\ f_y|\mathbf{x}\rangle \\ f_z|\mathbf{x}\rangle \end{pmatrix}, |\mathbf{x}\rangle = \begin{pmatrix} x \\ y \\ z \end{pmatrix}, \tag{10.8}$$

and

$$\nabla h|\mathbf{x}\rangle = \begin{pmatrix} h_x|\mathbf{x}\rangle \\ h_y|\mathbf{x}\rangle \\ h_z|\mathbf{x}\rangle \end{pmatrix}, \tag{10.9}$$

we see, by equating like components, that the vector equation

$$\nabla f|\mathbf{x}\rangle = \lambda \nabla h|\mathbf{x}\rangle, \tag{10.10}$$

is equivalent to the three scalar equations

$$\begin{aligned} f_x|\mathbf{x}\rangle &= \lambda h_x|\mathbf{x}\rangle, \\ f_y|\mathbf{x}\rangle &= \lambda h_y|\mathbf{x}\rangle, \\ f_z|\mathbf{x}\rangle &= \lambda h_z|\mathbf{x}\rangle. \end{aligned} \tag{10.11}$$

These scalar equations together with the constraint equation give a system of four equations to be solved for the four unknowns $x$, $y$, $z$, and $\lambda$.

### 10.3.2. Optimizing a Function Subject to Two Constraints

Some applications involve maximizing or minimizing an objective function $f$ subject to two or more constraints. Consider, for example, the problem of finding the extreme values of $f|\mathbf{x}\rangle, |\mathbf{x}\rangle = (x, y, z)^T$ subject to the two constraints,

$$g|\mathbf{x}\rangle = k \quad \text{and} \quad h|\mathbf{x}\rangle = l. \tag{10.12}$$

It can be shown that if $f$ has an extremum at $|\mathbf{x}^*\rangle = (a, b, c)^T$ subject to these constraints, then there are real numbers (Lagrange multipliers) $\lambda$ and $\mu$ such that

$$\nabla f|\mathbf{x}^*\rangle = \lambda \nabla g|\mathbf{x}^*\rangle + \mu \nabla h|\mathbf{x}^*\rangle. \tag{10.13}$$

Geometrically, we are looking for the extreme values of $f|\mathbf{x}\rangle$ on the curve of the intersection of the level surfaces $g|\mathbf{x}\rangle = k$ and $h|\mathbf{x}\rangle = l$. Condition (10.13) is a statement that at an extremum point $|\mathbf{x}^*\rangle$, the gradient of $f$ must lie in the plane determined by the gradient of $g$ and the gradient of $h$. (See Figure 10.6.) The vector (10.13) is equivalent to three scalar equations. When combined with the two constraint equations, this leads to a system of five equations that can be solved for the five unknowns $x$, $y$, $z$, $\lambda$, and $\mu$.

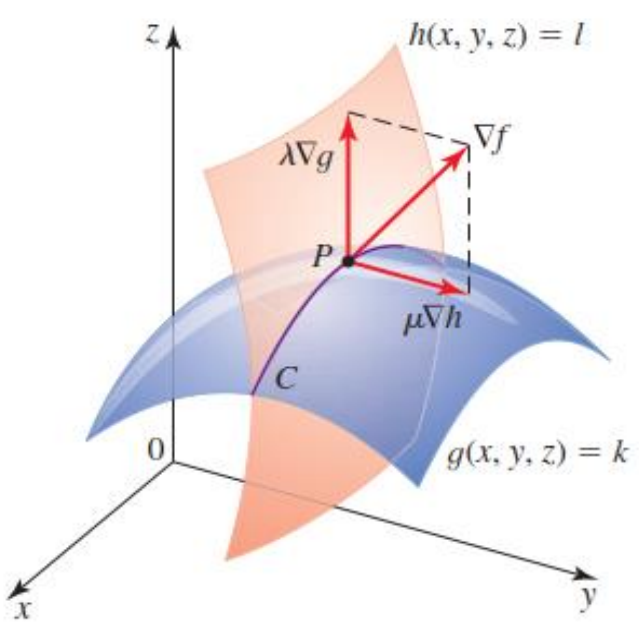

**Figure 10.6.** If $f$ has an extreme value at $P \equiv |\mathbf{x}^*\rangle \equiv (a, b, c)^T$, then $\nabla f|\mathbf{x}^*\rangle = \lambda \nabla g|\mathbf{x}^*\rangle + \mu \nabla h|\mathbf{x}^*\rangle$.

***Example 10.3***

Find the maximum and minimum values of the function $f|\mathbf{x}\rangle = 3x + 2y + 4z, |\mathbf{x}\rangle = (x, y, z)^T$ subject to the constraints $x - y + 2z = 1$ and $x^2 + y^2 = 4$.

***Solution***

Write the constraint equations in the form

$$h_1|\mathbf{x}\rangle = x - y + 2z - 1 = 0,$$
$$h_2|\mathbf{x}\rangle = x^2 + y^2 - 4 = 0.$$

Then the equation $\nabla f|\mathbf{x}\rangle = \lambda_1 \nabla h_1|\mathbf{x}\rangle + \lambda_2 \nabla h_2|\mathbf{x}\rangle$ becomes

$$\begin{pmatrix} 3 \\ 2 \\ 4 \end{pmatrix} = \lambda_1 \begin{pmatrix} 1 \\ -1 \\ 2 \end{pmatrix} + \lambda_2 \begin{pmatrix} 2x \\ 2y \\ 0 \end{pmatrix} = \begin{pmatrix} \lambda_1 + 2\lambda_2 x \\ -\lambda_1 + 2\lambda_2 y \\ 2\lambda_1 \end{pmatrix}.$$

Hence, we have the following system of five equations in the five variables, $x$, $y$, $z$, $\lambda_1$, and $\lambda_2$:

$$\lambda_1 + 2\lambda_2 x = 3,$$
$$-\lambda_1 + 2\lambda_2 y = 2,$$
$$2\lambda_1 = 4,$$
$$x - y + 2z = 1,$$
$$x^2 + y^2 = 4.$$

From the third equation, we have $\lambda_1 = 2$. Then

$$2 + 2\lambda_2 x = 3 \qquad \Rightarrow 2\lambda_2 x = 1,$$
$$-2 + 2\lambda_2 y = 2 \qquad \Rightarrow 2\lambda_2 y = 4.$$

The solution to the last two equations is

$$x = \frac{1}{2\lambda_2}, \quad y = \frac{2}{\lambda_2}.$$

Substituting these values of $x$ and $y$ into $x^2 + y^2 = 4$ yields

$$\left(\frac{1}{2\lambda_2}\right)^2 + \left(\frac{2}{\lambda_2}\right)^2 = 4 \qquad \Rightarrow \lambda_2^2 = \frac{17}{16}.$$

Therefore,

$$\lambda_2 = \pm\frac{\sqrt{17}}{4}.$$

So that

$$x = \pm\frac{2}{\sqrt{17}}, \qquad y = \pm\frac{8}{\sqrt{17}},$$

and the values of $z$ becomes $z = \frac{1}{2}\left(1 \pm \frac{6}{\sqrt{17}}\right)$. The value of $f$ at the point $\left(\frac{2}{\sqrt{17}}, \frac{8}{\sqrt{17}}, \frac{1}{2} + \frac{3}{\sqrt{17}}\right)^T$ is $2(1 + \sqrt{17})$, and the value of $f$ at the point $\left(-\frac{2}{\sqrt{17}}, -\frac{8}{\sqrt{17}}, \frac{1}{2} - \frac{3}{\sqrt{17}}\right)^T$ is $2(1 - \sqrt{17})$. Therefore, the maximum value of $f$ is $2(1 + \sqrt{17})$ and the minimum value of $f$ is $2(1 - \sqrt{17})$.

## 10.4 Non-Degenerate Constraint Qualification

The optimization problems with equality constraints can be formulated as:

$$\begin{aligned} &\text{Minimize} && f|\mathbf{x}\rangle, \\ &\text{Subject to} && h_1|\mathbf{x}\rangle = 0, \\ & && \vdots \\ & && h_m|\mathbf{x}\rangle = 0. \end{aligned} \tag{10.14}$$

This set of constraints defines a hypersurface in $\mathbb{R}^n$. Each function $h_i$ defines a surface $S_i = \{|\mathbf{x}\rangle: h_i|\mathbf{x}\rangle = 0\}$ of generally $n - 1$ dimensions in the space $\mathbb{R}^n$. This surface is smooth provided $h_i|\mathbf{x}\rangle \in C^1$. The $m$ constraints together defines a surface $S$, which is the intersection of the surfaces $S_1, \dots, S_m$; namely,

$$S = \{|\mathbf{x}\rangle: h_1|\mathbf{x}\rangle = 0\} \cap \cdots \cap \{|\mathbf{x}\rangle: h_m|\mathbf{x}\rangle = 0\}. \tag{10.15}$$

Using vector notation, we have

$$\begin{aligned} &\text{minimize} && f|\mathbf{x}\rangle, \\ &\text{subject to} && \mathbf{h}|\mathbf{x}\rangle = |\mathbf{0}\rangle. \end{aligned} \tag{10.16.1}$$

The vector-valued function $\mathbf{h}$ takes a point $|\mathbf{x}\rangle \in \mathbb{R}^n$ as input and produces the vector $\mathbf{h}|\mathbf{x}\rangle \in \mathbb{R}^m$ as output.

$$\mathbf{h}: \mathbb{R}^n \to \mathbb{R}^m, \quad \mathbf{h}|\mathbf{x}\rangle = \begin{pmatrix} h_1|\mathbf{x}\rangle \\ \vdots \\ h_m|\mathbf{x}\rangle \end{pmatrix}. \tag{10.16.2}$$

We assume that the function $\mathbf{h}$ is continuously differentiable, that is, $\mathbf{h} \in C^1$. Hence, we have

$$S = \{|\mathbf{x}\rangle : \mathbf{h}|\mathbf{x}\rangle = |\mathbf{0}\rangle\}. \tag{10.17}$$

**Definition (Feasible Set):** Any point satisfying the constraints is called a feasible point. The set of all feasible points

$$\{|\mathbf{x}\rangle \in \mathbb{R}^n : \mathbf{h}|\mathbf{x}\rangle = |\mathbf{0}\rangle\}, \tag{10.18}$$

is called the feasible set.

Remember, the Jacobian matrix of a vector-valued function of several variables is the matrix of all its first-order partial derivatives. Then the Jacobian matrix of $\mathbf{h}$ is defined to be an $m \times n$ matrix, denoted by $\mathbf{J}$, whose $(i,j)$th entry is $J_{ij} = \frac{\partial h_i}{\partial x_j}$, or explicitly

$$\mathbf{J} = \begin{pmatrix} \frac{\partial \mathbf{h}}{\partial x_1} & \cdots & \frac{\partial \mathbf{h}}{\partial x_n} \end{pmatrix} = \begin{pmatrix} \frac{\partial h_1}{\partial x_1} & \cdots & \frac{\partial h_1}{\partial x_n} \\ \vdots & \ddots & \vdots \\ \frac{\partial h_m}{\partial x_1} & \cdots & \frac{\partial h_m}{\partial x_n} \end{pmatrix} = \begin{pmatrix} \nabla h_1{}^T \\ \vdots \\ \nabla h_m{}^T \end{pmatrix}, \tag{10.19}$$

where $\nabla h_i{}^T$ is the transpose (row vector) of the gradient of the $i$ component. The Jacobian matrix, whose entries are functions of $|\mathbf{x}\rangle$, is denoted in various ways; common notations include $D\mathbf{h}$, $\mathbf{J}_\mathbf{h}$, and $\nabla\mathbf{h}$. Some authors define the Jacobian as the transpose of the form given above.

The Jacobian of a vector-valued function in several variables generalizes the gradient of a scalar-valued function in several variables, which in turn generalizes the derivative of a scalar-valued function of a single variable.

If $\mathbf{h}$ is differentiable at a point $|\mathbf{p}\rangle$ in $\mathbb{R}^n$, then its differential is represented by $\mathbf{J}_\mathbf{h}|\mathbf{p}\rangle$. In this case, the linear transformation represented by $\mathbf{J}_\mathbf{h}|\mathbf{p}\rangle \equiv \mathbf{J}_\mathbf{h}(\mathbf{p})$ is the best linear approximation of $\mathbf{h}$ near the point $|\mathbf{p}\rangle$, in a sense that

$$\mathbf{h}|\mathbf{x}\rangle \approx \mathbf{h}|\mathbf{p}\rangle + \mathbf{J}_\mathbf{h}(\mathbf{p})|\mathbf{x}-\mathbf{p}\rangle + \boldsymbol{O}(\|\mathbf{x}-\mathbf{p}\|), \tag{10.20}$$

where $\boldsymbol{O}(\|\mathbf{x}-\mathbf{p}\|)$ is a quantity that approaches zero much faster than the distance between $|\mathbf{x}\rangle$ and $|\mathbf{p}\rangle$ does as $|\mathbf{x}\rangle$ approaches $|\mathbf{p}\rangle$. This approximation specializes to the approximation of a scalar function of a single variable by its Taylor polynomial of degree one, namely

$$h(x) \approx h(p) + h'(p)(x-p) + O(x-p). \tag{10.21}$$

In this sense, the Jacobian may be regarded as a kind of "first-order derivative" of a vector-valued function of several variables. In particular, this means that the gradient of a scalar-valued function of several variables may, too, be regarded as its "first-order derivative".

For one constraint problem, consider the case where $\nabla h|\mathbf{x}^*\rangle = |\mathbf{0}\rangle$, or in other words, the point which minimizes $f|\mathbf{x}\rangle$ is also a critical point of $h|\mathbf{x}\rangle$. Remember our necessary condition for a minimum is $\nabla f|\mathbf{x}^*\rangle = \lambda\nabla h|\mathbf{x}^*\rangle$. Since $\nabla h|\mathbf{x}^*\rangle = |\mathbf{0}\rangle$, this implies that $\nabla f|\mathbf{x}^*\rangle = |\mathbf{0}\rangle$. However, this is the necessary condition for an unconstrained optimization problem, not a constrained one! In effect, when $\nabla h|\mathbf{x}^*\rangle = |\mathbf{0}\rangle$, the constraint is no longer considered in the problem. Hence, if any of the $\nabla h_i|\mathbf{x}^*\rangle$ is zero, then that constraint will not be considered in the analysis. Also, there will be a row of zeros in the Jacobian, so the Jacobian will not be full rank. The generalization of the condition that $\nabla h|\mathbf{x}^*\rangle \neq |\mathbf{0}\rangle$ for the case when $m = 1$ is that the Jacobian matrix must be full rank.

**Definition (Regular Point):** A point $|\mathbf{x}^*\rangle$ satisfying the constraints $h_1|\mathbf{x}^*\rangle = 0, \dots, h_m|\mathbf{x}^*\rangle = 0$ is said to be a regular point of the constraints if the gradient vectors $\nabla h_1|\mathbf{x}^*\rangle, \dots, \nabla h_m|\mathbf{x}^*\rangle$ are linearly independent.

**Definition (Non-Degenerate Constraint Qualification):** Let $D\mathbf{h}|\mathbf{x}^*\rangle$ be the Jacobian matrix of $\mathbf{h}$ at $|\mathbf{x}^*\rangle$, given by

$$D\mathbf{h}|\mathbf{x}^*\rangle = \begin{pmatrix} Dh_1|\mathbf{x}^*\rangle \\ \vdots \\ Dh_m|\mathbf{x}^*\rangle \end{pmatrix} = \begin{pmatrix} \nabla h_1(\mathbf{x}^*)^T \\ \vdots \\ \nabla h_m(\mathbf{x}^*)^T \end{pmatrix}. \quad (10.22)$$

Then, $|\mathbf{x}^*\rangle$ is regular if and only if rank $D\mathbf{h}|\mathbf{x}^*\rangle = m$, that is, the Jacobian matrix is of full rank.

**Remarks:**

- The definition states, in effect, that $|\mathbf{x}^*\rangle$ is a regular point of the constraints if it is a solution of $\mathbf{h}|\mathbf{x}\rangle = |\mathbf{0}\rangle$ and the Jacobian $D\mathbf{h}|\mathbf{x}^*\rangle$ has full row rank. The importance of a point $|\mathbf{x}^*\rangle$ being regular for a given set of equality constraints lies in the fact that a tangent plane of the hypersurface determined by the constraints at a regular point $|\mathbf{x}^*\rangle$ is well defined. Later in this chapter, the term 'tangent plane' will be used to express and describe important necessary as well as sufficient conditions for constrained optimization problems.
- Since $D\mathbf{h}|\mathbf{x}\rangle$ is a $m \times n$ matrix, it would not be possible for $|\mathbf{x}\rangle$ to be a regular point of the constraints if $m > n$. This leads to an upper bound for the number of independent equality constraints, i.e., $m \leq n$. Furthermore, if $m = n$, in many cases, the number of vectors $|\mathbf{x}\rangle$ that satisfy $\mathbf{h}|\mathbf{x}\rangle = |\mathbf{0}\rangle$ is finite and the optimization problem becomes a trivial one. For these reasons, we shall assume that $m < n$ throughout the rest of the lectures.

**From linear algebra, remember the following facts:**

- Given a linear map, $L: V \to W$ between two vector spaces $V$ and $W$, the kernel of $L$ (null space) is the vector subspace of all elements $|\mathbf{v}\rangle$ of $V$ such that $L|\mathbf{v}\rangle = |\mathbf{0}\rangle$, where $|\mathbf{0}\rangle$ denotes the zero vector in $W$.
- Consider a linear map represented as a $m \times n$ matrix $\mathbf{A}$ with coefficients in $\mathbb{R}$, that is operating on column vectors $|\mathbf{x}\rangle$ with $n$ components over $\mathbb{R}$. The kernel of this linear map is the set of solutions to the equation $\mathbf{A}|\mathbf{x}\rangle = |\mathbf{0}\rangle$. The dimension of the kernel of $\mathbf{A}$ is called the nullity of $\mathbf{A}$.
$$\text{Null}(\mathbf{A}) = \text{Ker}(\mathbf{A}) = \{|\mathbf{x}\rangle \in \mathbb{R}^n: \mathbf{A}|\mathbf{x}\rangle = |\mathbf{0}\rangle\}.$$
- The product $\mathbf{A}|\mathbf{x}\rangle$ can be written in terms of the dot product of vectors as follows:
$$\mathbf{A}|\mathbf{x}\rangle = \begin{pmatrix} \langle\mathbf{a}_1|\mathbf{x}\rangle \\ \vdots \\ \langle\mathbf{a}_m|\mathbf{x}\rangle \end{pmatrix}.$$
Here, $\langle\mathbf{a}_1|, \ldots, \langle\mathbf{a}_m|$ denote the rows of the matrix $\mathbf{A}$. It follows that $|\mathbf{x}\rangle$ is in the kernel of $\mathbf{A}$, if and only if $|\mathbf{x}\rangle$ is orthogonal to each of the row vectors of $\mathbf{A}$.
- The row space of a matrix $\mathbf{A}$ is the span of the row vectors of $\mathbf{A}$. The dimension of the row space of $\mathbf{A}$ is called the rank of $\mathbf{A}$, and the dimension of the kernel of $\mathbf{A}$ is called the nullity of $\mathbf{A}$. These quantities are related by the rank–nullity theorem
$$\text{nullity}(\mathbf{A}) + \text{rank}(\mathbf{A}) = n.$$
- So that the rank–nullity theorem can be restated as
$$\dim(\text{Ker } L) + \dim(\text{Im } L) = \dim(V).$$

**Definition (Dimension of the Surface S):** The set of equality constraints $h_1|\mathbf{x}\rangle = 0, \ldots, h_m|\mathbf{x}\rangle = 0$, $h_i: \mathbb{R}^n \to \mathbb{R}$, describes a surface

$$S = \{|\mathbf{x}\rangle \in \mathbb{R}^n: h_1|\mathbf{x}\rangle = 0, \ldots, h_m|\mathbf{x}\rangle = 0\}. \quad (10.23)$$

Assuming the points in $S$ are regular, the dimension of the surface $S$ is $n - m$.

***Example 10.4***

Let $n = 3$ and $m = 2$. Assuming regularity, the feasible set $S$ is a one-dimensional object (i.e., a curve in $\mathbb{R}^3$). For example, let

$$h_1|\mathbf{x}\rangle = x = 0,$$
$$h_2|\mathbf{x}\rangle = y - z^2 = 0.$$

In this case, $\nabla h_1|\mathbf{x}\rangle = (1,0,0)^T$, and $\nabla h_2|\mathbf{x}\rangle = (0,1,-2z)^T$. Hence, the vectors $\nabla h_1|\mathbf{x}\rangle$ and $\nabla h_2|\mathbf{x}\rangle$ are linearly independent in $\mathbb{R}^3$. Thus,

$$\dim S = \dim\{|\mathbf{x}\rangle: h_1|\mathbf{x}\rangle = 0, h_2|\mathbf{x}\rangle = 0\} = n - m = 1.$$

See Figure 10.7 for a graphical illustration.

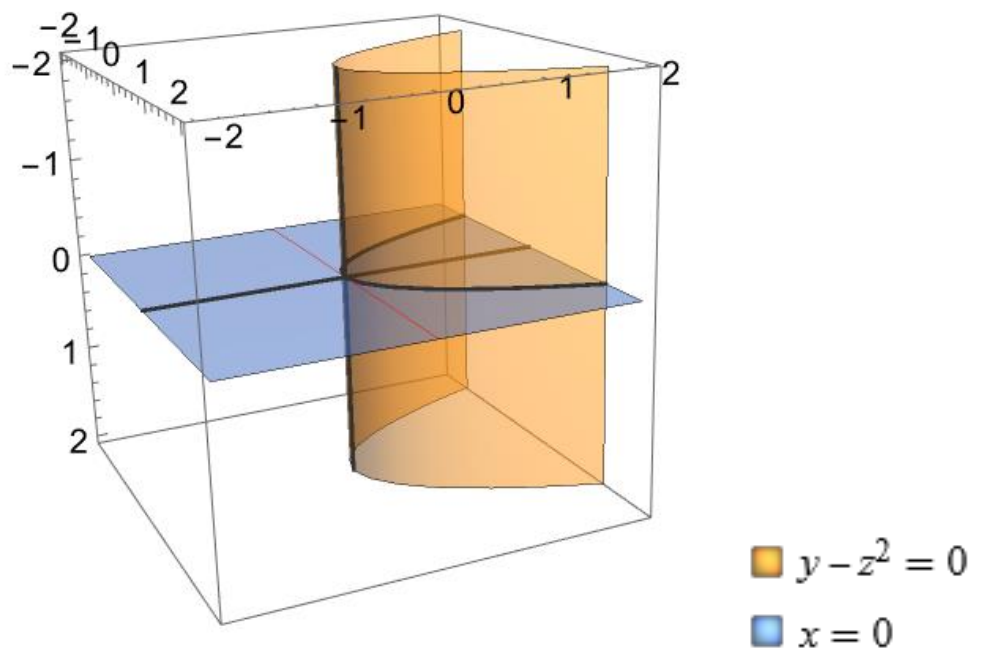

**Figure 10.7.** A one-dimensional surface in $\mathbb{R}^3$.

## 10.5 Tangent and Normal Spaces

Two important concepts are related to the gradients of the objective function, and the constraints of the optimization problem are the tangent plane and the normal plane [4-7]. In this section, we discuss these concepts in some detail.

**Definition (Curve $C$ on a Surface $S$):** A curve $C$ on a surface $S$ is a set of points $\{|\mathbf{x}(t)\rangle \in S: t \in (a,b)\}$, continuously parameterized by $t \in (a,b)$; that is, $|\mathbf{x}\rangle: (a,b) \to S$ is a continuous function.

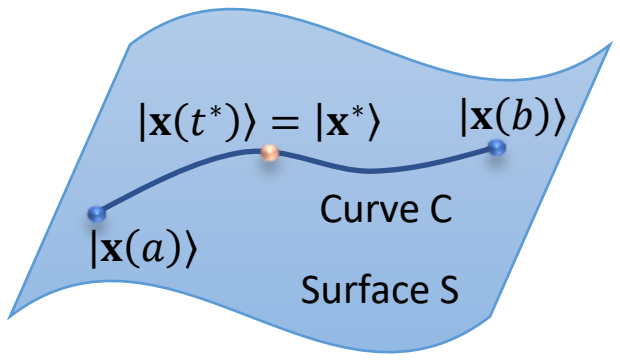

**Figure 10.8.** A curve on a surface.

A graphical illustration of the definition of a curve is given in Figure 10.8. The definition of a curve implies that all the points on the curve satisfy the equation describing the surface. The curve $C$ passes through a point $|\mathbf{x}^*\rangle$ if there exists $t^* \in (a,b)$ such that $|\mathbf{x}(t^*)\rangle = |\mathbf{x}^*\rangle$.

**Definition (Differentiable and Twice Differentiable Curve):** The curve $C = \{|\mathbf{x}(t)\rangle: t \in (a,b)\}$ is differentiable if

$$|\dot{\mathbf{x}}(t)\rangle = \frac{d}{dt}|\mathbf{x}(t)\rangle = \begin{pmatrix} \dot{x}_1(t) \\ \vdots \\ \dot{x}_n(t) \end{pmatrix}, \tag{10.24}$$

exists for all $t \in (a,b)$.
The curve $C = \{|\mathbf{x}(t)\rangle: t \in (a,b)\}$ is twice differentiable if

$$|\ddot{\mathbf{x}}(t)\rangle = \frac{d^2}{dt^2}|\mathbf{x}(t)\rangle = \begin{pmatrix} \ddot{x}_1(t) \\ \vdots \\ \ddot{x}_n(t) \end{pmatrix}, \tag{10.25}$$

exists for all $t \in (a,b)$.

Note that, the vector $|\dot{\mathbf{x}}(t^*)\rangle$ is tangent to the curve $C$ at $|\mathbf{x}^*\rangle$ (see Figure 10.9).

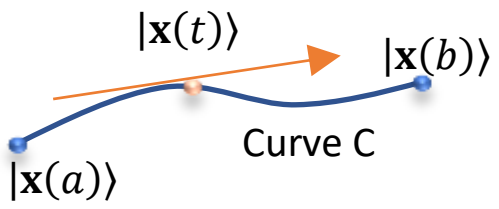

**Figure 10.9.** Geometric interpretation of the differentiability of a curve.

We are now ready to introduce the notions of tangent space. For this, recall the set

$$S = \{|\mathbf{x}\rangle \in \mathbb{R}^n : \mathbf{h}|\mathbf{x}\rangle = |\mathbf{0}\rangle\}, \tag{10.26}$$

where $\mathbf{h} \in C^1$. We think of $S$ as a surface in $\mathbb{R}^n$. Now consider all differentiable curves on $S$ passing through a point $|\mathbf{x}^*\rangle$. The tangent plane at $|\mathbf{x}^*\rangle$ is defined as the collection of the derivatives at $|\mathbf{x}^*\rangle$ of all these differentiable curves. Figure 10.10 illustrates the notion of a tangent plane.

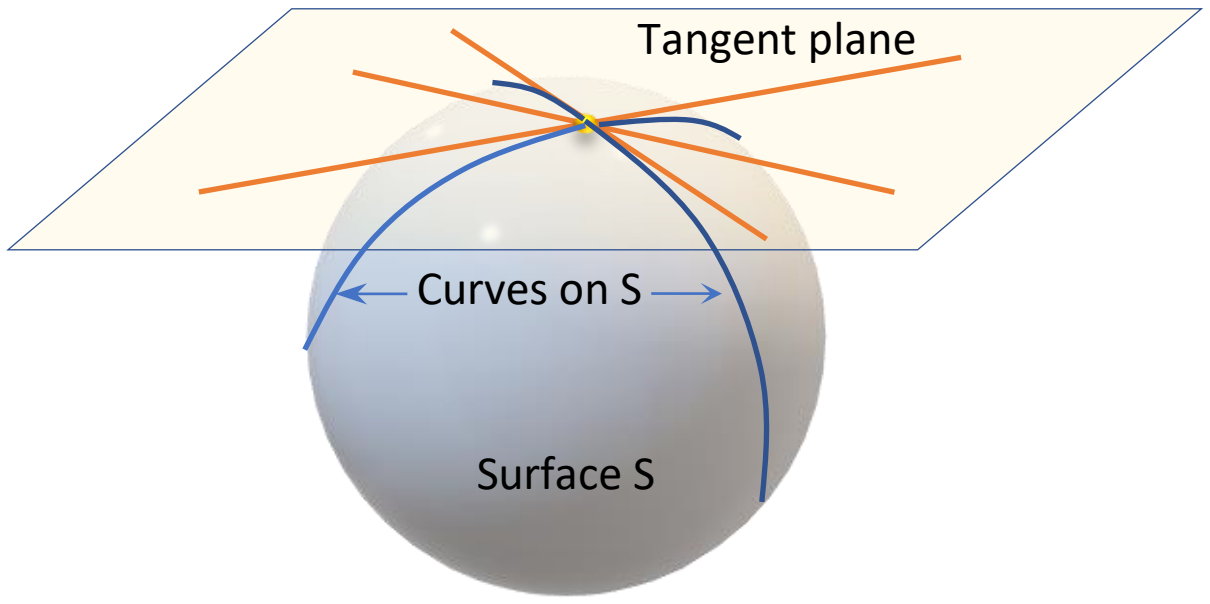

**Figure 10.10.** The tangent plane to the surface $S$ at the point $|\mathbf{x}^*\rangle$.

The tangent plane of a smooth function $f|\mathbf{x}\rangle$ at a given point $|\mathbf{x}^*\rangle$ can be defined as a hyperplane that passes through point $|\mathbf{x}^*\rangle$ with $\nabla f|\mathbf{x}^*\rangle$ as the normal. For example, for $n = 2$ the contours, tangent plane, and gradient of a smooth function are related to each other as illustrated in Figure 10.11. Following this idea, the tangent plane at point $|\mathbf{x}^*\rangle$ can be defined analytically as the set

$$T|\mathbf{x}^*\rangle = \{|\mathbf{x}\rangle : \langle \nabla f(\mathbf{x}^*)|\mathbf{x} - \mathbf{x}^*\rangle = 0\}. \tag{10.27}$$

In other words, a point $|\mathbf{x}\rangle$ lies on the tangent plane if the vector that connects $|\mathbf{x}^*\rangle$ to $|\mathbf{x}\rangle$ is orthogonal to the gradient $\nabla f|\mathbf{x}^*\rangle$, as can be seen in Figure 10.11.

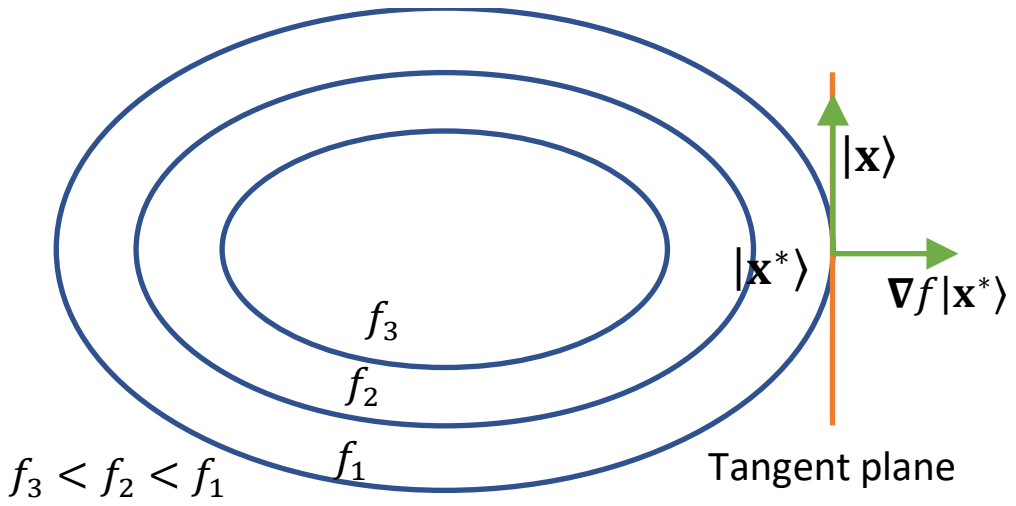

**Figure 10.11.** The gradient $\nabla f|\mathbf{x}^*\rangle$ is orthogonal to the tangent plane.

Ideally, we would like to express this tangent plane in terms of derivatives of constraint functions $h_i$ that define the surface. By using the Taylor series of the constraint function $h_i|\mathbf{x}\rangle$ at the feasible point $|\mathbf{x}^*\rangle$, we can write

$$\begin{aligned} h_i|\mathbf{x}^* + \mathbf{s}\rangle &= h_i|\mathbf{x}^*\rangle + \langle \nabla h_i(\mathbf{x}^*)|\mathbf{s}\rangle + O(\|\mathbf{s}\|) \\ &= \langle \nabla h_i(\mathbf{x}^*)|\mathbf{s}\rangle + O(\|\mathbf{s}\|), \end{aligned} \tag{10.28}$$

since $h_i|\mathbf{x}^*\rangle = 0$. It follows that $h_i|\mathbf{x}^* + \mathbf{s}\rangle = 0$ is equivalent to

$$\langle \nabla h_i(\mathbf{x}^*)|\mathbf{s}\rangle = 0 \quad \text{for } i = 1, \ldots, m. \tag{10.29}$$

In other words, $|\mathbf{s}\rangle$ is feasible if and only if it is orthogonal to the gradients of the constraint functions.

**Definition (Tangent Space):** The tangent space at a point $|\mathbf{x}^*\rangle$ on the surface $S = \{|\mathbf{x}\rangle \in \mathbb{R}^n : \mathbf{h}|\mathbf{x}\rangle = |\mathbf{0}\rangle\}$ is the set

$$T|\mathbf{x}^*\rangle = \{|\mathbf{y}\rangle : D\mathbf{h}(\mathbf{x}^*)|\mathbf{y}\rangle = |\mathbf{0}\rangle\}. \tag{10.30}$$

Notes:

1- The tangent space $T|\mathbf{x}^*\rangle$ is the null space of the matrix $D\mathbf{h}(\mathbf{x}^*)$, that is,

$$T|\mathbf{x}^*\rangle = \mathcal{N}\big(D\mathbf{h}(\mathbf{x}^*)\big). \tag{10.31}$$

The tangent space is, therefore, a subspace of $\mathbb{R}^n$.

2- Assuming $|\mathbf{x}^*\rangle$ is regular, the dimension of the tangent space is $n-m$, where $m$ is the number of equality constraints $h_i|\mathbf{x}^*\rangle = 0$.

***Example 10.5***

Let

$$S = \{|\mathbf{x}\rangle \in \mathbb{R}^3 : h_1|\mathbf{x}\rangle = x_1 = 0, h_2|\mathbf{x}\rangle = x_1 - x_2 = 0\}.$$

***Solution***

Then, $S$ is the $x_3$-axis in $\mathbb{R}^3$ (see Figure 10.12). We have

$$D\mathbf{h}(\mathbf{x}) = \begin{pmatrix} \nabla {h_1}^T \\ \nabla {h_2}^T \end{pmatrix} = \begin{pmatrix} 1 & 0 & 0 \\ 1 & -1 & 0 \end{pmatrix}.$$

Because $\nabla h_1$ and $\nabla h_2$ are linearly independent when evaluated at any $|\mathbf{x}\rangle \in S$, all the points of $S$ are regular. The tangent space at an arbitrary point of $S$ is

$$\begin{aligned} T|\mathbf{x}\rangle &= \{|\mathbf{y}\rangle : \langle \nabla h_1(\mathbf{x})|\mathbf{y}\rangle = 0, \langle \nabla h_2(\mathbf{x})|\mathbf{y}\rangle = 0\} \\ &= \left\{ |\mathbf{y}\rangle : \begin{pmatrix} 1 & 0 & 0 \\ 1 & -1 & 0 \end{pmatrix} \begin{pmatrix} y_1 \\ y_2 \\ y_3 \end{pmatrix} = |\mathbf{0}\rangle \right\} \\ &= \{(0,0,\alpha)^T : \alpha \in \mathbb{R}\} \\ &= \text{the } x_3\text{-axis in } \mathbb{R}^3. \end{aligned}$$

In this example, the tangent space $T|\mathbf{x}\rangle$ at any point $|\mathbf{x}\rangle \in S$ is a one-dimensional subspace of $\mathbb{R}^3$.

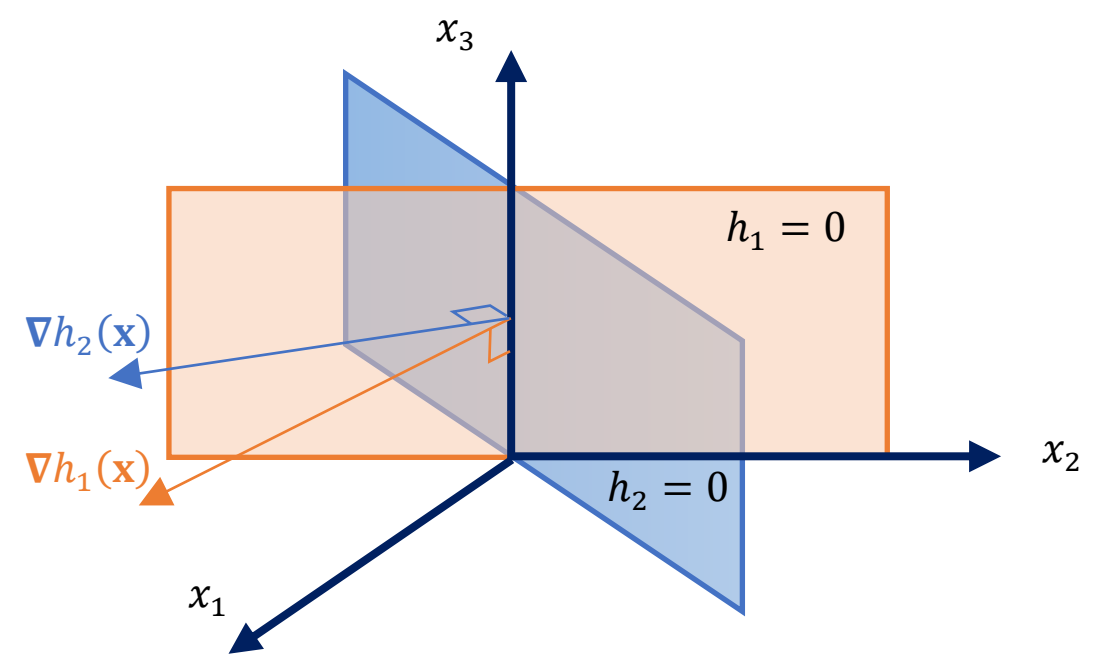

**Figure 10.12.** $S = \{|\mathbf{x}\rangle \in \mathbb{R}^3 : h_1|\mathbf{x}\rangle = x_1 = 0, h_2|\mathbf{x}\rangle = x_1 - x_2 = 0\}$.

**Theorem 10.2:** Suppose $|\mathbf{x}^*\rangle \in S$ is a regular point, and $T(\mathbf{x}^*)$ is the tangent space at $|\mathbf{x}^*\rangle$. Then, $|\mathbf{y}\rangle \in T(\mathbf{x}^*)$ if and only if there exists a differentiable curve in $S$ passing through $|\mathbf{x}^*\rangle$ with derivative $|\mathbf{y}\rangle$ at $|\mathbf{x}^*\rangle$.

**Proof:**

$\Leftarrow$: Suppose there exists a curve $\{|\mathbf{x}(t)\rangle : t \in (a,b)\}$ in $S$ such that $|\mathbf{x}(t^*)\rangle = |\mathbf{x}^*\rangle$ and $|\dot{\mathbf{x}}(t^*)\rangle = |\mathbf{y}\rangle$ for some $t^* \in (a,b)$. Then,

$$\mathbf{h}|\mathbf{x}(t)\rangle = |\mathbf{0}\rangle,$$

for all $t \in (a,b)$. If we differentiate the function $\mathbf{h}|\mathbf{x}(t)\rangle$ with respect to $t$ using the chain rule, we obtain

$$\frac{d}{dt}\mathbf{h}|\mathbf{x}(t)\rangle = D\mathbf{h}\big(\mathbf{x}(t)\big)|\dot{\mathbf{x}}(t)\rangle = |\mathbf{0}\rangle,$$

for all $t \in (a,b)$. Therefore, at $t^*$, we get

$$D\mathbf{h}(\mathbf{x}^*)|\mathbf{y}\rangle = |\mathbf{0}\rangle,$$

and hence $|\mathbf{y}\rangle \in T(\mathbf{x}^*)$.

$\Rightarrow$: we leave the opposite direction of the theorem as exercise.

■

**Definition (Normal Space):** The normal space $N(\mathbf{x}^*)$ at a point $|\mathbf{x}^*\rangle$ on the surface $S = \{|\mathbf{x}\rangle \in \mathbb{R}^n : \mathbf{h}|\mathbf{x}\rangle = |\mathbf{0}\rangle\}$ is the subspace of $\mathbb{R}^n$ spanned by the vectors $\nabla h_1|\mathbf{x}^*\rangle, \ldots, \nabla h_m|\mathbf{x}^*\rangle$, that is,

$$N(\mathbf{x}^*) = \text{span}\{\nabla h_1|\mathbf{x}^*\rangle, \ldots, \nabla h_m|\mathbf{x}^*\rangle\}$$

$$= \{|\mathbf{x}\rangle \in \mathbb{R}^n : |\mathbf{x}\rangle = \alpha_1 \nabla h_1 |\mathbf{x}^*\rangle + \cdots + \alpha_m \nabla h_m |\mathbf{x}^*\rangle, \alpha_1, \ldots, \alpha_m \in \mathbb{R}\}. \quad (10.32)$$

or

$$N(\mathbf{x}^*) = \{|\mathbf{x}\rangle \in \mathbb{R}^n : |\mathbf{x}\rangle = D\mathbf{h}(\mathbf{x}^*)^T |\mathbf{z}\rangle, |\mathbf{z}\rangle \in \mathbb{R}^m\}. \quad (10.33)$$

**Notes:**

- We can express the normal space $N(\mathbf{x}^*)$ as
$$N(\mathbf{x}^*) = \mathcal{R}(D\mathbf{h}(\mathbf{x}^*)^T), \quad (10.34)$$
that is, the range of the matrix $D\mathbf{h}(\mathbf{x}^*)^T$.
- Assuming $|\mathbf{x}^*\rangle$ is regular, the dimension of the normal space $N(\mathbf{x}^*)$ is $m$.

**Theorem 10.3:** We have
$$T(\mathbf{x}^*) = N(\mathbf{x}^*)^\perp \text{ and } T(\mathbf{x}^*)^\perp = N(\mathbf{x}^*). \quad (10.35)$$

**Proof:**

By definition of $T(\mathbf{x}^*)$, we may write

$$T(\mathbf{x}^*) = \{|\mathbf{y}\rangle \in \mathbb{R}^n : \langle \mathbf{x}|\mathbf{y}\rangle = 0 \text{ for all } |\mathbf{x}\rangle \in N(\mathbf{x}^*)\}.$$

Hence, by definition of $N(\mathbf{x}^*)$, we have $T(\mathbf{x}^*) = N(\mathbf{x}^*)^\perp$. It is easy to prove that $T(\mathbf{x}^*)^\perp = N(\mathbf{x}^*)$.

■

## 10.6 First-Order Necessary Conditions

The necessary conditions for a point $|\mathbf{x}^*\rangle$ to be a local minimizer are useful in two situations [4-7]:
(a) They can be used to exclude those points that do not satisfy the necessary conditions from the candidate points;
(b) They become sufficient conditions when the objective function in question is convex.
We now generalize the Lagrange theorem for the case when $f: \mathbb{R}^n \longrightarrow \mathbb{R}$ and $\mathbf{h}: \mathbb{R}^n \longrightarrow \mathbb{R}^m$, $m \le n$.

**Theorem 10.4: Lagrange Theorem (First-order necessary conditions).** Let $|\mathbf{x}^*\rangle$ be a local extremum point of $f: \mathbb{R}^n \longrightarrow \mathbb{R}$, subject to the constraints $\mathbf{h}|\mathbf{x}\rangle = |\mathbf{0}\rangle$, $\mathbf{h}: \mathbb{R}^n \longrightarrow \mathbb{R}^m$, $m \le n$. Assume that $|\mathbf{x}^*\rangle$ is a regular point. Then, there exists $|\boldsymbol{\lambda}^*\rangle \in \mathbb{R}^m$ such that

$$Df|\mathbf{x}^*\rangle + \langle \boldsymbol{\lambda}^*|(D\mathbf{h}|\mathbf{x}^*\rangle) = \langle \mathbf{0}|, \quad (10.36.1)$$
$$\nabla f|\mathbf{x}^*\rangle + (D\mathbf{h}|\mathbf{x}^*\rangle)^T|\boldsymbol{\lambda}^*\rangle = |\mathbf{0}\rangle, \quad (10.36.2)$$
$$\nabla f|\mathbf{x}^*\rangle + \langle D\mathbf{h}(\mathbf{x}^*)|\boldsymbol{\lambda}^*\rangle = |\mathbf{0}\rangle, \quad (10.36.3)$$
$$\nabla f|\mathbf{x}^*\rangle + \sum_{i=1}^{m} \lambda_i^* \nabla h_i |\mathbf{x}^*\rangle = |\mathbf{0}\rangle. \quad (10.36.4)$$

**Proof:**

We need to prove that

$$\nabla f|\mathbf{x}^*\rangle = -(D\mathbf{h}|\mathbf{x}^*\rangle)^T|\boldsymbol{\lambda}^*\rangle,$$

for some $|\boldsymbol{\lambda}^*\rangle \in \mathbb{R}^m$, that is, $\nabla f|\mathbf{x}^*\rangle \in N(\mathbf{x}^*)$. But, by Theorem 10.3, $N(\mathbf{x}^*) = T(\mathbf{x}^*)^\perp$. Therefore, it remains to show that $\nabla f|\mathbf{x}^*\rangle \in T(\mathbf{x}^*)^\perp$. Suppose

$$|\mathbf{y}\rangle \in T(\mathbf{x}^*).$$

Then, by Theorem 10.2, there exists a differentiable curve $\{|\mathbf{x}(t)\rangle : t \in (a,b)\}$ such that for all $t \in (a,b)$,

$$\mathbf{h}|\mathbf{x}(t)\rangle = |\mathbf{0}\rangle,$$

and there exists $t^* \in (a,b)$ satisfying

$$|\mathbf{x}(t^*)\rangle = |\mathbf{x}^*\rangle, \qquad |\dot{\mathbf{x}}(t^*)\rangle = |\mathbf{y}\rangle.$$

Now, consider the function $\phi(t) = f|\mathbf{x}(t)\rangle$. Note that $t^*$ is a local minimizer of $\phi(t)$. Hence, (using unconstrained local minimizers condition), we have

$$\frac{d}{dt}\phi(t^*) = 0,$$

and

$$\frac{d}{dt}\phi(t^*) = Df(\mathbf{x}^*)|\dot{\mathbf{x}}(t^*)\rangle$$
$$= Df(\mathbf{x}^*)|\mathbf{y}\rangle$$
$$= \langle \nabla f(\mathbf{x}^*)|\mathbf{y}\rangle$$
$$= 0.$$

So, all $|\mathbf{y}\rangle \in T(\mathbf{x}^*)$ satisfy

$$\langle \nabla f(\mathbf{x}^*)|\mathbf{y}\rangle = 0.$$

Hence, $\nabla f|\mathbf{x}^*\rangle$ is orthogonal to the tangent plane. that is

$$\nabla f|\mathbf{x}^*\rangle \in T(\mathbf{x}^*)^{\perp} = N(\mathbf{x}^*).$$

This implies that $\nabla f|\mathbf{x}^*\rangle$ is a linear combination of the gradients of $h_i$ at $|\mathbf{x}^*\rangle$. This completes the proof.

■

- Lagrange theorem states that if $|\mathbf{x}^*\rangle$ is an extremizer, then the gradient of the objective function $f$ can be expressed as a linear combination of the gradients of the constraints. The vector $|\boldsymbol{\lambda}^*\rangle$ is called the *Lagrange multiplier vector*, and its components are the *Lagrange multipliers*.
- It should be noted that the first-order necessary conditions $\nabla f|\mathbf{x}^*\rangle + \sum_{i=1}^{m} \lambda_i^* \nabla h_i|\mathbf{x}^*\rangle = |\mathbf{0}\rangle$ together with the constraints $\mathbf{h}|\mathbf{x}^*\rangle = |\mathbf{0}\rangle$ give a total of $n+m$ equations in the $n+m$ variables comprising $|\mathbf{x}^*\rangle, |\boldsymbol{\lambda}^*\rangle$.
- If $|\mathbf{x}^*\rangle$ is a local minimizer and $|\boldsymbol{\lambda}^*\rangle$ is the associated vector of Lagrange multipliers, the set $\{|\mathbf{x}^*\rangle, |\boldsymbol{\lambda}^*\rangle\}$ may be referred to as the minimizer set.
- The Lagrange condition is only necessary but not sufficient; that is, a point $|\mathbf{x}^*\rangle$ satisfying the above equations need not be an extremizer.
- Theorem 10.4 can be related to the first-order necessary conditions for a minimum for the case of unconstrained minimization as follows. If function $f|\mathbf{x}\rangle$ is minimized without constraints, we can consider the problem as the special case of (10.14) and (10.16) where the number of constraints is reduced to zero. In such a case, the condition of Theorem 10.4 becomes $\nabla f|\mathbf{x}^*\rangle = |\mathbf{0}\rangle$.

**Definition (Lagrangian Function):** The Lagrangian function $\mathcal{L}: \mathbb{R}^n \times \mathbb{R}^m \longrightarrow \mathbb{R}$ is defined as

$$\mathcal{L}(\mathbf{x}, \boldsymbol{\lambda}) = f|\mathbf{x}\rangle + \langle \boldsymbol{\lambda}|\mathbf{h}|\mathbf{x}\rangle. \quad (10.37)$$

The Lagrange condition for a local minimizer $|\mathbf{x}^*\rangle$ can be represented using the Lagrangian function as

$$D\mathcal{L}(\mathbf{x}^*, \boldsymbol{\lambda}^*) = \langle \mathbf{0}|, \quad (10.38)$$

for some $|\boldsymbol{\lambda}^*\rangle$. From (10.38), we can denote that the necessary condition in the Lagrange theorem is equivalent to the first-order necessary condition for unconstrained optimization applied to the Lagrangian function.
To explain this fact, let us denote the derivative of $\mathcal{L}$ with respect to $|\mathbf{x}\rangle$ and $|\boldsymbol{\lambda}\rangle$ as $D_x\mathcal{L}$, and $D_\lambda\mathcal{L}$ , respectively. Then,

$$D\mathcal{L}(\mathbf{x}, \boldsymbol{\lambda}) = \big(D_x\mathcal{L}(\mathbf{x}, \boldsymbol{\lambda}), D_\lambda\mathcal{L}(\mathbf{x}, \boldsymbol{\lambda})\big). \quad (10.39)$$

It is clear that $D_x\mathcal{L}(\mathbf{x}, \boldsymbol{\lambda}) = Df|\mathbf{x}\rangle + \langle \boldsymbol{\lambda}|D\mathbf{h}|\mathbf{x}\rangle$ and $D_\lambda\mathcal{L}(\mathbf{x}, \boldsymbol{\lambda}) = (\mathbf{h}|\mathbf{x}\rangle)^T$. Hence, the Lagrange theorem for a local minimizer $|\mathbf{x}^*\rangle$ can be stated as

$$D_x\mathcal{L}(\mathbf{x}^*, \boldsymbol{\lambda}^*) = \langle \mathbf{0}|,$$
$$D_\lambda\mathcal{L}(\mathbf{x}^*, \boldsymbol{\lambda}^*) = \langle \mathbf{0}|, \quad (10.40)$$

for some $|\boldsymbol{\lambda}^*\rangle$, which is equivalent to

$$D\mathcal{L}(\mathbf{x}^*, \boldsymbol{\lambda}^*) = \langle \mathbf{0}|. \quad (10.41)$$

So, in order to find possible extremizers, this entails solving the following equations:

$$D_x\mathcal{L}(\mathbf{x}, \boldsymbol{\lambda}) = \langle \mathbf{0}|,$$
$$D_\lambda\mathcal{L}(\mathbf{x}, \boldsymbol{\lambda}) = \langle \mathbf{0}|. \quad (10.42)$$

***Example 10.6***

$$\text{Minimize: } z = x_1 + x_2 + x_3,$$
$$\text{Subject to: } x_1^2 + x_2 = 3,$$
$$x_1 + 3x_2 + 2x_3 = 7.$$

***Solution***

The given program is equivalent to the unconstrained minimization of

$$z = \frac{1}{2}(x_1^2 + x_1 + 4),$$

which obviously has a solution. We may therefore apply the method of Lagrange multipliers to the original program standardized as

$$\text{minimize: } z = x_1 + x_2 + x_3,$$
$$\text{subject to: } x_1^2 + x_2 - 3 = 0,$$
$$x_1 + 3x_2 + 2x_3 - 7 = 0.$$

Here, $f|\mathbf{x}\rangle = x_1 + x_2 + x_3$, $|\mathbf{x}\rangle = (x_1, x_2, x_3)^T$,$n = 3$ (variables), $m = 2$ (constraints),

$$h_1|\mathbf{x}\rangle = x_1^2 + x_2 - 3,$$
$$h_2|\mathbf{x}\rangle = x_1 + 3x_2 + 2x_3 - 7.$$

The Lagrangian function is

$$\mathcal{L} = (x_1 + x_2 + x_3) - \lambda_1(x_1^2 + x_2 - 3) - \lambda_2(x_1 + 3x_2 + 2x_3 - 7),$$

and

$$\frac{\partial}{\partial x_1}\mathcal{L} = 1 - 2\lambda_1 x_1 - \lambda_2 = 0,$$
$$\frac{\partial}{\partial x_2}\mathcal{L} = 1 - \lambda_1 - 3\lambda_2 = 0,$$
$$\frac{\partial}{\partial x_3}\mathcal{L} = 1 - 2\lambda_2 = 0,$$
$$\frac{\partial}{\partial \lambda_1}\mathcal{L} = -(x_1^2 + x_2 - 3) = 0,$$
$$\frac{\partial}{\partial \lambda_2}\mathcal{L} = -(x_1 + 3x_2 + 2x_3 - 7) = 0.$$

Successively solving equations, we obtain $\lambda_2 = 0.5$, $\lambda_1 = 0.5$, $x_1 = 0.5$, $x_2 = 2.75$, and $x_3 = -0.875$, with

$$z = x_1 + x_2 + x_3 = 0.5 + 2.75 - 0.875 = 2.375.$$

Since the first partial derivatives of $f|\mathbf{x}\rangle$, $h_1|\mathbf{x}\rangle$, and $h_2|\mathbf{x}\rangle$ are all continuous, and since

$$\mathbf{J} = \begin{pmatrix} \frac{\partial h_1}{\partial x_1} & \frac{\partial h_1}{\partial x_2} & \frac{\partial h_1}{\partial x_3} \\ \frac{\partial h_2}{\partial x_1} & \frac{\partial h_2}{\partial x_2} & \frac{\partial h_2}{\partial x_3} \end{pmatrix} = \begin{pmatrix} 2x_1 & 1 & 0 \\ 1 & 3 & 2 \end{pmatrix},$$

is of rank 2 everywhere (the last two columns are linearly independent everywhere), either $x_1 = 0.5$, $x_2 = 2.75$, and $x_3 = -0.875$ is the optimal solution to the program, or no optimal solution exists. Checking feasible points in the region around $(0.5, 2.75, -0.875)^T$, we find that this point is indeed the location of a (global) minimum for the problem.

## 10.7 Second-Order Necessary and Sufficient Conditions

By an argument analogous to that used for the unconstrained case, we can also derive the corresponding second-order conditions for equality constrained problems [4-7]. We assume that $f: \mathbb{R}^n \to \mathbb{R}$ and $\mathbf{h}: \mathbb{R}^n \to \mathbb{R}^m$ are twice continuously differentiable, that is, $f, \mathbf{h} \in C^2$. Let

$$\mathcal{L}(\mathbf{x}, \boldsymbol{\lambda}) = f|\mathbf{x}\rangle + \langle \boldsymbol{\lambda}|\mathbf{h}|\mathbf{x}\rangle = f|\mathbf{x}\rangle + \lambda_1 h_1|\mathbf{x}\rangle + \cdots + \lambda_m h_m|\mathbf{x}\rangle, \tag{10.43}$$

be the Lagrangian function. Let $\boldsymbol{H}_{\mathcal{L}}(\mathbf{x}, \boldsymbol{\lambda})$ be the Hessian matrix of $\mathcal{L}(\mathbf{x}, \boldsymbol{\lambda})$ with respect to $|\mathbf{x}\rangle$, that is,

$$\boldsymbol{H}_{\mathcal{L}}(\mathbf{x}, \boldsymbol{\lambda}) = \boldsymbol{H}_f(\mathbf{x}) + \lambda_1 \boldsymbol{H}_1(\mathbf{x}) + \cdots + \lambda_m \boldsymbol{H}_m(\mathbf{x}), \tag{10.44}$$

where $\boldsymbol{H}_f(\mathbf{x})$ is the Hessian matrix of $f$ at $|\mathbf{x}\rangle$, and $\boldsymbol{H}_k(\mathbf{x})$ is the Hessian matrix of $h_k$ at $|\mathbf{x}\rangle$, $k = 1, \dots, m$, given by

$$\boldsymbol{H}_k(\mathbf{x}) = \begin{pmatrix} \frac{\partial^2}{\partial x_1^2} h_k|\mathbf{x}\rangle & \cdots & \frac{\partial^2 h_k}{\partial x_n \partial x_1} h_k|\mathbf{x}\rangle \\ \vdots & \ddots & \vdots \\ \frac{\partial^2}{\partial x_1 \partial x_n} h_k|\mathbf{x}\rangle & \cdots & \frac{\partial^2}{\partial x_n^2} h_k|\mathbf{x}\rangle \end{pmatrix}. \tag{10.45}$$

Let $\Sigma$ be:

$$\Sigma = \lambda_1 \boldsymbol{H}_1(\mathbf{x}) + \cdots + \lambda_m \boldsymbol{H}_m(\mathbf{x}). \tag{10.46}$$

Using the above notation, we can write

$$\boldsymbol{H}_{\mathcal{L}}(\mathbf{x},\boldsymbol{\lambda}) = \boldsymbol{H}_f|\mathbf{x}\rangle + \Sigma. \tag{10.47}$$

**Theorem 10.5 (Second-Order Necessary Conditions):** Let $|\mathbf{x}^*\rangle$ be a local minimizer of $f: \mathbb{R}^n \to \mathbb{R}$ subject to $\mathbf{h}|\mathbf{x}\rangle = |\mathbf{0}\rangle$, $\mathbf{h}: \mathbb{R}^n \to \mathbb{R}^m$, $m \le n$, and $f, \mathbf{h} \in C^2$. Suppose $|\mathbf{x}^*\rangle$ is regular. Then, there exists $|\boldsymbol{\lambda}^*\rangle \in \mathbb{R}^m$ such that

1. $Df|\mathbf{x}^*\rangle + \langle\boldsymbol{\lambda}^*|(D\mathbf{h}|\mathbf{x}^*\rangle) = \langle\mathbf{0}|$; and
2. for all $|\mathbf{y}\rangle \in T(\mathbf{x}^*)$, we have $\langle\mathbf{y}|\boldsymbol{H}_{\mathcal{L}}(\mathbf{x},\boldsymbol{\lambda})|\mathbf{y}\rangle \ge 0$.

**Proof:**

From Lagrange condition, Theorem 10.4, there is $\boldsymbol{\lambda}^* \in \mathbb{R}^m$ such that $Df|\mathbf{x}^*\rangle + \langle\boldsymbol{\lambda}^*|(D\mathbf{h}|\mathbf{x}^*\rangle) = \langle\mathbf{0}|$. It remains to prove the second part of the theorem. Suppose $|\mathbf{y}\rangle \in T(\mathbf{x}^*)$, that is, $|\mathbf{y}\rangle$ belongs to the tangent space to $S = \{|\mathbf{x}\rangle \in \mathbb{R}^n: \mathbf{h}|\mathbf{x}\rangle = |\mathbf{0}\rangle\}$ at $|\mathbf{x}^*\rangle$. Because $\mathbf{h} \in C^2$, following the argument of Theorem 10.2, there exists a twice differentiate curve $\{|\mathbf{x}(t)\rangle: t \in (a,b)\}$ on $S$ such that

$$|\mathbf{x}(t^*)\rangle = |\mathbf{x}^*\rangle,\ \ |\dot{\mathbf{x}}(t^*)\rangle = |\mathbf{y}\rangle,$$

for some $t^* \in (a,b)$. Observe that by assumption, $t^*$ is a local minimizer of the function $\phi(t) = f|\mathbf{x}(t)\rangle$. From the second-order necessary condition for unconstrained minimization, we obtain

$$\frac{d^2}{dt^2}\phi(t^*) \ge 0.$$

Using the following formula

$$\frac{d}{dt}\langle\mathbf{y}(t)|\mathbf{z}(t)\rangle = \left\langle\mathbf{z}(t)\middle|\frac{d}{dt}\mathbf{y}(t)\right\rangle + \left\langle\mathbf{y}(t)\middle|\frac{d}{dt}\mathbf{z}(t)\right\rangle,$$

the chain rule yields

$$\begin{aligned}\frac{d^2}{dt^2}\phi(t^*) &= \frac{d}{dt}\left[Df\big(\mathbf{x}(t^*)\big)|\dot{\mathbf{x}}(t^*)\rangle\right]\\ &= \big\langle\dot{\mathbf{x}}(t^*)\big|\boldsymbol{H}_f(\mathbf{x}^*)\big|\dot{\mathbf{x}}(t^*)\big\rangle + Df(\mathbf{x}^*)|\ddot{\mathbf{x}}(t^*)\rangle = \big\langle\mathbf{y}\big|\boldsymbol{H}_f(\mathbf{x}^*)\big|\mathbf{y}\big\rangle + Df(\mathbf{x}^*)|\ddot{\mathbf{x}}(t^*)\rangle \ge 0.\end{aligned}$$

Because $\mathbf{h}|\mathbf{x}(t)\rangle = |\mathbf{0}\rangle$ for all $t \in (a,b)$, we have

$$\frac{d^2}{dt^2}\langle\boldsymbol{\lambda}^*|\mathbf{h}|\mathbf{x}(t)\rangle = 0.$$

Thus, for all $t \in (a,b)$,

$$\begin{aligned}\frac{d^2}{dt^2}\langle\boldsymbol{\lambda}^*|\mathbf{h}|\mathbf{x}(t)\rangle &= \frac{d}{dt}\left[\left\langle\boldsymbol{\lambda}^*\middle|\frac{d}{dt}\mathbf{h}\big(\mathbf{x}(t)\big)\right\rangle\right]\\ &= \frac{d}{dt}\left[\sum_{k=1}^{m}\lambda_k^*\frac{d}{dt}h_k|\mathbf{x}(t)\rangle\right]\\ &= \frac{d}{dt}\left[\sum_{k=1}^{m}\lambda_k^* Dh_k\big(\mathbf{x}(t)\big)|\dot{\mathbf{x}}(t)\rangle\right]\\ &= \sum_{k=1}^{m}\lambda_k^*\frac{d}{dt}\left(Dh_k\big(\mathbf{x}(t)\big)|\dot{\mathbf{x}}(t)\rangle\right)\\ &= \sum_{k=1}^{m}\lambda_k^*\left[\big\langle\dot{\mathbf{x}}(t)\big|\boldsymbol{H}_k\big(\mathbf{x}(t)\big)\big|\dot{\mathbf{x}}(t)\big\rangle + Dh_k\big(\mathbf{x}(t)\big)|\ddot{\mathbf{x}}(t)\rangle\right]\\ &= \sum_{k=1}^{m}\left[\big\langle\dot{\mathbf{x}}(t)\big|\lambda_k^*\boldsymbol{H}_k\big(\mathbf{x}(t)\big)\big|\dot{\mathbf{x}}(t)\big\rangle + \lambda_k^* Dh_k\big(\mathbf{x}(t)\big)|\ddot{\mathbf{x}}(t)\rangle\right]\\ &= \big\langle\dot{\mathbf{x}}(t)\big|\textstyle\sum_{k=1}^{m}\lambda_k^*\boldsymbol{H}_k\big(\mathbf{x}(t)\big)\big|\dot{\mathbf{x}}(t)\big\rangle + \left(\displaystyle\sum_{k=1}^{m}\lambda_k^* Dh_k\big(\mathbf{x}(t)\big)\right)|\ddot{\mathbf{x}}(t)\rangle\\ &= \langle\dot{\mathbf{x}}(t)|\Sigma|\dot{\mathbf{x}}(t)\rangle + \langle\boldsymbol{\lambda}^*|D\mathbf{h}\big(\mathbf{x}(t)\big)|\ddot{\mathbf{x}}(t)\rangle = 0.\end{aligned}$$

In the case $t = t^*$, we have

$$\langle \mathbf{y}|\Sigma|\mathbf{y}\rangle + \langle \boldsymbol{\lambda}^*|D\mathbf{h}(\mathbf{x}^*)|\ddot{\mathbf{x}}(t^*)\rangle = 0.$$

Adding the above equation to the inequality

$$\langle \mathbf{y}|\boldsymbol{H}_f(\mathbf{x}^*)|\mathbf{y}\rangle + Df(\mathbf{x}^*)|\ddot{\mathbf{x}}(t^*)\rangle \geq 0,$$

yields

$$\langle \mathbf{y}|\left(\boldsymbol{H}_f(\mathbf{x}^*) + \Sigma\right)|\mathbf{y}\rangle + \left(Df(\mathbf{x}^*) + \langle \boldsymbol{\lambda}^*|D\mathbf{h}(\mathbf{x}^*)\right)|\ddot{\mathbf{x}}(t^*)\rangle \geq 0,$$

But, from Theorem 10.4, $Df|\mathbf{x}^*\rangle + \langle \boldsymbol{\lambda}^*|(D\mathbf{h}|\mathbf{x}^*\rangle) = \langle \mathbf{0}|$ . Subsequently,

$$\langle \mathbf{y}|\left(\boldsymbol{H}_f(\mathbf{x}^*) + \Sigma\right)|\mathbf{y}\rangle = \langle \mathbf{y}|\boldsymbol{H}_{\mathcal{L}}(\mathbf{x}^*, \boldsymbol{\lambda}^*)|\mathbf{y}\rangle \geq 0,$$

which complete the proof.

■

**Theorem 10.6 (Second-Order Sufficient Conditions):** Suppose $f, \mathbf{h} \in C^2$ and there exist a point $|\mathbf{x}^*\rangle \in \mathbb{R}^n$ and $|\boldsymbol{\lambda}^*\rangle \in \mathbb{R}^m$ such that

1. $Df|\mathbf{x}^*\rangle + \langle \boldsymbol{\lambda}^*|D\mathbf{h}|\mathbf{x}^*\rangle = \langle \mathbf{0}|$; and
2. for all $|\mathbf{y}\rangle \in T(\mathbf{x}^*)$, $|\mathbf{y}\rangle \neq |\mathbf{0}\rangle$, we have $\langle \mathbf{y}|\boldsymbol{H}_{\mathcal{L}}(\mathbf{x}^*, \boldsymbol{\lambda}^*)|\mathbf{y}\rangle > 0$.

Then, $\mathbf{x}^*$ is a strict local minimizer of $f$ subject to $\mathbf{h}|\mathbf{x}\rangle = |\mathbf{0}\rangle$.

From Theorem 10.6, if an $|\mathbf{x}^*\rangle$ satisfies the Lagrange condition, Theorem 10.4, and $\boldsymbol{H}_{\mathcal{L}}(\mathbf{x}^*, \boldsymbol{\lambda}^*)$ is positive definite on $T(\mathbf{x}^*)$, then $|\mathbf{x}^*\rangle$ is a strict local minimizer.

***Example 10.7***

Consider the problem

$$\text{Maximize } (x_1 - 1)^2 + (x_2 - 1)^2$$
$$\text{Subject to } x_1^2 + x_2^2 - 1 = 0.$$

***Solution***

The Lagrangian and subsequent first-order conditions would be

$$\mathcal{L}(x_1, x_2, \lambda) = (x_1 - 1)^2 + (x_2 - 1)^2 - \lambda(x_1^2 + x_2^2 - 1),$$

$$\boldsymbol{\nabla}_{\mathbf{x}}\mathcal{L}(x_1, x_2, \lambda) = \begin{pmatrix} 2x_1(1-\lambda) - 2 \\ 2x_2(1-\lambda) - 2 \end{pmatrix} = |\mathbf{0}\rangle,$$

From the two equations, we conclude $x_1 = x_2$, together with $x_1^2 + x_2^2 - 1 = 0$, we have two first-order stationary solutions $(x_1 = x_2 = \frac{1}{\sqrt{2}}, \lambda = 1 - \sqrt{2})$ and $(x_1 = x_2 = -\frac{1}{\sqrt{2}}, \lambda = 1 + \sqrt{2})$.

The Lagrangian Hessian matrix $\boldsymbol{H}_{\mathcal{L}}(\mathbf{x}, \boldsymbol{\lambda})$, at two $\lambda$s, becomes

$$\begin{pmatrix} 2(1-\lambda) & 0 \\ 0 & 2(1-\lambda) \end{pmatrix} \Rightarrow \begin{pmatrix} 2\sqrt{2} & 0 \\ 0 & 2\sqrt{2} \end{pmatrix} \text{ with } \left(\lambda = 1 - \sqrt{2}\right),$$
$$\Rightarrow \begin{pmatrix} -2\sqrt{2} & 0 \\ 0 & -2\sqrt{2} \end{pmatrix} \text{ with } \left(\lambda = 1 + \sqrt{2}\right),$$

where the first one is positive definite and the second negative definite, and they remain so in tangent subspace. Thus, $x_1 = x_2 = \frac{1}{\sqrt{2}}$ is a minimum and $x_1 = x_2 = -\frac{1}{\sqrt{2}}$ is a maximum.

## 10.8 Active and Inactive Inequality Constraints

Let us consider the simplest case — two variables and one inequality constraint [8]:

$$\text{Maximize } f|\mathbf{x}\rangle, \qquad (10.48)$$
$$\text{subject to } \phi|\mathbf{x}\rangle \leq b.$$

where $|\mathbf{x}\rangle = (x, y)^T$. In Figure 10.13, the red curve is the curve $\phi|\mathbf{x}\rangle = b$; the orange region is the constraint set $\phi|\mathbf{x}\rangle \leq b$. The blue lines are the contour lines of the function $f$. In Figure 10.13, one notes that the highest-level curve of $f$ which meets the constraint set meets it at the point $|\mathbf{p}\rangle$. Since $|\mathbf{p}\rangle$ lies on the boundary of the constraint set where $\phi|\mathbf{x}\rangle = b$, we say that the constraint is binding (or is active, effective, or tight) at $|\mathbf{p}\rangle$. The level set of $f$ and the level set of $\phi$ are tangent to each other at $|\mathbf{p}\rangle$. This means that $\nabla f|\mathbf{p}\rangle$ and $\nabla\phi|\mathbf{p}\rangle$ line up — point in the same direction or in opposite directions — and therefore that $\nabla f|\mathbf{p}\rangle$ is a multiple of $\nabla\phi|\mathbf{p}\rangle$. If we let $\lambda$ denote the multiplier, then $\nabla f|\mathbf{p}\rangle = \lambda\nabla\phi|\mathbf{p}\rangle$, or

$$\nabla f|\mathbf{p}\rangle - \lambda\nabla\phi|\mathbf{p}\rangle = |\mathbf{0}\rangle. \tag{10.49}$$

It is important to note that the gradient $\nabla f|\mathbf{p}\rangle$ points in the direction in which $f$ increases most rapidly at $|\mathbf{p}\rangle$. Moreover, $\nabla\phi|\mathbf{p}\rangle$ points to the set $\phi|\mathbf{x}\rangle \geq b$. Since $|\mathbf{p}\rangle$ maximizes $f$ on the set $\phi|\mathbf{x}\rangle \leq b$, the gradient $\nabla f|\mathbf{p}\rangle$ must point to the region where $\phi|\mathbf{x}\rangle \geq b$. Hence, two gradient vectors, $\nabla f|\mathbf{p}\rangle$ and $\nabla\phi|\mathbf{p}\rangle$, point in the same direction. Thus, if $\nabla f|\mathbf{p}\rangle = \lambda\nabla\phi|\mathbf{p}\rangle$, the multiplier $\lambda$ must be $\geq 0$.

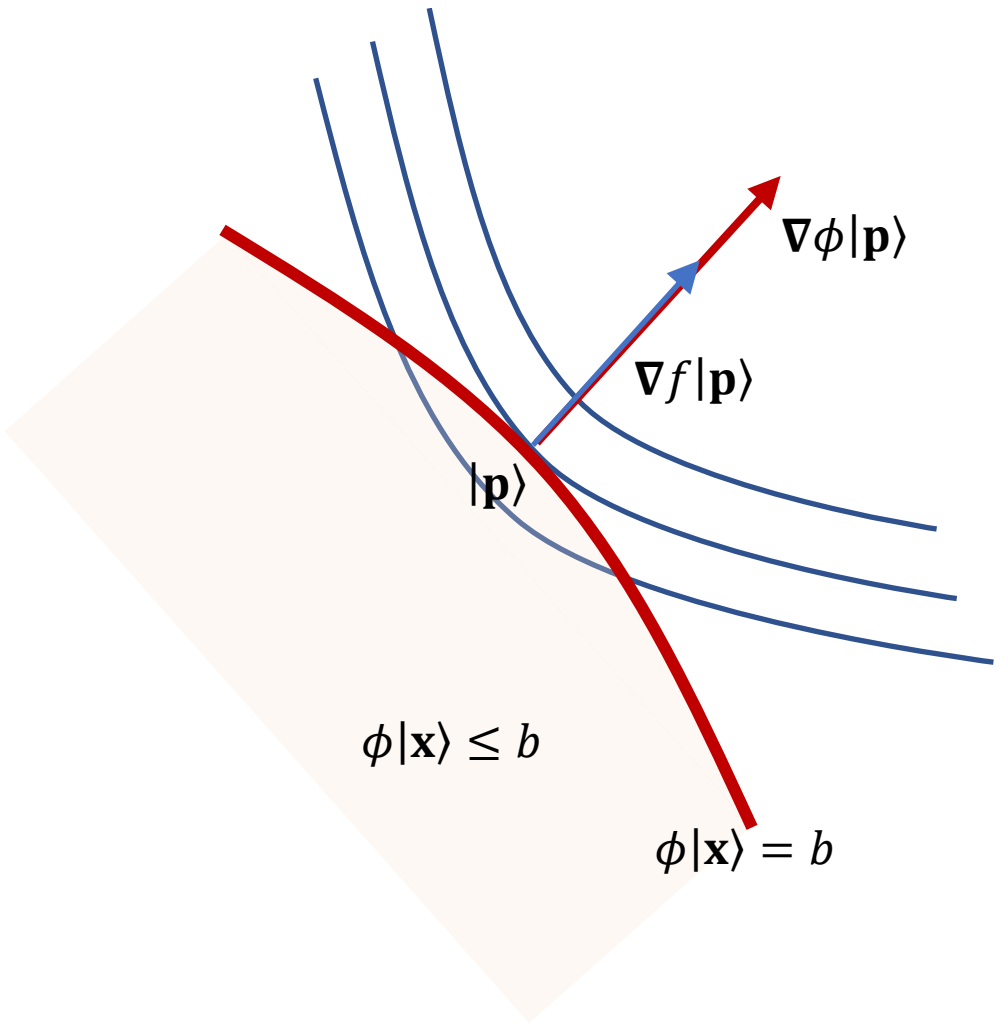

**Figure 10.13.** $\nabla f|\mathbf{p}\rangle$ and $\nabla\phi|\mathbf{p}\rangle$ in the same direction at the maximizer $|\mathbf{p}\rangle$.

Consequently, for the case in Figure 10.13, we still form the Lagrangian function

$$\mathcal{L}(x, y, \lambda) = f|\mathbf{x}\rangle - \lambda(\phi|\mathbf{x}\rangle - b), \tag{10.50}$$

using (10.49), we get

$$\frac{\partial\mathcal{L}}{\partial x} = \frac{\partial f}{\partial x} - \lambda\frac{\partial\phi}{\partial x} = 0 \quad \text{and} \quad \frac{\partial\mathcal{L}}{\partial y} = \frac{\partial f}{\partial y} - \lambda\frac{\partial\phi}{\partial y} = 0. \tag{10.51}$$

We require that the maximizer not be a critical point of the constraint function $\phi$.

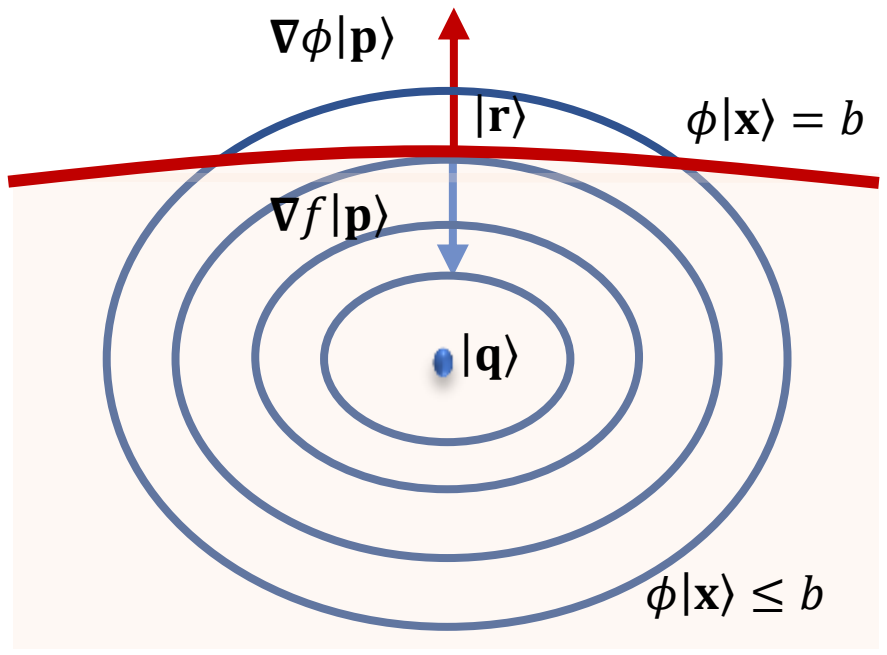

**Figure 10.14.** The constraint not binding.

Before we consider $\frac{\partial}{\partial\lambda}\mathcal{L}$, let us examine the following situation. Suppose that the maximum of $f$ on the constraint set $\phi|\mathbf{x}\rangle \leq b$ occurs not where $\phi|\mathbf{x}\rangle = b$ but at a point where $\phi|\mathbf{x}\rangle < b$. Figure 10.14 shows this case. The maximum of $f$ occurs at the point $|\mathbf{q}\rangle$ in the interior of the constraint set. There is a point $|\mathbf{r}\rangle$ on the level set $\phi|\mathbf{x}\rangle = b$ where this level set is tangent to a level set of $f$, but $\nabla f$ and $\nabla\phi$ point in opposite directions at $|\mathbf{r}\rangle$. In fact, one can increase the value of $f$ by moving further into the constraint set from $|\mathbf{r}\rangle$ until one reaches $|\mathbf{q}\rangle$. Since $|\mathbf{q}\rangle$ is in the interior of the constraint set, we say that the constraint is not binding (inactive, ineffective, loose) at $|\mathbf{q}\rangle$. Note that the point $|\mathbf{q}\rangle$ must be a local unconstrained max of $f$. Hence, we have

$$\frac{\partial}{\partial x}f|\mathbf{q}\rangle = 0 \quad \text{and} \quad \frac{\partial}{\partial y}f|\mathbf{q}\rangle = 0. \tag{10.52}$$

The derivatives of $\phi$ do not enter the criterion or the calculations at $|\mathbf{q}\rangle$. We can still use our Lagrangian

$$\mathcal{L}(x,y,\lambda) = f|\mathbf{x}\rangle - \lambda(\phi|\mathbf{x}\rangle - b), \tag{10.53}$$

and $\frac{\partial}{\partial x}\mathcal{L} = 0$, $\frac{\partial}{\partial y}\mathcal{L} = 0$, $\lambda = 0$. Setting $\lambda = 0$ causes the constraint function to drop out of the analysis; this is just what we want when the constraint is not binding at the max.

In summary, either the solution is on the boundary of the constraint set, that is, $\phi|\mathbf{x}\rangle - b = 0$ as in Figure 10.13. If so, the constraint is binding, there is a tangency at the solution, and the multiplier $\lambda$ must be $\geq 0$ . Or the solution is in the interior of the constraint set, as in Figure 10.14, in which case the multiplier $\lambda$ must be zero. If so, the constraint is not binding, and the solution is unaffected. A convenient way to summarize such a condition that one of two numbers be zero is to say that their product must be zero. Hence, the condition that either $\phi|\mathbf{x}\rangle - b = 0$ or $\lambda = 0$ can be represented as the following

$$\lambda(\phi|\mathbf{x}\rangle - b) = 0. \tag{10.54}$$

Since we do not know a priori whether or not the constraint will be binding at the maximizer, we cannot use the condition $\frac{\partial}{\partial\lambda}\mathcal{L} = 0$ that we used with equality constraints since this condition is equivalent to $\phi|\mathbf{x}\rangle - b = 0$. We will replace this statement with (10.54).

**Theorem 10.7.** Suppose that $f$ and $\phi$ are $C^1$ functions on $\mathbb{R}^n$ and that $|\mathbf{x}^*\rangle = (x^*, y^*)^T$ maximizes $f$ on the constraint set $\phi|\mathbf{x}\rangle \leq b$. If $\phi|\mathbf{x}^*\rangle = b$, suppose that

$$\frac{\partial}{\partial x}\phi|\mathbf{x}^*\rangle \neq 0 \quad \text{or} \quad \frac{\partial}{\partial y}\phi|\mathbf{x}^*\rangle \neq 0.$$

In any case, form the Lagrangian function

$$\mathcal{L}(x,y,\lambda) = f|\mathbf{x}\rangle - \lambda(\phi|\mathbf{x}\rangle - b).$$

Then, there is a multiplier $\lambda^*$ such that:

1. $\frac{\partial}{\partial x}\mathcal{L}(x^*, y^*, \lambda^*) = 0.$
2. $\frac{\partial}{\partial y}\mathcal{L}(x^*, y^*, \lambda^*) = 0.$
3. $\lambda^*(\phi(x^*, y^*) - b) = 0.$
4. $\lambda^* \geq 0.$
5. $\phi(x^*, y^*) \leq b.$

***Example 10.8***

$$\begin{aligned} &\text{maximize} && f|\mathbf{x}\rangle = xy, \qquad |\mathbf{x}\rangle = (x,y)^T \\ &\text{subject to} && \phi|\mathbf{x}\rangle = x^2 + y^2 \leq 1. \end{aligned}$$

***Solution***

The Lagrangian and conditions of Theorem 10.7 are

$$\begin{aligned} \mathcal{L}(x,y,\lambda) &= xy - \lambda(x^2 + y^2 - 1), \\ \frac{\partial}{\partial x}\mathcal{L} &= y - 2\lambda x = 0, \\ \frac{\partial}{\partial y}\mathcal{L} &= x - 2\lambda y = 0, \\ \lambda(x^2 + y^2 - 1) &= 0, \\ x^2 + y^2 &\leq 1, \\ \lambda &\geq 0. \end{aligned}$$

The $\frac{\partial}{\partial x}\mathcal{L} = \frac{\partial}{\partial y}\mathcal{L} = 0$, we have

$$\lambda = \frac{y}{2x} = \frac{x}{2y} \quad \text{or} \quad x^2 = y^2.$$

If $\lambda = 0$, then $x = y = 0$. This combination satisfies all the conditions, so it is a candidate for a solution.
If $\lambda \neq 0$, then the fourth equation becomes $x^2 + y^2 - 1 = 0$. Combining this with $x^2 = y^2$ , we find that $x^2 = y^2 = \frac{1}{2}$, or $x = \pm 1/\sqrt{2}$, $y = \pm 1/\sqrt{2}$. Combining these with the equation for $\lambda$ in $\lambda = \frac{y}{2x} = \frac{x}{2y}$, we find the following four candidates:

$$(x, y, \lambda) = \left(\frac{1}{\sqrt{2}}, \frac{1}{\sqrt{2}}, \frac{1}{\sqrt{2}}\right), \left(\frac{-1}{\sqrt{2}}, \frac{-1}{\sqrt{2}}, \frac{1}{\sqrt{2}}\right),$$

$$(x, y, \lambda) = \left(\frac{1}{\sqrt{2}}, \frac{-1}{\sqrt{2}}, \frac{-1}{\sqrt{2}}\right) = \left(\frac{-1}{\sqrt{2}}, \frac{1}{\sqrt{2}}, \frac{-1}{\sqrt{2}}\right).$$

The last two solutions are disregarded since the multiplier is negative (they are the solutions of minimizing $xy$ on the constraint $x^2 + y^2 \leq 1$). So, including $(0,0,0)$, there are three candidates which satisfy the conditions of Theorem 10.7. Plugging these three into the objective function, the two points with

$$x = y = \frac{1}{\sqrt{2}} \text{ and } x = y = -\frac{1}{\sqrt{2}}$$

are the solutions of our maximizing problem.

**Theorem 10.8:** Suppose that $f$, $\phi_1, \ldots, \phi_k$ are $C^1$ functions of $n$ variables. Suppose that $|\mathbf{x}^*\rangle \in \mathbb{R}^n$ is a local maximizer of $f$ by the $k$ inequalities

$$\phi_1|\mathbf{x}\rangle \leq b_1,$$
$$\ldots,$$
$$\phi_k|\mathbf{x}\rangle \leq b_k.$$

where $|\mathbf{x}\rangle = (x_1, \ldots, x_n)^T$. For ease of notation, assume that the first $k_0$ constraints are binding at $|\mathbf{x}^*\rangle$ and that the last $k - k_0$ constraints are not binding. Suppose that the following non-degenerate constraint qualification is satisfied at $|\mathbf{x}^*\rangle$.
The rank at $|\mathbf{x}^*\rangle$ of the Jacobian matrix of the binding constraints

$$\begin{pmatrix} \frac{\partial \phi_1|\mathbf{x}^*\rangle}{\partial x_1} & \cdots & \frac{\partial \phi_1|\mathbf{x}^*\rangle}{\partial x_n} \\ \vdots & \ddots & \vdots \\ \frac{\partial \phi_{k_0}|\mathbf{x}^*\rangle}{\partial x_1} & \cdots & \frac{\partial \phi_{k_0}|\mathbf{x}^*\rangle}{\partial x_n} \end{pmatrix},$$

is $k_0$—as large as it can be.
Form the Lagrangian

$$\mathcal{L}(x_1, \ldots x_n, \lambda_1, \ldots, \lambda_k) = f|\mathbf{x}\rangle - \lambda_1(\phi_1|\mathbf{x}\rangle - b_1) - \cdots - \lambda_k(\phi_k|\mathbf{x}\rangle - b_k).$$

Then, there exist multipliers $\lambda_1^*, \ldots, \lambda_k^*$ such that:

1. $\frac{\partial}{\partial x_1}\mathcal{L}(\mathbf{x}^*, \boldsymbol{\lambda}^*) = 0, \ldots, \frac{\partial}{\partial x_n}\mathcal{L}(\mathbf{x}^*, \boldsymbol{\lambda}^*) = 0$,
2. $\lambda_1^*(\phi_1|\mathbf{x}^*\rangle - b_1) = 0, \ldots, \lambda_k^*(\phi_k|\mathbf{x}^*\rangle - b_k) = 0$,
3. $\lambda_1^* \geq 0, \ldots, \lambda_k^* \geq 0$,
4. $\phi_1|\mathbf{x}^*\rangle \leq b_1, \ldots, \phi_k|\mathbf{x}^*\rangle \leq b_k$.

***Example 10.9***

$$\begin{aligned} \text{maximize} \quad & f|\mathbf{x}\rangle = xyz, \quad |\mathbf{x}\rangle = (x, y, z)^T \\ \text{subject to} \quad & \phi|\mathbf{x}\rangle = x + y + z \leq 1, \\ & x \geq 0, y \geq 0, z \geq 0. \end{aligned}$$

***Solution***

For consistent inequality constraints, we can write

$$-x \leq 0, \quad -y \leq 0, \quad -z \leq 0.$$

The Jacobian of the constraint functions is

$$\begin{pmatrix} 1 & 1 & 1 \\ -1 & 0 & 0 \\ 0 & -1 & 0 \\ 0 & 0 & -1 \end{pmatrix}.$$

Since its columns are linearly independent, it has rank three. Form the Lagrangian

$$\mathcal{L}(x,y,z,\lambda_1,\lambda_2,\lambda_3,\lambda_4) = xyz - \lambda_1(x+y+z-1) - \lambda_2(-x) - \lambda_3(-y) - \lambda_4(-z)$$
$$= xyz - \lambda_1(x+y+z-1) + \lambda_2 x + \lambda_3 y + \lambda_4 z.$$

The conditions of Theorem 10.8 are

$$\frac{\partial \mathcal{L}}{\partial x} = yz - \lambda_1 + \lambda_2 = 0 \quad \Rightarrow \quad \lambda_1 = yz + \lambda_2,$$
$$\frac{\partial \mathcal{L}}{\partial y} = xz - \lambda_1 + \lambda_3 = 0 \quad \Rightarrow \quad \lambda_1 = xz + \lambda_3,$$
$$\frac{\partial \mathcal{L}}{\partial z} = xy - \lambda_1 + \lambda_4 = 0 \quad \Rightarrow \quad \lambda_1 = xy + \lambda_4,$$
$$\lambda_1(x+y+z-1) = 0,$$
$$\lambda_2 x = 0,$$
$$\lambda_3 y = 0,$$
$$\lambda_4 z = 0,$$
$$\lambda_1 \geq 0,$$
$$\lambda_2 \geq 0,$$
$$\lambda_3 \geq 0,$$
$$\lambda_4 \geq 0,$$
$$x + y + z \leq 1,$$
$$x \geq 0,$$
$$y \geq 0,$$
$$z \geq 0.$$

We will look at two cases: $\lambda_1 = 0$ and $\lambda_1 > 0$.

If $\lambda_1 = 0$ in equation $\lambda_1 = yz + \lambda_2 = xz + \lambda_3 = xy + \lambda_4$, then because every variable in this equation is nonnegative,

$$yz = xz = xy = 0,$$
$$\lambda_1 = \lambda_2 = \lambda_3 = \lambda_4 = 0.$$

These equations lead to the (infinite) set of solution candidates in which two of the variables equal zero, and the third is any number in the interval [0,1]. In particular, the objective function equals zero for all $(x,y,z)$, which satisfy the above equations.

In the case $\lambda_1 > 0$. The solution of the constrained maximization problem is

$$x = y = z = \frac{1}{3}, \quad \lambda_1 = \frac{1}{9}, \quad f\left(\frac{1}{3},\frac{1}{3},\frac{1}{3}\right) = \frac{1}{27} > 0\ .$$

**Remark:**

Unlike the number of equality constraints, the number of inequality constraints, $m$, is not required to be less than $n$. For example, if we consider the case where all $\phi_j|\mathbf{x}\rangle$ for $1 \leq j \leq m$ are linear functions, then the constraints represent a polyhedron with $m$ facets, and the number of facets in such a polyhedron is obviously unlimited.

## 10.9 The Transformation to Equality Constraints

Now let us look to the solution of the problem in a different way. We have

$$\text{Minimize } f|\mathbf{x}\rangle, \qquad (10.55)$$
$$\text{subject to } \phi_j|\mathbf{x}\rangle \leq 0, j = 1,2,\ldots,m.$$

The inequality constraints in (10.55) can be transformed to equality constraints by adding slack variables. The idea is simple, by introducing positive slack variables $y_j$, $j = 1,\ldots,m$, such that the inequality constraints $\phi_j|\mathbf{x}\rangle \leq 0, j = 1,2,\ldots,m$ are transformed into equality constraints $G_j(\mathbf{x},\mathbf{y}) = 0, j = 1,\ldots,m$. Thus, the problem with inequality constraints has been replaced by a problem with equality constraints in which there are $n + m$ unknowns, ($x_i$ and $y_i$), and $m$ Lagrange multipliers $\lambda_i$.

Many choices of slack variables are available [2,6]; for instance, the positive variables $y_j, j = 1,\dots,m$ such that

$$G_j(\mathbf{x},\mathbf{y}) = \phi_j|\mathbf{x}\rangle + y_j = 0, j = 1,\dots,m. \tag{10.56}$$

***Example 10.10***

$$\begin{aligned} &\text{minimize} && f|\mathbf{x}\rangle = x_1^2 + x_2^2,\\ &\text{subject to} && x_1 \leq 1. \end{aligned}$$

***Solution***

The inequality constraint is

$$\phi|\mathbf{x}\rangle = x_1 - 1 \leq 0.$$

We can replace it with an equality constraint

$$G(\mathbf{x},\mathbf{y}) = x_1 - 1 + y_1 = 0, \;\; y_1 \geq 0.$$

The Lagrangian function becomes

$$\mathcal{L} = f + \lambda G = x_1^2 + x_2^2 + \lambda(x_1 - 1 + y_1).$$

At the optimum, the gradient of $\mathcal{L}$ is zero

$$\begin{pmatrix} 2x_1 + \lambda \\ 2x_2 \\ \lambda \end{pmatrix} = \begin{pmatrix} 0 \\ 0 \\ 0 \end{pmatrix}.$$

Hence, we have

$$\begin{aligned} 2x_1 + \lambda &= 0,\\ 2x_2 &= 0,\\ \lambda &= 0,\\ x_1 - 1 + y_1 &= 0, \end{aligned}$$

which gives $(x_1, x_2, y_1, \lambda) = (0,0,1,0)$. The condition $y_1 \geq 0$ is verified. The minimizer is inner to the domain.

The slack variable $y_j$ such that $y_j \geq 0$ can be simply replaced by the slack variable $y_j^2$. We have

$$G_j = \phi_j|\mathbf{x}\rangle + y_j^2 = 0, \qquad j = 1,2,\dots,m. \tag{10.57}$$

The problem now becomes

$$\begin{aligned} &\text{Minimize } f|\mathbf{x}\rangle,\\ &\text{subject to } G_j(\mathbf{x},\mathbf{y}) = g_j|\mathbf{x}\rangle + y_j^2 = 0, j = 1,2,\dots,m, \end{aligned} \tag{10.58}$$

where $|\mathbf{y}\rangle = (y_1, y_2, \dots, y_m)^T$ is the vector of slack variables.

***Example 10.11***

Consider again

$$\begin{aligned} &\text{minimize} && f|\mathbf{x}\rangle = x_1^2 + x_2^2,\\ &\text{subject to} && x_1 \leq 1. \end{aligned}$$

***Solution***

The inequality constraint is written as

$$\phi|\mathbf{x}\rangle = x_1 - 1 \leq 0,$$

which is replaced by an equality constraint

$$G(\mathbf{x},\mathbf{y}) = x_1 - 1 + y_1^2 = 0.$$

The Lagrangian function is thus equal to

$$\mathcal{L} = f + \lambda G = x_1^2 + x_2^2 + \lambda(x_1 - 1 + y_1^2).$$

At the optimum, the gradient of $\mathcal{L}$ is zero

$$\begin{pmatrix} 2x_1 + \lambda \\ 2x_2 \\ 2\lambda y_1 \end{pmatrix} = \begin{pmatrix} 0 \\ 0 \\ 0 \end{pmatrix}.$$

Thus, we get

$$2\lambda y_1 = 0,$$
$$x_1 - 1 + y_1^2 = 0,$$

which gives either $\lambda = 0$, or $y_1 = 0$.

- If $y_1 = 0$, we get an optimum at the limit of the domain $x_1 = 1$, $x_2 = 0$ with $\lambda = -2$. In this case, the function $f = 1$.
- If $\lambda = 0$, we get $x_1 = 0$, $x_2 = 0$, and $y_1 = \pm 1$. The solution is inside the domain, $f = 0$.
- Hence, the found minimum is thus $f = 0$ at point $(0,0)$.

Similar to the previously studied example, the optimum of the function $f|\mathbf{x}\rangle$, $|\mathbf{x}\rangle \in \mathbb{R}^n$ subject to $m$ inequality constraints $\phi_j|\mathbf{x}\rangle \le 0$ can be handled generally by selecting the slack variables in such a way that we are returned to a problem with $m$ equality constraints. Hence, the problem (10.55) can be solved conveniently by using the Lagrange multiplier method. For this, we construct the Lagrange function $\mathcal{L}$ as

$$\mathcal{L}(\mathbf{x}, \mathbf{y}, \boldsymbol{\lambda}) = f|\mathbf{x}\rangle + \sum_{j=1}^{m} \lambda_j G_j(\mathbf{x}, \mathbf{y}), \qquad (10.59)$$

where $|\boldsymbol{\lambda}\rangle = (\lambda_1, \lambda_2, \dots, \lambda_m)^T$ is the vector of Lagrange multipliers. In order to determine the stationary points of the Lagrange function, the following equations must be solved (necessary conditions):

$$\frac{\partial}{\partial x_i}\mathcal{L}(\mathbf{x}, \mathbf{y}, \boldsymbol{\lambda}) = 0, \qquad i = 1,2, \dots, n, \qquad (10.60.a)$$
$$\frac{\partial}{\partial \lambda_j}\mathcal{L}(\mathbf{x}, \mathbf{y}, \boldsymbol{\lambda}) = 0, \qquad j = 1,2, \dots, m, \qquad (10.60.b)$$
$$\frac{\partial}{\partial y_j}\mathcal{L}(\mathbf{x}, \mathbf{y}, \boldsymbol{\lambda}) = 0, \qquad j = 1,2, \dots, m, \qquad (10.60.c)$$

which corresponding to the following set of equations:

$$\frac{\partial}{\partial x_i} f|\mathbf{x}\rangle + \sum_{j=1}^{m} \lambda_j \frac{\partial \phi_j|\mathbf{x}\rangle}{\partial x_i} = 0, \qquad i = 1,2, \dots, n, \qquad (10.61.a)$$
$$\phi_j|\mathbf{x}\rangle + y_j^2 = 0, \qquad j = 1,2, \dots, m, \qquad (10.61.b)$$
$$2\lambda_j y_j = 0, \qquad j = 1,2, \dots, m. \qquad (10.61.c)$$

It can be seen that (10.60) or (10.61) represent $(n + 2m)$ equations in the $(n + 2m)$ unknowns, $|\mathbf{x}\rangle$, $|\boldsymbol{\lambda}\rangle$ and $|\mathbf{y}\rangle$. The solution of (10.60) or (10.61) thus gives the optimum solution vector $|\mathbf{x}^*\rangle$, the Lagrange multiplier vector, $|\boldsymbol{\lambda}^*\rangle$, and the slack variable vector, $|\mathbf{y}^*\rangle$. Equations (10.61.b) ensure that the constraints $\phi_i|\mathbf{x}\rangle \le 0$, $j = 1,2,\dots,m$, are satisfied, while (10.61.c) imply that either $\lambda_j = 0$ or $y_j = 0$.

- If $\lambda_j = 0$, it means that the $j$th constraint is not binding and can be ignored.
- If $y_j = 0$, it means that the constraint is binding ($\phi_j = 0$) at the optimum point.

Consider the separation of the constraints into two groups, $J_1$ and $J_2$, where $J_1 + J_2$ stands for the total set of constraints. Let the group $J_1$ represent the indices of the constraints that are binding at the optimum point and $J_2$ represent the indices of all the not binding constraints. Thus for $j \in J_1$, $y_j = 0$ (constraints are binding), for $j \in J_2$, $\lambda_j = 0$ (constraints are not binding), and (10.61.a) becomes

$$\frac{\partial f}{\partial x_i} + \sum_{j \in J_1} \lambda_j \frac{\partial \phi_j}{\partial x_i} = 0, \qquad i = 1,2, \dots, n. \qquad (10.62)$$

Similarly, (10.61.b) can be written as

$$\phi_j|\mathbf{x}\rangle = 0, \quad j \in J_1, \qquad (10.63.a)$$
$$\phi_j|\mathbf{x}\rangle + y_j^2 = 0, \quad j \in J_2. \qquad (10.63.b)$$

Equations (10.62) and (10.63) are $n+p+(m-p)=n+m$ equations in the $n+m$ unknowns $x_i$ $(i=1,2,\dots,n)$, $\lambda_j$ $(j\in J_1)$, and $y_j$ $(j\in J_2)$, where $p$ denotes the number of active constraints.

Assuming that the first $p$ constraints are binding, (10.62) can be expressed as

$$-\frac{\partial f}{\partial x_i}=\lambda_1\frac{\partial\phi_1}{\partial x_1}+\lambda_2\frac{\partial\phi_2}{\partial x_2}+\cdots+\lambda_p\frac{\partial\phi_p}{\partial x_p},\qquad i=1,2,\dots,n. \tag{10.64}$$

These equations can be written collectively as

$$-\boldsymbol{\nabla} f=\lambda_1\boldsymbol{\nabla}\phi_1+\lambda_2\boldsymbol{\nabla}\phi_2+\cdots+\lambda_p\boldsymbol{\nabla}\phi_p, \tag{10.65}$$

where $\boldsymbol{\nabla} f$ and $\boldsymbol{\nabla}\phi_j$ are the gradients of the objective function and the $j$th constraint, respectively:

$$\boldsymbol{\nabla} f=\begin{pmatrix}\frac{\partial f}{\partial x_1}\\ \frac{\partial f}{\partial x_2}\\ \vdots\\ \frac{\partial f}{\partial x_n}\end{pmatrix}\quad\text{and}\quad\boldsymbol{\nabla}\phi_j=\begin{pmatrix}\frac{\partial \phi_j}{\partial x_1}\\ \frac{\partial \phi_j}{\partial x_2}\\ \vdots\\ \frac{\partial \phi_j}{\partial x_n}\end{pmatrix}. \tag{10.66}$$

**Remarks:**

- According to (10.65), the negative of the gradient of the objective function can be written as a linear combination of the gradients of the active constraints at the optimum point.
- Moreover, we can demonstrate that in the case of a minimization problem, the values of $\lambda_j$ $(j\in J_1)$ must be positive.

  Let us assume, for the sake of simplicity, that the optimum point has only two active constraints ($p=2$) at the optimum point. Then (10.65) becomes

  $$-\boldsymbol{\nabla} f=\lambda_1\boldsymbol{\nabla}\phi_1+\lambda_2\boldsymbol{\nabla}\phi_2. \tag{10.67}$$

  Let $|\mathbf{s}\rangle$ be a feasible direction at the optimum point. Pre-multiplying (10.67) by $\langle\mathbf{s}|$ gives us

  $$-\langle\mathbf{s}|\boldsymbol{\nabla} f\rangle=\lambda_1\langle\mathbf{s}|\boldsymbol{\nabla}\phi_1\rangle+\lambda_2\langle\mathbf{s}|\boldsymbol{\nabla}\phi_2\rangle. \tag{10.68}$$

  Since $|\mathbf{s}\rangle$ is a feasible direction, it should satisfy the relations

  $$\langle\mathbf{s}|\boldsymbol{\nabla}\phi_1\rangle<0, \tag{10.69.a}$$
  $$\langle\mathbf{s}|\boldsymbol{\nabla}\phi_2\rangle<0. \tag{10.69.b}$$

  Thus if $\lambda_1>0$ and $\lambda_2>0$, the quantity $\langle\mathbf{s}|\boldsymbol{\nabla} f\rangle$ can always be seen to be positive. Since $\boldsymbol{\nabla} f$ represents the gradient direction, along which the value of the function increases at the maximum rate, $\langle\mathbf{s}|\boldsymbol{\nabla} f\rangle$ represents the component of the increment of $f$ along the direction $|\mathbf{s}\rangle$. If $\langle\mathbf{s}|\boldsymbol{\nabla} f\rangle>0$, the function value increases as we move along the direction $|\mathbf{s}\rangle$. As a result, if $\lambda_1$ and $\lambda_2$ are positive, we will be unable to identify any direction in the feasible domain along which the function value can be decreased further. Because the point at which (10.69) is valid is assumed to be optimum, $\lambda_1$ and $\lambda_2$ have to be positive. This reasoning can be extended to cases where there are more than two constraints active.

## 10.10 Karush-Kuhn-Tucker Condition

We consider the following problem:

$$\begin{aligned}&\text{minimize} && f|\mathbf{x}\rangle\\ &\text{subject to} && \mathbf{h}|\mathbf{x}\rangle=|\mathbf{0}\rangle,\\ & && \boldsymbol{\phi}|\mathbf{x}\rangle\le|\mathbf{0}\rangle,\end{aligned} \tag{10.70.a}$$

where,

$$|\mathbf{x}\rangle\in\mathbb{R}^n \tag{10.70.b}$$
$$f:\mathbb{R}^n\to\mathbb{R} \tag{10.70.c}$$

$$\mathbf{h}: \mathbb{R}^n \to \mathbb{R}^m, \mathbf{h}|\mathbf{x}\rangle = \begin{pmatrix} h_1|\mathbf{x}\rangle \\ \vdots \\ h_m|\mathbf{x}\rangle \end{pmatrix}, m \leq n \quad (10.70.d)$$

$$\boldsymbol{\phi}: \mathbb{R}^n \to \mathbb{R}^p, \boldsymbol{\phi}|\mathbf{x}\rangle = \begin{pmatrix} \phi_1|\mathbf{x}\rangle \\ \vdots \\ \phi_p|\mathbf{x}\rangle \end{pmatrix}. \quad (10.70.e)$$

For the above general problem, we adopt the following definitions.

**Definition (Active and Inactive Constraints):** An inequality constraint $\phi_j|\mathbf{x}\rangle \leq 0$ is said to be active at $|\mathbf{x}^*\rangle$ if $\phi_j|\mathbf{x}^*\rangle = 0$. It is inactive at $|\mathbf{x}^*\rangle$ if $\phi_j|\mathbf{x}^*\rangle < 0$.

By convention, we consider an equality constraint $h_i|\mathbf{x}\rangle = 0$ to be always active.

**Definition (Regular Point):** Let $|\mathbf{x}^*\rangle$ satisfy $\mathbf{h}|\mathbf{x}^*\rangle = |\mathbf{0}\rangle$, $\boldsymbol{\phi}|\mathbf{x}^*\rangle \leq |\mathbf{0}\rangle$, and let $J(\mathbf{x}^*)$ be the index set of active inequality constraints, that is,

$$J(\mathbf{x}^*) = \{j: \phi_j|\mathbf{x}^*\rangle = 0\} \quad (10.71)$$

Then, we say that $|\mathbf{x}^*\rangle$ is a regular point if the vectors

$$\nabla h_i|\mathbf{x}^*\rangle, \qquad \nabla\phi_i|\mathbf{x}^*\rangle, \qquad 1 \leq i \leq m, j \in J(\mathbf{x}^*) \quad (10.72)$$

are linearly independent.

For a point to be a local minimizer, we now demonstrate a first-order necessary condition. This condition is known as the Karush-Kuhn-Tucker (KKT) condition [9,10].

**Theorem 10.9 (Karush-Kuhn-Tucker (KKT) Theorem):** Let $f, \mathbf{h}, \boldsymbol{\phi} \in C^1$. Let $|\mathbf{x}^*\rangle$ be a regular point and a local minimizer for the problem of minimizing $f$ subject to $\mathbf{h}|\mathbf{x}\rangle = |\mathbf{0}\rangle, \boldsymbol{\phi}|\mathbf{x}\rangle \leq |\mathbf{0}\rangle$. Then, there exist $|\boldsymbol{\lambda}^*\rangle \in \mathbb{R}^m$ and $|\boldsymbol{\mu}^*\rangle \in \mathbb{R}^p$ such that
1. $|\boldsymbol{\mu}^*\rangle \geq |\mathbf{0}\rangle$;
2. $Df|\mathbf{x}^*\rangle + \langle\boldsymbol{\lambda}^*|(D\mathbf{h}|\mathbf{x}^*\rangle) + \langle\boldsymbol{\mu}^*|(D\boldsymbol{\phi}|\mathbf{x}^*\rangle) = \langle\mathbf{0}|$,
or $\nabla f|\mathbf{x}^*\rangle + (D\mathbf{h}|\mathbf{x}^*\rangle)^T|\boldsymbol{\lambda}^*\rangle + (D\boldsymbol{\phi}|\mathbf{x}^*\rangle)^T|\boldsymbol{\mu}^*\rangle = |\mathbf{0}\rangle$,
or $\nabla f|\mathbf{x}^*\rangle + \sum_{i=1}^{m}\lambda_i^*\nabla h_i|\mathbf{x}^*\rangle + \sum_{i=1}^{p}\mu_i^*\nabla\phi_i|\mathbf{x}^*\rangle = |\mathbf{0}\rangle$.
3. $\langle\boldsymbol{\mu}^*|\boldsymbol{\phi}(\mathbf{x}^*)\rangle = 0$.

The vector $|\boldsymbol{\lambda}^*\rangle$ and $|\boldsymbol{\mu}^*\rangle$ are called the Lagrange multiplier vector and the KKT multiplier vector, respectively. Additionally, we call their components Lagrange multipliers and KKT multipliers, respectively.

**Remark:**

1- Observe that $\mu_j^* \geq 0$ (by condition 1) and $\phi_j|\mathbf{x}^*\rangle \leq 0$. Therefore, the condition

$$\langle\boldsymbol{\mu}^*|\boldsymbol{\phi}(\mathbf{x}^*)\rangle = \mu_1^*\phi_1|\mathbf{x}^*\rangle + \cdots + \mu_p^*\phi_p|\mathbf{x}^*\rangle = 0 \quad (10.73)$$

implies that if $\phi_j|\mathbf{x}^*\rangle < 0$, then $\mu_j^* = 0$, that is, for all $j \notin J(\mathbf{x}^*)$, we have $\mu_j^* = 0$. In other words, the KKT multipliers $\mu_j^*$ corresponding to inactive constraints are zero. The other KKT multipliers, $\mu_i^*$, $i \in J(\mathbf{x}^*)$, are nonnegative; they may or may not be equal to zero.

2-

The KKT condition is composed of five parts two inequalities and three equations:
1. $|\boldsymbol{\mu}^*\rangle \geq |\mathbf{0}\rangle$;

2. $Df|\mathbf{x}^*\rangle + \langle\boldsymbol{\lambda}^*|(D\mathbf{h}|\mathbf{x}^*\rangle) + \langle\boldsymbol{\mu}^*|(D\boldsymbol{\phi}|\mathbf{x}^*\rangle) = \langle\mathbf{0}|$,
or $\nabla f|\mathbf{x}^*\rangle + (D\mathbf{h}|\mathbf{x}^*\rangle)^T|\boldsymbol{\lambda}^*\rangle + (D\boldsymbol{\phi}|\mathbf{x}^*\rangle)^T|\boldsymbol{\mu}^*\rangle = |\mathbf{0}\rangle$,
or $\nabla f|\mathbf{x}^*\rangle + \sum_{i=1}^{m}\lambda_i^*\nabla h_i|\mathbf{x}^*\rangle + \sum_{i=1}^{p}\mu_i^*\nabla\phi_i|\mathbf{x}^*\rangle = |\mathbf{0}\rangle$.

3. $\langle\boldsymbol{\mu}^*|\boldsymbol{\phi}(\mathbf{x}^*)\rangle = 0$;
4. $\mathbf{h}|\mathbf{x}^*\rangle = |\mathbf{0}\rangle$;
5. $\boldsymbol{\phi}|\mathbf{x}^*\rangle \leq |\mathbf{0}\rangle$.

When the inequality constraint has the form $\boldsymbol{\phi}(\mathbf{x}) \geq 0$, we can deduce the KKT condition. Specifically, consider the problem

$$\begin{aligned} &\text{minimize} && f|\mathbf{x}\rangle \\ &\text{subject to} && \mathbf{h}|\mathbf{x}\rangle = |\mathbf{0}\rangle, \\ &&& \boldsymbol{\phi}|\mathbf{x}\rangle \geq |\mathbf{0}\rangle, \end{aligned} \tag{10.74}$$

We multiply the inequality constraint function by $-1$ to obtain $-\boldsymbol{\phi}|\mathbf{x}\rangle \leq |\mathbf{0}\rangle$. Thus, the KKT condition for this case is $Df|\mathbf{x}^*\rangle + \langle\boldsymbol{\lambda}^*|(D\mathbf{h}|\mathbf{x}^*\rangle) + \langle\boldsymbol{\mu}^*|(D\boldsymbol{\phi}|\mathbf{x}^*\rangle) = \langle\mathbf{0}|$

1. $|\boldsymbol{\mu}^*\rangle \geq |\mathbf{0}\rangle$;
2. $Df|\mathbf{x}^*\rangle + \langle\boldsymbol{\lambda}^*|(D\mathbf{h}|\mathbf{x}^*\rangle) - \langle\boldsymbol{\mu}^*|(D\boldsymbol{\phi}|\mathbf{x}^*\rangle) = \langle\mathbf{0}|$;
3. $\langle\boldsymbol{\mu}^*|\boldsymbol{\phi}(\mathbf{x}^*)\rangle = 0$;
4. $\mathbf{h}|\mathbf{x}^*\rangle = |\mathbf{0}\rangle$;
5. $\boldsymbol{\phi}|\mathbf{x}^*\rangle \geq |\mathbf{0}\rangle$.

Changing the sign of $\boldsymbol{\mu}^*$, we obtain

1. $|\boldsymbol{\mu}^*\rangle \leq |\mathbf{0}\rangle$;
2. $Df|\mathbf{x}^*\rangle + \langle\boldsymbol{\lambda}^*|(D\mathbf{h}|\mathbf{x}^*\rangle) + \langle\boldsymbol{\mu}^*|(D\boldsymbol{\phi}|\mathbf{x}^*\rangle) = \langle\mathbf{0}|$;
3. $\langle\boldsymbol{\mu}^*|\boldsymbol{\phi}(\mathbf{x}^*)\rangle = 0$;
4. $\mathbf{h}|\mathbf{x}^*\rangle = |\mathbf{0}\rangle$;
5. $\boldsymbol{\phi}|\mathbf{x}^*\rangle \geq |\mathbf{0}\rangle$.

***Example 10.12***

$$\begin{aligned} &\text{minimize} && f|\mathbf{x}\rangle = x_1^2 + x_2^2 - 14x_1 - 6x_2 - 7 = (x_1-7)^2 + (x_2-3)^2 - 65 \\ &\text{subject to} && \phi_1|\mathbf{x}\rangle = x_1 + x_2 \leq 2 \\ &&& \phi_2|\mathbf{x}\rangle = x_1 + 2x_2 \leq 3 \end{aligned}$$

where $|\mathbf{x}\rangle = (x_1, x_2)^T$.

***Solution***

We can write the objective function and inequality constraints as

$$\begin{aligned} f|\mathbf{x}\rangle &= (x_1-7)^2 + (x_2-3)^2 - 65, \\ \phi_1|\mathbf{x}\rangle &= x_1 + x_2 - 2 \leq 0, \\ \phi_2|\mathbf{x}\rangle &= x_1 + 2x_2 - 3 \leq 0, \end{aligned}$$

with

$$\begin{aligned} \nabla f|\mathbf{x}\rangle &= (2(x_1-7), 2(x_2-3))^T, \\ \nabla\phi_1|\mathbf{x}\rangle &= (1,1)^T, \\ \nabla\phi_2|\mathbf{x}\rangle &= (1,2)^T, \end{aligned}$$

such that

$$\nabla f|\mathbf{x}^*\rangle + \sum_{i=1}^{m} \lambda_i^* \nabla h_i|\mathbf{x}^*\rangle + \sum_{i=1}^{p} \mu_i^* \nabla\phi_i|\mathbf{x}^*\rangle = |\mathbf{0}\rangle \Rightarrow \begin{pmatrix} 2(x_1-7) \\ 2(x_2-3) \end{pmatrix} + \mu_1\begin{pmatrix} 1 \\ 1 \end{pmatrix} + \mu_2\begin{pmatrix} 1 \\ 2 \end{pmatrix} = \begin{pmatrix} 0 \\ 0 \end{pmatrix}.$$

The Necessary (KKT) conditions are

$$\begin{aligned} 2x_1 - 14 + \mu_1 + \mu_2 &= 0, \\ 2x_2 - 6 + \mu_1 + 2\mu_2 &= 0, \\ \mu_1(x_1 + x_2 - 2) &= 0, \\ \mu_2(x_1 + 2x_2 - 3) &= 0, \\ \mu_1 &\geq 0, \\ \mu_2 &\geq 0, \end{aligned}$$

**Case 1:** $\mu_1 = 0,\ \mu_2 = 0,$

$$x_1 = \frac{14}{2} = 7, x_2 = \frac{6}{2} = 3. \text{ Violated constraints.}$$

**Case 2:** $\mu_1 > 0,\ \mu_2 = 0,$

$$2x_1 + \mu_1 = 14,$$
$$2x_2 + \mu_1 = 6,$$
$$x_1 + x_2 = 2.$$

The solution of the above system of equations is

$$x_1 = 3, x_2 = -1, \mu_1 = 8, \mu_2 = 0,$$

with $x_1 + 2x_2 = 3 - 2 < 3$ ⇒optimal solution

**Case 3:** $\mu_1 = 0,\ \mu_2 > 0$

$$2x_1 + \mu_2 = 14,$$
$$2x_2 + 2\mu_2 = 6,$$
$$x_1 + 2x_2 = 3,$$

The solution of the above system of equations is

$$x_1 = 5, x_2 = -1, \mu_1 = 0, \mu_2 = 4,$$

with $x_1 + x_2 = 5 - 1 = 4 > 2 \Rightarrow$ not optimal solution

**Case 4:** $\mu_1 > 0,\ \mu_2 > 0$

$$2x_1 + \mu_1 + \mu_2 = 14,$$
$$2x_2 + \mu_1 + 2\mu_2 = 6,$$
$$x_1 + x_2 = 2,$$
$$x_1 + 2x_2 = 3.$$

The solution of the above system of equations is

$$x_1 = 1, x_2 = 1, \mu_1 = 20, \mu_2 = -8 \Rightarrow \text{contradiction}$$

The feasible point satisfying the KKT condition is only a candidate for a minimizer. However, there is no guarantee that the point is indeed a minimizer, because the KKT condition is, in general, only necessary. A sufficient condition for a point to be a minimizer is given in the next section.

## 10.11 Second-Order Conditions

Just like in the case of equality constraints, we can also provide second-order necessary and sufficient conditions for extremum problems with inequality constraints [4]. In order to do this, we must define the following matrix:

$$\boldsymbol{H}_{\mathcal{L}}(\mathbf{x}, \boldsymbol{\lambda}, \boldsymbol{\mu}) = \boldsymbol{H}_f(\mathbf{x}) + \boldsymbol{\Sigma} + \acute{\boldsymbol{\Sigma}} \tag{10.75}$$

where $\boldsymbol{H}_f(\mathbf{x})$ is the Hessian matrix of $f$ at $|\mathbf{x}\rangle$, and the notation $\boldsymbol{\Sigma}$ represents

$$\boldsymbol{\Sigma} = \lambda_1 \boldsymbol{H}_1(\mathbf{x}) + \cdots + \lambda_m \boldsymbol{H}_m(\mathbf{x}), \tag{10.76}$$

where $\boldsymbol{H}_i(\mathbf{x})$ is the Hessian matrix of $h_i$ at $|\mathbf{x}\rangle$, $(i = 1, \dots, m)$. Similarly, the notation $\acute{\boldsymbol{\Sigma}}$ represents

$$\acute{\boldsymbol{\Sigma}} = \mu_1 \acute{\boldsymbol{H}}_1(\mathbf{x}) + \cdots + \mu_p \acute{\boldsymbol{H}}_p(\mathbf{x}) \tag{10.77}$$

where $\acute{\boldsymbol{H}}_i(\mathbf{x})$ is the Hessian of $\phi_i$ at $|\mathbf{x}\rangle$, $(i = 1, \dots, p)$ given by

$$\acute{\boldsymbol{H}}_i(\mathbf{x}) = \begin{pmatrix} \frac{\partial^2 \phi_i}{\partial x_1^2}|\mathbf{x}\rangle & \cdots & \frac{\partial^2 \phi_i}{\partial x_n \partial x_1}|\mathbf{x}\rangle \\ \vdots & & \vdots \\ \frac{\partial^2 \phi_i}{\partial x_1 \partial x_n}|\mathbf{x}\rangle & \cdots & \frac{\partial^2 \phi_i}{\partial x_n^2}|\mathbf{x}\rangle \end{pmatrix} \tag{10.78}$$

In the following theorem, we use

$$\hat{T}(\mathbf{x}^*) = \{|\mathbf{y}\rangle \in \mathbb{R}^n : D\mathbf{h}(\mathbf{x}^*)|\mathbf{y}\rangle = |\mathbf{0}\rangle, D\phi_i(\mathbf{x}^*)|\mathbf{y}\rangle = 0, j \in J(\mathbf{x}^*)\} \tag{10.79}$$

that is, the tangent space to the surface is defined by active constraints.

**Theorem 10.10 (Second-Order Necessary Conditions):** Let $|\mathbf{x}^*\rangle$ be a local minimizer of

$$f: \mathbb{R}^n \to \mathbb{R}$$

subject to

$$\mathbf{h}|\mathbf{x}\rangle = |\mathbf{0}\rangle,$$
$$\mathbf{\Phi}|\mathbf{x}\rangle \le |\mathbf{0}\rangle,$$
$$\mathbf{h}: \mathbb{R}^n \to \mathbb{R}^m, \quad m \le n,$$
$$\mathbf{g}: \mathbb{R}^n \to \mathbb{R}^p,$$

and $f$, $\mathbf{h}$, $\mathbf{\Phi} \in C^2$. Suppose $|\mathbf{x}^*\rangle$ is regular. Then, there exist $|\boldsymbol{\lambda}^*\rangle \in \mathbb{R}^m$ and $|\boldsymbol{\mu}^*\rangle \in \mathbb{R}^p$ such that:

1. $|\boldsymbol{\mu}^*\rangle \ge |\mathbf{0}\rangle$,
2. $Df|\mathbf{x}^*\rangle + \langle\boldsymbol{\lambda}^*|(D\mathbf{h}|\mathbf{x}^*\rangle) + \langle\boldsymbol{\mu}^*|(D\boldsymbol{\phi}|\mathbf{x}^*\rangle) = \langle\mathbf{0}|$,
   or $\nabla f|\mathbf{x}^*\rangle + (D\mathbf{h}|\mathbf{x}^*\rangle)^T|\boldsymbol{\lambda}^*\rangle + (D\boldsymbol{\phi}|\mathbf{x}^*\rangle)^T|\boldsymbol{\mu}^*\rangle = |\mathbf{0}\rangle$,
   or $\nabla f|\mathbf{x}^*\rangle + \sum_{i=1}^{m}\lambda_i^*\nabla h_i|\mathbf{x}^*\rangle + \sum_{i=1}^{p}\mu_i^*\nabla\phi_i|\mathbf{x}^*\rangle = |\mathbf{0}\rangle$.
3. $\langle\boldsymbol{\mu}^*|\boldsymbol{\phi}(\mathbf{x}^*)\rangle = 0$; and
4. For all $|\mathbf{y}\rangle \in T(\mathbf{x}^*)$ we have $\langle\mathbf{y}|\boldsymbol{H}_{\mathcal{L}}(\mathbf{x}^*, \boldsymbol{\lambda}^*, \boldsymbol{\mu}^*)|\mathbf{y}\rangle \ge 0$.

Let us define the following set:

$$\hat{T}(\mathbf{x}^*, \boldsymbol{\mu}^*) = \{|\mathbf{y}\rangle: D\mathbf{h}(\mathbf{x}^*)|\mathbf{y}\rangle = |\mathbf{0}\rangle, D\phi_i(\mathbf{x}^*)|\mathbf{y}\rangle = 0, i \in \bar{J}(\mathbf{x}^*, \boldsymbol{\mu}^*)\} \quad (10.80)$$

where $\bar{J}(\mathbf{x}^*, \boldsymbol{\mu}^*) = \{i: \phi_i|\mathbf{x}^*\rangle = 0, \mu_i^* > 0\}$. Note that $\bar{J}(\mathbf{x}^*, \boldsymbol{\mu}^*)$ is a subset of $J(\mathbf{x}^*)$, that is, $\bar{J}(\mathbf{x}^*, \boldsymbol{\mu}^*) \subset J(\mathbf{x}^*)$. This, in turn, implies that $T(\mathbf{x}^*)$ is a subset of $\hat{T}(\mathbf{x}^*, \boldsymbol{\mu}^*)$, that is, $T(\mathbf{x}^*) \subset \hat{T}(\mathbf{x}^*, \boldsymbol{\mu}^*)$.

**Theorem 10.11 (Second-Order Sufficient Conditions):** Suppose $f$, $\mathbf{h}$, $\boldsymbol{\phi} \in C^2$ and there exists a feasible point $|\mathbf{x}^*\rangle \in \mathbb{R}^n$ and vectors $|\boldsymbol{\lambda}^*\rangle \in \mathbb{R}^m$ and $|\boldsymbol{\mu}^*\rangle \in \mathbb{R}^p$, such that:

1. $|\boldsymbol{\mu}^*\rangle \ge |\mathbf{0}\rangle$,
2. $Df|\mathbf{x}^*\rangle + \langle\boldsymbol{\lambda}^*|(D\mathbf{h}|\mathbf{x}^*\rangle) + \langle\boldsymbol{\mu}^*|(D\boldsymbol{\phi}|\mathbf{x}^*\rangle) = \langle\mathbf{0}|$,
   or $\nabla f|\mathbf{x}^*\rangle + (D\mathbf{h}|\mathbf{x}^*\rangle)^T|\boldsymbol{\lambda}^*\rangle + (D\boldsymbol{\phi}|\mathbf{x}^*\rangle)^T|\boldsymbol{\mu}^*\rangle = |\mathbf{0}\rangle$,
   or $\nabla f|\mathbf{x}^*\rangle + \sum_{i=1}^{m}\lambda_i^*\nabla h_i|\mathbf{x}^*\rangle + \sum_{i=1}^{p}\mu_i^*\nabla\phi_i|\mathbf{x}^*\rangle = |\mathbf{0}\rangle$.
3. $\langle\boldsymbol{\mu}^*|\boldsymbol{\phi}(\mathbf{x}^*)\rangle = 0$; and
4. For all $|\mathbf{y}\rangle \in \hat{T}(\mathbf{x}^*, \boldsymbol{\mu}^*)$, $|\mathbf{y}\rangle \ne |\mathbf{0}\rangle$ we have $\langle\mathbf{y}|\boldsymbol{H}_{\mathcal{L}}(\mathbf{x}^*, \boldsymbol{\lambda}^*, \boldsymbol{\mu}^*)|\mathbf{y}\rangle > 0$.
5. Then, $|\mathbf{x}^*\rangle$ is a strict local minimizer of $f$ subject to $\mathbf{h}|\mathbf{x}\rangle = |\mathbf{0}\rangle$, $\boldsymbol{\phi}|\mathbf{x}\rangle \le |\mathbf{0}\rangle$.

***Example 10.13***

$$\begin{aligned} &\text{minimize} && f|\mathbf{x}\rangle = (x_1 - 1)^2 + x_2^2 \\ &\text{subject to} && \phi|\mathbf{x}\rangle = -x_1 + x_2^2 \ge 0 \end{aligned}$$

Suppose we want to verify whether $|\mathbf{x}^*\rangle = (0,0)^T$ is optimal.

***Solution***

$$\nabla f|\mathbf{x}\rangle = (2(x_1 - 1), 2x_2)^T,$$
$$\nabla\phi|\mathbf{x}\rangle = (-1, 2x_2)^T,$$
$$J = \{1\}.$$

Since $\nabla\phi|\mathbf{x}^*\rangle = (-1,0)^T$ is linearly independent, the constraint qualification is satisfied at $|\mathbf{x}^*\rangle$. The first-order KKT conditions are given by

$$2(x_1 - 1) + \mu_1 = 0,$$
$$2x_2 - 2x_2\mu_1 = 0,$$
$$\mu_1(-x_1 + x_2^2) = 0,$$
$$\mu_1 \ge 0.$$

Here $|\mathbf{x}^*\rangle = (0,0)^T$ and $\mu_1^* = 2$ satisfy the above conditions. Hence, $(\mathbf{x}^*, \mu_1^*)^T = (0,0,2)^T$ is KKT point and $|\mathbf{x}^*\rangle$ satisfies Theorem 10.9. Let us now apply the second-order necessary conditions to test whether $(0,0)^T$ is a local minimum to the problem. The first part of Theorem 10.10 is already satisfied, since $(\mathbf{x}^*, \mu_1^*)^T = (0,0,2)^T$ is a KKT point. To prove the second-order conditions, compute

$$\boldsymbol{H}_{\mathcal{L}}(\mathbf{x}, \mu) = \begin{pmatrix} 2 & 0 \\ 0 & 2 - 2\mu_1 \end{pmatrix}.$$

At $(\mathbf{x}^*, \mu^*)$,

$$\boldsymbol{H}_{\mathcal{L}}(\mathbf{x}^*, \mu^*) = \begin{pmatrix} 2 & 0 \\ 0 & -4 \end{pmatrix}.$$

We need to determine if

$$\left\langle \mathbf{y} \middle| \begin{pmatrix} 2 & 0 \\ 0 & -4 \end{pmatrix} \middle| \mathbf{y} \right\rangle \geq 0,$$

for all $|\mathbf{y}\rangle$ satisfying

$$D\phi_1(\mathbf{x}^*)|\mathbf{y}\rangle = 0 \quad \Rightarrow \quad (-1,0)\begin{pmatrix} y_1 \\ y_2 \end{pmatrix} = 0.$$

Simply, we need to consider only vectors $(y_1, y_2)^T \equiv (0, y_2)^T$ to satisfy $\left\langle \mathbf{y} \middle| \begin{pmatrix} 2 & 0 \\ 0 & -4 \end{pmatrix} \middle| \mathbf{y} \right\rangle \geq 0$. In other words,

$$(0, y_2)\begin{pmatrix} 2 & 0 \\ 0 & -4 \end{pmatrix}\begin{pmatrix} 0 \\ y_2 \end{pmatrix} = -4y_2^2 < 0 \quad \text{for all } y_2 \neq 0.$$

Thus, $|\mathbf{x}^*\rangle = (0,0)^T$ does not satisfy Theorem 10.10, so that is not a local minimum for the problem.

## 10.12 Mathematica Built-in Functions

The commands `FindMinimum` and `NMinimize`, and `FindMinimumPlot` can do constrained optimization for multivariable functions [11,12]. This section represents some examples.

### Mathematica Examples 10.1

```
Input     (* Find a local minimum,starting at x=7,subject to constraints 1<=x<=15: *)
          FindMinimum[
           {x Cos[x],1<=x<=15},
           {x,7}
           ]
Output    {-9.47729,{x->9.52933}}

Input     (* Find the minimum of a linear function,subject to linear and integer
          constraints: *)
          FindMinimum[
           {x+y,x+2 y>=3&&x>=0&&y>=0&&y∈Integers},
           {x,y}
           ]
Output    {2.,{x->0.,y->2}}

Input     (* Find a minimum of a function over a geometric region: *)
          FindMinimum[
           {x+y,{x,y}∈Disk[]},
           {x,y}
           ]
          Show[
           ContourPlot[x+y,{x,y}∈Disk[]],
           Graphics[{Red,PointSize[Large],Point[{x,y}/. Last[%]]}]
           ]
Output    {-1.41422,{x->-0.707108,y->-0.707108}}
Output
```

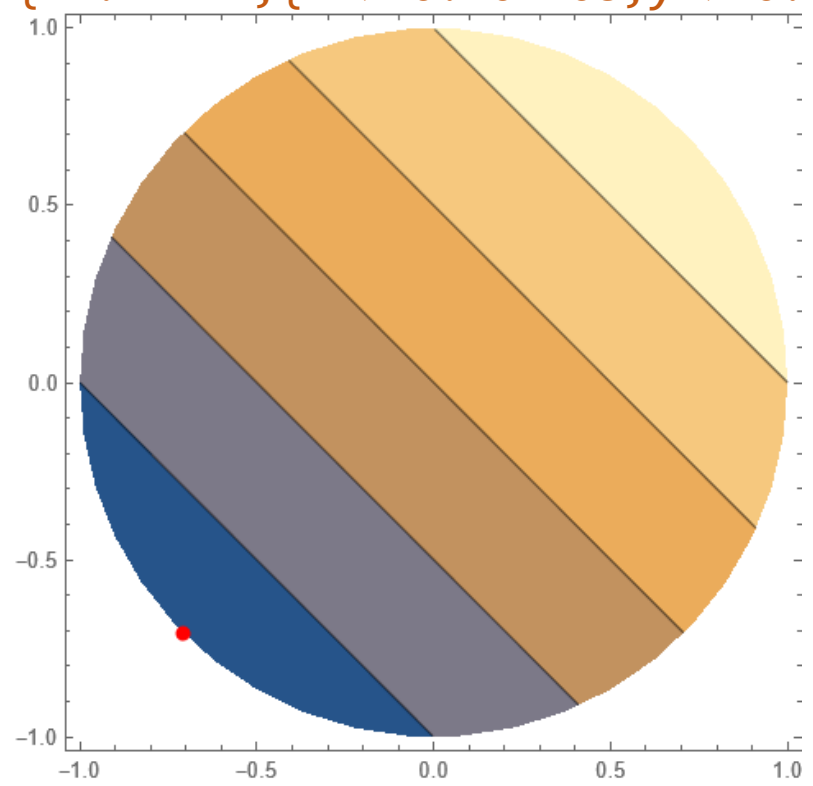

```
Input     (* Find the minimum distance between two regions: *)
          Subscript[ℛ,1]=Disk[];
          Subscript[ℛ,2]=InfiniteLine[{{-2,0},{0,2}}];
          FindMinimum[
           {(x-u)^2+(y-v)^2,{{x,y}∈Subscript[ℛ,1],{u,v}∈Subscript[ℛ,2]}},
           {x,y,u,v}
           ]
          Graphics[{{LightBlue,Subscript[ℛ,1]},{Green,Subscript[ℛ,2]},{Red,Point[{{x,y},{u
          ,v}}/. Last[%]]}}]
Output    {0.171573,{x->-0.707107,y->0.707107,u->-1.,v->1.}}
Output
```

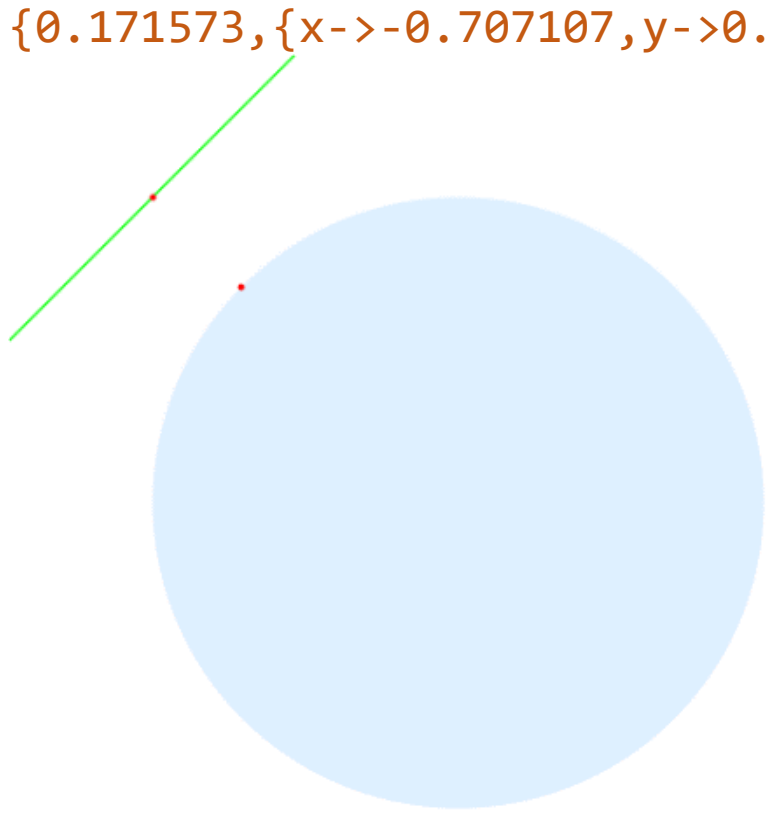

| | |
|---|---|
| $x \prec y$ or VectorLess[{x,y}] | yields True for vectors of length $n$ if xi<yi for all components $1 \le i \le n$. |
| $x \preccurlyeq y$ or VectorLessEqual[{x,y}] | yields True for vectors of length $n$ if xi≤yi for all components $1 \le i \le n$. |
| $x \succ y$ or VectorGreater[{x,y}] | yields True for vectors of length $n$ if xi>yi for all components $1 \le i \le n$. |
| $x \succcurlyeq y$ or VectorGreaterEqual[{x,y}] | yields True for vectors of length $n$ if xi≥yi for all components $1 \le i \le n$. |

### *Mathematica Examples 10.2*

```
Input     {1,2,3}≺{2,3,4}
Output    True

Input     {1,2,3}≺{2,3,3}
Output    False

Input     {1,2,3}≼{1,3,4}
Output    True

Input     {1,2,3}≼{1,2,1}
Output    False
```

| | |
|---|---|
| LinearOptimization[f,cons,vars] | finds values of variables vars that minimize the linear objective f subject to linear constraints cons. |
| LinearOptimization[c,{a,b}] | finds a real vector x that minimizes the linear objective $c.x$ subject to the linear inequality constraints $a.x + b \succcurlyeq 0$. |
| QuadraticOptimization[f,cons,vars] | finds values of variables vars that minimize the quadratic objective f subject to linear constraints cons. |
| QuadraticOptimization[{q,c},{a,b}] | finds a vector that minimizes the quadratic objective ½ $x.q.x + c.x$ subject to the linear inequality constraints $a.x + b \succcurlyeq 0$. |
| ConvexOptimization[f,cons,vars] | finds values of variables vars that minimize the convex objective function f subject to convex constraints cons. |
| LeastSquares[m,b] | finds an x that solves the linear least-squares problem for the matrix equation $m.x == b$. |
| LeastSquares[m,b] | gives a vector x that minimizes $Norm[m.x - b]$. |

**Linear optimization finds $x \in \mathbb{R}^n$ that solves the problem:**

| | |
|---|---|
| Minimize | $c.x$ |
| subject to constraints | $a.x + b \succcurlyeq 0, a_{eq}.x + b_{eq} = 0$ |
| where | $c \in \mathbb{R}^n, a \in \mathbb{R}^{m\times n}, b \in \mathbb{R}^m, a_{eq} \in \mathbb{R}^{k\times n}, b_{eq} \in \mathbb{R}^k$ |

**Quadratic optimization finds $x$ that solves the primal problem**

| | |
|---|---|
| Minimize | $½\, x.q.x + c.x$ |
| subject to constraints | $a.x + b \succcurlyeq 0, a_{eq}.x + b_{eq} = 0$ |
| where | $q \in S_+^n, c \in \mathbb{R}^n, a \in \mathbb{R}^{m\times n}, b \in \mathbb{R}^m, a_{eq} \in \mathbb{R}^{k\times n}, b_{eq} \in \mathbb{R}^k$ |
| | The space $S_+^n$ consists of symmetric positive semidefinite matrices. |

**Convex optimization finds $x \in \mathbb{R}^n$ that solves the following problem:**

| | |
|---|---|
| Minimize | $f_0(x)$ |
| subject to constraints | $f_j(x) \leq 0, j = 1, \ldots, k$ |
| where | $f_j$ are convex functions |

Convex optimization is global nonlinear optimization for convex functions with convex constraints. For convex problems, the global solution can be found.

***Mathematica Examples 10.3***

```
Input         (* Minimize x+y subject to the constraints x+2y>=3,x>=0,y>=0: *)
              res=LinearOptimization[
                 x+y,
                 {x+2 y>=3,x>=0,y>=0},
                 {x,y}
                 ]
Output        {x->0,y->3/2}

Input         (* Minimize x-y subject to the constraints x+y+z=1/2,x-2z=1 and 2x-y>=1: *)
              LinearOptimization[
               x-y,
               {x+y+z==1/2,x-2 z==1,2 x-y>=1},
               {x,y,z}
               ]
Output        {x->4/7,y->1/7,z->-(3/14)}

Input         (* Use VectorLessEqual to express several LessEqual inequality constraints at
              once: *)
              LinearOptimization[
               x+y,
               VectorLessEqual[{{x-2 y,-x+y},{3,2}}],
               {x,y}
               ]
              (* An equivalent form using scalar inequalities: *)
              LinearOptimization[
               x+y,
               {x-2 y<=3,-x+y<=2},
               {x,y}
               ]
Output        {x->-7,y->-5}
Output        {x->-7,y->-5}
```

```
Input     (* Specify the constraints using a combination of scalar and vector
          inequalities: *)
          LinearOptimization[
           x+y,
           {x+2 y>=3,-1≼{x,y}≼{2,1}},
           {x,y}
           ]
          (* An equivalent form using scalar inequalities: *)
          LinearOptimization[
           x+y,
           {x+2 y>=3,-1<=x<=2,-1<=y<=1},
           {x,y}
           ]
Output    {x->1,y->1}
Output    {x->1,y->1}

Input     (* Minimize 2 x^2 + 20 y^2 + 6 x y + 5 x subject to the constraint -x+y>=2: *)
          obj=2 x^2+20 y^2+6 x y+5 x;
          res=QuadraticOptimization[
            obj,
            -x+y>=2,
            {x,y}
            ]
Output    {x->-1.73214,y->0.267857}

Input     (* Minimize x^2 + y^2 subject to the equality constraint x+y=2 and the
          inequality constraints 1<=y<=2: *)
          res=QuadraticOptimization[
            x^2+y^2,
            {x+y==2,1<=y<=2},
            {x,y}
            ]
Output    {x->1.,y->1.}

Input     (* Define objective as 1/2 x.q.x+c.x and constraints as a.x+b ≽ 0 and a_eq.x+
          b_eq = 0: *)
          {q,c}={{{2,0},{0,2}},{0,0}};
          {a,b}={{{0,1},{0,-1}},{-1,2}};
          {Subscript[a,eq],Subscript[b,eq]}={{{1,1}},{-2}};
          (* Solve using matrix-vector inputs: *)
          QuadraticOptimization[{q,c},{a,b},{Subscript[a,eq],Subscript[b,eq]}]
Output    {1.,1.}

Input     (* Minimize x^2 + y^2 subject to the constraints x+y==1 and x^2+5y^2<1: *)
          ConvexOptimization[
           x^2+y^2,
           {x+y==1,x^2+5 y^2<1},
           {x,y}
           ]
Output    {x->0.666669,y->0.333332}

Input     (* Several linear inequality constraints can be expressed with
          VectorGreaterEqual: *)
          ConvexOptimization[
           x+y,
           VectorGreaterEqual[{{x+2 y,x},{3,-1}}],
           {x,y}
           ]
Output    {x->-1.,y->2.}

Input     (* Solve a simple least-squares problem: *)
```

```
              m = ({
                  {1, 1},
                  {1, 2},
                  {1, 3}
                 });
              b = {7, 7, 8};
              LeastSquares[m, b]
Output        {19/3,1/2}

Input         (*This finds a tuple that minimizes ‖m.{x,y} − b‖:*)
              Minimize[
               Norm[m.{x,y}-b],
               {x,y}
               ]
Output        {1/√6,{x->19/3,y->1/2}}

Input         (*Use LeastSquares to minimize ‖m.x − b‖:*)
              m = ({
                  {1, 2, 3},
                  {4, 5, 6},
                  {7, 8, 9}
                 });
              b = {2, -4, 2};
              LeastSquares[m, b]
Output        {0,0,0}
```